\documentclass[12pt]{article}%
\usepackage{e-jc}
\usepackage[all,cmtip]{xy}
\usepackage[colorlinks=true,citecolor=black,linkcolor=black,urlcolor=blue,
pdfauthor={Darij Grinberg, Tom Roby}, pdftitle={Iterative properties of
birational rowmotion}]{hyperref}
\usepackage{amsfonts}
\usepackage{amssymb}
\usepackage{amsmath}
\usepackage{framed}
\usepackage{comment}
\usepackage{needspace}
\usepackage{color}
\usepackage{amsthm}
\providecommand{\U}[1]{\protect\rule{.1in}{.1in}}
\theoremstyle{definition}
\newtheorem{theo}{Theorem}[section]
\newenvironment{theorem}[1][]
{\begin{theo}[#1]\begin{leftbar}}
{\end{leftbar}\end{theo}}
\newtheorem{lem}[theo]{Lemma}
\newenvironment{lemma}[1][]
{\begin{lem}[#1]\begin{leftbar}}
{\end{leftbar}\end{lem}}
\newtheorem{prop}[theo]{Proposition}
\newenvironment{proposition}[1][]
{\begin{prop}[#1]\begin{leftbar}}
{\end{leftbar}\end{prop}}
\newtheorem{defi}[theo]{Definition}
\newenvironment{definition}[1][]
{\begin{defi}[#1]\begin{leftbar}}
{\end{leftbar}\end{defi}}
\newtheorem{remk}[theo]{Remark}
\newenvironment{remark}[1][]
{\begin{remk}[#1]\begin{leftbar}}
{\end{leftbar}\end{remk}}
\newtheorem{coro}[theo]{Corollary}
\newenvironment{corollary}[1][]
{\begin{coro}[#1]\begin{leftbar}}
{\end{leftbar}\end{coro}}
\newtheorem{conv}[theo]{Convention}
\newenvironment{convention}[1][]
{\begin{conv}[#1]\begin{leftbar}}
{\end{leftbar}\end{conv}}
\newtheorem{conj}[theo]{Conjecture}
\newenvironment{conjecture}[1][]
{\begin{conj}[#1]\begin{leftbar}}
{\end{leftbar}\end{conj}}
\newtheorem{exmp}[theo]{Example}
\newenvironment{example}[1][]
{\begin{exmp}[#1]\begin{leftbar}}
{\end{leftbar}\end{exmp}}
\newenvironment{verlong}{}{}
\newenvironment{vershort}{}{}
\newenvironment{noncompile}{}{}
\excludecomment{verlong}
\includecomment{vershort}
\excludecomment{noncompile}
\def\volno{8.0}

\newcommand{\arxiv}[1]{\href{http://arxiv.org/abs/#1}{\texttt{arXiv:#1}}}
\begin{document}

\title{\textbf{Iterative properties of birational rowmotion}}
\author{Darij Grinberg\thanks{Supported by NSF grant 1001905.}\\{\small Department of Mathematics}\\[-0.8ex] {\small Massachusetts Institute of Technology}\\[-0.8ex] {\small Massachusetts, U.S.A.}\\{\small \texttt{darijgrinberg@gmail.com}}\\\phantom {!}
\and Tom Roby\thanks{Supported by NSF grant 1001905.}\\{\small Department of Mathematics}\\[-0.8ex] {\small University of Connecticut}\\[-0.8ex] {\small Connecticut, U.S.A.}\\{\small \texttt{tom.roby@uconn.edu}}}
\date{version {\volno} (\today). \\
\small{This is the arXiv version, not the published version. It has been updated in 2022 and 2026 to correct some typos. The published version has been abridged in several places and split into two papers.} }
\maketitle

\begin{abstract}
We study a birational map associated to any finite poset $P$. This map is a
far-reaching generalization (found by Einstein and Propp) of classical
rowmotion, which is a certain permutation of the set of order ideals of $P$.
Classical rowmotion has been studied by various authors (Fon-der-Flaass,
Cameron, Brouwer, Schrijver, Striker, Williams and many more) under different
guises (Striker-Williams promotion and Panyushev complementation are two
examples of maps equivalent to it). In contrast, birational rowmotion is new
and has yet to reveal several of its mysteries. In this paper, we prove that
birational rowmotion has order $p+q$ on the $\left(  p, q\right)  $-rectangle
poset (i.e., on the product of a $p$-element chain with a $q$-element chain);
we furthermore compute its orders on some triangle-shaped posets and on a
class of posets which we call ``skeletal'' (this class includes all graded
forests). In all cases mentioned, birational rowmotion turns out to have a
finite (and explicitly computable) order, a property it does not exhibit for
general finite posets (unlike classical rowmotion, which is a permutation of a
finite set). Our proof in the case of the rectangle poset uses an idea
introduced by Volkov (arXiv:hep-th/0606094) to prove the $AA$ case of the
Zamolodchikov periodicity conjecture; in fact, the finite order of birational
rowmotion on many posets can be considered an analogue to Zamolodchikov
periodicity. We comment on suspected, but so far enigmatic, connections to the
theory of root posets. We also make a digression to study classical rowmotion
on skeletal posets, since this case has seemingly been overlooked so far.


\end{abstract}

%

\tableofcontents

\begin{noncompile}
\vspace{1cm} \underline{\textbf{Remark on the level of detail}}

This paper comes in two versions: a regular one and a more detailed one. The
regular version is optimized for readability, leaving out the more
straightforward parts and technical arguments. The more detailed version has
many of them expanded (though not all, for the time being).

\begin{vershort}
This is the regular version of the paper. The more detailed one can be
obtained by replacing\newline\texttt{%
$\backslash$%
excludecomment\{verlong\}}\newline\texttt{%
$\backslash$%
includecomment\{vershort\}}\newline by\newline\texttt{%
$\backslash$%
includecomment\{verlong\}}\newline\texttt{%
$\backslash$%
excludecomment\{vershort\}}\newline in the preamble of the LaTeX sourcecode
and then compiling to PDF.
\end{vershort}

\begin{verlong}
This is the more detailed version of the paper. The two versions share the
same .tex file, with the only difference that there are two lines in the
preamble of the file which need to be modified in order to switch between the
short and the detailed version. Namely, these lines are \newline\texttt{%
$\backslash$%
excludecomment\{verlong\}}\newline\texttt{%
$\backslash$%
includecomment\{vershort\}}\newline for the short version and \newline\texttt{%
$\backslash$%
includecomment\{verlong\}}\newline\texttt{%
$\backslash$%
excludecomment\{vershort\}}\newline for the detailed one.
\end{verlong}
\end{noncompile}

\section*{Introduction}

The present paper had originally been intended as a companion paper to David
Einstein's and James Propp's work \cite{einstein-propp}, which introduced
piecewise-linear and birational rowmotion as extensions of the classical
concept of rowmotion on order ideals. While the present paper is
mathematically self-contained (and indeed gives some proofs on which
\cite{einstein-propp} relies), it provides only a modicum of motivation and
applications for the results it discusses. For the latter, the reader may
consult \cite{einstein-propp}.

Let $P$ be a finite poset, and $J\left(  P\right)  $ the set of the order
ideals\footnote{An \textit{order ideal} of a poset $P$ is a subset $S$ of $P$
such that every $s\in S$ and $p\in P$ with $p\leq s$ satisfy $p\in S$.} of
$P$. \textit{Rowmotion} is a classical map $J\left(  P\right)  \rightarrow
J\left(  P\right)  $ which can be defined in various ways, one of which is as
follows: For every $v\in P$, let $\mathbf{t}_{v}:J\left(  P\right)
\rightarrow J\left(  P\right)  $ be the map sending every order ideal $S\in
J\left(  P\right)  $ to $\left\{
\begin{array}
[c]{l}%
S\cup\left\{  v\right\}  \text{, if }v\notin S\text{ and }S\cup\left\{
v\right\}  \in J\left(  P\right)  ;\\
S\setminus\left\{  v\right\}  \text{, if }v\in S\text{ and }S\setminus\left\{
v\right\}  \in J\left(  P\right)  ;\\
S\text{, otherwise}%
\end{array}
\right.  $. These maps $\mathbf{t}_{v}$ are called \textit{classical
toggles}\footnote{or just \textit{toggles} in literature which doesn't occupy
itself with birational rowmotion}, since all they do is \textquotedblleft
toggle\textquotedblright\ an element into or out of an order ideal. Let
$\left(  v_{1},v_{2},...,v_{m}\right)  $ be a linear extension of $P$ (see
Definition \ref{def.linext} for the meaning of this). Then, (classical)
rowmotion is defined as the composition $\mathbf{t}_{v_{1}}\circ
\mathbf{t}_{v_{2}}\circ...\circ\mathbf{t}_{v_{m}}$ (which, as can be seen,
does not depend on the choice of the particular linear extension $\left(
v_{1},v_{2},...,v_{m}\right)  $). This rowmotion map has been studied from
various perspectives; in particular, it is isomorphic\footnote{By this,
we mean that there exists a bijection $\phi$ from $J\left(P\right)$ to the
set of all antichains of $P$ such that rowmotion is $\phi^{-1} \circ f
\circ \phi$.
} to the map $f$ of
Fon-der-Flaass \cite{fon-der-flaass}\footnote{Indeed, let $\mathcal{A}\left(
P\right)  $ denote the set of all antichains of $P$. Then, the map $J\left(
P\right)  \rightarrow\mathcal{A}\left(  P\right)  $ which sends every order
ideal $I\in J\left(  P\right)  $ to the antichain of the maximal elements of
$I$ is a bijection which intertwines rowmotion and Fon-der-Flaass' map $f$.},
the map $F^{-1}$ of Brouwer and Schrijver \cite{brouwer-schrijver}, and the
map $f^{-1}$ of Cameron and Fon-der-Flaass \cite{cameron-fon-der-flaass}%
\footnote{This time, the intertwining bijection from rowmotion to the map
$f^{-1}$ of \cite{cameron-fon-der-flaass} is given by mapping every order
ideal $I$ to its indicator function. This is a bijection from $J\left(
P\right)  $ to the set of Boolean monotonic functions $P\rightarrow\left\{
0,1\right\}  $.}. More recently, it has been studied (and christened
\textquotedblleft rowmotion\textquotedblright) in Striker and Williams
\cite{striker-williams}, where further sources and context are also given.
Since so much has already been said about this rowmotion map, we will only
briefly touch on its properties in Section \ref{sect.classicalr}, while most
of this paper will be spent studying a much more general construction.

Among the questions that have been posed about rowmotion, the most prevalent
was probably that of its order: While it clearly has finite order (being a
bijective map from the finite set $J\left(  P\right)  $ to itself), it turns
out to have a much smaller order than what one would naively expect when the
poset $P$ has certain ``special'' forms (e.g., a rectangle, a root poset, a
product of a rectangle with a $2$-chain, or -- apparently first considered in
this paper -- a forest). Most strikingly, when $P$ is the rectangle $\left[
p\right]  \times\left[  q\right]  $ (denoted $\operatorname*{Rect}\left(
p,q\right)  $ in Definition \ref{def.rect}), then the
\mbox{$\left(  p+q\right)  $-th power} of the rowmotion operator is the
identity map. This is proven in \cite[Theorem 3.6]{brouwer-schrijver} and
\cite[Theorem 2]{fon-der-flaass}\footnote{Another proof follows from two
observations made in \cite{propp-roby-full}: first, that the rowmotion
operator on the order ideals of the rectangle
$\left[p\right]\times\left[q\right]$ is equivalent to the operator named
$\Phi_A$ in \cite{propp-roby-full} (i.e., there is a bijection between
order ideals and antichains of $\left[p\right]\times\left[q\right]$
which intertwines these two
operators), and second, that the $\left(p+q\right)$-th power of this latter
operator $\Phi_A$ is the identity map (this is proven in
\cite[right after Proposition 26]{propp-roby-full}). This argument can also
be constructed from ideas given in \cite[\S 3.3.1]{propp-roby}.}. We will
(in Section
\ref{sect.classicalr}) give a simple algorithm to find the order of rowmotion
on graded forests and similar posets.

In \cite{einstein-propp}, David Einstein and James Propp have lifted the
rowmotion map from the set $J\left(  P\right)  $ of order ideals to the
progressively more general setups of:

\textbf{(a)} the order polytope $\mathcal{O}\left(  P\right)  $ of the poset
$P$ (as defined in \cite[Example 4.6.17]{stanley-ec1} or \cite[Definition
1.1]{stanley-polytopes}), and

\textbf{(b)} even more generally, the affine variety of $\mathbb{K}%
$-labellings of $P$ for $\mathbb{K}$ an arbitrary infinite field.

In case \textbf{(a)}, order ideals of $P$ are replaced by points in the order
polytope $\mathcal{O}\left(  P\right)  $, and the role of the map
$\mathbf{t}_{v}$ (for a given $v\in P$) is assumed by the map which reflects
the $v$-coordinate of a point in $\mathcal{O}\left(  P\right)  $ around the
midpoint of the interval of all values it could take without the point leaving
$\mathcal{O}\left(  P\right)  $ (while all other coordinates are considered
fixed). The operation of ``piecewise linear'' rowmotion (inspired by work of
Arkady Berenstein) is still defined as the composition of these reflection
maps in the same way as rowmotion is the composition of the toggles
$\mathbf{t}_{v}$. This ``piecewise linear'' rowmotion extends (interpolates,
even) classical rowmotion, as order ideals correspond to the vertices of the
order polytope $\mathcal{O}\left(  P\right)  $ (see \cite[Corollary
1.3]{stanley-polytopes}). We will not study case \textbf{(a)} here, since all
of the results we could find in this case can be obtained by tropicalization
from similar results for case \textbf{(b)}.

In case \textbf{(b)}, instead of order ideals of $P$ one considers maps from
the poset $\widehat{P}:=\left\{  0\right\}  \oplus P\oplus\left\{  1\right\}
$ (where $\oplus$ stands for the ordinal sum\footnote{More explicitly,
$\widehat{P}$ is the poset obtained by adding a new element $0$ to $P$, which
is set to be lower than every element of $P$, and adding a new element $1$ to
$P$, which is set to be higher than every element of $P$ (and $0$). We shall
repeat this definition in more formal terms in Definition \ref{def.Phat}.}) to
a given infinite field $\mathbb{K}$ (or, to speak more graphically, of all
labellings of the elements of $P$ by elements of $\mathbb{K}$, along with two
additional labels \textquotedblleft at the very bottom\textquotedblright\ and
\textquotedblleft at the very top\textquotedblright). The maps $\mathbf{t}%
_{v}$ are then replaced by certain birational maps which we call
\textit{birational }$v$\textit{-toggles} (Definition \ref{def.Tv}); the
resulting composition is called \textit{birational rowmotion} and denoted by
$R$. By a careful limiting procedure (the tropical limit), we can
\textquotedblleft degenerate\textquotedblright\ $R$ to the \textquotedblleft
piecewise linear\textquotedblright\ rowmotion of case \textbf{(a)}, and thus
it can be seen as an even higher generalization of classical rowmotion. We
refer to the body of this paper for precise definitions of these maps. Note
that birational $v$-toggles (but not birational rowmotion) in the case of a
rectangle poset have also appeared in \cite[(3.5)]{osz-geoRSK}, but
(apparently) have not been composed there in a way that yields birational rowmotion.

As in the case of classical rowmotion on $J\left(  P\right)  $, the most
interesting question is the order of this map $R$, which in general no longer
has an obvious reason to be finite (since the affine variety of $\mathbb{K}%
$-labellings is not a finite set like $J\left(  P\right)  $). Indeed, for some
posets $P$ this order is infinite. In this paper we will prove the following facts:

\begin{itemize}
\item Birational rowmotion (i.e., the map $R$) on any graded poset (in the
meaning of this word introduced in Definition \ref{def.graded.graded}) has a
very simple effect (namely, cyclic shifting) on the so-called ``w-tuple'' of a
labelling (a rather simple fingerprint of the labelling). This does not mean
$R$ itself has finite order (but turns out to be crucial in proving this in
several cases).

\item Birational rowmotion on graded forests and, slightly more generally,
skeletal posets (Definition \ref{def.skeletal}) has finite order (which can be
bounded from above by an iterative lcm, and also easily computed
algorithmically). Moreover, its order in these cases coincides with the order
of classical rowmotion (Section \ref{sect.classicalr}).

\item Birational rowmotion on a $p\times q$-rectangle has order $p+q$ and
satisfies a further symmetry property (Theorem \ref{thm.rect.antip.general}).
These results have originally been conjectured by James Propp and the second
author, and can be used as an alternative route to certain properties of
(Sch\"{u}tzenberger's) promotion map on semistandard Young tableaux.

\item Birational rowmotion on certain triangle-shaped posets (this is made
precise in Sections \ref{sect.tria}, \ref{sect.DeltaNabla},
\ref{sect.quartertri}) also has finite order (computed explicitly below). We
show this for three kinds of triangle-shaped posets (obtained by cutting the
$p\times p$-square in two along either of its two diagonals) and conjecture it
for a fourth (a quarter of a $p\times p$-square obtained by cutting it along
both diagonals).
\end{itemize}

The proof of the most difficult and fundamental case -- that of a $p\times
q$-rectangle -- is inspired by Volkov's proof of the ``rectangular''
(type-$AA$) Zamolodchikov conjecture \cite{volkov}, which uses a similar idea
of parametrizing (generic) $\mathbb{K}$-labellings by matrices (or tuples of
points in projective space). There is, of course, a striking similarity
between the fact itself and the Zamolodchikov conjecture; yet, we were not
able to reduce either result to the other.

Applications of the results of this paper (specifically Theorems
\ref{thm.rect.ord} and \ref{thm.rect.antip.general}) are found in
\cite{einstein-propp}. Further directions currently under study of the authors
are relations to the totally positive Grassmannian and generalizations to
further classes of posets.

An extended (12-page) abstract \cite{grinberg-roby} of this paper has been
published in the proceedings of the FPSAC 2014 conference.

For publication, this preprint has been split into two papers:
\begin{itemize}
\item
\textquotedblleft\textit{Iterative properties of birational rowmotion
I: generalities and skeletal posets}\textquotedblright\ (published in
\href{http://www.combinatorics.org/ojs/index.php/eljc/article/view/v23i1p33}{Volume 23, Issue 1 (2016)
of the Electronic Journal of Combinatorics}),
and
\item \textquotedblleft\textit{Iterative properties of
birational rowmotion II: rectangles and triangles}\textquotedblright
\ (published in
\href{http://www.combinatorics.org/ojs/index.php/eljc/article/view/v22i3p40}{Volume 22, Issue 3 (2015)
of the Electronic Journal of Combinatorics}).
\end{itemize}
These two papers are
somewhat less detailed than the preprint that you are currently reading;
they also (unlike this preprint) have undergone some stylistic changes
during the refereeing process.\footnote{However, all errors found by the
referees have been corrected both in the two papers and in this preprint.}

\subsection{Leitfaden}

The following Hasse diagram shows how the sections of this paper depend upon
each other.%
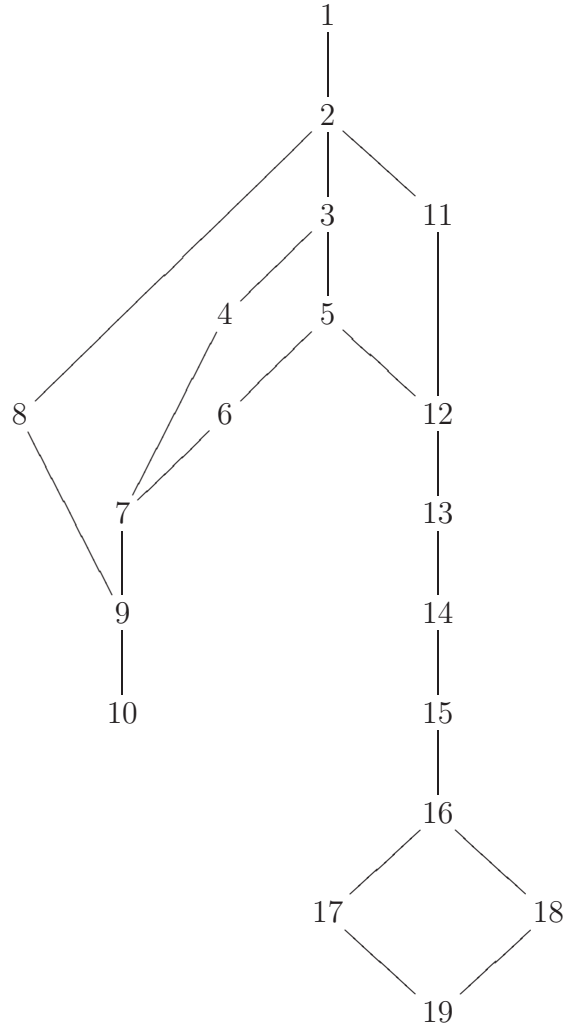
\[
\xymatrix{
& & & 1 \ar@{-}[d] & & \\
& & & 2 \ar@{-}[dddlll] \ar@{-}[d] \ar@{-}[dr] & & \\
& & & 3 \ar@{-}[dl] \ar@{-}[d] & 11 \ar@{-}[dd] & \\
& & 4 \ar@{-}[ddl] & 5 \ar@{-}[dl] \ar@{-}[dr] & & \\
8 \ar@{-}[ddr] & & 6 \ar@{-}[dl] & & 12 \ar@{-}[d] & \\
& 7 \ar@{-}[d] & & & 13 \ar@{-}[d] & \\
& 9 \ar@{-}[d] & & & 14 \ar@{-}[d] & \\
& 10 & & & 15 \ar@{-}[d] & \\
& & & & 16 \ar@{-}[dl] \ar@{-}[dr] & \\
& & & 17 \ar@{-}[dr] & & 18 \ar@{-}[dl] \\
& & & & 19 &
}
\]
A section $n$ depends substantially on a section $m$ if and only if $m>n$ in
the poset whose Hasse diagram is depicted above. Only substantial dependencies
are shown; dependencies upon definitions do not count as substantial (e.g.,
many sections depend on Definition \ref{def.ord}, but this does not make them
substantially dependent on Section \ref{sect.ord}), and dependencies which are
only used in proving inessential claims do not count (e.g., the proof of
Theorem \ref{thm.rect.ord} relies on Proposition \ref{prop.ord-projord} in
order to show that $\operatorname*{ord}\left(  R_{\operatorname*{Rect}\left(
p,q\right)  }\right)  =p+q$ rather than just $\operatorname*{ord}\left(
R_{\operatorname*{Rect}\left(  p,q\right)  }\right)  \mid p+q$, but since the
$\operatorname*{ord}\left(  R_{\operatorname*{Rect}\left(  p,q\right)
}\right)  \mid p+q$ statement is in our opinion the only important part of the
theorem, we do not count this as a dependency on Section \ref{sect.ord}).
Sections \ref{sect.negres} and \ref{sect.rootsys} are not shown.

No section of this paper depends on the Introduction.

\subsection{Acknowledgments}

When confronted with the (then open) problem of proving what is Theorem
\ref{thm.rect.ord} in this paper, Pavlo Pylyavskyy and Gregg Musiker suggested
reading \cite{volkov}. This suggestion proved highly useful and built the
cornerstone of this paper, without which the latter would have ended at its
``Skeletal posets'' section.

The notion of birational rowmotion is due to James Propp and Arkady
Berenstein. This paper owes James Propp also for a constant flow of
inspiration and useful suggestions.

David Einstein found errors in our computations, Bruce Sagan in
Proposition~\ref{prop.classical.for.ord.explicit}, and Hugh Thomas corrected
slips in the writing including an abuse of Zariski topology and some
accidental alternative history.

Nathan Williams noticed typos, too, and suggested a path connecting this
subject to the theory of minuscule posets (which we will not explore in this paper).

The first author came to know birational rowmotion in Alexander Postnikov's
combinatorics pre-seminar at MIT. Postnikov also suggested veins of further study.

Jessica Striker helped the first author understand some of the past work on
this subject, in particular the labyrinthine connections between the various
operators (rowmotion, Panyushev complementation, Striker-Williams promotion,
Sch\"{u}tzenberger promotion, etc.). The present paper explores merely one
corner of this labyrinth (the rowmotion corner).

We thank Dan Bump, Anne Schilling and the two referees of our FPSAC abstract
\cite{grinberg-roby} for further helpful comments. We also owe a number of
improvements in this paper to the suggestions of two anonymous EJC referees.

Both authors were partially supported by NSF grant \#1001905, and have
utilized the open-source CAS Sage (\cite{sage}, \cite{sage-combinat}) to
perform laborious computations. We thank Travis Scrimshaw, Fr\'{e}d\'{e}ric
Chapoton, Viviane Pons and Nathann Cohen for reviewing Sage patches relevant
to this project.

\section{Linear extensions of posets}

This first section serves to introduce some general notions concerning posets
and their linear extensions. In particular, we highlight that the set of
linear extensions of any finite poset is non-empty and connected by a simple
equivalence relation (Proposition \ref{prop.linext.transitive}). This will be
used in subsequent sections for defining the basic maps that we consider
throughout the paper.

Let us first get a basic convention out of the way:

\begin{convention}
We let $\mathbb{N}$ denote the set
$\left\{  0,1,2,...\right\}  $.
\end{convention}

We start by defining general notations related to posets:

\begin{definition}
Let $P$ be a poset. Let $u\in P$ and $v\in P$. In this definition, we will use
$\leq$, $<$, $\geq$ and $>$ to denote the lesser-or-equal relation, the lesser
relation, the greater-or-equal relation and the greater relation,
respectively, of the poset $P$.

\textbf{(a)} The elements $u$ and $v$ of $P$ are said to be
\textit{incomparable} if we have neither $u\leq v$ nor $u\geq v$.

\textbf{(b)} We write $u\lessdot v$ if we have $u<v$ and there is no $w\in P$
such that $u<w<v$. One often says that ``$u$ is covered by $v$'' to signify
that $u\lessdot v$.

\textbf{(c)} We write $u\gtrdot v$ if we have $u>v$ and there is no $w\in P$
such that $u>w>v$. (Thus, $u\gtrdot v$ holds if and only if $v\lessdot u$.)
One often says that ``$u$ covers $v$'' to signify that $u\gtrdot v$.

\textbf{(d)} An element $u$ of $P$ is called \textit{maximal} if every $v\in
P$ satisfying $v\geq u$ satisfies $v=u$. It is easy to see that every nonempty
finite poset has at least one maximal element.

\textbf{(e)} An element $u$ of $P$ is called \textit{minimal} if every $v\in
P$ satisfying $v\leq u$ satisfies $v=u$. It is easy to see that every nonempty
finite poset has at least one minimal element.

When any of these notations becomes ambiguous because the elements involved
belong to several different posets simultaneously, we will disambiguate it by
adding the words \textquotedblleft in $P$\textquotedblright\ (where $P$ is the
poset which we want to use).\footnotemark
\end{definition}

\begin{vershort}
\footnotetext{For instance, if $R$ denotes the poset $\mathbb{Z}$ endowed with
the reverse of its usual order, then we say (for instance) that
\textquotedblleft$1\lessdot0$ in $R$\textquotedblright\ rather than just
\textquotedblleft$1\lessdot0$\textquotedblright.}
\end{vershort}

\begin{verlong}
\footnotetext{For example, let $S$ be the poset defined as the set
$\mathbb{Z}$ endowed with the usual order. Further, let $R$ be the poset
defined as the set $\mathbb{N}$ with order relation given by each
$a\in\mathbb{N}$ being greater than $a+2$. Then, for example, we say
\textquotedblleft$2\lessdot0$ in $R$\textquotedblright\ to clarify that $2$ is
covered by $0$ with respect to the poset $R$, not the poset $S$. And similarly
we say \textquotedblleft$1$ and $0$ are incomparable in $R$\textquotedblright,
and \textquotedblleft$0$ is a minimal element of $R$\textquotedblright, etc.}
\end{verlong}

\begin{definition}
\label{def.linext}Let $P$ be a finite poset. A \textit{linear extension} of
$P$ will mean a list $\left(  v_{1},v_{2},...,v_{m}\right)  $ of the elements
of $P$ such that every element of $P$ occurs exactly once in this list, and
such that any $i\in\left\{  1,2,...,m\right\}  $ and $j\in\left\{
1,2,...,m\right\}  $ satisfying $v_{i}<v_{j}$ (where $<$ is the smaller
relation of $P$) must satisfy $i<j$.
\end{definition}

\begin{vershort}
A brief remark on this definition is in order. Stanley, in \cite[one paragraph
below the proof of Proposition 3.5.2]{stanley-ec1}, defines a linear extension
of a poset $P$ as an order-preserving bijection from $P$ to the chain
$\left\{  1,2,...,\left|  P\right|  \right\}  $; this is equivalent to our
definition (indeed, our linear extension $\left(  v_{1},v_{2},...,v_{m}%
\right)  $, whose length obviously is $m = \left|  P\right|  $, corresponds to
the bijection $P\rightarrow\left\{  1,2,...,\left|  P\right|  \right\}  $
which sends each $v_{i}$ to $i$). Another widespread definition of a linear
extension of $P$ is as a total order on $P$ compatible with the given order of
the poset $P$; this is equivalent to our definition as well (the total order
is the one defined by $v_{i}<v_{j}$ whenever $i<j$).
\end{vershort}

\begin{verlong}
\begin{convention}
Let $P$ be a poset. Let $u$ and $v$ be two elements of $P$. We say that
\textquotedblleft$u<v$ in $P$\textquotedblright\ if and only if $u$ is smaller
than $v$ with respect to the partial order on $P$. This notation can (and
will) be used even if we have \textbf{not} defined $<$ as the smaller relation
on $P$ (so that just saying \textquotedblleft$u<v$\textquotedblright\ without
mentioning $P$ would either make no sense, or mean something different than
saying \textquotedblleft$u<v$ in $P$\textquotedblright); therefore, it is
important to not leave out the \textquotedblleft in $P$\textquotedblright%
\ part. For example, if $R$ is the poset defined as the set $\mathbb{Z}$ with
order given by $a$ being smaller than $a-2$ for every $a\in\mathbb{Z}$, then
we have $5<1$ in $R$, but surely the statement \textquotedblleft%
$5<1$\textquotedblright\ doesn't hold without the \textquotedblleft in
$R$\textquotedblright\ context.

Similarly, \textquotedblleft$u\leq v$ in $P$\textquotedblright\ will mean that
$u$ is smaller or equal to $v$ with respect to the partial order on $P$, and
similarly the statements \textquotedblleft$u>v$ in $P$\textquotedblright\ and
\textquotedblleft$u\geq v$ in $P$\textquotedblright\ have to be understood.
For example, if $P$ and $Q$ are two posets which are identical as sets, but
the smaller relation of $P$ coincides with the greater relation of $Q$, then
saying \textquotedblleft$u<v$ in $P$\textquotedblright\ is equivalent to
saying \textquotedblleft$u>v$ in $Q$\textquotedblright.

We will occasionally leave out the words \textquotedblleft in $P$%
\textquotedblright\ from these notations when there really is no ambiguity
(i.e., when there is no way to interpret the statement \textquotedblleft%
$u<v$\textquotedblright\ other than as \textquotedblleft$u<v$ in
$P$\textquotedblright).
\end{convention}

\begin{remark}
\label{rmk.linext}The way linear extensions are defined in Definition
\ref{def.linext} is more or less the same as how linear extensions are
represented as lists in \cite[two paragraphs below Proposition 3.5.2]%
{stanley-ec1} (specifically, our list $\left(  v_{1},v_{2},...,v_{m}\right)  $
corresponds to what Stanley calls the \textquotedblleft permutation
$\sigma^{-1}\left(  1\right)  ,...,\sigma^{-1}\left(  p\right)  $%
\textquotedblright). But many other authors define the notion of a linear
extension differently from Definition \ref{def.linext}. Indeed, the usual way
to define a linear extension of a poset $P$ is by saying that a linear
extension of $P$ means a totally ordered set $\left(  P,\trianglelefteq
\right)  $ which has $P$ as its ground set and whose smaller-or-equal relation
$\trianglelefteq$ satisfies $p\trianglelefteq q$ for any two elements $p$ and
$q$ of $P$ satisfying $p\leq q$ in $P$. But when $P$ is a finite poset, this
definition is equivalent to the Definition \ref{def.linext} that we gave in
the sense that there is a bijection between:

\textbf{(a)} lists $\left(  v_{1},v_{2},...,v_{m}\right)  $ of the elements of
$P$ such that every element of $P$ occurs exactly once in this list, and such
that any $i\in\left\{  1,2,...,m\right\}  $ and $j\in\left\{
1,2,...,m\right\}  $ satisfying $v_{i}<v_{j}$ (where $<$ is the smaller
relation of $P$) must satisfy $i<j$,

and

\textbf{(b)} totally ordered sets $\left(  P,\trianglelefteq\right)  $ which
have $P$ as their ground set and whose smaller-or-equal relation
$\trianglelefteq$ satisfies $p\trianglelefteq q$ for any two elements $p$ and
$q$ of $P$ satisfying $p\leq q$ in $P$.

Indeed, given a totally ordered set $\left(  P,\trianglelefteq\right)  $ as in
\textbf{(b)}, we can construct a list $\left(  v_{1},v_{2},...,v_{m}\right)  $
as in \textbf{(a)} by setting $m=\left\vert P\right\vert $ and%
\[
\left(
\begin{array}
[c]{c}%
v_{i}=\left(  \text{the }i\text{-th smallest element of the totally ordered
set }\left(  P,\trianglelefteq\right)  \right) \\
\ \ \ \ \ \ \ \text{ for every }i\in\left\{  1,2,...,\left\vert P\right\vert
\right\}
\end{array}
\right)  .
\]
Conversely, given a list $\left(  v_{1},v_{2},...,v_{m}\right)  $ as in
\textbf{(a)}, we can obtain a totally ordered set $\left(  P,\trianglelefteq
\right)  $ as in \textbf{(b)} by letting any two elements $p$ and $q$ of $P$
satisfy $p\trianglelefteq q$ if and only if the indices $i$ and $j$ satisfying
$p=v_{i}$ and $q=v_{j}$ satisfy $i\leq j$. It is easy to see (and well-known)
that these ways to assign a list $\left(  v_{1},v_{2},...,v_{m}\right)  $ as
in \textbf{(a)} to a totally ordered set $\left(  P,\trianglelefteq\right)  $
as in \textbf{(b)} and vice versa are mutually inverse. Hence, lists $\left(
v_{1},v_{2},...,v_{m}\right)  $ as in \textbf{(a)} are in bijection to totally
ordered sets $\left(  P,\trianglelefteq\right)  $ as in \textbf{(b)}. So the
usual definition of linear extensions of $P$ is equivalent to our Definition
\ref{def.linext} (up to isomorphism).
\end{remark}
\end{verlong}

Notice that if $\left(  v_{1},v_{2},...,v_{m}\right)  $ is a linear extension
of a nonempty finite poset $P$, then $v_{1}$ is a minimal element of $P$ and
$v_{m}$ is a maximal element of $P$. The only linear extension of the empty
poset $\varnothing$ is the empty list $\left(  {}\right)  $.

\begin{theorem}
\label{thm.linext.ex}Let $P$ be a finite poset. Then, there exists a linear
extension of $P$.
\end{theorem}

\begin{vershort}
Theorem \ref{thm.linext.ex} is a well-known fact, and can be proven, e.g., by
induction over $\left\vert P\right\vert $ (with the induction step consisting
of splitting off a maximal element $u$ of $P$ and appending it to a linear
extension of the residual poset $P\setminus\left\{  u\right\}  $).
\end{vershort}

\begin{verlong}
Theorem \ref{thm.linext.ex} is a well-known fact, and can be proven, e.g.,
using the following basic facts:

\begin{theorem}
\label{thm.linext.induc}Let $P$ be a finite poset. Let $u$ be a maximal
element of $P$. Let $\left(  a_{1},a_{2},...,a_{\ell}\right)  $ be a linear
extension of the poset $P\setminus\left\{  u\right\}  $ (where $P\setminus
\left\{  u\right\}  $ is seen as a subposet of $P$ by restricting the smaller
relation of $P$ onto $P\setminus\left\{  u\right\}  $). Then, $\left(
a_{1},a_{2},...,a_{\ell},u\right)  $ is a linear extension of the poset $P$.
\end{theorem}

\begin{theorem}
\label{thm.linext.max}Let $P$ be a nonempty finite poset. Then, the poset $P$
has a maximal element.
\end{theorem}
\end{verlong}

The following proposition can be easily checked by the reader:

\begin{proposition}
\label{prop.linext.switch}Let $P$ be a finite poset. Let $\left(  v_{1}%
,v_{2},...,v_{m}\right)  $ be a linear extension of $P$. Let $i\in\left\{
1,2,...,m-1\right\}  $ be such that the elements $v_{i}$ and $v_{i+1}$ of $P$
are incomparable. Then, $\left(  v_{1},v_{2},...,v_{i-1},v_{i+1},v_{i}%
,v_{i+2},v_{i+3},...,v_{m}\right)  $ (this is the tuple obtained from the
tuple $\left(  v_{1},v_{2},...,v_{m}\right)  $ by interchanging the adjacent
entries $v_{i}$ and $v_{i+1}$) is a linear extension of $P$ as well.
\end{proposition}

\begin{definition}
\label{def.L(P)}Let $P$ be a finite poset. The set of all linear extensions of
$P$ will be called $\mathcal{L}\left(  P\right)  $. Thus, $\mathcal{L}\left(
P\right)  \neq\varnothing$ (by Theorem \ref{thm.linext.ex}).
\end{definition}

In our approach to birational rowmotion, we will use the following fact (which
is folklore and has applications in various contexts, including Young tableau theory):

\begin{proposition}
\label{prop.linext.transitive}Let $P$ be a finite poset. Let $\sim$ denote the
equivalence relation on $\mathcal{L}\left(  P\right)  $ generated by the
following requirement: For any linear extension $\left(  v_{1},v_{2}%
,...,v_{m}\right)  $ of $P$ and any $i\in\left\{  1,2,...,m-1\right\}  $ such
that the elements $v_{i}$ and $v_{i+1}$ of $P$ are incomparable, we set
$\left(  v_{1},v_{2},...,v_{m}\right)  \sim\left(  v_{1},v_{2},...,v_{i-1}%
,v_{i+1},v_{i},v_{i+2},v_{i+3},...,v_{m}\right)  $ (noting that $\left(
v_{1},v_{2},...,v_{i-1},v_{i+1},v_{i},v_{i+2},v_{i+3},...,v_{m}\right)  $ is
also a linear extension of $P$, because of Proposition
\ref{prop.linext.switch}). {\footnotemark\ } Then, any two elements of
$\mathcal{L}\left(  P\right)  $ are equivalent under the relation $\sim$.
\end{proposition}

\footnotetext{Here is a more formal way to restate this definition of $\sim$:
\par
We first introduce a binary relation $\equiv$ on the set $\mathcal{L}\left(
P\right)  $ as follows: If $\mathbf{v}$ and $\mathbf{w}$ are two linear
extensions of $P$, then we set $\mathbf{v}\equiv\mathbf{w}$ if and only if the
list $\mathbf{w}$ can be obtained from the list $\mathbf{v}$ by interchanging
two adjacent entries $v$ and $v^{\prime}$ which are incomparable in $P$. It is
clear that this binary relation $\equiv$ is symmetric. It is also clear that
for any linear extension $\left(  v_{1},v_{2},...,v_{m}\right)  $ of $P$ and
any $i\in\left\{  1,2,...,m-1\right\}  $ such that the elements $v_{i}$ and
$v_{i+1}$ of $P$ are incomparable, the list $\left(  v_{1},v_{2}%
,...,v_{i-1},v_{i+1},v_{i},v_{i+2},v_{i+3},...,v_{m}\right)  $ is also a
linear extension of $P$ (according to Proposition \ref{prop.linext.switch})
and satisfies $\left(  v_{1},v_{2},...,v_{m}\right)  \equiv\left(  v_{1}%
,v_{2},...,v_{i-1},v_{i+1},v_{i},v_{i+2},v_{i+3},...,v_{m}\right)  $. Now, we
define $\sim$ as the reflexive and transitive closure of the binary relation
$\equiv$. Then, $\sim$ is an equivalence relation on $\mathcal{L}\left(
P\right)  $.} This proposition is very basic (it generalizes the fact that the
symmetric group $S_{n}$ is generated by the adjacent-element transpositions)
and is classical, and proofs can be found in the literature. One proof is
in \cite[Proposition 4.1 (for the $\pi^{\prime}=\pi\tau_{j}$ case)]%
{ayyer-klee-schilling}; another is sketched in \cite[p. 79]{ruskey} and
presented in more detail in \cite[Lemma 1]{etienne}.\footnote{For the sake of
intellectual enrichment, let us outline yet another proof, which has been
suggested by Thomas McConville. This proof is geometric (it uses the theory of
hyperplane arrangements), and making it fully precise would require certain
topological technicalities which we shall not delve into; but the idea is
instructive and provides intuition. Namely, assume WLOG that $P=\left\{
1,2,\ldots,n\right\}  $ as sets. A \textit{listing} shall mean an $n$-tuple of
distinct elements of $P$. Thus, listings are in bijection with the
permutations of $\left\{  1,2,\ldots,n\right\}  $.
\par
Let $\mathcal{B}$ be the braid arrangement in $\mathbb{R}^{n}$ (that is, the
hyperplane arrangement formed by the hyperplanes $x_{i}=x_{j}$ for all $1\leq
i<j\leq n$). The chambers of $\mathcal{B}$ are known to be in a 1-to-1
correspondence with the listings: Namely, to any listing $\mathbf{p}=\left(
p_{1},p_{2},\ldots,p_{n}\right)  $ corresponds the chamber given by the
inequalities $x_{p_{1}}<x_{p_{2}}<\cdots<x_{p_{n}}$; we denote the latter
chamber by $C\left(  \mathbf{p}\right)  $.
\par
On the other hand, let $K_{P}$ be the open cone in $\mathbb{R}^{n}$ defined by
the inequalities $x_{i}<x_{j}$ for all pairs $\left(  i,j\right)  \in P^{2}$
satisfying $i<j$ in $P$. Then, a linear extension of $P$ is precisely a
listing $\mathbf{p}$ satisfying $C\left(  \mathbf{p}\right)  \subseteq K_{P}$.
\par
Let now $\mathbf{u}$ and $\mathbf{v}$ be two elements of $\mathcal{L}\left(
P\right)  $, that is, two linear extensions of $P$. Thus, $C\left(
\mathbf{u}\right)  \subseteq K_{P}$ and $C\left(  \mathbf{v}\right)  \subseteq
K_{P}$. Pick any two points $y\in C\left(  \mathbf{u}\right)  $ and $z\in
C\left(  \mathbf{v}\right)  $. Then, both $y$ and $z$ lie in $K_{P}$;
therefore, so does every point on the segment joining $y$ with $z$ (since
$K_{P}$ is convex). Thus, there exists a continuous path from $y$ to $z$
staying entirely inside $K_{P}$. By slightly deforming this path, we can
ensure that it never intersects more than one hyperplane of $\mathcal{B}$ at
the same point (at the expense of no longer being a straight line); if we do
this with care, then it still will remain inside $K_{P}$ (since $K_{P}$ is an
open set; here is where we are using some topology).
Consider this latter path. Let $C\left(  \mathbf{u}_{1}\right)
,C\left(  \mathbf{u}_{2}\right)  ,\ldots,C\left(  \mathbf{u}_{k}\right)  $ be
the chambers of $\mathcal{B}$ it traverses. Thus, all the listings
$\mathbf{u}_{1},\mathbf{u}_{2},\ldots,\mathbf{u}_{k}$ are linear extensions of
$P$ (since all the chambers $C\left(  \mathbf{u}_{1}\right)  ,C\left(
\mathbf{u}_{2}\right)  ,\ldots,C\left(  \mathbf{u}_{k}\right)  $ are contained
in $K_{P}$), thus belong to $\mathcal{L}\left(  P\right)  $. Moreover,
$C\left(  \mathbf{u}_{1}\right)  =C\left(  \mathbf{u}\right)  $ (since the
path starts at $y\in C\left(  \mathbf{u}\right)  $), so that $\mathbf{u}%
_{1}=\mathbf{u}$. Similarly, $\mathbf{u}_{k}=\mathbf{v}$. Moreover, for every
$i\in\left\{  1,2,\ldots,k-1\right\}  $, the chambers $C\left(  \mathbf{u}%
_{i}\right)  $ and $C\left(  \mathbf{u}_{i+1}\right)  $ are separated by
precisely one hyperplane. This can easily be translated as follows: For every
$i\in\left\{  1,2,\ldots,k-1\right\}  $, the listing $\mathbf{u}_{i+1}$ is
obtained from $\mathbf{u}_{i}$ by interchanging two adjacent entries. Hence,
for every $i\in\left\{  1,2,\ldots,k-1\right\}  $, we have $\mathbf{u}%
_{i+1}\sim\mathbf{u}_{i}$. Since $\sim$ is an equivalence relation, this shows
that $\mathbf{u}_{k}\sim\mathbf{u}_{1}$. In other words, $\mathbf{v}%
\sim\mathbf{u}$ (since $\mathbf{u}_{1}=\mathbf{u}$ and $\mathbf{u}%
_{k}=\mathbf{v}$), qed.}
In order to keep our paper self-contained, we will
prove it too. Our proof is based on the
following lemma (which is more or less a simple particular case of Proposition
\ref{prop.linext.transitive}):

\begin{lemma}
\label{lem.linext.transitive}Let $P$ be a finite poset. Define the equivalence
relation $\sim$ on $\mathcal{L}\left(  P\right)  $ as in Proposition
\ref{prop.linext.transitive}. Let $a_{1}$, $a_{2}$, $...$, $a_{k}$ be some
elements of $P$. Let $b_{1}$, $b_{2}$, $...$, $b_{\ell}$ be some further
elements of $P$. Let $u$ be a maximal element of $P$. Assume that $\left(
a_{1},a_{2},...,a_{k},u,b_{1},b_{2},...,b_{\ell}\right)  $ is a linear
extension of $P$. Then, $\left(  a_{1},a_{2},...,a_{k},b_{1},b_{2}%
,...,b_{\ell},u\right)  $ is a linear extension of $P$ satisfying $\left(
a_{1},a_{2},...,a_{k},u,b_{1},b_{2},...,b_{\ell}\right)  \sim\left(
a_{1},a_{2},...,a_{k},b_{1},b_{2},...,b_{\ell},u\right)  $.
\end{lemma}

\begin{proof}
[Proof of Lemma \ref{lem.linext.transitive} (sketched).]We will show that
every $i\in\left\{  0,1,...,\ell\right\}  $ satisfies the following assertion:%
\begin{equation}
\left(
\begin{array}
[c]{c}%
\text{The tuple }\left(  a_{1},a_{2},...,a_{k},b_{1},b_{2},...,b_{i}%
,u,b_{i+1},b_{i+2},...,b_{\ell}\right)  \text{ is a}\\
\text{linear extension of }P\text{ satisfying}\\
\left(  a_{1},a_{2},...,a_{k},u,b_{1},b_{2},...,b_{\ell}\right)  \sim\left(
a_{1},a_{2},...,a_{k},b_{1},b_{2},...,b_{i},u,b_{i+1},b_{i+2},...,b_{\ell
}\right)
\end{array}
\right)  . \label{pf.linext.transitive.lem1}%
\end{equation}

\textit{Proof of (\ref{pf.linext.transitive.lem1}):} We will prove
(\ref{pf.linext.transitive.lem1}) by induction over $i$:

\begin{vershort}
\textit{Induction base:} If $i=0$, then \newline$\left(  a_{1},a_{2}%
,...,a_{k},b_{1},b_{2},...,b_{i},u,b_{i+1},b_{i+2},...,b_{\ell}\right)
=\left(  a_{1},a_{2},...,a_{k},u,b_{1},b_{2},...,b_{\ell}\right)  $. Hence,
(\ref{pf.linext.transitive.lem1}) is a tautology for $i=0$, and the induction
base is done.
\end{vershort}

\begin{verlong}
\textit{Induction base:} If $i=0$, then \newline$\left(  a_{1},a_{2}%
,...,a_{k},b_{1},b_{2},...,b_{i},u,b_{i+1},b_{i+2},...,b_{\ell}\right)
=\left(  a_{1},a_{2},...,a_{k},u,b_{1},b_{2},...,b_{\ell}\right)  $. Hence, if
$i=0$, then $\left(  a_{1},a_{2},...,a_{k},b_{1},b_{2},...,b_{i}%
,u,b_{i+1},b_{i+2},...,b_{\ell}\right)  $ is a linear extension of $P$
(because we know that $\left(  a_{1},a_{2},...,a_{k},u,b_{1},b_{2}%
,...,b_{\ell}\right)  $ is a linear extension of $P$) and satisfies
\newline$\left(  a_{1},a_{2},...,a_{k},u,b_{1},b_{2},...,b_{\ell}\right)
\sim\left(  a_{1},a_{2},...,a_{k},b_{1},b_{2},...,b_{i},u,b_{i+1}%
,b_{i+2},...,b_{\ell}\right)  $ (because $\sim$ is an equivalence relation and
thus reflexive). In other words, if $i=0$, then
(\ref{pf.linext.transitive.lem1}) is true. This completes the induction base.
\end{verlong}

\textit{Induction step:} Let $I\in\left\{  1,2,...,\ell\right\}  $. Assume
that (\ref{pf.linext.transitive.lem1}) holds for $i=I-1$. We need to prove
that (\ref{pf.linext.transitive.lem1}) holds for $i=I$.

We have assumed that (\ref{pf.linext.transitive.lem1}) holds for $i=I-1$. In
other words, the tuple \newline$\left(  a_{1},a_{2},...,a_{k},b_{1}%
,b_{2},...,b_{I-1},u,b_{I-1+1},b_{I-1+2},...,b_{\ell}\right)  $ is a linear
extension of $P$ satisfying $\left(  a_{1},a_{2},...,a_{k},u,b_{1}%
,b_{2},...,b_{\ell}\right)  \sim\left(  a_{1},a_{2},...,a_{k},b_{1}%
,b_{2},...,b_{I-1},u,b_{I-1+1},b_{I-1+2},...,b_{\ell}\right)  $.

\begin{vershort}
Denote the smaller relation of $P$ by $<$. Since the tuple $\left(
a_{1},a_{2},...,a_{k},u,b_{1},b_{2},...,b_{\ell}\right)  $ is a linear
extension of $P$, we cannot have $u\geq b_{I}$ (because $u$ appears strictly
to the left of $b_{I}$ in this tuple). But we cannot have $u<b_{I}$ either
(since $u$ is a maximal element of $P$). Thus, $u$ and $b_{I}$ are incomparable.
\end{vershort}

\begin{verlong}
Denote the smaller relation of $P$ by $<$. Since the tuple $\left(
a_{1},a_{2},...,a_{k},u,b_{1},b_{2},...,b_{\ell}\right)  $ is a linear
extension of $P$, we cannot have $u>b_{I}$ (because $u$ appears to the left of
$b_{I}$ in this tuple), and we cannot have $u=b_{I}$ (because $u$ and $b_{I}$
appear in different places in this tuple). Thus, we cannot have $u\geq b_{I}%
$.\ \ \ \ \footnote{\textit{Proof.} Assume the contrary. Thus, we have $u\geq
b_{I}$. Since $u\neq b_{I}$ (because we cannot have $u=b_{I}$), this yields
$u>b_{I}$. But this contradicts the fact that we cannot have $u>b_{I}$. This
contradiction shows that our assumption was wrong, qed.} Also, $u$ is a
maximal element of $P$, so we cannot have $u<b_{I}$. Thus, we cannot have
$u\leq b_{I}$.\ \ \ \ \footnote{\textit{Proof.} Assume the contrary. Thus, we
have $u\leq b_{I}$. Since $u\neq b_{I}$ (because we cannot have $u=b_{I}$),
this yields $u<b_{I}$. But this contradicts the fact that we cannot have
$u<b_{I}$. This contradiction shows that our assumption was wrong, qed.}
Altogether, we thus know that $u$ and $b_{I}$ are incomparable (since we
cannot have $u\leq b_{I}$ and since we cannot have $u\geq b_{I}$).
\end{verlong}

\begin{vershort}
Now,%
\begin{align*}
&  \left(  a_{1},a_{2},...,a_{k},u,b_{1},b_{2},...,b_{\ell}\right) \\
&  \sim\left(  a_{1},a_{2},...,a_{k},b_{1},b_{2},...,b_{I-1},u,b_{I-1+1}%
,b_{I-1+2},...,b_{\ell}\right) \\
&  =\left(  a_{1},a_{2},...,a_{k},b_{1},b_{2},...,b_{I-1},u,b_{I}%
,b_{I+1},b_{I+2},...,b_{\ell}\right) \\
&  \sim\left(  a_{1},a_{2},...,a_{k},b_{1},b_{2},...,b_{I-1},b_{I}%
,u,b_{I+1},b_{I+2},...,b_{\ell}\right) \\
&  \ \ \ \ \ \ \ \ \ \ \left(  \text{by the definition of the relation }%
\sim\text{, since }u\text{ and }b_{I}\text{ are incomparable}\right) \\
&  =\left(  a_{1},a_{2},...,a_{k},b_{1},b_{2},...,b_{I},u,b_{I+1}%
,b_{I+2},...,b_{\ell}\right)  .
\end{align*}
The proof of this equivalence also shows that its right hand side is a linear
extension of $P$. Thus, (\ref{pf.linext.transitive.lem1}) holds for $i=I$.
This completes the induction step, whence (\ref{pf.linext.transitive.lem1}) is proven.
\end{vershort}

\begin{verlong}
We know that $\left(  a_{1},a_{2},...,a_{k},b_{1},b_{2},...,b_{I-1}%
,u,b_{I-1+1},b_{I-1+2},...,b_{\ell}\right)  $ is a linear extension of $P$.
Since%
\begin{align}
&  \left(  a_{1},a_{2},...,a_{k},b_{1},b_{2},...,b_{I-1},u,b_{I-1+1}%
,b_{I-1+2},...,b_{\ell}\right) \nonumber\\
&  =\left(  a_{1},a_{2},...,a_{k},b_{1},b_{2},...,b_{I-1},u,b_{I}%
,b_{I+1},...,b_{\ell}\right) \nonumber\\
&  =\left(  a_{1},a_{2},...,a_{k},b_{1},b_{2},...,b_{I-1},u,b_{I}%
,b_{I+1},b_{I+2},...,b_{\ell}\right)  , \label{pf.linext.transitive.lem1.1}%
\end{align}
this yields that $\left(  a_{1},a_{2},...,a_{k},b_{1},b_{2},...,b_{I-1}%
,u,b_{I},b_{I+1},b_{I+2},...,b_{\ell}\right)  $ is a linear extension of $P$.
Thus, $\left(  a_{1},a_{2},...,a_{k},b_{1},b_{2},...,b_{I-1},b_{I}%
,u,b_{I+1},b_{I+2},...,b_{\ell}\right)  $ (this is the tuple obtained from the
tuple $\left(  a_{1},a_{2},...,a_{k},b_{1},b_{2},...,b_{I-1},u,b_{I}%
,b_{I+1},b_{I+2},...,b_{\ell}\right)  $ by interchanging the adjacent entries
$u$ and $b_{I}$) is also a linear extension of $P$ (by the definition of the
relation $\sim$, since $u$ and $b_{I}$ are incomparable) and satisfies%
\begin{align}
&  \left(  a_{1},a_{2},...,a_{k},b_{1},b_{2},...,b_{I-1},u,b_{I}%
,b_{I+1},b_{I+2},...,b_{\ell}\right) \nonumber\\
&  \sim\left(  a_{1},a_{2},...,a_{k},b_{1},b_{2},...,b_{I-1},b_{I}%
,u,b_{I+1},b_{I+2},...,b_{\ell}\right)  \label{pf.linext.transitive.lem1.2}%
\end{align}
(again by the definition of the relation $\sim$, since $u$ and $b_{I}$ are incomparable).

Since $\left(  a_{1},a_{2},...,a_{k},b_{1},b_{2},...,b_{I-1},b_{I}%
,u,b_{I+1},b_{I+2},...,b_{\ell}\right)  $ is a linear extension of $P$, we
know that $\left(  a_{1},a_{2},...,a_{k},b_{1},b_{2},...,b_{I},u,b_{I+1}%
,b_{I+2},...,b_{\ell}\right)  $ is a linear extension of $P$ (because
clearly,
\begin{align*}
&  \left(  a_{1},a_{2},...,a_{k},b_{1},b_{2},...,b_{I-1},b_{I},u,b_{I+1}%
,b_{I+2},...,b_{\ell}\right) \\
&  = \left(  a_{1},a_{2},...,a_{k},b_{1},b_{2},...,b_{I},u,b_{I+1}%
,b_{I+2},...,b_{\ell}\right)
\end{align*}
).

Now,%
\begin{align*}
&  \left(  a_{1},a_{2},...,a_{k},u,b_{1},b_{2},...,b_{\ell}\right) \\
&  \sim\left(  a_{1},a_{2},...,a_{k},b_{1},b_{2},...,b_{I-1},u,b_{I-1+1}%
,b_{I-1+2},...,b_{\ell}\right) \\
&  =\left(  a_{1},a_{2},...,a_{k},b_{1},b_{2},...,b_{I-1},u,b_{I}%
,b_{I+1},b_{I+2},...,b_{\ell}\right)  \ \ \ \ \ \ \ \ \ \ \left(  \text{by
(\ref{pf.linext.transitive.lem1.1})}\right) \\
&  \sim\left(  a_{1},a_{2},...,a_{k},b_{1},b_{2},...,b_{I-1},b_{I}%
,u,b_{I+1},b_{I+2},...,b_{\ell}\right)  \ \ \ \ \ \ \ \ \ \ \left(  \text{by
(\ref{pf.linext.transitive.lem1.2})}\right) \\
&  =\left(  a_{1},a_{2},...,a_{k},b_{1},b_{2},...,b_{I},u,b_{I+1}%
,b_{I+2},...,b_{\ell}\right)  .
\end{align*}
We thus have shown that the tuple $\left(  a_{1},a_{2},...,a_{k},b_{1}%
,b_{2},...,b_{I},u,b_{I+1},b_{I+2},...,b_{\ell}\right)  $ is a linear
extension of $P$ satisfying \newline$\left(  a_{1},a_{2},...,a_{k}%
,u,b_{1},b_{2},...,b_{\ell}\right)  \sim\left(  a_{1},a_{2},...,a_{k}%
,b_{1},b_{2},...,b_{I},u,b_{I+1},b_{I+2},...,b_{\ell}\right)  $. In other
words, we have shown that (\ref{pf.linext.transitive.lem1}) holds for $i=I$.
This completes the induction step. The induction proof of
(\ref{pf.linext.transitive.lem1}) is thus finished.
\end{verlong}

\begin{vershort}
Lemma \ref{lem.linext.transitive} now follows by applying
(\ref{pf.linext.transitive.lem1}) to $i=\ell$. \qedhere
\end{vershort}

\begin{verlong}
Now we can apply (\ref{pf.linext.transitive.lem1}) to $i=\ell$. As a result,
we obtain that the tuple \newline$\left(  a_{1},a_{2},...,a_{k},b_{1}%
,b_{2},...,b_{\ell},u,b_{\ell+1},b_{\ell+2},...,b_{\ell}\right)  $ is a linear
extension of $P$ satisfying \newline$\left(  a_{1},a_{2},...,a_{k}%
,u,b_{1},b_{2},...,b_{\ell}\right)  \sim\left(  a_{1},a_{2},...,a_{k}%
,b_{1},b_{2},...,b_{\ell},u,b_{\ell+1},b_{\ell+2},...,b_{\ell}\right)  $.
Since \newline$\left(  a_{1},a_{2},...,a_{k},b_{1},b_{2},...,b_{\ell
},u,b_{\ell+1},b_{\ell+2},...,b_{\ell}\right)  =\left(  a_{1},a_{2}%
,...,a_{k},b_{1},b_{2},...,b_{\ell},u\right)  $, this sentence rewrites as
follows: The tuple $\left(  a_{1},a_{2},...,a_{k},b_{1},b_{2},...,b_{\ell
},u\right)  $ is a linear extension of $P$ satisfying $\left(  a_{1}%
,a_{2},...,a_{k},u,b_{1},b_{2},...,b_{\ell}\right)  \sim\left(  a_{1}%
,a_{2},...,a_{k},b_{1},b_{2},...,b_{\ell},u\right)  $. This proves Lemma
\ref{lem.linext.transitive}.
\end{verlong}
\end{proof}

\begin{vershort}
\begin{proof}
[Proof of Proposition \ref{prop.linext.transitive} (sketched).]We prove
Proposition \ref{prop.linext.transitive} by induction over $\left\vert
P\right\vert $. The induction base $\left\vert P\right\vert =0$ is trivial.
For the induction step, let $N$ be a positive integer. Assume that Proposition
\ref{prop.linext.transitive} is proven for all posets $P$ with $\left\vert
P\right\vert =N-1$. Now, let $P$ be a poset with $\left\vert P\right\vert =N$.

Let $\left(  v_{1},v_{2},...,v_{N}\right)  $ and $\left(  w_{1},w_{2}%
,...,w_{N}\right)  $ be two elements of $\mathcal{L}\left(  P\right)  $. We
are going to prove that $\left(  v_{1},v_{2},...,v_{N}\right)  \sim\left(
w_{1},w_{2},...,w_{N}\right)  $.

Let $u=v_{N}$. Then, $u$ is a maximal element of $P$ (since it comes last in
the linear extension $\left(  v_{1},v_{2},...,v_{N}\right)  $). Let $i$ be the
index satisfying $w_{i}=u$.

Consider the poset $P\setminus\left\{  u\right\}  $. This poset has size
$\left\vert P\setminus\left\{  u\right\}  \right\vert =\underbrace{\left\vert
P\right\vert }_{=N}-1=N-1$. Define a relation $\sim$ on $\mathcal{L}\left(
P\setminus\left\{  u\right\}  \right)  $ in the same way as the relation
$\sim$ on $\mathcal{L}\left(  P\right)  $ was defined. Recall that $u$ is a
maximal element of $P$. Hence,%
\begin{equation}
\left(
\begin{array}
[c]{c}%
\text{if }\left(  a_{1},a_{2},...,a_{N-1}\right)  \text{ is a linear extension
of }P\setminus\left\{  u\right\}  \text{, then}\\
\left(  a_{1},a_{2},...,a_{N-1},u\right)  \text{ is a linear extension of }P
\end{array}
\right)  . \label{pf.linext.transitive.0.short}%
\end{equation}
Moreover, just by recalling how the relations $\sim$ were defined, we can
easily see that%
\begin{equation}
\left(
\begin{array}
[c]{c}%
\text{if two linear extensions }\left(  a_{1},a_{2},...,a_{N-1}\right)  \text{
and }\left(  b_{1},b_{2},...,b_{N-1}\right)  \text{ of }P\setminus\left\{
u\right\} \\
\text{satisfy }\left(  a_{1},a_{2},...,a_{N-1}\right)  \sim\left(  b_{1}%
,b_{2},...,b_{N-1}\right)  \text{ in }\mathcal{L}\left(  P\setminus\left\{
u\right\}  \right)  \text{, then}\\
\left(  a_{1},a_{2},...,a_{N-1},u\right)  \text{ and }\left(  b_{1}%
,b_{2},...,b_{N-1},u\right)  \text{ are two linear extensions}\\
\text{of }P\text{ satisfying }\left(  a_{1},a_{2},...,a_{N-1},u\right)
\sim\left(  b_{1},b_{2},...,b_{N-1},u\right)  \text{ in }\mathcal{L}\left(
P\right)
\end{array}
\right)  \label{pf.linext.transitive.1.short}%
\end{equation}
(here, the fact that $\left(  a_{1},a_{2},...,a_{N-1},u\right)  $ and $\left(
b_{1},b_{2},...,b_{N-1},u\right)  $ are linear extensions of $P$ follows from
(\ref{pf.linext.transitive.0.short})).

It is rather clear that $\left(  v_{1},v_{2},...,v_{N-1}\right)  $ and
$\left(  w_{1},w_{2},...,w_{i-1},w_{i+1},w_{i+2},...,w_{N}\right)  $ are two
linear extensions of the poset $P\setminus\left\{  u\right\}  $ (since they
are obtained from the linear extensions $\left(  v_{1},v_{2},...,v_{N}\right)
$ and $\left(  w_{1},w_{2},...,w_{N}\right)  $ of $P$ by removing $u$). Since
we can apply Proposition \ref{prop.linext.transitive} to this poset
$P\setminus\left\{  u\right\}  $ in lieu of $P$ (by the induction hypothesis,
since $\left\vert P\setminus\left\{  u\right\}  \right\vert =N-1$), we thus
see that%
\[
\left(  v_{1},v_{2},...,v_{N-1}\right)  \sim\left(  w_{1},w_{2},...,w_{i-1}%
,w_{i+1},w_{i+2},...,w_{N}\right)
\]
in $\mathcal{L}\left(  P\setminus\left\{  u\right\}  \right)  $. By
(\ref{pf.linext.transitive.1.short}), this yields that $\left(  v_{1}%
,v_{2},...,v_{N-1},u\right)  $ and \newline$\left(  w_{1},w_{2},...,w_{i-1}%
,w_{i+1},w_{i+2},...,w_{N},u\right)  $ are two linear extensions of $P$
satisfying
\[
\left(  v_{1},v_{2},...,v_{N-1},u\right)  \sim\left(  w_{1},w_{2}%
,...,w_{i-1},w_{i+1},w_{i+2},...,w_{N},u\right)
\]
in $\mathcal{L}\left(  P\right)  $.

Now, we know that the tuple $\left(  w_{1},w_{2},...,w_{N}\right)  $ is a
linear extension of $P$. Since
\begin{align*}
&  \left(  w_{1},w_{2},...,w_{N}\right) \\
&  =\left(  w_{1},w_{2},...,w_{i-1},\underbrace{w_{i}}_{=u},w_{i+1}%
,w_{i+2},...,w_{N}\right)  =\left(  w_{1},w_{2},...,w_{i-1},u,w_{i+1}%
,w_{i+2},...,w_{N}\right)  ,
\end{align*}
this rewrites as follows: The tuple $\left(  w_{1},w_{2},...,w_{i-1}%
,u,w_{i+1},w_{i+2},...,w_{N}\right)  $ is a linear extension of $P$. Hence, we
can apply Lemma \ref{lem.linext.transitive} to $k=i-1$, $\ell=N-i$,
$a_{j}=w_{j}$ and $b_{j}=w_{i+j}$. As a result, we see that $\left(
w_{1},w_{2},...,w_{i-1},w_{i+1},w_{i+2},...,w_{N},u\right)  $ is a linear
extension of $P$ satisfying $\left(  w_{1},w_{2},...,w_{i-1},u,w_{i+1}%
,w_{i+2},...,w_{N}\right)  \sim\left(  w_{1},w_{2},...,w_{i-1},w_{i+1}%
,w_{i+2},...,w_{N},u\right)  $. Since the relation $\sim$ is symmetric
(because $\sim$ is an equivalence relation), this yields%
\[
\left(  w_{1},w_{2},...,w_{i-1},w_{i+1},w_{i+2},...,w_{N},u\right)
\sim\left(  w_{1},w_{2},...,w_{i-1},u,w_{i+1},w_{i+2},...,w_{N}\right)  .
\]

Altogether,%
\begin{align*}
\left(  v_{1},v_{2},...,v_{N}\right)   &  =\left(  v_{1},v_{2},...,v_{N-1}%
,\underbrace{v_{N}}_{=u}\right)  =\left(  v_{1},v_{2},...,v_{N-1},u\right) \\
&  \sim\left(  w_{1},w_{2},...,w_{i-1},w_{i+1},w_{i+2},...,w_{N},u\right) \\
&  \sim\left(  w_{1},w_{2},...,w_{i-1},\underbrace{u}_{=w_{i}},w_{i+1}%
,w_{i+2},...,w_{N}\right) \\
&  =\left(  w_{1},w_{2},...,w_{i-1},w_{i},w_{i+1},w_{i+2},...,w_{N}\right)
=\left(  w_{1},w_{2},...,w_{N}\right)  .
\end{align*}
We thus have shown that any two elements $\left(  v_{1},v_{2},...,v_{N}%
\right)  $ and $\left(  w_{1},w_{2},...,w_{N}\right)  $ of $\mathcal{L}\left(
P\right)  $ satisfy $\left(  v_{1},v_{2},...,v_{N}\right)  \sim\left(
w_{1},w_{2},...,w_{N}\right)  $. In other words, Proposition
\ref{prop.linext.transitive} is proven for $\left\vert P\right\vert =N$, so
the induction step is complete, and Proposition \ref{prop.linext.transitive}
is proven.
\end{proof}
\end{vershort}

\begin{verlong}
\begin{proof}
[Proof of Proposition \ref{prop.linext.transitive} (sketched).]We prove
Proposition \ref{prop.linext.transitive} by induction over $\left\vert
P\right\vert $. The induction base $\left\vert P\right\vert =0$ is trivial.
For the induction step, let $N$ be a positive integer. Assume that Proposition
\ref{prop.linext.transitive} is proven for all posets $P$ with $\left\vert
P\right\vert =N-1$. Now, let $P$ be a poset with $\left\vert P\right\vert =N$.

Let $\left(  v_{1},v_{2},...,v_{N}\right)  $ and $\left(  w_{1},w_{2}%
,...,w_{N}\right)  $ be two elements of $\mathcal{L}\left(  P\right)  $. We
are going to prove that $\left(  v_{1},v_{2},...,v_{N}\right)  \sim\left(
w_{1},w_{2},...,w_{N}\right)  $.

Let $u=v_{N}$. Then, $u$ is a maximal element of $P$ (since it comes last in
the linear extension $\left(  v_{1},v_{2},...,v_{N}\right)  $). Let $i$ be the
index satisfying $w_{i}=u$.

Consider the poset $P\setminus\left\{  u\right\}  $. This poset has size
$\left\vert P\setminus\left\{  u\right\}  \right\vert =\underbrace{\left\vert
P\right\vert }_{=N}-1=N-1$. Define a relation $\sim$ on $\mathcal{L}\left(
P\setminus\left\{  u\right\}  \right)  $ in the same way as the relation
$\sim$ on $\mathcal{L}\left(  P\right)  $ was defined. Recall that $u$ is a
maximal element of $P$. Hence,%
\begin{equation}
\left(
\begin{array}
[c]{c}%
\text{if }\left(  a_{1},a_{2},...,a_{N-1}\right)  \text{ is a linear extension
of }P\setminus\left\{  u\right\}  \text{, then}\\
\left(  a_{1},a_{2},...,a_{N-1},u\right)  \text{ is a linear extension of }P
\end{array}
\right)  \label{pf.linext.transitive.0}%
\end{equation}
(by Theorem \ref{thm.linext.induc}). Moreover, just by recalling how the
relations $\sim$ were defined, we can easily see that%
\begin{equation}
\left(
\begin{array}
[c]{c}%
\text{if two linear extensions }\left(  a_{1},a_{2},...,a_{N-1}\right)  \text{
and }\left(  b_{1},b_{2},...,b_{N-1}\right)  \text{ of }P\setminus\left\{
u\right\} \\
\text{satisfy }\left(  a_{1},a_{2},...,a_{N-1}\right)  \sim\left(  b_{1}%
,b_{2},...,b_{N-1}\right)  \text{ in }\mathcal{L}\left(  P\setminus\left\{
u\right\}  \right)  \text{, then}\\
\left(  a_{1},a_{2},...,a_{N-1},u\right)  \text{ and }\left(  b_{1}%
,b_{2},...,b_{N-1},u\right)  \text{ are two linear extensions}\\
\text{of }P\text{ satisfying }\left(  a_{1},a_{2},...,a_{N-1},u\right)
\sim\left(  b_{1},b_{2},...,b_{N-1},u\right)  \text{ in }\mathcal{L}\left(
P\right)
\end{array}
\right)  \label{pf.linext.transitive.1}%
\end{equation}
(here, the fact that $\left(  a_{1},a_{2},...,a_{N-1},u\right)  $ and $\left(
b_{1},b_{2},...,b_{N-1},u\right)  $ are linear extensions of $P$ follows from
(\ref{pf.linext.transitive.0})).

It is rather clear that $\left(  v_{1},v_{2},...,v_{N-1}\right)  $ and
$\left(  w_{1},w_{2},...,w_{i-1},w_{i+1},w_{i+2},...,w_{N}\right)  $ are two
linear extensions of the poset $P\setminus\left\{  u\right\}  $ (since they
are obtained from the linear extensions $\left(  v_{1},v_{2},...,v_{N}\right)
$ and $\left(  w_{1},w_{2},...,w_{N}\right)  $ of $P$ by removing $u$). Since
we can apply Proposition \ref{prop.linext.transitive} to this poset
$P\setminus\left\{  u\right\}  $ in lieu of $P$ (by the induction hypothesis,
since $\left\vert P\setminus\left\{  u\right\}  \right\vert =N-1$), we thus
see that%
\[
\left(  v_{1},v_{2},...,v_{N-1}\right)  \sim\left(  w_{1},w_{2},...,w_{i-1}%
,w_{i+1},w_{i+2},...,w_{N}\right)
\]
in $\mathcal{L}\left(  P\setminus\left\{  u\right\}  \right)  $. By
(\ref{pf.linext.transitive.1}), this yields that $\left(  v_{1},v_{2}%
,...,v_{N-1},u\right)  $ and \newline$\left(  w_{1},w_{2},...,w_{i-1}%
,w_{i+1},w_{i+2},...,w_{N},u\right)  $ are two linear extensions of $P$
satisfying
\[
\left(  v_{1},v_{2},...,v_{N-1},u\right)  \sim\left(  w_{1},w_{2}%
,...,w_{i-1},w_{i+1},w_{i+2},...,w_{N},u\right)
\]
in $\mathcal{L}\left(  P\right)  $.

Now, we know that the tuple $\left(  w_{1},w_{2},...,w_{N}\right)  $ is a
linear extension of $P$. Since
\begin{align*}
&  \left(  w_{1},w_{2},...,w_{N}\right) \\
&  =\left(  w_{1},w_{2},...,w_{i-1},\underbrace{w_{i}}_{=u},w_{i+1}%
,w_{i+2},...,w_{N}\right)  =\left(  w_{1},w_{2},...,w_{i-1},u,w_{i+1}%
,w_{i+2},...,w_{N}\right)  ,
\end{align*}
this rewrites as follows: The tuple $\left(  w_{1},w_{2},...,w_{i-1}%
,u,w_{i+1},w_{i+2},...,w_{N}\right)  $ is a linear extension of $P$. Hence, we
can apply Lemma \ref{lem.linext.transitive} to $k=i-1$, $\ell=N-i$,
$a_{j}=w_{j}$ and $b_{j}=w_{i+j}$. As a result, we obtain that $\left(
w_{1},w_{2},...,w_{i-1},w_{i+1},w_{i+2},...,w_{N},u\right)  $ is a linear
extension of $P$ satisfying \newline$\left(  w_{1},w_{2},...,w_{i-1}%
,u,w_{i+1},w_{i+2},...,w_{N}\right)  \sim\left(  w_{1},w_{2},...,w_{i-1}%
,w_{i+1},w_{i+2},...,w_{N},u\right)  $. Since the relation $\sim$ is symmetric
(because $\sim$ is an equivalence relation), this yields%
\[
\left(  w_{1},w_{2},...,w_{i-1},w_{i+1},w_{i+2},...,w_{N},u\right)
\sim\left(  w_{1},w_{2},...,w_{i-1},u,w_{i+1},w_{i+2},...,w_{N}\right)  .
\]

Altogether,%
\begin{align*}
\left(  v_{1},v_{2},...,v_{N}\right)   &  =\left(  v_{1},v_{2},...,v_{N-1}%
,\underbrace{v_{N}}_{=u}\right)  =\left(  v_{1},v_{2},...,v_{N-1},u\right) \\
&  \sim\left(  w_{1},w_{2},...,w_{i-1},w_{i+1},w_{i+2},...,w_{N},u\right) \\
&  \sim\left(  w_{1},w_{2},...,w_{i-1},\underbrace{u}_{=w_{i}},w_{i+1}%
,w_{i+2},...,w_{N}\right) \\
&  =\left(  w_{1},w_{2},...,w_{i-1},w_{i},w_{i+1},w_{i+2},...,w_{N}\right)
=\left(  w_{1},w_{2},...,w_{N}\right)  .
\end{align*}
We thus have shown that any two elements $\left(  v_{1},v_{2},...,v_{N}%
\right)  $ and $\left(  w_{1},w_{2},...,w_{N}\right)  $ of $\mathcal{L}\left(
P\right)  $ satisfy $\left(  v_{1},v_{2},...,v_{N}\right)  \sim\left(
w_{1},w_{2},...,w_{N}\right)  $. In other words, Proposition
\ref{prop.linext.transitive} is proven for $\left\vert P\right\vert =N$, so
the induction step is complete, and Proposition \ref{prop.linext.transitive}
is proven.
\end{proof}
\end{verlong}

\section{Birational rowmotion}

In this section, we introduce the basic objects whose nature we will
investigate: labellings of a finite poset $P$ (by elements of a field) and a
birational map between them called ``birational rowmotion''. This map
generalizes (in a certain sense) the notion of ordinary rowmotion on the set
$J\left(  P\right)  $ of order ideals of $P$ to the vastly more general
setting of field-valued labellings. We will discuss the technical concerns
raised by the definitions, and provide two examples and an alternative
description of birational rowmotion. A deeper study of birational rowmotion is
deferred to the following sections.

The concepts which we are going to define now go back to \cite{einstein-propp}
and earlier sources, and are often motivated there. The reader should be
warned that the notations used in \cite{einstein-propp} are not identical with
those used in the present paper (not to mention that \cite{einstein-propp} is
working over $\mathbb{R}_{+}$ rather than over fields as we do).

\begin{definition}
\label{def.Phat}Let $P$ be a poset. Then, $\widehat{P}$ will denote the poset
defined as follows: As a set, let $\widehat{P}$ be the disjoint union of the
set $P$ with the two-element set $\left\{  0,1\right\}  $. The
smaller-or-equal relation $\leq$ on $\widehat{P}$ will be given by%
\[
\left(  a\leq b\right)  \Longleftrightarrow\left(  \text{either }\left(  a\in
P\text{ and }b\in P\text{ and }a\leq b\text{ in }P\right)  \text{ or
}a=0\text{ or }b=1\right)
\]
\footnotemark. Here and in the following, we regard the canonical injection of
the set $P$ into the disjoint union $\widehat{P}$ as an inclusion; thus, $P$
becomes a subposet of $\widehat{P}$. In the terminology of Stanley's
\cite[section 3.2]{stanley-ec1}, this poset $\widehat{P}$ is the ordinal sum
$\left\{  0\right\}  \oplus P\oplus\left\{  1\right\}  $.
\end{definition}

\footnotetext{Here and in the following, the expression \textquotedblleft
either/or\textquotedblright\ always has a non-exclusive meaning. (Thus, in
particular, $0\leq1$ in $\widehat{P}$.)}

\begin{vershort}
\begin{convention}
Let $P$ be a finite poset, and let $u$ and $v$ be two elements of $P$. Then,
$u$ and $v$ are also elements of $\widehat{P}$ (since we are regarding $P$ as
a subposet of $\widehat{P}$). Thus, strictly speaking, statements like
\textquotedblleft$u<v$\textquotedblright\ or \textquotedblleft$u\lessdot
v$\textquotedblright\ are ambiguous because it is not clear whether they are
referring to the poset $P$ or to the poset $\widehat{P}$. However, this
ambiguity is irrelevant, because it is easily seen that the truth of each of
the statements \textquotedblleft$u<v$\textquotedblright, \textquotedblleft%
$u\leq v$\textquotedblright, \textquotedblleft$u>v$\textquotedblright,
\textquotedblleft$u\geq v$\textquotedblright, \textquotedblleft$u\lessdot
v$\textquotedblright, \textquotedblleft$u\gtrdot v$\textquotedblright\ and
\textquotedblleft$u$ and $v$ are incomparable\textquotedblright\ is
independent on whether it refers to the poset $P$ or to the poset
$\widehat{P}$. We are going to therefore omit mentioning the poset in these
statements, unless there are other reasons for us to do so.
\end{convention}
\end{vershort}

\begin{verlong}
Before we move on, we shall make some very simple observations about the
covering relation of the poset $\widehat{P}$ for a finite poset $P$:

\begin{proposition}
\label{prop.Phat.cover} Let $P$ be a finite poset. Consider the poset
$\widehat{P}$.

\textbf{(a)} We have $0\lessdot1$ in $\widehat{P}$ if and only if
$P=\varnothing$.

\textbf{(b)} Let $u\in P$. We have $0\lessdot u$ in $\widehat{P}$ if and only
if $u$ is a minimal element of $P$.

\textbf{(c)} Let $u\in P$. We have $u\lessdot1$ in $\widehat{P}$ if and only
if $u$ is a maximal element of $P$.

\textbf{(d)} Let $u\in P$ and $v\in P$. We have $u\lessdot v$ in $\widehat{P}$
if and only if $u\lessdot v$ in $P$.
\end{proposition}

\begin{corollary}
\label{cor.Phat.exist} Let $P$ be a finite poset. Consider the poset
$\widehat{P}$. Let $q\in\widehat{P}$.

\textbf{(a)} If $q\neq0$, then there exists a $u\in\widehat{P}$ such that
$u\lessdot q$ in $\widehat{P}$.

\textbf{(b)} If $q\neq1$, then there exists a $u\in\widehat{P}$ such that
$u\gtrdot q$ in $\widehat{P}$.

\textbf{(c)} If $q\in P$, then there exists a $u\in\widehat{P}$ such that
$u\lessdot q$ in $\widehat{P}$.

\textbf{(d)} If $q\in P$, then there exists a $u\in\widehat{P}$ such that
$u\gtrdot q$ in $\widehat{P}$.
\end{corollary}

\begin{convention}
Let $P$ be a finite poset, and let $u$ and $v$ be two elements of $P$. Then,
$u$ and $v$ are also elements of $\widehat{P}$ (since we are regarding $P$ as
a subposet of $\widehat{P}$). Thus, strictly speaking, statements like
\textquotedblleft$u<v$\textquotedblright\ or \textquotedblleft$u\lessdot
v$\textquotedblright\ are ambiguous because it is not clear whether they are
referring to the poset $P$ or to the poset $\widehat{P}$. However, this
ambiguity is irrelevant, because the statement $\left(  u<v\text{ in
}\widehat{P}\right)  $ is equivalent to the statement $\left(  u<v\text{ in
}P\right)  $ (by the construction of $\widehat{P}$), and because the statement
$\left(  u\lessdot v\text{ in }\widehat{P}\right)  $ is equivalent to the
statement $\left(  u\lessdot v\text{ in }P\right)  $ (by Proposition
\ref{prop.Phat.cover} \textbf{(d)}). This allows us to make statements of the
form \textquotedblleft$u<v$\textquotedblright, \textquotedblleft$u\leq
v$\textquotedblright, \textquotedblleft$u>v$\textquotedblright,
\textquotedblleft$u\geq v$\textquotedblright, \textquotedblleft$u\lessdot
v$\textquotedblright, \textquotedblleft$u\gtrdot v$\textquotedblright\ and
\textquotedblleft$u$ and $v$ are incomparable\textquotedblright\ without
worrying whether they refer to the poset $P$ or to the poset $\widehat{P}$. We
are going to make use of this in the following.
\end{convention}
\end{verlong}

\begin{definition}
\label{def.labelling}Let $P$ be a poset. Let $\mathbb{K}$ be a field. A
$\mathbb{K}$\textit{-labelling of }$P$ will mean a map $f:\widehat{P}%
\rightarrow\mathbb{K}$. Thus, $\mathbb{K}^{\widehat{P}}$ is the set of all
$\mathbb{K}$-labellings of $P$. If $f$ is a $\mathbb{K}$-labelling of $P$ and
$v$ is an element of $\widehat{P}$, then $f\left(  v\right)  $ will be called
the \textit{label of }$f$ \textit{at }$v$.
\end{definition}

\begin{definition}
In the following, whenever we are working with a field $\mathbb{K}$, we are
going to tacitly assume that $\mathbb{K}$ is either infinite or at least can
be enlarged when necessity arises. This assumption is needed in order to
clarify the notions of rational maps and generic elements of algebraic
varieties over $\mathbb{K}$. (We will not require $\mathbb{K}$ to be
algebraically closed.)

We will use the terminology of algebraic varieties and rational maps between
them, although the only algebraic varieties that we will be considering are
products of affine and projective spaces, as well as their open subsets. We
use the punctured arrow $\dashrightarrow$ to signify rational maps (i.e., a
rational map from a variety $U$ to a variety $V$ is called a rational map
$U\dashrightarrow V$). A rational map $U\dashrightarrow V$ is said to be
\textit{dominant} if its image is dense in $V$ (with respect to the Zariski topology).

The words \textquotedblleft generic\textquotedblright\ and \textquotedblleft
almost\textquotedblright\ will always refer to the Zariski topology. For
example, if $U$ is a finite set, then an assertion saying that some statement
holds \textquotedblleft for almost every point $p\in\mathbb{K}^{U}%
$\textquotedblright\ is supposed to mean that there is a Zariski-dense open
subset $D$ of $\mathbb{K}^{U}$ such that this statement holds for every point
$p\in D$. A \textquotedblleft generic\textquotedblright\ point on an algebraic
variety $V$ (for example, this can be a \textquotedblleft generic
matrix\textquotedblright\ when $V$ is a space of matrices, or a
\textquotedblleft generic $\mathbb{K}$-labelling of a poset $P$%
\textquotedblright\ when $V$ is the space of all $\mathbb{K}$-labellings of
$P$) means a point lying in some fixed Zariski-dense open subset $S$ of $V$;
the concrete definition of $S$ can usually be inferred from the context
(often, it will be the subset of $V$ on which everything we want to do with
our point is well-defined), but of course should never depend on the actual
point. (Note that one often has to read the whole proof in order to be able to
tell what this $S$ is. This is similar to the use of the \textquotedblleft for
$\epsilon$ small enough\textquotedblright\ wording in analysis, where it is
often not clear until the end of the proof how small exactly the $\epsilon$
needs to be.) We are sometimes going to abuse notation and say that an
equality holds \textquotedblleft for every point\textquotedblright\ instead of
\textquotedblleft for almost every point\textquotedblright\ when it is really
clear what the $S$ is. (For example, if we say that \textquotedblleft the
equality $\dfrac{x^{3}-y^{3}}{x-y}=x^{2}+xy+y^{2}$ holds for every
$x\in\mathbb{K}$ and $y\in\mathbb{K}$\textquotedblright, it is clear that $S$
has to be the set $\mathbb{K}^{2}\setminus\left\{  \left(  x,y\right)
\in\mathbb{K}^{2}\mid x=y\right\}  $, because the left hand side of the
equality makes no sense when $\left(  x,y\right)  $ is outside of this set.)
\end{definition}

\begin{remark}
Most statements that we make below work not only for fields, but also more
generally for semifields\footnotemark\ such as the semifield $\mathbb{Q}_{+}$
of positive rationals or the tropical semiring. Some (but not all!) statements
actually simplify when the underlying field is replaced by a semifield in
which no two nonzero elements add to zero (because in such cases, e.g., the
denominators in (\ref{def.Tv.def}) cannot become zero unless some labels of
$f$ are $0$). Thus, working with such semifields instead of fields would save
us the trouble of having things defined \textquotedblleft almost
everywhere\textquotedblright. Moreover, applying our results to the tropical
semifield would yield some of the statements about order polytopes made in
\cite{einstein-propp}. Nevertheless, we prefer to work with fields, for the
following reasons:

-- While most of our results can be formulated for semifields, not all of them
can (and sometimes, even when a result holds over semifields, its proof might
not work over semifields). In particular, Proposition \ref{prop.Grasp.GraspR}
makes no sense over semifields, because determinants involve subtraction.
Also, if we were to work in semifields which \textbf{do} contain two nonzero
elements summing up to zero, then we would still have the issue of zero
denominators, but we are not aware of a theoretical framework in the spirit of
Zariski topology for fields to reassure us in this case that these issues are negligible.

-- If an identity between subtraction-free rational expressions (such as
$\dfrac{x^{3}+y^{3}}{x+y}+3xy=\left(  x+y\right)  ^{2}$) holds over every
field (as long as the denominators involved are nonzero), then it must hold
over every semifield as well (again as long as the denominators involved are
nonzero), even if the identity has only been proven with the help of
subtraction (e.g., a proof of $\dfrac{x^{3}+y^{3}}{x+y}+3xy=\left(
x+y\right)  ^{2}$ over a field can begin by simplifying $\dfrac{x^{3}+y^{3}%
}{x+y}$ to $x^{2}-xy+y^{2}$, a technique not available over a semifield). This
is simply because every true identity between subtraction-free rational
expressions can be verified by multiplying by a common denominator (an
operation which does not introduce any subtractions) and comparing
coefficients. Since our main results (such as Theorem
\ref{thm.rect.antip.general}, or the $p+q\mid\operatorname*{ord}\left(
R_{\operatorname*{Rect}\left(  p,q\right)  }\right)  $ part of Theorem
\ref{thm.rect.ord}) can be construed as identities between subtraction-free
rational expressions, this yields that all these results hold over any
semifield (provided the denominators are nonzero) if they hold over every
field. So we are not losing any generality by restricting ourselves to
considering only fields.
\end{remark}

\footnotetext{The word \textquotedblleft\textit{semifield}\textquotedblright%
\ here means a commutative semiring in which each element other than $0$ has a
multiplicative inverse. (In contrast to other authors' conventions, our
semifields do have zeroes.) A \textit{semiring} is defined as a set with two
binary operations called \textquotedblleft addition\textquotedblright\ and
\textquotedblleft multiplication\textquotedblright\ and two elements $0$ and
$1$ which satisfies all axioms of a ring (in particular, it must be
associative and satisfy $0\cdot a=a\cdot0=0$ and $1\cdot a=a\cdot1=a$ for all
$a$) except for having additive inverses.}

\begin{definition}
\label{def.Tv}Let $P$ be a finite poset. Let $\mathbb{K}$ be a field. Let
$v\in P$. We define a rational map $T_{v}:\mathbb{K}^{\widehat{P}%
}\dashrightarrow\mathbb{K}^{\widehat{P}}$ by%
\begin{equation}
\left(  T_{v}f\right)  \left(  w\right)  =\left\{
\begin{array}
[c]{l}%
f\left(  w\right)  ,\ \ \ \ \ \ \ \ \ \ \text{if }w\neq v;\\
\dfrac{1}{f\left(  v\right)  }\cdot\dfrac{\sum\limits_{\substack{u\in
\widehat{P};\\u\lessdot v}}f\left(  u\right)  }{\sum\limits_{\substack{u\in
\widehat{P};\\u\gtrdot v}}\dfrac{1}{f\left(  u\right)  }}%
,\ \ \ \ \ \ \ \ \ \ \text{if }w=v
\end{array}
\right.  \ \ \ \ \ \ \ \ \ \ \text{for all }w\in\widehat{P} \label{def.Tv.def}%
\end{equation}
for all $f\in\mathbb{K}^{\widehat{P}}$. Note that this rational map $T_{v}$ is
well-defined, because the right-hand side of (\ref{def.Tv.def}) is
well-defined on a Zariski-dense open subset of $\mathbb{K}^{\widehat{P}}$.
(This follows from the fact that for every $v\in P$, there is at least one
$u\in\widehat{P}$ such that $u\gtrdot v$\ \ \ \ \footnotemark.)

This rational map $T_{v}$ is called the $v$\textit{-toggle}.
\end{definition}

\begin{vershort}
\footnotetext{Indeed, either there is at least one $u\in P$ such that
$u\gtrdot v$ in $P$ (and therefore also $u\gtrdot v$ in $\widehat{P}$), or
else $v$ is maximal in $P$ and then we have $1\gtrdot v$ in $\widehat{P}$.}
\end{vershort}

\begin{verlong}
\footnotetext{This follows from Corollary \ref{cor.Phat.exist} \textbf{(d)}
(applied to $q=v$).}
\end{verlong}

The map $T_{v}$ that we have just introduced (although defined over the
semifield $\mathbb{R}_{+}$ instead of our field $\mathbb{K}$) is called a
\textquotedblleft birational toggle operation\textquotedblright\ in
\cite{einstein-propp} (where it is denoted by $\phi_{i}$ with $i$ being a
number indexing the elements $v$ of $P$; however, the same notation is used
for the \textquotedblleft tropicalized\textquotedblright\ version of $T_{v}$).
As is clear from its definition, it only changes the label at the element $v$.

Note also the following almost trivial fact:

\begin{proposition}
\label{prop.Tv.invo}Let $P$ be a finite poset. Let $\mathbb{K}$ be a field.
Let $v\in P$. Then, the rational map $T_{v}$ is an involution, i.e., the map
$T_{v}^{2}$ is well-defined on a Zariski-dense open subset of $\mathbb{K}%
^{\widehat{P}}$ and satisfies $T_{v}^{2}=\operatorname*{id}$ on this subset.
\end{proposition}

We are calling this \textquotedblleft almost trivial\textquotedblright%
\ because one subtlety is easily overlooked: We have to check that the map
$T_{v}^{2}$ is well-defined on a Zariski-dense open subset of $\mathbb{K}%
^{\widehat{P}}$; this requires observing that for every $v\in P$, there exists
at least one $u\in\widehat{P}$ such that $u\lessdot v$.

\begin{verlong}
\begin{proof}
[Proof of Proposition \ref{prop.Tv.invo} (sketched).]Recall the following
simple fact: If $A$ is an affine variety, then a Zariski-dense open subset of
a Zariski-dense open subset of $A$ always is a Zariski-dense open subset of
$A$. Applying this to $A=\mathbb{K}^{\widehat{P}}$, we conclude that a
Zariski-dense open subset of a Zariski-dense open subset of $\mathbb{K}%
^{\widehat{P}}$ always is a Zariski-dense open subset of $\mathbb{K}%
^{\widehat{P}}$.

Define a subset $\mathfrak{Z}$ of $\mathbb{K}^{\widehat{P}}$ by%
\[
\mathfrak{Z}=\left\{  f\in\mathbb{K}^{\widehat{P}}\ \mid\ f\left(  v\right)
\neq0\right\}  .
\]
Clearly, $\mathfrak{Z}$ is a Zariski-dense open subset of $\mathbb{K}%
^{\widehat{P}}$.

Define a subset $\mathfrak{W}$ of $\mathfrak{Z}$ by%
\[
\mathfrak{W}=\left\{  f\in\mathfrak{Z}\ \mid\ f\left(  u\right)  \neq0\text{
for every }u\in\widehat{P}\text{ satisfying }u\gtrdot v\right\}  .
\]
Clearly, $\mathfrak{W}$ is a Zariski-dense open subset of $\mathfrak{Z}$
(since there are only finitely many $u\in\widehat{P}$ satisfying $u\gtrdot
v$). Since $\mathfrak{Z}$ (in turn) is a Zariski-dense open subset of
$\mathbb{K}^{\widehat{P}}$, this yields that $\mathfrak{W}$ is a Zariski-dense
open subset of $\mathbb{K}^{\widehat{P}}$ (since a Zariski-dense open subset
of a Zariski-dense open subset of $\mathbb{K}^{\widehat{P}}$ always is a
Zariski-dense open subset of $\mathbb{K}^{\widehat{P}}$).

For every $f\in\mathfrak{W}$, we have $f\left(  u\right)  \neq0$ for every
$u\in\widehat{P}$ satisfying $u\gtrdot v$ (by the definition of $\mathfrak{W}%
$). Hence, for every $f\in\mathfrak{W}$, the term $\dfrac{1}{f\left(
u\right)  }$ is well-defined for every $u\in\widehat{P}$ satisfying $u\gtrdot
v$. Hence, for every $f\in\mathfrak{W}$, the term $\sum\limits_{\substack{u\in
\widehat{P};\\u\gtrdot v}}\dfrac{1}{f\left(  u\right)  }$ is well-defined.
Therefore, we can define a subset $\mathfrak{V}$ of $\mathfrak{W}$ by%
\[
\mathfrak{V}=\left\{  f\in\mathfrak{W}\ \mid\ \sum\limits_{\substack{u\in
\widehat{P};\\u\gtrdot v}}\dfrac{1}{f\left(  u\right)  }\neq0\right\}  .
\]
Consider this subset $\mathfrak{V}$. Clearly, $\mathfrak{V}$ is an open subset
of $\mathfrak{W}$ (in the Zariski topology). Moreover, the sum $\sum
\limits_{\substack{u\in\widehat{P};\\u\gtrdot v}}\dfrac{1}{f\left(  u\right)
}$ is not an empty sum for $f\in\mathfrak{W}$ (since there exists a
$u\in\widehat{P}$ such that $u\gtrdot v$ (by Corollary \ref{cor.Phat.exist}
\textbf{(d)}, applied to $q=v$)). Hence, $\sum\limits_{\substack{u\in
\widehat{P};\\u\gtrdot v}}\dfrac{1}{f\left(  u\right)  }$ is not identically
$0$ as a rational function in $f\in\mathbb{K}^{\widehat{P}}$. As a
consequence, $\sum\limits_{\substack{u\in\widehat{P};\\u\gtrdot v}}\dfrac
{1}{f\left(  u\right)  }$ is not identically $0$ as a rational function in
$f\in\mathfrak{W}$ (because $\mathfrak{W}$ is Zariski-dense and open in
$\mathbb{K}^{\widehat{P}}$). Thus, $\mathfrak{V}$ is a nonempty subset of
$\mathfrak{W}$. Since $\mathfrak{V}$ is a nonempty open subset of
$\mathfrak{W}$, we conclude that $\mathfrak{V}$ is a Zariski-dense open subset
of $\mathfrak{W}$ (because any nonempty open subset of an algebraic variety is
a Zariski-dense open subset of this variety). Since $\mathfrak{W}$ (in turn)
is a Zariski-dense open subset of $\mathbb{K}^{\widehat{P}}$, this yields that
$\mathfrak{V}$ is a Zariski-dense open subset of $\mathbb{K}^{\widehat{P}}$
(since a Zariski-dense open subset of a Zariski-dense open subset of
$\mathbb{K}^{\widehat{P}}$ always is a Zariski-dense open subset of
$\mathbb{K}^{\widehat{P}}$).

Next, we define a subset $\mathfrak{U}$ of $\mathfrak{V}$ by%
\[
\mathfrak{U}=\left\{  f\in\mathfrak{V}\ \mid\ \sum\limits_{\substack{u\in
\widehat{P};\\u\lessdot v}}f\left(  u\right)  \neq0\right\}  .
\]
Consider this subset $\mathfrak{U}$. Clearly, $\mathfrak{U}$ is an open subset
of $\mathfrak{V}$ (in the Zariski topology). Moreover, the sum $\sum
\limits_{\substack{u\in\widehat{P};\\u\lessdot v}}f\left(  u\right)  $ is not
an empty sum for $f\in\mathfrak{W}$ (since there exists a $u\in\widehat{P}$
such that $u\gtrdot v$ (by Corollary \ref{cor.Phat.exist} \textbf{(c)},
applied to $q=v$)). Hence, $\sum\limits_{\substack{u\in\widehat{P};\\u\lessdot
v}}f\left(  u\right)  $ is not identically $0$ as a rational function in
$f\in\mathbb{K}^{\widehat{P}}$. As a consequence,$\sum\limits_{\substack{u\in
\widehat{P};\\u\lessdot v}}f\left(  u\right)  $ is not identically $0$ as a
rational function in $f\in\mathfrak{V}$ (because $\mathfrak{V}$ is
Zariski-dense and open in $\mathbb{K}^{\widehat{P}}$). Thus, $\mathfrak{U}$ is
a nonempty subset of $\mathfrak{V}$. Since $\mathfrak{U}$ is a nonempty open
subset of $\mathfrak{V}$, we conclude that $\mathfrak{U}$ is a Zariski-dense
open subset of $\mathfrak{V}$ (because any nonempty open subset of an
algebraic variety is a Zariski-dense open subset of this variety). Since
$\mathfrak{V}$ (in turn) is a Zariski-dense open subset of $\mathbb{K}%
^{\widehat{P}}$, this yields that $\mathfrak{U}$ is a Zariski-dense open
subset of $\mathbb{K}^{\widehat{P}}$ (since a Zariski-dense open subset of a
Zariski-dense open subset of $\mathbb{K}^{\widehat{P}}$ always is a
Zariski-dense open subset of $\mathbb{K}^{\widehat{P}}$).

Now, let $f\in\mathfrak{U}$. We are going to prove that $T_{v}^{2}f$ is
well-defined and satisfies $T_{v}^{2}f=f$.

We have $f\in\mathfrak{U}=\left\{  f\in\mathfrak{V}\ \mid\ \sum
\limits_{\substack{u\in\widehat{P};\\u\lessdot v}}f\left(  u\right)
\neq0\right\}  $. In other words, $f$ is an element of $\mathfrak{V}$ and
satisfies $\sum\limits_{\substack{u\in\widehat{P};\\u\lessdot v}}f\left(
u\right)  \neq0$.

We have $f\in\mathfrak{V}=\left\{  f\in\mathfrak{W}\ \mid\ \sum
\limits_{\substack{u\in\widehat{P};\\u\gtrdot v}}\dfrac{1}{f\left(  u\right)
}\neq0\right\}  $. In other words, $f$ is an element of $\mathfrak{W}$ and
satisfies $\sum\limits_{\substack{u\in\widehat{P};\\u\gtrdot v}}\dfrac
{1}{f\left(  u\right)  }\neq0$.

We have $f\in\mathfrak{W}=\left\{  f\in\mathfrak{Z}\ \mid\ f\left(  u\right)
\neq0\text{ for every }u\in\widehat{P}\text{ satisfying }u\gtrdot v\right\}
$. In other words, $f$ is an element of $\mathfrak{Z}$ and satisfies $\left(
f\left(  u\right)  \neq0\text{ for every }u\in\widehat{P}\text{ satisfying
}u\gtrdot v\right)  $.

We have $f\in\mathfrak{Z}=\left\{  f\in\mathbb{K}^{\widehat{P}}\ \mid
\ f\left(  v\right)  \neq0\right\}  $. In other words, $f$ is an element of
$\mathbb{K}^{\widehat{P}}$ and satisfies $f\left(  v\right)  \neq0$.

By the definition of $T_{v}f$, we have%
\[
\left(  T_{v}f\right)  \left(  w\right)  =\left\{
\begin{array}
[c]{l}%
f\left(  w\right)  ,\ \ \ \ \ \ \ \ \ \ \text{if }w\neq v;\\
\dfrac{1}{f\left(  v\right)  }\cdot\dfrac{\sum\limits_{\substack{u\in
\widehat{P};\\u\lessdot v}}f\left(  u\right)  }{\sum\limits_{\substack{u\in
\widehat{P};\\u\gtrdot v}}\dfrac{1}{f\left(  u\right)  }}%
,\ \ \ \ \ \ \ \ \ \ \text{if }w=v
\end{array}
\right.  \ \ \ \ \ \ \ \ \ \ \text{for every }w\in\widehat{P}.
\]
None of the denominators on the right hand side of this equation vanish
(because we have $f\left(  v\right)  \neq0$ and $\left(  f\left(  u\right)
\neq0\text{ for every }u\in\widehat{P}\text{ satisfying }u\gtrdot v\right)  $
and $\sum\limits_{\substack{u\in\widehat{P};\\u\gtrdot v}}\dfrac{1}{f\left(
u\right)  }\neq0$). Therefore, $\left(  T_{v}f\right)  \left(  w\right)  $ is
well-defined for every $w\in\widehat{P}$. In other words, $T_{v}f$ is
well-defined. Moreover, (\ref{def.Tv.def}) (applied to $w=v$) yields%
\begin{align}
\left(  T_{v}f\right)  \left(  v\right)   &  =\left\{
\begin{array}
[c]{l}%
f\left(  v\right)  ,\ \ \ \ \ \ \ \ \ \ \text{if }v\neq v;\\
\dfrac{1}{f\left(  v\right)  }\cdot\dfrac{\sum\limits_{\substack{u\in
\widehat{P};\\u\lessdot v}}f\left(  u\right)  }{\sum\limits_{\substack{u\in
\widehat{P};\\u\gtrdot v}}\dfrac{1}{f\left(  u\right)  }}%
,\ \ \ \ \ \ \ \ \ \ \text{if }v=v
\end{array}
\right. \nonumber\\
&  =\dfrac{1}{f\left(  v\right)  }\cdot\dfrac{\sum\limits_{\substack{u\in
\widehat{P};\\u\lessdot v}}f\left(  u\right)  }{\sum\limits_{\substack{u\in
\widehat{P};\\u\gtrdot v}}\dfrac{1}{f\left(  u\right)  }}%
\ \ \ \ \ \ \ \ \ \ \left(  \text{since }v=v\right) \label{pf.Tv.invo.0}\\
&  \neq0\ \ \ \ \ \ \ \ \ \ \left(  \text{since }\sum\limits_{\substack{u\in
\widehat{P};\\u\lessdot v}}f\left(  u\right)  \neq0\right)  .\nonumber
\end{align}
On the other hand, every $w\in\widehat{P}$ satisfying $w\neq v$ satisfies
\begin{equation}
\left(  T_{v}f\right)  \left(  w\right)  =f\left(  w\right)  .
\label{pf.Tv.invo.1}%
\end{equation}
\footnote{\textit{Proof of (\ref{pf.Tv.invo.1}):} Let $w\in\widehat{P}$ be
such that $w\neq v$. Then, (\ref{def.Tv.def}) yields%
\[
\left(  T_{v}f\right)  \left(  w\right)  =\left\{
\begin{array}
[c]{l}%
f\left(  w\right)  ,\ \ \ \ \ \ \ \ \ \ \text{if }w\neq v;\\
\dfrac{1}{f\left(  v\right)  }\cdot\dfrac{\sum\limits_{\substack{u\in
\widehat{P};\\u\lessdot v}}f\left(  u\right)  }{\sum\limits_{\substack{u\in
\widehat{P};\\u\gtrdot v}}\dfrac{1}{f\left(  u\right)  }}%
,\ \ \ \ \ \ \ \ \ \ \text{if }w=v
\end{array}
\right.  =f\left(  w\right)  \ \ \ \ \ \ \ \ \ \ \left(  \text{since }w\neq
v\right)  .
\]
This proves (\ref{pf.Tv.invo.1}).} Hence,%
\begin{equation}
\left(  T_{v}f\right)  \left(  u\right)  \neq0\text{ for every }%
u\in\widehat{P}\text{ satisfying }u\gtrdot v. \label{pf.Tv.invo.2}%
\end{equation}
\footnote{\textit{Proof of (\ref{pf.Tv.invo.2}):} Let $u\in\widehat{P}$ be
such that $u\gtrdot v$. We have $u>v$ in $\widehat{P}$ (since $u\gtrdot v$),
and thus $u\neq v$. Hence, (\ref{pf.Tv.invo.1}) (applied to $w=u$) yields
$\left(  T_{v}f\right)  \left(  u\right)  =f\left(  u\right)  \neq0$ (since we
know that $\left(  f\left(  u\right)  \neq0\text{ for every }u\in
\widehat{P}\text{ satisfying }u\gtrdot v\right)  $). This proves
(\ref{pf.Tv.invo.2}).} Also,%
\begin{equation}
\sum\limits_{\substack{u\in\widehat{P};\\u\gtrdot v}}\dfrac{1}{\left(
T_{v}f\right)  \left(  u\right)  }=\sum\limits_{\substack{u\in\widehat{P}%
;\\u\gtrdot v}}\dfrac{1}{f\left(  u\right)  } \label{pf.Tv.invo.3}%
\end{equation}
\footnote{\textit{Proof of (\ref{pf.Tv.invo.3}):} Every $u\in\widehat{P}$ such
that $u\gtrdot v$ must satisfy $u>v$ in $\widehat{P}$ (since it satisfies
$u\gtrdot v$). Hence, every $u\in\widehat{P}$ such that $u\gtrdot v$ must
satisfy $u\neq v$. Therefore, every $u\in\widehat{P}$ such that $u\gtrdot v$
must satisfy $\left(  T_{v}f\right)  \left(  u\right)  =f\left(  u\right)  $
(by (\ref{pf.Tv.invo.1}) (applied to $w=u$)). Thus,%
\[
\sum\limits_{\substack{u\in\widehat{P};\\u\gtrdot v}}\dfrac{1}{\left(
T_{v}f\right)  \left(  u\right)  }=\sum\limits_{\substack{u\in\widehat{P}%
;\\u\gtrdot v}}\dfrac{1}{f\left(  u\right)  },
\]
so that (\ref{pf.Tv.invo.3}) is proven.}. Finally,%
\begin{equation}
\sum\limits_{\substack{u\in\widehat{P};\\u\lessdot v}}\left(  T_{v}f\right)
\left(  u\right)  =\sum\limits_{\substack{u\in\widehat{P};\\u\lessdot
v}}f\left(  u\right)  \label{pf.Tv.invo.4}%
\end{equation}
\footnote{\textit{Proof of (\ref{pf.Tv.invo.4}):} Every $u\in\widehat{P}$ such
that $u\lessdot v$ must satisfy $u<v$ in $\widehat{P}$ (since it satisfies
$u\lessdot v$). Hence, every $u\in\widehat{P}$ such that $u\lessdot v$ must
satisfy $u\neq v$. Therefore, every $u\in\widehat{P}$ such that $u\lessdot v$
must satisfy $\left(  T_{v}f\right)  \left(  u\right)  =f\left(  u\right)  $
(by (\ref{pf.Tv.invo.1}) (applied to $w=u$)). Thus,%
\[
\sum\limits_{\substack{u\in\widehat{P};\\u\lessdot v}}\left(  T_{v}f\right)
\left(  u\right)  =\sum\limits_{\substack{u\in\widehat{P};\\u\lessdot
v}}f\left(  u\right)  ,
\]
so that (\ref{pf.Tv.invo.4}) is proven.}. Now, the definition of $T_{v}\left(
T_{v}f\right)  $ yields%
\[
\left(  T_{v}\left(  T_{v}f\right)  \right)  \left(  w\right)  =\left\{
\begin{array}
[c]{l}%
\left(  T_{v}f\right)  \left(  w\right)  ,\ \ \ \ \ \ \ \ \ \ \text{if }w\neq
v;\\
\dfrac{1}{\left(  T_{v}f\right)  \left(  v\right)  }\cdot\dfrac{\sum
\limits_{\substack{u\in\widehat{P};\\u\lessdot v}}\left(  T_{v}f\right)
\left(  u\right)  }{\sum\limits_{\substack{u\in\widehat{P};\\u\gtrdot
v}}\dfrac{1}{\left(  T_{v}f\right)  \left(  u\right)  }}%
,\ \ \ \ \ \ \ \ \ \ \text{if }w=v
\end{array}
\right.  \ \ \ \ \ \ \ \ \ \ \text{for all }w\in\widehat{P}.
\]
None of the denominators on the right hand side of this equation vanish
(because we have $\left(  T_{v}f\right)  \left(  v\right)  \neq0$ and $\left(
\left(  T_{v}f\right)  \left(  u\right)  \neq0\text{ for every }%
u\in\widehat{P}\text{ satisfying }u\gtrdot v\right)  $ and $\sum
\limits_{\substack{u\in\widehat{P};\\u\gtrdot v}}\dfrac{1}{\left(
T_{v}f\right)  \left(  u\right)  }=\sum\limits_{\substack{u\in\widehat{P}%
;\\u\gtrdot v}}\dfrac{1}{f\left(  u\right)  }\neq0$). Therefore, $\left(
T_{v}\left(  T_{v}f\right)  \right)  \left(  w\right)  $ is well-defined for
every $w\in\widehat{P}$. In other words, $T_{v}\left(  T_{v}f\right)  $ is
well-defined. In other words, $T_{v}^{2}f$ is well-defined (since
$T_{v}\left(  T_{v}f\right)  =\underbrace{\left(  T_{v}\circ T_{v}\right)
}_{=T_{v}^{2}}f=T_{v}^{2}f$). Moreover, (\ref{def.Tv.def}) (applied to $v$ and
$T_{v}f$ instead of $w$ and $f$) yields%
\begin{align*}
\left(  T_{v}\left(  T_{v}f\right)  \right)  \left(  v\right)   &  =\left\{
\begin{array}
[c]{l}%
\left(  T_{v}f\right)  \left(  w\right)  ,\ \ \ \ \ \ \ \ \ \ \text{if }v\neq
v;\\
\dfrac{1}{\left(  T_{v}f\right)  \left(  v\right)  }\cdot\dfrac{\sum
\limits_{\substack{u\in\widehat{P};\\u\lessdot v}}\left(  T_{v}f\right)
\left(  u\right)  }{\sum\limits_{\substack{u\in\widehat{P};\\u\gtrdot
v}}\dfrac{1}{\left(  T_{v}f\right)  \left(  u\right)  }}%
,\ \ \ \ \ \ \ \ \ \ \text{if }v=v
\end{array}
\right. \\
&  =\dfrac{1}{\left(  T_{v}f\right)  \left(  v\right)  }\cdot\dfrac
{\sum\limits_{\substack{u\in\widehat{P};\\u\lessdot v}}\left(  T_{v}f\right)
\left(  u\right)  }{\sum\limits_{\substack{u\in\widehat{P};\\u\gtrdot
v}}\dfrac{1}{\left(  T_{v}f\right)  \left(  u\right)  }}%
\ \ \ \ \ \ \ \ \ \ \left(  \text{since }v=v\right) \\
&  =\underbrace{\dfrac{1}{\left(  T_{v}f\right)  \left(  v\right)  }%
}_{\substack{=\dfrac{1}{\dfrac{1}{f\left(  v\right)  }\cdot\dfrac
{\sum\limits_{\substack{u\in\widehat{P};\\u\lessdot v}}f\left(  u\right)
}{\sum\limits_{\substack{u\in\widehat{P};\\u\gtrdot v}}\dfrac{1}{f\left(
u\right)  }}}\\\text{(by (\ref{pf.Tv.invo.0}))}}}\cdot\left(  \underbrace{\sum
\limits_{\substack{u\in\widehat{P};\\u\gtrdot v}}\dfrac{1}{\left(
T_{v}f\right)  \left(  u\right)  }}_{\substack{=\sum\limits_{\substack{u\in
\widehat{P};\\u\gtrdot v}}\dfrac{1}{f\left(  u\right)  }\\\text{(by
(\ref{pf.Tv.invo.3}))}}}\right)  ^{-1}\cdot\left(  \underbrace{\sum
\limits_{\substack{u\in\widehat{P};\\u\lessdot v}}\left(  T_{v}f\right)
\left(  u\right)  }_{\substack{=\sum\limits_{\substack{u\in\widehat{P}%
;\\u\lessdot v}}f\left(  u\right)  \\\text{(by (\ref{pf.Tv.invo.4}))}}}\right)
\\
&  =\dfrac{1}{\dfrac{1}{f\left(  v\right)  }\cdot\dfrac{\sum
\limits_{\substack{u\in\widehat{P};\\u\lessdot v}}f\left(  u\right)  }%
{\sum\limits_{\substack{u\in\widehat{P};\\u\gtrdot v}}\dfrac{1}{f\left(
u\right)  }}}\cdot\left(  \sum\limits_{\substack{u\in\widehat{P};\\u\gtrdot
v}}\dfrac{1}{f\left(  u\right)  }\right)  ^{-1}\cdot\left(  \sum
\limits_{\substack{u\in\widehat{P};\\u\lessdot v}}f\left(  u\right)  \right)
=f\left(  v\right)  .
\end{align*}
Now, every $w\in\widehat{P}$ satisfies%
\begin{equation}
\left(  T_{v}\left(  T_{v}f\right)  \right)  \left(  w\right)  =f\left(
w\right)  . \label{pf.Tv.invo.6}%
\end{equation}
\footnote{\textit{Proof of (\ref{pf.Tv.invo.6}):} Let $w\in\widehat{P}$. We
are going to prove that $\left(  T_{v}\left(  T_{v}f\right)  \right)  \left(
w\right)  =f\left(  w\right)  $.
\par
If $w=v$, then this follows from $\left(  T_{v}\left(  T_{v}f\right)  \right)
\left(  \underbrace{w}_{=v}\right)  =\left(  T_{v}\left(  T_{v}f\right)
\right)  \left(  v\right)  =f\left(  \underbrace{v}_{=w}\right)  =f\left(
w\right)  $. Hence, for the rest of the proof of $\left(  T_{v}\left(
T_{v}f\right)  \right)  \left(  w\right)  =f\left(  w\right)  $, we can WLOG
assume that we don't have $w=v$. Assume this.
\par
We don't have $w=v$. Thus, we have $w\neq v$. Now, (\ref{def.Tv.def}) (applied
to $T_{v}f$ instead of $f$) yields%
\begin{align*}
\left(  T_{v}\left(  T_{v}f\right)  \right)  \left(  w\right)   &  =\left\{
\begin{array}
[c]{l}%
\left(  T_{v}f\right)  \left(  w\right)  ,\ \ \ \ \ \ \ \ \ \ \text{if }w\neq
v;\\
\dfrac{1}{\left(  T_{v}f\right)  \left(  v\right)  }\cdot\dfrac{\sum
\limits_{\substack{u\in\widehat{P};\\u\lessdot v}}\left(  T_{v}f\right)
\left(  u\right)  }{\sum\limits_{\substack{u\in\widehat{P};\\u\gtrdot
v}}\dfrac{1}{\left(  T_{v}f\right)  \left(  u\right)  }}%
,\ \ \ \ \ \ \ \ \ \ \text{if }w=v
\end{array}
\right. \\
&  =\left(  T_{v}f\right)  \left(  w\right)  \ \ \ \ \ \ \ \ \ \ \left(
\text{since }w\neq v\right) \\
&  =f\left(  w\right)  \ \ \ \ \ \ \ \ \ \ \left(  \text{by
(\ref{pf.Tv.invo.1})}\right)  .
\end{align*}
This proves $\left(  T_{v}\left(  T_{v}f\right)  \right)  \left(  w\right)
=f\left(  w\right)  $. Thus, (\ref{pf.Tv.invo.6}) is proven.} In other words,
$T_{v}\left(  T_{v}f\right)  =f$. Now, $\underbrace{T_{v}^{2}}_{=T_{v}\circ
T_{v}}f=\left(  T_{v}\circ T_{v}\right)  f=T_{v}\left(  T_{v}f\right)  =f$.

Now, forget that we fixed $f$. We thus have shown that for every
$f\in\mathfrak{U}$, the $\mathbb{K}$-labelling $T_{v}^{2}f$ is well-defined
and satisfies $T_{v}^{2}f=f$. In other words, the map $T_{v}^{2}$ is
well-defined on the subset $\mathfrak{U}$ of $\mathbb{K}^{\widehat{P}}$ and
satisfies $T_{v}^{2}=\operatorname*{id}$ on this subset. Since $\mathfrak{U}$
is a Zariski-dense open subset of $\mathbb{K}^{\widehat{P}}$, this yields that
the map $T_{v}^{2}$ is well-defined on a Zariski-dense open subset of
$\mathbb{K}^{\widehat{P}}$ and satisfies $T_{v}^{2}=\operatorname*{id}$ on
this subset. This proves Proposition \ref{prop.Tv.invo}.
\end{proof}
\end{verlong}

Proposition \ref{prop.Tv.invo} yields the following:

\begin{corollary}
\label{cor.Tv.dom}Let $P$ be a finite poset. Let $\mathbb{K}$ be a field. Let
$v\in P$. Then, the map $T_{v}$ is a dominant rational map.
\end{corollary}

The reader should remember that dominant rational maps (unlike general
rational maps) can be composed, and their compositions are still dominant
rational maps. Of course, we are brushing aside subtleties like the fact that
dominant rational maps are defined only over infinite fields (unless we are
considering them in a sufficiently formal sense); as far as this paper is
concerned, it never hurts to extend the field $\mathbb{K}$ (say, by
introducing a new indeterminate), so when in doubt the reader can assume that
the field $\mathbb{K}$ is infinite.

The following proposition is trivially obtained by rewriting (\ref{def.Tv.def}%
); we are merely stating it for easier reference in proofs:

\begin{proposition}
\label{prop.Tv}Let $P$ be a finite poset. Let $\mathbb{K}$ be a field. Let
$v\in P$. For every $f\in\mathbb{K}^{\widehat{P}}$ for which $T_{v}f$ is
well-defined, we have:

\textbf{(a)} Every $w\in\widehat{P}$ such that $w\neq v$ satisfies $\left(
T_{v}f\right)  \left(  w\right)  =f\left(  w\right)  $.

\textbf{(b)} We have%
\[
\left(  T_{v}f\right)  \left(  v\right)  =\dfrac{1}{f\left(  v\right)  }%
\cdot\dfrac{\sum\limits_{\substack{u\in\widehat{P};\\u\lessdot v}}f\left(
u\right)  }{\sum\limits_{\substack{u\in\widehat{P};\\u\gtrdot v}}\dfrac
{1}{f\left(  u\right)  }}.
\]

\end{proposition}

It is very easy to check the following ``locality principle'':

\begin{proposition}
\label{prop.Tv.commute}Let $P$ be a finite poset. Let $\mathbb{K}$ be a field.
Let $v\in P$ and $w\in P$. Then, $T_{v}\circ T_{w}=T_{w}\circ T_{v}$, unless
we have either $v \lessdot w$ or $w \lessdot v$.
\end{proposition}

\begin{vershort}
\begin{proof}
[Proof of Proposition \ref{prop.Tv.commute} (sketched).]Assume that neither $v
\lessdot w$ nor $w \lessdot v$. Also, WLOG, assume that $v \neq w$, lest the
claim of the proposition be obvious.

The action of $T_{v}$ on a labelling of $P$ merely changes the label at $v$.
The new value depends on the label at $v$, on the labels at the elements
$u\in\widehat{P}$ satisfying $u\lessdot v$, and on the labels at the elements
$u\in\widehat{P}$ satisfying $u\gtrdot v$. A similar thing can be said about
the action of $T_{w}$. Since we have neither $v\lessdot w$ nor $w\lessdot v$
nor $v=w$, it thus becomes clear that the actions of $T_{v}$ and $T_{w}$ don't
interfere with each other, in the sense that the changes made by either of
them are the same no matter whether the other has been applied before it or
not. That is, $T_{v}\circ T_{w}=T_{w}\circ T_{v}$, so that Proposition
\ref{prop.Tv.commute} is proven.
\end{proof}
\end{vershort}

\begin{verlong}
\begin{proof}
[Proof of Proposition \ref{prop.Tv.commute} (sketched).]Assume that neither $v
\lessdot w$ nor $w \lessdot v$ holds. We need to show that $T_{v}\circ
T_{w}=T_{w}\circ T_{v}$. In other words, we need to show that $\left(
T_{v}\circ T_{w}\right)  f=\left(  T_{w}\circ T_{v}\right)  f$ for every
$f\in\mathbb{K}^{\widehat{P}}$ (for which both sides are defined). So let
$f\in\mathbb{K}^{\widehat{P}}$.

We can WLOG assume that $v\neq w$, because otherwise the claim is obvious.
Assume this. Recall that we don't have $v \lessdot w$; in other words, we
don't have $w \gtrdot v$.

We need to prove that $\left(  T_{v}\circ T_{w}\right)  f=\left(  T_{w}\circ
T_{v}\right)  f$. In order to prove this, it is enough to check that $\left(
\left(  T_{v}\circ T_{w}\right)  f\right)  \left(  p\right)  =\left(  \left(
T_{w}\circ T_{v}\right)  f\right)  \left(  p\right)  $ for every
$p\in\widehat{P}$. So let $p\in\widehat{P}$.

We need to check that $\left(  \left(  T_{v}\circ T_{w}\right)  f\right)
\left(  p\right)  =\left(  \left(  T_{w}\circ T_{v}\right)  f\right)  \left(
p\right)  $. If $p\notin\left\{  v,w\right\}  $, then neither of the maps
$T_{v}$ and $T_{w}$ changes the label at $p$ when acting on a labelling
(because of Proposition \ref{prop.Tv} \textbf{(a)}). Hence, if $p\notin%
\left\{  v,w\right\}  $, then the equality $\left(  \left(  T_{v}\circ
T_{w}\right)  f\right)  \left(  p\right)  =\left(  \left(  T_{w}\circ
T_{v}\right)  f\right)  \left(  p\right)  $ holds for obvious reasons (both
sides of it are $=f\left(  p\right)  $). Thus, we can WLOG assume that we
don't have $p\notin\left\{  v,w\right\}  $. Assume this. So $p\in\left\{
v,w\right\}  $. WLOG, say that $p=v$. Hence, $p=v\neq w$.

Since $p=v$, we have%
\begin{align}
\left(  \left(  T_{v}\circ T_{w}\right)  f\right)  \left(  p\right)   &
=\left(  \left(  T_{v}\circ T_{w}\right)  f\right)  \left(  v\right)  =\left(
T_{v}\left(  T_{w}f\right)  \right)  \left(  v\right) \nonumber\\
&  =\dfrac{1}{\left(  T_{w}f\right)  \left(  v\right)  }\cdot\dfrac
{\sum\limits_{\substack{u\in\widehat{P};\\u\lessdot v}}\left(  T_{w}f\right)
\left(  u\right)  }{\sum\limits_{\substack{u\in\widehat{P};\\u\gtrdot
v}}\dfrac{1}{\left(  T_{w}f\right)  \left(  u\right)  }}%
\ \ \ \ \ \ \ \ \ \ \left(  \text{by Proposition \ref{prop.Tv} \textbf{(b)}%
}\right)  . \label{pf.Tv.commute.1}%
\end{align}
But since $w\neq v$, we have $\left(  T_{w}f\right)  \left(  v\right)
=f\left(  v\right)  $ (by Proposition \ref{prop.Tv} \textbf{(a)}). Moreover,
for every $u\in\widehat{P}$ satisfying $u\lessdot v$, we have $u\neq w$
(because $u\lessdot v$, while we don't have $w\lessdot v$). Hence,%
\begin{equation}
\text{for every }u\in\widehat{P}\text{ satisfying }u\lessdot v\text{, we have
}\left(  T_{w}f\right)  \left(  u\right)  =f\left(  u\right)
\label{pf.Tv.commute.3}%
\end{equation}
(by Proposition \ref{prop.Tv} \textbf{(a)}). Furthermore, for every
$u\in\widehat{P}$ satisfying $u\gtrdot v$, we have $u\neq w$ (because
$u\gtrdot v$, while we don't have $w\gtrdot v$). Hence,%
\begin{equation}
\text{for every }u\in\widehat{P}\text{ satisfying }u\gtrdot v\text{, we have
}\left(  T_{w}f\right)  \left(  u\right)  =f\left(  u\right)
\label{pf.Tv.commute.4}%
\end{equation}
(by Proposition \ref{prop.Tv} \textbf{(a)}). In light of
(\ref{pf.Tv.commute.3}) and (\ref{pf.Tv.commute.4}) and the equality $\left(
T_{w}f\right)  \left(  v\right)  =f\left(  v\right)  $, we can rewrite
(\ref{pf.Tv.commute.1}) as%
\begin{equation}
\left(  \left(  T_{v}\circ T_{w}\right)  f\right)  \left(  p\right)
=\dfrac{1}{f\left(  v\right)  }\cdot\dfrac{\sum\limits_{\substack{u\in
\widehat{P};\\u\lessdot v}}f\left(  u\right)  }{\sum\limits_{\substack{u\in
\widehat{P};\\u\gtrdot v}}\dfrac{1}{f\left(  u\right)  }}.
\label{pf.Tv.commute.6}%
\end{equation}

On the other hand, $p=v$, so that%
\begin{align*}
\left(  \left(  T_{w}\circ T_{v}\right)  f\right)  \left(  p\right)   &
=\left(  \left(  T_{w}\circ T_{v}\right)  f\right)  \left(  v\right)  =\left(
T_{w}\left(  T_{v}f\right)  \right)  \left(  v\right)  =\left(  T_{v}f\right)
\left(  v\right) \\
&  \ \ \ \ \ \ \ \ \ \ \left(  \text{by Proposition \ref{prop.Tv}
\textbf{(a)}, since }w\neq v\right) \\
&  =\dfrac{1}{f\left(  v\right)  }\cdot\dfrac{\sum\limits_{\substack{u\in
\widehat{P};\\u \lessdot v}}f\left(  u\right)  }{\sum\limits_{\substack{u\in
\widehat{P};\\u \gtrdot v}}\dfrac{1}{f\left(  u\right)  }}%
\ \ \ \ \ \ \ \ \ \ \left(  \text{by Proposition \ref{prop.Tv} \textbf{(b)}%
}\right)  .
\end{align*}
Comparing this with (\ref{pf.Tv.commute.6}), we obtain $\left(  \left(
T_{v}\circ T_{w}\right)  f\right)  \left(  p\right)  =\left(  \left(
T_{w}\circ T_{v}\right)  f\right)  \left(  p\right)  $. This completes the
proof of Proposition \ref{prop.Tv.commute}.
\end{proof}
\end{verlong}

\begin{corollary}
\label{cor.Tv.commute}Let $P$ be a finite poset. Let $\mathbb{K}$ be a field.
Let $v$ and $w$ be two elements of $P$ which are incomparable. Then,
$T_{v}\circ T_{w}=T_{w}\circ T_{v}$.
\end{corollary}

\begin{vershort}
\begin{proof}
[\nopunct]This follows from Proposition \ref{prop.Tv.commute} because
incomparable elements never cover each other.
\end{proof}
\end{vershort}

\begin{verlong}
\begin{proof}
[Proof of Corollary \ref{cor.Tv.commute} (sketched).]We know that $v$ and $w$
are incomparable. Hence, we have neither $v \lessdot w$ nor $w \lessdot v$.
Thus, Proposition \ref{prop.Tv.commute} applies, and we obtain $T_{v}\circ
T_{w}=T_{w}\circ T_{v}$, qed.
\end{proof}
\end{verlong}

\begin{vershort}
Combining Corollary \ref{cor.Tv.commute} with Proposition
\ref{prop.linext.transitive}, we obtain:
\end{vershort}

\begin{corollary}
\label{cor.R.welldef}Let $P$ be a finite poset. Let $\mathbb{K}$ be a field.
Let $\left(  v_{1},v_{2},...,v_{m}\right)  $ be a linear extension of $P$.
Then, the dominant rational map $T_{v_{1}}\circ T_{v_{2}}\circ...\circ
T_{v_{m}}:\mathbb{K}^{\widehat{P}}\dashrightarrow\mathbb{K}^{\widehat{P}}$ is
well-defined and independent of the choice of the linear extension $\left(
v_{1},v_{2},...,v_{m}\right)  $.
\end{corollary}

\begin{verlong}
\begin{proof}
[Proof of Corollary \ref{cor.R.welldef} (sketched).]It is clear that the
dominant rational map $T_{v_{1}}\circ T_{v_{2}}\circ...\circ T_{v_{m}}$ is
well-defined (indeed, Corollary \ref{cor.Tv.dom} yields that each $T_{v_{i}}$
is dominant, and it is known that the composition of dominant rational maps is
a well-defined dominant rational map). We will now show that it is independent
of the choice of the linear extension $\left(  v_{1},v_{2},...,v_{m}\right)  $.

Unfix the linear extension $\left(  v_{1},v_{2},...,v_{m}\right)  $.

Consider the equivalence relation $\sim$ on $\mathcal{L}\left(  P\right)  $
introduced in Proposition \ref{prop.linext.transitive}. According to
Proposition \ref{prop.linext.transitive}, any two elements of $\mathcal{L}%
\left(  P\right)  $ are equivalent under the relation $\sim$.

Now, we notice the following fact: For any linear extension $\left(
v_{1},v_{2},...,v_{m}\right)  $ of $P$ and any $i\in\left\{
1,2,...,m-1\right\}  $ such that the elements $v_{i}$ and $v_{i+1}$ of $P$ are
incomparable, we have%
\[
T_{v_{1}}\circ T_{v_{2}}\circ...\circ T_{v_{m}}=T_{v_{1}}\circ T_{v_{2}}%
\circ...\circ T_{v_{i-1}}\circ T_{v_{i+1}}\circ T_{v_{i}}\circ T_{v_{i+2}%
}\circ T_{v_{i+3}}\circ...\circ T_{v_{m}}%
\]
(indeed, this follows from $T_{v_{i}}\circ T_{v_{i+1}}=T_{v_{i+1}}\circ
T_{v_{i}}$, which is a consequence of Corollary \ref{cor.Tv.commute}). Since
the equivalence relation $\sim$ is generated by the elementary relations
$\left(  v_{1},v_{2},...,v_{m}\right)  \sim\left(  v_{1},v_{2},...,v_{i-1}%
,v_{i+1},v_{i},v_{i+2},v_{i+3},...,v_{m}\right)  $, this yields that for any
two linear extensions $\left(  v_{1},v_{2},...,v_{m}\right)  $ and $\left(
w_{1},w_{2},...,w_{m}\right)  $ of $P$ satisfying $\left(  v_{1}%
,v_{2},...,v_{m}\right)  \sim\left(  w_{1},w_{2},...,w_{m}\right)  $, we have
$T_{v_{1}}\circ T_{v_{2}}\circ...\circ T_{v_{m}}=T_{w_{1}}\circ T_{w_{2}}%
\circ...\circ T_{w_{m}}$. But since \textbf{any} two linear extensions
$\left(  v_{1},v_{2},...,v_{m}\right)  $ and $\left(  w_{1},w_{2}%
,...,w_{m}\right)  $ of $P$ satisfy $\left(  v_{1},v_{2},...,v_{m}\right)
\sim\left(  w_{1},w_{2},...,w_{m}\right)  $ (because any two elements of
$\mathcal{L}\left(  P\right)  $ are equivalent under the relation $\sim$),
this can be rewritten as follows: For any two linear extensions $\left(
v_{1},v_{2},...,v_{m}\right)  $ and $\left(  w_{1},w_{2},...,w_{m}\right)  $
of $P$, we have $T_{v_{1}}\circ T_{v_{2}}\circ...\circ T_{v_{m}}=T_{w_{1}%
}\circ T_{w_{2}}\circ...\circ T_{w_{m}}$. In other words, the dominant
rational map $T_{v_{1}}\circ T_{v_{2}}\circ...\circ T_{v_{m}}$ (where $\left(
v_{1},v_{2},...,v_{m}\right)  $ is a linear extension of $P$) is independent
on the choice of the linear extension $\left(  v_{1},v_{2},...,v_{m}\right)
$. This completes the proof of Corollary \ref{cor.R.welldef}.
\end{proof}
\end{verlong}

\begin{definition}
\label{def.rm}Let $P$ be a finite poset. Let $\mathbb{K}$ be a field.
\textit{Birational rowmotion} is defined as the dominant rational map
$T_{v_{1}}\circ T_{v_{2}}\circ...\circ T_{v_{m}}:\mathbb{K}^{\widehat{P}%
}\dashrightarrow\mathbb{K}^{\widehat{P}}$, where $\left(  v_{1},v_{2}%
,...,v_{m}\right)  $ is a linear extension of $P$. This rational map is
well-defined (in particular, it does not depend on the linear extension
$\left(  v_{1},v_{2},...,v_{m}\right)  $ chosen) because of Corollary
\ref{cor.R.welldef} (and also because a linear extension of $P$ always exists;
this is Theorem \ref{thm.linext.ex}). This rational map will be denoted by $R$.
\end{definition}

The reason for the names ``birational toggle'' and ``birational rowmotion'' is
explained in the paper \cite{einstein-propp}, in which birational rowmotion
(again, defined over $\mathbb{R}_{+}$ rather than over $\mathbb{K}$) is
denoted (serendipitously from the standpoint of the second author of this
paper) by $\rho_{\mathcal{B}}$.

\newpage

\begin{example}
\label{ex.rowmotion.31}Let us demonstrate the effect of birational toggles and
birational rowmotion on a rather simple $4$-element poset. Namely, for this
example, we let $P$ be the poset $\left\{  p,q_{1},q_{2},q_{3}\right\}  $ with
order relation defined by setting $p<q_{i}$ for each $i\in\left\{
1,2,3\right\}  $. This poset has Hasse diagram
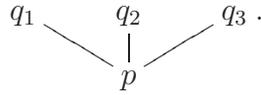
\[
\xymatrixrowsep{0.9pc}\xymatrixcolsep{2pc}\xymatrix{
q_1 \ar@{-}[rd] & q_2 \ar@{-}[d] & q_3 \ar@{-}[ld] \\
& p &
}.
\]
The extended poset $\widehat{P}$ has Hasse diagram
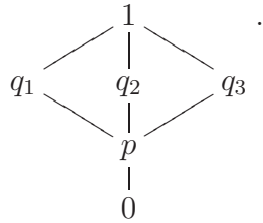
\[
\xymatrixrowsep{0.9pc}\xymatrixcolsep{2pc}\xymatrix{
& 1 \ar@{-}[rd] \ar@{-}[d] \ar@{-}[ld] & \\
q_1 \ar@{-}[rd] & q_2 \ar@{-}[d] & q_3 \ar@{-}[ld] \\
& p \ar@{-}[d] & \\
& 0 &
}.
\]
We can visualize a $\mathbb{K}$-labelling $f$ of $P$ by replacing, in the
Hasse diagram of $\widehat{P}$, each element $v\in\widehat{P}$ by the label
$f\left(  v\right)  $. Let $f$ be a $\mathbb{K}$-labelling sending $0$, $p$,
$q_{1}$, $q_{2}$, $q_{3}$, and $1$ to $a$, $w$, $x_{1}$, $x_{2}$, $x_{3}$, and
$b$, respectively (for some elements $a$, $b$, $w$, $x_{1}$, $x_{2}$, $x_{3}$
of $\mathbb{K}$); this $f$ is then visualized as follows:
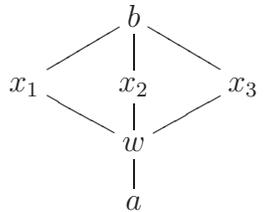
\[
\xymatrixrowsep{0.9pc}\xymatrixcolsep{2pc}\xymatrix{
& b \ar@{-}[rd] \ar@{-}[d] \ar@{-}[ld]  & \\
x_1 \ar@{-}[rd] & x_2 \ar@{-}[d] & x_3 \ar@{-}[ld] \\
& w \ar@{-}[d] & \\
& a &
}.
\]
Now, recall the definition of birational rowmotion $R$ on our poset $P$. Since
$\left(  p,q_{1},q_{2},q_{3}\right)  $ is a linear extension of $P$, we have
$R=T_{p}\circ T_{q_{1}}\circ T_{q_{2}}\circ T_{q_{3}}$. Let us track how this
transforms our labelling $f$:

\newpage We first apply $T_{q_{3}}$, obtaining%
\[
\xymatrixrowsep{0.9pc}\xymatrixcolsep{2pc}\xymatrix{
& & b \ar@{-}[rd] \ar@{-}[d] \ar@{-}[ld]  & \\
T_{q_3} f = \  & x_1 \ar@{-}[rd] & x_2 \ar@{-}[d] & {\color{red} \frac{bw}{x_3}} \ar@{-}[ld] \\
& & w \ar@{-}[d] & \\
& & a &
}
\]
(where we colored the label at $q_{3}$ red to signify that it is the label at
the element which got toggled). Indeed, the only label that changes under
$T_{q_{3}}$ is the one at $q_{3}$, and this label becomes%
\[
\left(  T_{q_{3}}f\right)  \left(  q_{3}\right)  =\dfrac{1}{f\left(
q_{3}\right)  }\cdot\dfrac{\sum\limits_{\substack{u\in\widehat{P};\\u\lessdot
q_{3}}}f\left(  u\right)  }{\sum\limits_{\substack{u\in\widehat{P};\\u\gtrdot
q_{3}}}\dfrac{1}{f\left(  u\right)  }}=\dfrac{1}{f\left(  q_{3}\right)  }%
\cdot\dfrac{f\left(  p\right)  }{\left(  \dfrac{1}{f\left(  1\right)
}\right)  }=\dfrac{1}{x_{3}}\cdot\dfrac{w}{\left(  \dfrac{1}{b}\right)
}=\dfrac{bw}{x_{3}}.
\]
Having applied $T_{q_{3}}$, we next apply $T_{q_{2}}$, obtaining%
\[
\xymatrixrowsep{0.9pc}\xymatrixcolsep{2pc}\xymatrix{
& & b \ar@{-}[rd] \ar@{-}[d] \ar@{-}[ld]  & \\
T_{q_2} T_{q_3} f = \  & x_1 \ar@{-}[rd] & {\color{red} \frac{bw}{x_2}} \ar@{-}[d] & \frac{bw}{x_3} \ar@{-}[ld] \\
& & w \ar@{-}[d] & \\
& & a &
}.
\]
Next, we apply $T_{q_{1}}$, obtaining%
\[
\xymatrixrowsep{0.9pc}\xymatrixcolsep{2pc}\xymatrix{
& & b \ar@{-}[rd] \ar@{-}[d] \ar@{-}[ld]  & \\
T_{q_1} T_{q_2} T_{q_3} f = \  & {\color{red} \frac{bw}{x_1}} \ar@{-}[rd] & \frac{bw}{x_2} \ar@{-}[d] & \frac{bw}{x_3} \ar@{-}[ld] \\
& & w \ar@{-}[d] & \\
& & a &
}.
\]
Finally, we apply $T_{p}$, resulting in
\[
\xymatrixrowsep{0.9pc}\xymatrixcolsep{2pc}\xymatrix{
& & b \ar@{-}[rd] \ar@{-}[d] \ar@{-}[ld]  & \\
T_p T_{q_1} T_{q_2} T_{q_3} f = \  & \frac{bw}{x_1} \ar@{-}[rd] & \frac{bw}{x_2} \ar@{-}[d] & \frac{bw}{x_3} \ar@{-}[ld] \\
& & {\color{red} \frac{ab}{x_1+x_2+x_3}} \ar@{-}[d] & \\
& & a &
},
\]
since the birational $p$-toggle $T_{p}$ has changed the label at $p$ to%
\begin{align*}
&  \left(  T_{p}T_{q_{1}}T_{q_{2}}T_{q_{3}}f\right)  \left(  p\right) \\
&  =\dfrac{1}{\left(  T_{q_{1}}T_{q_{2}}T_{q_{3}}f\right)  \left(  p\right)
}\cdot\dfrac{\sum\limits_{\substack{u\in\widehat{P};\\u\lessdot p}}\left(
T_{q_{1}}T_{q_{2}}T_{q_{3}}f\right)  \left(  u\right)  }{\sum
\limits_{\substack{u\in\widehat{P};\\u\gtrdot p}}\dfrac{1}{\left(  T_{q_{1}%
}T_{q_{2}}T_{q_{3}}f\right)  \left(  u\right)  }}\\
&  =\dfrac{1}{\left(  T_{q_{1}}T_{q_{2}}T_{q_{3}}f\right)  \left(  p\right)
}\cdot\dfrac{\left(  T_{q_{1}}T_{q_{2}}T_{q_{3}}f\right)  \left(  0\right)
}{\dfrac{1}{\left(  T_{q_{1}}T_{q_{2}}T_{q_{3}}f\right)  \left(  q_{1}\right)
}+\dfrac{1}{\left(  T_{q_{1}}T_{q_{2}}T_{q_{3}}f\right)  \left(  q_{2}\right)
}+\dfrac{1}{\left(  T_{q_{1}}T_{q_{2}}T_{q_{3}}f\right)  \left(  q_{3}\right)
}}\\
&  =\dfrac{1}{w}\cdot\dfrac{a}{\dfrac{1}{bw\diagup x_{1}}+\dfrac{1}{bw\diagup
x_{2}}+\dfrac{1}{bw\diagup x_{3}}}=\dfrac{ab}{x_{1}+x_{2}+x_{3}}.
\end{align*}
We thus have computed $Rf$ (since $R = T_{p} T_{q_{1}} T_{q_{2}} T_{q_{3}}$).
By repeating this procedure (or just substituting the labels of
$Rf$ obtained as variables), we can compute $R^{2}f$, $R^{3}f$ etc..
Specifically, we obtain%
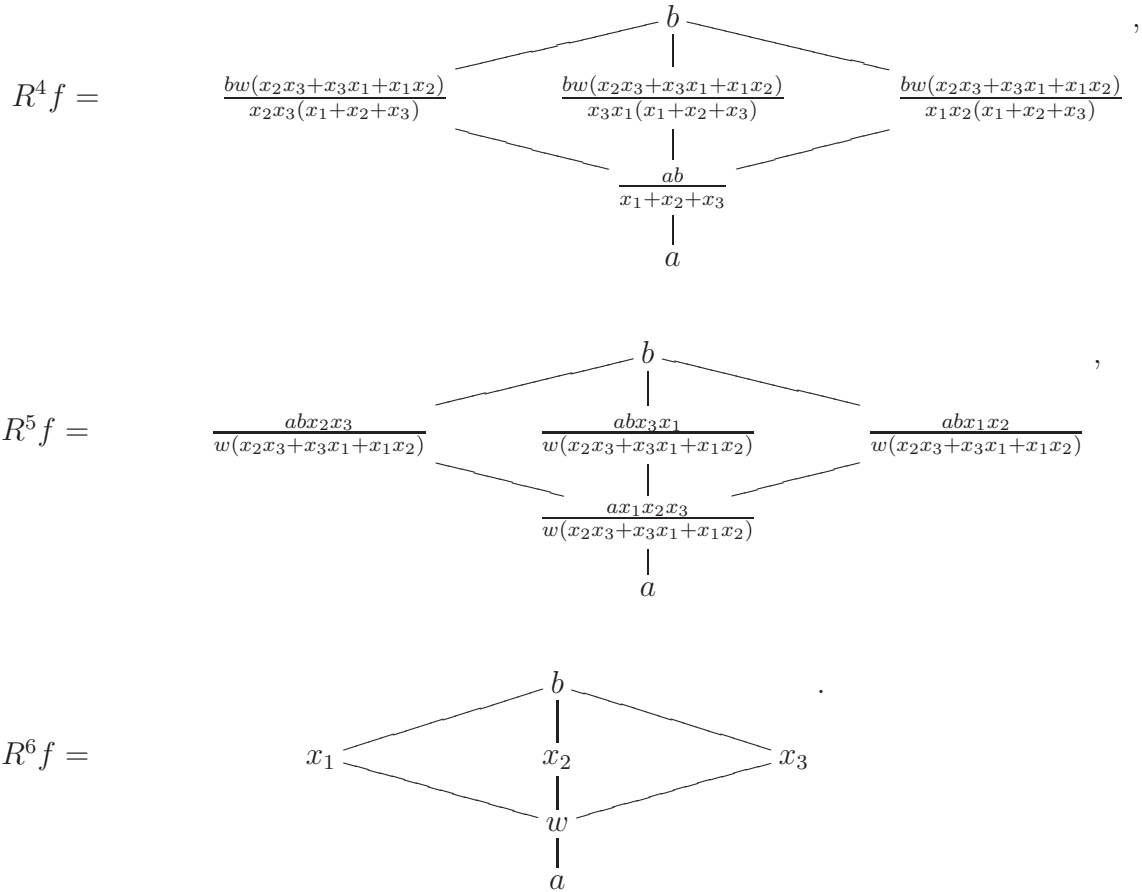
\begin{align*}
&
\xymatrixrowsep{0.9pc}\xymatrixcolsep{6pc}\xymatrix{ & & b \ar@{-}[rd] \ar@{-}[d] \ar@{-}[ld] & \\ R f = \  & \frac{bw}{x_1} \ar@{-}[rd] & \frac{bw}{x_2} \ar@{-}[d] & \frac{bw}{x_3} \ar@{-}[ld] \\ & & \frac{ab}{x_1+x_2+x_3} \ar@{-}[d] & \\ & & a & },\\
& \\
&
\xymatrixrowsep{0.9pc}\xymatrixcolsep{3pc}\xymatrix{ & & b \ar@{-}[rd] \ar@{-}[d] \ar@{-}[ld] & \\ R^2 f = \  & \frac{abx_1}{w\left(x_1+x_2+x_3\right)} \ar@{-}[rd] & \frac{abx_2}{w\left(x_1+x_2+x_3\right)} \ar@{-}[d] & \frac{abx_3}{w\left(x_1+x_2+x_3\right)} \ar@{-}[ld] \\ & & \frac{ax_1x_2x_3}{w\left(x_2x_3+x_3x_1+x_1x_2\right)} \ar@{-}[d] & \\ & & a & },\\
& \\
&
\xymatrixrowsep{0.9pc}\xymatrixcolsep{3pc}\xymatrix{ & & b \ar@{-}[rd] \ar@{-}[d] \ar@{-}[ld] & \\ R^3 f = \  & \frac{x_2x_3\left(x_1+x_2+x_3\right)}{x_2x_3+x_3x_1+x_1x_2} \ar@{-}[rd] & \frac{x_3x_1\left(x_1+x_2+x_3\right)}{x_2x_3+x_3x_1+x_1x_2} \ar@{-}[d] & \frac{x_1x_2\left(x_1+x_2+x_3\right)}{x_2x_3+x_3x_1+x_1x_2} \ar@{-}[ld] \\ & & w \ar@{-}[d] & \\ & & a & },\\
&
\end{align*}%
\begin{align*}
&
\ \xymatrixrowsep{0.9pc}\xymatrixcolsep{3pc}\xymatrix{ & & b \ar@{-}[rd] \ar@{-}[d] \ar@{-}[ld] & \\ R^4 f = \  & \frac{bw\left(x_2x_3+x_3x_1+x_1x_2\right)}{x_2x_3\left(x_1+x_2+x_3\right)} \ar@{-}[rd] & \frac{bw\left(x_2x_3+x_3x_1+x_1x_2\right)}{x_3x_1\left(x_1+x_2+x_3\right)} \ar@{-}[d] & \frac{bw\left(x_2x_3+x_3x_1+x_1x_2\right)}{x_1x_2\left(x_1+x_2+x_3\right)} \ar@{-}[ld] \\ & & \frac{ab}{x_1+x_2+x_3} \ar@{-}[d] & \\ & & a & },\\
& \\
&
\xymatrixrowsep{0.9pc}\xymatrixcolsep{3pc}\xymatrix{ & & b \ar@{-}[rd] \ar@{-}[d] \ar@{-}[ld] & \\ R^5 f = \  & \frac{abx_2x_3}{w\left(x_2x_3+x_3x_1+x_1x_2\right)} \ar@{-}[rd] & \frac{abx_3x_1}{w\left(x_2x_3+x_3x_1+x_1x_2\right)} \ar@{-}[d] & \frac{abx_1x_2}{w\left(x_2x_3+x_3x_1+x_1x_2\right)} \ar@{-}[ld] \\ & & \frac{ax_1x_2x_3}{w\left(x_2x_3+x_3x_1+x_1x_2\right)} \ar@{-}[d] & \\ & & a & },\\
& \\
&
\xymatrixrowsep{0.9pc}\xymatrixcolsep{6pc}\xymatrix{ & & b \ar@{-}[rd] \ar@{-}[d] \ar@{-}[ld] & \\ R^6 f = \  & x_1 \ar@{-}[rd] & x_2 \ar@{-}[d] & x_3 \ar@{-}[ld] \\ & & w \ar@{-}[d] & \\ & & a & }.
\end{align*}
There are several patterns here that catch the eye, some of which are related
to the very simple structure of $P$ and don't seem to generalize well.
However, the most striking observation here is that $R^{n}f=f$ for some
positive integer $n$ (namely, $n=6$ for this particular $P$). We will see in
Proposition \ref{prop.skeletal.ords} that this generalizes to a rather wide
class of posets, which we call \textquotedblleft skeletal
posets\textquotedblright\ (defined in Definition \ref{def.skeletal}), a class
of posets which contain (in particular) all graded forests such as our poset
$P$ here. (See Definition \ref{def.skeletal} for the definitions of the
concepts involved here.)
\end{example}

\begin{example}
\label{ex.rowmotion.2x2}Let us demonstrate the effect of birational toggles
and birational rowmotion on another $4$-element poset. Namely, for this
example, we let $P$ be the poset $\left\{  1,2\right\}  \times\left\{
1,2\right\}  $ with order relation defined by setting $\left(  i,k\right)
\leq\left(  i^{\prime},k^{\prime}\right)  $ if and only if $\left(  i\leq
i^{\prime}\text{ and }k\leq k^{\prime}\right)  $. This poset will later be
called the \textquotedblleft$2\times2$-rectangle\textquotedblright\ in
Definition \ref{def.rect}. It has Hasse diagram
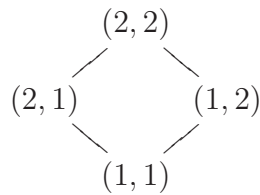
\[
\xymatrixrowsep{0.9pc}\xymatrixcolsep{0.20pc}\xymatrix{
& \left(2,2\right) \ar@{-}[rd] \ar@{-}[ld] & \\
\left(2,1\right) \ar@{-}[rd] & & \left(1,2\right) \ar@{-}[ld] \\
& \left(1,1\right) &
}.
\]
The extended poset $\widehat{P}$ has Hasse diagram
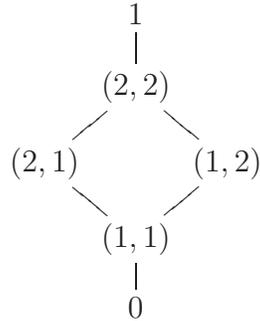
\[
\xymatrixrowsep{0.9pc}\xymatrixcolsep{0.20pc}\xymatrix{
& 1 \ar@{-}[d] & \\
& \left(2,2\right) \ar@{-}[rd] \ar@{-}[ld] & \\
\left(2,1\right) \ar@{-}[rd] & & \left(1,2\right) \ar@{-}[ld] \\
& \left(1,1\right) \ar@{-}[d] & \\
& 0 &
}.
\]
We can visualize a $\mathbb{K}$-labelling $f$ of $P$ by replacing, in the
Hasse diagram of $\widehat{P}$, each element $v\in\widehat{P}$ by the label
$f\left(  v\right)  $. Let $f$ be a $\mathbb{K}$-labelling sending $0$,
$\left(  1,1\right)  $, $\left(  1,2\right)  $, $\left(  2,1\right)  $,
$\left(  2,2\right)  $, and $1$ to $a$, $w$, $y$, $x$, $z$, and $b$,
respectively (for some elements $a$, $b$, $x$, $y$, $z$, $w$ of $\mathbb{K}$);
this $f$ is then visualized as follows:
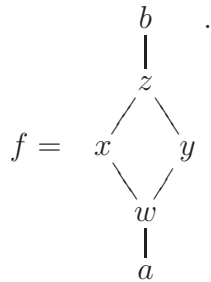
\[
\xymatrixrowsep{0.9pc}\xymatrixcolsep{0.20pc}\xymatrix{
& & b \ar@{-}[d] & \\
& & z \ar@{-}[rd] \ar@{-}[ld] & \\
f = \ & x \ar@{-}[rd] & & y \ar@{-}[ld] \\
& & w \ar@{-}[d] & \\
& & a &
}.
\]
Now, recall the definition of birational rowmotion $R$ on our poset $P$. Since
$\left(  \left(  1,1\right)  ,\left(  1,2\right)  ,\left(  2,1\right)
,\left(  2,2\right)  \right)  $ is a linear extension of $P$, we have
$R=T_{\left(  1,1\right)  }\circ T_{\left(  1,2\right)  }\circ T_{\left(
2,1\right)  }\circ T_{\left(  2,2\right)  }$. Let us track how this transforms
our labelling $f$:

\newpage We first apply $T_{\left(  2,2\right)  }$, obtaining%
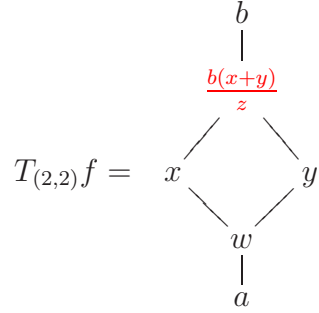
\[
\xymatrixrowsep{0.9pc}\xymatrixcolsep{0.20pc}\xymatrix{
& & b \ar@{-}[d] & \\
& & {\color{red} \frac{b\left(x+y\right)}{z}} \ar@{-}[rd] \ar@{-}[ld] & \\
T_{\left(  2,2\right)  } f = \ & x \ar@{-}[rd] & & y \ar@{-}[ld] \\
& & w \ar@{-}[d] & \\
& & a &
}
\]
(where we colored the label at $\left(  2,2\right)  $ red to signify that it
is the label at the element which got toggled). Indeed, the only label that
changes under $T_{\left(  2,2\right)  }$ is the one at $\left(  2,2\right)  $,
and this label becomes%
\begin{align*}
\left(  T_{\left(  2,2\right)  }f\right)  \left(  2,2\right)   &  =\dfrac
{1}{f\left(  \left(  2,2\right)  \right)  }\cdot\dfrac{\sum
\limits_{\substack{u\in\widehat{P};\\u\lessdot\left(  2,2\right)  }}f\left(
u\right)  }{\sum\limits_{\substack{u\in\widehat{P};\\u\gtrdot\left(
2,2\right)  }}\dfrac{1}{f\left(  u\right)  }}=\dfrac{1}{f\left(  \left(
2,2\right)  \right)  }\cdot\dfrac{f\left(  \left(  1,2\right)  \right)
+f\left(  \left(  2,1\right)  \right)  }{\left(  \dfrac{1}{f\left(  1\right)
}\right)  }\\
&  =\dfrac{1}{z}\cdot\dfrac{y+x}{\left(  \dfrac{1}{b}\right)  }=\dfrac
{b\left(  x+y\right)  }{z}.
\end{align*}
Having applied $T_{\left(  2,2\right)  }$, we next apply $T_{\left(
2,1\right)  }$, obtaining%
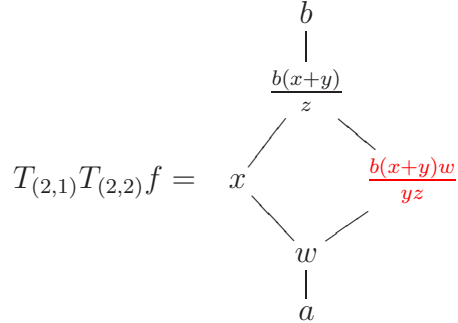
\[
\xymatrixrowsep{0.9pc}\xymatrixcolsep{0.20pc}\xymatrix{
& & b \ar@{-}[d] & \\
& & \frac{b\left(x+y\right)}{z} \ar@{-}[rd] \ar@{-}[ld] & \\
T_{\left(  2,1\right)  }T_{\left(  2,2\right)  }f
= \ & {\color{red} \frac{b\left(x+y\right)w}{xz}} \ar@{-}[rd] & & y \ar@{-}[ld] \\
& & w \ar@{-}[d] & \\
& & a &
}.
\]
Next, we apply $T_{\left(  1,2\right)  }$, obtaining%
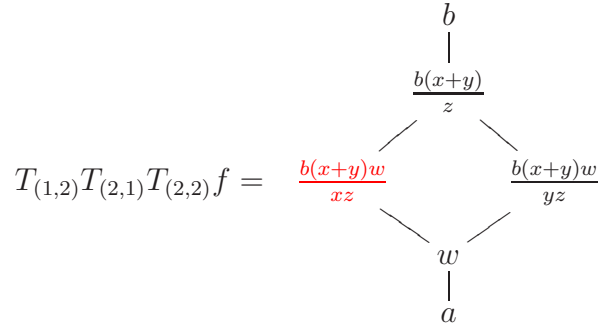
\[
\xymatrixrowsep{0.9pc}\xymatrixcolsep{0.20pc}\xymatrix{
& & b \ar@{-}[d] & \\
& & \frac{b\left(x+y\right)}{z} \ar@{-}[rd] \ar@{-}[ld] & \\
T_{\left(  1,2\right)  }T_{\left(  2,1\right)  }T_{\left(  2,2\right)
}f
= \ & {\frac{b\left(x+y\right)w}{xz}} \ar@{-}[rd] & & {\color{red} \frac{b\left(x+y\right)w}{yz}} \ar@{-}[ld] \\
& & w \ar@{-}[d] & \\
& & a &
}.
\]
Finally, we apply $T_{\left(  1,1\right)  }$, resulting in
\[
\xymatrixrowsep{0.9pc}\xymatrixcolsep{0.20pc}\xymatrix{
& & b \ar@{-}[d] & \\
& & \frac{b\left(x+y\right)}{z} \ar@{-}[rd] \ar@{-}[ld] & \\
Rf=T_{\left(  1,1\right)  }T_{\left(  1,2\right)  }T_{\left(  2,1\right)
}T_{\left(  2,2\right)  }f = \ & \frac{b\left(x+y\right)w}{xz} \ar@{-}[rd] & & \frac{b\left(x+y\right)w}{yz} \ar@{-}[ld] \\
& & {\color{red} \frac{ab}{z}} \ar@{-}[d] & \\
& & a &
}
\]
(after cancelling terms). We thus have computed $Rf$. By repeating this
procedure (or just substituting the labels of $Rf$ obtained as variables), we
can compute $R^{2}f$, $R^{3}f$ etc.. Specifically, we obtain%
\begin{align*}
&
\xymatrixrowsep{0.9pc}\xymatrixcolsep{0.20pc}\xymatrix{ & & b \ar@{-}[d] & \\ & & \frac{b\left(x+y\right)}{z} \ar@{-}[rd] \ar@{-}[ld] & \\ Rf = \ & \frac{b\left(x+y\right)w}{xz} \ar@{-}[rd] & & \frac{b\left(x+y\right)w}{yz} \ar@{-}[ld] \\ & & \frac{ab}{z} \ar@{-}[d] & \\ & & a & },\ \ \ \ \ \ \ \ \ \ \xymatrixrowsep{0.9pc}\xymatrixcolsep{0.20pc}\xymatrix{ & & b \ar@{-}[d] & \\ & & \frac{b\left(x+y\right)w}{xy} \ar@{-}[rd] \ar@{-}[ld] & \\ R^2 f = \ & \frac{ab}{y} \ar@{-}[rd] & & \frac{ab}{x} \ar@{-}[ld] \\ & & \frac{az}{x+y} \ar@{-}[d] & \\ & & a & },\\
& \\
&
\xymatrixrowsep{0.9pc}\xymatrixcolsep{0.20pc}\xymatrix{ & & b \ar@{-}[d] & \\ & & \frac{ab}{w} \ar@{-}[rd] \ar@{-}[ld] & \\ R^3 f = \ & \frac{ayz}{\left(x+y\right)w} \ar@{-}[rd] & & \frac{axz}{\left(x+y\right)w} \ar@{-}[ld] \\ & & \frac{axy}{\left(x+y\right)w} \ar@{-}[d] & \\ & & a & },\ \ \ \ \ \ \ \ \ \ \xymatrixrowsep{0.9pc}\xymatrixcolsep{0.20pc}\xymatrix{ & & b \ar@{-}[d] & \\ & & z \ar@{-}[rd] \ar@{-}[ld] & \\ R^4f = \ & x \ar@{-}[rd] & & y \ar@{-}[ld] \\ & & w \ar@{-}[d] & \\ & & a & }.
\end{align*}
There are two surprises here. First, it turns out that $R^{4}f=f$. This is not
obvious, but generalizes in at least two ways: On the one hand, our poset $P$
is a particular case of what we call a \textquotedblleft skeletal
poset\textquotedblright\ (Definition \ref{def.skeletal}), a class of posets
which all share the property (Proposition \ref{prop.skeletal.ords}) that
$R^{n}=\operatorname*{id}$ for some sufficiently high positive integer $n$
(which can be explicitly computed). On the other hand, our poset $P$ is a
particular case of rectangle posets, which turn out (Theorem
\ref{thm.rect.ord}) to satisfy $R^{p+q}=\operatorname*{id}$ with $p$ and $q$
being the side lengths (here, $2$ and $2$) of the rectangle. Second, on a more
subtle level, the rational functions appearing as labels in $Rf$, $R^{2}f$ and
$R^{3}f$ are not as \textquotedblleft wild\textquotedblright\ as one might
expect. The values $\left(  Rf\right)  \left(  \left(  1,1\right)  \right)  $,
$\left(  R^{2}f\right)  \left(  \left(  1,2\right)  \right)  $, $\left(
R^{2}f\right)  \left(  \left(  2,1\right)  \right)  $ and $\left(
R^{3}f\right)  \left(  \left(  2,2\right)  \right)  $ each have the form
$\dfrac{ab}{f\left(  v\right)  }$ for some $v\in P$. This is a
\textquotedblleft reciprocity\textquotedblright\ phenomenon which turns out to
generalize to arbitrary rectangles (Theorem \ref{thm.rect.antip.general}).

In the above calculation, we used the linear extension $\left(  \left(
1,1\right)  ,\left(  1,2\right)  ,\left(  2,1\right)  ,\left(  2,2\right)
\right)  $ of $P$ to compute $R$ as $T_{\left(  1,1\right)  }\circ T_{\left(
1,2\right)  }\circ T_{\left(  2,1\right)  }\circ T_{\left(  2,2\right)  }$. We
could have just as well used the linear extension $\left(  \left(  1,1\right)
,\left(  2,1\right)  ,\left(  1,2\right)  ,\left(  2,2\right)  \right)  $,
obtaining the same result. But we could not have used the list $\left(
\left(  1,1\right)  ,\left(  1,2\right)  ,\left(  2,2\right)  ,\left(
2,1\right)  \right)  $ (for example), since it is not a linear extension (and
indeed, the order of $T_{\left(  1,1\right)  }\circ T_{\left(  1,2\right)
}\circ T_{\left(  2,2\right)  }\circ T_{\left(  2,1\right)  }$ is infinite, as
follows from the results of \cite[\S 12.2]{einstein-propp}).
\end{example}

Let us state another proposition, which describes birational rowmotion implicitly:

\begin{proposition}
\label{prop.R.implicit}Let $P$ be a finite poset. Let $\mathbb{K}$ be a field.
Let $v\in P$. Let $f\in\mathbb{K}^{\widehat{P}}$. Then,%
\begin{equation}
\left(  Rf\right)  \left(  v\right)  =\dfrac{1}{f\left(  v\right)  }%
\cdot\dfrac{\sum\limits_{\substack{u\in\widehat{P};\\u\lessdot v}}f\left(
u\right)  }{\sum\limits_{\substack{u\in\widehat{P};\\u\gtrdot v}}\dfrac
{1}{\left(  Rf\right)  \left(  u\right)  }}. \label{prop.R.implicit.eq}%
\end{equation}

\end{proposition}

Here (and in statements further down this paper), we are taking the liberty to
leave assumptions such as \textquotedblleft Assume that $Rf$ is
well-defined\textquotedblright\ unsaid (for instance, such an assumption is
needed in Proposition \ref{prop.R.implicit}) because these assumptions are
satisfied when the parameters belong to some Zariski-dense open subset of
their domains.

\begin{vershort}
\begin{proof}
[Proof of Proposition \ref{prop.R.implicit} (sketched).]Fix a linear extension
$\left(  v_{1},v_{2},...,v_{m}\right)  $ of $P$. Recall that $R$ has been
defined as the composition $T_{v_{1}}\circ T_{v_{2}}\circ...\circ T_{v_{m}}$.
Hence, $Rf$ can be obtained from $f$ by traversing the linear extension
$\left(  v_{1},v_{2},...,v_{m}\right)  $ from right to left (thus starting
with the largest element $v_{m}$, then proceeding to $v_{m-1}$, etc.), and at
every step toggling the element being traversed. When an element $v$ is being
toggled, the elements $u\in\widehat{P}$ satisfying $u\lessdot v$ have not yet
been toggled (they are further left than $v$ in the linear extension), whereas
those satisfying $u\gtrdot v$ have been toggled already. Denoting the state of
the $\mathbb{K}$-labelling \text{before} the $v$-toggle by $g$, we see that
the state \textit{after} the $v$-toggle will be $T_{v}g$ with%
\begin{equation}
\left(  T_{v}g\right)  \left(  w\right)  =\left\{
\begin{array}
[c]{l}%
g\left(  w\right)  ,\ \ \ \ \ \ \ \ \ \ \text{if }w\neq v;\\
\dfrac{1}{g\left(  v\right)  }\cdot\dfrac{\sum\limits_{\substack{u\in
\widehat{P};\\u\lessdot v}}g\left(  u\right)  }{\sum\limits_{\substack{u\in
\widehat{P};\\u\gtrdot v}}\dfrac{1}{g\left(  u\right)  }}%
,\ \ \ \ \ \ \ \ \ \ \text{if }w=v
\end{array}
\right.  \ \ \ \ \ \ \ \ \ \ \text{for all }w\in\widehat{P}.
\label{pf.R.implicit.short.1}%
\end{equation}

But $g\left(  v\right)  =f\left(  v\right)  $ (since $v$ has not yet been
toggled at the time of $g$) and $\left(  T_{v}g\right)  \left(  v\right)
=\left(  Rf\right)  \left(  v\right)  $ (since $v$ has been toggled at the
time of $T_{v}g$, and is not going to be toggled ever again during the process
of computing $Rf$); moreover, all $u\in\widehat{P}$ satisfying $u\lessdot v$
satisfy $g\left(  u\right)  =f\left(  u\right)  $ (since these $u$ have not
yet been toggled), whereas all $u\in\widehat{P}$ satisfying $u\gtrdot v$
satisfy $g\left(  u\right)  =\left(  Rf\right)  \left(  u\right)  $ (since
these $u$ have already been toggled and will not be toggled ever again). Thus,
(\ref{pf.R.implicit.short.1}) (applied to $w=v$) transforms into
(\ref{prop.R.implicit.eq}). Proposition \ref{prop.R.implicit} is proven.
\end{proof}
\end{vershort}

\begin{verlong}
This proposition is, in fact, an easy consequence of the definition of $R$.
Loosely speaking, the main idea of the proof is that when you traverse a
linear extension $\left(  v_{1},v_{2},...,v_{m}\right)  $ of a poset $P$ from
right to left (thus starting with $v_{m}$, then proceeding to $v_{m-1}$,
etc.), at every step the elements $u\in\widehat{P}$ satisfying $u\lessdot v$
(where $v$ is the element you are currently visiting) lie ahead of you,
whereas those satisfying $u\gtrdot v$ are behind you. For the sake of
completeness, let us give the full proof:

\begin{proof}
[Proof of Proposition \ref{prop.R.implicit} (sketched).]Let $\left(
v_{1},v_{2},...,v_{m}\right)  $ be a linear extension of $P$. Let $i$ be the
index satisfying $v_{i}=v$. By the definition of birational rowmotion $R$, we
have $R=T_{v_{1}}\circ T_{v_{2}}\circ...\circ T_{v_{m}}$.

Let $A=T_{v_{i+1}}\circ T_{v_{i+2}}\circ...\circ T_{v_{m}}$ and $B=T_{v_{1}%
}\circ T_{v_{2}}\circ...\circ T_{v_{i-1}}$. Then,%
\[
R=T_{v_{1}}\circ T_{v_{2}}\circ...\circ T_{v_{m}}=\underbrace{T_{v_{1}}\circ
T_{v_{2}}\circ...\circ T_{v_{i-1}}}_{=B}\circ\underbrace{T_{v_{i}}}_{=T_{v}%
}\circ\underbrace{T_{v_{i+1}}\circ T_{v_{i+2}}\circ...\circ T_{v_{m}}}%
_{=A}=B\circ T_{v}\circ A.
\]

Now:

\begin{itemize}
\item Each of the maps $T_{v_{j}}$ with $j\neq i$ leaves the label at $v$
invariant when acting on a $\mathbb{K}$-labelling. Hence, each of the maps $B$
and $A$ leaves the label at $v$ invariant (since $B$ and $A$ are compositions
of maps $T_{v_{j}}$ with $j\neq i$). Thus, $\left(  B\left(  \left(
T_{v}\circ A\right)  f\right)  \right)  \left(  v\right)  =\left(  \left(
T_{v}\circ A\right)  f\right)  \left(  v\right)  $ and $\left(  Af\right)
\left(  v\right)  =f\left(  v\right)  $. Now, since $R=B\circ T_{v}\circ A$,
we have%
\begin{align}
\left(  Rf\right)  \left(  v\right)   &  =\left(  \left(  B\circ T_{v}\circ
A\right)  f\right)  \left(  v\right)  =\left(  B\left(  \left(  T_{v}\circ
A\right)  f\right)  \right)  \left(  v\right)  =\left(  \left(  T_{v}\circ
A\right)  f\right)  \left(  v\right) \nonumber\\
&  =\left(  T_{v}\left(  Af\right)  \right)  \left(  v\right) \nonumber\\
&  =\dfrac{1}{\left(  Af\right)  \left(  v\right)  }\cdot\dfrac{\sum
\limits_{\substack{u\in\widehat{P};\\u \lessdot v}}\left(  Af\right)  \left(
u\right)  }{\sum\limits_{\substack{u\in\widehat{P};\\u \gtrdot v}}\dfrac
{1}{\left(  Af\right)  \left(  u\right)  }}\ \ \ \ \ \ \ \ \ \ \left(
\text{by Proposition \ref{prop.Tv} \textbf{(b)}}\right) \nonumber\\
&  =\dfrac{1}{f\left(  v\right)  }\cdot\dfrac{\sum\limits_{\substack{u\in
\widehat{P};\\u \lessdot v}}\left(  Af\right)  \left(  u\right)  }%
{\sum\limits_{\substack{u\in\widehat{P};\\u \gtrdot v}}\dfrac{1}{\left(
Af\right)  \left(  u\right)  }} \label{pf.R.implicit.1}%
\end{align}
(since $\left(  Af\right)  \left(  v\right)  =f\left(  v\right)  $).

\item Let $u\in\widehat{P}$ be such that $u \lessdot v$. Then, $u<v=v_{i}$ in
$\widehat{P}$. Hence, $u$ is none of the elements $v_{i+1}$, $v_{i+2}$, $...$,
$v_{m}$ (because $\left(  v_{1},v_{2},...,v_{m}\right)  $ is a linear
extension of $P$). Thus, each of the maps $T_{v_{i+1}}$, $T_{v_{i+2}}$, $...$,
$T_{v_{m}}$ leaves the label at $u$ invariant when acting on a $\mathbb{K}%
$-labelling. Therefore, $A$ also leaves the label at $u$ invariant (since $A$
is a composition of these maps $T_{v_{i+1}}$, $T_{v_{i+2}}$, $...$, $T_{v_{m}%
}$). Hence, $\left(  Af\right)  \left(  u\right)  =f\left(  u\right)  $.

Forget that we fixed $u$. We have thus shown that%
\begin{equation}
\left(  Af\right)  \left(  u\right)  =f\left(  u\right)
\ \ \ \ \ \ \ \ \ \ \text{for every }u\in\widehat{P}\text{ such that } u
\lessdot v. \label{pf.R.implicit.2}%
\end{equation}

\item Let $u\in\widehat{P}$ be such that $u \gtrdot v$. Then, $u>v=v_{i}$ in
$\widehat{P}$. Hence, $u$ is none of the elements $v_{1}$, $v_{2}$, $...$,
$v_{i-1}$ (because $\left(  v_{1},v_{2},...,v_{m}\right)  $ is a linear
extension of $P$). Thus, each of the maps $T_{v_{1}}$, $T_{v_{2}}$, $...$,
$T_{v_{i-1}}$ leaves the label at $u$ invariant when acting on a $\mathbb{K}%
$-labelling. Therefore, $B$ also leaves the label at $u$ invariant (since $B$
is a composition of these maps $T_{v_{1}}$, $T_{v_{2}}$, $...$, $T_{v_{i-1}}%
$). Since $T_{v}$ also leaves the label at $u$ invariant (because $u\neq v$),
this yields that the composition $B\circ T_{v}$ also leaves the label at $u$
invariant. Hence, $\left(  \left(  B\circ T_{v}\right)  \left(  Af\right)
\right)  \left(  u\right)  =\left(  Af\right)  \left(  u\right)  $, so that%
\[
\left(  Af\right)  \left(  u\right)  =\left(  \left(  B\circ T_{v}\right)
\left(  Af\right)  \right)  \left(  u\right)  =\left(  \underbrace{\left(
B\circ T_{v}\circ A\right)  }_{=R}f\right)  \left(  u\right)  =\left(
Rf\right)  \left(  u\right)  .
\]

Forget that we fixed $u$. We thus have proven that%
\begin{equation}
\left(  Af\right)  \left(  u\right)  =\left(  Rf\right)  \left(  u\right)
\ \ \ \ \ \ \ \ \ \ \text{for every }u\in\widehat{P}\text{ such that } u
\gtrdot v. \label{pf.R.implicit.3}%
\end{equation}

\end{itemize}

Now, substituting (\ref{pf.R.implicit.2}) and (\ref{pf.R.implicit.3}) into
(\ref{pf.R.implicit.1}), we obtain%
\[
\left(  Rf\right)  \left(  v\right)  =\dfrac{1}{f\left(  v\right)  }%
\cdot\dfrac{\sum\limits_{\substack{u\in\widehat{P};\\u \lessdot v}}f\left(
u\right)  }{\sum\limits_{\substack{u\in\widehat{P};\\u \gtrdot v}}\dfrac
{1}{\left(  Rf\right)  \left(  u\right)  }}.
\]
This proves Proposition \ref{prop.R.implicit}.
\end{proof}
\end{verlong}

Here is a little triviality to complete the picture of Proposition
\ref{prop.R.implicit}:

\begin{proposition}
\label{prop.R.implicit.01}Let $P$ be a finite poset. Let $\mathbb{K}$ be a
field. Let $f\in\mathbb{K}^{\widehat{P}}$. Then, $\left(  Rf\right)  \left(
0\right)  =f\left(  0\right)  $ and $\left(  Rf\right)  \left(  1\right)
=f\left(  1\right)  $.
\end{proposition}

\begin{vershort}
\begin{proof}
[\nopunct]This is clear since no toggle changes the labels at $0$ and~$1$.
\end{proof}
\end{vershort}

\begin{verlong}
\begin{proof}
[Proof of Proposition \ref{prop.R.implicit.01} (sketched).]Let $\left(
v_{1},v_{2},...,v_{m}\right)  $ be a linear extension of $P$. By the
definition of birational rowmotion $R$, we have $R=T_{v_{1}}\circ T_{v_{2}%
}\circ...\circ T_{v_{m}}$. But each of the maps $T_{v_{j}}$ leaves the label
at $0$ invariant when acting on a $\mathbb{K}$-labelling (because $v_{j}\neq
0$). Therefore, the composition $T_{v_{1}}\circ T_{v_{2}}\circ...\circ
T_{v_{m}}$ of these maps also leaves the label at $0$ invariant when acting on
a $\mathbb{K}$-labelling. Thus, $\left(  \left(  T_{v_{1}}\circ T_{v_{2}}%
\circ...\circ T_{v_{m}}\right)  f\right)  \left(  0\right)  =f\left(
0\right)  $. Since $T_{v_{1}}\circ T_{v_{2}}\circ...\circ T_{v_{m}}=R$, this
rewrites as $\left(  Rf\right)  \left(  0\right)  =f\left(  0\right)  $.
Similarly, $\left(  Rf\right)  \left(  1\right)  =f\left(  1\right)  $. This
proves Proposition \ref{prop.R.implicit.01}.
\end{proof}
\end{verlong}

We will often use Proposition \ref{prop.R.implicit.01} tacitly. A trivial
corollary of Proposition \ref{prop.R.implicit.01} is:

\begin{corollary}
\label{cor.R.implicit.01}Let $P$ be a finite poset. Let $\mathbb{K}$ be a
field. Let $f\in\mathbb{K}^{\widehat{P}}$ and $\ell
\in\mathbb{N}$. Then, $\left(  R^{\ell}f\right)  \left(  0\right)  =f\left(
0\right)  $ and $\left(  R^{\ell}f\right)  \left(  1\right)  =f\left(
1\right)  $.
\end{corollary}

(Recall that $\mathbb{N}$ denotes the set
$\left\{  0,1,2,...\right\}  $ in this paper.)

\begin{verlong}
\begin{proof}
[Proof of Corollary \ref{cor.R.implicit.01} (sketched).]We will prove
Corollary \ref{cor.R.implicit.01} by induction over $\ell$.

\textit{Induction base:} We have $\left(  \underbrace{R^{0}}%
_{=\operatorname*{id}}f\right)  \left(  0\right)  =\underbrace{\left(
\operatorname*{id}f\right)  }_{=f}\left(  0\right)  =f\left(  0\right)  $ and
$\left(  \underbrace{R^{0}}_{=\operatorname*{id}}f\right)  \left(  1\right)
=\underbrace{\left(  \operatorname*{id}f\right)  }_{=f}\left(  1\right)
=f\left(  1\right)  $. Thus, Corollary \ref{cor.R.implicit.01} is proven for
$\ell=1$.

\textit{Induction step:} Let $L\in\mathbb{N}$. Assume that Corollary
\ref{cor.R.implicit.01} holds for $\ell=L$. We need to show that Corollary
\ref{cor.R.implicit.01} holds for $\ell=L+1$.

We have assumed that Corollary \ref{cor.R.implicit.01} holds for $\ell=L$. In
other words, we have $\left(  R^{L}f\right)  \left(  0\right)  =f\left(
0\right)  $ and $\left(  R^{L}f\right)  \left(  1\right)  =f\left(  1\right)
$. But Proposition \ref{prop.R.implicit.01} (applied to $R^{L}f$ instead of
$f$) yields $\left(  R\left(  R^{L}f\right)  \right)  \left(  0\right)
=\left(  R^{L}f\right)  \left(  0\right)  $ and $\left(  R\left(
R^{L}f\right)  \right)  \left(  1\right)  =\left(  R^{L}f\right)  \left(
1\right)  $. Now, $R^{L+1}=R\circ R^{L}$, so that $R^{L+1}f=\left(  R\circ
R^{L}\right)  f=R\left(  R^{L}f\right)  $, and hence%
\[
\underbrace{\left(  R^{L+1}f\right)  }_{=R\left(  R^{L}f\right)  }\left(
0\right)  =\left(  R\left(  R^{L}f\right)  \right)  \left(  0\right)  =\left(
R^{L}f\right)  \left(  0\right)  =f\left(  0\right)
\]
and
\[
\underbrace{\left(  R^{L+1}f\right)  }_{=R\left(  R^{L}f\right)  }\left(
1\right)  =\left(  R\left(  R^{L}f\right)  \right)  \left(  1\right)  =\left(
R^{L}f\right)  \left(  1\right)  =f\left(  1\right)  .
\]
So we have shown that we have $\left(  R^{L+1}f\right)  \left(  0\right)
=f\left(  0\right)  $ and $\left(  R^{L+1}f\right)  \left(  1\right)
=f\left(  1\right)  $. In other words, Corollary \ref{cor.R.implicit.01} holds
for $\ell=L+1$. Hence, the induction step is complete. Thus, Corollary
\ref{cor.R.implicit.01} is proven by induction.
\end{proof}
\end{verlong}

We will also need a converse of Propositions \ref{prop.R.implicit} and
\ref{prop.R.implicit.01}:

\begin{proposition}
\label{prop.R.implicit.converse}Let $P$ be a finite poset. Let $\mathbb{K}$ be
a field. Let $f\in\mathbb{K}^{\widehat{P}}$ and $g\in\mathbb{K}^{\widehat{P}}$
be such that $f\left(  0\right)  =g\left(  0\right)  $ and $f\left(  1\right)
=g\left(  1\right)  $. Assume that%
\begin{equation}
g\left(  v\right)  =\dfrac{1}{f\left(  v\right)  }\cdot\dfrac{\sum
\limits_{\substack{u\in\widehat{P};\\u\lessdot v}}f\left(  u\right)  }%
{\sum\limits_{\substack{u\in\widehat{P};\\u\gtrdot v}}\dfrac{1}{g\left(
u\right)  }}\ \ \ \ \ \ \ \ \ \ \text{for every }v\in P.
\label{prop.R.implicit.converse.eq}%
\end{equation}
(This means, in particular, that we assume that all denominators in
(\ref{prop.R.implicit.converse.eq}) are nonzero.) Then, $g=Rf$.\ \ \ \ \footnotemark
\end{proposition}

\footnotetext{More precisely, $Rf$ \textbf{is well-defined} and equal to $g$.}

\begin{vershort}
\begin{proof}
[Proof of Proposition \ref{prop.R.implicit.converse} (sketched).]It is clearly
enough to show that $g\left(  v\right)  =\left(  Rf\right)  \left(  v\right)
$ for every $v\in\widehat{P}$. Since this is clear for $v=0$ (since $g\left(
0\right)  =f\left(  0\right)  =\left(  Rf\right)  \left(  0\right)  $), we
only need to consider the case when $v\in\left\{  1\right\}  \cup P$. In this
case, we can prove $g\left(  v\right)  =\left(  Rf\right)  \left(  v\right)  $
by descending induction over $v$ -- that is, we assume as an induction
hypothesis that $g\left(  u\right)  =\left(  Rf\right)  \left(  u\right)  $
holds for all elements $u\in\left\{  1\right\}  \cup P$ which are greater than
$v$ in $\widehat{P}$. The induction base ($v=1$) is clear (just like $v=0$),
and the induction step follows by comparing (\ref{prop.R.implicit.eq}) with
(\ref{prop.R.implicit.converse.eq}). We leave the details (including a check
that $Rf$ is well-defined, which piggybacks on the induction) to the reader.
\end{proof}
\end{vershort}

\begin{verlong}
\begin{proof}
[Proof of Proposition \ref{prop.R.implicit.converse} (sketched).]We denote by
$<$ the smaller relation of the poset $P$. We also denote by $<$ the smaller
relation of the poset $\widehat{P}$. We thus have given two meanings to the
notation $<$, but these two meanings don't conflict (since the restriction of
the smaller relation of $\widehat{P}$ to $P$ is the smaller relation of $P$).

There exists a linear extension of $P$ (by Theorem \ref{thm.linext.ex}). Fix
such a linear extension $\left(  v_{1},v_{2},...,v_{m}\right)  $. Extend this
$m$-tuple $\left(  v_{1},v_{2},...,v_{m}\right)  $ to an $\left(  m+2\right)
$-tuple $\left(  v_{0},v_{1},v_{2},...,v_{m},v_{m+1}\right)  $ by defining
$v_{0}=0$ and $v_{m+1}=1$.

We know that $\left(  v_{1},v_{2},...,v_{m}\right)  $ is a linear extension of
$P$. In other words, $\left(  v_{1},v_{2},...,v_{m}\right)  $ is a list of the
elements of $P$ such that every element of $P$ occurs exactly once in this
list, and such that any $i\in\left\{  1,2,...,m\right\}  $ and $j\in\left\{
1,2,...,m\right\}  $ satisfying $v_{i}<v_{j}$ must satisfy $i<j$%
\ \ \ \ \footnote{This is because $\left(  v_{1},v_{2},...,v_{m}\right)  $ is
a linear extension of $P$ if and only if $\left(  v_{1},v_{2},...,v_{m}%
\right)  $ is a list of the elements of $P$ such that every element of $P$
occurs exactly once in this list, and such that any $i\in\left\{
1,2,...,m\right\}  $ and $j\in\left\{  1,2,...,m\right\}  $ satisfying
$v_{i}<v_{j}$ must satisfy $i<j$ (by the definition of a \textquotedblleft
linear extension\textquotedblright).}. In particular, $\left(  v_{1}%
,v_{2},...,v_{m}\right)  $ is a list of the elements of $P$ such that every
element of $P$ occurs exactly once in this list. Hence, $\left\{  v_{1}%
,v_{2},...,v_{m}\right\}  =P$. Now, $\left\{  v_{0},v_{1},v_{2},...,v_{m}%
,v_{m+1}\right\}  =\underbrace{\left\{  v_{1},v_{2},...,v_{m}\right\}  }%
_{=P}\cup\left\{  \underbrace{v_{0}}_{=0},\underbrace{v_{m+1}}_{=1}\right\}
=P\cup\left\{  0,1\right\}  =\widehat{P}$, qed.
\begin{equation}
\left\{  v_{0},v_{1},v_{2},...,v_{m},v_{m+1}\right\}  =\widehat{P}
\label{pf.R.implicit.converse.trueproof.Phat}%
\end{equation}
\footnote{\textit{Proof.} We know that $\left(  v_{1},v_{2},...,v_{m}\right)
$ is a linear extension of $P$. In other words, $\left(  v_{1},v_{2}%
,...,v_{m}\right)  $ is a list of the elements of $P$ such that every element
of $P$ occurs exactly once in this list, and such that any $i\in\left\{
1,2,...,m\right\}  $ and $j\in\left\{  1,2,...,m\right\}  $ satisfying
$v_{i}<v_{j}$ must satisfy $i<j$ (because $\left(  v_{1},v_{2},...,v_{m}%
\right)  $ is a linear extension of $P$ if and only if $\left(  v_{1}%
,v_{2},...,v_{m}\right)  $ is a list of the elements of $P$ such that every
element of $P$ occurs exactly once in this list, and such that any
$i\in\left\{  1,2,...,m\right\}  $ and $j\in\left\{  1,2,...,m\right\}  $
satisfying $v_{i}<v_{j}$ must satisfy $i<j$ (by the definition of a
\textquotedblleft linear extension\textquotedblright)). In particular,
$\left(  v_{1},v_{2},...,v_{m}\right)  $ is a list of the elements of $P$ such
that every element of $P$ occurs exactly once in this list. Hence, $\left\{
v_{1},v_{2},...,v_{m}\right\}  =P$. Now, $\left\{  v_{0},v_{1},v_{2}%
,...,v_{m},v_{m+1}\right\}  =\underbrace{\left\{  v_{1},v_{2},...,v_{m}%
\right\}  }_{=P}\cup\left\{  \underbrace{v_{0}}_{=0},\underbrace{v_{m+1}}%
_{=1}\right\}  =P\cup\left\{  0,1\right\}  =\widehat{P}$, qed.}.

Recall that any $i\in\left\{  1,2,...,m\right\}  $ and $j\in\left\{
1,2,...,m\right\}  $ satisfying $v_{i}<v_{j}$ must satisfy%
\begin{equation}
i<j. \label{pf.R.implicit.converse.trueproof.linextP}%
\end{equation}
Consequently, any $i\in\left\{  0,1,...,m+1\right\}  $ and $j\in\left\{
0,1,...,m+1\right\}  $ satisfying $v_{i}<v_{j}$ must satisfy
\begin{equation}
i<j. \label{pf.R.implicit.converse.trueproof.linextP01}%
\end{equation}
\footnote{\textit{Proof of (\ref{pf.R.implicit.converse.trueproof.linextP01}%
):} Let $i\in\left\{  0,1,...,m+1\right\}  $ and $j\in\left\{
0,1,...,m+1\right\}  $ be such that $v_{i}<v_{j}$. We need to prove that
$i<j$.
\par
We have $v_{i}<v_{j}$, so that $v_{i}\neq v_{j}$ and thus $i\neq j$.
\par
Assume (for the sake of contradiction) that $i=m+1$. Then, $v_{i}=v_{m+1}=1$
(by the definition of $v_{m+1}$). But the definition of $\widehat{P}$ yields
that $p\leq1$ for every $p\in\widehat{P}$. Applying this to $p=v_{j}$, this
yields $v_{j}\leq1$. Now, $v_{i}=1$, so that $1=v_{i}<v_{j}\leq1$, which is
absurd. This contradiction shows that our assumption (that $i=m+1$) was wrong.
Hence, we cannot have $i=m+1$. Thus, we have $i\neq m+1$.
\par
Assume that $i=0$. Then, $j\in\left\{  0,1,...,m+1\right\}  $, so that
$j\geq0$. Thus, $0\leq j$, so that $i=0\leq j$ and thus $i<j$ (since $i\neq
j$). Now, forget that we assumed that $i=0$. We thus have shown that $i<j$ if
$i=0$. Hence, our proof of $i<j$ is complete if $i=0$. Thus, for the rest of
this proof, we can WLOG assume that we don't have $i=0$. Assume this.
\par
We don't have $i=0$. Hence, we have $i\neq0$. Combined with $i\neq m+1$, this
yields $i\notin\left\{  0,m+1\right\}  $. Since $i\in\left\{
0,1,...,m+1\right\}  $ and $i\notin\left\{  0,m+1\right\}  $, we have
$i\in\left\{  0,1,...,m+1\right\}  \setminus\left\{  0,m+1\right\}  =\left\{
1,2,...,m\right\}  $.
\par
Assume (for the sake of contradiction) that $j=0$. Then, $v_{j}=v_{0}=0$ (by
the definition of $v_{0}$). But the definition of $\widehat{P}$ yields that
$0\leq p$ for every $p\in\widehat{P}$. Applying this to $p=v_{i}$, this yields
$0\leq v_{i}$. Now, $0\leq v_{i}<v_{j}=0$, which is absurd. This contradiction
shows that our assumption (that $j=0$) was wrong. Hence, we cannot have $j=0$.
Thus, we have $j\neq0$.
\par
Assume that $j=m+1$. Then, $i\in\left\{  0,1,...,m+1\right\}  $, so that
$i\leq m+1$. Thus, $i\leq m+1=j$ (since $j=m+1$), hence $i<j$ (since $i\neq
j$). Now, forget that we assumed that $j=m+1$. We thus have shown that $i<j$
if $j=m+1$. Hence, our proof of $i<j$ is complete if $j=m+1$. Thus, for the
rest of this proof, we can WLOG assume that we don't have $j=m+1$. Assume
this.
\par
We don't have $j=m+1$. Hence, we have $j\neq m+1$. Combined with $j\neq0$,
this yields $j\notin\left\{  0,m+1\right\}  $. Since $j\in\left\{
0,1,...,m+1\right\}  $ and $j\notin\left\{  0,m+1\right\}  $, we have
$j\in\left\{  0,1,...,m+1\right\}  \setminus\left\{  0,m+1\right\}  =\left\{
1,2,...,m\right\}  $.
\par
We now know that $i\in\left\{  1,2,...,m\right\}  $ and $j\in\left\{
1,2,...,m\right\}  $. Hence, (\ref{pf.R.implicit.converse.trueproof.linextP})
yields $i<j$. This completes the proof of $i<j$, and thus
(\ref{pf.R.implicit.converse.trueproof.linextP01}) is proven.}

For every $i\in\left\{  0,1,...,m\right\}  $, let us denote the rational map
$T_{v_{m-i+1}}\circ T_{v_{m-i+2}}\circ...\circ T_{v_{m}}:\mathbb{K}%
^{\widehat{P}}\dashrightarrow\mathbb{K}^{\widehat{P}}$ by $\xi_{i}$. Then, the
definition of $\xi_{m}$ yields%
\[
\xi_{m}=T_{v_{m-m+1}}\circ T_{v_{m-m+2}}\circ...\circ T_{v_{m}}=T_{v_{1}}\circ
T_{v_{2}}\circ...\circ T_{v_{m}}=R
\]
(since $R=T_{v_{1}}\circ T_{v_{2}}\circ...\circ T_{v_{m}}$ (by the definition
of $R$)). On the other hand, the definition of $\xi_{0}$ yields%
\[
\xi_{0}=T_{v_{m-0+1}}\circ T_{v_{m-0+2}}\circ...\circ T_{v_{m}}=T_{v_{m+1}%
}\circ T_{v_{m+2}}\circ...\circ T_{v_{m}}=\left(  \text{empty composition}%
\right)  =\operatorname*{id}.
\]

We are going to prove that for every $i\in\left\{  0,1,...,m\right\}  $, the
following assertion holds:

\textit{Assertion A:} The $\mathbb{K}$-labelling $\xi_{i}f\in\mathbb{K}%
^{\widehat{P}}$ is well-defined, and every $j\in\left\{  0,1,...,m+1\right\}
$ satisfies%
\[
\left(  \xi_{i}f\right)  \left(  v_{j}\right)  =\left\{
\begin{array}
[c]{c}%
g\left(  v_{j}\right)  ,\ \ \ \ \ \ \ \ \ \ \text{if }j>m-i;\\
f\left(  v_{j}\right)  ,\ \ \ \ \ \ \ \ \ \ \text{if }j\leq m-i
\end{array}
\right.  .
\]

In fact, we will prove Assertion A by induction over $i$:

\textit{Induction base:} Assertion A is satisfied for $i=0$%
\ \ \ \ \footnote{\textit{Proof.} Assume that $i=0$. Hence, $m-i=m-0=m$, so
that $m=m-i$. But $i=0$, and thus $\xi_{i}=\xi_{0}=\operatorname*{id}$, so
that $\xi_{i}f=\operatorname*{id}f=f$ is well-defined. Now, let $j\in\left\{
0,1,...,m+1\right\}  $. We are going to prove that
\begin{equation}
\left(  \xi_{i}f\right)  \left(  v_{j}\right)  =\left\{
\begin{array}
[c]{c}%
g\left(  v_{j}\right)  ,\ \ \ \ \ \ \ \ \ \ \text{if }j>m-i;\\
f\left(  v_{j}\right)  ,\ \ \ \ \ \ \ \ \ \ \text{if }j\leq m-i
\end{array}
\right.  . \label{pf.R.implicit.converse.trueproof.base.goal}%
\end{equation}
\par
Let us first assume that $j\leq m$. Then,%
\begin{equation}
\left\{
\begin{array}
[c]{c}%
g\left(  v_{j}\right)  ,\ \ \ \ \ \ \ \ \ \ \text{if }j>m-i;\\
f\left(  v_{j}\right)  ,\ \ \ \ \ \ \ \ \ \ \text{if }j\leq m-i
\end{array}
\right.  =f\left(  v_{j}\right)  \ \ \ \ \ \ \ \ \ \ \left(  \text{since
}j\leq m=m-i\right)  . \label{pf.R.implicit.converse.trueproof.base.1}%
\end{equation}
Now,%
\[
\underbrace{\left(  \xi_{i}f\right)  }_{=f}\left(  v_{j}\right)  =f\left(
v_{j}\right)  =\left\{
\begin{array}
[c]{c}%
g\left(  v_{j}\right)  ,\ \ \ \ \ \ \ \ \ \ \text{if }j>m-i;\\
f\left(  v_{j}\right)  ,\ \ \ \ \ \ \ \ \ \ \text{if }j\leq m-i
\end{array}
\right.
\]
(by (\ref{pf.R.implicit.converse.trueproof.base.1})). Thus,
(\ref{pf.R.implicit.converse.trueproof.base.goal}) is proven under the
assumption that $j\leq m$.
\par
Now, forget that we assumed that $j\leq m$. We thus have proven
(\ref{pf.R.implicit.converse.trueproof.base.goal}) under the assumption that
$j\leq m$. Hence, for the rest of the proof of
(\ref{pf.R.implicit.converse.trueproof.base.goal}), we can WLOG assume that we
don't have $j\leq m$. Assume this. We don't have $j\leq m$. Hence, we have
$j>m$. Since $j\in\left\{  0,1,...,m+1\right\}  $, this yields $j=m+1$. Thus,
$v_{j}=v_{m+1}=1$ (by the definition of $v_{j}$). Hence, $g\left(
\underbrace{v_{j}}_{=1}\right)  =g\left(  1\right)  =f\left(  1\right)  $
(because $f\left(  1\right)  =g\left(  1\right)  $). But%
\begin{align}
\left\{
\begin{array}
[c]{c}%
g\left(  v_{j}\right)  ,\ \ \ \ \ \ \ \ \ \ \text{if }j>m-i;\\
f\left(  v_{j}\right)  ,\ \ \ \ \ \ \ \ \ \ \text{if }j\leq m-i
\end{array}
\right.   &  =g\left(  v_{j}\right)  \ \ \ \ \ \ \ \ \ \ \left(  \text{since
}j=m+1>m=m-i\right) \label{pf.R.implicit.converse.trueproof.base.2}\\
&  =f\left(  1\right)  .\nonumber
\end{align}
Now,%
\[
\underbrace{\left(  \xi_{i}f\right)  }_{=f}\left(  v_{j}\right)  =f\left(
\underbrace{v_{j}}_{=1}\right)  =f\left(  1\right)  =\left\{
\begin{array}
[c]{c}%
g\left(  v_{j}\right)  ,\ \ \ \ \ \ \ \ \ \ \text{if }j>m-i;\\
f\left(  v_{j}\right)  ,\ \ \ \ \ \ \ \ \ \ \text{if }j\leq m-i
\end{array}
\right.
\]
(by (\ref{pf.R.implicit.converse.trueproof.base.2})). Hence,
(\ref{pf.R.implicit.converse.trueproof.base.goal}) is proven.
\par
Now, forget that we fixed $j$. We thus have shown that every $j\in\left\{
0,1,...,m+1\right\}  $ satisfies%
\[
\left(  \xi_{i}f\right)  \left(  v_{j}\right)  =\left\{
\begin{array}
[c]{c}%
g\left(  v_{j}\right)  ,\ \ \ \ \ \ \ \ \ \ \text{if }j>m-i;\\
f\left(  v_{j}\right)  ,\ \ \ \ \ \ \ \ \ \ \text{if }j\leq m-i
\end{array}
\right.  .
\]
Hence, Assertion A is proven, qed.}. Thus, the induction base is complete.

\textit{Induction step:} Let $I\in\left\{  0,1,...,m-1\right\}  $. Assume that
Assertion A is satisfied for $i=I$. We need to prove that Assertion A is
satisfied for $i=I+1$.

We know that Assertion A is satisfied for $i=I$. In other words, the
$\mathbb{K}$-labelling $\xi_{I}f\in\mathbb{K}^{\widehat{P}}$ is well-defined,
and every $j\in\left\{  0,1,...,m+1\right\}  $ satisfies%
\begin{equation}
\left(  \xi_{I}f\right)  \left(  v_{j}\right)  =\left\{
\begin{array}
[c]{c}%
g\left(  v_{j}\right)  ,\ \ \ \ \ \ \ \ \ \ \text{if }j>m-I;\\
f\left(  v_{j}\right)  ,\ \ \ \ \ \ \ \ \ \ \text{if }j\leq m-I
\end{array}
\right.  . \label{pf.R.implicit.converse.trueproof.step.1}%
\end{equation}

Notice that $I\in\left\{  0,1,...,m-1\right\}  $, so that $m-I\in\left\{
1,2,...,m\right\}  $. Thus, $v_{m-I}\in\left\{  v_{1},v_{2},...,v_{m}\right\}
=P$.

Let $v=v_{m-I}$. Then, $v=v_{m-I}\in P$ and $\xi_{I+1}=T_{v}\circ\xi_{I}%
$\ \ \ \ \footnote{\textit{Proof.} We have $\xi_{I}=T_{v_{m-I+1}}\circ
T_{v_{m-I+2}}\circ...\circ T_{v_{m}}$ (by the definition of $\xi_{I}$), so
that $T_{v_{m-I+1}}\circ T_{v_{m-I+2}}\circ...\circ T_{v_{m}}=\xi_{I}$. But
the definition of $\xi_{I+1}$ yields%
\begin{align*}
\xi_{I+1}  &  =T_{v_{m-\left(  I+1\right)  +1}}\circ T_{v_{m-\left(
I+1\right)  +2}}\circ...\circ T_{v_{m}}=T_{v_{m-I}}\circ T_{v_{m-I+1}}%
\circ...\circ T_{v_{m}}\\
&  =\underbrace{T_{v_{m-I}}}_{\substack{=T_{v}\\\text{(since }v_{m-I}%
=v\\\text{(since }v=v_{m-I}\text{))}}}\circ\underbrace{\left(  T_{v_{m-I+1}%
}\circ T_{v_{m-I+2}}\circ...\circ T_{v_{m}}\right)  }_{=\xi_{I}}=T_{v}\circ
\xi_{I},
\end{align*}
qed.}. We know that (\ref{prop.R.implicit.converse.eq}) holds, and thus the
right hand side of (\ref{prop.R.implicit.converse.eq}) is well-defined. In
other words,%
\begin{equation}
\dfrac{1}{f\left(  v\right)  }\cdot\dfrac{\sum\limits_{\substack{u\in
\widehat{P};\\u\lessdot v}}f\left(  u\right)  }{\sum\limits_{\substack{u\in
\widehat{P};\\u\gtrdot v}}\dfrac{1}{g\left(  u\right)  }}\text{ is
well-defined.} \label{pf.R.implicit.converse.trueproof.step.wd.use}%
\end{equation}
In particular, $\sum\limits_{\substack{u\in\widehat{P};\\u\gtrdot v}}\dfrac
{1}{g\left(  u\right)  }$ is well-defined. But%
\begin{equation}
\left(  \xi_{I}f\right)  \left(  u\right)  =f\left(  u\right)
\ \ \ \ \ \ \ \ \ \ \text{for every }u\in\widehat{P}\text{ satisfying
}u\lessdot v. \label{pf.R.implicit.converse.trueproof.step.2a}%
\end{equation}
\footnote{\textit{Proof of (\ref{pf.R.implicit.converse.trueproof.step.2a}):}
Let $u\in\widehat{P}$ be such that $u\lessdot v$. We have $u\in\widehat{P}%
=\left\{  v_{0},v_{1},v_{2},...,v_{m},v_{m+1}\right\}  $ (by
(\ref{pf.R.implicit.converse.trueproof.Phat})). Thus, there exists a
$j\in\left\{  0,1,...,m+1\right\}  $ such that $u=v_{j}$. Consider this $j$.
We have $u=v_{j}$, so that $v_{j}=u\lessdot v=v_{m-I}$, hence $v_{j}<v_{m-I}$.
Thus, $j<m-I$ (by (\ref{pf.R.implicit.converse.trueproof.linextP01}), applied
to $j$ and $m-I$ instead of $i$ and $j$). This entails $j\leq m-I$. Now,
(\ref{pf.R.implicit.converse.trueproof.step.1}) yields%
\[
\left(  \xi_{I}f\right)  \left(  v_{j}\right)  =\left\{
\begin{array}
[c]{c}%
g\left(  v_{j}\right)  ,\ \ \ \ \ \ \ \ \ \ \text{if }j>m-I;\\
f\left(  v_{j}\right)  ,\ \ \ \ \ \ \ \ \ \ \text{if }j\leq m-I
\end{array}
\right.  =f\left(  v_{j}\right)  \ \ \ \ \ \ \ \ \ \ \left(  \text{since
}j\leq m-I\right)  .
\]
Since $u=v_{j}$, this rewrites as $\left(  \xi_{I}f\right)  \left(  u\right)
=f\left(  u\right)  $. This proves
(\ref{pf.R.implicit.converse.trueproof.step.2a}).} Thus,
\begin{equation}
\sum\limits_{\substack{u\in\widehat{P};\\u\lessdot v}}\underbrace{\left(
\xi_{I}f\right)  \left(  u\right)  }_{\substack{=f\left(  u\right)
\\\text{(by (\ref{pf.R.implicit.converse.trueproof.step.2a}))}}}=\sum
\limits_{\substack{u\in\widehat{P};\\u\lessdot v}}f\left(  u\right)  .
\label{pf.R.implicit.converse.trueproof.step.2as}%
\end{equation}
Furthermore,%
\begin{equation}
\left(  \xi_{I}f\right)  \left(  u\right)  =g\left(  u\right)
\ \ \ \ \ \ \ \ \ \ \text{for every }u\in\widehat{P}\text{ satisfying
}u\gtrdot v. \label{pf.R.implicit.converse.trueproof.step.2b}%
\end{equation}
\footnote{\textit{Proof of (\ref{pf.R.implicit.converse.trueproof.step.2b}):}
Let $u\in\widehat{P}$ be such that $u\gtrdot v$. We have $u\in\widehat{P}%
=\left\{  v_{0},v_{1},v_{2},...,v_{m},v_{m+1}\right\}  $ (by
(\ref{pf.R.implicit.converse.trueproof.Phat})). Thus, there exists a
$j\in\left\{  0,1,...,m+1\right\}  $ such that $u=v_{j}$. Consider this $j$.
We have $u=v_{j}$, so that $v_{j}=u\gtrdot v=v_{m-I}$, hence $v_{j}>v_{m-I}$,
so that $v_{m-I}<v_{j}$. Thus, $m-I<j$ (by
(\ref{pf.R.implicit.converse.trueproof.linextP01}), applied to $m-I$ and $j$
instead of $i$ and $j$). In other words, $j>m-I$. Now,
(\ref{pf.R.implicit.converse.trueproof.step.1}) yields%
\[
\left(  \xi_{I}f\right)  \left(  v_{j}\right)  =\left\{
\begin{array}
[c]{c}%
g\left(  v_{j}\right)  ,\ \ \ \ \ \ \ \ \ \ \text{if }j>m-I;\\
f\left(  v_{j}\right)  ,\ \ \ \ \ \ \ \ \ \ \text{if }j\leq m-I
\end{array}
\right.  =g\left(  v_{j}\right)  \ \ \ \ \ \ \ \ \ \ \left(  \text{since
}j>m-I\right)  .
\]
Since $u=v_{j}$, this rewrites as $\left(  \xi_{I}f\right)  \left(  u\right)
=g\left(  u\right)  $. This proves
(\ref{pf.R.implicit.converse.trueproof.step.2b}).} Thus,%
\begin{equation}
\sum\limits_{\substack{u\in\widehat{P};\\u\gtrdot v}}\underbrace{\dfrac
{1}{\left(  \xi_{I}f\right)  \left(  u\right)  }}_{\substack{=\dfrac
{1}{g\left(  u\right)  }\\\text{(because
(\ref{pf.R.implicit.converse.trueproof.step.2b})}\\\text{yields }\left(
\xi_{I}f\right)  \left(  u\right)  =g\left(  u\right)  \text{)}}%
}=\sum\limits_{\substack{u\in\widehat{P};\\u\gtrdot v}}\dfrac{1}{g\left(
u\right)  } \label{pf.R.implicit.converse.trueproof.step.2bs}%
\end{equation}
(and this is well-defined because we know that $\sum\limits_{\substack{u\in
\widehat{P};\\u\gtrdot v}}\dfrac{1}{g\left(  u\right)  }$ is well-defined).
Finally,%
\begin{equation}
\left(  \xi_{I}f\right)  \left(  v\right)  =f\left(  v\right)  .
\label{pf.R.implicit.converse.trueproof.step.2v}%
\end{equation}
\footnote{\textit{Proof of (\ref{pf.R.implicit.converse.trueproof.step.2v}):}
The equality (\ref{pf.R.implicit.converse.trueproof.step.1}) (applied to
$j=m-I$) yields%
\[
\left(  \xi_{I}f\right)  \left(  v_{m-I}\right)  =\left\{
\begin{array}
[c]{c}%
g\left(  v_{m-I}\right)  ,\ \ \ \ \ \ \ \ \ \ \text{if }j>m-I;\\
f\left(  v_{m-I}\right)  ,\ \ \ \ \ \ \ \ \ \ \text{if }j\leq m-I
\end{array}
\right.  =f\left(  v_{m-I}\right)  \ \ \ \ \ \ \ \ \ \ \left(  \text{since
}m-I\leq m-I\right)  .
\]
Since $v=v_{m-I}$, this rewrites as $\left(  \xi_{I}f\right)  \left(
u\right)  =f\left(  u\right)  $. This proves
(\ref{pf.R.implicit.converse.trueproof.step.2v}).}

Now, we are going to prove that $T_{v}\left(  \xi_{I}f\right)  $ is
well-defined. In fact, let us recall that $T_{v}\left(  \xi_{I}f\right)  $ is
defined by%
\begin{equation}
\left(  T_{v}\left(  \xi_{I}f\right)  \right)  \left(  w\right)  =\left\{
\begin{array}
[c]{l}%
\left(  \xi_{I}f\right)  \left(  w\right)  ,\ \ \ \ \ \ \ \ \ \ \text{if
}w\neq v;\\
\dfrac{1}{\left(  \xi_{I}f\right)  \left(  v\right)  }\cdot\dfrac
{\sum\limits_{\substack{u\in\widehat{P};\\u\lessdot v}}\left(  \xi
_{I}f\right)  \left(  u\right)  }{\sum\limits_{\substack{u\in\widehat{P}%
;\\u\gtrdot v}}\dfrac{1}{\left(  \xi_{I}f\right)  \left(  u\right)  }%
},\ \ \ \ \ \ \ \ \ \ \text{if }w=v
\end{array}
\right.  \ \ \ \ \ \ \ \ \ \ \text{for all }w\in\widehat{P}.
\label{pf.R.implicit.converse.trueproof.step.wd1}%
\end{equation}
Hence, in order to prove that $T_{v}\left(  \xi_{I}f\right)  $ is
well-defined, we need to show that the right hand side of
(\ref{pf.R.implicit.converse.trueproof.step.wd1}) is well-defined for every
$w\in\widehat{P}$.

For every $w\in\widehat{P}$, the right hand side of
(\ref{pf.R.implicit.converse.trueproof.step.wd1}) is%
\begin{align*}
&  \left\{
\begin{array}
[c]{l}%
\left(  \xi_{I}f\right)  \left(  w\right)  ,\ \ \ \ \ \ \ \ \ \ \text{if
}w\neq v;\\
\dfrac{1}{\left(  \xi_{I}f\right)  \left(  v\right)  }\cdot\dfrac
{\sum\limits_{\substack{u\in\widehat{P};\\u\lessdot v}}\left(  \xi
_{I}f\right)  \left(  u\right)  }{\sum\limits_{\substack{u\in\widehat{P}%
;\\u\gtrdot v}}\dfrac{1}{\left(  \xi_{I}f\right)  \left(  u\right)  }%
},\ \ \ \ \ \ \ \ \ \ \text{if }w=v
\end{array}
\right. \\
&  =\left\{
\begin{array}
[c]{l}%
\left(  \xi_{I}f\right)  \left(  w\right)  ,\ \ \ \ \ \ \ \ \ \ \text{if
}w\neq v;\\
\dfrac{1}{f\left(  v\right)  }\cdot\dfrac{\sum\limits_{\substack{u\in
\widehat{P};\\u\lessdot v}}f\left(  u\right)  }{\sum\limits_{\substack{u\in
\widehat{P};\\u\gtrdot v}}\dfrac{1}{g\left(  u\right)  }}%
,\ \ \ \ \ \ \ \ \ \ \text{if }w=v
\end{array}
\right. \\
&  \ \ \ \ \ \ \ \ \ \ \left(  \text{by
(\ref{pf.R.implicit.converse.trueproof.step.2v}),
(\ref{pf.R.implicit.converse.trueproof.step.2as}) and
(\ref{pf.R.implicit.converse.trueproof.step.2bs})}\right)  .
\end{align*}
This is well-defined (since $\left(  \xi_{I}f\right)  \left(  w\right)  $ is
well-defined (because $\xi_{I}f$ is well-defined) and since $\dfrac
{1}{f\left(  v\right)  }\cdot\dfrac{\sum\limits_{\substack{u\in\widehat{P}%
;\\u\lessdot v}}f\left(  u\right)  }{\sum\limits_{\substack{u\in
\widehat{P};\\u\gtrdot v}}\dfrac{1}{g\left(  u\right)  }}$ is well-defined
(because of (\ref{pf.R.implicit.converse.trueproof.step.wd.use}))). Thus, we
have shown that the right hand side of
(\ref{pf.R.implicit.converse.trueproof.step.wd1}) is well-defined for every
$w\in\widehat{P}$. Since (\ref{pf.R.implicit.converse.trueproof.step.wd1}) is
the definition of $T_{v}\left(  \xi_{I}f\right)  $, this yields that
$T_{v}\left(  \xi_{I}f\right)  $ is well-defined. In other words, $\xi_{I+1}f$
is well-defined (because $\underbrace{\xi_{I+1}}_{=T_{v}\circ\xi_{I}}f=\left(
T_{v}\circ\xi_{I}\right)  f=T_{v}\left(  \xi_{I}f\right)  $).

We are now going to prove that every $j\in\left\{  0,1,...,m+1\right\}  $
satisfies%
\begin{equation}
\left(  \xi_{I+1}f\right)  \left(  v_{j}\right)  =\left\{
\begin{array}
[c]{c}%
g\left(  v_{j}\right)  ,\ \ \ \ \ \ \ \ \ \ \text{if }j>m-\left(  I+1\right)
;\\
f\left(  v_{j}\right)  ,\ \ \ \ \ \ \ \ \ \ \text{if }j\leq m-\left(
I+1\right)
\end{array}
\right.  . \label{pf.R.implicit.converse.trueproof.step.goal}%
\end{equation}

\textit{Proof of (\ref{pf.R.implicit.converse.trueproof.step.goal}):} Let
$j\in\left\{  0,1,...,m+1\right\}  $. We need to prove that
(\ref{pf.R.implicit.converse.trueproof.step.goal}) holds.

We distinguish between two cases:

\textit{Case 1:} We have $v_{j}\neq v$.

\textit{Case 2:} We have $v_{j}=v$.

Let us deal with Case 1 first. In this case, we have $v_{j}\neq v$. Hence,
$v_{j}\neq v=v_{m-I}$, so that $j\neq m-I$. Notice that $\underbrace{\xi
_{I+1}}_{=T_{v}\circ\xi_{I}}f=\left(  T_{v}\circ\xi_{I}\right)  f=T_{v}\left(
\xi_{I}f\right)  $. Hence,%
\begin{align}
\left(  \xi_{I+1}f\right)  \left(  v_{j}\right)   &  =\left(  T_{v}\left(
\xi_{I}f\right)  \right)  \left(  v_{j}\right) \nonumber\\
&  =\left\{
\begin{array}
[c]{l}%
\left(  \xi_{I}f\right)  \left(  v_{j}\right)  ,\ \ \ \ \ \ \ \ \ \ \text{if
}v_{j}\neq v;\\
\dfrac{1}{\left(  \xi_{I}f\right)  \left(  v\right)  }\cdot\dfrac
{\sum\limits_{\substack{u\in\widehat{P};\\u\lessdot v}}\left(  \xi
_{I}f\right)  \left(  u\right)  }{\sum\limits_{\substack{u\in\widehat{P}%
;\\u\gtrdot v}}\dfrac{1}{\left(  \xi_{I}f\right)  \left(  u\right)  }%
},\ \ \ \ \ \ \ \ \ \ \text{if }v_{j}=v
\end{array}
\right.  \ \ \ \ \ \ \ \ \ \ \left(  \text{by
(\ref{pf.R.implicit.converse.trueproof.step.wd1}), applied to }w=v_{j}\right)
\nonumber\\
&  =\left(  \xi_{I}f\right)  \left(  v_{j}\right)  \ \ \ \ \ \ \ \ \ \ \left(
\text{since }v_{j}\neq v\right) \nonumber\\
&  =\left\{
\begin{array}
[c]{c}%
g\left(  v_{j}\right)  ,\ \ \ \ \ \ \ \ \ \ \text{if }j>m-I;\\
f\left(  v_{j}\right)  ,\ \ \ \ \ \ \ \ \ \ \text{if }j\leq m-I
\end{array}
\right.  \ \ \ \ \ \ \ \ \ \ \left(  \text{by
(\ref{pf.R.implicit.converse.trueproof.step.1})}\right)  .
\label{pf.R.implicit.converse.trueproof.step.c1}%
\end{align}

Now, we have $j\neq m-I$. Hence, we must be in one of the following two subcases:

\textit{Subcase 1.1:} We have $j<m-I$.

\textit{Subcase 1.2:} We have $j>m-I$.

Let us first consider Subcase 1.1. In this subcase, we have $j<m-I$. Hence,
$j\leq m-I-1$ (since $j$ and $m-I$ are integers), so that $j\leq
m-I-1=m-\left(  I+1\right)  $. Also, $j\leq m-I$ (since $j<m-I$). Thus,
(\ref{pf.R.implicit.converse.trueproof.step.c1}) becomes%
\[
\left(  \xi_{I+1}f\right)  \left(  v_{j}\right)  =\left\{
\begin{array}
[c]{c}%
g\left(  v_{j}\right)  ,\ \ \ \ \ \ \ \ \ \ \text{if }j>m-I;\\
f\left(  v_{j}\right)  ,\ \ \ \ \ \ \ \ \ \ \text{if }j\leq m-I
\end{array}
\right.  =f\left(  v_{j}\right)  \ \ \ \ \ \ \ \ \ \ \left(  \text{since
}j\leq m-I\right)  .
\]
Compared with%
\[
\left\{
\begin{array}
[c]{c}%
g\left(  v_{j}\right)  ,\ \ \ \ \ \ \ \ \ \ \text{if }j>m-\left(  I+1\right)
;\\
f\left(  v_{j}\right)  ,\ \ \ \ \ \ \ \ \ \ \text{if }j\leq m-\left(
I+1\right)
\end{array}
\right.  =f\left(  v_{j}\right)  \ \ \ \ \ \ \ \ \ \ \left(  \text{since
}j\leq m-\left(  I+1\right)  \right)  ,
\]
this yields%
\[
\left(  \xi_{I+1}f\right)  \left(  v_{j}\right)  =\left\{
\begin{array}
[c]{c}%
g\left(  v_{j}\right)  ,\ \ \ \ \ \ \ \ \ \ \text{if }j>m-\left(  I+1\right)
;\\
f\left(  v_{j}\right)  ,\ \ \ \ \ \ \ \ \ \ \text{if }j\leq m-\left(
I+1\right)
\end{array}
\right.  .
\]
Hence, (\ref{pf.R.implicit.converse.trueproof.step.goal}) is proven in Subcase 1.1.

Let us now proceed to Subcase 1.2. In this subcase, we have $j>m-I$. Hence,
$j>m-I>m-I-1=m-\left(  I+1\right)  $. Now,
(\ref{pf.R.implicit.converse.trueproof.step.c1}) becomes%
\[
\left(  \xi_{I+1}f\right)  \left(  v_{j}\right)  =\left\{
\begin{array}
[c]{c}%
g\left(  v_{j}\right)  ,\ \ \ \ \ \ \ \ \ \ \text{if }j>m-I;\\
f\left(  v_{j}\right)  ,\ \ \ \ \ \ \ \ \ \ \text{if }j\leq m-I
\end{array}
\right.  =g\left(  v_{j}\right)  \ \ \ \ \ \ \ \ \ \ \left(  \text{since
}j>m-I\right)  .
\]
Compared with%
\[
\left\{
\begin{array}
[c]{c}%
g\left(  v_{j}\right)  ,\ \ \ \ \ \ \ \ \ \ \text{if }j>m-\left(  I+1\right)
;\\
f\left(  v_{j}\right)  ,\ \ \ \ \ \ \ \ \ \ \text{if }j\leq m-\left(
I+1\right)
\end{array}
\right.  =g\left(  v_{j}\right)  \ \ \ \ \ \ \ \ \ \ \left(  \text{since
}j>m-\left(  I+1\right)  \right)  ,
\]
this yields%
\[
\left(  \xi_{I+1}f\right)  \left(  v_{j}\right)  =\left\{
\begin{array}
[c]{c}%
g\left(  v_{j}\right)  ,\ \ \ \ \ \ \ \ \ \ \text{if }j>m-\left(  I+1\right)
;\\
f\left(  v_{j}\right)  ,\ \ \ \ \ \ \ \ \ \ \text{if }j\leq m-\left(
I+1\right)
\end{array}
\right.  .
\]
Hence, (\ref{pf.R.implicit.converse.trueproof.step.goal}) is proven in Subcase 1.2.

We have thus proven (\ref{pf.R.implicit.converse.trueproof.step.goal}) in each
of the two Subcases 1.1 and 1.2. Since these two Subcases cover all of Case 1,
this yields that (\ref{pf.R.implicit.converse.trueproof.step.goal}) always
holds in Case 1.

Now, let us consider Case 2. In this case, we have $v_{j}=v$. Hence,
$j=m-I$\ \ \ \ \footnote{\textit{Proof.} Assume (for the sake of
contradiction) that $j=0$. Then, $v_{j}=v_{0}=0$ (by the definition of $v_{0}%
$). But $v_{j}=v$, so that $v=v_{j}=0\notin P$, contradicting $v\in P$. This
contradiction shows that our assumption (that $j=0$) was wrong. Hence, we
cannot have $j=0$. We thus have $j\neq0$.
\par
Assume (for the sake of contradiction) that $j=m+1$. Then, $v_{j}=v_{m+1}=1$
(by the definition of $v_{m+1}$). But $v_{j}=v$, so that $v=v_{j}=m+1\notin
P$, contradicting $v\in P$. This contradiction shows that our assumption (that
$j=m+1$) was wrong. Hence, we cannot have $j=m+1$. We thus have $j\neq m+1$.
\par
Combining $j\neq0$ with $j\neq m+1$, we obtain $j\notin\left\{  0,m+1\right\}
$. Since $j\in\left\{  0,1,...,m+1\right\}  $ but $j\notin\left\{
0,m+1\right\}  $, we must have $j\in\left\{  0,1,...,m+1\right\}
\setminus\left\{  0,m+1\right\}  =\left\{  1,2,...,m\right\}  $.
\par
Now, recall that $\left(  v_{1},v_{2},...,v_{m}\right)  $ is a list of the
elements of $P$ such that every element of $P$ occurs exactly once in this
list. In particular, every element of $P$ occurs exactly once in the list
$\left(  v_{1},v_{2},...,v_{m}\right)  $. Hence, the element $v$ occurs
exactly once in the list $\left(  v_{1},v_{2},...,v_{m}\right)  $ (since $v$
is an element of $P$). In other words, there exists exactly one $k\in\left\{
1,2,...,m\right\}  $ such that $k=v$. In particular, there exists at most one
$k\in\left\{  1,2,...,m\right\}  $ such that $k=v$. In other words,%
\begin{equation}
\left(
\begin{array}
[c]{c}%
\text{if }k_{1}\text{ and }k_{2}\text{ are two elements }k\in\left\{
1,2,...,m\right\}  \text{ such that }k=v\text{,}\\
\text{then }k_{1}=k_{2}%
\end{array}
\right)  . \label{pf.R.implicit.converse.trueproof.step.k1k2}%
\end{equation}
\par
We know that $j\in\left\{  1,2,...,m\right\}  $ and $v_{j}=v$. In other words,
$j$ is an element $k\in\left\{  1,2,...,m\right\}  $ such that $k=v$. Also,
$m-I\in\left\{  1,2,...,m\right\}  $ and $v_{m-I}=v$. In other words, $m-I$ is
an element $k\in\left\{  1,2,...,m\right\}  $ such that $k=v$. Hence,
(\ref{pf.R.implicit.converse.trueproof.step.k1k2}) (applied to $k_{1}=j$ and
$k_{2}=m-I$) yields $j=m-I$, qed.}. From
(\ref{pf.R.implicit.converse.trueproof.step.1}), we obtain%
\[
\left(  \xi_{I}f\right)  \left(  v_{j}\right)  =\left\{
\begin{array}
[c]{c}%
g\left(  v_{j}\right)  ,\ \ \ \ \ \ \ \ \ \ \text{if }j>m-I;\\
f\left(  v_{j}\right)  ,\ \ \ \ \ \ \ \ \ \ \text{if }j\leq m-I
\end{array}
\right.  =f\left(  v_{j}\right)
\]
(since $j=m-I\leq m-I$). Notice that $\underbrace{\xi_{I+1}}_{=T_{v}\circ
\xi_{I}}f=\left(  T_{v}\circ\xi_{I}\right)  f=T_{v}\left(  \xi_{I}f\right)  $.
Hence,%
\begin{align*}
\left(  \xi_{I+1}f\right)  \left(  v_{j}\right)   &  =\left(  T_{v}\left(
\xi_{I}f\right)  \right)  \left(  v_{j}\right) \\
&  =\left\{
\begin{array}
[c]{l}%
\left(  \xi_{I}f\right)  \left(  v_{j}\right)  ,\ \ \ \ \ \ \ \ \ \ \text{if
}v_{j}\neq v;\\
\dfrac{1}{\left(  \xi_{I}f\right)  \left(  v\right)  }\cdot\dfrac
{\sum\limits_{\substack{u\in\widehat{P};\\u\lessdot v}}\left(  \xi
_{I}f\right)  \left(  u\right)  }{\sum\limits_{\substack{u\in\widehat{P}%
;\\u\gtrdot v}}\dfrac{1}{\left(  \xi_{I}f\right)  \left(  u\right)  }%
},\ \ \ \ \ \ \ \ \ \ \text{if }v_{j}=v
\end{array}
\right.  \ \ \ \ \ \ \ \ \ \ \left(  \text{by
(\ref{pf.R.implicit.converse.trueproof.step.wd1}), applied to }w=v_{j}\right)
\\
&  =\dfrac{1}{\left(  \xi_{I}f\right)  \left(  v\right)  }\cdot\dfrac
{\sum\limits_{\substack{u\in\widehat{P};\\u\lessdot v}}\left(  \xi
_{I}f\right)  \left(  u\right)  }{\sum\limits_{\substack{u\in\widehat{P}%
;\\u\gtrdot v}}\dfrac{1}{\left(  \xi_{I}f\right)  \left(  u\right)  }%
}\ \ \ \ \ \ \ \ \ \ \left(  \text{since }v_{j}=v\right) \\
&  =\dfrac{1}{f\left(  v\right)  }\cdot\dfrac{\sum\limits_{\substack{u\in
\widehat{P};\\u\lessdot v}}f\left(  u\right)  }{\sum\limits_{\substack{u\in
\widehat{P};\\u\gtrdot v}}\dfrac{1}{g\left(  u\right)  }}%
\ \ \ \ \ \ \ \ \ \ \left(  \text{by
(\ref{pf.R.implicit.converse.trueproof.step.2v}),
(\ref{pf.R.implicit.converse.trueproof.step.2as}) and
(\ref{pf.R.implicit.converse.trueproof.step.2bs})}\right) \\
&  =g\left(  v\right)  \ \ \ \ \ \ \ \ \ \ \left(  \text{by
(\ref{prop.R.implicit.converse.eq})}\right)  .
\end{align*}
Compared with%
\begin{align*}
&  \left\{
\begin{array}
[c]{c}%
g\left(  v_{j}\right)  ,\ \ \ \ \ \ \ \ \ \ \text{if }j>m-\left(  I+1\right)
;\\
f\left(  v_{j}\right)  ,\ \ \ \ \ \ \ \ \ \ \text{if }j\leq m-\left(
I+1\right)
\end{array}
\right. \\
&  =g\left(  \underbrace{v_{j}}_{=v}\right)  \ \ \ \ \ \ \ \ \ \ \left(
\text{since }j=m-I>m-I-1=m-\left(  I+1\right)  \right) \\
&  =g\left(  v\right)  ,
\end{align*}
this yields
\[
\left(  \xi_{I+1}f\right)  \left(  v_{j}\right)  =\left\{
\begin{array}
[c]{c}%
g\left(  v_{j}\right)  ,\ \ \ \ \ \ \ \ \ \ \text{if }j>m-\left(  I+1\right)
;\\
f\left(  v_{j}\right)  ,\ \ \ \ \ \ \ \ \ \ \text{if }j\leq m-\left(
I+1\right)
\end{array}
\right.  .
\]
Hence, (\ref{pf.R.implicit.converse.trueproof.step.goal}) is proven in Case 2.

We thus have proven (\ref{pf.R.implicit.converse.trueproof.step.goal}) in each
of the two Cases 1 and 2. Since these two Cases cover all possibilities, this
yields that (\ref{pf.R.implicit.converse.trueproof.step.goal}) always holds.
In other words, the proof of (\ref{pf.R.implicit.converse.trueproof.step.goal}%
) is complete.

Altogether, we thus have proven that the $\mathbb{K}$-labelling $\xi_{I+1}%
f\in\mathbb{K}^{\widehat{P}}$ is well-defined, and every $j\in\left\{
0,1,...,m+1\right\}  $ satisfies%
\[
\left(  \xi_{I+1}f\right)  \left(  v_{j}\right)  =\left\{
\begin{array}
[c]{c}%
g\left(  v_{j}\right)  ,\ \ \ \ \ \ \ \ \ \ \text{if }j>m-\left(  I+1\right)
;\\
f\left(  v_{j}\right)  ,\ \ \ \ \ \ \ \ \ \ \text{if }j\leq m-\left(
I+1\right)
\end{array}
\right.  .
\]
In other words, Assertion A is satisfied for $i=I+1$. This completes the
induction step. The induction proof of Assertion A is thus complete.

Now, we can apply Assertion A to $i=m$ (because Assertion A is proven for
every $i\in\left\{  0,1,...,m\right\}  $). As a result, we conclude that the
$\mathbb{K}$-labelling $\xi_{m}f\in\mathbb{K}^{\widehat{P}}$ is well-defined,
and that every $j\in\left\{  0,1,...,m+1\right\}  $ satisfies%
\begin{equation}
\left(  \xi_{m}f\right)  \left(  v_{j}\right)  =\left\{
\begin{array}
[c]{c}%
g\left(  v_{j}\right)  ,\ \ \ \ \ \ \ \ \ \ \text{if }j>m-m;\\
f\left(  v_{j}\right)  ,\ \ \ \ \ \ \ \ \ \ \text{if }j\leq m-m
\end{array}
\right.  . \label{pf.R.implicit.converse.trueproof.almostthere}%
\end{equation}

We know that the $\mathbb{K}$-labelling $\xi_{m}f\in\mathbb{K}^{\widehat{P}}$
is well-defined. In other words, $\mathbb{K}$-labelling $Rf\in\mathbb{K}%
^{\widehat{P}}$ is well-defined (since $\xi_{m}=R$). We shall now show that
$g=Rf$.

In fact, let $v\in\widehat{P}$ be arbitrary. We are going to prove that
$\left(  Rf\right)  \left(  v\right)  =g\left(  v\right)  $.

We have $v\in\widehat{P}=\left\{  v_{0},v_{1},v_{2},...,v_{m},v_{m+1}\right\}
$ (by (\ref{pf.R.implicit.converse.trueproof.Phat})). Hence, there exists a
$j\in\left\{  0,1,...,m+1\right\}  $ such that $v=v_{j}$. Consider this $j$.
If $j=0$, then $\left(  Rf\right)  \left(  v\right)  =g\left(  v\right)  $ is
obviously true\footnote{\textit{Proof.} Assume that $j=0$. Then, $v_{j}%
=v_{0}=0$ (by the definition of $v_{0}$), so that $v=v_{j}=0$. Hence, $\left(
Rf\right)  \left(  \underbrace{v}_{=0}\right)  =\left(  Rf\right)  \left(
0\right)  =f\left(  0\right)  $ (by Proposition \ref{prop.R.implicit.01}).
Compared with $g\left(  \underbrace{v}_{=0}\right)  =g\left(  0\right)
=f\left(  0\right)  $ (since $f\left(  0\right)  =g\left(  0\right)  $), this
yields $\left(  Rf\right)  \left(  v\right)  =g\left(  v\right)  $, qed.}.
Hence, for the rest of the proof of $\left(  Rf\right)  \left(  v\right)
=g\left(  v\right)  $, we can WLOG assume that we don't have $j=0$. Assume this.

We have $j\in\left\{  0,1,...,m+1\right\}  $ but we don't have $j=0$. Thus,
$j\in\left\{  1,2,...,m+1\right\}  $, so that $j>0=m-m$. Now,
(\ref{pf.R.implicit.converse.trueproof.almostthere}) becomes%
\[
\left(  \xi_{m}f\right)  \left(  v_{j}\right)  =\left\{
\begin{array}
[c]{c}%
g\left(  v_{j}\right)  ,\ \ \ \ \ \ \ \ \ \ \text{if }j>m-m;\\
f\left(  v_{j}\right)  ,\ \ \ \ \ \ \ \ \ \ \text{if }j\leq m-m
\end{array}
\right.  =g\left(  v_{j}\right)  \ \ \ \ \ \ \ \ \ \ \left(  \text{since
}j>m-m\right)  .
\]
Now, $\xi_{m}=R$ and $v=v_{j}$, so that $\left(  \underbrace{\xi_{m}}%
_{=R}f\right)  \left(  \underbrace{v_{j}}_{=v}\right)  =\left(  Rf\right)
\left(  v\right)  $. Hence,%
\[
\left(  Rf\right)  \left(  v\right)  =\left(  \xi_{m}f\right)  \left(
v_{j}\right)  =g\left(  \underbrace{v_{j}}_{=v}\right)  =g\left(  v\right)  .
\]
We thus have proven $\left(  Rf\right)  \left(  v\right)  =g\left(  v\right)
$.

Now, forget that we fixed $v$. We thus have shown that $\left(  Rf\right)
\left(  v\right)  =g\left(  v\right)  $ for every $v\in\widehat{P}$. In other
words, $Rf=g$, so that $g=Rf$. This proves Proposition
\ref{prop.R.implicit.converse}.
\end{proof}

\begin{noncompile}
The following alternative proof of Proposition \ref{prop.R.implicit.converse}
is incomplete, since it does not check that $Rf$ is well-defined. It is
archived here solely for the purpose of not throwing it away.

\begin{proof}
[Incomplete proof of Proposition \ref{prop.R.implicit.converse} (sketched).]%
Let us first notice that $\left(  Rf\right)  \left(  1\right)  =f\left(
1\right)  $ (by Proposition \ref{prop.R.implicit.01}), so that $\left(
Rf\right)  \left(  1\right)  =f\left(  1\right)  =g\left(  1\right)  $.
Similarly, $\left(  Rf\right)  \left(  0\right)  =g\left(  0\right)  $.

Let $>$ denote the greater relation of the poset $\widehat{P}$. For every
$v\in\widehat{P}$, let $N\left(  v\right)  $ be the nonnegative integer
$\left\vert \left\{  x\in\widehat{P}\ \mid\ x>v\right\}  \right\vert $ (this
is well-defined, because $\widehat{P}$ is finite). We are going to show that
for every $i\in\mathbb{N}$, we have%
\begin{equation}
g\left(  v\right)  =\left(  Rf\right)  \left(  v\right)
\ \ \ \ \ \ \ \ \ \ \text{for every }v\in\widehat{P}\text{ satisfying
}N\left(  v\right)  \leq i. \label{pf.R.implicit.converse.1}%
\end{equation}

\textit{Proof of (\ref{pf.R.implicit.converse.1}):} We will prove
(\ref{pf.R.implicit.converse.1}) by induction over $i$:

\textit{Induction base:} Let $v\in\widehat{P}$ be such that $N\left(
v\right)  \leq0$. Thus, $\left\vert \left\{  x\in\widehat{P}\ \mid
\ x>v\right\}  \right\vert =N\left(  v\right)  =0$, hence $\left\{
x\in\widehat{P}\ \mid\ x>v\right\}  =\varnothing$. But the element $1$ of
$\widehat{P}$ is greater than any other element of $\widehat{P}$ (by the
construction of $\widehat{P}$). Hence, $1\geq v$. But we don't have $1>v$
(since otherwise, we would have $1>v$ and thus $1\in\left\{  x\in
\widehat{P}\ \mid\ x>v\right\}  =\varnothing$, which is absurd). Since $1\geq
v$ but not $1>v$, we must have $1=v$, so that $g\left(  \underbrace{v}%
_{=1}\right)  =g\left(  1\right)  $. Compared with $\left(  Rf\right)  \left(
\underbrace{v}_{=1}\right)  =\left(  Rf\right)  \left(  1\right)  =g\left(
1\right)  $, this yields $g\left(  v\right)  =\left(  Rf\right)  \left(
v\right)  $.

Now, forget that we fixed $v$. We thus have proven that $g\left(  v\right)
=\left(  Rf\right)  \left(  v\right)  $ for every $v\in\widehat{P}$ satisfying
$N\left(  v\right)  \leq0$. In other words, (\ref{pf.R.implicit.converse.1})
holds for $i=0$. This completes the induction base.

\textit{Induction step:} Let $I\in\mathbb{N}$. Assume that
(\ref{pf.R.implicit.converse.1}) is proven for $i=I$. We need to prove that
(\ref{pf.R.implicit.converse.1}) holds for $i=I+1$.

We know that (\ref{pf.R.implicit.converse.1}) is proven for $i=I$. In other
words, we have%
\begin{equation}
g\left(  v\right)  =\left(  Rf\right)  \left(  v\right)
\ \ \ \ \ \ \ \ \ \ \text{for every }v\in\widehat{P}\text{ satisfying
}N\left(  v\right)  \leq I. \label{pf.R.implicit.converse.indass}%
\end{equation}

Now, let $v\in\widehat{P}$ be such that $N\left(  v\right)  \leq I+1$. We are
going to prove that $g\left(  v\right)  =\left(  Rf\right)  \left(  v\right)
$.

If $v=1$, then $g\left(  v\right)  =\left(  Rf\right)  \left(  v\right)  $ is
true (because $g\left(  1\right)  =\left(  Rf\right)  \left(  1\right)  $).
Hence, for the rest of the proof of $g\left(  v\right)  =\left(  Rf\right)
\left(  v\right)  $, we can WLOG assume that we don't have $v=1$. Assume this.

If $v=0$, then $g\left(  v\right)  =\left(  Rf\right)  \left(  v\right)  $ is
true (because $g\left(  0\right)  =\left(  Rf\right)  \left(  0\right)  $).
Hence, for the rest of the proof of $g\left(  v\right)  =\left(  Rf\right)
\left(  v\right)  $, we can WLOG assume that we don't have $v=0$. Assume this.

We have $v\in\widehat{P}$, but we don't have $v=0$, and we don't have $v=1$.
Thus, $v\in\widehat{P}\setminus\left\{  0,1\right\}  =P$. Hence, we can apply
(\ref{prop.R.implicit.converse.eq}).

Let $u\in\widehat{P}$ be such that $u\gtrdot v$. Then, $u>v$. Hence, $\left\{
x\in\widehat{P}\ \mid\ x>u\right\}  $ is a proper subset of $\left\{
x\in\widehat{P}\ \mid\ x>v\right\}  $\ \ \ \ \footnote{\textit{Proof.} Since
$u\in\widehat{P}$ and $u>v$, we have $u\in\left\{  x\in\widehat{P}%
\ \mid\ x>v\right\}  $. But we don't have $u>u$, and therefore it is not true
that $u\in\widehat{P}$ and $u>u$. Hence, $u\notin\left\{  x\in\widehat{P}%
\ \mid\ x>u\right\}  $. Comparing this with $u\in\left\{  x\in\widehat{P}%
\ \mid\ x>v\right\}  $, we conclude that $\left\{  x\in\widehat{P}%
\ \mid\ x>u\right\}  \neq\left\{  x\in\widehat{P}\ \mid\ x>v\right\}  $.
\par
Let $y\in\left\{  x\in\widehat{P}\ \mid\ x>u\right\}  $. Then, $y\in
\widehat{P}$ and $y>u$. Since $y\in\widehat{P}$ and $y>u>v$, we have
$y\in\left\{  x\in\widehat{P}\ \mid\ x>v\right\}  $.
\par
Now, forget that we fixed $y$. We have thus shown that $y\in\left\{
x\in\widehat{P}\ \mid\ x>v\right\}  $ for every $y\in\left\{  x\in
\widehat{P}\ \mid\ x>u\right\}  $. In other words, $\left\{  x\in
\widehat{P}\ \mid\ x>u\right\}  \subseteq\left\{  x\in\widehat{P}%
\ \mid\ x>v\right\}  $. Combined with $\left\{  x\in\widehat{P}\ \mid
\ x>u\right\}  \neq\left\{  x\in\widehat{P}\ \mid\ x>v\right\}  $, this yields
that $\left\{  x\in\widehat{P}\ \mid\ x>u\right\}  $ is a proper subset of
$\left\{  x\in\widehat{P}\ \mid\ x>v\right\}  $, qed.}. Thus, $\left\vert
\left\{  x\in\widehat{P}\ \mid\ x>u\right\}  \right\vert <\left\vert \left\{
x\in\widehat{P}\ \mid\ x>v\right\}  \right\vert $ (since the cardinalities of
these two sets are finite). Now, by the definition of $N\left(  u\right)  $,
we have%
\[
N\left(  u\right)  =\left\vert \left\{  x\in\widehat{P}\ \mid\ x>u\right\}
\right\vert <\left\vert \left\{  x\in\widehat{P}\ \mid\ x>v\right\}
\right\vert =N\left(  v\right)  \leq I+1.
\]
Since $N\left(  u\right)  $ and $I$ are integers, this yields $N\left(
u\right)  \leq I$. Hence, we can apply (\ref{pf.R.implicit.converse.indass})
to $u$ instead of $v$. We conclude that $g\left(  u\right)  =\left(
Rf\right)  \left(  u\right)  $.

Now, forget that we fixed $u$. We have thus shown that%
\[
g\left(  u\right)  =\left(  Rf\right)  \left(  u\right)
\ \ \ \ \ \ \ \ \ \ \text{for every }u\in\widehat{P}\text{ satisfying
}u\gtrdot v.
\]
Hence,
\begin{equation}
\sum\limits_{\substack{u\in\widehat{P};\\u\gtrdot v}}\dfrac{1}{g\left(
u\right)  }=\sum\limits_{\substack{u\in\widehat{P};\\u\gtrdot v}}\dfrac
{1}{\left(  Rf\right)  \left(  u\right)  }. \label{pf.R.implicit.converse.3}%
\end{equation}

Now, from (\ref{prop.R.implicit.converse.eq}), we see that%
\begin{align*}
g\left(  v\right)   &  =\dfrac{1}{f\left(  v\right)  }\cdot\dfrac
{\sum\limits_{\substack{u\in\widehat{P};\\u\lessdot v}}f\left(  u\right)
}{\sum\limits_{\substack{u\in\widehat{P};\\u\gtrdot v}}\dfrac{1}{g\left(
u\right)  }}=\dfrac{1}{f\left(  v\right)  }\cdot\dfrac{\sum
\limits_{\substack{u\in\widehat{P};\\u\lessdot v}}f\left(  u\right)  }%
{\sum\limits_{\substack{u\in\widehat{P};\\u\gtrdot v}}\dfrac{1}{\left(
Rf\right)  \left(  u\right)  }}\ \ \ \ \ \ \ \ \ \ \left(  \text{by
(\ref{pf.R.implicit.converse.3})}\right) \\
&  =\left(  Rf\right)  \left(  v\right)  \ \ \ \ \ \ \ \ \ \ \left(  \text{by
Proposition \ref{prop.R.implicit.eq}}\right)  .
\end{align*}

Now, forget that we fixed $v$. We thus have proven that $g\left(  v\right)
=\left(  Rf\right)  \left(  v\right)  $ for every $v\in\widehat{P}$ satisfying
$N\left(  v\right)  \leq I+1$. In other words, (\ref{pf.R.implicit.converse.1}%
) holds for $i=I+1$. This completes the induction step. The induction proof of
(\ref{pf.R.implicit.converse.1}) is thus complete.

So (\ref{pf.R.implicit.converse.1}) is proven now. Hence, every $v\in
\widehat{P}$ satisfies $g\left(  v\right)  =\left(  Rf\right)  \left(
v\right)  $ (by (\ref{pf.R.implicit.converse.1}), applied to $i=N\left(
v\right)  $). In other words, $g=Rf$. This proves Proposition
\ref{prop.R.implicit.converse}.
\end{proof}
\end{noncompile}
\end{verlong}

\begin{vershort}
As an aside, at this point we could give an alternative proof of Corollary
\ref{cor.R.welldef}, foregoing the use of Proposition
\ref{prop.linext.transitive}. In fact, the proofs of Propositions
\ref{prop.R.implicit}, \ref{prop.R.implicit.01} and
\ref{prop.R.implicit.converse} only used that $R$ is a composition $T_{v_{1}%
}\circ T_{v_{2}}\circ...\circ T_{v_{m}}$ for \textbf{some} linear extension
$\left(  v_{1},v_{2},...,v_{m}\right)  $ of $P$. Thus, starting with
\textbf{any} linear extension $\left(  v_{1},v_{2},...,v_{m}\right)  $ of $P$,
we could have defined $R$ as the composition $T_{v_{1}}\circ T_{v_{2}}%
\circ...\circ T_{v_{m}}$, and then used Propositions \ref{prop.R.implicit},
\ref{prop.R.implicit.01} and \ref{prop.R.implicit.converse} to characterize
the image $Rf$ of a $\mathbb{K}$-labelling $f$ under this map $R$ in a unique
way without reference to $\left(  v_{1},v_{2},...,v_{m}\right)  $, and thus
concluded that $R$ does not depend on $\left(  v_{1},v_{2},...,v_{m}\right)
$. The details of this derivation are left to the reader.

On a related note, Proposition \ref{prop.R.implicit}, Proposition
\ref{prop.R.implicit.01} and Proposition \ref{prop.R.implicit.converse}
combined can be used as an alternative definition of birational rowmotion $R$,
which works even when the poset $P$ fails to be finite, as long as for every
$v\in P$, there exist only finitely many $u\in P$ satisfying $u>v$ and there
exist only finitely many $u\in P$ satisfying $u\lessdot v$ (provided that some
technicalities arising from Zariski topology on infinite-dimensional spaces
are dealt with).\footnote{The asymmetry between the $>$ and $\lessdot$ signs
in this requirement is intentional. For instance, birational rowmotion can be
defined (but will not be invertible) for the poset $\left\{  0,-1,-2,-3,\ldots
\right\}  $ (with the usual order relation), but not for the poset $\left\{
0,1,2,3,\ldots\right\}  $ (again with the usual order relation).} We will not
dwell on this.
\end{vershort}

\begin{verlong}
As an aside, the proofs of Propositions \ref{prop.R.implicit},
\ref{prop.R.implicit.01} and \ref{prop.R.implicit.converse} allow us to give
an alternative proof of Corollary \ref{cor.R.welldef}, foregoing the use of
Proposition \ref{prop.linext.transitive}. This proof proceeds by noticing the following:

\begin{lemma}
\label{lem.R.implicit.elem}Let $P$ be a finite poset. Let $\left(  v_{1}%
,v_{2},...,v_{m}\right)  $ be a linear extension of $P$. Let $\mathbb{K}$ be a
field. Let $v\in P$. Let $f\in\mathbb{K}^{\widehat{P}}$. Then,%
\[
\left(  \left(  T_{v_{1}}\circ T_{v_{2}}\circ...\circ T_{v_{m}}\right)
f\right)  \left(  v\right)  =\dfrac{1}{f\left(  v\right)  }\cdot\dfrac
{\sum\limits_{\substack{u\in\widehat{P};\\u\lessdot v}}f\left(  u\right)
}{\sum\limits_{\substack{u\in\widehat{P};\\u\gtrdot v}}\dfrac{1}{\left(
\left(  T_{v_{1}}\circ T_{v_{2}}\circ...\circ T_{v_{m}}\right)  f\right)
\left(  u\right)  }}.
\]

\end{lemma}

\begin{proof}
[Proof of Lemma \ref{lem.R.implicit.elem} (sketched).]The proof of Lemma
\ref{lem.R.implicit.elem} is immediately obtained by replacing every
occurrence of \textquotedblleft$R$\textquotedblright\ by \textquotedblleft%
$T_{v_{1}}\circ T_{v_{2}}\circ...\circ T_{v_{m}}$\textquotedblright\ in the
proof of Proposition \ref{prop.R.implicit}.
\end{proof}

\begin{lemma}
\label{lem.R.implicit.01.elem}Let $P$ be a finite poset. Let $\left(
v_{1},v_{2},...,v_{m}\right)  $ be a linear extension of $P$. Let $\mathbb{K}$
be a field. Let $f\in\mathbb{K}^{\widehat{P}}$. Then, $\left(  \left(
T_{v_{1}}\circ T_{v_{2}}\circ...\circ T_{v_{m}}\right)  f\right)  \left(
0\right)  =f\left(  0\right)  $ and $\left(  \left(  T_{v_{1}}\circ T_{v_{2}%
}\circ...\circ T_{v_{m}}\right)  f\right)  \left(  1\right)  =f\left(
1\right)  $.
\end{lemma}

\begin{proof}
[Proof of Lemma \ref{lem.R.implicit.01.elem} (sketched).]The proof of Lemma
\ref{lem.R.implicit.01.elem} is immediately obtained by replacing every
occurrence of \textquotedblleft$R$\textquotedblright\ by \textquotedblleft%
$T_{v_{1}}\circ T_{v_{2}}\circ...\circ T_{v_{m}}$\textquotedblright\ in the
proof of Proposition \ref{prop.R.implicit.01}.
\end{proof}

\begin{lemma}
\label{lem.R.implicit.converse.elem}Let $P$ be a finite poset. Let $\left(
v_{1},v_{2},...,v_{m}\right)  $ be a linear extension of $P$. Let $\mathbb{K}$
be a field. Let $f\in\mathbb{K}^{\widehat{P}}$ and $g\in\mathbb{K}%
^{\widehat{P}}$ be such that $f\left(  0\right)  =g\left(  0\right)  $ and
$f\left(  1\right)  =g\left(  1\right)  $. Assume that%
\[
g\left(  v\right)  =\dfrac{1}{f\left(  v\right)  }\cdot\dfrac{\sum
\limits_{\substack{u\in\widehat{P};\\u\lessdot v}}f\left(  u\right)  }%
{\sum\limits_{\substack{u\in\widehat{P};\\u\gtrdot v}}\dfrac{1}{g\left(
u\right)  }}\ \ \ \ \ \ \ \ \ \ \text{for every }v\in P.
\]
Then, $g=\left(  T_{v_{1}}\circ T_{v_{2}}\circ...\circ T_{v_{m}}\right)  f$.
\end{lemma}

\begin{proof}
[Proof of Lemma \ref{lem.R.implicit.converse.elem} (sketched).]The proof of
Lemma \ref{lem.R.implicit.converse.elem} is immediately obtained by replacing
every occurrence of ``$R$'' by ``$T_{v_{1}}\circ T_{v_{2}}\circ...\circ
T_{v_{m}}$'' in the proof of Proposition \ref{prop.R.implicit.converse}.
\end{proof}

\begin{proof}
[Second proof of Corollary \ref{cor.R.welldef} (sketched).]Just as in the
above proof of Corollary \ref{cor.R.welldef}, we see that the dominant
rational map $T_{v_{1}}\circ T_{v_{2}}\circ...\circ T_{v_{m}}:\mathbb{K}%
^{\widehat{P}}\dashrightarrow\mathbb{K}^{\widehat{P}}$ is well-defined. We
only need to show that it is independent of the choice of the linear extension
$\left(  v_{1},v_{2},...,v_{m}\right)  $.

Let $\left(  w_{1},w_{2},...,w_{m}\right)  $ be a further linear extension of
$P$. Let $f\in\mathbb{K}^{\widehat{P}}$. Applying Lemma
\ref{lem.R.implicit.01.elem} to $\left(  w_{1},w_{2},...,w_{m}\right)  $
instead of $\left(  v_{1},v_{2},...,v_{n}\right)  $, we see that $\left(
\left(  T_{w_{1}}\circ T_{w_{2}}\circ...\circ T_{w_{m}}\right)  f\right)
\left(  0\right)  =f\left(  0\right)  $ and $\left(  \left(  T_{w_{1}}\circ
T_{w_{2}}\circ...\circ T_{w_{m}}\right)  f\right)  \left(  1\right)  =f\left(
1\right)  $. Moreover, applying Lemma \ref{lem.R.implicit.elem} to
\newline$\left(  w_{1},w_{2},...,w_{m}\right)  $ instead of $\left(
v_{1},v_{2},...,v_{n}\right)  $, we see that every $v\in P$ satisfies%
\[
\left(  \left(  T_{w_{1}}\circ T_{w_{2}}\circ...\circ T_{w_{m}}\right)
f\right)  \left(  v\right)  =\dfrac{1}{f\left(  v\right)  }\cdot\dfrac
{\sum\limits_{\substack{u\in\widehat{P};\\u\lessdot v}}f\left(  u\right)
}{\sum\limits_{\substack{u\in\widehat{P};\\u\gtrdot v}}\dfrac{1}{\left(
\left(  T_{w_{1}}\circ T_{w_{2}}\circ...\circ T_{w_{m}}\right)  f\right)
\left(  u\right)  }}.
\]
Hence, the condition of Lemma \ref{lem.R.implicit.converse.elem} is satisfied
if $g$ is set to be \newline$\left(  T_{w_{1}}\circ T_{w_{2}}\circ...\circ
T_{w_{m}}\right)  f$. Thus, applying Lemma \ref{lem.R.implicit.converse.elem}
to $g=\left(  T_{w_{1}}\circ T_{w_{2}}\circ...\circ T_{w_{m}}\right)  f$, we
see that $\left(  T_{w_{1}}\circ T_{w_{2}}\circ...\circ T_{w_{m}}\right)
f=\left(  T_{v_{1}}\circ T_{v_{2}}\circ...\circ T_{v_{m}}\right)  f$.

Now, forget that we fixed $f$. We thus have shown that $\left(  T_{w_{1}}\circ
T_{w_{2}}\circ...\circ T_{w_{m}}\right)  f=\left(  T_{v_{1}}\circ T_{v_{2}%
}\circ...\circ T_{v_{m}}\right)  f$ for every $f\in\mathbb{K}^{\widehat{P}}$.
In other words, $T_{w_{1}}\circ T_{w_{2}}\circ...\circ T_{w_{m}}=T_{v_{1}%
}\circ T_{v_{2}}\circ...\circ T_{v_{m}}$.

Now, the map $T_{w_{1}}\circ T_{w_{2}}\circ...\circ T_{w_{m}}$ is independent
of the choice of the linear extension $\left(  v_{1},v_{2},...,v_{m}\right)
$.\ \ \ \ \footnote{It depends on the choice of the linear extension $\left(
w_{1},w_{2},...,w_{m}\right)  $, but this latter extension can be chosen once
and for all and is unrelated to the choice of $\left(  v_{1},v_{2}%
,...,v_{m}\right)  $.} In other words, the map $T_{v_{1}}\circ T_{v_{2}}%
\circ...\circ T_{v_{m}}$ is independent of the choice of the linear extension
$\left(  v_{1},v_{2},...,v_{m}\right)  $ (since $T_{w_{1}}\circ T_{w_{2}}%
\circ...\circ T_{w_{m}}=T_{v_{1}}\circ T_{v_{2}}\circ...\circ T_{v_{m}}$).
This proves Corollary \ref{cor.R.welldef} again.
\end{proof}

On a related note, Proposition \ref{prop.R.implicit}, Proposition
\ref{prop.R.implicit.01} and Proposition \ref{prop.R.implicit.converse}
combined can be used as an alternative definition of birational rowmotion $R$,
which works even when the poset $P$ fails to be finite, as long as for every
$v\in P$, there exist only finitely many $u\in P$ satisfying $u>v$ and there
exist only finitely many $u\in P$ satisfying $u\lessdot v$%
.\ \ \ \ \footnote{The asymmetry between the $>$ and $\lessdot$ signs in this
requirement is intentional. For instance, birational rowmotion can be defined
(but will not be invertible) for the poset $\left\{  0,-1,-2,-3,\ldots
\right\}  $ (with the usual order relation), but not for the poset $\left\{
0,1,2,3,\ldots\right\}  $ (again with the usual order relation).} We will not
dwell on this.
\end{verlong}

Another general property of birational rowmotion concerns the question of what
happens if the birational toggles are composed not in the ``from top to
bottom'' order as in the definition of birational rowmotion, but the other way
round. It turns out that the result is the inverse of birational rowmotion:

\begin{proposition}
\label{prop.R.inverse}Let $P$ be a finite poset. Let $\mathbb{K}$ be a field.
Then, birational rowmotion $R$ is invertible (as a rational map). Its inverse
$R^{-1}$ is $T_{v_{m}}\circ T_{v_{m-1}}\circ...\circ T_{v_{1}}:\mathbb{K}%
^{\widehat{P}}\dashrightarrow\mathbb{K}^{\widehat{P}}$, where $\left(
v_{1},v_{2},...,v_{m}\right)  $ is a linear extension of $P$.
\end{proposition}

\begin{proof}
[Proof of Proposition \ref{prop.R.inverse} (sketched).]We know that $T_{w}$ is
an involution for every $w\in P$. Thus, in particular, for every $w\in P$, the
map $T_{w}$ is invertible and satisfies $T_{w}^{-1}=T_{w}$.

Let $\left(  v_{1},v_{2},...,v_{m}\right)  $ be a linear extension of $P$.
Then, $R=T_{v_{1}}\circ T_{v_{2}}\circ...\circ T_{v_{m}}$ (by the definition
of $R$), so that $R^{-1}=T_{v_{m}}^{-1}\circ T_{v_{m-1}}^{-1}\circ...\circ
T_{v_{1}}^{-1}$ (this makes sense since the map $T_{w}$ is invertible for
every $w\in P$). Since $T_{w}^{-1}=T_{w}$ for every $w\in P$, this simplifies
to $R^{-1}=T_{v_{m}}\circ T_{v_{m-1}}\circ...\circ T_{v_{1}}$. This proves
Proposition \ref{prop.R.inverse}.
\end{proof}

\section{Graded posets}

In this section, we restrict our attention to what we call \textit{graded
posets} (a notion that encompasses most of the posets we are interested in;
see Definition \ref{def.graded.n-graded}), and define (for this kind of
posets) a family of \textquotedblleft refined rowmotion\textquotedblright%
\ operators $R_{i}$ which toggle only the labels of the $i$-th degree of the
poset. These each turn out to be involutions, and their composition from top
to bottom degree is $R$ on the entire poset. We will later on use these
$R_{i}$ to get a better understanding of $R$ on graded posets.

Let us first introduce our notion of a graded poset:

\begin{definition}
\label{def.graded.n-graded}Let $P$ be a finite poset. Let $n$ be a nonnegative
integer. We say that the poset $P$ is $n$\textit{-graded} if there exists a
surjective map $\deg:P\rightarrow\left\{  1,2,...,n\right\}  $ such that the
following three assertions hold:

\textit{Assertion 1:} Any two elements $u$ and $v$ of $P$ such that $u \gtrdot
v$ satisfy $\deg u=\deg v+1$.

\textit{Assertion 2:} We have $\deg u=1$ for every minimal element $u$ of $P$.

\textit{Assertion 3:} We have $\deg v=n$ for every maximal element $v$ of $P$.
\end{definition}

Note that the word ``surjective'' in Definition \ref{def.graded.n-graded} is
almost superfluous: Indeed, whenever $P\neq\varnothing$, then any map
$\deg:P\rightarrow\left\{  1,2,...,n\right\}  $ satisfying the Assertions 1, 2
and 3 of Definition \ref{def.graded.n-graded} is automatically surjective
(this is easy to prove). But if $P=\varnothing$, such a map exists (vacuously)
for every $n$, whereas requiring surjectivity forced $n=0$.

\begin{example}
The poset $\left\{  1,2\right\}  \times\left\{  1,2\right\}  $ studied in
Example \ref{ex.rowmotion.2x2} is $3$-graded. The poset $P$ studied in Example
\ref{ex.rowmotion.31} is $2$-graded. The empty poset is $0$-graded, but not
$n$-graded for any positive $n$. A chain with $k$ elements is $k$-graded.
\end{example}

\begin{verlong}
\begin{proposition}
\label{prop.graded.unique} Let $P$ be a finite poset. Then, there exists at
most one pair $\left(  n,\deg\right)  $ consisting of a nonnegative integer
$n$ and a surjective map $\deg:P\rightarrow\left\{  1,2,...,n\right\}  $ such
that the Assertions 1, 2 and 3 of Definition \ref{def.graded.n-graded} hold.
\end{proposition}

We leave the proof of Proposition \ref{prop.graded.unique} to the reader (it
proceeds, as one could expect, by induction). We will use Proposition
\ref{prop.graded.unique} as a justification for calling $\deg$ the ``degree
map'' of $P$ without having to worry about whether there could be several such
maps. However, if one were to examine our arguments more closely, one would
see that none of them really relies on Proposition \ref{prop.graded.unique};
if we didn't know Proposition \ref{prop.graded.unique}, we could perform the
same reasoning with the only difference that we would have to define an
$n$-graded poset $P$ as a \textit{triple} of a poset $P$, a nonnegative
integer $n$ and a surjective map $\deg:P\rightarrow\left\{  1,2,...,n\right\}
$ such that the Assertions 1, 2 and 3 of Definition \ref{def.graded.n-graded}
hold, instead of defining it as a poset $P$ for which such $n$ and $\deg$ exist.

We should also mention that Proposition \ref{prop.graded.unique} would still
be valid if we would remove Assertion 3 from Definition
\ref{def.graded.n-graded}. However, this would change the notion of an
``$n$-graded poset'' to a weaker one, which we don't want.
\end{verlong}

\begin{vershort}
\begin{definition}
\label{def.graded.graded}Let $P$ be a finite poset. We say that the poset $P$
is \textit{graded} if there exists an $n\in\mathbb{N}$ such that $P$ is
$n$-graded. This $n$ is then called the \textit{height} of $P$.
\end{definition}
\end{vershort}

\begin{verlong}
\begin{corollary}
\label{cor.graded.unique} Let $P$ be a finite poset. There exists at most one
$n\in\mathbb{N}$ such that the poset $P$ is $n$-graded.
\end{corollary}

\begin{definition}
\label{def.graded.graded}
Let $P$ be a finite poset. We say that the poset $P$ is \textit{graded} if
there exists an $n\in\mathbb{N}$ such that $P$ is $n$-graded. This $n$ is
unique (because of Corollary \ref{cor.graded.unique}), and is called the
\textit{height} of $P$.
\end{definition}
\end{verlong}

The reader should be warned that the notion of a \textquotedblleft graded
poset\textquotedblright\ is not standard across literature; we have found at
least four non-equivalent definitions of this notion in different sources.

\begin{definition}
\label{def.graded.Phat}Let $n\in\mathbb{N}$. Let $P$ be an $n$-graded poset.
Then, there exists a surjective map $\deg:P\rightarrow\left\{
1,2,...,n\right\}  $ that satisfies the Assertions 1, 2 and 3 of Definition
\ref{def.graded.n-graded}. A moment of thought reveals that such a map $\deg$
is also uniquely determined by $P$\ \ \ \ \footnotemark. Thus, we will call
$\deg$ the \textit{degree map} of $P$.

Moreover, we extend this map $\deg$ to a map $\widehat{P}\rightarrow\left\{
0,1,...,n+1\right\}  $ by letting it map $0$ to $0$ and $1$ to $n+1$. This
extended map will also be denoted by $\deg$ and called the degree map. Notice
that this extended map $\deg$ still satisfies Assertion 1 of Definition
\ref{def.graded.n-graded} if $P$ is replaced by $\widehat{P}$ in that assertion.

For every $i\in\left\{  0,1,...,n+1\right\}  $, we will denote by
$\widehat{P}_{i}$ the subset $\deg^{-1}\left(  \left\{  i\right\}  \right)  $
of $\widehat{P}$. For every $v\in\widehat{P}$, the number $\deg v$ is called
the \textit{degree} of $v$.
\end{definition}

\begin{vershort}
\footnotetext{In fact, if $v\in P$, then it is easy to see that $\deg v$
equals the number of elements of any maximal chain in $P$ with highest element
$v$. This clearly determines $\deg v$ uniquely.}
\end{vershort}

\begin{verlong}
\footnotetext{In fact, Proposition \ref{prop.graded.unique} shows that the
pair $\left(  n,\deg\right)  $ is uniquely determined by $P$. As a
consequence, $\deg$ is uniquely determined by $P$.}
\end{verlong}

The notion of an \textquotedblleft$n$-graded poset\textquotedblright\ we just
defined is identical with the notion of a \textquotedblleft graded finite
poset of rank $n-1$\textquotedblright\ as defined in \cite[\S 3.1]%
{stanley-ec1}. The degree of an element $v$ of $P$ as defined in Definition
\ref{def.graded.Phat} is off by $1$ from the rank of $v$ in $P$ in the sense
of \cite[\S 3.1]{stanley-ec1}, but the degree $\deg v$ of an element $v$ of
$\widehat{P}$ equals its rank in $\widehat{P}$ in the sense of \cite[\S 3.1]%
{stanley-ec1}.

The way we extended the map $\deg:P\rightarrow\left\{  1,2,...,n\right\}  $ to
a map $\deg:\widehat{P}\rightarrow\left\{  0,1,...,n+1\right\}  $ in
Definition \ref{def.graded.Phat}, of course, was not arbitrary. In fact, it
was tailored to make the following true:

\begin{proposition}
\label{prop.graded.Phat.why}Let $n\in\mathbb{N}$. Let $P$ be an $n$-graded
poset. Let $u\in\widehat{P}$ and $v\in\widehat{P}$. Consider the map
$\deg:\widehat{P}\rightarrow\left\{  0,1,...,n+1\right\}  $ defined in
Definition \ref{def.graded.Phat}.

\textbf{(a)} If $u\lessdot v$ in $\widehat{P}$, then $\deg u=\deg v-1$.

\textbf{(b)} If $u<v$ in $\widehat{P}$, then $\deg u<\deg v$.

\textbf{(c)} If $u<v$ in $\widehat{P}$ and $\deg u=\deg v-1$, then $u\lessdot
v$ in $\widehat{P}$.

\textbf{(d)} If $u\neq v$ and $\deg u=\deg v$, then $u$ and $v$ are
incomparable in $\widehat{P}$.
\end{proposition}

\begin{proof}
[\nopunct]The rather simple proofs of these facts are left to the reader.
(Note that part \textbf{(a)} incorporates all three Assertions 1, 2 and 3 of
Definition \ref{def.graded.n-graded}.)
\end{proof}

In words, Proposition \ref{prop.graded.Phat.why} \textbf{(d)} states that any
two distinct elements of $\widehat{P}$ having the same degree are
incomparable. We will use this several times below.

One important observation is that any two distinct elements of a graded poset
having the same degree are incomparable. Hence:

\begin{corollary}
\label{cor.Ri.welldef}Let $n\in\mathbb{N}$. Let $\mathbb{K}$ be a field. Let
$P$ be an $n$-graded poset. Let $i\in\left\{  1,2,...,n\right\}  $. Let
$\left(  u_{1},u_{2},...,u_{k}\right)  $ be any list of the elements of
$\widehat{P}_{i}$ with every element of $\widehat{P}_{i}$ appearing exactly
once in the list. Then, the dominant rational map $T_{u_{1}}\circ T_{u_{2}%
}\circ...\circ T_{u_{k}}:\mathbb{K}^{\widehat{P}}\dashrightarrow
\mathbb{K}^{\widehat{P}}$ is well-defined and independent of the choice of the
list $\left(  u_{1},u_{2},...,u_{k}\right)  $.
\end{corollary}

\begin{proof}
[Proof of Corollary \ref{cor.Ri.welldef} (sketched).]This is analogous to the
proof of Corollary \ref{cor.R.welldef}, because any two distinct elements of
$\widehat{P}_{i}$ are incomparable. (In place of the set $\mathcal{L}\left(
P\right)  $ now serves the set of all lists of elements of $\widehat{P}_{i}$
(with every element of $\widehat{P}_{i}$ appearing exactly once in the list).
Any two elements of this latter set are equivalent under the relation $\sim$,
because any two adjacent elements in such a list of elements of $\widehat{P}%
_{i}$ are incomparable and can thus be switched.)
\end{proof}

\begin{definition}
\label{def.Ri}Let $n\in\mathbb{N}$. Let $\mathbb{K}$ be a field. Let $P$ be an
$n$-graded poset. Let $i\in\left\{  1,2,...,n\right\}  $. Then, let $R_{i}$
denote the dominant rational map $T_{u_{1}}\circ T_{u_{2}}\circ...\circ
T_{u_{k}}:\mathbb{K}^{\widehat{P}}\dashrightarrow\mathbb{K}^{\widehat{P}}$,
where $\left(  u_{1},u_{2},...,u_{k}\right)  $ is any list of the elements of
$\widehat{P}_{i}$ with every element of $\widehat{P}_{i}$ appearing exactly
once in the list. This map $T_{u_{1}}\circ T_{u_{2}}\circ...\circ T_{u_{k}}$
is well-defined (in particular, it does not depend on the list $\left(
u_{1},u_{2},...,u_{k}\right)  $) because of Corollary \ref{cor.Ri.welldef}.
\end{definition}

\begin{proposition}
\label{prop.Ri.R}Let $n\in\mathbb{N}$. Let $\mathbb{K}$ be a field. Let $P$ be
an $n$-graded poset. Then,%
\begin{equation}
R=R_{1}\circ R_{2}\circ...\circ R_{n}. \label{def.Ri.R}%
\end{equation}

\end{proposition}

\begin{proof}
[Proof of Proposition \ref{prop.Ri.R} (sketched).]For every $i\in\left\{
1,2,...,n\right\}  $, let $\left(  u_{1}^{\left[  i\right]  },u_{2}^{\left[
i\right]  },...,u_{k_{i}}^{\left[  i\right]  }\right)  $ be a list of the
elements of $\widehat{P}_{i}$ with every element of $\widehat{P}_{i}$
appearing exactly once in the list. Then, every $i\in\left\{
1,2,...,n\right\}  $ satisfies $R_{i}=T_{u_{1}^{\left[  i\right]  }}\circ
T_{u_{2}^{\left[  i\right]  }}\circ...\circ T_{u_{k_{i}}^{\left[  i\right]  }%
}$.

But any listing of the elements of $P$ in order of increasing degree is a
linear extension of $P$ (because any two distinct elements of a graded poset
having the same degree are incomparable). Thus,
\[
\left(  u_{1}^{\left[  1\right]  },u_{2}^{\left[  1\right]  },...,u_{k_{1}%
}^{\left[  1\right]  },\ \ \ u_{1}^{\left[  2\right]  },u_{2}^{\left[
2\right]  },...,u_{k_{2}}^{\left[  2\right]  },\ \ \ ...,\ \ \ u_{1}^{\left[
n\right]  },u_{2}^{\left[  n\right]  },...,u_{k_{n}}^{\left[  n\right]
}\right)
\]
is a linear extension of $P$. Thus, by the definition of $R$, we have%
\begin{align*}
R  &  =\left(  T_{u_{1}^{\left[  1\right]  }}\circ T_{u_{2}^{\left[  1\right]
}}\circ...\circ T_{u_{k_{1}}^{\left[  1\right]  }}\right)  \circ\left(
T_{u_{1}^{\left[  2\right]  }}\circ T_{u_{2}^{\left[  2\right]  }}%
\circ...\circ T_{u_{k_{2}}^{\left[  2\right]  }}\right)  \circ...\circ\left(
T_{u_{1}^{\left[  n\right]  }}\circ T_{u_{2}^{\left[  n\right]  }}%
\circ...\circ T_{u_{k_{n}}^{\left[  n\right]  }}\right) \\
&  =R_{1}\circ R_{2}\circ...\circ R_{n}%
\end{align*}
(since every $i\in\left\{  1,2,...,n\right\}  $ satisfies $T_{u_{1}^{\left[
i\right]  }}\circ T_{u_{2}^{\left[  i\right]  }}\circ...\circ T_{u_{k_{i}%
}^{\left[  i\right]  }}=R_{i}$). This proves Proposition \ref{prop.Ri.R}.
\end{proof}

We recall that birational rowmotion is a composition of toggle maps. As
Proposition \ref{prop.Ri.R} shows, the operators $R_{i}$ are an
\textquotedblleft intermediate\textquotedblright\ step between these toggle
maps and birational rowmotion as a whole, though they are defined only when
the poset $P$ is graded. They will be rather useful for us in our
understanding of birational rowmotion (and the condition on $P$ to be graded
doesn't prevent us from using them, since most of our results concern only
graded posets anyway).

\begin{proposition}
\label{prop.Ri.invo}Let $n\in\mathbb{N}$. Let $\mathbb{K}$ be a field. Let $P$
be an $n$-graded poset. Let $i\in\left\{  1,2,...,n\right\}  $. Then, $R_{i}$
is an involution (that is, $R_{i}^{2}=\operatorname*{id}$ on the set where
$R_{i}$ is defined).
\end{proposition}

\begin{vershort}
\begin{proof}
[Proof of Proposition \ref{prop.Ri.invo} (sketched).]We defined $R_{i}$ as the
composition $T_{u_{1}}\circ T_{u_{2}}\circ...\circ T_{u_{k}}$ of the toggles
$T_{u_{i}}$ where $\left(  u_{1},u_{2},...,u_{k}\right)  $ is any list of the
elements of $\widehat{P}_{i}$ with every element of $\widehat{P}_{i}$
appearing exactly once in the list. These toggles are involutions and commute
(the latter because any two distinct elements of $\widehat{P}_{i}$ are
incomparable, having the same degree in $P$). Since a composition of commuting
involutions is always an involution, this shows that $R_{i}$ is an involution, qed.
\end{proof}
\end{vershort}

\begin{verlong}
\begin{proof}
[Proof of Proposition \ref{prop.Ri.invo} (sketched).]Let $\left(  u_{1}%
,u_{2},...,u_{k}\right)  $ be any list of the elements of $\widehat{P}_{i}$
with every element of $\widehat{P}_{i}$ appearing exactly once in the list.
Then, for any two distinct elements $p$ and $q$ of $\left\{
1,2,...,k\right\}  $, the elements $u_{p}$ and $u_{q}$ of $\widehat{P}$ are
incomparable (since $u_{p}$ and $u_{q}$ are two distinct elements of
$\widehat{P}_{i}$, thus two distinct elements of $\widehat{P}$ having the same
degree $i$, but it is clear that any two distinct elements of $\widehat{P}$
having the same degree are incomparable). Hence, for any two distinct elements
$p$ and $q$ of $\left\{  1,2,...,k\right\}  $, we have $T_{u_{p}}\circ
T_{u_{q}}=T_{u_{q}}\circ T_{u_{p}}$ (by Corollary \ref{cor.Tv.commute}).
Therefore, the maps $T_{u_{1}}$, $T_{u_{2}}$, $...$, $T_{u_{k}}$ commute
pairwise. Hence:%
\begin{align*}
\left(  T_{u_{1}}\circ T_{u_{2}}\circ...\circ T_{u_{k}}\right)  ^{2}  &
=T_{u_{1}}^{2}\circ T_{u_{2}}^{2}\circ...\circ T_{u_{k}}^{2}%
=\operatorname*{id}\circ\operatorname*{id}\circ...\circ\operatorname*{id}\\
&  \ \ \ \ \ \ \ \ \ \ \left(
\begin{array}
[c]{c}%
\text{since Proposition \ref{prop.Tv.invo} yields that }T_{u_{p}}%
^{2}=\operatorname*{id}\\
\text{for every }p\in\left\{  1,2,...,k\right\}
\end{array}
\right) \\
&  =\operatorname*{id}.
\end{align*}
Since $T_{u_{1}}\circ T_{u_{2}}\circ...\circ T_{u_{k}}=R_{i}$ (by the
definition of $R_{i}$), this rewrites as $R_{i}^{2}=\operatorname*{id}$. This
proves Proposition \ref{prop.Ri.invo}.
\end{proof}
\end{verlong}

Similarly to Proposition \ref{prop.R.implicit}, we have:

\begin{proposition}
\label{prop.Ri.implicit}Let $n\in\mathbb{N}$. Let $P$ be an $n$-graded poset.
Let $i\in\left\{  1,2,...,n\right\}  $. Let $\mathbb{K}$ be a field. Let
$v\in\widehat{P}$. Let $f\in\mathbb{K}^{\widehat{P}}$.

\textbf{(a)} If $\deg v\neq i$, then $\left(  R_{i}f\right)  \left(  v\right)
=f\left(  v\right)  $.

\textbf{(b)} If $\deg v=i$, then
\begin{equation}
\left(  R_{i}f\right)  \left(  v\right)  =\dfrac{1}{f\left(  v\right)  }%
\cdot\dfrac{\sum\limits_{\substack{u\in\widehat{P};\\u\lessdot v}}f\left(
u\right)  }{\sum\limits_{\substack{u\in\widehat{P};\\u\gtrdot v}}\dfrac
{1}{f\left(  u\right)  }}. \label{prop.Ri.implicit.eq}%
\end{equation}

\end{proposition}

\begin{vershort}
\begin{proof}
[\nopunct]The proof of this proposition is very similar to that of Proposition
\ref{prop.R.implicit} and therefore left to the reader.
\end{proof}
\end{vershort}

\begin{verlong}
The proof of this proposition is very similar to that of Proposition
\ref{prop.R.implicit}:

\begin{proof}
[Proof of Proposition \ref{prop.Ri.implicit} (sketched).]Let $\left(
u_{1},u_{2},...,u_{k}\right)  $ be any list of the elements of $\widehat{P}%
_{i}$ with every element of $\widehat{P}_{i}$ appearing exactly once in the
list. By the definition of $R_{i}$, we have $R_{i}=T_{u_{1}}\circ T_{u_{2}%
}\circ...\circ T_{u_{k}}$.

\textbf{(b)} Assume that $\deg v=i$. Then, $v\in\widehat{P}_{i}$, so that $v$
is an entry of the list $\left(  u_{1},u_{2},...,u_{k}\right)  $. Let $j$ be
the index satisfying $u_{j}=v$.

Let $A=T_{u_{j+1}}\circ T_{u_{j+2}}\circ...\circ T_{u_{k}}$ and $B=T_{u_{1}%
}\circ T_{u_{2}}\circ...\circ T_{u_{j-1}}$. Then,%
\[
R_{i}=T_{u_{1}}\circ T_{u_{2}}\circ...\circ T_{u_{k}}=\underbrace{T_{u_{1}%
}\circ T_{u_{2}}\circ...\circ T_{u_{j-1}}}_{=B}\circ\underbrace{T_{u_{j}}%
}_{=T_{v}}\circ\underbrace{T_{u_{j+1}}\circ T_{u_{j+2}}\circ...\circ T_{v_{m}%
}}_{=A}=B\circ T_{v}\circ A.
\]

Now:

\begin{itemize}
\item Each of the maps $T_{u_{j}}$ with $j\neq k$ leaves the label at $v$
invariant when acting on a $\mathbb{K}$-labelling. Hence, each of the maps $B$
and $A$ leaves the label at $v$ invariant (since $B$ and $A$ are compositions
of maps $T_{u_{j}}$ with $j\neq k$). Thus, $\left(  B\left(  \left(
T_{v}\circ A\right)  f\right)  \right)  \left(  v\right)  =\left(  \left(
T_{v}\circ A\right)  f\right)  \left(  v\right)  $ and $\left(  Af\right)
\left(  v\right)  =f\left(  v\right)  $. Now, since $R_{i}=B\circ T_{v}\circ
A$, we have%
\begin{align}
\left(  R_{i}f\right)  \left(  v\right)   &  =\left(  \left(  B\circ
T_{v}\circ A\right)  f\right)  \left(  v\right)  =\left(  B\left(  \left(
T_{v}\circ A\right)  f\right)  \right)  \left(  v\right)  =\left(  \left(
T_{v}\circ A\right)  f\right)  \left(  v\right) \nonumber\\
&  =\left(  T_{v}\left(  Af\right)  \right)  \left(  v\right) \nonumber\\
&  =\dfrac{1}{\left(  Af\right)  \left(  v\right)  }\cdot\dfrac{\sum
\limits_{\substack{u\in\widehat{P};\\u \lessdot v}}\left(  Af\right)  \left(
u\right)  }{\sum\limits_{\substack{u\in\widehat{P};\\u \gtrdot v}}\dfrac
{1}{\left(  Af\right)  \left(  u\right)  }}\ \ \ \ \ \ \ \ \ \ \left(
\text{by Proposition \ref{prop.Tv} \textbf{(b)}}\right) \nonumber\\
&  =\dfrac{1}{f\left(  v\right)  }\cdot\dfrac{\sum\limits_{\substack{u\in
\widehat{P};\\u \lessdot v}}\left(  Af\right)  \left(  u\right)  }%
{\sum\limits_{\substack{u\in\widehat{P};\\u \gtrdot v}}\dfrac{1}{\left(
Af\right)  \left(  u\right)  }} \label{pf.Ri.implicit.1}%
\end{align}
(since $\left(  Af\right)  \left(  v\right)  =f\left(  v\right)  $).

\item Let $u\in\widehat{P}$ be such that $u\lessdot v$. Then, Proposition
\ref{prop.graded.Phat.why} \textbf{(a)} yields $\deg u=\underbrace{\deg
v}_{=i}-1=i-1\neq i$, so that $u\notin\widehat{P}_{i}$. Hence, $u$ is none of
the elements $u_{j+1}$, $u_{j+2}$, $...$, $u_{k}$ (for the simple reason that
all of these elements belong to $\widehat{P}_{i}$). Consequently, each of the
maps $T_{u_{j+1}}$, $T_{u_{j+2}}$, $...$, $T_{u_{k}}$ leaves the label at $u$
invariant when acting on a $\mathbb{K}$-labelling. Therefore, $A$ also leaves
the label at $u$ invariant (since $A$ is a composition of these maps
$T_{u_{j+1}}$, $T_{u_{j+2}}$, $...$, $T_{u_{k}}$). Hence, $\left(  Af\right)
\left(  u\right)  =f\left(  u\right)  $.

Forget that we fixed $u$. We have thus shown that%
\begin{equation}
\left(  Af\right)  \left(  u\right)  =f\left(  u\right)
\ \ \ \ \ \ \ \ \ \ \text{for every }u\in\widehat{P}\text{ such that } u
\lessdot v. \label{pf.Ri.implicit.2}%
\end{equation}

\item An argument completely analogous to the one we proved
(\ref{pf.Ri.implicit.2}) (except that we now apply Proposition
\ref{prop.graded.Phat.why} \textbf{(a)} to $v$ and $u$ instead of $u$ and $v$)
shows that%
\begin{equation}
\left(  Af\right)  \left(  u\right)  =f\left(  u\right)
\ \ \ \ \ \ \ \ \ \ \text{for every }u\in\widehat{P}\text{ such that }u\gtrdot
v. \label{pf.Ri.implicit.3}%
\end{equation}

\end{itemize}

Now, substituting (\ref{pf.Ri.implicit.2}) and (\ref{pf.Ri.implicit.3}) into
(\ref{pf.Ri.implicit.1}), we obtain%
\[
\left(  R_{i}f\right)  \left(  v\right)  =\dfrac{1}{f\left(  v\right)  }%
\cdot\dfrac{\sum\limits_{\substack{u\in\widehat{P};\\u \lessdot v}}f\left(
u\right)  }{\sum\limits_{\substack{u\in\widehat{P};\\u \gtrdot v}}\dfrac
{1}{f\left(  u\right)  }}.
\]
This proves Proposition \ref{prop.Ri.implicit} \textbf{(b)}.

\textbf{(a)} Assume that $\deg v\neq i$. Then, $v\notin\widehat{P}_{i}$, so
that $v$ is \textbf{not} an entry of the list $\left(  u_{1},u_{2}%
,...,u_{k}\right)  $. In other words, $v$ is none of the elements $u_{1}$,
$u_{2}$, $...$, $u_{k}$. But this yields that each of the maps $T_{u_{1}}$,
$T_{u_{2}}$, $...$, $T_{u_{k}}$ leaves the label at $v$ invariant when acting
on a $\mathbb{K}$-labelling. Hence, so does the composition $T_{u_{1}}\circ
T_{u_{2}}\circ...\circ T_{u_{k}}$ of these maps. Thus, $\left(  \left(
T_{u_{1}}\circ T_{u_{2}}\circ...\circ T_{u_{k}}\right)  f\right)  \left(
v\right)  =f\left(  v\right)  $. Since $T_{u_{1}}\circ T_{u_{2}}\circ...\circ
T_{u_{k}}=R_{i}$, this rewrites as $\left(  R_{i}f\right)  \left(  v\right)
=f\left(  v\right)  $. This proves Proposition \ref{prop.Ri.implicit}
\textbf{(a)}.
\end{proof}

While we have not made any precise statements about well-definedness in
Proposition \ref{prop.Ri.implicit}, we can actually show the following by a
more careful argument:

\begin{proposition}
\label{prop.Ri.implicit.wd}Let $n\in\mathbb{N}$. Let $P$ be an $n$-graded
poset. Let $i\in\left\{  1,2,...,n\right\}  $. Let $\mathbb{K}$ be a field.
Let $v\in\widehat{P}$. Let $f\in\mathbb{K}^{\widehat{P}}$. Assume that
$\dfrac{1}{f\left(  v\right)  }\cdot\dfrac{\sum\limits_{\substack{u\in
\widehat{P};\\u\lessdot v}}f\left(  u\right)  }{\sum\limits_{\substack{u\in
\widehat{P};\\u\gtrdot v}}\dfrac{1}{f\left(  u\right)  }}$ is well-defined for
every $v \in\widehat{P}_{i}$. Then, $R_{i} f$ is well-defined.
\end{proposition}

\begin{proof}
[Proof of Proposition \ref{prop.Ri.implicit.wd} (sketched).]Let $\left(
u_{1},u_{2},...,u_{k}\right)  $ be any list of the elements of $\widehat{P}%
_{i}$ with every element of $\widehat{P}_{i}$ appearing exactly once in the
list. By the definition of $R_{i}$, we have $R_{i}=T_{u_{1}}\circ T_{u_{2}%
}\circ...\circ T_{u_{k}}$.

For every $i\in\left\{  0,1,...,k\right\}  $, denote the rational map
$T_{u_{k-i+1}}\circ T_{u_{k-i+2}}\circ...\circ T_{u_{k}}:\mathbb{K}%
^{\widehat{P}}\dashrightarrow\mathbb{K}^{\widehat{P}}$ by $\rho_{i}$. The
definition of $\rho_{k}$ yields%
\[
\rho_{k}=T_{u_{k-k+1}}\circ T_{u_{k-k+2}}\circ...\circ T_{u_{k}}=T_{u_{1}%
}\circ T_{u_{2}}\circ...\circ T_{u_{k}}=R_{i}%
\]
(since $R_{i}=T_{u_{1}}\circ T_{u_{2}}\circ...\circ T_{u_{k}}$ (by the
definition of $R_i$)). On the other hand, the definition of $\rho_{0}$
yields%
\[
\rho_{0}=T_{u_{k-0+1}}\circ T_{u_{k-0+2}}\circ...\circ T_{u_{k}}=T_{u_{k+1}%
}\circ T_{u_{k+2}}\circ...\circ T_{u_{k}}=\left(  \text{empty composition}%
\right)  =\operatorname*{id}.
\]

We will now show that
\begin{equation}
\rho_{i}f\text{ is well-defined for every }i\in\left\{  0,1,...,k\right\}  .
\label{pf.Ri.implicit.wd.1}%
\end{equation}

\textit{Proof of (\ref{pf.Ri.implicit.wd.1}):} We will prove
(\ref{pf.Ri.implicit.wd.1}) by induction over $i$:

\textit{Induction base:} If $i=0$, then $\underbrace{\rho_{i}}_{=\rho
_{0}=\operatorname*{id}}f=\operatorname*{id}f=f$ is clearly well-defined. In
other words, (\ref{pf.Ri.implicit.wd.1}) is proven for $i=0$. This completes
the induction base.

\textit{Induction step:} Let $I\in\left\{  0,1,...,k-1\right\}  $. Assume that
(\ref{pf.Ri.implicit.wd.1}) holds for $i=I$. We now need to prove that
(\ref{pf.Ri.implicit.wd.1}) holds for $i=I+1$.

We know that (\ref{pf.Ri.implicit.wd.1}) holds for $i=I$. In other words,
$\rho_{I}f$ is well-defined.

Let $v=u_{k-I}$. Then, $v\in\widehat{P}_{i}$.\ \ \ \ \footnote{\textit{Proof.}
We know that $\left(  u_{1},u_{2},...,u_{k}\right)  $ is a list of the
elements of $\widehat{P}_{i}$. Thus, every entry of the list $\left(
u_{1},u_{2},...,u_{k}\right)  $ is an element of $\widehat{P}_{i}$. Since $v$
is an entry of the list $\left(  u_{1},u_{2},...,u_{k}\right)  $ (because
$v=u_{k-I}$), this yields that $v$ is an element of $\widehat{P}_{i}$. In
other words, $v\in\widehat{P}_{i}$, qed.} Hence, $\dfrac{1}{f\left(  v\right)
}\cdot\dfrac{\sum\limits_{\substack{u\in\widehat{P};\\u\lessdot v}}f\left(
u\right)  }{\sum\limits_{\substack{u\in\widehat{P};\\u\gtrdot v}}\dfrac
{1}{f\left(  u\right)  }}$ is well-defined (by the assumptions of Proposition
\ref{prop.Ri.implicit.wd}). In particular, $\sum\limits_{\substack{u\in
\widehat{P};\\u\gtrdot v}}\dfrac{1}{f\left(  u\right)  }$ is well-defined. On
the other hand, $\rho_{I+1}=T_{v}\circ\rho_{I}$%
\ \ \ \ \footnote{\textit{Proof.} We have $\rho_{I}=T_{u_{k-I+1}}\circ
T_{u_{k-I+2}}\circ...\circ T_{u_{k}}$ (by the definition of $\rho_{I}$), so
that $T_{u_{k-I+1}}\circ T_{u_{k-I+2}}\circ...\circ T_{u_{k}}=\rho_{I}$. But
the definition of $\rho_{I+1}$ yields%
\begin{align*}
\rho_{I+1}  &  =T_{u_{k-\left(  I+1\right)  +1}}\circ T_{u_{k-\left(
I+1\right)  +2}}\circ...\circ T_{u_{k}}=T_{u_{k-I}}\circ T_{u_{k-I+1}}%
\circ...\circ T_{u_{k}}\\
&  =\underbrace{T_{u_{k-I}}}_{\substack{=T_{v}\\\text{(since }u_{k-I}%
=v\\\text{(since }v=u_{k-I}\text{))}}}\circ\underbrace{\left(  T_{u_{k-I+1}%
}\circ T_{u_{k-I+2}}\circ...\circ T_{u_{k}}\right)  }_{=\rho_{I}}=T_{v}%
\circ\rho_{I},
\end{align*}
qed.}.

Notice that%
\begin{equation}
\left(  \rho_{I}f\right)  \left(  u\right)  =f\left(  u\right)
\ \ \ \ \ \ \ \ \ \ \text{for every }u\in\widehat{P}\text{ satisfying
}u\lessdot v. \label{pf.Ri.implicit.wd.2a}%
\end{equation}
\footnote{\textit{Proof of (\ref{pf.Ri.implicit.wd.2a}):} Let $u\in
\widehat{P}$ satisfy $u\lessdot v$. Then, Proposition
\ref{prop.graded.Phat.why} \textbf{(a)} yields $\deg u=\underbrace{\deg
v}_{=i}-1=i-1\neq i$.
\par
Notice that $\rho_{I}=T_{u_{k-I+1}}\circ T_{u_{k-I+2}}\circ...\circ T_{u_k}$ (by the
definition of $\rho_{I}$). Hence, $\rho_{I}$ is a composition of the maps
$T_{u_{k-I+1}}$, $T_{u_{k-I+2}}$, $...$, $T_{u_k}$.
\par
But $\left(  u_{1},u_{2},...,u_{k}\right)  $ is a list of elements of
$\widehat{P}_{i}$. Hence, $u_{\ell}$ is an element of $\widehat{P}_{i}$ for
every $\ell\in\left\{  1,2,...,k\right\}  $.
\par
Now, let $\ell\in\left\{  k-I+1,k-I+2,...,k\right\}  $. Then, $\ell\in\left\{
k-I+1,k-I+2,...,k\right\}  \subseteq\left\{  1,2,...,k\right\}  $, so that
$u_{\ell}$ is an element of $\widehat{P}_{i}$ (as we just have shown). Hence,
$u_{\ell}\in\widehat{P}_{i}=\deg^{-1}\left(  \left\{  i\right\}  \right)  $
(by the definition of $\widehat{P}_{i}$). In other words, $\deg\left(
u_{\ell}\right)  \in\left\{  i\right\}  $, so that $\deg\left(  u_{\ell
}\right)  =i$. Compared with $\deg u\neq i$, this yields $\deg\left(  u_{\ell
}\right)  \neq\deg u$, so that $u_{\ell}\neq u$. Thus, the map $T_{u_{\ell}}$
leaves the label at $u$ invariant when acting on a $\mathbb{K}$-labelling
(because of Proposition \ref{prop.Tv} \textbf{(a)}).
\par
Now, forget that we fixed $\ell$. We thus have shown that for every $\ell
\in\left\{  k-I+1,k-I+2,...,k\right\}  $, the map $T_{u_{\ell}}$ leaves the
label at $u$ invariant when acting on a $\mathbb{K}$-labelling. In other
words, each of the maps $T_{u_{k-I+1}}$, $T_{u_{k-I+2}}$, $...$, $T_{u_k}$ leaves the
label at $u$ invariant when acting on a $\mathbb{K}$-labelling. Hence, the map
$\rho_{I}$ also leaves the label at $u$ invariant when acting on a
$\mathbb{K}$-labelling (since $\rho_{I}$ is a composition of these maps
 $T_{u_{k-I+1}}$, $T_{u_{k-I+2}}$, $...$, $T_{u_k}$). In other words, $\left(  \rho
_{I}g\right)  \left(  u\right)  =g\left(  u\right)  $ for every $\mathbb{K}%
$-labelling $g\in\mathbb{K}^{\widehat{P}}$. Applying this to $g=f$, we obtain
$\left(  \rho_{I}f\right)  \left(  u\right)  =f\left(  u\right)  $. This
proves (\ref{pf.Ri.implicit.wd.2a}).} Thus,%
\begin{equation}
\sum\limits_{\substack{u\in\widehat{P};\\u\lessdot v}}\underbrace{\left(
\rho_{I}f\right)  \left(  u\right)  }_{\substack{=f\left(  u\right)
\\\text{(by (\ref{pf.Ri.implicit.wd.2a}))}}}=\sum\limits_{\substack{u\in
\widehat{P};\\u\lessdot v}}f\left(  u\right)  . \label{pf.Ri.implicit.wd.2as}%
\end{equation}
Also,%
\begin{equation}
\left(  \rho_{I}f\right)  \left(  u\right)  =f\left(  u\right)
\ \ \ \ \ \ \ \ \ \ \text{for every }u\in\widehat{P}\text{ satisfying
}u\gtrdot v. \label{pf.Ri.implicit.wd.2b}%
\end{equation}
\footnote{\textit{Proof of (\ref{pf.Ri.implicit.wd.2b}):} Let $u\in
\widehat{P}$ satisfy $u\gtrdot v$. Then, $v\lessdot u$. Thus, Proposition
\ref{prop.graded.Phat.why} \textbf{(a)} (applied to $v$ and $u$ instead of $u$
and $v$) yields $\deg v=\deg u-1$, so that $\deg u=\underbrace{\deg v}%
_{=i}+1=i+1\neq i$. From here on, we can prove $\left(  \rho_{I}f\right)
\left(  u\right)  =f\left(  u\right)  $ by arguing as in the proof of
(\ref{pf.Ri.implicit.wd.2a}). This proves (\ref{pf.Ri.implicit.wd.2b}).} Thus,%
\begin{equation}
\sum\limits_{\substack{u\in\widehat{P};\\u\gtrdot v}}\underbrace{\dfrac
{1}{\left(  \rho_{I}f\right)  \left(  u\right)  }}_{\substack{=\dfrac
{1}{f\left(  u\right)  }\\\text{(because (\ref{pf.Ri.implicit.wd.2b}%
)}\\\text{yields }\left(  \rho_{I}f\right)  \left(  u\right)  =f\left(
u\right)  \text{)}}}=\sum\limits_{\substack{u\in\widehat{P};\\u\gtrdot
v}}\dfrac{1}{f\left(  u\right)  } \label{pf.Ri.implicit.wd.2bs}%
\end{equation}
(and this is well-defined because we know that $\sum\limits_{\substack{u\in
\widehat{P};\\u\gtrdot v}}\dfrac{1}{f\left(  u\right)  }$ is well-defined).
Finally,%
\begin{equation}
\left(  \rho_{I}f\right)  \left(  v\right)  =f\left(  v\right)  .
\label{pf.Ri.implicit.wd.2v}%
\end{equation}
\footnote{\textit{Proof of (\ref{pf.Ri.implicit.wd.2v}):} Notice that
$\rho_{I}=T_{u_{k-I+1}}\circ T_{u_{k-I+2}}\circ...\circ T_{u_k}$ (by the definition of
$\rho_{I}$). Hence, $\rho_{I}$ is a composition of the maps
$T_{u_{k-I+1}}$, $T_{u_{k-I+2}}$, $...$, $T_{u_k}$.
\par
Let $\ell\in\left\{  k-I+1,k-I+2,...,k\right\}  $. Then, $\ell\geq k-I+1>k-I$,
so that $\ell\neq k-I$.
\par
Assume (for the sake of contradiction) that $u_{\ell}=v$. Now, recall that
$\left(  u_{1},u_{2},...,u_{k}\right)  $ is a list of the elements of
$\widehat{P}_{i}$ with every element of $\widehat{P}_{i}$ appearing exactly
once in the list. Hence, every element of $\widehat{P}_{i}$ appears exactly
once in the list $\left(  u_{1},u_{2},...,u_{k}\right)  $. In particular,
every element of $\widehat{P}_{i}$ appears at most once in the list $\left(
u_{1},u_{2},...,u_{k}\right)  $. In particular, $v$ appears at most once in
the list $\left(  u_{1},u_{2},...,u_{k}\right)  $ (since $v$ is an element of
$\widehat{P}_{i}$). In other words, if $\alpha$ and $\beta$ are two elements
of $\left\{  1,2,...,k\right\}  $ such that $u_{\alpha}=v$ and $u_{\beta}=v$,
then $\alpha=\beta$. Applying this to $\alpha=\ell$ and $\beta=k-I$, we obtain
$\ell=k-I$ (since $u_{\ell}=v$ and $u_{k-I}=v$). But this contradicts
$\ell\neq k-I$. This contradiction shows that our assumption (that $u_{\ell
}=v$) was false. Hence, we cannot have $u_{\ell}=v$. In other words, we have
$u_{\ell}\neq v$. Thus, the map $T_{u_{\ell}}$ leaves the label at $v$
invariant when acting on a $\mathbb{K}$-labelling (because of Proposition
\ref{prop.Tv} \textbf{(a)}).
\par
Now, forget that we fixed $\ell$. We thus have shown that for every $\ell
\in\left\{  k-I+1,k-I+2,...,k\right\}  $, the map $T_{u_{\ell}}$ leaves the
label at $v$ invariant when acting on a $\mathbb{K}$-labelling. In other
words, each of the maps $T_{u_{k-I+1}}$, $T_{u_{k-I+2}}$, $...$, $T_{u_k}$ leaves the
label at $v$ invariant when acting on a $\mathbb{K}$-labelling. Hence, the map
$\rho_{I}$ also leaves the label at $v$ invariant when acting on a
$\mathbb{K}$-labelling (since $\rho_{I}$ is a composition of these maps
$T_{u_{k-I+1}}$, $T_{u_{k-I+2}}$, $...$, $T_{u_k}$). In other words, $\left(  \rho
_{I}g\right)  \left(  v\right)  =g\left(  v\right)  $ for every $\mathbb{K}%
$-labelling $g\in\mathbb{K}^{\widehat{P}}$. Applying this to $g=f$, we obtain
$\left(  \rho_{I}f\right)  \left(  v\right)  =f\left(  v\right)  $. This
proves (\ref{pf.Ri.implicit.wd.2v}).}

Now, we are going to prove that $T_{v}\left(  \rho_{I}f\right)  $ is
well-defined. In fact, let us recall that $T_{v}\left(  \rho_{I}f\right)  $ is
defined by%
\begin{equation}
\left(  T_{v}\left(  \rho_{I}f\right)  \right)  \left(  w\right)  =\left\{
\begin{array}
[c]{l}%
\left(  \rho_{I}f\right)  \left(  w\right)  ,\ \ \ \ \ \ \ \ \ \ \text{if
}w\neq v;\\
\dfrac{1}{\left(  \rho_{I}f\right)  \left(  v\right)  }\cdot\dfrac
{\sum\limits_{\substack{u\in\widehat{P};\\u\lessdot v}}\left(  \rho
_{I}f\right)  \left(  u\right)  }{\sum\limits_{\substack{u\in\widehat{P}%
;\\u\gtrdot v}}\dfrac{1}{\left(  \rho_{I}f\right)  \left(  u\right)  }%
},\ \ \ \ \ \ \ \ \ \ \text{if }w=v
\end{array}
\right.  \ \ \ \ \ \ \ \ \ \ \text{for all }w\in\widehat{P}.
\label{pf.Ri.implicit.wd.step.1}%
\end{equation}
Hence, in order to prove that $T_{v}\left(  \rho_{I}f\right)  $ is
well-defined, we need to show that the right hand side of
(\ref{pf.Ri.implicit.wd.step.1}) is well-defined for every $w\in\widehat{P}$.

For every $w\in\widehat{P}$, the right hand side of
(\ref{pf.Ri.implicit.wd.step.1}) is%
\begin{align*}
&  \left\{
\begin{array}
[c]{l}%
\left(  \rho_{I}f\right)  \left(  w\right)  ,\ \ \ \ \ \ \ \ \ \ \text{if
}w\neq v;\\
\dfrac{1}{\left(  \rho_{I}f\right)  \left(  v\right)  }\cdot\dfrac
{\sum\limits_{\substack{u\in\widehat{P};\\u\lessdot v}}\left(  \rho
_{I}f\right)  \left(  u\right)  }{\sum\limits_{\substack{u\in\widehat{P}%
;\\u\gtrdot v}}\dfrac{1}{\left(  \rho_{I}f\right)  \left(  u\right)  }%
},\ \ \ \ \ \ \ \ \ \ \text{if }w=v
\end{array}
\right. \\
&  =\left\{
\begin{array}
[c]{l}%
\left(  \rho_{I}f\right)  \left(  w\right)  ,\ \ \ \ \ \ \ \ \ \ \text{if
}w\neq v;\\
\dfrac{1}{f\left(  v\right)  }\cdot\dfrac{\sum\limits_{\substack{u\in
\widehat{P};\\u\lessdot v}}f\left(  u\right)  }{\sum\limits_{\substack{u\in
\widehat{P};\\u\gtrdot v}}\dfrac{1}{f\left(  u\right)  }}%
,\ \ \ \ \ \ \ \ \ \ \text{if }w=v
\end{array}
\right. \\
&  \ \ \ \ \ \ \ \ \ \ \left(  \text{by (\ref{pf.Ri.implicit.wd.2v}),
(\ref{pf.Ri.implicit.wd.2as}) and (\ref{pf.Ri.implicit.wd.2bs})}\right)  .
\end{align*}
This is well-defined (since $\left(  \rho_{I}f\right)  \left(  w\right)  $ is
well-defined (because $\rho_{I}f$ is well-defined) and since $\dfrac
{1}{f\left(  v\right)  }\cdot\dfrac{\sum\limits_{\substack{u\in\widehat{P}%
;\\u\lessdot v}}f\left(  u\right)  }{\sum\limits_{\substack{u\in
\widehat{P};\\u\gtrdot v}}\dfrac{1}{f\left(  u\right)  }}$ is well-defined (by
the assumptions of Proposition \ref{prop.Ri.implicit.wd})). Thus, we have shown
that the right hand side of (\ref{pf.Ri.implicit.wd.step.1}) is well-defined
for every $w\in\widehat{P}$. Since (\ref{pf.Ri.implicit.wd.step.1}) is the
definition of $T_{v}\left(  \rho_{I}f\right)  $, this yields that
$T_{v}\left(  \rho_{I}f\right)  $ is well-defined. In other words, $\rho
_{I+1}f$ is well-defined (because $\underbrace{\rho_{I+1}}_{=T_{v}\circ
\rho_{I}}f=\left(  T_{v}\circ\rho_{I}\right)  f=T_{v}\left(  \rho_{I}f\right)
$). In other words, (\ref{pf.Ri.implicit.wd.1}) holds for $i=I+1$. This
completes the induction step. The induction proof of
(\ref{pf.Ri.implicit.wd.1}) is thus finished.

Now we have proven (\ref{pf.Ri.implicit.wd.1}). Hence, we can apply
(\ref{pf.Ri.implicit.wd.1}) to $i=k$. As a result, we see that $\rho_{k}f$ is
well-defined. In other words, $R_{i}f$ is well-defined (since $\rho_{k}=R_{i}%
$). This proves Proposition \ref{prop.Ri.implicit.wd}.
\end{proof}
\end{verlong}

Notice that using the proof of Proposition \ref{prop.Ri.implicit}, it is easy
to give an alternative proof of Corollary \ref{cor.Ri.welldef} (in the same
way as we saw that an alternative proof of Corollary \ref{cor.R.welldef} could
be given using the proofs of Propositions \ref{prop.R.implicit},
\ref{prop.R.implicit.01} and \ref{prop.R.implicit.converse}).

\begin{verlong}
Let us give an alternative proof of Proposition \ref{prop.Ri.R}, which has the
advantage of being generalizable to certain situations in which $P$ is allowed
to be infinite (although we don't study such situations):

\begin{proof}
[Second proof of Proposition \ref{prop.Ri.R} (sketched).]We denote by $<$ the
smaller relation of the poset $P$. We also denote by $<$ the smaller relation
of the poset $\widehat{P}$. We thus have given two meanings to the notation
$<$, but these two meanings don't conflict (since the restriction of the
smaller relation of $\widehat{P}$ to $P$ is the smaller relation of $P$).

For every $i\in\left\{  0,1,...,n-1\right\}  $, let us denote the rational map
$R_{n-i+1}\circ R_{n-i+2}\circ...\circ R_{n}:\mathbb{K}^{\widehat{P}%
}\dashrightarrow\mathbb{K}^{\widehat{P}}$ by $\rho_{i}$. Then, the definition
of $\rho_{n}$ yields%
\[
\rho_{n}=R_{n-n+1}\circ R_{n-n+2}\circ...\circ R_{n}=R_{1}\circ R_{2}%
\circ...\circ R_{n}.
\]
On the other hand, the definition of $\rho_{0}$ yields%
\[
\rho_{0}=R_{n-0+1}\circ R_{n-0+2}\circ...\circ R_{n}=R_{n+1}\circ R_{n+2}%
\circ...\circ R_{n}=\left(  \text{empty composition}\right)
=\operatorname*{id}.
\]

Let $f\in\mathbb{K}^{\widehat{P}}$ be a $\mathbb{K}$-labelling of $P$ for
which $Rf$ is well-defined.

We are going to prove that for every $i\in\left\{  0,1,...,n\right\}  $, the
following assertion holds:

\textit{Assertion A:} The $\mathbb{K}$-labelling $\rho_{i}f\in\mathbb{K}%
^{\widehat{P}}$ is well-defined, and every $v\in\widehat{P}$ satisfies%
\[
\left(  \rho_{i}f\right)  \left(  v\right)  =\left\{
\begin{array}
[c]{c}%
\left(  Rf\right)  \left(  v\right)  ,\ \ \ \ \ \ \ \ \ \ \text{if }\deg
v>n-i;\\
f\left(  v\right)  ,\ \ \ \ \ \ \ \ \ \ \text{if }\deg v\leq n-i
\end{array}
\right.  .
\]

In fact, we will prove Assertion A by induction over $i$:

\textit{Induction base:} Assertion A is satisfied for $i=0$%
\ \ \ \ \footnote{\textit{Proof.} Assume that $i=0$. Hence, $n-i=n-0=n$, so
that $n=n-i$. But $i=0$, and thus $\rho_{i}=\rho_{0}=\operatorname*{id} $, so
that $\rho_{i}f=\operatorname*{id} f=f$ is well-defined. Now, let
$v\in\widehat{P}$. We are going to prove that
\begin{equation}
\left(  \rho_{i}f\right)  \left(  v\right)  =\left\{
\begin{array}
[c]{c}%
\left(  Rf\right)  \left(  v\right)  ,\ \ \ \ \ \ \ \ \ \ \text{if }\deg
v>n-i;\\
f\left(  v\right)  ,\ \ \ \ \ \ \ \ \ \ \text{if }\deg v\leq n-i
\end{array}
\right.  . \label{pf.Ri.R.base.goal}%
\end{equation}
\par
Let us first assume that $v\neq1$. Then, $\deg v\leq n$ and%
\begin{equation}
\left\{
\begin{array}
[c]{c}%
\left(  Rf\right)  \left(  v\right)  ,\ \ \ \ \ \ \ \ \ \ \text{if }\deg
v>n-i;\\
f\left(  v\right)  ,\ \ \ \ \ \ \ \ \ \ \text{if }\deg v\leq n-i
\end{array}
\right.  =f\left(  v\right)  \ \ \ \ \ \ \ \ \ \ \left(  \text{since }\deg
v\leq n=n-i\right)  . \label{pf.Ri.R.base.1}%
\end{equation}
Now,%
\[
\underbrace{\left(  \rho_{i}f\right)  }_{=f}\left(  v\right)  =f\left(
v\right)  =\left\{
\begin{array}
[c]{c}%
\left(  Rf\right)  \left(  v\right)  ,\ \ \ \ \ \ \ \ \ \ \text{if }\deg
v>n-i;\\
f\left(  v\right)  ,\ \ \ \ \ \ \ \ \ \ \text{if }\deg v\leq n-i
\end{array}
\right.
\]
(by (\ref{pf.Ri.R.base.1})). Thus, (\ref{pf.Ri.R.base.goal}) is proven under
the assumption that $v\neq1$.
\par
Now, forget that we assumed that $v\neq1$. We thus have proven
(\ref{pf.Ri.R.base.goal}) under the assumption that $v\neq1$. Hence, for the
rest of the proof of (\ref{pf.Ri.R.base.goal}), we can WLOG assume that we
don't have $v\neq1$. Assume this. We don't have $v\neq1$. Hence, we have
$v=1$. Hence, $\left(  Rf\right)  \left(  \underbrace{v}_{=1}\right)  =\left(
Rf\right)  \left(  1\right)  =f\left(  1\right)  $ (by Proposition
\ref{prop.R.implicit.01}). But $\deg\underbrace{v}_{=1}=\deg1=n+1$ and%
\begin{align}
\left\{
\begin{array}
[c]{c}%
\left(  Rf\right)  \left(  v\right)  ,\ \ \ \ \ \ \ \ \ \ \text{if }\deg
v>n-i;\\
f\left(  v\right)  ,\ \ \ \ \ \ \ \ \ \ \text{if }\deg v\leq n-i
\end{array}
\right.   &  =\left(  Rf\right)  \left(  v\right)  \ \ \ \ \ \ \ \ \ \ \left(
\text{since }\deg v=n+1>n=n-i\right) \label{pf.Ri.R.base.2}\\
&  =f\left(  1\right)  .\nonumber
\end{align}
Now,%
\[
\underbrace{\left(  \rho_{i}f\right)  }_{=f}\left(  v\right)  =f\left(
\underbrace{v}_{=1}\right)  =f\left(  1\right)  =\left\{
\begin{array}
[c]{c}%
\left(  Rf\right)  \left(  v\right)  ,\ \ \ \ \ \ \ \ \ \ \text{if }\deg
v>n-i;\\
f\left(  v\right)  ,\ \ \ \ \ \ \ \ \ \ \text{if }\deg v\leq n-i
\end{array}
\right.
\]
(by (\ref{pf.Ri.R.base.2})). Hence, (\ref{pf.Ri.R.base.goal}) is proven.
\par
Now, forget that we fixed $v$. We thus have shown that every $v\in P$
satisfies%
\[
\left(  \rho_{i}f\right)  \left(  v\right)  =\left\{
\begin{array}
[c]{c}%
\left(  Rf\right)  \left(  v\right)  ,\ \ \ \ \ \ \ \ \ \ \text{if }\deg
v>n-i;\\
f\left(  v\right)  ,\ \ \ \ \ \ \ \ \ \ \text{if }\deg v\leq n-i
\end{array}
\right.  .
\]
Hence, Assertion A is proven, qed.}. Thus, the induction base is complete.

\textit{Induction step:} Let $I\in\left\{  0,1,...,n-1\right\}  $. Assume that
Assertion A is satisfied for $i=I$. We need to prove that Assertion A is
satisfied for $i=I+1$.

We know that Assertion A is satisfied for $i=I$. In other words, the
$\mathbb{K}$-labelling $\rho_{I}f\in\mathbb{K}^{\widehat{P}}$ is well-defined,
and every $v\in\widehat{P}$ satisfies%
\begin{equation}
\left(  \rho_{I}f\right)  \left(  v\right)  =\left\{
\begin{array}
[c]{c}%
\left(  Rf\right)  \left(  v\right)  ,\ \ \ \ \ \ \ \ \ \ \text{if }\deg
v>n-I;\\
f\left(  v\right)  ,\ \ \ \ \ \ \ \ \ \ \text{if }\deg v\leq n-I
\end{array}
\right.  . \label{pf.Ri.R.step.1}%
\end{equation}

Notice that $I\in\left\{  0,1,...,n-1\right\}  $, so that $n-I\in\left\{
1,2,...,n\right\}  $.

We have $\rho_{I+1}=R_{n-I}\circ\rho_{I}$\ \ \ \ \footnote{\textit{Proof.} We
have $\rho_{I}=R_{n-I+1}\circ R_{n-I+2}\circ...\circ R_{n}$ (by the definition
of $\rho_{I}$), so that $R_{n-I+1}\circ R_{n-I+2}\circ...\circ R_{n}=\rho_{I}%
$. But the definition of $\rho_{I+1}$ yields%
\begin{align*}
\rho_{I+1}  &  =R_{n-\left(  I+1\right)  +1}\circ R_{n-\left(  I+1\right)
+2}\circ...\circ R_{n}=R_{n-I}\circ R_{n-I+1}\circ...\circ R_{n}\\
&  =R_{n-I}\circ\underbrace{\left(  R_{n-I+1}\circ R_{n-I+2}\circ...\circ
R_{n}\right)  }_{=\rho_{I}}=R_{n-I}\circ\rho_{I},
\end{align*}
qed.}. We are going to show that $\rho_{I+1}f$ is well-defined.

In fact, let $v\in\widehat{P}_{n-I}$ be arbitrary. Then, $v\in\widehat{P}%
_{n-I}=\deg^{-1}\left(  \left\{  n-I\right\}  \right)  $ (by the definition of
$\widehat{P}_{n-I}$), so that $\deg v\in\left\{  n-I\right\}  $. Hence, $\deg
v=n-I$.

We have $v\in P$\ \ \ \ \footnote{\textit{Proof.} We have $\deg v=n-I\in
\left\{  1,2,...,n\right\}  $. In particular, $\deg v\neq0=\deg0$, so that
$v\neq0$. Also, $\deg v\in\left\{  1,2,...,n\right\}  $, so that $\deg v\neq
n+1=\deg1$, and thus $v\neq1$. Combined with $v\neq0$, this yields
$v\notin\left\{  0,1\right\}  $. Since $v\in\widehat{P}_{i}\subseteq
\widehat{P}$ but $v\notin\left\{  0,1\right\}  $, we have $v\in\widehat{P}%
\setminus\left\{  0,1\right\}  =P$, qed.}. Recall that $Rf$ is well-defined.
Proposition \ref{prop.Ri.implicit} yields%
\begin{equation}
\left(  Rf\right)  \left(  v\right)  =\dfrac{1}{f\left(  v\right)  }%
\cdot\dfrac{\sum\limits_{\substack{u\in\widehat{P};\\u\lessdot v}}f\left(
u\right)  }{\sum\limits_{\substack{u\in\widehat{P};\\u\gtrdot v}}\dfrac
{1}{\left(  Rf\right)  \left(  u\right)  }}. \label{pf.Ri.R.step.wd.1}%
\end{equation}
We have
\begin{equation}
f\left(  u\right)  =\left(  \rho_{I}f\right)  \left(  u\right)
\ \ \ \ \ \ \ \ \ \ \text{for every }u\in\widehat{P}\text{ such that
}u\lessdot v. \label{pf.Ri.R.step.wd.2a}%
\end{equation}
\footnote{\textit{Proof of (\ref{pf.Ri.R.step.wd.2a}):} Let $u\in\widehat{P}$
be such that $u\lessdot v$. Then, Proposition \ref{prop.graded.Phat.why}
\textbf{(a)} yields $\deg u=\deg v-1<\deg v=n-I\leq n-I$. But
(\ref{pf.Ri.R.step.1}) (applied to $u$ instead of $v$) yields
\[
\left(  \rho_{I}f\right)  \left(  u\right)  =\left\{
\begin{array}
[c]{c}%
\left(  Rf\right)  \left(  u\right)  ,\ \ \ \ \ \ \ \ \ \ \text{if }\deg
u>n-I;\\
f\left(  u\right)  ,\ \ \ \ \ \ \ \ \ \ \text{if }\deg u\leq n-I
\end{array}
\right.  =f\left(  u\right)  \ \ \ \ \ \ \ \ \ \ \left(  \text{since }\deg
u\leq n-I\right)  .
\]
Thus, (\ref{pf.Ri.R.step.wd.2a}) is proven.} Thus,%
\begin{equation}
\sum\limits_{\substack{u\in\widehat{P};\\u\lessdot v}}\underbrace{f\left(
u\right)  }_{\substack{=\left(  \rho_{I}f\right)  \left(  u\right)
\\\text{(by (\ref{pf.Ri.R.step.wd.2a}))}}}=\sum\limits_{\substack{u\in
\widehat{P};\\u\lessdot v}}\left(  \rho_{I}f\right)  \left(  u\right)  .
\label{pf.Ri.R.step.wd.2as}%
\end{equation}
Also,
\begin{equation}
\left(  Rf\right)  \left(  u\right)  =\left(  \rho_{I}f\right)  \left(
u\right)  \ \ \ \ \ \ \ \ \ \ \text{for every }u\in\widehat{P}\text{
satisfying }u\gtrdot v. \label{pf.Ri.R.step.wd.2b}%
\end{equation}
\footnote{\textit{Proof of (\ref{pf.Ri.R.step.wd.2b}):} Let $u\in\widehat{P}$
be such that $u\gtrdot v$. Then, $v\lessdot u$. Hence, Proposition
\ref{prop.graded.Phat.why} \textbf{(a)} (applied to $v$ and $u$ instead of $u$
and $v$) yields $\deg v=\deg u-1$, so that $\deg u=\deg v+1>\deg v=n-I$. But
(\ref{pf.Ri.R.step.1}) (applied to $u$ instead of $v$) yields
\[
\left(  \rho_{I}f\right)  \left(  u\right)  =\left\{
\begin{array}
[c]{c}%
\left(  Rf\right)  \left(  u\right)  ,\ \ \ \ \ \ \ \ \ \ \text{if }\deg
u>n-I;\\
f\left(  u\right)  ,\ \ \ \ \ \ \ \ \ \ \text{if }\deg u\leq n-I
\end{array}
\right.  =\left(  Rf\right)  \left(  u\right)  \ \ \ \ \ \ \ \ \ \ \left(
\text{since }\deg u>n-I\right)  .
\]
Thus, (\ref{pf.Ri.R.step.wd.2b}) is proven.} Hence,%
\begin{equation}
\sum\limits_{\substack{u\in\widehat{P};\\u\gtrdot v}}\underbrace{\dfrac
{1}{\left(  Rf\right)  \left(  u\right)  }}_{\substack{=\dfrac{1}{\left(
\rho_{I}f\right)  \left(  u\right)  }\\\text{(because
(\ref{pf.Ri.R.step.wd.2b})}\\\text{yields }\left(  Rf\right)  \left(
u\right)  =\left(  \rho_{I}f\right)  \left(  u\right)  \text{)}}%
}=\sum\limits_{\substack{u\in\widehat{P};\\u\gtrdot v}}\dfrac{1}{\left(
\rho_{I}f\right)  \left(  u\right)  }. \label{pf.Ri.R.step.wd.2bs}%
\end{equation}
Finally,%
\begin{equation}
f\left(  v\right)  =\left(  \rho_{I}f\right)  \left(  v\right)
\label{pf.Ri.R.step.wd.2v}%
\end{equation}
\footnote{\textit{Proof of (\ref{pf.Ri.R.step.wd.2v}):} Applying
(\ref{pf.Ri.R.step.1}), we obtain%
\[
\left(  \rho_{I}f\right)  \left(  v\right)  =\left\{
\begin{array}
[c]{c}%
\left(  Rf\right)  \left(  v\right)  ,\ \ \ \ \ \ \ \ \ \ \text{if }\deg
v>n-I;\\
f\left(  v\right)  ,\ \ \ \ \ \ \ \ \ \ \text{if }\deg v\leq n-I
\end{array}
\right.  =f\left(  v\right)  \ \ \ \ \ \ \ \ \ \ \left(  \text{since }\deg
v=n-I\leq n-I\right)  .
\]
This proves (\ref{pf.Ri.R.step.wd.2v}).}. Now, (\ref{pf.Ri.R.step.wd.1})
becomes%
\begin{equation}
\left(  Rf\right)  \left(  v\right)  =\dfrac{1}{f\left(  v\right)  }%
\cdot\dfrac{\sum\limits_{\substack{u\in\widehat{P};\\u\lessdot v}}f\left(
u\right)  }{\sum\limits_{\substack{u\in\widehat{P};\\u\gtrdot v}}\dfrac
{1}{\left(  Rf\right)  \left(  u\right)  }}=\dfrac{1}{\left(  \rho
_{I}f\right)  \left(  v\right)  }\cdot\dfrac{\sum\limits_{\substack{u\in
\widehat{P};\\u\lessdot v}}\left(  \rho_{I}f\right)  \left(  u\right)  }%
{\sum\limits_{\substack{u\in\widehat{P};\\u\gtrdot v}}\dfrac{1}{\left(
\rho_{I}f\right)  \left(  u\right)  }} \label{pf.Ri.R.step.Rfv}%
\end{equation}
(by (\ref{pf.Ri.R.step.wd.2v}), (\ref{pf.Ri.R.step.wd.2as}) and
(\ref{pf.Ri.R.step.wd.2bs})). In particular, $\dfrac{1}{\left(  \rho
_{I}f\right)  \left(  v\right)  }\cdot\dfrac{\sum\limits_{\substack{u\in
\widehat{P};\\u\lessdot v}}\left(  \rho_{I}f\right)  \left(  u\right)  }%
{\sum\limits_{\substack{u\in\widehat{P};\\u\gtrdot v}}\dfrac{1}{\left(
\rho_{I}f\right)  \left(  u\right)  }}$ is well-defined.

Now, forget that we fixed $v$. We thus have shown that $\dfrac{1}{\left(
\rho_{I}f\right)  \left(  v\right)  }\cdot\dfrac{\sum\limits_{\substack{u\in
\widehat{P};\\u\lessdot v}}\left(  \rho_{I}f\right)  \left(  u\right)  }%
{\sum\limits_{\substack{u\in\widehat{P};\\u\gtrdot v}}\dfrac{1}{\left(
\rho_{I}f\right)  \left(  u\right)  }}$ is well-defined for every
$v\in\widehat{P}_{n-I}$. Thus, Proposition \ref{prop.Ri.implicit.wd} (applied
to $n-I$ and $\rho_{I}f$ instead of $i$ and $f$) yields that $R_{n-I}\left(
\rho_{I}f\right)  $ is well-defined. In other words, $\rho_{I+1}f$ is
well-defined (since $\underbrace{\rho_{I+1}}_{=R_{n-I}\circ\rho_{I}}f=\left(
R_{n-I}\circ\rho_{I}\right)  f=R_{n-I}\left(  \rho_{I}f\right)  $).

We are now going to prove that every $v\in\widehat{P}$ satisfies%
\begin{equation}
\left(  \rho_{I+1}f\right)  \left(  v\right)  =\left\{
\begin{array}
[c]{c}%
\left(  Rf\right)  \left(  v\right)  ,\ \ \ \ \ \ \ \ \ \ \text{if }\deg
v>n-\left(  I+1\right)  ;\\
f\left(  v\right)  ,\ \ \ \ \ \ \ \ \ \ \text{if }\deg v\leq n-\left(
I+1\right)
\end{array}
\right.  . \label{pf.Ri.R.step.goal}%
\end{equation}

\textit{Proof of (\ref{pf.Ri.R.step.goal}):} Let $v\in\widehat{P}$. We need to
prove that (\ref{pf.Ri.R.step.goal}) holds.

We distinguish between two cases:

\textit{Case 1:} We have $\deg v\neq n-I$.

\textit{Case 2:} We have $\deg v=n-I$.

Let us deal with Case 1 first. In this case, we have $\deg v\neq n-I$. Hence,
Proposition \ref{prop.Ri.implicit} \textbf{(a)} (applied to $n-I$ and
$\rho_{I}f$ instead of $i$ and $f$) yields $\left(  R_{n-I}\left(  \rho
_{I}f\right)  \right)  \left(  v\right)  =\left(  \rho_{I}f\right)  \left(
v\right)  $. Now,%
\begin{align}
\left(  \underbrace{\rho_{I+1}}_{=R_{n-I}\circ\rho_{I}}f\right)  \left(
v\right)   &  =\underbrace{\left(  \left(  R_{n-I}\circ\rho_{I}\right)
\left(  f\right)  \right)  }_{=R_{n-I}\left(  \rho_{I}f\right)  }\left(
v\right)  =\left(  R_{n-I}\left(  \rho_{I}f\right)  \right)  \left(  v\right)
\nonumber\\
&  =\left(  \rho_{I}f\right)  \left(  v\right)  =\left\{
\begin{array}
[c]{c}%
\left(  Rf\right)  \left(  v\right)  ,\ \ \ \ \ \ \ \ \ \ \text{if }\deg
v>n-I;\\
f\left(  v\right)  ,\ \ \ \ \ \ \ \ \ \ \text{if }\deg v\leq n-I
\end{array}
\right.  \label{pf.Ri.R.step.c1}%
\end{align}
(by (\ref{pf.Ri.R.step.1})).

But we have $\deg v\neq n-I$. Hence, we must be in one of the following two subcases:

\textit{Subcase 1.1:} We have $\deg v<n-I$.

\textit{Subcase 1.2:} We have $\deg v>n-I$.

Let us first consider Subcase 1.1. In this subcase, we have $\deg v<n-I$.
Hence, $\deg v\leq n-I-1$ (since $\deg v$ and $n-I$ are integers), so that
$\deg v\leq n-I-1=n-\left(  I+1\right)  $. Also, $\deg v\leq n-I$ (since $\deg
v<n-I$). Thus, (\ref{pf.Ri.R.step.c1}) becomes%
\[
\left(  \rho_{I+1}f\right)  \left(  v\right)  =\left\{
\begin{array}
[c]{c}%
\left(  Rf\right)  \left(  v\right)  ,\ \ \ \ \ \ \ \ \ \ \text{if }\deg
v>n-I;\\
f\left(  v\right)  ,\ \ \ \ \ \ \ \ \ \ \text{if }\deg v\leq n-I
\end{array}
\right.  =f\left(  v\right)  \ \ \ \ \ \ \ \ \ \ \left(  \text{since }\deg
v<n-I\right)  .
\]
Compared with%
\[
\left\{
\begin{array}
[c]{c}%
\left(  Rf\right)  \left(  v\right)  ,\ \ \ \ \ \ \ \ \ \ \text{if }\deg
v>n-\left(  I+1\right)  ;\\
f\left(  v\right)  ,\ \ \ \ \ \ \ \ \ \ \text{if }\deg v\leq n-\left(
I+1\right)
\end{array}
\right.  =f\left(  v\right)  \ \ \ \ \ \ \ \ \ \ \left(  \text{since }\deg
v\leq n-\left(  I+1\right)  \right)  ,
\]
this yields%
\[
\left(  \rho_{I+1}f\right)  \left(  v\right)  =\left\{
\begin{array}
[c]{c}%
\left(  Rf\right)  \left(  v\right)  ,\ \ \ \ \ \ \ \ \ \ \text{if }\deg
v>n-\left(  I+1\right)  ;\\
f\left(  v\right)  ,\ \ \ \ \ \ \ \ \ \ \text{if }\deg v\leq n-\left(
I+1\right)
\end{array}
\right.  .
\]
Hence, (\ref{pf.Ri.R.step.goal}) is proven in Subcase 1.1.

Let us now proceed to Subcase 1.2. In this subcase, we have $\deg v>n-I$.
Hence, $\deg v>n-I>n-I-1=n-\left(  I+1\right)  $. Now, (\ref{pf.Ri.R.step.c1})
becomes%
\[
\left(  \rho_{I+1}f\right)  \left(  v\right)  =\left\{
\begin{array}
[c]{c}%
\left(  Rf\right)  \left(  v\right)  ,\ \ \ \ \ \ \ \ \ \ \text{if }\deg
v>n-I;\\
f\left(  v\right)  ,\ \ \ \ \ \ \ \ \ \ \text{if }\deg v\leq n-I
\end{array}
\right.  =\left(  Rf\right)  \left(  v\right)  \ \ \ \ \ \ \ \ \ \ \left(
\text{since }\deg v>n-I\right)  .
\]
Compared with%
\[
\left\{
\begin{array}
[c]{c}%
\left(  Rf\right)  \left(  v\right)  ,\ \ \ \ \ \ \ \ \ \ \text{if }\deg
v>n-\left(  I+1\right)  ;\\
f\left(  v\right)  ,\ \ \ \ \ \ \ \ \ \ \text{if }\deg v\leq n-\left(
I+1\right)
\end{array}
\right.  =\left(  Rf\right)  \left(  v\right)  \ \ \ \ \ \ \ \ \ \ \left(
\text{since }\deg v>n-\left(  I+1\right)  \right)  ,
\]
this yields%
\[
\left(  \rho_{I+1}f\right)  \left(  v\right)  =\left\{
\begin{array}
[c]{c}%
\left(  Rf\right)  \left(  v\right)  ,\ \ \ \ \ \ \ \ \ \ \text{if }\deg
v>n-\left(  I+1\right)  ;\\
f\left(  v\right)  ,\ \ \ \ \ \ \ \ \ \ \text{if }\deg v\leq n-\left(
I+1\right)
\end{array}
\right.  .
\]
Hence, (\ref{pf.Ri.R.step.goal}) is proven in Subcase 1.2.

We have thus proven (\ref{pf.Ri.R.step.goal}) in each of the two Subcases 1.1
and 1.2. Since these two Subcases cover all of Case 1, this yields that
(\ref{pf.Ri.R.step.goal}) always holds in Case 1.

Now, let us consider Case 2. In this case, we have $\deg v=n-I$. Hence,
Proposition \ref{prop.Ri.implicit} \textbf{(b)} (applied to $n-I$ and
$\rho_{I}f$ instead of $i$ and $f$) yields%
\begin{equation}
\left(  R_{n-I}\left(  \rho_{I}f\right)  \right)  \left(  v\right)  =\dfrac
{1}{\left(  \rho_{I}f\right)  \left(  v\right)  }\cdot\dfrac{\sum
\limits_{\substack{u\in\widehat{P};\\u\lessdot v}}\left(  \rho_{I}f\right)
\left(  u\right)  }{\sum\limits_{\substack{u\in\widehat{P};\\u\gtrdot
v}}\dfrac{1}{\left(  \rho_{I}f\right)  \left(  u\right)  }}=\left(  Rf\right)
\left(  v\right)  \label{pf.Ri.R.step.Rfv2}%
\end{equation}
(by (\ref{pf.Ri.R.step.Rfv})). But from (\ref{pf.Ri.R.step.1}), we obtain%
\[
\left(  \rho_{I}f\right)  \left(  v\right)  =\left\{
\begin{array}
[c]{c}%
\left(  Rf\right)  \left(  v\right)  ,\ \ \ \ \ \ \ \ \ \ \text{if }\deg
v>n-I;\\
f\left(  v\right)  ,\ \ \ \ \ \ \ \ \ \ \text{if }\deg v\leq n-I
\end{array}
\right.  =f\left(  v\right)
\]
(since $\deg v=n-I\leq n-I$). Notice that $\underbrace{\rho_{I+1}}%
_{=R_{n-I}\circ\rho_{I}}f=\left(  R_{n-I}\circ\rho_{I}\right)  f=R_{n-I}%
\left(  \rho_{I}f\right)  $. Hence,%
\[
\left(  \rho_{I+1}f\right)  \left(  v\right)  =\left(  R_{n-I}\left(  \rho
_{I}f\right)  \right)  \left(  v\right)  =\left(  Rf\right)  \left(  v\right)
\ \ \ \ \ \ \ \ \ \ \left(  \text{by (\ref{pf.Ri.R.step.Rfv2})}\right)  .
\]
Compared with%
\begin{align*}
&  \left\{
\begin{array}
[c]{c}%
\left(  Rf\right)  \left(  v\right)  ,\ \ \ \ \ \ \ \ \ \ \text{if }\deg
v>n-\left(  I+1\right)  ;\\
f\left(  v\right)  ,\ \ \ \ \ \ \ \ \ \ \text{if }\deg v\leq n-\left(
I+1\right)
\end{array}
\right. \\
&  =\left(  Rf\right)  \left(  v\right)  \ \ \ \ \ \ \ \ \ \ \left(
\text{since }\deg v=n-I>n-I-1=n-\left(  I+1\right)  \right)  ,
\end{align*}
this yields
\[
\left(  \rho_{I+1}f\right)  \left(  v\right)  =\left\{
\begin{array}
[c]{c}%
\left(  Rf\right)  \left(  v\right)  ,\ \ \ \ \ \ \ \ \ \ \text{if }\deg
v>n-\left(  I+1\right)  ;\\
f\left(  v\right)  ,\ \ \ \ \ \ \ \ \ \ \text{if }\deg v\leq n-\left(
I+1\right)
\end{array}
\right.  .
\]
Hence, (\ref{pf.Ri.R.step.goal}) is proven in Case 2.

We thus have proven (\ref{pf.Ri.R.step.goal}) in each of the two Cases 1 and
2. Since these two Cases cover all possibilities, this yields that
(\ref{pf.Ri.R.step.goal}) always holds. In other words, the proof of
(\ref{pf.Ri.R.step.goal}) is complete.

Altogether, we thus have proven that the $\mathbb{K}$-labelling $\rho
_{I+1}f\in\mathbb{K}^{\widehat{P}}$ is well-defined, and every $v\in
\widehat{P}$ satisfies%
\[
\left(  \rho_{I+1}f\right)  \left(  v\right)  =\left\{
\begin{array}
[c]{c}%
\left(  Rf\right)  \left(  v\right)  ,\ \ \ \ \ \ \ \ \ \ \text{if }\deg
v>n-\left(  I+1\right)  ;\\
f\left(  v\right)  ,\ \ \ \ \ \ \ \ \ \ \text{if }\deg v\leq n-\left(
I+1\right)
\end{array}
\right.  .
\]
In other words, Assertion A is satisfied for $i=I+1$. This completes the
induction step. The induction proof of Assertion A is thus complete.

Now, we can apply Assertion A to $i=n$ (because Assertion A is proven for
every $i\in\left\{  0,1,...,n\right\}  $). As a result, we conclude that the
$\mathbb{K}$-labelling $\rho_{n}f\in\mathbb{K}^{\widehat{P}}$ is well-defined,
and that every $v\in\widehat{P}$ satisfies%
\begin{equation}
\left(  \rho_{n}f\right)  \left(  v\right)  =\left\{
\begin{array}
[c]{c}%
\left(  Rf\right)  \left(  v\right)  ,\ \ \ \ \ \ \ \ \ \ \text{if }\deg
v>n-n;\\
f\left(  v\right)  ,\ \ \ \ \ \ \ \ \ \ \text{if }\deg v\leq n-n
\end{array}
\right.  . \label{pf.Ri.R.almostthere}%
\end{equation}

We know that the $\mathbb{K}$-labelling $\rho_{n}f\in\mathbb{K}^{\widehat{P}}$
is well-defined. We shall now show that $\rho_{n}f=Rf$.

In fact, let $v\in\widehat{P}$ be arbitrary. We are going to prove that
$\left(  \rho_{n}f\right)  \left(  v\right)  =\left(  Rf\right)  \left(
v\right)  $.

If $v=0$, then $\left(  \rho_{n}f\right)  \left(  v\right)  =\left(
Rf\right)  \left(  v\right)  $ can be checked
immediately\footnote{\textit{Proof.} Assume that $v=0$. Then, $\deg v=\deg
0=0$. Now, (\ref{pf.Ri.R.almostthere}) yields%
\begin{align*}
\left(  \rho_{n}f\right)  \left(  v\right)   &  =\left\{
\begin{array}
[c]{c}%
\left(  Rf\right)  \left(  v\right)  ,\ \ \ \ \ \ \ \ \ \ \text{if }\deg
v>n-n;\\
f\left(  v\right)  ,\ \ \ \ \ \ \ \ \ \ \text{if }\deg v\leq n-n
\end{array}
\right.  =f\left(  \underbrace{v}_{=0}\right)  \ \ \ \ \ \ \ \ \ \ \left(
\text{since }\deg v=0=n-n\leq n-n\right) \\
&  =f\left(  0\right)  =\left(  Rf\right)  \left(  0\right)
\ \ \ \ \ \ \ \ \ \ \left(  \text{since Proposition \ref{prop.R.implicit.01}
yields }\left(  Rf\right)  \left(  0\right)  =f\left(  0\right)  \right)  .
\end{align*}
This rewrites as $\left(  \rho_{n}f\right)  \left(  v\right)  =\left(
Rf\right)  \left(  v\right)  $ (since $v=0$), qed.}. Hence, for the rest of
the proof of $\left(  \rho_{n}f\right)  \left(  v\right)  =\left(  Rf\right)
\left(  v\right)  $, we can WLOG assume that we don't have $v=0$. Assume this.
Now, we have $\deg v>0$\ \ \ \ \footnote{\textit{Proof.} If $v\in P$, then
$\deg v\in\left\{  1,2,...,n\right\}  $ (because $\deg$ was defined to be a
map $P\rightarrow\left\{  1,2,...,n\right\}  $). Hence, if $v\in P$, then
$\deg v>0$. Thus, the proof of $\deg v>0$ is complete if $v\in P$. Thus, for
the rest of the proof of $\deg v>0$, we can WLOG assume that we don't have
$v\in P$. Assume this. We don't have $v\in P$. So we have $v\notin P$. Since
$v\in\widehat{P}$ but $v\notin P$, we have $v\in\widehat{P}\setminus
P=\left\{  0,1\right\}  $. But since we don't have $v=0$, this yields
$v\in\left\{  0,1\right\}  \setminus\left\{  0\right\}  =\left\{  1\right\}
$, so that $v=1$ and thus $\deg v=\deg1=\underbrace{n}_{\geq0}+1\geq1>0$,
qed.}. The equality (\ref{pf.Ri.R.almostthere}) yields%
\[
\left(  \rho_{n}f\right)  \left(  v\right)  =\left\{
\begin{array}
[c]{c}%
\left(  Rf\right)  \left(  v\right)  ,\ \ \ \ \ \ \ \ \ \ \text{if }\deg
v>n-n;\\
f\left(  v\right)  ,\ \ \ \ \ \ \ \ \ \ \text{if }\deg v\leq n-n
\end{array}
\right.  =\left(  Rf\right)  \left(  v\right)
\]
(since $\deg v>0=n-n$). This completes the proof of $\left(  \rho_{n}f\right)
\left(  v\right)  =\left(  Rf\right)  \left(  v\right)  $.

Now, forget that we fixed $v$. We thus have shown that $\left(  \rho
_{n}f\right)  \left(  v\right)  =\left(  Rf\right)  \left(  v\right)  $ for
every $v\in\widehat{P}$. In other words, $\rho_{n}f=Rf$.

Now, forget that we fixed $f$. We thus have shown that whenever $f\in
\mathbb{K}^{\widehat{P}}$ is a $\mathbb{K}$-labelling of $P$ for which $Rf$ is
well-defined, we have $\rho_{n}f=Rf$. In other words, $\rho_{n}=R$ (as
rational maps). Since $\rho_{n}=R_{1}\circ R_{2}\circ...\circ R_{n}$, this
rewrites as $R_{1}\circ R_{2}\circ...\circ R_{n}=R$. This proves Proposition
\ref{prop.Ri.R}.
\end{proof}
\end{verlong}

\section{w-tuples}

This section continues the study of birational rowmotion on graded posets by
introducing a ``fingerprint'' or ``checksum'' of a $\mathbb{K}$-labelling
called the w-tuple, defined by summing ratios of elements between successive
degrees (i.e., rows in the Hasse diagram). This w-tuple serves to extract some
information from a $\mathbb{K}$-labelling; we will later see how to make the
``rest'' of the labelling more manageable.

\begin{definition}
\label{def.wi}Let $n\in\mathbb{N}$. Let $\mathbb{K}$ be a field. Let $P$ be an
$n$-graded poset. Let $f\in\mathbb{K}^{\widehat{P}}$. Let $i\in\left\{
0,1,...,n\right\}  $. Then, $\mathbf{w}_{i}\left(  f\right)  $ will denote the
element of $\mathbb{K}$ defined by%
\[
\mathbf{w}_{i}\left(  f\right)  =\sum_{\substack{x\in\widehat{P}_{i}%
;\ y\in\widehat{P}_{i+1};\\y\gtrdot x}}\dfrac{f\left(  x\right)  }{f\left(
y\right)  }.
\]
(This element is not always defined, but is defined in the ``generic'' case
when $0\notin f\left(  \widehat{P}\right)  $.)
\end{definition}

Intuitively, one could think of $\mathbf{w}_{i}\left(  f\right)  $ as a kind
of \textquotedblleft checksum\textquotedblright\ for the labelling $f$ which
displays how much its labels at degree $i+1$ differ from those at degree $i$.
Of course, in general, the knowledge of $\mathbf{w}_{i}\left(  f\right)  $ for
all $i\in\left\{  0,1,...,n\right\}  $ is far from sufficient to reconstruct
the whole labelling $f$; however, in Definition \ref{def.hgeq}, we will
introduce the so-called homogenization of $f$, which will provide
\textquotedblleft complementary data\textquotedblright\ to these
$\mathbf{w}_{i}\left(  f\right)  $. As for now, let us show that the
$\mathbf{w}_{i}\left(  f\right)  $ behave in a rather simple way under the
maps $R$ and $R_{j}$.

\begin{definition}
Let $n\in\mathbb{N}$. Let $\mathbb{K}$ be a field. Let $P$ be an $n$-graded
poset. Let $f\in\mathbb{K}^{\widehat{P}}$. The $\left(  n+1\right)  $-tuple
$\left(  \mathbf{w}_{0}\left(  f\right)  ,\mathbf{w}_{1}\left(  f\right)
,...,\mathbf{w}_{n}\left(  f\right)  \right)  $ will be called the
\textit{w-tuple} of the $\mathbb{K}$-labelling $f$.
\end{definition}

It is easy to see:

\begin{proposition}
\label{prop.wi.Ri}Let $n\in\mathbb{N}$. Let $\mathbb{K}$ be a field. Let $P$
be an $n$-graded poset. Let $i\in\left\{  1,2,...,n\right\}  $. Then, every
$f\in\mathbb{K}^{\widehat{P}}$ satisfies%
\begin{align*}
&  \left(  \mathbf{w}_{0}\left(  R_{i}f\right)  ,\mathbf{w}_{1}\left(
R_{i}f\right)  ,...,\mathbf{w}_{n}\left(  R_{i}f\right)  \right) \\
&  =\left(  \mathbf{w}_{0}\left(  f\right)  ,\mathbf{w}_{1}\left(  f\right)
,...,\mathbf{w}_{i-2}\left(  f\right)  ,\mathbf{w}_{i}\left(  f\right)
,\mathbf{w}_{i-1}\left(  f\right)  ,\mathbf{w}_{i+1}\left(  f\right)
,\mathbf{w}_{i+2}\left(  f\right)  ,...,\mathbf{w}_{n}\left(  f\right)
\right)  .
\end{align*}

In other words, the map $R_{i}$ changes the w-tuple of a $\mathbb{K}%
$-labelling by interchanging its $\left(  i-1\right)  $-st entry with its
$i$-th entry (where the entries are labelled starting at $0$).
\end{proposition}

\begin{proof}
[Proof of Proposition \ref{prop.wi.Ri} (sketched).]%
\begin{verlong}
First we notice something simple: For each $x\in\widehat{P}_{i}$, we have that%
\begin{equation}
\left(  \text{every element }u\text{ of }\widehat{P}\text{ satisfying
}u\gtrdot x\text{ satisfies }u\in\widehat{P}_{i+1}\right)
\label{pf.wi.Ri.triv1}%
\end{equation}
and%
\begin{equation}
\left(  \text{every element }u\text{ of }\widehat{P}\text{ satisfying
}u\lessdot x\text{ satisfies }u\in\widehat{P}_{i-1}\right)  .
\label{pf.wi.Ri.triv2}%
\end{equation}
These facts follow from the definition of an $n$-graded poset (specifically,
from Assertion 1 in Definition \ref{def.graded.n-graded}).
\end{verlong}

Let $f\in\mathbb{K}^{\widehat{P}}$. We need to show that every $j\in\left\{
0,1,...,n\right\}  $ satisfies%
\begin{equation}
\mathbf{w}_{j}\left(  R_{i}f\right)  =\mathbf{w}_{\tau_{i}\left(  j\right)
}\left(  f\right)  , \label{pf.wi.Ri.1}%
\end{equation}
where $\tau_{i}$ is the permutation of the set $\left\{  0,1,...,n\right\}  $
which transposes $i-1$ with $i$ (while leaving all other elements of this set invariant).

\textit{Proof of (\ref{pf.wi.Ri.1}):} Let $j\in\left\{  0,1,...,n\right\}  $.
We distinguish between three cases:

\textit{Case 1:} We have $j=i$.

\textit{Case 2:} We have $j=i-1$.

\textit{Case 3:} We have $j\notin\left\{  i-1,i\right\}  $.

\begin{vershort}
Let us first consider Case 1. In this case, we have $j=i$. By the definition
of $\mathbf{w}_{i}\left(  R_{i}f\right)  $, we have%
\begin{align}
\mathbf{w}_{i}\left(  R_{i}f\right)   &  =\sum_{\substack{x\in\widehat{P}%
_{i};\ y\in\widehat{P}_{i+1};\\y\gtrdot x}}\dfrac{\left(  R_{i}f\right)
\left(  x\right)  }{\left(  R_{i}f\right)  \left(  y\right)  }=\sum
\limits_{x\in\widehat{P}_{i}}\left(  R_{i}f\right)  \left(  x\right)
\sum\limits_{\substack{y\in\widehat{P}_{i+1};\\y\gtrdot x}}\left(
\underbrace{\left(  R_{i}f\right)  \left(  y\right)  }_{\substack{=f\left(
y\right)  \\\text{(by Proposition \ref{prop.Ri.implicit} \textbf{(a)})}%
}}\right)  ^{-1}\nonumber\\
&  =\sum\limits_{x\in\widehat{P}_{i}}\left(  R_{i}f\right)  \left(  x\right)
\sum\limits_{\substack{y\in\widehat{P}_{i+1};\\y\gtrdot x}}\left(  f\left(
y\right)  \right)  ^{-1}. \label{pf.wi.Ri.c1.short.1}%
\end{align}
But every $x\in\widehat{P}_{i}$ satisfies%
\begin{align*}
\left(  R_{i}f\right)  \left(  x\right)   &  =\dfrac{1}{f\left(  x\right)
}\cdot\dfrac{\sum\limits_{\substack{u\in\widehat{P};\\u\lessdot x}}f\left(
u\right)  }{\sum\limits_{\substack{u\in\widehat{P};\\u\gtrdot x}}\dfrac
{1}{f\left(  u\right)  }}\ \ \ \ \ \ \ \ \ \ \left(  \text{by Proposition
\ref{prop.Ri.implicit} \textbf{(b)}}\right) \\
&  =\dfrac{1}{f\left(  x\right)  }\cdot\sum\limits_{\substack{u\in
\widehat{P};\\u\lessdot x}}f\left(  u\right)  \cdot\left(  \sum
\limits_{\substack{u\in\widehat{P};\\u\gtrdot x}}\left(  f\left(  u\right)
\right)  ^{-1}\right)  ^{-1}=\dfrac{1}{f\left(  x\right)  }\cdot
\sum\limits_{\substack{y\in\widehat{P};\\y\lessdot x}}f\left(  y\right)
\cdot\left(  \sum\limits_{\substack{y\in\widehat{P};\\y\gtrdot x}}\left(
f\left(  y\right)  \right)  ^{-1}\right)  ^{-1}\\
&  =\dfrac{1}{f\left(  x\right)  }\cdot\sum\limits_{\substack{y\in
\widehat{P}_{i-1};\\y\lessdot x}}f\left(  y\right)  \cdot\left(
\sum\limits_{\substack{y\in\widehat{P}_{i+1};\\y\gtrdot x}}\left(  f\left(
y\right)  \right)  ^{-1}\right)  ^{-1}%
\end{align*}
(here, we replaced $y\in\widehat{P}$ by $y\in\widehat{P}_{i-1}$ in the first
sum (because every $y\in\widehat{P}$ satisfying $y\lessdot x$ must belong to
$\widehat{P}_{i-1}$\ \ \ \ \footnote{since $x\in\widehat{P}_{i}$}) and we
replaced $y\in\widehat{P}$ by $y\in\widehat{P}_{i+1}$ in the second sum (for
similar reasons)) and thus%
\[
\left(  R_{i}f\right)  \left(  x\right)  \sum\limits_{\substack{y\in
\widehat{P}_{i+1};\\y\gtrdot x}}\left(  f\left(  y\right)  \right)
^{-1}=\dfrac{1}{f\left(  x\right)  }\cdot\sum\limits_{\substack{y\in
\widehat{P}_{i-1};\\y\lessdot x}}f\left(  y\right)  =\sum
\limits_{\substack{y\in\widehat{P}_{i-1};\\y\lessdot x}}\dfrac{f\left(
y\right)  }{f\left(  x\right)  }.
\]
Hence, (\ref{pf.wi.Ri.c1.short.1}) becomes%
\begin{align}
\mathbf{w}_{i}\left(  R_{i}f\right)   &  =\sum\limits_{x\in\widehat{P}_{i}%
}\underbrace{\left(  R_{i}f\right)  \left(  x\right)  \sum
\limits_{\substack{y\in\widehat{P}_{i+1};\\y\gtrdot x}}\left(  f\left(
y\right)  \right)  ^{-1}}_{=\sum\limits_{\substack{y\in\widehat{P}%
_{i-1};\\y\lessdot x}}\dfrac{f\left(  y\right)  }{f\left(  x\right)  }}%
=\sum\limits_{x\in\widehat{P}_{i}}\sum\limits_{\substack{y\in\widehat{P}%
_{i-1};\\y\lessdot x}}\dfrac{f\left(  y\right)  }{f\left(  x\right)  }%
=\sum\limits_{\substack{y\in\widehat{P}_{i-1};\ x\in\widehat{P}_{i};\\x\gtrdot
y}}\dfrac{f\left(  y\right)  }{f\left(  x\right)  }\nonumber\\
&  =\sum\limits_{\substack{x\in\widehat{P}_{i-1};\ y\in\widehat{P}%
_{i}\\y\gtrdot x}}\dfrac{f\left(  x\right)  }{f\left(  y\right)
}\ \ \ \ \ \ \ \ \ \ \left(  \text{here, we switched the indices in the
sum}\right) \nonumber\\
&  =\mathbf{w}_{i-1}\left(  f\right)  \ \ \ \ \ \ \ \ \ \ \left(  \text{by the
definition of }\mathbf{w}_{i-1}\left(  f\right)  \right)
\label{pf.wi.Ri.short.main}\\
&  =\mathbf{w}_{\tau_{i}\left(  i\right)  }\left(  f\right)  .\nonumber
\end{align}
In other words, (\ref{pf.wi.Ri.1}) holds for $j=i$. Thus, (\ref{pf.wi.Ri.1})
is proven in Case 1.
\end{vershort}

\begin{verlong}
Let us first consider Case 1. In this case, we have $j=i$. By the definition
of $\mathbf{w}_{i}\left(  R_{i}f\right)  $, we have%
\begin{align*}
\mathbf{w}_{i}\left(  R_{i}f\right)   &  =\sum_{\substack{x\in\widehat{P}%
_{i};\ y\in\widehat{P}_{i+1};\\y\gtrdot x}}\dfrac{\left(  R_{i}f\right)
\left(  x\right)  }{\left(  R_{i}f\right)  \left(  y\right)  }%
=\underbrace{\sum_{\substack{x\in\widehat{P}_{i};\ y\in\widehat{P}%
_{i+1};\\y\gtrdot x}}}_{=\sum\limits_{x\in\widehat{P}_{i}}\sum
\limits_{\substack{y\in\widehat{P}_{i+1};\\y\gtrdot x}}}\left(
\underbrace{\left(  R_{i}f\right)  \left(  y\right)  }_{\substack{=f\left(
y\right)  \\\text{(by Proposition \ref{prop.Ri.implicit} \textbf{(a)})}%
}}\right)  ^{-1}\cdot\left(  R_{i}f\right)  \left(  x\right) \\
&  =\sum\limits_{x\in\widehat{P}_{i}}\sum\limits_{\substack{y\in
\widehat{P}_{i+1};\\y\gtrdot x}}\left(  f\left(  y\right)  \right)  ^{-1}%
\cdot\left(  R_{i}f\right)  \left(  x\right) \\
&  =\sum\limits_{x\in\widehat{P}_{i}}\underbrace{\left(  \sum
\limits_{\substack{y\in\widehat{P}_{i+1};\\y\gtrdot x}}\left(  f\left(
y\right)  \right)  ^{-1}\right)  }_{\substack{=\sum\limits_{\substack{y\in
\widehat{P}_{i+1};\\y\gtrdot x}}\dfrac{1}{f\left(  y\right)  }=\sum
\limits_{\substack{u\in\widehat{P}_{i+1};\\u\gtrdot x}}\dfrac{1}{f\left(
u\right)  }\\=\sum\limits_{\substack{u\in\widehat{P};\\u\gtrdot x}}\dfrac
{1}{f\left(  u\right)  }\\\text{(here, we replaced the condition ``}%
u\in\widehat{P}_{i+1}\text{''}\\\text{by ``}u\in\widehat{P}\text{'' in the
sum, since every}\\\text{element }u\text{ of }\widehat{P}\text{ satisfying
}u\gtrdot x\text{ lies in }\widehat{P}_{i+1}\\\text{(by (\ref{pf.wi.Ri.triv1}%
)))}}}\underbrace{\left(  R_{i}f\right)  \left(  x\right)  }%
_{\substack{=\dfrac{1}{f\left(  x\right)  }\cdot\dfrac{\sum
\limits_{\substack{u\in\widehat{P};\\u\lessdot x}}f\left(  u\right)  }%
{\sum\limits_{\substack{u\in\widehat{P};\\u\gtrdot x}}\dfrac{1}{f\left(
u\right)  }}\\\text{(by Proposition \ref{prop.Ri.implicit} \textbf{(b)})}}}
\end{align*}%
\begin{align*}
&  =\sum\limits_{x\in\widehat{P}_{i}}\left(  \sum\limits_{\substack{u\in
\widehat{P};\\u\gtrdot x}}\dfrac{1}{f\left(  u\right)  }\right)  \dfrac
{1}{f\left(  x\right)  }\cdot\dfrac{\sum\limits_{\substack{u\in\widehat{P}%
;\\u\lessdot x}}f\left(  u\right)  }{\sum\limits_{\substack{u\in
\widehat{P};\\u\gtrdot x}}\dfrac{1}{f\left(  u\right)  }}\\
&  =\sum\limits_{x\in\widehat{P}_{i}}\dfrac{1}{f\left(  x\right)  }%
\cdot\underbrace{\left(  \sum\limits_{\substack{u\in\widehat{P};\\u\lessdot
x}}f\left(  u\right)  \right)  }_{\substack{=\sum\limits_{\substack{u\in
\widehat{P}_{i-1};\\u\lessdot x}}f\left(  u\right)  \\\text{(here, we replaced
the condition ``}u\in\widehat{P}\text{''}\\\text{by ``}u\in\widehat{P}%
_{i-1}\text{'' in the sum, since every}\\\text{element }u\text{ of
}\widehat{P}\text{ satisfying }u\lessdot x\\\text{satisfies }u\in
\widehat{P}_{i-1}\text{ (by (\ref{pf.wi.Ri.triv2})))}}}=\sum\limits_{x\in
\widehat{P}_{i}}\dfrac{1}{f\left(  x\right)  }\cdot\left(  \sum
\limits_{\substack{u\in\widehat{P}_{i-1};\\u\lessdot x}}f\left(  u\right)
\right) \\
&  =\underbrace{\sum\limits_{x\in\widehat{P}_{i}}\sum\limits_{\substack{u\in
\widehat{P}_{i-1};\\u\lessdot x}}}_{=\sum\limits_{\substack{u\in
\widehat{P}_{i-1};\ x\in\widehat{P}_{i};\\u\lessdot x}}=\sum
\limits_{\substack{u\in\widehat{P}_{i-1};\ x\in\widehat{P}_{i};\\x\gtrdot u}%
}}\underbrace{\dfrac{1}{f\left(  x\right)  }\cdot f\left(  u\right)
}_{=\dfrac{f\left(  u\right)  }{f\left(  x\right)  }}=\sum
\limits_{\substack{u\in\widehat{P}_{i-1};\ x\in\widehat{P}_{i};\\x\gtrdot
u}}\dfrac{f\left(  u\right)  }{f\left(  x\right)  }\\
&  =\sum\limits_{\substack{x\in\widehat{P}_{i-1};\ y\in\widehat{P}%
_{i};\\y\gtrdot x}}\dfrac{f\left(  x\right)  }{f\left(  y\right)
}\ \ \ \ \ \ \ \ \ \ \left(
\begin{array}
[c]{c}%
\text{here, we renamed the summation}\\
\text{indices }u\text{ and }x\text{ as }x\text{ and }y
\end{array}
\right)  .
\end{align*}
Compared with%
\[
\mathbf{w}_{i-1}\left(  f\right)  =\sum_{\substack{x\in\widehat{P}%
_{i-1};\ y\in\widehat{P}_{i};\\y\gtrdot x}}\dfrac{f\left(  x\right)
}{f\left(  y\right)  }\ \ \ \ \ \ \ \ \ \ \left(  \text{by the definition of
}\mathbf{w}_{i-1}\left(  f\right)  \right)  ,
\]
this yields
\begin{align}
\mathbf{w}_{i}\left(  R_{i}f\right)   &  =\mathbf{w}_{i-1}\left(  f\right)
\label{pf.wi.Ri.main}\\
&  =\mathbf{w}_{\tau_{i}\left(  i\right)  }\left(  f\right) \nonumber
\end{align}
(since $i-1=\tau_{i}\left(  i\right)  $). In other words, (\ref{pf.wi.Ri.1})
holds for $j=i$. Thus, (\ref{pf.wi.Ri.1}) is proven in Case 1.
\end{verlong}

\begin{vershort}
Let us now consider Case 2. In this case, $j=i-1$. Now, it can be shown that
$\mathbf{w}_{i-1}\left(  R_{i}f\right)  =\mathbf{w}_{i}\left(  f\right)  $.
This can be proven either in a similar way to how we proved $\mathbf{w}%
_{i}\left(  R_{i}f\right)  =\mathbf{w}_{i-1}\left(  f\right)  $ (the details
of this are left to the reader), or by noticing that%
\begin{align*}
\mathbf{w}_{i}\left(  f\right)   &  =\mathbf{w}_{i}\left(  R_{i}^{2}f\right)
\ \ \ \ \ \ \ \ \ \ \left(
\begin{array}
[c]{c}%
\text{since Proposition \ref{prop.Ri.invo} yields that }R_{i}^{2}%
=\operatorname*{id}\text{,}\\
\text{hence }\mathbf{w}_{i}\left(  R_{i}^{2}f\right)  =\mathbf{w}_{i}\left(
\operatorname*{id}f\right)  =\mathbf{w}_{i}\left(  f\right)
\end{array}
\right) \\
&  =\mathbf{w}_{i}\left(  R_{i}\left(  R_{i}f\right)  \right)  =\mathbf{w}%
_{i-1}\left(  R_{i}f\right)  \ \ \ \ \ \ \ \ \ \ \left(  \text{by
(\ref{pf.wi.Ri.short.main}), applied to }R_{i}f\text{ instead of }f\right)  .
\end{align*}
Either way, we end up knowing that $\mathbf{w}_{i-1}\left(  R_{i}f\right)
=\mathbf{w}_{i}\left(  f\right)  $. Thus, $\mathbf{w}_{i-1}\left(
R_{i}f\right)  =\mathbf{w}_{i}\left(  f\right)  =\mathbf{w}_{\tau_{i}\left(
i-1\right)  }\left(  f\right)  $. In other words, (\ref{pf.wi.Ri.1}) holds for
$j=i-1$. Thus, (\ref{pf.wi.Ri.1}) is proven in Case 2.
\end{vershort}

\begin{verlong}
Let us now consider Case 2. In this case, $j=i-1$. Now, it can be shown that
$\mathbf{w}_{i-1}\left(  R_{i}f\right)  =\mathbf{w}_{i}\left(  f\right)  $.
This can be proven either in a similar way to how we proved $\mathbf{w}%
_{i}\left(  R_{i}f\right)  =\mathbf{w}_{i-1}\left(  f\right)  $ (the details
of this are left to the reader), or by noticing that%
\begin{align*}
\mathbf{w}_{i}\left(  f\right)   &  =\mathbf{w}_{i}\left(  R_{i}^{2}f\right)
\ \ \ \ \ \ \ \ \ \ \left(
\begin{array}
[c]{c}%
\text{since Proposition \ref{prop.Ri.invo} yields that }R_{i}^{2}%
=\operatorname*{id}\text{,}\\
\text{hence }\mathbf{w}_{i}\left(  R_{i}^{2}f\right)  =\mathbf{w}_{i}\left(
\operatorname*{id}f\right)  =\mathbf{w}_{i}\left(  f\right)
\end{array}
\right) \\
&  =\mathbf{w}_{i}\left(  R_{i}\left(  R_{i}f\right)  \right)  =\mathbf{w}%
_{i-1}\left(  R_{i}f\right)  \ \ \ \ \ \ \ \ \ \ \left(  \text{by
(\ref{pf.wi.Ri.main}), applied to }R_{i}f\text{ instead of }f\right)  .
\end{align*}
Either way, we end up knowing that $\mathbf{w}_{i-1}\left(  R_{i}f\right)
=\mathbf{w}_{i}\left(  f\right)  $. Thus, $\mathbf{w}_{i-1}\left(
R_{i}f\right)  =\mathbf{w}_{i}\left(  f\right)  =\mathbf{w}_{\tau_{i}\left(
i-1\right)  }\left(  f\right)  $ (since $i=\tau_{i}\left(  i-1\right)  $). In
other words, (\ref{pf.wi.Ri.1}) holds for $j=i-1$. Thus, (\ref{pf.wi.Ri.1}) is
proven in Case 2.
\end{verlong}

Let us finally consider Case 3. In this case, $j\notin\left\{  i-1,i\right\}
$. Hence, $\tau_{i}\left(  j\right)  =j$. On the other hand, by the definition
of $\mathbf{w}_{j}\left(  R_{i}f\right)  $, we have%
\begin{align*}
\mathbf{w}_{j}\left(  R_{i}f\right)   &  =\sum_{\substack{x\in\widehat{P}%
_{j};\ y\in\widehat{P}_{j+1};\\y \gtrdot x}}\dfrac{\left(  R_{i}f\right)
\left(  x\right)  }{\left(  R_{i}f\right)  \left(  y\right)  }=\sum
_{\substack{x\in\widehat{P}_{j};\ y\in\widehat{P}_{j+1};\\y \gtrdot x}}\left(
\underbrace{\left(  R_{i}f\right)  \left(  y\right)  }_{\substack{=f\left(
y\right)  \\\text{(by Proposition \ref{prop.Ri.implicit} \textbf{(a)})}%
}}\right)  ^{-1}\cdot\underbrace{\left(  R_{i}f\right)  \left(  x\right)
}_{\substack{=f\left(  x\right)  \\\text{(by Proposition
\ref{prop.Ri.implicit} \textbf{(a)})}}}\\
&  =\sum_{\substack{x\in\widehat{P}_{j};\ y\in\widehat{P}_{j+1};\\y \gtrdot
x}}\left(  f\left(  y\right)  \right)  ^{-1}\cdot f\left(  x\right)
=\sum_{\substack{x\in\widehat{P}_{j};\ y\in\widehat{P}_{j+1};\\y \gtrdot x
}}\dfrac{f\left(  x\right)  }{f\left(  y\right)  }.
\end{align*}
Compared with $\mathbf{w}_{j}\left(  f\right)  =\sum\limits_{\substack{x\in
\widehat{P}_{j};\ y\in\widehat{P}_{j+1};\\y \gtrdot x}}\dfrac{f\left(
x\right)  }{f\left(  y\right)  }$ (by the definition of $\mathbf{w}_{j}\left(
f\right)  $), this yields $\mathbf{w}_{j}\left(  R_{i}f\right)  =\mathbf{w}%
_{j}\left(  f\right)  $. Since $j=\tau_{i}\left(  j\right)  $, this becomes
$\mathbf{w}_{j}\left(  R_{i}f\right)  =\mathbf{w}_{\tau_{i}\left(  j\right)
}\left(  f\right)  $. Hence, (\ref{pf.wi.Ri.1}) is proven in Case 3.

We have thus proven (\ref{pf.wi.Ri.1}) in each of the three possible cases 1,
2 and 3. This completes the proof of (\ref{pf.wi.Ri.1}) and thus of
Proposition \ref{prop.wi.Ri}.
\end{proof}

From Proposition \ref{prop.wi.Ri}, and (\ref{def.Ri.R}), we conclude:

\begin{proposition}
\label{prop.wi.R}Let $n\in\mathbb{N}$. Let $\mathbb{K}$ be a field. Let $P$ be
an $n$-graded poset. Then, every $f\in\mathbb{K}^{\widehat{P}}$ satisfies%
\[
\left(  \mathbf{w}_{0}\left(  Rf\right)  ,\mathbf{w}_{1}\left(  Rf\right)
,...,\mathbf{w}_{n}\left(  Rf\right)  \right)  =\left(  \mathbf{w}_{n}\left(
f\right)  ,\mathbf{w}_{0}\left(  f\right)  ,\mathbf{w}_{1}\left(  f\right)
,...,\mathbf{w}_{n-1}\left(  f\right)  \right)  .
\]

In other words, the map $R$ changes the w-tuple of a $\mathbb{K}$-labelling by
shifting it cyclically.
\end{proposition}

\begin{vershort}
\begin{proof}
[Proof of Proposition \ref{prop.wi.R} (sketched).]Proposition \ref{prop.Ri.R}
yields $R=R_{1}\circ R_{2}\circ...\circ R_{n}$. But for every $i\in\left\{
1,2,...,n\right\}  $, recall from Proposition \ref{prop.wi.Ri} that the map
$R_{i}$ changes the w-tuple of a $\mathbb{K}$-labelling by interchanging its
$\left(  i-1\right)  $-st entry with its $i$-th entry (where the entries are
labelled starting at $0$). Hence, the effect of the compound map $R=R_{1}\circ
R_{2}\circ...\circ R_{n}$ on the w-tuple is that of first interchanging the
$\left(  n-1\right)  $-st entry with the $n$-th entry, then interchanging the
$\left(  n-2\right)  $-st entry with the $\left(  n-1\right)  $-st entry, and
so on, through to finally interchanging the $0$-th entry with the $1$-st
entry. But this latter sequence of interchanges is equivalent to a cyclic
shift of the w-tuple\footnote{Indeed, the composition $\left(  0,1\right)
\circ\left(  1,2\right)  \circ...\circ\left(  n-1,n\right)  $ of
transpositions in the symmetric group on the set $\left\{  0,1,...,n\right\}
$ is the $\left(  n+1\right)  $-cycle $\left(  0,1,...,n\right)  $.}. Hence,
the map $R$ changes the w-tuple of a $\mathbb{K}$-labelling by shifting it
cyclically, qed.
\end{proof}
\end{vershort}

\begin{verlong}
\begin{proof}
[Proof of Proposition \ref{prop.wi.R} (sketched).]Let $\mathfrak{S}$ be the
symmetric group on the set $\left\{  0,1,...,n\right\}  $. For every
$i\in\left\{  1,2,...,n\right\}  $, let $\tau_{i}$ be the permutation of the
set $\left\{  0,1,...,n\right\}  $ which transposes $i-1$ with $i$ (while
leaving all other elements of this set invariant). Then, these $\tau_{i}$ are
elements of $\mathfrak{S}$ and satisfy%
\begin{equation}
\tau_{1}\circ\tau_{2}\circ...\circ\tau_{n}=\left(  0,1,...,n\right)
\label{pf.wi.R.cycle}%
\end{equation}
(where $\left(  0,1,...,n\right)  $ denotes the cyclic permutation of the set
$\left\{  0,1,...,n\right\}  $ which sends every $i$ to the remainder of $i+1$
modulo $n+1$). (Indeed, this is just a restatement of the known fact that
$s_{1}\circ s_{2}\circ...\circ s_{n}=\left(  1,2,...,n+1\right)  $ in the
symmetric group $S_{n+1}$.)

Now, let the symmetric group $\mathfrak{S}$ act on the $\left(  n+1\right)
$-tuples of elements of $\mathbb{K}$ in the obvious way:%
\[
\left.
\begin{array}
[c]{l}%
\sigma\left(  \alpha_{0},\alpha_{1},...,\alpha_{n}\right)  =\left(
\alpha_{\sigma^{-1}\left(  0\right)  },\alpha_{\sigma^{-1}\left(  1\right)
},...,\alpha_{\sigma^{-1}\left(  n\right)  }\right) \\
\ \ \ \ \ \ \ \ \ \ \text{for every }\sigma\in\mathfrak{S}\text{ and }\left(
\alpha_{0},\alpha_{1},...,\alpha_{n}\right)  \in\mathbb{K}^{n+1}%
\end{array}
\right.  .
\]

Define a rational map $\mathbf{w}:\mathbb{K}^{\widehat{P}}\dashrightarrow
\mathbb{K}^{n+1}$ by setting%
\begin{equation}
\mathbf{w}\left(  f\right)  =\left(  \mathbf{w}_{0}\left(  f\right)
,\mathbf{w}_{1}\left(  f\right)  ,...,\mathbf{w}_{n}\left(  f\right)  \right)
\ \ \ \ \ \ \ \ \ \ \text{for every }f\in\mathbb{K}^{\widehat{P}}.
\label{pf.wi.R.w}%
\end{equation}
(So $\mathbf{w}$ is the map sending every $\mathbb{K}$-labelling of $P$ to its
w-tuple.) Then, from Proposition \ref{prop.wi.Ri}, we obtain%
\begin{equation}
\mathbf{w}\circ R_{i}=\tau_{i}\circ\mathbf{w}\ \ \ \ \ \ \ \ \ \ \text{for
every }i\in\left\{  1,2,...,n\right\}  \label{pf.wi.R.step}%
\end{equation}
\footnote{\textit{Proof of (\ref{pf.wi.R.step}):} Let $i\in\left\{
1,2,...,n\right\}  $. Every $f\in\mathbb{K}^{\widehat{P}}$ satisfies%
\begin{align*}
\left(  \mathbf{w}\circ R_{i}\right)  \left(  f\right)   &  =\mathbf{w}\left(
R_{i}\left(  f\right)  \right)  =\left(  \mathbf{w}_{0}\left(  R_{i}f\right)
,\mathbf{w}_{1}\left(  R_{i}f\right)  ,...,\mathbf{w}_{n}\left(
R_{i}f\right)  \right)  \ \ \ \ \ \ \ \ \ \ \left(  \text{by the definition of
}\mathbf{w}\right) \\
&  =\left(  \mathbf{w}_{0}\left(  f\right)  ,\mathbf{w}_{1}\left(  f\right)
,...,\mathbf{w}_{i-2}\left(  f\right)  ,\mathbf{w}_{i}\left(  f\right)
,\mathbf{w}_{i-1}\left(  f\right)  ,\mathbf{w}_{i+1}\left(  f\right)
,\mathbf{w}_{i+2}\left(  f\right)  ,...,\mathbf{w}_{n}\left(  f\right)
\right) \\
&  =\tau_{i}\underbrace{\left(  \mathbf{w}_{0}\left(  f\right)  ,\mathbf{w}%
_{1}\left(  f\right)  ,...,\mathbf{w}_{n}\left(  f\right)  \right)
}_{\substack{=\mathbf{w}\left(  f\right)  \\\text{(by (\ref{pf.wi.R.w}))}%
}}=\tau_{i}\mathbf{w}\left(  f\right)  =\left(  \tau_{i}\circ\mathbf{w}%
\right)  \left(  f\right)  .
\end{align*}
Thus, $\mathbf{w}\circ R_{i}=\tau_{i}\circ\mathbf{w}$, qed.}. But Proposition
\ref{prop.Ri.R} yields $R=R_{1}\circ R_{2}\circ...\circ R_{n}$, so that%
\begin{align}
\mathbf{w}\circ R  &  =\mathbf{w}\circ R_{1}\circ R_{2}\circ...\circ
R_{n}\nonumber\\
&  =\underbrace{\mathbf{w}\circ R_{1}}_{\substack{=\tau_{1}\circ
\mathbf{w}\\\text{(by (\ref{pf.wi.R.step}), applied to }i=1\text{)}}}\circ
R_{2}\circ R_{3}\circ R_{4}\circ...\circ R_{n}\nonumber\\
&  =\tau_{1}\circ\underbrace{\mathbf{w}\circ R_{2}}_{\substack{=\tau_{2}%
\circ\mathbf{w}\\\text{(by (\ref{pf.wi.R.step}), applied to }i=2\text{)}%
}}\circ R_{3}\circ R_{4}\circ...\circ R_{n}\nonumber\\
&  =\tau_{1}\circ\tau_{2}\circ\mathbf{w}\circ R_{3}\circ R_{4}\circ...\circ
R_{n}\nonumber\\
&  =...\nonumber\\
&  =\tau_{1}\circ\tau_{2}\circ...\circ\tau_{n-1}\circ\underbrace{\mathbf{w}%
\circ R_{n}}_{\substack{=\tau_{n}\circ\mathbf{w}\\\text{(by
(\ref{pf.wi.R.step}), applied to }i=n\text{)}}}\nonumber\\
&  =\underbrace{\left(  \tau_{1}\circ\tau_{2}\circ...\circ\tau_{n}\right)
}_{\substack{=\left(  0,1,...,n\right)  \\\text{(by (\ref{pf.wi.R.cycle}))}%
}}\circ\mathbf{w}=\left(  0,1,...,n\right)  \circ\mathbf{w}.
\label{pf.wi.R.almost}%
\end{align}

Now, every $f\in\mathbb{K}^{\widehat{P}}$ satisfies $\mathbf{w}\left(
Rf\right)  =\left(  \mathbf{w}_{0}\left(  Rf\right)  ,\mathbf{w}_{1}\left(
Rf\right)  ,...,\mathbf{w}_{n}\left(  Rf\right)  \right)  $ (by the definition
of $\mathbf{w}\left(  Rf\right)  $), thus%
\begin{align*}
\left(  \mathbf{w}_{0}\left(  Rf\right)  ,\mathbf{w}_{1}\left(  Rf\right)
,...,\mathbf{w}_{n}\left(  Rf\right)  \right)   &  =\mathbf{w}\left(
Rf\right)  =\underbrace{\left(  \mathbf{w}\circ R\right)  }%
_{\substack{=\left(  0,1,...,n\right)  \circ\mathbf{w}\\\text{(by
(\ref{pf.wi.R.almost}))}}}f=\left(  \left(  0,1,...,n\right)  \circ
\mathbf{w}\right)  f\\
&  =\left(  0,1,...,n\right)  \underbrace{\left(  \mathbf{w}\left(  f\right)
\right)  }_{=\left(  \mathbf{w}_{0}\left(  f\right)  ,\mathbf{w}_{1}\left(
f\right)  ,...,\mathbf{w}_{n}\left(  f\right)  \right)  }\\
&  =\left(  0,1,...,n\right)  \left(  \mathbf{w}_{0}\left(  f\right)
,\mathbf{w}_{1}\left(  f\right)  ,...,\mathbf{w}_{n}\left(  f\right)  \right)
\\
&  =\left(  \mathbf{w}_{n}\left(  f\right)  ,\mathbf{w}_{0}\left(  f\right)
,\mathbf{w}_{1}\left(  f\right)  ,...,\mathbf{w}_{n-1}\left(  f\right)
\right)  .
\end{align*}
This proves Proposition \ref{prop.wi.R}.
\end{proof}
\end{verlong}

\begin{vershort}
As a consequence of Proposition \ref{prop.wi.R}, the map $R^{n+1}$ (for an
$n$-graded poset $P$) leaves the w-tuple of a $\mathbb{K}$-labelling fixed.
\end{vershort}

\begin{verlong}
An easy consequence of Proposition \ref{prop.wi.R} is the following fact:

\begin{corollary}
\label{cor.wi.Rn+1}Let $n\in\mathbb{N}$. Let $\mathbb{K}$ be a field. Let $P$
be an $n$-graded poset. Then, every $f\in\mathbb{K}^{\widehat{P}}$ satisfies%
\begin{equation}
\left(  \mathbf{w}_{0}\left(  R^{n+1}f\right)  ,\mathbf{w}_{1}\left(
R^{n+1}f\right)  ,...,\mathbf{w}_{n}\left(  R^{n+1}f\right)  \right)  =\left(
\mathbf{w}_{0}\left(  f\right)  ,\mathbf{w}_{1}\left(  f\right)
,...,\mathbf{w}_{n}\left(  f\right)  \right)  . \label{cor.wi.Rn+1.eq}%
\end{equation}

In other words, the map $R^{n+1}$ leaves the w-tuple of a $\mathbb{K}%
$-labelling fixed.
\end{corollary}

Roughly speaking, Corollary \ref{cor.wi.Rn+1} is proven by $n+1$-fold
application of Proposition \ref{prop.wi.R}. Here is a more formal proof:

\begin{proof}
[Proof of Corollary \ref{cor.wi.Rn+1} (sketched).]Consider the group
$\mathfrak{S}$ defined in the proof of Proposition \ref{prop.wi.R}, and its
action on the $\left(  n+1\right)  $-tuples of elements of $\mathbb{K}$
defined ibidem. Also, define the element $\left(  0,1,...,n\right)  $ of
$\mathfrak{S}$ just as in the proof of Proposition \ref{prop.wi.R}.

From Proposition \ref{prop.wi.R}, we know that the map $R$ changes the w-tuple
of a $\mathbb{K}$-labelling by shifting it cyclically. In other words, almost
every $\mathbb{K}$-labelling $f$ of $P$ satisfies%
\[
\left(  \text{the w-tuple of }Rf\right)  =\left(  0,1,...,n\right)  \left(
\text{the w-tuple of }f\right)  .
\]
From this, it is easy to see that for every $k\in\mathbb{N}$, almost every
$\mathbb{K}$-labelling $f$ of $P$ satisfies%
\[
\left(  \text{the w-tuple of }R^{k}f\right)  =\left(  0,1,...,n\right)
^{k}\left(  \text{the w-tuple of }f\right)
\]
(indeed, this can be proven by induction over $k$). Applying this to $k=n+1$,
we obtain%
\begin{align*}
\left(  \text{the w-tuple of }R^{n+1}f\right)   &  =\underbrace{\left(
0,1,...,n\right)  ^{n+1}}_{\substack{=\operatorname*{id}\\\text{(since the
permutation }\left(  0,1,...,n\right)  \in\mathfrak{S}\\\text{has order
}n+1\text{)}}}\left(  \text{the w-tuple of }f\right) \\
&  =\operatorname*{id}\left(  \text{the w-tuple of }f\right)  =\left(
\text{the w-tuple of }f\right) \\
&  =\left(  \mathbf{w}_{0}\left(  f\right)  ,\mathbf{w}_{1}\left(  f\right)
,...,\mathbf{w}_{n}\left(  f\right)  \right)
\end{align*}
(by the definition of the w-tuple). Compared with%
\[
\left(  \text{the w-tuple of }R^{n+1}f\right)  =\left(  \mathbf{w}_{0}\left(
R^{n+1}f\right)  ,\mathbf{w}_{1}\left(  R^{n+1}f\right)  ,...,\mathbf{w}%
_{n}\left(  R^{n+1}f\right)  \right)
\]
(by the definition of the w-tuple), this yields%
\[
\left(  \mathbf{w}_{0}\left(  R^{n+1}f\right)  ,\mathbf{w}_{1}\left(
R^{n+1}f\right)  ,...,\mathbf{w}_{n}\left(  R^{n+1}f\right)  \right)  =\left(
\mathbf{w}_{0}\left(  f\right)  ,\mathbf{w}_{1}\left(  f\right)
,...,\mathbf{w}_{n}\left(  f\right)  \right)  .
\]
This proves Corollary \ref{cor.wi.Rn+1}.
\end{proof}
\end{verlong}

\section{Graded rescaling of labellings}

In general, birational rowmotion $R$ has something that one might call an
\textquotedblleft avalanche effect\textquotedblright: If $f$ and $g$ are two
$\mathbb{K}$-labellings of a poset $P$ which differ from each other only in
their labels at one single element $v$, then the labellings $Rf$ and $Rg$ (in
general) differ at all elements covering $v$ and all elements beneath $v$, and
further applications of $R$ make the labellings even more different. Thus, a
change of just one label in a labelling will often \textquotedblleft
spread\textquotedblright\ through a large part of the poset when $R$ is
repeatedly applied; the effect of such a change is hard to track in general.
Thus, knowing the behavior of one particular $\mathbb{K}$-labelling $f$ under
$R$ does not help us at understanding the behaviors of $\mathbb{K}$-labellings
obtained from $f$ by changing labels at particular elements. However, if $P$
is a \textbf{graded} poset and we simultaneously multiply the labels at
\textbf{all elements of a given degree} in a given labelling of $P$ with a
given scalar, then the changes this causes to the behavior of the labelling
under $R$ are rather predictable. We are going to formalize this observation
in this section, proving some explicit formulas for how birational rowmotion
$R$ and its iterates react to such rescalings. These explicit formulas will be
subsumed into slick conclusions in Section \ref{sect.homogeneous}, where we
will introduce a notion of \emph{homogeneous equivalence} which formalizes the
idea of a \textquotedblleft labelling modulo scalar factors at each
degree\textquotedblright.

\begin{definition}
Let $\mathbb{K}$ be a field. Then, $\mathbb{K}^{\times}$ denotes the
multiplicative group of nonzero elements of $\mathbb{K}$.
\end{definition}

The following definition formalizes the idea of multiplying the labels at all
elements of a certain degree with one and the same scalar factor:

\begin{definition}
\label{def.bemol}Let $n\in\mathbb{N}$. Let $\mathbb{K}$ be a field. Let $P$ be
an $n$-graded poset. For every $\mathbb{K}$-labelling $f\in\mathbb{K}%
^{\widehat{P}}$ and any $\left(  n+2\right)  $-tuple $\left(  a_{0}%
,a_{1},...,a_{n+1}\right)  \in\left(  \mathbb{K}^{\times}\right)  ^{n+2}$, we
define a $\mathbb{K}$-labelling $\left(  a_{0},a_{1},...,a_{n+1}\right)  \flat
f\in\mathbb{K}^{\widehat{P}}$ by%
\[
\left(  \left(  a_{0},a_{1},...,a_{n+1}\right)  \flat f\right)  \left(
v\right)  =a_{\deg v}\cdot f\left(  v\right)  \ \ \ \ \ \ \ \ \ \ \text{for
every }v\in\widehat{P}.
\]

\end{definition}

\begin{vershort}
Straightforward application of this definition and that of $R_{i}$ shows:
\end{vershort}

\begin{proposition}
\label{prop.Ri.scalmult}Let $n\in\mathbb{N}$. Let $\mathbb{K}$ be a field. Let
$P$ be an $n$-graded poset. Let us use the notation introduced in Definition
\ref{def.bemol}.

Let $f\in\mathbb{K}^{\widehat{P}}$ be a $\mathbb{K}$-labelling. Let $\left(
a_{0},a_{1},...,a_{n+1}\right)  \in\left(  \mathbb{K}^{\times}\right)  ^{n+2}%
$. Let $i\in\left\{  1,2,...,n\right\}  $. Then,%
\begin{align*}
&  R_{i}\left(  \left(  a_{0},a_{1},...,a_{n+1}\right)  \flat f\right) \\
&  =\left(  a_{0},a_{1},...,a_{i-1},\dfrac{a_{i+1}a_{i-1}}{a_{i}}%
,a_{i+1},a_{i+2},...,a_{n+1}\right)  \flat\left(  R_{i}f\right)
\end{align*}
(provided that $R_{i}f$ is well-defined).
\end{proposition}

\begin{verlong}
\begin{proof}
[Proof of Proposition \ref{prop.Ri.scalmult} (sketched).]Let $v\in\widehat{P}%
$. We need to show that
\begin{align}
&  \left(  R_{i}\left(  \left(  a_{0},a_{1},...,a_{n+1}\right)  \flat
f\right)  \right)  \left(  v\right) \nonumber\\
&  =\left(  \left(  a_{0},a_{1},...,a_{i-1},\dfrac{a_{i+1}a_{i-1}}{a_{i}%
},a_{i+1},a_{i+2},...,a_{n+1}\right)  \flat\left(  R_{i}f\right)  \right)
\left(  v\right)  . \label{pf.Ri.scalmult.1}%
\end{align}
Once (\ref{pf.Ri.scalmult.1}) is proven, the proof of Proposition
\ref{prop.Ri.scalmult} will clearly be complete.

We distinguish between two cases:

\textit{Case 1:} We have $\deg v\neq i$.

\textit{Case 2:} We have $\deg v=i$.

Let us first consider Case 1. In this case, we have $\deg v\neq i$. Hence,%
\begin{align}
&  \left(  R_{i}\left(  \left(  a_{0},a_{1},...,a_{n+1}\right)  \flat
f\right)  \right)  \left(  v\right) \nonumber\\
&  =\left(  \left(  a_{0},a_{1},...,a_{n+1}\right)  \flat f\right)  \left(
v\right)  \ \ \ \ \ \ \ \ \ \ \left(  \text{by Proposition
\ref{prop.Ri.implicit} \textbf{(a)}}\right) \nonumber\\
&  =a_{\deg v}\cdot f\left(  v\right)  \ \ \ \ \ \ \ \ \ \ \left(  \text{by
the definition of }\left(  a_{0},a_{1},...,a_{n+1}\right)  \flat f\right)
\label{pf.Ri.scalmult.c1.1}%
\end{align}
and%
\begin{equation}
\left(  \left(  a_{0},a_{1},...,a_{i-1},\dfrac{a_{i+1}a_{i-1}}{a_{i}}%
,a_{i+1},a_{i+2},...,a_{n+1}\right)  \flat\left(  R_{i}f\right)  \right)
\left(  v\right)  =a_{\deg v}\cdot f\left(  v\right)
\label{pf.Ri.scalmult.c1.2}%
\end{equation}
(by the definition of $\left(  a_{0},a_{1},...,a_{i-1},\dfrac{a_{i+1}a_{i-1}%
}{a_{i}},a_{i+1},a_{i+2},...,a_{n+1}\right)  \flat\left(  R_{i}f\right)  $,
since $\deg v\neq i$). Comparing the two equalities (\ref{pf.Ri.scalmult.c1.1}%
) and (\ref{pf.Ri.scalmult.c1.2}), we obtain%
\begin{align*}
&  \left(  R_{i}\left(  \left(  a_{0},a_{1},...,a_{n+1}\right)  \flat
f\right)  \right)  \left(  v\right) \\
&  =\left(  \left(  a_{0},a_{1},...,a_{i-1},\dfrac{a_{i+1}a_{i-1}}{a_{i}%
},a_{i+1},a_{i+2},...,a_{n+1}\right)  \flat\left(  R_{i}f\right)  \right)
\left(  v\right)  .
\end{align*}
Thus, (\ref{pf.Ri.scalmult.1}) is proven in Case 1.

We can now proceed to Case 2. In this case, we have $\deg v=i$. Hence,%
\begin{align}
&  \left(  R_{i}\left(  \left(  a_{0},a_{1},...,a_{n+1}\right)  \flat
f\right)  \right)  \left(  v\right) \nonumber\\
&  =\dfrac{1}{\left(  \left(  a_{0},a_{1},...,a_{n+1}\right)  \flat f\right)
\left(  v\right)  }\cdot\dfrac{\sum\limits_{\substack{u\in\widehat{P}%
;\\u\lessdot v}}\left(  \left(  a_{0},a_{1},...,a_{n+1}\right)  \flat
f\right)  \left(  u\right)  }{\sum\limits_{\substack{u\in\widehat{P}%
;\\u\gtrdot v}}\dfrac{1}{\left(  \left(  a_{0},a_{1},...,a_{n+1}\right)  \flat
f\right)  \left(  u\right)  }}\label{pf.Ri.scalmult.c2.1}\\
&  \ \ \ \ \ \ \ \ \ \ \left(  \text{by Proposition \ref{prop.Ri.implicit}
\textbf{(b)}}\right)  .\nonumber
\end{align}

Now, we notice the following:

\begin{itemize}
\item For any element $u\in\widehat{P}$ such that $u\lessdot v$, we have
\begin{align*}
&  \left(  \left(  a_{0},a_{1},...,a_{n+1}\right)  \flat f\right)  \left(
u\right) \\
&  =a_{\deg u}\cdot f\left(  u\right)  \ \ \ \ \ \ \ \ \ \ \left(  \text{by
the definition of }\left(  a_{0},a_{1},...,a_{n+1}\right)  \flat f\right) \\
&  =a_{i-1}\cdot f\left(  u\right)  \ \ \ \ \ \ \ \ \ \ \left(
\begin{array}
[c]{c}%
\text{since }u\lessdot v\text{, so that}\\
\deg u=\underbrace{\deg v}_{=i}-1=i-1\text{, so that }a_{\deg u}=a_{i-1}%
\end{array}
\right)  .
\end{align*}
Hence,%
\begin{equation}
\sum\limits_{\substack{u\in\widehat{P};\\u\lessdot v}}\underbrace{\left(
\left(  a_{0},a_{1},...,a_{n+1}\right)  \flat f\right)  \left(  u\right)
}_{=a_{i-1}\cdot f\left(  u\right)  }=\sum\limits_{\substack{u\in
\widehat{P};\\u\lessdot v}}a_{i-1}\cdot f\left(  u\right)  =a_{i-1}\cdot
\sum\limits_{\substack{u\in\widehat{P};\\u\lessdot v}}f\left(  u\right)  .
\label{pf.Ri.scalmult.c2.2}%
\end{equation}

\item Similarly to (\ref{pf.Ri.scalmult.c2.2}), we can show that%
\begin{equation}
\sum\limits_{\substack{u\in\widehat{P};\\u \gtrdot v}}\dfrac{1}{\left(
\left(  a_{0},a_{1},...,a_{n+1}\right)  \flat f\right)  \left(  u\right)
}=\dfrac{1}{a_{i+1}}\cdot\sum\limits_{\substack{u\in\widehat{P};\\u \gtrdot
v}}\dfrac{1}{f\left(  u\right)  }. \label{pf.Ri.scalmult.c2.3}%
\end{equation}

\item We have%
\begin{align}
&  \left(  \left(  a_{0},a_{1},...,a_{n+1}\right)  \flat f\right)  \left(
v\right) \nonumber\\
&  =a_{\deg v}\cdot f\left(  v\right)  \ \ \ \ \ \ \ \ \ \ \left(  \text{by
the definition of }\left(  a_{0},a_{1},...,a_{n+1}\right)  \flat f\right)
\nonumber\\
&  =a_{i}\cdot f\left(  v\right)  \ \ \ \ \ \ \ \ \ \ \left(  \text{since
}\deg v=i\right)  . \label{pf.Ri.scalmult.c2.4}%
\end{align}

\end{itemize}

Now, substituting (\ref{pf.Ri.scalmult.c2.2}), (\ref{pf.Ri.scalmult.c2.3}) and
(\ref{pf.Ri.scalmult.c2.4}) into (\ref{pf.Ri.scalmult.c2.1}), we obtain%
\begin{align}
&  \left(  R_{i}\left(  \left(  a_{0},a_{1},...,a_{n+1}\right)  \flat
f\right)  \right)  \left(  v\right) \nonumber\\
&  =\dfrac{1}{a_{i}\cdot f\left(  v\right)  }\cdot\dfrac{a_{i-1}\cdot
\sum\limits_{\substack{u\in\widehat{P};\\u \lessdot v}}f\left(  u\right)
}{\dfrac{1}{a_{i+1}}\cdot\sum\limits_{\substack{u\in\widehat{P};\\u \gtrdot
v}}\dfrac{1}{f\left(  u\right)  }}=\dfrac{a_{i-1}a_{i+1}}{a_{i}}%
\cdot\underbrace{\dfrac{1}{f\left(  v\right)  }\cdot\dfrac{\sum
\limits_{\substack{u\in\widehat{P};\\u \lessdot v}}f\left(  u\right)  }%
{\sum\limits_{\substack{u\in\widehat{P};\\u \gtrdot v}}\dfrac{1}{f\left(
u\right)  }}}_{\substack{=\left(  R_{i}f\right)  \left(  v\right)  \\\text{(by
Proposition \ref{prop.Ri.implicit} \textbf{(b)})}}}\nonumber\\
&  =\dfrac{a_{i-1}a_{i+1}}{a_{i}}\cdot\left(  R_{i}f\right)  \left(  v\right)
. \label{pf.Ri.scalmult.c2.8}%
\end{align}

But by the definition of $\left(  a_{0},a_{1},...,a_{i-1},\dfrac
{a_{i+1}a_{i-1}}{a_{i}},a_{i+1},a_{i+2},...,a_{n+1}\right)  \flat\left(
R_{i}f\right)  $ and because of $\deg v=i$, we have%
\begin{align*}
&  \left(  \left(  a_{0},a_{1},...,a_{i-1},\dfrac{a_{i+1}a_{i-1}}{a_{i}%
},a_{i+1},a_{i+2},...,a_{n+1}\right)  \flat\left(  R_{i}f\right)  \right)
\left(  v\right) \\
&  =\dfrac{a_{i+1}a_{i-1}}{a_{i}}\cdot\left(  R_{i}f\right)  \left(  v\right)
.
\end{align*}
Comparing this with (\ref{pf.Ri.scalmult.c2.8}), we obtain%
\begin{align*}
&  \left(  R_{i}\left(  \left(  a_{0},a_{1},...,a_{n+1}\right)  \flat
f\right)  \right)  \left(  v\right) \\
&  =\left(  \left(  a_{0},a_{1},...,a_{i-1},\dfrac{a_{i+1}a_{i-1}}{a_{i}%
},a_{i+1},a_{i+2},...,a_{n+1}\right)  \flat\left(  R_{i}f\right)  \right)
\left(  v\right)  .
\end{align*}
Thus, (\ref{pf.Ri.scalmult.1}) is proven in Case 2.

We have now proven (\ref{pf.Ri.scalmult.1}) in both possible cases 1 and 2.
This completes the proof of (\ref{pf.Ri.scalmult.1}) and thus of Proposition
\ref{prop.Ri.scalmult}.
\end{proof}
\end{verlong}

A similar result can be obtained for $R$ instead of $R_{i}$:

\begin{proposition}
\label{prop.R.scalmult}Let $n\in\mathbb{N}$. Let $\mathbb{K}$ be a field. Let
$P$ be an $n$-graded poset. For every $\mathbb{K}$-labelling $f\in
\mathbb{K}^{\widehat{P}}$ and any $\left(  n+2\right)  $-tuple $\left(
a_{0},a_{1},...,a_{n+1}\right)  \in\left(  \mathbb{K}^{\times}\right)  ^{n+2}%
$, we define a $\mathbb{K}$-labelling $\left(  a_{0},a_{1},...,a_{n+1}\right)
\flat f\in\mathbb{K}^{\widehat{P}}$ as in Definition \ref{def.bemol}.

Let $f\in\mathbb{K}^{\widehat{P}}$ be a $\mathbb{K}$-labelling. Let $\left(
a_{0},a_{1},...,a_{n+1}\right)  \in\left(  \mathbb{K}^{\times}\right)  ^{n+2}%
$. Then,%
\[
R\left(  \left(  a_{0},a_{1},...,a_{n+1}\right)  \flat f\right)  =\left(
a_{0},ga_{0},ga_{1},...,ga_{n-1},a_{n+1}\right)  \flat\left(  Rf\right)  ,
\]
where $g=\dfrac{a_{n+1}}{a_{n}}$ (provided that $Rf$ is well-defined).
\end{proposition}

\begin{proof}
[Proof of Proposition \ref{prop.R.scalmult} (sketched).]Let $g=\dfrac{a_{n+1}%
}{a_{n}}$. We claim that every $j\in\left\{  1,2,...,n+1\right\}  $ satisfies%
\begin{align}
&  \left(  R_{j}\circ R_{j+1}\circ...\circ R_{n}\right)  \left(  \left(
a_{0},a_{1},...,a_{n+1}\right)  \flat f\right) \nonumber\\
&  =\left(  a_{0},a_{1},a_{2},...,a_{j-1},ga_{j-1},ga_{j},...,ga_{n-1}%
,a_{n+1}\right)  \flat\left(  \left(  R_{j}\circ R_{j+1}\circ...\circ
R_{n}\right)  f\right)  . \label{pf.R.scalmult.j}%
\end{align}

\begin{vershort}
Indeed, (\ref{pf.R.scalmult.j}) is easily verified by reverse induction over
$j$ (that is, induction over $n+1-j$), using Proposition
\ref{prop.Ri.scalmult} in the step. Now, applying (\ref{pf.R.scalmult.j}) to
$j=1$ and recalling that $R = R_{1} \circ R_{2} \circ... \circ R_{n}$, we
obtain the claim of Proposition \ref{prop.R.scalmult}.
\end{vershort}

\begin{verlong}
\textit{Proof of (\ref{pf.R.scalmult.j}):} We will prove
(\ref{pf.R.scalmult.j}) by reverse induction over $j$ (that is, induction over
$n+1-j$). The induction base case (the case $j=n+1$) is obvious. For the
induction step, let $i\in\left\{  1,2,...,n+1\right\}  $ be such that $i<n+1$.
Assume that (\ref{pf.R.scalmult.j}) holds for $j=i+1$, and let us show that
(\ref{pf.R.scalmult.j}) holds for $j=i$. We have\footnote{In this computation,
we will be tacitly using the fact that the $\left(  i+2\right)  $-nd entry of
the $\left(  n+2\right)  $-tuple $\left(  a_{0},ga_{0},ga_{1},...,ga_{n-1}%
,a_{n+1}\right)  $ is $ga_{i}$. This is clear if $i<n$; in the case when
$i=n$, it follows from the fact that the $\left(  n+2\right)  $-nd entry of
the $\left(  n+2\right)  $-tuple $\left(  a_{0},ga_{0},ga_{1},...,ga_{n-1}%
,a_{n+1}\right)  $ is $a_{n+1}=\underbrace{\dfrac{a_{n+1}}{a_{n}}}_{=g}%
a_{n}=ga_{n}$.}%
\begin{align*}
&  \left(  R_{i}\circ R_{i+1}\circ...\circ R_{n}\right)  \left(  \left(
a_{0},a_{1},...,a_{n+1}\right)  \flat f\right) \\
&  =R_{i}\left(  \underbrace{\left(  R_{i+1}\circ R_{i+2}\circ...\circ
R_{n}\right)  \left(  \left(  a_{0},a_{1},...,a_{n+1}\right)  \flat f\right)
}_{\substack{=\left(  a_{0},a_{1},a_{2},...,a_{i},ga_{i},ga_{i+1}%
,...,ga_{n-1},a_{n+1}\right)  \flat\left(  \left(  R_{i+1}\circ R_{i+2}%
\circ...\circ R_{n}\right)  f\right)  \\\text{(by (\ref{pf.R.scalmult.j}),
applied to }j=i+1\\\text{(since we assumed that (\ref{pf.R.scalmult.j}) holds
for }j=i+1\text{))}}}\right) \\
&  =R_{i}\left(  \left(  a_{0},a_{1},a_{2},...,a_{i},ga_{i},ga_{i+1}%
,...,ga_{n-1},a_{n+1}\right)  \flat\left(  \left(  R_{i+1}\circ R_{i+2}%
\circ...\circ R_{n}\right)  f\right)  \right) \\
&  =\left(  a_{0},a_{1},...,a_{i-1},\underbrace{\dfrac{ga_{i}a_{i-1}}{a_{i}}%
}_{=ga_{i-1}},ga_{i},ga_{i+1},...,ga_{n-1},a_{n+1}\right)  \flat
\underbrace{\left(  R_{i}\left(  \left(  R_{i+1}\circ R_{i+2}\circ...\circ
R_{n}\right)  f\right)  \right)  }_{\substack{=\left(  R_{i}\circ\left(
R_{i+1}\circ R_{i+2}\circ...\circ R_{n}\right)  \right)  f\\=\left(
R_{i}\circ R_{i+1}\circ...\circ R_{n}\right)  f}}\\
&  \ \ \ \ \ \ \ \ \ \ \left(
\begin{array}
[c]{c}%
\text{by Proposition \ref{prop.Ri.scalmult}, applied to }\left(  R_{i+1}\circ
R_{i+2}\circ...\circ R_{n}\right)  f\text{ and }\\
\left(  a_{0},a_{1},a_{2},...,a_{i},ga_{i},ga_{i+1},...,ga_{n-1}%
,a_{n+1}\right)  \text{ instead of }f\text{ and }\left(  a_{0},a_{1}%
,...,a_{n+1}\right)
\end{array}
\right) \\
&  =\left(  a_{0},a_{1},...,a_{i-1},ga_{i-1},ga_{i},ga_{i+1},...,ga_{n-1}%
,a_{n+1}\right)  \flat\left(  \left(  R_{i}\circ R_{i+1}\circ...\circ
R_{n}\right)  f\right) \\
&  =\left(  a_{0},a_{1},...,a_{i-1},ga_{i-1},ga_{i},...,ga_{n-1}%
,a_{n+1}\right)  \flat\left(  \left(  R_{i}\circ R_{i+1}\circ...\circ
R_{n}\right)  f\right)  .
\end{align*}
This proves (\ref{pf.R.scalmult.j}) for $j=i$, thus completing the induction.

Now that (\ref{pf.R.scalmult.j}) is proven, we can recall that $R=R_{1}\circ
R_{2}\circ...\circ R_{n}$ (by Proposition \ref{prop.Ri.R}), so that%
\begin{align*}
&  R\left(  \left(  a_{0},a_{1},...,a_{n+1}\right)  \flat f\right) \\
&  =\left(  R_{1}\circ R_{2}\circ...\circ R_{n}\right)  \left(  \left(
a_{0},a_{1},...,a_{n+1}\right)  \flat f\right) \\
&  =\left(  a_{0},a_{1},a_{2},...,a_{1-1},ga_{1-1},ga_{1},...,ga_{n-1}%
,a_{n+1}\right)  \flat\left(  \left(  R_{1}\circ R_{2}\circ...\circ
R_{n}\right)  f\right) \\
&  \ \ \ \ \ \ \ \ \ \ \left(  \text{by (\ref{pf.R.scalmult.j}), applied to
}j=1\right) \\
&  =\left(  a_{0},ga_{0},ga_{1},...,ga_{n-1},a_{n+1}\right)  \flat\left(
\underbrace{\left(  R_{1}\circ R_{2}\circ...\circ R_{n}\right)  }_{=R}f\right)
\\
&  =\left(  a_{0},ga_{0},ga_{1},...,ga_{n-1},a_{n+1}\right)  \flat\left(
Rf\right)  .
\end{align*}
This proves Proposition \ref{prop.R.scalmult}.
\end{verlong}
\end{proof}

We can go further and generalize Proposition \ref{prop.R.scalmult} to iterated
birational rowmotion:

\begin{proposition}
\label{prop.Rl.scalmult}Let $n\in\mathbb{N}$. Let $\mathbb{K}$ be a field. Let
$P$ be an $n$-graded poset. For every $\mathbb{K}$-labelling $f\in
\mathbb{K}^{\widehat{P}}$ and any $\left(  n+2\right)  $-tuple $\left(
a_{0},a_{1},...,a_{n+1}\right)  \in\left(  \mathbb{K}^{\times}\right)  ^{n+2}%
$, we define a $\mathbb{K}$-labelling $\left(  a_{0},a_{1},...,a_{n+1}\right)
\flat f\in\mathbb{K}^{\widehat{P}}$ as in Definition \ref{def.bemol}.

Let $\left(  a_{0},a_{1},...,a_{n+1}\right)  \in\left(  \mathbb{K}^{\times
}\right)  ^{n+2}$. For every $\ell\in\left\{  0,1,...,n+1\right\}  $ and
$k\in\left\{  0,1,...,n+1\right\}  $, define an element $\widehat{a}%
_{k}^{\left(  \ell\right)  }\in\mathbb{K}^{\times}$ by%
\[
\widehat{a}_{k}^{\left(  \ell\right)  }=\left\{
\begin{array}
[c]{c}%
\dfrac{a_{n+1}a_{k-\ell}}{a_{n+1-\ell}},\ \ \ \ \ \ \ \ \ \ \text{if }%
k\geq\ell;\\
\dfrac{a_{n+1+k-\ell}a_{0}}{a_{n+1-\ell}},\ \ \ \ \ \ \ \ \ \ \text{if }k<\ell
\end{array}
\right.  .
\]

Let $f\in\mathbb{K}^{\widehat{P}}$ be a $\mathbb{K}$-labelling. Then, every
$\ell\in\left\{  0,1,...,n+1\right\}  $ satisfies%
\[
R^{\ell}\left(  \left(  a_{0},a_{1},...,a_{n+1}\right)  \flat f\right)
=\left(  \widehat{a}_{0}^{\left(  \ell\right)  },\widehat{a}_{1}^{\left(
\ell\right)  },...,\widehat{a}_{n+1}^{\left(  \ell\right)  }\right)
\flat\left(  R^{\ell}f\right)
\]
(provided that $R^{\ell}f$ is well-defined).
\end{proposition}

\begin{example}
For this example, let $n=3$, and let $P$ be a $3$-graded poset. Then,%
\begin{align*}
\left(  \widehat{a}_{0}^{\left(  0\right)  },\widehat{a}_{1}^{\left(
0\right)  },\widehat{a}_{2}^{\left(  0\right)  },\widehat{a}_{3}^{\left(
0\right)  },\widehat{a}_{4}^{\left(  0\right)  }\right)    & =\left(
a_{0},a_{1},a_{2},a_{3},a_{4}\right)  ;\\
\left(  \widehat{a}_{0}^{\left(  1\right)  },\widehat{a}_{1}^{\left(
1\right)  },\widehat{a}_{2}^{\left(  1\right)  },\widehat{a}_{3}^{\left(
1\right)  },\widehat{a}_{4}^{\left(  1\right)  }\right)    & =\left(
a_{0},\dfrac{a_{4}a_{0}}{a_{3}},\dfrac{a_{4}a_{1}}{a_{3}},\dfrac{a_{4}a_{2}%
}{a_{3}},a_{4}\right)  ;\\
\left(  \widehat{a}_{0}^{\left(  2\right)  },\widehat{a}_{1}^{\left(
2\right)  },\widehat{a}_{2}^{\left(  2\right)  },\widehat{a}_{3}^{\left(
2\right)  },\widehat{a}_{4}^{\left(  2\right)  }\right)    & =\left(
a_{0},\dfrac{a_{3}a_{0}}{a_{2}},\dfrac{a_{4}a_{0}}{a_{2}},\dfrac{a_{4}a_{1}%
}{a_{2}},a_{4}\right)  ;\\
\left(  \widehat{a}_{0}^{\left(  3\right)  },\widehat{a}_{1}^{\left(
3\right)  },\widehat{a}_{2}^{\left(  3\right)  },\widehat{a}_{3}^{\left(
3\right)  },\widehat{a}_{4}^{\left(  3\right)  }\right)    & =\left(
a_{0},\dfrac{a_{2}a_{0}}{a_{1}},\dfrac{a_{3}a_{0}}{a_{1}},\dfrac{a_{4}a_{0}%
}{a_{1}},a_{4}\right)  ;\\
\left(  \widehat{a}_{0}^{\left(  4\right)  },\widehat{a}_{1}^{\left(
4\right)  },\widehat{a}_{2}^{\left(  4\right)  },\widehat{a}_{3}^{\left(
4\right)  },\widehat{a}_{4}^{\left(  4\right)  }\right)    & =\left(
a_{0},a_{1},a_{2},a_{3},a_{4}\right)  .
\end{align*}
More generally, we always have $\left(  \widehat{a}_{0}^{\left(  0\right)
},\widehat{a}_{1}^{\left(  0\right)  },...,\widehat{a}_{n+1}^{\left(
0\right)  }\right)  =\left(  a_{0},a_{1},...,a_{n+1}\right)  $ and $\left(
\widehat{a}_{0}^{\left(  n+1\right)  },\widehat{a}_{1}^{\left(  n+1\right)
},...,\widehat{a}_{n+1}^{\left(  n+1\right)  }\right)  =\left(  a_{0}%
,a_{1},...,a_{n+1}\right)  $ (as can be verified directly).
\end{example}

\begin{vershort}
\begin{proof}
[Proof of Proposition \ref{prop.Rl.scalmult} (sketched).]This proof is a
completely straightforward induction over $\ell$, with the base case being
trivial and the induction step relying on Proposition \ref{prop.R.scalmult}.
It is useful to notice that every $\ell\in\left\{  0,1,...,n+1\right\}  $ and
$k\in\left\{  0,1,...,n+1\right\}  $ satisfy
\[
\widehat{a}_{k}^{\left(  \ell\right)  }=\dfrac{a_{n+1+k-\ell}a_{0}%
}{a_{n+1-\ell}}\ \ \ \ \ \ \ \ \ \ \text{if }k\leq\ell
\]
to simplify the computations (this identity follows from the definition when
$k<\ell$ and can be easily checked for $k=\ell$).
\end{proof}
\end{vershort}

\begin{verlong}
\begin{proof}
[Proof of Proposition \ref{prop.Rl.scalmult} (sketched).]This proof is
straightforward (and this shouldn't come as a surprise). We will nevertheless
present it.

Let us first notice that every $\ell\in\left\{  0,1,...,n+1\right\}  $ and
$k\in\left\{  0,1,...,n+1\right\}  $ satisfy%
\begin{equation}
\widehat{a}_{k}^{\left(  \ell\right)  }=\dfrac{a_{n+1}a_{k-\ell}}{a_{n+1-\ell
}}\ \ \ \ \ \ \ \ \ \ \text{if }k\geq\ell\label{pf.Rl.scalmult.geq}%
\end{equation}
\footnote{This follows directly from the definition of $\widehat{a}%
_{k}^{\left(  \ell\right)  }$.} and%
\begin{equation}
\widehat{a}_{k}^{\left(  \ell\right)  }=\dfrac{a_{n+1+k-\ell}a_{0}%
}{a_{n+1-\ell}}\ \ \ \ \ \ \ \ \ \ \text{if }k\leq\ell
\label{pf.Rl.scalmult.leq}%
\end{equation}
\footnote{\textit{Proof of (\ref{pf.Rl.scalmult.leq}):} Let $\ell\in\left\{
0,1,...,n+1\right\}  $ and $k\in\left\{  0,1,...,n+1\right\}  $ be such that
$k\leq\ell$. If $k<\ell$, then (\ref{pf.Rl.scalmult.leq}) follows directly
from the definition of $\widehat{a}_{k}^{\left(  \ell\right)  }$. So WLOG
assume that we don't have $k<\ell$. Thus, $k=\ell$ (since $k\leq\ell$). Hence,%
\begin{align*}
\widehat{a}_{k}^{\left(  \ell\right)  }  &  =\widehat{a}_{\ell}^{\left(
\ell\right)  }=\left\{
\begin{array}
[c]{c}%
\dfrac{a_{n+1}a_{\ell-\ell}}{a_{n+1-\ell}},\ \ \ \ \ \ \ \ \ \ \text{if }%
\ell\geq\ell;\\
\dfrac{a_{n+1+\ell-\ell}a_{0}}{a_{n+1-\ell}},\ \ \ \ \ \ \ \ \ \ \text{if
}\ell<\ell
\end{array}
\right.  \ \ \ \ \ \ \ \ \ \ \left(  \text{by the definition of }%
\widehat{a}_{\ell}^{\left(  \ell\right)  }\right) \\
&  =\dfrac{a_{n+1}a_{\ell-\ell}}{a_{n+1-\ell}}\ \ \ \ \ \ \ \ \ \ \left(
\text{since }\ell\geq\ell\right) \\
&  =\dfrac{a_{n+1}a_{0}}{a_{n+1-\ell}}%
\end{align*}
and%
\begin{align*}
\dfrac{a_{n+1+k-\ell}a_{0}}{a_{n+1-\ell}}  &  =\dfrac{a_{n+1+\ell-\ell}a_{0}%
}{a_{n+1-\ell}}\ \ \ \ \ \ \ \ \ \ \left(  \text{since }k=\ell\right) \\
&  =\dfrac{a_{n+1}a_{0}}{a_{n+1-\ell}}.
\end{align*}
Thus, $\widehat{a}_{k}^{\left(  \ell\right)  }=\dfrac{a_{n+1}a_{0}%
}{a_{n+1-\ell}}=\dfrac{a_{n+1+k-\ell}a_{0}}{a_{n+1-\ell}}$. This proves
(\ref{pf.Rl.scalmult.leq}).}. In particular, every $\ell\in\left\{
0,1,...,n\right\}  $ satisfies%
\begin{equation}
\widehat{a}_{n}^{\left(  \ell\right)  }=\dfrac{a_{n+1}a_{n-\ell}}{a_{n+1-\ell
}}\ \ \ \ \ \ \ \ \ \ \left(  \text{by (\ref{pf.Rl.scalmult.geq}), applied to
}k=n\right)  . \label{pf.Rl.scalmult.n}%
\end{equation}
Furthermore, every $\ell\in\left\{  0,1,...,n+1\right\}  $ satisfies%
\begin{align}
\widehat{a}_{n+1}^{\left(  \ell\right)  }  &  =\dfrac{a_{n+1}a_{n+1-\ell}%
}{a_{n+1-\ell}}\ \ \ \ \ \ \ \ \ \ \left(  \text{by (\ref{pf.Rl.scalmult.geq}%
), applied to }k=n+1\right) \nonumber\\
&  =a_{n+1}. \label{pf.Rl.scalmult.n+1}%
\end{align}
Moreover, every $\ell\in\left\{  0,1,...,n+1\right\}  $ satisfies%
\begin{align}
\widehat{a}_{0}^{\left(  \ell\right)  }  &  =\dfrac{a_{n+1+0-\ell}a_{0}%
}{a_{n+1-\ell}}\ \ \ \ \ \ \ \ \ \ \left(  \text{by (\ref{pf.Rl.scalmult.leq}%
), applied to }k=0\right) \nonumber\\
&  =\dfrac{a_{n+1-\ell}a_{0}}{a_{n+1-\ell}}=a_{0}. \label{pf.Rl.scalmult.0}%
\end{align}

We will prove Proposition \ref{prop.Rl.scalmult} by induction over $\ell$. For
the induction base (the case $\ell=0$), it is enough to convince oneself that
$\widehat{a}_{k}^{\left(  0\right)  }=a_{k}$ for every $k\in\left\{
0,1,...,n+1\right\}  $, which is really obvious. Now, let us proceed to the
induction step.

Let $L\in\left\{  1,2,...,n+1\right\}  $. Assume that Proposition
\ref{prop.Rl.scalmult} holds for $\ell=L-1$. We need to prove that Proposition
\ref{prop.Rl.scalmult} holds for $\ell=L$.

We know that Proposition \ref{prop.Rl.scalmult} holds for $\ell=L-1$. In other
words,
\begin{equation}
R^{L-1}\left(  \left(  a_{0},a_{1},...,a_{n+1}\right)  \flat f\right)
=\left(  \widehat{a}_{0}^{\left(  L-1\right)  },\widehat{a}_{1}^{\left(
L-1\right)  },...,\widehat{a}_{n+1}^{\left(  L-1\right)  }\right)
\flat\left(  R^{L-1}f\right)  . \label{pf.Rl.scalmult.1}%
\end{equation}

Let $\widehat{g}=\dfrac{\widehat{a}_{n+1}^{\left(  L-1\right)  }}%
{\widehat{a}_{n}^{\left(  L-1\right)  }}$. Then,%
\begin{align*}
\widehat{g}  &  =\dfrac{\widehat{a}_{n+1}^{\left(  L-1\right)  }}%
{\widehat{a}_{n}^{\left(  L-1\right)  }}=\underbrace{\widehat{a}%
_{n+1}^{\left(  L-1\right)  }}_{\substack{=a_{n+1}\\\text{(by
(\ref{pf.Rl.scalmult.n+1}), applied}\\\text{to }L-1\text{ instead of }%
\ell\text{)}}}\diagup\underbrace{\widehat{a}_{n}^{\left(  L-1\right)  }%
}_{\substack{=\dfrac{a_{n+1}a_{n-\left(  L-1\right)  }}{a_{n+1-\left(
L-1\right)  }}\\\text{(by (\ref{pf.Rl.scalmult.n}), applied}\\\text{to
}L-1\text{ instead of }\ell\text{)}}}=a_{n+1}\diagup\dfrac{a_{n+1}a_{n-\left(
L-1\right)  }}{a_{n+1-\left(  L-1\right)  }}\\
&  =\dfrac{a_{n+1-\left(  L-1\right)  }}{a_{n-\left(  L-1\right)  }}.
\end{align*}

Now, $R^{L}=R\circ R^{L-1}$, so that%
\begin{align}
&  R^{L}\left(  \left(  a_{0},a_{1},...,a_{n+1}\right)  \flat f\right)
\nonumber\\
&  =\left(  R\circ R^{L-1}\right)  \left(  \left(  a_{0},a_{1},...,a_{n+1}%
\right)  \flat f\right)  =R\left(  \underbrace{R^{L-1}\left(  \left(
a_{0},a_{1},...,a_{n+1}\right)  \flat f\right)  }_{\substack{=\left(
\widehat{a}_{0}^{\left(  L-1\right)  },\widehat{a}_{1}^{\left(  L-1\right)
},...,\widehat{a}_{n+1}^{\left(  L-1\right)  }\right)  \flat\left(
R^{L-1}f\right)  \\\text{(by (\ref{pf.Rl.scalmult.1}))}}}\right) \nonumber\\
&  =R\left(  \left(  \widehat{a}_{0}^{\left(  L-1\right)  },\widehat{a}%
_{1}^{\left(  L-1\right)  },...,\widehat{a}_{n+1}^{\left(  L-1\right)
}\right)  \flat\left(  R^{L-1}f\right)  \right) \nonumber\\
&  =\left(  \widehat{a}_{0}^{\left(  L-1\right)  },\widehat{g}\widehat{a}%
_{0}^{\left(  L-1\right)  },\widehat{g}\widehat{a}_{1}^{\left(  L-1\right)
},...,\widehat{g}\widehat{a}_{n-1}^{\left(  L-1\right)  },\widehat{a}%
_{n+1}^{\left(  L-1\right)  }\right)  \flat\left(  \underbrace{R\left(
R^{L-1}f\right)  }_{=\left(  R\circ R^{L-1}\right)  f=R^{L}f}\right)
\nonumber\\
&  \ \ \ \ \ \ \ \ \ \ \left(
\begin{array}
[c]{c}%
\text{by Proposition \ref{prop.R.scalmult}, applied to } R^{L-1}f\text{, }\\
\left(  \widehat{a}_{0}^{\left(  L-1\right)  },\widehat{a}_{1}^{\left(
L-1\right)  },...,\widehat{a}_{n+1}^{\left(  L-1\right)  }\right)  \text{ and
}\widehat{g}\text{ instead of }f\text{, }\left(  a_{0},a_{1},...,a_{n+1}%
\right)  \text{ and }g
\end{array}
\right) \nonumber\\
&  =\left(  \widehat{a}_{0}^{\left(  L-1\right)  },\widehat{g}\widehat{a}%
_{0}^{\left(  L-1\right)  },\widehat{g}\widehat{a}_{1}^{\left(  L-1\right)
},...,\widehat{g}\widehat{a}_{n-1}^{\left(  L-1\right)  },\widehat{a}%
_{n+1}^{\left(  L-1\right)  }\right)  \flat\left(  R^{L}f\right)  .
\label{pf.Rl.scalmult.U}%
\end{align}

Now, let us notice that
\begin{equation}
\widehat{g}\widehat{a}_{k}^{\left(  L-1\right)  }=\widehat{a}_{k+1}^{\left(
L\right)  }\ \ \ \ \ \ \ \ \ \ \text{for every }k\in\left\{
0,1,...,n-1\right\}  . \label{pf.Rl.scalmult.main}%
\end{equation}
\footnote{\textit{Proof of (\ref{pf.Rl.scalmult.main}):} Let $k\in\left\{
0,1,...,n-1\right\}  $. We need to prove that $\widehat{g}\widehat{a}%
_{k}^{\left(  L-1\right)  }=\widehat{a}_{k+1}^{\left(  L\right)  }$. We
distinguish between two cases:
\par
\textit{Case 1:} We have $k\geq L-1$.
\par
\textit{Case 2:} We have $k<L-1$.
\par
Let us first consider Case 1. In this case, we have $k\geq L-1$. Hence,
(\ref{pf.Rl.scalmult.geq}) (applied to $\ell=L-1$) yields $\widehat{a}%
_{k}^{\left(  L-1\right)  }=\dfrac{a_{n+1}a_{k-\left(  L-1\right)  }%
}{a_{n+1-\left(  L-1\right)  }}$. Multiplying $\widehat{g}=\dfrac
{a_{n+1-\left(  L-1\right)  }}{a_{n-\left(  L-1\right)  }}$ with this
equality, we obtain%
\[
\widehat{g}\widehat{a}_{k}^{\left(  L-1\right)  }=\dfrac{a_{n+1-\left(
L-1\right)  }}{a_{n-\left(  L-1\right)  }}\cdot\dfrac{a_{n+1}a_{k-\left(
L-1\right)  }}{a_{n+1-\left(  L-1\right)  }}=\dfrac{a_{n+1}a_{k-\left(
L-1\right)  }}{a_{n-\left(  L-1\right)  }}=\dfrac{a_{n+1}a_{k+1-L}}{a_{n+1-L}%
}.
\]
On the other hand, $k\geq L-1$, so that $k+1\geq L$. Hence,
(\ref{pf.Rl.scalmult.geq}) (applied to $k+1$ and $L$ instead of $k$ and $\ell
$) yields $\widehat{a}_{k+1}^{\left(  L\right)  }=\dfrac{a_{n+1}a_{\left(
k+1\right)  -L}}{a_{n+1-L}}$. Thus, $\widehat{g}\widehat{a}_{k}^{\left(
L-1\right)  }=\dfrac{a_{n+1}a_{k+1-L}}{a_{n+1-L}}=\dfrac{a_{n+1}a_{\left(
k+1\right)  -L}}{a_{n+1-L}}=\widehat{a}_{k+1}^{\left(  L\right)  }$. Hence,
$\widehat{g}\widehat{a}_{k}^{\left(  L-1\right)  }=\widehat{a}_{k+1}^{\left(
L\right)  }$ is proven in Case 1.
\par
Let us now consider Case 2. In this case, we have $k<L-1$. Hence,
(\ref{pf.Rl.scalmult.leq}) (applied to $\ell=L-1$) yields $\widehat{a}%
_{k}^{\left(  L-1\right)  }=\dfrac{a_{n+1+k-\left(  L-1\right)  }a_{0}%
}{a_{n+1-\left(  L-1\right)  }}$. Multiplying $\widehat{g}=\dfrac
{a_{n+1-\left(  L-1\right)  }}{a_{n-\left(  L-1\right)  }}$ with this
equality, we obtain%
\[
\widehat{g}\widehat{a}_{k}^{\left(  L-1\right)  }=\dfrac{a_{n+1-\left(
L-1\right)  }}{a_{n-\left(  L-1\right)  }}\cdot\dfrac{a_{n+1+k-\left(
L-1\right)  }a_{0}}{a_{n+1-\left(  L-1\right)  }}=\dfrac{a_{n+1+k-\left(
L-1\right)  }a_{0}}{a_{n-\left(  L-1\right)  }}=\dfrac{a_{n+1+\left(
k+1\right)  -L}a_{0}}{a_{n+1-L}}.
\]
On the other hand, $k<L-1$, so that $k+1<L$. Hence, (\ref{pf.Rl.scalmult.leq})
(applied to $k+1$ and $L$ instead of $k$ and $\ell$) yields $\widehat{a}%
_{k+1}^{\left(  L\right)  }=\dfrac{a_{n+1+\left(  k+1\right)  -L}a_{0}%
}{a_{n+1-L}}$. Thus, $\widehat{g}\widehat{a}_{k}^{\left(  L-1\right)  }%
=\dfrac{a_{n+1+\left(  k+1\right)  -L}a_{0}}{a_{n+1-L}}=\widehat{a}%
_{k+1}^{\left(  L\right)  }$. Hence, $\widehat{g}\widehat{a}_{k}^{\left(
L-1\right)  }=\widehat{a}_{k+1}^{\left(  L\right)  }$ is proven in Case 2.
\par
Now, we have proven $\widehat{g}\widehat{a}_{k}^{\left(  L-1\right)
}=\widehat{a}_{k+1}^{\left(  L\right)  }$ in each of the two cases 1 and 2.
Since these two cases cover all possibilities, this yields that $\widehat{g}%
\widehat{a}_{k}^{\left(  L-1\right)  }=\widehat{a}_{k+1}^{\left(  L\right)  }$
always holds. This completes the proof of (\ref{pf.Rl.scalmult.main}).}
Furthermore,%
\begin{equation}
\widehat{a}_{0}^{\left(  L-1\right)  }=\widehat{a}_{0}^{\left(  L\right)  }
\label{pf.Rl.scalmult.main.0}%
\end{equation}
(because (\ref{pf.Rl.scalmult.0}) yields $\widehat{a}_{0}^{\left(  L-1\right)
}=a_{0}$ and $\widehat{a}_{0}^{\left(  L\right)  }=a_{0}$) and%
\begin{equation}
\widehat{a}_{n+1}^{\left(  L-1\right)  }=\widehat{a}_{n+1}^{\left(  L\right)
} \label{pf.Rl.scalmult.main.n+1}%
\end{equation}
(because (\ref{pf.Rl.scalmult.n+1}) yields $\widehat{a}_{n+1}^{\left(
L-1\right)  }=a_{n+1}$ and $\widehat{a}_{n+1}^{\left(  L\right)  }=a_{n+1}$).
Due to the equalities (\ref{pf.Rl.scalmult.main.0}),
(\ref{pf.Rl.scalmult.main}) and (\ref{pf.Rl.scalmult.main.n+1}), we have%
\begin{align*}
&  \left(  \widehat{a}_{0}^{\left(  L-1\right)  },\widehat{g}\widehat{a}%
_{0}^{\left(  L-1\right)  },\widehat{g}\widehat{a}_{1}^{\left(  L-1\right)
},...,\widehat{g}\widehat{a}_{n-1}^{\left(  L-1\right)  },\widehat{a}%
_{n+1}^{\left(  L-1\right)  }\right) \\
&  =\left(  \widehat{a}_{0}^{\left(  L\right)  },\widehat{a}_{1}^{\left(
L\right)  },\widehat{a}_{2}^{\left(  L\right)  },...,\widehat{a}_{n}^{\left(
L\right)  },\widehat{a}_{n+1}^{\left(  L\right)  }\right)  =\left(
\widehat{a}_{0}^{\left(  L\right)  },\widehat{a}_{1}^{\left(  L\right)
},...,\widehat{a}_{n+1}^{\left(  L\right)  }\right)  .
\end{align*}
Thus, (\ref{pf.Rl.scalmult.U}) becomes%
\begin{align*}
R^{L}\left(  \left(  a_{0},a_{1},...,a_{n+1}\right)  \flat f\right)   &
=\underbrace{\left(  \widehat{a}_{0}^{\left(  L-1\right)  },\widehat{g}%
\widehat{a}_{0}^{\left(  L-1\right)  },\widehat{g}\widehat{a}_{1}^{\left(
L-1\right)  },...,\widehat{g}\widehat{a}_{n-1}^{\left(  L-1\right)
},\widehat{a}_{n+1}^{\left(  L-1\right)  }\right)  }_{=\left(  \widehat{a}%
_{0}^{\left(  L\right)  },\widehat{a}_{1}^{\left(  L\right)  },...,\widehat{a}%
_{n+1}^{\left(  L\right)  }\right)  }\flat\left(  R^{L}f\right) \\
&  =\left(  \widehat{a}_{0}^{\left(  L\right)  },\widehat{a}_{1}^{\left(
L\right)  },...,\widehat{a}_{n+1}^{\left(  L\right)  }\right)  \flat\left(
R^{L}f\right)  .
\end{align*}
In other words, Proposition \ref{prop.Rl.scalmult} holds for $\ell=L$. This
completes the induction step, and thus Proposition \ref{prop.Rl.scalmult} is
proven by induction.
\end{proof}
\end{verlong}

As a consequence of Proposition \ref{prop.Rl.scalmult}, we notice a very
simple behavior of rescaled labellings under $R^{n+1}$ for an $n$-graded poset
$P$:

\begin{corollary}
\label{cor.Rl.scalmult}Let $n\in\mathbb{N}$. Let $\mathbb{K}$ be a field. Let
$P$ be an $n$-graded poset. For every $\mathbb{K}$-labelling $f\in
\mathbb{K}^{\widehat{P}}$ and any $\left(  n+2\right)  $-tuple $\left(
a_{0},a_{1},...,a_{n+1}\right)  \in\left(  \mathbb{K}^{\times}\right)  ^{n+2}%
$, we define a $\mathbb{K}$-labelling $\left(  a_{0},a_{1},...,a_{n+1}\right)
\flat f\in\mathbb{K}^{\widehat{P}}$ as in Definition \ref{def.bemol}.

Let $\left(  a_{0},a_{1},...,a_{n+1}\right)  \in\left(  \mathbb{K}^{\times
}\right)  ^{n+2}$. Let $f\in\mathbb{K}^{\widehat{P}}$ be a $\mathbb{K}%
$-labelling. Then,%
\[
R^{n+1}\left(  \left(  a_{0},a_{1},...,a_{n+1}\right)  \flat f\right)
=\left(  a_{0},a_{1},...,a_{n+1}\right)  \flat\left(  R^{n+1}f\right)
\]
(provided that $R^{n+1}f$ is well-defined).
\end{corollary}

\begin{vershort}
\begin{proof}
[\nopunct]We leave deriving Corollary \ref{cor.Rl.scalmult} from Proposition
\ref{prop.Rl.scalmult} to the reader.
\end{proof}
\end{vershort}

\begin{verlong}
\begin{proof}
[Proof of Corollary \ref{cor.Rl.scalmult} (sketched).]Let us use the notations
introduced in Proposition \ref{prop.Rl.scalmult}. It is easy to see that every
$\ell\in\left\{  0,1,...,n+1\right\}  $ and $k\in\left\{  0,1,...,n+1\right\}
$ satisfy (\ref{pf.Rl.scalmult.leq}) (in fact, this can be proven in the same
way as in the proof of Proposition \ref{prop.Rl.scalmult}). Thus, every
$k\in\left\{  0,1,...,n+1\right\}  $ satisfies%
\begin{align*}
\widehat{a}_{k}^{\left(  n+1\right)  }  &  =\dfrac{a_{n+1+k-\left(
n+1\right)  }a_{0}}{a_{n+1-\left(  n+1\right)  }}\ \ \ \ \ \ \ \ \ \ \left(
\text{by (\ref{pf.Rl.scalmult.leq}), applied to }\ell=n+1\text{ (since }k\leq
n+1\text{)}\right) \\
&  =\dfrac{a_{k}a_{0}}{a_{0}}\ \ \ \ \ \ \ \ \ \ \left(  \text{since
}n+1+k-\left(  n+1\right)  =k\text{ and }n+1-\left(  n+1\right)  =0\right) \\
&  =a_{k}.
\end{align*}
In other words, $\left(  \widehat{a}_{0}^{\left(  n+1\right)  },\widehat{a}%
_{1}^{\left(  n+1\right)  },...,\widehat{a}_{n+1}^{\left(  n+1\right)
}\right)  =\left(  a_{0},a_{1},...,a_{n+1}\right)  $. But Proposition
\ref{prop.Rl.scalmult} (applied to $\ell=n+1$) yields%
\begin{align*}
R^{n+1}\left(  \left(  a_{0},a_{1},...,a_{n+1}\right)  \flat f\right)   &
=\underbrace{\left(  \widehat{a}_{0}^{\left(  n+1\right)  },\widehat{a}%
_{1}^{\left(  n+1\right)  },...,\widehat{a}_{n+1}^{\left(  n+1\right)
}\right)  }_{=\left(  a_{0},a_{1},...,a_{n+1}\right)  }\flat\left(
R^{n+1}f\right) \\
&  =\left(  a_{0},a_{1},...,a_{n+1}\right)  \flat\left(  R^{n+1}f\right)  .
\end{align*}
This proves Corollary \ref{cor.Rl.scalmult}.
\end{proof}
\end{verlong}

Let us furthermore record how the rescaling of labels according to their
degree affects their w-tuples (as defined in Definition \ref{def.wi}):

\begin{proposition}
\label{prop.w.scalmult}Let $n\in\mathbb{N}$. Let $\mathbb{K}$ be a field. Let
$P$ be an $n$-graded poset. For every $\mathbb{K}$-labelling $f\in
\mathbb{K}^{\widehat{P}}$ and any $\left(  n+2\right)  $-tuple $\left(
a_{0},a_{1},...,a_{n+1}\right)  \in\left(  \mathbb{K}^{\times}\right)  ^{n+2}%
$, we define a $\mathbb{K}$-labelling $\left(  a_{0},a_{1},...,a_{n+1}\right)
\flat f\in\mathbb{K}^{\widehat{P}}$ as in Definition \ref{def.bemol}.

Let $f\in\mathbb{K}^{\widehat{P}}$ be a $\mathbb{K}$-labelling of $P$. Let
$\left(  a_{0},a_{1},...,a_{n+1}\right)  \in\left(  \mathbb{K}^{\times
}\right)  ^{n+2}$. Then, the w-tuple of the $\mathbb{K}$-labelling $\left(
a_{0},a_{1},...,a_{n+1}\right)  \flat f$ is%
\[
\left(  \dfrac{a_{0}}{a_{1}}\mathbf{w}_{0}\left(  f\right)  ,\dfrac{a_{1}%
}{a_{2}}\mathbf{w}_{1}\left(  f\right)  ,...,\dfrac{a_{n}}{a_{n+1}}%
\mathbf{w}_{n}\left(  f\right)  \right)  .
\]

\end{proposition}

\begin{vershort}
\begin{proof}
[\nopunct]Proposition \ref{prop.w.scalmult} follows by computation using just
the definitions of the notions involved.
\end{proof}
\end{vershort}

\begin{verlong}
\begin{proof}
[Proof of Proposition \ref{prop.w.scalmult} (sketched).]Let $g=\left(
a_{0},a_{1},...,a_{n+1}\right)  \flat f$. Then, every $j\in\left\{
0,1,...,n+1\right\}  $ and every $v\in\widehat{P}_{j}$ satisfy%
\begin{align}
g\left(  v\right)   &  =\left(  \left(  a_{0},a_{1},...,a_{n+1}\right)  \flat
f\right)  \left(  v\right)  =a_{\deg v}\cdot f\left(  v\right) \nonumber\\
&  \ \ \ \ \ \ \ \ \ \ \left(  \text{by the definition of }\left(  a_{0}%
,a_{1},...,a_{n+1}\right)  \flat f\right) \nonumber\\
&  =a_{j}\cdot f\left(  v\right)  \ \ \ \ \ \ \ \ \ \ \left(  \text{since
}v\in\widehat{P}_{j}\text{, and therefore }\deg v=j\right)  .
\label{pf.w.scalmult.1}%
\end{align}
Let $i\in\left\{  0,1,...,n\right\}  $. Then, by the definition of
$\mathbf{w}_{i}\left(  g\right)  $, we have%
\begin{align}
\mathbf{w}_{i}\left(  g\right)   &  =\sum_{\substack{x\in\widehat{P}%
_{i};\ y\in\widehat{P}_{i+1};\\y \gtrdot x}}\dfrac{g\left(  x\right)
}{g\left(  y\right)  }=\sum_{\substack{x\in\widehat{P}_{i};\ y\in
\widehat{P}_{i+1};\\y \gtrdot x}}\left(  \underbrace{g\left(  y\right)
}_{\substack{=a_{i+1}\cdot f\left(  y\right)  \\\text{(by
(\ref{pf.w.scalmult.1}), applied to}\\j=i+1\text{ and }v=y\text{)}}}\right)
^{-1}\cdot\underbrace{g\left(  x\right)  }_{\substack{=a_{i}\cdot f\left(
x\right)  \\\text{(by (\ref{pf.w.scalmult.1}), applied to}\\j=i\text{ and
}v=x\text{)}}}\nonumber\\
&  =\sum_{\substack{x\in\widehat{P}_{i};\ y\in\widehat{P}_{i+1};\\y \gtrdot
x}}\left(  a_{i+1}\cdot f\left(  y\right)  \right)  ^{-1}\cdot a_{i}\cdot
f\left(  x\right)  =\sum_{\substack{x\in\widehat{P}_{i};\ y\in\widehat{P}%
_{i+1};\\y \gtrdot x}}\dfrac{a_{i}}{a_{i+1}}\cdot\dfrac{f\left(  x\right)
}{f\left(  y\right)  }\nonumber\\
&  =\dfrac{a_{i}}{a_{i+1}}\cdot\underbrace{\sum_{\substack{x\in\widehat{P}%
_{i};\ y\in\widehat{P}_{i+1};\\y \gtrdot x}}\dfrac{f\left(  x\right)
}{f\left(  y\right)  }}_{=\mathbf{w}_{i}\left(  f\right)  }=\dfrac{a_{i}%
}{a_{i+1}}\mathbf{w}_{i}\left(  f\right)  . \label{pf.w.scalmult.2}%
\end{align}

Now, forget that we fixed $i$. Now, $\left(  a_{0},a_{1},...,a_{n+1}\right)
\flat f=g$. Thus, the w-tuple of the $\mathbb{K}$-labelling $\left(
a_{0},a_{1},...,a_{n+1}\right)  \flat f$ equals the w-tuple of $g$, and this
is the $\left(  n+1\right)  $-tuple%
\[
\left(  \mathbf{w}_{0}\left(  g\right)  ,\mathbf{w}_{1}\left(  g\right)
,...,\mathbf{w}_{n}\left(  g\right)  \right)  =\left(  \dfrac{a_{0}}{a_{1}%
}\mathbf{w}_{0}\left(  f\right)  ,\dfrac{a_{1}}{a_{2}}\mathbf{w}_{1}\left(
f\right)  ,...,\dfrac{a_{n}}{a_{n+1}}\mathbf{w}_{n}\left(  f\right)  \right)
\]
(because every $i\in\left\{  0,1,...,n\right\}  $ satisfies
(\ref{pf.w.scalmult.2})). This proves Proposition \ref{prop.w.scalmult}.
\end{proof}

Let us state a complete and utter triviality before we end this section:

\begin{lemma}
\label{lem.bemol.cancel} Let $n\in\mathbb{N}$. Let $\mathbb{K}$ be a field.
Let $P$ be an $n$-graded poset. For every $\mathbb{K}$-labelling
$f\in\mathbb{K}^{\widehat{P}}$ and any $\left(  n+2\right)  $-tuple $\left(
a_{0},a_{1},...,a_{n+1}\right)  \in\left(  \mathbb{K}^{\times}\right)  ^{n+2}%
$, we define a $\mathbb{K}$-labelling $\left(  a_{0},a_{1},...,a_{n+1}\right)
\flat f\in\mathbb{K}^{\widehat{P}}$ as in Definition \ref{def.bemol}.

Let $f\in\mathbb{K}^{\widehat{P}}$ and $g\in\mathbb{K}^{\widehat{P}}$ be two
$\mathbb{K}$-labellings of $P$. Let $\left(  a_{0},a_{1},...,a_{n+1}\right)
\in\left(  \mathbb{K}^{\times}\right)  ^{n+2}$. Assume that $\left(
a_{0},a_{1},...,a_{n+1}\right)  \flat f=\left(  a_{0},a_{1},...,a_{n+1}%
\right)  \flat g$. Then, $f=g$.
\end{lemma}

Roughly speaking, Lemma \ref{lem.bemol.cancel} states that if we have an
equality of the form $\left(  a_{0},a_{1},...,a_{n+1}\right)  \flat f=\left(
a_{0},a_{1},...,a_{n+1}\right)  \flat g$, we can \textquotedblleft
cancel\textquotedblright\ the \textquotedblleft$\left(  a_{0},a_{1}%
,...,a_{n+1}\right)  \flat$\textquotedblright\ from this equality, to obtain
$f=g$.

\begin{proof}
[Proof of Lemma \ref{lem.bemol.cancel}.]We know that $\left(  a_{0}%
,a_{1},...,a_{n+1}\right)  \in\left(  \mathbb{K}^{\times}\right)  ^{n+2}$.
Hence, $a_{i} \in\mathbb{K}^{\times}$ for every $i\in\left\{
0,1,...,n+1\right\}  $. Thus, $a_{i} \neq0$ for every $i\in\left\{
0,1,...,n+1\right\}  $.

Now, let $v\in\widehat{P}$. We have $a_{\deg v}\neq0$ (since $a_{i}\neq0$ for
every $i\in\left\{  0,1,...,n+1\right\}  $). Now,%
\begin{align*}
&  \underbrace{\left(  \left(  a_{0},a_{1},...,a_{n+1}\right)  \flat f\right)
}_{=\left(  a_{0},a_{1},...,a_{n+1}\right)  \flat g}\left(  v\right) \\
&  =\left(  \left(  a_{0},a_{1},...,a_{n+1}\right)  \flat g\right)  \left(
v\right)  =a_{\deg v}\cdot g\left(  v\right)  \ \ \ \ \ \ \ \ \ \ \left(
\text{by the definition of }\left(  a_{0},a_{1},...,a_{n+1}\right)  \flat
g\right)  .
\end{align*}
Compared with%
\begin{align*}
&  \left(  \left(  a_{0},a_{1},...,a_{n+1}\right)  \flat f\right)  \left(
v\right) \\
&  =a_{\deg v}\cdot f\left(  v\right)  \ \ \ \ \ \ \ \ \ \ \left(  \text{by
the definition of }\left(  a_{0},a_{1},...,a_{n+1}\right)  \flat f\right)  ,
\end{align*}
this yields $a_{\deg v}\cdot g\left(  v\right)  =a_{\deg v}\cdot f\left(
v\right)  $. We can cancel $a_{\deg v}$ from this equation (since $a_{\deg
v}\neq0$), and thus obtain $g\left(  v\right)  =f\left(  v\right)  $.

Now, forget that we fixed $v$. We thus have shown that $g\left(  v\right)
=f\left(  v\right)  $ for every $v\in\widehat{P}$. In other words, $g=f$.
Thus, $f=g$, so that Lemma \ref{lem.bemol.cancel} is proven.
\end{proof}
\end{verlong}

\section{\label{sect.homogeneous}Homogeneous labellings}

In the previous section, we have quantified how the rescaling of all labels at
a given degree affects a labelling (of a graded poset) under birational
rowmotion. In this section, we will introduce a notion of \textquotedblleft
homogeneous labellings\textquotedblright\ which (roughly speaking) are
\textquotedblleft labellings up to rescaling at a given
degree\textquotedblright\ in the same way as a point in a projective space can
be regarded as (roughly speaking) \textquotedblleft a point in the affine
space up to rescaling the coordinates\textquotedblright. To be precise, we
will need to restrict ourselves to considering only \textquotedblleft
zero-free\textquotedblright\ labellings (a Zariski-dense open subset of all
labellings) for the same reason as we need to exclude $0$ when defining a
projective space. Once done with the definitions, we will see that birational
rowmotion (and the maps $R_{i}$) can be defined on homogeneous labellings (it
is here that we will make use of the results of the previous section).

Let us begin with the definitions:

\begin{definition}
Let $\mathbb{K}$ be a field.

\textbf{(a)} For every $\mathbb{K}$-vector space $V$, let $\mathbb{P}\left(
V\right)  $ denote the projective space of $V$ (that is, the set of
equivalence classes of vectors in $V\setminus\left\{  0\right\}  $
modulo proportionality).

\textbf{(b)} For every $n\in\mathbb{N}$, we let $\mathbb{P}^{n}\left(
\mathbb{K}\right)  $ denote the projective space $\mathbb{P}\left(
\mathbb{K}^{n+1}\right)  $.
\end{definition}

\begin{definition}
\label{def.hgeq}Let $n\in\mathbb{N}$. Let $\mathbb{K}$ be a field. Let $P$ be
an $n$-graded poset.

\textbf{(a)} Denote by $\overline{\mathbb{K}^{\widehat{P}}}$ the product
$\prod\limits_{i=1}^{n}\mathbb{P}\left(  \mathbb{K}^{\widehat{P}_{i}}\right)
$ of projective spaces. Notice that the product is just a Cartesian product of
algebraic varieties, and a reader unfamiliar with algebraic geometry can just
regard it as a Cartesian product of sets.\footnotemark\ 

We have $\overline{\mathbb{K}^{\widehat{P}}}=\prod\limits_{i=1}^{n}%
\mathbb{P}\left(  \mathbb{K}^{\widehat{P}_{i}}\right)  \cong\prod
\limits_{i=1}^{n}\mathbb{P}^{\left\vert \widehat{P}_{i}\right\vert -1}\left(
\mathbb{K}\right)  $ (since every $i \in \left\{1, 2, ..., n\right\}$
satisfies
$\mathbb{P}\left(  \mathbb{K}^{\widehat{P}_{i}}\right)
\cong \mathbb{P}^{\left\vert \widehat{P}_{i}\right\vert -1}\left(
\mathbb{K}\right)$). We denote the elements of $\overline{\mathbb{K}%
^{\widehat{P}}}$ as \textit{homogeneous labellings}.

Notice that $\overline{\mathbb{K}^{\widehat{P}}}=\prod\limits_{i=1}%
^{n}\mathbb{P}\left(  \mathbb{K}^{\widehat{P}_{i}}\right)  \cong%
\prod\limits_{i=0}^{n+1}\mathbb{P}\left(  \mathbb{K}^{\widehat{P}_{i}}\right)
$ (as algebraic varieties). This is because $\mathbb{K}^{\widehat{P}_{0}}$ and
$\mathbb{K}^{\widehat{P}_{n+1}}$ are $1$-dimensional vector spaces
(since $\left|\widehat{P}_{0}\right| = 1$ and
$\left|\widehat{P}_{n+1}\right| = 1$), and thus
the projective spaces $\mathbb{P}\left(  \mathbb{K}^{\widehat{P}_{0}}\right)
$ and $\mathbb{P}\left(  \mathbb{K}^{\widehat{P}_{n+1}}\right)  $ each consist
of a single point.

\textbf{(b)} A $\mathbb{K}$-labelling $f\in\mathbb{K}^{\widehat{P}}$ is said
to be \textit{zero-free} if for every $i\in\left\{  0,1,...,n+1\right\}  $,
there exists some $v\in\widehat{P}_{i}$ satisfying $f\left(  v\right)  \neq0$.
(In other words, a $\mathbb{K}$-labelling $f\in\mathbb{K}^{\widehat{P}}$ is
said to be zero-free if there exists no $i\in\left\{  0,1,...,n+1\right\}  $
such that $f$ is identically $0$ on all elements of $\widehat{P}$ having
degree $i$.) Let $\mathbb{K}_{\neq0}^{\widehat{P}}$ be the set of all
zero-free $\mathbb{K}$-labellings. Clearly, this set $\mathbb{K}_{\neq
0}^{\widehat{P}}$ is a Zariski-dense open subset of $\mathbb{K}^{\widehat{P}}$.

\textbf{(c)} Identify the set $\mathbb{K}^{\widehat{P}}$ with $\prod
\limits_{i=0}^{n+1}\mathbb{K}^{\widehat{P}_{i}}$ in the obvious way (since
$\widehat{P}$, regarded as a set, is the disjoint union of the sets
$\widehat{P}_{i}$ over all $i\in\left\{  0,1,...,n+1\right\}  $).

Using the identifications $\mathbb{K}^{\widehat{P}} \cong
\prod\limits_{i=0}^{n+1}\mathbb{K}^{\widehat{P}_{i}}$ and
$\overline{\mathbb{K}^{\widehat{P}}}
\cong \prod\limits_{i=0}^{n+1}
\mathbb{P}\left(  \mathbb{K}^{\widehat{P}_{i}}\right)$, we now
define a
rational map $\pi:\mathbb{K}^{\widehat{P}}\dashrightarrow\overline
{\mathbb{K}^{\widehat{P}}}$ as the product of the canonical projections
$\mathbb{K}^{\widehat{P}_{i}}\dashrightarrow\mathbb{P}\left(  \mathbb{K}%
^{\widehat{P}_{i}}\right)  $ (which are defined everywhere outside of the
$\left\{  0\right\}  $ subsets) over all $i\in\left\{  0,1,...,n+1\right\}  $.
Notice that the domain of definition of this rational map $\pi$ is precisely
$\mathbb{K}_{\neq0}^{\widehat{P}}$. For every $f\in\mathbb{K}^{\widehat{P}}$,
we denote $\pi\left(  f\right)  $ as the \textit{homogenization} of the
$\mathbb{K}$-labelling $f$.

\textbf{(d)} Two zero-free $\mathbb{K}$-labellings $f\in\mathbb{K}%
^{\widehat{P}}$ and $g\in\mathbb{K}^{\widehat{P}}$ are said to be
\textit{homogeneously equivalent} if and only if they satisfy one of the
following equivalent conditions:

\textit{Condition 1:} For every $i\in\left\{  0,1,...,n+1\right\}  $ and any
two elements $x$ and $y$ of $\widehat{P}_{i}$, we have $\dfrac{f\left(
x\right)  }{f\left(  y\right)  }=\dfrac{g\left(  x\right)  }{g\left(
y\right)  }$.

\textit{Condition 2:} There exists an $\left(  n+2\right)  $-tuple $\left(
a_{0},a_{1},...,a_{n+1}\right)  \in\left(  \mathbb{K}^{\times}\right)  ^{n+2}$
such that every $x\in\widehat{P}$ satisfies $g\left(  x\right)  =a_{\deg
x}\cdot f\left(  x\right)  $.

\textit{Condition 3:} We have $\pi\left(  f\right)  =\pi\left(  g\right)  $.

(The equivalence between these three conditions is very easy to check. We will
never actually use Condition 1.)
\end{definition}

\footnotetext{The structure of algebraic variety will only be needed to define
the Zariski topology on $\overline{\mathbb{K}^{\widehat{P}}}$, which is more
or less obvious already (e.g., when we say that something holds ``for almost
every element $x$ of $\prod\limits_{i=1}^{n}\mathbb{P}\left(  \mathbb{K}%
^{\widehat{P}_{i}}\right)  $'', we could equivalently say that it holds ``for
$x=\operatorname*{proj}\left(  X\right)  $ for almost every element $X$ of
$\prod\limits_{i=1}^{n}\left(  \mathbb{K}^{\widehat{P}_{i}}\setminus\left\{
0\right\}  \right)  $'', where $\operatorname*{proj}$ is the canonical map
$\prod\limits_{i=1}^{n}\left(  \mathbb{K}^{\widehat{P}_{i}}\setminus\left\{
0\right\}  \right)  \rightarrow\prod\limits_{i=1}^{n}\mathbb{P}\left(
\mathbb{K}^{\widehat{P}_{i}}\right)  $ defined as the product of the
projections $\mathbb{K}^{\widehat{P}_{i}}\setminus\left\{  0\right\}
\rightarrow\mathbb{P}\left(  \mathbb{K}^{\widehat{P}_{i}}\right)  $).}

\begin{remark}
Clearly, homogeneous equivalence is an equivalence relation on the set
$\mathbb{K}_{\neq0}^{\widehat{P}}$ of all zero-free $\mathbb{K}$-labellings.
We can identify $\overline{\mathbb{K}^{\widehat{P}}}$ with the quotient of the
set $\mathbb{K}_{\neq0}^{\widehat{P}}$ modulo this relation. Then, $\pi$
becomes the canonical projection map $\mathbb{K}^{\widehat{P}}\dashrightarrow
\overline{\mathbb{K}^{\widehat{P}}}$.
\end{remark}

One remark about the notion \textquotedblleft zero-free\textquotedblright:
Being zero-free is a very weak condition on a $\mathbb{K}$-labelling (indeed
the zero-free $\mathbb{K}$-labellings form a Zariski-dense open subset of the
space of all $\mathbb{K}$-labellings), and the $\mathbb{K}$-labellings which
don't satisfy this condition are rather useless for us (if $f$ is a
$\mathbb{K}$-labelling which is not zero-free, then $R^{2}f$ is not
well-defined, and usually not even $Rf$ is well-defined). We are almost never
giving up any generality if we require a labelling to be zero-free.

\begin{remark}
Let $n\in\mathbb{N}$. Let $\mathbb{K}$ be a field. Let $P$ be an $n$-graded
poset. For every $\mathbb{K}$-labelling $f\in\mathbb{K}^{\widehat{P}}$ and any
$\left(  n+2\right)  $-tuple $\left(  a_{0},a_{1},...,a_{n+1}\right)
\in\left(  \mathbb{K}^{\times}\right)  ^{n+2}$, we define a $\mathbb{K}%
$-labelling $\left(  a_{0},a_{1},...,a_{n+1}\right)  \flat f\in\mathbb{K}%
^{\widehat{P}}$ as in Definition \ref{def.bemol}.

Let $f\in\mathbb{K}^{\widehat{P}}$ be a zero-free $\mathbb{K}$-labelling of
$P$. Let $\left(  a_{0},a_{1},...,a_{n+1}\right)  \in\left(  \mathbb{K}%
^{\times}\right)  ^{n+2}$. Then, the $\mathbb{K}$-labelling $\left(
a_{0},a_{1},...,a_{n+1}\right)  \flat f$ is also zero-free. (This follows
immediately from the definitions.)
\end{remark}

\begin{definition}
\label{def.hgeq.pii}Let $n\in\mathbb{N}$. Let $\mathbb{K}$ be a field. Let $P$
be an $n$-graded poset. For every zero-free $f\in\mathbb{K}^{\widehat{P}}$ and
every $i \in\left\{  1,2,...,n\right\}  $, the image of the restriction of $f
: \widehat{P} \to\mathbb{K}$ to $\widehat{P}_{i}$ under the canonical
projection $\mathbb{K}^{\widehat{P}_{i}}\dashrightarrow\mathbb{P}\left(
\mathbb{K}^{\widehat{P}_{i}}\right)  $ will be denoted by $\pi_{i}\left(
f\right)  $. This image $\pi_{i}\left(  f\right)  $ encodes the values of $f$
at the elements of $\widehat{P}$ of degree $i$ up to multiplying all these
values by a common nonzero scalar. Notice that
\begin{equation}
\pi\left(  f\right)  =\left(  \pi_{1}\left(  f\right)  ,\pi_{2}\left(
f\right)  ,...,\pi_{n}\left(  f\right)  \right)  \label{def.hgeq.pii.eq}%
\end{equation}
for every $f\in\mathbb{K}^{\widehat{P}}$. (Here, the right hand side of
(\ref{def.hgeq.pii.eq}) is regarded as an element of $\overline{\mathbb{K}%
^{\widehat{P}}}$ because it belongs to $\prod\limits_{i=1}^{n}\mathbb{P}%
\left(  \mathbb{K}^{\widehat{P}_{i}}\right)  =\overline{\mathbb{K}%
^{\widehat{P}}}$.)
\end{definition}

We are next going to see:

\begin{corollary}
\label{cor.hgRi}Let $n\in\mathbb{N}$. Let $\mathbb{K}$ be a field. Let $P$ be
an $n$-graded poset. Let $i\in\left\{  1,2,...,n\right\}  $. If $f\in
\mathbb{K}^{\widehat{P}}$ and $g\in\mathbb{K}^{\widehat{P}}$ are two
homogeneously equivalent zero-free $\mathbb{K}$-labellings, then $R_{i}f$ is
homogeneously equivalent to $R_{i}g$ (as long as $R_{i}f$ and $R_{i}g$ are zero-free).
\end{corollary}

\begin{corollary}
\label{cor.hgR}Let $n\in\mathbb{N}$. Let $\mathbb{K}$ be a field. Let $P$ be
an $n$-graded poset. If $f\in\mathbb{K}^{\widehat{P}}$ and $g\in
\mathbb{K}^{\widehat{P}}$ are two homogeneously equivalent zero-free
$\mathbb{K}$-labellings, then $Rf$ is homogeneously equivalent to $Rg$ (as
long as $Rf$ and $Rg$ are zero-free).
\end{corollary}

Notice that Corollary \ref{cor.hgRi} would not be valid if we were to replace
$R_{i}$ by a single toggle $T_{v}$! So the operators $R_{i}$ in some sense
combine the nice properties of $T_{v}$ (like being an involution, cf.
Proposition \ref{prop.Ri.invo}) with the nice properties of $R$ (like having
an easily describable action on w-tuples, cf. Proposition \ref{prop.wi.Ri},
and respecting homogeneous equivalence, cf. Corollary \ref{cor.hgRi}).

\begin{proof}
[Proof of Corollary \ref{cor.hgRi} (sketched).]Let $f\in\mathbb{K}%
^{\widehat{P}}$ and $g\in\mathbb{K}^{\widehat{P}}$ be two homogeneously
equivalent zero-free $\mathbb{K}$-labellings.

We know that $f$ and $g$ are homogeneously equivalent. By Condition 2 in
Definition \ref{def.hgeq} \textbf{(d)}, this means that there exists an
$\left(  n+2\right)  $-tuple $\left(  a_{0},a_{1},...,a_{n+1}\right)
\in\left(  \mathbb{K}^{\times}\right)  ^{n+2}$ such that every $x\in
\widehat{P}$ satisfies $g\left(  x\right)  =a_{\deg x}\cdot f\left(  x\right)
$. In other words, there exists an $\left(  n+2\right)  $-tuple $\left(
a_{0},a_{1},...,a_{n+1}\right)  \in\left(  \mathbb{K}^{\times}\right)  ^{n+2}$
such that
\[
g=\left(  a_{0},a_{1},...,a_{n+1}\right)  \flat f.
\]
Consider this $\left(  n+2\right)  $-tuple $\left(  a_{0},a_{1},...,a_{n+1}%
\right)  $. Since $g=\left(  a_{0},a_{1},...,a_{n+1}\right)  \flat f$, we
have
\begin{align*}
&  R_{i}g=R_{i}\left(  \left(  a_{0},a_{1},...,a_{n+1}\right)  \flat f\right)
\\
&  =\left(  a_{0},a_{1},...,a_{i-1},\dfrac{a_{i+1}a_{i-1}}{a_{i}}%
,a_{i+1},a_{i+2},...,a_{n+1}\right)  \flat\left(  R_{i}f\right)
\end{align*}
(by Proposition \ref{prop.Ri.scalmult}). Hence, there exists an $\left(
n+2\right)  $-tuple $\left(  b_{0},b_{1},...,b_{n+1}\right)  \in\left(
\mathbb{K}^{\times}\right)  ^{n+2}$ such that%
\[
R_{i}g=\left(  b_{0},b_{1},...,b_{n+1}\right)  \flat\left(  R_{i}f\right)
\]
(namely, $\left(  b_{0},b_{1},...,b_{n+1}\right)  =\left(  a_{0}%
,a_{1},...,a_{i-1},\dfrac{a_{i+1}a_{i-1}}{a_{i}},a_{i+1},a_{i+2}%
,...,a_{n+1}\right)  $). In other words, there exists an $\left(  n+2\right)
$-tuple $\left(  b_{0},b_{1},...,b_{n+1}\right)  \in\left(  \mathbb{K}%
^{\times}\right)  ^{n+2}$ such that every $x\in\widehat{P}$ satisfies $\left(
R_{i}g\right)  \left(  x\right)  =b_{\deg x}\cdot\left(  R_{i}f\right)
\left(  x\right)  $. But this is precisely Condition 2 in Definition
\ref{def.hgeq} \textbf{(d)}, stated for the labellings $R_{i}f$ and $R_{i}g$
instead of $f$ and $g$. Hence, $R_{i}f$ and $R_{i}g$ are homogeneously
equivalent. This proves Corollary \ref{cor.hgRi}.
\end{proof}

\begin{proof}
[Proof of Corollary \ref{cor.hgR} (sketched).]Let $f\in\mathbb{K}%
^{\widehat{P}}$ and $g\in\mathbb{K}^{\widehat{P}}$ be two homogeneously
equivalent zero-free $\mathbb{K}$-labellings. By iterative application of
Corollary \ref{cor.hgRi}, we then conclude that the $\mathbb{K}$-labellings
$\left(  R_{1}\circ R_{2}\circ...\circ R_{n}\right)  f$ and $\left(
R_{1}\circ R_{2}\circ...\circ R_{n}\right)  g$ are homogeneously equivalent
(if they are well-defined). Since $R_{1}\circ R_{2}\circ...\circ R_{n}=R$ (by
Proposition \ref{prop.Ri.R}), this rewrites as follows: The $\mathbb{K}%
$-labellings $Rf$ and $Rg$ are homogeneously equivalent. This proves Corollary
\ref{cor.hgR}.
\end{proof}

\begin{verlong}
Let us actually state the obvious:

\begin{proposition}
\label{prop.bemol.hgeq}Let $n\in\mathbb{N}$. Let $\mathbb{K}$ be a field. Let
$P$ be an $n$-graded poset. Let us use the notation introduced in Definition
\ref{def.bemol}.

\textbf{(a)} If $f\in\mathbb{K}^{\widehat{P}}$ is a zero-free $\mathbb{K}%
$-labelling, and $\left(  a_{0},a_{1},...,a_{n+1}\right)  \in\left(
\mathbb{K}^{\times}\right)  ^{n+2}$, then the $\mathbb{K}$-labellings $f$ and
$\left(  a_{0},a_{1},...,a_{n+1}\right)  \flat f$ are homogeneously equivalent.

\textbf{(b)} Let $f\in\mathbb{K}^{\widehat{P}}$ and $g\in\mathbb{K}%
^{\widehat{P}}$ be two zero-free $\mathbb{K}$-labellings. If $f$ and $g$ are
homogeneously equivalent, then there exists an $\left(  a_{0},a_{1}%
,...,a_{n+1}\right)  \in\left(  \mathbb{K}^{\times}\right)  ^{n+2}$ such that
$g=\left(  a_{0},a_{1},...,a_{n+1}\right)  \flat f$.
\end{proposition}

\begin{proof}
[Proof of Proposition \ref{prop.bemol.hgeq} (sketched).]\textbf{(a)} Let
$f\in\mathbb{K}^{\widehat{P}}$ be a zero-free $\mathbb{K}$-labelling. Let
$\left(  a_{0},a_{1},...,a_{n+1}\right)  \in\left(  \mathbb{K}^{\times
}\right)  ^{n+2}$. Let $g=\left(  a_{0},a_{1},...,a_{n+1}\right)  \flat f$.
Then, every $x\in\widehat{P}$ satisfies%
\[
g\left(  x\right)  =\left(  \left(  a_{0},a_{1},...,a_{n+1}\right)  \flat
f\right)  \left(  x\right)  =a_{\deg x}\cdot f\left(  x\right)
\]
(by the definition of $\left(  a_{0},a_{1},...,a_{n+1}\right)  \flat f$).
Hence, there exists an $\left(  n+2\right)  $-tuple $\left(  a_{0}%
,a_{1},...,a_{n+1}\right)  \in\left(  \mathbb{K}^{\times}\right)  ^{n+2}$ such
that every $x\in\widehat{P}$ satisfies $g\left(  x\right)  =a_{\deg x}\cdot
f\left(  x\right)  $. But this is exactly Condition 2 of Definition
\ref{def.hgeq} \textbf{(d)}. Since this condition is a necessary and
sufficient condition for $f$ and $g$ being homogeneously equivalent, this
yields that the $\mathbb{K}$-labellings $f$ and $g$ are homogeneously
equivalent. Since $g=\left(  a_{0},a_{1},...,a_{n+1}\right)  \flat f$, this
rewrites as follows: The $\mathbb{K}$-labellings $f$ and $\left(  a_{0}%
,a_{1},...,a_{n+1}\right)  \flat f$ are homogeneously equivalent. Proposition
\ref{prop.bemol.hgeq} \textbf{(a)} is thus proven.

\textbf{(b)} The proof of Proposition \ref{prop.bemol.hgeq} \textbf{(b)} is
obtained by arguing the proof of Proposition \ref{prop.bemol.hgeq}
\textbf{(a)} backwards.
\end{proof}
\end{verlong}

Let us introduce a general piece of notation:

\Needspace{14\baselineskip}

\begin{definition}
Let $S$ and $T$ be two sets. Let $\sim_{S}$ be an equivalence relation on the
set $S$, and let $\sim_{T}$ be an equivalence relation on the set $T$. Let
$\overline{S}$ be the quotient of the set $S$ modulo the equivalence relation
$\sim_{S}$, and let $\overline{T}$ be the quotient of the set $T$ modulo the
equivalence relation $\sim_{T}$. Let $\pi_{S}:S\rightarrow\overline{S}$ and
$\pi_{T}:T\rightarrow\overline{T}$ be the canonical projections of a set on
its quotient. Let $f:S\rightarrow T$ be a map. If $\overline{f}:\overline
{S}\rightarrow\overline{T}$ is a map for which the diagram%
\[
\xymatrix{
S \ar[r]^f \ar[d]_{\pi_S} & T \ar[d]^{\pi_T} \\
\overline{S} \ar[r]_{\overline{f}} & \overline{T}
}
\]
is commutative, then we say that \textquotedblleft the map $f$
\textit{descends} to the map $\overline{f}$\textquotedblright. It is easy to
see that there exists \textbf{at most one} map $\overline{f}:\overline
{S}\rightarrow\overline{T}$ such that the map $f$ descends to the map
$\overline{f}$ (for given $S$, $T$, $\sim_{S}$, $\sim_{T}$ and $f$). Moreover,
the existence of a map $\overline{f}:\overline{S}\rightarrow\overline{T}$ such
that the map $f$ descends to the map $\overline{f}$ is equivalent to the
statement that every two elements $x$ and $y$ of $S$ satisfying $x\sim_{S}y$
satisfy $f\left(  x\right)  \sim_{T}f\left(  y\right)  $.

The above statements are not literally true if we replace the map
$f:S\rightarrow T$ by a partial map $f:S\dashrightarrow T$. However, when $S$
and $T$ are two algebraic varieties and $\sim_{S}$ and $\sim_{T}$ are
algebraic equivalences (i.e., equivalence relations defined by polynomial
relations between coordinates of points) and $f:S\dashrightarrow T$ is a
\textbf{rational map}, then the above statements still are true (of course,
with $\overline{f}$ being a partial map).
\end{definition}

\Needspace{15\baselineskip}

\begin{definition}
\label{def.hgRi}Let $n\in\mathbb{N}$. Let $\mathbb{K}$ be a field. Let $P$ be
an $n$-graded poset. Let $i\in\left\{  1,2,...,n\right\}  $. Because of
Corollary \ref{cor.hgRi}, the rational map $R_{i}:\mathbb{K}^{\widehat{P}%
}\dashrightarrow\mathbb{K}^{\widehat{P}}$ descends (through the projection
$\pi:\mathbb{K}^{\widehat{P}}\dashrightarrow\overline{\mathbb{K}^{\widehat{P}%
}}$) to a partial map $\overline{\mathbb{K}^{\widehat{P}}}\dashrightarrow
\overline{\mathbb{K}^{\widehat{P}}}$. We denote this partial map
$\overline{\mathbb{K}^{\widehat{P}}}\dashrightarrow\overline{\mathbb{K}%
^{\widehat{P}}}$ by $\overline{R_{i}}$. Thus, the diagram%
\begin{equation}
\xymatrixcolsep{5pc}\xymatrix{ \mathbb K^{\widehat P} \ar@{-->}[r]^{R_i} \ar@{-->}[d]_-{\pi} & \mathbb K^{\widehat P} \ar@{-->}[d]^-{\pi} \\ \overline{\mathbb K^{\widehat P}} \ar@{-->}[r]_{\overline{R_i}} & \overline {\mathbb K^{\widehat P}} }
\label{def.hgRi.commut}%
\end{equation}
is commutative.
\end{definition}

\Needspace{15\baselineskip}

\begin{definition}
\label{def.hgR}Let $n\in\mathbb{N}$. Let $\mathbb{K}$ be a field. Let $P$ be
an $n$-graded poset. We define the partial map $\overline{R}:\overline
{\mathbb{K}^{\widehat{P}}}\dashrightarrow\overline{\mathbb{K}^{\widehat{P}}}$
by%
\[
\overline{R}=\overline{R_{1}}\circ\overline{R_{2}}\circ...\circ\overline
{R_{n}}.
\]
Then, the diagram%
\begin{equation}
\xymatrixcolsep{5pc}\xymatrix{ \mathbb K^{\widehat P} \ar@{-->}[r]^{R} \ar@{-->}[d]_-{\pi} & \mathbb K^{\widehat P} \ar@{-->}[d]^-{\pi} \\ \overline{\mathbb K^{\widehat P}} \ar@{-->}[r]_{\overline{R}} & \overline {\mathbb K^{\widehat P}} }
\label{def.hgR.commut}%
\end{equation}
is commutative\footnotemark. In other words, $\overline{R}$ is the partial map
$\overline{\mathbb{K}^{\widehat{P}}}\dashrightarrow\overline{\mathbb{K}%
^{\widehat{P}}}$ to which the partial map $R:\mathbb{K}^{\widehat{P}%
}\dashrightarrow\mathbb{K}^{\widehat{P}}$ descends (through the projection
$\pi:\mathbb{K}^{\widehat{P}}\dashrightarrow\overline{\mathbb{K}^{\widehat{P}%
}}$).
\end{definition}

\footnotetext{\textit{Proof.} We have $R=R_{1}\circ R_{2}\circ...\circ R_{n}$
and $\overline{R}=\overline{R_{1}}\circ\overline{R_{2}}\circ...\circ
\overline{R_{n}}$. Hence, the diagram (\ref{def.hgR.commut}) can be obtained
by stringing together the diagrams (\ref{def.hgRi.commut}) for all
$i\in\left\{  1,2,...,n\right\}  $ and then removing the ``interior edges''.
Therefore, the diagram (\ref{def.hgR.commut}) is commutative (since the
diagrams (\ref{def.hgRi.commut}) are commutative for all $i$), qed.}Next, we
formulate a result which says something to the extent of ``a zero-free
$\mathbb{K}$-labelling $f\in\mathbb{K}^{\widehat{P}}$ is almost always
uniquely determined by its w-tuple \newline$\left(  \mathbf{w}_{0}\left(
f\right)  ,\mathbf{w}_{1}\left(  f\right)  ,...,\mathbf{w}_{n}\left(
f\right)  \right)  $, its homogenization $\pi\left(  f\right)  $ and the value
$f\left(  0\right)  $''. The words ``almost always'' are required here because
otherwise the statement would be wrong; but they have to be made precise. Here
is the exact statement we want to make:

\begin{proposition}
\label{prop.reconstruct}Let $n\in\mathbb{N}$. Let $\mathbb{K}$ be a field. Let
$P$ be an $n$-graded poset. Let $f$ and $g$ be two zero-free $\mathbb{K}%
$-labellings in $\mathbb{K}^{\widehat{P}}$ such that $\left(  \mathbf{w}%
_{0}\left(  f\right)  ,\mathbf{w}_{1}\left(  f\right)  ,...,\mathbf{w}%
_{n}\left(  f\right)  \right)  =\left(  \mathbf{w}_{0}\left(  g\right)
,\mathbf{w}_{1}\left(  g\right)  ,...,\mathbf{w}_{n}\left(  g\right)  \right)
$ and such that no $i\in\left\{  0,1,...,n\right\}  $ satisfies $\mathbf{w}%
_{i}\left(  f\right)  =0$. Also assume that $\pi\left(  f\right)  =\pi\left(
g\right)  $ and $f\left(  0\right)  =g\left(  0\right)  $. Then, $f=g$.
\end{proposition}

Proposition \ref{prop.reconstruct} is easily proven by reconstructing $f$ and
$g$ \textquotedblleft bottom-up\textquotedblright\ along $\widehat{P}$.
Alternatively, we can prove Proposition \ref{prop.reconstruct} directly using
Proposition \ref{prop.w.scalmult}, as follows:

\begin{proof}
[Proof of Proposition \ref{prop.reconstruct} (sketched).]Since $\pi\left(
f\right)  =\pi\left(  g\right)  $, we know that $f$ and $g$ are homogeneously
equivalent. By Condition 2 in Definition \ref{def.hgeq} \textbf{(d)}, this
means that there exists an $\left(  n+2\right)  $-tuple $\left(  a_{0}%
,a_{1},...,a_{n+1}\right)  \in\left(  \mathbb{K}^{\times}\right)  ^{n+2}$ such
that every $x\in\widehat{P}$ satisfies $g\left(  x\right)  =a_{\deg x}\cdot
f\left(  x\right)  $. In other words, there exists an $\left(  n+2\right)
$-tuple $\left(  a_{0},a_{1},...,a_{n+1}\right)  \in\left(  \mathbb{K}%
^{\times}\right)  ^{n+2}$ such that
\[
g=\left(  a_{0},a_{1},...,a_{n+1}\right)  \flat f
\]
(where $\left(  a_{0},a_{1},...,a_{n+1}\right)  \flat f\in\mathbb{K}%
^{\widehat{P}}$ is defined as in Definition \ref{def.bemol}). Consider this
$\left(  n+2\right)  $-tuple $\left(  a_{0},a_{1},...,a_{n+1}\right)  $.

Since $g=\left(  a_{0},a_{1},...,a_{n+1}\right)  \flat f$, we know that
\begin{align*}
\left(  \text{the w-tuple of }g\right)   &  =\left(  \text{the w-tuple of
}\left(  a_{0},a_{1},...,a_{n+1}\right)  \flat f\right) \\
&  =\left(  \dfrac{a_{0}}{a_{1}}\mathbf{w}_{0}\left(  f\right)  ,\dfrac{a_{1}%
}{a_{2}}\mathbf{w}_{1}\left(  f\right)  ,...,\dfrac{a_{n}}{a_{n+1}}%
\mathbf{w}_{n}\left(  f\right)  \right)
\end{align*}
(by Proposition \ref{prop.w.scalmult}). Compared with%
\[
\left(  \text{the w-tuple of }g\right)  =\left(  \mathbf{w}_{0}\left(
g\right)  ,\mathbf{w}_{1}\left(  g\right)  ,...,\mathbf{w}_{n}\left(
g\right)  \right)  =\left(  \mathbf{w}_{0}\left(  f\right)  ,\mathbf{w}%
_{1}\left(  f\right)  ,...,\mathbf{w}_{n}\left(  f\right)  \right)  ,
\]
this yields%
\[
\left(  \dfrac{a_{0}}{a_{1}}\mathbf{w}_{0}\left(  f\right)  ,\dfrac{a_{1}%
}{a_{2}}\mathbf{w}_{1}\left(  f\right)  ,...,\dfrac{a_{n}}{a_{n+1}}%
\mathbf{w}_{n}\left(  f\right)  \right)  =\left(  \mathbf{w}_{0}\left(
f\right)  ,\mathbf{w}_{1}\left(  f\right)  ,...,\mathbf{w}_{n}\left(
f\right)  \right)  .
\]
In other words, $\dfrac{a_{i}}{a_{i+1}}\mathbf{w}_{i}\left(  f\right)
=\mathbf{w}_{i}\left(  f\right)  $ for every $i\in\left\{  0,1,...,n\right\}
$. Hence, $\dfrac{a_{i}}{a_{i+1}}=1$ for every $i\in\left\{
0,1,...,n\right\}  $ (here, we cancelled out $\mathbf{w}_{i}\left(  f\right)
$, because by assumption we don't have $\mathbf{w}_{i}\left(  f\right)  =0$).
In other words, $a_{i}=a_{i+1}$ for every $i\in\left\{  0,1,...,n\right\}  $.
Thus, $a_{0}=a_{1}=...=a_{n+1}$.

\begin{vershort}
But since $g=\left(  a_{0},a_{1},...,a_{n+1}\right)  \flat f$, we have
$g\left(  0\right)  =\left(  \left(  a_{0},a_{1},...,a_{n+1}\right)  \flat
f\right)  \left(  0\right)  =a_{\deg0}\cdot f\left(  0\right)  =a_{0}\cdot
f\left(  0\right)  $ (since $\deg0=0$), so that $f\left(  0\right)  =g\left(
0\right)  =a_{0}\cdot f\left(  0\right)  $. Since $f\left(  0\right)  \neq0$
(because $f$ is zero-free, and the only element of $\widehat{P}_{0}$ is $0$),
we can cancel $f\left(  0\right)  $ here and obtain $1=a_{0}$. In view of
this, $a_{0}=a_{1}=...=a_{n+1}$ becomes $a_{0}=a_{1}=...=a_{n+1}=1$. Thus,
$\left(  a_{0},a_{1},...,a_{n+1}\right)  =\left(  \underbrace{1,1,...,1}%
_{n+2\text{ times}}\right)  $, so that $g=\left(  a_{0},a_{1},...,a_{n+1}%
\right)  \flat f=\left(  \underbrace{1,1,...,1}_{n+2\text{ times}}\right)
\flat f=f$, proving Proposition \ref{prop.reconstruct}.
\end{vershort}

\begin{verlong}
But since $g=\left(  a_{0},a_{1},...,a_{n+1}\right)  \flat f$, we have
$g\left(  0\right)  =\left(  \left(  a_{0},a_{1},...,a_{n+1}\right)  \flat
f\right)  \left(  0\right)  =a_{\deg0}\cdot f\left(  0\right)  $ (by the
definition of $\left(  a_{0},a_{1},...,a_{n+1}\right)  \flat f$). Since
$\deg0=0$, this rewrites as $g\left(  0\right)  =a_{0}\cdot f\left(  0\right)
$. Thus, $f\left(  0\right)  =g\left(  0\right)  =a_{0}\cdot f\left(
0\right)  $. Since $f\left(  0\right)  \neq0$ (because $f$ is zero-free, and
the only element of $\widehat{P}_{0}$ is $0$), we can cancel $f\left(
0\right)  $ here and obtain $1=a_{0}$. In view of this, $a_{0}=a_{1}%
=...=a_{n+1}$ becomes $a_{0}=a_{1}=...=a_{n+1}=1$. Thus, $\left(  a_{0}%
,a_{1},...,a_{n+1}\right)  =\left(  \underbrace{1,1,...,1}_{n+2\text{ times}%
}\right)  $, so that%
\[
g=\underbrace{\left(  a_{0},a_{1},...,a_{n+1}\right)  }_{=\left(
\underbrace{1,1,...,1}_{n+2\text{ times}}\right)  }\flat f=\left(
\underbrace{1,1,...,1}_{n+2\text{ times}}\right)  \flat f=f.
\]
This proves Proposition \ref{prop.reconstruct}.
\end{verlong}
\end{proof}

\begin{definition}
Let $\mathbb{K}$ be a field. In the following, if $S$ is a finite set, and $q$
is an element of a projective space $\mathbb{P}\left(  \mathbb{K}^{S}\right)
$ of the free vector space with basis $S$, and $k$ is an integer, then $q^{k}$
will denote the element of $\mathbb{P}\left(  \mathbb{K}^{S}\right)  $
obtained by replacing every homogeneous coordinate of $q$ by its $k$-th power.
This is well-defined (and will mostly be used for $k=-1$). In particular, this
definition applies to $S=\left\{  1,2,...,n\right\}  $ for $n\in\mathbb{N}$
(in which case $\mathbb{K}^{S}=\mathbb{K}^{n}$).
\end{definition}

We can explicitly describe the action of the $\overline{R_{i}}$ when the
``structure of the poset $P$ between degrees $i-1$, $i$ and $i+1$'' is
particularly simple:

\begin{proposition}
\label{prop.hgRi.1}Let $n\in\mathbb{N}$. Let $\mathbb{K}$ be a field. Let $P$
be an $n$-graded poset. Fix $i\in\left\{  1,2,...,n\right\}  $. Assume that
every $u\in\widehat{P}_{i}$ and every $v\in\widehat{P}_{i+1}$ satisfy
$u\lessdot v$. Assume further that every $u\in\widehat{P}_{i-1}$ and every
$v\in\widehat{P}_{i}$ satisfy $u\lessdot v$. Let $f\in\mathbb{K}^{\widehat{P}%
}$. Then,%
\begin{align*}
&  \left(  \pi_{1}\left(  R_{i}f\right)  ,\pi_{2}\left(  R_{i}f\right)
,...,\pi_{n}\left(  R_{i}f\right)  \right) \\
&  =\left(  \pi_{1}\left(  f\right)  ,\pi_{2}\left(  f\right)  ,...,\pi
_{i-1}\left(  f\right)  ,\left(  \pi_{i}\left(  f\right)  \right)  ^{-1}%
,\pi_{i+1}\left(  f\right)  ,\pi_{i+2}\left(  f\right)  ,...,\pi_{n}\left(
f\right)  \right)  .
\end{align*}

\end{proposition}

From this proposition, we obtain two corollaries:

\begin{corollary}
\label{cor.hgRi.1}Let $n\in\mathbb{N}$. Let $\mathbb{K}$ be a field. Let $P$
be an $n$-graded poset. Fix $i\in\left\{  1,2,...,n\right\}  $. Assume that
every $u\in\widehat{P}_{i}$ and every $v\in\widehat{P}_{i+1}$ satisfy
$u\lessdot v$. Assume further that every $u\in\widehat{P}_{i-1}$ and every
$v\in\widehat{P}_{i}$ satisfy $u\lessdot v$. Let $\widetilde{f}=\left(
\widetilde{f}_{1},\widetilde{f}_{2},...,\widetilde{f}_{n}\right)  \in
\overline{\mathbb{K}^{\widehat{P}}}$. Then,%
\[
\overline{R_{i}}\left(  \widetilde{f}\right)  =\left(  \widetilde{f}%
_{1},\widetilde{f}_{2},...,\widetilde{f}_{i-1},\widetilde{f}_{i}%
^{-1},\widetilde{f}_{i+1},\widetilde{f}_{i+2},...,\widetilde{f}_{n}\right)  .
\]

\end{corollary}

\begin{corollary}
\label{cor.hgRi.2}Let $n\in\mathbb{N}$. Let $\mathbb{K}$ be a field. Let $P$
be an $n$-graded poset. Assume that, for every $i\in\left\{
1,2,...,n-1\right\}  $, every $u \in\widehat{P}_{i}$ and every $v
\in\widehat{P}_{i+1}$ satisfy $u \lessdot v$. Let $f\in\mathbb{K}%
^{\widehat{P}}$ be zero-free. Then,%
\[
\left(  \pi_{1}\left(  Rf\right)  ,\pi_{2}\left(  Rf\right)  ,...,\pi
_{n}\left(  Rf\right)  \right)  =\left(  \left(  \pi_{1}\left(  f\right)
\right)  ^{-1},\left(  \pi_{2}\left(  f\right)  \right)  ^{-1},...,\left(
\pi_{n}\left(  f\right)  \right)  ^{-1}\right)  .
\]

\end{corollary}

\section{\label{sect.ord}Order}

In this short section, we will relate the orders of the maps $R$ and
$\overline{R}$ for a graded poset $P$. The relation will later be used to gain
knowledge on both of these orders.

We begin by defining the order of a partial map:

\begin{definition}
\label{def.ord}Let $S$ be a set.

\textbf{(a)} If $\alpha$ and $\beta$ are two partial maps from the set $S$,
then we write \textquotedblleft$\alpha=\beta$\textquotedblright\ if and only
if every $s\in S$ for which both $\alpha\left(  s\right)  $ and $\beta\left(
s\right)  $ are well-defined satisfies $\alpha\left(  s\right)  =\beta\left(
s\right)  $. This is, per se, not a well-behaved notation (e.g., it is
possible that three partial maps $\alpha$, $\beta$ and $\gamma$ satisfy
$\alpha=\beta$ and $\beta=\gamma$ but not $\alpha=\gamma$). However, we are
going to use this notation for rational maps and their quotients (and, of
course, total maps) only; in all of these cases, the notation \textbf{is}
well-behaved (e.g., if $\alpha$, $\beta$ and $\gamma$ are three rational maps
satisfying $\alpha=\beta$ and $\beta=\gamma$, then $\alpha=\gamma$, because
the intersection of two Zariski-dense open subsets is Zariski-dense and open).

\textbf{(b)} The \textit{order} of a partial map $\varphi:S\dashrightarrow S$
is defined to be the smallest positive integer $k$ satisfying $\varphi
^{k}=\operatorname*{id}\nolimits_{S}$, if such a positive integer $k$ exists,
and $\infty$ otherwise. Here, we are disregarding the fact that $\varphi$ is
only a partial map; we will be working only with dominant rational maps and
their quotients (and total maps), so nothing will go wrong.

We denote the order of a partial map $\varphi:S\dashrightarrow S$ as
$\operatorname*{ord}\varphi$.
\end{definition}

\begin{convention}
In the following, we are going to occasionally make arithmetical statements
involving the symbol $\infty$. We declare that $0$ and $\infty$ are divisible
by $\infty$, but no positive integer is divisible by $\infty$. We further
declare that every positive integer (but not $0$) divides $\infty$. We set
$\operatorname{lcm}\left(  a,\infty\right)  $ and $\operatorname{lcm}\left(
\infty,a\right)  $ to mean $\infty$ whenever $a$ is a positive integer.
\end{convention}

As a consequence of Proposition \ref{prop.reconstruct}, we have:

\begin{proposition}
\label{prop.ord-projord}Let $n\in\mathbb{N}$. Let $\mathbb{K}$ be a field. Let
$P$ be an $n$-graded poset. Then, $\operatorname*{ord}R=\operatorname{lcm}%
\left(  n+1,\operatorname*{ord}\overline{R}\right)  $. (Recall that
$\operatorname{lcm}\left(  n+1,\infty\right)  $ is to be understood as
$\infty$.)
\end{proposition}

The proof of this boils down to considering the effect of $R$ on the w-tuple
\newline$\left(  \mathbf{w}_{0}\left(  f\right)  ,\mathbf{w}_{1}\left(
f\right)  ,...,\mathbf{w}_{n}\left(  f\right)  \right)  $ and on the
homogenization $\pi\left(  f\right)  $ of a $\mathbb{K}$-labelling $f$. The
effect on the w-tuple is a cyclic shift (by Proposition \ref{prop.wi.R}),
which has order $n+1$. The effect on the homogenization is $\overline{R}$. It
is now easy to see (invoking Proposition \ref{prop.reconstruct}) that the
order of $R$ is the $\operatorname{lcm}$ of the orders of these two actions.
Here are the details of this derivation:

\begin{proof}
[Proof of Proposition \ref{prop.ord-projord} (sketched).]%
\begin{verlong}
Consider the group $\mathfrak{S}$ defined in the proof of Proposition
\ref{prop.wi.R}, and its action on the $\left(  n+1\right)  $-tuples of
elements of $\mathbb{K}$ defined ibidem. Also, define the element $\left(
0,1,...,n\right)  $ of $\mathfrak{S}$ just as in the proof of Proposition
\ref{prop.wi.R}.
\end{verlong}

\textit{1st step:} The commutativity of the diagram (\ref{def.hgR.commut})
yields $\overline{R}\circ\pi=\pi\circ R$. Hence,
\begin{equation}
\text{every }\ell\in\mathbb{N}\text{ satisfies }\overline{R}^{\ell}\circ
\pi=\pi\circ R^{\ell} \label{pf.ord-projord.Rl}%
\end{equation}
(this is clear by induction over $\ell$). Thus, if some $\ell\in\mathbb{N}$
satisfies $R^{\ell}=\operatorname*{id}$, then it satisfies $\overline{R}%
^{\ell}=\operatorname*{id}$ as well\footnote{\textit{Proof.} Let $\ell
\in\mathbb{N}$ be such that $R^{\ell}=\operatorname*{id}$. Then, $\overline
{R}^{\ell}\circ\pi=\pi\circ\underbrace{R^{\ell}}_{=\operatorname*{id}}=\pi$.
Since $\pi$ is right-cancellable (since $\pi$ is surjective), this yields
$\overline{R}^{\ell}=\operatorname*{id}$, qed.}. Hence, $\operatorname*{ord}%
\overline{R}\mid\operatorname*{ord}R$ (recall that every positive integer
divides $\infty$, but only $0$ and $\infty$ are divisible by $\infty$). In
particular, if $\operatorname*{ord}\overline{R}=\infty$, then
$\operatorname*{ord}R=\infty$. Thus, Proposition \ref{prop.ord-projord} is
obvious in the case when $\operatorname*{ord}\overline{R}=\infty$. Hence, for
the rest of the proof of Proposition \ref{prop.ord-projord}, we can WLOG
assume that $\operatorname*{ord}\overline{R}\neq\infty$. Assume this.

\textit{2nd step:} Since $\operatorname*{ord}\overline{R}\neq\infty$, we know
that $\operatorname*{ord}\overline{R}$ is a positive integer. Let $m$ be this
positive integer. Then, $m=\operatorname*{ord}\overline{R}$, so that
$\overline{R}^{m}=\operatorname*{id}$.

Let $\ell=\operatorname{lcm}\left(  n+1,m\right)  $. Then, $n+1\mid\ell$ and
$m\mid\ell$. Since $\operatorname*{ord}\overline{R}=m\mid\ell$, we have
$\overline{R}^{\ell}=\operatorname*{id}$. But from (\ref{pf.ord-projord.Rl}),
we have $\pi\circ R^{\ell}=\underbrace{\overline{R}^{\ell}}%
_{=\operatorname*{id}}\circ\pi=\pi$.

We are now going to prove that $R^{\ell}=\operatorname*{id}$. In order to
prove this, it is clearly enough to show that almost every (in the sense of
Zariski topology) zero-free $\mathbb{K}$-labelling $f$ of $P$ satisfies
$R^{\ell}f=\operatorname*{id}f$ (because $R^{\ell}f=\operatorname*{id}f$ is a
polynomial identity in the labels of $f$). But it is easily shown that for
almost every (in the sense of Zariski topology) zero-free $\mathbb{K}%
$-labelling $f$ of $P$, the w-tuple $\left(  \mathbf{w}_{0}\left(  f\right)
,\mathbf{w}_{1}\left(  f\right)  ,...,\mathbf{w}_{n}\left(  f\right)  \right)
$ of $f$ consists of nonzero elements of $\mathbb{K}$.\ \ \ \ \footnotemark

\begin{vershort}
\footnotetext{\textit{Proof.} We will prove a slightly better result: Almost
every $f\in\mathbb{K}^{\widehat{P}}$ is a zero-free $\mathbb{K}$-labelling of
$P$ with the property that%
\begin{equation}
\left(  \mathbf{w}_{i}\left(  f\right)  \text{ is well-defined and nonzero for
every }i\in\left\{  0,1,...,n\right\}  \right)  .
\label{pf.ord-projord.short.step2.zari.1}%
\end{equation}
\par
In fact, the condition (\ref{pf.ord-projord.short.step2.zari.1}) on an
$f\in\mathbb{K}^{\widehat{P}}$ is a requirement saying that certain rational
expressions in the values of $f$ do not vanish (namely, the denominators in
$\mathbf{w}_{i}\left(  f\right)  $ and the sums $\mathbf{w}_{i}\left(
f\right)  $ themselves). If we can prove that none of these expressions is
identically zero, then it will follow that for almost every $f\in
\mathbb{K}^{\widehat{P}}$, none of these expressions vanishes (because there
are only finitely many expressions whose vanishing we are trying to avoid, and
the infiniteness of $\mathbb{K}$ allows us to avoid them all if none of them
is identically zero); thus (\ref{pf.ord-projord.short.step2.zari.1}) will
follow and we will be done. Hence, it remains to show that none of these
expressions is identically zero.
\par
Assume the contrary. Then, one of our rational expressions -- either a
denominator in one of the $\mathbf{w}_{i}\left(  f\right)  $, or one of the
sums $\mathbf{w}_{i}\left(  f\right)  $ -- identically vanishes. It must be
one of the sums $\mathbf{w}_{i}\left(  f\right)  $, since the denominators in
the $\mathbf{w}_{i}\left(  f\right)  $ cannot identically vanish (they are
simply values $f\left(  y\right)  $). So there exists some $i\in\left\{
0,1,...,n\right\}  $ such that every $\mathbb{K}$-labelling $f$ of $P$ (for
which $\mathbf{w}_{i}\left(  f\right)  $ is well-defined) satisfies
$\mathbf{w}_{i}\left(  f\right)  =0$. Consider this $i$. Notice that $i\leq n$
and thus $1\notin\widehat{P}_{i}$.
\par
We have
\[
0=\mathbf{w}_{i}\left(  f\right)  =\sum_{\substack{x\in\widehat{P}_{i}%
;\ y\in\widehat{P}_{i+1};\\y\gtrdot x}}\dfrac{f\left(  x\right)  }{f\left(
y\right)  }=\sum_{x\in\widehat{P}_{i}}f\left(  x\right)  \sum_{\substack{y\in
\widehat{P}_{i+1};\\y\gtrdot x}}\dfrac{1}{f\left(  y\right)  }.
\]
This forces the sum $\sum_{\substack{y\in\widehat{P}_{i+1};\\y\gtrdot
x}}\dfrac{1}{f\left(  y\right)  }$ to be identically $0$ for every
$x\in\widehat{P}_{i}$ (because these sums for different values of $x$ are
prevented from canceling each other by the completely independent $f\left(
x\right)  $ coefficients in front of them). Fix some $x\in\widehat{P}_{i}$
(such an $x$ clearly exists since $\deg:\widehat{P}\rightarrow\left\{
0,1,...,n+1\right\}  $ is surjective), and ponder what it means for the sum
$\sum_{\substack{y\in\widehat{P}_{i+1};\\y\gtrdot x}}\dfrac{1}{f\left(
y\right)  }$ to be identically $0$. It means that this sum is empty, i.e.,
that there exists no $y\in\widehat{P}_{i+1}$ satisfying $y\gtrdot x$. But this
can only happen when $x=1$, which is not the case in our situation (because
$x\in\widehat{P}_{i}$ and $1\notin\widehat{P}_{i}$). So we have obtained a
contradiction.}
\end{vershort}

\begin{verlong}
\footnotetext{\textit{Proof.} Let us first fix an $i\in\left\{
0,1,...,n\right\}  $. Recall that $\mathbf{w}_{i}\left(  f\right)
=\sum_{\substack{x\in\widehat{P}_{i};\ y\in\widehat{P}_{i+1};\\y\gtrdot
x}}\dfrac{f\left(  x\right)  }{f\left(  y\right)  }$ for any zero-free
$\mathbb{K}$-labelling $f$ of $P$. We will now construct a $\mathbb{K}%
$-labelling $f$ of $P$ for which $\mathbf{w}_{i}\left(  f\right)  $ is
well-defined and nonzero.
\par
The map $\deg:P\rightarrow\left\{  1,2,...,n\right\}  $ is surjective. Hence,
its extension $\deg:\widehat{P}\rightarrow\left\{  0,1,...,n+1\right\}  $ is
surjective as well (since $\deg0=0$ and $\deg1=n+1$). As a consequence, there
exists some $x_{0}\in\widehat{P}$ satisfying $\deg\left(  x_{0}\right)  =i$.
Pick such an $x_{0}$. We have $\deg\left(  x_{0}\right)  =i\neq n+1$ (since
$i\in\left\{  0,1,...,n\right\}  $), and thus $\deg\left(  x_{0}\right)  \neq
n+1=\deg1$. Hence, $x_{0}\neq1$.
\par
Corollary \ref{cor.Phat.exist} \textbf{(b)} (applied to $x_{0}$ instead of
$q$) now yields that there exists a $u\in\widehat{P}$ such that $u\gtrdot
x_{0}$ in $\widehat{P}$. Denote this $u$ by $y_{0}$. Thus, $y_{0}$ is an
element of $\widehat{P}$ and satisfies $y_{0}\gtrdot x_{0}$.
\par
We have $y_{0}\gtrdot x_{0}$, hence $x_{0}\lessdot y_{0}$ and thus
$\deg\left(  x_{0}\right)  =\deg\left(  y_{0}\right)  -1$ (by Proposition
\ref{prop.graded.Phat.why} \textbf{(a)}, applied to $u=x_{0}$ and $v=y_{0}$).
Thus, $\deg\left(  y_{0}\right)  =\underbrace{\deg\left(  x_{0}\right)  }%
_{=i}+1=i+1$, so that $y_{0}\in\deg^{-1}\left(  \left\{  i+1\right\}  \right)
=\widehat{P}_{i+1}$. Thus, there exists at least one $y\in\widehat{P}_{i+1}$
satisfying $y\gtrdot x_{0}$ (namely, $y=y_{0}$).
\par
Also, $\deg\left(  x_{0}\right)  =i$, so that $x_{0}\in\deg^{-1}\left(
\left\{  i\right\}  \right)  =\widehat{P}_{i}$.
\par
For any $\mathbb{K}$-labelling $f$ of $P$, let $\mathbf{u}\left(  f\right)  $
denote the element $\sum_{\substack{y\in\widehat{P}_{i+1};\\y\gtrdot x_{0}%
}}\dfrac{1}{f\left(  y\right)  }$ of $\mathbb{K}$. Also, for any $\mathbb{K}%
$-labelling $f$ of $P$, let $\mathbf{v}\left(  f\right)  $ denote the element
$\sum_{\substack{x\in\widehat{P}_{i};\\x\neq x_{0}}}\sum_{\substack{y\in
\widehat{P}_{i+1};\\y\gtrdot x}}\dfrac{f\left(  x\right)  }{f\left(  y\right)
}$ of $\mathbb{K}$. Then, every $\mathbb{K}$-labelling $f$ of $P$ satisfies%
\begin{align*}
\mathbf{w}_{i}\left(  f\right)   &  =\underbrace{\sum_{\substack{x\in
\widehat{P}_{i};\ y\in\widehat{P}_{i+1};\\y\gtrdot x}}}_{=\sum_{x\in
\widehat{P}_{i}}\sum_{\substack{y\in\widehat{P}_{i+1};\\y\gtrdot x}}}%
\dfrac{f\left(  x\right)  }{f\left(  y\right)  }=\sum_{x\in\widehat{P}_{i}%
}\sum_{\substack{y\in\widehat{P}_{i+1};\\y\gtrdot x}}\dfrac{f\left(  x\right)
}{f\left(  y\right)  }=\underbrace{\sum_{\substack{x\in\widehat{P}_{i};\\x\neq
x_{0}}}\sum_{\substack{y\in\widehat{P}_{i+1};\\y\gtrdot x}}\dfrac{f\left(
x\right)  }{f\left(  y\right)  }}_{\substack{=\mathbf{v}\left(  f\right)
\\\text{(by the definition of }\mathbf{v}\left(  f\right)  \text{)}%
}}+\underbrace{\sum_{\substack{y\in\widehat{P}_{i+1};\\y\gtrdot x_{0}}%
}\dfrac{f\left(  x_{0}\right)  }{f\left(  y\right)  }}_{=f\left(
x_{0}\right)  \cdot\sum_{\substack{y\in\widehat{P}_{i+1};\\y\gtrdot x_{0}%
}}\dfrac{1}{f\left(  y\right)  }}\\
&  \ \ \ \ \ \ \ \ \ \ \left(  \text{here, we have split off the addend for
}x=x_{0}\text{ from the outer sum}\right) \\
&  =\mathbf{v}\left(  f\right)  +f\left(  x_{0}\right)  \cdot\underbrace{\sum
_{\substack{y\in\widehat{P}_{i+1};\\y\gtrdot x_{0}}}\dfrac{1}{f\left(
y\right)  }}_{\substack{=\mathbf{u}\left(  f\right)  \\\text{(by the
definition of }\mathbf{u}\left(  f\right)  \text{)}}}=\mathbf{v}\left(
f\right)  +f\left(  x_{0}\right)  \cdot\mathbf{u}\left(  f\right)  .
\end{align*}
\par
We now construct a $\mathbb{K}$-labelling $f$ of $P$ as follows:
\par
\begin{itemize}
\item We choose the values $f\left(  y\right)  \in\mathbb{K}$ for all
$y\in\widehat{P}_{i+1}$ satisfying $y\gtrdot x_{0}$ in such a way that these
values $f\left(  y\right)  $ are nonzero and satisfy $\sum_{\substack{y\in
\widehat{P}_{i+1};\\y\gtrdot x_{0}}}\dfrac{1}{f\left(  y\right)  }\neq0$.
(This is clearly possible because there exists at least one $y\in
\widehat{P}_{i+1}$ satisfying $y\gtrdot x_{0}$.) As a result, $\mathbf{u}%
\left(  f\right)  =\sum_{\substack{y\in\widehat{P}_{i+1};\\y\gtrdot x_{0}%
}}\dfrac{1}{f\left(  y\right)  }$ is well-defined and nonzero.
\par
\item We choose the values $f\left(  y\right)  \in\mathbb{K}$ for all
remaining $y\in\widehat{P}$ in such a way that they are nonzero (but can
otherwise be arbitrary -- for example, we can set them to be $1$), with
only one exception: We do not yet decided what value
$f\left(  x_{0}\right)  $ will have. As a result, $\mathbf{v}%
\left(  f\right)  =\sum_{\substack{x\in\widehat{P}_{i};\\x\neq x_{0}}%
}\sum_{\substack{y\in\widehat{P}_{i+1};\\y\gtrdot x}}\dfrac{f\left(  x\right)
}{f\left(  y\right)  }$ is well-defined (because any $x\in\widehat{P}_{i}$
satisfying $x\neq x_{0}$ satisfies $x\neq x_{0}$, and because every
$y\in\widehat{P}_{i+1}$ satisfies $y\neq x_{0}$).
\par
\item Finally, it remains to choose $f\left(  x_{0}\right)  \in\mathbb{K}$. We
choose $f\left(  x_{0}\right)  $ in such a way that $\mathbf{v}\left(
f\right)  +f\left(  x_{0}\right)  \cdot\mathbf{u}\left(  f\right)  =1$. (This
is possible because $\mathbf{u}\left(  f\right)  $ and $\mathbf{v}\left(
f\right)  $ are two already well-defined elements of $\mathbb{K}$, and because
$\mathbf{u}\left(  f\right)  $ is nonzero.)
\end{itemize}
\par
Thus, we have chosen a value $f\left(  y\right)  \in\mathbb{K}$ for every
$y\in\widehat{P}$. In other words, we have constructed a $\mathbb{K}%
$-labelling $f$ of $P$. This $\mathbb{K}$-labelling has the property that
$\mathbf{w}_{i}\left(  f\right)  =\mathbf{v}\left(  f\right)  +f\left(
x_{0}\right)  \cdot\mathbf{u}\left(  f\right)  \neq0$, so that $\mathbf{w}%
_{i}\left(  f\right)  $ is well-defined and nonzero.
\par
Let us now forget that we defined $f$. We thus have constructed a $\mathbb{K}%
$-labelling $f$ of $P$ for which $\mathbf{w}_{i}\left(  f\right)  $ is
well-defined and nonzero. Hence, there exists a $\mathbb{K}$-labelling $f$ of
$P$ for which $\mathbf{w}_{i}\left(  f\right)  $ is well-defined and nonzero.
As a consequence, the set
\[
\left\{  f\text{ is a }\mathbb{K}\text{-labelling of }P\ \mid\ \mathbf{w}%
_{i}\left(  f\right)  \text{ is well-defined and nonzero}\right\}
\]
is nonempty. Therefore, this set is a dense open subset of $\mathbb{K}%
^{\widehat{P}}$ in the Zariski topology (because it is clearly an open subset
of $\mathbb{K}^{\widehat{P}}$ (since it is characterized by the nonvanishing
of certain rational expressions in the values of $f$)).
\par
Let us now forget that we fixed $i$. We thus have shown that for every
$i\in\left\{  0,1,...,n\right\}  $, the set
\[
\left\{  f\text{ is a }\mathbb{K}\text{-labelling of }P\ \mid\ \mathbf{w}%
_{i}\left(  f\right)  \text{ is well-defined and nonzero}\right\}
\]
is a dense open subset of $\mathbb{K}^{\widehat{P}}$ in the Zariski topology.
Hence, the intersection%
\[
\bigcap_{i\in\left\{  0,1,...,n\right\}  }\left\{  f\text{ is a }%
\mathbb{K}\text{-labelling of }P\ \mid\ \mathbf{w}_{i}\left(  f\right)  \text{
is well-defined and nonzero}\right\}
\]
of all these subsets is also a dense open subset of $\mathbb{K}^{\widehat{P}}$
in the Zariski topology (since the intersection of finitely many dense open
subsets in a topology must always be a dense open subset). Since
\begin{align*}
&  \bigcap_{i\in\left\{  0,1,...,n\right\}  }\left\{  f\text{ is a }%
\mathbb{K}\text{-labelling of }P\ \mid\ \mathbf{w}_{i}\left(  f\right)  \text{
is well-defined and nonzero}\right\} \\
&  =\left\{  f\text{ is a }\mathbb{K}\text{-labelling of }P\ \mid
\ \mathbf{w}_{i}\left(  f\right)  \text{ is well-defined and nonzero for every
}i\in\left\{  0,1,...,n\right\}  \right\}  ,
\end{align*}
this rewrites as follows: The set $\left\{  f\text{ is a }\mathbb{K}%
\text{-labelling of }P\ \mid\ \mathbf{w}_{i}\left(  f\right)  \text{ is
well-defined and nonzero for every }i\in\left\{  0,1,...,n\right\}  \right\}
$ is a dense open subset of $\mathbb{K}^{\widehat{P}}$ in the Zariski
topology.
\par
Now, the set%
\begin{align*}
&  \left\{  f\text{ is a zero-free $\mathbb{K}$-labelling of }P\text{\ $\mid
$\ }\mathbf{w}\text{$_{i}\left(  f\right)  $ is well-defined and nonzero for
every }i\in\left\{  0,1,...,n\right\}  \right\} \\
&  =\left\{  f\text{ is a }\mathbb{K}\text{-labelling of }P\ \mid
\ \mathbf{w}_{i}\left(  f\right)  \text{ is well-defined and nonzero for every
}i\in\left\{  0,1,...,n\right\}  \right\} \\
&  \ \ \ \ \ \ \ \ \ \ \cap\left\{  f\text{ is a }\mathbb{K}\text{-labelling
of }P\ \mid\ f\text{ is zero-free}\right\}
\end{align*}
is a dense open subset of $\mathbb{K}^{\widehat{P}}$ in the Zariski topology
(because it is the intersection of the sets $\left\{  f\text{ is a }%
\mathbb{K}\text{-labelling of }P\ \mid\ \mathbf{w}_{i}\left(  f\right)  \text{
is well-defined and nonzero for every }i\in\left\{  0,1,...,n\right\}
\right\}  $ and $\left\{  f\text{ is a }\mathbb{K}\text{-labelling of }%
P\ \mid\ f\text{ is zero-free}\right\}  $, each of which is a dense open
subset of $\mathbb{K}^{\widehat{P}}$). In other words, almost every
$f\in\mathbb{K}^{\widehat{P}}$ is a zero-free $\mathbb{K}$-labelling of $P$
with the property that $\left(  \mathbf{w}_{i}\left(  f\right)  \text{ is
well-defined and nonzero for every }i\in\left\{  0,1,...,n\right\}  \right)
$. As a consequence, almost every zero-free $\mathbb{K}$-labelling of $P$ has
the property that $\left(  \mathbf{w}_{i}\left(  f\right)  \text{ is
well-defined and nonzero for every }i\in\left\{  0,1,...,n\right\}  \right)
$. In other words, for almost every (in the sense of Zariski topology)
zero-free $\mathbb{K}$-labelling $f$ of $P$, the w-tuple $\left(
\mathbf{w}_{0}\left(  f\right)  ,\mathbf{w}_{1}\left(  f\right)
,...,\mathbf{w}_{n}\left(  f\right)  \right)  $ of $f$ consists of nonzero
elements of $\mathbb{K}$. Qed.}
\end{verlong}

Hence, in order to prove $R^{\ell}=\operatorname*{id}$, it is enough to show
that every zero-free $\mathbb{K}$-labelling $f$ of $P$ for which the w-tuple
$\left(  \mathbf{w}_{0}\left(  f\right)  ,\mathbf{w}_{1}\left(  f\right)
,...,\mathbf{w}_{n}\left(  f\right)  \right)  $ of $f$ consists of nonzero
elements of $\mathbb{K}$ satisfies $R^{\ell}f=\operatorname*{id}f$. This is
what we are going to do now.

So let $f$ be a zero-free $\mathbb{K}$-labelling of $P$ for which the w-tuple
$\left(  \mathbf{w}_{0}\left(  f\right)  ,\mathbf{w}_{1}\left(  f\right)
,...,\mathbf{w}_{n}\left(  f\right)  \right)  $ of $f$ consists of nonzero
elements of $\mathbb{K}$. We will prove that $R^{\ell}f=\operatorname*{id}f$.

\begin{vershort}
From Proposition \ref{prop.wi.R}, we know that the map $R$ changes the w-tuple
of a $\mathbb{K}$-labelling by shifting it cyclically. Hence, for every
$k\in\mathbb{N}$, the map $R^{k}$ changes the w-tuple of a $\mathbb{K}%
$-labelling by shifting it cyclically $k$ times. If this $k$ is divisible by
$n+1$, then this obviously means that the map $R^{k}$ preserves the w-tuple of
a $\mathbb{K}$-labelling (because the w-tuple has $n+1$ entries, and thus
shifting it cyclically for a multiple of $n+1$ times leaves it invariant).
Hence, the w-tuple of $f$ equals the w-tuple of $R^{\ell}f$. Recalling the
definition of a w-tuple, we can rewrite this as follows:%
\[
\left(  \mathbf{w}_{0}\left(  f\right)  ,\mathbf{w}_{1}\left(  f\right)
,...,\mathbf{w}_{n}\left(  f\right)  \right)  =\left(  \mathbf{w}_{0}\left(
R^{\ell}f\right)  ,\mathbf{w}_{1}\left(  R^{\ell}f\right)  ,...,\mathbf{w}%
_{n}\left(  R^{\ell}f\right)  \right)  .
\]

\end{vershort}

\begin{verlong}
From Proposition \ref{prop.wi.R}, we know that the map $R$ changes the w-tuple
of a $\mathbb{K}$-labelling by shifting it cyclically. In other words, almost
every $\mathbb{K}$-labelling $g$ of $P$ satisfies%
\[
\left(  \text{the w-tuple of }Rg\right)  =\left(  0,1,...,n\right)  \left(
\text{the w-tuple of }g\right)  .
\]
From this, it is easy to see that for every $k\in\mathbb{N}$, almost every
$\mathbb{K}$-labelling $g$ of $P$ satisfies%
\begin{equation}
\left(  \text{the w-tuple of }R^{k}g\right)  =\left(  0,1,...,n\right)
^{k}\left(  \text{the w-tuple of }g\right)  \label{pf.ord-projord.cycle}%
\end{equation}
(indeed, this can be proven by induction over $k$). Applying this to $g=f$ and
$k=\ell$, we obtain%
\begin{align*}
\left(  \text{the w-tuple of }R^{\ell}f\right)   &  =\underbrace{\left(
0,1,...,n\right)  ^{\ell}}_{\substack{=\operatorname*{id}\\\text{(since
}\left(  0,1,...,n\right)  \text{ has order }n+1\text{ in }\mathfrak{S}%
\text{,}\\\text{and since }n+1\mid\ell\text{)}}}\left(  \text{the w-tuple of
}f\right) \\
&  =\left(  \text{the w-tuple of }f\right)  =\left(  \mathbf{w}_{0}\left(
f\right)  ,\mathbf{w}_{1}\left(  f\right)  ,...,\mathbf{w}_{n}\left(
f\right)  \right)
\end{align*}
(by the definition of the w-tuple). Compared with%
\[
\left(  \text{the w-tuple of }R^{\ell}f\right)  =\left(  \mathbf{w}_{0}\left(
R^{\ell}f\right)  ,\mathbf{w}_{1}\left(  R^{\ell}f\right)  ,...,\mathbf{w}%
_{n}\left(  R^{\ell}f\right)  \right)
\]
(by the definition of the w-tuple), we obtain%
\[
\left(  \mathbf{w}_{0}\left(  f\right)  ,\mathbf{w}_{1}\left(  f\right)
,...,\mathbf{w}_{n}\left(  f\right)  \right)  =\left(  \mathbf{w}_{0}\left(
R^{\ell}f\right)  ,\mathbf{w}_{1}\left(  R^{\ell}f\right)  ,...,\mathbf{w}%
_{n}\left(  R^{\ell}f\right)  \right)  .
\]

\end{verlong}

Moreover, by assumption, the w-tuple $\left(  \mathbf{w}_{0}\left(  f\right)
,\mathbf{w}_{1}\left(  f\right)  ,...,\mathbf{w}_{n}\left(  f\right)  \right)
$ of $f$ consists of nonzero elements of $\mathbb{K}$. In other words, no
$i\in\left\{  0,1,...,n\right\}  $ satisfies $\mathbf{w}_{i}\left(  f\right)
=0$.

Furthermore $\pi\left(  R^{\ell}f\right)  =\underbrace{\left(  \pi\circ
R^{\ell}\right)  }_{=\pi}f=\pi\left(  f\right)  $.

Also, Corollary \ref{cor.R.implicit.01} (applied to $k=\ell$) yields $\left(
R^{\ell}f\right)  \left(  0\right)  =f\left(  0\right)  $.

We now can apply Proposition \ref{prop.reconstruct} to $g=R^{\ell}f$. As a
result, we obtain $R^{\ell}f=f$. In other words, $R^{\ell}f=\operatorname*{id}%
f$.

Now forget that we fixed $f$. We have thus shown that $R^{\ell}%
f=\operatorname*{id}f$ for every zero-free $\mathbb{K}$-labelling $f$ of $P$
for which the w-tuple $\left(  \mathbf{w}_{0}\left(  f\right)  ,\mathbf{w}%
_{1}\left(  f\right)  ,...,\mathbf{w}_{n}\left(  f\right)  \right)  $ of $f$
consists of nonzero elements of $\mathbb{K}$. Therefore, we have shown that
$R^{\ell}=\operatorname*{id}$ (by what we have said above). Thus,
$\operatorname*{ord}R\mid\ell=\operatorname{lcm}\left(  n+1,\underbrace{m}%
_{=\operatorname*{ord}\overline{R}}\right)  =\operatorname{lcm}\left(
n+1,\operatorname*{ord}\overline{R}\right)  $.

\textit{3rd step:} We now will show that $\operatorname{lcm}\left(
n+1,\operatorname*{ord}\overline{R}\right)  \mid\operatorname*{ord}R$.

In order to do that, we assume WLOG that $\operatorname*{ord}R\neq\infty$
(because otherwise, \newline$\operatorname{lcm}\left(  n+1,\operatorname*{ord}%
\overline{R}\right)  \mid\operatorname*{ord}R$ is obvious). Hence,
$\operatorname*{ord}R$ is a positive integer. Denote this positive integer by
$q$. So, $q=\operatorname*{ord}R$.

It is easy to see that for almost every (in the sense of Zariski topology)
zero-free $\mathbb{K}$-labelling $f$ of $P$, the entries of the w-tuple
$\left(  \mathbf{w}_{0}\left(  f\right)  ,\mathbf{w}_{1}\left(  f\right)
,...,\mathbf{w}_{n}\left(  f\right)  \right)  $ of $f$ are pairwise distinct.
Hence, there exists a zero-free $\mathbb{K}$-labelling $f$ of $P$ such that
the entries of the w-tuple $\left(  \mathbf{w}_{0}\left(  f\right)
,\mathbf{w}_{1}\left(  f\right)  ,...,\mathbf{w}_{n}\left(  f\right)  \right)
$ of $f$ are pairwise distinct and such that $R^{k}f$ is well-defined for all
$k\in\left\{  0,1,...,q\right\}  $. Consider such an $f$.

Since $q=\operatorname*{ord}R$, we have $R^{q}=\operatorname*{id}$, so that
$R^{q}f=f$.

\begin{vershort}
Recall once again (from the 2nd step) that for every $k\in\mathbb{N}$, the map
$R^{k}$ changes the w-tuple of a $\mathbb{K}$-labelling by shifting it
cyclically $k$ times. In particular, the map $R^{q}$ changes the w-tuple of
the $\mathbb{K}$-labelling $f$ by shifting it cyclically $q$ times. In other
words, the w-tuple of $R^{q}f$ is obtained from the w-tuple of $f$ by shifting
it cyclically $q$ times. Since $R^{q}f=f$, this rewrites as follows: The
w-tuple of $f$ is obtained from the w-tuple of $f$ by shifting it cyclically
$q$ times. In other words, the w-tuple of $f$ is fixed under a $q$-fold cyclic
shift. But since the w-tuple of $f$ is an $\left(  n+1\right)  $-tuple of
pairwise distinct entries, this can only happen if $n+1\mid q$. Hence, we have
$n+1\mid q$.
\end{vershort}

\begin{verlong}
Applying (\ref{pf.ord-projord.cycle}) to $g=f$ and $k=q$, we obtain%
\[
\left(  \text{the w-tuple of }R^{q}f\right)  =\left(  0,1,...,n\right)
^{q}\left(  \text{the w-tuple of }f\right)  .
\]
Since $R^{q}f=f$, this rewrites as
\[
\left(  \text{the w-tuple of }f\right)  =\left(  0,1,...,n\right)  ^{q}\left(
\text{the w-tuple of }f\right)  .
\]
Since $\left(  \text{the w-tuple of }f\right)  =\left(  \mathbf{w}_{0}\left(
f\right)  ,\mathbf{w}_{1}\left(  f\right)  ,...,\mathbf{w}_{n}\left(
f\right)  \right)  $, this rewrites as%
\begin{equation}
\left(  \mathbf{w}_{0}\left(  f\right)  ,\mathbf{w}_{1}\left(  f\right)
,...,\mathbf{w}_{n}\left(  f\right)  \right)  =\left(  0,1,...,n\right)
^{q}\left(  \mathbf{w}_{0}\left(  f\right)  ,\mathbf{w}_{1}\left(  f\right)
,...,\mathbf{w}_{n}\left(  f\right)  \right)  . \label{pf.ord-projord.3cycle}%
\end{equation}
But since the entries of the w-tuple $\left(  \mathbf{w}_{0}\left(  f\right)
,\mathbf{w}_{1}\left(  f\right)  ,...,\mathbf{w}_{n}\left(  f\right)  \right)
$ are pairwise distinct, the only permutation in $\mathfrak{S}$ which fixes
$\left(  \mathbf{w}_{0}\left(  f\right)  ,\mathbf{w}_{1}\left(  f\right)
,...,\mathbf{w}_{n}\left(  f\right)  \right)  $ is the identity permutation.
Hence, (\ref{pf.ord-projord.3cycle}) shows that $\left(  0,1,...,n\right)
^{q}$ is the identity permutation. Since the order of $\left(
0,1,...,n\right)  $ in $\mathfrak{S}$ is $n+1$, this yields that $n+1\mid q$.
\end{verlong}

Combining $n+1\mid q=\operatorname*{ord}R$ with $\operatorname*{ord}%
\overline{R}\mid\operatorname*{ord}R$, we obtain $\operatorname{lcm}\left(
n+1,\operatorname*{ord}\overline{R}\right)  \mid\operatorname*{ord}R$.
Combining this with $\operatorname*{ord}R\mid\operatorname{lcm}\left(
n+1,\operatorname*{ord}\overline{R}\right)  $, we obtain $\operatorname*{ord}%
R=\operatorname{lcm}\left(  n+1,\operatorname*{ord}\overline{R}\right)  $.
This proves Proposition \ref{prop.ord-projord}.
\end{proof}

\section{The opposite poset}

Before we move on to the first interesting class of posets for which we can
compute the order of birational rowmotion, let us prove an easy ``symmetry
property'' of birational rowmotion.

\begin{definition}
\label{def.op}Let $P$ be a poset. Then, $P^{\operatorname*{op}}$ will denote
the poset defined on the same ground set as $P$ but with the order relation
defined by%
\[
\left(  \left(  a<_{P^{\operatorname*{op}}}b\text{ if and only if }%
b<_{P}a\right)  \text{ for all }a\in P\text{ and }b\in P\right)
\]
(where $<_{P}$ denotes the smaller relation of the poset $P$, and where
$<_{P^{\operatorname*{op}}}$ denotes the smaller relation of the poset
$P^{\operatorname*{op}}$ which we are defining). The poset
$P^{\operatorname*{op}}$ is called the \textit{opposite poset} of $P$.
\end{definition}

Note that $P^{\operatorname*{op}}$ is called the \textit{dual} of the poset
$P$ in \cite{stanley-ec1}.

\begin{remark}
It is clear that $\left(  P^{\operatorname*{op}}\right)  ^{\operatorname*{op}%
}=P$ for any poset $P$. Also, if $n\in\mathbb{N}$, and if $P$ is an $n$-graded
poset, then $P^{\operatorname*{op}}$ is an $n$-graded poset.
\end{remark}

\begin{definition}
Let $P$ be a finite poset. Let $\mathbb{K}$ be a field. We denote the maps $R$
and $\overline{R}$ by $R_{P}$ and $\overline{R}_{P}$, respectively, so as to
make their dependence on $P$ explicit.
\end{definition}

We can now state a symmetry property of $\operatorname*{ord}R$ (as defined in
Definition \ref{def.ord}):

\begin{proposition}
\label{prop.op.ord}Let $P$ be a finite poset. Let $\mathbb{K}$ be a field.
Then, $\operatorname*{ord}\left(  R_{P^{\operatorname*{op}}}\right)
=\operatorname*{ord}\left(  R_{P}\right)  $ and $\operatorname*{ord}\left(
\overline{R}_{P^{\operatorname*{op}}}\right)  =\operatorname*{ord}\left(
\overline{R}_{P}\right)  $.
\end{proposition}

\begin{proof}
[Proof of Proposition \ref{prop.op.ord} (sketched).]Define a rational map
$\kappa:\mathbb{K}^{\widehat{P}}\dashrightarrow\mathbb{K}%
^{\widehat{P^{\operatorname*{op}}}}$ by%
\[
\left(  \kappa f\right)  \left(  w\right)  =\left\{
\begin{array}
[c]{c}%
\dfrac{1}{f\left(  w\right)  },\ \ \ \ \ \ \ \ \ \ \text{if }w\in P;\\
\dfrac{1}{f\left(  1\right)  },\ \ \ \ \ \ \ \ \ \ \text{if }w=0;\\
\dfrac{1}{f\left(  0\right)  },\ \ \ \ \ \ \ \ \ \ \text{if }w=1
\end{array}
\right.  \ \ \ \ \ \ \ \ \ \ \text{for every }w\in
\widehat{P^{\operatorname*{op}}}\text{ for every }f\in\mathbb{K}
^{\widehat{P}}.
\]
This map $\kappa$ is a birational map. (Its inverse map is defined in the same way.)

We claim that $\kappa\circ R_{P}=R_{P^{\operatorname*{op}}}^{-1}\circ\kappa$.

Indeed, it is easy to see (by computation) that every element $v\in P$
satisfies
\begin{equation}
\kappa\circ T_{v}=T_{v}\circ\kappa, \label{pf.op.ord.1}%
\end{equation}
where the $T_{v}$ on the left hand side is defined with respect to the poset
$P$, and the $T_{v}$ on the right hand side is defined with respect to the
poset $P^{\operatorname*{op}}$. Now, let $\left(  v_{1},v_{2},...,v_{m}%
\right)  $ be a linear extension of $P$. Then, $\left(  v_{m},v_{m-1}%
,...,v_{1}\right)  $ is a linear extension of $P^{\operatorname*{op}}$, so
that Proposition \ref{prop.R.inverse} (applied to $P^{\operatorname*{op}}$ and
$\left(  v_{m},v_{m-1},...,v_{1}\right)  $ instead of $P$ and $\left(
v_{1},v_{2},...,v_{m}\right)  $) yields that $R_{P^{\operatorname*{op}}}%
^{-1}=T_{v_{1}}\circ T_{v_{2}}\circ...\circ T_{v_{m}}:\mathbb{K}%
^{\widehat{P^{\operatorname*{op}}}}\dashrightarrow\mathbb{K}%
^{\widehat{P^{\operatorname*{op}}}}$. On the other hand, the definition of
$R_{P}$ yields $R_{P}=T_{v_{1}}\circ T_{v_{2}}\circ...\circ T_{v_{m}%
}:\mathbb{K}^{\widehat{P}}\dashrightarrow\mathbb{K}^{\widehat{P}}$. Now, using
(\ref{pf.op.ord.1}), it is easy to see that%
\[
\kappa\circ\left(  T_{v_{1}}\circ T_{v_{2}}\circ...\circ T_{v_{m}}\right)
=\left(  T_{v_{1}}\circ T_{v_{2}}\circ...\circ T_{v_{m}}\right)  \circ\kappa.
\]
Since the $T_{v_{1}}\circ T_{v_{2}}\circ...\circ T_{v_{m}}$ on the left hand
side equals $R_{P}$, and the $T_{v_{1}}\circ T_{v_{2}}\circ...\circ T_{v_{m}}$
on the right hand side equals $R_{P^{\operatorname*{op}}}^{-1}$, this rewrites
as $\kappa\circ R_{P}=R_{P^{\operatorname*{op}}}^{-1}\circ\kappa$. Since
$\kappa$ is a birational map, this shows that $R_{P}$ and
$R_{P^{\operatorname*{op}}}^{-1}$ are birationally equivalent, so that
$\operatorname*{ord}\left(  R_{P}\right)  =\operatorname*{ord}\left(
R_{P^{\operatorname*{op}}}^{-1}\right)  =\operatorname*{ord}\left(
R_{P^{\operatorname*{op}}}\right)  $. Since $\kappa$ commutes with
homogenization, we also obtain the birational equivalence of the maps
$\overline{R}_{P}$ and $\overline{R}_{P^{\operatorname*{op}}}^{-1}$, whence
$\operatorname*{ord}\left(  \overline{R}_{P}\right)  =\operatorname*{ord}%
\left(  \overline{R}_{P^{\operatorname*{op}}}^{-1}\right)
=\operatorname*{ord}\left(  \overline{R}_{P^{\operatorname*{op}}}\right)  $.
This proves Proposition \ref{prop.op.ord}.
\end{proof}

\section{Skeletal posets}

We will now introduce a class of posets which we call ``skeletal posets''.
Roughly speaking, these are graded posets built up inductively from the empty
poset by the operations of disjoint union (but only allowing disjoint union of
two $n$-graded posets for one and the same value of $n$) and ``grafting'' on
an antichain (generalizing the idea of grafting a tree on a new root). In
particular, all graded forests (oriented either away from the roots or towards
the roots) will belong to this class of posets, but also various other posets.
We begin by defining the notions involved:

\begin{definition}
Let $n\in\mathbb{N}$. Let $P$ and $Q$ be two $n$-graded posets. We denote by
$PQ$ the disjoint union of the posets $P$ and $Q$. (This disjoint union is
denoted by $P+Q$ in \cite[\S 3.2]{stanley-ec1}. Its poset structure is defined
in such a way that any element of $P$ and any element of $Q$ are incomparable,
while $P$ and $Q$ are subposets of $PQ$.) Clearly, $PQ$ is again an $n$-graded poset.
\end{definition}


\begin{definition}
Let $n\in\mathbb{N}$. Let $P$ be an $n$-graded poset. Let $k$ be a positive
integer. We denote by $B_{k}P$ the result of adding $k$ new elements to the
poset $P$, and declaring these $k$ elements to be smaller than each of the
elements of $P$ (but incomparable with each other). Clearly, $B_{k}P$ is an
$\left(  n+1\right)  $-graded poset.
\end{definition}

\begin{definition}
Let $n\in\mathbb{N}$. Let $P$ be an $n$-graded poset. Let $k$ be a positive
integer. We denote by $B_{k}^{\prime}P$ the result of adding $k$ new elements
to the poset $P$, and declaring these $k$ elements to be larger than each of
the elements of $P$ (but incomparable with each other). Clearly,
$B_{k}^{\prime}P$ is an $\left(  n+1\right)  $-graded poset.
\end{definition}

If $P$ is an $n$-graded poset and $k$ is a positive integer, then, in the
notations of Stanley (\cite[\S 3.2]{stanley-ec1}), we have $B_{k}P=A_{k}\oplus
P$ and $B_{k}^{\prime}P=P\oplus A_{k}$, where $A_{k}$ denotes the $k$-element antichain.

It is easy to see that $B_{k}P$ and $B_{k}^{\prime}P$ are ``symmetric''
notions with respect to taking the opposite poset:

\begin{proposition}
\label{prop.skeletal.op}Let $n\in\mathbb{N}$. Let $P$ be an $n$-graded poset.
Then, $B_{k}^{\prime}P=\left(  B_{k}\left(  P^{\operatorname*{op}}\right)
\right)  ^{\operatorname*{op}}$. (Here, we are using the notation introduced
in Definition \ref{def.op}.)
\end{proposition}

We now define the notion of a skeletal poset:

\begin{definition}
\label{def.skeletal}We define the class of \textit{skeletal posets}
inductively by means of the following axioms:

-- The empty poset is skeletal.

-- If $P$ is an $n$-graded skeletal poset and $k$ is a positive integer, then
the posets $B_{k}P$ and $B_{k}^{\prime}P$ are skeletal.

-- If $n$ is a nonnegative integer and $P$ and $Q$ are two $n$-graded skeletal
posets, then the poset $PQ$ is skeletal.

Notice that every skeletal poset is graded. Also, notice that every graded
rooted forest (made into a poset by having every node smaller than its
children) is a skeletal poset. (Indeed, every graded rooted forest can be
constructed from $\varnothing$ using merely the operations $P\mapsto B_{1}P$
and $\left(  P,Q\right)  \mapsto PQ$.) Also, every graded rooted arborescence
(i.e., the opposite poset of a graded rooted tree) is a skeletal poset (for a
similar reason).
\end{definition}

\begin{example}
The rooted forest
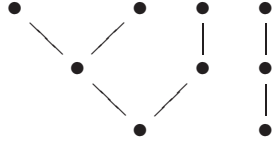
\[
\xymatrix@C=1em@R=1em{
\bullet\ar@{-}[dr] & & \bullet\ar@{-}[dl] & \bullet\ar@{-}[d] & \bullet\ar@{-}[d]
\\
& \bullet\ar@{-}[dr] & & \bullet\ar@{-}[dl] & \bullet\ar@{-}[d] \\
& & \bullet & & \bullet
}
\]
is skeletal, and in fact can be written as
$\left(
B_{1}\left(  \left(  B_{1}\left(  B_{2}\varnothing\right)  \right)  \left(
B_{1}\left(  B_{1}\varnothing\right)  \right)  \right)  \right)  \left(
B_{1}\left(  B_{1}\left(  B_{1}\varnothing\right)  \right)  \right)  $. (This
form of writing is not unique, since $B_{2}\varnothing=\left(  B_{1}%
\varnothing\right)  \left(  B_{1}\varnothing\right)  $.)

The tree
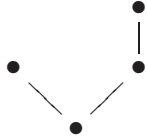
\[
\xymatrix@C=1em@R=1em{ & & \bullet\ar@{-}[d] \\ \bullet
\ar@{-}[dr] & & \bullet\ar@{-}[dl] \\ & \bullet& }
\]
can be written as
$B_{1}\left(  \left(  B_{1}\varnothing\right)  \left(  B_{1}\left(
B_{1}\varnothing\right)  \right)  \right)  $, but is \textbf{not} skeletal
because $B_{1}\varnothing$ and $B_{1}\left(  B_{1}\varnothing\right)  $ are
not $n$-graded with one and the same $n$.

The poset
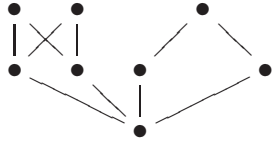
\[
\xymatrix@C=1em@R=1em{ \bullet\ar@{-}[dr] \ar@{-}[d] &
\bullet\ar@{-}[d] \ar@{-}[dl] & & \bullet\ar@{-}[dl]\ar@{-}[dr] &
\\ \bullet\ar@{-}[drr] & \bullet\ar@{-}[dr] & \bullet\ar@{-}[d] & &
\bullet\ar@{-}[dll] \\ & & \bullet& & }
\]
is neither a tree nor an
arborescence, but it has the form $B_{1}\left(  \left(  B_{2}\left(
B_{2}\varnothing\right)  \right)  \left(  B_{1}^{\prime}\left(  B_{2}%
\varnothing\right)  \right)  \right)  $ and is skeletal.
\end{example}

Our main result on skeletal posets is the following:

\begin{proposition}
\label{prop.skeletal.ords}Let $P$ be a skeletal poset. Let $\mathbb{K}$ be a
field. Then, $\operatorname*{ord}\left(  R_{P}\right)  $ and
$\operatorname*{ord}\left(  \overline{R}_{P}\right)  $ are finite.
\end{proposition}

In order to be able to prove this proposition, we first build up some
machinery for determining $\operatorname*{ord}\left(  R_{P}\right)  $ and
$\operatorname*{ord}\left(  \overline{R}_{P}\right)  $ given such orders in
smaller posets. Here is a very basic fact to get started:

\begin{proposition}
\label{prop.PQ.ord}Fix $n\in\mathbb{N}$. Let $P$ and $Q$ be two $n$-graded
posets. Let $\mathbb{K}$ be a field. Then, $\operatorname*{ord}\left(
R_{PQ}\right)  =\operatorname{lcm}\left(  \operatorname*{ord}\left(
R_{P}\right)  ,\operatorname*{ord}\left(  R_{Q}\right)  \right)  $.
\end{proposition}

\begin{vershort}
\begin{proof}
[Proof of Proposition \ref{prop.PQ.ord}.]The proof of this is as easy as it
looks: a $\mathbb{K}$-labelling of the disjoint union $PQ$ can be regarded as
a pair of a $\mathbb{K}$-labelling of $P$ and a $\mathbb{K}$-labelling of $Q$
(with identical labels at $0$ and $1$), and the map $R$ (as well as all
$R_{i}$) acts on these labellings independently.
\end{proof}
\end{vershort}

\begin{verlong}
\begin{proof}
[Proof of Proposition \ref{prop.PQ.ord}.]We WLOG assume that $P$ and $Q$ are
disjoint. Thus, we can regard $P$ and $Q$ as subposets of $PQ$. Let us do so.
Moreover, let us also regard $\widehat{P}$ as a subposet of $\widehat{PQ}$ by
using the inclusion $P\rightarrow PQ$ and identifying the element $0$ of
$\widehat{P}$ with the element $0$ of $\widehat{PQ}$ and identifying the
element $1$ of $\widehat{P}$ with the element $1$ of $\widehat{PQ}$.
Similarly, let us regard $\widehat{Q}$ as a subposet of $\widehat{PQ}$. Then,
if $f\in\mathbb{K}^{\widehat{PQ}}$ is any $\mathbb{K}$-labelling, then the
restrictions $f\mid_{\widehat{P}}$ and $f\mid_{\widehat{Q}}$ are $\mathbb{K}%
$-labellings in $\mathbb{K}^{\widehat{P}}$ and $\mathbb{K}^{\widehat{Q}}$, respectively.

Notice that if $v$ is any element of $P$, then%
\begin{equation}
\left(
\begin{array}
[c]{c}%
\text{any element }u\text{ of }\widehat{PQ}\text{ satisfying }u\lessdot
v\text{ in }\widehat{PQ}\\
\text{must be an element }u\text{ of }\widehat{P}\text{ satisfying }u\lessdot
v\text{ in }\widehat{P}%
\end{array}
\right)  . \label{pf.PQ.ord.1a}%
\end{equation}
\footnote{\textit{Proof of (\ref{pf.PQ.ord.1a}):} Let $v\in P$.
\par
Let $u$ be an element of $\widehat{PQ}$ satisfying $u\lessdot v$ in
$\widehat{PQ}$. Hence, $u<v$ in $\widehat{PQ}$. Thus, we cannot have $u=1$
(because otherwise, we would have $u=1$, and therefore $1=u<v$ in
$\widehat{PQ}$, which would contradict the fact that there exists no
$\alpha\in\widehat{PQ}$ such that $1<\alpha$ in $\widehat{PQ}$). In other
words, we must have $u\neq1$.
\par
Let us first show that $u\in\widehat{P}$. In fact, assume (for the sake of
contradiction) that we don't have $u\in\widehat{P}$. Thus, $u\notin%
\widehat{P}$, so that $u\in\widehat{PQ}\setminus\widehat{P}=Q$. By the
construction of $PQ$, we know that any element of $P$ and any element of $Q$
are incomparable in $PQ$. Hence, $v$ and $u$ are incomparable in $PQ$ (since
$v$ is an element of $P$ and since $u$ is an element of $Q$). Thus, $v$ and
$u$ are also incomparable in $\widehat{PQ}$ (because any two elements of $PQ$
which are incomparable in $PQ$ are also incomparable in $\widehat{PQ}$). This
contradicts the fact that $u<v$ in $\widehat{PQ}$. This contradiction shows
that our assumption (that we don't have $u\in\widehat{P}$) was wrong. Hence,
we have $u\in\widehat{P}$.
\par
Next, let us prove that $u<v$ in $\widehat{P}$. This is obvious if $u=0$
(because clearly, $0<v$ in $\widehat{P}$). Hence, for the rest of the proof
that $u<v$ in $P$, we can WLOG assume that we don't have $u=0$. Assume this.
We don't have $u=0$. So we have $u\neq0$. Combined with $u\neq1$, this yields
$u\notin\left\{  0,1\right\}  $. Since $u\in\widehat{P}$ and $u\notin\left\{
0,1\right\}  $, we have $u\in\widehat{P}\setminus\left\{  0,1\right\}  =P$.
Now, if $\alpha$ and $\beta$ are two elements of $PQ$ satisfying $\alpha
<\beta$ in $\widehat{PQ}$, then $\alpha<\beta$ also holds in $PQ$. Applying
this to $\alpha=u$ and $\beta=v$, we obtain that $u<v$ holds in $PQ$ (since
$u\in P\subseteq PQ$, $v\in P\subseteq PQ$ and since $u<v$ in $\widehat{PQ}$).
Recall now that if $\alpha$ and $\beta$ are two elements of $P$ satisfying
$\alpha<\beta$ in $PQ$, then $\alpha<\beta$ also holds in $P$. Applying this
to $\alpha=u$ and $\beta=v$, we obtain that $u<v$ holds in $P$ (since $u<v$ in
$PQ$). Thus, $u<v$ holds in $\widehat{P}$ as well (since $P$ is a subposet of
$\widehat{P}$). We have thus shown that $u<v$ in $\widehat{P}$.
\par
We now need to show that $u\lessdot v$ in $\widehat{P}$. Let $w$ be an element
of $\widehat{P}$ satisfying $u<w<v$ in $\widehat{P}$. Then, $u<w<v$ in
$\widehat{PQ}$ (since $u<w<v$ in $\widehat{P}$, and since $\widehat{P}$ is a
subposet of $\widehat{PQ}$). This contradicts the fact that $u\lessdot v$ in
$\widehat{PQ}$. Now, forget that we fixed $w$. We thus have proven a
contradiction for every element $w$ of $\widehat{P}$ satisfying $u<w<v$ in
$\widehat{P}$. Hence, there exists no element $w$ of $\widehat{P}$ satisfying
$u<w<v$ in $\widehat{P}$. Thus, $u\lessdot v$ in $\widehat{P}$ (since $u<v$ in
$\widehat{P}$). Altogether, $u\in\widehat{P}$ and $u\lessdot v$ in
$\widehat{P}$.
\par
Now, forget that we fixed $u$. We thus have shown that any element $u$ of
$\widehat{PQ}$ satisfying $u\lessdot v$ in $\widehat{PQ}$ must satisfy
$u\in\widehat{P}$ and $u\lessdot v$ in $\widehat{P}$. In other words, any
element $u$ of $\widehat{PQ}$ satisfying $u\lessdot v$ in $\widehat{PQ}$ must
be an element $u$ of $\widehat{P}$ satisfying $u\lessdot v$ in $\widehat{P}$.
This proves (\ref{pf.PQ.ord.1a}).} Also, if $v$ is an element of $P$, then%
\begin{equation}
\left(
\begin{array}
[c]{c}%
\text{any element }u\text{ of }\widehat{P}\text{ satisfying }u\lessdot v\text{
in }\widehat{P}\\
\text{must be an element }u\text{ of }\widehat{PQ}\text{ satisfying }u\lessdot
v\text{ in }\widehat{PQ}%
\end{array}
\right)  . \label{pf.PQ.ord.1b}%
\end{equation}
\footnote{\textit{Proof of (\ref{pf.PQ.ord.1b}):} Let $v\in P$.
\par
Let $u$ be an element of $\widehat{P}$ satisfying $u\lessdot v$ in
$\widehat{P}$.
\par
We know that $u\in\widehat{P}\subseteq\widehat{PQ}$ (since $\widehat{P}$ is a
subposet of $\widehat{PQ}$). Also, $u<v$ in $\widehat{P}$ (since $u\lessdot v$
in $\widehat{P}$), so that $u<v$ in $\widehat{PQ}$ (again because
$\widehat{P}$ is a subposet of $\widehat{PQ}$).
\par
Now, we are going to prove that $u\lessdot v$ in $\widehat{PQ}$. Let $w$ be an
element of $\widehat{PQ}$ satisfying $u<w<v$ in $\widehat{PQ}$.
\par
We cannot have $w=1$ (because otherwise, we would have $w=1$, and therefore
$1=w<v$ in $\widehat{PQ}$, which would contradict the fact that there exists
no $\alpha\in\widehat{PQ}$ such that $1<\alpha$ in $\widehat{PQ}$). In other
words, we must have $w\neq1$.
\par
Also, we cannot have $w=0$ (because otherwise, we would have $w=0$, and
therefore $u<w=0$ in $\widehat{PQ}$, which would contradict the fact that
there exists no $\alpha\in\widehat{PQ}$ such that $\alpha<0$ in $\widehat{PQ}%
$). In other words, we must have $w\neq0$. Combined with $w\neq1$, this yields
$w\notin\left\{  0,1\right\}  $. Hence, $w\in\widehat{PQ}\setminus\left\{
0,1\right\}  =PQ$.
\par
Now, assume (for the sake of contradiction) that we don't have $w\in P$. Then,
$w\notin P$. Since $w\in PQ$ and $w\notin P$, we have $w\in PQ\setminus P=Q$.
By the construction of $PQ$, we know that any element of $P$ and any element
of $Q$ are incomparable in $PQ$. Hence, $v$ and $w$ are incomparable in $PQ$
(since $v$ is an element of $P$ and since $w$ is an element of $Q$). Thus, $v$
and $w$ are also incomparable in $\widehat{PQ}$ (because any two elements of
$PQ$ which are incomparable in $PQ$ are also incomparable in $\widehat{PQ}$).
This contradicts the fact that $w<v$ in $\widehat{PQ}$. This contradiction
shows that our assumption (that we don't have $w\in P$) was wrong. Hence, we
have $w\in P$.
\par
Recall that if $\alpha$ and $\beta$ are two elements of $PQ$ satisfying
$\alpha<\beta$ in $\widehat{PQ}$, then $\alpha<\beta$ also holds in $PQ$.
Applying this to $\alpha=w$ and $\beta=v$, we obtain that $w<v$ holds in $PQ$
(since $w\in P\subseteq PQ$, $v\in P\subseteq PQ$ and since $w<v$ in
$\widehat{PQ}$). Recall now that if $\alpha$ and $\beta$ are two elements of
$P$ satisfying $\alpha<\beta$ in $PQ$, then $\alpha<\beta$ also holds in $P$.
Applying this to $\alpha=w$ and $\beta=v$, we obtain that $w<v$ holds in $P$
(since $w<v$ in $PQ$). Thus, $w<v$ holds in $\widehat{P}$ as well (since $P$
is a subposet of $\widehat{P}$).
\par
Let us next prove that $u<w$ in $\widehat{P}$. This is obvious if $u=0$
(because clearly, $0<w$ in $\widehat{P}$). Hence, for the rest of the proof
that $u<w$ in $\widehat{P}$, we can WLOG assume that we don't have $u=0$.
Assume this. We don't have $u=0$. So we have $u\neq0$. But we also don't have
$u=1$ (because $u<v\leq1$ in $\widehat{P}$). Hence, $u\neq1$. Since $u\neq0$
and $u\neq1$, we have $u\notin\left\{  0,1\right\}  $, so that $u\in
\widehat{P}\setminus\left\{  0,1\right\}  =P$. Now, recall that if $\alpha$
and $\beta$ are two elements of $PQ$ satisfying $\alpha<\beta$ in
$\widehat{PQ}$, then $\alpha<\beta$ also holds in $PQ$. Applying this to
$\alpha=u$ and $\beta=w$, we obtain that $u<w$ holds in $PQ$ (since $u\in
P\subseteq PQ$, $w\in P\subseteq PQ$ and since $u<w$ in $\widehat{PQ}$).
Recall now that if $\alpha$ and $\beta$ are two elements of $P$ satisfying
$\alpha<\beta$ in $PQ$, then $\alpha<\beta$ also holds in $P$. Applying this
to $\alpha=u$ and $\beta=w$, we obtain that $u<w$ holds in $P$ (since $u<w$ in
$PQ$). Thus, $u<w$ holds in $\widehat{P}$ as well (since $P$ is a subposet of
$\widehat{P}$). We thus have proven that $u<w$ in $\widehat{P}$.
\par
We now have $u<w<v$ in $\widehat{P}$. This contradicts the fact that
$u\lessdot v$ in $\widehat{P}$. Now, forget that we fixed $w$. We thus have
proven a contradiction for every element $w$ of $\widehat{PQ}$ satisfying
$u<w<v$ in $\widehat{PQ}$. Hence, there exists no element $w$ of $\widehat{PQ}$
satisfying $u<w<v$ in $\widehat{PQ}$. Thus, $u\lessdot v$ in $\widehat{PQ}$
(since $u<v$ in $\widehat{PQ}$). Altogether, $u\in\widehat{PQ}$ and $u\lessdot
v$ in $\widehat{PQ}$.
\par
Now, forget that we fixed $u$. We thus have shown that any element $u$ of
$\widehat{P}$ satisfying $u\lessdot v$ in $\widehat{P}$ must satisfy
$u\in\widehat{PQ}$ and $u\lessdot v$ in $\widehat{PQ}$. In other words, any
element $u$ of $\widehat{P}$ satisfying $u\lessdot v$ in $\widehat{P}$ must be
an element $u$ of $\widehat{PQ}$ satisfying $u\lessdot v$ in $\widehat{PQ}$.
This proves (\ref{pf.PQ.ord.1b}).} Now, if $v$ is an element of $P$, then%
\begin{equation}
\left(
\begin{array}
[c]{c}%
\text{the elements }u\text{ of }\widehat{PQ}\text{ satisfying }u\lessdot
v\text{ in }\widehat{PQ}\\
\text{are precisely the elements }u\text{ of }\widehat{P}\text{ satisfying
}u\lessdot v\text{ in }\widehat{P}%
\end{array}
\right)  . \label{pf.PQ.ord.1}%
\end{equation}
\footnote{\textit{Proof of (\ref{pf.PQ.ord.1}):} Let $v\in P$. Combining the
results (\ref{pf.PQ.ord.1a}) and (\ref{pf.PQ.ord.1b}), we see that the
elements $u$ of $\widehat{PQ}$ satisfying $u\lessdot v$ in $\widehat{PQ}$ are
precisely the elements $u$ of $\widehat{P}$ satisfying $u\lessdot v$ in
$\widehat{P}$. This proves (\ref{pf.PQ.ord.1}).} Similarly, if $v$ is an
element of $P$, then%
\begin{equation}
\left(
\begin{array}
[c]{c}%
\text{the elements }u\text{ of }\widehat{PQ}\text{ satisfying }u\gtrdot
v\text{ in }\widehat{PQ}\\
\text{are precisely the elements }u\text{ of }\widehat{P}\text{ satisfying
}u\gtrdot v\text{ in }\widehat{P}%
\end{array}
\right)  . \label{pf.PQ.ord.1o}%
\end{equation}

We have%
\begin{equation}
R_{P}\left(  f\mid_{\widehat{P}}\right)  =\left(  R_{PQ}f\right)
\mid_{\widehat{P}}\ \ \ \ \ \ \ \ \ \ \text{for every }f\in\mathbb{K}%
^{\widehat{PQ}}. \label{pf.PQ.ord.2}%
\end{equation}
\footnote{\textit{Proof of (\ref{pf.PQ.ord.2}):} Let $f\in\mathbb{K}%
^{\widehat{PQ}}$. Let $v\in P$. Then,%
\begin{equation}
\left(  \left(  R_{PQ}f\right)  \mid_{\widehat{P}}\right)  \left(  v\right)
=\left(  R_{PQ}f\right)  \left(  v\right)  =\dfrac{1}{f\left(  v\right)
}\cdot\dfrac{\sum\limits_{\substack{u\in\widehat{PQ};\\u\lessdot v}}f\left(
u\right)  }{\sum\limits_{\substack{u\in\widehat{PQ};\\u\gtrdot v}}\dfrac
{1}{\left(  R_{PQ}f\right)  \left(  u\right)  }} \label{pf.PQ.ord.2.pf.1}%
\end{equation}
(by Proposition \ref{prop.R.implicit}, applied to $PQ$ instead of $P$).
\par
Now, (\ref{pf.PQ.ord.1}) yields that the elements $u$ of $\widehat{PQ}$
satisfying $u\lessdot v$ are precisely the elements $u$ of $\widehat{P}$
satisfying $u\lessdot v$. Hence, we can replace the summation sign
$\sum\limits_{\substack{u\in\widehat{PQ};\\u\lessdot v}}$ by a $\sum
\limits_{\substack{u\in\widehat{P};\\u\lessdot v}}$ sign in $\sum
\limits_{\substack{u\in\widehat{PQ};\\u\lessdot v}}f\left(  u\right)  $.
Hence, $\sum\limits_{\substack{u\in\widehat{PQ};\\u\lessdot v}}f\left(
u\right)  =\sum\limits_{\substack{u\in\widehat{P};\\u\lessdot v}}f\left(
u\right)  $. Also, we can replace the summation sign $\sum
\limits_{\substack{u\in\widehat{PQ};\\u\gtrdot v}}$ by a $\sum
\limits_{\substack{u\in\widehat{P};\\u\gtrdot v}}$ sign in $\sum
\limits_{\substack{u\in\widehat{PQ};\\u\gtrdot v}}\dfrac{1}{\left(
R_{PQ}f\right)  \left(  u\right)  }$ (because (\ref{pf.PQ.ord.1o}) yields that
the elements $u$ of $\widehat{PQ}$ satisfying $u\gtrdot v$ are precisely the
elements $u$ of $\widehat{P}$ satisfying $u\gtrdot v$). Thus,
\[
\sum\limits_{\substack{u\in\widehat{PQ};\\u\gtrdot v}}\dfrac{1}{\left(
R_{PQ}f\right)  \left(  u\right)  }=\sum\limits_{\substack{u\in\widehat{P}%
;\\u\gtrdot v}}\dfrac{1}{\left(  R_{PQ}f\right)  \left(  u\right)  }%
=\sum\limits_{\substack{u\in\widehat{P};\\u\gtrdot v}}\dfrac{1}{\left(
\left(  R_{PQ}f\right)  \mid_{\widehat{P}}\right)  \left(  u\right)  }%
\]
(since every $u\in\widehat{P}$ satisfying $u\gtrdot v$ satisfies $\left(
R_{PQ}f\right)  \left(  u\right)  =\left(  \left(  R_{PQ}f\right)
\mid_{\widehat{P}}\right)  \left(  u\right)  $ (since $u\in\widehat{P}$)).
Now, (\ref{pf.PQ.ord.2.pf.1}) becomes%
\begin{equation}
\left(  \left(  R_{PQ}f\right)  \mid_{\widehat{P}}\right)  \left(  v\right)
=\dfrac{1}{f\left(  v\right)  }\cdot\dfrac{\sum\limits_{\substack{u\in
\widehat{PQ};\\u\lessdot v}}f\left(  u\right)  }{\sum\limits_{\substack{u\in
\widehat{PQ};\\u\gtrdot v}}\dfrac{1}{\left(  R_{PQ}f\right)  \left(  u\right)
}}=\dfrac{1}{f\left(  v\right)  }\cdot\dfrac{\sum\limits_{\substack{u\in
\widehat{P};\\u\lessdot v}}f\left(  u\right)  }{\sum\limits_{\substack{u\in
\widehat{P};\\u\gtrdot v}}\dfrac{1}{\left(  \left(  R_{PQ}f\right)
\mid_{\widehat{P}}\right)  \left(  u\right)  }} \label{pf.PQ.ord.2.pf.2}%
\end{equation}
(since $\sum\limits_{\substack{u\in\widehat{PQ};\\u\lessdot v}}f\left(
u\right)  =\sum\limits_{\substack{u\in\widehat{P};\\u\lessdot v}}f\left(
u\right)  $ and $\sum\limits_{\substack{u\in\widehat{PQ};\\u\gtrdot v}%
}\dfrac{1}{\left(  R_{PQ}f\right)  \left(  u\right)  }=\sum
\limits_{\substack{u\in\widehat{P};\\u\gtrdot v}}\dfrac{1}{\left(  \left(
R_{PQ}f\right)  \mid_{\widehat{P}}\right)  \left(  u\right)  }$).
\par
Now, forget that we fixed $v$. We thus have shown that (\ref{pf.PQ.ord.2.pf.1}%
) holds for every $v\in P$.
\par
Notice also that $\left(  \left(  R_{PQ}f\right)  \mid_{\widehat{P}}\right)
\left(  0\right)  =\left(  R_{PQ}f\right)  \left(  0\right)  =f\left(
0\right)  $ (by Proposition \ref{prop.R.implicit.01}, applied to $PQ$ instead
of $P$) and thus%
\[
\left(  \left(  R_{PQ}f\right)  \mid_{\widehat{P}}\right)  \left(  0\right)
=f\left(  0\right)  =\left(  f\mid_{\widehat{P}}\right)  \left(  0\right)  .
\]
Similarly, $\left(  \left(  R_{PQ}f\right)  \mid_{\widehat{P}}\right)  \left(
1\right)  =\left(  f\mid_{\widehat{P}}\right)  \left(  1\right)  $.
\par
So we know that $\left(  \left(  R_{PQ}f\right)  \mid_{\widehat{P}}\right)
\left(  0\right)  =\left(  f\mid_{\widehat{P}}\right)  \left(  0\right)  $ and
$\left(  \left(  R_{PQ}f\right)  \mid_{\widehat{P}}\right)  \left(  1\right)
=\left(  f\mid_{\widehat{P}}\right)  \left(  1\right)  $, and we also know
that (\ref{pf.PQ.ord.2.pf.2}) holds for every $v\in P$. Hence, Proposition
\ref{prop.R.implicit.converse} (applied to $\left(  R_{PQ}f\right)
\mid_{\widehat{P}}$ and $f\mid_{\widehat{P}}$ instead of $g$ and $f$) yields
$\left(  R_{PQ}f\right)  \mid_{\widehat{P}}=R_{P}\left(  f\mid_{\widehat{P}%
}\right)  $. In other words, $R_{P}\left(  f\mid_{\widehat{P}}\right)
=\left(  R_{PQ}f\right)  \mid_{\widehat{P}}$. This proves (\ref{pf.PQ.ord.2}%
).} Thus,%
\begin{equation}
R_{P}^{\ell}\left(  f\mid_{\widehat{P}}\right)  =\left(  R_{PQ}^{\ell
}f\right)  \mid_{\widehat{P}}\ \ \ \ \ \ \ \ \ \ \text{for every }%
f\in\mathbb{K}^{\widehat{PQ}}\text{ and }\ell\in\mathbb{N}.
\label{pf.PQ.ord.3}%
\end{equation}
\footnote{\textit{Proof of (\ref{pf.PQ.ord.3}):} We will prove
(\ref{pf.PQ.ord.3}) by induction over $\ell$.
\par
\textit{Induction base:} Every $f\in\mathbb{K}^{\widehat{PQ}}$ satisfies
$\underbrace{R_{P}^{0}}_{=\operatorname*{id}}\left(  f\mid_{\widehat{P}%
}\right)  =\operatorname*{id}\left(  f\mid_{\widehat{P}}\right)
=f\mid_{\widehat{P}}$ and $\left.  \left(  \underbrace{R_{PQ}^{0}%
}_{=\operatorname*{id}}f\right)  \mid_{\widehat{P}}\right.  =\left.  \left(
\operatorname*{id}f\right)  \mid_{\widehat{P}}\right.  =f\mid_{\widehat{P}}$.
Hence, every $f\in\mathbb{K}^{\widehat{PQ}}$ satisfies $R_{P}^{0}\left(
f\mid_{\widehat{P}}\right)  =\left.  f\mid_{\widehat{P}}\right.  =\left(
R_{PQ}^{0}f\right)  \mid_{\widehat{P}}$. Thus, (\ref{pf.PQ.ord.3}) holds for
$\ell=0$. This completes the induction base.
\par
\textit{Induction step:} Let $L\in\mathbb{N}$. Assume that (\ref{pf.PQ.ord.3})
holds for $\ell=L$. We must prove that (\ref{pf.PQ.ord.3}) also holds for
$\ell=L+1$.
\par
We know that (\ref{pf.PQ.ord.3}) holds for $\ell=L$. In other words,%
\begin{equation}
R_{P}^{L}\left(  f\mid_{\widehat{P}}\right)  =\left(  R_{PQ}^{L}f\right)
\mid_{\widehat{P}}\ \ \ \ \ \ \ \ \ \ \text{for every }f\in\mathbb{K}%
^{\widehat{PQ}}\text{.} \label{pf.PQ.ord.3.pf.1}%
\end{equation}
\par
Now, every $f\in\mathbb{K}^{\widehat{PQ}}$ satisfies%
\begin{align*}
\underbrace{R_{P}^{L+1}}_{=R_{P}\circ R_{P}^{L}}\left(  f\mid_{\widehat{P}%
}\right)   &  =\left(  R_{P}\circ R_{P}^{L}\right)  \left(  f\mid
_{\widehat{P}}\right)  =R_{P}\left(  \underbrace{R_{P}^{L}\left(
f\mid_{\widehat{P}}\right)  }_{\substack{=\left(  R_{PQ}^{L}f\right)
\mid_{\widehat{P}}\\\text{(by (\ref{pf.PQ.ord.3.pf.1}))}}}\right) \\
&  =R_{P}\left(  \left(  R_{PQ}^{L}f\right)  \mid_{\widehat{P}}\right)
=\left(  \underbrace{R_{PQ}\left(  R_{PQ}^{L}f\right)  }_{=\left(  R_{PQ}\circ
R_{PQ}^{L}\right)  f}\right)  \mid_{\widehat{P}}\\
&  \ \ \ \ \ \ \ \ \ \ \left(  \text{by (\ref{pf.PQ.ord.2}), applied to
}R_{PQ}^{L}f\text{ instead of }f\right) \\
&  =\left.  \left(  \underbrace{\left(  R_{PQ}\circ R_{PQ}^{L}\right)
}_{=R_{PQ}^{L+1}}f\right)  \mid_{\widehat{P}}\right.  =\left(  R_{PQ}%
^{L+1}f\right)  \mid_{\widehat{P}}.
\end{align*}
In other words, $R_{P}^{L+1}\left(  f\mid_{\widehat{P}}\right)  =\left(
R_{PQ}^{L+1}f\right)  \mid_{\widehat{P}}$ for every $f\in\mathbb{K}%
^{\widehat{PQ}}$. In other words, (\ref{pf.PQ.ord.3}) holds for $\ell=L+1$.
This completes the induction step. The induction proof of (\ref{pf.PQ.ord.3})
is thus complete.}

Let us now show that $\operatorname*{ord}\left(  R_{PQ}\right)  \mid
\operatorname{lcm}\left(  \operatorname*{ord}\left(  R_{P}\right)
,\operatorname*{ord}\left(  R_{Q}\right)  \right)  $.

\textit{Proof of }$\operatorname*{ord}\left(  R_{PQ}\right)  \mid
\operatorname{lcm}\left(  \operatorname*{ord}\left(  R_{P}\right)
,\operatorname*{ord}\left(  R_{Q}\right)  \right)  $\textit{:} First, we can
WLOG assume that $\operatorname*{ord}\left(  R_{P}\right)  $ and
$\operatorname*{ord}\left(  R_{Q}\right)  $ are finite, because otherwise the
claim $\operatorname*{ord}\left(  R_{PQ}\right)  \mid\operatorname{lcm}\left(
\operatorname*{ord}\left(  R_{P}\right)  ,\operatorname*{ord}\left(
R_{Q}\right)  \right)  $ is trivial\footnote{Indeed, unless both
$\operatorname*{ord}\left(  R_{P}\right)  $ and $\operatorname*{ord}\left(
R_{Q}\right)  $ are finite, we have $\operatorname{lcm}\left(
\operatorname*{ord}\left(  R_{P}\right)  ,\operatorname*{ord}\left(
R_{Q}\right)  \right)  =\infty$ (and we know that $\infty$ is divisible by any
positive integer).}. Assume this. Then, $\operatorname*{ord}\left(
R_{P}\right)  $ and $\operatorname*{ord}\left(  R_{Q}\right)  $ are finite, so
that $\operatorname{lcm}\left(  \operatorname*{ord}\left(  R_{P}\right)
,\operatorname*{ord}\left(  R_{Q}\right)  \right)  $ is finite. Denote
$\operatorname{lcm}\left(  \operatorname*{ord}\left(  R_{P}\right)
,\operatorname*{ord}\left(  R_{Q}\right)  \right)  $ by $\ell$. Then,
$\ell=\operatorname{lcm}\left(  \operatorname*{ord}\left(  R_{P}\right)
,\operatorname*{ord}\left(  R_{Q}\right)  \right)  $ is divisible by
$\operatorname*{ord}\left(  R_{P}\right)  $, and thus we have $R_{P}^{\ell
}=\operatorname*{id}$.

Now, let $f\in\mathbb{K}^{\widehat{PQ}}$ be arbitrary such that $R_{PQ}^{\ell
}f$ is well-defined. We have%
\begin{align*}
\left.  \left(  R_{PQ}^{\ell}f\right)  \mid_{\widehat{P}}\right.   &
=\underbrace{R_{P}^{\ell}}_{=\operatorname*{id}}\left(  f\mid_{\widehat{P}%
}\right)  \ \ \ \ \ \ \ \ \ \ \left(  \text{by (\ref{pf.PQ.ord.3})}\right) \\
&  =\operatorname*{id}\left(  f\mid_{\widehat{P}}\right)  =\left.
f\mid_{\widehat{P}}\right.  .
\end{align*}
Similarly, $\left.  \left(  R_{PQ}^{\ell}f\right)  \mid_{\widehat{Q}}\right.
=\left.  f\mid_{\widehat{Q}}\right.  $. Now, it is easy to see that
$R_{PQ}^{\ell}f=f$\ \ \ \ \footnote{\textit{Proof.} Let $v\in\widehat{PQ}$. We
are going to prove that $\left(  R_{PQ}^{\ell}f\right)  \left(  v\right)
=f\left(  v\right)  $.
\par
Indeed, if $v\in\widehat{P}$, then%
\[
\left(  R_{PQ}^{\ell}f\right)  \left(  v\right)  =\left(  \underbrace{\left(
R_{PQ}^{\ell}f\right)  \mid_{\widehat{P}}}_{=\left.  f\mid_{\widehat{P}%
}\right.  }\right)  \left(  v\right)  =\left(  f\mid_{\widehat{P}}\right)
\left(  v\right)  =f\left(  v\right)  .
\]
Hence, $\left(  R_{PQ}^{\ell}f\right)  \left(  v\right)  =f\left(  v\right)  $
is proven if $v\in\widehat{P}$. Thus, for the rest of the proof of $\left(
R_{PQ}^{\ell}f\right)  \left(  v\right)  =f\left(  v\right)  $, we can WLOG
assume that we don't have $v\in\widehat{P}$. Assume this.
\par
We don't have $v\in\widehat{P}$. Thus, we have $v\notin\widehat{P}$. Hence,
$v\in\widehat{PQ}\setminus\widehat{P}=Q\subseteq\widehat{Q}$. Thus,%
\[
\left(  R_{PQ}^{\ell}f\right)  \left(  v\right)  =\left(  \underbrace{\left(
R_{PQ}^{\ell}f\right)  \mid_{\widehat{Q}}}_{=\left.  f\mid_{\widehat{Q}%
}\right.  }\right)  \left(  v\right)  =\left(  f\mid_{\widehat{Q}}\right)
\left(  v\right)  =f\left(  v\right)  .
\]
Thus, $\left(  R_{PQ}^{\ell}f\right)  \left(  v\right)  =f\left(  v\right)  $
is proven.
\par
Now, forget that we fixed $v$. We thus have shown that $\left(  R_{PQ}^{\ell
}f\right)  \left(  v\right)  =f\left(  v\right)  $ for every $v\in
\widehat{PQ}$. In other words, $R_{PQ}^{\ell}f=f$, qed.}.

Now, forget that we fixed $f$. We thus have proven that $R_{PQ}^{\ell}f=f$ for
every $f\in\mathbb{K}^{\widehat{PQ}}$ for which $R_{PQ}^{\ell}f$ is
well-defined. In other words, $R_{PQ}^{\ell}=\operatorname*{id}$, so that
$\operatorname*{ord}\left(  R_{PQ}\right)  \mid\ell=\operatorname{lcm}\left(
\operatorname*{ord}\left(  R_{P}\right)  ,\operatorname*{ord}\left(
R_{Q}\right)  \right)  $.

Thus, we have shown that $\operatorname*{ord}\left(  R_{PQ}\right)
\mid\operatorname{lcm}\left(  \operatorname*{ord}\left(  R_{P}\right)
,\operatorname*{ord}\left(  R_{Q}\right)  \right)  $.

Next, we are going to prove that $\operatorname*{ord}\left(  R_{P}\right)
\mid\operatorname*{ord}\left(  R_{PQ}\right)  $.

\textit{Proof of }$\operatorname*{ord}\left(  R_{P}\right)  \mid
\operatorname*{ord}\left(  R_{PQ}\right)  $\textit{:} First, we can WLOG
assume that $\operatorname*{ord}\left(  R_{PQ}\right)  $ is finite, because
otherwise the claim $\operatorname*{ord}\left(  R_{P}\right)  \mid
\operatorname*{ord}\left(  R_{PQ}\right)  $ is trivial.

Let $f\in\mathbb{K}^{\widehat{P}}$ be any sufficiently generic (in the sense
of Zariski topology) $\mathbb{K}$-labelling of $P$. Clearly, there exists a
$\mathbb{K}$-labelling $F$ of $PQ$ which satisfies $F\mid_{\widehat{P}}=f$ and
is sufficiently generic for $R_{PQ}^{\operatorname*{ord}\left(  R_{PQ}\right)
}F$ to be well-defined\footnote{In fact, such a $\mathbb{K}$-labelling $F$ can
be obtained by extending $f$ from $\widehat{P}$ to $\widehat{PQ}$, choosing
the images of the elements of $\widehat{PQ}\setminus\widehat{P}=Q$
generically.}. Fix such a $\mathbb{K}$-labelling $F$. Then, (\ref{pf.PQ.ord.3}%
) (applied to $F$ and $\operatorname*{ord}\left(  R_{PQ}\right)  $ instead of
$f$ and $\ell$) yields%
\[
R_{P}^{\operatorname*{ord}\left(  R_{PQ}\right)  }\left(  F\mid_{\widehat{P}%
}\right)  =\left.  \left(  \underbrace{R_{PQ}^{\operatorname*{ord}\left(
R_{PQ}\right)  }}_{=\operatorname*{id}}F\right)  \mid_{\widehat{P}}\right.
=\left.  \left(  \operatorname*{id}F\right)  \mid_{\widehat{P}}\right.
=\left.  F\mid_{\widehat{P}}\right.  =f.
\]
Hence, $f=R_{P}^{\operatorname*{ord}\left(  R_{PQ}\right)  }%
\underbrace{\left(  F\mid_{\widehat{P}}\right)  }_{=f}=R_{P}%
^{\operatorname*{ord}\left(  R_{PQ}\right)  }f$. In other words,
$R_{P}^{\operatorname*{ord}\left(  R_{PQ}\right)  }f=f$.

Now, forget that we fixed $f$. We thus have proven that every sufficiently
generic $\mathbb{K}$-labelling $f$ of $P$ satisfies $R_{P}%
^{\operatorname*{ord}\left(  R_{PQ}\right)  }f=f$. In other words,
$R_{P}^{\operatorname*{ord}\left(  R_{PQ}\right)  }=\operatorname*{id}$.
Hence, $\operatorname*{ord}\left(  R_{P}\right)  \mid\operatorname*{ord}%
\left(  R_{PQ}\right)  $.

We thus have proven $\operatorname*{ord}\left(  R_{P}\right)  \mid
\operatorname*{ord}\left(  R_{PQ}\right)  $. Combining this with the statement
that $\operatorname*{ord}\left(  R_{Q}\right)  \mid\operatorname*{ord}\left(
R_{PQ}\right)  $ (whose proof is similar), we obtain $\operatorname{lcm}%
\left(  \operatorname*{ord}\left(  R_{P}\right)  ,\operatorname*{ord}\left(
R_{Q}\right)  \right)  \mid\operatorname*{ord}\left(  R_{PQ}\right)  $.
Combined with $\operatorname*{ord}\left(  R_{PQ}\right)  \mid
\operatorname{lcm}\left(  \operatorname*{ord}\left(  R_{P}\right)
,\operatorname*{ord}\left(  R_{Q}\right)  \right)  $, this yields
$\operatorname*{ord}\left(  R_{PQ}\right)  =\operatorname{lcm}\left(
\operatorname*{ord}\left(  R_{P}\right)  ,\operatorname*{ord}\left(
R_{Q}\right)  \right)  $. Thus, Proposition \ref{prop.PQ.ord} is proven.
\end{proof}
\end{verlong}

The analogue of Proposition \ref{prop.PQ.ord} with all $R$'s replaced by
$\overline{R}$'s is false. Instead, $\operatorname*{ord}\left(  \overline
{R}_{PQ}\right)  $ can be computed as follows:\footnote{The following
proposition is, in some sense, uninteresting, as it is a negative result (it
merely serves to convince one that $\operatorname*{ord}\left(  \overline
{R}_{PQ}\right)  $ is not lower than what is expected from Propositions
\ref{prop.ord-projord} and \ref{prop.PQ.ord}).}

\begin{proposition}
\label{prop.PQ.ord2}Fix $n\in\mathbb{N}$. Let $P$ and $Q$ be two $n$-graded
posets. Let $\mathbb{K}$ be a field. Then, $\operatorname*{ord}\left(
\overline{R}_{PQ}\right)  =\operatorname{lcm}\left(  \operatorname*{ord}%
\left(  R_{P}\right)  ,\operatorname*{ord}\left(  R_{Q}\right)  \right)  $.
\end{proposition}

\begin{proof}
[Proof of Proposition \ref{prop.PQ.ord2} (sketched).]Assume WLOG that $n\neq0$
(else, everything is obvious). Hence, $P$ and $Q$ are nonempty (being $n$-graded).

Proposition \ref{prop.ord-projord} yields $\operatorname*{ord}\left(
R_{PQ}\right)  =\operatorname{lcm}\left(  n+1,\operatorname*{ord}\left(
\overline{R}_{PQ}\right)  \right)  $.

WLOG assume that $\operatorname*{ord}\left(  R_{P}\right)  $ and
$\operatorname*{ord}\left(  R_{Q}\right)  $ are finite\footnote{Otherwise,
$\operatorname{lcm}\left(  \operatorname*{ord}\left(  R_{P}\right)
,\operatorname*{ord}\left(  R_{Q}\right)  \right)  $ is infinite, whence
$\operatorname*{ord}\left(  R_{PQ}\right)  $ is infinite (by Proposition
\ref{prop.PQ.ord}), whence $\operatorname*{ord}\left(  \overline{R}%
_{PQ}\right)  $ is infinite (because $\operatorname*{ord}\left(
R_{PQ}\right)  =\operatorname{lcm}\left(  n+1,\operatorname*{ord}\left(
\overline{R}_{PQ}\right)  \right)  $), whence Proposition \ref{prop.PQ.ord2}
is trivial.}. Then, Proposition \ref{prop.PQ.ord} shows that
$\operatorname*{ord}\left(  R_{PQ}\right)  =\operatorname{lcm}\left(
\operatorname*{ord}\left(  R_{P}\right)  ,\operatorname*{ord}\left(
R_{Q}\right)  \right)  $ is finite, so that $\operatorname*{ord}\left(
\overline{R}_{PQ}\right)  $ is finite (because $\operatorname*{ord}\left(
R_{PQ}\right)  =\operatorname{lcm}\left(  n+1,\operatorname*{ord}\left(
\overline{R}_{PQ}\right)  \right)  $). Let $\ell$ be $\operatorname*{ord}%
\left(  \overline{R}_{PQ}\right)  $. Then, $\ell$ is finite and satisfies
$\overline{R}_{PQ}^{\ell}=\operatorname*{id}$. We will show that $n+1\mid\ell$.

\begin{vershort}
For every $\mathbb{K}$-labelling $f$ of $PQ$ and every $i\in\left\{
0,1,...,n\right\}  $, define two elements $\mathbf{w}_{i}^{\left(  1\right)
}\left(  f\right)  $ and $\mathbf{w}_{i}^{\left(  2\right)  }\left(  f\right)
$ of $\mathbb{K}$ by%
\[
\mathbf{w}_{i}^{\left(  1\right)  }\left(  f\right)  =\sum_{\substack{x\in
\widehat{P}_{i};\ y\in\widehat{P}_{i+1};\\y\gtrdot x}}\dfrac{f\left(
x\right)  }{f\left(  y\right)  }\ \ \ \ \ \ \ \ \ \ \text{and}%
\ \ \ \ \ \ \ \ \ \ \mathbf{w}_{i}^{\left(  2\right)  }\left(  f\right)
=\sum_{\substack{x\in\widehat{Q}_{i};\ y\in\widehat{Q}_{i+1};\\y\gtrdot
x}}\dfrac{f\left(  x\right)  }{f\left(  y\right)  }%
\]
(where, of course, $\widehat{P}_{j}$ and $\widehat{Q}_{j}$ are embedded into
$\widehat{PQ}_{j}$ for every $j\in\left\{  0,1,...,n+1\right\}  $ in the
obvious way). These elements $\mathbf{w}_{i}^{\left(  1\right)  }\left(
f\right)  $ and $\mathbf{w}_{i}^{\left(  2\right)  }\left(  f\right)  $ are
defined not for every $f$, but for \textquotedblleft almost
every\textquotedblright\ $f$ in the sense of Zariski topology. We denote the
$\left(  n+1\right)  $-tuple
\[
\left(  \mathbf{w}_{0}^{\left(  1\right)  }\left(  f\right)  \diagup
\mathbf{w}_{0}^{\left(  2\right)  }\left(  f\right)  ,\ \mathbf{w}%
_{1}^{\left(  1\right)  }\left(  f\right)  \diagup\mathbf{w}_{1}^{\left(
2\right)  }\left(  f\right)  ,\ ...,\ \mathbf{w}_{n}^{\left(  1\right)
}\left(  f\right)  \diagup\mathbf{w}_{n}^{\left(  2\right)  }\left(  f\right)
\right)
\]
as the \textit{comparative w-tuple} of the labelling $f$. The advantage of
comparative w-tuples over usual w-tuples is the following fact: If $f$ and $g$
are two homogeneously equivalent $\mathbb{K}$-labellings of $PQ$, then%
\begin{equation}
\left(  \text{the comparative w-tuple of }f\right)  =\left(  \text{the
comparative w-tuple of }g\right)  . \label{pf.PQ.ord2.hgeq.short}%
\end{equation}
(This is easy to check and has no analogue for regular w-tuples.)
\end{vershort}

\begin{verlong}
For every $\mathbb{K}$-labelling $f$ of $PQ$ and every $i\in\left\{
0,1,...,n\right\}  $, define two elements $\mathbf{w}_{i}^{\left(  1\right)
}\left(  f\right)  $ and $\mathbf{w}_{i}^{\left(  2\right)  }\left(  f\right)
$ of $\mathbb{K}$ by%
\[
\mathbf{w}_{i}^{\left(  1\right)  }\left(  f\right)  =\sum_{\substack{x\in
\widehat{P}_{i};\ y\in\widehat{P}_{i+1};\\y\gtrdot x}}\dfrac{f\left(
x\right)  }{f\left(  y\right)  }\ \ \ \ \ \ \ \ \ \ \text{and}%
\ \ \ \ \ \ \ \ \ \ \mathbf{w}_{i}^{\left(  2\right)  }\left(  f\right)
=\sum_{\substack{x\in\widehat{Q}_{i};\ y\in\widehat{Q}_{i+1};\\y\gtrdot
x}}\dfrac{f\left(  x\right)  }{f\left(  y\right)  }%
\]
(where, of course, $\widehat{P}_{j}$ and $\widehat{Q}_{j}$ are embedded into
$\widehat{PQ}_{j}$ for every $j\in\left\{  0,1,...,n+1\right\}  $ in the
obvious way). These elements $\mathbf{w}_{i}^{\left(  1\right)  }\left(
f\right)  $ and $\mathbf{w}_{i}^{\left(  2\right)  }\left(  f\right)  $ are
defined not for every $f$, but for \textquotedblleft almost
every\textquotedblright\ $f$ in the sense of Zariski topology. We denote the
$\left(  n+1\right)  $-tuple
\[
\left(  \mathbf{w}_{0}^{\left(  1\right)  }\left(  f\right)  \diagup
\mathbf{w}_{0}^{\left(  2\right)  }\left(  f\right)  ,\ \mathbf{w}%
_{1}^{\left(  1\right)  }\left(  f\right)  \diagup\mathbf{w}_{1}^{\left(
2\right)  }\left(  f\right)  ,\ ...,\ \mathbf{w}_{n}^{\left(  1\right)
}\left(  f\right)  \diagup\mathbf{w}_{n}^{\left(  2\right)  }\left(  f\right)
\right)
\]
as the \textit{comparative w-tuple} of the labelling $f$. The advantage of
comparative w-tuples over usual w-tuples is the following fact: If $f$ and $g$
are two homogeneously equivalent $\mathbb{K}$-labellings of $PQ$, then%
\begin{equation}
\left(  \text{the comparative w-tuple of }f\right)  =\left(  \text{the
comparative w-tuple of }g\right)  . \label{pf.PQ.ord2.hgeq}%
\end{equation}
\footnote{\textit{Proof of (\ref{pf.PQ.ord2.hgeq}):} Let $f$ and $g$ be two
homogeneously equivalent zero-free $\mathbb{K}$-labellings of $PQ$. Then, just
as in the proof of Proposition \ref{prop.reconstruct}, we can show that there
exists an $\left(  n+2\right)  $-tuple $\left(  a_{0},a_{1},...,a_{n+1}%
\right)  \in\left(  \mathbb{K}^{\times}\right)  ^{n+1}$ such that $g=\left(
a_{0},a_{1},...,a_{n+1}\right)  \flat f$ (where $\left(  a_{0},a_{1}%
,...,a_{n+1}\right)  \flat f$ is defined as in Definition \ref{def.bemol}).
Consider this $\left(  n+2\right)  $-tuple $\left(  a_{0},a_{1},...,a_{n+1}%
\right)  $.
\par
In order to prove (\ref{pf.PQ.ord2.hgeq}), it is clearly enough to show that
every $i\in\left\{  0,1,...,n\right\}  $ satisfies $\mathbf{w}_{i}^{\left(
1\right)  }\left(  f\right)  \diagup\mathbf{w}_{i}^{\left(  2\right)  }\left(
f\right)  =\mathbf{w}_{i}^{\left(  1\right)  }\left(  g\right)  \diagup
\mathbf{w}_{i}^{\left(  2\right)  }\left(  g\right)  $. But this is easy,
since we can show (similarly to the proof of (\ref{pf.w.scalmult.2})) that
$\mathbf{w}_{i}^{\left(  1\right)  }\left(  g\right)  =\dfrac{a_{i}}{a_{i+1}%
}\mathbf{w}_{i}^{\left(  1\right)  }\left(  f\right)  $ and $\mathbf{w}%
_{i}^{\left(  2\right)  }\left(  g\right)  =\dfrac{a_{i}}{a_{i+1}}%
\mathbf{w}_{i}^{\left(  2\right)  }\left(  f\right)  $. This proves
(\ref{pf.PQ.ord2.hgeq}).} (This has no analogue for regular w-tuples.)
\end{verlong}

It is furthermore easy to see (in analogy to Proposition \ref{prop.wi.R}) that
the map $R_{PQ}$ changes the comparative w-tuple of a $\mathbb{K}$-labelling
by shifting it cyclically.

But it is also easy to see (the nonemptiness of $P$ and $Q$ must be used here)
that there exists some $f\in\mathbb{K}^{\widehat{PQ}}$ such that the ratios
$\mathbf{w}_{i}^{\left(  1\right)  }\left(  f\right)  \diagup\mathbf{w}%
_{i}^{\left(  2\right)  }\left(  f\right)  $ are well-defined and pairwise
distinct for all $i\in\left\{  0,1,...,n\right\}  $ and such that $R^{j}f$ is
well-defined for every $j\in\left\{  0,1,...,\ell\right\}  $. Consider such an
$f$. The ratios $\mathbf{w}_{i}^{\left(  1\right)  }\left(  f\right)
\diagup\mathbf{w}_{i}^{\left(  2\right)  }\left(  f\right)  $ are pairwise
distinct for all $i\in\left\{  0,1,...,n\right\}  $; that is, the comparative
w-tuple of $f$ contains no two equal entries.

\begin{vershort}
Since $\overline{R}_{PQ}^{\ell}=\operatorname*{id}$, we have $\overline
{R}_{PQ}^{\ell}\left(  \pi\left(  f\right)  \right)  =\pi\left(  f\right)  $.
The commutativity of the diagram (\ref{def.hgR.commut}) yields $\overline
{R}_{PQ}^{\ell}\circ\pi=\pi\circ R_{PQ}^{\ell}$. Now,%
\[
\pi\left(  f\right)  =\overline{R}_{PQ}^{\ell}\left(  \pi\left(  f\right)
\right)  =\left(  \underbrace{\overline{R}_{PQ}^{\ell}\circ\pi}_{=\pi\circ
R_{PQ}^{\ell}}\right)  \left(  f\right)  =\left(  \pi\circ R_{PQ}^{\ell
}\right)  \left(  f\right)  =\pi\left(  R_{PQ}^{\ell}f\right)  .
\]
In other words, the labellings $f$ and $R_{PQ}^{\ell}f$ are homogeneously
equivalent. Thus,%
\begin{equation}
\left(  \text{the comparative w-tuple of }f\right)  =\left(  \text{the
comparative w-tuple of }R_{PQ}^{\ell}f\right)  \label{pf.PQ.ord2.2.short}%
\end{equation}
(by (\ref{pf.PQ.ord2.hgeq.short})).
\end{vershort}

\begin{verlong}
Since $\overline{R}_{PQ}^{\ell}=\operatorname*{id}$, we have $\overline
{R}_{PQ}^{\ell}\left(  \pi\left(  f\right)  \right)  =\operatorname*{id}%
\left(  \pi\left(  f\right)  \right)  =\pi\left(  f\right)  $. The
commutativity of the diagram (\ref{def.hgR.commut}) yields $\overline{R}%
_{PQ}^{\ell}\circ\pi=\pi\circ R_{PQ}^{\ell}$. Now,%
\[
\pi\left(  f\right)  =\overline{R}_{PQ}^{\ell}\left(  \pi\left(  f\right)
\right)  =\left(  \underbrace{\overline{R}_{PQ}^{\ell}\circ\pi}_{=\pi\circ
R_{PQ}^{\ell}}\right)  \left(  f\right)  =\left(  \pi\circ R_{PQ}^{\ell
}\right)  \left(  f\right)  =\pi\left(  R_{PQ}^{\ell}f\right)  .
\]
In other words, the labellings $f$ and $R_{PQ}^{\ell}f$ are homogeneously
equivalent. Thus,%
\begin{equation}
\left(  \text{the comparative w-tuple of }f\right)  =\left(  \text{the
comparative w-tuple of }R_{PQ}^{\ell}f\right)  \label{pf.PQ.ord2.2}%
\end{equation}
(by (\ref{pf.PQ.ord2.hgeq})).
\end{verlong}

\begin{vershort}
Now, recall that the map $R_{PQ}$ changes the comparative w-tuple of a
$\mathbb{K}$-labelling by shifting it cyclically. Hence, for every
$k\in\mathbb{N}$, the map $R_{PQ}^{k}$ changes the comparative w-tuple of a
$\mathbb{K}$-labelling by shifting it cyclically $k$ times. Applying this to
the $\mathbb{K}$-labelling $f$ and to $k=\ell$, we see that the comparative
w-tuple of $R_{PQ}^{\ell}f$ is obtained from the comparative w-tuple of $f$ by
an $\ell$-fold cyclic shift. Due to (\ref{pf.PQ.ord2.2.short}), this rewrites
as follows: The comparative w-tuple of $f$ is obtained from the comparative
w-tuple of $f$ by an $\ell$-fold cyclic shift. In other words, the comparative
w-tuple of $f$ is invariant under an $\ell$-fold cyclic shift. But since the
comparative w-tuple of $f$ consists of $n+1$ pairwise distinct entries, this
is impossible unless $n+1\mid\ell$. Hence, we must have $n+1\mid\ell$.
\end{vershort}

\begin{verlong}
Consider the group $\mathfrak{S}$ defined in the proof of Proposition
\ref{prop.wi.R}, and its action on the $\left(  n+1\right)  $-tuples of
elements of $\mathbb{K}$ defined ibidem. Also, define the element $\left(
0,1,...,n\right)  $ of $\mathfrak{S}$ just as in the proof of Proposition
\ref{prop.wi.R}. We know that the map $R_{PQ}$ changes the comparative w-tuple
of a $\mathbb{K}$-labelling by shifting it cyclically. In other words, almost
every $\mathbb{K}$-labelling $g$ of $PQ$ satisfies%
\[
\left(  \text{the comparative w-tuple of }R_{PQ}g\right)  =\left(
0,1,...,n\right)  \left(  \text{the comparative w-tuple of }g\right)  .
\]
From this, it is easy to see that for every $k\in\mathbb{N}$, almost every
$\mathbb{K}$-labelling $g$ of $PQ$ satisfies%
\[
\left(  \text{the comparative w-tuple of }R_{PQ}^{k}g\right)  =\left(
0,1,...,n\right)  ^{k}\left(  \text{the comparative w-tuple of }g\right)
\]
(indeed, this can be proven by induction over $k$). Applying this to $g=f$ and
$k=\ell$, we obtain%
\[
\left(  \text{the comparative w-tuple of }R_{PQ}^{\ell}f\right)  =\left(
0,1,...,n\right)  ^{\ell}\left(  \text{the comparative w-tuple of }f\right)
.
\]
Compared with (\ref{pf.PQ.ord2.2}), this yields%
\[
\left(  0,1,...,n\right)  ^{\ell}\left(  \text{the comparative w-tuple of
}f\right)  =\left(  \text{the comparative w-tuple of }f\right)  .
\]
In other words, the comparative w-tuple of $f$ is fixed by the permutation
$\left(  0,1,...,n\right)  ^{\ell}\in\mathfrak{S}$. If $n+1\nmid\ell$, then
this can only happen if the comparative w-tuple of $f$ has two equal entries
(because if $n+1\nmid\ell$, then the permutation $\left(  0,1,...,n\right)
^{\ell}\in\mathfrak{S}$ is not the identity); but this is impossible (since
the comparative w-tuple of $f$ contains no two equal entries). Hence, we
cannot have $n+1\nmid\ell$. Thus, $n+1\mid\ell$.
\end{verlong}

Now,%
\[
\operatorname*{ord}\left(  R_{PQ}\right)  =\operatorname{lcm}\left(
n+1,\operatorname*{ord}\left(  \overline{R}_{PQ}\right)  \right)
=\operatorname*{ord}\left(  \overline{R}_{PQ}\right)
\]
(since $n+1\mid\ell=\operatorname*{ord}\left(  \overline{R}_{PQ}\right)  $).
Hence,%
\[
\operatorname*{ord}\left(  \overline{R}_{PQ}\right)  =\operatorname*{ord}%
\left(  R_{PQ}\right)  =\operatorname{lcm}\left(  \operatorname*{ord}\left(
R_{P}\right)  ,\operatorname*{ord}\left(  R_{Q}\right)  \right)  .
\]
This proves Proposition \ref{prop.PQ.ord2}.
\end{proof}

Now, let us track the effect of $B_{k}$ on the order of $\overline{R}$:

\begin{proposition}
\label{prop.Bk.ord}Let $n\in\mathbb{N}$. Let $P$ be an $n$-graded poset. Let
$\mathbb{K}$ be a field.

\textbf{(a)} We have $\operatorname*{ord}\left(  \overline{R}_{B_{1}P}\right)
=\operatorname*{ord}\left(  \overline{R}_{P}\right)  $.

\textbf{(b)} For every integer $k>1$, we have $\operatorname*{ord}\left(
\overline{R}_{B_{k}P}\right)  =\operatorname{lcm}\left(  2,\operatorname*{ord}%
\left(  \overline{R}_{P}\right)  \right)  $.
\end{proposition}

\begin{proof}
[Proof of Proposition \ref{prop.Bk.ord} (sketched).]We will be proving parts
\textbf{(a)} and \textbf{(b)} together. Let $k$ be a positive integer (this
has to be $1$ for proving part \textbf{(a)}). We need to prove that%
\begin{equation}
\operatorname*{ord}\left(  \overline{R}_{B_{k}P}\right)  =\left\{
\begin{array}
[c]{l}%
\operatorname{lcm}\left(  2,\operatorname*{ord}\left(  \overline{R}%
_{P}\right)  \right)  ,\ \ \ \ \ \ \ \ \ \ \text{if }k>1;\\
\operatorname*{ord}\left(  \overline{R}_{P}\right)
,\ \ \ \ \ \ \ \ \ \ \text{if }k=1
\end{array}
\right.  . \label{pf.Bk.ord.1}%
\end{equation}
Proving this clearly will prove both parts \textbf{(a)} and \textbf{(b)} of
Proposition \ref{prop.Bk.ord}.

\begin{vershort}
Let us make some conventions:

\begin{itemize}
\item For any $n$-tuple $\left(  \alpha_{1},\alpha_{2},...,\alpha_{n}\right)
$ and any object $\beta$, let $\beta\rightthreetimes\alpha$ denote the
$\left(  n+1\right)  $-tuple $\left(  \beta,\alpha_{1},\alpha_{2}%
,...,\alpha_{n}\right)  $.

\item We are going to identify $P$ with a subposet of $B_{k}P$ in the obvious
way. But of course, the degree map of $B_{k}P$ restricted to $P$ is not
identical with the degree map of $P$ (but rather differs from it by $1$), so
we will have to distinguish between \textquotedblleft degree in $P$%
\textquotedblright\ and \textquotedblleft degree in $B_{k}P$\textquotedblright%
. We identify the elements $0$ and $1$ of $\widehat{P}$ with the elements $0$
and $1$ of $\widehat{B_{k}P}$, respectively. Thus, $\widehat{P}$ becomes a
subposet of $\widehat{B_{k}P}$. However, it is not generally true that every
$u\lessdot v$ in $\widehat{P}$ must satisfy $u\lessdot v$ in $\widehat{B_{k}%
P}$.

\item We have a rational map $\pi:\mathbb{K}^{\widehat{P}}\dashrightarrow
\overline{\mathbb{K}^{\widehat{P}}}$ and a rational map $\pi:\mathbb{K}%
^{\widehat{B_{k}P}}\dashrightarrow\overline{\mathbb{K}^{\widehat{B_{k}P}}}$
denoted by the same letter. This is not problematic, because these two maps
can be distinguished by their different domains. We will also use the letter
$\pi$ to denote the rational map $\mathbb{K}^{k}\dashrightarrow\mathbb{P}%
\left(  \mathbb{K}^{k}\right)  $ obtained from the canonical projection
$\mathbb{K}^{k}\setminus\left\{  0\right\}  \rightarrow\mathbb{P}\left(
\mathbb{K}^{k}\right)  $ of the nonzero vectors in $\mathbb{K}^{k}$ onto the
projective space.
\end{itemize}

Now, we recall that the construction of $B_{k}P$ from $P$ involved adding $k$
new (pairwise incomparable) elements smaller than all existing elements of $P$
to the poset. This operation clearly raises the degree of every element of $P$
by $1$\ \ \ \ \footnote{In terms of the Hasse diagram, this can be regarded as
the $k$ new elements ``bumping up'' all existing elements of $P$ by $1$
degree.}, whereas the $k$ newly added elements all obtain degree $1$ in
$B_{k}P$. Formally speaking, this means that $\widehat{B_{k}P}_{i}%
=\widehat{P}_{i-1}$ for every $i\in\left\{  2,3,...,n+1\right\}  $, while
$\widehat{B_{k}P}_{1}$ is a $k$-element set. Moreover, for any $i\in\left\{
2,3,...,n+1\right\}  $, any $u\in\widehat{B_{k}P}_{i}=\widehat{P}_{i-1}$ and
any $v\in\widehat{B_{k}P}_{i+1}=\widehat{P}_{i}$, we have
\[
u\lessdot v\text{ in }\widehat{B_{k}P}\text{ if and only if }u\lessdot v\text{
in }\widehat{P}\text{.}%
\]
(This would not be true if we would allow $i=1$, $u\in\widehat{P}_{0}$ and
$v\in\widehat{P}_{1}$.)
\end{vershort}

\begin{verlong}
Let us make some conventions:

\begin{itemize}
\item For any $n$-tuple $\alpha$ and any object $\beta$, let $\beta
\rightthreetimes\alpha$ denote the $\left(  n+1\right)  $-tuple $\left(
\beta,\alpha_{1},\alpha_{2},...,\alpha_{n}\right)  $, where $\left(
\alpha_{1},\alpha_{2},...,\alpha_{n}\right)  =\alpha$. (In other words,
$\beta\rightthreetimes\alpha$ is obtained by prepending $\beta$ to $\alpha$.)

\item We are going to identify $P$ with a subposet of $B_{k}P$ in the obvious
way. But of course, the degree map of $B_{k}P$ restricted to $P$ is not
identical with the degree map of $P$ (but rather differs from it by $1$).
Therefore, we cannot denote both of these degree maps by $\deg$ (as we would
have normally done) without running into ambiguities. Instead, we shall denote
the degree map of $B_{k}P$ by $\deg_{B_{k}P}$, and we shall denote the degree
map of $P$ by $\deg_{P}$. We identify the elements $0$ and $1$ of
$\widehat{P}$ with the elements $0$ and $1$ of $\widehat{B_{k}P}$,
respectively. Thus, $\widehat{P}$ becomes a subposet of $\widehat{B_{k}P}$.
However, it is not generally true that every $u\in\widehat{P}$ and
$v\in\widehat{P}$ such that $u\lessdot v$ in $\widehat{P}$ must satisfy
$u\lessdot v$ in $\widehat{B_{k}P}$.

\item We have a rational map $\pi:\mathbb{K}^{\widehat{P}}\dashrightarrow
\overline{\mathbb{K}^{\widehat{P}}}$ and a rational map $\pi:\mathbb{K}%
^{\widehat{B_{k}P}}\dashrightarrow\overline{\mathbb{K}^{\widehat{B_{k}P}}}$
denoted by the same letter. This is not problematic, because these two maps
can be distinguished by their different domains. We will also use the letter
$\pi$ to denote the rational map $\mathbb{K}^{k}\dashrightarrow\mathbb{P}%
\left(  \mathbb{K}^{k}\right)  $ obtained from the canonical projection
$\mathbb{K}^{k}\setminus\left\{  0\right\}  \rightarrow\mathbb{P}\left(
\mathbb{K}^{k}\right)  $ of the nonzero vectors in $\mathbb{K}^{k}$ onto the
projective space.
\end{itemize}

Now, we recall that the construction of $B_{k}P$ from $P$ involved adding $k$
new (pairwise incomparable) elements smaller than all existing elements of $P$
to the poset. This operation clearly raises the degree of every element of $P$
by $1$\ \ \ \ \footnote{In terms of the Hasse diagram, this can be regarded as
the $k$ new elements ``bumping up'' all existing elements of $P$ by $1$
\par
degree.}, whereas the $k$ newly added elements all obtain degree $1$ in
$B_{k}P$. Formally speaking, this means that $\widehat{B_{k}P}_{i}%
=\widehat{P}_{i-1}$ for every $i\in\left\{  2,3,...,n+1\right\}  $, while
$\widehat{B_{k}P}_{1}$ is a $k$-element set. Moreover, for any $i\in\left\{
2,3,...,n+1\right\}  $, any $u\in\widehat{B_{k}P}_{i}=\widehat{P}_{i-1}$ and
any $v\in\widehat{B_{k}P}_{i+1}=\widehat{P}_{i}$, we have
\begin{equation}
u\lessdot v\text{ in }\widehat{B_{k}P}\text{ if and only if }u\lessdot v\text{
in }\widehat{P}\text{.} \label{pf.Bk.ord.covercover}%
\end{equation}
(This would not be true if we would allow $i=1$, $u\in\widehat{P}_{0}$ and
$v\in\widehat{P}_{1}$.)
\end{verlong}

We have $\mathbb{K}^{\widehat{B_{k}P}_{i}}=\mathbb{K}^{\widehat{P}_{i-1}}$ for
every $i\in\left\{  2,3,...,n+1\right\}  $ (since $\widehat{B_{k}P}%
_{i}=\widehat{P}_{i-1}$ for every $i\in\left\{  2,3,...,n+1\right\}  $),
whereas $\mathbb{K}^{\widehat{B_{k}P}_{1}}\cong\mathbb{K}^{k}$ (since
$\widehat{B_{k}P}_{1}$ is a $k$-element set). We will actually identify
$\mathbb{K}^{\widehat{B_{k}P}_{1}}$ with $\mathbb{K}^{k}$. Now,%
\begin{align}
\overline{\mathbb{K}^{\widehat{B_{k}P}}}  &  =\prod\limits_{i=1}%
^{n+1}\mathbb{P}\left(  \mathbb{K}^{\widehat{B_{k}P}_{i}}\right)
=\mathbb{P}\left(  \underbrace{\mathbb{K}^{\widehat{B_{k}P}_{1}}}%
_{=\mathbb{K}^{k}}\right)  \times\prod\limits_{i=2}^{n+1}\mathbb{P}\left(
\underbrace{\mathbb{K}^{\widehat{B_{k}P}_{i}}}_{=\mathbb{K}^{\widehat{P}%
_{i-1}}}\right) \nonumber\\
&  =\mathbb{P}\left(  \mathbb{K}^{k}\right)  \times\prod\limits_{i=2}%
^{n+1}\mathbb{P}\left(  \mathbb{K}^{\widehat{P}_{i-1}}\right)  =\mathbb{P}%
\left(  \mathbb{K}^{k}\right)  \times\underbrace{\prod\limits_{i=1}%
^{n}\mathbb{P}\left(  \mathbb{K}^{\widehat{P}_{i}}\right)  }_{=\overline
{\mathbb{K}^{\widehat{P}}}}=\mathbb{P}\left(  \mathbb{K}^{k}\right)
\times\overline{\mathbb{K}^{\widehat{P}}}. \label{pf.Bk.ord.decomp}%
\end{align}
Thus, the elements of $\overline{\mathbb{K}^{\widehat{B_{k}P}}}$ have the form
$\widetilde{p}\rightthreetimes\widetilde{g}$, where $\widetilde{p}%
\in\mathbb{P}\left(  \mathbb{K}^{k}\right)  $ and $\widetilde{g}\in
\overline{\mathbb{K}^{\widehat{P}}}$.

On the other hand, recall that $\widehat{P}$ is a subposet of $\widehat{B_{k}%
P}$. More precisely, $\widehat{P}$ is the subposet $\widehat{B_{k}P}%
\setminus\widehat{B_{k}P}_{1}$ of $\widehat{B_{k}P}$. Thus, we can define a
map $\Phi:\mathbb{K}^{k}\times\mathbb{K}^{\widehat{P}}\rightarrow
\mathbb{K}^{\widehat{B_{k}P}}$ by setting%
\[
\left(  \Phi\left(  p,g\right)  \right)  \left(  v\right)  =\left\{
\begin{array}
[c]{c}%
p\left(  v\right)  ,\ \ \ \ \ \ \ \ \ \ \text{if }v\in\widehat{B_{k}P}_{1};\\
g\left(  v\right)  ,\ \ \ \ \ \ \ \ \ \ \text{if }v\notin\widehat{B_{k}P}_{1}%
\end{array}
\right.  \ \ \ \ \ \ \ \ \ \ \text{for every }v\in\widehat{B_{k}P}%
\]
for every $\left(  p,g\right)  \in\mathbb{K}^{k}\times\mathbb{K}^{\widehat{P}%
}$. Here, the term $p\left(  v\right)  $ is to be understood by means of
regarding $p$ as an element of $\mathbb{K}^{\widehat{B_{k}P}_{1}}$ (since
$p\in\mathbb{K}^{k}=\mathbb{K}^{\widehat{B_{k}P}_{1}}$). Clearly, $\Phi$ is a
bijection. Moreover, it is easy to see that
\begin{equation}
\pi\left(  \Phi\left(  p,g\right)  \right)  =\pi\left(  p\right)
\rightthreetimes\pi\left(  g\right)  \ \ \ \ \ \ \ \ \ \ \text{for all }%
p\in\mathbb{K}^{k}\text{ and }g\in\mathbb{K}^{\widehat{P}}
\label{pf.Bk.ord.pitriv}%
\end{equation}
(where the $\pi$ on the left hand side is the map $\pi:\mathbb{K}%
^{\widehat{B_{k}P}}\dashrightarrow\overline{\mathbb{K}^{\widehat{B_{k}P}}}$,
whereas the $\pi$ in \textquotedblleft$\pi\left(  p\right)  $%
\textquotedblright\ is the map $\pi:\mathbb{K}^{k}\dashrightarrow
\mathbb{P}\left(  \mathbb{K}^{k}\right)  $, and the $\pi$ in \textquotedblleft%
$\pi\left(  g\right)  $\textquotedblright\ is the map $\pi:\mathbb{K}%
^{\widehat{P}}\dashrightarrow\overline{\mathbb{K}^{\widehat{P}}}$).

\begin{verlong}
\textit{Proof of (\ref{pf.Bk.ord.pitriv}):} Let $p\in\mathbb{K}^{k}$ and
$g\in\mathbb{K}^{\widehat{P}}$. From the definition of $\Phi\left(
p,g\right)  $, it is clear that $\pi_{1}\left(  \Phi\left(  p,g\right)
\right)  =\pi\left(  p\right)  $ and that every $i\in\left\{
2,3,...,n+1\right\}  $ satisfies $\pi_{i}\left(  \Phi\left(  p,g\right)
\right)  =\pi_{i-1}\left(  g\right)  $ (because every $i\in\left\{
2,3,...,n+1\right\}  $ satisfies $\widehat{B_{k}P}_{i}=\widehat{P}_{i-1}$).
But (\ref{def.hgeq.pii.eq}) (applied to $n+1$, $B_{k}P$ and $\Phi\left(
p,g\right)  $ instead of $n$, $P$ and $f$) yields%
\begin{align*}
\pi\left(  \Phi\left(  p,g\right)  \right)   &  =\left(  \pi_{1}\left(
\Phi\left(  p,g\right)  \right)  ,\pi_{2}\left(  \Phi\left(  p,g\right)
\right)  ,...,\pi_{n+1}\left(  \Phi\left(  p,g\right)  \right)  \right) \\
&  =\underbrace{\pi_{1}\left(  \Phi\left(  p,g\right)  \right)  }_{=\pi\left(
p\right)  }\rightthreetimes\underbrace{\left(  \pi_{2}\left(  \Phi\left(
p,g\right)  \right)  ,\pi_{3}\left(  \Phi\left(  p,g\right)  \right)
,...,\pi_{n+1}\left(  \Phi\left(  p,g\right)  \right)  \right)  }%
_{\substack{=\left(  \pi_{1}\left(  g\right)  ,\pi_{2}\left(  g\right)
,...,\pi_{n}\left(  g\right)  \right)  \\\text{(since every }i\in\left\{
2,3,...,n+1\right\}  \text{ satisfies }\pi_{i}\left(  \Phi\left(  p,g\right)
\right)  =\pi_{i-1}\left(  g\right)  \text{)}}}\\
&  \ \ \ \ \ \ \ \ \ \ \left(
\begin{array}
[c]{c}%
\text{by the definition of}\\
\pi_{1}\left(  \Phi\left(  p,g\right)  \right)  \rightthreetimes\left(
\pi_{2}\left(  \Phi\left(  p,g\right)  \right)  ,\pi_{3}\left(  \Phi\left(
p,g\right)  \right)  ,...,\pi_{n+1}\left(  \Phi\left(  p,g\right)  \right)
\right)
\end{array}
\right) \\
&  =\pi\left(  p\right)  \rightthreetimes\underbrace{\left(  \pi_{1}\left(
g\right)  ,\pi_{2}\left(  g\right)  ,...,\pi_{n}\left(  g\right)  \right)
}_{\substack{=\pi\left(  g\right)  \\\text{(by (\ref{def.hgeq.pii.eq}),
applied to }f=g\text{)}}}=\pi\left(  p\right)  \rightthreetimes\pi\left(
g\right)  .
\end{align*}
This proves (\ref{pf.Bk.ord.pitriv}).
\end{verlong}

Now, we claim that every $\widetilde{p}\in\mathbb{P}\left(  \mathbb{K}%
^{k}\right)  $ and $\widetilde{g}\in\overline{\mathbb{K}^{\widehat{P}}}$
satisfy
\begin{equation}
\left(\overline{R_{i}}\right)_{B_{k}P}\left(
\widetilde{p}\rightthreetimes\widetilde{g}%
\right)
= \widetilde{p} \rightthreetimes
\left(\overline{R_{i-1}}\right)_{P}
\left( \widetilde{g}\right)
 \ \ \ \ \ \ \ \ \ \ \text{for all }i\in\left\{
2,3,...,n+1\right\}
\label{pf.Bk.ord.Ri.big}
\end{equation}
and%
\begin{equation}
\left(\overline{R_{1}}\right)_{B_{k}P}\left(
\widetilde{p}\rightthreetimes\widetilde{g}%
\right)  =\widetilde{p}^{-1}\rightthreetimes\widetilde{g}.
\label{pf.Bk.ord.Ri.small}%
\end{equation}

\begin{vershort}
\textit{Proof of (\ref{pf.Bk.ord.Ri.big}) and (\ref{pf.Bk.ord.Ri.small}):} In
order to prove (\ref{pf.Bk.ord.Ri.big}), it is clearly enough to show that
every $p\in\mathbb{K}^{k}$ and $g\in\mathbb{K}^{\widehat{P}}$ satisfy%
\begin{equation}
\left(  R_{i}\right)  _{B_{k}P}\left(  p\rightthreetimes g\right)  \sim
p\rightthreetimes\left(  R_{i-1}\right)  _{P}\left(  g\right)
\ \ \ \ \ \ \ \ \ \ \text{for all }i\in\left\{  2,3,...,n+1\right\}  ,
\label{pf.Bk.ord.Ri.big.pf.short.1}%
\end{equation}
where the sign $\sim$ stands for homogeneous equivalence.

It is easy to prove the relation (\ref{pf.Bk.ord.Ri.big.pf.short.1}) for $i>2$
(because if $i>2$, then the elements of $\widehat{B_{k}P}$ having degrees
$i-1$, $i$ and $i+1$ are precisely the elements of $\widehat{P}$ having
degrees $i-2$, $i-1$ and $i$, and therefore toggling the elements of
$\widehat{B_{k}P}_{i}$ in $p\rightthreetimes g$ has precisely the same effect
as toggling the elements of $\widehat{P}_{i-1}$ in $g$ while leaving $p$
fixed, so that we even get the stronger assertion that $\left(  R_{i}\right)
_{B_{k}P}\left(  p\rightthreetimes g\right)  =p\rightthreetimes\left(
R_{i-1}\right)  _{P}\left(  g\right)  $). It is not much harder to check that
it also holds for $i=2$ (indeed, for $i=2$, the only difference between
toggling the elements of $\widehat{B_{k}P}_{i}$ in $p\rightthreetimes g$ and
toggling the elements of $\widehat{P}_{i-1}$ in $g$ while leaving $p$ fixed is
a scalar factor which is identical across all elements being toggled in either
poset\footnote{because every $u\in\widehat{B_{k}P}_{1}$ and every
$v\in\widehat{B_{k}P}_{2}$ satisfy $u\lessdot v$}; therefore the results are
the same up to homogeneous equivalence).

Finally, (\ref{pf.Bk.ord.Ri.small}) is trivial to check (e.g., using Corollary
\ref{cor.hgRi.1}).
\end{vershort}

\begin{verlong}
\textit{Proof of (\ref{pf.Bk.ord.Ri.big}):} Let $i\in\left\{
2,3,...,n+1\right\}  $. Let $\widetilde{p}\in\mathbb{P}\left(  \mathbb{K}%
^{k}\right)  $ and $\widetilde{g}\in\overline{\mathbb{K}^{\widehat{P}}}$.

Since $\widetilde{g}\in\overline{\mathbb{K}^{\widehat{P}}}$, there exists some
$g\in\mathbb{K}^{\widehat{P}}$ such that $\widetilde{g}=\pi\left(  g\right)
$. Consider this $g$.

Since $\widetilde{p}\in\mathbb{P}\left(  \mathbb{K}^{k}\right)  $, there
exists some nonzero $p\in\mathbb{K}^{k}$ such that $\widetilde{p}=\pi\left(
p\right)  $. Consider this $p$.

We will now show that $R_{i}\left(  \Phi\left(  p,g\right)  \right)  $ is
homogeneously equivalent to $\Phi\left(  p,R_{i-1}g\right)  $.

Let us first notice that $i$ and $i+1$ both belong to $\left\{
2,3,...,n+2\right\}  $. Hence, $\widehat{B_{k}P}_{i}=\widehat{P}_{i-1}$ and
$\widehat{B_{k}P}_{i+1}=\widehat{P}_{i}$.

Now, we distinguish between two cases:

\textit{Case 1:} We have $i>2$.

\textit{Case 2:} We have $i=2$.

Let us consider Case 1 first. In this case, $i>2$. Hence, $i-1\in\left\{
2,3,...,n+1\right\}  $, so that $\widehat{B_{k}P}_{i-1}=\widehat{P}_{i-2}$.

We are going to prove that
\begin{equation}
R_{i}\left(  \Phi\left(  p,g\right)  \right)  =\Phi\left(  p,R_{i-1}g\right)
. \label{pf.Bk.ord.Ri.big.c1.goal}%
\end{equation}

Indeed, let $v\in\widehat{B_{k}P}$. Then, we are going to prove that
\[
\left(  R_{i}\left(  \Phi\left(  p,g\right)  \right)  \right)  \left(
v\right)  =\left(  \Phi\left(  p,R_{i-1}g\right)  \right)  \left(  v\right)
.
\]

If $v\notin\widehat{B_{k}P}_{i}$, then $\left(  R_{i}\left(  \Phi\left(
p,g\right)  \right)  \right)  \left(  v\right)  =\left(  \Phi\left(
p,R_{i-1}g\right)  \right)  \left(  v\right)  $ is very easy to
check\footnote{\textit{Proof.} Assume that $v\notin\widehat{B_{k}P}_{i}$.
Then, $\deg_{B_{k}P}v\neq i$. Thus, Proposition \ref{prop.Ri.implicit}
\textbf{(a)} (applied to $B_{k}P$, $n+1$ and $\Phi\left(  p,g\right)  $
instead of $P$, $n$ and $f$) yields%
\begin{equation}
\left(  R_{i}\left(  \Phi\left(  p,g\right)  \right)  \right)  \left(
v\right)  =\left(  \Phi\left(  p,g\right)  \right)  \left(  v\right)
=\left\{
\begin{array}
[c]{c}%
p\left(  v\right)  ,\ \ \ \ \ \ \ \ \ \ \text{if }v\in\widehat{B_{k}P}_{1};\\
g\left(  v\right)  ,\ \ \ \ \ \ \ \ \ \ \text{if }v\notin\widehat{B_{k}P}_{1}%
\end{array}
\right.  \label{pf.Bk.ord.Ri.big.c1.goal.0}%
\end{equation}
(by the definition of $\Phi\left(  p,g\right)  $).
\par
On the other hand, $v\notin\widehat{B_{k}P}_{i}=\widehat{P}_{i-1}$. Hence, if
$v\notin\widehat{B_{k}P}_{1}$, then $\deg_{P}v\neq i-1$ (notice that $\deg
_{P}v$ is well-defined in this case, because $v\in P$ (since $v\notin%
\widehat{B_{k}P}_{1}$)). Hence, Proposition \ref{prop.Ri.implicit}
\textbf{(a)} (applied to $i-1$ and $g$ instead of $i$ and $f$) yields that, if
$v\notin\widehat{B_{k}P}_{1}$, then $\left(  R_{i-1}g\right)  \left(
v\right)  =g\left(  v\right)  $. But the definition of $\Phi\left(
p,R_{i-1}g\right)  $ yields%
\[
\left(  \Phi\left(  p,R_{i-1}g\right)  \right)  \left(  v\right)  =\left\{
\begin{array}
[c]{c}%
p\left(  v\right)  ,\ \ \ \ \ \ \ \ \ \ \text{if }v\in\widehat{B_{k}P}_{1};\\
\left(  R_{i-1}g\right)  \left(  v\right)  ,\ \ \ \ \ \ \ \ \ \ \text{if
}v\notin\widehat{B_{k}P}_{1}%
\end{array}
\right.  =\left\{
\begin{array}
[c]{c}%
p\left(  v\right)  ,\ \ \ \ \ \ \ \ \ \ \text{if }v\in\widehat{B_{k}P}_{1};\\
g\left(  v\right)  ,\ \ \ \ \ \ \ \ \ \ \text{if }v\notin\widehat{B_{k}P}_{1}%
\end{array}
\right.
\]
(because if $v\notin\widehat{B_{k}P}_{1}$, then $\left(  R_{i-1}g\right)
\left(  v\right)  =g\left(  v\right)  $). Compared with
(\ref{pf.Bk.ord.Ri.big.c1.goal.0}), this yields $\left(  R_{i}\left(
\Phi\left(  p,g\right)  \right)  \right)  \left(  v\right)  =\left(
\Phi\left(  p,R_{i-1}g\right)  \right)  \left(  v\right)  $, qed.}. Hence, for
the rest of the proof of $\left(  R_{i}\left(  \Phi\left(  p,g\right)
\right)  \right)  \left(  v\right)  =\left(  \Phi\left(  R_{i-1}p,g\right)
\right)  \left(  v\right)  $, we can WLOG assume that we don't have
$v\notin\widehat{B_{k}P}_{i}$. Assume this. Thus, $v\in\widehat{B_{k}P}_{i}$.
Hence, $\deg_{B_{k}P}v=i$. Thus, Proposition \ref{prop.Ri.implicit}
\textbf{(b)} (applied to $B_{k}P$, $n+1$ and $\Phi\left(  p,g\right)  $
instead of $P$, $n$ and $f$) yields%
\begin{equation}
\left(  R_{i}\left(  \Phi\left(  p,g\right)  \right)  \right)  \left(
v\right)  =\dfrac{1}{\left(  \Phi\left(  p,g\right)  \right)  \left(
v\right)  }\cdot\dfrac{\sum\limits_{\substack{u\in\widehat{B_{k}P};\\u\lessdot
v}}\left(  \Phi\left(  p,g\right)  \right)  \left(  u\right)  }{\sum
\limits_{\substack{u\in\widehat{B_{k}P};\\u\gtrdot v}}\dfrac{1}{\left(
\Phi\left(  p,g\right)  \right)  \left(  u\right)  }}.
\label{pf.Bk.ord.Ri.big.c1.goal.1}%
\end{equation}
We are now going to simplify the right hand side of this equality.

First of all, $\deg_{B_{k}P}v=i>2>1$, so that $v\notin\widehat{B_{k}P}_{1}$.
Hence, the definition of $\Phi\left(  p,g\right)  $ yields%
\begin{equation}
\left(  \Phi\left(  p,g\right)  \right)  \left(  v\right)  =\left\{
\begin{array}
[c]{c}%
p\left(  v\right)  ,\ \ \ \ \ \ \ \ \ \ \text{if }v\in\widehat{B_{k}P}_{1};\\
g\left(  v\right)  ,\ \ \ \ \ \ \ \ \ \ \text{if }v\notin\widehat{B_{k}P}_{1}%
\end{array}
\right.  =g\left(  v\right)  \label{pf.Bk.ord.Ri.big.c1.goal.2}%
\end{equation}
(since $v\notin\widehat{B_{k}P}_{1}$).

But recall that $\widehat{B_{k}P}_{i}=\widehat{P}_{i-1}$ and $\widehat{B_{k}%
P}_{i-1}=\widehat{P}_{i-2}$. Moreover, we have $i-1\in\left\{
2,3,...,n+1\right\}  $. Hence, any $u\in\widehat{B_{k}P}_{i-1}=\widehat{P}%
_{i-2}$ and any $v\in\widehat{B_{k}P}_{i}=\widehat{P}_{i+1}$ satisfies%
\[
u\lessdot v\text{ in }\widehat{B_{k}P}\text{ if and only if }u\lessdot v\text{
in }\widehat{P}%
\]
(by (\ref{pf.Bk.ord.covercover}), applied to $i-1$ instead of $i$). Thus, the
elements $u\in\widehat{B_{k}P}$ satisfying $u\lessdot v$ are precisely the
elements $u\in\widehat{P}$ satisfying $u\lessdot v$. For similar (but even
simpler) reasons, we can see that the elements $u\in\widehat{B_{k}P}$
satisfying $u\gtrdot v$ are precisely the elements $u\in\widehat{P}$
satisfying $u\gtrdot v$.

Every $u\in\widehat{B_{k}P}$ satisfying $u\lessdot v$ satisfies%
\begin{align*}
\left(  \Phi\left(  p,g\right)  \right)  \left(  u\right)   &  =\left\{
\begin{array}
[c]{c}%
p\left(  u\right)  ,\ \ \ \ \ \ \ \ \ \ \text{if }u\in\widehat{B_{k}P}_{1};\\
g\left(  u\right)  ,\ \ \ \ \ \ \ \ \ \ \text{if }u\notin\widehat{B_{k}P}_{1}%
\end{array}
\right.  \ \ \ \ \ \ \ \ \ \ \left(  \text{by the definition of }\Phi\left(
p,g\right)  \right) \\
&  =g\left(  u\right)  \ \ \ \ \ \ \ \ \ \ \left(
\begin{array}
[c]{c}%
\text{since }u\lessdot v\text{, hence }\deg_{B_{k}P}u=\underbrace{\deg
_{B_{k}P}v}_{=i>2}-1>2-1=1\text{,}\\
\text{so that }u\notin\widehat{B_{k}P}_{1}%
\end{array}
\right)  .
\end{align*}
Hence,%
\begin{equation}
\sum\limits_{\substack{u\in\widehat{B_{k}P};\\u\lessdot v}}\left(  \Phi\left(
p,g\right)  \right)  \left(  u\right)  =\sum\limits_{\substack{u\in
\widehat{B_{k}P};\\u\lessdot v}}g\left(  u\right)  =\sum
\limits_{\substack{u\in\widehat{P};\\u\lessdot v}}g\left(  u\right)
\label{pf.Bk.ord.Ri.big.c1.goal.3}%
\end{equation}
(since the elements $u\in\widehat{B_{k}P}$ satisfying $u\lessdot v$ are
precisely the elements $u\in\widehat{P}$ satisfying $u\lessdot v$).

Similarly,
\begin{equation}
\sum\limits_{\substack{u\in\widehat{B_{k}P};\\u\gtrdot v}}\dfrac{1}{\left(
\Phi\left(  p,g\right)  \right)  \left(  u\right)  }=\sum
\limits_{\substack{u\in\widehat{P};\\u\gtrdot v}}\dfrac{1}{g\left(  u\right)
} \label{pf.Bk.ord.Ri.big.c1.goal.4}%
\end{equation}
Substituting (\ref{pf.Bk.ord.Ri.big.c1.goal.2}),
(\ref{pf.Bk.ord.Ri.big.c1.goal.3}) and (\ref{pf.Bk.ord.Ri.big.c1.goal.4}) into
(\ref{pf.Bk.ord.Ri.big.c1.goal.1}), we obtain%
\begin{equation}
\left(  R_{i}\left(  \Phi\left(  p,g\right)  \right)  \right)  \left(
v\right)  =\dfrac{1}{g\left(  v\right)  }\cdot\dfrac{\sum
\limits_{\substack{u\in\widehat{P};\\u\lessdot v}}g\left(  u\right)  }%
{\sum\limits_{\substack{u\in\widehat{P};\\u\gtrdot v}}\dfrac{1}{g\left(
u\right)  }}. \label{pf.Bk.ord.Ri.big.c1.goal.5}%
\end{equation}
On the other hand, since $v\in\widehat{B_{k}P}_{i}=\widehat{P}_{i-1}$, we have
$\deg_{P}v=i-1$. Hence, Proposition \ref{prop.Ri.implicit} \textbf{(b)}
(applied to $i-1$ and $g$ instead of $i$ and $f$) yields%
\[
\left(  R_{i-1}g\right)  \left(  v\right)  =\dfrac{1}{g\left(  v\right)
}\cdot\dfrac{\sum\limits_{\substack{u\in\widehat{P};\\u\lessdot v}}g\left(
u\right)  }{\sum\limits_{\substack{u\in\widehat{P};\\u\gtrdot v}}\dfrac
{1}{g\left(  u\right)  }}.
\]
Compared with (\ref{pf.Bk.ord.Ri.big.c1.goal.5}), this yields%
\begin{equation}
\left(  R_{i-1}g\right)  \left(  v\right)  =\left(  R_{i}\left(  \Phi\left(
p,g\right)  \right)  \right)  \left(  v\right)  .
\label{pf.Bk.ord.Ri.big.c1.goal.6}%
\end{equation}
But by the definition of $\Phi\left(  p,R_{i-1}g\right)  $, we have%
\begin{align*}
\left(  \Phi\left(  p,R_{i-1}g\right)  \right)  \left(  v\right)   &
=\left\{
\begin{array}
[c]{c}%
p\left(  v\right)  ,\ \ \ \ \ \ \ \ \ \ \text{if }v\in\widehat{B_{k}P}_{1};\\
\left(  R_{i-1}g\right)  \left(  v\right)  ,\ \ \ \ \ \ \ \ \ \ \text{if
}v\notin\widehat{B_{k}P}_{1}%
\end{array}
\right. \\
&  =\left(  R_{i-1}g\right)  \left(  v\right)  \ \ \ \ \ \ \ \ \ \ \left(
\text{since }v\notin\widehat{B_{k}P}_{1}\right) \\
&  =\left(  R_{i}\left(  \Phi\left(  p,g\right)  \right)  \right)  \left(
v\right)  \ \ \ \ \ \ \ \ \ \ \left(  \text{by
(\ref{pf.Bk.ord.Ri.big.c1.goal.6})}\right)  .
\end{align*}
We thus have shown that $\left(  R_{i}\left(  \Phi\left(  p,g\right)  \right)
\right)  \left(  v\right)  =\left(  \Phi\left(  p,R_{i-1}g\right)  \right)
\left(  v\right)  $.

Now, forget that we fixed $v$. We have shown that every $v\in\widehat{B_{k}P}$
satisfies \newline$\left(  R_{i}\left(  \Phi\left(  p,g\right)  \right)
\right)  \left(  v\right)  =\left(  \Phi\left(  p,R_{i-1}g\right)  \right)
\left(  v\right)  $. In other words, $R_{i}\left(  \Phi\left(  p,g\right)
\right)  =\Phi\left(  p,R_{i-1}g\right)  $. Hence, $R_{i}\left(  \Phi\left(
p,g\right)  \right)  $ is homogeneously equivalent to $\Phi\left(
p,R_{i-1}g\right)  $.

This completes the proof that $R_{i}\left(  \Phi\left(  p,g\right)  \right)  $
is homogeneously equivalent to $\Phi\left(  p,R_{i-1}g\right)  $ in Case 1.

Let us now consider Case 2. In this case, $i=2$. Hence, $i-1=1$, so that
$\widehat{B_{k}P}_{i-1}=\widehat{B_{k}P}_{1}$.

We WLOG assume that the elements $g\left(  0\right)  $ and $\sum
\limits_{u\in\widehat{B_{k}P}_{1}}p\left(  u\right)  $ of $\mathbb{K}$ are
nonzero. (This assumption is really WLOG, because we are proving a polynomial
identity, and thus restriction to a Zariski-dense open subset can do no harm.)

Define an $\left(  n+3\right)  $-tuple $\left(  a_{0},a_{1},...,a_{n+2}%
\right)  \in\left(  \mathbb{K}^{\times}\right)  ^{n+3}$ by%
\[
a_{j}=\left\{
\begin{array}
[c]{c}%
1,\ \ \ \ \ \ \ \ \ \ \text{if }j\neq i;\\
\dfrac{1}{g\left(  0\right)  }\sum\limits_{u\in\widehat{B_{k}P}_{1}}p\left(
u\right)  ,\ \ \ \ \ \ \ \ \ \ \text{if }j=i
\end{array}
\right.  \ \ \ \ \ \ \ \ \ \ \text{for every }j\in\left\{
0,1,...,n+2\right\}  .
\]
Define an $h\in\mathbb{K}^{\widehat{B_{k}P}}$ by $h=\left(  a_{0}%
,a_{1},...,a_{n+2}\right)  \flat\left(  \Phi\left(  p,R_{i-1}g\right)
\right)  $ (where \newline$\left(  a_{0},a_{1},...,a_{n+2}\right)
\flat\left(  \Phi\left(  p,R_{i-1}g\right)  \right)  $ is to be understood
according to Definition \ref{def.bemol}).

We are going to prove that
\begin{equation}
R_{i}\left(  \Phi\left(  p,g\right)  \right)  =h.
\label{pf.Bk.ord.Ri.big.c2.goal}%
\end{equation}

Indeed, let $v\in\widehat{B_{k}P}$. Then, we are going to prove that $\left(
R_{i}\left(  \Phi\left(  p,g\right)  \right)  \right)  \left(  v\right)
=h\left(  v\right)  $.

If $v\notin\widehat{B_{k}P}_{i}$, then $\left(  R_{i}\left(  \Phi\left(
p,g\right)  \right)  \right)  \left(  v\right)  =h\left(  v\right)  $ is very
easy to check\footnote{\textit{Proof.} Assume that $v\notin\widehat{B_{k}%
P}_{i}$. Then, $\deg_{B_{k}P}v\neq i$. Thus, Proposition
\ref{prop.Ri.implicit} \textbf{(a)} (applied to $B_{k}P$, $n+1$ and
$\Phi\left(  p,g\right)  $ instead of $P$, $n$ and $f$) yields%
\begin{equation}
\left(  R_{i}\left(  \Phi\left(  p,g\right)  \right)  \right)  \left(
v\right)  =\left(  \Phi\left(  p,g\right)  \right)  \left(  v\right)
=\left\{
\begin{array}
[c]{c}%
p\left(  v\right)  ,\ \ \ \ \ \ \ \ \ \ \text{if }v\in\widehat{B_{k}P}_{1};\\
g\left(  v\right)  ,\ \ \ \ \ \ \ \ \ \ \text{if }v\notin\widehat{B_{k}P}_{1}%
\end{array}
\right.  \label{pf.Bk.ord.Ri.big.c2.goal.0}%
\end{equation}
(by the definition of $\Phi\left(  p,g\right)  $).
\par
On the other hand, $v\notin\widehat{B_{k}P}_{i}=\widehat{P}_{i-1}$. Hence, if
$v\notin\widehat{B_{k}P}_{1}$, then $\deg_{P}v\neq i-1$ (notice that $\deg
_{P}v$ is well-defined in this case, because $v\in P$ (since $v\notin%
\widehat{B_{k}P}_{1}$)). Hence, Proposition \ref{prop.Ri.implicit}
\textbf{(a)} (applied to $i-1$ and $g$ instead of $i$ and $f$) yields that, if
$v\notin\widehat{B_{k}P}_{1}$, then $\left(  R_{i-1}g\right)  \left(
v\right)  =g\left(  v\right)  $.
\par
But $h=\left(  a_{0},a_{1},...,a_{n+2}\right)  \flat\left(  \Phi\left(
p,R_{i-1}g\right)  \right)  $, so that
\begin{align*}
h\left(  v\right)   &  =\left(  \left(  a_{0},a_{1},...,a_{n+2}\right)
\flat\left(  \Phi\left(  p,R_{i-1}g\right)  \right)  \right)  \left(  v\right)
\\
&  =\underbrace{a_{\deg_{B_{k}P}v}}_{\substack{=1\\\text{(since }\deg_{B_{k}%
P}v\neq i\text{)}}}\cdot\left(  \Phi\left(  p,R_{i-1}g\right)  \right)
\left(  v\right)  \ \ \ \ \ \ \ \ \ \ \left(  \text{by the definition of
}\left(  a_{0},a_{1},...,a_{n+2}\right)  \flat\left(  \Phi\left(
p,R_{i-1}g\right)  \right)  \right) \\
&  =\left(  \Phi\left(  p,R_{i-1}g\right)  \right)  \left(  v\right)
=\left\{
\begin{array}
[c]{c}%
p\left(  v\right)  ,\ \ \ \ \ \ \ \ \ \ \text{if }v\in\widehat{B_{k}P}_{1};\\
\left(  R_{i-1}g\right)  \left(  v\right)  ,\ \ \ \ \ \ \ \ \ \ \text{if
}v\notin\widehat{B_{k}P}_{1}%
\end{array}
\right. \\
&  \ \ \ \ \ \ \ \ \ \ \left(  \text{by the definition of }\Phi\left(
p,R_{i-1}g\right)  \right) \\
&  =\left\{
\begin{array}
[c]{c}%
p\left(  v\right)  ,\ \ \ \ \ \ \ \ \ \ \text{if }v\in\widehat{B_{k}P}_{1};\\
g\left(  v\right)  ,\ \ \ \ \ \ \ \ \ \ \text{if }v\notin\widehat{B_{k}P}_{1}%
\end{array}
\right.
\end{align*}
(because if $v\notin\widehat{B_{k}P}_{1}$, then $\left(  R_{i-1}g\right)
\left(  v\right)  =g\left(  v\right)  $). Compared with
(\ref{pf.Bk.ord.Ri.big.c2.goal.0}), this yields $\left(  R_{i}\left(
\Phi\left(  p,g\right)  \right)  \right)  \left(  v\right)  =h\left(
v\right)  $, qed.}. Hence, for the rest of the proof of $\left(  R_{i}\left(
\Phi\left(  p,g\right)  \right)  \right)  \left(  v\right)  =h\left(
v\right)  $, we can WLOG assume that we don't have $v\notin\widehat{B_{k}%
P}_{i}$. Assume this. Thus, $v\in\widehat{B_{k}P}_{i}$. Hence, $\deg_{B_{k}%
P}v=i$. Thus, Proposition \ref{prop.Ri.implicit} \textbf{(b)} (applied to
$B_{k}P$, $n+1$ and $\Phi\left(  p,g\right)  $ instead of $P$, $n$ and $f$)
yields%
\begin{equation}
\left(  R_{i}\left(  \Phi\left(  p,g\right)  \right)  \right)  \left(
v\right)  =\dfrac{1}{\left(  \Phi\left(  p,g\right)  \right)  \left(
v\right)  }\cdot\dfrac{\sum\limits_{\substack{u\in\widehat{B_{k}P};\\u\lessdot
v}}\left(  \Phi\left(  p,g\right)  \right)  \left(  u\right)  }{\sum
\limits_{\substack{u\in\widehat{B_{k}P};\\u\gtrdot v}}\dfrac{1}{\left(
\Phi\left(  p,g\right)  \right)  \left(  u\right)  }}.
\label{pf.Bk.ord.Ri.big.c2.goal.1}%
\end{equation}
We are now going to simplify the right hand side of this equality.

First of all, we can prove the equalities (\ref{pf.Bk.ord.Ri.big.c1.goal.2})
and (\ref{pf.Bk.ord.Ri.big.c1.goal.4}) in the same way as in Case 1. But we
are no longer able to prove (\ref{pf.Bk.ord.Ri.big.c1.goal.3}). Instead, we
have the following argument:

Since $i=2$, we know that the elements $u\in\widehat{B_{k}P}$ satisfying
$u\lessdot v$ are precisely the $k$ elements of $\widehat{B_{k}P}_{1}$
(because the $k$ elements of $\widehat{B_{k}P}_{1}$ are smaller than any other
elements of $B_{k}P$ by the construction of $B_{k}P$). Hence,%
\begin{equation}
\sum\limits_{\substack{u\in\widehat{B_{k}P};\\u\lessdot v}}\left(  \Phi\left(
p,g\right)  \right)  \left(  u\right)  =\sum\limits_{u\in\widehat{B_{k}P}_{1}%
}\left(  \Phi\left(  p,g\right)  \right)  \left(  u\right)  =\sum
\limits_{u\in\widehat{B_{k}P}_{1}}p\left(  u\right)
\label{pf.Bk.ord.Ri.big.c2.goal.3}%
\end{equation}
(because every $u\in\widehat{B_{k}P}_{1}$ satisfies%
\begin{align*}
\left(  \Phi\left(  p,g\right)  \right)  \left(  u\right)   &  =\left\{
\begin{array}
[c]{c}%
p\left(  u\right)  ,\ \ \ \ \ \ \ \ \ \ \text{if }u\in\widehat{B_{k}P}_{1};\\
g\left(  u\right)  ,\ \ \ \ \ \ \ \ \ \ \text{if }u\notin\widehat{B_{k}P}_{1}%
\end{array}
\right.  \ \ \ \ \ \ \ \ \ \ \left(  \text{by the definition of }\Phi\left(
p,g\right)  \right) \\
&  =p\left(  u\right)  \ \ \ \ \ \ \ \ \ \ \left(  \text{since }%
u\in\widehat{B_{k}P}\right)
\end{align*}
).

Substituting (\ref{pf.Bk.ord.Ri.big.c1.goal.2}),
(\ref{pf.Bk.ord.Ri.big.c2.goal.3}) and (\ref{pf.Bk.ord.Ri.big.c1.goal.4}) into
(\ref{pf.Bk.ord.Ri.big.c2.goal.1}), we obtain%
\begin{equation}
\left(  R_{i}\left(  \Phi\left(  p,g\right)  \right)  \right)  \left(
v\right)  =\dfrac{1}{g\left(  v\right)  }\cdot\dfrac{\sum\limits_{u\in
\widehat{B_{k}P}_{1}}p\left(  u\right)  }{\sum\limits_{\substack{u\in
\widehat{P};\\u\gtrdot v}}\dfrac{1}{g\left(  u\right)  }}.
\label{pf.Bk.ord.Ri.big.c2.goal.5}%
\end{equation}
On the other hand, since $v\in\widehat{B_{k}P}_{i}=\widehat{P}_{i-1}$, we have
$\deg_{P}v=i-1$. Hence, Proposition \ref{prop.Ri.implicit} \textbf{(b)}
(applied to $i-1$ and $g$ instead of $i$ and $f$) yields%
\begin{equation}
\left(  R_{i-1}g\right)  \left(  v\right)  =\dfrac{1}{g\left(  v\right)
}\cdot\dfrac{\sum\limits_{\substack{u\in\widehat{P};\\u\lessdot v}}g\left(
u\right)  }{\sum\limits_{\substack{u\in\widehat{P};\\u\gtrdot v}}\dfrac
{1}{g\left(  u\right)  }}. \label{pf.Bk.ord.Ri.big.c2.goal.8}%
\end{equation}
Since $v\in\widehat{P}_{i-1}=\widehat{P}_{1}$ (because $i-1=1$), we know that
$v$ is a minimal element of $P$. Thus, there is only one element
$u\in\widehat{P}$ satisfying $u\lessdot v$, namely the element $0$. Hence,
$\sum\limits_{\substack{u\in\widehat{P};\\u\lessdot v}}g\left(  u\right)
=g\left(  0\right)  $. Therefore, (\ref{pf.Bk.ord.Ri.big.c2.goal.8})
simplifies to%
\begin{equation}
\left(  R_{i-1}g\right)  \left(  v\right)  =\dfrac{1}{g\left(  v\right)
}\cdot\dfrac{g\left(  0\right)  }{\sum\limits_{\substack{u\in\widehat{P}%
;\\u\gtrdot v}}\dfrac{1}{g\left(  u\right)  }}.
\label{pf.Bk.ord.Ri.big.c2.goal.9}%
\end{equation}
But $h=\left(  a_{0},a_{1},...,a_{n+2}\right)  \flat\left(  \Phi\left(
p,R_{i-1}g\right)  \right)  $, so that%
\begin{align*}
h\left(  v\right)   &  =\left(  \left(  a_{0},a_{1},...,a_{n+2}\right)
\flat\left(  \Phi\left(  p,R_{i-1}g\right)  \right)  \right)  \left(  v\right)
\\
&  =\underbrace{a_{\deg_{B_{k}P}v}}_{\substack{=a_{i}\\\text{(since }%
\deg_{B_{k}P}v=i\text{)}}}\cdot\underbrace{\left(  \Phi\left(  p,R_{i-1}%
g\right)  \right)  \left(  v\right)  }_{\substack{=\left\{
\begin{array}
[c]{c}%
p\left(  v\right)  ,\ \ \ \ \ \ \ \ \ \ \text{if }v\in\widehat{B_{k}P}_{1};\\
\left(  R_{i-1}g\right)  \left(  v\right)  ,\ \ \ \ \ \ \ \ \ \ \text{if
}v\notin\widehat{B_{k}P}_{1}%
\end{array}
\right.  \\\text{(by the definition of }\Phi\left(  p,R_{i-1}g\right)
\text{)}}}\\
&  =\underbrace{a_{i}}_{=\dfrac{1}{g\left(  0\right)  }\sum\limits_{u\in
\widehat{B_{k}P}_{1}}p\left(  u\right)  }\cdot\underbrace{\left\{
\begin{array}
[c]{c}%
p\left(  v\right)  ,\ \ \ \ \ \ \ \ \ \ \text{if }v\in\widehat{B_{k}P}_{1};\\
\left(  R_{i-1}g\right)  \left(  v\right)  ,\ \ \ \ \ \ \ \ \ \ \text{if
}v\notin\widehat{B_{k}P}_{1}%
\end{array}
\right.  }_{\substack{=\left(  R_{i-1}g\right)  \left(  v\right)
\\\text{(since }v\notin\widehat{B_{k}P}_{1}\text{)}}}\\
&  =\dfrac{1}{g\left(  0\right)  }\sum\limits_{u\in\widehat{B_{k}P}_{1}%
}p\left(  u\right)  \cdot\left(  R_{i-1}g\right)  \left(  v\right) \\
&  =\dfrac{1}{g\left(  0\right)  }\sum\limits_{u\in\widehat{B_{k}P}_{1}%
}p\left(  u\right)  \cdot\dfrac{1}{g\left(  v\right)  }\cdot\dfrac{g\left(
0\right)  }{\sum\limits_{\substack{u\in\widehat{P};\\u\gtrdot v}}\dfrac
{1}{g\left(  u\right)  }}\ \ \ \ \ \ \ \ \ \ \left(  \text{by
(\ref{pf.Bk.ord.Ri.big.c2.goal.9})}\right) \\
&  =\dfrac{1}{g\left(  v\right)  }\cdot\dfrac{\sum\limits_{u\in\widehat{B_{k}%
P}_{1}}p\left(  u\right)  }{\sum\limits_{\substack{u\in\widehat{P};\\u\gtrdot
v}}\dfrac{1}{g\left(  u\right)  }}=\left(  R_{i}\left(  \Phi\left(
p,g\right)  \right)  \right)  \left(  v\right)  \ \ \ \ \ \ \ \ \ \ \left(
\text{by (\ref{pf.Bk.ord.Ri.big.c2.goal.5})}\right)  .
\end{align*}
We thus have shown that $\left(  R_{i}\left(  \Phi\left(  p,g\right)  \right)
\right)  \left(  v\right)  =h\left(  v\right)  $.

Now, forget that we fixed $v$. We have shown that every $v\in\widehat{B_{k}P}$
satisfies \newline$\left(  R_{i}\left(  \Phi\left(  p,g\right)  \right)
\right)  \left(  v\right)  =h\left(  v\right)  $. In other words,
$R_{i}\left(  \Phi\left(  p,g\right)  \right)  =h$. Thus,%
\[
R_{i}\left(  \Phi\left(  p,g\right)  \right)  =h=\left(  a_{0},a_{1}%
,...,a_{n+2}\right)  \flat\left(  \Phi\left(  p,R_{i-1}g\right)  \right)  .
\]
Hence, $R_{i}\left(  \Phi\left(  p,g\right)  \right)  $ is homogeneously
equivalent to $\Phi\left(  p,R_{i-1}g\right)  $ (because \newline$\left(
a_{0},a_{1},...,a_{n+2}\right)  \flat\left(  \Phi\left(  p,R_{i-1}g\right)
\right)  $ is homogeneously equivalent to $\Phi\left(  p,R_{i-1}g\right)  $
because of Proposition \ref{prop.bemol.hgeq} \textbf{(a)}).

This completes the proof that $R_{i}\left(  \Phi\left(  p,g\right)  \right)  $
is homogeneously equivalent to $\Phi\left(  p,R_{i-1}g\right)  $ in Case 2.

Thus, we have shown that $R_{i}\left(  \Phi\left(  p,g\right)  \right)  $ is
homogeneously equivalent to $\Phi\left(  p,R_{i-1}g\right)  $ in both Cases 1
and 2. Since these two Cases cover all possibilities, this shows that
$R_{i}\left(  \Phi\left(  p,g\right)  \right)  $ is always homogeneously
equivalent to $\Phi\left(  p,R_{i-1}g\right)  $.

In other words, $\pi\left(  R_{i}\left(  \Phi\left(  p,g\right)  \right)
\right)  =\pi\left(  \Phi\left(  p,R_{i-1}g\right)  \right)  $. Since
\begin{align*}
\pi\left(  R_{i}\left(  \Phi\left(  p,g\right)  \right)  \right)   &
=\underbrace{\left(  \pi\circ R_{i}\right)  }_{\substack{=\overline{R_{i}%
}\circ\pi\\\text{(since the diagram}\\\text{(\ref{def.hgRi.commut}) is
commutative)}}}\left(  \Phi\left(  p,g\right)  \right)  =\left(
\overline{R_{i}}\circ\pi\right)  \left(  \Phi\left(  p,g\right)  \right) \\
&  =\overline{R_{i}}\left(  \underbrace{\pi\left(  \Phi\left(  p,g\right)
\right)  }_{\substack{=\pi\left(  p\right)  \rightthreetimes\pi\left(
g\right)  \\\text{(by (\ref{pf.Bk.ord.pitriv}))}}}\right)
=\underbrace{\overline{R_{i}}}_{=\left(\overline{R_{i}}\right)_{B_{k}P}}
\left(
\underbrace{\pi\left(  p\right)  }_{=\widetilde{p}}\rightthreetimes
\underbrace{\pi\left(  g\right)  }_{=\widetilde{g}}\right) \\
&  =\left(\overline{R_{i}}\right)_{B_{k}P}\left(
\widetilde{p}\rightthreetimes \widetilde{g}\right)
\end{align*}
and%
\begin{align*}
&  \pi\left(  \Phi\left(  p,R_{i-1}g\right)  \right) \\
&  =\underbrace{\pi\left(  p\right)  }_{=\widetilde{p}}\rightthreetimes
\underbrace{\pi\left(  R_{i-1}g\right)  }_{=\left(  \pi\circ R_{i-1}\right)
\left(  g\right)  }\ \ \ \ \ \ \ \ \ \ \left(  \text{by
(\ref{pf.Bk.ord.pitriv}), applied to }R_{i-1}g\text{ instead of }g\right) \\
&  =\widetilde{p}\rightthreetimes\underbrace{\left(  \pi\circ R_{i-1}\right)
}_{\substack{=\overline{R_{i-1}}\circ\pi\\\text{(since the diagram}%
\\\text{(\ref{def.hgRi.commut}) is commutative,}\\\text{for }i-1\text{ instead
of }i\text{)}}}\left(  g\right)  =\widetilde{p}\rightthreetimes
\underbrace{\left(  \overline{R_{i-1}}\circ\pi\right)  \left(  g\right)
}_{=\overline{R_{i-1}}\left(  \pi\left(  g\right)  \right)  }\\
&  =\widetilde{p}\rightthreetimes\underbrace{\overline{R_{i-1}}}%
_{=\left(\overline{R_{i-1}}\right)_{P}}
\left(  \underbrace{\pi\left(  g\right)
}_{=\widetilde{g}}\right)
= \widetilde{p} \rightthreetimes \left(\overline{R_{i-1}}\right)_{P}
\left(  \widetilde{g}\right)  ,
\end{align*}
this rewrites as $\left(\overline{R_{i}}\right)_{B_{k}P}\left(
\widetilde{p} \rightthreetimes\widetilde{g}\right)
=\widetilde{p}\rightthreetimes
\left(\overline{R_{i-1}}\right)_{P}\left(  \widetilde{g}\right)  $. Thus,
(\ref{pf.Bk.ord.Ri.big}) is proven.

\textit{Proof of (\ref{pf.Bk.ord.Ri.small}):} Let $\widetilde{p}\in
\mathbb{P}\left(  \mathbb{K}^{k}\right)  $ and $\widetilde{g}\in
\overline{\mathbb{K}^{\widehat{P}}}$.

Since $\widetilde{g}\in\overline{\mathbb{K}^{\widehat{P}}}$, there exists some
$g\in\mathbb{K}^{\widehat{P}}$ such that $\widetilde{g}=\pi\left(  g\right)
$. Consider this $g$.

Since $\widetilde{p}\in\mathbb{P}\left(  \mathbb{K}^{k}\right)  $, there
exists some nonzero $p\in\mathbb{K}^{k}$ such that $\widetilde{p}=\pi\left(
p\right)  $. Consider this $p$.

Let $p^{-1}$ denote the result of replacing every coordinate of $p\in
\mathbb{K}^{k}$ by its inverse. Then, $\left(  \pi\left(  p\right)  \right)
^{-1}=\pi\left(  p^{-1}\right)  $ (because $\left(  \pi\left(  p\right)
\right)  ^{-1}$ is the result of replacing every homogeneous coordinate of
$\pi\left(  p\right)  $ by its inverse, and because the homogeneous
coordinates of $\pi\left(  p\right)  $ are simply the coordinates of $p$).
Since $\pi\left(  p\right)  =\widetilde{p}$, this rewrites as follows:%
\begin{equation}
\widetilde{p}^{-1}=\pi\left(  p^{-1}\right)  . \label{pf.Bk.ord.Ri.small.1}%
\end{equation}

Notice that $p^{-1}\left(  v\right)  =\left(  p\left(  v\right)  \right)
^{-1}$ for every $v\in\widehat{B_{k}P}_{1}$ (by the definition of $p^{-1}$).

We will now show that $R_{1}\left(  \Phi\left(  p,g\right)  \right)  $ is
homogeneously equivalent to $\Phi\left(  p^{-1},g\right)  $.

Let us first notice that $2\in\left\{  2,3,...,n+2\right\}  $, so that
$\widehat{B_{k}P}_{2}=\widehat{P}_{1}$.

We WLOG assume that the elements $g\left(  0\right)  $ and $\sum
\limits_{u\in\widehat{P}_{1}}\dfrac{1}{g\left(  u\right)  }$ of $\mathbb{K}$
are nonzero. (This assumption is really WLOG, because we are proving a
polynomial identity, and thus restriction to a Zariski-dense and open subset
can do no harm.)

Define an $\left(  n+3\right)  $-tuple $\left(  a_{0},a_{1},...,a_{n+2}%
\right)  \in\left(  \mathbb{K}^{\times}\right)  ^{n+3}$ by%
\[
a_{j}=\left\{
\begin{array}
[c]{c}%
1,\ \ \ \ \ \ \ \ \ \ \text{if }j\neq1;\\
\dfrac{g\left(  0\right)  }{\sum\limits_{u\in\widehat{P}_{1}}\dfrac
{1}{g\left(  u\right)  }},\ \ \ \ \ \ \ \ \ \ \text{if }j=1
\end{array}
\right.  \ \ \ \ \ \ \ \ \ \ \text{for every }j\in\left\{
0,1,...,n+2\right\}  .
\]
Define an $h\in\mathbb{K}^{\widehat{B_{k}P}}$ by $h=\left(  a_{0}%
,a_{1},...,a_{n+2}\right)  \flat\left(  \Phi\left(  p^{-1},g\right)  \right)
$ (where \newline$\left(  a_{0},a_{1},...,a_{n+2}\right)  \flat\left(
\Phi\left(  p^{-1},g\right)  \right)  $ is to be understood according to
Definition \ref{def.bemol}).

We are going to prove that
\begin{equation}
R_{1}\left(  \Phi\left(  p,g\right)  \right)  =h.
\label{pf.Bk.ord.Ri.small.goal}%
\end{equation}

Indeed, let $v\in\widehat{B_{k}P}$. Then, we are going to prove that $\left(
R_{1}\left(  \Phi\left(  p,g\right)  \right)  \right)  \left(  v\right)
=h\left(  v\right)  $.

If $v\notin\widehat{B_{k}P}_{1}$, then $\left(  R_{1}\left(  \Phi\left(
p,g\right)  \right)  \right)  \left(  v\right)  =h\left(  v\right)  $ is very
easy to check\footnote{\textit{Proof.} Assume that $v\notin\widehat{B_{k}%
P}_{1}$. Then, $\deg_{B_{k}P}v\neq1$. Thus, Proposition \ref{prop.Ri.implicit}
\textbf{(a)} (applied to $B_{k}P$, $n+1$, $1$ and $\Phi\left(  p,g\right)  $
instead of $P$, $i$, $n$ and $f$) yields%
\begin{align}
\left(  R_{1}\left(  \Phi\left(  p,g\right)  \right)  \right)  \left(
v\right)   &  =\left(  \Phi\left(  p,g\right)  \right)  \left(  v\right)
=\left\{
\begin{array}
[c]{c}%
p\left(  v\right)  ,\ \ \ \ \ \ \ \ \ \ \text{if }v\in\widehat{B_{k}P}_{1};\\
g\left(  v\right)  ,\ \ \ \ \ \ \ \ \ \ \text{if }v\notin\widehat{B_{k}P}_{1}%
\end{array}
\right.  \ \ \ \ \ \ \ \ \ \ \left(  \text{by the definition of }\Phi\left(
p,g\right)  \right) \nonumber\\
&  =g\left(  v\right)  \ \ \ \ \ \ \ \ \ \ \left(  \text{since }%
v\notin\widehat{B_{k}P}_{1}\right)  . \label{pf.Bk.ord.Ri.small.goal.0}%
\end{align}
\par
But $h=\left(  a_{0},a_{1},...,a_{n+2}\right)  \flat\left(  \Phi\left(
p^{-1},g\right)  \right)  $, so that
\begin{align*}
h\left(  v\right)   &  =\left(  \left(  a_{0},a_{1},...,a_{n+2}\right)
\flat\left(  \Phi\left(  p^{-1},g\right)  \right)  \right)  \left(  v\right)
\\
&  =\underbrace{a_{\deg_{B_{k}P}v}}_{\substack{=1\\\text{(since }\deg_{B_{k}%
P}v\neq1\text{)}}}\cdot\left(  \Phi\left(  p^{-1},g\right)  \right)  \left(
v\right)  \ \ \ \ \ \ \ \ \ \ \left(  \text{by the definition of }\left(
a_{0},a_{1},...,a_{n+2}\right)  \flat\left(  \Phi\left(  p^{-1},g\right)
\right)  \right) \\
&  =\left(  \Phi\left(  p^{-1},g\right)  \right)  \left(  v\right)  =\left\{
\begin{array}
[c]{c}%
p^{-1}\left(  v\right)  ,\ \ \ \ \ \ \ \ \ \ \text{if }v\in\widehat{B_{k}%
P}_{1};\\
g\left(  v\right)  ,\ \ \ \ \ \ \ \ \ \ \text{if }v\notin\widehat{B_{k}P}_{1}%
\end{array}
\right. \\
&  \ \ \ \ \ \ \ \ \ \ \left(  \text{by the definition of }\Phi\left(
p^{-1},g\right)  \right) \\
&  =g\left(  v\right)  \ \ \ \ \ \ \ \ \ \ \left(  \text{since }%
v\notin\widehat{B_{k}P}_{1}\right)  .
\end{align*}
Compared with (\ref{pf.Bk.ord.Ri.small.goal.0}), this yields $\left(
R_{1}\left(  \Phi\left(  p,g\right)  \right)  \right)  \left(  v\right)
=h\left(  v\right)  $, qed.}. Hence, for the rest of the proof of $\left(
R_{1}\left(  \Phi\left(  p,g\right)  \right)  \right)  \left(  v\right)
=h\left(  v\right)  $, we can WLOG assume that we don't have $v\notin%
\widehat{B_{k}P}_{1}$. Assume this. Thus, $v\in\widehat{B_{k}P}_{1}$. Hence,
$\deg_{B_{k}P}v=1$. Thus, Proposition \ref{prop.Ri.implicit} \textbf{(b)}
(applied to $B_{k}P$, $n+1$, $1$ and $\Phi\left(  p,g\right)  $ instead of
$P$, $n$, $i$ and $f$) yields%
\begin{equation}
\left(  R_{1}\left(  \Phi\left(  p,g\right)  \right)  \right)  \left(
v\right)  =\dfrac{1}{\left(  \Phi\left(  p,g\right)  \right)  \left(
v\right)  }\cdot\dfrac{\sum\limits_{\substack{u\in\widehat{B_{k}P};\\u\lessdot
v}}\left(  \Phi\left(  p,g\right)  \right)  \left(  u\right)  }{\sum
\limits_{\substack{u\in\widehat{B_{k}P};\\u\gtrdot v}}\dfrac{1}{\left(
\Phi\left(  p,g\right)  \right)  \left(  u\right)  }}.
\label{pf.Bk.ord.Ri.small.goal.1}%
\end{equation}
We are now going to simplify the right hand side of this equality.

First of all, the definition of $\Phi\left(  p,g\right)  $ yields
\begin{equation}
\left(  \Phi\left(  p,g\right)  \right)  \left(  v\right)  =\left\{
\begin{array}
[c]{c}%
p\left(  v\right)  ,\ \ \ \ \ \ \ \ \ \ \text{if }v\in\widehat{B_{k}P}_{1};\\
g\left(  v\right)  ,\ \ \ \ \ \ \ \ \ \ \text{if }v\notin\widehat{B_{k}P}_{1}%
\end{array}
\right.  =p\left(  v\right)  \label{pf.Bk.ord.Ri.small.goal.2}%
\end{equation}
(since $v\in\widehat{B_{k}P}_{1}$).

Since $\deg_{B_{k}P}v=1$, the element $v$ is a minimal element of $B_{k}P$.
Thus, there exists only one element $u\in\widehat{B_{k}P}$ satisfying
$u\lessdot v$, and this element is $0$. Hence,%
\begin{equation}
\sum\limits_{\substack{u\in\widehat{B_{k}P};\\u\lessdot v}}\left(  \Phi\left(
p,g\right)  \right)  \left(  u\right)  =\left(  \Phi\left(  p,g\right)
\right)  \left(  0\right)  =\left\{
\begin{array}
[c]{c}%
p\left(  0\right)  ,\ \ \ \ \ \ \ \ \ \ \text{if }0\in\widehat{B_{k}P}_{1};\\
g\left(  0\right)  ,\ \ \ \ \ \ \ \ \ \ \text{if }0\notin\widehat{B_{k}P}_{1}%
\end{array}
\right.  =g\left(  0\right)  \label{pf.Bk.ord.Ri.small.goal.3}%
\end{equation}
(since $0\notin\widehat{B_{k}P}_{1}$).

Recall that the poset $B_{k}P$ was obtained from the poset $P$ by adding $k$
new elements, each of which is set to be smaller than all elements of $P$. One
of these $k$ new elements is $v$ (since $v\in\widehat{B_{k}P}_{1}$). Thus, the
elements $u\in\widehat{B_{k}P}$ satisfying $u\gtrdot v$ are precisely the
minimal elements of $P$. Since the minimal elements of $P$ are precisely the
elements of $\widehat{P}_{1}$ (because $P$ is $n$-graded), this yields that
the elements $u\in\widehat{B_{k}P}$ satisfying $u\gtrdot v$ are precisely the
elements of $\widehat{P}_{1}$. Thus,%
\[
\sum\limits_{\substack{u\in\widehat{B_{k}P};\\u\gtrdot v}}\dfrac{1}{\left(
\Phi\left(  p,g\right)  \right)  \left(  u\right)  }=\sum\limits_{u\in
\widehat{P}_{1}}\dfrac{1}{\left(  \Phi\left(  p,g\right)  \right)  \left(
u\right)  }.
\]
Since every $u\in\widehat{P}_{1}$ satisfies%
\begin{align*}
\left(  \Phi\left(  p,g\right)  \right)  \left(  u\right)   &  =\left\{
\begin{array}
[c]{c}%
p\left(  u\right)  ,\ \ \ \ \ \ \ \ \ \ \text{if }u\in\widehat{B_{k}P}_{1};\\
g\left(  u\right)  ,\ \ \ \ \ \ \ \ \ \ \text{if }u\notin\widehat{B_{k}P}_{1}%
\end{array}
\right.  \ \ \ \ \ \ \ \ \ \ \left(  \text{by the definition of }\Phi\left(
p,g\right)  \right) \\
&  =g\left(  u\right)  \ \ \ \ \ \ \ \ \ \ \left(
\begin{array}
[c]{c}%
\text{since }u\notin\widehat{B_{k}P}_{1}\text{ (because }u\in\widehat{P}%
_{1}=\widehat{B_{k}P}_{2}\text{, hence}\\
\deg_{B_{k}P}u=2\neq1\text{, hence }u\notin\widehat{B_{k}P}_{1}\text{)}%
\end{array}
\right)  ,
\end{align*}
this rewrites as%
\begin{equation}
\sum\limits_{\substack{u\in\widehat{B_{k}P};\\u\gtrdot v}}\dfrac{1}{\left(
\Phi\left(  p,g\right)  \right)  \left(  u\right)  }=\sum\limits_{u\in
\widehat{P}_{1}}\dfrac{1}{g\left(  u\right)  }.
\label{pf.Bk.ord.Ri.small.goal.4}%
\end{equation}

Substituting (\ref{pf.Bk.ord.Ri.small.goal.2}),
(\ref{pf.Bk.ord.Ri.small.goal.3}) and (\ref{pf.Bk.ord.Ri.small.goal.4}) into
(\ref{pf.Bk.ord.Ri.small.goal.1}), we obtain%
\[
\left(  R_{1}\left(  \Phi\left(  p,g\right)  \right)  \right)  \left(
v\right)  =\dfrac{1}{p\left(  v\right)  }\cdot\dfrac{g\left(  0\right)  }%
{\sum\limits_{u\in\widehat{P}_{1}}\dfrac{1}{g\left(  u\right)  }}.
\]
Compared with%
\begin{align*}
h\left(  v\right)   &  =\left(  \left(  a_{0},a_{1},...,a_{n+2}\right)
\flat\left(  \Phi\left(  p^{-1},g\right)  \right)  \right)  \left(  v\right)
\ \ \ \ \ \ \ \ \ \ \left(  \text{since }h=\left(  a_{0},a_{1},...,a_{n+2}%
\right)  \flat\left(  \Phi\left(  p^{-1},g\right)  \right)  \right) \\
&  =\underbrace{a_{\deg_{B_{k}P}v}}_{\substack{=a_{1}\\\text{(since }%
\deg_{B_{k}P}v=1\text{)}}}\cdot\underbrace{\left(  \Phi\left(  p^{-1}%
,g\right)  \right)  \left(  v\right)  }_{\substack{=\left\{
\begin{array}
[c]{c}%
p^{-1}\left(  v\right)  ,\ \ \ \ \ \ \ \ \ \ \text{if }v\in\widehat{B_{k}%
P}_{1};\\
g\left(  v\right)  ,\ \ \ \ \ \ \ \ \ \ \text{if }v\notin\widehat{B_{k}P}_{1}%
\end{array}
\right.  \\\text{(by the definition of }\Phi\left(  p^{-1},g\right)  \text{)}%
}}\\
&  =\underbrace{a_{1}}_{=\dfrac{g\left(  0\right)  }{\sum\limits_{u\in
\widehat{P}_{1}}\dfrac{1}{g\left(  u\right)  }}}\cdot\underbrace{\left\{
\begin{array}
[c]{c}%
p^{-1}\left(  v\right)  ,\ \ \ \ \ \ \ \ \ \ \text{if }v\in\widehat{B_{k}%
P}_{1};\\
g\left(  v\right)  ,\ \ \ \ \ \ \ \ \ \ \text{if }v\notin\widehat{B_{k}P}_{1}%
\end{array}
\right.  }_{\substack{=p^{-1}\left(  v\right)  \\\text{(since }v\in
\widehat{B_{k}P}_{1}\text{)}}}\\
&  =\dfrac{g\left(  0\right)  }{\sum\limits_{u\in\widehat{P}_{1}}\dfrac
{1}{g\left(  u\right)  }}\cdot\underbrace{p^{-1}\left(  v\right)  }_{=\left(
p\left(  v\right)  \right)  ^{-1}}=\dfrac{g\left(  0\right)  }{\sum
\limits_{u\in\widehat{P}_{1}}\dfrac{1}{g\left(  u\right)  }}\cdot\left(
p\left(  v\right)  \right)  ^{-1}=\dfrac{1}{p\left(  v\right)  }\cdot
\dfrac{g\left(  0\right)  }{\sum\limits_{u\in\widehat{P}_{1}}\dfrac
{1}{g\left(  u\right)  }},
\end{align*}
this yields $\left(  R_{1}\left(  \Phi\left(  p,g\right)  \right)  \right)
\left(  v\right)  =h\left(  v\right)  $.

Now, forget that we fixed $v$. We have shown that every $v\in\widehat{B_{k}P}$
satisfies \newline$\left(  R_{1}\left(  \Phi\left(  p,g\right)  \right)
\right)  \left(  v\right)  =h\left(  v\right)  $. In other words,
$R_{1}\left(  \Phi\left(  p,g\right)  \right)  =h$. Thus,%
\[
R_{1}\left(  \Phi\left(  p,g\right)  \right)  =h=\left(  a_{0},a_{1}%
,...,a_{n+2}\right)  \flat\left(  \Phi\left(  p^{-1},g\right)  \right)  .
\]
Hence, $R_{1}\left(  \Phi\left(  p,g\right)  \right)  $ is homogeneously
equivalent to $\Phi\left(  p^{-1},g\right)  $ (because \newline$\left(
a_{0},a_{1},...,a_{n+2}\right)  \flat\left(  \Phi\left(  p^{-1},g\right)
\right)  $ is homogeneously equivalent to $\Phi\left(  p^{-1},g\right)  $
(according to Proposition \ref{prop.bemol.hgeq} \textbf{(a)})). In other
words, $\pi\left(  R_{1}\left(  \Phi\left(  p,g\right)  \right)  \right)
=\pi\left(  \Phi\left(  p^{-1},g\right)  \right)  $. Since
\begin{align*}
\pi\left(  R_{1}\left(  \Phi\left(  p,g\right)  \right)  \right)   &
=\underbrace{\left(  \pi\circ R_{1}\right)  }_{\substack{=\overline{R_{1}%
}\circ\pi\\\text{(since the diagram}\\\text{(\ref{def.hgRi.commut}) is
commutative,}\\\text{for }i=1\text{)}}}\left(  \Phi\left(  p,g\right)
\right)  =\left(  \overline{R_{1}}\circ\pi\right)  \left(  \Phi\left(
p,g\right)  \right) \\
&  =\overline{R_{1}}\left(  \underbrace{\pi\left(  \Phi\left(  p,g\right)
\right)  }_{\substack{=\pi\left(  p\right)  \rightthreetimes\pi\left(
g\right)  \\\text{(by (\ref{pf.Bk.ord.pitriv}))}}}\right)
=\underbrace{\overline{R_{1}}}_{=\left(\overline{R_{1}}\right)_{B_{k}P}}
\left(
\underbrace{\pi\left(  p\right)  }_{=\widetilde{p}}\rightthreetimes
\underbrace{\pi\left(  g\right)  }_{=\widetilde{g}}\right) \\
&  =\left(\overline{R_{1}}\right)_{B_{k}P}
\left(  \widetilde{p}\rightthreetimes \widetilde{g}\right)
\end{align*}
and%
\begin{align*}
&  \pi\left(  \Phi\left(  p^{-1},g\right)  \right) \\
&  =\underbrace{\pi\left(  p^{-1}\right)  }_{\substack{=\widetilde{p}%
^{-1}\\\text{(by (\ref{pf.Bk.ord.Ri.small.1}))}}}\rightthreetimes
\underbrace{\pi\left(  g\right)  }_{=\widetilde{g}}\ \ \ \ \ \ \ \ \ \ \left(
\text{by (\ref{pf.Bk.ord.pitriv}), applied to }p^{-1}\text{ instead of
}p\right) \\
&  =\widetilde{p}^{-1}\rightthreetimes\widetilde{g},
\end{align*}
this rewrites as $\left(\overline{R_{1}}\right)_{B_{k}P}\left(
\widetilde{p}\rightthreetimes\widetilde{g}\right)
=\widetilde{p}^{-1}\rightthreetimes
\widetilde{g}$. Thus, (\ref{pf.Bk.ord.Ri.small}) is proven.
\end{verlong}

But recall that $\overline{R}=\overline{R_{1}}\circ\overline{R_{2}}%
\circ...\circ\overline{R_{n}}$ for any $n$-graded poset. Hence, $\overline
{R}_{B_{k}P}
=\left(\overline{R_{1}}\right)_{B_{k}P}
\circ\left(\overline{R_{2}}\right)_{B_{k}P}
\circ\left(\overline{R_{3}}\right)_{B_{k}P}\circ...
\circ\left(\overline{R_{n+1}}\right)_{B_{k}P}$
(because $B_{k}P$ is an $\left(  n+1\right)  $-graded poset) and $\overline
{R}_{P}
=
\left(\overline{R_{1}}\right)_{P}\circ
\left(\overline{R_{2}}\right)_{P}\circ...\circ
\left(\overline{R_{n}}\right)_{P}$ (because $P$ is an $n$-graded poset).
Because of these
equalities, and because of (\ref{pf.Bk.ord.Ri.big}) and
(\ref{pf.Bk.ord.Ri.small}), it is now easy to see that every $\widetilde{p}%
\in\mathbb{P}\left(  \mathbb{K}^{k}\right)  $ and $\widetilde{g}\in
\overline{\mathbb{K}^{\widehat{P}}}$ satisfy
\begin{equation}
\overline{R}_{B_{k}P}\left(  \widetilde{p}\rightthreetimes\widetilde{g}%
\right)  =\widetilde{p}^{-1}\rightthreetimes\overline{R}_{P}\left(
\widetilde{g}\right)  . \label{pf.Bk.ord.R}%
\end{equation}

\begin{verlong}
\textit{Proof of (\ref{pf.Bk.ord.R}):} Let $\widetilde{p}\in\mathbb{P}\left(
\mathbb{K}^{k}\right)  $ and $\widetilde{g}\in\overline{\mathbb{K}%
^{\widehat{P}}}$ be arbitrary. For every $m\in\left\{  1,2,...,n+2\right\}  $,
let $\mathfrak{A}_{m}$ denote the rational map
\[
\left(\overline{R_{m}}\right)_{B_{k}P}
\circ\left(\overline{R_{m+1}}\right)_{B_{k}P}\circ...
\circ\left(\overline{R_{n+1}}\right)_{B_{k}P}
:\overline{\mathbb{K}^{\widehat{B_{k}P}}}
\rightarrow\overline{\mathbb{K}^{\widehat{B_{k}P}}}.
\]
For every $m\in\left\{  1,2,...,n+1\right\}  $, let $\mathfrak{B}_{m}$ denote
the rational map%
\[
\left(\overline{R_{m}}\right)_{P}\circ
\left(\overline{R_{m+1}}\right)_{P}\circ...\circ
\left(\overline{R_{n}}\right)_{P}
:\overline{\mathbb{K}^{\widehat{P}}}\rightarrow\overline{\mathbb{K}%
^{\widehat{P}}}.
\]

It is easy to see that every $\ell\in\left\{  0,1,...,n\right\}  $ satisfies%
\begin{equation}
\mathfrak{A}_{n+2-\ell}\left(  \widetilde{p}\rightthreetimes\widetilde{g}%
\right)  =\widetilde{p}\rightthreetimes\mathfrak{B}_{n+1-\ell}\left(
\widetilde{g}\right)  . \label{pf.Bk.ord.R.pf.ind}%
\end{equation}
\footnote{\textit{Proof of (\ref{pf.Bk.ord.R.pf.ind}):} We will prove
(\ref{pf.Bk.ord.R.pf.ind}) by induction over $\ell$:
\par
\textit{Induction base:} We have
\begin{align*}
\mathfrak{A}_{n+2-0}  &  =\mathfrak{A}_{n+2}
=\left(\overline{R_{n+2}}\right)_{B_{k}P}
\circ\left(\overline{R_{n+2+1}}\right)_{B_{k}P}\circ...
\circ\left(\overline{R_{n+1}}\right)_{B_{k}P}
\ \ \ \ \ \ \ \ \ \ \left(  \text{by the definition of }\mathfrak{A}%
_{n+2}\right) \\
&  =\left(  \text{empty composition}\right)  =\operatorname*{id},
\end{align*}
so that $\mathfrak{A}_{n+2-0}\left(  \widetilde{p}\rightthreetimes
\widetilde{g}\right)  =\operatorname*{id}\left(  \widetilde{p}\rightthreetimes
\widetilde{g}\right)  =\widetilde{p}\rightthreetimes\widetilde{g}$. But%
\begin{align*}
\mathfrak{B}_{n+1-0}  &  =\mathfrak{B}_{n+1}
=
\left(\overline{R_{n+1}}\right)_{P}\circ
\left(\overline{R_{n+1+1}}\right)_{P}\circ...\circ
\left(\overline{R_{n}}\right)_{P}
\ \ \ \ \ \ \ \ \ \ \left(  \text{by the definition of }\mathfrak{B}%
_{n+1}\right) \\
&  =\left(  \text{empty composition}\right)  =\operatorname*{id},
\end{align*}
so that $\widetilde{p}\rightthreetimes\mathfrak{B}_{n+1-0}\left(
\widetilde{g}\right)  =\widetilde{p}\rightthreetimes\operatorname*{id}\left(
\widetilde{g}\right)  =\widetilde{p}\rightthreetimes\widetilde{g}$. Now,
$\mathfrak{A}_{n+2-0}\left(  \widetilde{p}\rightthreetimes\widetilde{g}%
\right)  =\widetilde{p}\rightthreetimes\widetilde{g}=\widetilde{p}%
\rightthreetimes\mathfrak{B}_{n+1-0}\left(  \widetilde{g}\right)  $. In other
words, (\ref{pf.Bk.ord.R.pf.ind}) holds for $\ell=0$. This completes the
induction base.
\par
\textit{Induction step:} Let $L\in\left\{  0,1,...,n\right\}  $ be such that
$L<n$. Assume that (\ref{pf.Bk.ord.R.pf.ind}) holds for $\ell=L$. We need to
prove that (\ref{pf.Bk.ord.R.pf.ind}) holds for $\ell=L+1$.
\par
Notice that $L\in\left\{  0,1,...,n-1\right\}  $ (since $L\in\left\{
0,1,...,n\right\}  $ and $L<n$). Hence, $n+1-L\in\left\{  2,3,...,n+1\right\}
$.
\par
Since (\ref{pf.Bk.ord.R.pf.ind}) holds for $\ell=L$, we have%
\begin{equation}
\mathfrak{A}_{n+2-L}\left(  \widetilde{p}\rightthreetimes\widetilde{g}\right)
=\widetilde{p}\rightthreetimes\mathfrak{B}_{n+1-L}\left(  \widetilde{g}%
\right)  . \label{pf.Bk.ord.R.pf.ind.pf.1}%
\end{equation}
\par
Now, the definition of $\mathfrak{A}_{n+2-L}$ yields $\mathfrak{A}%
_{n+2-L}
=\left(\overline{R_{n+2-L}}\right)_{B_{k}P}
\circ\left(\overline{R_{n+2-L+1}}\right)_{B_{k}P}\circ...
\circ\left(\overline{R_{n+1}}\right)_{B_{k}P}$. But we have
\begin{align*}
\mathfrak{A}_{n+2-\left(  L+1\right)  }  &  =\mathfrak{A}_{n+2-L-1}%
=\left(\overline{R_{n+2-L-1}}\right)_{B_{k}P}
\circ\left(\overline{R_{n+2-L-1+1}}\right)_{B_{k}P}
\circ...\circ\left(\overline{R_{n+1}}\right)_{B_{k}P}\\
&  \ \ \ \ \ \ \ \ \ \ \left(  \text{by the definition of }\mathfrak{A}%
_{n+2-L-1}\right) \\
&  =\left(\overline{R_{n+2-L-1}}\right)_{B_{k}P}
\circ\left(\overline{R_{n+2-L}}\right)_{B_{k}P}
\circ...\circ\left(\overline{R_{n+1}}\right)_{B_{k}P}\\
&  =\underbrace{\left(\overline{R_{n+2-L-1}}\right)_{B_{k}P}}_{
=\left(\overline{R_{n+1-L}}\right)_{B_{k}P}}\circ
\underbrace{\left(\overline{R_{n+2-L}}\right)_{B_{k}P}\circ
\left(\overline{R_{n+2-L+1}}\right)_{B_{k}P}\circ...
\circ\left(\overline{R_{n+1}}\right)_{B_{k}P}}_{=\mathfrak{A}%
_{n+2-L}}
=\left(\overline{R_{n+1-L}}\right)_{B_{k}P}\circ
\mathfrak{A}_{n+2-L},
\end{align*}
so that
\begin{align}
\mathfrak{A}_{n+2-\left(  L+1\right)  }\left(  \widetilde{p}\rightthreetimes
\widetilde{g}\right)
&  =\left( \left(\overline{R_{n+1-L}}\right)_{B_{k}P}%
\circ\mathfrak{A}_{n+2-L}\right)  \left(  \widetilde{p}\rightthreetimes
\widetilde{g}\right)
= \left(\overline{R_{n+1-L}}\right)_{B_{k}P}\left(
\underbrace{\mathfrak{A}_{n+2-L}\left(  \widetilde{p}\rightthreetimes
\widetilde{g}\right)  }_{\substack{=\widetilde{p}\rightthreetimes
\mathfrak{B}_{n+1-L}\left(  \widetilde{g}\right)  \\\text{(by
(\ref{pf.Bk.ord.R.pf.ind.pf.1}))}}}\right) \nonumber\\
&  =\left(\overline{R_{n+1-L}}\right)_{B_{k}P}\left(
\widetilde{p}\rightthreetimes
\mathfrak{B}_{n+1-L}\left(  \widetilde{g}\right)  \right)  =\widetilde{p}%
\rightthreetimes\underbrace{\left(\overline{R_{n+1-L-1}}\right)_{P}
}_{=\left(\overline{R_{n-L}}\right)_{P}}
\left(  \mathfrak{B}_{n+1-L}\left(  \widetilde{g}\right)  \right)
.\nonumber\\
&  \ \ \ \ \ \ \ \ \ \ \left(
\begin{array}
[c]{c}%
\text{by (\ref{pf.Bk.ord.Ri.big}) (applied to }n+1-L\text{ and }%
\mathfrak{B}_{n+1-L}\left(  \widetilde{g}\right)  \text{ instead of }i\text{
and }\widetilde{g}\text{)}\\
\text{(since }n+1-L\in\left\{  2,3,...,n+1\right\}  \text{)}%
\end{array}
\right) \nonumber\\
&  =\widetilde{p}\rightthreetimes
\left(\overline{R_{n-L}}\right)_{P}\left(  \mathfrak{B}%
_{n+1-L}\left(  \widetilde{g}\right)  \right)  .
\label{pf.Bk.ord.R.pf.ind.pf.3}%
\end{align}
\par
On the other hand, the definition of $\mathfrak{B}_{n+1-L}$ yields
$\mathfrak{B}_{n+1-L}
=
\left(\overline{R_{n+1-L}}\right)_{P}\circ
\left(\overline{R_{n+1-L+1}}\right)_{P}\circ...\circ
\left(\overline{R_{n}}\right)_{P}$. But
\begin{align*}
\mathfrak{B}_{n+1-\left(  L+1\right)  }  &  =\mathfrak{B}_{n+1-L-1}%
=
\left(\overline{R_{n+1-L-1}}\right)_{P}\circ
\left(\overline{R_{n+1-L-1+1}}\right)_{P}\circ ...\circ
\left(\overline{R_{n}}\right)_{P}\\
&  \ \ \ \ \ \ \ \ \ \ \left(  \text{by the definition of }\mathfrak{B}%
_{n+1-L-1}\right) \\
&  =
\left(\overline{R_{n+1-L-1}}\right)_{P}\circ
\left(\overline{R_{n+1-L}}\right)_{P}\circ ...\circ
\left(\overline{R_{n}}\right)_{P}\\
&  =\underbrace{\left(\overline{R_{n+1-L-1}}\right)_{P}}_{
 =\left(\overline{R_{n-L}}\right)_{P}} \circ
\underbrace{\left(\overline{R_{n+1-L}}\right)_{P}\circ
\left(\overline{R_{n+1-L+1}}\right)_{P} \circ...\circ
\left(\overline{R_{n}}\right)_{P}}_{=\mathfrak{B}_{n+1-L}}
= \left(\overline{R_{n-L}}\right)_{P}\circ\mathfrak{B}_{n+1-L}.
\end{align*}
Hence,%
\begin{equation}
\mathfrak{B}_{n+1-\left(  L+1\right)  }\left(  \widetilde{g}\right)  =\left(
\left(\overline{R_{n-L}}\right)_{P}\circ
\mathfrak{B}_{n+1-L}\right)  \left(  \widetilde{g}%
\right)
= \left(\overline{R_{n-L}}\right)_{P}\left(  \mathfrak{B}_{n+1-L}\left(
\widetilde{\mathfrak{g}}\right)  \right)  . \label{pf.Bk.ord.R.pf.ind.pf.4}%
\end{equation}
Now, (\ref{pf.Bk.ord.R.pf.ind.pf.3}) becomes%
\[
\mathfrak{A}_{n+2-\left(  L+1\right)  }\left(  \widetilde{p}\rightthreetimes
\widetilde{g}\right)  =\widetilde{p}\rightthreetimes
\underbrace{\left(\overline{R_{n-L}}\right)_{P}\left( 
 \mathfrak{B}_{n+1-L}\left(  \widetilde{g}\right)
\right)  }_{\substack{=\mathfrak{B}_{n+1-\left(  L+1\right)  }\left(
\widetilde{g}\right)  \\\text{(by (\ref{pf.Bk.ord.R.pf.ind.pf.4}))}%
}}=\widetilde{p}\rightthreetimes\mathfrak{B}_{n+1-\left(  L+1\right)  }\left(
\widetilde{g}\right)  .
\]
In other words, (\ref{pf.Bk.ord.R.pf.ind}) holds for $\ell=L+1$. This
completes the induction step. The induction proof of (\ref{pf.Bk.ord.R.pf.ind}%
) is thus finished.}

Now, we can apply (\ref{pf.Bk.ord.R.pf.ind}) to $\ell=n$. As a consequence, we
obtain $\mathfrak{A}_{n+2-n}\left(  \widetilde{p}\rightthreetimes
\widetilde{g}\right)  =\widetilde{p}\rightthreetimes\mathfrak{B}%
_{n+1-n}\left(  \widetilde{g}\right)  $. Since $n+2-n=2$ and $n+1-n=1$, this
simplifies to $\mathfrak{A}_{2}\left(  \widetilde{p}\rightthreetimes
\widetilde{g}\right)  =\widetilde{p}\rightthreetimes\mathfrak{B}_{1}\left(
\widetilde{g}\right)  $. Since%
\begin{align*}
\mathfrak{B}_{1}
&  =
\left(\overline{R_{1}}\right)_{P}\circ
\left(\overline{R_{2}}\right)_{P}\circ ...\circ
\left(\overline{R_{n}}\right)_{P}
\ \ \ \ \ \ \ \ \ \ \left(  \text{by the
definition of }\mathfrak{B}_{1}\right) \\
&  =\overline{R}_{P},
\end{align*}
this rewrites as
\begin{equation}
\mathfrak{A}_{2}\left(  \widetilde{p}\rightthreetimes\widetilde{g}\right)
=\widetilde{p}\rightthreetimes\overline{R}_{P}\left(  \widetilde{g}\right)  .
\label{pf.Bk.ord.R.pf.ind.pf.6}%
\end{equation}

Now, the definition of $\mathfrak{A}_{2}$ yields
\begin{equation}
\mathfrak{A}_{2}
= \left(\overline{R_{2}}\right)_{B_{k}P}
\circ\left(\overline{R_{3}}\right)_{B_{k}P}
\circ...\circ\left(\overline{R_{n+1}}\right)_{B_{k}P}.
\label{pf.Bk.ord.R.pf.ind.pf.7}
\end{equation}
But recall that%
\[
\overline{R}_{B_{k}P}
= \left(\overline{R_{1}}\right)_{B_{k}P}
\circ \underbrace{\left(\overline{R_{2}}\right)_{B_{k}P}
\circ\left(\overline{R_{3}}\right)_{B_{k}P}\circ...
\circ\left(\overline{R_{n+1}}\right)_{B_{k}P}}_{\substack{
=\mathfrak{A}_{2}\\\text{(by (\ref{pf.Bk.ord.R.pf.ind.pf.7}))}}}
=\left(\overline{R_{1}}\right)_{B_{k}P}\circ \mathfrak{A}_{2},
\]
so that%
\begin{align*}
\overline{R}_{B_{k}P}\left(  \widetilde{p}\rightthreetimes\widetilde{g}%
\right)   &  =\left(  \left(\overline{R_{1}}\right)_{B_{k}P}
\circ\mathfrak{A}_{2}\right)
\left(  \widetilde{p}\rightthreetimes\widetilde{g}\right)
= \left(\overline{R_{1}}\right)_{B_{k}P}\left(
\underbrace{\mathfrak{A}_{2}\left(  \widetilde{p}%
\rightthreetimes\widetilde{g}\right)  }_{\substack{=\widetilde{p}%
\rightthreetimes\overline{R}_{P}\left(  \widetilde{g}\right)  \\\text{(by
(\ref{pf.Bk.ord.R.pf.ind.pf.6}))}}}\right) \\
&  =\left(\overline{R_{1}}\right)_{B_{k}P}\left(
\widetilde{p}\rightthreetimes\overline{R}_{P}
\left(  \widetilde{g}\right)  \right)  =\widetilde{p}^{-1}%
\rightthreetimes\overline{R}_{P}\left(  \widetilde{g}\right)
\end{align*}
(by (\ref{pf.Bk.ord.Ri.small}), applied to $\overline{R}_{P}\left(
\widetilde{g}\right)  $ instead of $\widetilde{g}$). This proves
(\ref{pf.Bk.ord.R}).
\end{verlong}

Furthermore, every $\widetilde{p}\in\mathbb{P}\left(  \mathbb{K}^{k}\right)  $
and $\widetilde{g}\in\overline{\mathbb{K}^{\widehat{P}}}$ satisfy
\begin{equation}
\overline{R}_{B_{k}P}^{\ell}\left(  \widetilde{p}\rightthreetimes
\widetilde{g}\right)  =\widetilde{p}^{\left(  -1\right)  ^{\ell}%
}\rightthreetimes\overline{R}_{P}^{\ell}\left(  \widetilde{g}\right)
\ \ \ \ \ \ \ \ \ \ \text{for all }\ell\in\mathbb{N}. \label{pf.Bk.ord.Rl}%
\end{equation}
(This is proven by induction over $\ell$, using (\ref{pf.Bk.ord.R}).)

We know that the elements of $\overline{\mathbb{K}^{\widehat{B_{k}P}}}$ have
the form $\widetilde{p}\rightthreetimes\widetilde{g}$, where $\widetilde{p}%
\in\mathbb{P}\left(  \mathbb{K}^{k}\right)  $ and $\widetilde{g}\in
\overline{\mathbb{K}^{\widehat{P}}}$. Conversely, every element $\widetilde{p}%
\rightthreetimes\widetilde{g}$ with $\widetilde{p}\in\mathbb{P}\left(
\mathbb{K}^{k}\right)  $ and $\widetilde{g}\in\overline{\mathbb{K}%
^{\widehat{P}}}$ lies in $\overline{\mathbb{K}^{\widehat{B_{k}P}}}$. Hence,
for every $\ell\in\mathbb{N}$, we have the following equivalence of
assertions:%
\begin{align*}
&  \ \left(  \text{we have }\overline{R}_{B_{k}P}^{\ell}=\operatorname*{id}%
\right) \\
&  \Longleftrightarrow\ \left(  \text{every }\widetilde{p}\in\mathbb{P}\left(
\mathbb{K}^{k}\right)  \text{ and }\widetilde{g}\in\overline{\mathbb{K}%
^{\widehat{P}}}\text{ satisfy }\overline{R}_{B_{k}P}^{\ell}\left(
\widetilde{p}\rightthreetimes\widetilde{g}\right)  =\widetilde{p}%
\rightthreetimes\widetilde{g}\right) \\
&  \Longleftrightarrow\ \left(  \text{every }\widetilde{p}\in\mathbb{P}\left(
\mathbb{K}^{k}\right)  \text{ and }\widetilde{g}\in\overline{\mathbb{K}%
^{\widehat{P}}}\text{ satisfy }\widetilde{p}^{\left(  -1\right)  ^{\ell}%
}\rightthreetimes\overline{R}_{P}^{\ell}\left(  \widetilde{g}\right)
=\widetilde{p}\rightthreetimes\widetilde{g}\right) \\
&  \ \ \ \ \ \ \ \ \ \ \left(  \text{because of (\ref{pf.Bk.ord.Rl})}\right)
\\
&  \Longleftrightarrow\ \left(  \text{every }\widetilde{p}\in\mathbb{P}\left(
\mathbb{K}^{k}\right)  \text{ and }\widetilde{g}\in\overline{\mathbb{K}%
^{\widehat{P}}}\text{ satisfy }\widetilde{p}^{\left(  -1\right)  ^{\ell}%
}=\widetilde{p}\text{ and }\overline{R}_{P}^{\ell}\left(  \widetilde{g}%
\right)  =\widetilde{g}\right) \\
&  \Longleftrightarrow\ \left(  \underbrace{\text{every }\widetilde{p}%
\in\mathbb{P}\left(  \mathbb{K}^{k}\right)  \text{ satisfies }\widetilde{p}%
^{\left(  -1\right)  ^{\ell}}=\widetilde{p}}_{\text{this is equivalent to
}\left(  2\mid\ell\text{ if }k>1\right)  }\text{, and }\underbrace{\text{every
}\widetilde{g}\in\overline{\mathbb{K}^{\widehat{P}}}\text{ satisfies
}\overline{R}_{P}^{\ell}\left(  \widetilde{g}\right)  =\widetilde{g}%
}_{\text{this is equivalent to }\overline{R}_{P}^{\ell}=\operatorname*{id}%
}\right) \\
&  \ \ \ \ \ \ \ \ \ \ \left(  \text{since the sets }\mathbb{P}\left(
\mathbb{K}^{k}\right)  \text{ and }\overline{\mathbb{K}^{\widehat{P}}}\text{
are nonempty}\right) \\
&  \Longleftrightarrow\ \left(  \text{we have }\left(  2\mid\ell\text{ if
}k>1\right)  \text{ and }\underbrace{\overline{R}_{P}^{\ell}%
=\operatorname*{id}}_{\text{this is equivalent to }\operatorname*{ord}\left(
\overline{R}_{P}\right)  \mid\ell}\right) \\
&  \Longleftrightarrow\ \left(  \text{we have }\left(  2\mid\ell\text{ if
}k>1\right)  \text{ and }\operatorname*{ord}\left(  \overline{R}_{P}\right)
\mid\ell\right) \\
&  \Longleftrightarrow\ \left\{
\begin{array}
[c]{c}%
\left(  \text{we have }2\mid\ell\text{ and }\operatorname*{ord}\left(
\overline{R}_{P}\right)  \mid\ell\right)  ,\ \ \ \ \ \ \ \ \ \ \text{if
}k>1;\\
\left(  \text{we have }\operatorname*{ord}\left(  \overline{R}_{P}\right)
\mid\ell\right)  ,\ \ \ \ \ \ \ \ \ \ \text{if }k=1
\end{array}
\right. \\
&  \Longleftrightarrow\ \left\{
\begin{array}
[c]{c}%
\left(  \text{we have }\operatorname{lcm}\left(  2,\operatorname*{ord}\left(
\overline{R}_{P}\right)  \right)  \mid\ell\right)
,\ \ \ \ \ \ \ \ \ \ \text{if }k>1;\\
\left(  \text{we have }\operatorname*{ord}\left(  \overline{R}_{P}\right)
\mid\ell\right)  ,\ \ \ \ \ \ \ \ \ \ \text{if }k=1
\end{array}
\right. \\
&  \Longleftrightarrow\ \left(  \text{we have }\left\{
\begin{array}
[c]{l}%
\operatorname{lcm}\left(  2,\operatorname*{ord}\left(  \overline{R}%
_{P}\right)  \right)  ,\ \ \ \ \ \ \ \ \ \ \text{if }k>1;\\
\operatorname*{ord}\left(  \overline{R}_{P}\right)
,\ \ \ \ \ \ \ \ \ \ \text{if }k=1
\end{array}
\right.  \mid\ell\right)  .
\end{align*}
Hence, for every $\ell\in\mathbb{N}$, we have the following equivalence of
assertions:%
\begin{align*}
\left(  \text{we have }\operatorname*{ord}\left(  \overline{R}_{B_{k}%
P}\right)  \mid\ell\right)  \  &  \Longleftrightarrow\ \left(  \text{we have
}\overline{R}_{B_{k}P}^{\ell}=\operatorname*{id}\right) \\
&  \Longleftrightarrow\ \left(  \text{we have }\left\{
\begin{array}
[c]{l}%
\operatorname{lcm}\left(  2,\operatorname*{ord}\left(  \overline{R}%
_{P}\right)  \right)  ,\ \ \ \ \ \ \ \ \ \ \text{if }k>1;\\
\operatorname*{ord}\left(  \overline{R}_{P}\right)
,\ \ \ \ \ \ \ \ \ \ \text{if }k=1
\end{array}
\right.  \mid\ell\right)  .
\end{align*}
Consequently, $\operatorname*{ord}\left(  \overline{R}_{B_{k}P}\right)
=\left\{
\begin{array}
[c]{l}%
\operatorname{lcm}\left(  2,\operatorname*{ord}\left(  \overline{R}%
_{P}\right)  \right)  ,\ \ \ \ \ \ \ \ \ \ \text{if }k>1;\\
\operatorname*{ord}\left(  \overline{R}_{P}\right)
,\ \ \ \ \ \ \ \ \ \ \text{if }k=1
\end{array}
\right.  $. This is exactly what (\ref{pf.Bk.ord.1}) claims. Thus,
(\ref{pf.Bk.ord.1}) is proven, and with it Proposition \ref{prop.Bk.ord}.
\end{proof}

Here is an analogue of Proposition \ref{prop.Bk.ord}:

\begin{proposition}
\label{prop.B'k.ord}Let $n\in\mathbb{N}$. Let $P$ be an $n$-graded poset. Let
$\mathbb{K}$ be a field.

\textbf{(a)} We have $\operatorname*{ord}\left(  \overline{R}_{B_{1}^{\prime
}P}\right)  =\operatorname*{ord}\left(  \overline{R}_{P}\right)  $.

\textbf{(b)} For every integer $k>1$, we have $\operatorname*{ord}\left(
\overline{R}_{B_{k}^{\prime}P}\right)  =\operatorname{lcm}\left(
2,\operatorname*{ord}\left(  \overline{R}_{P}\right)  \right)  $.
\end{proposition}

\begin{proof}
[\nopunct]The proof of this is very similar (though not exactly identical) to
that of Proposition \ref{prop.Bk.ord}. Alternatively, it is easy to deduce
Proposition \ref{prop.B'k.ord} from Proposition \ref{prop.Bk.ord} using
Proposition \ref{prop.op.ord} and Proposition \ref{prop.skeletal.op}.
\end{proof}

Proposition \ref{prop.skeletal.ords} is easily shown by induction using
Propositions \ref{prop.PQ.ord}, \ref{prop.Bk.ord}, \ref{prop.B'k.ord} and
\ref{prop.ord-projord}. Moreover, using Propositions \ref{prop.PQ.ord},
\ref{prop.PQ.ord2}, \ref{prop.Bk.ord}, \ref{prop.B'k.ord} and
\ref{prop.ord-projord}, we can recursively compute (rather than just bound
from the above) the orders of $R_{P}$ and $\overline{R}_{P}$ for any skeletal
poset $P$ without doing any computations in $\mathbb{K}$. (This also shows
that the orders of $R_{P}$ and $\overline{R}_{P}$ don't depend on the base
field $\mathbb{K}$ as long as $\mathbb{K}$ is infinite and $P$ is skeletal.)

In the case of forests and trees we can also use this induction to establish a
concrete bound:

\begin{corollary}
\label{cor.for.ord}Let $n\in\mathbb{N}$. Let $P$ be an $n$-graded poset. Let
$\mathbb{K}$ be a field. Assume that $P$ is a rooted forest (made into a poset
by having every node smaller than its children).

\textbf{(a)} Then, $\operatorname*{ord}\left(  R_{P}\right)  \mid
\operatorname{lcm}\left(  1,2,...,n+1\right)  $.

\textbf{(b)} Moreover, if $P$ is a tree, then $\operatorname*{ord}\left(
\overline{R}_{P}\right)  \mid\operatorname{lcm}\left(  1,2,...,n\right)  $.
\end{corollary}

Corollary \ref{cor.for.ord} is also valid if we replace ``every node smaller
than its children'' by ``every node larger than its children'', and the proof
is exactly analogous.

\begin{vershort}
\begin{proof}
[Proof of Corollary \ref{cor.for.ord} (sketched).]\textbf{(a)} Corollary
\ref{cor.for.ord} \textbf{(a)} can be proven by strong induction over
$\left\vert P\right\vert $. Indeed, if $P$ is an $n$-graded poset and a rooted
forest, then we must be in one of the following three cases:

\textit{Case 1:} We have $P=\varnothing$.

\textit{Case 2:} The rooted forest $P$ is a tree.

\textit{Case 3:} The rooted forest $P$ is a disjoint union of more than one tree.

The validity of Corollary \ref{cor.for.ord} is trivial in Case 1, and in Case
3 it follows from the induction hypothesis using Proposition \ref{prop.PQ.ord}%
. In Case 2, we have $P=B_{1}Q$ for some rooted forest $Q$, which is
necessarily $\left(  n-1\right)  $-graded; thus, the induction hypothesis
(applied to $Q$ instead of $P$) yields $\operatorname*{ord}\left(
R_{Q}\right)  \mid\operatorname{lcm}\left(  1,2,...,\left(  n-1\right)
+1\right)  =\operatorname{lcm}\left(  1,2,...,n\right)  $, and we obtain%
\begin{align*}
\operatorname*{ord}\left(  \overline{R}_{P}\right)   &  =\operatorname*{ord}%
\left(  \overline{R}_{B_{1}Q}\right)  =\operatorname*{ord}\left(  \overline
{R}_{Q}\right)  \ \ \ \ \ \ \ \ \ \ \left(  \text{by Proposition
\ref{prop.Bk.ord} \textbf{(a)}}\right) \\
&  \mid\operatorname{lcm}\left(  1,2,...,n\right)
\end{align*}
and%
\begin{align*}
\operatorname*{ord}\left(  R_{P}\right)   &  =\operatorname{lcm}\left(
n+1,\underbrace{\operatorname*{ord}\left(  \overline{R}_{P}\right)  }%
_{\mid\operatorname{lcm}\left(  1,2,...,n\right)  }\right)
\ \ \ \ \ \ \ \ \ \ \left(  \text{by Proposition \ref{prop.ord-projord}%
}\right) \\
&  \mid\operatorname{lcm}\left(  n+1,\operatorname{lcm}\left(
1,2,...,n\right)  \right)  =\operatorname{lcm}\left(  1,2,...,n+1\right)  .
\end{align*}
Thus, the induction step is complete in each of the three Cases.

\textbf{(b)} If $P$ is a tree, then we must be in Case 2 of the above case
distinction, and thus we have $\operatorname*{ord}\left(  \overline{R}%
_{P}\right)  \mid\operatorname{lcm}\left(  1,2,...,n\right)  $ as shown above.
Corollary \ref{cor.for.ord} is therefore proven.
\end{proof}
\end{vershort}

\begin{verlong}
\begin{proof}
[Proof of Corollary \ref{cor.for.ord} (sketched).]\textbf{(a)} We shall prove
Corollary \ref{cor.for.ord} \textbf{(a)} by strong induction over $\left\vert
P\right\vert $. So we fix a $K\in\mathbb{N}$, and we assume that Corollary
\ref{cor.for.ord} \textbf{(a)} holds under the assumption that $\left\vert
P\right\vert <K$. We now need to show that Corollary \ref{cor.for.ord}
\textbf{(a)} also holds under the assumption that $\left\vert P\right\vert =K$.

We have assumed that Corollary \ref{cor.for.ord} \textbf{(a)} holds under the
assumption that $\left\vert P\right\vert <K$. In other words,%
\begin{equation}
\left(
\begin{array}
[c]{c}%
\text{every }n\in\mathbb{N}\text{, every }n\text{-graded poset }P\text{ and
every field }\mathbb{K}\text{ such that }P\text{ is a rooted}\\
\text{forest such that }\left\vert P\right\vert <K\text{ satisfy
}\operatorname*{ord}\left(  R_{P}\right)  \mid\operatorname{lcm}\left(
1,2,...,n+1\right)
\end{array}
\right)  . \label{pf.for.ord.1}%
\end{equation}

Now, let us show that Corollary \ref{cor.for.ord} \textbf{(a)} also holds
under the assumption that $\left\vert P\right\vert =K$. So let $n\in
\mathbb{N}$, let $\mathbb{K}$ be a field, and let $P$ be an $n$-graded poset
such that $P$ is a rooted forest and such that $\left\vert P\right\vert =K$.
Due to the structure of rooted forests, we must then be in one of the
following three cases:

\textit{Case 1:} We have $P=\varnothing$.

\textit{Case 2:} The rooted forest $P$ is a tree.

\textit{Case 3:} The rooted forest $P$ is a disjoint union of more than one tree.

Case 1 is trivial and left to the reader.

Let us now consider Case 2. In this case, the rooted forest $P$ is a tree.
Hence, as a poset, we have $P=B_{1}Q$ for some rooted forest $Q$. Consider
this $Q$. Since $P=B_{1}Q$, we have $\left\vert P\right\vert =\left\vert
B_{1}Q\right\vert =\left\vert Q\right\vert +1$, so that $\left\vert
Q\right\vert =\underbrace{\left\vert P\right\vert }_{=K}-1=K-1<K$. Moreover,
$Q$ is easily seen to be an $\left(  n-1\right)  $-graded poset (since
$B_{1}Q=P$ is an $n$-graded poset). Hence, we can apply (\ref{pf.for.ord.1})
to $Q$ instead of $P$. As a result, we obtain that $\operatorname*{ord}\left(
R_{Q}\right)  \mid\operatorname{lcm}\left(  1,2,...,\left(  n-1\right)
+1\right)  =\operatorname{lcm}\left(  1,2,...,n\right)  $.

Proposition \ref{prop.ord-projord} (applied to $Q$ instead of $P$) yields that
$\operatorname*{ord}\left(  R_{Q}\right)  =\operatorname{lcm}\left(
n+1,\operatorname*{ord}\left(  \overline{R}_{Q}\right)  \right)  $. Now,%
\[
\operatorname*{ord}\left(  \overline{R}_{Q}\right)  \mid\operatorname{lcm}%
\left(  n+1,\operatorname*{ord}\left(  \overline{R}_{Q}\right)  \right)
=\operatorname*{ord}\left(  R_{Q}\right)  \mid\operatorname{lcm}\left(
1,2,...,n\right)  .
\]

But Proposition \ref{prop.Bk.ord} \textbf{(a)} (applied to $Q$ instead of $P$)
yields $\operatorname*{ord}\left(  \overline{R}_{B_{1}Q}\right)
=\operatorname*{ord}\left(  \overline{R}_{Q}\right)  $.

Since $Q=B_{1}P$, we have%
\begin{equation}
\operatorname*{ord}\left(  \overline{R}_{P}\right)  =\operatorname*{ord}%
\left(  \overline{R}_{B_{1}Q}\right)  =\operatorname*{ord}\left(  \overline
{R}_{Q}\right)  \mid\operatorname{lcm}\left(  1,2,...,n\right)  .
\label{pf.for.ord.4}%
\end{equation}
Applying Proposition \ref{prop.ord-projord}, we obtain%
\begin{align*}
\operatorname*{ord}\left(  R_{P}\right)   &  =\operatorname{lcm}\left(
n+1,\underbrace{\operatorname*{ord}\left(  \overline{R}_{P}\right)
}_{\substack{\mid\operatorname{lcm}\left(  1,2,...,n\right)  \\\text{(by
(\ref{pf.for.ord.4}))}}}\right) \\
&  \mid\operatorname{lcm}\left(  n+1,\operatorname{lcm}\left(
1,2,...,n\right)  \right)  =\operatorname{lcm}\left(  1,2,...,n+1\right)  .
\end{align*}
We thus have proven Corollary \ref{cor.for.ord} \textbf{(a)} for our poset $P$
in Case 2.

Finally, let us consider Case 3. In this case, the rooted forest $P$ is a
disjoint union of more than one tree. Hence, the rooted forest $P$ is a
disjoint union of two \textbf{nonempty} rooted forests\footnote{These forests
don't necessarily have to be trees.}. In other words, there exist two nonempty
rooted forests $Q_{1}$ and $Q_{2}$ such that $P=Q_{1}Q_{2}$. Consider these
$Q_{1}$ and $Q_{2}$. Clearly, $\left\vert Q_{2}\right\vert >0$ (since $Q_{2}$
is nonempty).

Since $P=Q_{1}Q_{2}$, we have $\left\vert P\right\vert =\left\vert Q_{1}%
Q_{2}\right\vert =\left\vert Q_{1}\right\vert +\left\vert Q_{2}\right\vert
>\left\vert Q_{1}\right\vert $ (since $\left\vert Q_{2}\right\vert >0$), so
that $\left\vert Q_{1}\right\vert <\left\vert P\right\vert =K$. Also, $Q_{1}$
is an $n$-graded poset (this is easy to prove because $Q_{1}Q_{2}$ is an
$n$-graded poset). Hence, we can apply (\ref{pf.for.ord.1}) to $Q_{1}$ instead
of $P$. As a result, we obtain that $\operatorname*{ord}\left(  R_{Q_{1}%
}\right)  \mid\operatorname{lcm}\left(  1,2,...,n+1\right)  $. The same
argument (but with $Q_{1}$ and $Q_{2}$ switched) yields $\operatorname*{ord}%
\left(  R_{Q_{2}}\right)  \mid\operatorname{lcm}\left(  1,2,...,n+1\right)  $.
Now, $P=Q_{1}Q_{2}$, so that%
\begin{align*}
\operatorname*{ord}\left(  R_{P}\right)   &  =\operatorname*{ord}\left(
R_{Q_{1}Q_{2}}\right)  =\operatorname{lcm}\left(
\underbrace{\operatorname*{ord}\left(  R_{Q_{1}}\right)  }_{\mid
\operatorname{lcm}\left(  1,2,...,n+1\right)  }%
,\underbrace{\operatorname*{ord}\left(  R_{Q_{2}}\right)  }_{\mid
\operatorname{lcm}\left(  1,2,...,n+1\right)  }\right) \\
&  \ \ \ \ \ \ \ \ \ \ \left(  \text{by Proposition \ref{prop.PQ.ord}, applied
to }Q_{1}\text{ and }Q_{2}\text{ instead of }P\text{ and }Q\right) \\
&  \mid\operatorname{lcm}\left(  \operatorname{lcm}\left(  1,2,...,n+1\right)
,\operatorname{lcm}\left(  1,2,...,n+1\right)  \right)  =\operatorname{lcm}%
\left(  1,2,...,n+1\right)  .
\end{align*}
We thus have proven Corollary \ref{cor.for.ord} \textbf{(a)} for our poset $P$
in Case 3.

We have now proven Corollary \ref{cor.for.ord} \textbf{(a)} for our poset $P$
in each of the three Cases 1, 2 and 3. Since these three Cases cover all
possibilities, this yields that Corollary \ref{cor.for.ord} \textbf{(a)} holds
for our poset $P$ in every situation.

Now, forget that we fixed $P$. We thus have shown that if $n\in\mathbb{N}$, if
$\mathbb{K}$ is a field, and if $P$ is an $n$-graded poset such that $P$ is a
rooted forest and such that $\left\vert P\right\vert =K$, then Corollary
\ref{cor.for.ord} \textbf{(a)} holds for this poset $P$. In other words,
Corollary \ref{cor.for.ord} \textbf{(a)} holds under the assumption that
$\left\vert P\right\vert =K$. This completes the induction step. The induction
proof of Corollary \ref{cor.for.ord} \textbf{(a)} is thus complete.

\textbf{(b)} Now, assume that $P$ is a tree. Then, (\ref{pf.for.ord.4}) holds
(this can be proven just as in Case 2 of the proof of Corollary
\ref{cor.for.ord} \textbf{(a)}, but using Corollary \ref{cor.for.ord}
\textbf{(a)} instead of (\ref{pf.for.ord.1})). This proves Corollary
\ref{cor.for.ord} \textbf{(b)}.

[Some more details might be useful in this proof?]
\end{proof}
\end{verlong}

\section{\label{sect.classicalr}Interlude: Classical rowmotion on skeletal
posets}

The above results concerning birational rowmotion on skeletal posets suggest
the question of what can be said about \textbf{classical} rowmotion (on the
set of order ideals) on this class of posets. Indeed, while the classical
rowmotion map (as opposed to the birational one) has been the object of
several studies (e.g., \cite{striker-williams} and
\cite{cameron-fon-der-flaass}), it seems that this rather simple case has
never been explicitly covered. Let us therefore go on a tangent to bridge this
gap and derive the counterparts of Propositions \ref{prop.Bk.ord} and
\ref{prop.skeletal.ords} and Corollary \ref{cor.for.ord} for classical
rowmotion. Nothing of what we do in this Section \ref{sect.classicalr} will be
relevant to later sections, so this section can be skipped.

First, we define the maps involved.

\begin{definition}
Let $P$ be a poset.

\textbf{(a)} An \textit{order ideal} of $P$ means a subset $S$ of $P$ such
that every $s\in S$ and $p\in P$ with $p\leq s$ satisfy $p\in S$.

\textbf{(b)} The set of all order ideals of $P$ will be denoted by $J\left(
P\right)  $.
\end{definition}

Here is the definition of (classical) toggles on order ideals (an analogue of
Definition \ref{def.Tv}):

\begin{definition}
\label{def.classical.tv}Let $P$ be a finite poset. Let $v\in P$. Define a map
$\mathbf{t}_{v}:J\left(  P\right)  \rightarrow J\left(  P\right)  $ by%
\[
\mathbf{t}_{v}\left(  S\right)  =\left\{
\begin{array}
[c]{l}%
S\cup\left\{  v\right\}  \text{, if }v\notin S\text{ and }S\cup\left\{
v\right\}  \in J\left(  P\right)  ;\\
S\setminus\left\{  v\right\}  \text{, if }v\in S\text{ and }S\setminus\left\{
v\right\}  \in J\left(  P\right)  ;\\
S\text{, otherwise}%
\end{array}
\right.  \ \ \ \ \ \ \ \ \ \ \text{for every }S\in J\left(  P\right)  .
\]
(This is clearly well-defined.) This map $\mathbf{t}_{v}$ will be called the
\textit{classical }$v$\textit{-toggle}.
\end{definition}

We can rewrite this definition in more \textquotedblleft
local\textquotedblright\ terms, by replacing the conditions \textquotedblleft%
$S\cup\left\{  v\right\}  \in J\left(  P\right)  $\textquotedblright\ and
\textquotedblleft$S\setminus\left\{  v\right\}  \in J\left(  P\right)
$\textquotedblright\ by the respectively equivalent conditions
\textquotedblleft every element $u\in P$ satisfying $u\lessdot v$ lies in
$S$\textquotedblright\ and \textquotedblleft no element $u\in P$ satisfying
$u\gtrdot v$ lies in $S$\textquotedblright\ (in fact, the equivalence of these
conditions is easily seen). Hence, we obtain the following analogue to
Proposition \ref{prop.Tv}:

\begin{proposition}
\label{prop.classical.Tv}Let $P$ be a finite poset. Let $v\in P$. For every
$S\in J\left(  P\right)  $, we have:

\textbf{(a)} If $w$ is an element of $P$ such that $w\neq v$, then we have
$w\in\mathbf{t}_{v}\left(  S\right)  $ if and only if $w\in S$.

\textbf{(b)} We have $v\in\mathbf{t}_{v}\left(  S\right)  $ if and only if
\begin{align*}
&  \left(  v\in S\text{ and not }\left(  \text{no element }u\in P\text{
satisfying }u\gtrdot v\text{ lies in }S\right)  \right) \\
&  \text{or }\left(  v\notin S\text{ and }\left(  \text{every element }u\in
P\text{ satisfying }u\lessdot v\text{ lies in }S\right)  \right)  .
\end{align*}

\end{proposition}

While the complicated logical statement in Proposition \ref{prop.classical.Tv}
\textbf{(b)} can be simplified, the form we have stated it in exhibits its
similarity to Proposition \ref{prop.Tv} particularly well. This, in fact, is
more than a similarity: If we allow $\mathbb{K}$ to be a semifield rather than
a field, we can regard the classical $v$-toggle $\mathbf{t}_{v}$ as a
restriction of the birational toggle $T_{v}$ (when $\mathbb{K}$ is chosen
appropriately)\footnote{Here are the details: Let $\operatorname*{Trop}%
\mathbb{Z}$ be the tropical semiring over $\mathbb{Z}$, that is, the semiring
obtained by endowing the set $\mathbb{Z}\cup\left\{  -\infty\right\}  $ with
the binary operation $\left(  a,b\right)  \mapsto\max\left\{  a,b\right\}  $
as \textquotedblleft addition\textquotedblright\ and the binary operation
$\left(  a,b\right)  \mapsto a+b$ as \textquotedblleft
multiplication\textquotedblright\ (where the usual rules for sums involving
$-\infty$ apply). Then, $\operatorname*{Trop}\mathbb{Z}$ is a semifield, with
$\left(  a,b\right)  \mapsto a-b$ serving as \textquotedblleft
subtraction\textquotedblright, with $-\infty$ serving as \textquotedblleft
zero\textquotedblright\ and with the integer $0$ serving as \textquotedblleft
one\textquotedblright. Now, to every order ideal $S\in J\left(  P\right)  $,
we can assign a $\left(  \operatorname*{Trop}\mathbb{Z}\right)  $-labelling
$\operatorname*{tlab}S\in\left(  \operatorname*{Trop}\mathbb{Z}\right)
^{\widehat{P}}$, defined by%
\[
\left(  \operatorname*{tlab}S\right)  \left(  v\right)  =\left\{
\begin{array}
[c]{c}%
1,\text{ if }v\notin S\cup\left\{  0\right\}  ;\\
0,\ \text{if }v\in S\cup\left\{  0\right\}
\end{array}
\right.  .
\]
This yields a map $\operatorname*{tlab}:J\left(  P\right)  \rightarrow\left(
\operatorname*{Trop}\mathbb{Z}\right)  ^{\widehat{P}}$, obviously injective.
This map $\operatorname*{tlab}$ satisfies $T_{v}\circ\operatorname*{tlab}%
=\operatorname*{tlab}\circ\mathbf{t}_{v}$ for every $v\in P$. This allows us
to regard the classical toggles $\mathbf{t}_{v}$ as restrictions of the
birational toggles $T_{v}$, if we consider this map $\operatorname*{tlab}$ as
an inclusion. This reasoning goes back to Einstein and Propp
\cite{einstein-propp}.}. Hence, some theorems about birational toggles can be
used to derive analogous theorems about classical toggles\footnote{For
example, we could derive Proposition \ref{prop.classical.tv.commute} from
Proposition \ref{prop.Tv.commute} using this tactic. However, we could not
derive (say) Proposition \ref{prop.classical.skeletal.ords} from Proposition
\ref{prop.skeletal.ords} this way, because the order of a restriction of a
permutation could be a proper divisor of the order of the permutation.}. We
will not use this tactic in the following, because often it will be easier to
study the classical $v$-toggles on their own. However, many of the properties
of classical toggles (and classical rowmotion) that we are going to discuss
will have proofs that are parallel to the proofs of the analogous results
about birational toggles. We will omit these proofs when the analogy is
glaring enough.

We have the following easily-verified analogues of Proposition
\ref{prop.Tv.invo}, Proposition \ref{prop.Tv.commute} and Corollary
\ref{cor.R.welldef}:

\begin{proposition}
\label{prop.classical.tv.invo}Let $P$ be a finite poset. Let $v\in P$. Then,
the map $\mathbf{t}_{v}$ is an involution on $J\left(  P\right)  $ (that is,
we have $\mathbf{t}_{v}^{2}=\operatorname*{id}$).
\end{proposition}

\begin{proposition}
\label{prop.classical.tv.commute}Let $P$ be a finite poset. Let $v\in P$ and
$w\in P$. Then, $\mathbf{t}_{v}\circ\mathbf{t}_{w}=\mathbf{t}_{w}%
\circ\mathbf{t}_{v}$, unless we have either $v\lessdot w$ or $w\lessdot v$.
\end{proposition}

\begin{corollary}
\label{cor.classical.r.welldef}Let $P$ be a finite poset. Let $\left(
v_{1},v_{2},...,v_{m}\right)  $ be a linear extension of $P$. Then, the map
$\mathbf{t}_{v_{1}}\circ\mathbf{t}_{v_{2}}\circ...\circ\mathbf{t}_{v_{m}%
}:J\left(  P\right)  \rightarrow J\left(  P\right)  $ is well-defined and
independent of the choice of the linear extension $\left(  v_{1}%
,v_{2},...,v_{m}\right)  $.
\end{corollary}

The three results above are observations made on \cite[page 546]%
{cameron-fon-der-flaass} (in somewhat different notation).

Two convenient advantages of the classical setup are that we don't have to
worry about denominators becoming zero, so our maps are actual maps rather
than partial maps, and that we don't have to pass to the poset $\widehat{P}$.

We can now define rowmotion in analogy to Definition \ref{def.rm}:

\begin{definition}
\label{def.classical.rm}Let $P$ be a finite poset. \textit{Classical
rowmotion} (simply called \textquotedblleft rowmotion\textquotedblright\ in
existing literature) is defined as the map $\mathbf{t}_{v_{1}}\circ
\mathbf{t}_{v_{2}}\circ...\circ\mathbf{t}_{v_{m}}:J\left(  P\right)
\rightarrow J\left(  P\right)  $, where $\left(  v_{1},v_{2},...,v_{m}\right)
$ is a linear extension of $P$. This map is well-defined (in particular, it
does not depend on the linear extension $\left(  v_{1},v_{2},...,v_{m}\right)
$ chosen) because of Corollary \ref{cor.classical.r.welldef} (and also because
of the fact that a linear extension of $P$ exists; this is Theorem
\ref{thm.linext.ex}). This map will be denoted by $\mathbf{r}$.
\end{definition}

To highlight the similarities between the classical and birational cases, let
us state the analogue of Proposition \ref{prop.R.implicit}:

\begin{proposition}
\label{prop.classical.r.implicit}Let $P$ be a finite poset. Let $v\in P$. Let
$S\in J\left(  P\right)  $. Then, $v\in\mathbf{r}\left(  S\right)  $ holds if
and only if the following two conditions hold:

\textit{Condition 1:} Every $u\in P$ satisfying $u\lessdot v$ belongs to $S$.

\textit{Condition 2:} Either $v\notin S$, or there exists an $u\in
\mathbf{r}\left(  S\right)  $ satisfying $u\gtrdot v$. (Recall that the
expression \textquotedblleft either/or\textquotedblright\ is meant non-exclusively.)
\end{proposition}

This proposition is easily seen to be equivalent to the following well-known
equivalent description of rowmotion (\cite[Lemma 1]{cameron-fon-der-flaass},
translated into our notation):

\begin{proposition}
\label{prop.classical.r.maxmin}Let $P$ be a finite poset. Let $S\in J\left(
P\right)  $. Then, the maximal elements of $\mathbf{r}\left(  S\right)  $ are
precisely the minimal elements of $P\setminus S$.
\end{proposition}

We record the analogue of Proposition \ref{prop.R.implicit.converse}:

\begin{proposition}
\label{prop.classical.r.implicit.converse}Let $P$ be a finite poset. Let $S$
and $T$ be two order ideals of $P$. Assume that for every $v\in P$, the
relation $v\in T$ holds if and only if Conditions 1 and 2 of Proposition
\ref{prop.classical.r.implicit} hold with $\mathbf{r}\left(  S\right)  $
replaced by $T$. Then, $T=\mathbf{r}\left(  S\right)  $.
\end{proposition}

In analogy to Proposition \ref{prop.R.inverse}, we have:

\begin{proposition}
\label{prop.classical.r.inverse}Let $P$ be a finite poset. Then, classical
rowmotion $\mathbf{r}$ is invertible. Its inverse $\mathbf{r}^{-1}$ is
$\mathbf{t}_{v_{m}}\circ\mathbf{t}_{v_{m-1}}\circ...\circ\mathbf{t}_{v_{1}%
}:J\left(  P\right)  \rightarrow J\left(  P\right)  $, where $\left(
v_{1},v_{2},...,v_{m}\right)  $ is a linear extension of $P$.
\end{proposition}

We can study graded posets again. In analogy to Corollary \ref{cor.Ri.welldef}%
, Definition \ref{def.Ri}, Proposition \ref{prop.Ri.R} and Proposition
\ref{prop.Ri.invo}, we have:

\begin{corollary}
\label{cor.classical.ri.welldef}Let $n\in\mathbb{N}$. Let $P$ be an $n$-graded
poset. Let $i\in\left\{  1,2,...,n\right\}  $. Let $\left(  u_{1}%
,u_{2},...,u_{k}\right)  $ be any list of the elements of $\widehat{P}_{i}$
with every element of $\widehat{P}_{i}$ appearing exactly once in the list.
(Note that $\widehat{P}_{i}$ is simply $\left\{  v\in P\ \mid\ \deg
v=i\right\}  $, because $i$ equals neither $0$ nor $n+1$. We are using the
notation $\widehat{P}_{i}$ despite not working with $\widehat{P}$ merely to
stress some analogies.) Then, the map $\mathbf{t}_{u_{1}}\circ\mathbf{t}%
_{u_{2}}\circ...\circ\mathbf{t}_{u_{k}}:J\left(  P\right)  \rightarrow
J\left(  P\right)  $ is well-defined and independent of the choice of the list
$\left(  u_{1},u_{2},...,u_{k}\right)  $.
\end{corollary}

\begin{definition}
\label{def.classical.ri}Let $n\in\mathbb{N}$. Let $P$ be an $n$-graded poset.
Let $i\in\left\{  1,2,...,n\right\}  $. Then, let $\mathbf{r}_{i}$ denote the
map $\mathbf{t}_{u_{1}}\circ\mathbf{t}_{u_{2}}\circ...\circ\mathbf{t}_{u_{k}%
}:J\left(  P\right)  \rightarrow J\left(  P\right)  $, where $\left(
u_{1},u_{2},...,u_{k}\right)  $ is any list of the elements of $\widehat{P}%
_{i}$ with every element of $\widehat{P}_{i}$ appearing exactly once in the
list. This map $\mathbf{t}_{u_{1}}\circ\mathbf{t}_{u_{2}}\circ...\circ
\mathbf{t}_{u_{k}}$ is well-defined (in particular, it does not depend on the
list $\left(  u_{1},u_{2},...,u_{k}\right)  $) because of Corollary
\ref{cor.classical.ri.welldef}.
\end{definition}

\begin{proposition}
\label{prop.classical.ri.r}Let $n\in\mathbb{N}$. Let $P$ be an $n$-graded
poset. Then,%
\[
\mathbf{r}=\mathbf{r}_{1}\circ\mathbf{r}_{2}\circ...\circ\mathbf{r}_{n}.
\]

\end{proposition}

\begin{proposition}
\label{prop.classical.ri.invo}Let $n\in\mathbb{N}$. Let $P$ be an $n$-graded
poset. Let $i\in\left\{  1,2,...,n\right\}  $. Then, $\mathbf{r}_{i}$ is an
involution on $J\left(  P\right)  $ (that is, $\mathbf{r}_{i}^{2}%
=\operatorname*{id}$).
\end{proposition}

A parody of w-tuples can also be defined. The following is analogous to
Definition \ref{def.wi}:

\begin{definition}
\label{def.classical.wi}Let $n\in\mathbb{N}$. Let $P$ be an $n$-graded poset.
Let $S\in J\left(  P\right)  $. Let $i\in\left\{  0,1,...,n\right\}  $. Then,
$\mathbf{w}_{i}\left(  S\right)  $ will denote the integer%
\[
\left\{
\begin{array}
[c]{l}%
1,\text{ if }P_{i}\subseteq S\text{ and }P_{i+1}\cap S=\varnothing\\
0,\text{ otherwise}%
\end{array}
\right.  .
\]
Here, we are using the notation $P_{j}$ for the subset $\deg^{-1}\left(
\left\{  j\right\}  \right)  $ of $P$; this subset is empty if $j=0$ and also
empty if $j=n+1$.
\end{definition}

Analogues of Proposition \ref{prop.wi.Ri} and Proposition \ref{prop.wi.R} are
easily found:

\begin{proposition}
\label{prop.classical.wi.ri}Let $n\in\mathbb{N}$. Let $P$ be an $n$-graded
poset. Let $i\in\left\{  1,2,...,n\right\}  $. Then, every $S\in J\left(
P\right)  $ satisfies%
\begin{align*}
&  \left(  \mathbf{w}_{0}\left(  \mathbf{r}_{i}\left(  S\right)  \right)
,\mathbf{w}_{1}\left(  \mathbf{r}_{i}\left(  S\right)  \right)
,...,\mathbf{w}_{n}\left(  \mathbf{r}_{i}\left(  S\right)  \right)  \right) \\
&  =\left(  \mathbf{w}_{0}\left(  S\right)  ,\mathbf{w}_{1}\left(  S\right)
,...,\mathbf{w}_{i-2}\left(  S\right)  ,\mathbf{w}_{i}\left(  S\right)
,\mathbf{w}_{i-1}\left(  S\right)  ,\mathbf{w}_{i+1}\left(  S\right)
,\mathbf{w}_{i+2}\left(  S\right)  ,...,\mathbf{w}_{n}\left(  S\right)
\right)  .
\end{align*}

\end{proposition}

\begin{proposition}
\label{prop.classical.wi.r}Let $n\in\mathbb{N}$. Let $P$ be an $n$-graded
poset. Then, every $S\in J\left(  P\right)  $ satisfies%
\[
\left(  \mathbf{w}_{0}\left(  \mathbf{r}\left(  S\right)  \right)
,\mathbf{w}_{1}\left(  \mathbf{r}\left(  S\right)  \right)  ,...,\mathbf{w}%
_{n}\left(  \mathbf{r}\left(  S\right)  \right)  \right)  =\left(
\mathbf{w}_{n}\left(  S\right)  ,\mathbf{w}_{0}\left(  S\right)
,\mathbf{w}_{1}\left(  S\right)  ,...,\mathbf{w}_{n-1}\left(  S\right)
\right)  .
\]

\end{proposition}

\begin{verlong}
Here is an analogue of Corollary \ref{cor.wi.Rn+1}:

\begin{corollary}
\label{cor.classical.wi.rn+1}Let $n\in\mathbb{N}$. Let $P$ be an $n$-graded
poset. Then, every $S\in J\left(  P\right)  $ satisfies%
\[
\left(  \mathbf{w}_{0}\left(  \mathbf{r}^{n+1}\left(  S\right)  \right)
,\mathbf{w}_{1}\left(  \mathbf{r}^{n+1}\left(  S\right)  \right)
,...,\mathbf{w}_{n}\left(  \mathbf{r}^{n+1}\left(  S\right)  \right)  \right)
=\left(  \mathbf{w}_{0}\left(  S\right)  ,\mathbf{w}_{1}\left(  S\right)
,...,\mathbf{w}_{n}\left(  S\right)  \right)  .
\]

\end{corollary}
\end{verlong}

However, the $\left(  n+1\right)  $-tuple $\left(  \mathbf{w}_{0}\left(
S\right)  ,\mathbf{w}_{1}\left(  S\right)  ,...,\mathbf{w}_{n}\left(
S\right)  \right)  $ obtained from an order ideal $S$ is not particularly
informative. In fact, it is $\left(  0,0,...,0\right)  $ for ``most'' order
ideals; here is what this means precisely:

\begin{definition}
\label{def.classical.level}Let $n\in\mathbb{N}$. Let $P$ be an $n$-graded
poset. An order ideal of $P$ is said to be \textit{level} if and only if it
has the form $P_{1}\cup P_{2}\cup...\cup P_{i}$ for some $i\in\left\{
0,1,...,n\right\}  $.
\end{definition}

Easy properties of level order ideals are:

\begin{proposition}
\label{prop.classical.level}Let $n\in\mathbb{N}$. Let $P$ be an $n$-graded poset.

\textbf{(a)} There exist precisely $n+1$ level order ideals of $P$, and those
form an orbit under classical rowmotion $\mathbf{r}$. Namely, one has%
\[
\mathbf{r}\left(  P_{1}\cup P_{2}\cup...\cup P_{i}\right)  =\left\{
\begin{array}
[c]{l}%
P_{1}\cup P_{2}\cup...\cup P_{i+1},\ \ \ \ \ \ \ \ \ \ \text{if }i<n;\\
\varnothing,\ \ \ \ \ \ \ \ \ \ \text{if }i=n
\end{array}
\right.  .
\]

\textbf{(b)} If $S\in J\left(  P\right)  $, then $\left(  \mathbf{w}%
_{0}\left(  S\right)  ,\mathbf{w}_{1}\left(  S\right)  ,...,\mathbf{w}%
_{n}\left(  S\right)  \right)  =\left(  0,0,...,0\right)  $ unless $S$ is level.
\end{proposition}

Now, we can define an (arguably toylike, but, as we will see, rather useful)
analogue of homogeneous equivalence. In somewhat questionable analogy with
Definition \ref{def.hgeq}, we set:

\begin{definition}
\label{def.classical.hgeq}Let $n\in\mathbb{N}$. Let $P$ be an $n$-graded poset.

Two order ideals $S$ and $T$ of $P$ are said to be \textit{homogeneously
equivalent} if and only if either both $S$ and $T$ are level or we have $S=T$.
Clearly, being homogeneously equivalent is an equivalence relation. Let
$\overline{J\left(  P\right)  }$ denote the set of equivalence classes of
elements of $J\left(  P\right)  $ modulo this relation. Let $\pi$ denote the
canonical projection $J\left(  P\right)  \rightarrow\overline{J\left(
P\right)  }$. (We distinguish this map $\pi$ from the map $\pi$ defined in
Definition \ref{def.hgeq} by the fact that they act on different objects.)
\end{definition}

The following analogue of Corollary \ref{cor.hgR} is almost trivial:

\begin{corollary}
\label{cor.classical.hgr}Let $n\in\mathbb{N}$. Let $P$ be an $n$-graded poset.
If $S$ and $T$ are two homogeneously equivalent order ideals of $P$, then
$\mathbf{r}\left(  S\right)  $ is homogeneously equivalent to $\mathbf{r}%
\left(  T\right)  $.
\end{corollary}

(An analogue of Corollary \ref{cor.hgRi} exists as well.) We also have the
following analogue of Proposition \ref{prop.reconstruct}:

\begin{proposition}
\label{prop.classical.reconstruct}Let $n\in\mathbb{N}$. Let $P$ be an
$n$-graded poset. Let $S$ and $T$ be two order ideals of $P$ such that
$\left(  \mathbf{w}_{0}\left(  S\right)  ,\mathbf{w}_{1}\left(  S\right)
,...,\mathbf{w}_{n}\left(  S\right)  \right)  =\left(  \mathbf{w}_{0}\left(
T\right)  ,\mathbf{w}_{1}\left(  T\right)  ,...,\mathbf{w}_{n}\left(
T\right)  \right)  $. Also assume that $\pi\left(  S\right)  =\pi\left(
T\right)  $. Then, $S=T$.
\end{proposition}

We can furthermore state analogues of Definitions \ref{def.hgRi} and
\ref{def.hgR}:

\Needspace{15\baselineskip}

\begin{definition}
\label{def.classical.hgri}Let $n\in\mathbb{N}$. Let $P$ be an $n$-graded
poset. Let $i\in\left\{  1,2,...,n\right\}  $. The map $\mathbf{r}%
_{i}:J\left(  P\right)  \rightarrow J\left(  P\right)  $ descends (through the
projection $\pi:J\left(  P\right)  \rightarrow\overline{J\left(  P\right)  }$)
to a map $\overline{J\left(  P\right)  }\rightarrow\overline{J\left(
P\right)  }$. We denote this map $\overline{J\left(  P\right)  }%
\rightarrow\overline{J\left(  P\right)  }$ by $\overline{\mathbf{r}_{i}}$.
Thus, the diagram%
\[
\xymatrixcolsep{5pc}\xymatrix{
J\left(P\right) \ar[r]^{\mathbf r_i} \ar[d]_-{\pi} & J\left(P\right)
\ar[d]^-{\pi} \\
\overline{J\left(P\right)} \ar[r]_{\overline{\mathbf r_i}}
& \overline{J\left(P\right)}
}
\]
is commutative.
\end{definition}

\Needspace{15\baselineskip}

\begin{definition}
\label{def.classical.hgr}Let $n\in\mathbb{N}$. Let $P$ be an $n$-graded poset.
We define the map $\overline{\mathbf{r}}:\overline{J\left(  P\right)
}\rightarrow\overline{J\left(  P\right)  }$ by%
\[
\overline{\mathbf{r}}=\overline{\mathbf{r}_{1}}\circ\overline{\mathbf{r}_{2}%
}\circ...\circ\overline{\mathbf{r}_{n}}.
\]
Then, the diagram%
\begin{equation}
\xymatrixcolsep{5pc}\xymatrix{ J\left(P\right) \ar[r]^{\mathbf r} \ar[d]_-{\pi} & J\left(P\right) \ar[d]^-{\pi} \\ \overline{J\left(P\right)} \ar[r]_{\overline{\mathbf r}} & \overline {J\left(P\right)} }
\label{def.classical.hgr.commute}%
\end{equation}
is commutative. In other words, $\overline{\mathbf{r}}$ is the map
$\overline{J\left(  P\right)  }\rightarrow\overline{J\left(  P\right)  }$ to
which the map $\mathbf{r}:J\left(  P\right)  \rightarrow J\left(  P\right)  $
descends (through the projection $\pi:J\left(  P\right)  \rightarrow
\overline{J\left(  P\right)  }$).
\end{definition}

It might seem that the map $\overline{\mathbf{r}}$ is not worth considering,
since its cycle structure differs from the cycle structure of $\mathbf{r}$
only in the collapsing of an $\left(  n+1\right)  $-cycle (the one formed by
all level order ideals) to a point. However, triviality in combinatorics does
not preclude usefulness, and we will employ the ``projective'' version
$\overline{\mathbf{r}}$ of classical rowmotion as a stirrup in determining the
order of classical rowmotion $\mathbf{r}$ on skeletal posets.

We have the following simple relation between the orders of $\mathbf{r}$ and
$\overline{\mathbf{r}}$:

\begin{proposition}
\label{prop.classical.ord-projord}Let $n\in\mathbb{N}$. Let $P$ be an
$n$-graded poset. Then, $\operatorname*{ord}\mathbf{r}=\operatorname{lcm}%
\left(  n+1,\operatorname*{ord}\overline{\mathbf{r}}\right)  $.
\end{proposition}

Notice that our convention to define $\operatorname{lcm}\left(  n+1,\infty
\right)  $ as $\infty$ is irrelevant for Proposition
\ref{prop.classical.ord-projord}: In fact, in the situation of Proposition
\ref{prop.classical.ord-projord}, both $\operatorname*{ord}\mathbf{r}$ and
$\operatorname*{ord}\overline{\mathbf{r}}$ are clearly (finite) positive
integers\footnote{Indeed, the maps $\mathbf{r}$ and $\overline{\mathbf{r}}$
are permutations of finite sets (namely, of the sets $J\left(  P\right)  $ and
$\overline{J\left(  P\right)  }$) and thus have finite orders.}.

\begin{proof}
[Proof of Proposition \ref{prop.classical.ord-projord} (sketched).]We know
that $\mathbf{r}$ is an invertible map $J\left(  P\right)  \rightarrow
J\left(  P\right)  $, thus a permutation of the finite set $J\left(  P\right)
$. Hence, $\operatorname*{ord}\mathbf{r}$ is the $\operatorname{lcm}$ of the
lengths of the cycles of this permutation $\mathbf{r}$. Similarly,
$\operatorname*{ord}\overline{\mathbf{r}}$ is the $\operatorname{lcm}$ of the
lengths of the cycles of the permutation $\overline{\mathbf{r}}$ of the finite
set $\overline{J\left(  P\right)  }$.

Let $Z_{1}$, $Z_{2}$, $...$, $Z_{k}$ be the cycles of the permutation
$\mathbf{r}$ of $J\left(  P\right)  $. We assume WLOG that $Z_{1}$ is the
cycle consisting of the $n+1$ level order ideals (because we know that they
form a cycle). Thus, $\left\vert Z_{1}\right\vert =n+1$. Since
$\operatorname*{ord}\mathbf{r}$ is the $\operatorname{lcm}$ of the lengths of
the cycles of the permutation $\mathbf{r}$, we have $\operatorname*{ord}%
\mathbf{r}=\operatorname{lcm}\left(  \left\vert Z_{1}\right\vert ,\left\vert
Z_{2}\right\vert ,...,\left\vert Z_{k}\right\vert \right)  $.

Now, let us recall that $\overline{J\left(  P\right)  }$ is the quotient of
$J\left(  P\right)  $ modulo homogeneous equivalence. But homogeneous
equivalence merely identifies the $n+1$ level order ideals, while all other
elements of $J\left(  P\right)  $ are still pairwise non-equivalent. Hence,
the cycles of the permutation $\overline{\mathbf{r}}$ of $\overline{J\left(
P\right)  }$ are $\pi\left(  Z_{1}\right)  $, $\pi\left(  Z_{2}\right)  $,
$...$, $\pi\left(  Z_{k}\right)  $, and while $\pi\left(  Z_{2}\right)  $,
$\pi\left(  Z_{3}\right)  $, $...$, $\pi\left(  Z_{k}\right)  $ are isomorphic
to $Z_{2}$, $Z_{3}$, $...$, $Z_{k}$, respectively, the first cycle $\pi\left(
Z_{1}\right)  $ (being the projection of the cycle of the level order ideals)
now has length $1$. Now, $\operatorname*{ord}\overline{\mathbf{r}}$ is the
$\operatorname{lcm}$ of the lengths of the cycles of the permutation
$\overline{\mathbf{r}}$ of the finite set $\overline{J\left(  P\right)  }$.
Since these cycles are $\pi\left(  Z_{1}\right)  $, $\pi\left(  Z_{2}\right)
$, $...$, $\pi\left(  Z_{k}\right)  $, this yields%
\begin{align*}
\operatorname*{ord}\overline{\mathbf{r}}  &  =\operatorname{lcm}\left(
\left\vert \pi\left(  Z_{1}\right)  \right\vert ,\left\vert \pi\left(
Z_{2}\right)  \right\vert ,...,\left\vert \pi\left(  Z_{k}\right)  \right\vert
\right)  =\operatorname{lcm}\left(  \underbrace{\left\vert \pi\left(
Z_{1}\right)  \right\vert }_{=1},\left\vert \pi\left(  Z_{2}\right)
\right\vert ,\left\vert \pi\left(  Z_{3}\right)  \right\vert ,...,\left\vert
\pi\left(  Z_{k}\right)  \right\vert \right) \\
&  =\operatorname{lcm}\left(  1,\left\vert \pi\left(  Z_{2}\right)
\right\vert ,\left\vert \pi\left(  Z_{3}\right)  \right\vert ,...,\left\vert
\pi\left(  Z_{k}\right)  \right\vert \right)  =\operatorname{lcm}\left(
\left\vert \pi\left(  Z_{2}\right)  \right\vert ,\left\vert \pi\left(
Z_{3}\right)  \right\vert ,...,\left\vert \pi\left(  Z_{k}\right)  \right\vert
\right) \\
&  =\operatorname{lcm}\left(  \left\vert Z_{2}\right\vert ,\left\vert
Z_{3}\right\vert ,...,\left\vert Z_{k}\right\vert \right)
\ \ \ \ \ \ \ \ \ \ \left(
\begin{array}
[c]{c}%
\text{since }\pi\left(  Z_{2}\right)  \text{, }\pi\left(  Z_{3}\right)
\text{, }...\text{, }\pi\left(  Z_{k}\right)  \text{ are}\\
\text{isomorphic to }Z_{2}\text{, }Z_{3}\text{, }...\text{, }Z_{k}\text{,
respectively}%
\end{array}
\right)  .
\end{align*}
Now,%
\begin{align*}
\operatorname*{ord}\mathbf{r}  &  =\operatorname{lcm}\left(  \left\vert
Z_{1}\right\vert ,\left\vert Z_{2}\right\vert ,...,\left\vert Z_{k}\right\vert
\right)  =\operatorname{lcm}\left(  \underbrace{\left\vert Z_{1}\right\vert
}_{=n+1},\underbrace{\operatorname{lcm}\left(  \left\vert Z_{2}\right\vert
,\left\vert Z_{3}\right\vert ,...,\left\vert Z_{k}\right\vert \right)
}_{=\operatorname*{ord}\overline{\mathbf{r}}}\right) \\
&  =\operatorname{lcm}\left(  n+1,\operatorname*{ord}\overline{\mathbf{r}%
}\right)  .
\end{align*}
This proves Proposition \ref{prop.classical.ord-projord}.
\end{proof}

\begin{verlong}
Here is another way to formulate the above proof of Proposition
\ref{prop.classical.ord-projord}. Namely, we shall derive this proposition
from a simple combinatorial lemma:

\begin{lemma}
\label{lem.prop.classical.ord-projord.lem} Let $G$ and $H$ be finite sets, and
let $\pi:G\rightarrow H$ be a surjective map. Let $\sigma$ be a permutation of
$G$, and let $\overline{\sigma}$ be a permutation of $H$ such that
$\overline{\sigma}\circ\pi=\pi\circ\sigma$. Let $L$ be a cycle of the
permutation $\sigma$. Assume that for any two elements $S$ and $T$ of $G$, we
have $\pi\left(  S\right)  =\pi\left(  T\right)  $ if and only if either both
$S$ and $T$ belong to $L$ or we have $S=T$. Then,%
\[
\operatorname*{ord}\sigma=\operatorname{lcm}\left(  \operatorname*{ord}%
\overline{\sigma},\left\vert L\right\vert \right)  .
\]

\end{lemma}

\begin{proof}
[Proof of Lemma \ref{lem.prop.classical.ord-projord.lem} (sketched).]We know
that the order of a permutation is the $\operatorname{lcm}$ of the lengths of
its cycles. Thus,%
\[
\operatorname*{ord}\sigma=\operatorname{lcm}\left\{  \left\vert C\right\vert
\ \mid\ C\text{ is a cycle of }\sigma\right\}  \ \ \ \ \ \ \ \ \ \ \text{and}%
\ \ \ \ \ \ \ \ \ \ \operatorname*{ord}\overline{\sigma}=\operatorname{lcm}%
\left\{  \left\vert D\right\vert \ \mid\ D\text{ is a cycle of }%
\overline{\sigma}\right\}  .
\]

But $\overline{\sigma}\circ\pi=\pi\circ\sigma$. Thus, there is a surjective
map%
\begin{align*}
\left\{  \text{the cycles of }\sigma\right\}   &  \rightarrow\left\{
\text{the cycles of }\overline{\sigma}\right\}  ,\\
C  &  \mapsto\pi\left(  C\right)  .
\end{align*}
\footnote{This map is also bijective in our case, but we do not need this.}
Hence,%
\[
\left\{  \left\vert D\right\vert \ \mid\ D\text{ is a cycle of }%
\overline{\sigma}\right\}  =\left\{  \left\vert \pi\left(  C\right)
\right\vert \ \mid\ C\text{ is a cycle of }\sigma\right\}  .
\]

Now, recall that the map $\pi$ sends all elements of $L$ to one and the same
element of $H$ (since $\pi\left(  S\right)  =\pi\left(  T\right)  $ if and
only if either both $S$ and $T$ belong to $L$ or we have $S=T$). Thus,
$\left\vert \pi\left(  L\right)  \right\vert =1$.

But the map $\pi$ is injective on $G\setminus L$ (since $\pi\left(  S\right)
=\pi\left(  T\right)  $ if and only if either both $S$ and $T$ belong to $L$
or we have $S=T$), and therefore, for every cycle $C$ of $\sigma$ other than
$L$, we have that $\pi\left(  C\right)  $ is a cycle of $\overline{\sigma}$
with the same length as $C$. Thus,%
\begin{equation}
\left\vert \pi\left(  C\right)  \right\vert =\left\vert C\right\vert
\ \ \ \ \ \ \ \ \ \ \text{for every cycle }\sigma\text{ of }C\text{ other than
}L. \label{pf.lem.prop.classical.ord-projord.lem.1}%
\end{equation}
Hence,%
\begin{align*}
&  \left\{  \left\vert \pi\left(  C\right)  \right\vert \ \mid\ C\text{ is a
cycle of }\sigma\right\} \\
&  =\left\{  \underbrace{\left\vert \pi\left(  L\right)  \right\vert }%
_{=1}\right\}  \cup\left\{  \underbrace{\left\vert \pi\left(  C\right)
\right\vert }_{\substack{=\left\vert C\right\vert \\\text{(by
(\ref{pf.lem.prop.classical.ord-projord.lem.1}))}}}\ \mid\ C\text{ is a cycle
of }\sigma\text{ other than }L\right\} \\
&  =\left\{  1\right\}  \cup\left\{  \left\vert C\right\vert \ \mid\ C\text{
is a cycle of }\sigma\text{ other than }L\right\}  .
\end{align*}
Now,%
\begin{align}
\operatorname*{ord}\overline{\sigma}  &  =\operatorname{lcm}%
\underbrace{\left\{  \left\vert D\right\vert \ \mid\ D\text{ is a cycle of
}\overline{\sigma}\right\}  }_{=\left\{  1\right\}  \cup\left\{  \left\vert
C\right\vert \ \mid\ C\text{ is a cycle of }\sigma\text{ other than
}L\right\}  }\nonumber\\
&  =\operatorname{lcm}\left(  \left\{  1\right\}  \cup\left\{  \left\vert
C\right\vert \ \mid\ C\text{ is a cycle of }\sigma\text{ other than
}L\right\}  \right) \nonumber\\
&  =\operatorname{lcm}\left(  1,\operatorname{lcm}\left\{  \left\vert
C\right\vert \ \mid\ C\text{ is a cycle of }\sigma\text{ other than
}L\right\}  \right) \nonumber\\
&  =\operatorname{lcm}\left\{  \left\vert C\right\vert \ \mid\ C\text{ is a
cycle of }\sigma\text{ other than }L\right\}  .
\label{pf.lem.prop.classical.ord-projord.lem.2}%
\end{align}
Now,%
\begin{align*}
\operatorname*{ord}\sigma &  =\operatorname{lcm}\underbrace{\left\{
\left\vert C\right\vert \ \mid\ C\text{ is a cycle of }\sigma\right\}
}_{\substack{=\left\{  \left\vert L\right\vert \right\}  \cup\left\{
\left\vert C\right\vert \ \mid\ C\text{ is a cycle of }\sigma\text{ other than
}L\right\}  \\\text{(since }L\text{ is a cycle of }\sigma\text{)}}}\\
&  =\operatorname{lcm}\left(  \left\{  \left\vert L\right\vert \right\}
\cup\left\{  \left\vert C\right\vert \ \mid\ C\text{ is a cycle of }%
\sigma\text{ other than }L\right\}  \right) \\
&  =\operatorname{lcm}\left(  \left\vert L\right\vert
,\underbrace{\operatorname{lcm}\left\{  \left\vert C\right\vert \ \mid
\ C\text{ is a cycle of }\sigma\text{ other than }L\right\}  }%
_{\substack{=\operatorname*{ord}\overline{\sigma}\\\text{(by
(\ref{pf.lem.prop.classical.ord-projord.lem.2}))}}}\right)
=\operatorname{lcm}\left(  \left\vert L\right\vert ,\operatorname*{ord}%
\overline{\sigma}\right) \\
&= \operatorname{lcm}\left(  \operatorname*{ord}%
\overline{\sigma},\left\vert L\right\vert \right) .
\end{align*}
This proves Lemma \ref{lem.prop.classical.ord-projord.lem}.
\end{proof}

\begin{proof}
[Second proof of Proposition \ref{prop.classical.ord-projord} (sketched).]The
map $\pi:J\left(  P\right)  \rightarrow\overline{J\left(  P\right)  }$ is
surjective. The permutations $\mathbf{r}$ and $\overline{\mathbf{r}}$ of
$J\left(  P\right)  $ and $\overline{J\left(  P\right)  }$ satisfy
$\overline{\mathbf{r}}\circ\pi=\pi\circ\mathbf{r}$ (since the diagram
(\ref{def.classical.hgr.commute}) commutes). Let $L$ denote the set of all
level order ideals in $J\left(  P\right)  $. Then, $\left\vert L\right\vert
=n+1$ (by Proposition \ref{prop.classical.level}). The set $L$ is an orbit
under $\mathbf{r}$ (again by Proposition \ref{prop.classical.level}), that is,
a cycle of $\mathbf{r}$. For any two elements $S$ and $T$ of $J\left(
P\right)  $, we have $\pi\left(  S\right)  =\pi\left(  T\right)  $ if and only
if $S$ and $T$ are homogeneously equivalent (because $\pi$ is the canonical
projection $J\left(  P\right)  \rightarrow\overline{J\left(  P\right)  }$),
which is equivalent to saying that either both $S$ and $T$ are level (that is,
both $S$ and $T$ belong to $L$) or we have $S=T$. Hence, all the assumptions
of Lemma \ref{lem.prop.classical.ord-projord.lem} are satisfied for
$G=J\left(  P\right)  $, $H=\overline{J\left(  P\right)  }$, $\sigma
=\mathbf{r}$ and $\overline{\sigma}=\overline{\mathbf{r}}$. Therefore,
applying Lemma \ref{lem.prop.classical.ord-projord.lem} to this setting, we
obtain $\operatorname*{ord}\mathbf{r}=\operatorname{lcm}\left(
\underbrace{\left\vert L\right\vert }_{=n+1},\operatorname*{ord}%
\overline{\mathbf{r}}\right)  =\operatorname*{ord}\mathbf{r}%
=\operatorname{lcm}\left(  n+1,\operatorname*{ord}\overline{\mathbf{r}%
}\right)  $.
\end{proof}
\end{verlong}

Our goal is to make a statement about the order of classical rowmotion on
skeletal posets. Of course, the finiteness of these orders is obvious in this
case, because $J\left(  P\right)  $ is a finite set. However, we can make
stronger claims:

\begin{proposition}
\label{prop.classical.skeletal.ords}Let $P$ be a skeletal poset. Let
$\mathbb{K}$ be a field. Then, $\operatorname*{ord}\left(  R_{P}\right)
=\operatorname*{ord}\left(  \mathbf{r}_{P}\right)  $ and $\operatorname*{ord}%
\left(  \overline{R}_{P}\right)  =\operatorname*{ord}\left(  \overline
{\mathbf{r}}_{P}\right)  $.
\end{proposition}

Here, we are using the following convention:

\begin{definition}
\label{def.classical.conv.rP}Let $P$ be a finite poset. We denote the maps
$\mathbf{r}$ and $\overline{\mathbf{r}}$ by $\mathbf{r}_{P}$ and
$\overline{\mathbf{r}}_{P}$, respectively, so as to make their dependence on
$P$ explicit.
\end{definition}

Proposition \ref{prop.classical.skeletal.ords} yields (in particular) that the
order of classical rowmotion coincides with the order of birational rowmotion
(whatever the base field) for skeletal posets. This was conjectured by James
Propp (private communication) for the case of $P$ a tree. We are going to
prove Proposition \ref{prop.classical.skeletal.ords} by exhibiting further
analogies between classical and birational rowmotion. First of all, the
following proposition is just as trivial as its birational counterpart
(Proposition \ref{prop.PQ.ord}):

\begin{proposition}
\label{prop.classical.PQ.ord}Let $n\in\mathbb{N}$. Let $P$ and $Q$ be two
$n$-graded posets. Then, $\operatorname*{ord}\left(  \mathbf{r}_{PQ}\right)
=\operatorname{lcm}\left(  \operatorname*{ord}\left(  \mathbf{r}_{P}\right)
,\operatorname*{ord}\left(  \mathbf{r}_{Q}\right)  \right)  $.
\end{proposition}

We can show a simple counterpart of this proposition for $\operatorname*{ord}%
\left(  \overline{\mathbf{r}}_{PQ}\right)  $ (but still with
$\operatorname*{ord}\left(  \mathbf{r}_{P}\right)  $ and $\operatorname*{ord}%
\left(  \mathbf{r}_{Q}\right)  $ on the right hand side!):

\begin{proposition}
\label{prop.classical.PQ.ord2}Let $n\in\mathbb{N}$. Let $P$ and $Q$ be
two\textbf{ }$n$-graded posets. Then, $\operatorname*{ord}\left(
\overline{\mathbf{r}}_{PQ}\right)  =\operatorname{lcm}\left(
\operatorname*{ord}\left(  \mathbf{r}_{P}\right)  ,\operatorname*{ord}\left(
\mathbf{r}_{Q}\right)  \right)  $.
\end{proposition}

\begin{proof}
[Proof of Proposition \ref{prop.classical.PQ.ord2} (sketched).]WLOG, assume
that $n\neq0$ (else, the statement is trivial). Hence, $P$ and $Q$ are nonempty.

Consider the order ideal $P$ of $PQ$. Then, one can easily see (by induction)
that every $i\in\left\{  0,1,...,n+1\right\}  $ satisfies%
\begin{align*}
&  \mathbf{r}_{PQ}^{i}\left(  P\right) \\
&  =\left(  \left\{
\begin{array}
[c]{l}%
P_{1}\cup P_{2}\cup...\cup P_{i-1},\ \ \ \ \ \ \ \ \ \ \text{if }i>0;\\
P,\ \ \ \ \ \ \ \ \ \ \text{if }i=0
\end{array}
\right.  \right)  \cup\left(  \left\{
\begin{array}
[c]{l}%
Q_{1}\cup Q_{2}\cup...\cup Q_{i},\ \ \ \ \ \ \ \ \ \ \text{if }i\leq n;\\
\varnothing,\ \ \ \ \ \ \ \ \ \ \text{if }i=n+1
\end{array}
\right.  \right)  .
\end{align*}
From this, it follows that the smallest positive integer $k$ satisfying
$\mathbf{r}_{PQ}^{k}\left(  P\right)  =P$ is $n+1$. Since $P$ is not level
(as an order ideal of $PQ$),
this does not change under applying $\pi$; that is, the smallest positive
integer $k$ satisfying $\overline{\mathbf{r}}_{PQ}^{k}\left(  P\right)  =P$ is
still $n+1$. Hence, $n+1\mid\operatorname*{ord}\left(  \overline{\mathbf{r}%
}_{PQ}\right)  $. But Proposition \ref{prop.classical.ord-projord} (applied to
$PQ$ instead of $P$) yields%
\[
\operatorname*{ord}\left(  \mathbf{r}_{PQ}\right)  =\operatorname{lcm}\left(
n+1,\operatorname*{ord}\left(  \overline{\mathbf{r}}_{PQ}\right)  \right)
=\operatorname*{ord}\left(  \overline{\mathbf{r}}_{PQ}\right)
\]
(since $n+1\mid\operatorname*{ord}\left(  \overline{\mathbf{r}}_{PQ}\right)
$), so that%
\[
\operatorname*{ord}\left(  \overline{\mathbf{r}}_{PQ}\right)
=\operatorname*{ord}\left(  \mathbf{r}_{PQ}\right)  =\operatorname{lcm}\left(
\operatorname*{ord}\left(  \mathbf{r}_{P}\right)  ,\operatorname*{ord}\left(
\mathbf{r}_{Q}\right)  \right)
\]
(by Proposition \ref{prop.classical.PQ.ord}). This proves Proposition
\ref{prop.classical.PQ.ord2}.
\end{proof}

More interesting is the analogue of Proposition \ref{prop.Bk.ord}:

\begin{proposition}
\label{prop.classical.Bk.ord}Let $n\in\mathbb{N}$. Let $P$ be an $n$-graded poset.

\textbf{(a)} We have $\operatorname*{ord}\left(  \overline{\mathbf{r}}%
_{B_{1}P}\right)  =\operatorname*{ord}\left(  \overline{\mathbf{r}}%
_{P}\right)  $.

\textbf{(b)} For every integer $k>1$, we have $\operatorname*{ord}\left(
\overline{\mathbf{r}}_{B_{k}P}\right)  =\operatorname{lcm}\left(
2,\operatorname*{ord}\left(  \overline{\mathbf{r}}_{P}\right)  \right)  $.
\end{proposition}

\begin{proof}
[Proof of Proposition \ref{prop.classical.Bk.ord} (sketched).]We will be
proving parts \textbf{(a)} and \textbf{(b)} together. Let $k$ be a positive
integer (this has to be $1$ for proving part \textbf{(a)}). We need to prove
that%
\begin{equation}
\operatorname*{ord}\left(  \overline{\mathbf{r}}_{B_{k}P}\right)  =\left\{
\begin{array}
[c]{l}%
\operatorname{lcm}\left(  2,\operatorname*{ord}\left(  \overline{\mathbf{r}%
}_{P}\right)  \right)  ,\ \ \ \ \ \ \ \ \ \ \text{if }k>1;\\
\operatorname*{ord}\left(  \overline{\mathbf{r}}_{P}\right)
,\ \ \ \ \ \ \ \ \ \ \text{if }k=1
\end{array}
\right.  . \label{pf.classical.Bk.ord.1}%
\end{equation}
Proving this clearly will prove both parts \textbf{(a)} and \textbf{(b)} of
Proposition \ref{prop.classical.Bk.ord}.

Notice that $B_{k}P$ is an $\left(  n+1\right)  $-graded poset. For every
$\ell\in\left\{  1,2,...,n+1\right\}  $, let $\left(  B_{k}P\right)  _{\ell}$
be the subset $\deg^{-1}\left(  \left\{  \ell\right\}  \right)  $ of $B_{k}P$.
Thus, $\left(  B_{k}P\right)  _{\ell}=\left\{  v\in B_{k}P\ \mid\ \deg
v=\ell\right\}  $. In particular, $\left(  B_{k}P\right)  _{1}$ is the set of
all minimal elements of $B_{k}P$, so that $\left(  B_{k}P\right)  _{1}$ is an
antichain of size $k$ (by the construction of $B_{k}P$). We also have
\begin{align}
v>w
\ \ \ \ \ \ \ \ \ \ \text{for every }
w\in\left(  B_{k}P\right)  _{1}\text{ and }v\in P .
\label{pf.classical.Bk.ord.wlessv}
\end{align}

For every graded poset $Q$, the map $\overline{\mathbf{r}}_{Q}$ is an
invertible map $\overline{J\left(  Q\right)  }\rightarrow\overline{J\left(
Q\right)  }$, that is, a permutation of the finite set $\overline{J\left(
Q\right)  }$. Hence, its order $\operatorname*{ord}\left(  \overline
{\mathbf{r}}_{Q}\right)  $ is the $\operatorname{lcm}$ of the lengths of the
orbits of this map $\overline{\mathbf{r}}_{Q}$. We are going to compare the
orbits of the maps $\overline{\mathbf{r}}_{B_{k}P}$ and $\overline{\mathbf{r}%
}_{P}$.

Define a map $\phi:J\left(  P\right)  \rightarrow J\left(  B_{k}P\right)  $ by%
\[
\phi\left(  S\right)  =\left(  B_{k}P\right)  _{1}\cup
S\ \ \ \ \ \ \ \ \ \ \text{for every }S\in J\left(  P\right)  .
\]
It is easy to see that this map $\phi$ is well-defined (that is, $\left(
B_{k}P\right)  _{1}\cup S$ is an order ideal of $B_{k}P$ for every $S\in
J\left(  P\right)  $), and that it sends level order ideals of $P$ to level
order ideals of $B_{k}P$. Hence, it preserves homogeneous equivalence, so that
it induces a map $\overline{J\left(  P\right)  }\rightarrow\overline{J\left(
B_{k}P\right)  }$. Denote this map $\overline{J\left(  P\right)  }%
\rightarrow\overline{J\left(  B_{k}P\right)  }$ by $\overline{\phi}$. Thus,
$\overline{\phi}\circ\pi=\pi\circ\phi$.

It is moreover easy to see that $\overline{\mathbf{r}}_{B_{k}P}\circ
\overline{\phi}=\overline{\phi}\circ\overline{\mathbf{r}}_{P}$%
\ \ \ \ \footnote{\textit{Proof.} In order to prove this, it is enough to show
that for every $S\in J\left(  P\right)  $, the order ideals $\left(
\mathbf{r}_{B_{k}P}\circ\phi\right)  \left(  S\right)  $ and $\left(
\phi\circ\mathbf{r}_{P}\right)  \left(  S\right)  $ are homogeneously
equivalent. This is clear in the case when $S$ is level (because both $\left(
\mathbf{r}_{B_{k}P}\circ\phi\right)  \left(  S\right)  $ and $\left(
\phi\circ\mathbf{r}_{P}\right)  \left(  S\right)  $ are level in this case),
so let us WLOG assume that $S$ is not level. Then, we can actually show that
$\left(  \mathbf{r}_{B_{k}P}\circ\phi\right)  \left(  S\right)  $ and $\left(
\phi\circ\mathbf{r}_{P}\right)  \left(  S\right)  $ are identical. Indeed, it
is easy to see that:
\par
\begin{itemize}
\item for every $T\in J\left(  P\right)  $ and every $v\in P$, we have
$\left(  \mathbf{t}_{v}\circ\phi\right)  \left(  T\right)  =\left(  \phi
\circ\mathbf{t}_{v}\right)  \left(  T\right)  $;
\par
\item for every nonempty $T\in J\left(  P\right)  $ and every $w\in\left(
B_{k}P\right)  _{1}$, we have $\left(  \mathbf{t}_{w}\circ\phi\right)  \left(
T\right)  =\phi\left(  T\right)  $.
\end{itemize}
\par
Using these facts, and the definition of classical rowmotion as a composition
of classical toggle maps $\mathbf{t}_{v}$, we can then easily see that
$\left(  \mathbf{r}_{B_{k}P}\circ\phi\right)  \left(  S\right)  =\left(
\phi\circ\mathbf{r}_{P}\right)  \left(  S\right)  $. This completes the proof
of $\overline{\mathbf{r}}_{B_{k}P}\circ\overline{\phi}=\overline{\phi}%
\circ\overline{\mathbf{r}}_{P}$.}. Hence, the subset $\overline{\phi}\left(
\overline{J\left(  P\right)  }\right)  $ is closed under application of the
map $\overline{\mathbf{r}}_{B_{k}P}$.

The map $\overline{\phi}$ also is injective (this is very easy to see again,
since the only order ideals of $P$ which are mapped to level order ideals by
$\phi$ are themselves level). Thus, $\operatorname*{ord}\left(  \overline
{\mathbf{r}}_{B_{k}P}\mid_{\overline{\phi}\left(  \overline{J\left(  P\right)
}\right)  }\right)  =\operatorname*{ord}\left(  \overline{\mathbf{r}}%
_{P}\right)  $ (because the injectivity of $\overline{\phi}$ allows us to
identify $\overline{J\left(  P\right)  }$ with $\overline{\phi}\left(
\overline{J\left(  P\right)  }\right)  $ along the map $\overline{\phi}$, and
then the equality $\overline{\mathbf{r}}_{B_{k}P}\circ\overline{\phi
}=\overline{\phi}\circ\overline{\mathbf{r}}_{P}$ rewrites as $\overline
{\mathbf{r}}_{B_{k}P}\mid_{\overline{\phi}\left(  \overline{J\left(  P\right)
}\right)  }=\overline{\mathbf{r}}_{P}$, so that $\operatorname*{ord}\left(
\overline{\mathbf{r}}_{B_{k}P}\mid_{\overline{\phi}\left(  \overline{J\left(
P\right)  }\right)  }\right)  =\operatorname*{ord}\left(  \overline
{\mathbf{r}}_{P}\right)  $).

Let $H$ be the set of all nonempty proper subsets of $\left(  B_{k}P\right)
_{1}$. It is clear that $H\subseteq J\left(  B_{k}P\right)  $. Notice that
$H=\varnothing$ if $k=1$. Every $T\in H$ satisfies%
\[
\mathbf{r}_{B_{k}P}\left(  T\right)  =\left(  B_{k}P\right)  _{1}\setminus T
\]
(this is easy to see from any definition of classical rowmotion, or from
Proposition \ref{prop.classical.r.implicit}). Hence, the set $H$ is closed
under application of the map $\mathbf{r}_{B_{k}P}$, and this map
$\mathbf{r}_{B_{k}P}$ maps every element of $H$ to its complement in $\left(
B_{k}P\right)  _{1}$. In particular, this shows that
\[
\operatorname*{ord}%
\left(  \mathbf{r}_{B_{k}P}\mid_{H}\right)  =\left\{
\begin{array}
[c]{l}%
2,\ \ \ \ \ \ \ \ \ \ \text{if }k>1;\\
1,\ \ \ \ \ \ \ \ \ \ \text{if }k=1
\end{array}
\right.  .
\]

We now use the map $\pi$ to identify the set $H$ with its projection
$\pi\left(  H\right)  $ under $\pi$ (this is allowed because $\pi$ is
injective on $H$). This identification entails $\left.  \overline{\mathbf{r}%
}_{B_{k}P}\mid_{H}\right.  =\mathbf{r}_{B_{k}P}\mid_{H}$. In particular, the
set $H$ is closed under application of the map $\overline{\mathbf{r}}_{B_{k}%
P}$.

However, it is easy to see that%
\begin{equation}
J\left(  B_{k}P\right)  =\left\{  \varnothing\right\}  \cup H\cup\phi\left(
J\left(  P\right)  \right)  .\label{pf.prop.classical.Bk.ord.union1}%
\end{equation}

[\textit{Proof of (\ref{pf.prop.classical.Bk.ord.union1}):} Clearly, the three
sets $\left\{  \varnothing\right\}  $, $H$ and $\phi\left(  J\left(  P\right)
\right)  $ are subsets of $J\left(  B_{k}P\right)  $. Thus, their union
$\left\{  \varnothing\right\}  \cup H\cup\phi\left(  J\left(  P\right)
\right)  $ is a subset of $J\left(  B_{k}P\right)  $ as well. It thus remains
to prove that it is not a proper subset.

Assume the contrary. Thus, there exists an order ideal $T\in J\left(
B_{k}P\right)  $ that is not contained in the union $\left\{  \varnothing
\right\}  \cup H\cup\phi\left(  J\left(  P\right)  \right)  $. Consider this
$T$.

We note that every subset of $\left(  B_{k}P\right)  _{1}$ is contained in the
union $\left\{  \varnothing\right\}  \cup H\cup\phi\left(  J\left(  P\right)
\right)  $ (indeed, the empty subset $\varnothing$ is contained in $\left\{
\varnothing\right\}  $; the full subset $\left(  B_{k}P\right)  _{1}$ is
contained in $\phi\left(  J\left(  P\right)  \right)  $ because it can be
written as $\phi\left(  \varnothing\right)  $ for $\varnothing\in J\left(
P\right)  $; and all remaining subsets of $\left(  B_{k}P\right)  _{1}$ are
contained in $H$ by the definition of $H$). Thus, $T$ cannot be a subset of
$\left(  B_{k}P\right)  _{1}$ (since $T$ is not contained in this union).
Hence, $T$ contains some element $v\in B_{k}P$ of degree $>1$. This element
$v$ belongs to $P$, and thus is larger than every element of $\left(
B_{k}P\right)  _{1}$ (by the definition of (\ref{pf.classical.Bk.ord.wlessv}%
)). Therefore, an order ideal of $B_{k}P$ that contains $v$ must necessarily
contain $\left(  B_{k}P\right)  _{1}$ as a subset. Thus, $T$ must contain
$\left(  B_{k}P\right)  _{1}$ as a subset (since $T$ contains $v$). This shows
that $T=\left(  B_{k}P\right)  _{1}\cup\left(  T\setminus\left(
B_{k}P\right)  _{1}\right)  $. Moreover, since $T$ is an order ideal of
$B_{k}P$, we can easily see that $T\setminus\left(  B_{k}P\right)  _{1}$ is an
order ideal of $P$, and satisfies $\phi\left(  T\setminus\left(
B_{k}P\right)  _{1}\right)  =\left(  B_{k}P\right)  _{1}\cup\left(
T\setminus\left(  B_{k}P\right)  _{1}\right)  =T$. Hence, $T=\phi\left(
T^{\prime}\right)  $ for some $T^{\prime}\in J\left(  P\right)  $ (namely, for
$T^{\prime}=T\setminus\left(  B_{k}P\right)  _{1}$). Therefore, $T\in
\phi\left(  J\left(  P\right)  \right)  \subseteq\left\{  \varnothing\right\}
\cup H\cup\phi\left(  J\left(  P\right)  \right)  $. This contradicts the fact
that $T$ is not contained in the union $\left\{  \varnothing\right\}  \cup
H\cup\phi\left(  J\left(  P\right)  \right)  $. This contradiction shows that
our assumption was false. This proves (\ref{pf.prop.classical.Bk.ord.union1}).]

Now, applying the projection $\pi:J\left(  B_{k}P\right)  \rightarrow
\overline{J\left(  B_{k}P\right)  }$ to both sides of
(\ref{pf.prop.classical.Bk.ord.union1}), we obtain%
\begin{align*}
\overline{J\left(  B_{k}P\right)  }  & =\pi\left(  \left\{  \varnothing
\right\}  \cup H\cup\phi\left(  J\left(  P\right)  \right)  \right)
=\underbrace{\pi\left(  \left\{  \varnothing\right\}  \right)  }%
_{\substack{=\left\{  \pi\left(  \varnothing\right)  \right\}  \subseteq
\pi\left(  \phi\left(  J\left(  P\right)  \right)  \right)  \\\text{(since
}\pi\left(  \varnothing\right)  \in\pi\left(  \phi\left(  J\left(  P\right)
\right)  \right)  \\\text{(because }\varnothing=\phi\left(  \varnothing
\right)  \in\phi\left(  J\left(  P\right)  \right)  \text{))}}}\cup
\underbrace{\pi\left(  H\right)  }_{=H}\cup\pi\left(  \phi\left(  J\left(
P\right)  \right)  \right)  \\
& \subseteq\pi\left(  \phi\left(  J\left(  P\right)  \right)  \right)  \cup
H\cup\pi\left(  \phi\left(  J\left(  P\right)  \right)  \right)
=H\cup\underbrace{\pi\left(  \phi\left(  J\left(  P\right)  \right)  \right)
}_{\substack{=\overline{\phi}\left(  \overline{J\left(  P\right)  }\right)
\\\text{(by the definition of }\overline{\phi}\text{)}}}=H\cup\overline{\phi
}\left(  \overline{J\left(  P\right)  }\right)  .
\end{align*}

In other words, the set $\overline{J\left(  B_{k}P\right)  }$ is the union of
the two subsets $H$ and $\overline{\phi}\left(  \overline{J\left(  P\right)
}\right)  $. Moreover, these two subsets are disjoint\footnote{\textit{Proof.}
Let us compare their elements:
\par
\begin{itemize}
\item Each element of $H$ is a nonempty proper subset of $\left(
B_{k}P\right)  _{1}$.
\par
\item Each element of $\overline{\phi}\left(  \overline{J\left(  P\right)
}\right)  =\pi\left(  \phi\left(  J\left(  P\right)  \right)  \right)  $ comes
(via the map $\pi$) from a set of the form $\phi\left(  S\right)  =\left(
B_{k}P\right)  _{1}\cup S$ with $S\in J\left(  P\right)  $. Any such set
clearly contains $\left(  B_{k}P\right)  _{1}$ as a subset.
\end{itemize}
\par
Thus, an element of $H$ could equal an element of $\overline{\phi}\left(
\overline{J\left(  P\right)  }\right)  $ only if the map $\pi$ would equate
these two. However, the map $\pi$ equates level order ideals only; thus, $\pi$
does not equate $H$ to any other ideal (since $H$ is not level). Hence, an
element of $H$ cannot equal an element of $\overline{\phi}\left(
\overline{J\left(  P\right)  }\right)  $. In other words, the sets $H$ and
$\overline{\phi}\left(  \overline{J\left(  P\right)  }\right)  $ are
disjoint.}, and each of them is closed under application of the map
$\overline{\mathbf{r}}_{B_{k}P}$. Hence,%
\begin{align*}
\operatorname*{ord}\left(  \overline{\mathbf{r}}_{B_{k}P}\right)   &
=\operatorname{lcm}\left(  \operatorname*{ord}\left(  \underbrace{\overline
{\mathbf{r}}_{B_{k}P}\mid_{H}}_{=\mathbf{r}_{B_{k}P}\mid_{H}}\right)
,\operatorname*{ord}\left(  \overline{\mathbf{r}}_{B_{k}P}\mid_{\overline
{\phi}\left(  \overline{J\left(  P\right)  }\right)  }\right)  \right)  \\
&  =\operatorname{lcm}\left(  \underbrace{\operatorname*{ord}\left(
\mathbf{r}_{B_{k}P}\mid_{H}\right)  }_{=\left\{
\begin{array}
[c]{l}%
2,\ \ \ \ \ \ \ \ \ \ \text{if }k>1;\\
1,\ \ \ \ \ \ \ \ \ \ \text{if }k=1
\end{array}
\right.  },\underbrace{\operatorname*{ord}\left(  \overline{\mathbf{r}}%
_{B_{k}P}\mid_{\overline{\phi}\left(  \overline{J\left(  P\right)  }\right)
}\right)  }_{=\operatorname*{ord}\left(  \overline{\mathbf{r}}_{P}\right)
}\right)  \\
&  =\operatorname{lcm}\left(  \left\{
\begin{array}
[c]{l}%
2,\ \ \ \ \ \ \ \ \ \ \text{if }k>1;\\
1,\ \ \ \ \ \ \ \ \ \ \text{if }k=1
\end{array}
\right.  ,\operatorname*{ord}\left(  \overline{\mathbf{r}}_{P}\right)
\right)  \\
&  =\left\{
\begin{array}
[c]{l}%
\operatorname{lcm}\left(  2,\operatorname*{ord}\left(  \overline{\mathbf{r}%
}_{P}\right)  \right)  ,\ \ \ \ \ \ \ \ \ \ \text{if }k>1;\\
\operatorname*{ord}\left(  \overline{\mathbf{r}}_{P}\right)
,\ \ \ \ \ \ \ \ \ \ \text{if }k=1
\end{array}
\right.  .
\end{align*}
This proves (\ref{pf.classical.Bk.ord.1}). Thus, the proof of Proposition
\ref{prop.classical.Bk.ord} is complete.
\end{proof}

We can also formulate an analogue of Proposition \ref{prop.B'k.ord}:

\begin{proposition}
\label{prop.classical.B'k.ord}Let $n\in\mathbb{N}$. Let $P$ be an $n$-graded poset.

\textbf{(a)} We have $\operatorname*{ord}\left(  \overline{\mathbf{r}}%
_{B_{1}^{\prime}P}\right)  =\operatorname*{ord}\left(  \overline{\mathbf{r}%
}_{P}\right)  $.

\textbf{(b)} For every integer $k>1$, we have $\operatorname*{ord}\left(
\overline{\mathbf{r}}_{B_{k}^{\prime}P}\right)  =\operatorname{lcm}\left(
2,\operatorname*{ord}\left(  \overline{\mathbf{r}}_{P}\right)  \right)  $.
\end{proposition}

\begin{proof}
[\nopunct]The proof of this is fairly similar to that of Proposition
\ref{prop.classical.Bk.ord}.
\end{proof}

We can now prove Proposition \ref{prop.classical.skeletal.ords}:

\begin{proof}
[Proof of Proposition \ref{prop.classical.skeletal.ords} (sketched).]For any
skeletal poset $T$, we can compute $\operatorname*{ord}\left(  R_{T}\right)  $
and $\operatorname*{ord}\left(  \overline{R}_{T}\right)  $ inductively using
Proposition \ref{prop.PQ.ord}, Proposition \ref{prop.PQ.ord2}, Proposition
\ref{prop.Bk.ord} and Proposition \ref{prop.B'k.ord} (and the fact that
$\operatorname*{ord}\left(  R_{\varnothing}\right)  =1$ and
$\operatorname*{ord}\left(  \overline{R}_{\varnothing}\right)  =1$). More precisely:

\begin{itemize}
\item If $T$ is the empty poset $\varnothing$, then $\operatorname*{ord}%
\left(  R_{T}\right)  =\operatorname*{ord}\left(  R_{\varnothing}\right)  =1$
and $\operatorname*{ord}\left(  \overline{R}_{T}\right)  =\operatorname*{ord}%
\left(  \overline{R}_{\varnothing}\right)  =1$.

\item If $T$ has the form $B_{k}P$ for some $n$-graded skeletal poset $P$ and
some positive integer $k$, then Proposition \ref{prop.Bk.ord} yields
\[
\operatorname*{ord}\left(  \overline{R}_{T}\right)  =\operatorname*{ord}%
\left(  \overline{R}_{B_{k}P}\right)  =\left\{
\begin{array}
[c]{l}%
\operatorname{lcm}\left(  2,\operatorname*{ord}\left(  \overline{R}%
_{P}\right)  \right)  ,\ \ \ \ \ \ \ \ \ \ \text{if }k>1;\\
\operatorname*{ord}\left(  \overline{R}_{P}\right)
,\ \ \ \ \ \ \ \ \ \ \text{if }k=1
\end{array}
\right.  ,
\]
and Proposition \ref{prop.ord-projord} yields $\operatorname*{ord}\left(
R_{T}\right)  =\operatorname{lcm}\left(  n+1,\operatorname*{ord}\left(
\overline{R}_{T}\right)  \right)  $.

\item Analogously one can compute $\operatorname*{ord}\left(  R_{T}\right)  $
and $\operatorname*{ord}\left(  \overline{R}_{T}\right)  $ if $T$ has the form
$B_{k}^{\prime}P$.

\item If $T$ has the form $PQ$ for two WLOG nonempty $n$-graded skeletal
posets $P$ and $Q$, then Proposition \ref{prop.PQ.ord} yields
$\operatorname*{ord}\left(  R_{PQ}\right)  =\operatorname{lcm}\left(
\operatorname*{ord}\left(  R_{P}\right)  ,\operatorname*{ord}\left(
R_{Q}\right)  \right)  $, and Proposition \ref{prop.PQ.ord2} yields
$\operatorname*{ord}\left(  \overline{R}_{PQ}\right)  =\operatorname{lcm}%
\left(  \operatorname*{ord}\left(  R_{P}\right)  ,\operatorname*{ord}\left(
R_{Q}\right)  \right)  $.
\end{itemize}

This gives an algorithm for inductively computing $\operatorname*{ord}\left(
R_{T}\right)  $ and $\operatorname*{ord}\left(  \overline{R}_{T}\right)  $ for
a skeletal poset $T$. Using Proposition \ref{prop.classical.PQ.ord},
Proposition \ref{prop.classical.PQ.ord2}, Proposition
\ref{prop.classical.Bk.ord} and Proposition \ref{prop.classical.B'k.ord} (and
the fact that $\operatorname*{ord}\left(  \mathbf{r}_{\varnothing}\right)  =1$
and $\operatorname*{ord}\left(  \overline{\mathbf{r}}_{\varnothing}\right)
=1$) instead, we could similarly obtain an algorithm for inductively computing
$\operatorname*{ord}\left(  \mathbf{r}_{T}\right)  $ and $\operatorname*{ord}%
\left(  \overline{\mathbf{r}}_{T}\right)  $ for a skeletal poset $T$. And
these two algorithms are \textbf{the same}, because of the direct analogy
between the propositions that are used in the first algorithm and those used
in the second one. Therefore, $\operatorname*{ord}\left(  R_{P}\right)
=\operatorname*{ord}\left(  \mathbf{r}_{P}\right)  $ and $\operatorname*{ord}%
\left(  \overline{R}_{P}\right)  =\operatorname*{ord}\left(  \overline
{\mathbf{r}}_{P}\right)  $. This proves Proposition
\ref{prop.classical.skeletal.ords}.
\end{proof}

Proposition \ref{prop.classical.skeletal.ords} does not generalize to
arbitrary graded posets. Counterexamples to such a generalization can be found
in Section \ref{sect.negres}.

Finally, in analogy to Corollary \ref{cor.for.ord}, we can now show:

\begin{corollary}
\label{cor.classical.for.ord}Let $n\in\mathbb{N}$. Let $P$ be an $n$-graded
poset. Assume that $P$ is a rooted forest (made into a poset by having every
node smaller than its children).

\textbf{(a)} Then, $\operatorname*{ord}\left(  \mathbf{r}_{P}\right)
\mid\operatorname{lcm}\left(  1,2,...,n+1\right)  $.

\textbf{(b)} Moreover, if $P$ is a tree, then $\operatorname*{ord}\left(
\overline{\mathbf{r}}_{P}\right)  \mid\operatorname{lcm}\left(
1,2,...,n\right)  $.
\end{corollary}

Corollary \ref{cor.classical.for.ord} is also valid if we replace
\textquotedblleft every node smaller than its children\textquotedblright\ by
\textquotedblleft every node larger than its children\textquotedblright, and
the proof is exactly analogous.

Let us notice that the algorithm described in the proof of Proposition
\ref{prop.classical.skeletal.ords} can be turned into an explicit formula (not
just an upper bound as in Corollary \ref{cor.classical.for.ord}):\footnote{%
Update 2022: We thank Bruce Sagan for finding a typo in an earlier version
of this formula.}

\begin{proposition}
\label{prop.classical.for.ord.explicit}Let $n\in\mathbb{N}$. Let $P$ be an
$n$-graded poset. Assume that $P$ is a rooted forest (made into a poset by
having every node smaller than its children). Notice that $\left\vert
\widehat{P}_{i}\right\vert \leq\left\vert \widehat{P}_{i+1}\right\vert $ for
every $i\in\left\{  0,1,...,n-1\right\}  $ (where $\widehat{P}_{i}$ and
$\widehat{P}_{i+1}$ are defined as in Definition \ref{def.graded.Phat}). Then,%
\[
\operatorname*{ord}\left(  \overline{\mathbf{r}}_{P}\right)
=\operatorname{lcm}\left\{  n+1-i\ \mid\ i\in\left\{  0,1,...,n-1\right\}
;\ \left\vert \widehat{P}_{i}\right\vert <\left\vert \widehat{P}%
_{i+1}\right\vert \right\}  .
\]
(Of course, $\operatorname*{ord}\left(  \mathbf{r}_{P}\right)  $ can now be
computed by $\operatorname*{ord}\left(  \mathbf{r}_{P}\right)
=\operatorname{lcm}\left(  n+1,\operatorname*{ord}\left(  \overline
{\mathbf{r}}_{P}\right)  \right)  $.)
\end{proposition}

The same property therefore holds for birational rowmotion $R_{P}$ and its
homogeneous version $\overline{R}_{P}$.

\begin{proof}[Proof of Proposition \ref{prop.classical.for.ord.explicit} (sketched).]
We define $\mathcal{N}_{P}$ to be the set \newline
$\left\{  n+1-i\ \mid\ i\in\left\{
0,1,...,n-1\right\}  ;\ \left\vert \widehat{P}_{i}\right\vert <\left\vert
\widehat{P}_{i+1}\right\vert \right\}  $. Thus, we must prove that
$\operatorname*{ord}\left(  \overline{\mathbf{r}}_{P}\right)
=\operatorname{lcm}\left(  \mathcal{N}_{P}\right)  $.

We proceed by strong induction on $\left\vert P\right\vert $. So we fix an
$n$-graded poset $P$ for some $n\in\mathbb{N}$, and we set out to prove the
equality
\begin{equation}
\operatorname*{ord}\left(  \overline{\mathbf{r}}_{P}\right)
=\operatorname{lcm}\left(  \mathcal{N}_{P}\right)  ,
\label{pf.prop.classical.for.ord.explicit.goal}%
\end{equation}
assuming (as the induction hypothesis) that the analogous equality
$\operatorname*{ord}\left(  \overline{\mathbf{r}}_{P^{\prime}}\right)
=\operatorname{lcm}\left(  \mathcal{N}_{P^{\prime}}\right)  $ has been proved
for all $n^{\prime}\in\mathbb{N}$ and all $n^{\prime}$-graded posets
$P^{\prime}$ with $\left\vert P^{\prime}\right\vert <\left\vert P\right\vert $.

We are in one of the following three cases:

\textit{Case 0:} The poset $P$ has no minimal elements.

\textit{Case 1:} The poset $P$ has exactly one minimal element.

\textit{Case 2:} The poset $P$ has more than one minimal element.

Case 0 is easy: In this case, we must have $P=\varnothing$ (since every
nonempty finite poset has at least one minimal element), and thus it is easy
to see that both $\operatorname*{ord}\left(  \overline{\mathbf{r}}_{P}\right)
$ and $\operatorname{lcm}\left(  \mathcal{N}_{P}\right)  $ equal $1$ (since
$\overline{\mathbf{r}}_{P}$ is the identity map on the $1$-element set
$\overline{J\left(  P\right)  }=\left\{  \pi\left(  \varnothing\right)
\right\}  $, while $\mathcal{N}_{P}$ is an empty set and thus has
$\operatorname{lcm}$ equal to $1$). Thus,
(\ref{pf.prop.classical.for.ord.explicit.goal}) is proved in Case 0.

Let us now consider Case 1. In this case, the poset $P$ has exactly one
minimal element. Thus, $P=B_{1}Q$ for some poset $Q$. Consider this $Q$. This
poset $Q$ is $\left(  n-1\right)  $-graded (since $B_{1}Q=P$ is $n$-graded)
and is a rooted forest (since $B_{1}Q=P$ is a rooted forest). Since we
furthermore have $\left\vert Q\right\vert <\left\vert P\right\vert $ (because
$P=B_{1}Q$), we can apply our induction hypothesis to $n^{\prime}=n-1$ and
$P^{\prime}=Q$. We thus obtain%
\[
\operatorname*{ord}\left(  \overline{\mathbf{r}}_{Q}\right)
=\operatorname{lcm}\left(  \mathcal{N}_{Q}\right)  .
\]
On the other hand, from $P=B_{1}Q$, we obtain
\begin{align}
\operatorname*{ord}\left(  \overline{\mathbf{r}}_{P}\right)   &
=\operatorname*{ord}\left(  \overline{\mathbf{r}}_{B_{1}Q}\right)
=\operatorname*{ord}\left(  \overline{\mathbf{r}}_{Q}\right)
\ \ \ \ \ \ \ \ \ \ \left(  \text{by Proposition \ref{prop.classical.Bk.ord}
\textbf{(a)}}\right) \nonumber\\
&  =\operatorname{lcm}\left(  \mathcal{N}_{Q}\right)  .
\label{pf.prop.classical.for.ord.explicit.c1.0}%
\end{align}

For each $i\in\left\{  1,2,\ldots,n+1\right\}  $, there is a canonical
bijection $\widehat{P}_{i}\cong\widehat{Q}_{i-1}$ (since $P=B_{1}Q$); thus, we
have
\begin{equation}
\left\vert \widehat{P}_{i}\right\vert =\left\vert \widehat{Q}_{i-1}\right\vert
\ \ \ \ \ \ \ \ \ \ \text{for each }i\in\left\{  1,2,\ldots,n+1\right\}  .
\label{pf.prop.classical.for.ord.explicit.c1.1}%
\end{equation}
Moreover, $\left\vert \widehat{P}_{1}\right\vert =1$ (since $P=B_{1}Q$). Also,
$\left\vert \widehat{P}_{0}\right\vert =1$ (since $\widehat{P}_{0}=\left\{
0\right\}  $). Thus, $\left\vert \widehat{P}_{0}\right\vert =1=\left\vert
\widehat{P}_{1}\right\vert $. Hence, we don't have $\left\vert \widehat{P}%
_{0}\right\vert <\left\vert \widehat{P}_{1}\right\vert $. Therefore, $0$ is
\textbf{not} an $i\in\left\{  0,1,...,n-1\right\}  $ satisfying $\left\vert
\widehat{P}_{i}\right\vert <\left\vert \widehat{P}_{i+1}\right\vert $.

Now, the definition of $\mathcal{N}_{P}$ yields%
\begin{align*}
\mathcal{N}_{P}  &  =\left\{  n+1-i\ \mid\ i\in\left\{  0,1,...,n-1\right\}
;\ \left\vert \widehat{P}_{i}\right\vert <\left\vert \widehat{P}%
_{i+1}\right\vert \right\} \\
&  =\left\{  n+1-i\ \mid\ i\in\left\{  1,2,...,n-1\right\}
;\ \underbrace{\left\vert \widehat{P}_{i}\right\vert }_{\substack{=\left\vert
\widehat{Q}_{i-1}\right\vert \\\text{(by
(\ref{pf.prop.classical.for.ord.explicit.c1.1}))}}}<\underbrace{\left\vert
\widehat{P}_{i+1}\right\vert }_{\substack{=\left\vert \widehat{Q}%
_{i}\right\vert \\\text{(by (\ref{pf.prop.classical.for.ord.explicit.c1.1}))}%
}}\right\} \\
&  \ \ \ \ \ \ \ \ \ \ \ \ \ \ \ \ \ \ \ \ \left(
\begin{array}
[c]{c}%
\text{since we don't have }\left\vert \widehat{P}_{0}\right\vert <\left\vert
\widehat{P}_{1}\right\vert \text{, so that our set}\\
\text{contains no element corresponding to }i=0
\end{array}
\right) \\
&  =\left\{  n+1-i\ \mid\ i\in\left\{  1,2,...,n-1\right\}  ;\ \left\vert
\widehat{Q}_{i-1}\right\vert <\left\vert \widehat{Q}_{i}\right\vert \right\}
\\
&  =\left\{  \underbrace{n+1-\left(  i+1\right)  }_{=n-i}\ \mid\ i\in\left\{
0,1,\ldots,n-2\right\}  ;\ \left\vert \widehat{Q}_{i}\right\vert <\left\vert
\widehat{Q}_{i+1}\right\vert \right\} \\
&  \ \ \ \ \ \ \ \ \ \ \ \ \ \ \ \ \ \ \ \ \left(  \text{here, we have
substituted }i+1\text{ for the index }i\right) \\
&  =\left\{  n-i\ \mid\ i\in\left\{  0,1,\ldots,n-2\right\}  ;\ \left\vert
\widehat{Q}_{i}\right\vert <\left\vert \widehat{Q}_{i+1}\right\vert \right\}
\\
&  =\mathcal{N}_{Q}%
\end{align*}
(by the definition of $\mathcal{N}_{Q}$, since $Q$ is $\left(  n-1\right)
$-graded). In view of this, we can rewrite
(\ref{pf.prop.classical.for.ord.explicit.c1.0}) as $\operatorname*{ord}\left(
\overline{\mathbf{r}}_{P}\right)  =\operatorname{lcm}\left(  \mathcal{N}%
_{P}\right)  $. This proves (\ref{pf.prop.classical.for.ord.explicit.goal}) in
Case 1.

Let us now consider Case 2. In this case, the poset $P$ has more than one
minimal element. Hence, $P$ cannot be a rooted tree. Since $P$ is a rooted
forest, we thus conclude that $P$ consists of more than one rooted tree.
Therefore, $P$ is a disjoint union of more than one nonempty poset. In other
words, $P=QR$ for two nonempty posets $Q$ and $R$ (note that each of $Q$ and
$R$ can itself be a nontrivial disjoint union). Consider these $Q$ and $R$.
Since the poset $QR=P$ is $n$-graded, we see that the poset $Q$ is $n$-graded
(because any minimal element of $Q$ is a minimal element of $QR$, whereas any
maximal element of $Q$ is a maximal element of $QR$). Likewise, the poset $R$
is $n$-graded. Hence, Proposition \ref{prop.classical.PQ.ord2} yields
$\operatorname*{ord}\left(  \overline{\mathbf{r}}_{QR}\right)
=\operatorname{lcm}\left(  \operatorname*{ord}\left(  \mathbf{r}_{Q}\right)
,\operatorname*{ord}\left(  \mathbf{r}_{R}\right)  \right)  $. In view of
$P=QR$, we can rewrite this as%
\begin{equation}
\operatorname*{ord}\left(  \overline{\mathbf{r}}_{P}\right)
=\operatorname{lcm}\left(  \operatorname*{ord}\left(  \mathbf{r}_{Q}\right)
,\operatorname*{ord}\left(  \mathbf{r}_{R}\right)  \right)  .
\label{pf.prop.classical.for.ord.explicit.c2.0}%
\end{equation}

However, since $Q$ is $n$-graded, we have $\operatorname*{ord}\left(
\mathbf{r}_{Q}\right)  =\operatorname{lcm}\left(  n+1,\operatorname*{ord}%
\left(  \overline{\mathbf{r}}_{Q}\right)  \right)  $ (by Proposition
\ref{prop.classical.ord-projord}). Similarly, $\operatorname*{ord}\left(
\mathbf{r}_{R}\right)  =\operatorname{lcm}\left(  n+1,\operatorname*{ord}%
\left(  \overline{\mathbf{r}}_{R}\right)  \right)  $. Using of these two
equalities, we can rewrite (\ref{pf.prop.classical.for.ord.explicit.c2.0}) as%
\begin{align}
\operatorname*{ord}\left(  \overline{\mathbf{r}}_{P}\right)   &
=\operatorname{lcm}\left(  \operatorname{lcm}\left(  n+1,\operatorname*{ord}%
\left(  \overline{\mathbf{r}}_{Q}\right)  \right)  ,\operatorname{lcm}\left(
n+1,\operatorname*{ord}\left(  \overline{\mathbf{r}}_{R}\right)  \right)
\right) \nonumber\\
&  =\operatorname{lcm}\left(  n+1,\operatorname*{ord}\left(  \overline
{\mathbf{r}}_{Q}\right)  ,\operatorname*{ord}\left(  \overline{\mathbf{r}}%
_{R}\right)  \right)  . \label{pf.prop.classical.for.ord.explicit.c2.1}%
\end{align}

From $P=QR$, we obtain $\left\vert P\right\vert =\left\vert QR\right\vert
=\left\vert Q\right\vert +\left\vert R\right\vert >\left\vert Q\right\vert $
(since $R$ is nonempty). Hence, we can apply our induction hypothesis to
$P^{\prime}=Q$ and $n^{\prime}=n$. We thus obtain $\operatorname*{ord}\left(
\overline{\mathbf{r}}_{Q}\right)  =\operatorname{lcm}\left(  \mathcal{N}%
_{Q}\right)  $. Similarly, $\operatorname*{ord}\left(  \overline{\mathbf{r}%
}_{R}\right)  =\operatorname{lcm}\left(  \mathcal{N}_{R}\right)  $. Using
these two equalities, we can rewrite
(\ref{pf.prop.classical.for.ord.explicit.c2.1}) as%
\begin{align}
\operatorname*{ord}\left(  \overline{\mathbf{r}}_{P}\right)   &
=\operatorname{lcm}\left(  n+1,\operatorname{lcm}\left(  \mathcal{N}%
_{Q}\right)  ,\operatorname{lcm}\left(  \mathcal{N}_{R}\right)  \right)
\nonumber\\
&  =\operatorname{lcm}\left(  \left\{  n+1\right\}  \cup\mathcal{N}_{Q}%
\cup\mathcal{N}_{R}\right)  . \label{pf.prop.classical.for.ord.explicit.c2.2}%
\end{align}

On the other hand, from $P=QR$, we see that%
\begin{equation}
\left\vert \widehat{P}_{i}\right\vert =\left\vert \widehat{Q}_{i}\right\vert
+\left\vert \widehat{R}_{i}\right\vert \ \ \ \ \ \ \ \ \ \ \text{for each
}i\in\left\{  1,2,\ldots,n\right\}  .
\label{pf.prop.classical.for.ord.explicit.c2.Pi}%
\end{equation}
In particular, $\left\vert \widehat{P}_{1}\right\vert =\underbrace{\left\vert
\widehat{Q}_{1}\right\vert }_{\geq1}+\underbrace{\left\vert \widehat{R}%
_{1}\right\vert }_{\geq1}\geq1+1=2>1$. Thus, $\left\vert \widehat{P}%
_{0}\right\vert =1<\left\vert \widehat{P}_{1}\right\vert $.

Now, let $i\in\left\{  1,2,\ldots,n-1\right\}  $ be arbitrary. Since $Q$ is a
rooted forest, we can easily see that $\left\vert \widehat{Q}_{i}\right\vert
\leq\left\vert \widehat{Q}_{i+1}\right\vert $\ \ \ \ \footnote{\textit{Proof.}
The poset $Q$ is $n$-graded; thus, any maximal element of $Q$ has degree $n$.
Any element $v\in\widehat{Q}_{i}$ has degree $\deg v=i\leq n-1<n$, and thus
cannot be maximal (by the preceding sentence). Thus, any element
$v\in\widehat{Q}_{i}$ has a child (viewed as a vertex of the forest $Q$).
Moreover, any child of $v$ must cover $v$ in the poset $Q$, and thus must have
degree $i+1$ (since $v$ has degree $i$). In other words, any child of $v$ must
belong to $\widehat{Q}_{i+1}$.
\par
Thus, we can define a map $\psi:\widehat{Q}_{i}\rightarrow\widehat{Q}_{i+1}$
as follows: For any $v\in\widehat{Q}_{i}$, we randomly pick a child of $v$,
and we let $\psi\left(  v\right)  $ be this child. This map $\psi$ is
injective (since no two distinct parents can have a common child in a rooted
forest). Thus, we have found an injective map from $\widehat{Q}_{i}$ to
$\widehat{Q}_{i+1}$ (namely, $\psi$). This shows that $\left\vert
\widehat{Q}_{i}\right\vert \leq\left\vert \widehat{Q}_{i+1}\right\vert $.}.
Similarly, $\left\vert \widehat{R}_{i}\right\vert \leq\left\vert
\widehat{R}_{i+1}\right\vert $. By summing these two inequalities, we obtain
the inequality%
\begin{equation}
\left\vert \widehat{Q}_{i}\right\vert +\left\vert \widehat{R}_{i}\right\vert
\leq\left\vert \widehat{Q}_{i+1}\right\vert +\left\vert \widehat{R}%
_{i+1}\right\vert . \label{pf.prop.classical.for.ord.explicit.c2.4}%
\end{equation}
Clearly, this inequality (\ref{pf.prop.classical.for.ord.explicit.c2.4}) is
strict if and only if any one of the two inequalities $\left\vert
\widehat{Q}_{i}\right\vert \leq\left\vert \widehat{Q}_{i+1}\right\vert $ and
$\left\vert \widehat{R}_{i}\right\vert \leq\left\vert \widehat{R}%
_{i+1}\right\vert $ is strict. However, the inequality
(\ref{pf.prop.classical.for.ord.explicit.c2.4}) can be rewritten as%
\begin{equation}
\left\vert \widehat{P}_{i}\right\vert \leq\left\vert \widehat{P}%
_{i+1}\right\vert \label{pf.prop.classical.for.ord.explicit.c2.5}%
\end{equation}
(since (\ref{pf.prop.classical.for.ord.explicit.c2.Pi}) yields $\left\vert
\widehat{P}_{i}\right\vert =\left\vert \widehat{Q}_{i}\right\vert +\left\vert
\widehat{R}_{i}\right\vert $ and $\left\vert \widehat{P}_{i+1}\right\vert
=\left\vert \widehat{Q}_{i+1}\right\vert +\left\vert \widehat{R}%
_{i+1}\right\vert $). Thus, we conclude that the inequality
(\ref{pf.prop.classical.for.ord.explicit.c2.5}) is strict if and only if any
one of the two inequalities $\left\vert \widehat{Q}_{i}\right\vert
\leq\left\vert \widehat{Q}_{i+1}\right\vert $ and $\left\vert \widehat{R}%
_{i}\right\vert \leq\left\vert \widehat{R}_{i+1}\right\vert $ is strict. In
other words, the logical equivalence%
\begin{equation}
\left(  \left\vert \widehat{P}_{i}\right\vert <\left\vert \widehat{P}%
_{i+1}\right\vert \right)  \ \Longleftrightarrow\ \left(  \left\vert
\widehat{Q}_{i}\right\vert <\left\vert \widehat{Q}_{i+1}\right\vert \text{ or
}\left\vert \widehat{R}_{i}\right\vert <\left\vert \widehat{R}_{i+1}%
\right\vert \right)  \label{pf.prop.classical.for.ord.explicit.c2.6}%
\end{equation}
holds.

Forget that we fixed $i$. We thus have proved the equivalence
(\ref{pf.prop.classical.for.ord.explicit.c2.6}) for each $i\in\left\{
1,2,\ldots,n-1\right\}  $.

Now, the definition of $\mathcal{N}_{P}$ yields%
\begin{align*}
\mathcal{N}_{P} &  =\left\{  n+1-i\ \mid\ i\in\left\{  0,1,...,n-1\right\}
;\ \left\vert \widehat{P}_{i}\right\vert <\left\vert \widehat{P}%
_{i+1}\right\vert \right\}  \\
&  =\left\{  n+1\right\}  \cup\left\{  n+1-i\ \mid\ i\in\left\{
1,2,...,n-1\right\}  ;\ \underbrace{\left\vert \widehat{P}_{i}\right\vert
<\left\vert \widehat{P}_{i+1}\right\vert }_{\substack{\Longleftrightarrow
\ \left(  \left\vert \widehat{Q}_{i}\right\vert <\left\vert \widehat{Q}%
_{i+1}\right\vert \text{ or }\left\vert \widehat{R}_{i}\right\vert <\left\vert
\widehat{R}_{i+1}\right\vert \right)  \\\text{(by
(\ref{pf.prop.classical.for.ord.explicit.c2.6}))}}}\right\}  \\
&  \ \ \ \ \ \ \ \ \ \ \ \ \ \ \ \ \ \ \ \ \left(  \text{since }\left\vert
\widehat{P}_{0}\right\vert <\left\vert \widehat{P}_{1}\right\vert \right)  \\
&  =\left\{  n+1\right\}  \cup\underbrace{\left\{  n+1-i\ \mid\ i\in\left\{
1,2,...,n-1\right\}  ;\ \left(  \left\vert \widehat{Q}_{i}\right\vert
<\left\vert \widehat{Q}_{i+1}\right\vert \text{ or }\left\vert \widehat{R}%
_{i}\right\vert <\left\vert \widehat{R}_{i+1}\right\vert \right)  \right\}
}_{=\left\{  n+1-i\ \mid\ i\in\left\{  1,2,...,n-1\right\}  ;\ \left\vert
\widehat{Q}_{i}\right\vert <\left\vert \widehat{Q}_{i+1}\right\vert \right\}
\cup\left\{  n+1-i\ \mid\ i\in\left\{  1,2,...,n-1\right\}  ;\ \left\vert
\widehat{R}_{i}\right\vert <\left\vert \widehat{R}_{i+1}\right\vert \right\}
}\\
&  =\left\{  n+1\right\}  \cup\left\{  n+1-i\ \mid\ i\in\left\{
1,2,...,n-1\right\}  ;\ \left\vert \widehat{Q}_{i}\right\vert <\left\vert
\widehat{Q}_{i+1}\right\vert \right\}  \\
&  \ \ \ \ \ \ \ \ \ \ \cup\left\{  n+1-i\ \mid\ i\in\left\{
1,2,...,n-1\right\}  ;\ \left\vert \widehat{R}_{i}\right\vert <\left\vert
\widehat{R}_{i+1}\right\vert \right\}  \\
&  =\left\{  n+1\right\}  \cup\underbrace{\left\{  n+1-i\ \mid\ i\in\left\{
0,1,...,n-1\right\}  ;\ \left\vert \widehat{Q}_{i}\right\vert <\left\vert
\widehat{Q}_{i+1}\right\vert \right\}  }_{\substack{=\mathcal{N}%
_{Q}\\\text{(by the definition of }\mathcal{N}_{Q}\text{, since }Q\text{ is
}n\text{-graded)}}}\\
&  \ \ \ \ \ \ \ \ \ \ \cup\underbrace{\left\{  n+1-i\ \mid\ i\in\left\{
0,1,...,n-1\right\}  ;\ \left\vert \widehat{R}_{i}\right\vert <\left\vert
\widehat{R}_{i+1}\right\vert \right\}  }_{\substack{=\mathcal{N}%
_{R}\\\text{(by the definition of }\mathcal{N}_{R}\text{, since }R\text{ is
}n\text{-graded)}}}\\
&  \ \ \ \ \ \ \ \ \ \ \ \ \ \ \ \ \ \ \ \ \left(
\begin{array}
[c]{c}%
\text{here, we have extended the indexing set for the two sets}\\
\left\{  n+1-i\ \mid\ i\in\left\{  1,2,...,n-1\right\}  ;\ \left\vert
\widehat{Q}_{i}\right\vert <\left\vert \widehat{Q}_{i+1}\right\vert \right\}
\\
\text{and }\left\{  n+1-i\ \mid\ i\in\left\{  0,1,...,n-1\right\}
;\ \left\vert \widehat{R}_{i}\right\vert <\left\vert \widehat{R}%
_{i+1}\right\vert \right\}  \\
\text{from }\left\{  1,2,\ldots,n-1\right\}  \text{ to }\left\{
0,1,\ldots,n-1\right\}  \text{;}\\
\text{this can potentially insert the element }n+1\text{ into these two
sets,}\\
\text{but ultimately does not affect their union with }\left\{  n+1\right\}
\text{,}\\
\text{since the set }\left\{  n+1\right\}  \text{ already contains }n+1\text{
anyway}%
\end{array}
\right)  \\
&  =\left\{  n+1\right\}  \cup\mathcal{N}_{Q}\cup\mathcal{N}_{R}.
\end{align*}
Hence,%
\[
\operatorname{lcm}\left(  \mathcal{N}_{P}\right)  =\operatorname{lcm}\left(
\left\{  n+1\right\}  \cup\mathcal{N}_{Q}\cup\mathcal{N}_{R}\right)  .
\]
Comparing this with (\ref{pf.prop.classical.for.ord.explicit.c2.2}), we obtain
$\operatorname*{ord}\left(  \overline{\mathbf{r}}_{P}\right)
=\operatorname{lcm}\left(  \mathcal{N}_{P}\right)  $. Hence,
(\ref{pf.prop.classical.for.ord.explicit.goal}) is proved in Case 2.

We have now proved the equality (\ref{pf.prop.classical.for.ord.explicit.goal}%
) in all three Cases 0, 1 and 2. Thus,
(\ref{pf.prop.classical.for.ord.explicit.goal}) always holds. This completes
the induction step, and thus Proposition \ref{prop.classical.for.ord.explicit}
is proven.
\end{proof}

\section{The rectangle: statements of the results}

\begin{definition}
\label{def.rect}Let $p$ and $q$ be two positive integers. The $p\times
q$\textit{-rectangle} will denote the poset $\left\{  1,2,...,p\right\}
\times\left\{  1,2,...,q\right\}  $ with order defined as follows: For two
elements $\left(  i,k\right)  $ and $\left(  i^{\prime},k^{\prime}\right)  $
of $\left\{  1,2,...,p\right\}  \times\left\{  1,2,...,q\right\}  $, we set
$\left(  i,k\right)  \leq\left(  i^{\prime},k^{\prime}\right)  $ if and only
if $\left(  i\leq i^{\prime}\text{ and }k\leq k^{\prime}\right)  $.
\end{definition}

\Needspace{9\baselineskip}

\begin{example}
Here is the Hasse diagram of the $2\times3$-rectangle:%
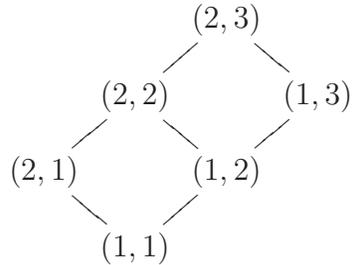
\[
\xymatrixrowsep{0.9pc}\xymatrixcolsep{0.20pc}\xymatrix{
& & \left(2,3\right) \ar@{-}[rd] \ar@{-}[ld] & \\
& \left(2,2\right) \ar@{-}[rd] \ar@{-}[ld] & & \left(1,3\right) \ar@{-}[ld]\\
\left(2,1\right) \ar@{-}[rd] & & \left(1,2\right) \ar@{-}[ld] & \\
& \left(1,1\right) & &
}.
\]

\end{example}

\begin{remark}
Let $p$ and $q$ be positive integers. The $p\times q$-rectangle is denoted by
$\left[  p\right]  \times\left[  q\right]  $ in the papers
\cite{striker-williams}, \cite{einstein-propp}, \cite{propp-roby}
and \cite{propp-roby-full}.
\end{remark}

\begin{remark}
\label{rmk.rect.cover}Let $p$ and $q$ be two positive integers. Let
$\operatorname*{Rect}\left(  p,q\right)  $ denote the $p\times q$-rectangle.

\textbf{(a)} The $p\times q$-rectangle is a $\left(  p+q-1\right)  $-graded
poset, with $\deg\left(  \left(  i,k\right)  \right)  =i+k-1$ for all $\left(
i,k\right)  \in\operatorname*{Rect}\left(  p,q\right)  $.

\textbf{(b)} Let $\left(  i,k\right)  $ and $\left(  i^{\prime},k^{\prime
}\right)  $ be two elements of $\operatorname*{Rect}\left(  p,q\right)  $.
Then, $\left(  i,k\right)  \lessdot\left(  i^{\prime},k^{\prime}\right)  $ if
and only if either $\left(  i^{\prime}=i\text{ and }k^{\prime}=k+1\right)  $
or $\left(  k^{\prime}=k\text{ and }i^{\prime}=i+1\right)  $.
\end{remark}

We are going to use Remark \ref{rmk.rect.cover} without explicit mention.

The following theorem was conjectured by James Propp and the second author:

\begin{theorem}
\label{thm.rect.ord}Let $\operatorname*{Rect}\left(  p,q\right)  $ denote the
$p\times q$-rectangle. Let $\mathbb{K}$ be a field. Then, $\operatorname*{ord}%
\left(  R_{\operatorname*{Rect}\left(  p,q\right)  }\right)  =p+q$.
\end{theorem}

This is a birational analogue (and, using the reasoning of
\cite{einstein-propp}, generalization) of the classical fact (appearing in
\cite[Theorem 3.1]{striker-williams} and \cite[Theorem 2]{fon-der-flaass})
that $\operatorname*{ord}\left(  \mathbf{r}_{\operatorname*{Rect}\left(
p,q\right)  }\right)  =p+q$ (using the notations of Definition
\ref{def.classical.rm} and Definition \ref{def.classical.conv.rP}).

Notice that Proposition \ref{prop.ord-projord} yields that $p+q\mid
\operatorname*{ord}\left(  R_{\operatorname*{Rect}\left(  p,q\right)
}\right)  $, so all that needs to be proven in order to verify Theorem
\ref{thm.rect.ord} is showing that $R_{\operatorname*{Rect}\left(  p,q\right)
}^{p+q}=\operatorname*{id}$.

Notice also that in the case when $p\leq2$ and $q\leq2$, Theorem
\ref{thm.rect.ord} follows rather easily from Propositions \ref{prop.Bk.ord}
\textbf{(a)}, \ref{prop.B'k.ord} \textbf{(a)} and \ref{prop.ord-projord}
(because $\operatorname*{Rect}\left(  p,q\right)  $ is a skeletal poset in
this case), but this approach does not generalize to any interesting cases.

\begin{remark}
\label{rmk.rect.jdt}Theorem \ref{thm.rect.ord} generalizes a well-known
property of promotion on semistandard Young tableaux of rectangular shape,
albeit not in an obvious way. Here are some details (which a reader
unacquainted with Young tableaux can freely skip):

Let $N$ be a nonnegative integer, and let $\lambda$ be a partition. Let
$\operatorname*{SSYT}\nolimits_{N}\lambda$ denote the set of all semistandard
Young tableaux of shape $\lambda$ whose entries are all $\leq N$. One can
define a map $\operatorname*{Pro}:\operatorname*{SSYT}\nolimits_{N}%
\lambda\rightarrow\operatorname*{SSYT}\nolimits_{N}\lambda$ called
\textit{jeu-de-taquin promotion} (or \textit{Sch\"{u}tzenberger promotion}, or
simply \textit{promotion} when no ambiguities can arise); see \cite[\S 5.1]%
{russell} for a precise definition. (The definition in \cite[\S 2]{rhoades} is
different -- it defines the \textit{inverse} of this map. Conventions differ.)
This map has some interesting properties already for arbitrary $\lambda$, but
the most interesting situation is that of $\lambda$ being a rectangular
partition (i.e., a partition all of whose nonzero parts are equal). In this
situation, a folklore theorem states that $\operatorname*{Pro}\nolimits^{N}%
=\operatorname*{id}$. (The particular case of this theorem when
$\operatorname*{Pro}$ is applied only to \textbf{standard} Young tableaux is
well-known -- see, e.g., \cite[Theorem 4.4]{haiman} --, but the only proof of
the general theorem that we were able to find in literature is Rhoades's --
\cite[Corollary 5.6]{rhoades} --, which makes use of Kazhdan-Lusztig theory.)

Theorem \ref{thm.rect.ord} can be used to give an alternative proof of this
$\operatorname*{Pro}\nolimits^{N}=\operatorname*{id}$ theorem. See a future
version of \cite{einstein-propp} (or, for the time being, \cite[\S 2, pp.
4--5]{einstein-propp-fpsac}) for how this works.

Note that \cite[\S 5.1]{russell}, \cite[\S 2]{rhoades} and
\cite{einstein-propp} give three different definitions of promotion. The
definitions in \cite[\S 5.1]{russell} and in \cite{einstein-propp} are
equivalent, while the definition in \cite[\S 2]{rhoades} defines the inverse
of the map defined in the other two sources. Unfortunately, we were unable to
find the proofs of these facts in existing literature; they are claimed in
\cite[Propositions 2.5 and 2.6]{kirillov-berenstein}, and can be proven using
the concept of tableau switching \cite[Definition 2.2.1]{leeuwen-lrr}.
\end{remark}

Besides Theorem \ref{thm.rect.ord}, we can also state some kind of symmetry
property of birational rowmotion on the $p\times q$-rectangle (referred to as
the \textquotedblleft pairing property\textquotedblright\ in
\cite{einstein-propp}), which was also conjectured by Propp and the second author:

\begin{theorem}
\label{thm.rect.antip.general}Let $\operatorname*{Rect}\left(  p,q\right)  $
denote the $p\times q$-rectangle. Let $\mathbb{K}$ be a field. Let
$f\in\mathbb{K}^{\widehat{\operatorname*{Rect}\left(  p,q\right)  }}$. Assume
that $R_{\operatorname*{Rect}\left(  p,q\right)  }^{\ell}f$ is well-defined
for every $\ell\in\left\{  0,1,...,i+k-1\right\}  $. Let $\left(  i,k\right)
\in\operatorname*{Rect}\left(  p,q\right)  $. Then,%
\[
f\left(  \left(  p+1-i,q+1-k\right)  \right)  =\dfrac{f\left(  0\right)
f\left(  1\right)  }{\left(  R_{\operatorname*{Rect}\left(  p,q\right)
}^{i+k-1}f\right)  \left(  \left(  i,k\right)  \right)  }.
\]

\end{theorem}

This Theorem generalizes the \textquotedblleft reciprocity
phenomenon\textquotedblright\ observed on the $2\times2$-rectangle in Example
\ref{ex.rowmotion.2x2}.

\begin{remark}
\label{rmk.rect.ord.application}While Theorem \ref{thm.rect.ord} only makes a
statement about $R_{\operatorname*{Rect}\left(  p,q\right)  }$, it can be used
(in combination with Proposition \ref{prop.Bk.ord} and others) to derive upper
bounds on the orders of $R_{P}$ and $\overline{R}_{P}$ for some other posets
$P$. Here is an example: Let $\mathbb{K}$ be a field. For the duration of this
remark, let us denote the poset $\operatorname*{Rect}\left(  2,3\right)
\setminus\left\{  \left(  1,1\right)  ,\left(  2,3\right)  \right\}  $ by $N$.
(The Hasse diagram of this poset has the rather simple form
$\xymatrixrowsep{0.9pc}\xymatrixcolsep{0.20pc}
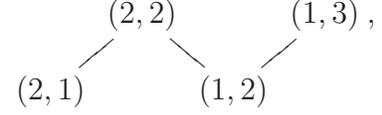
\xymatrix{
& \left(2,2\right) \ar@{-}[rd] \ar@{-}[ld] & & \left(1,3\right) \ar@{-}[ld]\\
\left(2,1\right) & & \left(1,2\right) & \\
}$, which explains why we have chosen to call it $N$ here.) Then,
$\operatorname*{ord}\left(  R_{N}\right)  \mid15$ and $\operatorname*{ord}%
\left(  \overline{R}_{N}\right)  \mid5$. This can be proven as follows: We
have $\operatorname*{Rect}\left(  2,3\right)  \cong B_{1}\left(  B_{1}%
^{\prime}N\right)  $ and therefore%
\begin{align*}
\operatorname*{ord}\left(  \overline{R}_{\operatorname*{Rect}\left(
2,3\right)  }\right)   &  =\operatorname*{ord}\left(  \overline{R}%
_{B_{1}\left(  B_{1}^{\prime}N\right)  }\right)  =\operatorname*{ord}\left(
\overline{R}_{B_{1}^{\prime}N}\right)  \ \ \ \ \ \ \ \ \ \ \left(  \text{by
Proposition \ref{prop.Bk.ord} \textbf{(a)}}\right) \\
&  =\operatorname*{ord}\left(  \overline{R}_{N}\right)
\ \ \ \ \ \ \ \ \ \ \left(  \text{by Proposition \ref{prop.B'k.ord}
\textbf{(a)}}\right)  ,
\end{align*}
so that
\begin{align*}
\operatorname*{ord}\left(  \overline{R}_{N}\right)   &  =\operatorname*{ord}%
\left(  \overline{R}_{\operatorname*{Rect}\left(  2,3\right)  }\right) \\
&  \mid\operatorname{lcm}\left(  4+1,\operatorname*{ord}\left(  \overline
{R}_{\operatorname*{Rect}\left(  2,3\right)  }\right)  \right)
=\operatorname*{ord}\left(  R_{\operatorname*{Rect}\left(  2,3\right)
}\right) \\
&  \ \ \ \ \ \ \ \ \ \ \left(  \text{by Proposition \ref{prop.ord-projord},
since }\operatorname*{Rect}\left(  2,3\right)  \text{ is }4\text{-graded}%
\right) \\
&  =2+3\ \ \ \ \ \ \ \ \ \ \left(  \text{by Theorem \ref{thm.rect.ord}}\right)
\\
&  =5
\end{align*}
and thus%
\begin{align*}
\operatorname*{ord}\left(  R_{N}\right)   &  =\operatorname{lcm}\left(
\underbrace{2+1}_{=3},\underbrace{\operatorname*{ord}\left(  \overline{R}%
_{N}\right)  }_{\mid5}\right)  \ \ \ \ \ \ \ \ \ \ \left(  \text{by
Proposition \ref{prop.ord-projord}, since }N\text{ is }2\text{-graded}\right)
\\
&  \mid\operatorname{lcm}\left(  3,5\right)  =15.
\end{align*}
It can actually be shown that $\operatorname*{ord}\left(  R_{N}\right)  =15$
and $\operatorname*{ord}\left(  \overline{R}_{N}\right)  =5$ by direct computation.

In the same vein it can be shown that $\operatorname*{ord}\left(  \overline
{R}_{\operatorname*{Rect}\left(  p,q\right)  \setminus\left\{  \left(
1,1\right)  ,\left(  p,q\right)  \right\}  }\right)  \mid p+q$ and
$\operatorname*{ord}\left(  R_{\operatorname*{Rect}\left(  p,q\right)
\setminus\left\{  \left(  1,1\right)  ,\left(  p,q\right)  \right\}  }\right)
\mid\operatorname{lcm}\left(  p+q-2,p+q\right)  $ for any integers $p>1$ and
$q>1$. This doesn't, however, generalize to arbitrary posets obtained by
removing some ranks from $\operatorname*{Rect}\left(  p,q\right)  $ (indeed,
sometimes birational rowmotion doesn't even have finite order on such posets,
cf. Section \ref{sect.negres}).
\end{remark}

\section{Reduced labellings}

The proof that we give for Theorem \ref{thm.rect.ord} and Theorem
\ref{thm.rect.antip.general} is largely inspired by the proof of
Zamolodchikov's conjecture in case $AA$ given by Volkov in \cite{volkov}%
\footnote{``Case $AA$'' refers to the Cartesian product of the Dynkin diagrams
of two type-$A$ root systems. This, of course, is a rectangle, just as in our
Theorem \ref{thm.rect.ord}.}. This is not very surprising because the orbit of
a $\mathbb{K}$-labelling under birational rowmotion appears superficially
similar to a solution of a $Y$-system of type $AA$. Yet we do not see a way to
derive Theorem \ref{thm.rect.ord} from Zamolodchikov's conjecture or vice
versa. (It should be noticed that Zamolodchikov's Y-system has an obvious
``reducibility property'', namely consisting of two decoupled subsystems -- a
property at least not obviously satisfied in the case of birational rowmotion.)

The first step towards our proof of Theorem \ref{thm.rect.ord} is to restrict
attention to so-called \textit{reduced labellings}. Let us define these first:

\begin{definition}
Let $\operatorname*{Rect}\left(  p,q\right)  $ denote the $p\times
q$-rectangle. Let $\mathbb{K}$ be a field. A labelling $f\in\mathbb{K}%
^{\widehat{\operatorname*{Rect}\left(  p,q\right)  }}$ is said to be
\textit{reduced} if $f\left(  0\right)  =f\left(  1\right)  =1$. The set of
all reduced labellings $f\in\mathbb{K}^{\widehat{\operatorname*{Rect}\left(
p,q\right)  }}$ will be identified with $\mathbb{K}^{\operatorname*{Rect}%
\left(  p,q\right)  }$ in the obvious way.

Note that fixing the values of $f\left(  0\right)  $ and $f\left(  1\right)  $
like this makes $f$ ``less generic'', but still the operator
$R_{\operatorname*{Rect}\left(  p,q\right)  }$ restricts to a rational map
from the variety of all reduced labellings $f\in\mathbb{K}%
^{\widehat{\operatorname*{Rect}\left(  p,q\right)  }}$ to itself. (This is
because the operator $R_{\operatorname*{Rect}\left(  p,q\right)  }$ does not
change the values at $0$ and $1$, and does not degenerate from setting
$f\left(  0\right)  =f\left(  1\right)  =1$.)
\end{definition}

Reduced labellings are not much less general than arbitrary labellings: In
fact, every zero-free $\mathbb{K}$-labelling $f$ of a graded poset $P$ is
homogeneously equivalent to a reduced labelling. Thus, many results can be
proven for all labellings by means of proving them for reduced labellings
first, and then extending them to general labellings by means of homogeneous
equivalence.\footnote{A slightly different way to reduce the case of a general
labelling to that of a reduced one is taken in \cite[\S 4]{einstein-propp}.}
We will use this tactic in our proof of Theorem \ref{thm.rect.ord}. Here is
how this works:

\begin{proposition}
\label{prop.rect.reduce}Let $\operatorname*{Rect}\left(  p,q\right)  $ denote
the $p\times q$-rectangle. Let $\mathbb{K}$ be a field. Assume that almost
every (in the Zariski sense) reduced labelling $f\in\mathbb{K}%
^{\widehat{\operatorname*{Rect}\left(  p,q\right)  }}$ satisfies
$R_{\operatorname*{Rect}\left(  p,q\right)  }^{p+q}f=f$. Then,
$\operatorname*{ord}\left(  R_{\operatorname*{Rect}\left(  p,q\right)
}\right)  =p+q$.
\end{proposition}

\begin{vershort}
\begin{proof}
[Proof of Proposition \ref{prop.rect.reduce} (sketched).]Let $g\in
\mathbb{K}^{\widehat{\operatorname*{Rect}\left(  p,q\right)  }}$ be any
$\mathbb{K}$-labelling of $\operatorname*{Rect}\left(  p,q\right)  $ which is
sufficiently generic for $R_{\operatorname*{Rect}\left(  p,q\right)  }^{p+q}g$
to be well-defined.

We use the notation of Definition \ref{def.bemol}. Recall that
$\operatorname*{Rect}\left(  p,q\right)  $ is a $\left(  p+q-1\right)
$-graded poset. We can easily find a $\left(  p+q+1\right)  $-tuple $\left(
a_{0},a_{1},...,a_{p+q}\right)  \in\left(  \mathbb{K}^{\times}\right)
^{p+q+1}$ such that $\left(  a_{0},a_{1},...,a_{p+q}\right)  \flat g$ is a
reduced $\mathbb{K}$-labelling (in fact, set $a_{0}=\dfrac{1}{g\left(
0\right)  }$ and $a_{p+q}=\dfrac{1}{g\left(  1\right)  }$, and choose all
other $a_{i}$ arbitrarily). Corollary \ref{cor.Rl.scalmult} (applied to
$p+q-1$, $\operatorname*{Rect}\left(  p,q\right)  $ and $g$ instead of $n$,
$P$ and $f$) then yields%
\begin{equation}
R_{\operatorname*{Rect}\left(  p,q\right)  }^{p+q}\left(  \left(  a_{0}%
,a_{1},...,a_{p+q}\right)  \flat g\right)  =\left(  a_{0},a_{1},...,a_{p+q}%
\right)  \flat\left(  R_{\operatorname*{Rect}\left(  p,q\right)  }%
^{p+q}g\right)  . \label{pf.rect.reduce.short.2}%
\end{equation}

We have assumed that almost every (in the Zariski sense) reduced labelling
$f\in\mathbb{K}^{\widehat{\operatorname*{Rect}\left(  p,q\right)  }}$
satisfies $R_{\operatorname*{Rect}\left(  p,q\right)  }^{p+q}f=f$. Thus, every
reduced labelling $f\in\mathbb{K}^{\widehat{\operatorname*{Rect}\left(
p,q\right)  }}$ for which $R_{\operatorname*{Rect}\left(  p,q\right)  }%
^{p+q}f$ is well-defined satisfies $R_{\operatorname*{Rect}\left(  p,q\right)
}^{p+q}f=f$ (because $R_{\operatorname*{Rect}\left(  p,q\right)  }^{p+q}f=f$
can be written as an equality between rational functions in the labels of $f$,
and thus it must hold everywhere if it holds on a Zariski-dense open subset).
Applying this to $f=\left(  a_{0},a_{1},...,a_{p+q}\right)  \flat g$, we
obtain that $R_{\operatorname*{Rect}\left(  p,q\right)  }^{p+q}\left(  \left(
a_{0},a_{1},...,a_{p+q}\right)  \flat g\right)  =\left(  a_{0},a_{1}%
,...,a_{p+q}\right)  \flat g$. Thus,%
\begin{align}
\left(  a_{0},a_{1},...,a_{p+q}\right)  \flat g  &  =R_{\operatorname*{Rect}%
\left(  p,q\right)  }^{p+q}\left(  \left(  a_{0},a_{1},...,a_{p+q}\right)
\flat g\right) \nonumber\\
&  =\left(  a_{0},a_{1},...,a_{p+q}\right)  \flat\left(
R_{\operatorname*{Rect}\left(  p,q\right)  }^{p+q}g\right)
\ \ \ \ \ \ \ \ \ \ \left(  \text{by (\ref{pf.rect.reduce.short.2})}\right)  .
\label{pf.rect.reduce.short.3}%
\end{align}
We can cancel the \textquotedblleft$\left(  a_{0},a_{1},...,a_{p+q}\right)
\flat$\textquotedblright\ from both sides of this equality (because all the
$a_{i}$ are nonzero), and thus obtain $g=R_{\operatorname*{Rect}\left(
p,q\right)  }^{p+q}g$.

Now, forget that we fixed $g$. We thus have proven that
$g=R_{\operatorname*{Rect}\left(  p,q\right)  }^{p+q}g$ holds for every
$\mathbb{K}$-labelling $g\in\mathbb{K}^{\widehat{\operatorname*{Rect}\left(
p,q\right)  }}$ of $\operatorname*{Rect}\left(  p,q\right)  $ which is
sufficiently generic for $R_{\operatorname*{Rect}\left(  p,q\right)  }^{p+q}g$
to be well-defined. In other words, $R_{\operatorname*{Rect}\left(
p,q\right)  }^{p+q}=\operatorname*{id}$ as partial maps. Hence,
$\operatorname*{ord}\left(  R_{\operatorname*{Rect}\left(  p,q\right)
}\right)  \mid p+q$.

On the other hand, Proposition \ref{prop.ord-projord} (applied to
$P=\operatorname*{Rect}\left(  p,q\right)  $ and $n=p+q-1$) yields
$\operatorname*{ord}\left(  R_{\operatorname*{Rect}\left(  p,q\right)
}\right)  =\operatorname{lcm}\left(  \left(  p+q-1\right)
+1,\operatorname*{ord}\left(  \overline{R}_{\operatorname*{Rect}\left(
p,q\right)  }\right)  \right)  $. Hence, $\operatorname*{ord}\left(
R_{\operatorname*{Rect}\left(  p,q\right)  }\right)  $ is divisible by
$\left(  p+q-1\right)  +1=p+q$. Combined with $\operatorname*{ord}\left(
R_{\operatorname*{Rect}\left(  p,q\right)  }\right)  \mid p+q$, this yields
$\operatorname*{ord}\left(  R_{\operatorname*{Rect}\left(  p,q\right)
}\right)  =p+q$. This proves Proposition \ref{prop.rect.reduce}.
\end{proof}
\end{vershort}

\begin{verlong}
\begin{proof}
[First proof of Proposition \ref{prop.rect.reduce} (sketched).]Let
$g\in\mathbb{K}^{\widehat{\operatorname*{Rect}\left(  p,q\right)  }}$ be any
generic labelling. Then, there exists a generic reduced labelling
$f\in\mathbb{K}^{\widehat{\operatorname*{Rect}\left(  p,q\right)  }}$ such
that $f$ is homogeneously equivalent to $g$ (namely, one can obtain such an
$f$ by simply replacing the labels of $g$ at $0$ and $1$ by $1$, while leaving
all the other labels unchanged). For this $f$, we know that
$R_{\operatorname*{Rect}\left(  p,q\right)  }^{p+q}f=f$ (by assumption, since
$f$ is reduced). Since $f$ is homogeneously equivalent to $g$, we have
$\pi\left(  f\right)  =\pi\left(  g\right)  $.

Now, let $\pi$ be the canonical projection map $\mathbb{K}^{\widehat{P}%
}\dashrightarrow\overline{\mathbb{K}^{\widehat{P}}}$ as defined in Definition
\ref{def.hgeq}. Since the diagram (\ref{def.hgR.commut}) commutes, we have
$\overline{R}_{\operatorname*{Rect}\left(  p,q\right)  }^{p+q}\left(
\pi\left(  f\right)  \right)  =\pi\left(  \underbrace{R_{\operatorname*{Rect}%
\left(  p,q\right)  }^{p+q}f}_{=f}\right)  =\pi\left(  f\right)  $. But since
$\pi\left(  f\right)  =\pi\left(  g\right) $, this becomes $\overline
{R}_{\operatorname*{Rect}\left(  p,q\right)  }^{p+q}\left(  \pi\left(
g\right)  \right)  =\pi\left(  g\right)  $.

We thus have shown that every generic labelling $g\in\mathbb{K}%
^{\widehat{\operatorname*{Rect}\left(  p,q\right)  }}$ satisfies $\overline
{R}_{\operatorname*{Rect}\left(  p,q\right)  }^{p+q}\left(  \pi\left(
g\right)  \right)  =\pi\left(  g\right)  $. Since $\pi$ is surjective, this
yields $\overline{R}_{\operatorname*{Rect}\left(  p,q\right)  }^{p+q}%
=\operatorname*{id}$, so that $\operatorname*{ord}\left(  \overline
{R}_{\operatorname*{Rect}\left(  p,q\right)  }\right)  $ is finite and divides
$p+q$.

Proposition \ref{prop.ord-projord} (applied to $P=\operatorname*{Rect}\left(
p,q\right)  $ and $n=p+q-1$) yields
\begin{align*}
\operatorname*{ord}\left(  R_{\operatorname*{Rect}\left(  p,q\right)
}\right)   &  =\operatorname{lcm}\left(  \left(  p+q-1\right)
+1,\operatorname*{ord}\left(  \overline{R}_{\operatorname*{Rect}\left(
p,q\right)  }\right)  \right) \\
&  =\operatorname{lcm}\left(  p+q,\operatorname*{ord}\left(  \overline
{R}_{\operatorname*{Rect}\left(  p,q\right)  }\right)  \right)  =p+q
\end{align*}
(since $\operatorname*{ord}\left(  \overline{R}_{\operatorname*{Rect}\left(
p,q\right)  }\right)  $ divides $p+q$). This proves Proposition
\ref{prop.rect.reduce}.
\end{proof}

\begin{proof}
[Second proof of Proposition \ref{prop.rect.reduce} (sketched).]Let
$g\in\mathbb{K}^{\widehat{\operatorname*{Rect}\left(  p,q\right)  }}$ be any
$\mathbb{K}$-labelling of $\operatorname*{Rect}\left(  p,q\right)  $ which is
sufficiently generic for $R_{\operatorname*{Rect}\left(  p,q\right)  }^{p+q}g$
to be well-defined.

For any $n\in\mathbb{N}$, for any $n$-graded poset $P$, for any $\mathbb{K}%
$-labelling $f\in\mathbb{K}^{\widehat{P}}$ and for any $\left(  n+2\right)
$-tuple $\left(  a_{0},a_{1},...,a_{n+1}\right)  \in\left(  \mathbb{K}%
^{\times}\right)  ^{n+2}$, we define a $\mathbb{K}$-labelling $\left(
a_{0},a_{1},...,a_{n+1}\right)  \flat f\in\mathbb{K}^{\widehat{P}}$ as in
Definition \ref{def.bemol}. Since $\operatorname*{Rect}\left(  p,q\right)  $
is a $\left(  p+q-1\right)  $-graded poset, we therefore can define a
$\mathbb{K}$-labelling $\left(  a_{0},a_{1},...,a_{p+q}\right)  \flat
g\in\mathbb{K}^{\widehat{\operatorname*{Rect}\left(  p,q\right)  }}$ for every
$\left(  a_{0},a_{1},...,a_{p+q}\right)  \in\left(  \mathbb{K}^{\times
}\right)  ^{p+q+1}$.

Now, define an $\left(  a_{0},a_{1},...,a_{p+q}\right)  \in\left(
\mathbb{K}^{\times}\right)  ^{p+q+1}$ by
\[
\left(  a_{i}=\left\{
\begin{array}
[c]{c}%
\dfrac{1}{g\left(  0\right)  },\ \ \ \ \ \ \ \ \ \ \text{if }i=0;\\
\dfrac{1}{g\left(  1\right)  },\ \ \ \ \ \ \ \ \ \ \text{if }i=p+q;\\
1,\ \ \ \ \ \ \ \ \ \ \text{if }i\in\left\{  1,2,...,p+q-1\right\}
\end{array}
\right.  \ \ \ \ \ \ \ \ \ \ \text{for every }i\in\left\{
0,1,...,p+q\right\}  \right)  .
\]
Then, $\left(  a_{0},a_{1},...,a_{p+q}\right)  \flat g$ is a reduced
$\mathbb{K}$-labelling\footnote{\textit{Proof.} By the definition of $a_{0}$,
we have $a_{0}=\left\{
\begin{array}
[c]{c}%
\dfrac{1}{g\left(  0\right)  },\ \ \ \ \ \ \ \ \ \ \text{if }0=0;\\
\dfrac{1}{g\left(  1\right)  },\ \ \ \ \ \ \ \ \ \ \text{if }0=p+q;\\
1,\ \ \ \ \ \ \ \ \ \ \text{if }0\in\left\{  1,2,...,p+q-1\right\}
\end{array}
\right.  =\dfrac{1}{g\left(  0\right)  }$ (since $0=0$). By the definition of
$\left(  a_{0},a_{1},...,a_{p+q}\right)  \flat g$, we have%
\begin{align*}
\left(  \left(  a_{0},a_{1},...,a_{p+q}\right)  \flat g\right)  \left(
0\right)   &  =a_{\deg0}\cdot g\left(  0\right)  =\underbrace{a_{0}}%
_{=\dfrac{1}{g\left(  0\right)  }}\cdot g\left(  0\right)
\ \ \ \ \ \ \ \ \ \ \left(  \text{since }\deg0=0\right) \\
&  =\dfrac{1}{g\left(  0\right)  }\cdot g\left(  0\right)  =1.
\end{align*}
\par
By the definition of $a_{p+q}$, we have $a_{p+q}=\left\{
\begin{array}
[c]{c}%
\dfrac{1}{g\left(  0\right)  },\ \ \ \ \ \ \ \ \ \ \text{if }p+q=0;\\
\dfrac{1}{g\left(  1\right)  },\ \ \ \ \ \ \ \ \ \ \text{if }p+q=p+q;\\
1,\ \ \ \ \ \ \ \ \ \ \text{if }p+q\in\left\{  1,2,...,p+q-1\right\}
\end{array}
\right.  =\dfrac{1}{g\left(  1\right)  }$ (since $p+q=p+q$). By the definition
of $\left(  a_{0},a_{1},...,a_{p+q}\right)  \flat g$, we have%
\begin{align*}
\left(  \left(  a_{0},a_{1},...,a_{p+q}\right)  \flat g\right)  \left(
1\right)   &  =a_{\deg1}\cdot g\left(  1\right)  =\underbrace{a_{p+q}%
}_{=\dfrac{1}{g\left(  1\right)  }}\cdot g\left(  1\right)
\ \ \ \ \ \ \ \ \ \ \left(  \text{since }\deg1=p+q\right) \\
&  =\dfrac{1}{g\left(  1\right)  }\cdot g\left(  1\right)  =1.
\end{align*}
Thus, $\left(  \left(  a_{0},a_{1},...,a_{p+q}\right)  \flat g\right)  \left(
0\right)  =1$ and $\left(  \left(  a_{0},a_{1},...,a_{p+q}\right)  \flat
g\right)  \left(  1\right)  =1$. In other words, $\left(  \left(  a_{0}%
,a_{1},...,a_{p+q}\right)  \flat g\right)  \left(  0\right)  =\left(  \left(
a_{0},a_{1},...,a_{p+q}\right)  \flat g\right)  \left(  1\right)  =1$. In
other words, $\left(  a_{0},a_{1},...,a_{p+q}\right)  \flat g$ is a reduced
$\mathbb{K}$-labelling, qed.}. Corollary \ref{cor.Rl.scalmult} (applied to
$p+q-1$, $\operatorname*{Rect}\left(  p,q\right)  $ and $g$ instead of $n$,
$P$ and $f$) yields%
\[
R_{\operatorname*{Rect}\left(  p,q\right)  }^{\left(  p+q-1\right)  +1}\left(
\left(  a_{0},a_{1},...,a_{\left(  p+q-1\right)  +1}\right)  \flat g\right)
=\left(  a_{0},a_{1},...,a_{\left(  p+q-1\right)  +1}\right)  \flat\left(
R_{\operatorname*{Rect}\left(  p,q\right)  }^{\left(  p+q-1\right)
+1}g\right)  .
\]
Since $\left(  p+q-1\right)  +1=p+q$, this simplifies to%
\begin{equation}
R_{\operatorname*{Rect}\left(  p,q\right)  }^{p+q}\left(  \left(  a_{0}%
,a_{1},...,a_{p+q}\right)  \flat g\right)  =\left(  a_{0},a_{1},...,a_{p+q}%
\right)  \flat\left(  R_{\operatorname*{Rect}\left(  p,q\right)  }%
^{p+q}g\right)  . \label{pf.rect.reduce.2}%
\end{equation}
In particular, $R_{\operatorname*{Rect}\left(  p,q\right)  }^{p+q}\left(
\left(  a_{0},a_{1},...,a_{p+q}\right)  \flat g\right)  $ is well-defined.

We have assumed that almost every (in the Zariski sense) reduced labelling
$f\in\mathbb{K}^{\widehat{\operatorname*{Rect}\left(  p,q\right)  }}$
satisfies $R_{\operatorname*{Rect}\left(  p,q\right)  }^{p+q}f=f$. Thus, every
reduced labelling $f\in\mathbb{K}^{\widehat{\operatorname*{Rect}\left(
p,q\right)  }}$ for which $R_{\operatorname*{Rect}\left(  p,q\right)  }%
^{p+q}f$ is well-defined satisfies $R_{\operatorname*{Rect}\left(  p,q\right)
}^{p+q}f=f$ (because $R_{\operatorname*{Rect}\left(  p,q\right)  }^{p+q}f=f$
can be written as an equality between rational functions in the labels of $f$,
and thus it must hold everywhere if it holds on a Zariski-dense open subset).
Applying this to $f=\left(  a_{0},a_{1},...,a_{p+q}\right)  \flat g$, we
obtain that $R_{\operatorname*{Rect}\left(  p,q\right)  }^{p+q}\left(  \left(
a_{0},a_{1},...,a_{p+q}\right)  \flat g\right)  =\left(  a_{0},a_{1}%
,...,a_{p+q}\right)  \flat g$. Thus,%
\begin{align}
\left(  a_{0},a_{1},...,a_{p+q}\right)  \flat g  &  =R_{\operatorname*{Rect}%
\left(  p,q\right)  }^{p+q}\left(  \left(  a_{0},a_{1},...,a_{p+q}\right)
\flat g\right) \nonumber\\
&  =\left(  a_{0},a_{1},...,a_{p+q}\right)  \flat\left(
R_{\operatorname*{Rect}\left(  p,q\right)  }^{p+q}g\right)
\ \ \ \ \ \ \ \ \ \ \left(  \text{by (\ref{pf.rect.reduce.2})}\right)  .
\label{pf.rect.reduce.3}%
\end{align}

Now, let $v\in\widehat{\operatorname*{Rect}\left(  p,q\right)  }$. Then,%
\begin{align*}
&  \underbrace{\left(  \left(  a_{0},a_{1},...,a_{p+q}\right)  \flat g\right)
}_{\substack{=\left(  a_{0},a_{1},...,a_{p+q}\right)  \flat\left(
R_{\operatorname*{Rect}\left(  p,q\right)  }^{p+q}g\right)  \\\text{(by
(\ref{pf.rect.reduce.3}))}}}\left(  v\right) \\
&  =\left(  \left(  a_{0},a_{1},...,a_{p+q}\right)  \flat\left(
R_{\operatorname*{Rect}\left(  p,q\right)  }^{p+q}g\right)  \right)  \left(
v\right) \\
&  =a_{\deg v}\cdot\left(  R_{\operatorname*{Rect}\left(  p,q\right)  }%
^{p+q}g\right)  \left(  v\right)  \ \ \ \ \ \ \ \ \ \ \left(  \text{by the
definition of }\left(  a_{0},a_{1},...,a_{p+q}\right)  \flat\left(
R_{\operatorname*{Rect}\left(  p,q\right)  }^{p+q}g\right)  \right)  .
\end{align*}
In other words,%
\begin{align*}
&  a_{\deg v}\cdot\left(  R_{\operatorname*{Rect}\left(  p,q\right)  }%
^{p+q}g\right)  \left(  v\right) \\
&  =\left(  \left(  a_{0},a_{1},...,a_{p+q}\right)  \flat g\right)  \left(
v\right)  =a_{\deg v}\cdot g\left(  v\right)  \ \ \ \ \ \ \ \ \ \ \left(
\text{by the definition of }\left(  a_{0},a_{1},...,a_{p+q}\right)  \flat
g\right)  .
\end{align*}
We can divide this equality by $a_{\deg v}$ (since $a_{\deg v}\neq0$ (since
$a_{\deg v}\in\mathbb{K}^{\times}$)). As a result, we obtain $\left(
R_{\operatorname*{Rect}\left(  p,q\right)  }^{p+q}g\right)  \left(  v\right)
=g\left(  v\right)  $.

Now, forget that we fixed $v$. We thus have shown that every $v\in
\widehat{\operatorname*{Rect}\left(  p,q\right)  }$ satisfies $\left(
R_{\operatorname*{Rect}\left(  p,q\right)  }^{p+q}g\right)  \left(  v\right)
=g\left(  v\right)  $. In other words, $R_{\operatorname*{Rect}\left(
p,q\right)  }^{p+q}g=g$.

Now, forget that we fixed $g$. We thus have proven that
$R_{\operatorname*{Rect}\left(  p,q\right)  }^{p+q}g=g$ holds for every
$\mathbb{K}$-labelling $g\in\mathbb{K}^{\widehat{\operatorname*{Rect}\left(
p,q\right)  }}$ of $\operatorname*{Rect}\left(  p,q\right)  $ which is
sufficiently generic for $R_{\operatorname*{Rect}\left(  p,q\right)  }^{p+q}g$
to be well-defined. In other words, $R_{\operatorname*{Rect}\left(
p,q\right)  }^{p+q}=\operatorname*{id}$ as partial maps. Hence,
$\operatorname*{ord}\left(  R_{\operatorname*{Rect}\left(  p,q\right)
}\right)  \mid p+q$.

On the other hand, Proposition \ref{prop.ord-projord} (applied to
$P=\operatorname*{Rect}\left(  p,q\right)  $ and $n=p+q-1$) yields
$\operatorname*{ord}\left(  R_{\operatorname*{Rect}\left(  p,q\right)
}\right)  =\operatorname{lcm}\left(  \left(  p+q-1\right)
+1,\operatorname*{ord}\left(  \overline{R}_{\operatorname*{Rect}\left(
p,q\right)  }\right)  \right)  $. Hence, $\operatorname*{ord}\left(
R_{\operatorname*{Rect}\left(  p,q\right)  }\right)  $ is divisible by
$\left(  p+q-1\right)  +1=p+q$. Combined with $\operatorname*{ord}\left(
R_{\operatorname*{Rect}\left(  p,q\right)  }\right)  \mid p+q$, this yields
$\operatorname*{ord}\left(  R_{\operatorname*{Rect}\left(  p,q\right)
}\right)  =p+q$. This proves Proposition \ref{prop.rect.reduce}.
\end{proof}
\end{verlong}

Let us also formulate the particular case of Theorem
\ref{thm.rect.antip.general} for reduced labellings:

\begin{theorem}
\label{thm.rect.antip}Let $\operatorname*{Rect}\left(  p,q\right)  $ denote
the $p\times q$-rectangle. Let $\mathbb{K}$ be a field. Let $f\in
\mathbb{K}^{\widehat{\operatorname*{Rect}\left(  p,q\right)  }}$ be reduced.
Assume that $R_{\operatorname*{Rect}\left(  p,q\right)  }^{\ell}f$ is
well-defined for every $\ell\in\left\{  0,1,...,i+k-1\right\}  $. Let $\left(
i,k\right)  \in\operatorname*{Rect}\left(  p,q\right)  $. Then,%
\[
f\left(  \left(  p+1-i,q+1-k\right)  \right)  =\dfrac{1}{\left(
R_{\operatorname*{Rect}\left(  p,q\right)  }^{i+k-1}f\right)  \left(  \left(
i,k\right)  \right)  }.
\]

\end{theorem}

We will prove this before we prove the general form (Theorem
\ref{thm.rect.antip.general}), and in fact we are going to derive Theorem
\ref{thm.rect.antip.general} from its particular case, Theorem
\ref{thm.rect.antip}. We are not going to encumber this section with the
derivation; its details can be found in Section \ref{sect.rect.finish}.

\section{The Grassmannian parametrization: statements}

In this section, we are going to introduce the main actor in our proof of
Theorem \ref{thm.rect.ord}: an assignment of a reduced $\mathbb{K}$-labelling
of $\operatorname*{Rect}\left(  p,q\right)  $, denoted $\operatorname*{Grasp}%
\nolimits_{j}A$, to any integer $j$ and almost any matrix $A\in\mathbb{K}%
^{p\times\left(  p+q\right)  }$ (Definition \ref{def.Grasp}). This assignment
will give us a family of $\mathbb{K}$-labellings of $\operatorname*{Rect}%
\left(  p,q\right)  $ which is large enough to cover almost all reduced
$\mathbb{K}$-labellings of $\operatorname*{Rect}\left(  p,q\right)  $ (this is
formalized in Proposition \ref{prop.Grasp.generic}), while at the same time
the construction of this assignment makes it easy to track the behavior of the
$\mathbb{K}$-labellings in this family through multiple iterations of
birational rowmotion. Indeed, we will see that birational rowmotion has a very
simple effect on the reduced $\mathbb{K}$-labelling $\operatorname*{Grasp}%
\nolimits_{j}A$ (Proposition \ref{prop.Grasp.GraspR}).

\begin{definition}
\label{def.minors}Let $\mathbb{K}$ be a commutative ring. Let $A\in
\mathbb{K}^{u\times v}$ be a $u\times v$-matrix for some nonnegative integers
$u$ and $v$. (This means, at least in this paper, a matrix with $u$ rows and
$v$ columns.)

\textbf{(a)} For every $i\in\left\{  1,2,...,v\right\}  $, let $A_{i}$ denote
the $i$-th column of $A$.

\textbf{(b)} Moreover, we extend this definition to all $i\in\mathbb{Z}$ as
follows: For every $i\in\mathbb{Z}$, let%
\[
A_{i}=\left(  -1\right)  ^{\left(  u-1\right)  \left(  i-i^{\prime}\right)
\diagup v}\cdot A_{i^{\prime}},
\]
where $i^{\prime}$ is the element of $\left\{  1,2,...,v\right\}  $ which is
congruent to $i$ modulo $v$. (Thus, $A_{v+i}=\left(  -1\right)  ^{u-1}A_{i}$
for every $i\in\mathbb{Z}$. Consequently, the sequence $\left(  A_{i}\right)
_{i\in\mathbb{Z}}$ is periodic with period dividing $2v$, and if $u$ is odd,
the period also divides $v$.)

\textbf{(c)} For any four integers $a$, $b$, $c$ and $d$ satisfying $a\leq b$
and $c\leq d$, we let $A\left[  a:b\mid c:d\right] $ be the matrix whose
columns (from left to right) are $A_{a}$, $A_{a+1}$, $...$, $A_{b-1}$, $A_{c}%
$, $A_{c+1}$, $...$, $A_{d-1}$. (This matrix has $b-a+d-c$
columns.)\footnotemark\ When $b-a+d-c=u$,
this matrix $A\left[  a:b\mid c:d\right]  $ is a square matrix (with $u$ rows
and $u$ columns), and thus has a determinant $\det\left(  A\left[  a:b\mid
c:d\right]  \right) $.

\textbf{(d)} We extend the definition of $\det\left(  A\left[  a:b\mid
c:d\right]  \right)  $ to encompass the case when $b=a-1$ or $d=c-1$, by
setting $\det\left(  A\left[  a:b\mid c:d\right]  \right)  =0$ in this case
(although the matrix $A\left[  a:b\mid c:d\right]  $ itself is not defined in
this case).
\end{definition}

\footnotetext{Some remarks on this matrix $A\left[  a:b\mid c:d\right] $ are
appropriate at this point.
\par
\textbf{1.} We notice that we allow the case $a=b$. In this case, obviously,
the columns of the matrix $A\left[  a:b\mid c:d\right] $ are $A_{c}$,
$A_{c+1}$, $...$, $A_{d-1}$, so the precise value of $a=b$ does not matter.
Similarly, the case $c=d$ is allowed.
\par
\textbf{2.} The matrix $A\left[  a:b\mid c:d\right] $ is not always a
submatrix of $A$. Its columns are columns of $A$ multiplied with $1$ or $-1$;
they can appear several times and need not appear in the same order as they
appear in $A$.}

The reader should be warned that, for $\det\left(  A\left[  a:b\mid
c:d\right]  \right) $ to be defined, we need $b-a+d-c=u$ (not just
$b-a+d-c\equiv u\operatorname{mod}v$, despite the apparent periodicity in the
construction of the matrix $A$.)

\begin{example}
If $A$ is the $2\times3$-matrix $\left(
\begin{array}
[c]{ccc}%
3 & 5 & 7\\
4 & 1 & 9
\end{array}
\right)  $, then Definition \ref{def.minors} \textbf{(b)} yields (for
instance) $A_{5}=\left(  -1\right)  ^{\left(  2-1\right)  \left(  5-2\right)
\diagup3}\cdot A_{2}=-A_{2}=-\left(
\begin{array}
[c]{c}%
5\\
1
\end{array}
\right)  =\left(
\begin{array}
[c]{c}%
-5\\
-1
\end{array}
\right)  $ and $A_{-4}=\left(  -1\right)  ^{\left(  2-1\right)  \left(
\left(  -4\right)  -2\right)  \diagup3}\cdot A_{2}=A_{2}=\left(
\begin{array}
[c]{c}%
5\\
1
\end{array}
\right)  $.

If $A$ is the $3\times2$-matrix $\left(
\begin{array}
[c]{cc}%
1 & 2\\
3 & 2\\
-5 & 4
\end{array}
\right)  $, then Definition \ref{def.minors} \textbf{(b)} yields (for
instance) $A_{0}=\left(  -1\right)  ^{\left(  3-1\right)  \left(  0-2\right)
\diagup2}\cdot A_{2}=A_{2}=\left(
\begin{array}
[c]{c}%
2\\
2\\
4
\end{array}
\right)  $.
\end{example}

\begin{remark}
Some parts of Definition \ref{def.minors} might look accidental and haphazard;
here are some motivations and aide-memoires:

The choice of sign in Definition \ref{def.minors} \textbf{(b)} is not only the
``right'' one for what we are going to do below, but also naturally appears in
\cite[Remark 3.3]{postnikov}. It guarantees, among other things, that if
$A\in\mathbb{R}^{u\times v}$ is totally nonnegative, then the matrix having
columns $A_{1+i}$, $A_{2+i}$, $...$, $A_{v+i}$ is totally nonnegative for
every $i\in\mathbb{Z}$.

The notation $A\left[  a:b\mid c:d\right]  $ in Definition \ref{def.minors}
\textbf{(c)} borrows from Python's notation $\left[  x:y\right]  $ for taking
indices from the interval $\left\{  x,x+1,...,y-1\right\}  $.

The convention to define $\det\left(  A\left[  a:b\mid c:d\right]  \right) $
as $0$ in Definition \ref{def.minors} \textbf{(d)} can be motivated using
exterior algebra as follows: If we identify $\wedge^{u}\left(  \mathbb{K}%
^{u}\right)  $ with $\mathbb{K}$ by equating with $1\in\mathbb{K}$ the wedge
product $e_{1}\wedge e_{2}\wedge...\wedge e_{u}$ of the standard basis
vectors, then $\det\left(  A\left[  a:b\mid c:d\right]  \right)  =A_{a}\wedge
A_{a+1}\wedge...\wedge A_{b-1}\wedge A_{c}\wedge A_{c+1}\wedge...\wedge
A_{d-1}$; this belongs to the product of $\wedge^{b-a}\left(  \mathbb{K}%
^{u}\right)  $ with $\wedge^{d-c}\left(  \mathbb{K}^{u}\right)  $ in
$\wedge^{u}\left(  \mathbb{K}^{u}\right)  $. If $b=a-1$, then this product is
$0$ (since $\wedge^{b-a}\left(  \mathbb{K}^{u}\right)  =\wedge^{-1}\left(
\mathbb{K}^{u}\right)  =0$), so $\det\left(  A\left[  a:b\mid c:d\right]
\right)  $ has to be $0$ in this case.
\end{remark}

Before we go any further, we make several straightforward observations about
the notations we have just introduced.

\begin{proposition}
\label{prop.minors.0}Let $\mathbb{K}$ be a field. Let $A\in\mathbb{K}^{u\times
v}$ be a $u\times v$-matrix for some nonnegative integers $u$ and $v$. Let
$a$, $b$, $c$ and $d$ be four integers satisfying $a\leq b$ and $c\leq d$ and
$b-a+d-c=u$. Assume that some element of the interval $\left\{
a,a+1,...,b-1\right\}  $ is congruent to some element of the interval
$\left\{  c,c+1,...,d-1\right\}  $ modulo $v$. Then, $\det\left(  A\left[
a:b\mid c:d\right]  \right)  =0$.
\end{proposition}

\begin{proof}
[Proof of Proposition \ref{prop.minors.0}.]The assumption yields that the
matrix $A\left[  a:b\mid c:d\right]  $ has two columns which are proportional
to each other by a factor of $\pm1$. Hence, this matrix has determinant $0$.
\end{proof}

\begin{proposition}
\label{prop.minors.antisymm}Let $\mathbb{K}$ be a field. Let $A\in
\mathbb{K}^{u\times v}$ be a $u\times v$-matrix for some nonnegative integers
$u$ and $v$. Let $a$, $b$, $c$ and $d$ be four integers satisfying $a\leq b$
and $c\leq d$ and $b-a+d-c=u$. Then,%
\[
\det\left(  A\left[  a:b\mid c:d\right]  \right)  =\left(  -1\right)
^{\left(  b-a\right)  \left(  d-c\right)  }\det\left(  A\left[  c:d\mid
a:b\right]  \right)  .
\]

\end{proposition}

\begin{proof}
[Proof of Proposition \ref{prop.minors.antisymm}.]This follows from the fact
that permuting the columns of a matrix multiplies its determinant by the sign
of the corresponding permutation.
\end{proof}

\begin{proposition}
\label{prop.minors.complete}Let $\mathbb{K}$ be a field. Let $A\in
\mathbb{K}^{u\times v}$ be a $u\times v$-matrix for some nonnegative integers
$u$ and $v$. Let $a$, $b_{1}$, $b_{2}$ and $c$ be four integers satisfying
$a\leq b_{1}\leq c$ and $a\leq b_{2}\leq c$. Then,%
\[
A\left[  a:b_{1}\mid b_{1}:c\right]  =A\left[  a:b_{2}\mid b_{2}:c\right]  .
\]

\end{proposition}

\begin{proof}
[Proof of Proposition \ref{prop.minors.complete}.]Both matrices $A\left[
a:b_{1}\mid b_{1}:c\right]  $ and $A\left[  a:b_{2}\mid b_{2}:c\right]  $ are
simply the matrix with columns $A_{a}$, $A_{a+1}$, $...$, $A_{c-1}$.
\end{proof}

\begin{proposition}
\label{prop.minors.trivial}Let $\mathbb{K}$ be a field. Let $A\in
\mathbb{K}^{u\times v}$ be a $u\times v$-matrix for some nonnegative integers
$u$ and $v$. Let $c$ and $d$ be two integers satisfying $c\leq d$. Then:

\textbf{(a)} Any integers $a_{1}$ and $a_{2}$ satisfy%
\[
A\left[  a_{1}:a_{1}\mid c:d\right]  =A\left[  a_{2}:a_{2}\mid c:d\right]  .
\]

\textbf{(b)} Any integers $a_{1}$ and $a_{2}$ satisfy%
\[
A\left[  c:d\mid a_{1}:a_{1}\right]  =A\left[  c:d\mid a_{2}:a_{2}\right]  .
\]

\textbf{(c)} If $a$ and $b$ are integers satisfying $c\leq b\leq d$, then%
\[
A\left[  c:b\mid b:d\right]  =A\left[  c:d\mid a:a\right]  .
\]

\end{proposition}

\begin{vershort}
\begin{proof}
[Proof of Proposition \ref{prop.minors.trivial}.]All six matrices in the above
equalities are simply the matrix with columns $A_{c}$, $A_{c+1}$, $...$,
$A_{d-1}$.
\end{proof}
\end{vershort}

\begin{verlong}
\begin{proof}
[Proof of Proposition \ref{prop.minors.trivial}.]\textbf{(a)} Both matrices
$A\left[  a_{1}:a_{1}\mid c:d\right]  $ and $A\left[  a_{2}:a_{2}\mid
c:d\right]  $ are simply the matrix with columns $A_{c}$, $A_{c+1}$, $...$,
$A_{d-1}$.

\textbf{(b)} Analogous.

\textbf{(c)} Similar.
\end{proof}
\end{verlong}

\begin{proposition}
\label{prop.minors.period}Let $\mathbb{K}$ be a field. Let $A\in
\mathbb{K}^{u\times v}$ be a $u\times v$-matrix for some nonnegative integers
$u$ and $v$. Let $a$, $b$, $c$ and $d$ be four integers satisfying $a\leq b$
and $c\leq d$ and $b-a+d-c=u$.

\textbf{(a)} We have%
\[
\det\left(  A\left[  v+a:v+b\mid v+c:v+d\right]  \right)  =\det\left(
A\left[  a:b\mid c:d\right]  \right)  .
\]

\textbf{(b)} We have%
\[
\det\left(  A\left[  a:b\mid v+c:v+d\right]  \right)  =\left(  -1\right)
^{\left(  u-1\right)  \left(  d-c\right)  }\det\left(  A\left[  a:b\mid
c:d\right]  \right)  .
\]

\textbf{(c)} We have%
\[
\det\left(  A\left[  a:b\mid v+c:v+d\right]  \right)  =\det\left(  A\left[
c:d\mid a:b\right]  \right)  .
\]

\end{proposition}

\begin{vershort}
\begin{proof}
[Proof of Proposition \ref{prop.minors.period} (sketched).]Nothing about this
is anything more than trivial. Part \textbf{(a)} and \textbf{(b)} follow from
the fact that $A_{v+i}=\left(  -1\right)  ^{u-1}A_{i}$ for every
$i\in\mathbb{Z}$ (which is owed to Definition \ref{def.minors} \textbf{(b)})
and the multilinearity of the determinant. The proof of part \textbf{(c)}
additionally uses Proposition \ref{prop.minors.antisymm} and a careful sign
computation (notice that $\left(  -1\right)  ^{\left(  d-c-1\right)  \left(
d-c\right)  }=1$ because $\left(  d-c-1\right)  \left(  d-c\right)  $ is even,
no matter what the parities of $c$ and $d$ are). All details can be easily
filled in by the reader.
\end{proof}
\end{vershort}

\begin{verlong}
\begin{proof}
[Proof of Proposition \ref{prop.minors.period} (sketched).]\textbf{(a)} By the
definition of \newline$A\left[  v+a:v+b\mid v+c:v+d\right]  $, we know that
$A\left[  v+a:v+b\mid v+c:v+d\right]  $ is the matrix whose columns (from left
to right) are $A_{v+a}$, $A_{v+a+1}$, $...$, $A_{v+b-1}$, $A_{v+c}$,
$A_{v+c+1}$, $...$, $A_{v+d-1}$. Meanwhile, $A\left[  a:b\mid c:d\right]  $ is
the matrix whose columns (from left to right) are $A_{a}$, $A_{a+1}$, $...$,
$A_{b-1}$, $A_{c}$, $A_{c+1}$, $...$, $A_{d-1}$. Hence, each column of the
matrix $A\left[  v+a:v+b\mid v+c:v+d\right]  $ equals $\left(  -1\right)
^{u-1}$ times the corresponding column of the matrix $A\left[  a:b\mid
c:d\right]  $ (because $A_{v+i}=\left(  -1\right)  ^{u-1}A_{i}$ for every
$i\in\mathbb{Z}$ (as follows from Definition \ref{def.minors} \textbf{(b)})).
Consequently,
\[
A\left[  v+a:v+b\mid v+c:v+d\right]  =\left(  -1\right)  ^{u-1}\cdot A\left[
a:b\mid c:d\right]  ,
\]
so that%
\begin{align*}
&  \det\left(  A\left[  v+a:v+b\mid v+c:v+d\right]  \right) \\
&  =\det\left(  \left(  -1\right)  ^{u-1}\cdot A\left[  a:b\mid c:d\right]
\right) \\
&  =\underbrace{\left(  \left(  -1\right)  ^{u-1}\right)  ^{u}}%
_{\substack{=\left(  -1\right)  ^{\left(  u-1\right)  u}=1\\\text{(since
}\left(  u-1\right)  u\text{ is even)}}}\cdot\det\left(  A\left[  a:b\mid
c:d\right]  \right)  =\det\left(  A\left[  a:b\mid c:d\right]  \right)  .
\end{align*}
This proves Proposition \ref{prop.minors.period} \textbf{(a)}.

\textbf{(b)} By the definition of $A\left[  a:b\mid v+c:v+d\right]  $, we know
that $A\left[  a:b\mid v+c:v+d\right]  $ is the matrix whose columns (from
left to right) are $A_{a}$, $A_{a+1}$, $...$, $A_{b-1}$, $A_{v+c}$,
$A_{v+c+1}$, $...$, $A_{v+d-1}$. Since $A_{v+i}=\left(  -1\right)  ^{u-1}%
A_{i}$ for every $i\in\mathbb{Z}$ (by Definition \ref{def.minors}
\textbf{(b)}), this rewrites as follows: $A\left[  a:b\mid v+c:v+d\right]  $
is the matrix whose columns (from left to right) are $A_{a}$, $A_{a+1}$,
$...$, $A_{b-1}$, $\left(  -1\right)  ^{u-1}A_{c}$, $\left(  -1\right)
^{u-1}A_{c+1}$, $...$, $\left(  -1\right)  ^{u-1}A_{d-1}$. Hence, the matrix
$A\left[  a:b\mid v+c:v+d\right]  $ is obtained from the matrix $A\left[
a:b\mid c:d\right]  $ by multiplying each of its last $d-c$ columns with
$\left(  -1\right)  ^{u-1}$ (since $A\left[  a:b\mid c:d\right]  $ is the
matrix whose columns (from left to right) are $A_{a}$, $A_{a+1}$, $...$,
$A_{b-1}$, $A_{c}$, $A_{c+1}$, $...$, $A_{d-1}$). Therefore, by the
multilinearity of the determinant, we have%
\begin{align*}
&  \det\left(  A\left[  a:b\mid v+c:v+d\right]  \right) \\
&  =\underbrace{\left(  \left(  -1\right)  ^{u-1}\right)  ^{d-c}}_{=\left(
-1\right)  ^{\left(  u-1\right)  \left(  d-c\right)  }}\det\left(  A\left[
a:b\mid c:d\right]  \right)  =\left(  -1\right)  ^{\left(  u-1\right)  \left(
d-c\right)  }\det\left(  A\left[  a:b\mid c:d\right]  \right)  .
\end{align*}
This proves Proposition \ref{prop.minors.period} \textbf{(b)}.

\textbf{(c)} By Proposition \ref{prop.minors.period} \textbf{(b)}, we have%
\begin{align*}
&  \det\left(  A\left[  a:b\mid v+c:v+d\right]  \right) \\
&  =\left(  -1\right)  ^{\left(  u-1\right)  \left(  d-c\right)
}\underbrace{\det\left(  A\left[  a:b\mid c:d\right]  \right)  }%
_{\substack{=\left(  -1\right)  ^{\left(  b-a\right)  \left(  d-c\right)
}\det\left(  A\left[  c:d\mid a:b\right]  \right)  \\\text{(by Proposition
\ref{prop.minors.antisymm})}}}\\
&  =\left(  -1\right)  ^{\left(  u-1\right)  \left(  d-c\right)  }\left(
-1\right)  ^{\left(  b-a\right)  \left(  d-c\right)  }\det\left(  A\left[
c:d\mid a:b\right]  \right)  .
\end{align*}
Since
\begin{align*}
&  \left(  -1\right)  ^{\left(  u-1\right)  \left(  d-c\right)  }%
\underbrace{\left(  -1\right)  ^{\left(  b-a\right)  \left(  d-c\right)  }%
}_{=\left(  \left(  -1\right)  ^{-1}\right)  ^{\left(  b-a\right)  \left(
d-c\right)  }=\left(  -1\right)  ^{-\left(  b-a\right)  \left(  d-c\right)  }%
}\\
&  =\left(  -1\right)  ^{\left(  u-1\right)  \left(  d-c\right)  }\left(
-1\right)  ^{-\left(  b-a\right)  \left(  d-c\right)  }=\left(  -1\right)
^{\left(  u-1\right)  \left(  d-c\right)  -\left(  b-a\right)  \left(
d-c\right)  }=\left(  -1\right)  ^{\left(  u-1-\left(  b-a\right)  \right)
\left(  d-c\right)  }\\
&  =\left(  -1\right)  ^{\left(  d-c-1\right)  \left(  d-c\right)
}\ \ \ \ \ \ \ \ \ \ \left(  \text{since }u-1-\left(  b-a\right)  =d-c-1\text{
(because }b-a+d-c=u\text{)}\right) \\
&  =1\ \ \ \ \ \ \ \ \ \ \left(
\begin{array}
[c]{cc}%
\text{since }\left(  -1\right)  ^{\left(  r-1\right)  r}=1\text{ for every
}r\in\mathbb{Z} & \\
\text{ (since }\left(  r-1\right)  r\text{ is even for every }r\in
\mathbb{Z}\text{)} &
\end{array}
\right)  ,
\end{align*}
this simplifies to%
\begin{align*}
&  \det\left(  A\left[  a:b\mid v+c:v+d\right]  \right) \\
&  =\underbrace{\left(  -1\right)  ^{\left(  u-1\right)  \left(  d-c\right)
}\left(  -1\right)  ^{\left(  b-a\right)  \left(  d-c\right)  }}_{=1}%
\det\left(  A\left[  c:d\mid a:b\right]  \right)  =\det\left(  A\left[
c:d\mid a:b\right]  \right)  .
\end{align*}
This proves Proposition \ref{prop.minors.period} \textbf{(c)}.
\end{proof}
\end{verlong}

\begin{definition}
\label{def.Grasp}Let $\mathbb{K}$ be a field. Let $p$ and $q$ be two positive
integers. Let $A\in\mathbb{K}^{p\times\left(  p+q\right)  }$. Let
$j\in\mathbb{Z}$.

\textbf{(a)} We define a map $\operatorname*{Grasp}\nolimits_{j}A\in
\mathbb{K}^{\operatorname*{Rect}\left(  p,q\right)  }$ by%
\begin{align}
\left(  \operatorname*{Grasp}\nolimits_{j}A\right)  \left(  \left(
i,k\right)  \right)   &  =\dfrac{\det\left(  A\left[  j+1:j+i\mid
j+i+k-1:j+p+k\right]  \right)  }{\det\left(  A\left[  j:j+i\mid
j+i+k:j+p+k\right]  \right)  }\label{def.Grasp.def}\\
&  \ \ \ \ \ \ \ \ \ \ \left.  \text{for every }\left(  i,k\right)
\in\operatorname*{Rect}\left(  p,q\right)  =\left\{  1,2,...,p\right\}
\times\left\{  1,2,...,q\right\}  \right. \nonumber
\end{align}
(this is well-defined when the matrix $A$ is sufficiently generic (in the
sense of Zariski topology), since the matrix $A\left[  j:j+i\mid
j+i+k:j+p+k\right]  $ is obtained by picking $p$ distinct columns out of $A$,
some possibly multiplied with $\left(  -1\right)  ^{u-1}$). This map
$\operatorname*{Grasp}\nolimits_{j}A$ will be considered as a reduced
$\mathbb{K}$-labelling of $\operatorname*{Rect}\left(  p,q\right)  $ (since we
are identifying the set of all reduced labellings $f\in\mathbb{K}%
^{\widehat{\operatorname*{Rect}\left(  p,q\right)  }}$ with $\mathbb{K}%
^{\operatorname*{Rect}\left(  p,q\right)  }$).

\textbf{(b)} It will be handy to extend the map $\operatorname*{Grasp}%
\nolimits_{j}A$ to a slightly larger domain by blindly following
(\ref{def.Grasp.def}) (and using Definition \ref{def.minors} \textbf{(d)}),
accepting the fact that outside $\left\{  1,2,...,p\right\}  \times\left\{
1,2,...,q\right\}  $ its values can be \textquotedblleft
infinity\textquotedblright:%
\begin{align*}
\left(  \operatorname*{Grasp}\nolimits_{j}A\right)  \left(  \left(
0,k\right)  \right)   &  =0\ \ \ \ \ \ \ \ \ \ \text{for all }k\in\left\{
1,2,...,q\right\}  ;\\
\left(  \operatorname*{Grasp}\nolimits_{j}A\right)  \left(  \left(
p+1,k\right)  \right)   &  =\infty\ \ \ \ \ \ \ \ \ \ \text{for all }%
k\in\left\{  1,2,...,q\right\}  ;\\
\left(  \operatorname*{Grasp}\nolimits_{j}A\right)  \left(  \left(
i,0\right)  \right)   &  =0\ \ \ \ \ \ \ \ \ \ \text{for all }i\in\left\{
1,2,...,p\right\}  ;\\
\left(  \operatorname*{Grasp}\nolimits_{j}A\right)  \left(  \left(
i,q+1\right)  \right)   &  =\infty\ \ \ \ \ \ \ \ \ \ \text{for all }%
i\in\left\{  1,2,...,p\right\}  .
\end{align*}
(We treat $\infty$ as a symbol with the properties $\dfrac{1}{0} = \infty$ and
$\dfrac{1}{\infty} = 0$.)
\end{definition}

The notation \textquotedblleft$\operatorname*{Grasp}$\textquotedblright%
\ harkens back to \textquotedblleft Grassmannian
parametrization\textquotedblright, as we will later parametrize (generic)
reduced labellings on $\operatorname*{Rect}\left(  p,q\right)  $ by matrices
via this map $\operatorname*{Grasp}\nolimits_{0}$. The reason for the word
\textquotedblleft Grassmannian\textquotedblright\ is that, while we have
defined $\operatorname*{Grasp}\nolimits_{j}$ as a rational map from the matrix
space $\mathbb{K}^{p\times\left(  p+q\right)  }$, it actually is not defined
outside of the Zariski-dense open subset $\mathbb{K}_{\operatorname*{rk}%
=p}^{p\times\left(  p+q\right)  }$ of $\mathbb{K}^{p\times\left(  p+q\right)
}$ formed by all matrices whose rank is $p$, and on that subset $\mathbb{K}%
_{\operatorname*{rk}=p}^{p\times\left(  p+q\right)  }$ it factors through the
quotient of $\mathbb{K}_{\operatorname*{rk}=p}^{p\times\left(  p+q\right)  }$
by the left multiplication action of $\operatorname*{GL}\nolimits_{p}%
\mathbb{K}$ (because it is easy to see that $\operatorname*{Grasp}%
\nolimits_{j}A$ is invariant under row transformations of $A$); this quotient
is a well-known avatar of the Grassmannian.

The formula (\ref{def.Grasp.def}) is inspired by the $Y_{ijk}$ of Volkov's
\cite{volkov}; similar expressions (in a different context) also appear in
\cite[Theorem 4.21]{kirillov-intro}.

\begin{example}
If $p=2$, $q=2$ and $A=\left(
\begin{array}
[c]{cccc}%
a_{11} & a_{12} & a_{13} & a_{14}\\
a_{21} & a_{22} & a_{23} & a_{24}%
\end{array}
\right)  $, then%
\[
\left(  \operatorname*{Grasp}\nolimits_{0}A\right)  \left(  \left(
1,1\right)  \right)  =\dfrac{\det\left(  A\left[  1:1\mid1:3\right]  \right)
}{\det\left(  A\left[  0:1\mid2:3\right]  \right)  }=\dfrac{\det\left(
\begin{array}
[c]{cc}%
a_{11} & a_{12}\\
a_{21} & a_{22}%
\end{array}
\right)  }{\det\left(
\begin{array}
[c]{cc}%
-a_{14} & a_{12}\\
-a_{24} & a_{22}%
\end{array}
\right)  }=\dfrac{a_{11}a_{22}-a_{12}a_{21}}{a_{12}a_{24}-a_{14}a_{22}}%
\]
and%
\[
\left(  \operatorname*{Grasp}\nolimits_{1}A\right)  \left(  \left(
1,2\right)  \right)  =\dfrac{\det\left(  A\left[  2:2\mid3:5\right]  \right)
}{\det\left(  A\left[  1:2\mid4:5\right]  \right)  }=\dfrac{\det\left(
\begin{array}
[c]{cc}%
a_{13} & a_{14}\\
a_{23} & a_{24}%
\end{array}
\right)  }{\det\left(
\begin{array}
[c]{cc}%
a_{11} & a_{14}\\
a_{21} & a_{24}%
\end{array}
\right)  }=\dfrac{a_{13}a_{24}-a_{14}a_{23}}{a_{11}a_{24}-a_{14}a_{21}}.
\]

\end{example}

We will see more examples of values of $\operatorname*{Grasp}\nolimits_{0}A$
in Example \ref{ex.Grasp.generic}.

\begin{proposition}
\label{prop.Grasp.period}Let $\mathbb{K}$ be a field. Let $p$ and $q$ be two
positive integers. Let $A\in\mathbb{K}^{p\times\left(  p+q\right)  }$ be a
matrix. Then, $\operatorname*{Grasp}\nolimits_{j}A=\operatorname*{Grasp}%
\nolimits_{p+q+j}A$ for every $j\in\mathbb{Z}$ (provided that $A$ is
sufficiently generic in the sense of Zariski topology for
$\operatorname*{Grasp}\nolimits_{j} A$ to be well-defined).
\end{proposition}

\begin{proof}
[Proof of Proposition \ref{prop.Grasp.period} (sketched).]We need to show
that
\[
\left(  \operatorname*{Grasp}\nolimits_{j}A\right)  \left(  \left(
i,k\right)  \right)  =\left(  \operatorname*{Grasp}\nolimits_{p+q+j}A\right)
\left(  \left(  i,k\right)  \right)
\]
for every $\left(  i,k\right)  \in\left\{  1,2,...,p\right\}  \times\left\{
1,2,...,q\right\}  $. But we have%
\begin{align*}
&  A\left[  p+q+j:p+q+j+i\mid p+q+j+i+k:p+q+j+p+k\right] \\
&  =A\left[  j:j+i\mid j+i+k:j+p+k\right]
\end{align*}
(by Proposition \ref{prop.minors.period} \textbf{(a)}, applied to $u=p$,
$v=p+q$, $a=j$, $b=j+i$, $c=j+i+k$ and $d=j+p+k$) and%
\begin{align*}
&  A\left[  p+q+j+1:p+q+j+i\mid p+q+j+i+k-1:p+q+j+p+k\right] \\
&  =A\left[  j+1:j+i\mid j+i+k-1:j+p+k\right]
\end{align*}
(by Proposition \ref{prop.minors.period} \textbf{(a)}, applied to $u=p$,
$v=p+q$, $a=j+1$, $b=j+i$, $c=j+i+k-1$ and $d=j+p+k$). Using these equalities,
we immediately obtain $\left(  \operatorname*{Grasp}\nolimits_{j}A\right)
\left(  \left(  i,k\right)  \right)  =\left(  \operatorname*{Grasp}%
\nolimits_{p+q+j}A\right)  \left(  \left(  i,k\right)  \right)  $ from the
definition of $\operatorname*{Grasp}\nolimits_{j}A$. Proposition
\ref{prop.Grasp.period} is proven.
\end{proof}

\begin{proposition}
\label{prop.Grasp.antipode}Let $\mathbb{K}$ be a field. Let $p$ and $q$ be two
positive integers. Let $A\in\mathbb{K}^{p\times\left(  p+q\right)  }$ be a
matrix. Let $\left(  i,k\right)  \in\operatorname*{Rect}\left(  p,q\right)  $
and $j\in\mathbb{Z}$. Then,%
\[
\left(  \operatorname*{Grasp}\nolimits_{j}A\right)  \left(  \left(
i,k\right)  \right)  =\dfrac{1}{\left(  \operatorname*{Grasp}%
\nolimits_{j+i+k-1}A\right)  \left(  \left(  p+1-i,q+1-k\right)  \right)  }%
\]
(provided that $A$ is sufficiently generic in the sense of Zariski topology
for $\left(  \operatorname*{Grasp}\nolimits_{j}A\right)  \left(  \left(
i,k\right)  \right)  $ and $\left(  \operatorname*{Grasp}\nolimits_{j+i+k-1}%
A\right)  \left(  \left(  p+1-i,q+1-k\right)  \right)  $ to be well-defined).
\end{proposition}

\begin{vershort}
\begin{proof}
The proof of Proposition \ref{prop.Grasp.antipode} is completely
straightforward: one merely needs to expand the definitions of $\left(
\operatorname*{Grasp}\nolimits_{j}A\right)  \left(  \left(  i,k\right)
\right)  $ and $\left(  \operatorname*{Grasp}\nolimits_{j+i+k-1}A\right)
\left(  \left(  p+1-i,q+1-k\right)  \right)  $ and to apply Proposition
\ref{prop.minors.period} \textbf{(c)} twice.
\end{proof}
\end{vershort}

\begin{verlong}
The following proof is completely straightforward:

\begin{proof}
[Proof of Proposition \ref{prop.Grasp.antipode} (sketched).]We have $\left(
i,k\right)  \in\operatorname*{Rect}\left(  p,q\right)  =\left\{
1,2,...,p\right\}  \times\left\{  1,2,...,q\right\}  $, so that $i\in\left\{
1,2,...,p\right\}  $ and $k\in\left\{  1,2,...,q\right\}  $.

Let $i^{\prime}=p+1-i$ and $k^{\prime}=q+1-k$. Let $j^{\prime}=j+i+k-1$. We
have
\begin{align*}
j^{\prime}+1  &  =j+i+k\ \ \ \ \ \ \ \ \ \ \left(  \text{since }j^{\prime
}=j+i+k-1\right)  ;\\
\underbrace{j^{\prime}}_{=j+i+k-1}+\underbrace{i^{\prime}}_{=p+1-i}  &
=\left(  j+i+k-1\right)  +\left(  p+1-i\right)  =j+k+p;\\
\underbrace{j^{\prime}+i^{\prime}}_{=j+k+p}+\underbrace{k^{\prime}}%
_{=q+1-k}-1  &  =\left(  j+k+p\right)  +\left(  q+1-k\right)  -1=p+q+j;\\
j^{\prime}+i^{\prime}+k^{\prime}  &  =\underbrace{\left(  j^{\prime}%
+i^{\prime}+k^{\prime}-1\right)  }_{=p+q+j}+1=p+q+j+1;\\
\underbrace{j^{\prime}}_{=j+i+k-1}+p+\underbrace{k^{\prime}}_{=q+1-k}  &
=\left(  j+i+k-1\right)  +p+\left(  q+1-k\right)  =p+q+j+i.
\end{align*}

From Proposition \ref{prop.minors.period} \textbf{(c)} (applied to $u=p$,
$v=p+q$, $a=j+i+k$, $b=j+k+p$, $c=j$ and $d=j+i$), we obtain%
\begin{align}
&  \det\left(  A\left[  j+i+k:j+k+p\mid p+q+j:p+q+j+i\right]  \right)
\nonumber\\
&  =\det\left(  A\left[  j:j+i\mid j+i+k:j+k+p\right]  \right)  .
\label{pf.Grasp.antipode.1}%
\end{align}
From Proposition \ref{prop.minors.period} \textbf{(c)} (applied to $u=p$,
$v=p+q$, $a=j+i+k-1$, $b=j+k+p$, $c=j+1$ and $d=j+i$), we obtain%
\begin{align}
&  \det\left(  A\left[  j+i+k-1:j+k+p\mid p+q+j+1:p+q+j+i\right]  \right)
\nonumber\\
&  =\det\left(  A\left[  j+1:j+i\mid j+i+k-1:j+k+p\right]  \right)  .
\label{pf.Grasp.antipode.2}%
\end{align}

Since $i^{\prime}=p+1-i\in\left\{  1,2,...,p\right\}  $ (because $i\in\left\{
1,2,...,p\right\}  $) and $k^{\prime}=q+1-k\in\left\{  1,2,...,q\right\}  $
(because $k\in\left\{  1,2,...,q\right\}  $), we have $\left(  i^{\prime
},k^{\prime}\right)  \in\left\{  1,2,...,p\right\}  \times\left\{
1,2,...,q\right\}  =\operatorname*{Rect}\left(  p,q\right)  $. Hence, by the
definition of $\left(  \operatorname*{Grasp}\nolimits_{j^{\prime}}A\right)
\left(  \left(  i^{\prime},k^{\prime}\right)  \right)  $, we have%
\begin{align*}
&  \left(  \operatorname*{Grasp}\nolimits_{j^{\prime}}A\right)  \left(
\left(  i^{\prime},k^{\prime}\right)  \right) \\
&  =\dfrac{\det\left(  A\left[  j^{\prime}+1:j^{\prime}+i^{\prime}\mid
j^{\prime}+i^{\prime}+k^{\prime}-1:j^{\prime}+p+k^{\prime}\right]  \right)
}{\det\left(  A\left[  j^{\prime}:j^{\prime}+i^{\prime}\mid j^{\prime
}+i^{\prime}+k^{\prime}:j^{\prime}+p+k^{\prime}\right]  \right)  }\\
&  =\dfrac{\det\left(  A\left[  j+i+k:j+k+p\mid p+q+j:p+q+j+i\right]  \right)
}{\det\left(  A\left[  j+i+k-1:j+k+p\mid p+q+j+1:p+q+j+i\right]  \right)  }\\
&  \ \ \ \ \ \ \ \ \ \ \left(
\begin{array}
[c]{c}%
\text{since }j^{\prime}+1=j+i+k\text{, }j^{\prime}=j+i+k-1\text{, }j^{\prime
}+i^{\prime}=j+k+p\text{,}\\
j^{\prime}+i^{\prime}+k^{\prime}-1=p+q+j\text{, }j^{\prime}+i^{\prime
}+k^{\prime}=p+q+j+1\\
\text{and }j^{\prime}+p+k^{\prime}=p+q+j+i
\end{array}
\right) \\
&  =\dfrac{\det\left(  A\left[  j:j+i\mid j+i+k:j+k+p\right]  \right)  }%
{\det\left(  A\left[  j+1:j+i\mid j+i+k-1:j+k+p\right]  \right)
}\ \ \ \ \ \ \ \ \ \ \left(  \text{by (\ref{pf.Grasp.antipode.1}) and
(\ref{pf.Grasp.antipode.2})}\right) \\
&  =\dfrac{\det\left(  A\left[  j:j+i\mid j+i+k:j+p+k\right]  \right)  }%
{\det\left(  A\left[  j+1:j+i\mid j+i+k-1:j+p+k\right]  \right)  }=\dfrac
{1}{\left(  \operatorname*{Grasp}\nolimits_{j}A\right)  \left(  \left(
i,k\right)  \right)  }\\
&  \ \ \ \ \ \ \ \ \ \ \left(  \text{this is clear by taking the inverse of
(\ref{def.Grasp.def})}\right)  .
\end{align*}
Compared with $\left(  \operatorname*{Grasp}\nolimits_{j^{\prime}}A\right)
\left(  \left(  i^{\prime},k^{\prime}\right)  \right)  =\left(
\operatorname*{Grasp}\nolimits_{j+i+k-1}A\right)  \left(  \left(
p+1-i,q+1-k\right)  \right)  $ (because $i^{\prime}=p+1-i$ and $k^{\prime
}=q+1-k$ and $j^{\prime}=j+i+k-1$), this yields%
\[
\left(  \operatorname*{Grasp}\nolimits_{j+i+k-1}A\right)  \left(  \left(
p+1-i,q+1-k\right)  \right)  =\dfrac{1}{\left(  \operatorname*{Grasp}%
\nolimits_{j}A\right)  \left(  \left(  i,k\right)  \right)  }.
\]
This proves Proposition \ref{prop.Grasp.antipode}.
\end{proof}
\end{verlong}

Now, let us state the two facts which will combine to a proof of Theorem
\ref{thm.rect.ord}:

\begin{proposition}
\label{prop.Grasp.GraspR}Let $\mathbb{K}$ be a field. Let $p$ and $q$ be two
positive integers. Let $A\in\mathbb{K}^{p\times\left(  p+q\right)  }$ be a
matrix. Let $j\in\mathbb{Z}$. Then,%
\[
\operatorname*{Grasp}\nolimits_{j}A=R_{\operatorname*{Rect}\left(  p,q\right)
}\left(  \operatorname*{Grasp}\nolimits_{j+1}A\right)
\]
(provided that $A$ is sufficiently generic in the sense of Zariski topology
for the two sides of this equality to be well-defined).
\end{proposition}

\begin{proposition}
\label{prop.Grasp.generic}Let $\mathbb{K}$ be a field. Let $p$ and $q$ be two
positive integers. For almost every (in the Zariski sense) $f\in
\mathbb{K}^{\operatorname*{Rect}\left(  p,q\right)  }$, there exists a matrix
$A\in\mathbb{K}^{p\times\left(  p+q\right)  }$ satisfying
$f=\operatorname*{Grasp}\nolimits_{0}A$.
\end{proposition}

Once these propositions are proven, Theorems \ref{thm.rect.ord},
\ref{thm.rect.antip} and \ref{thm.rect.antip.general} will be rather easy to
obtain. We delay the details of this until Section \ref{sect.rect.finish}.

\section{The Pl\"{u}cker-Ptolemy relation}

This section is devoted to proving Proposition \ref{prop.Grasp.GraspR}. Before
we proceed to the proof, we will need some fundamental identities concerning
determinants of matrices. Our main tool is the following fact, which we call
the \textit{Pl\"{u}cker-Ptolemy relation}:

\begin{theorem}
\label{thm.pluecker.ptolemy}Let $\mathbb{K}$ be a field. Let $A\in
\mathbb{K}^{u\times v}$ be a $u\times v$-matrix for some nonnegative integers
$u$ and $v$. Let $a$, $b$, $c$ and $d$ be four integers satisfying $a\leq b+1$
and $c\leq d+1$ and $b-a+d-c=u-2$. Then,%
\begin{align*}
&  \det\left(  A\left[  a-1:b\mid c:d+1\right]  \right)  \cdot\det\left(
A\left[  a:b+1\mid c-1:d\right]  \right) \\
&  +\det\left(  A\left[  a:b\mid c-1:d+1\right]  \right)  \cdot\det\left(
A\left[  a-1:b+1\mid c:d\right]  \right) \\
&  =\det\left(  A\left[  a-1:b\mid c-1:d\right]  \right)  \cdot\det\left(
A\left[  a:b+1\mid c:d+1\right]  \right)  .
\end{align*}

\end{theorem}

Notice that the special case of this theorem for $v=u+2$, $a=2$, $b=p$,
$c=p+2$ and $d=p+q$ is the following lemma:

\begin{lemma}
\textit{\label{lem.pluecker.ptolemy}}Let $\mathbb{K}$ be a field. Let
$u\in\mathbb{N}$. Let $B\in\mathbb{K}^{u\times\left(  u+2\right)  }$ be a
$u\times\left(  u+2\right)  $-matrix. Let $p$ and $q$ be two integers $\geq2$
satisfying $p+q=u+2$. Then,%
\begin{align}
&  \det\left(  B\left[  1:p\mid p+2:p+q+1\right]  \right)  \cdot\det\left(
B\left[  2:p+1\mid p+1:p+q\right]  \right) \nonumber\\
&  +\det\left(  B\left[  2:p\mid p+1:p+q+1\right]  \right)  \cdot\det\left(
B\left[  1:p+1\mid p+2:p+q\right]  \right) \nonumber\\
&  =\det\left(  B\left[  1:p\mid p+1:p+q\right]  \right)  \cdot\det\left(
B\left[  2:p+1\mid p+2:p+q+1\right]  \right)  .
\label{lem.pluecker.ptolemy.eq}%
\end{align}

\end{lemma}

\begin{proof}
[Proof of Theorem \ref{thm.pluecker.ptolemy} (sketched).]If $a=b-1$ or
$c=d-1$, then Theorem \ref{thm.pluecker.ptolemy} degenerates to a triviality
(namely, $0+0=0$). Hence, for the rest of this proof, we assume WLOG that
neither $a=b-1$ nor $c=d-1$. Hence, $a\leq b$ and $c\leq d$.

Now, Theorem \ref{thm.pluecker.ptolemy} follows from the Pl\"{u}cker relations
(see, e.g., \cite[(QR)]{kleiman-laksov}) applied to the $u\times\left(
u+2\right)  $-matrix $A\left[  a-1:b+1\mid c-1:d+1\right]  $. But let us show
an alternative proof of Theorem \ref{thm.pluecker.ptolemy} which avoids the
use of the Pl\"{u}cker relations:

Let $p=b-a+2$ and $q=d-c+2$. Then, $p\geq2$, $q\geq2$ and $p+q=u+2$.

Let $B$ be the matrix whose columns (from left to right) are $A_{a-1}$,
$A_{a}$, $...$, $A_{b}$, $A_{c-1}$, $A_{c}$, $...$, $A_{d}$. Then, $B$ is a
$u\times\left(  u+2\right)  $-matrix and satisfies%
\begin{align*}
A\left[  a-1:b\mid c:d+1\right]   &  =B\left[  1:p\mid p+2:p+q+1\right]  ;\\
A\left[  a:b+1\mid c-1:d\right]   &  =B\left[  2:p+1\mid p+1:p+q\right]  ;\\
A\left[  a:b\mid c-1:d+1\right]   &  =B\left[  2:p\mid p+1:p+q+1\right]  ;\\
A\left[  a-1:b+1\mid c:d\right]   &  =B\left[  1:p+1\mid p+2:p+q\right]  ;\\
A\left[  a-1:b\mid c-1:d\right]   &  =B\left[  1:p\mid p+1:p+q\right]  ;\\
A\left[  a:b+1\mid c:d+1\right]   &  =B\left[  2:p+1\mid p+2:p+q+1\right]  .
\end{align*}
Hence, the equality that we have to prove, namely%
\begin{align*}
&  \det\left(  A\left[  a-1:b\mid c:d+1\right]  \right)  \cdot\det\left(
A\left[  a:b+1\mid c-1:d\right]  \right) \\
&  +\det\left(  A\left[  a:b\mid c-1:d+1\right]  \right)  \cdot\det\left(
A\left[  a-1:b+1\mid c:d\right]  \right) \\
&  =\det\left(  A\left[  a-1:b\mid c-1:d\right]  \right)  \cdot\det\left(
A\left[  a:b+1\mid c:d+1\right]  \right)  ,
\end{align*}
rewrites as
\begin{align*}
&  \det\left(  B\left[  1:p\mid p+2:p+q+1\right]  \right)  \cdot\det\left(
B\left[  2:p+1\mid p+1:p+q\right]  \right) \\
&  +\det\left(  B\left[  2:p\mid p+1:p+q+1\right]  \right)  \cdot\det\left(
B\left[  1:p+1\mid p+2:p+q\right]  \right) \\
&  =\det\left(  B\left[  1:p\mid p+1:p+q\right]  \right)  \cdot\det\left(
B\left[  2:p+1\mid p+2:p+q+1\right]  \right)  .
\end{align*}
But this equality follows from Lemma \ref{lem.pluecker.ptolemy}. Hence, in
order to complete the proof of Theorem \ref{thm.pluecker.ptolemy}, we only
need to verify Lemma \ref{lem.pluecker.ptolemy}.
\end{proof}

\begin{vershort}
\begin{proof}
[Proof of Lemma \ref{lem.pluecker.ptolemy} (sketched).]Let $\left(
e_{1},e_{2},...,e_{u}\right)  $ be the standard basis of the $\mathbb{K}%
$-vector space $\mathbb{K}^{u}$.

Let $\alpha$ and $\beta$ be the $\left(  p-1\right)  $-st entries of the
columns $B_{1}$ and $B_{p+q}$ of $B$. Let $\gamma$ and $\delta$ be the $p$-th
entries of the columns $B_{1}$ and $B_{p+q}$ of $B$.

We need to prove (\ref{lem.pluecker.ptolemy.eq}). Since
(\ref{lem.pluecker.ptolemy.eq}) is a polynomial identity in the entries of
$B$, let us WLOG assume that the columns $B_{2}$, $B_{3}$, $...$, $B_{p+q-1}$
of $B$ (these are the middle $u$ among the altogether $u+2=p+q$ columns of
$B$) are linearly independent (since $u$ vectors in $\mathbb{K}^{u}$ in
general position are linearly independent). Then, by applying row
transformations to the matrix $B$, we can transform these columns into the
basis vectors $e_{1}$, $e_{2}$, $...$, $e_{u}$ of $\mathbb{K}^{u}$. Since the
equality (\ref{lem.pluecker.ptolemy.eq}) is preserved under row
transformations of $B$ (indeed, row transformations of $B$ amount to row
transformations of all six matrices appearing in
(\ref{lem.pluecker.ptolemy.eq}), and thus their only effect on the equality
(\ref{lem.pluecker.ptolemy.eq}) is to multiply the six determinants appearing
in (\ref{lem.pluecker.ptolemy.eq}) by certain scalar factors, but these scalar
factors are all equal and thus don't affect the validity of the equality), we
can therefore WLOG assume that the columns $B_{2}$, $B_{3}$, $...$,
$B_{p+q-1}$ of $B$ \textbf{are} the basis vectors $e_{1}$, $e_{2}$, $...$,
$e_{u}$ of $\mathbb{K}^{u}$. The matrix $B$ then looks as follows:%
\[
\left(
\begin{array}
[c]{cccccccccccc}%
\ast & 1 & 0 & \cdots & 0 & 0 & 0 & 0 & 0 & \cdots & 0 & \ast\\
\ast & 0 & 1 & \cdots & 0 & 0 & 0 & 0 & 0 & \cdots & 0 & \ast\\
\vdots & \vdots & \vdots & \ddots & \vdots & \vdots & \vdots & \vdots & \vdots
& \ddots & \vdots & \ast\\
\ast & 0 & 0 & \cdots & 1 & 0 & 0 & 0 & 0 & \cdots & 0 & \ast\\
\ast & 0 & 0 & \cdots & 0 & 1 & 0 & 0 & 0 & \cdots & 0 & \ast\\
\alpha & 0 & 0 & \cdots & 0 & 0 & 1 & 0 & 0 & \cdots & 0 & \beta\\
\gamma & 0 & 0 & \cdots & 0 & 0 & 0 & 1 & 0 & \cdots & 0 & \delta\\
\ast & 0 & 0 & \cdots & 0 & 0 & 0 & 0 & 1 & \cdots & 0 & \ast\\
\vdots & \vdots & \vdots & \ddots & \vdots & \vdots & \vdots & \vdots & \vdots
& \ddots & \vdots & \vdots\\
\ast & 0 & 0 & \cdots & 0 & 0 & 0 & 0 & 0 & \cdots & 1 & \ast
\end{array}
\right)  ,
\]
where asterisks ($\ast$) signify entries which we are not concerned with.

Now, there is a method to simplify the determinant of a matrix if some columns
of this matrix are known to belong to the standard basis $\left(  e_{1}%
,e_{2},...,e_{u}\right)  $. Indeed, such a matrix can first be brought to a
block-triangular form by permuting columns (which affects the determinant by
$\left(  -1\right)  ^{\sigma}$, with $\sigma$ being the sign of the
permutation used), and then its determinant can be evaluated using the fact
that the determinant of a block-triangular matrix is the product of the
determinants of its diagonal blocks. Applying this method to each of the six
matrices appearing in (\ref{lem.pluecker.ptolemy.eq}), we obtain%
\begin{align*}
\det\left(  B\left[  1:p\mid p+2:p+q+1\right]  \right)   &  =\left(
-1\right)  ^{p+q}\left(  \alpha\delta-\beta\gamma\right)  ;\\
\det\left(  B\left[  2:p+1\mid p+1:p+q\right]  \right)   &  =1;\\
\det\left(  B\left[  2:p\mid p+1:p+q+1\right]  \right)   &  =\left(
-1\right)  ^{q-1}\beta;\\
\det\left(  B\left[  1:p+1\mid p+2:p+q\right]  \right)   &  =\left(
-1\right)  ^{p-1}\gamma;\\
\det\left(  B\left[  1:p\mid p+1:p+q\right]  \right)   &  =\left(  -1\right)
^{p-2}\alpha;\\
\det\left(  B\left[  2:p+1\mid p+2:p+q+1\right]  \right)   &  =\left(
-1\right)  ^{q-2}\delta.
\end{align*}
Hence, (\ref{lem.pluecker.ptolemy.eq}) rewrites as%
\[
\left(  -1\right)  ^{p+q}\left(  \alpha\delta-\beta\gamma\right)
\cdot1+\left(  -1\right)  ^{q-1}\beta\cdot\left(  -1\right)  ^{p-1}%
\gamma=\left(  -1\right)  ^{p-2}\alpha\cdot\left(  -1\right)  ^{q-2}\delta.
\]
Upon cancelling the signs, this simplifies to $\left(  \alpha\delta
-\beta\gamma\right)  +\beta\gamma=\alpha\delta$, which is trivially true. Thus
we have proven (\ref{lem.pluecker.ptolemy.eq}). Hence, Lemma
\ref{lem.pluecker.ptolemy} is proven.
\end{proof}

\begin{remark}
Instead of transforming the middle $p+q$ columns of the matrix $B$ to the
standard basis vectors $e_{1}$, $e_{2}$, $...$, $e_{u}$ of $\mathbb{K}^{u}$ as
we did in the proof of Lemma \ref{lem.pluecker.ptolemy}, we could have
transformed the first and last columns of $B$ into the two last standard basis
vectors $e_{u-1}$ and $e_{u}$. The resulting identity would have been
Dodgson's condensation identity (which appears, e.g., in
\cite[\textit{(Alice)}]{zeilberger-twotime}), applied to the matrix formed by
the remaining $u$ columns of $B$ and after some interchange of rows and columns.
\end{remark}
\end{vershort}

\begin{verlong}
We will sketch two proofs of Lemma \ref{lem.pluecker.ptolemy}.

\begin{proof}
[First proof of Lemma \ref{lem.pluecker.ptolemy} (sketched).]Let $\left(
e_{1},e_{2},...,e_{u}\right)  $ be the standard basis of the $\mathbb{K}%
$-vector space $\mathbb{K}^{u}$.

Let $\alpha$ and $\beta$ be the $\left(  p-1\right)  $-st entries of the
columns $B_{1}$ and $B_{p+q}$ of $B$. Let $\gamma$ and $\delta$ be the $p$-th
entries of the columns $B_{1}$ and $B_{p+q}$ of $B$.

We need to prove (\ref{lem.pluecker.ptolemy.eq}). Since
(\ref{lem.pluecker.ptolemy.eq}) is a polynomial identity in the entries of
$B$, we can WLOG assume that the columns $B_{2}$, $B_{3}$, $...$, $B_{p+q-1}$
of $B$ (these are the middle $u$ among the altogether $u+2=p+q$ columns of
$B$) are linearly independent (since $u$ vectors in $\mathbb{K}^{u}$ in
general position are linearly independent). Assume this. Then, by applying row
transformations to the matrix $B$, we can transform these columns into the
basis vectors $e_{1}$, $e_{2}$, $...$, $e_{u}$ of $\mathbb{K}^{u}$. Since the
equality (\ref{lem.pluecker.ptolemy.eq}) is preserved under row
transformations of $B$ (indeed, row transformations of $B$ amount to row
transformations of all six matrices appearing in
(\ref{lem.pluecker.ptolemy.eq}), and thus their only effect on the equality
(\ref{lem.pluecker.ptolemy.eq}) is to multiply the six determinants appearing
in (\ref{lem.pluecker.ptolemy.eq}) by certain scalar factors, but these scalar
factors are all equal and thus don't affect the validity of the equality), we
can therefore WLOG assume that the columns $B_{2}$, $B_{3}$, $...$,
$B_{p+q-1}$ of $B$ \textbf{are} the basis vectors $e_{1}$, $e_{2}$, $...$,
$e_{u}$ of $\mathbb{K}^{u}$. Assume this. The matrix $B$ then looks as
follows:%
\[
\left(
\begin{array}
[c]{cccccccccccc}%
\ast & 1 & 0 & \cdots & 0 & 0 & 0 & 0 & 0 & \cdots & 0 & \ast\\
\ast & 0 & 1 & \cdots & 0 & 0 & 0 & 0 & 0 & \cdots & 0 & \ast\\
\vdots & \vdots & \vdots & \ddots & \vdots & \vdots & \vdots & \vdots & \vdots
& \ddots & \vdots & \ast\\
\ast & 0 & 0 & \cdots & 1 & 0 & 0 & 0 & 0 & \cdots & 0 & \ast\\
\ast & 0 & 0 & \cdots & 0 & 1 & 0 & 0 & 0 & \cdots & 0 & \ast\\
\alpha & 0 & 0 & \cdots & 0 & 0 & 1 & 0 & 0 & \cdots & 0 & \beta\\
\gamma & 0 & 0 & \cdots & 0 & 0 & 0 & 1 & 0 & \cdots & 0 & \delta\\
\ast & 0 & 0 & \cdots & 0 & 0 & 0 & 0 & 1 & \cdots & 0 & \ast\\
\vdots & \vdots & \vdots & \ddots & \vdots & \vdots & \vdots & \vdots & \vdots
& \ddots & \vdots & \vdots\\
\ast & 0 & 0 & \cdots & 0 & 0 & 0 & 0 & 0 & \cdots & 1 & \ast
\end{array}
\right)  ,
\]
where asterisks ($\ast$) signify entries which we are not concerned with.

Now, there is a method to simplify the determinant of a matrix if some columns
of this matrix are known to belong to the standard basis $\left(  e_{1}%
,e_{2},...,e_{u}\right)  $. Indeed, such a matrix can first be brought to a
block-triangular form by permuting columns (which affects the determinant by
$\left(  -1\right)  ^{\sigma}$, with $\sigma$ being the sign of the
permutation used), and then its determinant can be evaluated using the fact
that the determinant of a block-triangular matrix is the product of the
determinants of its diagonal blocks. Applying this method to each of the six
matrices appearing in (\ref{lem.pluecker.ptolemy.eq}), we obtain%
\begin{align*}
\det\left(  B\left[  1:p\mid p+2:p+q+1\right]  \right)   &  =\left(
-1\right)  ^{p+q}\left(  \alpha\delta-\beta\gamma\right)  ;\\
\det\left(  B\left[  2:p+1\mid p+1:p+q\right]  \right)   &  =1;\\
\det\left(  B\left[  2:p\mid p+1:p+q+1\right]  \right)   &  =\left(
-1\right)  ^{q-1}\beta;\\
\det\left(  B\left[  1:p+1\mid p+2:p+q\right]  \right)   &  =\left(
-1\right)  ^{p-1}\gamma;\\
\det\left(  B\left[  1:p\mid p+1:p+q\right]  \right)   &  =\left(  -1\right)
^{p-2}\alpha;\\
\det\left(  B\left[  2:p+1\mid p+2:p+q+1\right]  \right)   &  =\left(
-1\right)  ^{q-2}\delta.
\end{align*}
Hence, (\ref{lem.pluecker.ptolemy.eq}) rewrites as%
\[
\left(  -1\right)  ^{p+q}\left(  \alpha\delta-\beta\gamma\right)
\cdot1+\left(  -1\right)  ^{q-1}\beta\cdot\left(  -1\right)  ^{p-1}%
\gamma=\left(  -1\right)  ^{p-2}\alpha\cdot\left(  -1\right)  ^{q-2}\delta.
\]
Upon cancelling the signs, this simplifies to $\left(  \alpha\delta
-\beta\gamma\right)  +\beta\gamma=\alpha\delta$, which is trivially true. Thus
we have proven (\ref{lem.pluecker.ptolemy.eq}). Hence, Lemma
\ref{lem.pluecker.ptolemy} is proven.
\end{proof}

\begin{remark}
Instead of transforming the middle $p+q$ columns of the matrix $B$ to the
standard basis vectors $e_{1}$, $e_{2}$, $...$, $e_{u}$ of $\mathbb{K}^{u}$ as
we did in the first proof of Lemma \ref{lem.pluecker.ptolemy}, we could have
transformed the first and last columns of $B$ into the two last standard basis
vectors $e_{u-1}$ and $e_{u}$. The resulting identity would have been
Dodgson's condensation identity (which appears, e.g., in
\cite[\textit{(Alice)}]{zeilberger-twotime}), applied to the matrix formed by
the remaining $u$ columns of $B$ and after some interchange of rows and columns.
\end{remark}

Before we begin with the second proof of Lemma \ref{lem.pluecker.ptolemy}, let
us state the Pl\"{u}cker relations (at least in one of their forms):
\end{verlong}

\begin{verlong}
\begin{proposition}
\label{prop.pluecker.pluecker} Let $\mathbb{K}$ be a commutative ring. Let $u$
be a positive integer. Let $A\in\mathbb{K}^{u\times\left(  u-1\right)  }$ be a
$u\times\left(  u-1\right)  $-matrix. Let $w_{1}$, $w_{2}$, $...$, $w_{u+1}$
be $u+1$ column vectors in $\mathbb{K}^{u}$. Then,
\begin{equation}
\sum\limits_{r=1}^{u+1}\left(  -1\right)  ^{r}\det\left(  A\mid w_{r}\right)
\cdot\det\left(  w_{1}\mid w_{2}\mid...\mid\widehat{w_{r}}\mid...\mid
w_{u+1}\right)  =0. \label{prop.pluecker.pluecker.eq}%
\end{equation}
Here, we are using the conventions introduced in Convention \ref{conv.attach}
and Convention \ref{conv.widehat} below, as well as the following convention:
For every $\ell\in\mathbb{N}$, we regard every column vector $v\in
\mathbb{K}^{\ell}$ as an $\ell\times1$-matrix. (For example, the column vector
$w_{r}$ in the expression $\left(  A\mid w_{r}\right)  $ in
(\ref{prop.pluecker.pluecker.eq}) has to be understood as a $u\times1$-matrix;
this ensures that $\left(  A\mid w_{r}\right)  $ in
(\ref{prop.pluecker.pluecker.eq}) is a well-defined $u\times u$-matrix, and
similarly $\left(  w_{1}\mid w_{2}\mid...\mid\widehat{w_{r}}\mid...\mid
w_{u+1}\right)  $ in (\ref{prop.pluecker.pluecker.eq}) is a well-defined
$u\times u$-matrix.)
\end{proposition}

\begin{convention}
\label{conv.widehat} Let $k$ be a nonnegative integer. Let $p_{1}$, $p_{2}$,
$...$, $p_{k}$ be $k$ arbitrary objects. Let $i$ be an element of $\left\{
1,2,...,k\right\}  $. Then, the notation \textquotedblleft$p_{1}$, $p_{2}$,
$...$, $\widehat{p_{i}}$, $...$, $p_{k}$\textquotedblright\ is a shorthand for
\textquotedblleft$p_{1}$, $p_{2}$, $...$, $p_{i-1}$, $p_{i+1}$, $p_{i+2}$,
$...$, $p_{k}$\textquotedblright\ (that is, all elements $p_{1}$, $p_{2}$,
$...$, $p_{k}$ in this order but with the $i$-th element $p_{i}$ excluded). In
other words, when we list the $k$ objects $p_{1}$, $p_{2}$, $...$, $p_{k}$ and
put a hat over the $i$-th element $\widehat{p_{i}}$, the hat signifies that we
are meaning to exclude the $i$-th element from the list.

Consequently, the notation $\left(  p_{1},p_{2},...,\widehat{p_{i}}%
,...,p_{k}\right)  $ denotes the $\left(  k-1\right)  $-tuple $\left(
p_{1},p_{2},...,p_{i-1},p_{i+1},p_{i+2},...,p_{k}\right)  $. Similarly, if
$p_{1}$, $p_{2}$, $...$, $p_{i-1}$, $p_{i+1}$, $p_{i+2}$, $...$, $p_{k}$ are
matrices with the same number of rows, then the notation $\left(  p_{1}\mid
p_{2}\mid...\mid\widehat{p_{i}}\mid...\mid p_{k}\right)  $ denotes the matrix
$\left(  p_{1}\mid p_{2}\mid...\mid p_{i-1}\mid p_{i+1}\mid p_{i+2}\mid...\mid
p_{k}\right)  $.
\end{convention}

(Convention \ref{conv.attach} can be found in Section \ref{sect.dominance}.)

Proposition \ref{prop.pluecker.pluecker} is easily seen to be a consequence of
\cite[(QR)]{kleiman-laksov} in the case when $\mathbb{K}$ is a field and the
vectors $w_{1}$, $w_{2}$, $...$, $w_{u+1}$ span the $\mathbb{K}$-vector space
$\mathbb{K}^{u}$. It is moreover easy to deduce its general validity from the
fact that it holds in this case (by Zariski density, since it is a polynomial
identity in the entries of $A$ and in the coordinates of the vectors $w_{i}$).
However, we will give a proof of Proposition \ref{prop.pluecker.pluecker} here
in order to keep this paper self-contained.\footnote{Update (2022): A nicer
proof of Proposition \ref{prop.pluecker.pluecker} can be found in Theorem
6.150 of \href{https://arxiv.org/abs/2008.09862v2}{Darij Grinberg,
\textit{Notes on the combinatorial fundamentals of algebra}
(arXiv:2008.09862v2)}. (Apply the latter theorem to the $u\times\left(
u+1\right)  $-matrix $B=\left(  w_{1}\mid w_{2}\mid...\mid w_{u+1}\right)
$.)}

Before we prove Proposition \ref{prop.pluecker.pluecker}, let us recall two
known facts:

\begin{proposition}
\label{prop.cofactor.laplace}Let $\mathbb{K}$ be a commutative ring. Let $n$
be a positive integer. Let $B\in\mathbb{K}^{n\times n}$. Then,%
\begin{align*}
\det B  &  =\sum_{\ell=1}^{n}\left(  \text{the }\left(  1,\ell\right)
\text{-th entry of the matrix }B\right) \\
&  \ \ \ \ \ \ \ \ \ \ \ \ \ \ \ \ \ \ \ \ \cdot\left(  \text{the }\left(
1,\ell\right)  \text{-th cofactor of the matrix }B\right)  .
\end{align*}

\end{proposition}

\begin{proof}
[Proof of Proposition \ref{prop.cofactor.laplace} (sketched).]Proposition
\ref{prop.cofactor.laplace} is simply the formula for computing the
determinant $\det B$ of the matrix $B$ using Laplace expansion with respect to
the first row.
\end{proof}

\begin{proposition}
\label{prop.cofactor.laplace.lastcol}Let $\mathbb{K}$ be a commutative ring.
Let $n$ be a positive integer. Let $B\in\mathbb{K}^{n\times n}$. Then,%
\begin{align*}
\det B  &  =\sum_{\ell=1}^{n}\left(  \text{the }\left(  \ell,n\right)
\text{-th entry of the matrix }B\right) \\
&  \ \ \ \ \ \ \ \ \ \ \ \ \ \ \ \ \ \ \ \ \cdot\left(  \text{the }\left(
\ell,n\right)  \text{-th cofactor of the matrix }B\right)  .
\end{align*}

\end{proposition}

\begin{proof}
[Proof of Proposition \ref{prop.cofactor.laplace.lastcol} (sketched).]%
Proposition \ref{prop.cofactor.laplace.lastcol} is simply the formula for
computing the determinant $\det B$ of the matrix $B$ using Laplace expansion
with respect to the $n$-th column.
\end{proof}

\begin{proof}
[Proof of Proposition \ref{prop.pluecker.pluecker}]In the following, we are
going to use Convention \ref{conv.cols} and Convention \ref{conv.rows}. We are
also going to use some results from Section \ref{sect.dominance}; this will
not constitute circular reasoning because these results don't rely on any of
the claims that we are currently proving.

For every $r\in\left\{  1,2,...,u+1\right\}  $ and $\ell\in\left\{
1,2,...,u\right\}  $, we have%
\begin{align}
&  \left(  \text{the }\left(  \ell,u\right)  \text{-th cofactor of the matrix
}\left(  A\mid w_{r}\right)  \right) \nonumber\\
&  =\left(  -1\right)  ^{u-\ell}\cdot\det\left(  \operatorname*{rows}%
\nolimits_{1,2,...,\widehat{\ell},...,u}A\right)  .
\label{pf.pluecker.pluecker.1}%
\end{align}
\footnote{\textit{Proof of (\ref{pf.pluecker.pluecker.1}):} Let $r\in\left\{
1,2,...,u+1\right\}  $ and $\ell\in\left\{  1,2,...,u\right\}  $. By the
definition of a cofactor, it is clear that the $\left(  \ell,u\right)  $-th
cofactor of the matrix $\left(  A\mid w_{r}\right)  $ is $\left(  -1\right)
^{u-\ell}$ times the determinant of the matrix obtained by removing row $\ell$
and column $u$ from the matrix $\left(  A\mid w_{r}\right)  $. In other words,%
\begin{align*}
&  \left(  \text{the }\left(  \ell,u\right)  \text{-th cofactor of the matrix
}\left(  A\mid w_{r}\right)  \right) \\
&  =\left(  -1\right)  ^{u-\ell}\cdot\det\underbrace{\left(  \text{the matrix
obtained by removing row }\ell\text{ and column }u\text{ from the matrix
}\left(  A\mid w_{r}\right)  \right)  }_{=\operatorname*{rows}%
\nolimits_{1,2,...,\widehat{\ell},...,u}\left(  \operatorname*{cols}%
\nolimits_{1,2,...,u-1}\left(  A\mid w_{r}\right)  \right)  }\\
&  =\left(  -1\right)  ^{u-\ell}\cdot\det\left(  \operatorname*{rows}%
\nolimits_{1,2,...,\widehat{\ell},...,u}\left(
\underbrace{\operatorname*{cols}\nolimits_{1,2,...,u-1}\left(  A\mid
w_{r}\right)  }_{\substack{=A\\\text{(by Corollary \ref{cor.cols.attach1a},
applied}\\\text{to }u-1\text{, }1\text{, }A\text{ and }w_{r}\text{ instead of
}v_{1}\text{, }v_{2}\text{, }A_{1}\text{ and }A_{2}\text{)}}}\right)  \right)
\\
&  =\left(  -1\right)  ^{u-\ell}\cdot\det\left(  \operatorname*{rows}%
\nolimits_{1,2,...,\widehat{\ell},...,u}A\right)  .
\end{align*}
This proves (\ref{pf.pluecker.pluecker.1}).}

Let $W\in\mathbb{K}^{u\times\left(  u+1\right)  }$ be the matrix whose columns
(from left to right) are $w_{1}$, $w_{2}$, $...$, $w_{u+1}$. Let
$w_{1}^{\prime}$, $w_{2}^{\prime}$, $...$, $w_{u}^{\prime}$ be the rows of
this matrix $W$ (from top to bottom). For every $\left(  i,j\right)
\in\left\{  1,2,...,u\right\}  \times\left\{  1,2,...,u+1\right\}  $, let
$w_{i,j}$ denote the $\left(  i,j\right)  $-th entry of the matrix $W$.

For every $r\in\left\{  1,2,...,u+1\right\}  $ and $\ell\in\left\{
1,2,...,u\right\}  $, we have%
\begin{equation}
\left(  \text{the }\left(  \ell,u\right)  \text{-th entry of the matrix
}\left(  A\mid w_{r}\right)  \right)  =w_{\ell,r}.
\label{pf.pluecker.pluecker.2}%
\end{equation}
\footnote{\textit{Proof of (\ref{pf.pluecker.pluecker.2}):} Let $r\in\left\{
1,2,...,u+1\right\}  $ and $\ell\in\left\{  1,2,...,u\right\}  $. We have
\begin{equation}
\left(  A\mid w_{r}\right)  _{u}=\left(  \text{the }u\text{-th column of the
matrix }\left(  A\mid w_{r}\right)  \right)
\label{pf.pluecker.pluecker.2.pf.1}%
\end{equation}
(since $\left(  A\mid w_{r}\right)  $ is a $u\times u$-matrix, and since
$1\leq u\leq u$). But the columns of the matrix $W$ (from left to right) are
$w_{1}$, $w_{2}$, $...$, $w_{u+1}$. Hence, $w_{r}$ is the $r$-th column of the
matrix $W$. In other words,%
\begin{equation}
\left(  \text{the column vector }w_{r}\right)  =\left(  \text{the }r\text{-th
column of the matrix }W\right)  . \label{pf.pluecker.pluecker.2.pf.2}%
\end{equation}
\par
On the other hand, the rows of the matrix $W$ (from top to bottom) are
$w_{1}^{\prime}$, $w_{2}^{\prime}$, $...$, $w_{u}^{\prime}$. Hence, the $\ell
$-th row of the matrix $W$ is $w_{\ell}^{\prime}$. In other words,%
\begin{equation}
\left(  \text{the }\ell\text{-th row of the matrix }W\right)  =w_{\ell
}^{\prime}. \label{pf.pluecker.pluecker.2.pf.3}%
\end{equation}
\par
Now,%
\begin{align*}
&  \left(  \text{the }\left(  \ell,u\right)  \text{-th entry of the matrix
}\left(  A\mid w_{r}\right)  \right) \\
&  =\left(  \text{the }\ell\text{-th entry of }\underbrace{\text{the
}u\text{-th column of the matrix }\left(  A\mid w_{r}\right)  }%
_{\substack{=\left(  A\mid w_{r}\right)  _{u}\\\text{(by
(\ref{pf.pluecker.pluecker.2.pf.1}))}}}\right) \\
&  =\left(  \text{the }\ell\text{-th entry of }\left(  A\mid w_{r}\right)
_{u}\right) \\
&  =\left(  \text{the }\ell\text{-th entry of }\left(  A\mid w_{r}\right)
_{\left(  u-1\right)  +1}\right)  \ \ \ \ \ \ \ \ \ \ \left(  \text{since
}u=\left(  u-1\right)  +1\right) \\
&  =\left(  \text{the }\left(  \ell,1\right)  \text{-th entry of the matrix
}w_{r}\right) \\
&  \ \ \ \ \ \ \ \ \ \ \left(
\begin{array}
[c]{c}%
\text{by Proposition \ref{prop.attach2.entry}, applied to }u-1\text{,
}1\text{, }A\text{, }w_{r}\text{, }1\text{ and }\ell\\
\text{instead of }v_{1}\text{, }v_{2}\text{, }A_{1}\text{, }A_{2}\text{, }%
\ell\text{ and }i
\end{array}
\right) \\
&  =\left(  \text{the }\ell\text{-th entry of }\underbrace{\text{the column
vector }w_{r}}_{\substack{=\left(  \text{the }r\text{-th column of the matrix
}W\right)  \\\text{(by (\ref{pf.pluecker.pluecker.2.pf.2}))}}}\right) \\
&  \ \ \ \ \ \ \ \ \ \ \left(  \text{because of how we are identifying column
vectors in }\mathbb{K}^{u}\text{ with }u\times1\text{-matrices}\right) \\
&  =\left(  \text{the }\ell\text{-th entry of the }r\text{-th column of the
matrix }W\right)  =\left(  \text{the }\left(  \ell,r\right)  \text{-th entry
of the matrix }W\right) \\
&  =w_{\ell,r}\ \ \ \ \ \ \ \ \ \ \left(  \text{since }w_{\ell,r}\text{ is the
}\left(  \ell,r\right)  \text{-th entry of the matrix }W\text{ (by the
definition of }w_{\ell,r}\text{)}\right)  .
\end{align*}
\par
This proves (\ref{pf.pluecker.pluecker.2}).}

Now, for every $r\in\left\{  1,2,...,u+1\right\}  $, we have%
\begin{align*}
&  \det\left(  A\mid w_{r}\right) \\
&  =\sum_{\ell=1}^{u}\underbrace{\left(  \text{the }\left(  \ell,u\right)
\text{-th entry of the matrix }\left(  A\mid w_{r}\right)  \right)
}_{\substack{=w_{\ell,r}\\\text{(by (\ref{pf.pluecker.pluecker.2}))}}}\\
&  \ \ \ \ \ \ \ \ \ \ \ \ \ \ \ \ \ \ \ \ \cdot\underbrace{\left(  \text{the
}\left(  \ell,u\right)  \text{-th cofactor of the matrix }\left(  A\mid
w_{r}\right)  \right)  }_{\substack{=\left(  -1\right)  ^{u-\ell}\cdot
\det\left(  \operatorname*{rows}\nolimits_{1,2,...,\widehat{\ell}%
,...,u}A\right)  \\\text{(by (\ref{pf.pluecker.pluecker.1}))}}}\\
&  \ \ \ \ \ \ \ \ \ \ \left(  \text{by Proposition
\ref{prop.cofactor.laplace.lastcol}, applied to }u\text{ and }\left(  A\mid
w_{r}\right)  \text{ instead of }n\text{ and }B\right) \\
&  =\sum_{\ell=1}^{u}w_{\ell,r}\cdot\left(  -1\right)  ^{u-\ell}\cdot
\det\left(  \operatorname*{rows}\nolimits_{1,2,...,\widehat{\ell}%
,...,u}A\right) \\
&  =\sum_{k=1}^{u}w_{k,r}\cdot\left(  -1\right)  ^{u-k}\cdot\det\left(
\operatorname*{rows}\nolimits_{1,2,...,\widehat{k},...,u}A\right)
\end{align*}
(here, we renamed the index $\ell$ as $k$ in the sum). Therefore,%
\begin{align}
&  \sum\limits_{r=1}^{u+1}\left(  -1\right)  ^{r}\underbrace{\det\left(  A\mid
w_{r}\right)  }_{=\sum\limits_{k=1}^{u}w_{k,r}\cdot\left(  -1\right)
^{u-k}\cdot\det\left(  \operatorname*{rows}\nolimits_{1,2,...,\widehat{k}%
,...,u}A\right)  }\cdot\det\left(  w_{1}\mid w_{2}\mid...\mid\widehat{w_{r}%
}\mid...\mid w_{u+1}\right) \nonumber\\
&  =\sum\limits_{r=1}^{u+1}\left(  -1\right)  ^{r}\left(  \sum\limits_{k=1}%
^{u}w_{k,r}\cdot\left(  -1\right)  ^{u-k}\cdot\det\left(  \operatorname*{rows}%
\nolimits_{1,2,...,\widehat{k},...,u}A\right)  \right) \nonumber\\
&  \ \ \ \ \ \ \ \ \ \ \cdot\det\left(  w_{1}\mid w_{2}\mid...\mid
\widehat{w_{r}}\mid...\mid w_{u+1}\right) \nonumber\\
&  =\underbrace{\sum\limits_{r=1}^{u+1}\sum\limits_{k=1}^{u}}_{=\sum
\limits_{k=1}^{u}\sum\limits_{r=1}^{u+1}}\left(  -1\right)  ^{r}w_{k,r}%
\cdot\left(  -1\right)  ^{u-k}\cdot\det\left(  \operatorname*{rows}%
\nolimits_{1,2,...,\widehat{k},...,u}A\right)  \cdot\det\left(  w_{1}\mid
w_{2}\mid...\mid\widehat{w_{r}}\mid...\mid w_{u+1}\right) \nonumber\\
&  =\sum\limits_{k=1}^{u}\sum\limits_{r=1}^{u+1}\left(  -1\right)  ^{r}%
w_{k,r}\cdot\left(  -1\right)  ^{u-k}\cdot\det\left(  \operatorname*{rows}%
\nolimits_{1,2,...,\widehat{k},...,u}A\right)  \cdot\det\left(  w_{1}\mid
w_{2}\mid...\mid\widehat{w_{r}}\mid...\mid w_{u+1}\right) \nonumber\\
&  =\sum\limits_{k=1}^{u}\left(  -1\right)  ^{u-k}\cdot\det\left(
\operatorname*{rows}\nolimits_{1,2,...,\widehat{k},...,u}A\right) \nonumber\\
&  \ \ \ \ \ \ \ \ \ \ \cdot\left(  \sum\limits_{r=1}^{u+1}\left(  -1\right)
^{r}w_{k,r}\cdot\det\left(  w_{1}\mid w_{2}\mid...\mid\widehat{w_{r}}%
\mid...\mid w_{u+1}\right)  \right)  . \label{pf.pluecker.pluecker.5}%
\end{align}

We will now show that every $k\in\left\{  1,2,...,u\right\}  $ satisfies%
\begin{equation}
\sum\limits_{r=1}^{u+1}\left(  -1\right)  ^{r}w_{k,r}\cdot\det\left(
w_{1}\mid w_{2}\mid...\mid\widehat{w_{r}}\mid...\mid w_{u+1}\right)  =0.
\label{pf.pluecker.pluecker.6}%
\end{equation}

\textit{Proof of (\ref{pf.pluecker.pluecker.6}):} Fix $k\in\left\{
1,2,...,u\right\}  $. Let $V\in\mathbb{K}^{\left(  u+1\right)  \times\left(
u+1\right)  }$ be the square matrix whose rows (from top to bottom) are
$w_{k}^{\prime}$, $w_{1}^{\prime}$, $w_{2}^{\prime}$, $...$, $w_{u}^{\prime}$.
Then, the row vector $w_{k}^{\prime}$ appears twice as a row in the matrix $V$
(one time as the topmost row $w_{k}^{\prime}$, and one more time among the
remaining rows $w_{1}^{\prime}$, $w_{2}^{\prime}$, $...$, $w_{u}^{\prime}$).
Hence, there exists a row vector appearing twice as a row in the matrix $V$.
In other words, the matrix $V$ has two equal rows. Since any square matrix
which has two equal rows must have determinant zero, this yields that the
matrix $V$ has determinant zero. In other words, $\det V=0$.

But if $p$ and $q$ are any nonnegative integers, if $B$ is a $p\times
q$-matrix, and if $i_{1}$, $i_{2}$, $...$, $i_{k}$ are some elements of
$\left\{  1,2,...,p\right\}  $, then $\operatorname*{rows}\nolimits_{i_{1}%
,i_{2},...,i_{k}}B$ is the matrix whose rows (from top to bottom) are the rows
labelled $i_{1}$, $i_{2}$, $...$, $i_{k}$ of the matrix $B$%
.\ \ \ \ \footnote{This is how $\operatorname*{rows}\nolimits_{i_{1}%
,i_{2},...,i_{k}}B$ is defined in Convention \ref{conv.rows}.} Applying this
to $p=u+1$, $q=u+1$, $B=V$ and $\left(  i_{1},i_{2},...,i_{k}\right)  =\left(
2,3,...,u+1\right)  $, we conclude that $\operatorname*{rows}%
\nolimits_{2,3,...,u+1}V$ is the matrix whose rows (from top to bottom) are
the rows labelled $2$, $3$, $...$, $u+1$ of the matrix $V$. Since the rows
labelled $2$, $3$, $...$, $u+1$ of the matrix $V$ are $w_{1}^{\prime}$,
$w_{2}^{\prime}$, $...$, $w_{u}^{\prime}$ (because the rows of the matrix $V$
(from top to bottom) are $w_{k}^{\prime}$, $w_{1}^{\prime}$, $w_{2}^{\prime}$,
$...$, $w_{u}^{\prime}$), this rewrites as follows: $\operatorname*{rows}%
\nolimits_{2,3,...,u+1}V$ is the matrix whose rows (from top to bottom) are
$w_{1}^{\prime}$, $w_{2}^{\prime}$, $...$, $w_{u}^{\prime}$. But since the
matrix whose rows (from top to bottom) are $w_{1}^{\prime}$, $w_{2}^{\prime}$,
$...$, $w_{u}^{\prime}$ is the matrix $W$ (because $w_{1}^{\prime}$,
$w_{2}^{\prime}$, $...$, $w_{u}^{\prime}$ are the rows of this matrix $W$
(from top to bottom)), this rewrites as follows: $\operatorname*{rows}%
\nolimits_{2,3,...,u+1}V$ is the matrix $W$. In other words,
$\operatorname*{rows}\nolimits_{2,3,...,u+1}V=W$.

Let us make a trivial observation: If $z_{1}$, $z_{2}$, $...$, $z_{u}$ are $u$
vectors in $\mathbb{K}^{u}$, then
\begin{equation}
\left(
\begin{array}
[c]{c}%
\left(  z_{1}\mid z_{2}\mid...\mid z_{u}\right)  \text{ is the }u\times
u\text{-matrix whose columns}\\
\text{(from left to right) are }z_{1}\text{, }z_{2}\text{, }...\text{, }z_{u}%
\end{array}
\right)  . \label{pf.pluecker.pluecker.6.pf.0}%
\end{equation}
\footnote{\textit{Proof of (\ref{pf.pluecker.pluecker.6.pf.0}):} Let $z_{1}$,
$z_{2}$, $...$, $z_{u}$ be $u$ vectors in $\mathbb{K}^{u}$. Then, $z_{i}$ is a
$u\times1$-matrix for every $i\in\left\{  1,2,...,u\right\}  $ (because
$z_{i}\in\mathbb{K}^{u}$ for every $i\in\left\{  1,2,...,u\right\}  $). Hence,
$\left(  z_{1}\mid z_{2}\mid...\mid z_{u}\right)  $ is a $u\times\left(
\underbrace{1+1+...+1}_{u\text{ times }1}\right)  $-matrix. Since
$\underbrace{1+1+...+1}_{u\text{ times }1}=u$, this rewrites as follows:
$\left(  z_{1}\mid z_{2}\mid...\mid z_{u}\right)  $ is a $u\times u$-matrix.
\par
Now, (\ref{conv.attach.formal}) (applied to $u$, $u$, $1$, $z_{m}$ and $1$
instead of $k$, $p$, $q_{m}$, $A_{m}$ and $j$) yields that for every
$i\in\left\{  1,2,...,u\right\}  $, we have%
\begin{align}
&  \left(  \text{the }\left(  \underbrace{\left(  1+1+...+1\right)
}_{i-1\text{ times }1}+1\right)  \text{-th column of the matrix }\left(
z_{1}\mid z_{2}\mid...\mid z_{u}\right)  \right) \nonumber\\
&  =\left(  \text{the }1\text{-st column of the matrix }z_{i}\right)  .
\label{pf.pluecker.pluecker.6.pf.0.pf.1}%
\end{align}
\par
Now, fix $i\in\left\{  1,2,...,u\right\}  $. Then,
$\underbrace{\underbrace{\left(  1+1+...+1\right)  }_{i-1\text{ times }1}%
}_{=i-1}+1=\left(  i-1\right)  +1=i$, so that%
\begin{align*}
&  \left(  \text{the }\left(  \underbrace{\left(  1+1+...+1\right)
}_{i-1\text{ times }1}+1\right)  \text{-th column of the matrix }\left(
z_{1}\mid z_{2}\mid...\mid z_{u}\right)  \right) \\
&  =\left(  \text{the }i\text{-th column of the matrix }\left(  z_{1}\mid
z_{2}\mid...\mid z_{u}\right)  \right)  .
\end{align*}
Compared with (\ref{pf.pluecker.pluecker.6.pf.0.pf.1}), this becomes%
\begin{align*}
&  \left(  \text{the }i\text{-th column of the matrix }\left(  z_{1}\mid
z_{2}\mid...\mid z_{u}\right)  \right) \\
&  =\left(  \text{the }1\text{-st column of the matrix }z_{i}\right)  =z_{i}%
\end{align*}
(because of how we identify column vectors in $\mathbb{K}^{u}$ with $u\times
1$-matrices).
\par
Now, forget that we fixed $i$. We thus have shown that every $i\in\left\{
1,2,...,u\right\}  $ satisfies $\left(  \text{the }i\text{-th column of the
matrix }\left(  z_{1}\mid z_{2}\mid...\mid z_{u}\right)  \right)  =z_{i}$.
Since the matrix $\left(  z_{1}\mid z_{2}\mid...\mid z_{u}\right)  $ has
altogether $u$ columns (because $\left(  z_{1}\mid z_{2}\mid...\mid
z_{u}\right)  $ is a $u\times u$-matrix), this yields that the columns of the
matrix $\left(  z_{1}\mid z_{2}\mid...\mid z_{u}\right)  $ (from left to
right) are $z_{1}$, $z_{2}$, $...$, $z_{u}$. In other words, $\left(
z_{1}\mid z_{2}\mid...\mid z_{u}\right)  $ is the $u\times u$-matrix whose
columns (from left to right) are $z_{1}$, $z_{2}$, $...$, $z_{u}$. This proves
(\ref{pf.pluecker.pluecker.6.pf.0}).}

For every $\ell\in\left\{  1,2,...,u+1\right\}  $, we have%
\begin{equation}
\left(  \text{the }\left(  1,\ell\right)  \text{-th entry of the matrix
}V\right)  =w_{k,\ell} \label{pf.pluecker.pluecker.6.pf.1}%
\end{equation}
\footnote{\textit{Proof of (\ref{pf.pluecker.pluecker.6.pf.1}):} Let $\ell
\in\left\{  1,2,...,u+1\right\}  $. Recall that $V$ is the square matrix whose
rows (from top to bottom) are $w_{k}^{\prime}$, $w_{1}^{\prime}$,
$w_{2}^{\prime}$, $...$, $w_{u}^{\prime}$. Hence, the $1$-st row of the matrix
$V$ is $w_{k}^{\prime}$. In other words,%
\begin{equation}
\left(  \text{the }1\text{-st row of the matrix }V\right)  =w_{k}^{\prime}.
\label{pf.pluecker.pluecker.6.pf.1.pf.1}%
\end{equation}
\par
But $w_{1}^{\prime}$, $w_{2}^{\prime}$, $...$, $w_{u}^{\prime}$ are the rows
of the matrix $W$ (from top to bottom). Hence, the $k$-th row of the matrix
$W$ is $w_{k}^{\prime}$. In other words,%
\begin{equation}
\left(  \text{the }k\text{-th row of the matrix }W\right)  =w_{k}^{\prime}.
\label{pf.pluecker.pluecker.6.pf.1.pf.2}%
\end{equation}
\par
But by the definition of $w_{k,\ell}$, we know that $w_{k,\ell}$ is the
$\left(  k,\ell\right)  $-th entry of the matrix $W$. In other words,%
\begin{align*}
w_{k,\ell}  &  =\left(  \text{the }\left(  k,\ell\right)  \text{-th entry of
the matrix }W\right) \\
&  =\left(  \text{the }\ell\text{-th entry of }\underbrace{\text{the
}k\text{-th row of the matrix }W}_{\substack{=w_{k}^{\prime}\\\text{(by
(\ref{pf.pluecker.pluecker.6.pf.1.pf.2}))}}}\right) \\
&  =\left(  \text{the }\ell\text{-th entry of }w_{k}^{\prime}\right)  .
\end{align*}
Compared with%
\begin{align*}
&  \left(  \text{the }\left(  1,\ell\right)  \text{-th entry of the matrix
}V\right) \\
&  =\left(  \text{the }\ell\text{-th entry of }\underbrace{\text{the
}1\text{-st row of the matrix }V}_{\substack{=w_{k}^{\prime}\\\text{(by
(\ref{pf.pluecker.pluecker.6.pf.1.pf.1}))}}}\right) \\
&  =\left(  \text{the }\ell\text{-th entry of }w_{k}^{\prime}\right)  ,
\end{align*}
this yields $\left(  \text{the }\left(  1,\ell\right)  \text{-th entry of the
matrix }V\right)  =w_{k,\ell}$. This proves (\ref{pf.pluecker.pluecker.6.pf.1}%
).}. For every $r\in\left\{  1,2,...,u+1\right\}  $, we have%
\begin{align}
&  \left(  \text{the }\left(  1,r\right)  \text{-th cofactor of the matrix
}V\right) \nonumber\\
&  =\left(  -1\right)  ^{r-1}\det\left(  w_{1}\mid w_{2}\mid...\mid
\widehat{w_{r}}\mid...\mid w_{u+1}\right)  .
\label{pf.pluecker.pluecker.6.pf.2}%
\end{align}
\footnote{\textit{Proof of (\ref{pf.pluecker.pluecker.6.pf.2}):} Let
$r\in\left\{  1,2,...,u+1\right\}  $. By the definition of a cofactor, it is
clear that the $\left(  1,r\right)  $-th cofactor of the matrix $V$ is
$\left(  -1\right)  ^{r-1}$ times the determinant of the matrix obtained by
removing row $1$ and column $r$ from the matrix $V$. In other words,%
\begin{align}
&  \left(  \text{the }\left(  1,r\right)  \text{-th cofactor of the matrix
}V\right) \nonumber\\
&  =\left(  -1\right)  ^{r-1}\cdot\det\underbrace{\left(  \text{the matrix
obtained by removing row }1\text{ and column }r\text{ from the matrix
}V\right)  }_{=\operatorname*{cols}\nolimits_{1,2,...,\widehat{r}%
,...,u+1}\left(  \operatorname*{rows}\nolimits_{2,3,...,u+1}V\right)
}\nonumber\\
&  =\left(  -1\right)  ^{r-1}\cdot\det\left(  \operatorname*{cols}%
\nolimits_{1,2,...,\widehat{r},...,u+1}\underbrace{\left(
\operatorname*{rows}\nolimits_{2,3,...,u+1}V\right)  }_{=W}\right) \nonumber\\
&  =\left(  -1\right)  ^{r-1}\cdot\det\left(  \operatorname*{cols}%
\nolimits_{1,2,...,\widehat{r},...,u+1}W\right)  .
\label{pf.pluecker.pluecker.6.pf.2.pf.1}%
\end{align}
\par
Now, recall that if $p$ and $q$ are any nonnegative integers, if $B$ is a
$p\times q$-matrix, and if $j_{1}$, $j_{2}$, $...$, $j_{\ell}$ are some
elements of $\left\{  1,2,...,q\right\}  $, then $\operatorname*{cols}%
\nolimits_{j_{1},j_{2},...,j_{\ell}}B$ is the matrix whose columns (from left
to right) are the columns labelled $j_{1}$, $j_{2}$, $...$, $j_{\ell}$ of the
matrix $B$ (because this is how $\operatorname*{cols}\nolimits_{j_{1}%
,j_{2},...,j_{\ell}}B$ was defined in Convention \ref{conv.cols}). Applying
this to $u$, $u+1$, $W$ and $\left(  1,2,...,\widehat{r},...,u+1\right)  $
instead of $p$, $q$, $B$ and $\left(  j_{1},j_{2},...,j_{\ell}\right)  $, we
conclude that $\operatorname*{cols}\nolimits_{1,2,...,\widehat{r},...,u+1}W$
is the matrix whose columns (from left to right) are the columns labelled $1$,
$2$, $...$, $\widehat{r}$, $...$, $u+1$ of the matrix $W$. Since the columns
labelled $1$, $2$, $...$, $\widehat{r}$, $...$, $u+1$ of the matrix $W$ are
$w_{1}$, $w_{2}$, $...$, $\widehat{w_{r}}$, $...$, $w_{u+1}$ (because the
columns of the matrix $W$ (from left to right) are $w_{1}$, $w_{2}$, $...$,
$w_{u+1}$), this rewrites as follows: $\operatorname*{cols}%
\nolimits_{1,2,...,\widehat{r},...,u+1}W$ is the matrix whose columns (from
left to right) are $w_{1}$, $w_{2}$, $...$, $\widehat{w_{r}}$, $...$,
$w_{u+1}$. In other words, $\operatorname*{cols}\nolimits_{1,2,...,\widehat{r}%
,...,u+1}W$ is the $u\times u$-matrix whose columns (from left to right) are
$w_{1}$, $w_{2}$, $...$, $\widehat{w_{r}}$, $...$, $w_{u+1}$.
\par
But we can apply (\ref{pf.pluecker.pluecker.6.pf.0}) to $\left(  z_{1}%
,z_{2},...,z_{u}\right)  =\left(  w_{1},w_{2},...,\widehat{w_{r}}%
,...,w_{u+1}\right)  $. As a consequence, we conclude that $\left(  w_{1}\mid
w_{2}\mid...\mid\widehat{w_{r}}\mid...\mid w_{u+1}\right)  $ is the $u\times
u$-matrix whose columns (from left to right) are $w_{1}$, $w_{2}$, $...$,
$\widehat{w_{r}}$, $...$, $w_{u+1}$. In other words, $\left(  w_{1}\mid
w_{2}\mid...\mid\widehat{w_{r}}\mid...\mid w_{u+1}\right)  $ is
$\operatorname*{cols}\nolimits_{1,2,...,\widehat{r},...,u+1}W$ (because
$\operatorname*{cols}\nolimits_{1,2,...,\widehat{r},...,u+1}W$ is the $u\times
u$-matrix whose columns (from left to right) are $w_{1}$, $w_{2}$, $...$,
$\widehat{w_{r}}$, $...$, $w_{u+1}$). In other words,%
\begin{equation}
\left(  w_{1}\mid w_{2}\mid...\mid\widehat{w_{r}}\mid...\mid w_{u+1}\right)
=\operatorname*{cols}\nolimits_{1,2,...,\widehat{r},...,u+1}W.
\label{pf.pluecker.pluecker.6.pf.2.pf.2}%
\end{equation}
\par
Now, (\ref{pf.pluecker.pluecker.6.pf.2.pf.1}) becomes%
\begin{align*}
\left(  \text{the }\left(  1,r\right)  \text{-th cofactor of the matrix
}V\right)   &  =\left(  -1\right)  ^{r-1}\cdot\det\left(
\underbrace{\operatorname*{cols}\nolimits_{1,2,...,\widehat{r},...,u+1}%
W}_{\substack{=\left(  w_{1}\mid w_{2}\mid...\mid\widehat{w_{r}}\mid...\mid
w_{u+1}\right)  \\\text{(by (\ref{pf.pluecker.pluecker.6.pf.2.pf.2}))}%
}}\right) \\
&  =\left(  -1\right)  ^{r-1}\det\left(  w_{1}\mid w_{2}\mid...\mid
\widehat{w_{r}}\mid...\mid w_{u+1}\right)  .
\end{align*}
This proves (\ref{pf.pluecker.pluecker.6.pf.2}).}

Now, recall that $\det V=0$. Hence,%
\begin{align*}
0  &  =\det V\\
&  =\sum_{\ell=1}^{u+1}\underbrace{\left(  \text{the }\left(  1,\ell\right)
\text{-th entry of the matrix }V\right)  }_{\substack{=w_{k,\ell}\\\text{(by
(\ref{pf.pluecker.pluecker.6.pf.1}))}}}\\
&  \ \ \ \ \ \ \ \ \ \ \ \ \ \ \ \ \ \ \ \ \cdot\left(  \text{the }\left(
1,\ell\right)  \text{-th cofactor of the matrix }V\right) \\
&  \ \ \ \ \ \ \ \ \ \ \left(
\begin{array}
[c]{c}%
\text{by Proposition \ref{prop.cofactor.laplace}, applied to }u+1\text{ and
}V\\
\text{instead of }n\text{ and }B
\end{array}
\right) \\
&  =\sum_{\ell=1}^{u+1}w_{k,\ell}\cdot\left(  \text{the }\left(
1,\ell\right)  \text{-th cofactor of the matrix }V\right) \\
&  =\sum_{r=1}^{u+1}w_{k,r}\cdot\underbrace{\left(  \text{the }\left(
1,r\right)  \text{-th cofactor of the matrix }V\right)  }_{\substack{=\left(
-1\right)  ^{r-1}\det\left(  w_{1}\mid w_{2}\mid...\mid\widehat{w_{r}}%
\mid...\mid w_{u+1}\right)  \\\text{(by (\ref{pf.pluecker.pluecker.6.pf.2}))}%
}}\\
&  \ \ \ \ \ \ \ \ \ \ \left(  \text{here, we renamed the summation index
}\ell\text{ as }r\right) \\
&  =\sum_{r=1}^{u+1}w_{k,r}\cdot\underbrace{\left(  -1\right)  ^{r-1}%
}_{=-\left(  -1\right)  ^{r}}\det\left(  w_{1}\mid w_{2}\mid...\mid
\widehat{w_{r}}\mid...\mid w_{u+1}\right) \\
&  =\sum_{r=1}^{u+1}w_{k,r}\cdot\left(  -\left(  -1\right)  ^{r}\right)
\det\left(  w_{1}\mid w_{2}\mid...\mid\widehat{w_{r}}\mid...\mid
w_{u+1}\right) \\
&  =-\sum_{r=1}^{u+1}\left(  -1\right)  ^{r}w_{k,r}\cdot\det\left(  w_{1}\mid
w_{2}\mid...\mid\widehat{w_{r}}\mid...\mid w_{u+1}\right)  .
\end{align*}
Adding $\sum\limits_{r=1}^{u+1}\left(  -1\right)  ^{r}w_{k,r}\cdot\det\left(
w_{1}\mid w_{2}\mid...\mid\widehat{w_{r}}\mid...\mid w_{u+1}\right)  $ to this
equality, we obtain%
\[
\sum\limits_{r=1}^{u+1}\left(  -1\right)  ^{r}w_{k,r}\cdot\det\left(
w_{1}\mid w_{2}\mid...\mid\widehat{w_{r}}\mid...\mid w_{u+1}\right)  =0.
\]
This proves (\ref{pf.pluecker.pluecker.6}).

Now, (\ref{pf.pluecker.pluecker.5}) becomes%
\begin{align*}
&  \sum\limits_{r=1}^{u+1}\left(  -1\right)  ^{r}\det\left(  A\mid
w_{r}\right)  \cdot\det\left(  w_{1}\mid w_{2}\mid...\mid\widehat{w_{r}}%
\mid...\mid w_{u+1}\right) \\
&  =\sum\limits_{k=1}^{u}\left(  -1\right)  ^{u-k}\cdot\det\left(
\operatorname*{rows}\nolimits_{1,2,...,\widehat{k},...,u}A\right) \\
&  \ \ \ \ \ \ \ \ \ \ \cdot\underbrace{\left(  \sum\limits_{r=1}^{u+1}\left(
-1\right)  ^{r}w_{k,r}\cdot\det\left(  w_{1}\mid w_{2}\mid...\mid
\widehat{w_{r}}\mid...\mid w_{u+1}\right)  \right)  }_{\substack{=0\\\text{(by
(\ref{pf.pluecker.pluecker.6}))}}}\\
&  =\sum\limits_{k=1}^{u}\left(  -1\right)  ^{u-k}\cdot\det\left(
\operatorname*{rows}\nolimits_{1,2,...,\widehat{k},...,u}A\right)  \cdot0=0.
\end{align*}
This proves Proposition \ref{prop.pluecker.pluecker}.
\end{proof}

\begin{proof}
[Second proof of Lemma \ref{lem.pluecker.ptolemy} (sketched).]We have $p\geq2$
and $q\geq2$, so that $p+q\geq2+2=4$. Thus, $u+2=p+q\geq4$, and therefore
$u\geq2>0$. Hence, $u$ is a positive integer. Since $q\geq2$, we have $2\leq
q$, so that $p+2\leq p+q=u+2$, and thus $p\leq u$, so that $p+1\leq u+1$.
Hence, $\left\{  1,p,p+1\right\}  \subseteq\left\{  1,2,...,u+1\right\}  $.

Let $A\in\mathbb{K}^{u\times\left(  u-1\right)  }$ be the $u\times\left(
u-1\right)  $-matrix $B\left[  2:p\mid p+2:p+q+1\right]  $. Then, $A$ is the
matrix whose columns (from left to right) are $B_{2}$, $B_{3}$, $...$,
$B_{p-1}$, $B_{p+2}$, $B_{p+3}$, $...$, $B_{p+q}$ (by the definition of
$B\left[  2:p\mid p+2:p+q+1\right]  $). In other words, the columns of $A$
(from left to right) are $B_{2}$, $B_{3}$, $...$, $B_{p-1}$, $B_{p+2}$,
$B_{p+3}$, $...$, $B_{p+q}$.

Define $u+1$ column vectors $w_{1}$, $w_{2}$, $...$, $w_{u+1}$ in
$\mathbb{K}^{u}$ by $\left(  w_{i}=B_{i}\text{ for every }i\in\left\{
1,2,...,u+1\right\}  \right)  $.

We know (from Proposition \ref{prop.pluecker.pluecker}) that
(\ref{prop.pluecker.pluecker.eq}) holds. Now, we are going to simplify the
left hand side of (\ref{prop.pluecker.pluecker.eq}).

We notice that every $r\in\left\{  1,2,...,u+1\right\}  $ satisfying
$r\notin\left\{  1,p,p+1\right\}  $ satisfies%
\begin{equation}
\det\left(  A\mid w_{r}\right)  =0. \label{pf.pluecker.ptolemy.2ndproof.1}%
\end{equation}
\footnote{\textit{Proof of (\ref{pf.pluecker.ptolemy.2ndproof.1}):} Let
$r\in\left\{  1,2,...,u+1\right\}  $ be such that $r\notin\left\{
1,p,p+1\right\}  $. Then, $r\leq u+1$ (since $r\in\left\{
1,2,...,u+1\right\}  $), so that $r\leq u+1\leq u+2=p+q$. Hence, $r\in\left\{
1,2,...,p+q\right\}  $ and $r\notin\left\{  1,p,p+1\right\}  $. In other
words, $r$ is one of the integers $2$, $3$, $...$, $p-1$, $p+2$, $p+3$, $...$,
$p+q$. Thus, the column vector $B_{r}$ is one of the column vectors $B_{2}$,
$B_{3}$, $...$, $B_{p-1}$, $B_{p+2}$, $B_{p+3}$, $...$, $B_{p+q}$.
\par
Now, recall that $A$ is the matrix whose columns (from left to right) are
$B_{2}$, $B_{3}$, $...$, $B_{p-1}$, $B_{p+2}$, $B_{p+3}$, $...$, $B_{p+q}$.
Hence, the columns of the matrix $A$ (from left to right) are the column
vectors $B_{2}$, $B_{3}$, $...$, $B_{p-1}$, $B_{p+2}$, $B_{p+3}$, $...$,
$B_{p+q}$. Hence, $B_{r}$ is one of the columns of the matrix $A$ (since
$B_{r}$ is one of the column vectors $B_{2}$, $B_{3}$, $...$, $B_{p-1}$,
$B_{p+2}$, $B_{p+3}$, $...$, $B_{p+q}$). But $w_{r}=B_{r}$ (by the definition
of $w_{r}$).
\par
Now, consider the matrix $\left(  A\mid w_{r}\right)  $. The column $B_{r}$
appears twice in this matrix: one time as a column of the matrix $A$ (since
$B_{r}$ is one of the columns of the matrix $A$), and one more time as the
single column of the $u\times1$-matrix $w_{r}$ (since $w_{r}=B_{r}$). Hence,
the matrix $\left(  A\mid w_{r}\right)  $ has a column appearing twice. In
other words, the matrix $\left(  A\mid w_{r}\right)  $ has two equal columns.
Therefore, the matrix $\left(  A\mid w_{r}\right)  $ has zero determinant
(because any matrix which has two equal columns must have zero determinant).
In other words, $\det\left(  A\mid w_{r}\right)  =0$. This proves
(\ref{pf.pluecker.ptolemy.2ndproof.1}).}

But from (\ref{prop.pluecker.pluecker.eq}), we obtain%
\begin{align*}
0  &  =\underbrace{\sum\limits_{r=1}^{u+1}}_{=\sum\limits_{r\in\left\{
1,2,...,u+1\right\}  }}\left(  -1\right)  ^{r}\det\left(  A\mid w_{r}\right)
\cdot\det\left(  w_{1}\mid w_{2}\mid...\mid\widehat{w_{r}}\mid...\mid
w_{u+1}\right) \\
&  =\sum\limits_{r\in\left\{  1,2,...,u+1\right\}  }\left(  -1\right)
^{r}\det\left(  A\mid w_{r}\right)  \cdot\det\left(  w_{1}\mid w_{2}%
\mid...\mid\widehat{w_{r}}\mid...\mid w_{u+1}\right) \\
&  =\underbrace{\sum\limits_{\substack{r\in\left\{  1,2,...,u+1\right\}
;\\r\in\left\{  1,p,p+1\right\}  }}}_{\substack{=\sum\limits_{r\in\left\{
1,p,p+1\right\}  }\\\text{(since }\left\{  1,p,p+1\right\}  \subseteq\left\{
1,2,...,u+1\right\}  \text{)}}}\left(  -1\right)  ^{r}\det\left(  A\mid
w_{r}\right)  \cdot\det\left(  w_{1}\mid w_{2}\mid...\mid\widehat{w_{r}}%
\mid...\mid w_{u+1}\right) \\
&  \ \ \ \ \ \ \ \ \ \ +\sum\limits_{\substack{r\in\left\{
1,2,...,u+1\right\}  ;\\r\notin\left\{  1,p,p+1\right\}  }}\left(  -1\right)
^{r}\underbrace{\det\left(  A\mid w_{r}\right)  }_{\substack{=0\\\text{(by
(\ref{pf.pluecker.ptolemy.2ndproof.1}))}}}\cdot\det\left(  w_{1}\mid w_{2}%
\mid...\mid\widehat{w_{r}}\mid...\mid w_{u+1}\right) \\
&  =\sum\limits_{r\in\left\{  1,p,p+1\right\}  }\left(  -1\right)  ^{r}%
\det\left(  A\mid w_{r}\right)  \cdot\det\left(  w_{1}\mid w_{2}\mid
...\mid\widehat{w_{r}}\mid...\mid w_{u+1}\right) \\
&  \ \ \ \ \ \ \ \ \ \ +\underbrace{\sum\limits_{\substack{r\in\left\{
1,2,...,u+1\right\}  ;\\r\notin\left\{  1,p,p+1\right\}  }}\left(  -1\right)
^{r}0\cdot\det\left(  w_{1}\mid w_{2}\mid...\mid\widehat{w_{r}}\mid...\mid
w_{u+1}\right)  }_{=0}\\
&  =\sum\limits_{r\in\left\{  1,p,p+1\right\}  }\left(  -1\right)  ^{r}%
\det\left(  A\mid w_{r}\right)  \cdot\det\left(  w_{1}\mid w_{2}\mid
...\mid\widehat{w_{r}}\mid...\mid w_{u+1}\right)
\end{align*}%
\begin{align}
&  =\left(  -1\right)  ^{1}\det\left(  A\mid w_{1}\right)  \cdot\det\left(
w_{1}\mid w_{2}\mid...\mid\widehat{w_{1}}\mid...\mid w_{u+1}\right)
\nonumber\\
&  \ \ \ \ \ \ \ \ \ \ +\left(  -1\right)  ^{p}\det\left(  A\mid w_{p}\right)
\cdot\det\left(  w_{1}\mid w_{2}\mid...\mid\widehat{w_{p}}\mid...\mid
w_{u+1}\right) \nonumber\\
&  \ \ \ \ \ \ \ \ \ \ +\left(  -1\right)  ^{p+1}\det\left(  A\mid
w_{p+1}\right)  \cdot\det\left(  w_{1}\mid w_{2}\mid...\mid\widehat{w_{p+1}%
}\mid...\mid w_{u+1}\right) \label{pf.pluecker.ptolemy.2ndproof.2}\\
&  \ \ \ \ \ \ \ \ \ \ \left(  \text{since }1\text{, }p\text{ and }p+1\text{
are pairwise distinct (because }p\geq2\text{)}\right)  .\nonumber
\end{align}
We shall now simplify the determinants appearing on the right hand side of
this equality.

Notice first that $p+q=u+2$, so that $p+q-1=u+2-1=u+1$. Hence, $\left(
-1\right)  ^{p+q-1}=\left(  -1\right)  ^{u+1}=-\left(  -1\right)  ^{u}$. Also,
$p+q=u+2$, so that $\left(  -1\right)  ^{p+q}=\left(  -1\right)
^{u+2}=\left(  -1\right)  ^{u}\underbrace{\left(  -1\right)  ^{2}}%
_{=1}=\left(  -1\right)  ^{u}$.

We have%
\begin{equation}
\det\left(  A\mid w_{1}\right)  =\left(  -1\right)  ^{u-1}\det\left(  B\left[
1:p\mid p+2:p+q+1\right]  \right)  . \label{pf.pluecker.ptolemy.2ndproof.3a1}%
\end{equation}
\footnote{\textit{Proof of (\ref{pf.pluecker.ptolemy.2ndproof.3a1}):} Let us
first notice that $w_{1}$ is a single column vector. Hence, the columns of
$w_{1}$ (from left to right) are $w_{1}$. In other words, the columns of
$w_{1}$ (from left to right) are $B_{1}$ (since $w_{1}=B_{1}$ (by the
definition of $w_{1}$)).
\par
We know (by the definition of $\left(  A\mid w_{1}\right)  $) that the matrix
$\left(  A\mid w_{1}\right)  $ is the matrix whose columns (from left to
right) are all the columns of $A$ (from left to right) and then all the
columns of $w_{1}$ (from left to right). Since the columns of $A$ (from left
to right) are $B_{2}$, $B_{3}$, $...$, $B_{p-1}$, $B_{p+2}$, $B_{p+3}$, $...$,
$B_{p+q}$, whereas the columns of $w_{1}$ (from left to right) are $B_{1}$,
this rewrites as follows: The matrix $\left(  A\mid w_{1}\right)  $ is the
matrix whose columns (from left to right) are $B_{2}$, $B_{3}$, $...$,
$B_{p-1}$, $B_{p+2}$, $B_{p+3}$, $...$, $B_{p+q}$, $B_{1}$.
\par
Let $\sigma$ be the permutation in $S_{u}$ which sends every $i\in\left\{
1,2,...,u\right\}  $ to $\left\{
\begin{array}
[c]{c}%
i+1,\text{ if }i<u;\\
1,\text{ if }i=u
\end{array}
\right.  $. Then, $\sigma$ sends the numbers $1, 2, \ldots, u$ to the
numbers $2, 3, \ldots, u, 1$, respectively. Thus,
$\sigma$ is a cycle of length $u$, so that $\left(
-1\right)  ^{\sigma}=\left(  -1\right)  ^{u-1}$.
\par
On the other hand, the matrix $B\left[  1:p\mid p+2:p+q+1\right]  $ is the
matrix whose columns (from left to right) are $B_{1}$, $B_{2}$, $...$,
$B_{p-1}$, $B_{p+2}$, $B_{p+3}$, $...$, $B_{p+q}$ (by the definition of
$B\left[  1:p\mid p+2:p+q+1\right]  $). Hence, permuting the columns of the
matrix $B\left[  1:p\mid p+2:p+q+1\right]  $ according to the permutation
$\sigma$ yields the matrix whose columns (from left to right) are $B_{2}$,
$B_{3}$, $...$, $B_{p-1}$, $B_{p+2}$, $B_{p+3}$, $...$, $B_{p+q}$, $B_{1}$
(since $\sigma$ sends the numbers $1, 2, \ldots, u$ to the
numbers $2, 3, \ldots, u, 1$, respectively). In
other words, permuting the columns of the matrix $B\left[  1:p\mid
p+2:p+q+1\right]  $ according to the permutation $\sigma$ yields the matrix
$\left(  A\mid w_{1}\right)  $ (since $\left(  A\mid w_{1}\right)  $ is the
matrix whose columns (from left to right) are $B_{2}$, $B_{3}$, $...$,
$B_{p-1}$, $B_{p+2}$, $B_{p+3}$, $...$, $B_{p+q}$, $B_{1}$). In other words,
the matrix $\left(  A\mid w_{1}\right)  $ is obtained from the matrix
$B\left[  1:p\mid p+2:p+q+1\right]  $ by permuting the columns according to
the permutation $\sigma$. Since permuting the columns of a square matrix
always multiplies its determinant by the sign of the permutation used, we
therefore have%
\begin{align*}
\det\left(  A\mid w_{1}\right)   &  =\underbrace{\left(  -1\right)  ^{\sigma}%
}_{=\left(  -1\right)  ^{u-1}}\det\left(  B\left[  1:p\mid p+2:p+q+1\right]
\right) \\
&  =\left(  -1\right)  ^{u-1}\det\left(  B\left[  1:p\mid p+2:p+q+1\right]
\right)  .
\end{align*}
This proves (\ref{pf.pluecker.ptolemy.2ndproof.3a1}).} Furthermore,%
\begin{equation}
\det\left(  w_{1}\mid w_{2}\mid...\mid\widehat{w_{1}}\mid...\mid
w_{u+1}\right)  =\det\left(  B\left[  2:p+1\mid p+1:p+q\right]  \right)  .
\label{pf.pluecker.ptolemy.2ndproof.3a2}%
\end{equation}
\footnote{\textit{Proof of (\ref{pf.pluecker.ptolemy.2ndproof.3a2}):} We know
(by the definition of $\left(  w_{2}\mid w_{3}\mid...\mid w_{u+1}\right)  $)
that the matrix $\left(  w_{2}\mid w_{3}\mid...\mid w_{u+1}\right)  $ is the
matrix whose columns (from left to right) are all the columns of $w_{2}$ (from
left to right), then all the columns of $w_{3}$ (from left to right), etc.,
and finally all the columns of $w_{u+1}$ (from left to right). In other words,
$\left(  w_{2}\mid w_{3}\mid...\mid w_{u+1}\right)  $ is the matrix whose
columns (from left to right) are $w_{2}$, $w_{3}$, $...$, $w_{u+1}$ (because
for any $i\in\left\{  2,3,...,u+1\right\}  $, all the columns of $w_{i}$ (from
left to right) are just $w_{i}$ (since $w_{i}$ is a single column vector)). In
other words, $\left(  w_{2}\mid w_{3}\mid...\mid w_{u+1}\right)  $ is the
matrix whose columns (from left to right) are $B_{2}$, $B_{3}$, $...$,
$B_{u+1}$ (since $w_{i}=B_{i}$ for every $i\in\left\{  2,3,...,u+1\right\}  $
(by the definition of $w_{i}$)).
\par
On the other hand, the matrix $B\left[  2:p+1\mid p+1:p+q\right]  $ is the
matrix whose columns (from left to right) are $B_{2}$, $B_{3}$, $...$, $B_{p}%
$, $B_{p+1}$, $B_{p+2}$, $...$, $B_{p+q-1}$ (by the definition of $B\left[
2:p+1\mid p+1:p+q\right]  $). In other words, the matrix $B\left[  2:p+1\mid
p+1:p+q\right]  $ is the matrix whose columns (from left to right) are $B_{2}%
$, $B_{3}$, $...$, $B_{p+q-1}$. In other words, the matrix $B\left[  2:p+1\mid
p+1:p+q\right]  $ is the matrix whose columns (from left to right) are $B_{2}%
$, $B_{3}$, $...$, $B_{u+1}$ (since $p+q-1=u+1$ (since $p+q=u+2$)). In other
words, the matrix $B\left[  2:p+1\mid p+1:p+q\right]  $ is the matrix $\left(
w_{2}\mid w_{3}\mid...\mid w_{u+1}\right)  $ (since $\left(  w_{2}\mid
w_{3}\mid...\mid w_{u+1}\right)  $ is the matrix whose columns (from left to
right) are $B_{2}$, $B_{3}$, $...$, $B_{u+1}$). Hence,%
\[
B\left[  2:p+1\mid p+1:p+q\right]  =\left(  w_{2}\mid w_{3}\mid...\mid
w_{u+1}\right)  .
\]
But clearly,%
\begin{align*}
\left(  w_{1}\mid w_{2}\mid...\mid\widehat{w_{1}}\mid...\mid w_{u+1}\right)
&  =\left(  w_{2}\mid w_{3}\mid...\mid w_{u+1}\right) \\
&  =B\left[  2:p+1\mid p+1:p+q\right]  .
\end{align*}
Taking determinants on both sides of this equation, we obtain
(\ref{pf.pluecker.ptolemy.2ndproof.3a2}).} Next, we notice that%
\begin{equation}
\det\left(  A\mid w_{p}\right)  =\left(  -1\right)  ^{q-1}\det\left(  B\left[
2:p+1\mid p+2:p+q+1\right]  \right)  .
\label{pf.pluecker.ptolemy.2ndproof.3b1}%
\end{equation}
\footnote{\textit{Proof of (\ref{pf.pluecker.ptolemy.2ndproof.3b1}):} Let us
first notice that $w_{p}$ is a single column vector. Hence, all the columns of
$w_{p}$ (from left to right) are $w_{p}$. In other words, all the columns of
$w_{p}$ (from left to right) are $B_{p}$ (since $w_{p}=B_{p}$ (by the
definition of $w_{p}$)).
\par
We know (by the definition of $\left(  A\mid w_{p}\right)  $) that the matrix
$\left(  A\mid w_{p}\right)  $ is the matrix whose columns (from left to
right) are all the columns of $A$ (from left to right) and then all the
columns of $w_{p}$ (from left to right). Since the columns of $A$ (from left
to right) are $B_{2}$, $B_{3}$, $...$, $B_{p-1}$, $B_{p+2}$, $B_{p+3}$, $...$,
$B_{p+q}$, whereas the columns of $w_{p}$ (from left to right) are $B_{p}$,
this rewrites as follows: The matrix $\left(  A\mid w_{p}\right)  $ is the
matrix whose columns (from left to right) are $B_{2}$, $B_{3}$, $...$,
$B_{p-1}$, $B_{p+2}$, $B_{p+3}$, $...$, $B_{p+q}$, $B_{p}$.
\par
Let $\sigma$ be the permutation in $S_{u}$ which sends every $i\in\left\{
1,2,...,u\right\}  $ to $\left\{
\begin{array}
[c]{c}%
i,\text{ if }i\leq p-2;\\
i+1,\text{ if }p-1 \leq i < u;\\
p-1,\text{ if }i=u
\end{array}
\right.  $. Then, $\sigma$ sends the numbers $1, 2, \ldots, u$ to the
numbers $1, 2, \ldots, p-2, p, p+1, \ldots, u, p-1$, respectively. Hence,
$\sigma$ is a cycle of length $u-\left(  p-1\right)
+1=\underbrace{u+2}_{=p+q}-p=p+q-p=q$, so that $\left(  -1\right)  ^{\sigma
}=\left(  -1\right)  ^{q-1}$.
\par
On the other hand, the matrix $B\left[  2:p+1\mid p+2:p+q+1\right]  $ is the
matrix whose columns (from left to right) are $B_{2}$, $B_{3}$, $...$, $B_{p}%
$, $B_{p+2}$, $B_{p+3}$, $...$, $B_{p+q}$ (by the definition of $B\left[
2:p+1\mid p+2:p+q+1\right]  $). Hence, permuting the columns of the matrix
$B\left[  2:p+1\mid p+2:p+q+1\right]  $ according to the permutation $\sigma$
yields the matrix whose columns (from left to right) are $B_{2}$, $B_{3}$,
$...$, $B_{p-1}$, $B_{p+2}$, $B_{p+3}$, $...$, $B_{p+q}$, $B_{p}$
(since $\sigma$ sends the numbers $1, 2, \ldots, u$ to the
numbers $1, 2, \ldots, p-2, p, p+1, \ldots, u, p-1$, respectively). In other
words, permuting the columns of the matrix $B\left[  2:p+1\mid
p+2:p+q+1\right]  $ according to the permutation $\sigma$ yields the matrix
$\left(  A\mid w_{p}\right)  $ (since $\left(  A\mid w_{p}\right)  $ is the
matrix whose columns (from left to right) are $B_{2}$, $B_{3}$, $...$,
$B_{p-1}$, $B_{p+2}$, $B_{p+3}$, $...$, $B_{p+q}$, $B_{p}$). In other words,
the matrix $\left(  A\mid w_{p}\right)  $ is obtained from the matrix
$B\left[  2:p+1\mid p+2:p+q+1\right]  $ by permuting the columns according to
the permutation $\sigma$. Since permuting the columns of a square matrix
always multiplies its determinant by the sign of the permutation used, we
therefore have%
\begin{align*}
\det\left(  A\mid w_{p}\right)   &  =\underbrace{\left(  -1\right)  ^{\sigma}%
}_{=\left(  -1\right)  ^{q-1}}\det\left(  B\left[  2:p+1\mid p+2:p+q+1\right]
\right) \\
&  =\left(  -1\right)  ^{q-1}\det\left(  B\left[  2:p+1\mid p+2:p+q+1\right]
\right)  .
\end{align*}
This proves (\ref{pf.pluecker.ptolemy.2ndproof.3b1}).} Also,%
\begin{equation}
\det\left(  w_{1}\mid w_{2}\mid...\mid\widehat{w_{p}}\mid...\mid
w_{u+1}\right)  =\det\left(  B\left[  1:p\mid p+1:p+q\right]  \right)  .
\label{pf.pluecker.ptolemy.2ndproof.3b2}%
\end{equation}
\footnote{\textit{Proof of (\ref{pf.pluecker.ptolemy.2ndproof.3b2}):} We know
(by the definition of $\left(  w_{1}\mid w_{2}\mid...\mid w_{p-1}\mid
w_{p+1}\mid w_{p+2}\mid...\mid w_{u+1}\right)  $) that the matrix $\left(
w_{1}\mid w_{2}\mid...\mid w_{p-1}\mid w_{p+1}\mid w_{p+2}\mid...\mid
w_{u+1}\right)  $ is the matrix whose columns (from left to right) are all the
columns of $w_{1}$ (from left to right), then all the columns of $w_{2}$ (from
left to right), etc., then all the columns of $w_{p-1}$ (from left to right),
then all the columns of $w_{p+1}$ (from left to right), then all the columns
of $w_{p+2}$ (from left to right), etc., and finally all the columns of
$w_{u+1}$ (from left to right). In other words, $\left(  w_{1}\mid w_{2}%
\mid...\mid w_{p-1}\mid w_{p+1}\mid w_{p+2}\mid...\mid w_{u+1}\right)  $ is
the matrix whose columns (from left to right) are $w_{1}$, $w_{2}$, $...$,
$w_{p-1}$, $w_{p+1}$, $w_{p+2}$, $...$, $w_{u+1}$ (because for any
$i\in\left\{  1,2,...,u+1\right\}  \setminus\left\{  p\right\}  $, all the
columns of $w_{i}$ (from left to right) are just $w_{i}$ (since $w_{i}$ is a
single column vector)). In other words, $\left(  w_{1}\mid w_{2}\mid...\mid
w_{p-1}\mid w_{p+1}\mid w_{p+2}\mid...\mid w_{u+1}\right)  $ is the matrix
whose columns (from left to right) are $B_{1}$, $B_{2}$, $...$, $B_{p-1}$,
$B_{p+1}$, $B_{p+2}$, $...$, $B_{u+1}$ (since $w_{i}=B_{i}$ for every
$i\in\left\{  1,2,...,u+1\right\}  \setminus\left\{  p\right\}  $ (by the
definition of $w_{i}$)).
\par
On the other hand, the matrix $B\left[  1:p\mid p+1:p+q\right]  $ is the
matrix whose columns (from left to right) are $B_{1}$, $B_{2}$, $...$,
$B_{p-1}$, $B_{p+1}$, $B_{p+2}$, $...$, $B_{p+q-1}$ (by the definition of
$B\left[  1:p\mid p+1:p+q\right]  $). In other words, the matrix $B\left[
1:p\mid p+1:p+q\right]  $ is the matrix whose columns (from left to right) are
$B_{1}$, $B_{2}$, $...$, $B_{p-1}$, $B_{p+1}$, $B_{p+2}$, $...$, $B_{u+1}$
(since $p+q-1=u+1$ (since $p+q=u+2$)). In other words, the matrix $B\left[
1:p\mid p+1:p+q\right]  $ is the matrix $\left(  w_{1}\mid w_{2}\mid...\mid
w_{p-1}\mid w_{p+1}\mid w_{p+2}\mid...\mid w_{u+1}\right)  $ (since $\left(
w_{1}\mid w_{2}\mid...\mid w_{p-1}\mid w_{p+1}\mid w_{p+2}\mid...\mid
w_{u+1}\right)  $ is the matrix whose columns (from left to right) are $B_{1}%
$, $B_{2}$, $...$, $B_{p-1}$, $B_{p+1}$, $B_{p+2}$, $...$, $B_{u+1}$). Hence,%
\[
B\left[  1:p\mid p+1:p+q\right]  =\left(  w_{1}\mid w_{2}\mid...\mid
w_{p-1}\mid w_{p+1}\mid w_{p+2}\mid...\mid w_{u+1}\right)  .
\]
But clearly,%
\begin{align*}
\left(  w_{1}\mid w_{2}\mid...\mid\widehat{w_{p}}\mid...\mid w_{u+1}\right)
&  =\left(  w_{1}\mid w_{2}\mid...\mid w_{p-1}\mid w_{p+1}\mid w_{p+2}%
\mid...\mid w_{u+1}\right) \\
&  =B\left[  1:p\mid p+1:p+q\right]  .
\end{align*}
Taking determinants on both sides of this equation, we obtain
(\ref{pf.pluecker.ptolemy.2ndproof.3b2}).} Furthermore,%
\begin{equation}
\det\left(  A\mid w_{p+1}\right)  =\left(  -1\right)  ^{q-1}\det\left(
B\left[  2:p\mid p+1:p+q+1\right]  \right)  .
\label{pf.pluecker.ptolemy.2ndproof.3c1}%
\end{equation}
\footnote{\textit{Proof of (\ref{pf.pluecker.ptolemy.2ndproof.3c1}):} Let us
first notice that $w_{p+1}$ is a single column vector. Hence, all the columns
of $w_{p+1}$ (from left to right) are $w_{p+1}$. In other words, all the
columns of $w_{p+1}$ (from left to right) are $B_{p+1}$ (since $w_{p+1}%
=B_{p+1}$ (by the definition of $w_{p+1}$)).
\par
We know (by the definition of $\left(  A\mid w_{p+1}\right)  $) that the
matrix $\left(  A\mid w_{p+1}\right)  $ is the matrix whose columns (from left
to right) are all the columns of $A$ (from left to right) and then all the
columns of $w_{p+1}$ (from left to right). Since the columns of $A$ (from left
to right) are $B_{2}$, $B_{3}$, $...$, $B_{p-1}$, $B_{p+2}$, $B_{p+3}$, $...$,
$B_{p+q}$, whereas the columns of $w_{p+1}$ (from left to right) are $B_{p+1}%
$, this rewrites as follows: The matrix $\left(  A\mid w_{p+1}\right)  $ is
the matrix whose columns (from left to right) are $B_{2}$, $B_{3}$, $...$,
$B_{p-1}$, $B_{p+2}$, $B_{p+3}$, $...$, $B_{p+q}$, $B_{p+1}$.
\par
Let $\sigma$ be the permutation in $S_{u}$ which sends every $i\in\left\{
1,2,...,u\right\}  $ to $\left\{
\begin{array}
[c]{c}%
i,\text{ if }i\leq p-2;\\
i+1,\text{ if }p-1 \leq i < u;\\
p-1,\text{ if }i=u
\end{array}
\right.  $. Then, $\sigma$ sends the numbers $1, 2, \ldots, u$ to the
numbers $1, 2, \ldots, p-2, p, p+1, \ldots, u, p-1$, respectively. Hence,
$\sigma$ is a cycle of length $u-\left(  p-1\right)
+1=\underbrace{u+2}_{=p+q}-p=p+q-p=q$, so that $\left(  -1\right)  ^{\sigma
}=\left(  -1\right)  ^{q-1}$.
\par
On the other hand, the matrix $B\left[  2:p\mid p+1:p+q+1\right]  $ is the
matrix whose columns (from left to right) are $B_{2}$, $B_{3}$, $...$,
$B_{p-1}$, $B_{p+1}$, $B_{p+2}$, $...$, $B_{p+q}$ (by the definition of
$B\left[  2:p\mid p+1:p+q+1\right]  $). Hence, permuting the columns of the
matrix $B\left[  2:p\mid p+1:p+q+1\right]  $ according to the permutation
$\sigma$ yields the matrix whose columns (from left to right) are $B_{2}$,
$B_{3}$, $...$, $B_{p-1}$, $B_{p+2}$, $B_{p+3}$, $...$, $B_{p+q}$, $B_{p+1}$
(since $\sigma$ sends the numbers $1, 2, \ldots, u$ to the
numbers $1, 2, \ldots, p-2, p, p+1, \ldots, u, p-1$, respectively).
In other words, permuting the columns of the matrix $B\left[  2:p\mid
p+1:p+q+1\right]  $ according to the permutation $\sigma$ yields the matrix
$\left(  A\mid w_{p+1}\right)  $ (since $\left(  A\mid w_{p+1}\right)  $ is
the matrix whose columns (from left to right) are $B_{2}$, $B_{3}$, $...$,
$B_{p-1}$, $B_{p+2}$, $B_{p+3}$, $...$, $B_{p+q}$, $B_{p+1}$). In other words,
the matrix $\left(  A\mid w_{p+1}\right)  $ is obtained from the matrix
$B\left[  2:p\mid p+1:p+q+1\right]  $ by permuting the columns according to
the permutation $\sigma$. Since permuting the columns of a square matrix
always multiplies its determinant by the sign of the permutation used, we
therefore have%
\begin{align*}
\det\left(  A\mid w_{p+1}\right)   &  =\underbrace{\left(  -1\right)
^{\sigma}}_{=\left(  -1\right)  ^{q-1}}\det\left(  B\left[  2:p\mid
p+1:p+q+1\right]  \right) \\
&  =\left(  -1\right)  ^{q-1}\det\left(  B\left[  2:p\mid p+1:p+q+1\right]
\right)  .
\end{align*}
This proves (\ref{pf.pluecker.ptolemy.2ndproof.3c1}).} Also,%
\begin{equation}
\det\left(  w_{1}\mid w_{2}\mid...\mid\widehat{w_{p+1}}\mid...\mid
w_{u+1}\right)  =\det\left(  B\left[  1:p+1\mid p+2:p+q\right]  \right)  .
\label{pf.pluecker.ptolemy.2ndproof.3c2}%
\end{equation}
\footnote{\textit{Proof of (\ref{pf.pluecker.ptolemy.2ndproof.3c2}):} We know
(by the definition of $\left(  w_{1}\mid w_{2}\mid...\mid w_{p}\mid
w_{p+2}\mid w_{p+3}\mid...\mid w_{u+1}\right)  $) that the matrix $\left(
w_{1}\mid w_{2}\mid...\mid w_{p}\mid w_{p+2}\mid w_{p+3}\mid...\mid
w_{u+1}\right)  $ is the matrix whose columns (from left to right) are all the
columns of $w_{1}$ (from left to right), then all the columns of $w_{2}$ (from
left to right), etc., then all the columns of $w_{p}$ (from left to right),
then all the columns of $w_{p+2}$ (from left to right), then all the columns
of $w_{p+3}$ (from left to right), etc., and finally all the columns of
$w_{u+1}$ (from left to right). In other words, $\left(  w_{1}\mid w_{2}%
\mid...\mid w_{p}\mid w_{p+2}\mid w_{p+3}\mid...\mid w_{u+1}\right)  $ is the
matrix whose columns (from left to right) are $w_{1}$, $w_{2}$, $...$, $w_{p}%
$, $w_{p+2}$, $w_{p+3}$, $...$, $w_{u+1}$ (because for any $i\in\left\{
1,2,...,u+1\right\}  \setminus\left\{  p+1\right\}  $, all the columns of
$w_{i}$ (from left to right) are just $w_{i}$ (since $w_{i}$ is a single
column vector)). In other words, $\left(  w_{1}\mid w_{2}\mid...\mid w_{p}\mid
w_{p+2}\mid w_{p+3}\mid...\mid w_{u+1}\right)  $ is the matrix whose columns
(from left to right) are $B_{1}$, $B_{2}$, $...$, $B_{p}$, $B_{p+2}$,
$B_{p+3}$, $...$, $B_{u+1}$ (since $w_{i}=B_{i}$ for every $i\in\left\{
1,2,...,u+1\right\}  \setminus\left\{  p+1\right\}  $ (by the definition of
$w_{i}$)).
\par
On the other hand, the matrix $B\left[  1:p+1\mid p+2:p+q\right]  $ is the
matrix whose columns (from left to right) are $B_{1}$, $B_{2}$, $...$, $B_{p}%
$, $B_{p+2}$, $B_{p+3}$, $...$, $B_{p+q-1}$ (by the definition of $B\left[
1:p+1\mid p+2:p+q\right]  $). In other words, the matrix $B\left[  1:p+1\mid
p+2:p+q\right]  $ is the matrix whose columns (from left to right) are $B_{1}%
$, $B_{2}$, $...$, $B_{p}$, $B_{p+2}$, $B_{p+3}$, $...$, $B_{u+1}$ (since
$p+q-1=u+1$ (since $p+q=u+2$)). In other words, the matrix $B\left[  1:p+1\mid
p+2:p+q\right]  $ is the matrix $\left(  w_{1}\mid w_{2}\mid...\mid w_{p}\mid
w_{p+2}\mid w_{p+3}\mid...\mid w_{u+1}\right)  $ (since $\left(  w_{1}\mid
w_{2}\mid...\mid w_{p}\mid w_{p+2}\mid w_{p+3}\mid...\mid w_{u+1}\right)  $ is
the matrix whose columns (from left to right) are $B_{1}$, $B_{2}$, $...$,
$B_{p}$, $B_{p+2}$, $B_{p+3}$, $...$, $B_{u+1}$). Hence,%
\[
B\left[  1:p+1\mid p+2:p+q\right]  =\left(  w_{1}\mid w_{2}\mid...\mid
w_{p}\mid w_{p+2}\mid w_{p+3}\mid...\mid w_{u+1}\right)  .
\]
But clearly,%
\begin{align*}
\left(  w_{1}\mid w_{2}\mid...\mid\widehat{w_{p+1}}\mid...\mid w_{u+1}\right)
&  =\left(  w_{1}\mid w_{2}\mid...\mid w_{p}\mid w_{p+2}\mid w_{p+3}%
\mid...\mid w_{u+1}\right) \\
&  =B\left[  1:p+1\mid p+2:p+q\right]  .
\end{align*}
Taking determinants on both sides of this equation, we obtain
(\ref{pf.pluecker.ptolemy.2ndproof.3c2}).}

Now, (\ref{pf.pluecker.ptolemy.2ndproof.2}) becomes%
\begin{align*}
0  &  =\left(  -1\right)  ^{1}\underbrace{\det\left(  A\mid w_{1}\right)
}_{\substack{=\left(  -1\right)  ^{u-1}\det\left(  B\left[  1:p\mid
p+2:p+q+1\right]  \right)  \\\text{(by (\ref{pf.pluecker.ptolemy.2ndproof.3a1}%
))}}}\cdot\underbrace{\det\left(  w_{1}\mid w_{2}\mid...\mid\widehat{w_{1}%
}\mid...\mid w_{u+1}\right)  }_{\substack{=\det\left(  B\left[  2:p+1\mid
p+1:p+q\right]  \right)  \\\text{(by (\ref{pf.pluecker.ptolemy.2ndproof.3a2}%
))}}}\\
&  \ \ \ \ \ \ \ \ \ \ +\left(  -1\right)  ^{p}\underbrace{\det\left(  A\mid
w_{p}\right)  }_{\substack{=\left(  -1\right)  ^{q-1}\det\left(  B\left[
2:p+1\mid p+2:p+q+1\right]  \right)  \\\text{(by
(\ref{pf.pluecker.ptolemy.2ndproof.3b1}))}}}\cdot\underbrace{\det\left(
w_{1}\mid w_{2}\mid...\mid\widehat{w_{p}}\mid...\mid w_{u+1}\right)
}_{\substack{=\det\left(  B\left[  1:p\mid p+1:p+q\right]  \right)
\\\text{(by (\ref{pf.pluecker.ptolemy.2ndproof.3b2}))}}}\\
&  \ \ \ \ \ \ \ \ \ \ +\left(  -1\right)  ^{p+1}\underbrace{\det\left(  A\mid
w_{p+1}\right)  }_{\substack{=\left(  -1\right)  ^{q-1}\det\left(  B\left[
2:p\mid p+1:p+q+1\right]  \right)  \\\text{(by
(\ref{pf.pluecker.ptolemy.2ndproof.3c1}))}}}\cdot\underbrace{\det\left(
w_{1}\mid w_{2}\mid...\mid\widehat{w_{p+1}}\mid...\mid w_{u+1}\right)
}_{\substack{=\det\left(  B\left[  1:p+1\mid p+2:p+q\right]  \right)
\\\text{(by (\ref{pf.pluecker.ptolemy.2ndproof.3c2}))}}}\\
&  =\underbrace{\left(  -1\right)  ^{1}\left(  -1\right)  ^{u-1}}_{=\left(
-1\right)  ^{1+\left(  u-1\right)  }=\left(  -1\right)  ^{u}}\det\left(
B\left[  1:p\mid p+2:p+q+1\right]  \right) \\
&  \ \ \ \ \ \ \ \ \ \ \ \ \ \ \ \ \ \ \ \ \cdot\det\left(  B\left[  2:p+1\mid
p+1:p+q\right]  \right) \\
&  \ \ \ \ \ \ \ \ \ \ +\underbrace{\left(  -1\right)  ^{p}\left(  -1\right)
^{q-1}}_{=\left(  -1\right)  ^{p+\left(  q-1\right)  }=\left(  -1\right)
^{p+q-1}=-\left(  -1\right)  ^{u}}\det\left(  B\left[  2:p+1\mid
p+2:p+q+1\right]  \right) \\
&  \ \ \ \ \ \ \ \ \ \ \ \ \ \ \ \ \ \ \ \ \cdot\det\left(  B\left[  1:p\mid
p+1:p+q\right]  \right) \\
&  \ \ \ \ \ \ \ \ \ \ +\underbrace{\left(  -1\right)  ^{p+1}\left(
-1\right)  ^{q-1}}_{=\left(  -1\right)  ^{\left(  p+1\right)  +\left(
q-1\right)  }=\left(  -1\right)  ^{p+q}=\left(  -1\right)  ^{u}}\det\left(
B\left[  2:p\mid p+1:p+q+1\right]  \right) \\
&  \ \ \ \ \ \ \ \ \ \ \ \ \ \ \ \ \ \ \ \ \cdot\det\left(  B\left[  1:p+1\mid
p+2:p+q\right]  \right) \\
&  =\left(  -1\right)  ^{u}\det\left(  B\left[  1:p\mid p+2:p+q+1\right]
\right)  \cdot\det\left(  B\left[  2:p+1\mid p+1:p+q\right]  \right) \\
&  \ \ \ \ \ \ \ \ \ \ -\left(  -1\right)  ^{u}\det\left(  B\left[  2:p+1\mid
p+2:p+q+1\right]  \right)  \cdot\det\left(  B\left[  1:p\mid p+1:p+q\right]
\right) \\
&  \ \ \ \ \ \ \ \ \ \ +\left(  -1\right)  ^{u}\det\left(  B\left[  2:p\mid
p+1:p+q+1\right]  \right)  \cdot\det\left(  B\left[  1:p+1\mid p+2:p+q\right]
\right)  .
\end{align*}
Dividing this equality by $\left(  -1\right)  ^{u}$, we obtain%
\begin{align*}
0  &  =\det\left(  B\left[  1:p\mid p+2:p+q+1\right]  \right)  \cdot
\det\left(  B\left[  2:p+1\mid p+1:p+q\right]  \right) \\
&  \ \ \ \ \ \ \ \ \ \ -\det\left(  B\left[  2:p+1\mid p+2:p+q+1\right]
\right)  \cdot\det\left(  B\left[  1:p\mid p+1:p+q\right]  \right) \\
&  \ \ \ \ \ \ \ \ \ \ +\det\left(  B\left[  2:p\mid p+1:p+q+1\right]
\right)  \cdot\det\left(  B\left[  1:p+1\mid p+2:p+q\right]  \right)  .
\end{align*}
Adding $\det\left(  B\left[  2:p+1\mid p+2:p+q+1\right]  \right)  \cdot
\det\left(  B\left[  1:p\mid p+1:p+q\right]  \right)  $ to this equality, we
obtain%
\begin{align*}
&  \det\left(  B\left[  2:p+1\mid p+2:p+q+1\right]  \right)  \cdot\det\left(
B\left[  1:p\mid p+1:p+q\right]  \right) \\
&  =\det\left(  B\left[  1:p\mid p+2:p+q+1\right]  \right)  \cdot\det\left(
B\left[  2:p+1\mid p+1:p+q\right]  \right) \\
&  \ \ \ \ \ \ \ \ \ \ +\det\left(  B\left[  2:p\mid p+1:p+q+1\right]
\right)  \cdot\det\left(  B\left[  1:p+1\mid p+2:p+q\right]  \right)  .
\end{align*}
Hence,%
\begin{align*}
&  \det\left(  B\left[  1:p\mid p+2:p+q+1\right]  \right)  \cdot\det\left(
B\left[  2:p+1\mid p+1:p+q\right]  \right) \\
&  \ \ \ \ \ \ \ \ \ \ +\det\left(  B\left[  2:p\mid p+1:p+q+1\right]
\right)  \cdot\det\left(  B\left[  1:p+1\mid p+2:p+q\right]  \right) \\
&  =\det\left(  B\left[  2:p+1\mid p+2:p+q+1\right]  \right)  \cdot\det\left(
B\left[  1:p\mid p+1:p+q\right]  \right) \\
&  =\det\left(  B\left[  1:p\mid p+1:p+q\right]  \right)  \cdot\det\left(
B\left[  2:p+1\mid p+2:p+q+1\right]  \right)  .
\end{align*}
This proves Lemma \ref{lem.pluecker.ptolemy}.
\end{proof}
\end{verlong}

\begin{proof}
[Proof of Proposition \ref{prop.Grasp.GraspR}.]Let $f=\operatorname*{Grasp}%
\nolimits_{j+1}A$ and $g=\operatorname*{Grasp}\nolimits_{j}A$.

Clearly, $f\left(  0\right)  =1=g\left(  0\right)  $ and $f\left(  1\right)
=1=g\left(  1\right)  $.

We want to show that $\operatorname*{Grasp}\nolimits_{j}%
A=R_{\operatorname*{Rect}\left(  p,q\right)  }\left(  \operatorname*{Grasp}%
\nolimits_{j+1}A\right)  $. In other words, we want to show that
$g=R_{\operatorname*{Rect}\left(  p,q\right)  }\left(  f\right)  $ (because
$g=\operatorname*{Grasp}\nolimits_{j}A$ and $f=\operatorname*{Grasp}%
\nolimits_{j+1}A$). According to Proposition \ref{prop.R.implicit.converse}
(applied to $P=\operatorname*{Rect}\left(  p,q\right)  $), this will follow
once we can show that%
\begin{equation}
g\left(  v\right)  =\dfrac{1}{f\left(  v\right)  }\cdot\dfrac{\sum
\limits_{\substack{u\in\widehat{\operatorname*{Rect}\left(  p,q\right)
};\\u\lessdot v}}f\left(  u\right)  }{\sum\limits_{\substack{u\in
\widehat{\operatorname*{Rect}\left(  p,q\right)  };\\u\gtrdot v}}\dfrac
{1}{g\left(  u\right)  }}\ \ \ \ \ \ \ \ \ \ \text{for every }v\in
\operatorname*{Rect}\left(  p,q\right)  . \label{pf.Grasp.GraspR.goal}%
\end{equation}

So let $v\in\operatorname*{Rect}\left(  p,q\right)  $. Thus, $v=\left(
i,k\right)  $ for some $i\in\left\{  1,2,...,p\right\}  $ and $k\in\left\{
1,2,...,q\right\}  $. Consider these $i$ and $k$. We must prove
(\ref{pf.Grasp.GraspR.goal}).

We are clearly in one of the following four cases:

\textit{Case 1:} We have $v\neq\left(  1,1\right)  $ and $v\neq\left(
p,q\right)  $.

\textit{Case 2:} We have $v=\left(  1,1\right)  $ and $v\neq\left(
p,q\right)  $.

\textit{Case 3:} We have $v\neq\left(  1,1\right)  $ and $v=\left(
p,q\right)  $.

\textit{Case 4:} We have $v=\left(  1,1\right)  $ and $v=\left(  p,q\right)  $.

Let us consider Case 1 first. In this case, we have $v\neq\left(  1,1\right)
$ and $v\neq\left(  p,q\right)  $. As a consequence, all elements
$u\in\widehat{\operatorname*{Rect}\left(  p,q\right)  }$ satisfying $u\lessdot
v$ belong to $\operatorname*{Rect}\left(  p,q\right)  $, and the same holds
for all $u\in\widehat{\operatorname*{Rect}\left(  p,q\right)  }$ satisfying
$u\gtrdot v$.

Due to the specific form of the poset $\operatorname*{Rect}\left(  p,q\right)
$, there are at most two elements $u$ of $\widehat{\operatorname*{Rect}\left(
p,q\right)  }$ satisfying $u\lessdot v$, namely $\left(  i,k-1\right)  $
(which exists only if $k\neq1$) and $\left(  i-1,k\right)  $ (which exists
only if $i\neq1$). Hence, the sum $\sum\limits_{\substack{u\in
\widehat{\operatorname*{Rect}\left(  p,q\right)  };\\u\lessdot v}}f\left(
u\right)  $ takes one of the three forms $f\left(  \left(  i,k-1\right)
\right)  +f\left(  \left(  i-1,k\right)  \right)  $, $f\left(  \left(
i,k-1\right)  \right)  $ and $f\left(  \left(  i-1,k\right)  \right)  $. Due
to Definition \ref{def.Grasp} \textbf{(b)}, all of these three forms can be
rewritten uniformly as $f\left(  \left(  i,k-1\right)  \right)  +f\left(
\left(  i-1,k\right)  \right)  $ (because if $\left(  i,k-1\right)
\notin\operatorname*{Rect}\left(  p,q\right)  $ then Definition
\ref{def.Grasp} \textbf{(b)} guarantees that $f\left(  \left(  i,k-1\right)
\right)  =0$, and similarly $f\left(  \left(  i-1,k\right)  \right)  =0$ if
$\left(  i-1,k\right)  \notin\operatorname*{Rect}\left(  p,q\right)  $). So we
have%
\begin{equation}
\sum\limits_{\substack{u\in\widehat{\operatorname*{Rect}\left(  p,q\right)
};\\u\lessdot v}}f\left(  u\right)  =f\left(  \left(  i,k-1\right)  \right)
+f\left(  \left(  i-1,k\right)  \right)  . \label{pf.Grasp.GraspR.f}%
\end{equation}
Similarly,
\begin{equation}
\sum\limits_{\substack{u\in\widehat{\operatorname*{Rect}\left(  p,q\right)
};\\u\gtrdot v}}\dfrac{1}{g\left(  u\right)  }=\dfrac{1}{g\left(  \left(
i,k+1\right)  \right)  }+\dfrac{1}{g\left(  \left(  i+1,k\right)  \right)  },
\label{pf.Grasp.GraspR.g}%
\end{equation}
where we set $\dfrac{1}{\infty}=0$ as usual.

But $f=\operatorname*{Grasp}\nolimits_{j+1}A$. Hence,%
\begin{align*}
&  f\left(  \left(  i,k-1\right)  \right) \\
&  =\left(  \operatorname*{Grasp}\nolimits_{j+1}A\right)  \left(  \left(
i,k-1\right)  \right) \\
&  =\dfrac{\det\left(  A\left[  \left(  j+1\right)  +1:\left(  j+1\right)
+i\mid\left(  j+1\right)  +i+\left(  k-1\right)  -1:\left(  j+1\right)
+p+\left(  k-1\right)  \right]  \right)  }{\det\left(  A\left[  j+1:\left(
j+1\right)  +i\mid\left(  j+1\right)  +i+\left(  k-1\right)  :\left(
j+1\right)  +p+\left(  k-1\right)  \right]  \right)  }\\
&  \ \ \ \ \ \ \ \ \ \ \left(  \text{by the definition of }%
\operatorname*{Grasp}\nolimits_{j+1}A\right) \\
&  =\dfrac{\det\left(  A\left[  j+2:j+i+1\mid j+i+k-1:j+p+k\right]  \right)
}{\det\left(  A\left[  j+1:j+i+1\mid j+i+k:j+p+k\right]  \right)  }%
\end{align*}
and
\begin{align*}
&  f\left(  \left(  i-1,k\right)  \right) \\
&  =\left(  \operatorname*{Grasp}\nolimits_{j+1}A\right)  \left(  \left(
i-1,k\right)  \right) \\
&  =\dfrac{\det\left(  A\left[  \left(  j+1\right)  +1:\left(  j+1\right)
+\left(  i-1\right)  \mid\left(  j+1\right)  +\left(  i-1\right)  +k-1:\left(
j+1\right)  +p+k\right]  \right)  }{\det\left(  A\left[  j+1:\left(
j+1\right)  +\left(  i-1\right)  \mid\left(  j+1\right)  +\left(  i-1\right)
+k:\left(  j+1\right)  +p+k\right]  \right)  }\\
&  \ \ \ \ \ \ \ \ \ \ \left(  \text{by the definition of }%
\operatorname*{Grasp}\nolimits_{j+1}A\right) \\
&  =\dfrac{\det\left(  A\left[  j+2:j+i\mid j+i+k-1:j+p+k+1\right]  \right)
}{\det\left(  A\left[  j+1:j+i\mid j+i+k:j+p+k+1\right]  \right)  }.
\end{align*}
Due to these two equalities, (\ref{pf.Grasp.GraspR.f}) becomes%
\begin{align}
&  \sum\limits_{\substack{u\in\widehat{\operatorname*{Rect}\left(  p,q\right)
};\\u \lessdot v}}f\left(  u\right) \nonumber\\
&  =\dfrac{\det\left(  A\left[  j+2:j+i+1\mid j+i+k-1:j+p+k\right]  \right)
}{\det\left(  A\left[  j+1:j+i+1\mid j+i+k:j+p+k\right]  \right)  }\nonumber\\
&  \ \ \ \ \ \ \ \ \ \ +\dfrac{\det\left(  A\left[  j+2:j+i\mid
j+i+k-1:j+p+k+1\right]  \right)  }{\det\left(  A\left[  j+1:j+i\mid
j+i+k:j+p+k+1\right]  \right)  }\nonumber\\
&  =\left(  \det\left(  A\left[  j+1:j+i+1\mid j+i+k:j+p+k\right]  \right)
\right)  ^{-1}\nonumber\\
&  \ \ \ \ \ \ \ \ \ \ \cdot\left(  \det\left(  A\left[  j+1:j+i\mid
j+i+k:j+p+k+1\right]  \right)  \right)  ^{-1}\nonumber\\
&  \ \ \ \ \ \ \ \ \ \ \cdot\left(  \det\left(  A\left[  j+1:j+i\mid
j+i+k:j+p+k+1\right]  \right)  \right. \nonumber\\
&  \ \ \ \ \ \ \ \ \ \ \left.  \ \ \ \ \ \ \ \ \ \ \cdot\det\left(  A\left[
j+2:j+i+1\mid j+i+k-1:j+p+k\right]  \right)  \right. \nonumber\\
&  \ \ \ \ \ \ \ \ \ \ \left.  +\det\left(  A\left[  j+2:j+i\mid
j+i+k-1:j+p+k+1\right]  \right)  \right. \nonumber\\
&  \ \ \ \ \ \ \ \ \ \ \left.  \ \ \ \ \ \ \ \ \ \ \cdot\det\left(  A\left[
j+1:j+i+1\mid j+i+k:j+p+k\right]  \right)  \right) \nonumber\\
&  =\left(  \det\left(  A\left[  j+1:j+i+1\mid j+i+k:j+p+k\right]  \right)
\right)  ^{-1}\nonumber\\
&  \ \ \ \ \ \ \ \ \ \ \cdot\left(  \det\left(  A\left[  j+1:j+i\mid
j+i+k:j+p+k+1\right]  \right)  \right)  ^{-1}\nonumber\\
&  \ \ \ \ \ \ \ \ \ \ \cdot\det\left(  A\left[  j+1:j+i\mid
j+i+k-1:j+p+k\right]  \right) \nonumber\\
&  \ \ \ \ \ \ \ \ \ \ \cdot\det\left(  A\left[  j+2:j+i+1\mid
j+i+k:j+p+k+1\right]  \right)  \label{pf.Grasp.GraspR.side1}%
\end{align}
(because applying Theorem \ref{thm.pluecker.ptolemy} to $a=j+2$, $b=j+i$,
$c=j+i+k$ and $d=j+p+k$ yields%
\begin{align*}
&  \det\left(  A\left[  j+1:j+i\mid j+i+k:j+p+k+1\right]  \right) \\
&  \ \ \ \ \ \ \ \ \ \ \cdot\det\left(  A\left[  j+2:j+i+1\mid
j+i+k-1:j+p+k\right]  \right) \\
&  +\det\left(  A\left[  j+2:j+i\mid j+i+k-1:j+p+k+1\right]  \right) \\
&  \ \ \ \ \ \ \ \ \ \ \cdot\det\left(  A\left[  j+1:j+i+1\mid
j+i+k:j+p+k\right]  \right) \\
&  =\det\left(  A\left[  j+1:j+i\mid j+i+k-1:j+p+k\right]  \right) \\
&  \ \ \ \ \ \ \ \ \ \ \cdot\det\left(  A\left[  j+2:j+i+1\mid
j+i+k:j+p+k+1\right]  \right)
\end{align*}
).

On the other hand, $g=\operatorname*{Grasp}\nolimits_{j}A$, so that%
\begin{align*}
&  g\left(  \left(  i,k+1\right)  \right) \\
&  =\left(  \operatorname*{Grasp}\nolimits_{j}A\right)  \left(  \left(
i,k+1\right)  \right)  =\dfrac{\det\left(  A\left[  j+1:j+i\mid j+i+\left(
k+1\right)  -1:j+p+\left(  k+1\right)  \right]  \right)  }{\det\left(
A\left[  j:j+i\mid j+i+\left(  k+1\right)  :j+p+\left(  k+1\right)  \right]
\right)  }\\
&  \ \ \ \ \ \ \ \ \ \ \left(  \text{by the definition of }%
\operatorname*{Grasp}\nolimits_{j}A\right) \\
&  =\dfrac{\det\left(  A\left[  j+1:j+i\mid j+i+k:j+p+k+1\right]  \right)
}{\det\left(  A\left[  j:j+i\mid j+i+k+1:j+p+k+1\right]  \right)  }%
\end{align*}
and therefore%
\begin{equation}
\dfrac{1}{g\left(  \left(  i,k+1\right)  \right)  }=\dfrac{\det\left(
A\left[  j:j+i\mid j+i+k+1:j+p+k+1\right]  \right)  }{\det\left(  A\left[
j+1:j+i\mid j+i+k:j+p+k+1\right]  \right)  }. \label{pf.Grasp.GraspR.g.1}%
\end{equation}
Also, from $g=\operatorname*{Grasp}\nolimits_{j}A$, we obtain%
\begin{align*}
&  g\left(  \left(  i+1,k\right)  \right) \\
&  =\left(  \operatorname*{Grasp}\nolimits_{j}A\right)  \left(  \left(
i-1,k\right)  \right)  =\dfrac{\det\left(  A\left[  j+1:j+\left(  i+1\right)
\mid j+\left(  i+1\right)  +k-1:j+p+k\right]  \right)  }{\det\left(  A\left[
j:j+\left(  i+1\right)  \mid j+\left(  i+1\right)  +k:j+p+k\right]  \right)
}\\
&  \ \ \ \ \ \ \ \ \ \ \left(  \text{by the definition of }%
\operatorname*{Grasp}\nolimits_{j}A\right) \\
&  =\dfrac{\det\left(  A\left[  j+1:j+i+1\mid j+i+k:j+p+k\right]  \right)
}{\det\left(  A\left[  j:j+i+1\mid j+i+k+1:j+p+k\right]  \right)  },
\end{align*}
so that%
\begin{equation}
\dfrac{1}{g\left(  \left(  i+1,k\right)  \right)  }=\dfrac{\det\left(
A\left[  j:j+i+1\mid j+i+k+1:j+p+k\right]  \right)  }{\det\left(  A\left[
j+1:j+i+1\mid j+i+k:j+p+k\right]  \right)  }. \label{pf.Grasp.GraspR.g.2}%
\end{equation}
Due to (\ref{pf.Grasp.GraspR.g.1}) and (\ref{pf.Grasp.GraspR.g.2}), the
equality (\ref{pf.Grasp.GraspR.g}) becomes%
\begin{align}
&  \sum\limits_{\substack{u\in\widehat{\operatorname*{Rect}\left(  p,q\right)
};\\u\gtrdot v}}\dfrac{1}{g\left(  u\right)  }\nonumber\\
&  =\dfrac{\det\left(  A\left[  j:j+i\mid j+i+k+1:j+p+k+1\right]  \right)
}{\det\left(  A\left[  j+1:j+i\mid j+i+k:j+p+k+1\right]  \right)  }\nonumber\\
&  \ \ \ \ \ \ \ \ \ \ +\dfrac{\det\left(  A\left[  j:j+i+1\mid
j+i+k+1:j+p+k\right]  \right)  }{\det\left(  A\left[  j+1:j+i+1\mid
j+i+k:j+p+k\right]  \right)  }\nonumber\\
&  =\left(  \det\left(  A\left[  j+1:j+i\mid j+i+k:j+p+k+1\right]  \right)
\right)  ^{-1}\nonumber\\
&  \ \ \ \ \ \ \ \ \ \ \cdot\left(  \det\left(  A\left[  j+1:j+i+1\mid
j+i+k:j+p+k\right]  \right)  \right)  ^{-1}\nonumber\\
&  \ \ \ \ \ \ \ \ \ \ \cdot\left(  \det\left(  A\left[  j:j+i\mid
j+i+k+1:j+p+k+1\right]  \right)  \right. \nonumber\\
&  \ \ \ \ \ \ \ \ \ \ \left.  \ \ \ \ \ \ \ \ \ \ \cdot\det\left(  A\left[
j+1:j+i+1\mid j+i+k:j+p+k\right]  \right)  \right. \nonumber\\
&  \ \ \ \ \ \ \ \ \ \ \left.  +\det\left(  A\left[  j+1:j+i\mid
j+i+k:j+p+k+1\right]  \right)  \right. \nonumber\\
&  \ \ \ \ \ \ \ \ \ \ \left.  \ \ \ \ \ \ \ \ \ \ \cdot\det\left(  A\left[
j:j+i+1\mid j+i+k+1:j+p+k\right]  \right)  \right) \nonumber\\
&  =\left(  \det\left(  A\left[  j+1:j+i\mid j+i+k:j+p+k+1\right]  \right)
\right)  ^{-1}\nonumber\\
&  \ \ \ \ \ \ \ \ \ \ \cdot\left(  \det\left(  A\left[  j+1:j+i+1\mid
j+i+k:j+p+k\right]  \right)  \right)  ^{-1}\nonumber\\
&  \ \ \ \ \ \ \ \ \ \ \cdot\det\left(  A\left[  j:j+i\mid j+i+k:j+p+k\right]
\right) \nonumber\\
&  \ \ \ \ \ \ \ \ \ \ \cdot\det\left(  A\left[  j+1:j+i+1\mid
j+i+k+1:j+p+k+1\right]  \right)  \label{pf.Grasp.GraspR.side2}%
\end{align}
(because applying Theorem \ref{thm.pluecker.ptolemy} to $a=j+1$, $b=j+i$,
$c=j+i+k+1$ and $d=j+p+k$ yields%
\begin{align*}
&  \det\left(  A\left[  j:j+i\mid j+i+k+1:j+p+k+1\right]  \right) \\
&  \ \ \ \ \ \ \ \ \ \ \cdot\det\left(  A\left[  j+1:j+i+1\mid
j+i+k:j+p+k\right]  \right) \\
&  +\det\left(  A\left[  j+1:j+i\mid j+i+k:j+p+k+1\right]  \right) \\
&  \ \ \ \ \ \ \ \ \ \ \cdot\det\left(  A\left[  j:j+i+1\mid
j+i+k+1:j+p+k\right]  \right) \\
&  =\det\left(  A\left[  j:j+i\mid j+i+k:j+p+k\right]  \right) \\
&  \ \ \ \ \ \ \ \ \ \ \cdot\det\left(  A\left[  j+1:j+i+1\mid
j+i+k+1:j+p+k+1\right]  \right)
\end{align*}
).

Since $v=\left(  i,k\right)  $ and $g=\operatorname*{Grasp}\nolimits_{j}A$, we
have%
\begin{align}
&  g\left(  v\right) \nonumber\\
&  =\left(  \operatorname*{Grasp}\nolimits_{j}A\right)  \left(  \left(
i,k\right)  \right)  =\dfrac{\det\left(  A\left[  j+1:j+i\mid
j+i+k-1:j+p+k\right]  \right)  }{\det\left(  A\left[  j:j+i\mid
j+i+k:j+p+k\right]  \right)  }\label{pf.Grasp.GraspR.side3}\\
&  \ \ \ \ \ \ \ \ \ \ \left(  \text{by the definition of }%
\operatorname*{Grasp}\nolimits_{j}A\right)  .\nonumber
\end{align}

Since $v=\left(  i,k\right)  $ and $f=\operatorname*{Grasp}\nolimits_{j+1}A$,
we have%
\begin{align}
&  f\left(  v\right) \nonumber\\
&  =\left(  \operatorname*{Grasp}\nolimits_{j+1}A\right)  \left(  \left(
i,k\right)  \right) \nonumber\\
&  =\dfrac{\det\left(  A\left[  \left(  j+1\right)  +1:\left(  j+1\right)
+i\mid\left(  j+1\right)  +i+k-1:\left(  j+1\right)  +p+k\right]  \right)
}{\det\left(  A\left[  j+1:\left(  j+1\right)  +i\mid\left(  j+1\right)
+i+k:\left(  j+1\right)  +p+k\right]  \right)  }\nonumber\\
&  \ \ \ \ \ \ \ \ \ \ \left(  \text{by the definition of }%
\operatorname*{Grasp}\nolimits_{j+1}A\right) \nonumber\\
&  =\dfrac{\det\left(  A\left[  j+2:j+i+1\mid j+i+k:j+p+k+1\right]  \right)
}{\det\left(  A\left[  j+1:j+i+1\mid j+i+k+1:j+p+k+1\right]  \right)  }.
\label{pf.Grasp.GraspR.side4}%
\end{align}

Now, we can rewrite the terms $\sum\limits_{\substack{u\in
\widehat{\operatorname*{Rect}\left(  p,q\right)  };\\u \lessdot v}}f\left(
u\right)  $, $\sum\limits_{\substack{u\in\widehat{\operatorname*{Rect}\left(
p,q\right)  };\\u \gtrdot v}}\dfrac{1}{g\left(  u\right)  }$, $g\left(
v\right)  $ and $f\left(  v\right)  $ in (\ref{pf.Grasp.GraspR.goal}) using
the equalities (\ref{pf.Grasp.GraspR.side1}), (\ref{pf.Grasp.GraspR.side2}),
(\ref{pf.Grasp.GraspR.side3}) and (\ref{pf.Grasp.GraspR.side4}), respectively.
The resulting equation is a tautology because all determinants cancel out
(this can be checked by the reader). Hence, (\ref{pf.Grasp.GraspR.goal}) is
proven in Case 1.

Let us now consider Case 3. In this case, we have $v\neq\left(  1,1\right)  $
and $v=\left(  p,q\right)  $. Hence, (\ref{pf.Grasp.GraspR.side1}),
(\ref{pf.Grasp.GraspR.side3}) and (\ref{pf.Grasp.GraspR.side4}) are still
valid, whereas (\ref{pf.Grasp.GraspR.side2}) gets superseded by the simpler
equality%
\begin{equation}
\sum\limits_{\substack{u\in\widehat{\operatorname*{Rect}\left(  p,q\right)
};\\u\gtrdot v}}\dfrac{1}{g\left(  u\right)  }=\dfrac{1}{g\left(  1\right)
}=\dfrac{1}{1}=1. \label{pf.Grasp.GraspR.side2a}%
\end{equation}

Now, we can rewrite the terms $\sum\limits_{\substack{u\in
\widehat{\operatorname*{Rect}\left(  p,q\right)  };\\u\lessdot v}}f\left(
u\right)  $, $\sum\limits_{\substack{u\in\widehat{\operatorname*{Rect}\left(
p,q\right)  };\\u\gtrdot v}}\dfrac{1}{g\left(  u\right)  }$, $g\left(
v\right)  $ and $f\left(  v\right)  $ in (\ref{pf.Grasp.GraspR.goal}) using
the equalities (\ref{pf.Grasp.GraspR.side1}), (\ref{pf.Grasp.GraspR.side2a}),
(\ref{pf.Grasp.GraspR.side3}) and (\ref{pf.Grasp.GraspR.side4}), respectively.
The resulting equation (after multiplying through with all denominators and
cancelling terms appearing on both sides) simplifies to%
\begin{align*}
&  \det\left(  A\left[  j+1:j+i+1\mid j+i+k:j+p+k\right]  \right) \\
&  \ \ \ \ \ \ \ \ \ \ \cdot\det\left(  A\left[  j+1:j+i\mid
j+i+k:j+p+k+1\right]  \right) \\
&  =\det\left(  A\left[  j+1:j+i+1\mid j+i+k+1:j+p+k+1\right]  \right) \\
&  \ \ \ \ \ \ \ \ \ \ \cdot\det\left(  A\left[  j:j+i\mid j+i+k:j+p+k\right]
\right)  .
\end{align*}
Since $i=p$ and $k=q$ (because $\left(  i,k\right)  =v=\left(  p,q\right)  $),
this rewrites as%
\begin{align*}
&  \det\left(  A\left[  j+1:j+p+1\mid j+p+q:j+p+q\right]  \right) \\
&  \ \ \ \ \ \ \ \ \ \ \cdot\det\left(  A\left[  j+1:j+p\mid
j+p+q:j+p+q+1\right]  \right) \\
&  =\det\left(  A\left[  j+1:j+p+1\mid j+p+q+1:j+p+q+1\right]  \right) \\
&  \ \ \ \ \ \ \ \ \ \ \cdot\det\left(  A\left[  j:j+p\mid j+p+q:j+p+q\right]
\right)  .
\end{align*}
But this follows from%
\begin{align*}
&  \det\left(  A\left[  j+1:j+p+1\mid j+p+q:j+p+q\right]  \right) \\
&  =\det\left(  A\left[  j+1:j+p+1\mid j+p+q+1:j+p+q+1\right]  \right)
\end{align*}
(which is clear from Proposition \ref{prop.minors.trivial} \textbf{(b)}) and
\begin{align*}
&  \det\left(  A\left[  j+1:j+p\mid j+p+q:j+p+q+1\right]  \right) \\
&  =\det\left(  A\left[  j:j+p\mid j+p+q:j+p+q\right]  \right)
\end{align*}
(which can be easily proven\footnote{\textit{Proof.} We have%
\begin{align*}
&  \det\left(  A\left[  j+1:j+p\mid j+p+q:j+p+q+1\right]  \right) \\
&  =\det\left(  A\left[  j+1:j+p\mid p+q+j:p+q+j+1\right]  \right)
=\det\left(  \underbrace{A\left[  j:j+1\mid j+1:j+p\right]  }%
_{\substack{=A\left[  j:j+p\mid j+p+q:j+p+q\right]  \\\text{(by Proposition
\ref{prop.minors.trivial} \textbf{(c)})}}}\right) \\
&  \ \ \ \ \ \ \ \ \ \ \left(  \text{by Proposition \ref{prop.minors.period}
\textbf{(c)}, applied to }u=p\text{, }v=p+q\text{, }a=j+1\text{, }b=j+p\text{,
}c=j\text{ and }d=j+1\right) \\
&  =\det\left(  A\left[  j:j+p\mid j+p+q:j+p+q\right]  \right)  ,
\end{align*}
qed.}). Thus, (\ref{pf.Grasp.GraspR.goal}) is proven in Case 3.

Let us next consider Case 2. In this case, we have $v=\left(  1,1\right)  $
and $v\neq\left(  p,q\right)  $. Hence, (\ref{pf.Grasp.GraspR.side2}),
(\ref{pf.Grasp.GraspR.side3}) and (\ref{pf.Grasp.GraspR.side4}) are still
valid, whereas (\ref{pf.Grasp.GraspR.side1}) gets superseded by the simpler
equality%
\begin{equation}
\sum\limits_{\substack{u\in\widehat{\operatorname*{Rect}\left(  p,q\right)
};\\u\lessdot v}}f\left(  u\right)  =f\left(  0\right)  =1.
\label{pf.Grasp.GraspR.side1a}%
\end{equation}
Now, we can rewrite the terms $\sum\limits_{\substack{u\in
\widehat{\operatorname*{Rect}\left(  p,q\right)  };\\u\lessdot v}}f\left(
u\right)  $, $\sum\limits_{\substack{u\in\widehat{\operatorname*{Rect}\left(
p,q\right)  };\\u\gtrdot v}}\dfrac{1}{g\left(  u\right)  }$, $g\left(
v\right)  $ and $f\left(  v\right)  $ in (\ref{pf.Grasp.GraspR.goal}) using
the equalities (\ref{pf.Grasp.GraspR.side1a}), (\ref{pf.Grasp.GraspR.side2}),
(\ref{pf.Grasp.GraspR.side3}) and (\ref{pf.Grasp.GraspR.side4}), respectively.
The resulting equation (after multiplying through with all denominators and
cancelling terms appearing on both sides) simplifies to
\begin{align*}
&  \det\left(  A\left[  j+1:j+i\mid j+i+k-1:j+p+k\right]  \right) \\
&  \ \ \ \ \ \ \ \ \ \ \cdot\det\left(  A\left[  j+2:j+i+1\mid
j+i+k:j+p+k+1\right]  \right) \\
&  =\det\left(  A\left[  j+1:j+i+1\mid j+i+k:j+p+k\right]  \right) \\
&  \ \ \ \ \ \ \ \ \ \ \cdot\det\left(  A\left[  j+1:j+i\mid
j+i+k:j+p+k+1\right]  \right)  .
\end{align*}
Since $i=1$ and $k=1$ (because $\left(  i,k\right)  =v=\left(  1,1\right)  $),
this rewrites as%
\begin{align*}
&  \det\left(  A\left[  j+1:j+1\mid j+1+1-1:j+p+1\right]  \right) \\
&  \ \ \ \ \ \ \ \ \ \ \cdot\det\left(  A\left[  j+2:j+1+1\mid
j+1+1:j+p+1+1\right]  \right) \\
&  =\det\left(  A\left[  j+1:j+1+1\mid j+1+1:j+p+1\right]  \right) \\
&  \ \ \ \ \ \ \ \ \ \ \cdot\det\left(  A\left[  j+1:j+1\mid
j+1+1:j+p+1+1\right]  \right)  .
\end{align*}
In other words, this rewrites as%
\begin{align*}
&  \det\left(  A\left[  j+1:j+1\mid j+1:j+p+1\right]  \right) \\
&  \ \ \ \ \ \ \ \ \ \ \cdot\det\left(  A\left[  j+2:j+2\mid j+2:j+p+2\right]
\right) \\
&  =\det\left(  A\left[  j+1:j+2\mid j+2:j+p+1\right]  \right) \\
&  \ \ \ \ \ \ \ \ \ \ \cdot\det\left(  A\left[  j+1:j+1\mid j+2:j+p+2\right]
\right)  .
\end{align*}
But this trivially follows from%
\[
\det\left(  A\left[  j+1:j+1\mid j+1:j+p+1\right]  \right)  =\det\left(
A\left[  j+1:j+2\mid j+2:j+p+1\right]  \right)
\]
(this is because of Proposition \ref{prop.minors.complete}) and%
\[
\det\left(  A\left[  j+2:j+2\mid j+2:j+p+2\right]  \right)  =\det\left(
A\left[  j+1:j+1\mid j+2:j+p+2\right]  \right)
\]
(this is because of Proposition \ref{prop.minors.trivial} \textbf{(a)}). This
proves (\ref{pf.Grasp.GraspR.goal}) in Case 2.

We have now proven (\ref{pf.Grasp.GraspR.goal}) in each of the Cases 1, 2 and
3. We leave the proof in Case 4 to the reader (this case is completely
straightforward, since it has $\left(  p,q\right)  =v=\left(  1,1\right)  $).
Thus, we now know that (\ref{pf.Grasp.GraspR.goal}) holds in each of the four
Cases 1, 2, 3 and 4. Since these four Cases cover all possibilities, this
yields that (\ref{pf.Grasp.GraspR.goal}) always holds. As we have seen, this
completes the proof of Proposition \ref{prop.Grasp.GraspR}.
\end{proof}

A remark seems in order, about why we paid so much attention to the
``degenerate'' Cases 2, 3 and 4. Indeed, only in Cases 3 and 4 have we used
the fact that the sequence $\left(  A_{n}\right)  _{n\in\mathbb{Z}}$ is
``$\left(  p+q\right)  $-periodic up to sign'' rather than just an arbitrary
sequence of length-$p$ column vectors. Had we left out these seemingly
straightforward cases, it would have seemed that the proof showed a result too
good to be true (because it is rather clear that the periodicity in the
definition of $A_{n}$ for general $n\in\mathbb{Z}$ is needed).

\section{\label{sect.dominance}Dominance of the Grassmannian parametrization}

\begin{verlong}
It remains to verify Proposition \ref{prop.Grasp.generic}. Before we do so, a
piece of notation that we should have introduced long ago:

\begin{definition}
\label{def.minors.simple}Let $\mathbb{K}$ be a commutative ring. Let
$A\in\mathbb{K}^{u\times v}$ be a $u\times v$-matrix for some nonnegative
integers $u$ and $v$.

For any integers $a$ and $b$ satisfying $a\leq b$, we let $A\left[
a:b\right]  $ be the matrix whose columns (from left to right) are $A_{a}$,
$A_{a+1}$, $...$, $A_{b-1}$. If $b-a=u$, then this matrix is square and thus
has a determinant $\det\left(  A\left[  a:b\right]  \right)  $.
\end{definition}

As usual, let us state a completely trivial observation:

\begin{proposition}
\label{prop.minors.simple.transi}Let $\mathbb{K}$ be a field. Let
$A\in\mathbb{K}^{u\times v}$ be a $u\times v$-matrix for some nonnegative
integers $u$ and $v$. Let $a$, $b$ and $c$ be integers satisfying $a\leq b\leq
c$. Then, $A\left[  a:c\right]  =A\left[  a:b\mid b:c\right]  $.
\end{proposition}
\end{verlong}

Let us show an example before we start proving Proposition
\ref{prop.Grasp.generic}.

\begin{example}
\label{ex.Grasp.generic}For this example, let $p=2$ and $q=2$. Let
$f\in\mathbb{K}^{\widehat{\operatorname*{Rect}\left(  2,2\right)  }}$ be a
generic reduced labelling. We want to construct a matrix $A\in\mathbb{K}%
^{2\times\left(  2+2\right)  }$ satisfying $f=\operatorname*{Grasp}%
\nolimits_{0}A$.

Clearly, the condition $f=\operatorname*{Grasp}\nolimits_{0}A$ imposes $4$
equations on the entries of $A$ (one for every element of
$\operatorname*{Rect}\left(  2,2\right)  $). Since the matrix $A$ we want to
find has a total of $8$ entries, we are therefore trying to solve an
underdetermined system. However, we can get rid of the superfluous freedom if
we additionally try to ensure that our matrix $A$ has the form $\left(
I_{p}\mid B\right)  $ for some $B\in\mathbb{K}^{2\times2}$ (where $\left(
I_{p}\mid B\right)  $ means the matrix obtained from the $p\times p$ identity
matrix $I_{p}$ by attaching the matrix $B$ to it on the right). Let us do this
now. So we are looking for a matrix $B\in\mathbb{K}^{2\times2}$ satisfying
$f=\operatorname*{Grasp}\nolimits_{0}\left(  I_{p}\mid B\right)  $. This puts
$4$ conditions on $4$ unknowns. Write $B=\left(
\begin{array}
[c]{cc}%
x & y\\
z & w
\end{array}
\right)  $. Then, $\left(  I_{p}\mid B\right)  =\left(
\begin{array}
[c]{cccc}%
1 & 0 & x & y\\
0 & 1 & z & w
\end{array}
\right)  $. Now,%
\begin{align*}
\left(  \operatorname*{Grasp}\nolimits_{0}\left(  I_{p}\mid B\right)  \right)
\left(  \left(  1,1\right)  \right)   &  =\dfrac{\det\left(  \left(  I_{p}\mid
B\right)  \left[  1:1\mid1:3\right]  \right)  }{\det\left(  \left(  I_{p}\mid
B\right)  \left[  0:1\mid2:3\right]  \right)  }=\dfrac{\det\left(
\begin{array}
[c]{cc}%
1 & 0\\
0 & 1
\end{array}
\right)  }{\det\left(
\begin{array}
[c]{cc}%
-y & 0\\
-w & 1
\end{array}
\right)  }=\dfrac{-1}{y};\\
\left(  \operatorname*{Grasp}\nolimits_{0}\left(  I_{p}\mid B\right)  \right)
\left(  \left(  1,2\right)  \right)   &  =\dfrac{\det\left(  \left(  I_{p}\mid
B\right)  \left[  1:1\mid2:4\right]  \right)  }{\det\left(  \left(  I_{p}\mid
B\right)  \left[  0:1\mid3:4\right]  \right)  }=\dfrac{\det\left(
\begin{array}
[c]{cc}%
0 & x\\
1 & z
\end{array}
\right)  }{\det\left(
\begin{array}
[c]{cc}%
-y & x\\
-w & z
\end{array}
\right)  }=\dfrac{-x}{wx-yz};\\
\left(  \operatorname*{Grasp}\nolimits_{0}\left(  I_{p}\mid B\right)  \right)
\left(  \left(  2,1\right)  \right)   &  =\dfrac{\det\left(  \left(  I_{p}\mid
B\right)  \left[  1:2\mid2:3\right]  \right)  }{\det\left(  \left(  I_{p}\mid
B\right)  \left[  0:2\mid3:3\right]  \right)  }=\dfrac{\det\left(
\begin{array}
[c]{cc}%
1 & 0\\
0 & 1
\end{array}
\right)  }{\det\left(
\begin{array}
[c]{cc}%
-y & 1\\
-w & 0
\end{array}
\right)  }=\dfrac{1}{w};\\
\left(  \operatorname*{Grasp}\nolimits_{0}\left(  I_{p}\mid B\right)  \right)
\left(  \left(  2,2\right)  \right)   &  =\dfrac{\det\left(  \left(  I_{p}\mid
B\right)  \left[  1:2\mid3:4\right]  \right)  }{\det\left(  \left(  I_{p}\mid
B\right)  \left[  0:2\mid4:4\right]  \right)  }=\dfrac{\det\left(
\begin{array}
[c]{cc}%
1 & x\\
0 & z
\end{array}
\right)  }{\det\left(
\begin{array}
[c]{cc}%
-y & 1\\
-w & 0
\end{array}
\right)  }=\dfrac{z}{w}.
\end{align*}
The requirement $f=\operatorname*{Grasp}\nolimits_{0}\left(  I_{p}\mid
B\right)  $ therefore translates into the system%
\[
\left\{
\begin{array}
[c]{lcl}%
f\left(  \left(  1,1\right)  \right)  & = & \dfrac{-1}{y};\\
f\left(  \left(  1,2\right)  \right)  & = & \dfrac{-x}{wx-yz};\\
f\left(  \left(  2,1\right)  \right)  & = & \dfrac{1}{w};\\[9pt]%
f\left(  \left(  2,2\right)  \right)  & = & \dfrac{z}{w}%
\end{array}
\right.  .
\]
This system can be solved by elimination: First, compute $w$ using $f\left(
\left(  2,1\right)  \right)  =\dfrac{1}{w}$, obtaining $w=\dfrac{1}{f\left(
\left(  2,1\right)  \right)  }$; then, compute $y$ using $f\left(  \left(
1,1\right)  \right)  =\dfrac{-1}{y}$, obtaining $y=\dfrac{-1}{f\left(  \left(
1,1\right)  \right)  }$; then, compute $z$ using $f\left(  \left(  2,2\right)
\right)  =\dfrac{z}{w}$ and the already eliminated $w$, obtaining
$z=\dfrac{f\left(  \left(  2,2\right)  \right)  }{f\left(  \left(  2,1\right)
\right)  }$; finally, compute $x$ using $f\left(  \left(  1,2\right)  \right)
=\dfrac{-x}{wx-yz}$ and the already eliminated $w,y,z$, obtaining
$x=\dfrac{-f\left(  \left(  1,2\right)  \right)  f\left(  \left(  2,2\right)
\right)  }{\left(  f\left(  \left(  1,2\right)  \right)  +f\left(  \left(
2,1\right)  \right)  \right)  f\left(  \left(  1,1\right)  \right)  }$. While
the denominators in these fractions can vanish, leading to underdetermination
or unsolvability, this will not happen for \textbf{generic} $f$.

This approach to solving $f=\operatorname*{Grasp}\nolimits_{0}A$ generalizes
to arbitrary $p$ and $q$, and motivates the following proof.
\end{example}

\begin{vershort}
We are now going to outline the proof of Proposition \ref{prop.Grasp.generic}.
As shouldn't be surprising after Example \ref{ex.Grasp.generic}, the
underlying idea of the proof is the following: For any fixed $f\in
\mathbb{K}^{\operatorname*{Rect}\left(  p,q\right)  }$, the equation
$f=\operatorname*{Grasp}\nolimits_{0}A$ (with $A$ an unknown matrix in
$\mathbb{K}^{p\times\left(  p+q\right)  }$) can be considered as a system of
$pq$ equations on $p\left(  p+q\right)  $ unknowns (the entries of $A$). While
this system is usually underdetermined, we can restrict the entries of $A$ by
requiring that the leftmost $p$ columns of $A$ form the $p\times p$ identity
matrix. Upon this restriction, we are left with $pq$ unknowns only, and for
$f$ sufficiently generic, the resulting system will be uniquely solvable by
\textquotedblleft triangular elimination\textquotedblright\ (i.e., there is an
equation containing only one unknown; then, when this unknown is eliminated,
the resulting system again contains an equation with only one unknown, and
once this one is eliminated, one gets a further system containing an equation
with only one unknown, and so forth) -- like a triangular system of linear
equations with nonzero entries on the diagonal, but without the linearity. Of
course, this is not a complete proof because the applicability of
\textquotedblleft triangular elimination\textquotedblright\ has to be proven,
not merely claimed. We are only going to sketch the ideas of this proof,
leaving all straightforward details to the reader to fill in. For the sake of
clarity, we are going to word the argument using algebraic properties of
families of rational functions instead of using the algorithmic nature of
\textquotedblleft triangular elimination\textquotedblright\ (similarly to how
most applications of linear algebra use the language of bases of vector spaces
rather than talk about the process of solving systems by Gaussian
elimination). While this clarity comes at the cost of a slight disconnect from
the motivation of the proof, we hope that the reader will still see how the
wind blows.

We first introduce some notation to capture the essence of ``triangular
elimination'' without having to talk about actually moving around variables in equations:

\begin{definition}
\label{def.algebraic.triangularity.short}Let $\mathbb{F}$ be a field. Let
$\mathbf{P}$ be a finite set.

\textbf{(a)} Let $x_{\mathbf{p}}$ be a new symbol for every $\mathbf{p}%
\in\mathbf{P}$. We will denote by $\mathbb{F}\left(  x_{\mathbf{P}}\right)  $
the field of rational functions over $\mathbb{F}$ in the indeterminates
$x_{\mathbf{p}}$ with $\mathbf{p}$ ranging over all elements of $\mathbf{P}$
(hence altogether $\left\vert \mathbf{P}\right\vert $ indeterminates). We also
will denote by $\mathbb{F}\left[  x_{\mathbf{P}}\right]  $ the ring of
polynomials over $\mathbb{F}$ in the indeterminates $x_{\mathbf{p}}$ with
$\mathbf{p}$ ranging over all elements of $\mathbf{P}$. (Thus, $\mathbb{F}%
\left(  x_{\mathbf{P}}\right)  =\mathbb{F}\left(  x_{\mathbf{p}_{1}%
},x_{\mathbf{p}_{2}},...,x_{\mathbf{p}_{n}}\right)  $ and $\mathbb{F}\left[
x_{\mathbf{P}}\right]  =\mathbb{F}\left[  x_{\mathbf{p}_{1}},x_{\mathbf{p}%
_{2}},...,x_{\mathbf{p}_{n}}\right]  $ if $\mathbf{P}$ is written in the form
$\mathbf{P}=\left\{  \mathbf{p}_{1},\mathbf{p}_{2},...,\mathbf{p}_{n}\right\}
$.) The symbols $x_{\mathbf{p}}$ are understood to be distinct, and are used
as commuting indeterminates. We regard $\mathbb{F}\left[  x_{\mathbf{P}%
}\right]  $ as a subring of $\mathbb{F}\left(  x_{\mathbf{P}}\right)  $, and
$\mathbb{F}\left(  x_{\mathbf{P}}\right)  $ as the field of quotients of
$\mathbb{F}\left[  x_{\mathbf{P}}\right]  $.

\textbf{(b)} If $\mathbf{Q}$ is a subset of $\mathbf{P}$, then $\mathbb{F}%
\left(  x_{\mathbf{Q}}\right)  $ can be canonically embedded into
$\mathbb{F}\left(  x_{\mathbf{P}}\right)  $, and $\mathbb{F}\left[
x_{\mathbf{Q}}\right]  $ can be canonically embedded into $\mathbb{F}\left[
x_{\mathbf{P}}\right]  $. We regard these embeddings as inclusions.

\textbf{(c)} Let $\mathbb{K}$ be a field extension of $\mathbb{F}$. Let $f$ be
an element of $\mathbb{F}\left(  x_{\mathbf{P}}\right)  $. If $\left(
a_{\mathbf{p}}\right)  _{\mathbf{p}\in\mathbf{P}}\in\mathbb{K}^{\mathbf{P}}$
is a family of elements of $\mathbb{K}$ indexed by elements of $\mathbf{P}$,
then we let $f\left(  \left(  a_{\mathbf{p}}\right)  _{\mathbf{p}\in
\mathbf{P}}\right)  $ denote the element of $\mathbb{K}$ obtained by
substituting $a_{\mathbf{p}}$ for $x_{\mathbf{p}}$ for each $\mathbf{p}%
\in\mathbf{P}$ in the rational function $f$. This $f\left(  \left(
a_{\mathbf{p}}\right)  _{\mathbf{p}\in\mathbf{P}}\right)  $ is defined only if
the substitution does not render the denominator equal to $0$. If $\mathbb{K}$
is infinite, this shows that $f\left(  \left(  a_{\mathbf{p}}\right)
_{\mathbf{p}\in\mathbf{P}}\right)  $ is defined for almost all $\left(
a_{\mathbf{p}}\right)  _{\mathbf{p}\in\mathbf{P}}\in\mathbb{K}^{\mathbf{P}}$
(with respect to the Zariski topology).

\textbf{(d)} Let $\mathbf{P}$ now be a finite totally ordered set, and let
$\vartriangleleft$ be the smaller relation of $\mathbf{P}$. For every
$\mathbf{p}\in\mathbf{P}$, let $\mathbf{p}\Downarrow$ denote the subset
$\left\{  \mathbf{v}\in\mathbf{P}\ \mid\ \mathbf{v}\vartriangleleft
\mathbf{p}\right\}  $ of $\mathbf{P}$. For every $\mathbf{p}\in\mathbf{P}$,
let $Q_{\mathbf{p}}$ be an element of $\mathbb{F}\left(  x_{\mathbf{P}%
}\right)  $.

We say that the family $\left(  Q_{\mathbf{p}}\right)  _{\mathbf{p}%
\in\mathbf{P}}$ is $\mathbf{P}$\textit{-triangular} if and only if the
following condition holds:

\textit{Algebraic triangularity condition:} For every $\mathbf{p}\in
\mathbf{P}$, there exist elements $\alpha_{\mathbf{p}}$, $\beta_{\mathbf{p}}$,
$\gamma_{\mathbf{p}}$, $\delta_{\mathbf{p}}$ of $\mathbb{F}\left(
x_{\mathbf{p}\Downarrow}\right)  $ such that $\alpha_{\mathbf{p}}%
\delta_{\mathbf{p}}-\beta_{\mathbf{p}}\gamma_{\mathbf{p}}\neq0$ and
$Q_{\mathbf{p}}=\dfrac{\alpha_{\mathbf{p}}x_{\mathbf{p}}+\beta_{\mathbf{p}}%
}{\gamma_{\mathbf{p}}x_{\mathbf{p}}+\delta_{\mathbf{p}}}$.\ \ \ \ \footnotemark
\end{definition}

\footnotetext{Notice that the fraction $\dfrac{\alpha_{\mathbf{p}%
}x_{\mathbf{p}}+\beta_{\mathbf{p}}}{\gamma_{\mathbf{p}}x_{\mathbf{p}}%
+\delta_{\mathbf{p}}}$ is well-defined for any four elements $\alpha
_{\mathbf{p}}$, $\beta_{\mathbf{p}}$, $\gamma_{\mathbf{p}}$, $\delta
_{\mathbf{p}}$ of $\mathbb{F}\left(  x_{\mathbf{p}\Downarrow}\right)  $ such
that $\alpha_{\mathbf{p}}\delta_{\mathbf{p}}-\beta_{\mathbf{p}}\gamma
_{\mathbf{p}}\neq0$. (Indeed, $\gamma_{\mathbf{p}}x_{\mathbf{p}}%
+\delta_{\mathbf{p}}\neq0$ in this case, as can easily be checked.)}

We will use $\mathbf{P}$-triangularity via the following fact:

\begin{lemma}
\label{lem.algebraic.triangularity.short}Let $\mathbb{F}$ be a field. Let
$\mathbf{P}$ be a finite totally ordered set. For every $\mathbf{p}%
\in\mathbf{P}$, let $Q_{\mathbf{p}}$ be an element of $\mathbb{F}\left(
x_{\mathbf{P}}\right)  $. Assume that $\left(  Q_{\mathbf{p}}\right)
_{\mathbf{p}\in\mathbf{P}}$ is a $\mathbf{P}$-triangular family. Then:

\textbf{(a)} The family $\left(  Q_{\mathbf{p}}\right)  _{\mathbf{p}%
\in\mathbf{P}}\in\left(  \mathbb{F}\left(  x_{\mathbf{P}}\right)  \right)
^{\mathbf{P}}$ is algebraically independent (over $\mathbb{F}$).

\textbf{(b)} There exists a $\mathbf{P}$-triangular family $\left(
R_{\mathbf{p}}\right)  _{\mathbf{p}\in\mathbf{P}}\in\left(  \mathbb{F}\left(
x_{\mathbf{P}}\right)  \right)  ^{\mathbf{P}}$ such that every $\mathbf{q}%
\in\mathbf{P}$ satisfies $Q_{\mathbf{q}}\left(  \left(  R_{\mathbf{p}}\right)
_{\mathbf{p}\in\mathbf{P}}\right)  =x_{\mathbf{q}}$.
\end{lemma}

\begin{proof}
[Proof of Lemma \ref{lem.algebraic.triangularity.short} (sketched).]As in the
definition of $\mathbf{P}$-triangularity, we let $\mathbf{p}\Downarrow$ denote
the subset $\left\{  \mathbf{v}\in\mathbf{P}\ \mid\ \mathbf{v}\vartriangleleft
\mathbf{p}\right\}  $ of $\mathbf{P}$ for every $\mathbf{p}\in\mathbf{P}$.

\textbf{(a)} Assume that the family $\left(  Q_{\mathbf{p}}\right)
_{\mathbf{p}\in\mathbf{P}}\in\left(  \mathbb{F}\left(  x_{\mathbf{P}}\right)
\right)  ^{\mathbf{P}}$ is not algebraically independent (over $\mathbb{F}$).
Then, some nonzero polynomial $P\in\mathbb{F}\left[  x_{\mathbf{P}}\right]  $
satisfies $P\left(  \left(  Q_{\mathbf{p}}\right)  _{\mathbf{p}\in\mathbf{P}%
}\right)  =0$. Fix such a $P$, and let $\mathbf{u}$ be the maximal (with
respect to the order on $\mathbf{P}$) element of $\mathbf{P}$ such that
$x_{\mathbf{u}}$ appears in $P$ (meaning that the degree of $P$ with respect
to the variable $x_{\mathbf{u}}$ is $>0$). Then, $P$ can be construed as a
non-constant polynomial in the variable $x_{\mathbf{u}}$ over the ring
$\mathbb{F}\left[  x_{\mathbf{u}\Downarrow}\right]  $. Hence, $P\left(
\left(  Q_{\mathbf{p}}\right)  _{\mathbf{p}\in\mathbf{P}}\right)  =0$ shows
that $Q_{\mathbf{u}}$ is algebraic over the subfield of $\mathbb{F}\left(
x_{\mathbf{P}}\right)  $ generated by the elements $Q_{\mathbf{v}}$ for
$\mathbf{v}\in\left.  \mathbf{u}\Downarrow\right.  $.

Now, notice that every $\mathbf{v}\in\left.  \mathbf{u}\Downarrow\right.  $
satisfies $Q_{\mathbf{v}}\in\mathbb{F}\left(  x_{\mathbf{u}\Downarrow}\right)
$\ \ \ \ \footnote{\textit{Proof.} Let $\mathbf{v}\in\left.  \mathbf{u}%
\Downarrow\right.  $. Then, $\mathbf{v}\vartriangleleft\mathbf{u}$, so that
$\left.  \mathbf{v}\Downarrow\right.  \subseteq\left.  \mathbf{u}%
\Downarrow\right.  $, hence $\mathbb{F}\left(  x_{\mathbf{v}\Downarrow
}\right)  \subseteq\mathbb{F}\left(  x_{\mathbf{u}\Downarrow}\right)  $.
\par
By the algebraic triangularity condition, we know that there exist elements
$\alpha_{\mathbf{v}}$, $\beta_{\mathbf{v}}$, $\gamma_{\mathbf{v}}$,
$\delta_{\mathbf{v}}$ of $\mathbb{F}\left(  x_{\mathbf{v}\Downarrow}\right)  $
such that $\alpha_{\mathbf{v}}\delta_{\mathbf{v}}-\beta_{\mathbf{v}}%
\gamma_{\mathbf{v}}\neq0$ and $Q_{\mathbf{v}}=\dfrac{\alpha_{\mathbf{v}%
}x_{\mathbf{v}}+\beta_{\mathbf{v}}}{\gamma_{\mathbf{v}}x_{\mathbf{v}}%
+\delta_{\mathbf{v}}}$. These elements $\alpha_{\mathbf{v}}$, $\beta
_{\mathbf{v}}$, $\gamma_{\mathbf{v}}$, $\delta_{\mathbf{v}}$ belong to
$\mathbb{F}\left(  x_{\mathbf{u}\Downarrow}\right)  $ (by virtue of lying in
$\mathbb{F}\left(  x_{\mathbf{v}\Downarrow}\right)  \subseteq\mathbb{F}\left(
x_{\mathbf{u}\Downarrow}\right)  $), and so does $x_{\mathbf{v}}$ (since
$\mathbf{v}\in\left.  \mathbf{u}\Downarrow\right.  $). Hence, the fraction
$\dfrac{\alpha_{\mathbf{v}}x_{\mathbf{v}}+\beta_{\mathbf{v}}}{\gamma
_{\mathbf{v}}x_{\mathbf{v}}+\delta_{\mathbf{v}}}$ also lies in $\mathbb{F}%
\left(  x_{\mathbf{u}\Downarrow}\right)  $. Since this fraction is
$Q_{\mathbf{v}}$, we thus have shown $Q_{\mathbf{v}}\in\mathbb{F}\left(
x_{\mathbf{u}\Downarrow}\right)  $, qed.}. Hence, the subfield of
$\mathbb{F}\left(  x_{\mathbf{P}}\right)  $ generated by the elements
$Q_{\mathbf{v}}$ for $\mathbf{v}\in\left.  \mathbf{u}\Downarrow\right.  $ is a
subfield of $\mathbb{F}\left(  x_{\mathbf{u}\Downarrow}\right)  $. Since
$Q_{\mathbf{u}}$ is algebraic over the former field, we thus conclude that
$Q_{\mathbf{u}}$ is \textquotedblleft all the more\textquotedblright%
\ algebraic over the latter field. But by the algebraic triangularity
condition, there exist elements $\alpha_{\mathbf{u}}$, $\beta_{\mathbf{u}}$,
$\gamma_{\mathbf{u}}$, $\delta_{\mathbf{u}}$ of $\mathbb{F}\left(
x_{\mathbf{u}\Downarrow}\right)  $ such that $\alpha_{\mathbf{u}}%
\delta_{\mathbf{u}}-\beta_{\mathbf{u}}\gamma_{\mathbf{u}}\neq0$ and
$Q_{\mathbf{u}}=\dfrac{\alpha_{\mathbf{u}}x_{\mathbf{u}}+\beta_{\mathbf{u}}%
}{\gamma_{\mathbf{u}}x_{\mathbf{u}}+\delta_{\mathbf{u}}}$. We can easily solve
the equation $Q_{\mathbf{u}}=\dfrac{\alpha_{\mathbf{u}}x_{\mathbf{u}}%
+\beta_{\mathbf{u}}}{\gamma_{\mathbf{u}}x_{\mathbf{u}}+\delta_{\mathbf{u}}}$
for $x_{\mathbf{u}}$ and obtain $x_{\mathbf{u}}=\dfrac{Q_{\mathbf{u}}%
\delta_{\mathbf{u}}-\beta_{\mathbf{u}}}{\alpha_{\mathbf{u}}-Q_{\mathbf{u}%
}\gamma_{\mathbf{u}}}$ (and the denominator here does not vanish because of
$\alpha_{\mathbf{u}}\delta_{\mathbf{u}}-\beta_{\mathbf{u}}\gamma_{\mathbf{u}%
}\neq0$). Therefore, $x_{\mathbf{u}}$ is algebraic over the field
$\mathbb{F}\left(  x_{\mathbf{u}\Downarrow}\right)  $ (because we know
$Q_{\mathbf{u}}$ to be algebraic over this field, whereas $\alpha_{\mathbf{u}%
}$, $\beta_{\mathbf{u}}$, $\gamma_{\mathbf{u}}$, $\delta_{\mathbf{u}}$ lie in
that field). But this is absurd since $\mathbf{u}\notin\left.  \mathbf{u}%
\Downarrow\right.  $. This contradiction shows that our assumption was wrong,
and Lemma \ref{lem.algebraic.triangularity.short} \textbf{(a)} is proven.

\textbf{(b)} We will construct the required family $\left(  R_{\mathbf{p}%
}\right)  _{\mathbf{p}\in\mathbf{P}}\in\left(  \mathbb{F}\left(
x_{\mathbf{P}}\right)  \right)  ^{\mathbf{P}}$ by induction. Of course, this
is trivial if $\mathbf{P}=\varnothing$, so let us assume that $\mathbf{P}$ is
nonempty. Let $\mathbf{m}$ be the maximum element of $\mathbf{P}$, and let us
assume that we have already constructed a $\mathbf{P}\setminus\left\{
\mathbf{m}\right\}  $-triangular family $\left(  R_{\mathbf{p}}\right)
_{\mathbf{p}\in\mathbf{P}\setminus\left\{  \mathbf{m}\right\}  }\in\left(
\mathbb{F}\left(  x_{\mathbf{P}\setminus\left\{  \mathbf{m}\right\}  }\right)
\right)  ^{\mathbf{P}\setminus\left\{  \mathbf{m}\right\}  }$ such that every
$\mathbf{q}\in\mathbf{P}\setminus\left\{  \mathbf{m}\right\}  $ satisfies
$Q_{\mathbf{q}}\left(  \left(  R_{\mathbf{p}}\right)  _{\mathbf{p}%
\in\mathbf{P}\setminus\left\{  \mathbf{m}\right\}  }\right)  =x_{\mathbf{q}}$.
We now only need to find an element $R_{\mathbf{m}}\in\mathbb{F}\left(
x_{\mathbf{P}}\right)  $ such that the resulting family $\left(
R_{\mathbf{p}}\right)  _{\mathbf{p}\in\mathbf{P}}\in\left(  \mathbb{F}\left(
x_{\mathbf{P}}\right)  \right)  ^{\mathbf{P}}$ will be $\mathbf{P}$-triangular
and satisfy $Q_{\mathbf{m}}\left(  \left(  R_{\mathbf{p}}\right)
_{\mathbf{p}\in\mathbf{P}}\right)  =x_{\mathbf{m}}$.

Since $\mathbf{m}$ is maximum, we have $\left.  \mathbf{m}\Downarrow\right.
=\mathbf{P}\setminus\left\{  \mathbf{m}\right\}  $.

We know that the family $\left(  R_{\mathbf{p}}\right)  _{\mathbf{p}%
\in\mathbf{P}\setminus\left\{  \mathbf{m}\right\}  }$ is $\mathbf{P}%
\setminus\left\{  \mathbf{m}\right\}  $-triangular. Hence, Lemma
\ref{lem.algebraic.triangularity.short} \textbf{(a)} (applied to this family)
yields that the family $\left(  R_{\mathbf{p}}\right)  _{\mathbf{p}%
\in\mathbf{P}\setminus\left\{  \mathbf{m}\right\}  }$ is algebraically
independent. This yields that it can be substituted into any rational function
in $\mathbb{F}\left(  x_{\mathbf{P}\setminus\left\{  \mathbf{m}\right\}
}\right)  $ (without running the risk of denominators becoming $0$).

The family $\left(  Q_{\mathbf{p}}\right)  _{\mathbf{p}\in\mathbf{P}}$ is
$\mathbf{P}$-triangular, so that (by the algebraic triangularity condition)
there exist elements $\alpha_{\mathbf{m}}$, $\beta_{\mathbf{m}}$,
$\gamma_{\mathbf{m}}$, $\delta_{\mathbf{m}}$ of $\mathbb{F}\left(
x_{\mathbf{m}\Downarrow}\right)  $ such that $\alpha_{\mathbf{m}}%
\delta_{\mathbf{m}}-\beta_{\mathbf{m}}\gamma_{\mathbf{m}}\neq0$ and
$Q_{\mathbf{m}}=\dfrac{\alpha_{\mathbf{m}}x_{\mathbf{m}}+\beta_{\mathbf{m}}%
}{\gamma_{\mathbf{m}}x_{\mathbf{m}}+\delta_{\mathbf{m}}}$. Consider these
$\alpha_{\mathbf{m}}$, $\beta_{\mathbf{m}}$, $\gamma_{\mathbf{m}}$,
$\delta_{\mathbf{m}}$. Now, define four elements $\alpha_{\mathbf{m}}^{\prime
}$, $\beta_{\mathbf{m}}^{\prime}$, $\gamma_{\mathbf{m}}^{\prime}$,
$\delta_{\mathbf{m}}^{\prime}$ of $\mathbb{F}\left(  x_{\mathbf{P}%
\setminus\left\{  \mathbf{m}\right\}  }\right)  $ by%
\begin{align*}
\alpha_{\mathbf{m}}^{\prime}  &  =\delta_{\mathbf{m}}\left(  \left(
R_{\mathbf{p}}\right)  _{\mathbf{p}\in\mathbf{P}\setminus\left\{
\mathbf{m}\right\}  }\right)  ,\ \ \ \ \ \ \ \ \ \ \beta_{\mathbf{m}}^{\prime
}=-\beta_{\mathbf{m}}\left(  \left(  R_{\mathbf{p}}\right)  _{\mathbf{p}%
\in\mathbf{P}\setminus\left\{  \mathbf{m}\right\}  }\right)  ,\\
\gamma_{\mathbf{m}}^{\prime}  &  =-\gamma_{\mathbf{m}}\left(  \left(
R_{\mathbf{p}}\right)  _{\mathbf{p}\in\mathbf{P}\setminus\left\{
\mathbf{m}\right\}  }\right)  ,\ \ \ \ \ \ \ \ \ \ \delta_{\mathbf{m}}%
^{\prime}=\alpha_{\mathbf{m}}\left(  \left(  R_{\mathbf{p}}\right)
_{\mathbf{p}\in\mathbf{P}\setminus\left\{  \mathbf{m}\right\}  }\right)  .
\end{align*}
Note that these are well-defined (because $\alpha_{\mathbf{m}}$,
$\beta_{\mathbf{m}}$, $\gamma_{\mathbf{m}}$, $\delta_{\mathbf{m}}$ belong to
$\mathbb{F}\left(  x_{\mathbf{m}\Downarrow}\right)  =\mathbb{F}\left(
x_{\mathbf{P}\setminus\left\{  \mathbf{m}\right\}  }\right)  $ and because the
family $\left(  R_{\mathbf{p}}\right)  _{\mathbf{p}\in\mathbf{P}%
\setminus\left\{  \mathbf{m}\right\}  }$ is algebraically independent) and
belong to $\mathbb{F}\left(  x_{\mathbf{m}\Downarrow}\right)  $ (since
$\mathbf{P}\setminus\left\{  \mathbf{m}\right\}  =\left.  \mathbf{m}%
\Downarrow\right.  $). They furthermore satisfy%
\begin{align*}
&  \alpha_{\mathbf{m}}^{\prime}\delta_{\mathbf{m}}^{\prime}-\beta_{\mathbf{m}%
}^{\prime}\gamma_{\mathbf{m}}^{\prime}\\
&  =\delta_{\mathbf{m}}\left(  \left(  R_{\mathbf{p}}\right)  _{\mathbf{p}%
\in\mathbf{P}\setminus\left\{  \mathbf{m}\right\}  }\right)  \cdot
\alpha_{\mathbf{m}}\left(  \left(  R_{\mathbf{p}}\right)  _{\mathbf{p}%
\in\mathbf{P}\setminus\left\{  \mathbf{m}\right\}  }\right)  -\left(
-\beta_{\mathbf{m}}\left(  \left(  R_{\mathbf{p}}\right)  _{\mathbf{p}%
\in\mathbf{P}\setminus\left\{  \mathbf{m}\right\}  }\right)  \right)  \left(
-\gamma_{\mathbf{m}}\left(  \left(  R_{\mathbf{p}}\right)  _{\mathbf{p}%
\in\mathbf{P}\setminus\left\{  \mathbf{m}\right\}  }\right)  \right) \\
&  =\underbrace{\left(  \delta_{\mathbf{m}}\alpha_{\mathbf{m}}-\left(
-\beta_{\mathbf{m}}\right)  \left(  -\gamma_{\mathbf{m}}\right)  \right)
}_{=\alpha_{\mathbf{m}}\delta_{\mathbf{m}}-\beta_{\mathbf{m}}\gamma
_{\mathbf{m}}\neq0}\left(  \left(  R_{\mathbf{p}}\right)  _{\mathbf{p}%
\in\mathbf{P}\setminus\left\{  \mathbf{m}\right\}  }\right)  \neq0
\end{align*}
(since $\left(  R_{\mathbf{p}}\right)  _{\mathbf{p}\in\mathbf{P}%
\setminus\left\{  \mathbf{m}\right\}  }$ is algebraically independent). Let us
now define $R_{\mathbf{m}}=\dfrac{\alpha_{\mathbf{m}}^{\prime}x_{\mathbf{m}%
}+\beta_{\mathbf{m}}^{\prime}}{\gamma_{\mathbf{m}}^{\prime}x_{\mathbf{m}%
}+\delta_{\mathbf{m}}^{\prime}}$. (This is easily seen to be well-defined
because $\alpha_{\mathbf{m}}^{\prime}\delta_{\mathbf{m}}^{\prime}%
-\beta_{\mathbf{m}}^{\prime}\gamma_{\mathbf{m}}^{\prime}\neq0$ entails
$\left(  \gamma_{\mathbf{m}}^{\prime},\delta_{\mathbf{m}}^{\prime}\right)
\neq\left(  0,0\right)  $.) Since the family $\left(  R_{\mathbf{p}}\right)
_{\mathbf{p}\in\mathbf{P}\setminus\left\{  \mathbf{m}\right\}  }$ is already
$\mathbf{P}\setminus\left\{  \mathbf{m}\right\}  $-triangular, and because of
the fact that $\alpha_{\mathbf{m}}^{\prime}$, $\beta_{\mathbf{m}}^{\prime}$,
$\gamma_{\mathbf{m}}^{\prime}$, $\delta_{\mathbf{m}}^{\prime}$ are elements of
$\mathbb{F}\left(  x_{\mathbf{m}\Downarrow}\right)  $ satisfying
$\alpha_{\mathbf{m}}^{\prime}\delta_{\mathbf{m}}^{\prime}-\beta_{\mathbf{m}%
}^{\prime}\gamma_{\mathbf{m}}^{\prime}\neq0$ and $R_{\mathbf{m}}=\dfrac
{\alpha_{\mathbf{m}}^{\prime}x_{\mathbf{m}}+\beta_{\mathbf{m}}^{\prime}%
}{\gamma_{\mathbf{m}}^{\prime}x_{\mathbf{m}}+\delta_{\mathbf{m}}^{\prime}}$,
we see that the family $\left(  R_{\mathbf{p}}\right)  _{\mathbf{p}%
\in\mathbf{P}}\in\left(  \mathbb{F}\left(  x_{\mathbf{P}}\right)  \right)
^{\mathbf{P}}$ is $\mathbf{P}$-triangular. We are now going to prove that
$Q_{\mathbf{m}}\left(  \left(  R_{\mathbf{p}}\right)  _{\mathbf{p}%
\in\mathbf{P}}\right)  =x_{\mathbf{m}}$, and then we will be done.

Since $Q_{\mathbf{m}}=\dfrac{\alpha_{\mathbf{m}}x_{\mathbf{m}}+\beta
_{\mathbf{m}}}{\gamma_{\mathbf{m}}x_{\mathbf{m}}+\delta_{\mathbf{m}}}$, we
have
\begin{equation}
Q_{\mathbf{m}}\left(  \left(  R_{\mathbf{p}}\right)  _{\mathbf{p}\in
\mathbf{P}}\right)  =\dfrac{\alpha_{\mathbf{m}}\left(  \left(  R_{\mathbf{p}%
}\right)  _{\mathbf{p}\in\mathbf{P}}\right)  R_{\mathbf{m}}+\beta_{\mathbf{m}%
}\left(  \left(  R_{\mathbf{p}}\right)  _{\mathbf{p}\in\mathbf{P}}\right)
}{\gamma_{\mathbf{m}}\left(  \left(  R_{\mathbf{p}}\right)  _{\mathbf{p}%
\in\mathbf{P}}\right)  R_{\mathbf{m}}+\delta_{\mathbf{m}}\left(  \left(
R_{\mathbf{p}}\right)  _{\mathbf{p}\in\mathbf{P}}\right)  }.
\label{pf.algebraic.triangularity.short.3}%
\end{equation}
But $\alpha_{\mathbf{m}}\in\mathbb{F}\left(  x_{\mathbf{P}\setminus\left\{
\mathbf{m}\right\}  }\right)  $, so that the variable $x_{\mathbf{m}}$ does
not appear in $\alpha_{\mathbf{m}}$ at all. Hence, $\alpha_{\mathbf{m}}\left(
\left(  R_{\mathbf{p}}\right)  _{\mathbf{p}\in\mathbf{P}}\right)
=\alpha_{\mathbf{m}}\left(  \left(  R_{\mathbf{p}}\right)  _{\mathbf{p}%
\in\mathbf{P}\setminus\left\{  \mathbf{m}\right\}  }\right)  =\delta
_{\mathbf{m}}^{\prime}$. Using this and the similarly proven equalities
$\beta_{\mathbf{m}}\left(  \left(  R_{\mathbf{p}}\right)  _{\mathbf{p}%
\in\mathbf{P}}\right)  =-\beta_{\mathbf{m}}^{\prime}$, $\gamma_{\mathbf{m}%
}\left(  \left(  R_{\mathbf{p}}\right)  _{\mathbf{p}\in\mathbf{P}}\right)
=-\gamma_{\mathbf{m}}^{\prime}$ and $\delta_{\mathbf{m}}\left(  \left(
R_{\mathbf{p}}\right)  _{\mathbf{p}\in\mathbf{P}}\right)  =\alpha_{\mathbf{m}%
}^{\prime}$, we can rewrite the equality
(\ref{pf.algebraic.triangularity.short.3}) as%
\[
Q_{\mathbf{m}}\left(  \left(  R_{\mathbf{p}}\right)  _{\mathbf{p}\in
\mathbf{P}}\right)  =\dfrac{\delta_{\mathbf{m}}^{\prime}R_{\mathbf{m}}%
-\beta_{\mathbf{m}}^{\prime}}{-\gamma_{\mathbf{m}}^{\prime}R_{\mathbf{m}%
}+\alpha_{\mathbf{m}}^{\prime}}.
\]
But the right hand side of this equality simplifies to $x_{\mathbf{m}}$ if we
recall that $R_{\mathbf{m}}=\dfrac{\alpha_{\mathbf{m}}^{\prime}x_{\mathbf{m}%
}+\beta_{\mathbf{m}}^{\prime}}{\gamma_{\mathbf{m}}^{\prime}x_{\mathbf{m}%
}+\delta_{\mathbf{m}}^{\prime}}$ (the proof of this is mechanical, using no
properties of $\alpha_{\mathbf{m}}^{\prime}$, $\beta_{\mathbf{m}}^{\prime}$,
$\gamma_{\mathbf{m}}^{\prime}$, $\delta_{\mathbf{m}}^{\prime}$ and
$x_{\mathbf{m}}$ other than lying in a field). Hence, we have shown that
$Q_{\mathbf{m}}\left(  \left(  R_{\mathbf{p}}\right)  _{\mathbf{p}%
\in\mathbf{P}}\right)  =x_{\mathbf{m}}$. As explained above, this completes
the (inductive) proof of Lemma \ref{lem.algebraic.triangularity.short}
\textbf{(b)}.
\end{proof}

We now can proceed to the proof of Proposition \ref{prop.Grasp.generic}:

\begin{proof}
[Proof of Proposition \ref{prop.Grasp.generic} (sketched).]Let $\mathbb{F}$ be
the prime field of $\mathbb{K}$. (This means either $\mathbb{Q}$ or
$\mathbb{F}_{p}$ depending on the characteristic of $\mathbb{K}$.) In the
following, the word ``algebraically independent'' will always mean
``algebraically independent over $\mathbb{F}$'' (rather than over $\mathbb{K}$
or over $\mathbb{Z}$).

Let $\mathbf{P}$ be a totally ordered set such that
\[
\mathbf{P}=\left\{  1,2,...,p\right\}  \times\left\{  1,2,...,q\right\}
\text{ as sets,}%
\]
and such that%
\[
\left(  i,k\right)  \trianglelefteq\left(  i^{\prime},k^{\prime}\right)
\text{ for all }\left(  i,k\right)  \in\mathbf{P}\text{ and }\left(
i^{\prime},k^{\prime}\right)  \in\mathbf{P}\text{ satisfying }\left(  i\geq
i^{\prime}\text{ and }k\leq k^{\prime}\right)  ,
\]
where $\trianglelefteq$ denotes the smaller-or-equal relation of $\mathbf{P}$.
Such a $\mathbf{P}$ clearly exists (in fact, there usually exist several such
$\mathbf{P}$, and it doesn't matter which of them we choose). We denote the
smaller relation of $\mathbf{P}$ by $\vartriangleleft$. We will later see what
this total order is good for (intuitively, it is an order in which the
variables can be eliminated; in other words, it makes our system behave like a
triangular matrix rather than like a triangular matrix with permuted columns),
but for now let us notice that it is generally not compatible with
$\operatorname*{Rect}\left(  p,q\right)  $.

Let $Z:\left\{  1,2,...,q\right\}  \rightarrow\left\{  1,2,...,q\right\}  $
denote the map which sends every $k\in\left\{  1,2,...,q-1\right\}  $ to $k+1$
and sends $q$ to $1$. Thus, $Z$ is a permutation in the symmetric group
$S_{q}$, and can be written in cycle notation as $\left(  1,2,...,q\right)  $.

Consider the field $\mathbb{F}\left(  x_{\mathbf{P}}\right)  $ and the ring
$\mathbb{F}\left[  x_{\mathbf{P}}\right]  $ defined as in Definition
\ref{def.algebraic.triangularity.short}.

Recall that we need to prove Proposition \ref{prop.Grasp.generic}. In other
words, we need to show that for almost every $f\in\mathbb{K}%
^{\operatorname*{Rect}\left(  p,q\right)  }$, there exists a matrix
$A\in\mathbb{K}^{p\times\left(  p+q\right)  }$ satisfying
$f=\operatorname*{Grasp}\nolimits_{0}A$.

In order to prove this, it is enough to show that there exists a matrix
$\widetilde{D}\in\left(  \mathbb{F}\left(  x_{\mathbf{P}}\right)  \right)
^{p\times\left(  p+q\right)  }$ satisfying%
\begin{equation}
x_{\mathbf{p}}=\left(  \operatorname*{Grasp}\nolimits_{0}\widetilde{D}\right)
\left(  \mathbf{p}\right)  \ \ \ \ \ \ \ \ \ \ \text{for every }\mathbf{p}%
\in\mathbf{P}\text{.} \label{pf.Grasp.generic.short.reduce-to-rational}%
\end{equation}
Indeed, once the existence of such a matrix $\widetilde{D}$ is proven, we will
be able to obtain a matrix $A\in\mathbb{K}^{p\times\left(  p+q\right)  }$
satisfying $f=\operatorname*{Grasp}\nolimits_{0}A$ for almost every
$f\in\mathbb{K}^{\operatorname*{Rect}\left(  p,q\right)  }$ simply by
substituting $f\left(  \mathbf{p}\right)  $ for every $x_{\mathbf{p}}$ in all
entries of the matrix $\widetilde{D}$\ \ \ \ \footnote{Indeed, this matrix $A$
(obtained by substitution of $f\left(  \mathbf{p}\right)  $ for $x_{\mathbf{p}%
}$) will be well-defined for almost every $f\in\mathbb{K}%
^{\operatorname*{Rect}\left(  p,q\right)  }$ (the ``almost'' is due to the
possibility of some denominators becoming $0$), and will satisfy $f\left(
\mathbf{p}\right)  =\left(  \operatorname*{Grasp}\nolimits_{0}A\right)
\left(  \mathbf{p}\right)  $ for every $\mathbf{p}\in\mathbf{P}$ (because
$\widetilde{D}$ satisfies (\ref{pf.Grasp.generic.short.reduce-to-rational})),
that is, $f=\operatorname*{Grasp}\nolimits_{0}A$.}. Hence, all we need to show
is the existence of a matrix $\widetilde{D}\in\left(  \mathbb{F}\left(
x_{\mathbf{P}}\right)  \right)  ^{p\times\left(  p+q\right)  }$ satisfying
(\ref{pf.Grasp.generic.short.reduce-to-rational}).

Define a matrix $C\in\left(  \mathbb{F}\left[  x_{\mathbf{P}}\right]  \right)
^{p\times q}$ by
\[
C=\left(  x_{\left(  i,Z\left(  k\right)  \right)  }\right)  _{1\leq i\leq
p,\ 1\leq k\leq q}.
\]
This is simply a matrix whose entries are all the indeterminates
$x_{\mathbf{p}}$ of the polynomial ring $\mathbb{F}\left[  x_{\mathbf{P}%
}\right]  $, albeit in a strange order. (The order, again, is tailored to make
the ``triangularity'' argument work nicely. This matrix $C$ is not going to be
directly related to the $\widetilde{D}$ we will construct, but will be used in
its construction.)

For every $\left(  i,k\right)  \in\mathbf{P}$, define an element
$\mathfrak{N}_{\left(  i,k\right)  }\in\mathbb{F}\left[  x_{\mathbf{P}%
}\right]  $ by%
\begin{equation}
\mathfrak{N}_{\left(  i,k\right)  }=\det\left(  \left(  I_{p}\mid C\right)
\left[  1:i\mid i+k-1:p+k\right]  \right)  .
\label{lem.Grasp.generic.short.Ndef}%
\end{equation}

For every $\left(  i,k\right)  \in\mathbf{P}$, define an element
$\mathfrak{D}_{\left(  i,k\right)  }\in\mathbb{F}\left[  x_{\mathbf{P}%
}\right]  $ by%
\begin{equation}
\mathfrak{D}_{\left(  i,k\right)  }=\det\left(  \left(  I_{p}\mid C\right)
\left[  0:i\mid i+k:p+k\right]  \right)  .
\label{lem.Grasp.generic.short.Ddef}%
\end{equation}

Our plan from here is the following:

\textit{Step 1:} We will find alternate expressions for the polynomials
$\mathfrak{N}_{\left(  i,k\right)  }$ and $\mathfrak{D}_{\left(  i,k\right)
}$ which will give us a better idea of what variables occur in these polynomials.

\textit{Step 2:} We will show that $\mathfrak{N}_{\left(  i,k\right)  }$ and
$\mathfrak{D}_{\left(  i,k\right)  }$ are nonzero for all $\left(  i,k\right)
\in\mathbf{P}$.

\textit{Step 3:} We will define a $Q_{\mathbf{p}}\in\mathbb{F}\left(
x_{\mathbf{P}}\right)  $ for every $\mathbf{p}\in\mathbf{P}$ by $Q_{\mathbf{p}%
}=\dfrac{\mathfrak{N}_{\mathbf{p}}}{\mathfrak{D}_{\mathbf{p}}}$, and we will
show that $Q_{\mathbf{p}}=\left(  \operatorname*{Grasp}\nolimits_{0}\left(
I_{p}\mid C\right)  \right)  \left(  \mathbf{p}\right)  $.

\textit{Step 4:} We will prove that the family $\left(  Q_{\mathbf{p}}\right)
_{\mathbf{p}\in\mathbf{P}}\in\left(  \mathbb{F}\left(  x_{\mathbf{P}}\right)
\right)  ^{\mathbf{P}}$ is $\mathbf{P}$-triangular.

\textit{Step 5:} We will use Lemma \ref{lem.algebraic.triangularity.short}
\textbf{(b)} and the result of Step 4 to find a matrix $\widetilde{D}%
\in\left(  \mathbb{F}\left(  x_{\mathbf{P}}\right)  \right)  ^{p\times\left(
p+q\right)  }$ satisfying (\ref{pf.Grasp.generic.short.reduce-to-rational}).

Let us now go into detail on each specific step (although we won't take that
detail very far).

\textit{Details of Step 1:} Let us introduce three more pieces of notation
pertaining to matrices:

\begin{itemize}
\item If $\ell\in\mathbb{N}$, and if $A_{1}$, $A_{2}$, $...$, $A_{k}$ are
several matrices with $\ell$ rows each, then $\left(  A_{1}\mid A_{2}%
\mid...\mid A_{k}\right)  $ will denote the matrix obtained by starting with
an (empty) $\ell\times0$-matrix, then attaching the matrix $A_{1}$ to it on
the right, then attaching the matrix $A_{2}$ to the result on the right, etc.,
and finally attaching the matrix $A_{k}$ to the result on the right. For
example, if $p$ is a nonnegative integer, and $B$ is a matrix with $p$ rows,
then $\left(  I_{p}\mid B\right)  $ means the matrix obtained from the
$p\times p$ identity matrix $I_{p}$ by attaching the matrix $B$ to it on the
right. (As a concrete example, $\left(  I_{2}\mid\left(
\begin{array}
[c]{cc}%
1 & -2\\
3 & 0
\end{array}
\right)  \right)  =\left(
\begin{array}
[c]{cccc}%
1 & 0 & 1 & -2\\
0 & 1 & 3 & 0
\end{array}
\right)  $.)

\item If $\ell\in\mathbb{N}$, if $B$ is a matrix with $\ell$ rows, and if
$i_{1}$, $i_{2}$, $...$, $i_{k}$ are some elements of $\left\{  1,2,...,\ell
\right\}  $, then $\operatorname*{rows}\nolimits_{i_{1},i_{2},...,i_{k}}B$
will denote the matrix whose rows (from top to bottom) are the rows labelled
$i_{1}$, $i_{2}$, $...$, $i_{k}$ of the matrix $B$.

\item If $u$ and $v$ are two nonnegative integers, and $A$ is a $u\times
v$-matrix, then, for any two integers $a$ and $b$ satisfying $a\leq b$, we let
$A\left[  a:b\right]  $ be the matrix whose columns (from left to right) are
$A_{a}$, $A_{a+1}$, $...$, $A_{b-1}$. This is a natural extension of the
notation introduced in Definition \ref{def.minors} \textbf{(c)} (or, rather,
the latter notation is a natural extension of the definition we just made) and
has the obvious property that if $a$, $b$ and $c$ are integers satisfying
$a\leq b\leq c$, then $A\left[  a:c\right]  =A\left[  a:b\mid b:c\right]  $.
\end{itemize}

We will use without proof a standard fact about determinants:

\begin{itemize}
\item Given a commutative ring $\mathbb{L}$, two nonnegative integers $a$ and
$b$ satisfying $a\geq b$, and a matrix $U\in\mathbb{L}^{a\times b}$, we have%
\begin{equation}
\det\left(  \left(
\begin{array}
[c]{c}%
I_{a-b}\\
0_{b\times\left(  a-b\right)  }%
\end{array}
\right)  \mid U\right)  =\det\left(  \operatorname*{rows}%
\nolimits_{a-b+1,a-b+2,...,a}U\right)
\label{pf.Grasp.generic.short.step1.block1}%
\end{equation}
and%
\begin{equation}
\det\left(  \left(
\begin{array}
[c]{c}%
0_{b\times\left(  a-b\right)  }\\
I_{a-b}%
\end{array}
\right)  \mid U\right)  =\left(  -1\right)  ^{b\left(  a-b\right)  }%
\det\left(  \operatorname*{rows}\nolimits_{1,2,...,b}U\right)  .
\label{pf.Grasp.generic.short.step1.block2}%
\end{equation}
(Here, $0_{u\times v}$ denotes the $u\times v$ zero matrix for all
$u\in\mathbb{N}$ and $v\in\mathbb{N}$, and $\left(
\begin{array}
[c]{c}%
I_{a-b}\\
0_{b\times\left(  a-b\right)  }%
\end{array}
\right)  $ and $\left(
\begin{array}
[c]{c}%
0_{b\times\left(  a-b\right)  }\\
I_{a-b}%
\end{array}
\right)  $ are to be read as block matrices.)
\end{itemize}

Now,
\[
\left(  I_{p}\mid C\right)  \left[  1:i\mid i+k-1:p+k\right]  =\left(  \left(
\begin{array}
[c]{c}%
I_{i-1}\\
0_{\left(  p-\left(  i-1\right)  \right)  \times\left(  i-1\right)  }%
\end{array}
\right)  \ \mid\ \left(  I_{p}\mid C\right)  \left[  i+k-1:p+k\right]
\right)  ,
\]
so that%
\begin{align*}
&  \det\left(  \left(  I_{p}\mid C\right)  \left[  1:i\mid i+k-1:p+k\right]
\right) \\
&  =\det\left(  \left(
\begin{array}
[c]{c}%
I_{i-1}\\
0_{\left(  p-\left(  i-1\right)  \right)  \times\left(  i-1\right)  }%
\end{array}
\right)  \ \mid\ \left(  I_{p}\mid C\right)  \left[  i+k-1:p+k\right]  \right)
\\
&  =\det\left(  \operatorname*{rows}\nolimits_{i,i+1,...,p}\left(  \left(
I_{p}\mid C\right)  \left[  i+k-1:p+k\right]  \right)  \right)
\end{align*}
(by (\ref{pf.Grasp.generic.short.step1.block1})). Thus,%
\begin{align}
\mathfrak{N}_{\left(  i,k\right)  }  &  =\det\left(  \left(  I_{p}\mid
C\right)  \left[  1:i\mid i+k-1:p+k\right]  \right) \nonumber\\
&  =\det\left(  \operatorname*{rows}\nolimits_{i,i+1,...,p}\left(  \left(
I_{p}\mid C\right)  \left[  i+k-1:p+k\right]  \right)  \right)  .
\label{pf.Grasp.generic.short.step1.N}%
\end{align}

Also,%
\begin{align*}
&  \left(  I_{p}\mid C\right)  \left[  0:i\mid i+k:p+k\right] \\
&  =\left(  \underbrace{\left(  I_{p}\mid C\right)  _{0}}_{\substack{=\left(
-1\right)  ^{p-1}C_{q}\\\text{(by Definition \ref{def.minors} \textbf{(b)})}%
}}\ \mid\ \left(
\begin{array}
[c]{c}%
I_{i-1}\\
0_{\left(  p-\left(  i-1\right)  \right)  \times\left(  i-1\right)  }%
\end{array}
\right)  \ \mid\ \left(  I_{p}\mid C\right)  \left[  i+k:p+k\right]  \right)
\\
&  =\left(  \left(  -1\right)  ^{p-1}C_{q}\ \mid\ \left(
\begin{array}
[c]{c}%
I_{i-1}\\
0_{\left(  p-\left(  i-1\right)  \right)  \times\left(  i-1\right)  }%
\end{array}
\right)  \ \mid\ \left(  I_{p}\mid C\right)  \left[  i+k:p+k\right]  \right)
,
\end{align*}
whence%
\begin{align*}
&  \det\left(  \left(  I_{p}\mid C\right)  \left[  0:i\mid i+k:p+k\right]
\right) \\
&  =\det\left(  \left(  -1\right)  ^{p-1}C_{q}\ \mid\ \left(
\begin{array}
[c]{c}%
I_{i-1}\\
0_{\left(  p-\left(  i-1\right)  \right)  \times\left(  i-1\right)  }%
\end{array}
\right)  \ \mid\ \left(  I_{p}\mid C\right)  \left[  i+k:p+k\right]  \right)
\\
&  =\left(  -1\right)  ^{p-1}\det\left(  C_{q}\ \mid\ \left(
\begin{array}
[c]{c}%
I_{i-1}\\
0_{\left(  p-\left(  i-1\right)  \right)  \times\left(  i-1\right)  }%
\end{array}
\right)  \ \mid\ \left(  I_{p}\mid C\right)  \left[  i+k:p+k\right]  \right)
\\
&  =\underbrace{\left(  -1\right)  ^{p-1}\left(  -1\right)  ^{i-1}}_{=\left(
-1\right)  ^{p-i}}\det\left(  \left(
\begin{array}
[c]{c}%
I_{i-1}\\
0_{\left(  p-\left(  i-1\right)  \right)  \times\left(  i-1\right)  }%
\end{array}
\right)  \ \mid\ C_{q}\ \mid\ \left(  I_{p}\mid C\right)  \left[
i+k:p+k\right]  \right) \\
&  \ \ \ \ \ \ \ \ \ \ \left(
\begin{array}
[c]{c}%
\text{since permuting the columns of a matrix multiplies the}\\
\text{determinant by the sign of the permutation}%
\end{array}
\right) \\
&  =\left(  -1\right)  ^{p-i}\det\left(  \left(
\begin{array}
[c]{c}%
I_{i-1}\\
0_{\left(  p-\left(  i-1\right)  \right)  \times\left(  i-1\right)  }%
\end{array}
\right)  \ \mid\ C_{q}\ \mid\ \left(  I_{p}\mid C\right)  \left[
i+k:p+k\right]  \right) \\
&  =\left(  -1\right)  ^{p-i}\det\left(  \operatorname*{rows}%
\nolimits_{i,i+1,...,p}\left(  C_{q}\ \mid\ \left(  I_{p}\mid C\right)
\left[  i+k:p+k\right]  \right)  \right)
\end{align*}
(by (\ref{pf.Grasp.generic.short.step1.block1})). Thus,%
\begin{align}
\mathfrak{D}_{\left(  i,k\right)  }  &  =\det\left(  \left(  I_{p}\mid
C\right)  \left[  0:i\mid i+k:p+k\right]  \right) \nonumber\\
&  =\left(  -1\right)  ^{p-i}\det\left(  \operatorname*{rows}%
\nolimits_{i,i+1,...,p}\left(  C_{q}\ \mid\ \left(  I_{p}\mid C\right)
\left[  i+k:p+k\right]  \right)  \right)  .
\label{pf.Grasp.generic.short.step1.D}%
\end{align}

We have thus found alternative formulas (\ref{pf.Grasp.generic.short.step1.N})
and (\ref{pf.Grasp.generic.short.step1.D}) for $\mathfrak{N}_{\left(
i,k\right)  }$ and $\mathfrak{D}_{\left(  i,k\right)  }$. While not shorter
than the definitions, these formulas involve smaller matrices (unless $i=1$)
and are more useful in understanding the monomials appearing in $\mathfrak{N}%
_{\left(  i,k\right)  }$ and $\mathfrak{D}_{\left(  i,k\right)  }$.

\textit{Details of Step 2:} We claim that $\mathfrak{N}_{\left(  i,k\right)
}$ and $\mathfrak{D}_{\left(  i,k\right)  }$ are nonzero for all $\left(
i,k\right)  \in\mathbf{P}$.

\textit{Proof.} Let $\left(  i,k\right)  \in\mathbf{P}$. Let us first check
that $\mathfrak{N}_{\left(  i,k\right)  }$ is nonzero.

There are, in fact, many ways to do this. Here is probably the shortest one:
Assume the contrary, i.e., assume that $\mathfrak{N}_{\left(  i,k\right)  }%
=0$. Then, every matrix $G\in\mathbb{F}^{p\times\left(  p+q\right)  }$
satisfies $\det\left(  G\left[  1:i\mid i+k-1:p+k\right]  \right)
=0$\ \ \ \ \footnote{\textit{Proof.} Let $\widetilde{\mathbb{F}}$ be a field
extension of $\mathbb{F}$ such that $\left\vert \widetilde{\mathbb{F}%
}\right\vert =\infty$. (We need this to make sense of Zariski density
arguments.) We are going to prove that every matrix $G\in\widetilde{\mathbb{F}%
}^{p\times\left(  p+q\right)  }$ satisfies $\det\left(  G\left[  1:i\mid
i+k-1:p+k\right]  \right)  =0$; this will clearly imply the same claim for
$G\in\mathbb{F}^{p\times\left(  p+q\right)  }$.
\par
Let $G\in\widetilde{\mathbb{F}}^{p\times\left(  p+q\right)  }$. We want to
prove that $\det\left(  G\left[  1:i\mid i+k-1:p+k\right]  \right)  =0$. Since
this is a polynomial identity in the entries of $G$, we can WLOG assume that
$G$ is generic enough that the first $p$ columns of $G$ are linearly
independent (since this just restricts $G$ to a Zariski-dense open subset of
$\widetilde{\mathbb{F}}^{p\times\left(  p+q\right)  }$). Assume this. Then, we
can write $G$ in the form $\left(  U\mid V\right)  $, with $U$ being the
matrix formed by the first $p$ columns of $G$, and $V$ being the matrix formed
by the remaining $q$ columns. Since the first $p$ columns of $G$ are linearly
independent, the matrix $U$ is invertible.
\par
Left multiplication by $U^{-1}$ acts on matrices column by column. This yields%
\[
U^{-1}\cdot\left(  G\left[  1:i\mid i+k-1:p+k\right]  \right)  =\left(
U^{-1}G\right)  \left[  1:i\mid i+k-1:p+k\right]  .
\]
Also, $U^{-1} \underbrace{G}_{=\left(  U\mid V\right)  } = U^{-1}\left(  U\mid
V\right)  = \left(  U^{-1}U \mid U^{-1}V\right)  = \left(  I_{p} \mid
U^{-1}V\right)  $.
\par
Now, we have $\mathfrak{N}_{\left(  i,k\right)  }=0$. Since $\mathfrak{N}%
_{\left(  i,k\right)  }=\det\left(  \left(  I_{p}\mid C\right)  \left[
1:i\mid i+k-1:p+k\right]  \right)  $, this yields that $\det\left(  \left(
I_{p}\mid C\right)  \left[  1:i\mid i+k-1:p+k\right]  \right)  =0$. But the
matrix $C$ is, in some sense, the ``most generic matrix'': namely, the entries
of the matrix $C$ are pairwise distinct commuting indeterminates, and
therefore we can obtain any other matrix (over a commutative $\mathbb{F}%
$-algebra) from $C$ by substituting the corresponding values for the
indeterminates. In particular, we can make a substitution that turns $C$ into
$U^{-1}V$. Thus, from $\det\left(  \left(  I_{p}\mid C\right)  \left[  1:i\mid
i+k-1:p+k\right]  \right)  =0$, we obtain $\det\left(  \left(  I_{p}\mid
U^{-1}V\right)  \left[  1:i\mid i+k-1:p+k\right]  \right)  =0$.
\par
Now,%
\begin{align*}
&  \left(  \det U\right)  ^{-1}\cdot\det\left(  G\left[  1:i\mid
i+k-1:p+k\right]  \right) \\
&  =\det\left(  \underbrace{U^{-1}\cdot\left(  G\left[  1:i\mid
i+k-1:p+k\right]  \right)  }_{=\left(  U^{-1}G\right)  \left[  1:i\mid
i+k-1:p+k\right]  }\right)  =\det\left(  \left(  \underbrace{U^{-1}%
G}_{=\left(  I_{p}\mid U^{-1}V\right)  }\right)  \left[  1:i\mid
i+k-1:p+k\right]  \right) \\
&  =\det\left(  \left(  I_{p}\mid U^{-1}V\right)  \left[  1:i\mid
i+k-1:p+k\right]  \right)  =0.
\end{align*}
Multiplying this with $\det U$ (which is nonzero since $U$ is invertible), we
obtain $\det\left(  G\left[  1:i\mid i+k-1:p+k\right]  \right)  =0$, qed.}.
But this is absurd, because we can pick $G$ to have the $p$ columns labelled
$1$, $2$, $...$, $i-1$, $i+k-1$, $i+k$, $...$, $p+k-1$ linearly independent.
This contradiction shows that our assumption was wrong. Hence, $\mathfrak{N}%
_{\left(  i,k\right)  }$ is nonzero. Similarly, $\mathfrak{D}_{\left(
i,k\right)  }$ is nonzero.

\textit{Details of Step 3:} Define a $Q_{\mathbf{p}}\in\mathbb{F}\left(
x_{\mathbf{P}}\right)  $ for every $\mathbf{p}\in\mathbf{P}$ by $Q_{\mathbf{p}%
}=\dfrac{\mathfrak{N}_{\mathbf{p}}}{\mathfrak{D}_{\mathbf{p}}}$. This is
well-defined because Step 2 has shown that $\mathfrak{D}_{\mathbf{p}}$ is
nonzero. Moreover, it is easy to see that every $\left(  i,k\right)
\in\mathbf{P}$ satisfies%
\[
Q_{\left(  i,k\right)  }=\left(  \operatorname*{Grasp}\nolimits_{0}\left(
I_{p}\mid C\right)  \right)  \left(  \left(  i,k\right)  \right)  .
\]
\footnote{Indeed, the definition of $\operatorname*{Grasp}\nolimits_{0}\left(
I_{p}\mid C\right)  $ yields%
\[
\left(  \operatorname*{Grasp}\nolimits_{0}\left(  I_{p}\mid C\right)  \right)
\left(  \left(  i,k\right)  \right)  =\dfrac{\det\left(  \left(  I_{p}\mid
C\right)  \left[  1:i\mid i+k-1:p+k\right]  \right)  }{\det\left(  \left(
I_{p}\mid C\right)  \left[  0:i\mid i+k:p+k\right]  \right)  }=\dfrac
{\mathfrak{N}_{\left(  i,k\right)  }}{\mathfrak{D}_{\left(  i,k\right)  }}%
\]
(by (\ref{lem.Grasp.generic.short.Ndef}) and
(\ref{lem.Grasp.generic.short.Ddef})).} In other words, every $\mathbf{p}%
\in\mathbf{P}$ satisfies%
\begin{equation}
Q_{\mathbf{p}}=\left(  \operatorname*{Grasp}\nolimits_{0}\left(  I_{p}\mid
C\right)  \right)  \left(  \mathbf{p}\right)  .
\label{pf.Grasp.generic.short.step3}%
\end{equation}

\textit{Details of Step 4:} We are now going to prove that the family $\left(
Q_{\mathbf{p}}\right)  _{\mathbf{p}\in\mathbf{P}}\in\left(  \mathbb{F}\left(
x_{\mathbf{P}}\right)  \right)  ^{\mathbf{P}}$ is $\mathbf{P}$-triangular.

By the definition of $\mathbf{P}$-triangularity, this requires showing that
for every $\mathbf{p}\in\mathbf{P}$, there exist elements $\alpha_{\mathbf{p}%
}$, $\beta_{\mathbf{p}}$, $\gamma_{\mathbf{p}}$, $\delta_{\mathbf{p}}$ of
$\mathbb{F}\left(  x_{\mathbf{p}\Downarrow}\right)  $ such that $\alpha
_{\mathbf{p}}\delta_{\mathbf{p}}-\beta_{\mathbf{p}}\gamma_{\mathbf{p}}\neq0$
and $Q_{\mathbf{p}}=\dfrac{\alpha_{\mathbf{p}}x_{\mathbf{p}}+\beta
_{\mathbf{p}}}{\gamma_{\mathbf{p}}x_{\mathbf{p}}+\delta_{\mathbf{p}}}$ (where
$\mathbf{p}\Downarrow$ is defined as in Definition
\ref{def.algebraic.triangularity.short} \textbf{(d)}). So fix $\mathbf{p}%
\in\mathbf{P}$. Write $\mathbf{p}$ in the form $\mathbf{p}=\left(  i,k\right)
$.

We will actually do something slightly better than we need. We will find
elements $\alpha_{\mathbf{p}}$, $\beta_{\mathbf{p}}$, $\gamma_{\mathbf{p}}$,
$\delta_{\mathbf{p}}$ of $\mathbb{F}\left[  x_{\mathbf{p}\Downarrow}\right]  $
(not just of $\mathbb{F}\left(  x_{\mathbf{p}\Downarrow}\right)  $) such that
$\alpha_{\mathbf{p}}\delta_{\mathbf{p}}-\beta_{\mathbf{p}}\gamma_{\mathbf{p}%
}\neq0$ and $\mathfrak{N}_{\mathbf{p}}=\alpha_{\mathbf{p}}x_{\mathbf{p}}%
+\beta_{\mathbf{p}}$ and $\mathfrak{D}_{\mathbf{p}}=\gamma_{\mathbf{p}%
}x_{\mathbf{p}}+\delta_{\mathbf{p}}$. (Of course, the conditions
$\mathfrak{N}_{\mathbf{p}}=\alpha_{\mathbf{p}}x_{\mathbf{p}}+\beta
_{\mathbf{p}}$ and $\mathfrak{D}_{\mathbf{p}}=\gamma_{\mathbf{p}}%
x_{\mathbf{p}}+\delta_{\mathbf{p}}$ combined imply $Q_{\mathbf{p}}%
=\dfrac{\alpha_{\mathbf{p}}x_{\mathbf{p}}+\beta_{\mathbf{p}}}{\gamma
_{\mathbf{p}}x_{\mathbf{p}}+\delta_{\mathbf{p}}}$, hence the yearned-for
$\mathbf{P}$-triangularity.)

Let us first deal with two ``boundary'' cases: the case when $k=1$, and the
case when $k\neq1$ but $i=p$.

The case when $k=1$ is very easy. In fact, in this case, it is easy to prove
that $\mathfrak{N}_{\mathbf{p}}=1$ (using
(\ref{pf.Grasp.generic.short.step1.N})) and that $\mathfrak{D}_{\mathbf{p}%
}=\left(  -1\right)  ^{i+p}x_{\mathbf{p}}$ (using
(\ref{pf.Grasp.generic.short.step1.D})). Consequently, we can take
$\alpha_{\mathbf{p}}=0$, $\beta_{\mathbf{p}}=1$, $\gamma_{\mathbf{p}}=\left(
-1\right)  ^{i+p}$ and $\delta_{\mathbf{p}}=0$, and it is clear that all three
requirements $\alpha_{\mathbf{p}}\delta_{\mathbf{p}}-\beta_{\mathbf{p}}%
\gamma_{\mathbf{p}}\neq0$ and $\mathfrak{N}_{\mathbf{p}}=\alpha_{\mathbf{p}%
}x_{\mathbf{p}}+\beta_{\mathbf{p}}$ and $\mathfrak{D}_{\mathbf{p}}%
=\gamma_{\mathbf{p}}x_{\mathbf{p}}+\delta_{\mathbf{p}}$ are satisfied.

The case when $k\neq1$ but $i=p$ is not much harder. In this case,
(\ref{pf.Grasp.generic.short.step1.N}) simplifies to $\mathfrak{N}%
_{\mathbf{p}}=x_{\mathbf{p}}$, and (\ref{pf.Grasp.generic.short.step1.D})
simplifies to $\mathfrak{D}_{\mathbf{p}}=x_{\left(  p,1\right)  }$. Hence, we
can take $\alpha_{\mathbf{p}}=1$, $\beta_{\mathbf{p}}=0$, $\gamma_{\mathbf{p}%
}=0$ and $\delta_{\mathbf{p}}=x_{\left(  p,1\right)  }$ to achieve
$\alpha_{\mathbf{p}}\delta_{\mathbf{p}}-\beta_{\mathbf{p}}\gamma_{\mathbf{p}%
}\neq0$ and $\mathfrak{N}_{\mathbf{p}}=\alpha_{\mathbf{p}}x_{\mathbf{p}}%
+\beta_{\mathbf{p}}$ and $\mathfrak{D}_{\mathbf{p}}=\gamma_{\mathbf{p}%
}x_{\mathbf{p}}+\delta_{\mathbf{p}}$. Note that this choice of $\delta
_{\mathbf{p}}$ is legitimate because $x_{\left(  p,1\right)  }$ does lie in
$\mathbb{F}\left[  x_{\mathbf{p}\Downarrow}\right]  $ (since $\left(
p,1\right)  \in\left.  \mathbf{p}\Downarrow\right.  $).

Now that these two cases are settled, let us deal with the remaining case. So
we have neither $k=1$ nor $i=p$.

Consider the matrix $\operatorname*{rows}\nolimits_{i,i+1,...,p}\left(
\left(  I_{p}\mid C\right)  \left[  i+k-1:p+k\right]  \right)  $ (this matrix
appears on the right hand side of (\ref{pf.Grasp.generic.short.step1.N})).
Each entry of this matrix comes either from the matrix $I_{p}$ or from the
matrix $C$. If it comes from $I_{p}$, it clearly lies in $\mathbb{F}\left[
x_{\mathbf{p}\Downarrow}\right]  $. If it comes from $C$, it has the form
$x_{\mathbf{q}}$ for some $\mathbf{q}\in\mathbf{P}$, and this $\mathbf{q}$
belongs to $\left.  \mathbf{p}\Downarrow\right.  $ unless the entry is the
$\left(  1,p-i+1\right)  $-th entry. Therefore, each entry of the matrix
$\left(  I_{p}\mid C\right)  \left[  i+k-1:p+k\right]  $ apart from the
$\left(  1,p-i+1\right)  $-th entry lies in $\mathbb{F}\left[  x_{\mathbf{p}%
\Downarrow}\right]  $, whereas the $\left(  1,p-i+1\right)  $-th entry is
$x_{\mathbf{p}}$. Hence, if we use Laplace expansion with respect to the first
row to compute the determinant of this matrix, we obtain a formula of the form%
\begin{align*}
&  \det\left(  \operatorname*{rows}\nolimits_{i,i+1,...,p}\left(  \left(
I_{p}\mid C\right)  \left[  i+k-1:p+k\right]  \right)  \right) \\
&  =x_{\mathbf{p}}\cdot\left(  \text{some polynomial in entries lying in
}\mathbb{F}\left[  x_{\mathbf{p}\Downarrow}\right]  \right) \\
&  \ \ \ \ \ \ \ \ \ \ +\left(  \text{more polynomials in entries lying in
}\mathbb{F}\left[  x_{\mathbf{p}\Downarrow}\right]  \right) \\
&  \in\mathbb{F}\left[  x_{\mathbf{p}\Downarrow}\right]  \cdot x_{\mathbf{p}%
}+\mathbb{F}\left[  x_{\mathbf{p}\Downarrow}\right]  .
\end{align*}
In other words, there exist elements $\alpha_{\mathbf{p}}$ and $\beta
_{\mathbf{p}}$ of $\mathbb{F}\left[  x_{\mathbf{p}\Downarrow}\right]  $ such
that \newline$\det\left(  \operatorname*{rows}\nolimits_{i,i+1,...,p}\left(
\left(  I_{p}\mid C\right)  \left[  i+k-1:p+k\right]  \right)  \right)
=\alpha_{\mathbf{p}}x_{\mathbf{p}}+\beta_{\mathbf{p}}$. Consider these
$\alpha_{\mathbf{p}}$ and $\beta_{\mathbf{p}}$. We have%
\begin{align}
\mathfrak{N}_{\mathbf{p}}  &  =\mathfrak{N}_{\left(  i,k\right)  }=\det\left(
\operatorname*{rows}\nolimits_{i,i+1,...,p}\left(  \left(  I_{p}\mid C\right)
\left[  i+k-1:p+k\right]  \right)  \right)  \ \ \ \ \ \ \ \ \ \ \left(
\text{by (\ref{pf.Grasp.generic.short.step1.N})}\right)
\label{pf.Grasp.generic.short.step4.N0}\\
&  =\alpha_{\mathbf{p}}x_{\mathbf{p}}+\beta_{\mathbf{p}}.
\label{pf.Grasp.generic.short.step4.N}%
\end{align}

We can similarly deal with the matrix $\operatorname*{rows}%
\nolimits_{i,i+1,...,p}\left(  C_{q}\ \mid\ \left(  I_{p}\mid C\right)
\left[  i+k:p+k\right]  \right)  $ which appears on the right hand side of
(\ref{pf.Grasp.generic.short.step1.D}). Again, each entry of this matrix apart
from the $\left(  1,p-i+1\right)  $-th entry lies in $\mathbb{F}\left[
x_{\mathbf{p}\Downarrow}\right]  $, whereas the $\left(  1,p-i+1\right)  $-th
entry is $x_{\mathbf{p}}$. Using Laplace expansion again, we thus see that
\[
\det\left(  \operatorname*{rows}\nolimits_{i,i+1,...,p}\left(  C_{q}%
\ \mid\ \left(  I_{p}\mid C\right)  \left[  i+k:p+k\right]  \right)  \right)
\in\mathbb{F}\left[  x_{\mathbf{p}\Downarrow}\right]  \cdot x_{\mathbf{p}%
}+\mathbb{F}\left[  x_{\mathbf{p}\Downarrow}\right]  ,
\]
so that%
\[
\left(  -1\right)  ^{p-i}\det\left(  \operatorname*{rows}%
\nolimits_{i,i+1,...,p}\left(  C_{q}\ \mid\ \left(  I_{p}\mid C\right)
\left[  i+k:p+k\right]  \right)  \right)  \in\mathbb{F}\left[  x_{\mathbf{p}%
\Downarrow}\right]  \cdot x_{\mathbf{p}}+\mathbb{F}\left[  x_{\mathbf{p}%
\Downarrow}\right]  .
\]
Hence, there exist elements $\gamma_{\mathbf{p}}$ and $\delta_{\mathbf{p}}$ of
$\mathbb{F}\left[  x_{\mathbf{p}\Downarrow}\right]  $ such that \newline%
$\left(  -1\right)  ^{p-i}\det\left(  \operatorname*{rows}%
\nolimits_{i,i+1,...,p}\left(  C_{q}\ \mid\ \left(  I_{p}\mid C\right)
\left[  i+k:p+k\right]  \right)  \right)  =\gamma_{\mathbf{p}}x_{\mathbf{p}%
}+\delta_{\mathbf{p}}$. Consider these $\gamma_{\mathbf{p}}$ and
$\delta_{\mathbf{p}}$. We have%
\begin{align}
\mathfrak{D}_{\mathbf{p}}  &  =\mathfrak{D}_{\left(  i,k\right)  }=\left(
-1\right)  ^{p-i}\det\left(  \operatorname*{rows}\nolimits_{i,i+1,...,p}%
\left(  C_{q}\ \mid\ \left(  I_{p}\mid C\right)  \left[  i+k:p+k\right]
\right)  \right)  \ \ \ \ \ \ \ \ \ \ \left(  \text{by
(\ref{pf.Grasp.generic.short.step1.D})}\right)
\label{pf.Grasp.generic.short.step4.D0}\\
&  =\gamma_{\mathbf{p}}x_{\mathbf{p}}+\delta_{\mathbf{p}}.\nonumber
\end{align}

We thus have found elements $\alpha_{\mathbf{p}}$, $\beta_{\mathbf{p}}$,
$\gamma_{\mathbf{p}}$, $\delta_{\mathbf{p}}$ of $\mathbb{F}\left[
x_{\mathbf{p}\Downarrow}\right]  $ satisfying $\mathfrak{N}_{\mathbf{p}%
}=\alpha_{\mathbf{p}}x_{\mathbf{p}}+\beta_{\mathbf{p}}$ and $\mathfrak{D}%
_{\mathbf{p}}=\gamma_{\mathbf{p}}x_{\mathbf{p}}+\delta_{\mathbf{p}}$. In order
to finish the proof of $\mathbf{P}$-triangularity, we only need to show that
$\alpha_{\mathbf{p}}\delta_{\mathbf{p}}-\beta_{\mathbf{p}}\gamma_{\mathbf{p}%
}\neq0$.

In order to achieve this goal, we notice that
\[
\alpha_{\mathbf{p}}\underbrace{\mathfrak{D}_{\mathbf{p}}}_{=\gamma
_{\mathbf{p}}x_{\mathbf{p}}+\delta_{\mathbf{p}}}-\underbrace{\mathfrak{N}%
_{\mathbf{p}}}_{=\alpha_{\mathbf{p}}x_{\mathbf{p}}+\beta_{\mathbf{p}}}%
\gamma_{\mathbf{p}}=\alpha_{\mathbf{p}}\left(  \gamma_{\mathbf{p}%
}x_{\mathbf{p}}+\delta_{\mathbf{p}}\right)  -\left(  \alpha_{\mathbf{p}%
}x_{\mathbf{p}}+\beta_{\mathbf{p}}\right)  \gamma_{\mathbf{p}}=\alpha
_{\mathbf{p}}\delta_{\mathbf{p}}-\beta_{\mathbf{p}}\gamma_{\mathbf{p}}.
\]
Hence, proving $\alpha_{\mathbf{p}}\delta_{\mathbf{p}}-\beta_{\mathbf{p}%
}\gamma_{\mathbf{p}}\neq0$ is equivalent to proving $\alpha_{\mathbf{p}%
}\mathfrak{D}_{\mathbf{p}}-\mathfrak{N}_{\mathbf{p}}\gamma_{\mathbf{p}}\neq0$.
It is the latter that we are going to do, because $\alpha_{\mathbf{p}}$,
$\mathfrak{D}_{\mathbf{p}}$, $\mathfrak{N}_{\mathbf{p}}$ and $\gamma
_{\mathbf{p}}$ are easier to get our hands on than $\beta_{\mathbf{p}}$ and
$\delta_{\mathbf{p}}$.

So we need to prove that $\alpha_{\mathbf{p}}\mathfrak{D}_{\mathbf{p}%
}-\mathfrak{N}_{\mathbf{p}}\gamma_{\mathbf{p}}\neq0$. To do so, we look back
at our proof of
\[
\det\left(  \operatorname*{rows}\nolimits_{i,i+1,...,p}\left(  \left(
I_{p}\mid C\right)  \left[  i+k-1:p+k\right]  \right)  \right)  \in
\mathbb{F}\left[  x_{\mathbf{p}\Downarrow}\right]  \cdot x_{\mathbf{p}%
}+\mathbb{F}\left[  x_{\mathbf{p}\Downarrow}\right]  .
\]
This proof proceeded by applying Laplace expansion with respect to the first
row to the matrix $\operatorname*{rows}\nolimits_{i,i+1,...,p}\left(  \left(
I_{p}\mid C\right)  \left[  i+k-1:p+k\right]  \right)  $. The only term
involving $x_{\mathbf{p}}$ was%
\[
x_{\mathbf{p}}\cdot\left(  \text{some polynomial in entries lying in
}\mathbb{F}\left[  x_{\mathbf{p}\Downarrow}\right]  \right)  .
\]
Recalling the statement of Laplace expansion, we notice that \textquotedblleft
some polynomial in entries lying in $\mathbb{F}\left[  x_{\mathbf{p}%
\Downarrow}\right]  $\textquotedblright\ in this term is actually the $\left(
1,p-i+1\right)  $-th cofactor of the matrix $\operatorname*{rows}%
\nolimits_{i,i+1,...,p}\left(  \left(  I_{p}\mid C\right)  \left[
i+k-1:p+k\right]  \right)  $. Hence,%
\begin{align}
\alpha_{\mathbf{p}}  &  =\left(  \text{the }\left(  1,p-i+1\right)  \text{-th
cofactor of }\operatorname*{rows}\nolimits_{i,i+1,...,p}\left(  \left(
I_{p}\mid C\right)  \left[  i+k-1:p+k\right]  \right)  \right) \nonumber\\
&  =\left(  -1\right)  ^{p-i}\cdot\det\left(  \operatorname*{rows}%
\nolimits_{i+1,i+2,...,p}\left(  \left(  I_{p}\mid C\right)  \left[
i+k-1:p+k-1\right]  \right)  \right)  .
\label{pf.Grasp.generic.short.step4.alpha}%
\end{align}
Similarly,%
\begin{equation}
\gamma_{\mathbf{p}}=\det\left(  \operatorname*{rows}\nolimits_{i+1,i+2,...,p}%
\left(  C_{q}\ \mid\ \left(  I_{p}\mid C\right)  \left[  i+k:p+k-1\right]
\right)  \right)  \label{pf.Grasp.generic.short.step4.gamma}%
\end{equation}
(note that we lost the sign $\left(  -1\right)  ^{p-i}$ from
(\ref{pf.Grasp.generic.short.step4.D0}) since it got cancelled against the
$\left(  -1\right)  ^{p-i}$ arising from the definition of
a cofactor).

Now, recall that we have neither $k=1$ nor $i=p$. Hence, $\left(
i+1,k-1\right)  $ also belongs to $\mathbf{P}$, so we can apply
(\ref{pf.Grasp.generic.short.step1.N}) to $\left(  i+1,k-1\right)  $ in lieu
of $\left(  i,k\right)  $, and obtain%
\[
\mathfrak{N}_{\left(  i+1,k-1\right)  }=\det\left(  \operatorname*{rows}%
\nolimits_{i+1,i+2,...,p}\left(  \left(  I_{p}\mid C\right)  \left[
i+k-1:p+k-1\right]  \right)  \right)  .
\]
In light of this, (\ref{pf.Grasp.generic.short.step4.alpha}) becomes%
\[
\alpha_{\mathbf{p}}=\left(  -1\right)  ^{p-i}\cdot\mathfrak{N}_{\left(
i+1,k-1\right)  }.
\]
Similarly, we can apply (\ref{pf.Grasp.generic.short.step1.D}) to $\left(
i+1,k-1\right)  $ in lieu of $\left(  i,k\right)  $, and use this to rewrite
(\ref{pf.Grasp.generic.short.step4.gamma}) as%
\[
\gamma_{\mathbf{p}}=\left(  -1\right)  ^{p-\left(  i+1\right)  }%
\cdot\mathfrak{D}_{\left(  i+1,k-1\right)  }.
\]
Hence,%
\begin{align*}
&  \underbrace{\alpha_{\mathbf{p}}}_{=\left(  -1\right)  ^{p-i}\cdot
\mathfrak{N}_{\left(  i+1,k-1\right)  }}\mathfrak{D}_{\mathbf{p}}%
-\mathfrak{N}_{\mathbf{p}}\underbrace{\gamma_{\mathbf{p}}}_{=\left(
-1\right)  ^{p-\left(  i+1\right)  }\cdot\mathfrak{D}_{\left(  i+1,k-1\right)
}}\\
&  =\left(  -1\right)  ^{p-i}\cdot\mathfrak{N}_{\left(  i+1,k-1\right)  }%
\cdot\mathfrak{D}_{\mathbf{p}}-\mathfrak{N}_{\mathbf{p}}\cdot
\underbrace{\left(  -1\right)  ^{p-\left(  i+1\right)  }}_{=-\left(
-1\right)  ^{p-i}}\cdot\mathfrak{D}_{\left(  i+1,k-1\right)  }\\
&  =\left(  -1\right)  ^{p-i}\cdot\left(  \mathfrak{N}_{\left(
i+1,k-1\right)  }\mathfrak{D}_{\mathbf{p}}+\mathfrak{N}_{\mathbf{p}%
}\mathfrak{D}_{\left(  i+1,k-1\right)  }\right)  .
\end{align*}
Thus, we can shift our goal from proving $\alpha_{\mathbf{p}}\mathfrak{D}%
_{\mathbf{p}}-\mathfrak{N}_{\mathbf{p}}\gamma_{\mathbf{p}}\neq0$ to proving
$\mathfrak{N}_{\left(  i+1,k-1\right)  }\mathfrak{D}_{\mathbf{p}}%
+\mathfrak{N}_{\mathbf{p}}\mathfrak{D}_{\left(  i+1,k-1\right)  }\neq0$.

But this turns out to be surprisingly simple: Since $\mathbf{p}=\left(
i,k\right)  $, we have%
\begin{align}
&  \mathfrak{N}_{\left(  i+1,k-1\right)  }\mathfrak{D}_{\mathbf{p}%
}+\mathfrak{N}_{\mathbf{p}}\mathfrak{D}_{\left(  i+1,k-1\right)  }\nonumber\\
&  =\mathfrak{N}_{\left(  i+1,k-1\right)  }\mathfrak{D}_{\left(  i,k\right)
}+\mathfrak{N}_{\left(  i,k\right)  }\mathfrak{D}_{\left(  i+1,k-1\right)
}=\mathfrak{D}_{\left(  i,k\right)  }\cdot\mathfrak{N}_{\left(
i+1,k-1\right)  }+\mathfrak{N}_{\left(  i,k\right)  }\cdot\mathfrak{D}%
_{\left(  i+1,k-1\right)  }\nonumber\\
&  =\det\left(  \left(  I_{p}\mid C\right)  \left[  0:i\mid i+k:p+k\right]
\right)  \cdot\det\left(  \left(  I_{p}\mid C\right)  \left[  1:i+1\mid
i+k-1:p+k-1\right]  \right) \nonumber\\
&  \ \ \ \ \ \ \ \ \ \ +\det\left(  \left(  I_{p}\mid C\right)  \left[
1:i\mid i+k-1:p+k\right]  \right)  \cdot\det\left(  \left(  I_{p}\mid
C\right)  \left[  0:i+1\mid i+k:p+k-1\right]  \right) \nonumber\\
&  \ \ \ \ \ \ \ \ \ \ \left(
\begin{array}
[c]{c}%
\text{here, we just substituted }\mathfrak{D}_{\left(  i,k\right)  }\text{,
}\mathfrak{N}_{\left(  i+1,k-1\right)  }\text{, }\mathfrak{N}_{\left(
i,k\right)  }\text{ and }\mathfrak{D}_{\left(  i+1,k-1\right)  }\\
\text{by their definitions}%
\end{array}
\right) \nonumber\\
&  =\det\left(  \left(  I_{p}\mid C\right)  \left[  0:i\mid
i+k-1:p+k-1\right]  \right)  \cdot\det\left(  \left(  I_{p}\mid C\right)
\left[  1:i+1\mid i+k:p+k\right]  \right)
\label{pf.Grasp.generic.short.step4.pluck1}%
\end{align}
(by Theorem \ref{thm.pluecker.ptolemy}, applied to $p$, $p+q$, $\left(
I_{p}\mid C\right)  $, $1$, $i$, $i+k$ and $p+k-1$ instead of $u$, $v$, $A$,
$a$, $b$, $c$ and $d$). On the other hand, $\left(  i,k-1\right)  $ and
$\left(  i+1,k\right)  $ also belong to $\mathbf{P}$ and satisfy%
\[
\mathfrak{D}_{\left(  i,k-1\right)  }=\det\left(  \left(  I_{p}\mid C\right)
\left[  0:i\mid i+k-1:p+k-1\right]  \right)
\]
and
\[
\mathfrak{N}_{\left(  i+1,k\right)  }=\det\left(  \left(  I_{p}\mid C\right)
\left[  1:i+1\mid i+k:p+k\right]  \right)
\]
(by the respective definitions of $\mathfrak{D}_{\left(  i,k-1\right)  }$ and
$\mathfrak{N}_{\left(  i+1,k\right)  }$). Hence,
(\ref{pf.Grasp.generic.short.step4.pluck1}) becomes%
\begin{align*}
&  \mathfrak{N}_{\left(  i+1,k-1\right)  }\mathfrak{D}_{\mathbf{p}%
}+\mathfrak{N}_{\mathbf{p}}\mathfrak{D}_{\left(  i+1,k-1\right)  }\\
&  =\underbrace{\det\left(  \left(  I_{p}\mid C\right)  \left[  0:i\mid
i+k-1:p+k-1\right]  \right)  }_{=\mathfrak{D}_{\left(  i,k-1\right)  }}%
\cdot\underbrace{\det\left(  \left(  I_{p}\mid C\right)  \left[  1:i+1\mid
i+k:p+k\right]  \right)  }_{=\mathfrak{N}_{\left(  i+1,k\right)  }}\\
&  =\mathfrak{D}_{\left(  i,k-1\right)  }\cdot\mathfrak{N}_{\left(
i+1,k\right)  }\neq0
\end{align*}
(since the result of Step 2 shows that $\mathfrak{D}_{\left(  i,k-1\right)  }$
and $\mathfrak{N}_{\left(  i+1,k\right)  }$ are nonzero). This finishes our
proof of $\mathfrak{N}_{\left(  i+1,k-1\right)  }\mathfrak{D}_{\mathbf{p}%
}+\mathfrak{N}_{\mathbf{p}}\mathfrak{D}_{\left(  i+1,k-1\right)  }\neq0$, thus
also of $\alpha_{\mathbf{p}}\mathfrak{D}_{\mathbf{p}}-\mathfrak{N}%
_{\mathbf{p}}\gamma_{\mathbf{p}}\neq0$, hence also of $\alpha_{\mathbf{p}%
}\delta_{\mathbf{p}}-\beta_{\mathbf{p}}\gamma_{\mathbf{p}}\neq0$, and
ultimately of the $\mathbf{P}$-triangularity of the family $\left(
Q_{\mathbf{p}}\right)  _{\mathbf{p}\in\mathbf{P}}$.

\textit{Details of Step 5:} Recall that our goal is to prove the existence of
a matrix $\widetilde{D}\in\left(  \mathbb{F}\left(  x_{\mathbf{P}}\right)
\right)  ^{p\times\left(  p+q\right)  }$ satisfying
(\ref{pf.Grasp.generic.short.reduce-to-rational}).

Since Step 4, we know that the family $\left(  Q_{\mathbf{p}}\right)
_{\mathbf{p}\in\mathbf{P}}\in\left(  \mathbb{F}\left(  x_{\mathbf{P}}\right)
\right)  ^{\mathbf{P}}$ is $\mathbf{P}$-triangular. Hence, Lemma
\ref{lem.algebraic.triangularity.short} \textbf{(b)} shows that there exists a
$\mathbf{P}$-triangular family $\left(  R_{\mathbf{p}}\right)  _{\mathbf{p}%
\in\mathbf{P}}\in\left(  \mathbb{F}\left(  x_{\mathbf{P}}\right)  \right)
^{\mathbf{P}}$ such that every $\mathbf{q}\in\mathbf{P}$ satisfies
$Q_{\mathbf{q}}\left(  \left(  R_{\mathbf{p}}\right)  _{\mathbf{p}%
\in\mathbf{P}}\right)  =x_{\mathbf{q}}$. Consider this $\left(  R_{\mathbf{p}%
}\right)  _{\mathbf{p}\in\mathbf{P}}$. Applying Lemma
\ref{lem.algebraic.triangularity.short} \textbf{(a)} to this family $\left(
R_{\mathbf{p}}\right)  _{\mathbf{p}\in\mathbf{P}}$, we conclude that $\left(
R_{\mathbf{p}}\right)  _{\mathbf{p}\in\mathbf{P}}$ is algebraically independent.

In Step 3, we have shown that $Q_{\mathbf{p}}=\left(  \operatorname*{Grasp}%
\nolimits_{0}\left(  I_{p}\mid C\right)  \right)  \left(  \mathbf{p}\right)  $
for every $\mathbf{p}\in\mathbf{P}$. Renaming $\mathbf{p}$ as $\mathbf{q}$, we
rewrite this as follows:%
\begin{equation}
Q_{\mathbf{q}}=\left(  \operatorname*{Grasp}\nolimits_{0}\left(  I_{p}\mid
C\right)  \right)  \left(  \mathbf{q}\right)  \ \ \ \ \ \ \ \ \ \ \text{for
every }\mathbf{q}\in\mathbf{P}.
\label{pf.algebraic.triangularity.short.step5.1}%
\end{equation}
Now, let $\widetilde{C}\in\left(  \mathbb{F}\left(  x_{\mathbf{P}}\right)
\right)  ^{p\times\left(  p+q\right)  }$ denote the matrix obtained from the
matrix $C\in\left(  \mathbb{F}\left[  x_{\mathbf{P}}\right]  \right)
^{p\times\left(  p+q\right)  }$ by substituting $\left(  R_{\mathbf{p}%
}\right)  _{\mathbf{p}\in\mathbf{P}}$ for the variables $\left(
x_{\mathbf{p}}\right)  _{\mathbf{p}\in\mathbf{P}}$. Since
(\ref{pf.algebraic.triangularity.short.step5.1}) is an identity between
rational functions in the variables $\left(  x_{\mathbf{p}}\right)
_{\mathbf{p}\in\mathbf{P}}$, we thus can substitute $\left(  R_{\mathbf{p}%
}\right)  _{\mathbf{p}\in\mathbf{P}}$ for the variables $\left(
x_{\mathbf{p}}\right)  _{\mathbf{p}\in\mathbf{P}}$ in
(\ref{pf.algebraic.triangularity.short.step5.1})\footnote{The substitution
does not suffer from vanishing denominators because $\left(  R_{\mathbf{p}%
}\right)  _{\mathbf{p}\in\mathbf{P}}$ is algebraically independent.}, and
obtain%
\[
Q_{\mathbf{q}}\left(  \left(  R_{\mathbf{p}}\right)  _{\mathbf{p}\in
\mathbf{P}}\right)  =\left(  \operatorname*{Grasp}\nolimits_{0}\left(
I_{p}\mid\widetilde{C}\right)  \right)  \left(  \mathbf{q}\right)
\ \ \ \ \ \ \ \ \ \ \text{for every }\mathbf{q}\in\mathbf{P}%
\]
(since this substitution takes the matrix $C$ to $\widetilde{C}$). But since
$Q_{\mathbf{q}}\left(  \left(  R_{\mathbf{p}}\right)  _{\mathbf{p}%
\in\mathbf{P}}\right)  =x_{\mathbf{q}}$ for every $\mathbf{q}\in\mathbf{P}$,
this rewrites as%
\[
x_{\mathbf{q}}=\left(  \operatorname*{Grasp}\nolimits_{0}\left(  I_{p}%
\mid\widetilde{C}\right)  \right)  \left(  \mathbf{q}\right)
\ \ \ \ \ \ \ \ \ \ \text{for every }\mathbf{q}\in\mathbf{P}.
\]
Upon renaming $\mathbf{q}$ as $\mathbf{p}$ again, this becomes%
\[
x_{\mathbf{p}}=\left(  \operatorname*{Grasp}\nolimits_{0}\left(  I_{p}%
\mid\widetilde{C}\right)  \right)  \left(  \mathbf{p}\right)
\ \ \ \ \ \ \ \ \ \ \text{for every }\mathbf{p}\in\mathbf{P}.
\]
Hence, there exists a matrix $\widetilde{D}\in\left(  \mathbb{F}\left(
x_{\mathbf{P}}\right)  \right)  ^{p\times\left(  p+q\right)  }$ satisfying
(\ref{pf.Grasp.generic.short.reduce-to-rational}) (namely, $\widetilde{D}%
=\left(  I_{p}\mid\widetilde{C}\right)  $). Thus, as we know, Proposition
\ref{prop.Grasp.generic} is proven.
\end{proof}
\end{vershort}

\begin{verlong}
Let us prepare for the proof by introducing some more notation.\footnote{The
reason why Definition \ref{def.algebraic.triangularity} begins with ``Let
$\mathbb{F}$ be a field.'' rather than ``Let $\mathbb{K}$ be a field.'' is
that we will later apply its statement to a field different from what we call
$\mathbb{K}$ in the context of birational rowmotion.}

\begin{definition}
\label{def.algebraic.triangularity}Let $\mathbb{F}$ be a field. Let
$\mathbf{P}$ be a finite set.

\textbf{(a)} Let $x_{\mathbf{p}}$ be a new symbol for every $\mathbf{p}%
\in\mathbf{P}$. We will denote by $\mathbb{F}\left(  x_{\mathbf{P}}\right)  $
the field of rational functions over $\mathbb{F}$ in the indeterminates
$x_{\mathbf{p}}$ with $\mathbf{p}$ ranging over all elements of $\mathbf{P}$
(hence altogether $\left\vert \mathbf{P}\right\vert $ indeterminates). (Thus,
$\mathbb{F}\left(  x_{\mathbf{P}}\right)  =\mathbb{F}\left(  x_{\mathbf{p}%
_{1}},x_{\mathbf{p}_{2}},...,x_{\mathbf{p}_{n}}\right)  $ if $\mathbf{P}$ is
written in the form $\mathbf{P}=\left\{  \mathbf{p}_{1},\mathbf{p}%
_{2},...,\mathbf{p}_{n}\right\}  $.) The symbols $x_{\mathbf{p}}$ are
understood to be distinct, and are used as commuting indeterminates.

\textbf{(b)} If $\mathbf{Q}$ is a subset of $\mathbf{P}$, then $\mathbb{F}%
\left(  x_{\mathbf{Q}}\right)  $ can be canonically embedded into
$\mathbb{F}\left(  x_{\mathbf{P}}\right)  $. We regard this embedding as an inclusion.

\textbf{(c)} Let $\mathbb{K}$ be a field extension of $\mathbb{F}$. Let $f$ be
an element of $\mathbb{F}\left(  x_{\mathbf{P}}\right)  $. If $\left(
a_{\mathbf{p}}\right)  _{\mathbf{p}\in\mathbf{P}}\in\mathbb{K}^{\mathbf{P}}$
is a family of elements of $\mathbb{K}$ indexed by elements of $\mathbf{P}$,
then we let $f\left(  \left(  a_{\mathbf{p}}\right)  _{\mathbf{p}\in
\mathbf{P}}\right)  $ denote the element of $\mathbb{K}$ obtained by
substituting $a_{\mathbf{p}}$ for $x_{\mathbf{p}}$ for each $\mathbf{p}%
\in\mathbf{P}$ in the rational function $f$. This $f\left(  \left(
a_{\mathbf{p}}\right)  _{\mathbf{p}\in\mathbf{P}}\right)  $ is defined only if
the substitution does not render the denominator equal to $0$. If $\mathbb{K}$
is infinite, this shows that $f\left(  \left(  a_{\mathbf{p}}\right)
_{\mathbf{p}\in\mathbf{P}}\right)  $ is defined for almost all $\left(
a_{\mathbf{p}}\right)  _{\mathbf{p}\in\mathbf{P}}\in\mathbb{K}^{\mathbf{P}}$
(with respect to the Zariski topology).
\end{definition}

We are now going to study automorphisms of the $\mathbb{F}$-algebra
$\mathbb{F}\left(  x_{\mathbf{P}}\right)  $, and more generally, $\mathbb{F}%
$-algebra homomorphisms from $\mathbb{F}\left(  x_{\mathbf{P}}\right)  $:

\begin{definition}
\label{def.algebraic.triangularity.ev}Let $\mathbb{F}$ be a field. Let
$\mathbb{K}$ be a field extension of $\mathbb{F}$. Let $\mathbf{P}$ be a
finite set. Let $\left(  Q_{\mathbf{p}}\right)  _{\mathbf{p}\in\mathbf{P}}%
\in\mathbb{K}^{\mathbf{P}}$ be an algebraically independent (over $\mathbb{F}%
$) family of elements of $\mathbb{K}$ indexed by elements of $\mathbf{P}$.
Then, for every $f\in\mathbb{F}\left(  x_{\mathbf{P}}\right)  $, the element
$f\left(  \left(  Q_{\mathbf{p}}\right)  _{\mathbf{p}\in\mathbf{P}}\right)  $
of $\mathbb{K}$ is well-defined\footnotemark. Hence, we can define a map
$\mathbb{F}\left(  x_{\mathbf{P}}\right)  \rightarrow\mathbb{K}$ by sending
every $f\in\mathbb{F}\left(  x_{\mathbf{P}}\right)  $ to $f\left(  \left(
Q_{\mathbf{p}}\right)  _{\mathbf{p}\in\mathbf{P}}\right)  $. Denote this map
by $\operatorname*{rev}\nolimits_{\left(  Q_{\mathbf{p}}\right)
_{\mathbf{p}\in\mathbf{P}}}$. It is easy to see that $\operatorname*{rev}%
\nolimits_{\left(  Q_{\mathbf{p}}\right)  _{\mathbf{p}\in\mathbf{P}}}$ is an
$\mathbb{F}$-algebra homomorphism and a field homomorphism from $\mathbb{F}%
\left(  x_{\mathbf{P}}\right)  $ to $\mathbb{K}$.
\end{definition}

\footnotetext{\textit{Proof.} Let $f\in\mathbb{F}\left(  x_{\mathbf{P}%
}\right)  $. Denote by $\mathbb{F}\left[  x_{\mathbf{P}}\right]  $ the
polynomial ring over $\mathbb{F}$ in the indeterminates $x_{\mathbf{p}}$ with
$\mathbf{p}$ ranging over all elements of $\mathbf{P}$. Then, $\mathbb{F}%
\left(  x_{\mathbf{P}}\right)  $ is the quotient field of $\mathbb{F}\left[
x_{\mathbf{P}}\right]  $. Hence, $f\in\mathbb{F}\left(  x_{\mathbf{P}}\right)
$ can be written as a quotient $\dfrac{q}{r}$ of two elements $q$ and $r$ of
$\mathbb{F}\left[  x_{\mathbf{P}}\right]  $ with $r\neq0$. Consider these $q$
and $r$. Since $\left(  Q_{\mathbf{p}}\right)  _{\mathbf{p}\in\mathbf{P}}$ is
algebraically independent, we have $r\left(  \left(  Q_{\mathbf{p}}\right)
_{\mathbf{p}\in\mathbf{P}}\right)  \neq0$ (because $r\neq0$). Hence,
$\dfrac{q\left(  \left(  Q_{\mathbf{p}}\right)  _{\mathbf{p}\in\mathbf{P}%
}\right)  }{r\left(  \left(  Q_{\mathbf{p}}\right)  _{\mathbf{p}\in\mathbf{P}%
}\right)  }$ is well-defined. But since $f=\dfrac{q}{r}$, we have $f\left(
\left(  Q_{\mathbf{p}}\right)  _{\mathbf{p}\in\mathbf{P}}\right)  =\dfrac
{q}{r}\left(  \left(  Q_{\mathbf{p}}\right)  _{\mathbf{p}\in\mathbf{P}%
}\right)  =\dfrac{q\left(  \left(  Q_{\mathbf{p}}\right)  _{\mathbf{p}%
\in\mathbf{P}}\right)  }{r\left(  \left(  Q_{\mathbf{p}}\right)
_{\mathbf{p}\in\mathbf{P}}\right)  }$. Thus, $f\left(  \left(  Q_{\mathbf{p}%
}\right)  _{\mathbf{p}\in\mathbf{P}}\right)  $ is well-defined (since
$\dfrac{q\left(  \left(  Q_{\mathbf{p}}\right)  _{\mathbf{p}\in\mathbf{P}%
}\right)  }{r\left(  \left(  Q_{\mathbf{p}}\right)  _{\mathbf{p}\in\mathbf{P}%
}\right)  }$ is well-defined), qed.} Note that the name ``rev'' in the
notation $\operatorname*{rev}\nolimits_{\left(  Q_{\mathbf{p}}\right)
_{\mathbf{p}\in\mathbf{P}}}$ introduced in Definition
\ref{def.algebraic.triangularity.ev} has been chosen to signify ``\textbf{r}%
ational \textbf{ev}aluation'', as the homomorphism $\operatorname*{rev}%
\nolimits_{\left(  Q_{\mathbf{p}}\right)  _{\mathbf{p}\in\mathbf{P}}}$ is a
rational-functions analogue of the well-known evaluation homomorphisms from
polynomial rings.

Of course, nothing stops us from applying Definition
\ref{def.algebraic.triangularity.ev} to the particular situation when
$\mathbb{K}=\mathbb{F}\left(  x_{\mathbf{P}}\right)  $. In this situation, the
resulting homomorphism $\operatorname*{rev}\nolimits_{\left(  Q_{\mathbf{p}%
}\right)  _{\mathbf{p}\in\mathbf{P}}}$ will be a field endomorphism of
$\mathbb{F}\left(  x_{\mathbf{P}}\right)  $. We are going to study a certain
(even more special) situation in which this homomorphism turns out to be
invertible. Namely, this will be the situation in which the family $\left(
Q_{\mathbf{p}}\right)  _{\mathbf{p}\in\mathbf{P}}$ is, in some sense,
``triangular'' with respect to the variables $x_{\mathbf{p}}$. We formalize
this notion in the following definition:

\begin{definition}
\label{def.algebraic.triangularity.triang}Let $\mathbb{F}$ be a field. Let
$\mathbf{P}$ be a finite totally ordered set, and let $\vartriangleleft$ be
the smaller relation of $\mathbf{P}$. For every $\mathbf{p}\in\mathbf{P}$, let
$\mathbf{p}\Downarrow$ denote the subset $\left\{  \mathbf{v}\in
\mathbf{P}\ \mid\ \mathbf{v}\vartriangleleft\mathbf{p}\right\}  $ of
$\mathbf{P}$. For every $\mathbf{p}\in\mathbf{P}$, let $Q_{\mathbf{p}}$ be an
element of $\mathbb{F}\left(  x_{\mathbf{P}}\right)  $.

We say that the family $\left(  Q_{\mathbf{p}}\right)  _{\mathbf{p}%
\in\mathbf{P}}$ is $\mathbf{P}$\textit{-triangular} if and only if the
following condition holds:

\textit{Algebraic triangularity condition:} For every $\mathbf{p}\in
\mathbf{P}$, there exist elements $\alpha_{\mathbf{p}}$, $\beta_{\mathbf{p}}$,
$\gamma_{\mathbf{p}}$, $\delta_{\mathbf{p}}$ of $\mathbb{F}\left(
x_{\mathbf{p}\Downarrow}\right)  $ such that $\alpha_{\mathbf{p}}%
\delta_{\mathbf{p}}-\beta_{\mathbf{p}}\gamma_{\mathbf{p}}\neq0$ and
$Q_{\mathbf{p}}=\dfrac{\alpha_{\mathbf{p}}x_{\mathbf{p}}+\beta_{\mathbf{p}}%
}{\gamma_{\mathbf{p}}x_{\mathbf{p}}+\delta_{\mathbf{p}}}$.\ \ \ \ \footnotemark
\end{definition}

\footnotetext{Notice that the fraction $\dfrac{\alpha_{\mathbf{p}%
}x_{\mathbf{p}}+\beta_{\mathbf{p}}}{\gamma_{\mathbf{p}}x_{\mathbf{p}}%
+\delta_{\mathbf{p}}}$ is well-defined for any four elements $\alpha
_{\mathbf{p}}$, $\beta_{\mathbf{p}}$, $\gamma_{\mathbf{p}}$, $\delta
_{\mathbf{p}}$ of $\mathbb{F}\left(  x_{\mathbf{p}\Downarrow}\right)  $ such
that $\alpha_{\mathbf{p}}\delta_{\mathbf{p}}-\beta_{\mathbf{p}}\gamma
_{\mathbf{p}}\neq0$. Indeed, let $\alpha_{\mathbf{p}}$, $\beta_{\mathbf{p}}$,
$\gamma_{\mathbf{p}}$, $\delta_{\mathbf{p}}$ be four elements of
$\mathbb{F}\left(  x_{\mathbf{p}\Downarrow}\right)  $ such that $\alpha
_{\mathbf{p}}\delta_{\mathbf{p}}-\beta_{\mathbf{p}}\gamma_{\mathbf{p}}\neq0$.
Then, $\left(  \gamma_{\mathbf{p}},\delta_{\mathbf{p}}\right)  \neq\left(
0,0\right)  $ (since $\alpha_{\mathbf{p}}\delta_{\mathbf{p}}-\beta
_{\mathbf{p}}\gamma_{\mathbf{p}}\neq0$). But we can regard $\gamma
_{\mathbf{p}}x_{\mathbf{p}}+\delta_{\mathbf{p}}$ as a polynomial in the
variable $x_{\mathbf{p}}$ over the ring $\mathbb{F}\left(  x_{\mathbf{p}%
\Downarrow}\right)  $ (since $\gamma_{\mathbf{p}}$ and $\delta_{\mathbf{p}}$
both lie in this ring $\mathbb{F}\left(  x_{\mathbf{p}\Downarrow}\right)  $,
whereas $x_{\mathbf{p}}$ is an indeterminate not appearing in this ring). This
polynomial has coefficients $\gamma_{\mathbf{p}}$ and $\delta_{\mathbf{p}}$,
and these coefficients are not both zero (since $\left(  \gamma_{\mathbf{p}%
},\delta_{\mathbf{p}}\right)  \neq\left(  0,0\right)  $). Hence, the
polynomial $\gamma_{\mathbf{p}}x_{\mathbf{p}}+\delta_{\mathbf{p}}$ is nonzero.
Thus, the fraction $\dfrac{\alpha_{\mathbf{p}}x_{\mathbf{p}}+\beta
_{\mathbf{p}}}{\gamma_{\mathbf{p}}x_{\mathbf{p}}+\delta_{\mathbf{p}}}$ is
well-defined, qed.}The next lemma is similar to the fact that triangular
matrices with invertible elements along the diagonal are invertible:

\begin{lemma}
\label{lem.algebraic.triangularity}Let $\mathbb{F}$ be a field. Let
$\mathbf{P}$ be a finite totally ordered set. For every $\mathbf{p}%
\in\mathbf{P}$, let $Q_{\mathbf{p}}$ be an element of $\mathbb{F}\left(
x_{\mathbf{P}}\right)  $. Assume that $\left(  Q_{\mathbf{p}}\right)
_{\mathbf{p}\in\mathbf{P}}$ is a $\mathbf{P}$-triangular family. Then:

\textbf{(a)} The family $\left(  Q_{\mathbf{p}}\right)  _{\mathbf{p}%
\in\mathbf{P}}\in\left(  \mathbb{F}\left(  x_{\mathbf{P}}\right)  \right)
^{\mathbf{P}}$ is algebraically independent (over $\mathbb{F}$). Hence, the
$\mathbb{F}$-algebra homomorphism $\operatorname*{rev}\nolimits_{\left(
Q_{\mathbf{p}}\right)  _{\mathbf{p}\in\mathbf{P}}}$ from $\mathbb{F}\left(
x_{\mathbf{P}}\right)  $ to $\mathbb{F}\left(  x_{\mathbf{P}}\right)  $ is well-defined.

\textbf{(b)} There exists a $\mathbf{P}$-triangular family $\left(
R_{\mathbf{p}}\right)  _{\mathbf{p}\in\mathbf{P}}\in\left(  \mathbb{F}\left(
x_{\mathbf{P}}\right)  \right)  ^{\mathbf{P}}$ such that the maps
$\operatorname*{rev}\nolimits_{\left(  Q_{\mathbf{p}}\right)  _{\mathbf{p}%
\in\mathbf{P}}}$ and $\operatorname*{rev}\nolimits_{\left(  R_{\mathbf{p}%
}\right)  _{\mathbf{p}\in\mathbf{P}}}$ are mutually inverse.
\end{lemma}

Before we prove this lemma, let us state two simpler facts:

\begin{lemma}
\label{lem.algebraic.triangularity.id}Let $\mathbb{F}$ be a field. Let
$\mathbb{K}$ be a field extension of $\mathbb{F}$. Let $\mathbf{P}$ be a
finite set. Let $\mathfrak{a}$ and $\mathfrak{b}$ be two $\mathbb{F}$-algebra
homomorphisms from $\mathbb{F}\left(  x_{\mathbf{P}}\right)  $ to $\mathbb{K}%
$. Assume that $\mathfrak{a}\left(  x_{\mathbf{p}}\right)  =\mathfrak{b}%
\left(  x_{\mathbf{p}}\right)  $ for every $\mathbf{p}\in\mathbf{P}$. Then,
$\mathfrak{a}=\mathfrak{b}$.
\end{lemma}

Lemma \ref{lem.algebraic.triangularity.id} can be very easily proven using the
universal property of the polynomial ring (in this case, the polynomial ring
over $\mathbb{F}$ in the indeterminates $x_{\mathbf{p}}$ with $\mathbf{p}$
ranging over all elements of $\mathbf{P}$) and the fact that every algebra
homomorphism from the quotient field of a ring $A$ to a field is uniquely
determined by its restriction to $A$ (this fact is applied to $A$ being the
polynomial ring just mentioned). Here is the full proof, for the sake of completeness:

\begin{proof}
[Proof of Lemma \ref{lem.algebraic.triangularity.id} (sketched).]Let
$\mathbb{F}\left[  x_{\mathbf{P}}\right]  $ denote the polynomial ring over
$\mathbb{F}$ in the indeterminates $x_{\mathbf{p}}$ with $\mathbf{p}$ ranging
over all elements of $\mathbf{P}$. We identify $\mathbb{F}\left[
x_{\mathbf{P}}\right]  $ with a subring of $\mathbb{F}\left(  x_{\mathbf{P}%
}\right)  $. Then, $\mathbb{F}\left(  x_{\mathbf{P}}\right)  $ is the quotient
field of $\mathbb{F}\left[  x_{\mathbf{P}}\right]  $. Let $f\in\mathbb{F}%
\left(  x_{\mathbf{P}}\right)  $. Then, we can write $f$ in the form
$f=\dfrac{u}{v}$ for some $u\in\mathbb{F}\left[  x_{\mathbf{P}}\right]  $ and
some $v\in\mathbb{F}\left[  x_{\mathbf{P}}\right]  $ satisfying $v\neq0$
(since $\mathbb{F}\left(  x_{\mathbf{P}}\right)  $ is the quotient field of
$\mathbb{F}\left[  x_{\mathbf{P}}\right]  $). Consider these $u$ and $v$.
Since $f=\dfrac{u}{v}$, we have $vf=u$. Since $v\neq0$, we know that $v$ is
invertible in $\mathbb{F}\left(  x_{\mathbf{P}}\right)  $. That is, the
multiplicative inverse $v^{-1}$ of $v$ is well-defined. Since $\mathfrak{a}$
is an $\mathbb{F}$-algebra homomorphism, we have $\mathfrak{a}\left(
v\right)  \cdot\mathfrak{a}\left(  v^{-1}\right)  =\mathfrak{a}\left(
\underbrace{vv^{-1}}_{=1}\right)  =\mathfrak{a}\left(  1\right)  =1$ (again
since $\mathfrak{a}$ is an $\mathbb{F}$-algebra homomorphism). Hence,
$\mathfrak{a}\left(  v\right)  $ is invertible. Now, since $\mathfrak{a}$ is
an $\mathbb{F}$-algebra homomorphism, we have $\mathfrak{a}\left(  v\right)
\cdot\mathfrak{a}\left(  f\right)  =\mathfrak{a}\left(  vf\right)
=\mathfrak{a}\left(  u\right)  $ (since $vf=u$). Thus, $\mathfrak{a}\left(
f\right)  =\dfrac{\mathfrak{a}\left(  u\right)  }{\mathfrak{a}\left(
v\right)  }$ (since $\mathfrak{a}\left(  v\right)  $ is invertible). The same
argument for $\mathfrak{b}$ instead of $\mathfrak{a}$ shows that
$\mathfrak{b}\left(  f\right)  =\dfrac{\mathfrak{b}\left(  u\right)
}{\mathfrak{b}\left(  v\right)  }$. Now, recall the universal property of the
polynomial ring $\mathbb{F}\left[  x_{\mathbf{P}}\right]  $. It states that
for any commutative $\mathbb{F}$-algebra $A$ and any family $\left(
\eta_{\mathbf{p}}\right)  _{\mathbf{p}\in\mathbf{P}}$ of elements of $A$,
there is a unique $\mathbb{F}$-algebra homomorphism $\varphi:\mathbb{F}\left[
x_{\mathbf{P}}\right]  \rightarrow A$ such that every $\mathbf{p}\in
\mathbf{P}$ satisfies $\varphi\left(  x_{\mathbf{p}}\right)  =\eta
_{\mathbf{p}}$. Applied to $A=\mathbb{K}$ and $\left(  \eta_{\mathbf{p}%
}\right)  _{\mathbf{p}\in\mathbf{P}}=\left(  \mathfrak{a}\left(
x_{\mathbf{p}}\right)  \right)  _{\mathbf{p}\in\mathbf{P}}$, we conclude the
following: There is a unique $\mathbb{F}$-algebra homomorphism $\varphi
:\mathbb{F}\left[  x_{\mathbf{P}}\right]  \rightarrow\mathbb{K}$ such that
every $\mathbf{p}\in\mathbf{P}$ satisfies $\varphi\left(  x_{\mathbf{p}%
}\right)  =\mathfrak{a}\left(  x_{\mathbf{p}}\right)  $. In particular, if
$\varphi_{1}$ and $\varphi_{2}$ are two $\mathbb{F}$-algebra homomorphisms
$\varphi:\mathbb{F}\left[  x_{\mathbf{P}}\right]  \rightarrow\mathbb{K}$ such
that every $\mathbf{p}\in\mathbf{P}$ satisfies $\varphi\left(  x_{\mathbf{p}%
}\right)  =\mathfrak{a}\left(  x_{\mathbf{p}}\right)  $, then $\varphi
_{1}=\varphi_{2}$. But we can apply this to $\varphi_{1}=\mathfrak{a}%
\mid_{\mathbb{F}\left[  x_{\mathbf{P}}\right]  }$ and $\varphi_{2}%
=\mathfrak{b}\mid_{\mathbb{F}\left[  x_{\mathbf{P}}\right]  }$ (since every
$\mathbf{p}\in\mathbf{P}$ satisfies $\left(  \mathfrak{a}\mid_{\mathbb{F}%
\left[  x_{\mathbf{P}}\right]  }\right)  \left(  x_{\mathbf{p}}\right)
=\mathfrak{a}\left(  x_{\mathbf{p}}\right)  $ and since every $\mathbf{p}%
\in\mathbf{P}$ satisfies $\left(  \mathfrak{b}\mid_{\mathbb{F}\left[
x_{\mathbf{P}}\right]  }\right)  \left(  x_{\mathbf{p}}\right)  =\mathfrak{b}%
\left(  x_{\mathbf{p}}\right)  =\mathfrak{a}\left(  x_{\mathbf{p}}\right)  $),
and thus conclude that $\mathfrak{a}\mid_{\mathbb{F}\left[  x_{\mathbf{P}%
}\right]  }=\mathfrak{b}\mid_{\mathbb{F}\left[  x_{\mathbf{P}}\right]  }$.
Since $u\in\mathbb{F}\left[  x_{\mathbf{P}}\right]  $, we have $\mathfrak{a}%
\left(  u\right)  =\underbrace{\left(  \mathfrak{a}\mid_{\mathbb{F}\left[
x_{\mathbf{P}}\right]  }\right)  }_{=\mathfrak{b}\mid_{\mathbb{F}\left[
x_{\mathbf{P}}\right]  }}\left(  u\right)  =\left(  \mathfrak{b}%
\mid_{\mathbb{F}\left[  x_{\mathbf{P}}\right]  }\right)  \left(  u\right)
=\mathfrak{b}\left(  u\right)  $. Similarly, we obtain $\mathfrak{a}\left(
v\right)  =\mathfrak{b}\left(  v\right)  $ by applying the same reasoning to
$v$ instead of $u$. Now,
\begin{align*}
\mathfrak{a}\left(  f\right)   &  =\dfrac{\mathfrak{a}\left(  u\right)
}{\mathfrak{a}\left(  v\right)  }=\dfrac{\mathfrak{b}\left(  u\right)
}{\mathfrak{b}\left(  v\right)  }\ \ \ \ \ \ \ \ \ \ \left(  \text{since
}\mathfrak{a}\left(  u\right)  =\mathfrak{b}\left(  u\right)  \text{ and
}\mathfrak{a}\left(  v\right)  =\mathfrak{b}\left(  v\right)  \right) \\
&  =\mathfrak{b}\left(  f\right)  .
\end{align*}
Now, forget that we fixed $f$. We thus have shown that every $f\in
\mathbb{F}\left(  x_{\mathbf{P}}\right)  $ satisfies $\mathfrak{a}\left(
f\right)  =\mathfrak{b}\left(  f\right)  $. In other words, $\mathfrak{a}%
=\mathfrak{b}$. This proves Lemma \ref{lem.algebraic.triangularity.id}.
\end{proof}

Our next auxiliary lemma shows how algebraic independency of a family can
sometimes be deduced from that of a subfamily:

\begin{lemma}
\label{lem.algebraic.triangularity.step}Let $\mathbb{F}$ be a field. Let
$\mathbf{P}$ be a finite set, and let $\mathbf{q}\in\mathbf{P}$. Let $\left(
Q_{\mathbf{p}}\right)  _{\mathbf{p}\in\mathbf{P}}\in\left(  \mathbb{F}\left(
x_{\mathbf{P}}\right)  \right)  ^{\mathbf{P}}$ be a family of elements of
$\mathbb{F}\left(  x_{\mathbf{P}}\right)  $ indexed by elements of
$\mathbf{P}$. Assume that $Q_{\mathbf{p}}\in\mathbb{F}\left(  x_{\mathbf{P}%
\setminus\left\{  \mathbf{q}\right\}  }\right)  $ for every $\mathbf{p}%
\in\mathbf{P}\setminus\left\{  \mathbf{q}\right\}  $. Furthermore, assume that
the family $\left(  Q_{\mathbf{p}}\right)  _{\mathbf{p}\in\mathbf{P}%
\setminus\left\{  \mathbf{q}\right\}  }$ is algebraically independent (over
$\mathbb{F}$). Let $\alpha$, $\beta$, $\gamma$ and $\delta$ be four elements
of $\mathbb{F}\left(  x_{\mathbf{P}\setminus\left\{  \mathbf{q}\right\}
}\right)  $ such that $\alpha\delta-\beta\gamma\neq0$. Finally, assume that
$Q_{\mathbf{q}}=\dfrac{\alpha x_{\mathbf{q}}+\beta}{\gamma x_{\mathbf{q}%
}+\delta}$. Then:

\textbf{(a)} We have $\alpha-Q_{\mathbf{q}}\gamma\neq0$ and $x_{\mathbf{q}%
}=\dfrac{Q_{\mathbf{q}}\delta-\beta}{\alpha-Q_{\mathbf{q}}\gamma}$.

\textbf{(b)} The family $\left(  Q_{\mathbf{p}}\right)  _{\mathbf{p}%
\in\mathbf{P}}$ is algebraically independent (over $\mathbb{F}$).
\end{lemma}

\begin{proof}
[Proof of Lemma \ref{lem.algebraic.triangularity.step} (sketched).]For every
subset $\mathbf{S}$ of $\mathbf{P}$, let $\mathbb{F}\left[  x_{\mathbf{S}%
}\right]  $ denote the polynomial ring over $\mathbb{F}$ in the indeterminates
$x_{\mathbf{p}}$ with $\mathbf{p}$ ranging over all elements of $\mathbf{S}$.
We identify $\mathbb{F}\left[  x_{\mathbf{S}}\right]  $ with a subring of
$\mathbb{F}\left[  x_{\mathbf{P}}\right]  $ for every subset $\mathbf{S}$ of
$\mathbf{P}$, and we identify $\mathbb{F}\left[  x_{\mathbf{P}}\right]  $ with
a subring of $\mathbb{F}\left(  x_{\mathbf{P}}\right)  $.

Let $h\in\mathbb{F}\left[  x_{\mathbf{P}}\right]  $ be such that $h\left(
\left(  Q_{\mathbf{p}}\right)  _{\mathbf{p}\in\mathbf{P}}\right)  =0$. We are
now going to prove that $h=0$. (This will yield that $\left(  Q_{\mathbf{p}%
}\right)  _{\mathbf{p}\in\mathbf{P}}$ is algebraically independent.)

Let us assume (for the sake of contradiction) that $h\neq0$.

We know that $\left(  Q_{\mathbf{p}}\right)  _{\mathbf{p}\in\mathbf{P}%
\setminus\left\{  \mathbf{q}\right\}  }\in\left(  \mathbb{F}\left(
x_{\mathbf{P}\setminus\left\{  \mathbf{q}\right\}  }\right)  \right)
^{\mathbf{P}\setminus\left\{  \mathbf{q}\right\}  }$ (since $Q_{\mathbf{p}}%
\in\mathbb{F}\left(  x_{\mathbf{P}\setminus\left\{  \mathbf{q}\right\}
}\right)  $ for every $\mathbf{p}\in\mathbf{P}\setminus\left\{  \mathbf{q}%
\right\}  $).

From $h\left(  \left(  Q_{\mathbf{p}}\right)  _{\mathbf{p}\in\mathbf{P}%
}\right)  =0$ and the fact that the family $\left(  Q_{\mathbf{p}}\right)
_{\mathbf{p}\in\mathbf{P}\setminus\left\{  \mathbf{q}\right\}  }\in\left(
\mathbb{F}\left(  x_{\mathbf{P}\setminus\left\{  \mathbf{q}\right\}  }\right)
\right)  ^{\mathbf{P}\setminus\left\{  \mathbf{q}\right\}  }$ is algebraically
independent, it is easy to see that the element $Q_{\mathbf{q}}\in
\mathbb{F}\left(  x_{\mathbf{P}}\right)  $ is algebraic over the subfield
$\mathbb{F}\left(  x_{\mathbf{P}\setminus\left\{  \mathbf{q}\right\}
}\right)  $\ \ \ \ \footnote{\textit{Proof.} Let $\widetilde{h}$ be the image
of the polynomial $h\in\mathbb{F}\left[  x_{\mathbf{P}}\right]  $ under the
canonical ring isomorphism $\mathbb{F}\left[  x_{\mathbf{P}}\right]
\rightarrow\left(  \mathbb{F}\left[  x_{\mathbf{P}\setminus\left\{
\mathbf{q}\right\}  }\right]  \right)  \left[  x_{\mathbf{q}}\right]  $. Then,
$\widetilde{h}$ is a polynomial in the variable $x_{\mathbf{q}}$ over the ring
$\mathbb{F}\left[  x_{\mathbf{P}\setminus\left\{  \mathbf{q}\right\}
}\right]  $. Let $\widehat{h}$ be the result of replacing every coefficient
$c$ of the polynomial $\widetilde{h}$ by $c\left(  \left(  Q_{\mathbf{p}%
}\right)  _{\mathbf{p}\in\mathbf{P}\setminus\left\{  \mathbf{q}\right\}
}\right)  $. Then, $\widehat{h}$ is a polynomial in the variable
$x_{\mathbf{q}}$ over the field $\mathbb{F}\left(  x_{\mathbf{P}%
\setminus\left\{  \mathbf{q}\right\}  }\right)  $ (because $Q_{\mathbf{p}}%
\in\mathbb{F}\left(  x_{\mathbf{P}\setminus\left\{  \mathbf{q}\right\}
}\right)  $ for every $\mathbf{p}\in\mathbf{P}\setminus\left\{  \mathbf{q}%
\right\}  $). In other words, $\widehat{h}\in\left(  \mathbb{F}\left(
x_{\mathbf{P}\setminus\left\{  \mathbf{q}\right\}  }\right)  \right)  \left[
x_{\mathbf{q}}\right]  $. Moreover, $h\left(  \left(  Q_{\mathbf{p}}\right)
_{\mathbf{p}\in\mathbf{P}}\right)  =\widehat{h}\left(  Q_{\mathbf{q}}\right)
$ (because the polynomial $h\left(  \left(  Q_{\mathbf{p}}\right)
_{\mathbf{p}\in\mathbf{P}}\right)  $ is the result of substituting
$Q_{\mathbf{p}}$ for $x_{\mathbf{p}}$ for every $\mathbf{p}\in\mathbf{P}$ in
the polynomial $h$, whereas $\widehat{h}\left(  Q_{\mathbf{q}}\right)  $ is
the result of first writing $h$ as a univariate polynomial in the variable
$x_{\mathbf{q}}$ over the ring $\mathbb{F}\left[  x_{\mathbf{P}\setminus
\left\{  \mathbf{q}\right\}  }\right]  $, then substituting $Q_{\mathbf{p}}$
for $x_{\mathbf{p}}$ for every $\mathbf{p}\in\mathbf{P}\setminus\left\{
\mathbf{q}\right\}  $ in each coefficient of this univariate polynomial, and
finally substituting $Q_{\mathbf{q}}$ for the variable $x_{\mathbf{q}}$; but
this clearly leads to the same result as just substituting $Q_{\mathbf{p}}$
for $x_{\mathbf{p}}$ for every $\mathbf{p}\in\mathbf{P}$ right away). Thus,
$\widehat{h}\left(  Q_{\mathbf{q}}\right)  =h\left(  \left(  Q_{\mathbf{p}%
}\right)  _{\mathbf{p}\in\mathbf{P}}\right)  =0$.
\par
But $\widetilde{h}$ is the image of $h$ under an isomorphism. Hence,
$\widetilde{h}\neq0$ (since $h=0$). Hence, at least one coefficient $c$ of
$\widetilde{h}$ satisfies $c\neq0$. This coefficient $c$ then satisfies
$c\left(  \left(  Q_{\mathbf{p}}\right)  _{\mathbf{p}\in\mathbf{P}%
\setminus\left\{  \mathbf{q}\right\}  }\right)  \neq0$ (because the family
$\left(  Q_{\mathbf{p}}\right)  _{\mathbf{p}\in\mathbf{P}\setminus\left\{
\mathbf{q}\right\}  }$ is algebraically independent). But $c\left(  \left(
Q_{\mathbf{p}}\right)  _{\mathbf{p}\in\mathbf{P}\setminus\left\{
\mathbf{q}\right\}  }\right)  $ is a coefficient of the polynomial
$\widehat{h}$ (by the construction of $\widehat{h}$, since $c$ is a
coefficient of $\widetilde{h}$). Since $c\left(  \left(  Q_{\mathbf{p}%
}\right)  _{\mathbf{p}\in\mathbf{P}\setminus\left\{  \mathbf{q}\right\}
}\right)  \neq0$, this yields that the polynomial $\widehat{h}$ has a nonzero
coefficient, thus is nonzero.
\par
Now, we know that $\widehat{h}\left(  Q_{\mathbf{q}}\right)  =0$. In other
words, $Q_{\mathbf{q}}$ is a root of the polynomial $\widehat{h}\in\left(
\mathbb{F}\left(  x_{\mathbf{P}\setminus\left\{  \mathbf{q}\right\}  }\right)
\right)  \left[  x_{\mathbf{q}}\right]  $. Since this polynomial $\widehat{h}$
is nonzero, this yields that $Q_{\mathbf{q}}$ is algebraic over the field
$\mathbb{F}\left(  x_{\mathbf{P}\setminus\left\{  \mathbf{q}\right\}
}\right)  $, qed.}.

Since $Q_{\mathbf{q}}=\dfrac{\alpha x_{\mathbf{q}}+\beta}{\gamma
x_{\mathbf{q}}+\delta}$, we have $Q_{\mathbf{q}}\left(  \gamma x_{\mathbf{q}%
}+\delta\right)  =\alpha x_{\mathbf{q}}+\beta$. In other words, $Q_{\mathbf{q}%
}\gamma x_{\mathbf{q}}+Q_{\mathbf{q}}\delta=\alpha x_{\mathbf{q}}+\beta$.
Thus,%
\begin{equation}
Q_{\mathbf{q}}\delta-\beta=\alpha x_{\mathbf{q}}-Q_{\mathbf{q}}\gamma
x_{\mathbf{q}}=\left(  \alpha-Q_{\mathbf{q}}\gamma\right)  x_{\mathbf{q}}.
\label{pf.algebraic.triangularity.step.1}%
\end{equation}
From here (and from $\alpha\delta-\beta\gamma\neq0$), it is easy to see that
$\alpha-Q_{\mathbf{q}}\gamma\neq0$\ \ \ \ \footnote{\textit{Proof.} Assume the
contrary. Thus, $\alpha-Q_{\mathbf{q}}\gamma=0$. In other words,
$\alpha=Q_{\mathbf{q}}\gamma$. But (\ref{pf.algebraic.triangularity.step.1})
yields $Q_{\mathbf{q}}\delta-\beta=\underbrace{\left(  \alpha-Q_{\mathbf{q}%
}\gamma\right)  }_{=0}x_{\mathbf{q}}=0$, hence $\beta=Q_{\mathbf{q}}\delta$.
Now, $\underbrace{\alpha}_{=Q_{\mathbf{q}}\gamma}\delta-\underbrace{\beta
}_{=Q_{\mathbf{q}}\delta}\gamma=Q_{\mathbf{q}}\gamma\delta-Q_{\mathbf{q}%
}\delta\gamma=0$, contradicting $\alpha\delta-\beta\gamma\neq0$. This
contradiction shows that our assumption was wrong, qed.}. Hence, we can divide
(\ref{pf.algebraic.triangularity.step.1}) by $\alpha-Q_{\mathbf{q}}\gamma$,
and obtain%
\[
x_{\mathbf{q}}=\dfrac{Q_{\mathbf{q}}\delta-\beta}{\alpha-Q_{\mathbf{q}}\gamma
}\in\left(  \mathbb{F}\left(  x_{\mathbf{P}\setminus\left\{  \mathbf{q}%
\right\}  }\right)  \right)  \left(  Q_{\mathbf{q}}\right)
\]
(because $\alpha$, $\beta$, $\gamma$, $\delta$ are elements of $\mathbb{F}%
\left(  x_{\mathbf{P}\setminus\left\{  \mathbf{q}\right\}  }\right)  $). Since
$Q_{\mathbf{q}}$ is algebraic over $\mathbb{F}\left(  x_{\mathbf{P}%
\setminus\left\{  \mathbf{q}\right\}  }\right)  $, this yields that
$x_{\mathbf{q}}$ is algebraic over $\mathbb{F}\left(  x_{\mathbf{P}%
\setminus\left\{  \mathbf{q}\right\}  }\right)  $. But this is absurd, since
$\mathbf{q}\notin\mathbf{P}\setminus\left\{  \mathbf{q}\right\}  $ and since
the $x_{\mathbf{p}}$ for all $\mathbf{p}\in\mathbf{P}$ are independent
indeterminates. Hence, we have obtained a contradiction from assuming that
$h\neq0$. The assumption must therefore have been false. Thus, we can't have
$h\neq0$, so that we have $h=0$.

Now, forget that we fixed $h$. We thus have shown that $h=0$ for every
$h\in\mathbb{F}\left[  x_{\mathbf{P}}\right]  $ such that $h\left(  \left(
Q_{\mathbf{p}}\right)  _{\mathbf{p}\in\mathbf{P}}\right)  =0$. In other words,
the family $\left(  Q_{\mathbf{p}}\right)  _{\mathbf{p}\in\mathbf{P}}%
\in\left(  \mathbb{F}\left(  x_{\mathbf{P}}\right)  \right)  ^{\mathbf{P}}$ is
algebraically independent. This proves Lemma
\ref{lem.algebraic.triangularity.step} \textbf{(b)}. We also have shown that
$\alpha-Q_{\mathbf{q}}\gamma\neq0$ and $x_{\mathbf{q}}=\dfrac{Q_{\mathbf{q}%
}\delta-\beta}{\alpha-Q_{\mathbf{q}}\gamma}$, so that Lemma
\ref{lem.algebraic.triangularity.step} \textbf{(a)} is proven as well.
\end{proof}

The next lemma provides a way to check $\mathbf{P}$-triangularity of a family
using a subfamily:

\begin{lemma}
\label{lem.algebraic.triangularity.ind}Let $\mathbb{F}$ be a field. Let
$\mathbf{P}$ be a nonempty finite totally ordered set, and let
$\vartriangleleft$ be the smaller relation of $\mathbf{P}$. For every
$\mathbf{p}\in\mathbf{P}$, let $\mathbf{p}\Downarrow$ denote the subset
$\left\{  \mathbf{v}\in\mathbf{P}\ \mid\ \mathbf{v}\vartriangleleft
\mathbf{p}\right\}  $ of $\mathbf{P}$. Let $\mathbf{q}$ be the maximal element
of $\mathbf{P}$.

Let $\left(  Q_{\mathbf{p}}\right)  _{\mathbf{p}\in\mathbf{P}}\in\left(
\mathbb{F}\left(  x_{\mathbf{P}}\right)  \right)  ^{\mathbf{P}}$ be a family
of elements of $\mathbb{F}\left(  x_{\mathbf{P}}\right)  $ indexed by elements
of $\mathbf{P}$. Assume that every $\mathbf{p}\in\mathbf{P}\setminus\left\{
\mathbf{q}\right\}  $ satisfies $Q_{\mathbf{p}}\in\mathbb{F}\left(
x_{\mathbf{P}\setminus\left\{  \mathbf{q}\right\}  }\right)  $. Moreover,
assume that $\left(  Q_{\mathbf{p}}\right)  _{\mathbf{p}\in\mathbf{P}%
\setminus\left\{  \mathbf{q}\right\}  }\in\left(  \mathbb{F}\left(
x_{\mathbf{P}\setminus\left\{  \mathbf{q}\right\}  }\right)  \right)
^{\mathbf{P}\setminus\left\{  \mathbf{q}\right\}  }$ is a $\mathbf{P}%
\setminus\left\{  \mathbf{q}\right\}  $-triangular family.

Let $\alpha$, $\beta$, $\gamma$, $\delta$ be elements of $\mathbb{F}\left(
x_{\mathbf{P}\setminus\left\{  \mathbf{q}\right\}  }\right)  $ such that
$\alpha\delta-\beta\gamma\neq0$ and $Q_{\mathbf{q}}=\dfrac{\alpha
x_{\mathbf{q}}+\beta}{\gamma x_{\mathbf{q}}+\delta}$. Then, the family
$\left(  Q_{\mathbf{p}}\right)  _{\mathbf{p}\in\mathbf{P}}\in\left(
\mathbb{F}\left(  x_{\mathbf{P}}\right)  \right)  ^{\mathbf{P}}$ is
$\mathbf{P}$-triangular.
\end{lemma}

\begin{proof}
[Proof of Lemma \ref{lem.algebraic.triangularity.ind} (sketched).]Clearly,
$\mathbf{P}\setminus\left\{  \mathbf{q}\right\}  $ is a subposet of
$\mathbf{P}$. If $\mathbf{p}$ is an element of $\mathbf{P}\setminus\left\{
\mathbf{q}\right\}  $, then the subset $\left\{  \mathbf{v}\in\mathbf{P}%
\ \mid\ \mathbf{v}\vartriangleleft\mathbf{p}\right\}  $ of $\mathbf{P}$ equals
the subset $\left\{  \mathbf{v}\in\mathbf{P}\setminus\left\{  \mathbf{q}%
\right\}  \ \mid\ \mathbf{v}\vartriangleleft\mathbf{p}\right\}  $ of
$\mathbf{P}\setminus\left\{  \mathbf{q}\right\}  $ (since $\mathbf{q}$ is
maximal and therefore not less than $\mathbf{p}$). In other words, if
$\mathbf{p}$ is an element of $\mathbf{P}\setminus\left\{  \mathbf{q}\right\}
$, then the set $\mathbf{p}\Downarrow$ does not depend on whether $\mathbf{p}$
is regarded as an element of $\mathbf{P}$ or as an element of $\mathbf{P}%
\setminus\left\{  \mathbf{q}\right\}  $.

We know that the family $\left(  Q_{\mathbf{p}}\right)  _{\mathbf{p}%
\in\mathbf{P}\setminus\left\{  \mathbf{q}\right\}  }$ is $\mathbf{P}%
\setminus\left\{  \mathbf{q}\right\}  $-triangular. By the definition of
``$\mathbf{P}\setminus\left\{  \mathbf{q}\right\}  $-triangular'', this means
that the algebraic triangularity condition\footnote{This condition was
described in Definition \ref{def.algebraic.triangularity.triang}.} with
$\mathbf{P}$ replaced by $\mathbf{P}\setminus\left\{  \mathbf{q}\right\}  $
holds. In other words,
\begin{equation}
\left(
\begin{array}
[c]{c}%
\text{for every }\mathbf{p}\in\mathbf{P}\setminus\left\{  \mathbf{q}\right\}
\text{, there exist elements }\alpha_{\mathbf{p}}\text{, }\beta_{\mathbf{p}%
}\text{, }\gamma_{\mathbf{p}}\text{, }\delta_{\mathbf{p}}\text{ of }%
\mathbb{F}\left(  x_{\mathbf{p}\Downarrow}\right) \\
\text{such that }\alpha_{\mathbf{p}}\delta_{\mathbf{p}}-\beta_{\mathbf{p}%
}\gamma_{\mathbf{p}}\neq0\text{ and }Q_{\mathbf{p}}=\dfrac{\alpha_{\mathbf{p}%
}x_{\mathbf{p}}+\beta_{\mathbf{p}}}{\gamma_{\mathbf{p}}x_{\mathbf{p}}%
+\delta_{\mathbf{p}}}%
\end{array}
\right)  . \label{pf.algebraic.triangularity.ind.1}%
\end{equation}
(Note that the notation $\mathbf{p}\Downarrow$ is unambiguous here, because
(as we said) the set $\mathbf{p}\Downarrow$ does not depend on whether
$\mathbf{p}$ is regarded as an element of $\mathbf{P}$ or as an element of
$\mathbf{P}\setminus\left\{  \mathbf{q}\right\}  $.)

Now, we are going to prove that the algebraic triangularity condition (without
any replacements) holds (for the family $\left(  Q_{\mathbf{p}}\right)
_{\mathbf{p}\in\mathbf{P}}$). In order to do so, we need to prove that for
every $\mathbf{p}\in\mathbf{P}$, there exist elements $\alpha_{\mathbf{p}}$,
$\beta_{\mathbf{p}}$, $\gamma_{\mathbf{p}}$, $\delta_{\mathbf{p}}$ of
$\mathbb{F}\left(  x_{\mathbf{p}\Downarrow}\right)  $ such that $\alpha
_{\mathbf{p}}\delta_{\mathbf{p}}-\beta_{\mathbf{p}}\gamma_{\mathbf{p}}\neq0$
and $Q_{\mathbf{p}}=\dfrac{\alpha_{\mathbf{p}}x_{\mathbf{p}}+\beta
_{\mathbf{p}}}{\gamma_{\mathbf{p}}x_{\mathbf{p}}+\delta_{\mathbf{p}}}$. Since
we already know that this holds for each $\mathbf{p}\in\mathbf{P}%
\setminus\left\{  \mathbf{q}\right\}  $ (because of
(\ref{pf.algebraic.triangularity.ind.1})), we only need to prove it for those
$\mathbf{p}\in\mathbf{P}$ which don't satisfy $\mathbf{p}\in\mathbf{P}%
\setminus\left\{  \mathbf{q}\right\}  $. In other words, we need to prove that
for every $\mathbf{p}\in\mathbf{P}$ which doesn't satisfy $\mathbf{p}%
\in\mathbf{P}\setminus\left\{  \mathbf{q}\right\}  $, there exist elements
$\alpha_{\mathbf{p}}$, $\beta_{\mathbf{p}}$, $\gamma_{\mathbf{p}}$,
$\delta_{\mathbf{p}}$ of $\mathbb{F}\left(  x_{\mathbf{p}\Downarrow}\right)  $
such that $\alpha_{\mathbf{p}}\delta_{\mathbf{p}}-\beta_{\mathbf{p}}%
\gamma_{\mathbf{p}}\neq0$ and $Q_{\mathbf{p}}=\dfrac{\alpha_{\mathbf{p}%
}x_{\mathbf{p}}+\beta_{\mathbf{p}}}{\gamma_{\mathbf{p}}x_{\mathbf{p}}%
+\delta_{\mathbf{p}}}$. Since the only $\mathbf{p}\in\mathbf{P}$ which doesn't
satisfy $\mathbf{p}\in\mathbf{P}\setminus\left\{  \mathbf{q}\right\}  $ is
$\mathbf{q}$, this rewrites as follows: We need to prove that for
$\mathbf{p}=\mathbf{q}$, there exist elements $\alpha_{\mathbf{p}}$,
$\beta_{\mathbf{p}}$, $\gamma_{\mathbf{p}}$, $\delta_{\mathbf{p}}$ of
$\mathbb{F}\left(  x_{\mathbf{p}\Downarrow}\right)  $ such that $\alpha
_{\mathbf{p}}\delta_{\mathbf{p}}-\beta_{\mathbf{p}}\gamma_{\mathbf{p}}\neq0$
and $Q_{\mathbf{p}}=\dfrac{\alpha_{\mathbf{p}}x_{\mathbf{p}}+\beta
_{\mathbf{p}}}{\gamma_{\mathbf{p}}x_{\mathbf{p}}+\delta_{\mathbf{p}}}$.

So let $\mathbf{p}=\mathbf{q}$. Then, $Q_{\mathbf{p}}=Q_{\mathbf{q}}$ and
$\left.  \mathbf{p}\Downarrow\right.  =\left.  \mathbf{q}\Downarrow\right.
=\left\{  \mathbf{v}\in\mathbf{P}\ \mid\ \mathbf{v}\vartriangleleft
\mathbf{q}\right\}  =\mathbf{P}\setminus\left\{  \mathbf{q}\right\}  $ (since
$\mathbf{q}$ is the maximal element of the totally ordered set $\mathbf{P}$),
so that $\mathbb{F}\left(  x_{\mathbf{p}\Downarrow}\right)  =\mathbb{F}\left(
x_{\mathbf{P}\setminus\left\{  \mathbf{q}\right\}  }\right)  $. But we know
that $\alpha$, $\beta$, $\gamma$, $\delta$ belong to $\mathbb{F}\left(
x_{\mathbf{P}\setminus\left\{  \mathbf{q}\right\}  }\right)  =\mathbb{F}%
\left(  x_{\mathbf{p}\Downarrow}\right)  $ and satisfy $\alpha\delta
-\beta\gamma\neq0$ and $Q_{\mathbf{p}}=Q_{\mathbf{q}}=\dfrac{\alpha
x_{\mathbf{q}}+\beta}{\gamma x_{\mathbf{q}}+\delta}=\dfrac{\alpha
x_{\mathbf{p}}+\beta}{\gamma x_{\mathbf{p}}+\delta}$ (since $\mathbf{q}%
=\mathbf{p}$). Hence, there exist elements $\alpha_{\mathbf{p}}$,
$\beta_{\mathbf{p}}$, $\gamma_{\mathbf{p}}$, $\delta_{\mathbf{p}}$ of
$\mathbb{F}\left(  x_{\mathbf{p}\Downarrow}\right)  $ such that $\alpha
_{\mathbf{p}}\delta_{\mathbf{p}}-\beta_{\mathbf{p}}\gamma_{\mathbf{p}}\neq0$
and $Q_{\mathbf{p}}=\dfrac{\alpha_{\mathbf{p}}x_{\mathbf{p}}+\beta
_{\mathbf{p}}}{\gamma_{\mathbf{p}}x_{\mathbf{p}}+\delta_{\mathbf{p}}}$
(namely, $\alpha$, $\beta$, $\gamma$, $\delta$). Hence, we have proven that
the algebraic triangularity condition (without any replacements) holds (for
the family $\left(  Q_{\mathbf{p}}\right)  _{\mathbf{p}\in\mathbf{P}}$). In
other words, the family $\left(  Q_{\mathbf{p}}\right)  _{\mathbf{p}%
\in\mathbf{P}}\in\left(  \mathbb{F}\left(  x_{\mathbf{P}}\right)  \right)
^{\mathbf{P}}$ is $\mathbf{P}$-triangular (by the definition of ``$\mathbf{P}%
$-triangular''). This proves Lemma \ref{lem.algebraic.triangularity.ind}.
\end{proof}

\begin{proof}
[Proof of Lemma \ref{lem.algebraic.triangularity} (sketched).]\textbf{(a)} We
will prove Lemma \ref{lem.algebraic.triangularity} \textbf{(a)} by induction
over $\left\vert \mathbf{P}\right\vert $:

\textit{Induction base:} If $\left\vert \mathbf{P}\right\vert =0$, then
$\mathbf{P}=\varnothing$. Hence, if $\left\vert \mathbf{P}\right\vert =0$,
then Lemma \ref{lem.algebraic.triangularity} \textbf{(a)} is obviously true
(because any empty family in a field is algebraically independent). Hence, the
induction base is complete.

\textit{Induction step:} Let $N\in\mathbb{N}$. Assume that Lemma
\ref{lem.algebraic.triangularity} \textbf{(a)} is proven in the case when
$\left\vert \mathbf{P}\right\vert =N$. We need to show that Lemma
\ref{lem.algebraic.triangularity} \textbf{(a)} holds in the case when
$\left\vert \mathbf{P}\right\vert =N+1$.

So assume that $\mathbf{P}$ is a finite totally ordered set satisfying
$\left\vert \mathbf{P}\right\vert =N+1$. Let $Q_{\mathbf{p}}$ be as in Lemma
\ref{lem.algebraic.triangularity}. Define the notations $\vartriangleleft$ and
$\mathbf{p}\Downarrow$ (for $\mathbf{p}\in\mathbf{P}$) as in Definition
\ref{def.algebraic.triangularity.triang}.

We know that the family $\left(  Q_{\mathbf{p}}\right)  _{\mathbf{p}%
\in\mathbf{P}}$ is $\mathbf{P}$-triangular, thus satisfies the algebraic
triangularity condition\footnote{This condition was described in Definition
\ref{def.algebraic.triangularity.triang}.}.

Let $\mathbf{q}$ be the maximal element of $\mathbf{P}$ (this $\mathbf{q}$
exists, since $\left\vert \mathbf{P}\right\vert =N+1\geq1$). Then, every
$\mathbf{p}\in\mathbf{P}\setminus\left\{  \mathbf{q}\right\}  $ satisfies
$Q_{\mathbf{p}}\in\mathbb{F}\left(  x_{\mathbf{P}\setminus\left\{
\mathbf{q}\right\}  }\right)  $.\ \ \ \ \footnote{\textit{Proof.} Let
$\mathbf{p}\in\mathbf{P}\setminus\left\{  \mathbf{q}\right\}  $. By the
algebraic triangularity condition, there exist elements $\alpha_{\mathbf{p}}$,
$\beta_{\mathbf{p}}$, $\gamma_{\mathbf{p}}$, $\delta_{\mathbf{p}}$ of
$\mathbb{F}\left(  x_{\mathbf{p}\Downarrow}\right)  $ such that $\alpha
_{\mathbf{p}}\delta_{\mathbf{p}}-\beta_{\mathbf{p}}\gamma_{\mathbf{p}}\neq0$
and $Q_{\mathbf{p}}=\dfrac{\alpha_{\mathbf{p}}x_{\mathbf{p}}+\beta
_{\mathbf{p}}}{\gamma_{\mathbf{p}}x_{\mathbf{p}}+\delta_{\mathbf{p}}}$.
Consider these elements $\alpha_{\mathbf{p}}$, $\beta_{\mathbf{p}}$,
$\gamma_{\mathbf{p}}$, $\delta_{\mathbf{p}}$. Since $\mathbf{q}$ is a maximal
element of $\mathbf{P}$, we have $\left.  \mathbf{p}\Downarrow\right.
\subseteq\mathbf{P}\setminus\left\{  \mathbf{q}\right\}  $, so that
$\mathbb{F}\left(  x_{\mathbf{p}\Downarrow}\right)  \subseteq\mathbb{F}\left(
x_{\mathbf{P}\setminus\left\{  \mathbf{q}\right\}  }\right)  $. Thus,
$\alpha_{\mathbf{p}}$, $\beta_{\mathbf{p}}$, $\gamma_{\mathbf{p}}$,
$\delta_{\mathbf{p}}$ are elements of $\mathbb{F}\left(  x_{\mathbf{P}%
\setminus\left\{  \mathbf{q}\right\}  }\right)  $ (since they are elements of
$\mathbb{F}\left(  x_{\mathbf{p}\Downarrow}\right)  $). Since $x_{\mathbf{p}}$
is also an element of $\mathbb{F}\left(  x_{\mathbf{P}\setminus\left\{
\mathbf{q}\right\}  }\right)  $ (because $\mathbf{p}\in\mathbf{P}%
\setminus\left\{  \mathbf{q}\right\}  $), this yields that $\dfrac
{\alpha_{\mathbf{p}}x_{\mathbf{p}}+\beta_{\mathbf{p}}}{\gamma_{\mathbf{p}%
}x_{\mathbf{p}}+\delta_{\mathbf{p}}}$ is an element of $\mathbb{F}\left(
x_{\mathbf{P}\setminus\left\{  \mathbf{q}\right\}  }\right)  $. In other
words, $Q_{\mathbf{p}}$ is an element of $\mathbb{F}\left(  x_{\mathbf{P}%
\setminus\left\{  \mathbf{q}\right\}  }\right)  $ (since $Q_{\mathbf{p}%
}=\dfrac{\alpha_{\mathbf{p}}x_{\mathbf{p}}+\beta_{\mathbf{p}}}{\gamma
_{\mathbf{p}}x_{\mathbf{p}}+\delta_{\mathbf{p}}}$), qed.}

Since $\mathbf{q}\in\mathbf{P}$, it is clear that $\mathbf{P}\setminus\left\{
\mathbf{q}\right\}  $ is a finite totally ordered set satisfying $\left\vert
\mathbf{P}\setminus\left\{  \mathbf{q}\right\}  \right\vert
=\underbrace{\left\vert \mathbf{P}\right\vert }_{=N+1}-1=N+1-1=N$. If
$\mathbf{p}$ is an element of $\mathbf{P}\setminus\left\{  \mathbf{q}\right\}
$, then the subset $\left\{  \mathbf{v}\in\mathbf{P}\ \mid\ \mathbf{v}%
\vartriangleleft\mathbf{p}\right\}  $ of $\mathbf{P}$ equals the subset
$\left\{  \mathbf{v}\in\mathbf{P}\setminus\left\{  \mathbf{q}\right\}
\ \mid\ \mathbf{v}\vartriangleleft\mathbf{p}\right\}  $ of $\mathbf{P}%
\setminus\left\{  \mathbf{q}\right\}  $ (since $\mathbf{q}$ is maximal and
therefore not less than $\mathbf{p}$). In other words, if $\mathbf{p}$ is an
element of $\mathbf{P}\setminus\left\{  \mathbf{q}\right\}  $, then the set
$\mathbf{p}\Downarrow$ does not depend on whether $\mathbf{p}$ is regarded as
an element of $\mathbf{P}$ or as an element of $\mathbf{P}\setminus\left\{
\mathbf{q}\right\}  $. Consequently, the algebraic triangularity condition
holds for $\mathbf{P}\setminus\left\{  \mathbf{q}\right\}  $ instead of
$\mathbf{P}$ (but with the same elements $Q_{\mathbf{p}}$) (because it holds
for $\mathbf{P}$). In other words, the family $\left(  Q_{\mathbf{p}}\right)
_{\mathbf{p}\in\mathbf{P}\setminus\left\{  \mathbf{q}\right\}  }$ is
$\mathbf{P}\setminus\left\{  \mathbf{q}\right\}  $-triangular.

We thus know the following facts:

\begin{itemize}
\item The finite totally ordered set $\mathbf{P}\setminus\left\{
\mathbf{q}\right\}  $ satisfies $\left\vert \mathbf{P}\setminus\left\{
\mathbf{q}\right\}  \right\vert =N$.

\item For every $\mathbf{p}\in\mathbf{P}\setminus\left\{  \mathbf{q}\right\}
$, the rational function $Q_{\mathbf{p}}$ is an element of $\mathbb{F}\left(
x_{\mathbf{P}\setminus\left\{  \mathbf{q}\right\}  }\right)  $.

\item The family $\left(  Q_{\mathbf{p}}\right)  _{\mathbf{p}\in
\mathbf{P}\setminus\left\{  \mathbf{q}\right\}  }$ is $\mathbf{P}%
\setminus\left\{  \mathbf{q}\right\}  $-triangular.
\end{itemize}

Due to these facts, we can apply Lemma \ref{lem.algebraic.triangularity}
\textbf{(a)} to $\mathbf{P}\setminus\left\{  \mathbf{q}\right\}  $ instead of
$\mathbf{P}$ (because we assumed that Lemma \ref{lem.algebraic.triangularity}
\textbf{(a)} is proven in the case when $\left\vert \mathbf{P}\right\vert
=N$). As a consequence, we obtain that the family $\left(  Q_{\mathbf{p}%
}\right)  _{\mathbf{p}\in\mathbf{P}\setminus\left\{  \mathbf{q}\right\}  }%
\in\left(  \mathbb{F}\left(  x_{\mathbf{P}\setminus\left\{  \mathbf{q}%
\right\}  }\right)  \right)  ^{\mathbf{P}\setminus\left\{  \mathbf{q}\right\}
}$ is algebraically independent.

But recall that $\left(  Q_{\mathbf{p}}\right)  _{\mathbf{p}\in\mathbf{P}}$ is
a $\mathbf{P}$-triangular family. Hence, it satisfies the algebraic
triangularity condition. Applying this algebraic triangularity condition to
$\mathbf{p}=\mathbf{q}$, we see that there exist elements $\alpha_{\mathbf{q}%
}$, $\beta_{\mathbf{q}}$, $\gamma_{\mathbf{q}}$, $\delta_{\mathbf{q}}$ of
$\mathbb{F}\left(  x_{\mathbf{q}\Downarrow}\right)  $ such that $\alpha
_{\mathbf{q}}\delta_{\mathbf{q}}-\beta_{\mathbf{q}}\gamma_{\mathbf{q}}\neq0$
and $Q_{\mathbf{q}}=\dfrac{\alpha_{\mathbf{q}}x_{\mathbf{q}}+\beta
_{\mathbf{q}}}{\gamma_{\mathbf{q}}x_{\mathbf{q}}+\delta_{\mathbf{q}}}$.
Consider these elements $\alpha_{\mathbf{q}}$, $\beta_{\mathbf{q}}$,
$\gamma_{\mathbf{q}}$, $\delta_{\mathbf{q}}$.

Since $\left.  \mathbf{q}\Downarrow\right.  \subseteq\mathbf{P}\setminus
\left\{  \mathbf{q}\right\}  $, we have $\mathbb{F}\left(  x_{\mathbf{q}%
\Downarrow}\right)  \subseteq\mathbb{F}\left(  x_{\mathbf{P}\setminus\left\{
\mathbf{q}\right\}  }\right)  $. Thus, $\alpha_{\mathbf{q}}$, $\beta
_{\mathbf{q}}$, $\gamma_{\mathbf{q}}$, $\delta_{\mathbf{q}}$ are elements of
$\mathbb{F}\left(  x_{\mathbf{P}\setminus\left\{  \mathbf{q}\right\}
}\right)  $ (since $\alpha_{\mathbf{q}}$, $\beta_{\mathbf{q}}$, $\gamma
_{\mathbf{q}}$, $\delta_{\mathbf{q}}$ are elements of $\mathbb{F}\left(
x_{\mathbf{q}\Downarrow}\right)  $).

We can now apply Lemma \ref{lem.algebraic.triangularity.step} \textbf{(b)} to
$\alpha_{\mathbf{q}}$, $\beta_{\mathbf{q}}$, $\gamma_{\mathbf{q}}$,
$\delta_{\mathbf{q}}$ instead of $\alpha$, $\beta$, $\gamma$, $\delta$. We
conclude that the family $\left(  Q_{\mathbf{p}}\right)  _{\mathbf{p}%
\in\mathbf{P}}$ is algebraically independent. Hence, Lemma
\ref{lem.algebraic.triangularity} \textbf{(a)} holds for our totally ordered
set $\mathbf{P}$. We thus have shown that Lemma
\ref{lem.algebraic.triangularity} \textbf{(a)} holds in the case when
$\left\vert \mathbf{P}\right\vert =N+1$. This finishes the induction step, and
so Lemma \ref{lem.algebraic.triangularity} \textbf{(a)} is proven by induction.

\textbf{(b)} We will prove Lemma \ref{lem.algebraic.triangularity}
\textbf{(b)} by induction over $\left\vert \mathbf{P}\right\vert $:

\textit{Induction base:} If $\mathbf{P}=\varnothing$, then Lemma
\ref{lem.algebraic.triangularity} \textbf{(b)} is obviously true (because the
empty family can be taken as $\left(  R_{\mathbf{p}}\right)  _{\mathbf{p}%
\in\mathbf{P}}$). Hence, the induction base is complete.

\textit{Induction step:} Let $N\in\mathbb{N}$. Assume that Lemma
\ref{lem.algebraic.triangularity} \textbf{(b)} is proven in the case when
$\left\vert \mathbf{P}\right\vert =N$. We need to show that Lemma
\ref{lem.algebraic.triangularity} \textbf{(b)} holds in the case when
$\left\vert \mathbf{P}\right\vert =N+1$.

So assume that $\mathbf{P}$ is a finite totally ordered set satisfying
$\left\vert \mathbf{P}\right\vert =N+1$. Let $Q_{\mathbf{p}}$ be as in Lemma
\ref{lem.algebraic.triangularity}. Define the notations $\vartriangleleft$ and
$\mathbf{p}\Downarrow$ (for $\mathbf{p}\in\mathbf{P}$) as in Definition
\ref{def.algebraic.triangularity.triang}.

We know that the family $\left(  Q_{\mathbf{p}}\right)  _{\mathbf{p}%
\in\mathbf{P}}$ is $\mathbf{P}$-triangular. Thus, it satisfies the algebraic
triangularity condition\footnote{This condition was described in Definition
\ref{def.algebraic.triangularity.triang}.}.

Let $\mathbf{q}$ be the maximal element of $\mathbf{P}$ (this $\mathbf{q}$
exists, since $\left\vert \mathbf{P}\right\vert =N+1\geq1$). Then, every
$\mathbf{p}\in\mathbf{P}\setminus\left\{  \mathbf{q}\right\}  $ satisfies
$Q_{\mathbf{p}}\in\mathbb{F}\left(  x_{\mathbf{P}\setminus\left\{
\mathbf{q}\right\}  }\right)  $. (This is proven in the same way as in the
above proof of Lemma \ref{lem.algebraic.triangularity} \textbf{(a)}.)

Since $\mathbf{q}\in\mathbf{P}$, it is clear that $\mathbf{P}\setminus\left\{
\mathbf{q}\right\}  $ is a finite totally ordered set satisfying $\left\vert
\mathbf{P}\setminus\left\{  \mathbf{q}\right\}  \right\vert
=\underbrace{\left\vert \mathbf{P}\right\vert }_{=N+1}-1=N+1-1=N$. If
$\mathbf{p}$ is an element of $\mathbf{P}\setminus\left\{  \mathbf{q}\right\}
$, then the subset $\left\{  \mathbf{v}\in\mathbf{P}\ \mid\ \mathbf{v}%
\vartriangleleft\mathbf{p}\right\}  $ of $\mathbf{P}$ equals the subset
$\left\{  \mathbf{v}\in\mathbf{P}\setminus\left\{  \mathbf{q}\right\}
\ \mid\ \mathbf{v}\vartriangleleft\mathbf{p}\right\}  $ of $\mathbf{P}%
\setminus\left\{  \mathbf{q}\right\}  $ (since $\mathbf{q}$ is maximal and
therefore not less than $\mathbf{p}$). In other words, if $\mathbf{p}$ is an
element of $\mathbf{P}\setminus\left\{  \mathbf{q}\right\}  $, then the set
$\mathbf{p}\Downarrow$ does not depend on whether $\mathbf{p}$ is regarded as
an element of $\mathbf{P}$ or as an element of $\mathbf{P}\setminus\left\{
\mathbf{q}\right\}  $. Consequently, the algebraic triangularity condition
holds for $\mathbf{P}\setminus\left\{  \mathbf{q}\right\}  $ instead of
$\mathbf{P}$ (but with the same elements $Q_{\mathbf{p}}$) (because it holds
for $\mathbf{P}$). In other words, the family $\left(  Q_{\mathbf{p}}\right)
_{\mathbf{p}\in\mathbf{P}\setminus\left\{  \mathbf{q}\right\}  }$ is
$\mathbf{P}\setminus\left\{  \mathbf{q}\right\}  $-triangular. Thus, the
family $\left(  Q_{\mathbf{p}}\right)  _{\mathbf{p}\in\mathbf{P}%
\setminus\left\{  \mathbf{q}\right\}  }$ is algebraically independent (by
Lemma \ref{lem.algebraic.triangularity} \textbf{(a)}, applied to
$\mathbf{P}\setminus\left\{  \mathbf{q}\right\}  $ instead of $\mathbf{P}$).

Altogether, we thus know the following facts:

\begin{itemize}
\item The finite totally ordered set $\mathbf{P}\setminus\left\{
\mathbf{q}\right\}  $ satisfies $\left\vert \mathbf{P}\setminus\left\{
\mathbf{q}\right\}  \right\vert =N$.

\item For every $\mathbf{p}\in\mathbf{P}\setminus\left\{  \mathbf{q}\right\}
$, the rational function $Q_{\mathbf{p}}$ is an element of $\mathbb{F}\left(
x_{\mathbf{P}\setminus\left\{  \mathbf{q}\right\}  }\right)  $.

\item The family $\left(  Q_{\mathbf{p}}\right)  _{\mathbf{p}\in
\mathbf{P}\setminus\left\{  \mathbf{q}\right\}  }$ is $\mathbf{P}%
\setminus\left\{  \mathbf{q}\right\}  $-triangular.
\end{itemize}

Due to these facts, we can apply Lemma \ref{lem.algebraic.triangularity}
\textbf{(b)} to $\mathbf{P}\setminus\left\{  \mathbf{q}\right\}  $ instead of
$\mathbf{P}$ (because we assumed that Lemma \ref{lem.algebraic.triangularity}
\textbf{(b)} is proven in the case when $\left\vert \mathbf{P}\right\vert
=N$). As a consequence, we obtain that there exists a $\mathbf{P}%
\setminus\left\{  \mathbf{q}\right\}  $-triangular family $\left(
R_{\mathbf{p}}\right)  _{\mathbf{p}\in\mathbf{P}\setminus\left\{
\mathbf{q}\right\}  }\in\left(  \mathbb{F}\left(  x_{\mathbf{P}\setminus
\left\{  \mathbf{q}\right\}  }\right)  \right)  ^{\mathbf{P}\setminus\left\{
\mathbf{q}\right\}  }$ such that the maps $\operatorname*{rev}%
\nolimits_{\left(  Q_{\mathbf{p}}\right)  _{\mathbf{p}\in\mathbf{P}%
\setminus\left\{  \mathbf{q}\right\}  }}$ and $\operatorname*{rev}%
\nolimits_{\left(  R_{\mathbf{p}}\right)  _{\mathbf{p}\in\mathbf{P}%
\setminus\left\{  \mathbf{q}\right\}  }}$ are mutually inverse. Consider this
$\mathbf{P}\setminus\left\{  \mathbf{q}\right\}  $-triangular family $\left(
R_{\mathbf{p}}\right)  _{\mathbf{p}\in\mathbf{P}\setminus\left\{
\mathbf{q}\right\}  }\in\left(  \mathbb{F}\left(  x_{\mathbf{P}\setminus
\left\{  \mathbf{q}\right\}  }\right)  \right)  ^{\mathbf{P}\setminus\left\{
\mathbf{q}\right\}  }$.

We are going to extend this family to a family $\left(  R_{\mathbf{p}}\right)
_{\mathbf{p}\in\mathbf{P}}\in\left(  \mathbb{F}\left(  x_{\mathbf{P}}\right)
\right)  ^{\mathbf{P}}$ by defining a new element $R_{\mathbf{q}}\in
\mathbb{F}\left(  x_{\mathbf{P}}\right)  $.

But before we define $R_{\mathbf{q}}$, let us make some more observations.

First of all, every $\mathbf{p}\in\mathbf{P}\setminus\left\{  \mathbf{q}%
\right\}  $ satisfies $R_{\mathbf{p}}\in\mathbb{F}\left(  x_{\mathbf{P}%
\setminus\left\{  \mathbf{q}\right\}  }\right)  $. (This is simply because
$\left(  R_{\mathbf{p}}\right)  _{\mathbf{p}\in\mathbf{P}\setminus\left\{
\mathbf{q}\right\}  }\in\left(  \mathbb{F}\left(  x_{\mathbf{P}\setminus
\left\{  \mathbf{q}\right\}  }\right)  \right)  ^{\mathbf{P}\setminus\left\{
\mathbf{q}\right\}  }$.)

Recall that $\left(  Q_{\mathbf{p}}\right)  _{\mathbf{p}\in\mathbf{P}}$ is a
$\mathbf{P}$-triangular family. Hence, it satisfies the algebraic
triangularity condition (described in Definition
\ref{def.algebraic.triangularity.triang}). Applying this algebraic
triangularity condition to $\mathbf{p}=\mathbf{q}$, we see that there exist
elements $\alpha_{\mathbf{q}}$, $\beta_{\mathbf{q}}$, $\gamma_{\mathbf{q}}$,
$\delta_{\mathbf{q}}$ of $\mathbb{F}\left(  x_{\mathbf{q}\Downarrow}\right)  $
such that $\alpha_{\mathbf{q}}\delta_{\mathbf{q}}-\beta_{\mathbf{q}}%
\gamma_{\mathbf{q}}\neq0$ and $Q_{\mathbf{q}}=\dfrac{\alpha_{\mathbf{q}%
}x_{\mathbf{q}}+\beta_{\mathbf{q}}}{\gamma_{\mathbf{q}}x_{\mathbf{q}}%
+\delta_{\mathbf{q}}}$. Denote these elements $\alpha_{\mathbf{q}}$,
$\beta_{\mathbf{q}}$, $\gamma_{\mathbf{q}}$, $\delta_{\mathbf{q}}$ by $\alpha
$, $\beta$, $\gamma$, $\delta$. Thus, $\alpha$, $\beta$, $\gamma$, $\delta$
are elements of $\mathbb{F}\left(  x_{\mathbf{q}\Downarrow}\right)  $ such
that
\[
\alpha\delta-\beta\gamma\neq0\ \ \ \ \ \ \ \ \ \ \text{and}%
\ \ \ \ \ \ \ \ \ \ Q_{\mathbf{q}}=\dfrac{\alpha x_{\mathbf{q}}+\beta}{\gamma
x_{\mathbf{q}}+\delta}.
\]

Since $\left.  \mathbf{q}\Downarrow\right.  \subseteq\mathbf{P}\setminus
\left\{  \mathbf{q}\right\}  $, we have $\mathbb{F}\left(  x_{\mathbf{q}%
\Downarrow}\right)  \subseteq\mathbb{F}\left(  x_{\mathbf{P}\setminus\left\{
\mathbf{q}\right\}  }\right)  $. Thus, $\alpha$, $\beta$, $\gamma$, $\delta$
are elements of $\mathbb{F}\left(  x_{\mathbf{P}\setminus\left\{
\mathbf{q}\right\}  }\right)  $ (since $\alpha$, $\beta$, $\gamma$, $\delta$
are elements of $\mathbb{F}\left(  x_{\mathbf{q}\Downarrow}\right)  $).

Since the family $\left(  R_{\mathbf{p}}\right)  _{\mathbf{p}\in
\mathbf{P}\setminus\left\{  \mathbf{q}\right\}  }\in\left(  \mathbb{F}\left(
x_{\mathbf{P}\setminus\left\{  \mathbf{q}\right\}  }\right)  \right)
^{\mathbf{P}\setminus\left\{  \mathbf{q}\right\}  }$ is $\mathbf{P}%
\setminus\left\{  \mathbf{q}\right\}  $-triangular, we can apply Lemma
\ref{lem.algebraic.triangularity} \textbf{(a)} to $\mathbf{P}\setminus\left\{
\mathbf{q}\right\}  $ and $R_{\mathbf{p}}$ instead of $\mathbf{P}$ and
$Q_{\mathbf{p}}$. We thus conclude that the family $\left(  R_{\mathbf{p}%
}\right)  _{\mathbf{p}\in\mathbf{P}\setminus\left\{  \mathbf{q}\right\}  }%
\in\left(  \mathbb{F}\left(  x_{\mathbf{P}\setminus\left\{  \mathbf{q}%
\right\}  }\right)  \right)  ^{\mathbf{P}\setminus\left\{  \mathbf{q}\right\}
}$ is algebraically independent. Hence, the $\mathbb{F}$-algebra homomorphism
$\operatorname*{rev}\nolimits_{\left(  R_{\mathbf{p}}\right)  _{\mathbf{p}%
\in\mathbf{P}\setminus\left\{  \mathbf{q}\right\}  }}$ from $\mathbb{F}\left(
x_{\mathbf{P}\setminus\left\{  \mathbf{q}\right\}  }\right)  $ to
$\mathbb{F}\left(  x_{\mathbf{P}\setminus\left\{  \mathbf{q}\right\}
}\right)  $ is well-defined. Denote this $\mathbb{F}$-algebra homomorphism by
$\rho$. Then,
\[
\rho=\operatorname*{rev}\nolimits_{\left(  R_{\mathbf{p}}\right)
_{\mathbf{p}\in\mathbf{P}\setminus\left\{  \mathbf{q}\right\}  }}%
\]
is a ring homomorphism from $\mathbb{F}\left(  x_{\mathbf{P}\setminus\left\{
\mathbf{q}\right\}  }\right)  $ to $\mathbb{F}\left(  x_{\mathbf{P}%
\setminus\left\{  \mathbf{q}\right\}  }\right)  $, hence a field homomorphism
from $\mathbb{F}\left(  x_{\mathbf{P}\setminus\left\{  \mathbf{q}\right\}
}\right)  $ to $\mathbb{F}\left(  x_{\mathbf{P}\setminus\left\{
\mathbf{q}\right\}  }\right)  $ (because any ring homomorphism from a field to
a field is a field homomorphism, and because $\mathbb{F}\left(  x_{\mathbf{P}%
\setminus\left\{  \mathbf{q}\right\}  }\right)  $ is a field). Consequently,
$\rho$ is injective (because any field homomorphism is injective).\footnote{Of
course, this also follows from the fact that the maps $\operatorname*{rev}%
\nolimits_{\left(  Q_{\mathbf{p}}\right)  _{\mathbf{p}\in\mathbf{P}%
\setminus\left\{  \mathbf{q}\right\}  }}$ and $\operatorname*{rev}%
\nolimits_{\left(  R_{\mathbf{p}}\right)  _{\mathbf{p}\in\mathbf{P}%
\setminus\left\{  \mathbf{q}\right\}  }}=\rho$ are mutually inverse.}

Since $\alpha$, $\beta$, $\gamma$, $\delta$ are elements of $\mathbb{F}\left(
x_{\mathbf{P}\setminus\left\{  \mathbf{q}\right\}  }\right)  $, the values
$\rho\left(  \alpha\right)  $, $\rho\left(  \beta\right)  $, $\rho\left(
\gamma\right)  $ and $\rho\left(  \delta\right)  $ are well-defined elements
of $\mathbb{F}\left(  x_{\mathbf{P}\setminus\left\{  \mathbf{q}\right\}
}\right)  $. It is easy to see that $\rho\left(  \alpha\right)  -x_{\mathbf{q}%
}\rho\left(  \gamma\right)  \neq0$\ \ \ \ \footnote{\textit{Proof.} Assume the
contrary. Thus, $\rho\left(  \alpha\right)  -x_{\mathbf{q}}\rho\left(
\gamma\right)  =0$. Hence, $\rho\left(  \alpha\right)  =x_{\mathbf{q}}%
\rho\left(  \gamma\right)  $.
\par
If $\gamma\neq0$, then $\rho\left(  \gamma\right)  \neq0$ (since $\rho$ is
injective). Hence, if $\gamma\neq0$, then we can divide $\rho\left(
\alpha\right)  =x_{\mathbf{q}}\rho\left(  \gamma\right)  $ by $\rho\left(
\gamma\right)  $ and obtain $\dfrac{\rho\left(  \alpha\right)  }{\rho\left(
\gamma\right)  }=x_{\mathbf{q}}$. Thus, if $\gamma\neq0$, then%
\begin{align*}
x_{\mathbf{q}}  &  =\dfrac{\rho\left(  \alpha\right)  }{\rho\left(
\gamma\right)  }=\rho\left(  \dfrac{\alpha}{\gamma}\right)
\ \ \ \ \ \ \ \ \ \ \left(  \text{since }\rho\text{ is a field homomorphism}%
\right) \\
&  \in\mathbb{F}\left(  x_{\mathbf{P}\setminus\left\{  \mathbf{q}\right\}
}\right)  \ \ \ \ \ \ \ \ \ \ \left(  \text{since }\mathbb{F}\left(
x_{\mathbf{P}\setminus\left\{  \mathbf{q}\right\}  }\right)  \text{ is the
target of }\rho\right)  ,
\end{align*}
which is absurd (because the indeterminate $x_{\mathbf{q}}$ is foreign to the
field $\mathbb{F}\left(  x_{\mathbf{P}\setminus\left\{  \mathbf{q}\right\}
}\right)  $). Thus, if $\gamma\neq0$, then we have obtained a contradiction.
Consequently, we cannot have $\gamma\neq0$. Hence, $\gamma=0$, so that
$\alpha\delta-\beta\underbrace{\gamma}_{=0}=\alpha\delta$. Thus, $\alpha
\delta=\alpha\delta-\beta\gamma\neq0$, so that $\alpha\neq0$ and thus
$\rho\left(  \alpha\right)  \neq0$ (since $\rho$ is injective). But this
contradicts $\rho\left(  \alpha\right)  =x_{\mathbf{q}}\rho\left(
\underbrace{\gamma}_{=0}\right)  =x_{\mathbf{q}}\underbrace{\rho\left(
0\right)  }_{=0}=0$. This contradiction shows that our assumption was wrong,
qed.}. Hence, the fraction $\dfrac{\rho\left(  \delta\right)  x_{\mathbf{q}%
}-\rho\left(  \beta\right)  }{\rho\left(  \alpha\right)  -x_{\mathbf{q}}%
\rho\left(  \gamma\right)  }$ is a well-defined element of $\mathbb{F}\left(
x_{\mathbf{P}}\right)  $.

Now, recall our family $\left(  R_{\mathbf{p}}\right)  _{\mathbf{p}%
\in\mathbf{P}\setminus\left\{  \mathbf{q}\right\}  }\in\left(  \mathbb{F}%
\left(  x_{\mathbf{P}\setminus\left\{  \mathbf{q}\right\}  }\right)  \right)
^{\mathbf{P}\setminus\left\{  \mathbf{q}\right\}  }\subseteq\left(
\mathbb{F}\left(  x_{\mathbf{P}}\right)  \right)  ^{\mathbf{P}\setminus
\left\{  \mathbf{q}\right\}  }$. Since $\mathbf{q}\notin\mathbf{P}%
\setminus\left\{  \mathbf{q}\right\}  $, an element $R_{\mathbf{q}}$ is not
yet defined. We define an element $R_{\mathbf{q}}\in\mathbb{F}\left(
x_{\mathbf{P}}\right)  $ by%
\[
R_{\mathbf{q}}=\dfrac{\rho\left(  \delta\right)  x_{\mathbf{q}}-\rho\left(
\beta\right)  }{\rho\left(  \alpha\right)  -x_{\mathbf{q}}\rho\left(
\gamma\right)  }.
\]
Having defined $R_{\mathbf{q}}$, we thus have extended our family $\left(
R_{\mathbf{p}}\right)  _{\mathbf{p}\in\mathbf{P}\setminus\left\{
\mathbf{q}\right\}  }\in\left(  \mathbb{F}\left(  x_{\mathbf{P}}\right)
\right)  ^{\mathbf{P}\setminus\left\{  \mathbf{q}\right\}  }$ to a family
$\left(  R_{\mathbf{p}}\right)  _{\mathbf{p}\in\mathbf{P}}\in\left(
\mathbb{F}\left(  x_{\mathbf{P}}\right)  \right)  ^{\mathbf{P}}$. We also
have
\[
R_{\mathbf{q}}=\dfrac{\rho\left(  \delta\right)  x_{\mathbf{q}}-\rho\left(
\beta\right)  }{\rho\left(  \alpha\right)  -x_{\mathbf{q}}\rho\left(
\gamma\right)  }=\dfrac{\rho\left(  \delta\right)  x_{\mathbf{q}}+\left(
-\rho\left(  \beta\right)  \right)  }{\left(  -\rho\left(  \gamma\right)
\right)  x_{\mathbf{q}}+\rho\left(  \alpha\right)  }%
\]
and
\begin{align*}
\rho\left(  \delta\right)  \rho\left(  \alpha\right)  -\left(  -\rho\left(
\beta\right)  \right)  \left(  -\rho\left(  \gamma\right)  \right)   &
=\rho\left(  \alpha\right)  \rho\left(  \delta\right)  -\rho\left(
\beta\right)  \rho\left(  \gamma\right)  =\rho\left(  \alpha\delta-\beta
\gamma\right) \\
&  \ \ \ \ \ \ \ \ \ \ \left(  \text{since }\rho\text{ is an }\mathbb{F}%
\text{-algebra homomorphism}\right) \\
&  \neq0\ \ \ \ \ \ \ \ \ \ \left(  \text{since }\rho\text{ is injective and
since }\alpha\delta-\beta\gamma\neq0\right)  .
\end{align*}

Now, applying Lemma \ref{lem.algebraic.triangularity.ind} to $R_{\mathbf{p}}$,
$\rho\left(  \delta\right)  $, $-\rho\left(  \beta\right)  $, $-\rho\left(
\gamma\right)  $ and $\rho\left(  \alpha\right)  $ instead of $Q_{\mathbf{q}}%
$, $\alpha$, $\beta$, $\gamma$ and $\delta$, we conclude that the family
$\left(  R_{\mathbf{p}}\right)  _{\mathbf{p}\in\mathbf{P}}\in\left(
\mathbb{F}\left(  x_{\mathbf{P}}\right)  \right)  ^{\mathbf{P}}$ is
$\mathbf{P}$-triangular.

Hence, we can apply Lemma \ref{lem.algebraic.triangularity} \textbf{(a)} to
$\left(  R_{\mathbf{p}}\right)  _{\mathbf{p}\in\mathbf{P}}$ instead of
$\left(  Q_{\mathbf{p}}\right)  _{\mathbf{p}\in\mathbf{P}}$. We conclude that
the family $\left(  R_{\mathbf{p}}\right)  _{\mathbf{p}\in\mathbf{P}}%
\in\left(  \mathbb{F}\left(  x_{\mathbf{P}}\right)  \right)  ^{\mathbf{P}}$ is
algebraically independent. Hence, the $\mathbb{F}$-algebra homomorphism
$\operatorname*{rev}\nolimits_{\left(  R_{\mathbf{p}}\right)  _{\mathbf{p}%
\in\mathbf{P}}}$ from $\mathbb{F}\left(  x_{\mathbf{P}}\right)  $ to
$\mathbb{F}\left(  x_{\mathbf{P}}\right)  $ is well-defined.

It now remains to show that the maps $\operatorname*{rev}\nolimits_{\left(
Q_{\mathbf{p}}\right)  _{\mathbf{p}\in\mathbf{P}}}$ and $\operatorname*{rev}%
\nolimits_{\left(  R_{\mathbf{p}}\right)  _{\mathbf{p}\in\mathbf{P}}}$ are
mutually inverse. Indeed, define two $\mathbb{F}$-algebra homomorphisms
$\mathfrak{a}$ and $\mathfrak{b}$ from $\mathbb{F}\left(  x_{\mathbf{P}%
}\right)  $ to $\mathbb{F}\left(  x_{\mathbf{P}}\right)  $ by
\[
\mathfrak{a}=\operatorname*{rev}\nolimits_{\left(  Q_{\mathbf{p}}\right)
_{\mathbf{p}\in\mathbf{P}}}\circ\operatorname*{rev}\nolimits_{\left(
R_{\mathbf{p}}\right)  _{\mathbf{p}\in\mathbf{P}}}%
\ \ \ \ \ \ \ \ \ \ \text{and}\ \ \ \ \ \ \ \ \ \ \mathfrak{b}%
=\operatorname*{id}.
\]
(These are indeed $\mathbb{F}$-algebra homomorphisms because
$\operatorname*{rev}\nolimits_{\left(  Q_{\mathbf{p}}\right)  _{\mathbf{p}%
\in\mathbf{P}}}$ and $\operatorname*{rev}\nolimits_{\left(  R_{\mathbf{p}%
}\right)  _{\mathbf{p}\in\mathbf{P}}}$ are $\mathbb{F}$-algebra
homomorphisms.) We are going to prove that $\mathfrak{a}=\mathfrak{b}$.

First of all, notice that
\begin{equation}
\operatorname*{rev}\nolimits_{\left(  R_{\mathbf{p}}\right)  _{\mathbf{p}%
\in\mathbf{P}}}\left(  f\right)  =\rho\left(  f\right)
\ \ \ \ \ \ \ \ \ \ \text{for every }f\in\mathbb{F}\left(  x_{\mathbf{P}%
\setminus\left\{  \mathbf{q}\right\}  }\right)  .
\label{pf.algebraic.triangularity.5.pf1}%
\end{equation}
\footnote{\textit{Proof of (\ref{pf.algebraic.triangularity.5.pf1}):} Let
$f\in\mathbb{F}\left(  x_{\mathbf{P}\setminus\left\{  \mathbf{q}\right\}
}\right)  $. Then, the variable $x_{\mathbf{q}}$ does not appear in $f$. We
thus have $f\left(  \left(  R_{\mathbf{p}}\right)  _{\mathbf{p}\in\mathbf{P}%
}\right)  =f\left(  \left(  R_{\mathbf{p}}\right)  _{\mathbf{p}\in
\mathbf{P}\setminus\left\{  \mathbf{q}\right\}  }\right)  $ (by means of our
identification of $\mathbb{F}\left(  x_{\mathbf{P}\setminus\left\{
\mathbf{q}\right\}  }\right)  $ with a subfield of $\mathbb{F}\left(
x_{\mathbf{P}}\right)  $). But $\operatorname*{rev}\nolimits_{\left(
R_{\mathbf{p}}\right)  _{\mathbf{p}\in\mathbf{P}}}\left(  f\right)  =f\left(
\left(  R_{\mathbf{p}}\right)  _{\mathbf{p}\in\mathbf{P}}\right)  $ (by the
definition of $\operatorname*{rev}\nolimits_{\left(  R_{\mathbf{p}}\right)
_{\mathbf{p}\in\mathbf{P}}}$). Since $\rho=\operatorname*{rev}%
\nolimits_{\left(  R_{\mathbf{p}}\right)  _{\mathbf{p}\in\mathbf{P}%
\setminus\left\{  \mathbf{q}\right\}  }}$, we have $\rho\left(  f\right)
=\operatorname*{rev}\nolimits_{\left(  R_{\mathbf{p}}\right)  _{\mathbf{p}%
\in\mathbf{P}\setminus\left\{  \mathbf{q}\right\}  }}\left(  f\right)
=f\left(  \left(  R_{\mathbf{p}}\right)  _{\mathbf{p}\in\mathbf{P}%
\setminus\left\{  \mathbf{q}\right\}  }\right)  $ (by the definition of
$\operatorname*{rev}\nolimits_{\left(  R_{\mathbf{p}}\right)  _{\mathbf{p}%
\in\mathbf{P}\setminus\left\{  \mathbf{q}\right\}  }}$). Thus,%
\[
\operatorname*{rev}\nolimits_{\left(  R_{\mathbf{p}}\right)  _{\mathbf{p}%
\in\mathbf{P}}}\left(  f\right)  =f\left(  \left(  R_{\mathbf{p}}\right)
_{\mathbf{p}\in\mathbf{P}}\right)  =f\left(  \left(  R_{\mathbf{p}}\right)
_{\mathbf{p}\in\mathbf{P}\setminus\left\{  \mathbf{q}\right\}  }\right)
=\rho\left(  f\right)  .
\]
This proves (\ref{pf.algebraic.triangularity.5.pf1}).} Moreover,%
\begin{equation}
\mathfrak{a}\left(  f\right)  =f\ \ \ \ \ \ \ \ \ \ \text{for every }%
f\in\mathbb{F}\left(  x_{\mathbf{P}\setminus\left\{  \mathbf{q}\right\}
}\right)  . \label{pf.algebraic.triangularity.5.pf2}%
\end{equation}
\footnote{\textit{Proof of (\ref{pf.algebraic.triangularity.5.pf2}):} Let
$f\in\mathbb{F}\left(  x_{\mathbf{P}\setminus\left\{  \mathbf{q}\right\}
}\right)  $. Let $g=\rho\left(  f\right)  $. Then, $g\in\mathbb{F}\left(
x_{\mathbf{P}\setminus\left\{  \mathbf{q}\right\}  }\right)  $ (since the
target of $\rho$ is $\mathbb{F}\left(  x_{\mathbf{P}\setminus\left\{
\mathbf{q}\right\}  }\right)  $). Hence, the variable $x_{\mathbf{q}}$ does
not appear in $g$. We thus have $g\left(  \left(  Q_{\mathbf{p}}\right)
_{\mathbf{p}\in\mathbf{P}}\right)  =g\left(  \left(  Q_{\mathbf{p}}\right)
_{\mathbf{p}\in\mathbf{P}\setminus\left\{  \mathbf{q}\right\}  }\right)  $ (by
means of our identification of $\mathbb{F}\left(  x_{\mathbf{P}\setminus
\left\{  \mathbf{q}\right\}  }\right)  $ with a subfield of $\mathbb{F}\left(
x_{\mathbf{P}}\right)  $). But $\operatorname*{rev}\nolimits_{\left(
Q_{\mathbf{p}}\right)  _{\mathbf{p}\in\mathbf{P}}}\left(  g\right)  =g\left(
\left(  Q_{\mathbf{p}}\right)  _{\mathbf{p}\in\mathbf{P}}\right)  $ (by the
definition of $\operatorname*{rev}\nolimits_{\left(  Q_{\mathbf{p}}\right)
_{\mathbf{p}\in\mathbf{P}}}$) and $\operatorname*{rev}\nolimits_{\left(
Q_{\mathbf{p}}\right)  _{\mathbf{p}\in\mathbf{P}\setminus\left\{
\mathbf{q}\right\}  }}\left(  g\right)  =g\left(  \left(  Q_{\mathbf{p}%
}\right)  _{\mathbf{p}\in\mathbf{P}\setminus\left\{  \mathbf{q}\right\}
}\right)  $ (by the definition of $\operatorname*{rev}\nolimits_{\left(
Q_{\mathbf{p}}\right)  _{\mathbf{p}\in\mathbf{P}\setminus\left\{
\mathbf{q}\right\}  }}$). Now, since $\mathfrak{a}=\operatorname*{rev}%
\nolimits_{\left(  Q_{\mathbf{p}}\right)  _{\mathbf{p}\in\mathbf{P}}}%
\circ\operatorname*{rev}\nolimits_{\left(  R_{\mathbf{p}}\right)
_{\mathbf{p}\in\mathbf{P}}}$, we have%
\begin{align*}
\mathfrak{a}\left(  f\right)   &  =\left(  \operatorname*{rev}%
\nolimits_{\left(  Q_{\mathbf{p}}\right)  _{\mathbf{p}\in\mathbf{P}}}%
\circ\operatorname*{rev}\nolimits_{\left(  R_{\mathbf{p}}\right)
_{\mathbf{p}\in\mathbf{P}}}\right)  \left(  f\right)  =\operatorname*{rev}%
\nolimits_{\left(  Q_{\mathbf{p}}\right)  _{\mathbf{p}\in\mathbf{P}}}\left(
\underbrace{\operatorname*{rev}\nolimits_{\left(  R_{\mathbf{p}}\right)
_{\mathbf{p}\in\mathbf{P}}}\left(  f\right)  }_{\substack{=\rho\left(
f\right)  \\\text{(by (\ref{pf.algebraic.triangularity.5.pf1}))}}}\right)
=\operatorname*{rev}\nolimits_{\left(  Q_{\mathbf{p}}\right)  _{\mathbf{p}%
\in\mathbf{P}}}\left(  \underbrace{\rho\left(  f\right)  }_{=g}\right) \\
&  =\operatorname*{rev}\nolimits_{\left(  Q_{\mathbf{p}}\right)
_{\mathbf{p}\in\mathbf{P}}}\left(  g\right)  =g\left(  \left(  Q_{\mathbf{p}%
}\right)  _{\mathbf{p}\in\mathbf{P}}\right)  =g\left(  \left(  Q_{\mathbf{p}%
}\right)  _{\mathbf{p}\in\mathbf{P}\setminus\left\{  \mathbf{q}\right\}
}\right)  =\operatorname*{rev}\nolimits_{\left(  Q_{\mathbf{p}}\right)
_{\mathbf{p}\in\mathbf{P}\setminus\left\{  \mathbf{q}\right\}  }}\left(
\underbrace{g}_{=\rho\left(  f\right)  }\right) \\
&  =\operatorname*{rev}\nolimits_{\left(  Q_{\mathbf{p}}\right)
_{\mathbf{p}\in\mathbf{P}\setminus\left\{  \mathbf{q}\right\}  }}\left(
\rho\left(  f\right)  \right)  =\left(  \operatorname*{rev}\nolimits_{\left(
Q_{\mathbf{p}}\right)  _{\mathbf{p}\in\mathbf{P}\setminus\left\{
\mathbf{q}\right\}  }}\circ\underbrace{\rho}_{=\operatorname*{rev}%
\nolimits_{\left(  R_{\mathbf{p}}\right)  _{\mathbf{p}\in\mathbf{P}%
\setminus\left\{  \mathbf{q}\right\}  }}}\right)  \left(  f\right) \\
&  =\underbrace{\left(  \operatorname*{rev}\nolimits_{\left(  Q_{\mathbf{p}%
}\right)  _{\mathbf{p}\in\mathbf{P}\setminus\left\{  \mathbf{q}\right\}  }%
}\circ\operatorname*{rev}\nolimits_{\left(  R_{\mathbf{p}}\right)
_{\mathbf{p}\in\mathbf{P}\setminus\left\{  \mathbf{q}\right\}  }}\right)
}_{\substack{=\operatorname*{id}\\\text{(since the maps }\operatorname*{rev}%
\nolimits_{\left(  Q_{\mathbf{p}}\right)  _{\mathbf{p}\in\mathbf{P}%
\setminus\left\{  \mathbf{q}\right\}  }}\text{ and}\\\operatorname*{rev}%
\nolimits_{\left(  R_{\mathbf{p}}\right)  _{\mathbf{p}\in\mathbf{P}%
\setminus\left\{  \mathbf{q}\right\}  }}\text{ are mutually inverse) }%
}}\left(  f\right)  =\operatorname*{id}\left(  f\right)  =f,
\end{align*}
so that (\ref{pf.algebraic.triangularity.5.pf2}) is proven.} Furthermore,%
\begin{equation}
\operatorname*{rev}\nolimits_{\left(  Q_{\mathbf{p}}\right)  _{\mathbf{p}%
\in\mathbf{P}}}\left(  \rho\left(  f\right)  \right)
=f\ \ \ \ \ \ \ \ \ \ \text{for every }f\in\mathbb{F}\left(  x_{\mathbf{P}%
\setminus\left\{  \mathbf{q}\right\}  }\right)  .
\label{pf.algebraic.triangularity.5.pf3}%
\end{equation}
\footnote{\textit{Proof of (\ref{pf.algebraic.triangularity.5.pf3}):} Let
$f\in\mathbb{F}\left(  x_{\mathbf{P}\setminus\left\{  \mathbf{q}\right\}
}\right)  $. Since $\rho\left(  f\right)  =\operatorname*{rev}%
\nolimits_{\left(  R_{\mathbf{p}}\right)  _{\mathbf{p}\in\mathbf{P}}}\left(
f\right)  $ (by (\ref{pf.algebraic.triangularity.5.pf1})), we have%
\[
\operatorname*{rev}\nolimits_{\left(  Q_{\mathbf{p}}\right)  _{\mathbf{p}%
\in\mathbf{P}}}\left(  \rho\left(  f\right)  \right)  =\operatorname*{rev}%
\nolimits_{\left(  Q_{\mathbf{p}}\right)  _{\mathbf{p}\in\mathbf{P}}}\left(
\operatorname*{rev}\nolimits_{\left(  R_{\mathbf{p}}\right)  _{\mathbf{p}%
\in\mathbf{P}}}\left(  f\right)  \right)  =\underbrace{\left(
\operatorname*{rev}\nolimits_{\left(  Q_{\mathbf{p}}\right)  _{\mathbf{p}%
\in\mathbf{P}}}\circ\operatorname*{rev}\nolimits_{\left(  R_{\mathbf{p}%
}\right)  _{\mathbf{p}\in\mathbf{P}}}\right)  }_{=\mathfrak{a}}\left(
f\right)  =\mathfrak{a}\left(  f\right)  =f
\]
(by (\ref{pf.algebraic.triangularity.5.pf2})). This proves
(\ref{pf.algebraic.triangularity.5.pf3}).} Applying
(\ref{pf.algebraic.triangularity.5.pf3}) to $f=\alpha$, we obtain
$\operatorname*{rev}\nolimits_{\left(  Q_{\mathbf{p}}\right)  _{\mathbf{p}%
\in\mathbf{P}}}\left(  \rho\left(  \alpha\right)  \right)  =\alpha$ (since
$\alpha\in\mathbb{F}\left(  x_{\mathbf{P}\setminus\left\{  \mathbf{q}\right\}
}\right)  $). Similarly, $\operatorname*{rev}\nolimits_{\left(  Q_{\mathbf{p}%
}\right)  _{\mathbf{p}\in\mathbf{P}}}\left(  \rho\left(  \beta\right)
\right)  =\beta$, $\operatorname*{rev}\nolimits_{\left(  Q_{\mathbf{p}%
}\right)  _{\mathbf{p}\in\mathbf{P}}}\left(  \rho\left(  \gamma\right)
\right)  =\gamma$ and $\operatorname*{rev}\nolimits_{\left(  Q_{\mathbf{p}%
}\right)  _{\mathbf{p}\in\mathbf{P}}}\left(  \rho\left(  \delta\right)
\right)  =\delta$.

Now, by the definition of $\operatorname*{rev}\nolimits_{\left(
R_{\mathbf{p}}\right)  _{\mathbf{p}\in\mathbf{P}}}$, we have%
\begin{equation}
\operatorname*{rev}\nolimits_{\left(  R_{\mathbf{p}}\right)  _{\mathbf{p}%
\in\mathbf{P}}}\left(  x_{\mathbf{q}}\right)  =x_{\mathbf{q}}\left(  \left(
R_{\mathbf{p}}\right)  _{\mathbf{p}\in\mathbf{P}}\right)  =R_{\mathbf{q}%
}=\dfrac{\rho\left(  \delta\right)  x_{\mathbf{q}}-\rho\left(  \beta\right)
}{\rho\left(  \alpha\right)  -x_{\mathbf{q}}\rho\left(  \gamma\right)  }.
\label{pf.algebraic.triangularity.5.pf5}%
\end{equation}
By the definition of $\operatorname*{rev}\nolimits_{\left(  Q_{\mathbf{p}%
}\right)  _{\mathbf{p}\in\mathbf{P}}}$, we have%
\begin{equation}
\operatorname*{rev}\nolimits_{\left(  Q_{\mathbf{p}}\right)  _{\mathbf{p}%
\in\mathbf{P}}}\left(  x_{\mathbf{q}}\right)  =x_{\mathbf{q}}\left(  \left(
Q_{\mathbf{p}}\right)  _{\mathbf{p}\in\mathbf{P}}\right)  =Q_{\mathbf{q}}.
\label{pf.algebraic.triangularity.5.pf6}%
\end{equation}

Since $\mathfrak{a}=\operatorname*{rev}\nolimits_{\left(  Q_{\mathbf{p}%
}\right)  _{\mathbf{p}\in\mathbf{P}}}\circ\operatorname*{rev}%
\nolimits_{\left(  R_{\mathbf{p}}\right)  _{\mathbf{p}\in\mathbf{P}}}$, we
have%
\begin{align}
\mathfrak{a}\left(  x_{\mathbf{q}}\right)   &  =\left(  \operatorname*{rev}%
\nolimits_{\left(  Q_{\mathbf{p}}\right)  _{\mathbf{p}\in\mathbf{P}}}%
\circ\operatorname*{rev}\nolimits_{\left(  R_{\mathbf{p}}\right)
_{\mathbf{p}\in\mathbf{P}}}\right)  \left(  x_{\mathbf{q}}\right)
=\operatorname*{rev}\nolimits_{\left(  Q_{\mathbf{p}}\right)  _{\mathbf{p}%
\in\mathbf{P}}}\left(  \underbrace{\operatorname*{rev}\nolimits_{\left(
R_{\mathbf{p}}\right)  _{\mathbf{p}\in\mathbf{P}}}\left(  x_{\mathbf{q}%
}\right)  }_{\substack{=\dfrac{\rho\left(  \delta\right)  x_{\mathbf{q}}%
-\rho\left(  \beta\right)  }{\rho\left(  \alpha\right)  -x_{\mathbf{q}}%
\rho\left(  \gamma\right)  }\\\text{(by
(\ref{pf.algebraic.triangularity.5.pf5}))}}}\right) \nonumber\\
&  =\operatorname*{rev}\nolimits_{\left(  Q_{\mathbf{p}}\right)
_{\mathbf{p}\in\mathbf{P}}}\left(  \dfrac{\rho\left(  \delta\right)
x_{\mathbf{q}}-\rho\left(  \beta\right)  }{\rho\left(  \alpha\right)
-x_{\mathbf{q}}\rho\left(  \gamma\right)  }\right)  =\dfrac
{\operatorname*{rev}\nolimits_{\left(  Q_{\mathbf{p}}\right)  _{\mathbf{p}%
\in\mathbf{P}}}\left(  \rho\left(  \delta\right)  \right)  \operatorname*{rev}%
\nolimits_{\left(  Q_{\mathbf{p}}\right)  _{\mathbf{p}\in\mathbf{P}}}\left(
x_{\mathbf{q}}\right)  -\operatorname*{rev}\nolimits_{\left(  Q_{\mathbf{p}%
}\right)  _{\mathbf{p}\in\mathbf{P}}}\left(  \rho\left(  \beta\right)
\right)  }{\operatorname*{rev}\nolimits_{\left(  Q_{\mathbf{p}}\right)
_{\mathbf{p}\in\mathbf{P}}}\left(  \rho\left(  \alpha\right)  \right)
-\operatorname*{rev}\nolimits_{\left(  Q_{\mathbf{p}}\right)  _{\mathbf{p}%
\in\mathbf{P}}}\left(  x_{\mathbf{q}}\right)  \operatorname*{rev}%
\nolimits_{\left(  Q_{\mathbf{p}}\right)  _{\mathbf{p}\in\mathbf{P}}}\left(
\rho\left(  \gamma\right)  \right)  }\nonumber\\
&  \ \ \ \ \ \ \ \ \ \ \left(  \text{since }\operatorname*{rev}%
\nolimits_{\left(  Q_{\mathbf{p}}\right)  _{\mathbf{p}\in\mathbf{P}}}\text{ is
a field homomorphism}\right) \nonumber\\
&  =\dfrac{\delta Q_{\mathbf{q}}-\beta}{\alpha-Q_{\mathbf{q}}\gamma
}\nonumber\\
&  \ \ \ \ \ \ \ \ \ \ \left(
\begin{array}
[c]{c}%
\text{since }\operatorname*{rev}\nolimits_{\left(  Q_{\mathbf{p}}\right)
_{\mathbf{p}\in\mathbf{P}}}\left(  \rho\left(  \alpha\right)  \right)
=\alpha\text{, }\operatorname*{rev}\nolimits_{\left(  Q_{\mathbf{p}}\right)
_{\mathbf{p}\in\mathbf{P}}}\left(  \rho\left(  \beta\right)  \right)
=\beta\text{, }\operatorname*{rev}\nolimits_{\left(  Q_{\mathbf{p}}\right)
_{\mathbf{p}\in\mathbf{P}}}\left(  \rho\left(  \gamma\right)  \right)
=\gamma\text{,}\\
\operatorname*{rev}\nolimits_{\left(  Q_{\mathbf{p}}\right)  _{\mathbf{p}%
\in\mathbf{P}}}\left(  \rho\left(  \delta\right)  \right)  =\delta\text{ and
}\operatorname*{rev}\nolimits_{\left(  Q_{\mathbf{p}}\right)  _{\mathbf{p}%
\in\mathbf{P}}}\left(  x_{\mathbf{q}}\right)  =Q_{\mathbf{q}}%
\end{array}
\right) \nonumber\\
&  =\dfrac{Q_{\mathbf{q}}\delta-\beta}{\alpha-Q_{\mathbf{q}}\gamma
}=x_{\mathbf{q}}\label{pf.algebraic.triangularity.5.pf8}\\
&  \ \ \ \ \ \ \ \ \ \ \left(  \text{since Lemma
\ref{lem.algebraic.triangularity.step} \textbf{(a)} yields }\alpha
-Q_{\mathbf{q}}\gamma\neq0\text{ and }x_{\mathbf{q}}=\dfrac{Q_{\mathbf{q}%
}\delta-\beta}{\alpha-Q_{\mathbf{q}}\gamma}\right)  .\nonumber
\end{align}

Now, it is easy to see that%
\begin{equation}
\mathfrak{a}\left(  x_{\mathbf{p}}\right)  =\mathfrak{b}\left(  x_{\mathbf{p}%
}\right)  \ \ \ \ \ \ \ \ \ \ \text{for every }\mathbf{p}\in\mathbf{P}.
\label{pf.algebraic.triangularity.5.pf9}%
\end{equation}
\footnote{\textit{Proof of (\ref{pf.algebraic.triangularity.5.pf9}):} Let
$\mathbf{p}\in\mathbf{P}$. Notice that $\mathfrak{b}=\operatorname*{id}$, thus
$\mathfrak{b}\left(  x_{\mathbf{p}}\right)  =x_{\mathbf{p}}$.
\par
If $\mathbf{p}=\mathbf{q}$, then (\ref{pf.algebraic.triangularity.5.pf9}) is
clearly true (because if $\mathbf{p}=\mathbf{q}$, then%
\begin{align*}
\mathfrak{a}\left(  x_{\mathbf{p}}\right)   &  =\mathfrak{a}\left(
x_{\mathbf{q}}\right)  =x_{\mathbf{q}}\ \ \ \ \ \ \ \ \ \ \left(  \text{by
(\ref{pf.algebraic.triangularity.5.pf8})}\right) \\
&  =x_{\mathbf{p}}\ \ \ \ \ \ \ \ \ \ \left(  \text{since }\mathbf{q}%
=\mathbf{p}\right) \\
&  =\mathfrak{b}\left(  x_{\mathbf{p}}\right)
\end{align*}
). Hence, for the rest of the proof of (\ref{pf.algebraic.triangularity.5.pf9}%
), we can WLOG assume that we don't have $\mathbf{p}=\mathbf{q}$. Assume this.
\par
We have $\mathbf{p}\in\mathbf{P}\setminus\left\{  \mathbf{q}\right\}  $ (since
$\mathbf{p}\in\mathbf{P}$ but we don't have $\mathbf{p}=\mathbf{q}$). Thus,
$x_{\mathbf{p}}\in\mathbb{F}\left(  x_{\mathbf{P}\setminus\left\{
\mathbf{q}\right\}  }\right)  $, so that we can apply
(\ref{pf.algebraic.triangularity.5.pf2}) to $f=x_{\mathbf{p}}$ and obtain
$\mathfrak{a}\left(  x_{\mathbf{p}}\right)  =x_{\mathbf{p}}=\mathfrak{b}%
\left(  x_{\mathbf{p}}\right)  $. This proves
(\ref{pf.algebraic.triangularity.5.pf9}).} Thus, Lemma
\ref{lem.algebraic.triangularity.id} (applied to $\mathbb{K}=\mathbb{F}\left(
x_{\mathbf{P}}\right)  $) yields that $\mathfrak{a}=\mathfrak{b}$. Since
$\mathfrak{a}=\operatorname*{rev}\nolimits_{\left(  Q_{\mathbf{p}}\right)
_{\mathbf{p}\in\mathbf{P}}}\circ\operatorname*{rev}\nolimits_{\left(
R_{\mathbf{p}}\right)  _{\mathbf{p}\in\mathbf{P}}}$ and $\mathfrak{b}%
=\operatorname*{id}$, this rewrites as
\[
\operatorname*{rev}\nolimits_{\left(  Q_{\mathbf{p}}\right)  _{\mathbf{p}%
\in\mathbf{P}}}\circ\operatorname*{rev}\nolimits_{\left(  R_{\mathbf{p}%
}\right)  _{\mathbf{p}\in\mathbf{P}}}=\operatorname*{id}.
\]
Hence,%
\begin{equation}
\underbrace{\operatorname*{rev}\nolimits_{\left(  Q_{\mathbf{p}}\right)
_{\mathbf{p}\in\mathbf{P}}}\circ\operatorname*{rev}\nolimits_{\left(
R_{\mathbf{p}}\right)  _{\mathbf{p}\in\mathbf{P}}}}_{=\operatorname*{id}}%
\circ\operatorname*{rev}\nolimits_{\left(  Q_{\mathbf{p}}\right)
_{\mathbf{p}\in\mathbf{P}}}=\operatorname*{rev}\nolimits_{\left(
Q_{\mathbf{p}}\right)  _{\mathbf{p}\in\mathbf{P}}}.
\label{pf.algebraic.triangularity.6}%
\end{equation}
But $\operatorname*{rev}\nolimits_{\left(  Q_{\mathbf{p}}\right)
_{\mathbf{p}\in\mathbf{P}}}$ is a field homomorphism, and thus injective
(since every field homomorphism is injective), therefore left-cancellable.
Thus, we can cancel $\operatorname*{rev}\nolimits_{\left(  Q_{\mathbf{p}%
}\right)  _{\mathbf{p}\in\mathbf{P}}}$ on the left in the equality
(\ref{pf.algebraic.triangularity.6}), and obtain $\operatorname*{rev}%
\nolimits_{\left(  R_{\mathbf{p}}\right)  _{\mathbf{p}\in\mathbf{P}}}%
\circ\operatorname*{rev}\nolimits_{\left(  Q_{\mathbf{p}}\right)
_{\mathbf{p}\in\mathbf{P}}}=\operatorname*{id}$. Combining this with
$\operatorname*{rev}\nolimits_{\left(  Q_{\mathbf{p}}\right)  _{\mathbf{p}%
\in\mathbf{P}}}\circ\operatorname*{rev}\nolimits_{\left(  R_{\mathbf{p}%
}\right)  _{\mathbf{p}\in\mathbf{P}}}=\operatorname*{id}$, we conclude that
the maps $\operatorname*{rev}\nolimits_{\left(  Q_{\mathbf{p}}\right)
_{\mathbf{p}\in\mathbf{P}}}$ and $\operatorname*{rev}\nolimits_{\left(
R_{\mathbf{p}}\right)  _{\mathbf{p}\in\mathbf{P}}}$ are mutually inverse.

We thus have constructed a $\mathbf{P}$-triangular family $\left(
R_{\mathbf{p}}\right)  _{\mathbf{p}\in\mathbf{P}}\in\left(  \mathbb{F}\left(
x_{\mathbf{P}}\right)  \right)  ^{\mathbf{P}}$ such that the maps
$\operatorname*{rev}\nolimits_{\left(  Q_{\mathbf{p}}\right)  _{\mathbf{p}%
\in\mathbf{P}}}$ and $\operatorname*{rev}\nolimits_{\left(  R_{\mathbf{p}%
}\right)  _{\mathbf{p}\in\mathbf{P}}}$ are mutually inverse. Thus, there
exists a $\mathbf{P}$-triangular family $\left(  R_{\mathbf{p}}\right)
_{\mathbf{p}\in\mathbf{P}}\in\left(  \mathbb{F}\left(  x_{\mathbf{P}}\right)
\right)  ^{\mathbf{P}}$ such that the maps $\operatorname*{rev}%
\nolimits_{\left(  Q_{\mathbf{p}}\right)  _{\mathbf{p}\in\mathbf{P}}}$ and
$\operatorname*{rev}\nolimits_{\left(  R_{\mathbf{p}}\right)  _{\mathbf{p}%
\in\mathbf{P}}}$ are mutually inverse. In other words, Lemma
\ref{lem.algebraic.triangularity} \textbf{(b)} holds for our totally ordered
set $\mathbf{P}$. We thus have shown that Lemma
\ref{lem.algebraic.triangularity} \textbf{(b)} holds in the case when
$\left\vert \mathbf{P}\right\vert =N+1$. This finishes the induction step, and
so Lemma \ref{lem.algebraic.triangularity} \textbf{(b)} is proven by induction.
\end{proof}

We next state a rather simple fact about substitution into rational functions:

\begin{lemma}
\label{lem.algebraic.triangularity.substass}Let $\mathbb{F}$ be a field. Let
$\mathbb{K}$ be a field extension of $\mathbb{F}$. Let $\mathbf{P}$ be a
finite set. Let $\left(  A_{\mathbf{p}}\right)  _{\mathbf{p}\in\mathbf{P}}%
\in\mathbb{K}^{\mathbf{P}}$ be a family of elements of $\mathbb{K}$ indexed by
elements of $\mathbf{P}$. Let $\left(  B_{\mathbf{r}}\right)  _{\mathbf{r}%
\in\mathbf{P}}\in\left(  \mathbb{F}\left(  x_{\mathbf{P}}\right)  \right)
^{\mathbf{P}}$ be a family of elements of $\mathbb{F}\left(  x_{\mathbf{P}%
}\right)  $ indexed by elements of $\mathbf{P}$. Let $f\in\mathbb{F}\left(
x_{\mathbf{P}}\right)  $ be arbitrary.

Assume that the element $B_{\mathbf{r}}\left(  \left(  A_{\mathbf{p}}\right)
_{\mathbf{p}\in\mathbf{P}}\right)  $ of $\mathbb{K}$ is well-defined (that is,
substituting $\left(  A_{\mathbf{p}}\right)  _{\mathbf{p}\in\mathbf{P}}$ for
the variables into the rational function $B_{\mathbf{r}}$ does not render the
denominator $0$) for every $\mathbf{r}\in\mathbf{P}$.

Assume further that the element $f\left(  \left(  B_{\mathbf{r}}\left(
\left(  A_{\mathbf{p}}\right)  _{\mathbf{p}\in\mathbf{P}}\right)  \right)
_{\mathbf{r}\in\mathbf{P}}\right)  $ of $\mathbb{K}$ is well-defined (that is,
substituting $\left(  B_{\mathbf{r}}\left(  \left(  A_{\mathbf{p}}\right)
_{\mathbf{p}\in\mathbf{P}}\right)  \right)  _{\mathbf{r}\in\mathbf{P}}$ for
the variables into the rational function $f$ does not render the denominator
$0$).

Then, the element $f\left(  \left(  B_{\mathbf{r}}\right)  _{\mathbf{r}%
\in\mathbf{P}}\right)  $ of $\mathbb{F}\left(  x_{\mathbf{P}}\right)  $ and
the element $\left(  f\left(  \left(  B_{\mathbf{r}}\right)  _{\mathbf{r}%
\in\mathbf{P}}\right)  \right)  \left(  \left(  A_{\mathbf{p}}\right)
_{\mathbf{p}\in\mathbf{P}}\right)  $ of $\mathbb{K}$ are also well-defined,
and satisfy%
\begin{equation}
\left(  f\left(  \left(  B_{\mathbf{r}}\right)  _{\mathbf{r}\in\mathbf{P}%
}\right)  \right)  \left(  \left(  A_{\mathbf{p}}\right)  _{\mathbf{p}%
\in\mathbf{P}}\right)  =f\left(  \left(  B_{\mathbf{r}}\left(  \left(
A_{\mathbf{p}}\right)  _{\mathbf{p}\in\mathbf{P}}\right)  \right)
_{\mathbf{r}\in\mathbf{P}}\right)  .
\label{lem.algebraic.triangularity.substass.eq}%
\end{equation}

\end{lemma}

Despite the bulky notation, Lemma \ref{lem.algebraic.triangularity.substass}
is a rather evident fact. For example, for $\left\vert \mathbf{P}\right\vert
=1$, it boils down to the fact that any two univariate rational functions $f$
and $B$ and any $A\in\mathbb{K}$ satisfy $\left(  f\left(  B\right)  \right)
\left(  A\right)  =f\left(  B\left(  A\right)  \right)  $, provided that
$f\left(  B\right)  $ and $\left(  f\left(  B\right)  \right)  \left(
A\right)  $ are well-defined. The formal proof is rather straightforward (it
proceeds by reduction first to the case of $f$ being a polynomial, then to the
case of $f$ being a monomial). Here are the details of the proof, purely for
the sake of completeness:

\begin{proof}
[Proof of Lemma \ref{lem.algebraic.triangularity.substass} (sketched).]Let
$\mathbb{F}\left[  x_{\mathbf{P}}\right]  $ denote the polynomial ring over
$\mathbb{F}$ in the indeterminates $x_{\mathbf{p}}$ with $\mathbf{p}$ ranging
over all elements of $\mathbf{P}$. We identify $\mathbb{F}\left[
x_{\mathbf{P}}\right]  $ with a subring of $\mathbb{F}\left(  x_{\mathbf{P}%
}\right)  $. Clearly, substitution of variables into a polynomial is a
particular case of substitution of variables into a rational function. Hence,
for every polynomial $g\in\mathbb{F}\left[  x_{\mathbf{P}}\right]  $, for
every field extension $\mathbb{L}$ of $\mathbb{K}$ and for every family
$\left(  a_{\mathbf{p}}\right)  _{\mathbf{p}\in\mathbf{P}}\in\mathbb{L}%
^{\mathbf{P}}$, the value $g\left(  \left(  a_{\mathbf{p}}\right)
_{\mathbf{p}\in\mathbf{P}}\right)  $ is the element of $\mathbb{L}$ obtained
by substituting $a_{\mathbf{p}}$ for $x_{\mathbf{p}}$ for each $\mathbf{p}%
\in\mathbf{P}$ in the polynomial $g$. (This yields, in particular, that
\begin{equation}
\text{this value }g\left(  \left(  a_{\mathbf{p}}\right)  _{\mathbf{p}%
\in\mathbf{P}}\right)  \text{ is always well-defined when }g\in\mathbb{F}%
\left[  x_{\mathbf{P}}\right]  ,
\label{pf.algebraic.triangularity.substass.wd}%
\end{equation}
because in a polynomial there are no denominators that can vanish upon
substitution.) We are going to use the word \textquotedblleft
monomial\textquotedblright\ to mean a monomial without coefficient. Now, it is
easy to see that for every monomial $m$ in $\mathbb{F}\left[  x_{\mathbf{P}%
}\right]  $, we have%
\begin{equation}
\left(  \left(  m\left(  \left(  B_{\mathbf{r}}\right)  _{\mathbf{r}%
\in\mathbf{P}}\right)  \right)  \left(  \left(  A_{\mathbf{p}}\right)
_{\mathbf{p}\in\mathbf{P}}\right)  \text{ is well-defined}\right)
\label{pf.algebraic.triangularity.substass.1a}%
\end{equation}
and%
\begin{equation}
\left(  m\left(  \left(  B_{\mathbf{r}}\right)  _{\mathbf{r}\in\mathbf{P}%
}\right)  \right)  \left(  \left(  A_{\mathbf{p}}\right)  _{\mathbf{p}%
\in\mathbf{P}}\right)  =m\left(  \left(  B_{\mathbf{r}}\left(  \left(
A_{\mathbf{p}}\right)  _{\mathbf{p}\in\mathbf{P}}\right)  \right)
_{\mathbf{r}\in\mathbf{P}}\right)  .
\label{pf.algebraic.triangularity.substass.1b}%
\end{equation}
\footnote{\textit{Proof of (\ref{pf.algebraic.triangularity.substass.1a}) and
(\ref{pf.algebraic.triangularity.substass.1b}):} Let $m$ be a monomial in
$\mathbb{F}\left[  x_{\mathbf{P}}\right]  $. Then, we can write $m$ in the
form $m=\prod\limits_{\mathbf{q}\in\mathbf{P}}x_{\mathbf{q}}^{n_{\mathbf{q}}}$
for some family $\left(  n_{\mathbf{q}}\right)  _{\mathbf{q}\in\mathbf{P}}%
\in\mathbb{N}^{\mathbf{P}}$. Consider this family $\left(  n_{\mathbf{q}%
}\right)  _{\mathbf{q}\in\mathbf{P}}$.
\par
Notice that $m\left(  \left(  B_{\mathbf{r}}\right)  _{\mathbf{r}\in
\mathbf{P}}\right)  $ and $m\left(  \left(  B_{\mathbf{r}}\left(  \left(
A_{\mathbf{p}}\right)  _{\mathbf{p}\in\mathbf{P}}\right)  \right)
_{\mathbf{r}\in\mathbf{P}}\right)  $ are well-defined (by
(\ref{pf.algebraic.triangularity.substass.wd}) (since $m\in\mathbb{F}\left[
x_{\mathbf{P}}\right]  $)).
\par
Since $m=\prod\limits_{\mathbf{q}\in\mathbf{P}}x_{\mathbf{q}}^{n_{\mathbf{q}}%
}$, we have%
\[
m\left(  \left(  B_{\mathbf{r}}\right)  _{\mathbf{r}\in\mathbf{P}}\right)
=\left(  \prod\limits_{\mathbf{q}\in\mathbf{P}}x_{\mathbf{q}}^{n_{\mathbf{q}}%
}\right)  \left(  \left(  B_{\mathbf{r}}\right)  _{\mathbf{r}\in\mathbf{P}%
}\right)  =\prod\limits_{\mathbf{q}\in\mathbf{P}}B_{\mathbf{q}}^{n_{\mathbf{q}%
}},
\]
so that%
\begin{equation}
\underbrace{\left(  m\left(  \left(  B_{\mathbf{r}}\right)  _{\mathbf{r}%
\in\mathbf{P}}\right)  \right)  }_{=\prod\limits_{\mathbf{q}\in\mathbf{P}%
}B_{\mathbf{q}}^{n_{\mathbf{q}}}}\left(  \left(  A_{\mathbf{p}}\right)
_{\mathbf{p}\in\mathbf{P}}\right)  =\left(  \prod\limits_{\mathbf{q}%
\in\mathbf{P}}B_{\mathbf{q}}^{n_{\mathbf{q}}}\right)  \left(  \left(
A_{\mathbf{p}}\right)  _{\mathbf{p}\in\mathbf{P}}\right)  =\prod
\limits_{\mathbf{q}\in\mathbf{P}}\left(  B_{\mathbf{q}}\left(  \left(
A_{\mathbf{p}}\right)  _{\mathbf{p}\in\mathbf{P}}\right)  \right)
^{n_{\mathbf{q}}}. \label{pf.algebraic.triangularity.substass.1.pf.a}%
\end{equation}
This computation, in particular, shows that $\left(  m\left(  \left(
B_{\mathbf{r}}\right)  _{\mathbf{r}\in\mathbf{P}}\right)  \right)  \left(
\left(  A_{\mathbf{p}}\right)  _{\mathbf{p}\in\mathbf{P}}\right)  $ is
well-defined. Thus, (\ref{pf.algebraic.triangularity.substass.1a}) is proven.
\par
But since $m=\prod\limits_{\mathbf{q}\in\mathbf{P}}x_{\mathbf{q}%
}^{n_{\mathbf{q}}}$, we have%
\[
m\left(  \left(  B_{\mathbf{r}}\left(  \left(  A_{\mathbf{p}}\right)
_{\mathbf{p}\in\mathbf{P}}\right)  \right)  _{\mathbf{r}\in\mathbf{P}}\right)
=\left(  \prod\limits_{\mathbf{q}\in\mathbf{P}}x_{\mathbf{q}}^{n_{\mathbf{q}}%
}\right)  \left(  \left(  B_{\mathbf{r}}\left(  \left(  A_{\mathbf{p}}\right)
_{\mathbf{p}\in\mathbf{P}}\right)  \right)  _{\mathbf{r}\in\mathbf{P}}\right)
=\prod\limits_{\mathbf{q}\in\mathbf{P}}\left(  B_{\mathbf{q}}\left(  \left(
A_{\mathbf{p}}\right)  _{\mathbf{p}\in\mathbf{P}}\right)  \right)
^{n_{\mathbf{q}}}.
\]
Compared with (\ref{pf.algebraic.triangularity.substass.1.pf.a}), this yields
$\left(  m\left(  \left(  B_{\mathbf{r}}\right)  _{\mathbf{r}\in\mathbf{P}%
}\right)  \right)  \left(  \left(  A_{\mathbf{p}}\right)  _{\mathbf{p}%
\in\mathbf{P}}\right)  =m\left(  \left(  B_{\mathbf{r}}\left(  \left(
A_{\mathbf{p}}\right)  _{\mathbf{p}\in\mathbf{P}}\right)  \right)
_{\mathbf{r}\in\mathbf{P}}\right)  $. Thus,
(\ref{pf.algebraic.triangularity.substass.1b}) is proven.} Using these facts,
it is easy to see the following: For every polynomial $g\in\mathbb{F}\left[
x_{\mathbf{P}}\right]  $, we have%
\begin{equation}
\left(  \left(  g\left(  \left(  B_{\mathbf{r}}\right)  _{\mathbf{r}%
\in\mathbf{P}}\right)  \right)  \left(  \left(  A_{\mathbf{p}}\right)
_{\mathbf{p}\in\mathbf{P}}\right)  \text{ is well-defined}\right)
\label{pf.algebraic.triangularity.substass.2a}%
\end{equation}
and%
\begin{equation}
\left(  g\left(  \left(  B_{\mathbf{r}}\right)  _{\mathbf{r}\in\mathbf{P}%
}\right)  \right)  \left(  \left(  A_{\mathbf{p}}\right)  _{\mathbf{p}%
\in\mathbf{P}}\right)  =g\left(  \left(  B_{\mathbf{r}}\left(  \left(
A_{\mathbf{p}}\right)  _{\mathbf{p}\in\mathbf{P}}\right)  \right)
_{\mathbf{r}\in\mathbf{P}}\right)  .
\label{pf.algebraic.triangularity.substass.2b}%
\end{equation}
\footnote{\textit{Proof of (\ref{pf.algebraic.triangularity.substass.2a}) and
(\ref{pf.algebraic.triangularity.substass.2b}):} Let $g\in\mathbb{F}\left[
x_{\mathbf{P}}\right]  $ be a polynomial. Then, $g$ is a $\mathbb{F}$-linear
combination of monomials in $\mathbb{F}\left[  x_{\mathbf{P}}\right]  $ (since
every polynomial is an $\mathbb{F}$-linear combination of monomials). In other
words, there exists some $n\in\mathbb{N}$, some elements $\lambda_{1}$,
$\lambda_{2}$, $...$, $\lambda_{n}$ of $\mathbb{F}$ and some monomials $m_{1}%
$, $m_{2}$, $...$, $m_{n}$ in $\mathbb{F}\left[  x_{\mathbf{P}}\right]  $ such
that $g=\sum\limits_{i=1}^{n}\lambda_{i}m_{i}$. Consider this $n$, these
$\lambda_{1}$, $\lambda_{2}$, $...$, $\lambda_{n}$ and these $m_{1}$, $m_{2}$,
$...$, $m_{n}$.
\par
Notice that $g\left(  \left(  B_{\mathbf{r}}\right)  _{\mathbf{r}\in
\mathbf{P}}\right)  $ and $g\left(  \left(  B_{\mathbf{r}}\left(  \left(
A_{\mathbf{p}}\right)  _{\mathbf{p}\in\mathbf{P}}\right)  \right)
_{\mathbf{r}\in\mathbf{P}}\right)  $ are well-defined (by
(\ref{pf.algebraic.triangularity.substass.wd}) (since $g\in\mathbb{F}\left[
x_{\mathbf{P}}\right]  $)).
\par
For every $i\in\left\{  1,2,...,n\right\}  $, the element $\left(
m_{i}\left(  \left(  B_{\mathbf{r}}\right)  _{\mathbf{r}\in\mathbf{P}}\right)
\right)  \left(  \left(  A_{\mathbf{p}}\right)  _{\mathbf{p}\in\mathbf{P}%
}\right)  $ is well-defined (by (\ref{pf.algebraic.triangularity.substass.1a}%
), applied to $m=m_{i}$).
\par
Since $g=\sum\limits_{i=1}^{n}\lambda_{i}m_{i}$, we have%
\[
g\left(  \left(  B_{\mathbf{r}}\right)  _{\mathbf{r}\in\mathbf{P}}\right)
=\left(  \sum\limits_{i=1}^{n}\lambda_{i}m_{i}\right)  \left(  \left(
B_{\mathbf{r}}\right)  _{\mathbf{r}\in\mathbf{P}}\right)  =\sum\limits_{i=1}%
^{n}\lambda_{i}m_{i}\left(  \left(  B_{\mathbf{r}}\right)  _{\mathbf{r}%
\in\mathbf{P}}\right)  ,
\]
so that%
\begin{align}
\underbrace{\left(  g\left(  \left(  B_{\mathbf{r}}\right)  _{\mathbf{r}%
\in\mathbf{P}}\right)  \right)  }_{=\sum\limits_{i=1}^{n}\lambda_{i}%
m_{i}\left(  \left(  B_{\mathbf{r}}\right)  _{\mathbf{r}\in\mathbf{P}}\right)
}\left(  \left(  A_{\mathbf{p}}\right)  _{\mathbf{p}\in\mathbf{P}}\right)   &
=\left(  \sum\limits_{i=1}^{n}\lambda_{i}m_{i}\left(  \left(  B_{\mathbf{r}%
}\right)  _{\mathbf{r}\in\mathbf{P}}\right)  \right)  \left(  \left(
A_{\mathbf{p}}\right)  _{\mathbf{p}\in\mathbf{P}}\right) \nonumber\\
&  =\sum\limits_{i=1}^{n}\lambda_{i}\underbrace{\left(  m_{i}\left(  \left(
B_{\mathbf{r}}\right)  _{\mathbf{r}\in\mathbf{P}}\right)  \right)  \left(
\left(  A_{\mathbf{p}}\right)  _{\mathbf{p}\in\mathbf{P}}\right)
}_{\substack{=m_{i}\left(  \left(  B_{\mathbf{r}}\left(  \left(
A_{\mathbf{p}}\right)  _{\mathbf{p}\in\mathbf{P}}\right)  \right)
_{\mathbf{r}\in\mathbf{P}}\right)  \\\text{(by
(\ref{pf.algebraic.triangularity.substass.1b}), applied to }m=m_{i}\text{)}%
}}\nonumber\\
&  =\sum\limits_{i=1}^{n}\lambda_{i}m_{i}\left(  \left(  B_{\mathbf{r}}\left(
\left(  A_{\mathbf{p}}\right)  _{\mathbf{p}\in\mathbf{P}}\right)  \right)
_{\mathbf{r}\in\mathbf{P}}\right)  .
\label{pf.algebraic.triangularity.substass.2.pf.a}%
\end{align}
This computation, in particular, shows that $\left(  g\left(  \left(
B_{\mathbf{r}}\right)  _{\mathbf{r}\in\mathbf{P}}\right)  \right)  \left(
\left(  A_{\mathbf{p}}\right)  _{\mathbf{p}\in\mathbf{P}}\right)  $ is
well-defined. Thus, (\ref{pf.algebraic.triangularity.substass.2a}) is proven.
\par
But since $g=\sum\limits_{i=1}^{n}\lambda_{i}m_{i}$, we have%
\begin{align*}
g\left(  \left(  B_{\mathbf{r}}\left(  \left(  A_{\mathbf{p}}\right)
_{\mathbf{p}\in\mathbf{P}}\right)  \right)  _{\mathbf{r}\in\mathbf{P}}\right)
&  =\left(  \sum\limits_{i=1}^{n}\lambda_{i}m_{i}\right)  \left(  \left(
B_{\mathbf{r}}\left(  \left(  A_{\mathbf{p}}\right)  _{\mathbf{p}\in
\mathbf{P}}\right)  \right)  _{\mathbf{r}\in\mathbf{P}}\right) \\
&  =\sum\limits_{i=1}^{n}\lambda_{i}m_{i}\left(  \left(  B_{\mathbf{r}}\left(
\left(  A_{\mathbf{p}}\right)  _{\mathbf{p}\in\mathbf{P}}\right)  \right)
_{\mathbf{r}\in\mathbf{P}}\right)  .
\end{align*}
Compared with (\ref{pf.algebraic.triangularity.substass.2.pf.a}), this yields
$\left(  g\left(  \left(  B_{\mathbf{r}}\right)  _{\mathbf{r}\in\mathbf{P}%
}\right)  \right)  \left(  \left(  A_{\mathbf{p}}\right)  _{\mathbf{p}%
\in\mathbf{P}}\right)  =g\left(  \left(  B_{\mathbf{r}}\left(  \left(
A_{\mathbf{p}}\right)  _{\mathbf{p}\in\mathbf{P}}\right)  \right)
_{\mathbf{r}\in\mathbf{P}}\right)  $. Thus,
(\ref{pf.algebraic.triangularity.substass.2b}) is proven.} Now, we know that
the element $f\left(  \left(  B_{\mathbf{r}}\left(  \left(  A_{\mathbf{p}%
}\right)  _{\mathbf{p}\in\mathbf{P}}\right)  \right)  _{\mathbf{r}%
\in\mathbf{P}}\right)  $ is well-defined. In other words, substituting
$B_{\mathbf{r}}\left(  \left(  A_{\mathbf{p}}\right)  _{\mathbf{p}%
\in\mathbf{P}}\right)  $ for $x_{\mathbf{r}}$ for every $\mathbf{r}%
\in\mathbf{P}$ into $f$ does not render the denominator zero. In other words,
there exists a polynomial $g\in\mathbb{F}\left[  x_{\mathbf{P}}\right]  $ and
a nonzero polynomial $h\in\mathbb{F}\left[  x_{\mathbf{P}}\right]  $ such that
$f=\dfrac{g}{h}$ and $h\left(  \left(  B_{\mathbf{r}}\left(  \left(
A_{\mathbf{p}}\right)  _{\mathbf{p}\in\mathbf{P}}\right)  \right)
_{\mathbf{r}\in\mathbf{P}}\right)  \neq0$. Consider these $g$ and $h$.
Applying (\ref{pf.algebraic.triangularity.substass.2b}) to $h$ instead of $g$,
we obtain%
\begin{equation}
\left(  h\left(  \left(  B_{\mathbf{r}}\right)  _{\mathbf{r}\in\mathbf{P}%
}\right)  \right)  \left(  \left(  A_{\mathbf{p}}\right)  _{\mathbf{p}%
\in\mathbf{P}}\right)  =h\left(  \left(  B_{\mathbf{r}}\left(  \left(
A_{\mathbf{p}}\right)  _{\mathbf{p}\in\mathbf{P}}\right)  \right)
_{\mathbf{r}\in\mathbf{P}}\right)  .
\label{pf.algebraic.triangularity.substass.3}%
\end{equation}
Since $h\left(  \left(  B_{\mathbf{r}}\left(  \left(  A_{\mathbf{p}}\right)
_{\mathbf{p}\in\mathbf{P}}\right)  \right)  _{\mathbf{r}\in\mathbf{P}}\right)
\neq0$, this yields%
\begin{equation}
\left(  h\left(  \left(  B_{\mathbf{r}}\right)  _{\mathbf{r}\in\mathbf{P}%
}\right)  \right)  \left(  \left(  A_{\mathbf{p}}\right)  _{\mathbf{p}%
\in\mathbf{P}}\right)  =h\left(  \left(  B_{\mathbf{r}}\left(  \left(
A_{\mathbf{p}}\right)  _{\mathbf{p}\in\mathbf{P}}\right)  \right)
_{\mathbf{r}\in\mathbf{P}}\right)  \neq0.
\label{pf.algebraic.triangularity.substass.3o}%
\end{equation}
Now, $h\left(  \left(  B_{\mathbf{r}}\right)  _{\mathbf{r}\in\mathbf{P}%
}\right)  \neq0$ (because otherwise, we would have $h\left(  \left(
B_{\mathbf{r}}\right)  _{\mathbf{r}\in\mathbf{P}}\right)  =0$, so that
$\underbrace{\left(  h\left(  \left(  B_{\mathbf{r}}\right)  _{\mathbf{r}%
\in\mathbf{P}}\right)  \right)  }_{=0}\left(  \left(  A_{\mathbf{p}}\right)
_{\mathbf{p}\in\mathbf{P}}\right)  =0\left(  \left(  A_{\mathbf{p}}\right)
_{\mathbf{p}\in\mathbf{P}}\right)  =0$, contradicting
(\ref{pf.algebraic.triangularity.substass.3o})). Thus, the substitution of
$B_{\mathbf{r}}$ for $x_{\mathbf{r}}$ for every $\mathbf{r}\in\mathbf{P}$ into
the rational function $\dfrac{g}{h}$ does not render the denominator zero.
Hence, the value $\dfrac{g}{h}\left(  \left(  B_{\mathbf{r}}\right)
_{\mathbf{r}\in\mathbf{P}}\right)  $ is well-defined. Since $\dfrac{g}{h}=f$,
this rewrites as follows: The value $f\left(  \left(  B_{\mathbf{r}}\right)
_{\mathbf{r}\in\mathbf{P}}\right)  $ is well-defined. Furthermore,
$f=\dfrac{g}{h}$, so that $f\left(  \left(  B_{\mathbf{r}}\right)
_{\mathbf{r}\in\mathbf{P}}\right)  =\dfrac{g}{h}\left(  \left(  B_{\mathbf{r}%
}\right)  _{\mathbf{r}\in\mathbf{P}}\right)  =\dfrac{g\left(  \left(
B_{\mathbf{r}}\right)  _{\mathbf{r}\in\mathbf{P}}\right)  }{h\left(  \left(
B_{\mathbf{r}}\right)  _{\mathbf{r}\in\mathbf{P}}\right)  }$. The substitution
of $A_{\mathbf{p}}$ for $x_{\mathbf{p}}$ for every $\mathbf{p}\in\mathbf{P}$
into the rational function $\dfrac{g\left(  \left(  B_{\mathbf{r}}\right)
_{\mathbf{r}\in\mathbf{P}}\right)  }{h\left(  \left(  B_{\mathbf{r}}\right)
_{\mathbf{r}\in\mathbf{P}}\right)  }$ does not render the denominator zero
(because (\ref{pf.algebraic.triangularity.substass.3o}) yields $\left(
h\left(  \left(  B_{\mathbf{r}}\right)  _{\mathbf{r}\in\mathbf{P}}\right)
\right)  \left(  \left(  A_{\mathbf{p}}\right)  _{\mathbf{p}\in\mathbf{P}%
}\right)  \neq0$). Hence, the value $\dfrac{g\left(  \left(  B_{\mathbf{r}%
}\right)  _{\mathbf{r}\in\mathbf{P}}\right)  }{h\left(  \left(  B_{\mathbf{r}%
}\right)  _{\mathbf{r}\in\mathbf{P}}\right)  }\left(  \left(  A_{\mathbf{p}%
}\right)  _{\mathbf{p}\in\mathbf{P}}\right)  $ is well-defined. Since
$\dfrac{g\left(  \left(  B_{\mathbf{r}}\right)  _{\mathbf{r}\in\mathbf{P}%
}\right)  }{h\left(  \left(  B_{\mathbf{r}}\right)  _{\mathbf{r}\in\mathbf{P}%
}\right)  }=f\left(  \left(  B_{\mathbf{r}}\right)  _{\mathbf{r}\in\mathbf{P}%
}\right)  $, this rewrites as follows: The value $\left(  f\left(  \left(
B_{\mathbf{r}}\right)  _{\mathbf{r}\in\mathbf{P}}\right)  \right)  \left(
\left(  A_{\mathbf{p}}\right)  _{\mathbf{p}\in\mathbf{P}}\right)  $ is
well-defined. Since both sides of the equality
(\ref{pf.algebraic.triangularity.substass.3}) are nonzero (because $h\left(
\left(  B_{\mathbf{r}}\left(  \left(  A_{\mathbf{p}}\right)  _{\mathbf{p}%
\in\mathbf{P}}\right)  \right)  _{\mathbf{r}\in\mathbf{P}}\right)  \neq0$), we
can divide the equality (\ref{pf.algebraic.triangularity.substass.2b}) by
(\ref{pf.algebraic.triangularity.substass.3}). As a result, we obtain%
\begin{equation}
\dfrac{\left(  g\left(  \left(  B_{\mathbf{r}}\right)  _{\mathbf{r}%
\in\mathbf{P}}\right)  \right)  \left(  \left(  A_{\mathbf{p}}\right)
_{\mathbf{p}\in\mathbf{P}}\right)  }{\left(  h\left(  \left(  B_{\mathbf{r}%
}\right)  _{\mathbf{r}\in\mathbf{P}}\right)  \right)  \left(  \left(
A_{\mathbf{p}}\right)  _{\mathbf{p}\in\mathbf{P}}\right)  }=\dfrac{g\left(
\left(  B_{\mathbf{r}}\left(  \left(  A_{\mathbf{p}}\right)  _{\mathbf{p}%
\in\mathbf{P}}\right)  \right)  _{\mathbf{r}\in\mathbf{P}}\right)  }{h\left(
\left(  B_{\mathbf{r}}\left(  \left(  A_{\mathbf{p}}\right)  _{\mathbf{p}%
\in\mathbf{P}}\right)  \right)  _{\mathbf{r}\in\mathbf{P}}\right)  }.
\label{pf.algebraic.triangularity.substass.4}%
\end{equation}
Since $f=\dfrac{g}{h}$, we have%
\begin{align*}
f\left(  \left(  B_{\mathbf{r}}\left(  \left(  A_{\mathbf{p}}\right)
_{\mathbf{p}\in\mathbf{P}}\right)  \right)  _{\mathbf{r}\in\mathbf{P}}\right)
&  =\dfrac{g}{h}\left(  \left(  B_{\mathbf{r}}\left(  \left(  A_{\mathbf{p}%
}\right)  _{\mathbf{p}\in\mathbf{P}}\right)  \right)  _{\mathbf{r}%
\in\mathbf{P}}\right)  =\dfrac{g\left(  \left(  B_{\mathbf{r}}\left(  \left(
A_{\mathbf{p}}\right)  _{\mathbf{p}\in\mathbf{P}}\right)  \right)
_{\mathbf{r}\in\mathbf{P}}\right)  }{h\left(  \left(  B_{\mathbf{r}}\left(
\left(  A_{\mathbf{p}}\right)  _{\mathbf{p}\in\mathbf{P}}\right)  \right)
_{\mathbf{r}\in\mathbf{P}}\right)  }\\
&  =\dfrac{\left(  g\left(  \left(  B_{\mathbf{r}}\right)  _{\mathbf{r}%
\in\mathbf{P}}\right)  \right)  \left(  \left(  A_{\mathbf{p}}\right)
_{\mathbf{p}\in\mathbf{P}}\right)  }{\left(  h\left(  \left(  B_{\mathbf{r}%
}\right)  _{\mathbf{r}\in\mathbf{P}}\right)  \right)  \left(  \left(
A_{\mathbf{p}}\right)  _{\mathbf{p}\in\mathbf{P}}\right)  }%
\ \ \ \ \ \ \ \ \ \ \left(  \text{by
(\ref{pf.algebraic.triangularity.substass.4})}\right)  .
\end{align*}
Compared with%
\begin{align*}
\underbrace{\left(  f\left(  \left(  B_{\mathbf{r}}\right)  _{\mathbf{r}%
\in\mathbf{P}}\right)  \right)  }_{=\dfrac{g\left(  \left(  B_{\mathbf{r}%
}\right)  _{\mathbf{r}\in\mathbf{P}}\right)  }{h\left(  \left(  B_{\mathbf{r}%
}\right)  _{\mathbf{r}\in\mathbf{P}}\right)  }}\left(  \left(  A_{\mathbf{p}%
}\right)  _{\mathbf{p}\in\mathbf{P}}\right)   &  =\dfrac{g\left(  \left(
B_{\mathbf{r}}\right)  _{\mathbf{r}\in\mathbf{P}}\right)  }{h\left(  \left(
B_{\mathbf{r}}\right)  _{\mathbf{r}\in\mathbf{P}}\right)  }\left(  \left(
A_{\mathbf{p}}\right)  _{\mathbf{p}\in\mathbf{P}}\right) \\
&  =\dfrac{\left(  g\left(  \left(  B_{\mathbf{r}}\right)  _{\mathbf{r}%
\in\mathbf{P}}\right)  \right)  \left(  \left(  A_{\mathbf{p}}\right)
_{\mathbf{p}\in\mathbf{P}}\right)  }{\left(  h\left(  \left(  B_{\mathbf{r}%
}\right)  _{\mathbf{r}\in\mathbf{P}}\right)  \right)  \left(  \left(
A_{\mathbf{p}}\right)  _{\mathbf{p}\in\mathbf{P}}\right)  },
\end{align*}
this yields $\left(  f\left(  \left(  B_{\mathbf{r}}\right)  _{\mathbf{r}%
\in\mathbf{P}}\right)  \right)  \left(  \left(  A_{\mathbf{p}}\right)
_{\mathbf{p}\in\mathbf{P}}\right)  =f\left(  \left(  B_{\mathbf{r}}\left(
\left(  A_{\mathbf{p}}\right)  _{\mathbf{p}\in\mathbf{P}}\right)  \right)
_{\mathbf{r}\in\mathbf{P}}\right)  $. This finishes the proof of Lemma
\ref{lem.algebraic.triangularity.substass}.
\end{proof}

We make three further conventions which will be in place for the rest of
Section \ref{sect.dominance}:

\begin{convention}
\label{conv.attach}Let $p\in\mathbb{N}$. If $A_{1}$, $A_{2}$, $...$, $A_{k}$
are several matrices with $p$ rows each, then $\left(  A_{1}\mid A_{2}%
\mid...\mid A_{k}\right)  $ will denote the matrix obtained by starting with
an (empty) $p\times0$-matrix, then attaching the matrix $A_{1}$ to it on the
right, then attaching the matrix $A_{2}$ to the result on the right, etc., and
finally attaching the matrix $A_{k}$ to the result on the right. In other
words, the matrix $\left(  A_{1}\mid A_{2}\mid...\mid A_{k}\right)  $ is
defined as the matrix whose columns (from left to right) are all the columns
of $A_{1}$ (from left to right), then all the columns of $A_{2}$ (from left to
right), etc., and finally all the columns of $A_{k}$ (from left to right).
Hence, if $A_{i}$ is a $p\times q_{i}$-matrix for every $i\in\left\{
1,2,...,k\right\}  $, then, for every $i\in\left\{  1,2,...,k\right\}  $ and
every $j\in\left\{  1,2,...,q_{i}\right\}  $, we have%
\begin{align}
&  \left(  \text{the }\left(  \left(  q_{1}+q_{2}+...+q_{i-1}\right)
+j\right)  \text{-th column of the matrix }\left(  A_{1}\mid A_{2}\mid...\mid
A_{k}\right)  \right) \nonumber\\
&  =\left(  \text{the }j\text{-th column of the matrix }A_{i}\right)  .
\label{conv.attach.formal}%
\end{align}

For example, if $p$ is a nonnegative integer, and $B$ is a matrix with $p$
rows, then $\left(  I_{p}\mid B\right)  $ means the matrix obtained from the
$p\times p$ identity matrix $I_{p}$ by attaching the matrix $B$ to it on the
right. (As a concrete example, $\left(  I_{2}\mid\left(
\begin{array}
[c]{cc}%
1 & -2\\
3 & 0
\end{array}
\right)  \right)  =\left(
\begin{array}
[c]{cccc}%
1 & 0 & 1 & -2\\
0 & 1 & 3 & 0
\end{array}
\right)  $.)
\end{convention}

\begin{convention}
\label{conv.rows} Let $p\in\mathbb{N}$ and $q\in\mathbb{N}$. Let $B$ be a
$p\times q$-matrix, and let $i_{1}$, $i_{2}$, $...$, $i_{k}$ be some elements
of $\left\{  1,2,...,p\right\}  $. Then, $\operatorname*{rows}\nolimits_{i_{1}%
,i_{2},...,i_{k}}B$ will denote the matrix whose rows (from top to bottom) are
the rows labelled $i_{1}$, $i_{2}$, $...$, $i_{k}$ of the matrix $B$. In other
words, if we write the matrix $B$ in the form $\left(  b_{u,v}\right)  _{1\leq
u\leq p,\ 1\leq v\leq q}$, then%
\begin{equation}
\operatorname*{rows}\nolimits_{i_{1},i_{2},...,i_{k}}B=\left(  b_{i_{u}%
,v}\right)  _{1\leq u\leq k,\ 1\leq v\leq q}. \label{conv.rows.formal}%
\end{equation}

\end{convention}

\begin{convention}
\label{conv.cols} Let $p\in\mathbb{N}$ and $q\in\mathbb{N}$. Let $B$ be a
$p\times q$-matrix, and let $j_{1}$, $j_{2}$, $...$, $j_{\ell}$ be some
elements of $\left\{  1,2,...,q\right\}  $. Then, $\operatorname*{cols}%
\nolimits_{j_{1},j_{2},...,j_{\ell}}B$ will denote the matrix whose columns
(from left to right) are the columns labelled $j_{1}$, $j_{2}$, $...$,
$j_{\ell}$ of the matrix $B$. In other words, if we write the matrix $B$ in
the form $\left(  b_{u,v}\right)  _{1\leq u\leq p,\ 1\leq v\leq q}$, then%
\begin{equation}
\operatorname*{cols}\nolimits_{j_{1},j_{2},...,j_{\ell}}B=\left(  b_{u,j_{v}%
}\right)  _{1\leq u\leq p,\ 1\leq v\leq\ell}. \label{conv.cols.formal}%
\end{equation}

\end{convention}

We notice some trivial properties of matrices:

\begin{proposition}
\label{prop.attach.cols}Let $\mathbb{K}$ be a commutative ring. Let $u$ and
$v$ be nonnegative integers. Let $A\in\mathbb{K}^{u\times v}$. Let $a$ and $b$
be integers such that $1\leq a\leq b\leq v+1$. Then,%
\[
A\left[  a:b\right]  =\operatorname*{cols}\nolimits_{a,a+1,...,b-1}A.
\]

\end{proposition}

\begin{proposition}
\label{prop.cols.empty}Let $\mathbb{K}$ be a commutative ring. Let $u$ and $v$
be nonnegative integers. Let $A\in\mathbb{K}^{u\times v}$. Let $a$ be an
integer. Then, $A\left[  a:a\right]  =\left(  \text{empty matrix}\right)  $.
\end{proposition}

\begin{proposition}
\label{prop.attach.matrices}Let $\mathbb{K}$ be a commutative ring. Let $u$
and $v$ be nonnegative integers. Let $A\in\mathbb{K}^{u\times v}$. Let $a$,
$b$, $c$ and $d$ be integers satisfying $a\leq b$ and $c\leq d$. Then,%
\[
A\left[  a:b\mid c:d\right]  =\left(  A\left[  a:b\right]  \mid A\left[
c:d\right]  \right)  .
\]
(Here, we are using the notations from Definition \ref{def.minors} and
Definition \ref{def.minors.simple}.)
\end{proposition}

\begin{proposition}
\label{prop.attach.rows}Let $\mathbb{K}$ be a commutative ring. Let $u$,
$v_{1}$ and $v_{2}$ be nonnegative integers. Let $A_{1}\in\mathbb{K}^{u\times
v_{1}}$ and $A_{2}\in\mathbb{K}^{u\times v_{2}}$. Let $i_{1}$, $i_{2}$, $...$,
$i_{k}$ be some elements of $\left\{  1,2,...,u\right\}  $. Then,%
\[
\operatorname*{rows}\nolimits_{i_{1},i_{2},...,i_{k}}\left(  A_{1}\mid
A_{2}\right)  =\left(  \operatorname*{rows}\nolimits_{i_{1},i_{2},...,i_{k}%
}\left(  A_{1}\right)  \mid\operatorname*{rows}\nolimits_{i_{1},i_{2}%
,...,i_{k}}\left(  A_{2}\right)  \right)  .
\]

\end{proposition}

\begin{proposition}
\label{prop.rows.entry}Let us regard column vectors of size $\ell$ as
$\ell\times1$-matrices for every $\ell\in\mathbb{N}$. Let us also regard
scalars as $1\times1$-matrices.

Let $\ell\in\mathbb{N}$. Let $\mathbb{K}$ be a commutative ring. Let
$U\in\mathbb{K}^{\ell}$ be a column vector. Let $i\in\left\{  1,2,...,\ell
\right\}  $. Then,%
\[
\operatorname*{rows}\nolimits_{i}U=\left(  \text{the }i\text{-th entry of
}U\right)  .
\]

\end{proposition}

\begin{proposition}
\label{prop.rows.cols}Let $\mathbb{K}$ be a commutative ring. Let $p$ and $q$
be nonnegative integers. Let $B\in\mathbb{K}^{p\times q}$. Let $i_{1}$,
$i_{2}$, $...$, $i_{k}$ be some elements of $\left\{  1,2,...,p\right\}  $.
Let $j_{1}$, $j_{2}$, $...$, $j_{\ell}$ be some elements of $\left\{
1,2,...,q\right\}  $. Then,%
\[
\operatorname*{cols}\nolimits_{j_{1},j_{2},...,j_{\ell}}\left(
\operatorname*{rows}\nolimits_{i_{1},i_{2},...,i_{k}}B\right)
=\operatorname*{rows}\nolimits_{i_{1},i_{2},...,i_{k}}\left(
\operatorname*{cols}\nolimits_{j_{1},j_{2},...,j_{\ell}}B\right)  .
\]

\end{proposition}

\begin{proposition}
\label{prop.cols.interval}Let $\mathbb{K}$ be a commutative ring. Let $p$ and
$q$ be nonnegative integers. Let $A\in\mathbb{K}^{p\times q}$. Let $a$ and $b$
be integers such that $a<b$. Then,%
\[
\operatorname*{cols}\nolimits_{1,2,...,b-a-1}\left(  A\left[  a:b\right]
\right)  =A\left[  a:b-1\right]  .
\]

\end{proposition}

\begin{proposition}
\label{prop.rows.rows}Let $\mathbb{K}$ be a commutative ring. Let $p$ and $q$
be nonnegative integers. Let $B\in\mathbb{K}^{p\times q}$. Let $i_{1}$,
$i_{2}$, $...$, $i_{k}$ be some elements of $\left\{  1,2,...,p\right\}  $
such that $k>0$. Then,%
\[
\operatorname*{rows}\nolimits_{2,3,...,k}\left(  \operatorname*{rows}%
\nolimits_{i_{1},i_{2},...,i_{k}}B\right)  =\operatorname*{rows}%
\nolimits_{i_{2},i_{3},...,i_{k}}B.
\]

\end{proposition}

\begin{proposition}
\label{prop.rows.row}Let $\mathbb{K}$ be a commutative ring. Let $p$ and $q$
be nonnegative integers. Let $B\in\mathbb{K}^{p\times q}$. Let $i_{1}$,
$i_{2}$, $...$, $i_{k}$ be some elements of $\left\{  1,2,...,p\right\}  $.
Let $j\in\left\{  1,2,...,k\right\}  $. Then,%
\[
\operatorname*{rows}\nolimits_{j}\left(  \operatorname*{rows}\nolimits_{i_{1}%
,i_{2},...,i_{k}}B\right)  =\operatorname*{rows}\nolimits_{i_{j}}B.
\]

\end{proposition}

\begin{proposition}
\label{prop.cols.attach1}Let $\mathbb{K}$ be a commutative ring. Let $u$,
$v_{1}$ and $v_{2}$ be nonnegative integers. Let $A_{1}\in\mathbb{K}^{u\times
v_{1}}$ and $A_{2}\in\mathbb{K}^{u\times v_{2}}$. Let $\ell\in\mathbb{N}$ be
such that $v_{1}\leq\ell\leq v_{1}+v_{2}$. Then,%
\[
\operatorname*{cols}\nolimits_{1,2,...,\ell}\left(  A_{1}\mid A_{2}\right)
=\left(  A_{1}\mid\operatorname*{cols}\nolimits_{1,2,...,\ell-v_{1}}\left(
A_{2}\right)  \right)  .
\]

\end{proposition}

\begin{corollary}
\label{cor.cols.attach1a}Let $\mathbb{K}$ be a commutative ring. Let $u$,
$v_{1}$ and $v_{2}$ be nonnegative integers. Let $A_{1}\in\mathbb{K}^{u\times
v_{1}}$ and $A_{2}\in\mathbb{K}^{u\times v_{2}}$. Then,%
\[
\operatorname*{cols}\nolimits_{1,2,...,v_{1}}\left(  A_{1}\mid A_{2}\right)
=A_{1}.
\]

\end{corollary}

\begin{proposition}
\label{prop.attach2.entry}Let $\mathbb{K}$ be a commutative ring. Let $u$,
$v_{1}$ and $v_{2}$ be nonnegative integers. Let $A_{1}\in\mathbb{K}^{u\times
v_{1}}$ and $A_{2}\in\mathbb{K}^{u\times v_{2}}$. Let $\ell\in\left\{
1,2,...,v_{2}\right\}  $. Let $i\in\left\{  1,2,...,u\right\}  $. Then, the
vector $\left(  A_{1}\mid A_{2}\right)  _{v_{1}+\ell}$ satisfies%
\[
\left(  \text{the }i\text{-th entry of }\left(  A_{1}\mid A_{2}\right)
_{v_{1}+\ell}\right)  =\left(  \text{the }\left(  i,\ell\right)  \text{-th
entry of the matrix }A_{2}\right)  .
\]

\end{proposition}

\begin{corollary}
\label{cor.cutoffI}Let $\mathbb{K}$ be a commutative ring. Let $u$, $v_{1}$
and $v_{2}$ be nonnegative integers. Let $A_{1}\in\mathbb{K}^{u\times v_{1}}$
and $A_{2}\in\mathbb{K}^{u\times v_{2}}$.

\textbf{(a)} Every $r\in\left\{  0,1,...,v_{1}\right\}  $ satisfies $\left(
A_{1}\mid A_{2}\right)  \left[  r+1:v_{1}+1\right]  =\operatorname*{cols}%
\nolimits_{r+1,r+2,...,v_{1}}\left(  A_{1}\right)  $.

\textbf{(b)} Every $i\in\left\{  1,2,...,v_{1}+1\right\}  $ satisfies $\left(
A_{1}\mid A_{2}\right)  \left[  1:i\right]  =\operatorname*{cols}%
\nolimits_{1,2,...,i-1}\left(  A_{1}\right)  $.
\end{corollary}

\begin{corollary}
\label{cor.cofactor.1}Let $\mathbb{K}$ be a commutative ring. Let $p$ and $q$
be nonnegative integers. Let $A\in\mathbb{K}^{p\times q}$. Let $a$ and $b$ be
integers such that $a<b$ and $b-a\leq p$. Let $i=p-\left(  b-a\right)  +1$.
Then,%
\begin{align*}
&  \left(  \text{the }\left(  1,p-i+1\right)  \text{-th cofactor of the matrix
}\operatorname*{rows}\nolimits_{i,i+1,...,p}\left(  A\left[  a:b\right]
\right)  \right) \\
&  =\left(  -1\right)  ^{p-i}\det\left(  \operatorname*{rows}%
\nolimits_{i+1,i+2,...,p}\left(  A\left[  a:b-1\right]  \right)  \right)  .
\end{align*}

\end{corollary}

\begin{corollary}
\label{cor.cofactor.2}Let $\mathbb{K}$ be a commutative ring. Let $p$ and $q$
be nonnegative integers. Let $A\in\mathbb{K}^{p\times q}$. Let $a$ and $b$ be
integers such that $a<b$ and $b-a\leq p$. Let $i=p-\left(  b-a\right)  $. Let
$\xi$ be a column vector of size $p$. We regard column vectors of size $\ell$
as $\ell\times1$-matrices for every $\ell\in\mathbb{N}$. Then,%
\begin{align*}
&  \left(  \text{the }\left(  1,p-i+1\right)  \text{-th cofactor of the matrix
}\operatorname*{rows}\nolimits_{i,i+1,...,p}\left(  \xi\mid A\left[
a:b\right]  \right)  \right) \\
&  =\left(  -1\right)  ^{p-i}\det\left(  \operatorname*{rows}%
\nolimits_{i+1,i+2,...,p}\left(  \xi\mid A\left[  a:b-1\right]  \right)
\right)  .
\end{align*}

\end{corollary}

And here is a simple property of determinants of attached matrices:

\begin{lemma}
\label{lem.attach.matrices}Let $\mathbb{K}$ be a commutative ring. Let $p$ and
$k$ be nonnegative integers such that $p\geq k$. Let $U\in\mathbb{K}^{p\times
k}$ be an $p\times k$-matrix. For any $u\in\mathbb{N}$ and $v\in\mathbb{N}$,
let $0_{u\times v}$ denote the $u\times v$ zero matrix.

\textbf{(a)} We have $\det\left(  \left(
\begin{array}
[c]{c}%
I_{p-k}\\
0_{k\times\left(  p-k\right)  }%
\end{array}
\right)  \mid U\right)  =\det\left(  \operatorname*{rows}%
\nolimits_{p-k+1,p-k+2,...,p}U\right)  $ (where $\left(
\begin{array}
[c]{c}%
I_{p-k}\\
0_{k\times\left(  p-k\right)  }%
\end{array}
\right)  $ is to be understood as a block matrix).

\textbf{(b)} We have $\det\left(  \left(
\begin{array}
[c]{c}%
0_{k\times\left(  p-k\right)  }\\
I_{p-k}%
\end{array}
\right)  \mid U\right)  =\left(  -1\right)  ^{k\left(  p-k\right)  }%
\det\left(  \operatorname*{rows}\nolimits_{1,2,...,k}U\right)  $ (where
$\left(
\begin{array}
[c]{c}%
0_{k\times\left(  p-k\right)  }\\
I_{p-k}%
\end{array}
\right)  $ is to be understood as a block matrix).

\textbf{(c)} Let $u_{1}$, $u_{2}$, $...$, $u_{p-k}$ be $p-k$ pairwise distinct
elements of $\left\{  1,2,...,p\right\}  $. For every $j\in\left\{
1,2,...,p\right\}  $, let $e_{j}$ be the vector in $\mathbb{K}^{p}$ whose
$j$-th coordinate is $1$ and whose other coordinates are all $0$. Let
$P\in\mathbb{K}^{p\times\left(  p-k\right)  }$ be the $p\times\left(
p-k\right)  $-matrix whose columns are $e_{u_{1}}$, $e_{u_{2}}$, $...$,
$e_{u_{p-k}}$ from left to right. Let $v_{1}$, $v_{2}$, $...$, $v_{k}$ be the
elements of the set $\left\{  1,2,...,p\right\}  \setminus\left\{  u_{1}%
,u_{2},...,u_{p-k}\right\}  $ listed in any order (but each element only
appearing once). Then, $\det\left(  P\mid U\right)  =\pm\det\left(
\operatorname*{rows}\nolimits_{v_{1},v_{2},...,v_{k}}U\right)  $.
\end{lemma}

All three statements of Lemma \ref{lem.attach.matrices} are standard facts
from the theory of determinants. (They can be all deduced from the fact that
the determinant of a block-triangular matrix is the product of the
determinants of its diagonal blocks. The reduction requires permuting the rows.)

We now show some statements of computational nature:

\begin{lemma}
\label{lem.Grasp.generic.calc} Let $\mathbb{F}$ be a field.

Let $p$ and $q$ be two positive integers. Let $\mathbf{P}$ be a totally
ordered set such that%
\[
\mathbf{P}=\left\{  1,2,...,p\right\}  \times\left\{  1,2,...,q\right\}
\text{ as sets.}%
\]
Let $\vartriangleleft$ denote the smaller relation of $\mathbf{P}$, and let
$\trianglelefteq$ denote the smaller-or-equal relation of $\mathbf{P}$. Assume
that%
\begin{equation}
\left(  i,k\right)  \trianglelefteq\left(  i^{\prime},k^{\prime}\right)
\text{ for all }\left(  i,k\right)  \in\mathbf{P}\text{ and }\left(
i^{\prime},k^{\prime}\right)  \in\mathbf{P}\text{ satisfying }\left(  i\geq
i^{\prime}\text{ and }k\leq k^{\prime}\right)  .
\label{lem.Grasp.generic.calc.ass}%
\end{equation}

Let $Z:\left\{  1,2,...,q\right\}  \rightarrow\left\{  1,2,...,q\right\}  $
denote the map which sends every $k\in\left\{  1,2,...,q-1\right\}  $ to $k+1$
and sends $q$ to $1$. Thus, $Z$ is a permutation in the symmetric group
$S_{q}$, and can be written in cycle notation as $\left(  1,2,...,q\right)  $.

For every subset $\mathbf{S}$ of $\mathbf{P}$, let us denote by $\mathbb{F}%
\left[  x_{\mathbf{S}}\right]  $ the polynomial ring over $\mathbb{F}$ in the
indeterminates $x_{\mathbf{p}}$ with $\mathbf{p}$ ranging over all elements of
$\mathbf{S}$. Clearly, there is a canonical embedding $\mathbb{F}\left[
x_{\mathbf{S}}\right]  \rightarrow\mathbb{F}\left[  x_{\mathbf{P}}\right]  $
for every subset $\mathbf{S}$ of $\mathbf{P}$. We will regard this embedding
as an inclusion.

Define a family $\left(  y_{\mathbf{p}}\right)  _{\mathbf{p}\in\mathbf{P}}%
\in\left(  \mathbb{F}\left[  x_{\mathbf{P}}\right]  \right)  ^{\mathbf{P}}$ of
elements of $\mathbb{F}\left[  x_{\mathbf{P}}\right]  $ by setting
\[
y_{\left(  i,k\right)  }=x_{\left(  i,Z\left(  k\right)  \right)
}\ \ \ \ \ \ \ \ \ \ \text{for all }\left(  i,k\right)  \in\mathbf{P}.
\]

Define a matrix $C\in\left(  \mathbb{F}\left[  x_{\mathbf{P}}\right]  \right)
^{p\times q}$ by
\[
C=\left(  y_{\left(  i,k\right)  }\right)  _{1\leq i\leq p,\ 1\leq k\leq q}.
\]

For every $\left(  i,k\right)  \in\mathbf{P}$, define an element
$\mathfrak{N}_{\left(  i,k\right)  }\in\mathbb{F}\left[  x_{\mathbf{P}%
}\right]  $ by%
\begin{equation}
\mathfrak{N}_{\left(  i,k\right)  }=\det\left(  \left(  I_{p}\mid C\right)
\left[  1:i\mid i+k-1:p+k\right]  \right)  .
\label{lem.Grasp.generic.independency.Ndef}%
\end{equation}

For every $\left(  i,k\right)  \in\mathbf{P}$, define an element
$\mathfrak{D}_{\left(  i,k\right)  }\in\mathbb{F}\left[  x_{\mathbf{P}%
}\right]  $ by%
\begin{equation}
\mathfrak{D}_{\left(  i,k\right)  }=\det\left(  \left(  I_{p}\mid C\right)
\left[  0:i\mid i+k:p+k\right]  \right)  .
\label{lem.Grasp.generic.independency.Ddef}%
\end{equation}

For every $\mathbf{p}\in\mathbf{P}$, let $\mathbf{p}\Downarrow$ denote the
subset $\left\{  \mathbf{v}\in\mathbf{P}\ \mid\ \mathbf{v}\vartriangleleft
\mathbf{p}\right\}  $ of $\mathbf{P}$.

\textbf{(a)} Every $\left(  r,s\right)  \in\mathbf{P}$ satisfies%
\begin{equation}
\mathfrak{N}_{\left(  r,s\right)  }=\det\left(  \operatorname*{rows}%
\nolimits_{r,r+1,...,p}\left(  \left(  I_{p}\mid C\right)  \left[
r+s-1:p+s\right]  \right)  \right)  \label{lem.Grasp.generic.calc.a.N}%
\end{equation}
and%
\begin{equation}
\mathfrak{D}_{\left(  r,s\right)  }=\left(  -1\right)  ^{r-1}\det\left(
\operatorname*{rows}\nolimits_{r,r+1,...,p}\left(  \left(  -1\right)
^{p-1}C_{q}\ \mid\ \left(  I_{p}\mid C\right)  \left[  r+s:p+s\right]
\right)  \right)  . \label{lem.Grasp.generic.calc.a.D}%
\end{equation}
(Here, we treat column vectors of length $\ell$ as $\ell\times1$-matrices for
every $\ell\in\mathbb{N}$. Thus, $\left(  -1\right)  ^{p-1}C_{q}$ is
considered as a $p\times1$-matrix.)

\textbf{(b)} Every $\left(  i,k\right)  \in\mathbf{P}$ satisfying $i\neq p$
and $k\neq1$ satisfies%
\begin{equation}
\mathfrak{N}_{\left(  i,k\right)  }\in\left(  -1\right)  ^{p-i}\cdot
\mathfrak{N}_{\left(  i+1,k-1\right)  }\cdot x_{\left(  i,k\right)
}+\mathbb{F}\left[  x_{\left(  i,k\right)  \Downarrow}\right]
\label{lem.Grasp.generic.calc.b.N}%
\end{equation}
and%
\begin{equation}
\mathfrak{D}_{\left(  i,k\right)  }\in\left(  -1\right)  ^{p-i+1}%
\cdot\mathfrak{D}_{\left(  i+1,k-1\right)  }\cdot x_{\left(  i,k\right)
}+\mathbb{F}\left[  x_{\left(  i,k\right)  \Downarrow}\right]
\label{lem.Grasp.generic.calc.b.D}%
\end{equation}
and%
\begin{equation}
\mathfrak{N}_{\left(  i+1,k-1\right)  }\in\mathbb{F}\left[  x_{\left(
i,k\right)  \Downarrow}\right]  \label{lem.Grasp.generic.calc.b.N2}%
\end{equation}
and%
\begin{equation}
\mathfrak{D}_{\left(  i+1,k-1\right)  }\in\mathbb{F}\left[  x_{\left(
i,k\right)  \Downarrow}\right]  . \label{lem.Grasp.generic.calc.b.D2}%
\end{equation}

\textbf{(c)} Every $\left(  i,k\right)  \in\mathbf{P}$ satisfying $i\neq p$
and $k\neq1$ satisfies%
\begin{equation}
\mathfrak{N}_{\left(  i,k\right)  }\mathfrak{D}_{\left(  i+1,k-1\right)
}+\mathfrak{D}_{\left(  i,k\right)  }\mathfrak{N}_{\left(  i+1,k-1\right)
}=\mathfrak{D}_{\left(  i,k-1\right)  }\mathfrak{N}_{\left(  i+1,k\right)  }.
\label{lem.Grasp.generic.calc.c}%
\end{equation}

\textbf{(d)} Every $i\in\left\{  1,2,...,p\right\}  $ satisfies%
\begin{equation}
\mathfrak{N}_{\left(  i,1\right)  }=1 \label{lem.Grasp.generic.calc.d.N}%
\end{equation}
and%
\begin{equation}
\mathfrak{D}_{\left(  i,1\right)  }=\left(  -1\right)  ^{i+p}x_{\left(
i,1\right)  }. \label{lem.Grasp.generic.calc.d.D}%
\end{equation}

\textbf{(e)} Every $k\in\left\{  2,3,...,q\right\}  $ satisfies%
\begin{equation}
\mathfrak{N}_{\left(  p,k\right)  }=x_{\left(  p,k\right)  }
\label{lem.Grasp.generic.calc.e.N}%
\end{equation}
and%
\begin{equation}
\mathfrak{D}_{\left(  p,k\right)  }=x_{\left(  p,1\right)  }.
\label{lem.Grasp.generic.calc.e.D}%
\end{equation}

\textbf{(f)} Every $\mathbf{p}\in\mathbf{P}$ satisfies%
\begin{equation}
\mathfrak{N}_{\mathbf{p}}\neq0 \label{lem.Grasp.generic.calc.f.N}%
\end{equation}
and%
\begin{equation}
\mathfrak{D}_{\mathbf{p}}\neq0. \label{lem.Grasp.generic.calc.f.D}%
\end{equation}

\textbf{(g)} For every $\mathbf{p}\in\mathbf{P}$, there exist elements
$\alpha_{\mathbf{p}}$, $\beta_{\mathbf{p}}$, $\gamma_{\mathbf{p}}$ and
$\delta_{\mathbf{p}}$ of $\mathbb{F}\left[  x_{\mathbf{p}\Downarrow}\right]  $
satisfying $\alpha_{\mathbf{p}}\delta_{\mathbf{p}}-\beta_{\mathbf{p}}%
\gamma_{\mathbf{p}}\neq0$, $\mathfrak{N}_{\mathbf{p}}=\alpha_{\mathbf{p}%
}x_{\mathbf{p}}+\beta_{\mathbf{p}}$ and $\mathfrak{D}_{\mathbf{p}}%
=\gamma_{\mathbf{p}}x_{\mathbf{p}}+\delta_{\mathbf{p}}$.
\end{lemma}

\begin{proof}
[Proof of Lemma \ref{lem.Grasp.generic.calc} (sketched).]Recall that
$C=\left(  y_{\left(  i,k\right)  }\right)  _{1\leq i\leq p,\ 1\leq k\leq q}$.
Hence, every $\left(  i,k\right)  \in\left\{  1,2,...,p\right\}
\times\left\{  1,2,...,q\right\}  $ satisfies%
\begin{equation}
\left(  \text{the }\left(  i,k\right)  \text{-th entry of the matrix
}C\right)  =y_{\left(  i,k\right)  }=x_{\left(  i,Z\left(  k\right)  \right)
} \label{pf.Grasp.generic.calc.entry}%
\end{equation}
(by the definition of $y_{\left(  i,k\right)  }$).

For any $u\in\mathbb{N}$ and $v\in\mathbb{N}$, let $0_{u\times v}$ denote the
$u\times v$ zero matrix.

\textbf{(a)} We treat column vectors of length $\ell$ as $\ell\times
1$-matrices. Thus, every $u\in\mathbb{N}$ and $a\in\mathbb{N}$ and every
matrix $A\in\mathbb{F}^{p\times u}$ satisfy $A\left[  a:a+1\right]  =A_{a}$.
Moreover, treating column vectors as matrices allows us to make sense of terms
like $\operatorname*{rows}\nolimits_{i,i+1,...,p}A$ when $A$ is a vector in
$\mathbb{F}^{p}$ and $i$ is an element of $\left\{  1,2,...,p+1\right\}  $.
For example, $\operatorname*{rows}\nolimits_{1,2}\left(
\begin{array}
[c]{c}%
1\\
0\\
0
\end{array}
\right)  =\left(
\begin{array}
[c]{c}%
1\\
0
\end{array}
\right)  $ if $p=3$.

Let $\left(  r,s\right)  \in\mathbf{P}$. Since $\left(  r,s\right)
\in\mathbf{P}=\left\{  1,2,...,p\right\}  \times\left\{  1,2,...,q\right\}  $,
we have $r\in\left\{  1,2,...,p\right\}  $ and $s\in\left\{
1,2,...,q\right\}  $.

We have $r\in\left\{  1,2,...,p\right\}  $, so that $r\leq p\leq p+1$. Since
$r\in\left\{  1,2,...,p\right\}  \subseteq\left\{  1,2,...,p+1\right\}  $, we
can apply Corollary \ref{cor.cutoffI} \textbf{(b)} to $\mathbb{F}\left[
x_{\mathbf{P}}\right]  $, $p$, $p$, $q$, $I_{p}$, $C$ and $r$ instead of
$\mathbb{K}$, $u$, $v_{1}$, $v_{2}$, $A_{1}$, $A_{2}$ and $i$. As a result, we
obtain%
\begin{equation}
\left(  I_{p}\mid C\right)  \left[  1:r\right]  =\operatorname*{cols}%
\nolimits_{1,2,...,r-1}\left(  I_{p}\right)  =\left(
\begin{array}
[c]{c}%
I_{r-1}\\
0_{\left(  p-r+1\right)  \times\left(  r-1\right)  }%
\end{array}
\right)  \label{pf.Grasp.generic.calc.a.Ip1r}%
\end{equation}
(where $\left(
\begin{array}
[c]{c}%
I_{r-1}\\
0_{\left(  p-r+1\right)  \times\left(  r-1\right)  }%
\end{array}
\right)  $ is to be understood as a block matrix). By Proposition
\ref{prop.attach.matrices} (applied to $\mathbb{F}\left[  x_{\mathbf{P}%
}\right]  $, $p$, $p+q$, $\left(  I_{p}\mid C\right)  $, $1$, $r$, $r+s-1$ and
$p+s$ instead of $\mathbb{K}$, $u$, $v$, $A$, $a$, $b$, $c$ and $d$), we have%
\begin{align*}
&  \left(  I_{p}\mid C\right)  \left[  1:r\mid r+s-1:p+s\right] \\
&  =\left(  \underbrace{\left(  I_{p}\mid C\right)  \left[  1:r\right]
}_{\substack{=\left(
\begin{array}
[c]{c}%
I_{r-1}\\
0_{\left(  p-r+1\right)  \times\left(  r-1\right)  }%
\end{array}
\right)  \\\text{(by (\ref{pf.Grasp.generic.calc.a.Ip1r}))}}}\ \mid\ \left(
I_{p}\mid C\right)  \left[  r+s-1:p+s\right]  \right)  .\\
&  =\left(  \left(
\begin{array}
[c]{c}%
I_{r-1}\\
0_{\left(  p-r+1\right)  \times\left(  r-1\right)  }%
\end{array}
\right)  \ \mid\ \left(  I_{p}\mid C\right)  \left[  r+s-1:p+s\right]
\right)  .
\end{align*}
Hence,%
\begin{align}
&  \det\left(  \left(  I_{p}\mid C\right)  \left[  1:r\mid r+s-1:p+s\right]
\right) \nonumber\\
&  =\det\left(  \left(
\begin{array}
[c]{c}%
I_{r-1}\\
0_{\left(  p-r+1\right)  \times\left(  r-1\right)  }%
\end{array}
\right)  \ \mid\ \left(  I_{p}\mid C\right)  \left[  r+s-1:p+s\right]  \right)
\nonumber\\
&  =\det\left(  \operatorname*{rows}\nolimits_{r,r+1,...,p}\left(  \left(
I_{p}\mid C\right)  \left[  r+s-1:p+s\right]  \right)  \right)
\label{pf.Grasp.generic.calc.a.simplN1}%
\end{align}
(by Lemma \ref{lem.attach.matrices} \textbf{(a)}, applied to $p-r+1$,
$\mathbb{F}\left[  x_{\mathbf{P}}\right]  $ and $\left(  I_{p}\mid C\right)
\left[  r+s-1:p+s\right]  $ instead of $k$, $\mathbb{K}$ and $U$). Now,%
\begin{align*}
\mathfrak{N}_{\left(  r,s\right)  }  &  =\det\left(  \left(  I_{p}\mid
C\right)  \left[  1:r\mid r+s-1:p+s\right]  \right)
\ \ \ \ \ \ \ \ \ \ \left(  \text{by the definition of }\mathfrak{N}_{\left(
r,s\right)  }\right) \\
&  =\det\left(  \operatorname*{rows}\nolimits_{r,r+1,...,p}\left(  \left(
I_{p}\mid C\right)  \left[  r+s-1:p+s\right]  \right)  \right)
\end{align*}
(by (\ref{pf.Grasp.generic.calc.a.simplN1})). This proves
(\ref{lem.Grasp.generic.calc.a.N}).

On the other hand, by Proposition \ref{prop.minors.simple.transi}, we have
\begin{align}
\left(  I_{p}\mid C\right)  \left[  0:r\right]   &  =\left(  I_{p}\mid
C\right)  \left[  0:1\mid1:r\right] \nonumber\\
&  =\left(  \underbrace{\left(  I_{p}\mid C\right)  \left[  0:1\right]
}_{\substack{=\left(  I_{p}\mid C\right)  _{0}=\left(  -1\right)
^{p-1}\left(  I_{p}\mid C\right)  _{p+q}\\\text{(by Definition
\ref{def.minors} \textbf{(b)})}}}\mid\underbrace{\left(  I_{p}\mid C\right)
\left[  1:r\right]  }_{\substack{=\left(
\begin{array}
[c]{c}%
I_{r-1}\\
0_{\left(  p-r+1\right)  \times\left(  r-1\right)  }%
\end{array}
\right)  \\\text{(by (\ref{pf.Grasp.generic.calc.a.Ip1r}))}}}\right)
\nonumber\\
&  \ \ \ \ \ \ \ \ \ \ \left(
\begin{array}
[c]{c}%
\text{by Proposition \ref{prop.attach.matrices}, applied to } \mathbb{F}%
\left[  x_{\mathbf{P}}\right]  \text{, }p\text{, } p+q\text{, }\\
\left(  I_{p}\mid C\right)  \text{, }0\text{, }1\text{, }1\text{ and }r\text{
instead of }\mathbb{K}\text{, }u\text{, }v\text{, }A\text{, }a\text{,
}b\text{, }c\text{ and }d
\end{array}
\right) \nonumber\\
&  =\left(  \left(  -1\right)  ^{p-1}\underbrace{\left(  I_{p}\mid C\right)
_{p+q}}_{\substack{=C_{q}\\\text{(since }q>0\text{ and }q\leq q\text{)}%
}}\ \mid\ \left(
\begin{array}
[c]{c}%
I_{r-1}\\
0_{\left(  p-r+1\right)  \times\left(  r-1\right)  }%
\end{array}
\right)  \right) \nonumber\\
&  =\left(  \left(  -1\right)  ^{p-1}C_{q}\ \mid\ \left(
\begin{array}
[c]{c}%
I_{r-1}\\
0_{\left(  p-r+1\right)  \times\left(  r-1\right)  }%
\end{array}
\right)  \right)  . \label{pf.Grasp.generic.calc.a.Dstep}%
\end{align}
Now, by Proposition \ref{prop.attach.matrices} (applied to $\mathbb{F}\left[
x_{\mathbf{P}}\right]  $, $p$, $p+q$, $\left(  I_{p}\mid C\right)  $, $0$,
$r$, $r+s$ and $p+s$ instead of $\mathbb{K}$, $u$, $v$, $A$, $a$, $b$, $c$ and
$d$), we have%
\begin{align*}
&  \left(  I_{p}\mid C\right)  \left[  0:r\mid r+s:p+s\right] \\
&  =\left(  \left(  I_{p}\mid C\right)  \left[  0:r\right]  \ \mid\ \left(
I_{p}\mid C\right)  \left[  r+s:p+s\right]  \right) \\
&  =\left(  \left(  \left(  -1\right)  ^{p-1}C_{q}\ \mid\ \left(
\begin{array}
[c]{c}%
I_{r-1}\\
0_{\left(  p-r+1\right)  \times\left(  r-1\right)  }%
\end{array}
\right)  \right)  \ \mid\ \left(  I_{p}\mid C\right)  \left[  r+s:p+s\right]
\right) \\
&  \ \ \ \ \ \ \ \ \ \ \left(  \text{by (\ref{pf.Grasp.generic.calc.a.Dstep}%
)}\right) \\
&  =\left(  \left(  -1\right)  ^{p-1}C_{q}\ \mid\ \left(
\begin{array}
[c]{c}%
I_{r-1}\\
0_{\left(  p-r+1\right)  \times\left(  r-1\right)  }%
\end{array}
\right)  \ \mid\ \left(  I_{p}\mid C\right)  \left[  r+s:p+s\right]  \right)
.
\end{align*}
Thus,%
\begin{align}
&  \det\left(  \left(  I_{p}\mid C\right)  \left[  0:r\mid r+s:p+s\right]
\right) \nonumber\\
&  =\det\left(  \left(  -1\right)  ^{p-1}C_{q}\ \mid\ \left(
\begin{array}
[c]{c}%
I_{r-1}\\
0_{\left(  p-r+1\right)  \times\left(  r-1\right)  }%
\end{array}
\right)  \ \mid\ \left(  I_{p}\mid C\right)  \left[  r+s:p+s\right]  \right)
\nonumber\\
&  =\left(  -1\right)  ^{r-1}\det\left(  \left(
\begin{array}
[c]{c}%
I_{r-1}\\
0_{\left(  p-r+1\right)  \times\left(  r-1\right)  }%
\end{array}
\right)  \ \mid\ \left(  -1\right)  ^{p-1}C_{q}\ \mid\ \left(  I_{p}\mid
C\right)  \left[  r+s:p+s\right]  \right) \nonumber\\
&  \ \ \ \ \ \ \ \ \ \ \left(
\begin{array}
[c]{c}%
\text{since permuting the columns of a matrix merely multiplies the}\\
\text{sign of its determinant by the sign of the permutation}%
\end{array}
\right) \nonumber\\
&  =\left(  -1\right)  ^{r-1}\underbrace{\det\left(  \left(
\begin{array}
[c]{c}%
I_{r-1}\\
0_{\left(  p-r+1\right)  \times\left(  r-1\right)  }%
\end{array}
\right)  \ \mid\ \left(  \left(  -1\right)  ^{p-1}C_{q}\ \mid\ \left(
I_{p}\mid C\right)  \left[  r+s:p+s\right]  \right)  \right)  }%
_{\substack{=\det\left(  \operatorname*{rows}\nolimits_{r,r+1,...,p}\left(
\left(  -1\right)  ^{p-1}C_{q}\ \mid\ \left(  I_{p}\mid C\right)  \left[
r+s:p+s\right]  \right)  \right)  \\\text{(by Lemma \ref{lem.attach.matrices}
\textbf{(a)}, applied to }p-r+1\text{, }\mathbb{F}\left[  x_{\mathbf{P}%
}\right]  \text{ and}\\\left(  \left(  -1\right)  ^{p-1}C_{q}\ \mid\ \left(
I_{p}\mid C\right)  \left[  r+s:p+s\right]  \right)  \text{ instead of
}k\text{, }\mathbb{K}\text{ and }U\text{)}}}\nonumber\\
&  =\left(  -1\right)  ^{r-1}\det\left(  \operatorname*{rows}%
\nolimits_{r,r+1,...,p}\left(  \left(  -1\right)  ^{p-1}C_{q}\ \mid\ \left(
I_{p}\mid C\right)  \left[  r+s:p+s\right]  \right)  \right)  .
\label{pf.Grasp.generic.calc.a.simplD1}%
\end{align}
Now, by the definition of $\mathfrak{D}_{\left(  r,s\right)  }$, we have
\begin{align*}
\mathfrak{D}_{\left(  r,s\right)  }  &  =\det\left(  \left(  I_{p}\mid
C\right)  \left[  0:r\mid r+s:p+s\right]  \right) \\
&  =\left(  -1\right)  ^{r-1}\det\left(  \operatorname*{rows}%
\nolimits_{r,r+1,...,p}\left(  \left(  -1\right)  ^{p-1}C_{q}\ \mid\ \left(
I_{p}\mid C\right)  \left[  r+s:p+s\right]  \right)  \right)
\end{align*}
(by (\ref{pf.Grasp.generic.calc.a.simplD1})). This proves
(\ref{lem.Grasp.generic.calc.a.D}). Thus, Lemma \ref{lem.Grasp.generic.calc}
\textbf{(a)} is proven.

\textbf{(b)} Let $L:\left\{  1,2,...,q\right\}  \rightarrow\left\{
1,2,...,q\right\}  $ denote the map which sends every $\ell\in\left\{
2,3,...,q\right\}  $ to $\ell-1$ and sends $1$ to $q$. Then, $L$ is a
permutation in the symmetric group $S_{q}$, and actually is the inverse of the
permutation $Z$. Hence, $L\circ Z=\operatorname*{id}$ and $Z\circ
L=\operatorname*{id}$.

Notice that every $\left(  i,k\right)  \in\mathbf{P}$ satisfies%
\begin{equation}
x_{\left(  i,k\right)  }=y_{\left(  i,L\left(  k\right)  \right)  }.
\label{pf.Grasp.generic.calc.b.leftmove}%
\end{equation}
\footnote{\textit{Proof of (\ref{pf.Grasp.generic.calc.b.leftmove}):} Let
$\left(  i,k\right)  \in\mathbf{P}$. Then, $\left(  i,k\right)  \in
\mathbf{P}=\left\{  1,2,...,p\right\}  \times\left\{  1,2,...,q\right\}  $, so
that $\left(  i,L\left(  k\right)  \right)  \in\left\{  1,2,...,p\right\}
\times\left\{  1,2,...,q\right\}  =\mathbf{P}$. Hence, by the definition of
$y_{\left(  i,L\left(  k\right)  \right)  }$, we have $y_{\left(  i,L\left(
k\right)  \right)  }=x_{\left(  i,Z\left(  L\left(  k\right)  \right)
\right)  }$. But since $Z\left(  L\left(  k\right)  \right)
=\underbrace{\left(  Z\circ L\right)  }_{=\operatorname*{id}}\left(  k\right)
=\operatorname*{id}\left(  k\right)  =k$, this rewrites as $y_{\left(
i,L\left(  k\right)  \right)  }=x_{\left(  i,k\right)  }$. This proves
(\ref{pf.Grasp.generic.calc.b.leftmove}).}

Now, let $\left(  i,k\right)  \in\mathbf{P}$ be such that $i\neq p$ and
$k\neq1$.

We have $\left(  i,k\right)  \in\mathbf{P}=\left\{  1,2,...,p\right\}
\times\left\{  1,2,...,q\right\}  $. Hence, $i\in\left\{  1,2,...,p\right\}  $
and $k\in\left\{  1,2,...,q\right\}  $. Since $i\in\left\{  1,2,...,p\right\}
$ and $i\neq p$, we have $i\in\left\{  1,2,...,p-1\right\}  $, so that $1\leq
i\leq p-1$. Since $k\in\left\{  1,2,...,q\right\}  $ and $k\neq1$, we have
$k\in\left\{  2,3,...,q\right\}  $. Hence, $2\leq k\leq q$. Also, $k\geq2$, so
that $k-1\geq1$. Since $i\in\left\{  1,2,...,p-1\right\}  $, we have
$i+1\in\left\{  2,3,...,p\right\}  $.

Since $k\in\left\{  2,3,...,q\right\}  $, we have $k-1\in\left\{
1,2,...,q-1\right\}  $, so that $Z\left(  k-1\right)  =\left(  k-1\right)  +1$
(by the definition of $Z$). Hence, $Z\left(  k-1\right)  =\left(  k-1\right)
+1=k$. Also, from $k-1\in\left\{  1,2,...,q-1\right\}  $, it follows that
$k-1\geq1>0$, so that $p+k-1>p$. Also, $p-\underbrace{i}_{\leq p-1}+1\geq
p-\left(  p-1\right)  +1=2\geq1$.

Set $\mathbf{p}=\left(  i,k\right)  $. Clearly, $\mathbf{p}=\left(
i,k\right)  \in\mathbf{P}$. Now, let $\mathbf{N}$ denote the $\left(
p-i+1\right)  \times\left(  p-i+1\right)  $-matrix $\operatorname*{rows}%
\nolimits_{i,i+1,...,p}\left(  \left(  I_{p}\mid C\right)  \left[
i+k-1:p+k\right]  \right)  $. Then,%
\begin{equation}
\text{the }\left(  1,p-i+1\right)  \text{-th entry of the matrix }%
\mathbf{N}\text{ is }x_{\mathbf{p}} \label{pf.Grasp.generic.calc.b.Ndet2}%
\end{equation}
\footnote{\textit{Proof of (\ref{pf.Grasp.generic.calc.b.Ndet2}):} Notice that
$1\leq p+k-1\leq p+q$ (since $1\leq p=p+\underbrace{0}_{\leq k-1}\leq p+k-1$
and $p+\underbrace{k}_{\leq q}-1\leq p+q-1\leq p+q$). Thus,
\begin{align*}
\left(  I_{p}\mid C\right)  _{p+k-1}  &  =\left(  \text{the }\left(
p+k-1\right)  \text{-th column of the matrix }\left(  I_{p}\mid C\right)
\right) \\
&  =\left(  \text{the }\left(  k-1\right)  \text{-th column of the matrix
}C\right) \\
&  \ \ \ \ \ \ \ \ \ \ \left(  \text{by the definition of }\left(  I_{p}\mid
C\right)  \text{, because }p+k-1>p\right)  .
\end{align*}
But since $\mathbf{N}=\operatorname*{rows}\nolimits_{i,i+1,...,p}\left(
\left(  I_{p}\mid C\right)  \left[  i+k-1:p+k\right]  \right)  $, we have%
\begin{align*}
&  \left(  \text{the }\left(  1,p-i+1\right)  \text{-th entry of the matrix
}\mathbf{N}\right) \\
&  =\left(  \text{the }\left(  1,p-i+1\right)  \text{-th entry of the matrix
}\operatorname*{rows}\nolimits_{i,i+1,...,p}\left(  \left(  I_{p}\mid
C\right)  \left[  i+k-1:p+k\right]  \right)  \right) \\
&  =\left(  \text{the }\left(  i,p-i+1\right)  \text{-th entry of the matrix
}\left(  I_{p}\mid C\right)  \left[  i+k-1:p+k\right]  \right) \\
&  \ \ \ \ \ \ \ \ \ \ \left(  \text{by the definition of }%
\operatorname*{rows}\nolimits_{i,i+1,...,p}\left(  \left(  I_{p}\mid C\right)
\left[  i+k-1:p+k\right]  \right)  \right) \\
&  =\left(  \text{the }i\text{-th entry of }\underbrace{\text{the }\left(
p-i+1\right)  \text{-th column of the matrix }\left(  I_{p}\mid C\right)
\left[  i+k-1:p+k\right]  }_{\substack{=\left(  I_{p}\mid C\right)  _{\left(
i+k-1\right)  +\left(  p-i+1\right)  -1}=\left(  I_{p}\mid C\right)
_{p+k-1}=\left(  \text{the }\left(  k-1\right)  \text{-th column of the matrix
}C\right)  \\\text{(as proven above)}}}\right) \\
&  =\left(  \text{the }i\text{-th entry of the }\left(  k-1\right)  \text{-th
column of the matrix }C\right) \\
&  =\left(  \text{the }\left(  i,k-1\right)  \text{-th entry of the matrix
}C\right) \\
&  =x_{\left(  i,Z\left(  k-1\right)  \right)  }\ \ \ \ \ \ \ \ \ \ \left(
\text{by (\ref{pf.Grasp.generic.calc.entry}), applied to }\left(
i,k-1\right)  \text{ instead of }\left(  i,k\right)  \right) \\
&  =x_{\mathbf{p}}\ \ \ \ \ \ \ \ \ \ \left(  \text{because }\left(
i,\underbrace{Z\left(  k-1\right)  }_{=k}\right)  =\left(  i,k\right)
=\mathbf{p}\right)  .
\end{align*}
This proves (\ref{pf.Grasp.generic.calc.b.Ndet2}).}, while%
\begin{equation}
\text{all entries of the matrix }\mathbf{N}\text{ other than the }\left(
1,p-i+1\right)  \text{-th entry belong to the ring }\mathbb{F}\left[
x_{\mathbf{p}\Downarrow}\right]  \label{pf.Grasp.generic.calc.b.Ndet3}%
\end{equation}
\footnote{\textit{Proof of (\ref{pf.Grasp.generic.calc.b.Ndet3}):} We need to
show that for every $\left(  u,v\right)  \in\left\{  1,2,...,p-i+1\right\}
\times\left\{  1,2,...,p-i+1\right\}  $ satisfying $\left(  u,v\right)
\neq\left(  1,p-i+1\right)  $, the $\left(  u,v\right)  $-th entry of the
matrix $\mathbf{N}$ belongs to the ring $\mathbb{F}\left[  x_{\mathbf{p}%
\Downarrow}\right]  $.
\par
So let $\left(  u,v\right)  \in\left\{  1,2,...,p-i+1\right\}  \times\left\{
1,2,...,p-i+1\right\}  $ be such that $\left(  u,v\right)  \neq\left(
1,p-i+1\right)  $. Let $\eta$ be the $\left(  u,v\right)  $-th entry of the
matrix $\mathbf{N}$. We are going to show that $\eta\in\mathbb{F}\left[
x_{\mathbf{p}\Downarrow}\right]  $.
\par
Since $\left(  u,v\right)  \in\left\{  1,2,...,p-i+1\right\}  \times\left\{
1,2,...,p-i+1\right\}  $, we have $u\in\left\{  1,2,...,p-i+1\right\}  $ and
$v\in\left\{  1,2,...,p-i+1\right\}  $. Thus, $1\leq u\leq p-i+1$ and $1\leq
v\leq p-i+1$.
\par
Recall that $\eta$ is the $\left(  u,v\right)  $-th entry of the matrix
$\mathbf{N}$. Thus,
\begin{align}
\eta &  =\left(  \text{the }\left(  u,v\right)  \text{-th entry of the matrix
}\mathbf{N}\right) \nonumber\\
&  =\left(  \text{the }\left(  u,v\right)  \text{-th entry of the matrix
}\operatorname*{rows}\nolimits_{i,i+1,...,p}\left(  \left(  I_{p}\mid
C\right)  \left[  i+k-1:p+k\right]  \right)  \right) \nonumber\\
&  \ \ \ \ \ \ \ \ \ \ \left(  \text{since }\mathbf{N}=\operatorname*{rows}%
\nolimits_{i,i+1,...,p}\left(  \left(  I_{p}\mid C\right)  \left[
i+k-1:p+k\right]  \right)  \right) \nonumber\\
&  =\left(  \text{the }\left(  i+u-1,v\right)  \text{-th entry of the matrix
}\left(  I_{p}\mid C\right)  \left[  i+k-1:p+k\right]  \right) \nonumber\\
&  \ \ \ \ \ \ \ \ \ \ \left(  \text{by the definition of the matrix
}\operatorname*{rows}\nolimits_{i,i+1,...,p}\left(  \left(  I_{p}\mid
C\right)  \left[  i+k-1:p+k\right]  \right)  \right) \nonumber\\
&  =\left(  \text{the }\left(  i+u-1\right)  \text{-th entry of }%
\underbrace{\text{the }v\text{-th column of the matrix }\left(  I_{p}\mid
C\right)  \left[  i+k-1:p+k\right]  }_{\substack{=\left(  I_{p}\mid C\right)
_{\left(  i+k-1\right)  +v-1}\\\text{(by the definition of }\left(  I_{p}\mid
C\right)  \left[  i+k-1:p+k\right]  \text{)}}}\right) \nonumber\\
&  =\left(  \text{the }\left(  i+u-1\right)  \text{-th entry of }\left(
I_{p}\mid C\right)  _{\left(  i+k-1\right)  +v-1}\right) \nonumber\\
&  =\left(  \text{the }\left(  i+u-1\right)  \text{-th entry of }\left(
I_{p}\mid C\right)  _{i+k-2+v}\right)  . \label{pf.Grasp.generic.calc.b.Nl1}%
\end{align}
\par
Notice that $i+k-2+v=\left(  \underbrace{i}_{\geq1}-1\right)
+\underbrace{\left(  k-1\right)  }_{\geq1}+\underbrace{v}_{\geq0}\geq\left(
1-1\right)  +1+0=1$. Also, $i+\underbrace{k}_{\leq q}-2+\underbrace{v}_{\leq
p-i+1}\leq i+q-2+p-i+1=p+q-1<p+q$.
\par
Since $i+k-2+v\geq1$ and $i+k-2+v\leq p+q$, it is clear that
\[
\left(  I_{p}\mid C\right)  _{i+k-2+v}=\left(  \text{the }\left(
i+k-2+v\right)  \text{-th column of the matrix }\left(  I_{p}\mid C\right)
\right)  .
\]
\par
If $i+k-2+v\leq p$, then we thus have%
\begin{align*}
\left(  I_{p}\mid C\right)  _{i+k-2+v}  &  =\left(  \text{the }\left(
i+k-2+v\right)  \text{-th column of the matrix }\left(  I_{p}\mid C\right)
\right) \\
&  =\left(  \text{the }\left(  i+k-2+v\right)  \text{-th column of the matrix
}I_{p}\right)
\end{align*}
(by the definition of $\left(  I_{p}\mid C\right)  $, since $1\leq i+k-2+v\leq
p$). Hence, if $i+k-2+v\leq p$, then (\ref{pf.Grasp.generic.calc.b.Nl1})
becomes%
\begin{align*}
\eta &  =\left(  \text{the }\left(  i+u-1\right)  \text{-th entry of
}\underbrace{\left(  I_{p}\mid C\right)  _{i+k-2+v}}_{=\left(  \text{the
}\left(  i+k-2+v\right)  \text{-th column of the matrix }I_{p}\right)
}\right) \\
&  =\left(  \text{the }\left(  i+u-1\right)  \text{-th entry of the }\left(
i+k-2+v\right)  \text{-th column of the matrix }I_{p}\right) \\
&  \in\mathbb{F}\ \ \ \ \ \ \ \ \ \ \left(  \text{since the matrix }%
I_{p}\text{ is defined over }\mathbb{F}\right) \\
&  \subseteq\mathbb{F}\left[  x_{\mathbf{p}\Downarrow}\right]  .
\end{align*}
Thus, if $i+k-2+v\leq p$, then $\eta\in\mathbb{F}\left[  x_{\mathbf{p}%
\Downarrow}\right]  $ is proven. Hence, for the rest of the proof of $\eta
\in\mathbb{F}\left[  x_{\mathbf{p}\Downarrow}\right]  $, we can WLOG assume
that we don't have $i+k-2+v\leq p$. Assume this.
\par
We don't have $i+k-2+v\leq p$. Thus, $i+k-2+v>p$. Hence, $p<i+k-2+v\leq p+q$.
Now,%
\begin{align*}
\left(  I_{p}\mid C\right)  _{i+k-2+v}  &  =\left(  \text{the }\left(
i+k-2+v\right)  \text{-th column of the matrix }\left(  I_{p}\mid C\right)
\right) \\
&  =\left(  \text{the }\left(  i+k-2+v-p\right)  \text{-th column of the
matrix }C\right)
\end{align*}
(by the definition of $\left(  I_{p}\mid C\right)  $, since $p<i+k-2+v\leq
p+q$). Thus, (\ref{pf.Grasp.generic.calc.b.Nl1}) becomes
\begin{align}
\eta &  =\left(  \text{the }\left(  i+u-1\right)  \text{-th entry of
}\underbrace{\left(  I_{p}\mid C\right)  _{i+k-2+v}}_{=\left(  \text{the
}\left(  i+k-2+v-p\right)  \text{-th column of the matrix }C\right)  }\right)
\nonumber\\
&  =\left(  \text{the }\left(  i+u-1\right)  \text{-th entry of the }\left(
i+k-2+v-p\right)  \text{-th column of the matrix }C\right) \nonumber\\
&  =\left(  \text{the }\left(  i+u-1,i+k-2+v-p\right)  \text{-th entry of the
matrix }C\right) \nonumber\\
&  =x_{\left(  i+u-1,Z\left(  i+k-2+v-p\right)  \right)  }
\label{pf.Grasp.generic.calc.b.Nl4}%
\end{align}
(by (\ref{pf.Grasp.generic.calc.entry}), applied to $\left(
i+u-1,i+k-2+v-p\right)  $ instead of $\left(  i,k\right)  $). But since
$i+k-2+v<p+q$, we have $i+k-2+v-p<q$, so that $i+k-2+v-p\in\left\{
1,2,...,q-1\right\}  $ (because we also have $\underbrace{i+k-2+v}%
_{>p}-p>p-p=0$). Hence, (by the definition of $Z$) we have $Z\left(
i+k-2+v-p\right)  =\left(  i+k-2+v-p\right)  +1=i+k-1+v-p$. Thus,
(\ref{pf.Grasp.generic.calc.b.Nl4}) becomes%
\begin{equation}
\eta=x_{\left(  i+u-1,Z\left(  i+k-2+v-p\right)  \right)  }=x_{\left(
i+u-1,i+k-1+v-p\right)  }\ \ \ \ \ \ \ \ \ \ \left(  \text{since }Z\left(
i+k-2+v-p\right)  =i+k-1+v-p\right)  . \label{pf.Grasp.generic.calc.b.Nl5}%
\end{equation}
\par
But we have $i+\underbrace{u}_{\geq1}-1\geq i+1-1=i$ and $i+k-1+v-p=k+\left(
i-1+\underbrace{v}_{\leq p-i+1}-p\right)  \leq k+\left(  i-1+\left(
p-i+1\right)  -p\right)  =k$.
\par
Let us assume (for the sake of contradiction) that $\left(
i+u-1,i+k-1+v-p\right)  =\left(  i,k\right)  $. Then, $i+u-1=i$ and
$i+k-1+v-p=k$. Hence, $u-1=\underbrace{\left(  i+u-1\right)  }_{=i}-i=i-i=0$,
so that $u=1$. Also, subtracting $k$ from the equality $i+k-1+v-p=k$, we
obtain $i-1+v-p=0$, so that $v=p-i+1$. Combined with $u=1$, this yields
$\left(  u,v\right)  =\left(  1,p-i+1\right)  $. This contradicts $\left(
u,v\right)  \neq\left(  1,p-i+1\right)  $. This contradiction shows that our
assumption (that $\left(  i+u-1,i+k-1+v-p\right)  =\left(  i,k\right)  $) was
false. In other words, we don't have $\left(  i+u-1,i+k-1+v-p\right)  =\left(
i,k\right)  $.
\par
We have $i+k-1+v-p=Z\left(  i+k-2+v-p\right)  \in\left\{  1,2,...,q\right\}  $
(because the target of the map $Z$ is $\left\{  1,2,...,q\right\}  $). Also,
$i+u-1\geq i\geq1$ and $i+\underbrace{u}_{\leq p-i+1}-1\leq i+\left(
p-i+1\right)  -1=p$. Hence, $1\leq i+u-1\leq p$, so that $i+u-1\in\left\{
1,2,...,p\right\}  $. Combined with $i+k-1+v-p\in\left\{  1,2,...,q\right\}
$, this yields%
\[
\left(  i+u-1,i+k-1+v-p\right)  \in\left\{  1,2,...,p\right\}  \times\left\{
1,2,...,q\right\}  =\mathbf{P}.
\]
\par
We have $\left(  i+u-1,i+k-1+v-p\right)  \trianglelefteq\left(  i,k\right)  $
(by (\ref{lem.Grasp.generic.calc.ass}), applied to $\left(
i+u-1,i+k-1+v-p\right)  $ and $\left(  i,k\right)  $ instead of $\left(
i,k\right)  $ and $\left(  i^{\prime},k^{\prime}\right)  $), because
$i+u-1\geq i$ and $i+k-1+v-p\leq k$. Combined with the fact that we don't
have $\left(  i+u-1,i+k-1+v-p\right)  =\left(  i,k\right)  $, this shows that
$\left(  i+u-1,i+k-1+v-p\right)  \vartriangleleft\left(  i,k\right)  $. Since
$\left(  i,k\right)  =\mathbf{p}$, this rewrites as $\left(
i+u-1,i+k-1+v-p\right)  \vartriangleleft\mathbf{p}$.
\par
Since $\left(  i+u-1,i+k-1+v-p\right)  \in\mathbf{P}$ and $\left(
i+u-1,i+k-1+v-p\right)  \vartriangleleft\mathbf{p}$, we have%
\[
\left(  i+u-1,i+k-1+v-p\right)  \in\left\{  \mathbf{v}\in\mathbf{P}%
\ \mid\ \mathbf{v}\vartriangleleft\mathbf{p}\right\}  =\left.  \mathbf{p}%
\Downarrow\right.  ,
\]
so that $x_{\left(  i+u-1,i+k-1+v-p\right)  }\in\mathbb{F}\left[
x_{\mathbf{p}\Downarrow}\right]  $. Thus, (\ref{pf.Grasp.generic.calc.b.Nl5})
yields $\eta=x_{\left(  i+u-1,i+k-1+v-p\right)  }\in\mathbb{F}\left[
x_{\mathbf{p}\Downarrow}\right]  $.
\par
We thus have proven that $\eta\in\mathbb{F}\left[  x_{\mathbf{p}\Downarrow
}\right]  $. In other words, $\eta$ belongs to $\mathbb{F}\left[
x_{\mathbf{p}\Downarrow}\right]  $. Since $\eta$ is the $\left(  u,v\right)
$-th entry of the matrix $\mathbf{N}$, this rewrites as follows: The $\left(
u,v\right)  $-th entry of the matrix $\mathbf{N}$ belongs to $\mathbb{F}%
\left[  x_{\mathbf{p}\Downarrow}\right]  $. This proves
(\ref{pf.Grasp.generic.calc.b.Ndet3}).}. Now, Proposition
\ref{prop.cofactor.laplace} (applied to $\mathbb{F}\left[  x_{\mathbf{P}%
}\right]  $, $\mathbf{N}$ and $p-i+1$ instead of $\mathbb{K}$, $B$ and $n$)
yields%
\begin{align}
\det\mathbf{N}  &  =\sum_{\ell=1}^{p-i+1}\left(  \text{the }\left(
1,\ell\right)  \text{-th entry of the matrix }\mathbf{N}\right) \nonumber\\
&  \ \ \ \ \ \ \ \ \ \ \ \ \ \ \ \ \ \ \ \ \cdot\left(  \text{the }\left(
1,\ell\right)  \text{-th cofactor of the matrix }\mathbf{N}\right)  .
\label{pf.Grasp.generic.calc.b.Ndet4}%
\end{align}
But for every $\ell\in\left\{  1,2,...,p-i+1\right\}  $, we have%
\begin{equation}
\left(  \text{the }\left(  1,\ell\right)  \text{-th cofactor of the matrix
}\mathbf{N}\right)  \in\mathbb{F}\left[  x_{\mathbf{p}\Downarrow}\right]  .
\label{pf.Grasp.generic.calc.b.Ndet5}%
\end{equation}
\footnote{\textit{Proof of (\ref{pf.Grasp.generic.calc.b.Ndet5}):} Let
$\ell\in\left\{  1,2,...,p-i+1\right\}  $. By the definition of a cofactor, we
have%
\begin{align*}
&  \left(  \text{the }\left(  1,\ell\right)  \text{-th cofactor of the matrix
}\mathbf{N}\right) \\
&  =\left(  -1\right)  ^{\ell-1}\det\left(  \text{the matrix obtained by
deleting row }1\text{ and column }\ell\text{ from }\mathbf{N}\right)  .
\end{align*}
But consider the matrix obtained by deleting row $1$ and column $\ell$ from
$\mathbf{N}$. All entries of this matrix are entries of $\mathbf{N}$ other
than the $\left(  1,p-i+1\right)  $-th entry (because the $\left(
1,p-i+1\right)  $-th entry of $\mathbf{N}$ lies in row $1$ of $\mathbf{N}$,
which is deleted in the new matrix), and hence belong to $\mathbb{F}\left[
x_{\mathbf{p}\Downarrow}\right]  $ (according to
(\ref{pf.Grasp.generic.calc.b.Ndet3})). Hence, the determinant of this matrix
belongs to $\mathbb{F}\left[  x_{\mathbf{p}\Downarrow}\right]  $ as well. In
other words, $\det\left(  \text{the matrix obtained by deleting row }1\text{
and column }\ell\text{ from }\mathbf{N}\right)  \in\mathbb{F}\left[
x_{\mathbf{p}\Downarrow}\right]  $. Thus,%
\begin{align*}
&  \left(  \text{the }\left(  1,\ell\right)  \text{-th cofactor of the matrix
}\mathbf{N}\right) \\
&  =\left(  -1\right)  ^{\ell-1}\underbrace{\det\left(  \text{the matrix
obtained by deleting row }1\text{ and column }\ell\text{ from }\mathbf{N}%
\right)  }_{\in\mathbb{F}\left[  x_{\mathbf{p}\Downarrow}\right]  }\\
&  \in\left(  -1\right)  ^{\ell-1}\mathbb{F}\left[  x_{\mathbf{p}\Downarrow
}\right]  \subseteq\mathbb{F}\left[  x_{\mathbf{p}\Downarrow}\right]  .
\end{align*}
This proves (\ref{pf.Grasp.generic.calc.b.Ndet5}).} Thus,
(\ref{pf.Grasp.generic.calc.b.Ndet4}) becomes%
\begin{align}
\det\mathbf{N}  &  =\sum_{\ell=1}^{p-i+1}\left(  \text{the }\left(
1,\ell\right)  \text{-th entry of the matrix }\mathbf{N}\right) \nonumber\\
&  \ \ \ \ \ \ \ \ \ \ \ \ \ \ \ \ \ \ \ \ \cdot\left(  \text{the }\left(
1,\ell\right)  \text{-th cofactor of the matrix }\mathbf{N}\right) \nonumber\\
&  =\sum_{\ell=1}^{p-i}\underbrace{\left(  \text{the }\left(  1,\ell\right)
\text{-th entry of the matrix }\mathbf{N}\right)  }_{\substack{\in
\mathbb{F}\left[  x_{\mathbf{p}\Downarrow}\right]  \\\text{(by
(\ref{pf.Grasp.generic.calc.b.Ndet3}) (because the}\\\left(  1,\ell\right)
\text{-th entry of the matrix }\mathbf{N}\text{ is an entry of the matrix
}\mathbf{N}\text{ other than}\\\text{the }\left(  1,p-i+1\right)  \text{-th
entry (since }\ell\neq p-i+1\text{)))}}}\nonumber\\
&  \ \ \ \ \ \ \ \ \ \ \ \ \ \ \ \ \ \ \ \ \cdot\underbrace{\left(  \text{the
}\left(  1,\ell\right)  \text{-th cofactor of the matrix }\mathbf{N}\right)
}_{\substack{\in\mathbb{F}\left[  x_{\mathbf{p}\Downarrow}\right]  \\\text{(by
(\ref{pf.Grasp.generic.calc.b.Ndet5}))}}}\nonumber\\
&  \ \ \ \ \ \ \ \ \ \ +\underbrace{\left(  \text{the }\left(  1,p-i+1\right)
\text{-th entry of the matrix }\mathbf{N}\right)  }_{\substack{=x_{\mathbf{p}%
}\\\text{(by (\ref{pf.Grasp.generic.calc.b.Ndet2}))}}}\nonumber\\
&  \ \ \ \ \ \ \ \ \ \ \ \ \ \ \ \ \ \ \ \ \cdot\left(  \text{the }\left(
1,p-i+1\right)  \text{-th cofactor of the matrix }\mathbf{N}\right)
\nonumber\\
&  \in\underbrace{\sum_{\ell=1}^{p-i}\mathbb{F}\left[  x_{\mathbf{p}%
\Downarrow}\right]  \cdot\mathbb{F}\left[  x_{\mathbf{p}\Downarrow}\right]
}_{\substack{\subseteq\mathbb{F}\left[  x_{\mathbf{p}\Downarrow}\right]
\\\text{(since }\mathbb{F}\left[  x_{\mathbf{p}\Downarrow}\right]  \text{ is
an }\mathbb{F}\text{-algebra)}}}\nonumber\\
&  \ \ \ \ \ \ \ \ \ \ +x_{\mathbf{p}}\cdot\left(  \text{the }\left(
1,p-i+1\right)  \text{-th cofactor of the matrix }\mathbf{N}\right)
\nonumber\\
&  \subseteq\mathbb{F}\left[  x_{\mathbf{p}\Downarrow}\right]  +x_{\mathbf{p}%
}\cdot\left(  \text{the }\left(  1,p-i+1\right)  \text{-th cofactor of the
matrix }\mathbf{N}\right)  . \label{pf.Grasp.generic.calc.b.Ndet6}%
\end{align}

But $\mathbf{N}=\operatorname*{rows}\nolimits_{i,i+1,...,p}\left(  \left(
I_{p}\mid C\right)  \left[  i+k-1:p+k\right]  \right)  $, so that%
\begin{align}
&  \left(  \text{the }\left(  1,p-i+1\right)  \text{-th cofactor of the matrix
}\mathbf{N}\right) \nonumber\\
&  =\left(  \text{the }\left(  1,p-i+1\right)  \text{-th cofactor of the
matrix }\operatorname*{rows}\nolimits_{i,i+1,...,p}\left(  \left(  I_{p}\mid
C\right)  \left[  i+k-1:p+k\right]  \right)  \right) \nonumber\\
&  =\left(  -1\right)  ^{p-i}\det\left(  \operatorname*{rows}%
\nolimits_{i+1,i+2,...,p}\left(  \left(  I_{p}\mid C\right)  \left[
i+k-1:p+k-1\right]  \right)  \right) \label{pf.Grasp.generic.calc.b.Ndet8}\\
&  \ \ \ \ \ \ \ \ \ \ \left(
\begin{array}
[c]{c}%
\text{by Corollary \ref{cor.cofactor.1}, applied to }\mathbb{F}\left[
x_{\mathbf{P}}\right]  \text{, }p+q\text{, }\left(  I_{p}\mid C\right)
\text{, }i+k-1\text{ and }p+k\\
\text{instead of }\mathbb{K}\text{, }q\text{, }A\text{, }a\text{ and }b\text{
(since }\underbrace{i}_{\leq p-1<p}+k-1<p+k-1<p+k\text{,}\\
\left(  p+k\right)  -\left(  i+k-1\right)  =p-\underbrace{i}_{\geq1}+1\leq
p-1+1=p\\
\text{and }i=p-\left(  \left(  p+k\right)  -\left(  i+k-1\right)  \right)
+1\text{)}%
\end{array}
\right)  .\nonumber
\end{align}
Multiplying this equality with $\left(  -1\right)  ^{p-i}$, we obtain%
\begin{align}
&  \left(  -1\right)  ^{p-i}\left(  \text{the }\left(  1,p-i+1\right)
\text{-th cofactor of the matrix }\mathbf{N}\right) \nonumber\\
&  =\underbrace{\left(  -1\right)  ^{p-i}\left(  -1\right)  ^{p-i}%
}_{\substack{=\left(  -1\right)  ^{\left(  p-i\right)  +\left(  p-i\right)
}=1\\\text{(since }\left(  p-i\right)  +\left(  p-i\right)  =2\left(
p-i\right)  \text{ is even)}}}\det\left(  \operatorname*{rows}%
\nolimits_{i+1,i+2,...,p}\left(  \left(  I_{p}\mid C\right)  \left[
i+k-1:p+k-1\right]  \right)  \right) \nonumber\\
&  =\det\left(  \operatorname*{rows}\nolimits_{i+1,i+2,...,p}\left(  \left(
I_{p}\mid C\right)  \left[  i+k-1:p+k-1\right]  \right)  \right)  .
\label{pf.Grasp.generic.calc.b.Ndet8a}%
\end{align}

On the other hand, combining $i+1\in\left\{  2,3,...,p\right\}  \subseteq
\left\{  1,2,...,p\right\}  $ with $k-1\in\left\{  1,2,...,q-1\right\}
\subseteq\left\{  1,2,...,q\right\}  $, we obtain $\left(  i+1,k-1\right)
\in\left\{  1,2,...,p\right\}  \times\left\{  1,2,...,q\right\}  =\mathbf{P}$.
Applying (\ref{lem.Grasp.generic.calc.a.N}) to $\left(  r,s\right)  =\left(
i+1,k-1\right)  $, we obtain%
\begin{align}
\mathfrak{N}_{\left(  i+1,k-1\right)  }  &  =\det\left(  \operatorname*{rows}%
\nolimits_{i+1,i+2,...,p}\left(  \left(  I_{p}\mid C\right)  \left[
\underbrace{\left(  i+1\right)  +\left(  k-1\right)  }_{=i+k}%
-1:\underbrace{p+\left(  k-1\right)  }_{=p+k-1}\right]  \right)  \right)
\nonumber\\
&  =\det\left(  \operatorname*{rows}\nolimits_{i+1,i+2,...,p}\left(  \left(
I_{p}\mid C\right)  \left[  i+k-1:p+k-1\right]  \right)  \right) \nonumber\\
&  =\left(  -1\right)  ^{p-i}\underbrace{\left(  \text{the }\left(
1,p-i+1\right)  \text{-th cofactor of the matrix }\mathbf{N}\right)
}_{\substack{\in\mathbb{F}\left[  x_{\mathbf{p}\Downarrow}\right]  \\\text{(by
(\ref{pf.Grasp.generic.calc.b.Ndet5}), applied to }\ell=p-i+1\text{)}%
}}\ \ \ \ \ \ \ \ \ \ \left(  \text{by (\ref{pf.Grasp.generic.calc.b.Ndet8a}%
)}\right) \label{pf.Grasp.generic.calc.b.Ndet9}\\
&  \in\left(  -1\right)  ^{p-i}\mathbb{F}\left[  x_{\mathbf{p}\Downarrow
}\right]  \subseteq\mathbb{F}\left[  x_{\mathbf{p}\Downarrow}\right]
=\mathbb{F}\left[  x_{\left(  i,k\right)  \Downarrow}\right]
\ \ \ \ \ \ \ \ \ \ \left(  \text{since }\mathbf{p}=\left(  i,k\right)
\right)  .\nonumber
\end{align}
Thus, (\ref{lem.Grasp.generic.calc.b.N2}) is proven. Also, multiplying
(\ref{pf.Grasp.generic.calc.b.Ndet9}) with $\left(  -1\right)  ^{p-i}$, we
obtain%
\begin{align}
&  \left(  -1\right)  ^{p-i}\mathfrak{N}_{\left(  i+1,k-1\right)  }\nonumber\\
&  =\underbrace{\left(  -1\right)  ^{p-i}\left(  -1\right)  ^{p-i}%
}_{\substack{=\left(  -1\right)  ^{\left(  p-i\right)  +\left(  p-i\right)
}=1\\\text{(since }\left(  p-i\right)  +\left(  p-i\right)  =2\left(
p-i\right)  \text{ is even)}}}\left(  \text{the }\left(  1,p-i+1\right)
\text{-th cofactor of the matrix }\mathbf{N}\right) \nonumber\\
&  =\left(  \text{the }\left(  1,p-i+1\right)  \text{-th cofactor of the
matrix }\mathbf{N}\right)  . \label{pf.Grasp.generic.calc.b.Ndet9a}%
\end{align}

Now, (\ref{pf.Grasp.generic.calc.b.Ndet6}) becomes%
\begin{align}
&  \det\mathbf{N}\nonumber\\
&  \in\mathbb{F}\left[  x_{\mathbf{p}\Downarrow}\right]  +x_{\mathbf{p}}%
\cdot\underbrace{\left(  \text{the }\left(  1,p-i+1\right)  \text{-th cofactor
of the matrix }\mathbf{N}\right)  }_{\substack{=\left(  -1\right)
^{p-i}\mathfrak{N}_{\left(  i+1,k-1\right)  }\\\text{(by
(\ref{pf.Grasp.generic.calc.b.Ndet9a}))}}}\nonumber\\
&  =\mathbb{F}\left[  x_{\mathbf{p}\Downarrow}\right]  +x_{\mathbf{p}}%
\cdot\left(  -1\right)  ^{p-i}\mathfrak{N}_{\left(  i+1,k-1\right)
}=\mathbb{F}\left[  x_{\mathbf{p}\Downarrow}\right]  +\left(  -1\right)
^{p-i}\cdot\mathfrak{N}_{\left(  i+1,k-1\right)  }\cdot x_{\mathbf{p}%
}\nonumber\\
&  =\left(  -1\right)  ^{p-i}\cdot\mathfrak{N}_{\left(  i+1,k-1\right)  }\cdot
x_{\mathbf{p}}+\mathbb{F}\left[  x_{\mathbf{p}\Downarrow}\right]  =\left(
-1\right)  ^{p-i}\cdot\mathfrak{N}_{\left(  i+1,k-1\right)  }\cdot x_{\left(
i,k\right)  }+\mathbb{F}\left[  x_{\left(  i,k\right)  \Downarrow}\right]
\label{pf.Grasp.generic.calc.b.Ndet38}%
\end{align}
(since $\mathbf{p}=\left(  i,k\right)  $). Now,
(\ref{lem.Grasp.generic.calc.a.N}) (applied to $\left(  r,s\right)  =\left(
i,k\right)  $) yields%
\begin{align*}
\mathfrak{N}_{\left(  i,k\right)  }  &  =\det\left(
\underbrace{\operatorname*{rows}\nolimits_{i,i+1,...,p}\left(  \left(
I_{p}\mid C\right)  \left[  i+k-1:p+k\right]  \right)  }_{=\mathbf{N}}\right)
\\
&  =\det\mathbf{N}\in\left(  -1\right)  ^{p-i}\cdot\mathfrak{N}_{\left(
i+1,k-1\right)  }\cdot x_{\left(  i,k\right)  }+\mathbb{F}\left[  x_{\left(
i,k\right)  \Downarrow}\right]
\end{align*}
(by (\ref{pf.Grasp.generic.calc.b.Ndet38})). Thus,
(\ref{lem.Grasp.generic.calc.b.N}) is proven.

On the other hand, let $\mathbf{M}$ be the $\left(  p-i+1\right)
\times\left(  p-i+1\right)  $-matrix \newline$\operatorname*{rows}%
\nolimits_{i,i+1,...,p}\left(  \left(  -1\right)  ^{p-1}C_{q}\ \mid\ \left(
I_{p}\mid C\right)  \left[  i+k:p+k\right]  \right)  $. Then,%
\begin{equation}
\text{the }\left(  1,p-i+1\right)  \text{-th entry of the matrix }%
\mathbf{M}\text{ is }x_{\mathbf{p}} \label{pf.Grasp.generic.calc.b.Ddet2}%
\end{equation}
\footnote{\textit{Proof of (\ref{pf.Grasp.generic.calc.b.Ddet2}):} Notice that
$1\leq p+k-1\leq p+q$ (since $1\leq p=p+\underbrace{0}_{\leq k-1}\leq p+k-1$
and $p+\underbrace{k}_{\leq q}-1\leq p+q-1\leq p+q$). Thus,
\begin{align*}
\left(  I_{p}\mid C\right)  _{p+k-1}  &  =\left(  \text{the }\left(
p+k-1\right)  \text{-th column of the matrix }\left(  I_{p}\mid C\right)
\right) \\
&  =\left(  \text{the }\left(  k-1\right)  \text{-th column of the matrix
}C\right) \\
&  \ \ \ \ \ \ \ \ \ \ \left(  \text{by the definition of }\left(  I_{p}\mid
C\right)  \text{, because }p+k-1>p\right)  .
\end{align*}
Also, $p-\underbrace{i}_{\leq p-1<p}+1>p-p+1=1$. But since $\mathbf{M}%
=\operatorname*{rows}\nolimits_{i,i+1,...,p}\left(  \left(  -1\right)
^{p-1}C_{q}\ \mid\ \left(  I_{p}\mid C\right)  \left[  i+k:p+k\right]
\right)  $, we have%
\begin{align*}
&  \left(  \text{the }\left(  1,p-i+1\right)  \text{-th entry of the matrix
}\mathbf{M}\right) \\
&  =\left(  \text{the }\left(  1,p-i+1\right)  \text{-th entry of the matrix
}\operatorname*{rows}\nolimits_{i,i+1,...,p}\left(  \left(  -1\right)
^{p-1}C_{q}\ \mid\ \left(  I_{p}\mid C\right)  \left[  i+k:p+k\right]
\right)  \right) \\
&  =\left(  \text{the }\left(  i,p-i+1\right)  \text{-th entry of the matrix
}\left(  \left(  -1\right)  ^{p-1}C_{q}\ \mid\ \left(  I_{p}\mid C\right)
\left[  i+k:p+k\right]  \right)  \right) \\
&  \ \ \ \ \ \ \ \ \ \ \left(  \text{by the definition of }%
\operatorname*{rows}\nolimits_{i,i+1,...,p}\left(  \left(  -1\right)
^{p-1}C_{q}\ \mid\ \left(  I_{p}\mid C\right)  \left[  i+k:p+k\right]
\right)  \right) \\
&  =\left(  \text{the }i\text{-th entry of }\underbrace{\text{the }\left(
p-i+1\right)  \text{-th column of the matrix }\left(  \left(  -1\right)
^{p-1}C_{q}\ \mid\ \left(  I_{p}\mid C\right)  \left[  i+k:p+k\right]
\right)  }_{\substack{=\left(  \text{the }\left(  p-i\right)  \text{-th column
of the matrix }\left(  I_{p}\mid C\right)  \left[  i+k:p+k\right]  \right)
\\\text{(by the definition of }\left(  \left(  -1\right)  ^{p-1}C_{q}%
\ \mid\ \left(  I_{p}\mid C\right)  \left[  i+k:p+k\right]  \right)  \text{,
since }p-i+1>1\text{)}}}\right) \\
&  =\left(  \text{the }i\text{-th entry of }\underbrace{\text{the }\left(
p-i\right)  \text{-th column of the matrix }\left(  I_{p}\mid C\right)
\left[  i+k:p+k\right]  }_{\substack{=\left(  I_{p}\mid C\right)  _{\left(
i+k\right)  +\left(  p-i\right)  -1}=\left(  I_{p}\mid C\right)
_{p+k-1}=\left(  \text{the }\left(  k-1\right)  \text{-th column of the matrix
}C\right)  \\\text{(as proven above)}}}\right) \\
&  =\left(  \text{the }i\text{-th entry of the }\left(  k-1\right)  \text{-th
column of the matrix }C\right) \\
&  =\left(  \text{the }\left(  i,k-1\right)  \text{-th entry of the matrix
}C\right) \\
&  =x_{\left(  i,Z\left(  k-1\right)  \right)  }\ \ \ \ \ \ \ \ \ \ \left(
\text{by (\ref{pf.Grasp.generic.calc.entry}), applied to }\left(
i,k-1\right)  \text{ instead of }\left(  i,k\right)  \right) \\
&  =x_{\mathbf{p}}\ \ \ \ \ \ \ \ \ \ \left(  \text{because }\left(
i,\underbrace{Z\left(  k-1\right)  }_{=k}\right)  =\left(  i,k\right)
=\mathbf{p}\right)  .
\end{align*}
This proves (\ref{pf.Grasp.generic.calc.b.Ddet2}).}, while%
\begin{equation}
\text{all entries of the matrix }\mathbf{M}\text{ other than the }\left(
1,p-i+1\right)  \text{-th entry belong to the ring }\mathbb{F}\left[
x_{\mathbf{p}\Downarrow}\right]  \label{pf.Grasp.generic.calc.b.Ddet3}%
\end{equation}
\footnote{\textit{Proof of (\ref{pf.Grasp.generic.calc.b.Ddet3}):} We need to
show that for every $\left(  u,v\right)  \in\left\{  1,2,...,p-i+1\right\}
\times\left\{  1,2,...,p-i+1\right\}  $ satisfying $\left(  u,v\right)
\neq\left(  1,p-i+1\right)  $, the $\left(  u,v\right)  $-th entry of the
matrix $\mathbf{M}$ belongs to the ring $\mathbb{F}\left[  x_{\mathbf{p}%
\Downarrow}\right]  $.
\par
So let $\left(  u,v\right)  \in\left\{  1,2,...,p-i+1\right\}  \times\left\{
1,2,...,p-i+1\right\}  $ be such that $\left(  u,v\right)  \neq\left(
1,p-i+1\right)  $. Let $\eta$ be the $\left(  u,v\right)  $-th entry of the
matrix $\mathbf{M}$. We are going to show that $\eta\in\mathbb{F}\left[
x_{\mathbf{p}\Downarrow}\right]  $.
\par
Since $\left(  u,v\right)  \in\left\{  1,2,...,p-i+1\right\}  \times\left\{
1,2,...,p-i+1\right\}  $, we have $u\in\left\{  1,2,...,p-i+1\right\}  $ and
$v\in\left\{  1,2,...,p-i+1\right\}  $. Thus, $1\leq u\leq p-i+1$ and $1\leq
v\leq p-i+1$.
\par
Since $u\leq p-i+1$, we have $i+\underbrace{u}_{\leq p-i+1}-1\leq i+\left(
p-i+1\right)  -1=p$. Also, $i+\underbrace{u}_{\geq1}-1\geq i+u-1\geq i\geq1$.
Thus, $1\leq i+u-1\leq p$, so that $i+u-1\in\left\{  1,2,...,p\right\}  $.
Combined with $1\in\left\{  1,2,...,q\right\}  $, this yields $\left(
i+u-1,1\right)  \in\left\{  1,2,...,p\right\}  \times\left\{
1,2,...,q\right\}  =\mathbf{P}$. Moreover, since $i+\underbrace{u}_{\geq
1}-1\geq i+1-1=i$ and $1\leq k$, we have $\left(  i+u-1,1\right)
\trianglelefteq\left(  i,k\right)  $ (by (\ref{lem.Grasp.generic.calc.ass}),
applied to $\left(  i+u-1,1\right)  $ and $\left(  i,k\right)  $ instead of
$\left(  i,k\right)  $ and $\left(  i^{\prime},k^{\prime}\right)  $). Also,
$1<2\leq k$, so that $1\neq k$ and thus $\left(  i+u-1,1\right)  \neq\left(
i,k\right)  $. Combined with $\left(  i+u-1,1\right)  \trianglelefteq\left(
i,k\right)  $, this yields $\left(  i+u-1,1\right)  \vartriangleleft\left(
i,k\right)  $. Since $\left(  i+u-1,1\right)  \in\mathbf{P}$ and $\left(
i+u-1,1\right)  \vartriangleleft\left(  i,k\right)  =\mathbf{p}$, we have%
\[
\left(  i+u-1,1\right)  \in\left\{  \mathbf{v}\in\mathbf{P}\ \mid
\ \mathbf{v}\vartriangleleft\mathbf{p}\right\}  =\left.  \mathbf{p}%
\Downarrow\right.  .
\]
\par
Recall that $\eta$ is the $\left(  u,v\right)  $-th entry of the matrix
$\mathbf{M}$. Thus,
\begin{align}
\eta &  =\left(  \text{the }\left(  u,v\right)  \text{-th entry of the matrix
}\mathbf{M}\right) \nonumber\\
&  =\left(  \text{the }\left(  u,v\right)  \text{-th entry of the matrix
}\operatorname*{rows}\nolimits_{i,i+1,...,p}\left(  \left(  -1\right)
^{p-1}C_{q}\ \mid\ \left(  I_{p}\mid C\right)  \left[  i+k:p+k\right]
\right)  \right) \nonumber\\
&  \ \ \ \ \ \ \ \ \ \ \left(  \text{since }\mathbf{M}=\operatorname*{rows}%
\nolimits_{i,i+1,...,p}\left(  \left(  -1\right)  ^{p-1}C_{q}\ \mid\ \left(
I_{p}\mid C\right)  \left[  i+k:p+k\right]  \right)  \right) \nonumber\\
&  =\left(  \text{the }\left(  i+u-1,v\right)  \text{-th entry of the matrix
}\left(  \left(  -1\right)  ^{p-1}C_{q}\ \mid\ \left(  I_{p}\mid C\right)
\left[  i+k:p+k\right]  \right)  \right)  \label{pf.Grasp.generic.calc.b.Dl1}%
\\
&  \ \ \ \ \ \ \ \ \ \ \left(  \text{by the definition of the matrix
}\operatorname*{rows}\nolimits_{i,i+1,...,p}\left(  \left(  -1\right)
^{p-1}C_{q}\ \mid\ \left(  I_{p}\mid C\right)  \left[  i+k:p+k\right]
\right)  \right)  .\nonumber
\end{align}
\par
If $v=1$, then this becomes%
\begin{align*}
\eta &  =\left(  \text{the }\left(  i+u-1,v\right)  \text{-th entry of the
matrix }\left(  \left(  -1\right)  ^{p-1}C_{q}\ \mid\ \left(  I_{p}\mid
C\right)  \left[  i+k:p+k\right]  \right)  \right) \\
&  =\left(  \text{the }\left(  i+u-1,1\right)  \text{-th entry of the matrix
}\left(  \left(  -1\right)  ^{p-1}C_{q}\ \mid\ \left(  I_{p}\mid C\right)
\left[  i+k:p+k\right]  \right)  \right) \\
&  \ \ \ \ \ \ \ \ \ \ \left(  \text{since }v=1\right) \\
&  =\left(  \text{the }\left(  i+u-1\right)  \text{-th entry of }%
\underbrace{\text{the }1\text{-st column of the matrix }\left(  \left(
-1\right)  ^{p-1}C_{q}\ \mid\ \left(  I_{p}\mid C\right)  \left[
i+k:p+k\right]  \right)  }_{\substack{=\left(  \text{the }1\text{-st column of
the matrix }\left(  -1\right)  ^{p-1}C_{q}\right)  \\\text{(by the definition
of }\left(  \left(  -1\right)  ^{p-1}C_{q}\ \mid\ \left(  I_{p}\mid C\right)
\left[  i+k:p+k\right]  \right)  \text{)}}}\right) \\
&  =\left(  \text{the }\left(  i+u-1\right)  \text{-th entry of }%
\underbrace{\text{the }1\text{-st column of the matrix }\left(  -1\right)
^{p-1}C_{q}}_{\substack{=\left(  \text{the vector }\left(  -1\right)
^{p-1}C_{q}\right)  \\\text{(since the matrix }\left(  -1\right)  ^{p-1}%
C_{q}\text{ is just a column vector)}}}\right) \\
&  =\left(  \text{the }\left(  i+u-1\right)  \text{-th entry of the vector
}\left(  -1\right)  ^{p-1}C_{q}\right) \\
&  =\left(  -1\right)  ^{p-1}\cdot\left(  \text{the }\left(  i+u-1\right)
\text{-th entry of the vector }\underbrace{C_{q}}_{\substack{=\left(
\text{the }q\text{-th column of the matrix }C\right)  \\\text{(since }1\leq
q\leq q\text{)}}}\right) \\
&  =\left(  -1\right)  ^{p-1}\cdot\underbrace{\left(  \text{the }\left(
i+u-1\right)  \text{-th entry of the }q\text{-th column of the matrix
}C\right)  }_{=\left(  \text{the }\left(  i+u-1,q\right)  \text{-th entry of
the matrix }C\right)  }\\
&  =\left(  -1\right)  ^{p-1}\cdot\underbrace{\left(  \text{the }\left(
i+u-1,q\right)  \text{-th entry of the matrix }C\right)  }%
_{\substack{=x_{\left(  i+u-1,Z\left(  q\right)  \right)  }\\\text{(by
(\ref{pf.Grasp.generic.calc.entry}), applied to }\left(  i+u-1,q\right)
\text{ instead of }\left(  i,k\right)  \text{)}}}\\
&  =\left(  -1\right)  ^{p-1}\cdot x_{\left(  i+u-1,Z\left(  q\right)
\right)  }=\left(  -1\right)  ^{p-1}\cdot\underbrace{x_{\left(
i+u-1,1\right)  }}_{\substack{\in\mathbb{F}\left[  x_{\mathbf{p}\Downarrow
}\right]  \\\text{(since }\left(  i+u-1,1\right)  \in\left.  \mathbf{p}%
\Downarrow\right.  \text{)}}}\ \ \ \ \ \ \ \ \ \ \left(  \text{since }Z\left(
q\right)  =1\text{ (by the definition of }Z\text{)}\right) \\
&  \in\left(  -1\right)  ^{p-1}\cdot\mathbb{F}\left[  x_{\mathbf{p}\Downarrow
}\right]  \subseteq\mathbb{F}\left[  x_{\mathbf{p}\Downarrow}\right]  .
\end{align*}
Hence, if $v=1$, then $\eta\in\mathbb{F}\left[  x_{\mathbf{p}\Downarrow
}\right]  $ is proven. Therefore, for the rest of the proof of $\eta
\in\mathbb{F}\left[  x_{\mathbf{p}\Downarrow}\right]  $, we can WLOG assume
that we don't have $v=1$. Assume this.
\par
We don't have $v=1$, but we have $v\in\left\{  1,2,...,p-i+1\right\}  $.
Hence, $v\in\left\{  2,3,...,p-i+1\right\}  $. Now,
(\ref{pf.Grasp.generic.calc.b.Dl1}) becomes%
\begin{align}
\eta &  =\left(  \text{the }\left(  i+u-1,v\right)  \text{-th entry of the
matrix }\left(  \left(  -1\right)  ^{p-1}C_{q}\ \mid\ \left(  I_{p}\mid
C\right)  \left[  i+k:p+k\right]  \right)  \right) \nonumber\\
&  =\left(  \text{the }\left(  i+u-1\right)  \text{-th entry of }%
\underbrace{\text{the }v\text{-th column of the matrix }\left(  \left(
-1\right)  ^{p-1}C_{q}\ \mid\ \left(  I_{p}\mid C\right)  \left[
i+k:p+k\right]  \right)  }_{\substack{=\left(  \text{the }\left(  v-1\right)
\text{-th column of the matrix }\left(  I_{p}\mid C\right)  \left[
i+k:p+k\right]  \right)  \\\text{(by the definition of }\left(  \left(
-1\right)  ^{p-1}C_{q}\ \mid\ \left(  I_{p}\mid C\right)  \left[
i+k:p+k\right]  \right)  \text{, since }v\in\left\{  2,3,...,p-i+1\right\}
\text{)}}}\right) \nonumber\\
&  =\left(  \text{the }\left(  i+u-1\right)  \text{-th entry of }%
\underbrace{\text{the }\left(  v-1\right)  \text{-th column of the matrix
}\left(  I_{p}\mid C\right)  \left[  i+k:p+k\right]  }_{\substack{=\left(
I_{p}\mid C\right)  _{\left(  i+k\right)  +\left(  v-1\right)  -1}\\\text{(by
the definition of }\left(  I_{p}\mid C\right)  \left[  i+k:p+k\right]
\text{)}}}\right) \nonumber\\
&  =\left(  \text{the }\left(  i+u-1\right)  \text{-th entry of }\left(
I_{p}\mid C\right)  _{\left(  i+k\right)  +\left(  v-1\right)  -1}\right)
\nonumber\\
&  =\left(  \text{the }\left(  i+u-1\right)  \text{-th entry of }\left(
I_{p}\mid C\right)  _{i+k-2+v}\right)  . \label{pf.Grasp.generic.calc.b.Dl2}%
\end{align}
\par
Notice that $i+k-2+v=\left(  \underbrace{i}_{\geq1}-1\right)
+\underbrace{\left(  k-1\right)  }_{\geq1}+\underbrace{v}_{\geq0}\geq\left(
1-1\right)  +1+0=1$. Also, $i+\underbrace{k}_{\leq q}-2+\underbrace{v}_{\leq
p-i+1}\leq i+q-2+p-i+1=p+q-1<p+q$.
\par
Since $i+k-2+v\geq1$ and $i+k-2+v\leq p+q$, it is clear that
\[
\left(  I_{p}\mid C\right)  _{i+k-2+v}=\left(  \text{the }\left(
i+k-2+v\right)  \text{-th column of the matrix }\left(  I_{p}\mid C\right)
\right)  .
\]
\par
If $i+k-2+v\leq p$, then we thus have%
\begin{align*}
\left(  I_{p}\mid C\right)  _{i+k-2+v}  &  =\left(  \text{the }\left(
i+k-2+v\right)  \text{-th column of the matrix }\left(  I_{p}\mid C\right)
\right) \\
&  =\left(  \text{the }\left(  i+k-2+v\right)  \text{-th column of the matrix
}I_{p}\right)
\end{align*}
(by the definition of $\left(  I_{p}\mid C\right)  $, since $1\leq i+k-2+v\leq
p$). Hence, if $i+k-2+v\leq p$, then (\ref{pf.Grasp.generic.calc.b.Dl2})
becomes%
\begin{align*}
\eta &  =\left(  \text{the }\left(  i+u-1\right)  \text{-th entry of
}\underbrace{\left(  I_{p}\mid C\right)  _{i+k-2+v}}_{=\left(  \text{the
}\left(  i+k-2+v\right)  \text{-th column of the matrix }I_{p}\right)
}\right) \\
&  =\left(  \text{the }\left(  i+u-1\right)  \text{-th entry of the }\left(
i+k-2+v\right)  \text{-th column of the matrix }I_{p}\right) \\
&  \in\mathbb{F}\ \ \ \ \ \ \ \ \ \ \left(  \text{since the matrix }%
I_{p}\text{ is defined over }\mathbb{F}\right) \\
&  \subseteq\mathbb{F}\left[  x_{\mathbf{p}\Downarrow}\right]  .
\end{align*}
Thus, if $i+k-2+v\leq p$, then $\eta\in\mathbb{F}\left[  x_{\mathbf{p}%
\Downarrow}\right]  $ is proven. Hence, for the rest of the proof of $\eta
\in\mathbb{F}\left[  x_{\mathbf{p}\Downarrow}\right]  $, we can WLOG assume
that we don't have $i+k-2+v\leq p$. Assume this.
\par
We don't have $i+k-2+v\leq p$. Thus, $i+k-2+v>p$. Hence, $p<i+k-2+v\leq p+q$.
Now,%
\begin{align*}
\left(  I_{p}\mid C\right)  _{i+k-2+v}  &  =\left(  \text{the }\left(
i+k-2+v\right)  \text{-th column of the matrix }\left(  I_{p}\mid C\right)
\right) \\
&  =\left(  \text{the }\left(  i+k-2+v-p\right)  \text{-th column of the
matrix }C\right)
\end{align*}
(by the definition of $\left(  I_{p}\mid C\right)  $, since $p<i+k-2+v\leq
p+q$). Thus, (\ref{pf.Grasp.generic.calc.b.Dl2}) becomes
\begin{align}
\eta &  =\left(  \text{the }\left(  i+u-1\right)  \text{-th entry of
}\underbrace{\left(  I_{p}\mid C\right)  _{i+k-2+v}}_{=\left(  \text{the
}\left(  i+k-2+v-p\right)  \text{-th column of the matrix }C\right)  }\right)
\nonumber\\
&  =\left(  \text{the }\left(  i+u-1\right)  \text{-th entry of the }\left(
i+k-2+v-p\right)  \text{-th column of the matrix }C\right) \nonumber\\
&  =\left(  \text{the }\left(  i+u-1,i+k-2+v-p\right)  \text{-th entry of the
matrix }C\right) \nonumber\\
&  =x_{\left(  i+u-1,Z\left(  i+k-2+v-p\right)  \right)  }
\label{pf.Grasp.generic.calc.b.Dl4}%
\end{align}
(by (\ref{pf.Grasp.generic.calc.entry}), applied to $\left(
i+u-1,i+k-2+v-p\right)  $ instead of $\left(  i,k\right)  $). But since
$i+k-2+v<p+q$, we have $i+k-2+v-p<q$, so that $i+k-2+v-p\in\left\{
1,2,...,q-1\right\}  $ (because we also have $\underbrace{i+k-2+v}%
_{>p}-p>p-p=0$). Hence, (by the definition of $Z$) we have $Z\left(
i+k-2+v-p\right)  =\left(  i+k-2+v-p\right)  +1=i+k-1+v-p$. Thus,
(\ref{pf.Grasp.generic.calc.b.Dl4}) becomes%
\begin{equation}
\eta=x_{\left(  i+u-1,Z\left(  i+k-2+v-p\right)  \right)  }=x_{\left(
i+u-1,i+k-1+v-p\right)  }\ \ \ \ \ \ \ \ \ \ \left(  \text{since }Z\left(
i+k-2+v-p\right)  =i+k-1+v-p\right)  . \label{pf.Grasp.generic.calc.b.Dl5}%
\end{equation}
\par
But we have $i+\underbrace{u}_{\geq1}-1\geq i+1-1=i$ and $i+k-1+v-p=k+\left(
i-1+\underbrace{v}_{\leq p-i+1}-p\right)  \leq k+\left(  i-1+\left(
p-i+1\right)  -p\right)  =k$.
\par
Let us assume (for the sake of contradiction) that $\left(
i+u-1,i+k-1+v-p\right)  =\left(  i,k\right)  $. Then, $i+u-1=i$ and
$i+k-1+v-p=k$. Hence, $u-1=\underbrace{\left(  i+u-1\right)  }_{=i}-i=i-i=0$,
so that $u=1$. Also, subtracting $k$ from the equality $i+k-1+v-p=k$, we
obtain $i-1+v-p=0$, so that $v=p-i+1$. Combined with $u=1$, this yields
$\left(  u,v\right)  =\left(  1,p-i+1\right)  $. This contradicts $\left(
u,v\right)  \neq\left(  1,p-i+1\right)  $. This contradiction shows that our
assumption (that $\left(  i+u-1,i+k-1+v-p\right)  =\left(  i,k\right)  $) was
false. In other words, we don't have $\left(  i+u-1,i+k-1+v-p\right)  =\left(
i,k\right)  $.
\par
We have $i+k-1+v-p=Z\left(  i+k-2+v-p\right)  \in\left\{  1,2,...,q\right\}  $
(because the target of the map $Z$ is $\left\{  1,2,...,q\right\}  $). Also,
$i+u-1\geq i\geq1$ and $i+\underbrace{u}_{\leq p-i+1}-1\leq i+\left(
p-i+1\right)  -1=p$. Hence, $1\leq i+u-1\leq p$, so that $i+u-1\in\left\{
1,2,...,p\right\}  $. Combined with $i+k-1+v-p\in\left\{  1,2,...,q\right\}
$, this yields%
\[
\left(  i+u-1,i+k-1+v-p\right)  \in\left\{  1,2,...,p\right\}  \times\left\{
1,2,...,q\right\}  =\mathbf{P}.
\]
\par
The element $\left(  i+u-1,i+k-1+v-p\right)  $ is smaller or equal to $\left(
i,k\right)  $ in the poset $\mathbf{P}$ (because $i+u-1\geq i$ and
$i+k-1+v-p\leq k$). Combined with the fact that we don't have $\left(
i+u-1,i+k-1+v-p\right)  =\left(  i,k\right)  $, this shows that the element
$\left(  i+u-1,i+k-1+v-p\right)  $ is (strictly) smaller than $\left(
i,k\right)  $ in the poset $\mathbf{P}$. In other words, $\left(
i+u-1,i+k-1+v-p\right)  \vartriangleleft\left(  i,k\right)  $. Since $\left(
i,k\right)  =\mathbf{p}$, this rewrites as $\left(  i+u-1,i+k-1+v-p\right)
\vartriangleleft\mathbf{p}$.
\par
Since $\left(  i+u-1,i+k-1+v-p\right)  \in\mathbf{P}$ and $\left(
i+u-1,i+k-1+v-p\right)  \vartriangleleft\mathbf{p}$, we have%
\[
\left(  i+u-1,i+k-1+v-p\right)  \in\left\{  \mathbf{v}\in\mathbf{P}%
\ \mid\ \mathbf{v}\vartriangleleft\mathbf{p}\right\}  =\left.  \mathbf{p}%
\Downarrow\right.  ,
\]
so that $x_{\left(  i+u-1,i+k-1+v-p\right)  }\in\mathbb{F}\left[
x_{\mathbf{p}\Downarrow}\right]  $. Thus, (\ref{pf.Grasp.generic.calc.b.Dl5})
yields $\eta=x_{\left(  i+u-1,i+k-1+v-p\right)  }\in\mathbb{F}\left[
x_{\mathbf{p}\Downarrow}\right]  $.
\par
We thus have proven that $\eta\in\mathbb{F}\left[  x_{\mathbf{p}\Downarrow
}\right]  $. In other words, $\eta$ belongs to $\mathbb{F}\left[
x_{\mathbf{p}\Downarrow}\right]  $. Since $\eta$ is the $\left(  u,v\right)
$-th entry of the matrix $\mathbf{M}$, this rewrites as follows: The $\left(
u,v\right)  $-th entry of the matrix $\mathbf{M}$ belongs to $\mathbb{F}%
\left[  x_{\mathbf{p}\Downarrow}\right]  $. This proves
(\ref{pf.Grasp.generic.calc.b.Ddet3}).}. Now, Proposition
\ref{prop.cofactor.laplace} (applied to $\mathbb{F}\left[  x_{\mathbf{P}%
}\right]  $, $\mathbf{M}$ and $p-i+1$ instead of $\mathbb{K}$, $B$ and $n$)
yields%
\begin{align}
\det\mathbf{M}  &  =\sum_{\ell=1}^{p-i+1}\left(  \text{the }\left(
1,\ell\right)  \text{-th entry of the matrix }\mathbf{M}\right) \nonumber\\
&  \ \ \ \ \ \ \ \ \ \ \ \ \ \ \ \ \ \ \ \ \cdot\left(  \text{the }\left(
1,\ell\right)  \text{-th cofactor of the matrix }\mathbf{M}\right)  .
\label{pf.Grasp.generic.calc.b.Ddet4}%
\end{align}
But for every $\ell\in\left\{  1,2,...,p-i+1\right\}  $, we have%
\begin{equation}
\left(  \text{the }\left(  1,\ell\right)  \text{-th cofactor of the matrix
}\mathbf{M}\right)  \in\mathbb{F}\left[  x_{\mathbf{p}\Downarrow}\right]  .
\label{pf.Grasp.generic.calc.b.Ddet5}%
\end{equation}
\footnote{\textit{Proof of (\ref{pf.Grasp.generic.calc.b.Ddet5}):} Let
$\ell\in\left\{  1,2,...,p-i+1\right\}  $. By the definition of a cofactor, we
have%
\begin{align*}
&  \left(  \text{the }\left(  1,\ell\right)  \text{-th cofactor of the matrix
}\mathbf{M}\right) \\
&  =\left(  -1\right)  ^{\ell-1}\det\left(  \text{the matrix obtained by
deleting row }1\text{ and column }\ell\text{ from }\mathbf{M}\right)  .
\end{align*}
But consider the matrix obtained by deleting row $1$ and column $\ell$ from
$\mathbf{M}$. All entries of this matrix are entries of $\mathbf{M}$ other
than the $\left(  1,p-i+1\right)  $-th entry (because the $\left(
1,p-i+1\right)  $-th entry of $\mathbf{M}$ lies in row $1$ of $\mathbf{M}$,
which is deleted in the new matrix), and hence belong to $\mathbb{F}\left[
x_{\mathbf{p}\Downarrow}\right]  $ (according to
(\ref{pf.Grasp.generic.calc.b.Ddet3})). Hence, the determinant of this matrix
belongs to $\mathbb{F}\left[  x_{\mathbf{p}\Downarrow}\right]  $ as well. In
other words, $\det\left(  \text{the matrix obtained by deleting row }1\text{
and column }\ell\text{ from }\mathbf{M}\right)  \in\mathbb{F}\left[
x_{\mathbf{p}\Downarrow}\right]  $. Thus,%
\begin{align*}
&  \left(  \text{the }\left(  1,\ell\right)  \text{-th cofactor of the matrix
}\mathbf{M}\right) \\
&  =\left(  -1\right)  ^{\ell-1}\underbrace{\det\left(  \text{the matrix
obtained by deleting row }1\text{ and column }\ell\text{ from }\mathbf{M}%
\right)  }_{\in\mathbb{F}\left[  x_{\mathbf{p}\Downarrow}\right]  }\\
&  \in\left(  -1\right)  ^{\ell-1}\mathbb{F}\left[  x_{\mathbf{p}\Downarrow
}\right]  \subseteq\mathbb{F}\left[  x_{\mathbf{p}\Downarrow}\right]  .
\end{align*}
This proves (\ref{pf.Grasp.generic.calc.b.Ddet5}).} Thus,
(\ref{pf.Grasp.generic.calc.b.Ddet4}) becomes%
\begin{align}
\det\mathbf{M}  &  =\sum_{\ell=1}^{p-i+1}\left(  \text{the }\left(
1,\ell\right)  \text{-th entry of the matrix }\mathbf{M}\right) \nonumber\\
&  \ \ \ \ \ \ \ \ \ \ \ \ \ \ \ \ \ \ \ \ \cdot\left(  \text{the }\left(
1,\ell\right)  \text{-th cofactor of the matrix }\mathbf{M}\right) \nonumber\\
&  =\sum_{\ell=1}^{p-i}\underbrace{\left(  \text{the }\left(  1,\ell\right)
\text{-th entry of the matrix }\mathbf{M}\right)  }_{\substack{\in
\mathbb{F}\left[  x_{\mathbf{p}\Downarrow}\right]  \\\text{(by
(\ref{pf.Grasp.generic.calc.b.Ddet3}) (because the}\\\left(  1,\ell\right)
\text{-th entry of the matrix }\mathbf{M}\text{ is an entry of the matrix
}\mathbf{M}\text{ other than}\\\text{the }\left(  1,p-i+1\right)  \text{-th
entry (since }\ell\neq p-i+1\text{)))}}}\nonumber\\
&  \ \ \ \ \ \ \ \ \ \ \ \ \ \ \ \ \ \ \ \ \cdot\underbrace{\left(  \text{the
}\left(  1,\ell\right)  \text{-th cofactor of the matrix }\mathbf{M}\right)
}_{\substack{\in\mathbb{F}\left[  x_{\mathbf{p}\Downarrow}\right]  \\\text{(by
(\ref{pf.Grasp.generic.calc.b.Ddet5}))}}}\nonumber\\
&  \ \ \ \ \ \ \ \ \ \ +\underbrace{\left(  \text{the }\left(  1,p-i+1\right)
\text{-th entry of the matrix }\mathbf{M}\right)  }_{\substack{=x_{\mathbf{p}%
}\\\text{(by (\ref{pf.Grasp.generic.calc.b.Ddet2}))}}}\nonumber\\
&  \ \ \ \ \ \ \ \ \ \ \ \ \ \ \ \ \ \ \ \ \cdot\left(  \text{the }\left(
1,p-i+1\right)  \text{-th cofactor of the matrix }\mathbf{M}\right)
\nonumber\\
&  \in\underbrace{\sum_{\ell=1}^{p-i}\mathbb{F}\left[  x_{\mathbf{p}%
\Downarrow}\right]  \cdot\mathbb{F}\left[  x_{\mathbf{p}\Downarrow}\right]
}_{\substack{\subseteq\mathbb{F}\left[  x_{\mathbf{p}\Downarrow}\right]
\\\text{(since }\mathbb{F}\left[  x_{\mathbf{p}\Downarrow}\right]  \text{ is
an }\mathbb{F}\text{-algebra)}}}\nonumber\\
&  \ \ \ \ \ \ \ \ \ \ +x_{\mathbf{p}}\cdot\left(  \text{the }\left(
1,p-i+1\right)  \text{-th cofactor of the matrix }\mathbf{M}\right)
\nonumber\\
&  \subseteq\mathbb{F}\left[  x_{\mathbf{p}\Downarrow}\right]  +x_{\mathbf{p}%
}\cdot\left(  \text{the }\left(  1,p-i+1\right)  \text{-th cofactor of the
matrix }\mathbf{M}\right)  . \label{pf.Grasp.generic.calc.b.Ddet6}%
\end{align}

But $\mathbf{M}=\operatorname*{rows}\nolimits_{i,i+1,...,p}\left(  \left(
-1\right)  ^{p-1}C_{q}\ \mid\ \left(  I_{p}\mid C\right)  \left[
i+k:p+k\right]  \right)  $, so that%
\begin{align}
&  \left(  \text{the }\left(  1,p-i+1\right)  \text{-th cofactor of the matrix
}\mathbf{M}\right) \nonumber\\
&  =\left(  \text{the }\left(  1,p-i+1\right)  \text{-th cofactor of the
matrix }\operatorname*{rows}\nolimits_{i,i+1,...,p}\left(  \left(  -1\right)
^{p-1}C_{q}\ \mid\ \left(  I_{p}\mid C\right)  \left[  i+k:p+k\right]
\right)  \right) \nonumber\\
&  =\left(  -1\right)  ^{p-i}\det\left(  \operatorname*{rows}%
\nolimits_{i+1,i+2,...,p}\left(  \left(  -1\right)  ^{p-1}C_{q}\ \mid\ \left(
I_{p}\mid C\right)  \left[  i+k:p+k-1\right]  \right)  \right)
\label{pf.Grasp.generic.calc.b.Ddet8}\\
&  \ \ \ \ \ \ \ \ \ \ \left(
\begin{array}
[c]{c}%
\text{by Corollary \ref{cor.cofactor.2}, applied to }\mathbb{F}\left[
x_{\mathbf{P}}\right]  \text{, }p+q\text{, }\left(  I_{p}\mid C\right)
\text{, }\left(  -1\right)  ^{p-1}C_{q}\text{, }i+k\text{ and }p+k\\
\text{instead of }\mathbb{K}\text{, }q\text{, }A\text{, }\xi\text{, }a\text{
and }b\text{ (since }\underbrace{i}_{\leq p-1<p}+k<p+k\text{,}\\
\left(  p+k\right)  -\left(  i+k\right)  =p-\underbrace{i}_{\geq1}\leq
p-1<p\text{ and }i=p-\left(  \left(  p+k\right)  -\left(  i+k\right)  \right)
\text{)}%
\end{array}
\right)  .\nonumber
\end{align}
Multiplying this equality with $\left(  -1\right)  ^{p-i}$, we obtain%
\begin{align}
&  \left(  -1\right)  ^{p-i}\left(  \text{the }\left(  1,p-i+1\right)
\text{-th cofactor of the matrix }\mathbf{M}\right) \nonumber\\
&  =\underbrace{\left(  -1\right)  ^{p-i}\left(  -1\right)  ^{p-i}%
}_{\substack{=\left(  -1\right)  ^{\left(  p-i\right)  +\left(  p-i\right)
}=1\\\text{(since }\left(  p-i\right)  +\left(  p-i\right)  =2\left(
p-i\right)  \\\text{is even)}}}\det\left(  \operatorname*{rows}%
\nolimits_{i+1,i+2,...,p}\left(  \left(  -1\right)  ^{p-1}C_{q}\ \mid\ \left(
I_{p}\mid C\right)  \left[  i+k:p+k-1\right]  \right)  \right) \nonumber\\
&  =\det\left(  \operatorname*{rows}\nolimits_{i+1,i+2,...,p}\left(  \left(
-1\right)  ^{p-1}C_{q}\ \mid\ \left(  I_{p}\mid C\right)  \left[
i+k:p+k-1\right]  \right)  \right)  . \label{pf.Grasp.generic.calc.b.Ddet8a}%
\end{align}

On the other hand, combining $i+1\in\left\{  2,3,...,p\right\}  \subseteq
\left\{  1,2,...,p\right\}  $ with $k-1\in\left\{  1,2,...,q-1\right\}
\subseteq\left\{  1,2,...,q\right\}  $, we obtain $\left(  i+1,k-1\right)
\in\left\{  1,2,...,p\right\}  \times\left\{  1,2,...,q\right\}  =\mathbf{P}$.
Applying (\ref{lem.Grasp.generic.calc.a.D}) to $\left(  r,s\right)  =\left(
i+1,k-1\right)  $, we obtain%
\begin{align}
&  \mathfrak{D}_{\left(  i+1,k-1\right)  }\nonumber\\
&  =\underbrace{\left(  -1\right)  ^{\left(  i+1\right)  -1}}_{=\left(
-1\right)  ^{i}}\det\left(  \operatorname*{rows}\nolimits_{i+1,i+2,...,p}%
\left(  \left(  -1\right)  ^{p-1}C_{q}\ \mid\ \left(  I_{p}\mid C\right)
\left[  \underbrace{\left(  i+1\right)  +\left(  k-1\right)  }_{=i+k}%
:\underbrace{p+\left(  k-1\right)  }_{=p+k-1}\right]  \right)  \right)
\nonumber\\
&  =\left(  -1\right)  ^{i}\underbrace{\det\left(  \operatorname*{rows}%
\nolimits_{i+1,i+2,...,p}\left(  \left(  -1\right)  ^{p-1}C_{q}\ \mid\ \left(
I_{p}\mid C\right)  \left[  i+k:p+k-1\right]  \right)  \right)  }%
_{\substack{=\left(  -1\right)  ^{p-i}\left(  \text{the }\left(
1,p-i+1\right)  \text{-th cofactor of the matrix }\mathbf{M}\right)
\\\text{(by (\ref{pf.Grasp.generic.calc.b.Ddet8a}))}}}\nonumber\\
&  =\underbrace{\left(  -1\right)  ^{i}\left(  -1\right)  ^{p-i}}_{=\left(
-1\right)  ^{i+\left(  p-i\right)  }=\left(  -1\right)  ^{p}}\left(  \text{the
}\left(  1,p-i+1\right)  \text{-th cofactor of the matrix }\mathbf{M}\right)
\nonumber\\
&  =\left(  -1\right)  ^{p}\underbrace{\left(  \text{the }\left(
1,p-i+1\right)  \text{-th cofactor of the matrix }\mathbf{M}\right)
}_{\substack{\in\mathbb{F}\left[  x_{\mathbf{p}\Downarrow}\right]  \\\text{(by
(\ref{pf.Grasp.generic.calc.b.Ddet5}), applied to }\ell=p-i+1\text{)}%
}}\label{pf.Grasp.generic.calc.b.Ddet9}\\
&  \in\left(  -1\right)  ^{p}\mathbb{F}\left[  x_{\mathbf{p}\Downarrow
}\right]  \subseteq\mathbb{F}\left[  x_{\mathbf{p}\Downarrow}\right]
=\mathbb{F}\left[  x_{\left(  i,k\right)  \Downarrow}\right]
\ \ \ \ \ \ \ \ \ \ \left(  \text{since }\mathbf{p}=\left(  i,k\right)
\right)  .\nonumber
\end{align}
Thus, (\ref{lem.Grasp.generic.calc.b.D2}) is proven. Also, dividing the
equality (\ref{pf.Grasp.generic.calc.b.Ddet9}) by $\left(  -1\right)  ^{p}$,
we obtain%
\[
\dfrac{1}{\left(  -1\right)  ^{p}}\mathfrak{D}_{\left(  i+1,k-1\right)
}=\left(  \text{the }\left(  1,p-i+1\right)  \text{-th cofactor of the matrix
}\mathbf{M}\right)  .
\]
Since $\dfrac{1}{\left(  -1\right)  ^{p}}=\left(  -1\right)  ^{p}$, this
simplifies to
\begin{equation}
\left(  -1\right)  ^{p}\mathfrak{D}_{\left(  i+1,k-1\right)  }=\left(
\text{the }\left(  1,p-i+1\right)  \text{-th cofactor of the matrix
}\mathbf{M}\right)  . \label{pf.Grasp.generic.calc.b.Ddet10}%
\end{equation}

Now, (\ref{pf.Grasp.generic.calc.b.Ddet6}) becomes%
\begin{align}
&  \det\mathbf{M}\nonumber\\
&  \in\mathbb{F}\left[  x_{\mathbf{p}\Downarrow}\right]  +x_{\mathbf{p}}%
\cdot\underbrace{\left(  \text{the }\left(  1,p-i+1\right)  \text{-th cofactor
of the matrix }\mathbf{M}\right)  }_{\substack{=\left(  -1\right)
^{p}\mathfrak{D}_{\left(  i+1,k-1\right)  }\\\text{(by
(\ref{pf.Grasp.generic.calc.b.Ddet10}))}}}\nonumber\\
&  =\mathbb{F}\left[  x_{\mathbf{p}\Downarrow}\right]  +x_{\mathbf{p}}%
\cdot\left(  -1\right)  ^{p}\mathfrak{D}_{\left(  i+1,k-1\right)  }%
=\mathbb{F}\left[  x_{\mathbf{p}\Downarrow}\right]  +\left(  -1\right)
^{p}\cdot\mathfrak{D}_{\left(  i+1,k-1\right)  }\cdot x_{\mathbf{p}}=\left(
-1\right)  ^{p}\cdot\mathfrak{D}_{\left(  i+1,k-1\right)  }\cdot
x_{\mathbf{p}}+\mathbb{F}\left[  x_{\mathbf{p}\Downarrow}\right] \nonumber\\
&  =\left(  -1\right)  ^{p}\cdot\mathfrak{D}_{\left(  i+1,k-1\right)  }\cdot
x_{\left(  i,k\right)  }+\mathbb{F}\left[  x_{\left(  i,k\right)  \Downarrow
}\right]  \label{pf.Grasp.generic.calc.b.D38}%
\end{align}
(since $\mathbf{p}=\left(  i,k\right)  $). But applying
(\ref{lem.Grasp.generic.calc.a.D}) to $\left(  r,s\right)  =\left(
i,k\right)  $, we obtain%
\begin{align*}
\mathfrak{D}_{\left(  i,k\right)  }  &  =\left(  -1\right)  ^{i-1}\det\left(
\underbrace{\operatorname*{rows}\nolimits_{i,i+1,...,p}\left(  \left(
-1\right)  ^{p-1}C_{q}\ \mid\ \left(  I_{p}\mid C\right)  \left[
i+k:p+k\right]  \right)  }_{=\mathbf{M}}\right) \\
&  =\left(  -1\right)  ^{i-1}\det\mathbf{M}\in\left(  -1\right)  ^{i-1}\left(
\left(  -1\right)  ^{p}\cdot\mathfrak{D}_{\left(  i+1,k-1\right)  }\cdot
x_{\left(  i,k\right)  }+\mathbb{F}\left[  x_{\left(  i,k\right)  \Downarrow
}\right]  \right)  \ \ \ \ \ \ \ \ \ \ \left(  \text{by
(\ref{pf.Grasp.generic.calc.b.D38})}\right) \\
&  \subseteq\underbrace{\left(  -1\right)  ^{i-1}\left(  -1\right)  ^{p}%
}_{\substack{=\left(  -1\right)  ^{i-1+p}=\left(  -1\right)  ^{p-i+1}%
\\\text{(since }i-1+p\equiv p-i+1\operatorname{mod}2\text{)}}}\cdot
\mathfrak{D}_{\left(  i+1,k-1\right)  }\cdot x_{\left(  i,k\right)
}+\underbrace{\left(  -1\right)  ^{i-1}\mathbb{F}\left[  x_{\left(
i,k\right)  \Downarrow}\right]  }_{\subseteq\mathbb{F}\left[  x_{\left(
i,k\right)  \Downarrow}\right]  }\\
&  \subseteq\left(  -1\right)  ^{p-i+1}\cdot\mathfrak{D}_{\left(
i+1,k-1\right)  }\cdot x_{\left(  i,k\right)  }+\mathbb{F}\left[  x_{\left(
i,k\right)  \Downarrow}\right]  .
\end{align*}
Thus, (\ref{lem.Grasp.generic.calc.b.D}) is proven.

Thus, the proof of Lemma \ref{lem.Grasp.generic.calc} \textbf{(b)} is complete.

\textbf{(c)} Let $\left(  i,k\right)  \in\mathbf{P}$ be such that $i\neq p$
and $k\neq1$.

We have $\left(  i,k\right)  \in\mathbf{P}=\left\{  1,2,...,p\right\}
\times\left\{  1,2,...,q\right\}  $. Hence, $i\in\left\{  1,2,...,p\right\}  $
and $k\in\left\{  1,2,...,q\right\}  $. Since $i\in\left\{  1,2,...,p\right\}
$ and $i\neq p$, we have $i\in\left\{  1,2,...,p-1\right\}  $, so that $1\leq
i\leq p-1$. Since $k\in\left\{  1,2,...,q\right\}  $ and $k\neq1$, we have
$k\in\left\{  2,3,...,q\right\}  $. Hence, $2\leq k\leq q$. Also, $k\geq2$, so
that $k-1\geq1$. Since $i\in\left\{  1,2,...,p-1\right\}  $, we have
$i+1\in\left\{  2,3,...,p\right\}  $.

We have $1\leq i\leq i+1$ and $\underbrace{i}_{\leq p-1}+k\leq\left(
p-1\right)  +k=p+k-1\leq\left(  p+k-1\right)  +1$ and $i-1+\left(
p+k-1\right)  -\left(  i+k\right)  =p-2$. Therefore, we can apply Theorem
\ref{thm.pluecker.ptolemy} to $p$, $p+q$, $\left(  I_{p}\mid C\right)  $, $1$,
$i$, $i+k$ and $p+k-1$ instead of $u$, $v$, $A$, $a$, $b$, $c$ and $d$. As a
result, we obtain%
\begin{align*}
&  \det\left(  \left(  I_{p}\mid C\right)  \left[  1-1:i\mid i+k:\left(
p+k-1\right)  +1\right]  \right) \\
&  \ \ \ \ \ \ \ \ \ \ \ \ \ \ \ \ \ \ \ \ \cdot\det\left(  \left(  I_{p}\mid
C\right)  \left[  1:i+1\mid\left(  i+k\right)  -1:p+k-1\right]  \right) \\
&  \ \ \ \ \ \ \ \ \ \ +\det\left(  \left(  I_{p}\mid C\right)  \left[
1:i\mid\left(  i+k\right)  -1:\left(  p+k-1\right)  +1\right]  \right) \\
&  \ \ \ \ \ \ \ \ \ \ \ \ \ \ \ \ \ \ \ \ \cdot\det\left(  \left(  I_{p}\mid
C\right)  \left[  1-1:i+1\mid i+k:p+k-1\right]  \right) \\
&  =\det\left(  \left(  I_{p}\mid C\right)  \left[  1-1:i\mid\left(
i+k\right)  -1:p+k-1\right]  \right) \\
&  \ \ \ \ \ \ \ \ \ \ \ \ \ \ \ \ \ \ \ \ \cdot\det\left(  \left(  I_{p}\mid
C\right)  \left[  1:i+1\mid i+k:\left(  p+k-1\right)  +1\right]  \right)  .
\end{align*}
Since $1-1=0$ and $\left(  p+k-1\right)  +1=p+k$, this rewrites as follows:%
\begin{align}
&  \det\left(  \left(  I_{p}\mid C\right)  \left[  0:i\mid i+k:p+k\right]
\right)  \cdot\det\left(  \left(  I_{p}\mid C\right)  \left[  1:i+1\mid\left(
i+k\right)  -1:p+k-1\right]  \right) \nonumber\\
&  \ \ \ \ \ \ \ \ \ \ +\det\left(  \left(  I_{p}\mid C\right)  \left[
1:i\mid\left(  i+k\right)  -1:p+k\right]  \right)  \cdot\det\left(  \left(
I_{p}\mid C\right)  \left[  0:i+1\mid i+k:p+k-1\right]  \right) \nonumber\\
&  =\det\left(  \left(  I_{p}\mid C\right)  \left[  0:i\mid\left(  i+k\right)
-1:p+k-1\right]  \right)  \cdot\det\left(  \left(  I_{p}\mid C\right)  \left[
1:i+1\mid i+k:p+k\right]  \right)  . \label{pf.Grasp.generic.calc.plu}%
\end{align}

But combining $i+1\in\left\{  2,3,...,p\right\}  \subseteq\left\{
1,2,...,p\right\}  $ with $k-1\in\left\{  1,2,...,q-1\right\}  \subseteq
\left\{  1,2,...,q\right\}  $, we obtain $\left(  i+1,k-1\right)  \in\left\{
1,2,...,p\right\}  \times\left\{  1,2,...,q\right\}  =\mathbf{P}$. Thus,
applying (\ref{lem.Grasp.generic.independency.Ndef}) to $\left(  r,s\right)
=\left(  i+1,k-1\right)  $, we obtain%
\begin{align}
\mathfrak{N}_{\left(  i+1,k-1\right)  }  &  =\det\left(  \left(  I_{p}\mid
C\right)  \left[  1:i+1\mid\underbrace{\left(  i+1\right)  +\left(
k-1\right)  -1}_{=\left(  i+k\right)  -1}:\underbrace{p+\left(  k-1\right)
}_{=p+k-1}\right]  \right) \nonumber\\
&  =\det\left(  \left(  I_{p}\mid C\right)  \left[  1:i+1\mid\left(
i+k\right)  -1:p+k-1\right]  \right)  . \label{pf.Grasp.generic.calc.c.2}%
\end{align}
Also, applying (\ref{lem.Grasp.generic.independency.Ddef}) to $\left(
r,s\right)  =\left(  i+1,k-1\right)  $, we obtain%
\begin{align}
\mathfrak{D}_{\left(  i+1,k-1\right)  }  &  =\det\left(  \left(  I_{p}\mid
C\right)  \left[  0:i+1\mid\underbrace{\left(  i+1\right)  +\left(
k-1\right)  }_{=i+k}:\underbrace{p+\left(  k-1\right)  }_{=p+k-1}\right]
\right) \nonumber\\
&  =\det\left(  \left(  I_{p}\mid C\right)  \left[  0:i+1\mid
i+k:p+k-1\right]  \right)  . \label{pf.Grasp.generic.calc.c.3}%
\end{align}
On the other hand, combining $i\in\left\{  1,2,...,p\right\}  $ with
$k-1\in\left\{  1,2,...,q\right\}  $, we obtain $\left(  i,k-1\right)
\in\left\{  1,2,...,p\right\}  \times\left\{  1,2,...,q\right\}  =\mathbf{P}$.
Thus, applying (\ref{lem.Grasp.generic.independency.Ddef}) to $\left(
r,s\right)  =\left(  i,k-1\right)  $, we obtain%
\begin{align}
\mathfrak{D}_{\left(  i,k-1\right)  }  &  =\det\left(  \left(  I_{p}\mid
C\right)  \left[  0:i\mid\underbrace{i+\left(  k-1\right)  }_{=\left(
i+k\right)  -1}:\underbrace{p+\left(  k-1\right)  }_{=p+k-1}\right]  \right)
\nonumber\\
&  =\det\left(  \left(  I_{p}\mid C\right)  \left[  0:i\mid\left(  i+k\right)
-1:p+k-1\right]  \right)  . \label{pf.Grasp.generic.calc.c.4}%
\end{align}
Finally, combining $i+1\in\left\{  1,2,...,p\right\}  $ with $k\in\left\{
1,2,...,q\right\}  $, we obtain $\left(  i+1,k\right)  \in\left\{
1,2,...,p\right\}  \times\left\{  1,2,...,q\right\}  =\mathbf{P}$. Thus,
applying (\ref{lem.Grasp.generic.independency.Ndef}) to $\left(  r,s\right)
=\left(  i+1,k\right)  $, we obtain%
\begin{align}
\mathfrak{N}_{\left(  i+1,k\right)  }  &  =\det\left(  \left(  I_{p}\mid
C\right)  \left[  1:i+1\mid\underbrace{\left(  i+1\right)  +k-1}%
_{=i+k}:p+k\right]  \right) \nonumber\\
&  =\det\left(  \left(  I_{p}\mid C\right)  \left[  1:i+1\mid i+k:p+k\right]
\right)  . \label{pf.Grasp.generic.calc.c.5}%
\end{align}

Now,
\begin{align*}
&  \mathfrak{N}_{\left(  i,k\right)  }\mathfrak{D}_{\left(  i+1,k-1\right)
}+\mathfrak{D}_{\left(  i,k\right)  }\mathfrak{N}_{\left(  i+1,k-1\right)  }\\
&  =\underbrace{\mathfrak{D}_{\left(  i,k\right)  }}_{\substack{=\det\left(
\left(  I_{p}\mid C\right)  \left[  0:i\mid i+k:p+k\right]  \right)
\\\text{(by (\ref{lem.Grasp.generic.independency.Ddef}))}}%
}\underbrace{\mathfrak{N}_{\left(  i+1,k-1\right)  }}_{\substack{=\det\left(
\left(  I_{p}\mid C\right)  \left[  1:i+1\mid\left(  i+k\right)
-1:p+k-1\right]  \right)  \\\text{(by (\ref{pf.Grasp.generic.calc.c.2}))}}}\\
&  \ \ \ \ \ \ \ \ \ \ +\underbrace{\mathfrak{N}_{\left(  i,k\right)  }%
}_{\substack{=\det\left(  \left(  I_{p}\mid C\right)  \left[  1:i\mid
i+k-1:p+k\right]  \right)  \\\text{(by
(\ref{lem.Grasp.generic.independency.Ndef}))}}}\underbrace{\mathfrak{D}%
_{\left(  i+1,k-1\right)  }}_{\substack{=\det\left(  \left(  I_{p}\mid
C\right)  \left[  0:i+1\mid i+k:p+k-1\right]  \right)  \\\text{(by
(\ref{pf.Grasp.generic.calc.c.3}))}}}\\
&  =\det\left(  \left(  I_{p}\mid C\right)  \left[  0:i\mid i+k:p+k\right]
\right)  \cdot\det\left(  \left(  I_{p}\mid C\right)  \left[  1:i+1\mid\left(
i+k\right)  -1:p+k-1\right]  \right) \\
&  \ \ \ \ \ \ \ \ \ \ +\det\left(  \left(  I_{p}\mid C\right)  \left[
1:i\mid i+k-1:p+k\right]  \right)  \cdot\det\left(  \left(  I_{p}\mid
C\right)  \left[  0:i+1\mid i+k:p+k-1\right]  \right) \\
&  =\underbrace{\det\left(  \left(  I_{p}\mid C\right)  \left[  0:i\mid\left(
i+k\right)  -1:p+k-1\right]  \right)  }_{\substack{=\mathfrak{D}_{\left(
i,k-1\right)  }\\\text{(by (\ref{pf.Grasp.generic.calc.c.4}))}}}\cdot
\underbrace{\det\left(  \left(  I_{p}\mid C\right)  \left[  1:i+1\mid
i+k:p+k\right]  \right)  }_{\substack{=\mathfrak{N}_{\left(  i+1,k\right)
}\\\text{(by (\ref{pf.Grasp.generic.calc.c.5}))}}}\\
&  \ \ \ \ \ \ \ \ \ \ \left(  \text{by (\ref{pf.Grasp.generic.calc.plu}%
)}\right) \\
&  =\mathfrak{D}_{\left(  i,k-1\right)  }\mathfrak{N}_{\left(  i+1,k\right)
}.
\end{align*}
This proves Lemma \ref{lem.Grasp.generic.calc} \textbf{(c)}.

\textbf{(d)} In the following, we are going to regard column vectors of size
$\ell$ as $\ell\times1$-matrices for every $\ell\in\mathbb{N}$. We will also
regard scalars as $1\times1$-matrices.

Let $i\in\left\{  1,2,...,p\right\}  $ be arbitrary. Then, $i-1\in\left\{
0,1,...,p-1\right\}  \subseteq\left\{  0,1,...,p\right\}  $, and%
\begin{align*}
&  \left(  I_{p}\mid C\right)  \left[  \underbrace{i+1-1}_{=\left(
i-1\right)  +1}:p+1\right] \\
&  =\left(  I_{p}\mid C\right)  \left[  \left(  i-1\right)  +1:p+1\right]
=\operatorname*{cols}\nolimits_{\left(  i-1\right)  +1,\left(  i-1\right)
+2,...,p}\left(  I_{p}\right) \\
&  \ \ \ \ \ \ \ \ \ \ \left(
\begin{array}
[c]{c}%
\text{by Corollary \ref{cor.cutoffI} \textbf{(a)} (applied to }\mathbb{F}%
\left[  x_{\mathbf{P}}\right]  \text{, }p\text{, }p\text{, }q\text{, }%
I_{p}\text{, }C\text{ and }i-1\\
\text{instead of }\mathbb{K}\text{, }u\text{, }v_{1}\text{, }v_{2}\text{,
}A_{1}\text{, }A_{2}\text{ and }r\text{)}%
\end{array}
\right) \\
&  =\operatorname*{cols}\nolimits_{i,i+1,...,p}\left(  I_{p}\right)  =\left(
\begin{array}
[c]{c}%
0_{\left(  i-1\right)  \times\left(  p-i+1\right)  }\\
I_{p-i+1}%
\end{array}
\right)
\end{align*}
(where $\left(
\begin{array}
[c]{c}%
0_{\left(  i-1\right)  \times\left(  p-i+1\right)  }\\
I_{p-i+1}%
\end{array}
\right)  $ is to be understood as a block matrix). By
(\ref{lem.Grasp.generic.calc.a.N}) (applied to $\left(  r,s\right)  =\left(
i,1\right)  $), we have
\begin{align*}
\mathfrak{N}_{\left(  i,1\right)  }  &  =\det\left(  \operatorname*{rows}%
\nolimits_{i,i+1,...,p}\left(  \underbrace{\left(  I_{p}\mid C\right)  \left[
i+1-1:p+1\right]  }_{=\left(
\begin{array}
[c]{c}%
0_{\left(  i-1\right)  \times\left(  p-i+1\right)  }\\
I_{p-i+1}%
\end{array}
\right)  }\right)  \right) \\
&  =\det\left(  \underbrace{\operatorname*{rows}\nolimits_{i,i+1,...,p}\left(
\begin{array}
[c]{c}%
0_{\left(  i-1\right)  \times\left(  p-i+1\right)  }\\
I_{p-i+1}%
\end{array}
\right)  }_{=I_{p-i+1}}\right)  =\det\left(  I_{p-i+1}\right)  =1.
\end{align*}
This proves (\ref{lem.Grasp.generic.calc.e.N}).

Also, we have $i\in\left\{  1,2,...,p\right\}  \subseteq\left\{
0,1,...,p\right\}  $, and%
\begin{align*}
&  \left(  I_{p}\mid C\right)  \left[  i+1:p+1\right]  =\operatorname*{cols}%
\nolimits_{i+1,i+2,...,p}\left(  I_{p}\right) \\
&  \ \ \ \ \ \ \ \ \ \ \left(
\begin{array}
[c]{c}%
\text{by Corollary \ref{cor.cutoffI} \textbf{(a)} (applied to }\mathbb{F}%
\left[  x_{\mathbf{P}}\right]  \text{, }p\text{, }p\text{, }q\text{, }%
I_{p}\text{, }C\text{ and }i\\
\text{instead of }\mathbb{K}\text{, }u\text{, }v_{1}\text{, }v_{2}\text{,
}A_{1}\text{, }A_{2}\text{ and }r\text{)}%
\end{array}
\right) \\
&  =\left(
\begin{array}
[c]{c}%
0_{i\times\left(  p-i\right)  }\\
I_{p-i}%
\end{array}
\right)
\end{align*}
(where $\left(
\begin{array}
[c]{c}%
0_{i\times\left(  p-i\right)  }\\
I_{p-i}%
\end{array}
\right)  $ is to be understood as a block matrix). Thus,%
\begin{align*}
&  \operatorname*{rows}\nolimits_{i,i+1,...,p}\left(  \left(  -1\right)
^{p-1}C_{q}\ \mid\ \underbrace{\left(  I_{p}\mid C\right)  \left[
i+1:p+1\right]  }_{=\left(
\begin{array}
[c]{c}%
0_{i\times\left(  p-i\right)  }\\
I_{p-i}%
\end{array}
\right)  }\right) \\
&  =\operatorname*{rows}\nolimits_{i,i+1,...,p}\left(  \left(  -1\right)
^{p-1}C_{q}\ \mid\ \left(
\begin{array}
[c]{c}%
0_{i\times\left(  p-i\right)  }\\
I_{p-i}%
\end{array}
\right)  \right) \\
&  =\left(  \operatorname*{rows}\nolimits_{i,i+1,...,p}\left(  \left(
-1\right)  ^{p-1}C_{q}\right)  \ \mid\ \underbrace{\operatorname*{rows}%
\nolimits_{i,i+1,...,p}\left(
\begin{array}
[c]{c}%
0_{i\times\left(  p-i\right)  }\\
I_{p-i}%
\end{array}
\right)  }_{=\left(
\begin{array}
[c]{c}%
0_{1\times\left(  p-i\right)  }\\
I_{p-i}%
\end{array}
\right)  }\right) \\
&  \ \ \ \ \ \ \ \ \ \ \left(
\begin{array}
[c]{c}%
\text{by Proposition \ref{prop.attach.rows}, applied to }\mathbb{F}\left[
x_{\mathbf{P}}\right]  \text{, }p\text{, }1\text{, }p-i\text{, }\left(
-1\right)  ^{p-1}C_{q}\text{, }\left(
\begin{array}
[c]{c}%
0_{i\times\left(  p-i\right)  }\\
I_{p-i}%
\end{array}
\right)  \text{,}\\
p-i+1\text{ and }\left(  i,i+1,...,p\right)  \text{ instead of }%
\mathbb{K}\text{, }u\text{, }v_{1}\text{, }v_{2}\text{, }A_{1}\text{, }%
A_{2}\text{, }k\text{ and }\left(  i_{1},i_{2},...,i_{k}\right)
\end{array}
\right) \\
&  =\left(  \operatorname*{rows}\nolimits_{i,i+1,...,p}\left(  \left(
-1\right)  ^{p-1}C_{q}\right)  \ \mid\ \left(
\begin{array}
[c]{c}%
0_{1\times\left(  p-i\right)  }\\
I_{p-i}%
\end{array}
\right)  \right)  ,
\end{align*}
so that%
\begin{align*}
&  \det\left(  \underbrace{\operatorname*{rows}\nolimits_{i,i+1,...,p}\left(
\left(  -1\right)  ^{p-1}C_{q}\ \mid\ \left(  I_{p}\mid C\right)  \left[
i+1:p+1\right]  \right)  }_{=\left(  \operatorname*{rows}%
\nolimits_{i,i+1,...,p}\left(  \left(  -1\right)  ^{p-1}C_{q}\right)
\ \mid\ \left(
\begin{array}
[c]{c}%
0_{1\times\left(  p-i\right)  }\\
I_{p-i}%
\end{array}
\right)  \right)  }\right) \\
&  =\det\left(  \operatorname*{rows}\nolimits_{i,i+1,...,p}\left(  \left(
-1\right)  ^{p-1}C_{q}\right)  \ \mid\ \left(
\begin{array}
[c]{c}%
0_{1\times\left(  p-i\right)  }\\
I_{p-i}%
\end{array}
\right)  \right) \\
&  =\left(  -1\right)  ^{p-i}\det\left(  \left(
\begin{array}
[c]{c}%
0_{1\times\left(  p-i\right)  }\\
I_{p-i}%
\end{array}
\right)  \ \mid\ \operatorname*{rows}\nolimits_{i,i+1,...,p}\left(  \left(
-1\right)  ^{p-1}C_{q}\right)  \right) \\
&  \ \ \ \ \ \ \ \ \ \ \left(
\begin{array}
[c]{c}%
\text{since permuting the columns of a matrix multiplies its determinant}\\
\text{by the sign of the permutation}%
\end{array}
\right) \\
&  =\left(  -1\right)  ^{\left(  p-i+1\right)  -1}\underbrace{\det\left(
\left(
\begin{array}
[c]{c}%
0_{1\times\left(  \left(  p-i+1\right)  -1\right)  }\\
I_{\left(  p-i+1\right)  -1}%
\end{array}
\right)  \ \mid\ \operatorname*{rows}\nolimits_{i,i+1,...,p}\left(  \left(
-1\right)  ^{p-1}C_{q}\right)  \right)  }_{\substack{=\left(  -1\right)
^{1\left(  \left(  p-i+1\right)  -1\right)  }\det\left(  \operatorname*{rows}%
\nolimits_{1}\left(  \operatorname*{rows}\nolimits_{i,i+1,...,p}\left(
\left(  -1\right)  ^{p-1}C_{q}\right)  \right)  \right)  \\\text{(by Lemma
\ref{lem.attach.matrices} \textbf{(b)}, applied to}\\\mathbb{F}\left[
x_{\mathbf{P}}\right]  \text{, }p-i+1\text{, }1\text{ and }%
\operatorname*{rows}\nolimits_{i,i+1,...,p}\left(  \left(  -1\right)
^{p-1}C_{q}\right)  \\\text{instead of }\mathbb{K}\text{, }p\text{, }k\text{
and }U\text{)}}}\\
&  \ \ \ \ \ \ \ \ \ \ \left(  \text{since }p-i=\left(  p-i+1\right)
-1\right) \\
&  =\underbrace{\left(  -1\right)  ^{\left(  p-i+1\right)  -1}\left(
-1\right)  ^{1\left(  \left(  p-i+1\right)  -1\right)  }}_{\substack{=\left(
-1\right)  ^{p-i}\left(  -1\right)  ^{p-i}=\left(  -1\right)  ^{\left(
p-i\right)  +\left(  p-i\right)  }=1\\\text{(since }\left(  p-i\right)
+\left(  p-i\right)  =2\left(  p-i\right)  \text{ is even)}}}\det\left(
\underbrace{\operatorname*{rows}\nolimits_{1}\left(  \operatorname*{rows}%
\nolimits_{i,i+1,...,p}\left(  \left(  -1\right)  ^{p-1}C_{q}\right)  \right)
}_{\substack{=\operatorname*{rows}\nolimits_{i}\left(  \left(  -1\right)
^{p-1}C_{q}\right)  \\\text{(by Proposition \ref{prop.rows.row}, applied to
}\mathbb{F}\left[  x_{\mathbf{P}}\right]  \text{, }1\text{, }\left(
-1\right)  ^{p-1}C_{q}\text{,}\\p-i+1\text{, }\left(  i,i+1,...,p\right)
\text{ and }1\text{ instead of }\mathbb{K}\text{, }q\text{, }B\text{,
}k\text{, }\left(  i_{1},i_{2},...,i_{k}\right)  \text{ and }j\text{)}%
}}\right)
\end{align*}%
\begin{align*}
&  =\det\left(  \underbrace{\operatorname*{rows}\nolimits_{i}\left(  \left(
-1\right)  ^{p-1}C_{q}\right)  }_{\substack{=\left(  \text{the }i\text{-th
entry of }\left(  -1\right)  ^{p-1}C_{q}\right)  \\\text{(by Proposition
\ref{prop.rows.entry}, applied to }\mathbb{F}\left[  x_{\mathbf{P}}\right]
\text{, }p\\\text{and }\left(  -1\right)  ^{p-1}C_{q}\text{ instead of
}\mathbb{K}\text{, }\ell\text{ and }U\text{)}}}\right) \\
&  =\det\left(  \text{the }i\text{-th entry of }\left(  -1\right)  ^{p-1}%
C_{q}\right)  =\left(  \text{the }i\text{-th entry of }\left(  -1\right)
^{p-1}C_{q}\right) \\
&  \ \ \ \ \ \ \ \ \ \ \left(
\begin{array}
[c]{c}%
\text{since the determinant of a scalar regarded as a }1\times1\text{-matrix}%
\\
\text{is simply this scalar}%
\end{array}
\right) \\
&  =\left(  -1\right)  ^{p-1}\left(  \text{the }i\text{-th entry of
}\underbrace{C_{q}}_{\substack{=\left(  \text{the }q\text{-th column of the
matrix }C\right)  \\\text{(since }1\leq q\leq q\text{)}}}\right) \\
&  =\left(  -1\right)  ^{p-1}\underbrace{\left(  \text{the }i\text{-th entry
of the }q\text{-th column of the matrix }C\right)  }_{\substack{=\left(
\text{the }\left(  i,q\right)  \text{-th entry of the matrix }C\right)
=x_{\left(  i,Z\left(  q\right)  \right)  }\\\text{(by
(\ref{pf.Grasp.generic.calc.entry}), applied to }\left(  i,q\right)
\text{ instead of }\left(  i,k\right)  \text{)}}}\\
&  =\left(  -1\right)  ^{p-1}x_{\left(  i,Z\left(  q\right)  \right)
}=\left(  -1\right)  ^{p-1}x_{\left(  i,1\right)  }\ \ \ \ \ \ \ \ \ \ \left(
\text{since }Z\left(  q\right)  =1\text{ (by the definition of }%
Z\text{)}\right)  .
\end{align*}
By (\ref{lem.Grasp.generic.calc.a.D}) (applied to $\left(  r,s\right)
=\left(  i,1\right)  $), we have
\begin{align*}
\mathfrak{D}_{\left(  i,1\right)  }  &  =\left(  -1\right)  ^{i-1}%
\underbrace{\det\left(  \operatorname*{rows}\nolimits_{i,i+1,...,p}\left(
\left(  -1\right)  ^{p-1}C_{q}\ \mid\ \left(  I_{p}\mid C\right)  \left[
i+1:p+1\right]  \right)  \right)  }_{=\left(  -1\right)  ^{p-1}x_{\left(
i,1\right)  }}\\
&  =\underbrace{\left(  -1\right)  ^{i-1}\left(  -1\right)  ^{p-1}%
}_{\substack{=\left(  -1\right)  ^{\left(  i-1\right)  +\left(  p-1\right)
}=\left(  -1\right)  ^{i+p}\\\text{(since }\left(  i-1\right)  +\left(
p-1\right)  \equiv i+p\operatorname{mod}2\text{)}}}x_{\left(  i,1\right)
}=\left(  -1\right)  ^{i+p}x_{\left(  i,1\right)  }.
\end{align*}
This proves (\ref{lem.Grasp.generic.calc.d.D}). Thus, the proof of Lemma
\ref{lem.Grasp.generic.calc} \textbf{(d)} is complete.

\textbf{(e)} In the following, we are going to regard column vectors of size
$\ell$ as $\ell\times1$-matrices for every $\ell\in\mathbb{N}$. We will also
regard scalars as $1\times1$-matrices.

Let $k\in\left\{  2,3,...,q\right\}  $. Then, $k-1\in\left\{
1,2,...,q-1\right\}  \subseteq\left\{  1,2,...,q\right\}  $. Since
$k-1\in\left\{  1,2,...,q-1\right\}  $, we have $Z\left(  k-1\right)  =\left(
k-1\right)  +1$ (by the definition of $Z$).

By the definition of $\left(  I_{p}\mid C\right)  \left[  p+k-1:p+k\right]  $,
it is clear that the matrix \newline$\left(  I_{p}\mid C\right)  \left[
p+k-1:p+k\right]  $ is the matrix whose column vectors (from left to right)
are $\left(  I_{p}\mid C\right)  _{p+k-1}$ (yes, there is only one column
vector). Thus, $\left(  I_{p}\mid C\right)  \left[  p+k-1:p+k\right]  $ is a
matrix having only one column vector, and this column vector is $\left(
I_{p}\mid C\right)  _{p+k-1}$. Since we are identifying column vectors of size
$\ell$ with $\ell\times1$-matrices for every $\ell\in\mathbb{N}$, we can
rewrite this as follows:%
\[
\left(  I_{p}\mid C\right)  \left[  p+k-1:p+k\right]  =\left(  I_{p}\mid
C\right)  _{p+k-1}.
\]
Hence,%
\begin{align}
&  \operatorname*{rows}\nolimits_{p}\left(  \left(  I_{p}\mid C\right)
\left[  p+k-1:p+k\right]  \right) \nonumber\\
&  =\operatorname*{rows}\nolimits_{p}\left(  \left(  I_{p}\mid C\right)
_{p+k-1}\right)  =\left(  \text{the }p\text{-th entry of }\left(  I_{p}\mid
C\right)  _{p+k-1}\right) \nonumber\\
&  \ \ \ \ \ \ \ \ \ \ \left(
\begin{array}
[c]{c}%
\text{by Proposition \ref{prop.rows.entry}, applied to }\mathbb{F}\left[
x_{\mathbf{P}}\right]  \text{, }p\text{, }p\text{ and }\left(  I_{p}\mid
C\right)  _{p+k-1}\\
\text{instead of }\mathbb{K}\text{, }\ell\text{, }i\text{ and }U
\end{array}
\right) \nonumber\\
&  =\left(  \text{the }\left(  p,k-1\right)  \text{-th entry of the matrix
}C\right) \nonumber\\
&  \ \ \ \ \ \ \ \ \ \ \left(
\begin{array}
[c]{c}%
\text{by Proposition \ref{prop.attach2.entry}, applied to }p\text{, }p\text{,
}q\text{, }I_{p}\text{, }C\text{, }k-1\text{ and }p\\
\text{instead of }u\text{, }v_{1}\text{, }v_{2}\text{, }A_{1}\text{, }%
A_{2}\text{, }\ell\text{ and }i
\end{array}
\right) \nonumber\\
&  =x_{\left(  p,Z\left(  k-1\right)  \right)  }\ \ \ \ \ \ \ \ \ \ \left(
\text{by (\ref{pf.Grasp.generic.calc.entry}), applied to }\left(
p,k-1\right)  \text{ instead of }\left(  i,k\right)  \right) \nonumber\\
&  =x_{\left(  p,k\right)  }\ \ \ \ \ \ \ \ \ \ \left(  \text{since }Z\left(
k-1\right)  =\left(  k-1\right)  +1=k\right)  .
\label{pf.Grasp.generic.calc.e.1}%
\end{align}
Applying (\ref{lem.Grasp.generic.calc.a.N}) to $\left(  r,s\right)  =\left(
p,k\right)  $, we obtain%
\begin{align*}
\mathfrak{N}_{\left(  p,k\right)  }  &  =\det\left(  \operatorname*{rows}%
\nolimits_{p,p+1,...,p}\left(  \left(  I_{p}\mid C\right)  \left[
p+k-1:p+k\right]  \right)  \right) \\
&  =\det\left(  \underbrace{\operatorname*{rows}\nolimits_{p}\left(  \left(
I_{p}\mid C\right)  \left[  p+k-1:p+k\right]  \right)  }%
_{\substack{=x_{\left(  p,k\right)  }\\\text{(by
(\ref{pf.Grasp.generic.calc.e.1}))}}}\right)  \ \ \ \ \ \ \ \ \ \ \left(
\text{since }\left(  p,p+1,...,p\right)  =\left(  p\right)  \right) \\
&  =\det\left(  x_{\left(  p,k\right)  }\right)  =x_{\left(  p,k\right)  }%
\end{align*}
(since the determinant of a $1\times1$-matrix is simply the unique entry of
this matrix). Thus, (\ref{lem.Grasp.generic.calc.e.N}) is proven.

Next, we apply (\ref{lem.Grasp.generic.calc.a.D}) to $\left(  r,s\right)
=\left(  p,k\right)  $. As a result, we obtain
\begin{align*}
\mathfrak{D}_{\left(  p,k\right)  }  &  =\left(  -1\right)  ^{p-1}\det\left(
\operatorname*{rows}\nolimits_{p,p+1,...,p}\left(  \left(  -1\right)
^{p-1}C_{q}\ \mid\ \underbrace{\left(  I_{p}\mid C\right)  \left[
p+s:p+s\right]  }_{\substack{=\left(  \text{empty matrix}\right)  \\\text{(by
Proposition \ref{prop.cols.empty}, applied to}\\\mathbb{F}\left[
x_{\mathbf{P}}\right]  \text{, }p\text{, }p+q\text{, }\left(  I_{p}\mid
C\right)  \text{ and }p+s\\\text{instead of }\mathbb{K}\text{, }u\text{,
}v\text{, }A\text{ and }a\text{)}}}\right)  \right) \\
&  =\left(  -1\right)  ^{p-1}\det\left(  \operatorname*{rows}%
\nolimits_{p,p+1,...,p}\underbrace{\left(  \left(  -1\right)  ^{p-1}%
C_{q}\ \mid\ \left(  \text{empty matrix}\right)  \right)  }_{=\left(
-1\right)  ^{p-1}C_{q}}\right) \\
&  =\left(  -1\right)  ^{p-1}\det\left(  \operatorname*{rows}%
\nolimits_{p,p+1,...,p}\left(  \left(  -1\right)  ^{p-1}C_{q}\right)  \right)
\\
&  =\left(  -1\right)  ^{p-1}\det\left(  \underbrace{\operatorname*{rows}%
\nolimits_{p}\left(  \left(  -1\right)  ^{p-1}C_{q}\right)  }%
_{\substack{=\left(  \text{the }p\text{-th entry of }\left(  -1\right)
^{p-1}C_{q}\right)  \\\text{(by Proposition \ref{prop.rows.entry}, applied to
}\mathbb{F}\left[  x_{\mathbf{P}}\right]  \text{, }p\text{, }\left(
-1\right)  ^{p-1}C_{q}\\\text{and }p\text{ instead of }\mathbb{K}\text{, }%
\ell\text{, }U\text{ and }i\text{)}}}\right) \\
&  \ \ \ \ \ \ \ \ \ \ \left(  \text{since }\left(  p,p+1,...,p\right)
=\left(  p\right)  \right) \\
&  =\left(  -1\right)  ^{p-1}\underbrace{\det\left(  \text{the }p\text{-th
entry of }\left(  -1\right)  ^{p-1}C_{q}\right)  }_{\substack{=\left(
\text{the }p\text{-th entry of }\left(  -1\right)  ^{p-1}C_{q}\right)
\\\text{(since the determinant of a }1\times1\text{-matrix is}\\\text{simply
the unique entry of this matrix)}}}\\
&  =\left(  -1\right)  ^{p-1}\underbrace{\left(  \text{the }p\text{-th entry
of }\left(  -1\right)  ^{p-1}C_{q}\right)  }_{=\left(  -1\right)
^{p-1}\left(  \text{the }p\text{-th entry of }C_{q}\right)  }\\
&  =\underbrace{\left(  -1\right)  ^{p-1}\left(  -1\right)  ^{p-1}%
}_{\substack{=\left(  -1\right)  ^{\left(  p-1\right)  +\left(  p-1\right)
}=1\\\text{(since }\left(  p-1\right)  +\left(  p-1\right)  =2\left(
p-1\right)  \text{ is even)}}}\left(  \text{the }p\text{-th entry of
}\underbrace{C_{q}}_{\substack{=\left(  \text{the }q\text{-th column of
}C\right)  \\\text{(since }1\leq q\leq q\text{)}}}\right) \\
&  =\left(  \text{the }p\text{-th entry of the }q\text{-th column of
}C\right)  =\left(  \text{the }\left(  p,q\right)  \text{-th entry of
}C\right) \\
&  =x_{\left(  p,Z\left(  q\right)  \right)  }\ \ \ \ \ \ \ \ \ \ \left(
\text{by (\ref{pf.Grasp.generic.calc.entry}), applied to }\left(  p,q\right)
\text{ instead of }\left(  i,k\right)  \right) \\
&  =x_{\left(  p,1\right)  }\ \ \ \ \ \ \ \ \ \ \left(  \text{since the
definition of }Z\text{ yields }Z\left(  q\right)  =1\right)  .
\end{align*}
This proves (\ref{lem.Grasp.generic.calc.e.D}). Thus, the proof of Lemma
\ref{lem.Grasp.generic.calc} \textbf{(e)} is complete.

\textbf{(f)} We will prove that every $\left(  i,k\right)  \in\mathbf{P}$
satisfies%
\begin{equation}
\mathfrak{N}_{\left(  i,k\right)  }\neq0\ \ \ \ \ \ \ \ \ \ \text{and}%
\ \ \ \ \ \ \ \ \ \ \mathfrak{D}_{\left(  i,k\right)  }\neq0.
\label{pf.Grasp.generic.calc.f.st}%
\end{equation}

First, let us notice that every $\left(  i,k\right)  \in\mathbf{P}$ satisfies
$k+p-i\in\mathbb{N}$\ \ \ \ \footnote{\textit{Proof.} Let $\left(  i,k\right)
\in\mathbf{P}$. Then, $\left(  i,k\right)  \in\mathbf{P}=\left\{
1,2,...,p\right\}  \times\left\{  1,2,...,q\right\}  $, so that $i\in\left\{
1,2,...,p\right\}  $ and $k\in\left\{  1,2,...,q\right\}  $. Since
$i\in\left\{  1,2,...,p\right\}  $, we have $1\leq i\leq p$. Since
$k\in\left\{  1,2,...,q\right\}  $, we have $1\leq k\leq q$. Now,
$\underbrace{k}_{\geq1}+p-\underbrace{i}_{\leq p}\geq1+p-p=1\geq0$, so that
$k+p-i\in\mathbb{N}$, qed.}.

We are going to prove (\ref{pf.Grasp.generic.calc.f.st}) by strong induction
over $k+p-i$\ \ \ \ \footnote{This kind of induction is legitimate as a proof
because $k+p-i\in\mathbb{N}$ for every $\left(  i,k\right)  \in\mathbf{P}$.}:

\textit{Induction step:}\footnote{A strong induction does not require an
induction base.} Let $N\in\mathbb{N}$. Let us assume that
(\ref{pf.Grasp.generic.calc.f.st}) holds for every $\left(  i,k\right)
\in\mathbf{P}$ satisfying $k+p-i<N$. We now need to prove that
(\ref{pf.Grasp.generic.calc.f.st}) holds for every $\left(  i,k\right)
\in\mathbf{P}$ satisfying $k+p-i=N$.

We have assumed that (\ref{pf.Grasp.generic.calc.f.st}) holds for every
$\left(  i,k\right)  \in\mathbf{P}$ satisfying $k+p-i<N$. In other words,%
\begin{equation}
\text{we have }\mathfrak{N}_{\left(  i,k\right)  }\neq0\text{ and
}\mathfrak{D}_{\left(  i,k\right)  }\neq0\text{ for every }\left(  i,k\right)
\in\mathbf{P}\text{ satisfying }k+p-i<N.
\label{pf.Grasp.generic.calc.f.indass}%
\end{equation}

Now, let $\left(  i,k\right)  \in\mathbf{P}$ be such that $k+p-i=N$. Since%
\[
\left(  i,k\right)  \in\mathbf{P}=\left\{  1,2,...,p\right\}  \times\left\{
1,2,...,q\right\}  ,
\]
we have $i\in\left\{  1,2,...,p\right\}  $ and $k\in\left\{
1,2,...,q\right\}  $.

We distinguish between three cases:

\textit{Case 1:} We have $k=1$.

\textit{Case 2:} We have $i=p$ but not $k=1$.

\textit{Case 3:} We have neither $k=1$ nor $i=p$.

Let us first consider Case 1. In this case, we have $k=1$. Thus,%
\begin{align*}
\mathfrak{N}_{\left(  i,k\right)  }  &  =\mathfrak{N}_{\left(  i,1\right)
}=1\ \ \ \ \ \ \ \ \ \ \left(  \text{by (\ref{lem.Grasp.generic.calc.d.N}%
)}\right) \\
&  \neq0
\end{align*}
and%
\begin{align*}
\mathfrak{D}_{\left(  i,k\right)  }  &  =\mathfrak{D}_{\left(  i,1\right)
}\ \ \ \ \ \ \ \ \ \ \left(  \text{since }k=1\right) \\
&  =\left(  -1\right)  ^{i+p}\underbrace{x_{\left(  i,1\right)  }}_{\neq
0}\ \ \ \ \ \ \ \ \ \ \left(  \text{by (\ref{lem.Grasp.generic.calc.d.D}%
)}\right) \\
&  \neq0.
\end{align*}
In other words, $\mathfrak{N}_{\left(  i,k\right)  }\neq0$ and $\mathfrak{D}%
_{\left(  i,k\right)  }\neq0$ is proven in Case 1.

Let us now consider Case 2. In this case, we have $i=p$ but not $k=1$. Then,
$k\in\left\{  2,3,...,q\right\}  $ (since we have $k\in\left\{
1,2,...,q\right\}  $ but not $k=1$). Since $i=p$, we have%
\begin{align*}
\mathfrak{N}_{\left(  i,k\right)  }  &  =\mathfrak{N}_{\left(  p,k\right)
}=x_{\left(  p,k\right)  }\ \ \ \ \ \ \ \ \ \ \left(  \text{by
(\ref{lem.Grasp.generic.calc.e.N})}\right) \\
&  \neq0
\end{align*}
and%
\begin{align*}
\mathfrak{D}_{\left(  i,k\right)  }  &  =\mathfrak{D}_{\left(  p,k\right)
}\ \ \ \ \ \ \ \ \ \ \left(  \text{since }i=p\right) \\
&  =x_{\left(  p,1\right)  }\ \ \ \ \ \ \ \ \ \ \left(  \text{by
(\ref{lem.Grasp.generic.calc.e.D})}\right) \\
&  \neq0.
\end{align*}
In other words, $\mathfrak{N}_{\left(  i,k\right)  }\neq0$ and $\mathfrak{D}%
_{\left(  i,k\right)  }\neq0$ is proven in Case 2.

Let us finally consider Case 3. In this case, we have neither $k=1$ nor $i=p$.
Then, $i\in\left\{  1,2,...,p-1\right\}  $ (since we have $i\in\left\{
1,2,...,p\right\}  $ but not $i=p$), and $k\in\left\{  2,3,...,q\right\}  $
(since we have $k\in\left\{  1,2,...,q\right\}  $ but not $k=1$). Since
$i\in\left\{  1,2,...,p-1\right\}  $, we have $i+1\in\left\{
2,3,...,p\right\}  \subseteq\left\{  1,2,...,p\right\}  $. Since $k\in\left\{
2,3,...,q\right\}  $, we have $k-1\in\left\{  1,2,...,q-1\right\}
\subseteq\left\{  1,2,...,q\right\}  $. Combining $i+1\in\left\{
1,2,...,p\right\}  $ with $k-1\in\left\{  1,2,...,q\right\}  $, we obtain
$\left(  i+1,k-1\right)  \in\left\{  1,2,...,p\right\}  \times\left\{
1,2,...,q\right\}  =\mathbf{P}$.

Since $\left(  i+1,k-1\right)  \in\mathbf{P}$ and $\left(  k-1\right)
+p-\left(  i+1\right)  =\underbrace{\left(  k+p-i\right)  }_{=N}-2=N-2<N$, it
is clear that we can apply (\ref{pf.Grasp.generic.calc.f.indass}) to $\left(
i+1,k-1\right)  $ instead of $\left(  i,k\right)  $. As a result, we obtain%
\[
\mathfrak{N}_{\left(  i+1,k-1\right)  }\neq0\text{ and }\mathfrak{D}_{\left(
i+1,k-1\right)  }\neq0.
\]

Now, for every subset $\mathbf{S}$ of $\mathbf{P}$, let $\mathbb{F}\left(
x_{\mathbf{S}}\right)  $ denote the field of rational functions over
$\mathbb{F}$ in the indeterminates $x_{\mathbf{p}}$ with $\mathbf{p}$ ranging
over all elements of $\mathbf{S}$ (hence altogether $\left\vert \mathbf{S}%
\right\vert $ indeterminates). Then, $\mathbb{F}\left(  x_{\mathbf{S}}\right)
$ is the fraction field of the polynomial ring $\mathbb{F}\left[
x_{\mathbf{S}}\right]  $.

Notice that every $\mathbf{q}\in\mathbf{P}$ satisfies $\mathbf{q}\notin\left.
\mathbf{q}\Downarrow\right.  $. Applied to $\mathbf{q}=\left(  i,k\right)  $,
this yields $\left(  i,k\right)  \notin\left.  \left(  i,k\right)
\Downarrow\right.  $. Hence, the variable $x_{\left(  i,k\right)  }$ is
foreign to the field $\mathbb{F}\left(  x_{\left(  i,k\right)  \Downarrow
}\right)  $. Thus, $x_{\left(  i,k\right)  }\notin\mathbb{F}\left(  x_{\left(
i,k\right)  \Downarrow}\right)  $.

We have $\mathfrak{N}_{\left(  i+1,k-1\right)  }\in\mathbb{F}\left[
x_{\left(  i,k\right)  \Downarrow}\right]  $ (by
(\ref{lem.Grasp.generic.calc.b.N2})). Hence, $\left(  -1\right)  ^{p-i}%
\cdot\mathfrak{N}_{\left(  i+1,k-1\right)  }\in\left(  -1\right)  ^{p-i}%
\cdot\mathbb{F}\left[  x_{\left(  i,k\right)  \Downarrow}\right]
\subseteq\mathbb{F}\left[  x_{\left(  i,k\right)  \Downarrow}\right]  $.
Moreover, $\left(  -1\right)  ^{p-i}\cdot\mathfrak{N}_{\left(  i+1,k-1\right)
}\neq0$ (since $\left(  -1\right)  ^{p-i}\neq0$ and $\mathfrak{N}_{\left(
i+1,k-1\right)  }\neq0$ and since $\mathbb{F}\left[  x_{\left(  i,k\right)
\Downarrow}\right]  $ is an integral domain).

Assume (for the sake of contradiction) that $\mathfrak{N}_{\left(  i,k\right)
}=0$. Then,%
\[
0=\mathfrak{N}_{\left(  i,k\right)  }\in\left(  -1\right)  ^{p-i}%
\cdot\mathfrak{N}_{\left(  i+1,k-1\right)  }\cdot x_{\left(  i,k\right)
}+\mathbb{F}\left[  x_{\left(  i,k\right)  \Downarrow}\right]
\]
(by (\ref{lem.Grasp.generic.calc.b.N})). In other words, there exists an
$\mathfrak{U}\in\mathbb{F}\left[  x_{\left(  i,k\right)  \Downarrow}\right]  $
such that
\[
0=\left(  -1\right)  ^{p-i}\cdot\mathfrak{N}_{\left(  i+1,k-1\right)  }\cdot
x_{\left(  i,k\right)  }+\mathfrak{U}.
\]
Consider this $\mathfrak{U}$.

We have $0=\left(  -1\right)  ^{p-i}\cdot\mathfrak{N}_{\left(  i+1,k-1\right)
}\cdot x_{\left(  i,k\right)  }+\mathfrak{U}$, hence $\mathfrak{U}=-\left(
-1\right)  ^{p-i}\cdot\mathfrak{N}_{\left(  i+1,k-1\right)  }\cdot x_{\left(
i,k\right)  }$. Dividing this equality by $\left(  -1\right)  ^{p-i}%
\cdot\mathfrak{N}_{\left(  i+1,k-1\right)  }$ (this is clearly allowed in
$\mathbb{F}\left(  x_{\mathbf{P}}\right)  $ since $\left(  -1\right)
^{p-i}\cdot\mathfrak{N}_{\left(  i+1,k-1\right)  }\neq0$), we obtain%
\[
\dfrac{\mathfrak{U}}{\left(  -1\right)  ^{p-i}\cdot\mathfrak{N}_{\left(
i+1,k-1\right)  }}=-x_{\left(  i,k\right)  }.
\]
Hence,%
\[
-x_{\left(  i,k\right)  }=\dfrac{\mathfrak{U}}{\left(  -1\right)  ^{p-i}%
\cdot\mathfrak{N}_{\left(  i+1,k-1\right)  }}\in\mathbb{F}\left(  x_{\left(
i,k\right)  \Downarrow}\right)
\]
(since $\mathfrak{U}$ and $\left(  -1\right)  ^{p-i}\cdot\mathfrak{N}_{\left(
i+1,k-1\right)  }$ both belong to $\mathbb{F}\left[  x_{\left(  i,k\right)
\Downarrow}\right]  $). Consequently, $x_{\left(  i,k\right)  }%
=-\underbrace{\left(  -x_{\left(  i,k\right)  }\right)  }_{\in\mathbb{F}%
\left(  x_{\left(  i,k\right)  \Downarrow}\right)  }\in-\mathbb{F}\left(
x_{\left(  i,k\right)  \Downarrow}\right)  \subseteq\mathbb{F}\left(
x_{\left(  i,k\right)  \Downarrow}\right)  $, contradicting $x_{\left(
i,k\right)  }\notin\mathbb{F}\left(  x_{\left(  i,k\right)  \Downarrow
}\right)  $. This contradiction shows that our assumption (that $\mathfrak{N}%
_{\left(  i,k\right)  }=0$) is wrong. Hence, we don't have $\mathfrak{N}%
_{\left(  i,k\right)  }=0$. Thus, $\mathfrak{N}_{\left(  i,k\right)  }\neq0$.

We have $\mathfrak{D}_{\left(  i+1,k-1\right)  }\in\mathbb{F}\left[
x_{\left(  i,k\right)  \Downarrow}\right]  $ (by
(\ref{lem.Grasp.generic.calc.b.D2})). Hence, $\left(  -1\right)  ^{p-i+1}%
\cdot\mathfrak{D}_{\left(  i+1,k-1\right)  }\in\left(  -1\right)
^{p-i+1}\cdot\mathbb{F}\left[  x_{\left(  i,k\right)  \Downarrow}\right]
\subseteq\mathbb{F}\left[  x_{\left(  i,k\right)  \Downarrow}\right]  $.
Moreover, $\left(  -1\right)  ^{p-i+1}\cdot\mathfrak{D}_{\left(
i+1,k-1\right)  }\neq0$ (since $\left(  -1\right)  ^{p-i+1}\neq0$ and
$\mathfrak{D}_{\left(  i+1,k-1\right)  }\neq0$ and since $\mathbb{F}\left[
x_{\left(  i,k\right)  \Downarrow}\right]  $ is an integral domain).

Assume (for the sake of contradiction) that $\mathfrak{D}_{\left(  i,k\right)
}=0$. Then,%
\[
0=\mathfrak{D}_{\left(  i,k\right)  }\in\left(  -1\right)  ^{p-i+1}%
\cdot\mathfrak{D}_{\left(  i+1,k-1\right)  }\cdot x_{\left(  i,k\right)
}+\mathbb{F}\left[  x_{\left(  i,k\right)  \Downarrow}\right]
\]
(by (\ref{lem.Grasp.generic.calc.b.D})). In other words, there exists an
$\mathfrak{V}\in\mathbb{F}\left[  x_{\left(  i,k\right)  \Downarrow}\right]  $
such that
\[
0=\left(  -1\right)  ^{p-i+1}\cdot\mathfrak{D}_{\left(  i+1,k-1\right)  }\cdot
x_{\left(  i,k\right)  }+\mathfrak{V}.
\]
Consider this $\mathfrak{V}$.

We have $0=\left(  -1\right)  ^{p-i+1}\cdot\mathfrak{D}_{\left(
i+1,k-1\right)  }\cdot x_{\left(  i,k\right)  }+\mathfrak{V}$, hence
$\mathfrak{V}=-\left(  -1\right)  ^{p-i+1}\cdot\mathfrak{D}_{\left(
i+1,k-1\right)  }\cdot x_{\left(  i,k\right)  }$. Dividing this equality by
$\left(  -1\right)  ^{p-i+1}\cdot\mathfrak{D}_{\left(  i+1,k-1\right)  }$
(this is clearly allowed in $\mathbb{F}\left(  x_{\mathbf{P}}\right)  $ since
$\left(  -1\right)  ^{p-i+1}\cdot\mathfrak{D}_{\left(  i+1,k-1\right)  }\neq
0$), we obtain%
\[
\dfrac{\mathfrak{V}}{\left(  -1\right)  ^{p-i+1}\cdot\mathfrak{D}_{\left(
i+1,k-1\right)  }}=-x_{\left(  i,k\right)  }.
\]
Hence,%
\[
-x_{\left(  i,k\right)  }=\dfrac{\mathfrak{V}}{\left(  -1\right)
^{p-i+1}\cdot\mathfrak{D}_{\left(  i+1,k-1\right)  }}\in\mathbb{F}\left(
x_{\left(  i,k\right)  \Downarrow}\right)
\]
(since $\mathfrak{V}$ and $\left(  -1\right)  ^{p-i+1}\cdot\mathfrak{D}%
_{\left(  i+1,k-1\right)  }$ both belong to $\mathbb{F}\left[  x_{\left(
i,k\right)  \Downarrow}\right]  $). Consequently, $x_{\left(  i,k\right)
}=-\underbrace{\left(  -x_{\left(  i,k\right)  }\right)  }_{\in\mathbb{F}%
\left(  x_{\left(  i,k\right)  \Downarrow}\right)  }\in-\mathbb{F}\left(
x_{\left(  i,k\right)  \Downarrow}\right)  \subseteq\mathbb{F}\left(
x_{\left(  i,k\right)  \Downarrow}\right)  $, contradicting $x_{\left(
i,k\right)  }\notin\mathbb{F}\left(  x_{\left(  i,k\right)  \Downarrow
}\right)  $. This contradiction shows that our assumption (that $\mathfrak{D}%
_{\left(  i,k\right)  }=0$) is wrong. Hence, we don't have $\mathfrak{D}%
_{\left(  i,k\right)  }=0$. Thus, $\mathfrak{D}_{\left(  i,k\right)  }\neq0$.

We now have proven that $\mathfrak{N}_{\left(  i,k\right)  }\neq0$ and
$\mathfrak{D}_{\left(  i,k\right)  }\neq0$ hold in Case 3.

Altogether, we have shown that $\mathfrak{N}_{\left(  i,k\right)  }\neq0$ and
$\mathfrak{D}_{\left(  i,k\right)  }\neq0$ hold in each of the three Cases 1,
2 and 3. Since these three Cases cover all possibilities, this yields that
$\mathfrak{N}_{\left(  i,k\right)  }\neq0$ and $\mathfrak{D}_{\left(
i,k\right)  }\neq0$ always hold.

Now forget that we fixed $\left(  i,k\right)  $. We have thus shown that
$\mathfrak{N}_{\left(  i,k\right)  }\neq0$ and $\mathfrak{D}_{\left(
i,k\right)  }\neq0$ hold for every $\left(  i,k\right)  \in\mathbf{P}$
satisfying $k+p-i=N$. In other words, (\ref{pf.Grasp.generic.calc.f.st}) holds
for every $\left(  i,k\right)  \in\mathbf{P}$ satisfying $k+p-i=N$. This
completes the induction step. The induction proof of
(\ref{pf.Grasp.generic.calc.f.st}) is thus complete.

Now, let $\mathbf{p}\in\mathbf{P}$ be arbitrary. Then, we can write
$\mathbf{p}$ in the form $\mathbf{p}=\left(  i,k\right)  $ for some
$i\in\left\{  1,2,...,p\right\}  $ and some $k\in\left\{  1,2,...,q\right\}  $
(since $\mathbf{p}\in\mathbf{P}=\left\{  1,2,...,p\right\}  \times\left\{
1,2,...,q\right\}  $). Consider these $i$ and $k$. Since $\mathbf{p}=\left(
i,k\right)  $, we have $\mathfrak{N}_{\mathbf{p}}=\mathfrak{N}_{\left(
i,k\right)  }\neq0$ (by (\ref{pf.Grasp.generic.calc.f.st})) and $\mathfrak{D}%
_{\mathbf{p}}=\mathfrak{D}_{\left(  i,k\right)  }\neq0$ (by
(\ref{pf.Grasp.generic.calc.f.st})). This proves Lemma
\ref{lem.Grasp.generic.calc} \textbf{(f)}.

\textbf{(g)} Let $\mathbf{p}\in\mathbf{P}$. Then, $\mathbf{p}\in
\mathbf{P}=\left\{  1,2,...,p\right\}  \times\left\{  1,2,...,q\right\}  $.
Hence, we can write $\mathbf{p}$ in the form $\mathbf{p}=\left(  i,k\right)  $
for some $i\in\left\{  1,2,...,p\right\}  $ and $k\in\left\{
1,2,...,q\right\}  $. Consider these $i$ and $k$.

We distinguish between three cases:

\textit{Case 1:} We have $k=1$.

\textit{Case 2:} We have $i=p$ but not $k=1$.

\textit{Case 3:} We have neither $k=1$ nor $i=p$.

Let us consider Case 1 first. In this case, we have $k=1$. Thus,
$\mathbf{p}=\left(  i,\underbrace{k}_{=1}\right)  =\left(  i,1\right)  $, so
that%
\begin{align*}
\mathfrak{N}_{\mathbf{p}}  &  =\mathfrak{N}_{\left(  i,1\right)
}=1\ \ \ \ \ \ \ \ \ \ \left(  \text{by (\ref{lem.Grasp.generic.calc.d.N}%
)}\right) \\
&  =0x_{\mathbf{p}}+1\ \ \ \ \ \ \ \ \ \ \left(  \text{since }%
\underbrace{0x_{\mathbf{p}}}_{=0}+1=1\right)
\end{align*}
and%
\begin{align*}
\mathfrak{D}_{\mathbf{p}}  &  =\mathfrak{D}_{\left(  i,1\right)
}\ \ \ \ \ \ \ \ \ \ \left(  \text{since }\mathbf{p}=\left(  i,1\right)
\right) \\
&  =\left(  -1\right)  ^{i+p}x_{\left(  i,1\right)  }%
\ \ \ \ \ \ \ \ \ \ \left(  \text{by (\ref{lem.Grasp.generic.calc.d.D}%
)}\right) \\
&  =\left(  -1\right)  ^{i+p}x_{\mathbf{p}}\ \ \ \ \ \ \ \ \ \ \left(
\text{since }\left(  i,1\right)  =\mathbf{p}\right) \\
&  =\left(  -1\right)  ^{i+p}x_{\mathbf{p}}+0.
\end{align*}
Also, it is clear that $0$, $1$, $\left(  -1\right)  ^{i+p}$ and $0$ are
elements of $\mathbb{F}\left[  x_{\mathbf{p}\Downarrow}\right]  $ and satisfy
$0\cdot0-1\cdot\left(  -1\right)  ^{i+p}=-\left(  -1\right)  ^{i+p}\neq0$.
Hence, there exist elements $\alpha_{\mathbf{p}}$, $\beta_{\mathbf{p}}$,
$\gamma_{\mathbf{p}}$ and $\delta_{\mathbf{p}}$ of $\mathbb{F}\left[
x_{\mathbf{p}\Downarrow}\right]  $ satisfying $\alpha_{\mathbf{p}}%
\delta_{\mathbf{p}}-\beta_{\mathbf{p}}\gamma_{\mathbf{p}}\neq0$,
$\mathfrak{N}_{\mathbf{p}}=\alpha_{\mathbf{p}}x_{\mathbf{p}}+\beta
_{\mathbf{p}}$ and $\mathfrak{D}_{\mathbf{p}}=\gamma_{\mathbf{p}}%
x_{\mathbf{p}}+\delta_{\mathbf{p}}$ (namely, $0$, $1$, $\left(  -1\right)
^{i+p}$ and $0$). In other words, Lemma \ref{lem.Grasp.generic.calc}
\textbf{(g)} is proven in Case 1.

Let us now consider Case 2. In this case, we have $i=p$ but not $k=1$. Then,
$k\in\left\{  2,3,...,q\right\}  $ (since we have $k\in\left\{
1,2,...,q\right\}  $ but not $k=1$). Now, $\mathbf{p}=\left(  \underbrace{i}%
_{=p},k\right)  =\left(  p,k\right)  $. It is easily seen that $\left(
p,1\right)  \in\left.  \mathbf{p}\Downarrow\right.  $%
\ \ \ \ \footnote{\textit{Proof.} Clearly, $\left(  p,1\right)  \in\left\{
1,2,...,p\right\}  \times\left\{  1,2,...,q\right\}  =\mathbf{P}$ and $\left(
p,k\right)  =\mathbf{p}\in\mathbf{P}$. Moreover, $p\geq p$ and $1\leq k$.
Hence, $\left(  p,1\right)  \trianglelefteq\left(  p,k\right)  $ (by
(\ref{lem.Grasp.generic.calc.ass}), applied to $\left(  p,1\right)  $ and
$\left(  p,k\right)  $ instead of $\left(  i,k\right)  $ and $\left(
i^{\prime},k^{\prime}\right)  $). We also have $k\neq1$ (since we don't have
$k=1$), so that $1\neq k$ and thus $\left(  p,1\right)  \neq\left(
p,k\right)  $. Combining this with $\left(  p,1\right)  \trianglelefteq\left(
p,k\right)  $, we obtain $\left(  p,1\right)  \vartriangleleft\left(
p,k\right)  $. Since $\left(  p,k\right)  =\mathbf{p}$, this rewrites as
$\left(  p,1\right)  \vartriangleleft\mathbf{p}$. Hence, $\left(  p,1\right)
\in\left\{  \mathbf{v}\in\mathbf{P}\mid\mathbf{v}\vartriangleleft
\mathbf{p}\right\}  =\left.  \mathbf{p}\Downarrow\right.  $ (because the
definition of $\mathbf{p}\Downarrow$ yields $\left.  \mathbf{p}\Downarrow
\right.  =\left\{  \mathbf{v}\in\mathbf{P}\ \mid\ \mathbf{v}\vartriangleleft
\mathbf{p}\right\}  $), qed.}. Hence, $x_{\left(  p,1\right)  }\in
\mathbb{F}\left[  x_{\mathbf{p}\Downarrow}\right]  $. Now,%
\begin{align*}
\mathfrak{N}_{\mathbf{p}}  &  =\mathfrak{N}_{\left(  p,k\right)
}\ \ \ \ \ \ \ \ \ \ \left(  \text{since }\mathbf{p}=\left(  p,k\right)
\right) \\
&  =x_{\left(  p,k\right)  }\ \ \ \ \ \ \ \ \ \ \left(  \text{by
(\ref{lem.Grasp.generic.calc.e.N})}\right) \\
&  =x_{\mathbf{p}}\ \ \ \ \ \ \ \ \ \ \left(  \text{since }\left(  p,k\right)
=\mathbf{p}\right) \\
&  =1x_{\mathbf{p}}=1x_{\mathbf{p}}+0
\end{align*}
and%
\begin{align*}
\mathfrak{D}_{\mathbf{p}}  &  =\mathfrak{D}_{\left(  p,k\right)
}\ \ \ \ \ \ \ \ \ \ \left(  \text{since }\mathbf{p}=\left(  p,k\right)
\right) \\
&  =x_{\left(  p,1\right)  }\ \ \ \ \ \ \ \ \ \ \left(  \text{by
(\ref{lem.Grasp.generic.calc.e.D})}\right) \\
&  =0x_{\mathbf{p}}+x_{\left(  p,1\right)  }\ \ \ \ \ \ \ \ \ \ \left(
\text{since }\underbrace{0x_{\mathbf{p}}}_{=0}+x_{\left(  p,1\right)
}=x_{\left(  p,1\right)  }\right)  .
\end{align*}
Now, $1$, $0$, $0$ and $x_{\left(  p,1\right)  }$ are elements of
$\mathbb{F}\left[  x_{\mathbf{p}\Downarrow}\right]  $ (indeed, we have proven
above that $x_{\left(  p,1\right)  }\in\mathbb{F}\left[  x_{\mathbf{p}%
\Downarrow}\right]  $) and satisfy $1\cdot x_{\left(  p,1\right)  }%
-0\cdot0=x_{\left(  p,1\right)  }\neq0$. Hence, there exist elements
$\alpha_{\mathbf{p}}$, $\beta_{\mathbf{p}}$, $\gamma_{\mathbf{p}}$ and
$\delta_{\mathbf{p}}$ of $\mathbb{F}\left[  x_{\mathbf{p}\Downarrow}\right]  $
satisfying $\alpha_{\mathbf{p}}\delta_{\mathbf{p}}-\beta_{\mathbf{p}}%
\gamma_{\mathbf{p}}\neq0$, $\mathfrak{N}_{\mathbf{p}}=\alpha_{\mathbf{p}%
}x_{\mathbf{p}}+\beta_{\mathbf{p}}$ and $\mathfrak{D}_{\mathbf{p}}%
=\gamma_{\mathbf{p}}x_{\mathbf{p}}+\delta_{\mathbf{p}}$ (namely, $1$, $0$, $0$
and $x_{\left(  p,1\right)  }$). In other words, Lemma
\ref{lem.Grasp.generic.calc} \textbf{(g)} is proven in Case 2.

Let us now consider Case 3. In this case, we have neither $k=1$ nor $i=p$.
Then, $i\in\left\{  1,2,...,p-1\right\}  $ (since we have $i\in\left\{
1,2,...,p\right\}  $ but not $i=p$), and $k\in\left\{  2,3,...,q\right\}  $
(since we have $k\in\left\{  1,2,...,q\right\}  $ but not $k=1$). Since
$i\in\left\{  1,2,...,p-1\right\}  $, we have $i+1\in\left\{
2,3,...,p\right\}  \subseteq\left\{  1,2,...,p\right\}  $. Since $k\in\left\{
2,3,...,q\right\}  $, we have $k-1\in\left\{  1,2,...,q-1\right\}
\subseteq\left\{  1,2,...,q\right\}  $.

Combining $i+1\in\left\{  1,2,...,p\right\}  $ with $k\in\left\{
1,2,...,q\right\}  $, we obtain $\left(  i+1,k\right)  \in\left\{
1,2,...,p\right\}  \times\left\{  1,2,...,q\right\}  =\mathbf{P}$. Hence, we
can apply (\ref{lem.Grasp.generic.calc.f.N}) to $\left(  i+1,k\right)  $
instead of $\mathbf{p}$. As a result, we obtain $\mathfrak{N}_{\left(
i+1,k\right)  }\neq0$.

Combining $i\in\left\{  1,2,...,p\right\}  $ with $k-1\in\left\{
1,2,...,q\right\}  $, we obtain $\left(  i,k-1\right)  \in\left\{
1,2,...,p\right\}  \times\left\{  1,2,...,q\right\}  =\mathbf{P}$. Hence, we
can apply (\ref{lem.Grasp.generic.calc.f.D}) to $\left(  i,k-1\right)  $
instead of $\mathbf{p}$. As a result, we obtain $\mathfrak{D}_{\left(
i,k-1\right)  }\neq0$.

Since $\mathfrak{D}_{\left(  i,k-1\right)  }\neq0$ and $\mathfrak{N}_{\left(
i+1,k\right)  }\neq0$, we have $\underbrace{\mathfrak{D}_{\left(
i,k-1\right)  }}_{\neq0}\underbrace{\mathfrak{N}_{\left(  i+1,k\right)  }%
}_{\neq0}\neq0$ (since $\mathbb{F}\left[  x_{\mathbf{P}}\right]  $ is an
integral domain).

We have%
\[
\mathfrak{N}_{\left(  i,k\right)  }\in\left(  -1\right)  ^{p-i}\cdot
\mathfrak{N}_{\left(  i+1,k-1\right)  }\cdot x_{\left(  i,k\right)
}+\mathbb{F}\left[  x_{\left(  i,k\right)  \Downarrow}\right]
\]
(by (\ref{lem.Grasp.generic.calc.b.N})). In other words, there exists an
$\mathfrak{U}\in\mathbb{F}\left[  x_{\left(  i,k\right)  \Downarrow}\right]  $
such that
\begin{equation}
\mathfrak{N}_{\left(  i,k\right)  }=\left(  -1\right)  ^{p-i}\cdot
\mathfrak{N}_{\left(  i+1,k-1\right)  }\cdot x_{\left(  i,k\right)
}+\mathfrak{U}. \label{pf.Grasp.generic.calc.g.N1}%
\end{equation}
Consider this $\mathfrak{U}$. Since $\left(  i,k\right)  =\mathbf{p}$, the
equality (\ref{pf.Grasp.generic.calc.g.N1}) becomes%
\begin{equation}
\mathfrak{N}_{\mathbf{p}}=\left(  -1\right)  ^{p-i}\cdot\mathfrak{N}_{\left(
i+1,k-1\right)  }\cdot x_{\mathbf{p}}+\mathfrak{U}.
\label{pf.Grasp.generic.calc.g.N2}%
\end{equation}

Furthermore,%
\[
\mathfrak{D}_{\left(  i,k\right)  }\in\left(  -1\right)  ^{p-i+1}%
\cdot\mathfrak{D}_{\left(  i+1,k-1\right)  }\cdot x_{\left(  i,k\right)
}+\mathbb{F}\left[  x_{\left(  i,k\right)  \Downarrow}\right]
\]
(by (\ref{lem.Grasp.generic.calc.b.D})). In other words, there exists an
$\mathfrak{V}\in\mathbb{F}\left[  x_{\left(  i,k\right)  \Downarrow}\right]  $
such that
\begin{equation}
\mathfrak{D}_{\left(  i,k\right)  }=\left(  -1\right)  ^{p-i+1}\cdot
\mathfrak{D}_{\left(  i+1,k-1\right)  }\cdot x_{\left(  i,k\right)
}+\mathfrak{V}. \label{pf.Grasp.generic.calc.g.D1}%
\end{equation}
Consider this $\mathfrak{V}$. Since $\left(  i,k\right)  =\mathbf{p}$, the
equality (\ref{pf.Grasp.generic.calc.g.D1}) becomes%
\begin{equation}
\mathfrak{D}_{\mathbf{p}}=\left(  -1\right)  ^{p-i+1}\cdot\mathfrak{D}%
_{\left(  i+1,k-1\right)  }\cdot x_{\mathbf{p}}+\mathfrak{V}.
\label{pf.Grasp.generic.calc.g.D2}%
\end{equation}

From (\ref{lem.Grasp.generic.calc.b.N2}), we have $\mathfrak{N}_{\left(
i+1,k-1\right)  }\in\mathbb{F}\left[  x_{\left(  i,k\right)  \Downarrow
}\right]  $, so that $\left(  -1\right)  ^{p-i}\cdot\mathfrak{N}_{\left(
i+1,k-1\right)  }\in\left(  -1\right)  ^{p-i}\cdot\mathbb{F}\left[  x_{\left(
i,k\right)  \Downarrow}\right]  \subseteq\mathbb{F}\left[  x_{\left(
i,k\right)  \Downarrow}\right]  $.

From (\ref{lem.Grasp.generic.calc.b.D2}), we have $\mathfrak{D}_{\left(
i+1,k-1\right)  }\in\mathbb{F}\left[  x_{\left(  i,k\right)  \Downarrow
}\right]  $, so that $\left(  -1\right)  ^{p-i+1}\cdot\mathfrak{D}_{\left(
i+1,k-1\right)  }\in\left(  -1\right)  ^{p-i+1}\cdot\mathbb{F}\left[
x_{\left(  i,k\right)  \Downarrow}\right]  \subseteq\mathbb{F}\left[
x_{\left(  i,k\right)  \Downarrow}\right]  $.

We now have%
\begin{align}
&  \left(  \left(  -1\right)  ^{p-i}\cdot\mathfrak{N}_{\left(  i+1,k-1\right)
}\right)  \cdot\underbrace{\mathfrak{V}}_{\substack{=\mathfrak{D}_{\left(
i,k\right)  }-\left(  -1\right)  ^{p-i+1}\cdot\mathfrak{D}_{\left(
i+1,k-1\right)  }\cdot x_{\left(  i,k\right)  }\\\text{(as a consequence of
(\ref{pf.Grasp.generic.calc.g.D1}))}}}\nonumber\\
&  \ \ \ \ \ \ \ \ \ \ -\underbrace{\mathfrak{U}}_{\substack{=\mathfrak{N}%
_{\left(  i,k\right)  }-\left(  -1\right)  ^{p-i}\cdot\mathfrak{N}_{\left(
i+1,k-1\right)  }\cdot x_{\left(  i,k\right)  }\\\text{(as a consequence of
(\ref{pf.Grasp.generic.calc.g.N1}))}}}\cdot\left(  \left(  -1\right)
^{p-i+1}\cdot\mathfrak{D}_{\left(  i+1,k-1\right)  }\right) \nonumber\\
&  =\left(  \left(  -1\right)  ^{p-i}\cdot\mathfrak{N}_{\left(
i+1,k-1\right)  }\right)  \cdot\left(  \mathfrak{D}_{\left(  i,k\right)
}-\left(  -1\right)  ^{p-i+1}\cdot\mathfrak{D}_{\left(  i+1,k-1\right)  }\cdot
x_{\left(  i,k\right)  }\right) \nonumber\\
&  \ \ \ \ \ \ \ \ \ \ -\left(  \mathfrak{N}_{\left(  i,k\right)  }-\left(
-1\right)  ^{p-i}\cdot\mathfrak{N}_{\left(  i+1,k-1\right)  }\cdot x_{\left(
i,k\right)  }\right)  \cdot\left(  \left(  -1\right)  ^{p-i+1}\cdot
\mathfrak{D}_{\left(  i+1,k-1\right)  }\right) \nonumber\\
&  =\left(  \left(  -1\right)  ^{p-i}\cdot\mathfrak{N}_{\left(
i+1,k-1\right)  }\cdot\mathfrak{D}_{\left(  i,k\right)  }-\left(  -1\right)
^{p-i}\cdot\mathfrak{N}_{\left(  i+1,k-1\right)  }\cdot\left(  -1\right)
^{p-i+1}\cdot\mathfrak{D}_{\left(  i+1,k-1\right)  }\cdot x_{\left(
i,k\right)  }\right) \nonumber\\
&  \ \ \ \ \ \ \ \ \ \ -\left(  \mathfrak{N}_{\left(  i,k\right)  }%
\cdot\left(  -1\right)  ^{p-i+1}\cdot\mathfrak{D}_{\left(  i+1,k-1\right)
}-\left(  -1\right)  ^{p-i}\cdot\mathfrak{N}_{\left(  i+1,k-1\right)  }\cdot
x_{\left(  i,k\right)  }\cdot\left(  -1\right)  ^{p-i+1}\cdot\mathfrak{D}%
_{\left(  i+1,k-1\right)  }\right) \nonumber\\
&  =\left(  \left(  -1\right)  ^{p-i}\cdot\mathfrak{N}_{\left(
i+1,k-1\right)  }\cdot\mathfrak{D}_{\left(  i,k\right)  }-\left(  -1\right)
^{p-i}\left(  -1\right)  ^{p-i+1}\cdot\mathfrak{N}_{\left(  i+1,k-1\right)
}\cdot\mathfrak{D}_{\left(  i+1,k-1\right)  }\cdot x_{\left(  i,k\right)
}\right) \nonumber\\
&  \ \ \ \ \ \ \ \ \ \ -\left(  \left(  -1\right)  ^{p-i+1}\cdot
\mathfrak{D}_{\left(  i+1,k-1\right)  }\cdot\mathfrak{N}_{\left(  i,k\right)
}-\left(  -1\right)  ^{p-i}\left(  -1\right)  ^{p-i+1}\cdot\mathfrak{N}%
_{\left(  i+1,k-1\right)  }\cdot\mathfrak{D}_{\left(  i+1,k-1\right)  }\cdot
x_{\left(  i,k\right)  }\right) \nonumber\\
&  =\left(  -1\right)  ^{p-i}\cdot\mathfrak{N}_{\left(  i+1,k-1\right)  }%
\cdot\mathfrak{D}_{\left(  i,k\right)  }-\underbrace{\left(  -1\right)
^{p-i+1}}_{=\left(  -1\right)  ^{p-i}\left(  -1\right)  ^{1}=\left(
-1\right)  ^{p-i}\left(  -1\right)  =-\left(  -1\right)  ^{p-i}}%
\cdot\mathfrak{D}_{\left(  i+1,k-1\right)  }\cdot\mathfrak{N}_{\left(
i,k\right)  }\nonumber\\
&  =\left(  -1\right)  ^{p-i}\cdot\mathfrak{N}_{\left(  i+1,k-1\right)  }%
\cdot\mathfrak{D}_{\left(  i,k\right)  }-\left(  -\left(  -1\right)
^{p-i}\right)  \cdot\mathfrak{D}_{\left(  i+1,k-1\right)  }\cdot
\mathfrak{N}_{\left(  i,k\right)  }\nonumber\\
&  =\left(  -1\right)  ^{p-i}\cdot\mathfrak{N}_{\left(  i+1,k-1\right)  }%
\cdot\mathfrak{D}_{\left(  i,k\right)  }+\left(  -1\right)  ^{p-i}%
\cdot\mathfrak{D}_{\left(  i+1,k-1\right)  }\cdot\mathfrak{N}_{\left(
i,k\right)  }\nonumber\\
&  =\left(  -1\right)  ^{p-i}\cdot\underbrace{\left(  \mathfrak{N}_{\left(
i+1,k-1\right)  }\cdot\mathfrak{D}_{\left(  i,k\right)  }+\mathfrak{D}%
_{\left(  i+1,k-1\right)  }\cdot\mathfrak{N}_{\left(  i,k\right)  }\right)
}_{\substack{=\mathfrak{D}_{\left(  i,k\right)  }\mathfrak{N}_{\left(
i+1,k-1\right)  }+\mathfrak{N}_{\left(  i,k\right)  }\mathfrak{D}_{\left(
i+1,k-1\right)  }\\=\mathfrak{N}_{\left(  i,k\right)  }\mathfrak{D}_{\left(
i+1,k-1\right)  }+\mathfrak{D}_{\left(  i,k\right)  }\mathfrak{N}_{\left(
i+1,k-1\right)  }\\=\mathfrak{D}_{\left(  i,k-1\right)  }\mathfrak{N}_{\left(
i+1,k\right)  }\\\text{(by Lemma \ref{lem.Grasp.generic.calc} \textbf{(c)})}%
}}\nonumber\\
&  =\left(  -1\right)  ^{p-i}\cdot\underbrace{\mathfrak{D}_{\left(
i,k-1\right)  }\mathfrak{N}_{\left(  i+1,k\right)  }}_{\neq0}\neq0.
\label{pf.Grasp.generic.calc.g.nonzero}%
\end{align}
Thus we know that the four elements $\left(  -1\right)  ^{p-i}\cdot
\mathfrak{N}_{\left(  i+1,k-1\right)  }$, $\mathfrak{U}$, $\left(  -1\right)
^{p-i+1}\cdot\mathfrak{D}_{\left(  i+1,k-1\right)  }$ and $\mathfrak{V}$ are
elements of $\mathbb{F}\left[  x_{\left(  i,k\right)  \Downarrow}\right]  $
and satisfy (\ref{pf.Grasp.generic.calc.g.N2}),
(\ref{pf.Grasp.generic.calc.g.D2}) and (\ref{pf.Grasp.generic.calc.g.nonzero}%
). Since $\left(  i,k\right)  =\mathbf{p}$, this rewrites as follows: The four
elements $\left(  -1\right)  ^{p-i}\cdot\mathfrak{N}_{\left(  i+1,k-1\right)
}$, $\mathfrak{U}$, $\left(  -1\right)  ^{p-i+1}\cdot\mathfrak{D}_{\left(
i+1,k-1\right)  }$ and $\mathfrak{V}$ are elements of $\mathbb{F}\left[
x_{\mathbf{p}\Downarrow}\right]  $ and satisfy
(\ref{pf.Grasp.generic.calc.g.N2}), (\ref{pf.Grasp.generic.calc.g.D2}) and
(\ref{pf.Grasp.generic.calc.g.nonzero}). Hence, there exist elements
$\alpha_{\mathbf{p}}$, $\beta_{\mathbf{p}}$, $\gamma_{\mathbf{p}}$ and
$\delta_{\mathbf{p}}$ of $\mathbb{F}\left[  x_{\mathbf{p}\Downarrow}\right]  $
satisfying $\alpha_{\mathbf{p}}\delta_{\mathbf{p}}-\beta_{\mathbf{p}}%
\gamma_{\mathbf{p}}\neq0$, $\mathfrak{N}_{\mathbf{p}}=\alpha_{\mathbf{p}%
}x_{\mathbf{p}}+\beta_{\mathbf{p}}$ and $\mathfrak{D}_{\mathbf{p}}%
=\gamma_{\mathbf{p}}x_{\mathbf{p}}+\delta_{\mathbf{p}}$ (namely, $\left(
-1\right)  ^{p-i}\cdot\mathfrak{N}_{\left(  i+1,k-1\right)  }$, $\mathfrak{U}%
$, $\left(  -1\right)  ^{p-i+1}\cdot\mathfrak{D}_{\left(  i+1,k-1\right)  }$
and $\mathfrak{V}$). In other words, Lemma \ref{lem.Grasp.generic.calc}
\textbf{(g)} is proven in Case 3.

We have thus proven Lemma \ref{lem.Grasp.generic.calc} \textbf{(g)} in each of
the three Cases 1, 2 and 3. Since these three Cases cover all possibilities,
this yields that Lemma \ref{lem.Grasp.generic.calc} \textbf{(g)} always holds.
The proof of Lemma \ref{lem.Grasp.generic.calc} \textbf{(g)} is thus complete.
\end{proof}

The next lemma is going to allow us to apply the above results to the proof of
Proposition \ref{prop.Grasp.generic}.

\begin{lemma}
\label{lem.Grasp.generic.independency} Let $\mathbb{F}$ be a field.

Let $p$ and $q$ be two positive integers. Let $\mathbf{P}$ be a totally
ordered set such that%
\[
\mathbf{P}=\left\{  1,2,...,p\right\}  \times\left\{  1,2,...,q\right\}
\text{ as sets.}%
\]
Let $\vartriangleleft$ denote the smaller relation of $\mathbf{P}$, and let
$\trianglelefteq$ denote the smaller-or-equal relation of $\mathbf{P}$. Assume
that%
\[
\left(  i,k\right)  \trianglelefteq\left(  i^{\prime},k^{\prime}\right)
\text{ for all }\left(  i,k\right)  \in\mathbf{P}\text{ and }\left(
i^{\prime},k^{\prime}\right)  \in\mathbf{P}\text{ satisfying }\left(  i\geq
i^{\prime}\text{ and }k\leq k^{\prime}\right)  .
\]

Let $Z:\left\{  1,2,...,q\right\}  \rightarrow\left\{  1,2,...,q\right\}  $
denote the map which sends every $k\in\left\{  1,2,...,q-1\right\}  $ to $k+1$
and sends $q$ to $1$. Thus, $Z$ is a permutation in the symmetric group
$S_{q}$, and can be written in cycle notation as $\left(  1,2,...,q\right)  $.

Denote by $\mathbb{F}\left[  x_{\mathbf{P}}\right]  $ the polynomial ring over
$\mathbb{F}$ in the indeterminates $x_{\mathbf{p}}$ with $\mathbf{p}$ ranging
over all elements of $\mathbf{P}$. Let us regard $\mathbb{F}\left[
x_{\mathbf{P}}\right]  $ as a subring of $\mathbb{F}\left(  x_{\mathbf{P}%
}\right)  $. (Recall that the field $\mathbb{F}\left(  x_{\mathbf{P}}\right)
$ is defined as in Definition \ref{def.algebraic.triangularity} \textbf{(a)}.)
It is clear that $\mathbb{F}\left(  x_{\mathbf{P}}\right)  $ is then the
quotient field of $\mathbb{F}\left[  x_{\mathbf{P}}\right]  $.

Define a family $\left(  y_{\mathbf{p}}\right)  _{\mathbf{p}\in\mathbf{P}}%
\in\left(  \mathbb{F}\left[  x_{\mathbf{P}}\right]  \right)  ^{\mathbf{P}}$ of
elements of $\mathbb{F}\left[  x_{\mathbf{P}}\right]  $ by setting
\[
y_{\left(  i,k\right)  }=x_{\left(  i,Z\left(  k\right)  \right)
}\ \ \ \ \ \ \ \ \ \ \text{for all }\left(  i,k\right)  \in\mathbf{P}.
\]

Define a matrix $C\in\left(  \mathbb{F}\left[  x_{\mathbf{P}}\right]  \right)
^{p\times q}$ by
\[
C=\left(  y_{\left(  i,k\right)  }\right)  _{1\leq i\leq p,\ 1\leq k\leq q}.
\]

Let $\mathbb{L}$ be the field $\mathbb{F}\left(  x_{\mathbf{P}}\right)  $.
Notice that $C\in\left(  \mathbb{F}\left[  x_{\mathbf{P}}\right]  \right)
^{p\times q}\subseteq\mathbb{L}^{p\times q}$.

For every $\left(  i,k\right)  \in\mathbf{P}$, define two elements
$\mathfrak{N}_{\left(  i,k\right)  }\in\mathbb{F}\left[  x_{\mathbf{P}%
}\right]  $ and $\mathfrak{D}_{\left(  i,k\right)  }\in\mathbb{F}\left[
x_{\mathbf{P}}\right]  $ as in Lemma \ref{lem.Grasp.generic.calc}.

\textbf{(a)} For every $\mathbf{p}\in\mathbf{P}$, the rational function
$\left(  \operatorname*{Grasp}\nolimits_{0}\left(  I_{p}\mid C\right)
\right)  \left(  \mathbf{p}\right)  \in\mathbb{F}\left(  x_{\mathbf{P}%
}\right)  =\mathbb{L}$ is well-defined and satisfies%
\begin{equation}
\left(  \operatorname*{Grasp}\nolimits_{0}\left(  I_{p}\mid C\right)  \right)
\left(  \mathbf{p}\right)  =\dfrac{\mathfrak{N}_{\mathbf{p}}}{\mathfrak{D}%
_{\mathbf{p}}}. \label{lem.Grasp.generic.independency.a}%
\end{equation}

\textbf{(b)} Define a family $\left(  Q_{\mathbf{p}}\right)  _{\mathbf{p}%
\in\mathbf{P}}\in\mathbb{L}^{\mathbf{P}}$ of elements of the field
$\mathbb{L}$ by setting
\[
Q_{\mathbf{p}}=\left(  \operatorname*{Grasp}\nolimits_{0}\left(  I_{p}\mid
C\right)  \right)  \left(  \mathbf{p}\right)  \ \ \ \ \ \ \ \ \ \ \text{for
all }\mathbf{p}\in\mathbf{P}.
\]
Then, the family $\left(  Q_{\mathbf{p}}\right)  _{\mathbf{p}\in\mathbf{P}}$
is $\mathbf{P}$-triangular.

\textbf{(c)} Let $\mathbb{M}$ be a field extension of $\mathbb{F}$. Let
$\left(  d_{\mathbf{r}}\right)  _{\mathbf{r}\in\mathbf{P}}\in\mathbb{M}%
^{\mathbf{P}}$ be any family. Let $D$ denote the matrix $\left(  d_{\left(
i,Z\left(  k\right)  \right)  }\right)  _{1\leq i\leq p,\ 1\leq k\leq q}%
\in\mathbb{M}^{p\times q}$. Then, for every $\left(  i,k\right)  \in
\mathbf{P}$, we have%
\begin{equation}
\mathfrak{N}_{\left(  i,k\right)  }\left(  \left(  d_{\mathbf{r}}\right)
_{\mathbf{r}\in\mathbf{P}}\right)  =\det\left(  \left(  I_{p}\mid D\right)
\left[  1:i\mid i+k-1:p+k\right]  \right)
\label{lem.Grasp.generic.independency.cN}%
\end{equation}
and%
\begin{equation}
\mathfrak{D}_{\left(  i,k\right)  }\left(  \left(  d_{\mathbf{r}}\right)
_{\mathbf{r}\in\mathbf{P}}\right)  =\det\left(  \left(  I_{p}\mid D\right)
\left[  0:i\mid i+k:p+k\right]  \right)  .
\label{lem.Grasp.generic.independency.cD}%
\end{equation}

\textbf{(d)} Consider the family $\left(  Q_{\mathbf{p}}\right)
_{\mathbf{p}\in\mathbf{P}}$ defined in part \textbf{(b)}. Let $\mathbb{M}$,
$\left(  d_{\mathbf{r}}\right)  _{\mathbf{r}\in\mathbf{P}}$ and $D$ be as in
Lemma \ref{lem.Grasp.generic.independency} \textbf{(c)}. Let $\mathbf{p}%
\in\mathbf{P}$ be such that $\mathfrak{D}_{\mathbf{p}}\left(  \left(
d_{\mathbf{r}}\right)  _{\mathbf{r}\in\mathbf{P}}\right)  \neq0$. Then, the
values $Q_{\mathbf{p}}\left(  \left(  d_{\mathbf{r}}\right)  _{\mathbf{r}%
\in\mathbf{P}}\right)  $ and $\left(  \operatorname*{Grasp}\nolimits_{0}%
\left(  I_{p}\mid D\right)  \right)  \left(  \mathbf{p}\right)  $ are
well-defined and satisfy%
\begin{equation}
Q_{\mathbf{p}}\left(  \left(  d_{\mathbf{r}}\right)  _{\mathbf{r}\in
\mathbf{P}}\right)  =\left(  \operatorname*{Grasp}\nolimits_{0}\left(
I_{p}\mid D\right)  \right)  \left(  \mathbf{p}\right)  .
\label{lem.Grasp.generic.independency.d}%
\end{equation}

\end{lemma}

Before we come to the proof of this lemma, let us see how this lemma yields
Proposition \ref{prop.Grasp.generic}:

\begin{proof}
[Proof of Proposition \ref{prop.Grasp.generic} using Lemma
\ref{lem.Grasp.generic.independency} (sketched).]Let $\mathbb{F}$ be the prime
field of $\mathbb{K}$. (Here, the \textit{prime field} of a field $\mathbb{K}$
denotes the smallest field contained in $\mathbb{K}$. This is either
$\mathbb{Q}$ or a finite field.) In the following, the words ``algebraically
independent'' will always mean ``algebraically independent over $\mathbb{F}$''
(rather than over $\mathbb{K}$ or over $\mathbb{Z}$). Clearly, $\mathbb{K}$ is
a field extension of $\mathbb{F}$.

Now, let $\widetilde{\mathbf{P}}$ be the poset $\left\{  1,2,...,p\right\}
\times\left\{  1,2,...,q\right\}  $ with the order defined by%
\[
\left(  \left(  i,k\right)  \text{ is smaller or equal to }\left(  i^{\prime
},k^{\prime}\right)  \text{ if and only if }\left(  i\geq i^{\prime}\text{ and
}k\leq k^{\prime}\right)  \right)  .
\]
\footnote{This poset differs from $\operatorname*{Rect}\left(  p,q\right)  $
(in general), although it has the same ground set.} Fix any total order
$\mathbf{P}$ on the set $\left\{  1,2,...,p\right\}  \times\left\{
1,2,...,q\right\}  $ which is compatible with the order of
$\widetilde{\mathbf{P}}$. (Such a total order clearly exists, because for any
finite poset $\mathbf{Q}$ there exists a total order on the set $\mathbf{Q}$
compatible with the order of $\mathbf{Q}$.) Thus, $\mathbf{P}$ is a totally
ordered set such that
\[
\mathbf{P}=\left\{  1,2,...,p\right\}  \times\left\{  1,2,...,q\right\}
\text{ as sets.}%
\]
We denote the smaller relation of this totally ordered set $\mathbf{P}$ by
$\vartriangleleft$, and we denote the smaller-or-equal relation of
$\mathbf{P}$ by $\trianglelefteq$. Then,%
\[
\left(  i,k\right)  \trianglelefteq\left(  i^{\prime},k^{\prime}\right)
\text{ for all }\left(  i,k\right)  \in\mathbf{P}\text{ and }\left(
i^{\prime},k^{\prime}\right)  \in\mathbf{P}\text{ satisfying }\left(  i\geq
i^{\prime}\text{ and }k\leq k^{\prime}\right)  .
\]
(because the total order $\mathbf{P}$ was chosen to be compatible with the
order of $\widetilde{\mathbf{P}}$).

Notice that%
\[
\mathbf{P}=\left\{  1,2,...,p\right\}  \times\left\{  1,2,...,q\right\}
=\operatorname*{Rect}\left(  p,q\right)  \text{ as sets}%
\]
(but not generally as posets).

We define a map $Z$, a field $\mathbb{L}$, a family $\left(  y_{\mathbf{p}%
}\right)  _{\mathbf{p}\in\mathbf{P}}$, a matrix $C$, a ring $\mathbb{F}\left[
x_{\mathbf{P}}\right]  $ and a family $\left(  Q_{\mathbf{p}}\right)
_{\mathbf{p}\in\mathbf{P}}$ as in Lemma \ref{lem.Grasp.generic.independency}.
Also, for every $\left(  i,k\right)  \in\mathbf{P}$, we define polynomials
$\mathfrak{N}_{\left(  i,k\right)  }$ and $\mathfrak{D}_{\left(  i,k\right)
}$ as in Lemma \ref{lem.Grasp.generic.calc}. We regard $\mathbb{F}\left[
x_{\mathbf{P}}\right]  $ as a subring of $\mathbb{F}\left(  x_{\mathbf{P}%
}\right)  $.

Then, Lemma \ref{lem.Grasp.generic.independency} \textbf{(b)} yields that the
family $\left(  Q_{\mathbf{p}}\right)  _{\mathbf{p}\in\mathbf{P}}$ is
$\mathbf{P}$-triangular. Hence, Lemma \ref{lem.algebraic.triangularity}
\textbf{(a)} yields that the family $\left(  Q_{\mathbf{p}}\right)
_{\mathbf{p}\in\mathbf{P}}\in\left(  \mathbb{F}\left(  x_{\mathbf{P}}\right)
\right)  ^{\mathbf{P}}$ is algebraically independent (over $\mathbb{F}$).
Hence, the $\mathbb{F}$-algebra homomorphism $\operatorname*{rev}%
\nolimits_{\left(  Q_{\mathbf{p}}\right)  _{\mathbf{p}\in\mathbf{P}}}$ from
$\mathbb{F}\left(  x_{\mathbf{P}}\right)  $ to $\mathbb{F}\left(
x_{\mathbf{P}}\right)  $ is well-defined. Furthermore, Lemma
\ref{lem.algebraic.triangularity} \textbf{(b)} yields that there exists a
$\mathbf{P}$-triangular family $\left(  R_{\mathbf{p}}\right)  _{\mathbf{p}%
\in\mathbf{P}}\in\left(  \mathbb{F}\left(  x_{\mathbf{P}}\right)  \right)
^{\mathbf{P}}$ such that the maps $\operatorname*{rev}\nolimits_{\left(
Q_{\mathbf{p}}\right)  _{\mathbf{p}\in\mathbf{P}}}$ and $\operatorname*{rev}%
\nolimits_{\left(  R_{\mathbf{p}}\right)  _{\mathbf{p}\in\mathbf{P}}}$ are
mutually inverse. Consider this family $\left(  R_{\mathbf{p}}\right)
_{\mathbf{p}\in\mathbf{P}}$. This family $\left(  R_{\mathbf{p}}\right)
_{\mathbf{p}\in\mathbf{P}}$ is $\mathbf{P}$-triangular and therefore
algebraically independent (by Lemma \ref{lem.algebraic.triangularity}
\textbf{(a)}, applied to $\left(  R_{\mathbf{p}}\right)  _{\mathbf{p}%
\in\mathbf{P}}$ instead of $\left(  Q_{\mathbf{p}}\right)  _{\mathbf{p}%
\in\mathbf{P}}$).

But (\ref{lem.Grasp.generic.calc.f.D}) says that every $\mathbf{p}%
\in\mathbf{P}$ satisfies $\mathfrak{D}_{\mathbf{p}}\neq0$. Renaming
$\mathbf{p}$ as $\mathbf{q}$ in this result, we obtain the following: Every
$\mathbf{q}\in\mathbf{P}$ satisfies $\mathfrak{D}_{\mathbf{q}}\neq0$. Hence,
every $\mathbf{q}\in\mathbf{P}$ satisfies $\mathfrak{D}_{\mathbf{q}}\left(
\left(  R_{\mathbf{p}}\right)  _{\mathbf{p}\in\mathbf{P}}\right)  \neq0$
(because the family $\left(  R_{\mathbf{p}}\right)  _{\mathbf{p}\in\mathbf{P}%
}$ is algebraically independent). Renaming the indices $\mathbf{p}$ and
$\mathbf{q}$ as $\mathbf{r}$ and $\mathbf{v}$ in this result, we obtain the
following: Every $\mathbf{v}\in\mathbf{P}$ satisfies%
\begin{equation}
\mathfrak{D}_{\mathbf{v}}\left(  \left(  R_{\mathbf{r}}\right)  _{\mathbf{r}%
\in\mathbf{P}}\right)  \neq0. \label{pf.Grasp.generic.Dvneq0}%
\end{equation}

Now, we notice that $\mathbf{P}$ is finite, and that every $\mathbf{v}%
\in\mathbf{P}$ satisfies $\mathfrak{D}_{\mathbf{v}}\left(  \left(
R_{\mathbf{r}}\right)  _{\mathbf{r}\in\mathbf{P}}\right)  \neq0$. Hence,
almost every\footnote{As usual, the words ``almost every'' are to be
understood in the sense of Zariski topology.} $f\in\mathbb{K}^{\mathbf{P}}$
satisfies%
\begin{equation}
\left(  R_{\mathbf{r}}\left(  \left(  f\left(  \mathbf{p}\right)  \right)
_{\mathbf{p}\in\mathbf{P}}\right)  \text{ is well-defined for all }%
\mathbf{r}\in\mathbf{P}\right)  \label{pf.Grasp.generic.genass1}%
\end{equation}
and%
\begin{equation}
\left(  \left(  \mathfrak{D}_{\mathbf{v}}\left(  \left(  R_{\mathbf{r}%
}\right)  _{\mathbf{r}\in\mathbf{P}}\right)  \right)  \left(  \left(  f\left(
\mathbf{p}\right)  \right)  _{\mathbf{p}\in\mathbf{P}}\right)  \neq0\text{ for
all }\mathbf{v}\in\mathbf{P}\right)  . \label{pf.Grasp.generic.genass2}%
\end{equation}

Fix an $f\in\mathbb{K}^{\mathbf{P}}$ satisfying
(\ref{pf.Grasp.generic.genass1}) and (\ref{pf.Grasp.generic.genass2}). Define
a matrix $B\in\mathbb{K}^{p\times q}$ by%
\[
B=\left(  R_{\left(  i,Z\left(  k\right)  \right)  }\left(  \left(  f\left(
\mathbf{p}\right)  \right)  _{\mathbf{p}\in\mathbf{P}}\right)  \right)
_{1\leq i\leq p,\ 1\leq k\leq q}.
\]

Fix $\mathbf{v}\in\mathbf{P}$. We have%
\begin{align*}
\operatorname*{rev}\nolimits_{\left(  Q_{\mathbf{p}}\right)  _{\mathbf{p}%
\in\mathbf{P}}}\left(  x_{\mathbf{v}}\right)   &  =x_{\mathbf{v}}\left(
\left(  Q_{\mathbf{p}}\right)  _{\mathbf{p}\in\mathbf{P}}\right)
\ \ \ \ \ \ \ \ \ \ \left(  \text{by the definition of }\operatorname*{rev}%
\nolimits_{\left(  Q_{\mathbf{p}}\right)  _{\mathbf{p}\in\mathbf{P}}}\right)
\\
&  =Q_{\mathbf{v}}.
\end{align*}
But we have $\operatorname*{rev}\nolimits_{\left(  R_{\mathbf{p}}\right)
_{\mathbf{p}\in\mathbf{P}}}\circ\operatorname*{rev}\nolimits_{\left(
Q_{\mathbf{p}}\right)  _{\mathbf{p}\in\mathbf{P}}}=\operatorname*{id}$ (since
the maps $\operatorname*{rev}\nolimits_{\left(  R_{\mathbf{p}}\right)
_{\mathbf{p}\in\mathbf{P}}}$ and $\operatorname*{rev}\nolimits_{\left(
Q_{\mathbf{p}}\right)  _{\mathbf{p}\in\mathbf{P}}}$ are mutually inverse).
Hence, $\left(  \operatorname*{rev}\nolimits_{\left(  R_{\mathbf{p}}\right)
_{\mathbf{p}\in\mathbf{P}}}\circ\operatorname*{rev}\nolimits_{\left(
Q_{\mathbf{p}}\right)  _{\mathbf{p}\in\mathbf{P}}}\right)  \left(
x_{\mathbf{v}}\right)  =\operatorname*{id}\left(  x_{\mathbf{v}}\right)
=x_{\mathbf{v}}$, so that%
\begin{align}
x_{\mathbf{v}}  &  =\left(  \operatorname*{rev}\nolimits_{\left(
R_{\mathbf{p}}\right)  _{\mathbf{p}\in\mathbf{P}}}\circ\operatorname*{rev}%
\nolimits_{\left(  Q_{\mathbf{p}}\right)  _{\mathbf{p}\in\mathbf{P}}}\right)
\left(  x_{\mathbf{v}}\right)  =\operatorname*{rev}\nolimits_{\left(
R_{\mathbf{p}}\right)  _{\mathbf{p}\in\mathbf{P}}}\left(
\underbrace{\operatorname*{rev}\nolimits_{\left(  Q_{\mathbf{p}}\right)
_{\mathbf{p}\in\mathbf{P}}}\left(  x_{\mathbf{v}}\right)  }_{=Q_{\mathbf{v}}%
}\right) \nonumber\\
&  =\operatorname*{rev}\nolimits_{\left(  R_{\mathbf{p}}\right)
_{\mathbf{p}\in\mathbf{P}}}\left(  Q_{\mathbf{v}}\right)  =Q_{\mathbf{v}%
}\left(  \left(  R_{\mathbf{p}}\right)  _{\mathbf{p}\in\mathbf{P}}\right)
\ \ \ \ \ \ \ \ \ \ \left(  \text{by the definition of }\operatorname*{rev}%
\nolimits_{\left(  R_{\mathbf{p}}\right)  _{\mathbf{p}\in\mathbf{P}}}\right)
\nonumber\\
&  =Q_{\mathbf{v}}\left(  \left(  R_{\mathbf{r}}\right)  _{\mathbf{r}%
\in\mathbf{P}}\right)  \label{pf.Grasp.generic.identity}%
\end{align}
(here, we renamed the index $\mathbf{p}$ as $\mathbf{r}$).

Recall that $R_{\mathbf{r}}\left(  \left(  f\left(  \mathbf{p}\right)
\right)  _{\mathbf{p}\in\mathbf{P}}\right)  $ is well-defined for all
$\mathbf{r}\in\mathbf{P}$ (by (\ref{pf.Grasp.generic.genass1})). Therefore,
$\mathfrak{D}_{\mathbf{v}}\left(  \left(  R_{\mathbf{r}}\left(  \left(
f\left(  \mathbf{p}\right)  \right)  _{\mathbf{p}\in\mathbf{P}}\right)
\right)  _{\mathbf{r}\in\mathbf{P}}\right)  $ is well-defined (since
$\mathfrak{D}_{\mathbf{v}}$ is a polynomial). Hence, we can apply Lemma
\ref{lem.algebraic.triangularity.substass} to $\left(  f\left(  \mathbf{p}%
\right)  \right)  _{\mathbf{p}\in\mathbf{P}}$, $\left(  R_{\mathbf{r}}\right)
_{\mathbf{r}\in\mathbf{P}}$ and $\mathfrak{D}_{\mathbf{v}}$ instead of
$\left(  A_{\mathbf{p}}\right)  _{\mathbf{p}\in\mathbf{P}}$, $\left(
B_{\mathbf{r}}\right)  _{\mathbf{r}\in\mathbf{P}}$ and $f$. As a consequence,
we obtain that the element $\mathfrak{D}_{\mathbf{v}}\left(  \left(
R_{\mathbf{r}}\right)  _{\mathbf{r}\in\mathbf{P}}\right)  $ of $\mathbb{F}%
\left(  x_{\mathbf{P}}\right)  $ and the element $\left(  \mathfrak{D}%
_{\mathbf{v}}\left(  \left(  R_{\mathbf{r}}\right)  _{\mathbf{r}\in\mathbf{P}%
}\right)  \right)  \left(  \left(  f\left(  \mathbf{p}\right)  \right)
_{\mathbf{p}\in\mathbf{P}}\right)  $ of $\mathbb{K}$ are also well-defined,
and satisfy%
\[
\left(  \mathfrak{D}_{\mathbf{v}}\left(  \left(  R_{\mathbf{r}}\right)
_{\mathbf{r}\in\mathbf{P}}\right)  \right)  \left(  \left(  f\left(
\mathbf{p}\right)  \right)  _{\mathbf{p}\in\mathbf{P}}\right)  =\mathfrak{D}%
_{\mathbf{v}}\left(  \left(  R_{\mathbf{r}}\left(  \left(  f\left(
\mathbf{p}\right)  \right)  _{\mathbf{p}\in\mathbf{P}}\right)  \right)
_{\mathbf{r}\in\mathbf{P}}\right)  .
\]
Hence, $\mathfrak{D}_{\mathbf{v}}\left(  \left(  R_{\mathbf{r}}\left(  \left(
f\left(  \mathbf{p}\right)  \right)  _{\mathbf{p}\in\mathbf{P}}\right)
\right)  _{\mathbf{r}\in\mathbf{P}}\right)  =\left(  \mathfrak{D}_{\mathbf{v}%
}\left(  \left(  R_{\mathbf{r}}\right)  _{\mathbf{r}\in\mathbf{P}}\right)
\right)  \left(  \left(  f\left(  \mathbf{p}\right)  \right)  _{\mathbf{p}%
\in\mathbf{P}}\right)  \neq0$ (by (\ref{pf.Grasp.generic.genass2})).
Therefore, Lemma \ref{lem.Grasp.generic.independency} \textbf{(d)} (applied to
$\mathbb{K}$, $\mathbf{v}$, $\left(  R_{\mathbf{r}}\left(  \left(  f\left(
\mathbf{p}\right)  \right)  _{\mathbf{p}\in\mathbf{P}}\right)  \right)
_{\mathbf{r}\in\mathbf{P}}$ and $B$ instead of $\mathbb{M}$, $\mathbf{p}$ and
$\left(  d_{\mathbf{r}}\right)  _{\mathbf{r}\in\mathbf{P}}$ and $D$) yields
that the values $Q_{\mathbf{v}}\left(  \left(  R_{\mathbf{r}}\left(  \left(
f\left(  \mathbf{p}\right)  \right)  _{\mathbf{p}\in\mathbf{P}}\right)
\right)  _{\mathbf{r}\in\mathbf{P}}\right)  $ and $\left(
\operatorname*{Grasp}\nolimits_{0}\left(  I_{p}\mid B\right)  \right)  \left(
\mathbf{v}\right)  $ are well-defined and satisfy%
\begin{equation}
Q_{\mathbf{v}}\left(  \left(  R_{\mathbf{r}}\left(  \left(  f\left(
\mathbf{p}\right)  \right)  _{\mathbf{p}\in\mathbf{P}}\right)  \right)
_{\mathbf{r}\in\mathbf{P}}\right)  =\left(  \operatorname*{Grasp}%
\nolimits_{0}\left(  I_{p}\mid B\right)  \right)  \left(  \mathbf{v}\right)
\label{pf.Grasp.generic.step1}%
\end{equation}
(because $B=\left(  R_{\left(  i,Z\left(  k\right)  \right)  }\left(  \left(
f\left(  \mathbf{p}\right)  \right)  _{\mathbf{p}\in\mathbf{P}}\right)
\right)  _{1\leq i\leq p,\ 1\leq k\leq q}$).

Now, recall that $R_{\mathbf{r}}\left(  \left(  f\left(  \mathbf{p}\right)
\right)  _{\mathbf{p}\in\mathbf{P}}\right)  $ is well-defined for all
$\mathbf{r}\in\mathbf{P}$ (by (\ref{pf.Grasp.generic.genass1})), and that the
value $Q_{\mathbf{v}}\left(  \left(  R_{\mathbf{r}}\left(  \left(  f\left(
\mathbf{p}\right)  \right)  _{\mathbf{p}\in\mathbf{P}}\right)  \right)
_{\mathbf{r}\in\mathbf{P}}\right)  $ is well-defined (as we now know). Hence,
we can apply Lemma \ref{lem.algebraic.triangularity.substass} to $\left(
f\left(  \mathbf{p}\right)  \right)  _{\mathbf{p}\in\mathbf{P}}$, $\left(
R_{\mathbf{r}}\right)  _{\mathbf{r}\in\mathbf{P}}$ and $Q_{\mathbf{v}}$
instead of $\left(  A_{\mathbf{p}}\right)  _{\mathbf{p}\in\mathbf{P}}$,
$\left(  B_{\mathbf{r}}\right)  _{\mathbf{r}\in\mathbf{P}}$ and $f$. As a
consequence, we obtain that the element $Q_{\mathbf{v}}\left(  \left(
R_{\mathbf{r}}\right)  _{\mathbf{r}\in\mathbf{P}}\right)  $ of $\mathbb{F}%
\left(  x_{\mathbf{P}}\right)  $ and the element $\left(  Q_{\mathbf{v}%
}\left(  \left(  R_{\mathbf{r}}\right)  _{\mathbf{r}\in\mathbf{P}}\right)
\right)  \left(  \left(  f\left(  \mathbf{p}\right)  \right)  _{\mathbf{p}%
\in\mathbf{P}}\right)  $ of $\mathbb{K}$ are also well-defined, and satisfy%
\[
\left(  Q_{\mathbf{v}}\left(  \left(  R_{\mathbf{r}}\right)  _{\mathbf{r}%
\in\mathbf{P}}\right)  \right)  \left(  \left(  f\left(  \mathbf{p}\right)
\right)  _{\mathbf{p}\in\mathbf{P}}\right)  =Q_{\mathbf{v}}\left(  \left(
R_{\mathbf{r}}\left(  \left(  f\left(  \mathbf{p}\right)  \right)
_{\mathbf{p}\in\mathbf{P}}\right)  \right)  _{\mathbf{r}\in\mathbf{P}}\right)
.
\]
Thus,%
\begin{align*}
Q_{\mathbf{v}}\left(  \left(  R_{\mathbf{r}}\left(  \left(  f\left(
\mathbf{p}\right)  \right)  _{\mathbf{p}\in\mathbf{P}}\right)  \right)
_{\mathbf{r}\in\mathbf{P}}\right)   &  =\left(  \underbrace{Q_{\mathbf{v}%
}\left(  \left(  R_{\mathbf{r}}\right)  _{\mathbf{r}\in\mathbf{P}}\right)
}_{\substack{=x_{\mathbf{v}}\\\text{(by (\ref{pf.Grasp.generic.identity}))}%
}}\right)  \left(  \left(  f\left(  \mathbf{p}\right)  \right)  _{\mathbf{p}%
\in\mathbf{P}}\right) \\
&  =x_{\mathbf{v}}\left(  \left(  f\left(  \mathbf{p}\right)  \right)
_{\mathbf{p}\in\mathbf{P}}\right)  =f\left(  \mathbf{v}\right)  .
\end{align*}
Compared with (\ref{pf.Grasp.generic.step1}), this yields%
\[
f\left(  \mathbf{v}\right)  =\left(  \operatorname*{Grasp}\nolimits_{0}\left(
I_{p}\mid B\right)  \right)  \left(  \mathbf{v}\right)  .
\]

Now, forget that we fixed $\mathbf{v}$. We thus have shown that $f\left(
\mathbf{v}\right)  =\left(  \operatorname*{Grasp}\nolimits_{0}\left(
I_{p}\mid B\right)  \right)  \left(  \mathbf{v}\right)  $ for every
$\mathbf{v}\in\mathbf{P}$. In other words, $f=\operatorname*{Grasp}%
\nolimits_{0}\left(  I_{p}\mid B\right)  $. Hence, there exists a matrix
$A\in\mathbb{K}^{p\times\left(  p+q\right)  }$ satisfying
$f=\operatorname*{Grasp}\nolimits_{0}A$ (namely, $A=\operatorname*{Grasp}%
\nolimits_{0}\left(  I_{p}\mid B\right)  $).

Now, forget that we fixed $f$. We have thus shown that for every
$f\in\mathbb{K}^{\mathbf{P}}$ satisfying (\ref{pf.Grasp.generic.genass1}) and
(\ref{pf.Grasp.generic.genass2}), there exists a matrix $A\in\mathbb{K}%
^{p\times\left(  p+q\right)  }$ satisfying $f=\operatorname*{Grasp}%
\nolimits_{0}A$. Since almost every $f\in\mathbb{K}^{\mathbf{P}}$ satisfies
(\ref{pf.Grasp.generic.genass1}) and (\ref{pf.Grasp.generic.genass2}), we have
therefore proven that for almost every $f\in\mathbb{K}^{\mathbf{P}}$, there
exists a matrix $A\in\mathbb{K}^{p\times\left(  p+q\right)  }$ satisfying
$f=\operatorname*{Grasp}\nolimits_{0}A$. Since $\mathbf{P}%
=\operatorname*{Rect}\left(  p,q\right)  $ as sets, this rewrites as follows:
For almost every $f\in\mathbb{K}^{\operatorname*{Rect}\left(  p,q\right)  }$,
there exists a matrix $A\in\mathbb{K}^{p\times\left(  p+q\right)  }$
satisfying $f=\operatorname*{Grasp}\nolimits_{0}A$. In other words,
Proposition \ref{prop.Grasp.generic} is proven.
\end{proof}

It now remains to prove Lemma \ref{lem.Grasp.generic.independency}.

\begin{proof}
[Proof of Lemma \ref{lem.Grasp.generic.independency} (sketched).]Let us first
introduce some notation:

For every subset $\mathbf{S}$ of $\mathbf{P}$, let $\mathbb{F}\left[
x_{\mathbf{S}}\right]  $ denote the polynomial ring over $\mathbb{F}$ in the
indeterminates $x_{\mathbf{p}}$ with $\mathbf{p}$ ranging over all elements of
$\mathbf{S}$. We identify $\mathbb{F}\left[  x_{\mathbf{S}}\right]  $ with a
subring of $\mathbb{F}\left[  x_{\mathbf{P}}\right]  $ for every subset
$\mathbf{S}$ of $\mathbf{P}$, and we identify $\mathbb{F}\left[
x_{\mathbf{P}}\right]  $ with a subring of $\mathbb{F}\left(  x_{\mathbf{P}%
}\right)  $.

For every $\mathbf{p}\in\mathbf{P}$, let $\mathbf{p}\Downarrow$ denote the
subset $\left\{  \mathbf{v}\in\mathbf{P}\ \mid\ \mathbf{v}\vartriangleleft
\mathbf{p}\right\}  $ of $\mathbf{P}$.

\textbf{(c)} \footnote{Here is a brief sketch of how we are going to prove
Lemma \ref{lem.Grasp.generic.independency} \textbf{(c)} (the actual proof in
all its detail is given in the body of the text and does not rely on this
sketch):
\par
We have $C=\left(  y_{\left(  i,k\right)  }\right)  _{1\leq
i\leq p,\ 1\leq k\leq q}=\left(  x_{\left(  i,Z\left(  k\right)  \right)
}\right)  _{1\leq i\leq p,\ 1\leq k\leq q}$ (because $y_{\left(  i,k\right)
}=x_{\left(  i,Z\left(  k\right)  \right)  }$ for all $i$ and $k$). Hence, if
we substitute the family $\left(  d_{\mathbf{r}}\right)  _{\mathbf{r}%
\in\mathbf{P}}$ for $\left(  x_{\mathbf{r}}\right)  _{\mathbf{r}\in\mathbf{P}%
}$, then the matrix $C$ becomes $\left(  d_{\left(  i,Z\left(  k\right)
\right)  }\right)  _{1\leq i\leq p,\ 1\leq k\leq q}=D$.
\par
Now, let $\left(  i,k\right)  \in\mathbf{P}$. The formula
(\ref{lem.Grasp.generic.independency.Ndef}) is a polynomial identity in the
variables $x_{\mathbf{r}}$. Substituting the family $\left(  d_{\mathbf{r}%
}\right)  _{\mathbf{r}\in\mathbf{P}}$ for $\left(  x_{\mathbf{r}}\right)
_{\mathbf{r}\in\mathbf{P}}$ in this identity, we obtain%
\[
\mathfrak{N}_{\left(  i,k\right)  }\left(  \left(  d_{\mathbf{r}}\right)
_{\mathbf{r}\in\mathbf{P}}\right)  =\det\left(  \left(  I_{p}\mid D\right)
\left[  1:i\mid i+k-1:p+k\right]  \right)
\]
(because we know that $C$ becomes $D$ under this substitution). Thus,
(\ref{lem.Grasp.generic.independency.cN}) is proven. Similarly,
(\ref{lem.Grasp.generic.independency.cD}) follows, so that Lemma
\ref{lem.Grasp.generic.independency} \textbf{(c)} is proven. We give a
detailed version of this argument in the actual body of the text.} By the
universal property of the polynomial ring $\mathbb{F}\left[  x_{\mathbf{P}%
}\right]  $, we have the following: For every commutative $\mathbb{F}$-algebra
$\mathbb{A}$ and every family $\left(  \eta_{\mathbf{r}}\right)
_{\mathbf{r}\in\mathbf{P}}\in\mathbb{A}^{\mathbf{P}}$ of elements of
$\mathbb{A}$, there exists a unique $\mathbb{F}$-algebra homomorphism
$\psi:\mathbb{F}\left[  x_{\mathbf{P}}\right]  \rightarrow\mathbb{A}$ such
that $\left(  \psi\left(  x_{\mathbf{r}}\right)  =\eta_{\mathbf{r}}\text{ for
every }\mathbf{r}\in\mathbf{P}\right)  $. This $\mathbb{F}$-algebra
homomorphism is denoted by $\operatorname*{ev}\nolimits_{\left(
\eta_{\mathbf{r}}\right)  _{\mathbf{r}\in\mathbf{P}}}$ and is called the
\textit{evaluation map at }$\left(  \eta_{\mathbf{r}}\right)  _{\mathbf{r}%
\in\mathbf{P}}$. It satisfies%
\begin{equation}
\operatorname*{ev}\nolimits_{\left(  \eta_{\mathbf{r}}\right)  _{\mathbf{r}%
\in\mathbf{P}}}\left(  f\right)  =f\left(  \left(  \eta_{\mathbf{r}}\right)
_{\mathbf{r}\in\mathbf{P}}\right)  \ \ \ \ \ \ \ \ \ \ \text{for every }%
f\in\mathbb{F}\left[  x_{\mathbf{P}}\right]  .
\label{pf.Grasp.generic.independency.d.ev}%
\end{equation}

Thus, the family $\left(  d_{\mathbf{r}}\right)  _{\mathbf{r}\in\mathbf{P}}%
\in\mathbb{M}^{\mathbf{P}}$ induces an $\mathbb{F}$-algebra homomorphism
$\operatorname*{ev}\nolimits_{\left(  d_{\mathbf{r}}\right)  _{\mathbf{r}%
\in\mathbf{P}}}:\mathbb{F}\left[  x_{\mathbf{P}}\right]  \rightarrow
\mathbb{M}$. Denote this homomorphism $\left(  d_{\mathbf{r}}\right)
_{\mathbf{r}\in\mathbf{P}}$ by $\mathfrak{d}$. Thus,%
\[
\mathfrak{d}=\operatorname*{ev}\nolimits_{\left(  d_{\mathbf{r}}\right)
_{\mathbf{r}\in\mathbf{P}}}.
\]
Hence, for every $f\in\mathbb{F}\left[  x_{\mathbf{P}}\right]  $, we have%
\begin{equation}
\mathfrak{d}\left(  f\right)  =\operatorname*{ev}\nolimits_{\left(
d_{\mathbf{r}}\right)  _{\mathbf{r}\in\mathbf{P}}}\left(  f\right)  =f\left(
\left(  d_{\mathbf{r}}\right)  _{\mathbf{r}\in\mathbf{P}}\right)
\label{pf.Grasp.generic.independency.d.d}%
\end{equation}
(by (\ref{pf.Grasp.generic.independency.d.ev}), applied to $\left(
\eta_{\mathbf{r}}\right)  _{\mathbf{r}\in\mathbf{P}}=\left(  d_{\mathbf{r}%
}\right)  _{\mathbf{r}\in\mathbf{P}}$).

We introduce a further notation: If $\mathbb{A}$ and $\mathbb{B}$ are two
$\mathbb{F}$-algebras, if $\varphi:\mathbb{A}\rightarrow\mathbb{B}$ is an
$\mathbb{F}$-algebra homomorphism, and if $u$ and $v$ are two nonnegative
integers, then $\varphi^{u\times v}$ will denote the $\mathbb{F}$-module
homomorphism $\mathbb{A}^{u\times v}\rightarrow\mathbb{B}^{u\times v}$
canonically induced by $\varphi$ on the groups of $u\times v$-matrices.
(Explicitly, this homomorphism $\varphi^{u\times v}$ takes any matrix in
$\mathbb{A}^{u\times v}$ and applies $\varphi$ to each of its entries.)

It is clear that the operation of attaching two matrices to each other is
canonical with respect to the ground ring. By this we mean that if
$\mathbb{A}$ and $\mathbb{B}$ are two $\mathbb{F}$-algebras, if $\varphi
:\mathbb{A}\rightarrow\mathbb{B}$ is an $\mathbb{F}$-algebra homomorphism, if
$u$, $v_{1}$ and $v_{2}$ are three nonnegative integers, and if $A_{1}%
\in\mathbb{A}^{u\times v_{1}}$ and $A_{2}\in\mathbb{A}^{u\times v_{2}}$ are
two matrices, then%
\begin{equation}
\varphi^{u\times\left(  v_{1}+v_{2}\right)  }\left(  A_{1}\mid A_{2}\right)
=\left(  \varphi^{u\times v_{1}}\left(  A_{1}\right)  \mid\varphi^{u\times
v_{2}}\left(  A_{2}\right)  \right)  .
\label{pf.Grasp.generic.independency.d.attach}%
\end{equation}

It is also clear that the definition of $A\left[  a:b\mid c:d\right]  $ in
Definition \ref{def.minors} is canonical with respect to the ground ring. By
this we mean that if $\mathbb{A}$ and $\mathbb{B}$ are two $\mathbb{F}%
$-algebras, if $\varphi:\mathbb{A}\rightarrow\mathbb{B}$ is an $\mathbb{F}%
$-algebra homomorphism, if $u$ and $v$ are nonnegative integers, if
$A\in\mathbb{A}^{u\times v}$ is a matrix, and if $a$, $b$, $c$ and $d$ are
four integers satisfying $a\leq b$ and $c\leq d$, then%
\begin{equation}
\varphi^{u\times\left(  b-a+d-c\right)  }\left(  A\left[  a:b\mid c:d\right]
\right)  =\left(  \varphi^{u\times v}\left(  A\right)  \right)  \left[
a:b\mid c:d\right]  . \label{pf.Grasp.generic.independency.d.interval}%
\end{equation}

Furthermore, the definition of the determinant of a square matrix is canonical
with respect to the ground ring. By this we mean that if $\mathbb{A}$ and
$\mathbb{B}$ are two commutative $\mathbb{F}$-algebras, if
$\varphi : \mathbb{A} \rightarrow\mathbb{B}$ is an
$\mathbb{F}$-algebra homomorphism, if $u$ is a
nonnegative integer, and if $A\in\mathbb{A}^{u\times u}$ is a square matrix,
then%
\begin{equation}
\varphi\left(  \det A\right)  =\det\left(  \varphi^{u\times u}\left(
A\right)  \right)  . \label{pf.Grasp.generic.independency.d.det}%
\end{equation}

Finally, the definition of the identity matrix is canonical with respect to
the ground ring. By this we mean that if $\mathbb{A}$ and $\mathbb{B}$ are two
$\mathbb{F}$-algebras, if $\varphi:\mathbb{A}\rightarrow\mathbb{B}$ is an
$\mathbb{F}$-algebra homomorphism, and if $u$ is a nonnegative integer, then%
\begin{equation}
\varphi^{u\times u}\left(  I_{u}\right)  =I_{u}.
\label{pf.Grasp.generic.independency.d.Iu}%
\end{equation}

Since%
\[
C=\left(  \underbrace{y_{\left(  i,k\right)  }}_{\substack{=x_{\left(
i,Z\left(  k\right)  \right)  }\\\text{(by the definition}\\\text{of
}y_{\left(  i,k\right)  }\text{)}}}\right)  _{1\leq i\leq p,\ 1\leq k\leq
q}=\left(  x_{\left(  i,Z\left(  k\right)  \right)  }\right)  _{1\leq i\leq
p,\ 1\leq k\leq q},
\]
we have%
\[
\mathfrak{d}^{p\times q}\left(  C\right)  =\mathfrak{d}^{p\times q}\left(
\left(  x_{\left(  i,Z\left(  k\right)  \right)  }\right)  _{1\leq i\leq
p,\ 1\leq k\leq q}\right)  =\left(  \mathfrak{d}\left(  x_{\left(  i,Z\left(
k\right)  \right)  }\right)  \right)  _{1\leq i\leq p,\ 1\leq k\leq q}%
\]
(by the definition of $\mathfrak{d}^{p\times q}$). Since every $i\in\left\{
1,2,...,p\right\}  $ and $k\in\left\{  1,2,...,q\right\}  $ satisfy%
\begin{align*}
\mathfrak{d}\left(  x_{\left(  i,Z\left(  k\right)  \right)  }\right)   &
=x_{\left(  i,Z\left(  k\right)  \right)  }\left(  \left(  d_{\mathbf{r}%
}\right)  _{\mathbf{r}\in\mathbf{P}}\right)  \ \ \ \ \ \ \ \ \ \ \left(
\text{by (\ref{pf.Grasp.generic.independency.d.d}), applied to }f=x_{\left(
i,Z\left(  k\right)  \right)  }\right) \\
&  =d_{\left(  i,Z\left(  k\right)  \right)  },
\end{align*}
this becomes%
\begin{equation}
\mathfrak{d}^{p\times q}\left(  C\right)  =\left(  \underbrace{\mathfrak{d}%
\left(  x_{\left(  i,Z\left(  k\right)  \right)  }\right)  }_{=d_{\left(
i,Z\left(  k\right)  \right)  }}\right)  _{1\leq i\leq p,\ 1\leq k\leq
q}=\left(  d_{\left(  i,Z\left(  k\right)  \right)  }\right)  _{1\leq i\leq
p,\ 1\leq k\leq q}=D. \label{pf.Grasp.generic.independency.d.dC}%
\end{equation}

Now, let $\left(  i,k\right)  \in\mathbf{P}$. Applying
(\ref{pf.Grasp.generic.independency.d.d}) to $f=\mathfrak{N}_{\left(
i,k\right)  }$, we obtain $\mathfrak{d}\left(  \mathfrak{N}_{\left(
i,k\right)  }\right)  =\mathfrak{N}_{\left(  i,k\right)  }\left(  \left(
d_{\mathbf{r}}\right)  _{\mathbf{r}\in\mathbf{P}}\right)  $. Hence,%
\begin{align*}
\mathfrak{N}_{\left(  i,k\right)  }\left(  \left(  d_{\mathbf{r}}\right)
_{\mathbf{r}\in\mathbf{P}}\right)   &  =\mathfrak{d}\left(  \mathfrak{N}%
_{\left(  i,k\right)  }\right)  =\mathfrak{d}\left(  \det\left(  \left(
I_{p}\mid C\right)  \left[  1:i\mid i+k-1:p+k\right]  \right)  \right) \\
&  \ \ \ \ \ \ \ \ \ \ \left(  \text{by
(\ref{lem.Grasp.generic.independency.Ndef})}\right) \\
&  =\det\left(  \underbrace{\mathfrak{d}^{p\times p}}_{\substack{=\mathfrak{d}%
^{p\times\left(  i-1+\left(  p+k\right)  -\left(  i+k-1\right)  \right)
}\\\text{(since }p=i-1+\left(  p+k\right)  -\left(  i+k-1\right)  \text{)}%
}}\left(  \left(  I_{p}\mid C\right)  \left[  1:i\mid i+k-1:p+k\right]
\right)  \right) \\
&  \ \ \ \ \ \ \ \ \ \ \left(
\begin{array}
[c]{c}%
\text{by (\ref{pf.Grasp.generic.independency.d.det}), applied to }%
\mathbb{F}\left[  x_{\mathbf{P}}\right]  \text{, }\mathbb{M}\text{,
}\mathfrak{d}\text{, }p\text{ and}\\
\left(  I_{p}\mid C\right)  \left[  1:i\mid i+k-1:p+k\right]  \text{ instead
of }\mathbb{A}\text{, }\mathbb{B}\text{, }\varphi\text{, }u\text{ and }A
\end{array}
\right) \\
&  =\det\left(  \mathfrak{d}^{p\times\left(  i-1+\left(  p+k\right)  -\left(
i+k-1\right)  \right)  }\left(  \left(  I_{p}\mid C\right)  \left[  1:i\mid
i+k-1:p+k\right]  \right)  \right)  .
\end{align*}
Since%
\begin{align*}
&  \mathfrak{d}^{p\times\left(  i-1+\left(  p+k\right)  -\left(  i+k-1\right)
\right)  }\left(  \left(  I_{p}\mid C\right)  \left[  1:i\mid
i+k-1:p+k\right]  \right) \\
&  =\left(  \underbrace{\mathfrak{d}^{p\times\left(  p+q\right)  }\left(
I_{p}\mid C\right)  }_{\substack{=\left(  \mathfrak{d}^{p\times p}\left(
I_{p}\right)  \mid\mathfrak{d}^{p\times q}\left(  C\right)  \right)
\\\text{(by (\ref{pf.Grasp.generic.independency.d.attach}), applied
to}\\\mathbb{F}\left[  x_{\mathbf{P}}\right]  \text{, }\mathbb{M}\text{,
}\mathfrak{d}\text{, }p\text{, }p\text{, }q\text{, }I_{p}\text{ and }C\text{
instead of}\\\mathbb{A}\text{, }\mathbb{B}\text{, }\varphi\text{, }u\text{,
}v_{1}\text{, }v_{2}\text{, }A_{1}\text{ and }A_{2}\text{)}}}\right)  \left[
1:i\mid i+k-1:p+k\right] \\
&  \ \ \ \ \ \ \ \ \ \ \left(
\begin{array}
[c]{c}%
\text{by (\ref{pf.Grasp.generic.independency.d.interval}), applied to
}\mathbb{F}\left[  x_{\mathbf{P}}\right]  \text{, }\mathbb{M}\text{,
}\mathfrak{d}\text{, }p\text{, }p+q\text{,}\\
1\text{, }i\text{, }i+k-1\text{, }p+k\text{ and }\left(  I_{p}\mid C\right)
\text{ instead of }\mathbb{A}\text{, }\mathbb{B}\text{, }\varphi\text{,
}u\text{, }v\text{, }a\text{, }b\text{, }c\text{, }d\text{ and }A
\end{array}
\right) \\
&  =\left(  \underbrace{\mathfrak{d}^{p\times p}\left(  I_{p}\right)
}_{\substack{=I_{p}\\\text{(by (\ref{pf.Grasp.generic.independency.d.Iu}),
applied to}\\\mathbb{F}\left[  x_{\mathbf{P}}\right]  \text{, }\mathbb{M}%
\text{, }\mathfrak{d}\text{, and }p\text{ instead of}\\\mathbb{A}\text{,
}\mathbb{B}\text{, }\varphi\text{, and }u\text{)}}}\mid
\underbrace{\mathfrak{d}^{p\times q}\left(  C\right)  }%
_{\substack{=D\\\text{(by (\ref{pf.Grasp.generic.independency.d.dC}))}%
}}\right)  \left[  1:i\mid i+k-1:p+k\right] \\
&  =\left(  I_{p}\mid D\right)  \left[  1:i\mid i+k-1:p+k\right]  ,
\end{align*}
this becomes%
\begin{align*}
\mathfrak{N}_{\left(  i,k\right)  }\left(  \left(  d_{\mathbf{r}}\right)
_{\mathbf{r}\in\mathbf{P}}\right)   &  =\det\left(  \underbrace{\mathfrak{d}%
^{p\times\left(  i-1+\left(  p+k\right)  -\left(  i-k-1\right)  \right)
}\left(  \left(  I_{p}\mid C\right)  \left[  1:i\mid i+k-1:p+k\right]
\right)  }_{=\left(  I_{p}\mid D\right)  \left[  1:i\mid i+k-1:p+k\right]
}\right) \\
&  =\det\left(  \left(  I_{p}\mid D\right)  \left[  1:i\mid i+k-1:p+k\right]
\right)  .
\end{align*}
This proves (\ref{lem.Grasp.generic.independency.cN}).

On the other hand, applying (\ref{pf.Grasp.generic.independency.d.d}) to
$f=\mathfrak{D}_{\left(  i,k\right)  }$, we obtain $\mathfrak{d}\left(
\mathfrak{D}_{\left(  i,k\right)  }\right)  =\mathfrak{D}_{\left(  i,k\right)
}\left(  \left(  d_{\mathbf{r}}\right)  _{\mathbf{r}\in\mathbf{P}}\right)  $.
Hence,%
\begin{align*}
\mathfrak{D}_{\left(  i,k\right)  }\left(  \left(  d_{\mathbf{r}}\right)
_{\mathbf{r}\in\mathbf{P}}\right)   &  =\mathfrak{d}\left(  \mathfrak{D}%
_{\left(  i,k\right)  }\right)  =\mathfrak{d}\left(  \det\left(  \left(
I_{p}\mid C\right)  \left[  0:i\mid i+k:p+k\right]  \right)  \right) \\
&  \ \ \ \ \ \ \ \ \ \ \left(  \text{by
(\ref{lem.Grasp.generic.independency.Ddef})}\right) \\
&  =\det\left(  \underbrace{\mathfrak{d}^{p\times p}}_{\substack{=\mathfrak{d}%
^{p\times\left(  i-0+\left(  p+k\right)  -\left(  i+k\right)  \right)
}\\\text{(since }p=i-0+\left(  p+k\right)  -\left(  i+k\right)  \text{)}%
}}\left(  \left(  I_{p}\mid C\right)  \left[  0:i\mid i+k:p+k\right]  \right)
\right) \\
&  \ \ \ \ \ \ \ \ \ \ \left(
\begin{array}
[c]{c}%
\text{by (\ref{pf.Grasp.generic.independency.d.det}), applied to }%
\mathbb{F}\left[  x_{\mathbf{P}}\right]  \text{, }\mathbb{M}\text{,
}\mathfrak{d}\text{, }p\text{ and}\\
\left(  I_{p}\mid C\right)  \left[  0:i\mid i+k:p+k\right]  \text{ instead of
}\mathbb{A}\text{, }\mathbb{B}\text{, }\varphi\text{, }u\text{ and }A
\end{array}
\right) \\
&  =\det\left(  \mathfrak{d}^{p\times\left(  i-0+\left(  p+k\right)  -\left(
i+k\right)  \right)  }\left(  \left(  I_{p}\mid C\right)  \left[  0:i\mid
i+k:p+k\right]  \right)  \right)  .
\end{align*}
Since%
\begin{align*}
&  \mathfrak{d}^{p\times\left(  i-0+\left(  p+k\right)  -\left(  i+k\right)
\right)  }\left(  \left(  I_{p}\mid C\right)  \left[  0:i\mid i+k:p+k\right]
\right) \\
&  =\left(  \underbrace{\mathfrak{d}^{p\times\left(  p+q\right)  }\left(
I_{p}\mid C\right)  }_{\substack{=\left(  \mathfrak{d}^{p\times p}\left(
I_{p}\right)  \mid\mathfrak{d}^{p\times q}\left(  C\right)  \right)
\\\text{(by (\ref{pf.Grasp.generic.independency.d.attach}), applied
to}\\\mathbb{F}\left[  x_{\mathbf{P}}\right]  \text{, }\mathbb{M}\text{,
}\mathfrak{d}\text{, }p\text{, }p\text{, }q\text{, }I_{p}\text{ and }C\text{
instead of}\\\mathbb{A}\text{, }\mathbb{B}\text{, }\varphi\text{, }u\text{,
}v_{1}\text{, }v_{2}\text{, }A_{1}\text{ and }A_{2}\text{)}}}\right)  \left[
0:i\mid i+k:p+k\right] \\
&  \ \ \ \ \ \ \ \ \ \ \left(
\begin{array}
[c]{c}%
\text{by (\ref{pf.Grasp.generic.independency.d.interval}), applied to
}\mathbb{F}\left[  x_{\mathbf{P}}\right]  \text{, }\mathbb{M}\text{,
}\mathfrak{d}\text{, }p\text{, }p+q\text{,}\\
0\text{, }i\text{, }i+k\text{, }p+k\text{ and }\left(  I_{p}\mid C\right)
\text{ instead of }\mathbb{A}\text{, }\mathbb{B}\text{, }\varphi\text{,
}u\text{, }v\text{, }a\text{, }b\text{, }c\text{, }d\text{ and }A
\end{array}
\right) \\
&  =\left(  \underbrace{\mathfrak{d}^{p\times p}\left(  I_{p}\right)
}_{\substack{=I_{p}\\\text{(by (\ref{pf.Grasp.generic.independency.d.Iu}),
applied to}\\\mathbb{F}\left[  x_{\mathbf{P}}\right]  \text{, }\mathbb{M}%
\text{, }\mathfrak{d}\text{, and }p\text{ instead of}\\\mathbb{A}\text{,
}\mathbb{B}\text{, }\varphi\text{, and }u\text{)}}}\mid
\underbrace{\mathfrak{d}^{p\times q}\left(  C\right)  }%
_{\substack{=D\\\text{(by (\ref{pf.Grasp.generic.independency.d.dC}))}%
}}\right)  \left[  0:i\mid i+k:p+k\right] \\
&  =\left(  I_{p}\mid D\right)  \left[  0:i\mid i+k:p+k\right]  ,
\end{align*}
this becomes%
\begin{align*}
\mathfrak{D}_{\left(  i,k\right)  }\left(  \left(  d_{\mathbf{r}}\right)
_{\mathbf{r}\in\mathbf{P}}\right)   &  =\det\left(  \underbrace{\mathfrak{d}%
^{p\times\left(  i-0+\left(  p+k\right)  -\left(  i+k\right)  \right)
}\left(  \left(  I_{p}\mid C\right)  \left[  0:i\mid i+k:p+k\right]  \right)
}_{=\left(  I_{p}\mid D\right)  \left[  0:i\mid i+k:p+k\right]  }\right) \\
&  =\det\left(  \left(  I_{p}\mid D\right)  \left[  0:i\mid i+k:p+k\right]
\right)  .
\end{align*}
This proves (\ref{lem.Grasp.generic.independency.cD}). This proves Lemma
\ref{lem.Grasp.generic.independency} \textbf{(c)}.

\textbf{(a)} Let $\mathbf{p}\in\mathbf{P}$. Since $\mathbf{p}\in
\mathbf{P}=\left\{  1,2,...,p\right\}  \times\left\{  1,2,...,q\right\}  $, we
can write $\mathbf{p}$ in the form $\mathbf{p}=\left(  i,k\right)  $ for some
$i\in\left\{  1,2,...,p\right\}  $ and $k\in\left\{  1,2,...,q\right\}  $.
Consider these$\ i$ and $k$.

Let us check that $\left(  \operatorname*{Grasp}\nolimits_{0}\left(  I_{p}\mid
C\right)  \right)  \left(  \left(  i,k\right)  \right)  $ is well-defined. In
fact, by the definition of $\left(  \operatorname*{Grasp}\nolimits_{0}\left(
I_{p}\mid C\right)  \right)  \left(  \left(  i,k\right)  \right)  $, we have%
\[
\left(  \operatorname*{Grasp}\nolimits_{0}\left(  I_{p}\mid C\right)  \right)
\left(  \left(  i,k\right)  \right)  =\dfrac{\det\left(  \left(  I_{p}\mid
C\right)  \left[  0+1:0+i\mid0+i+k-1:0+p+k\right]  \right)  }{\det\left(
\left(  I_{p}\mid C\right)  \left[  0+0:0+i\mid0+i+k:0+p+k\right]  \right)
}.
\]
Hence, the element $\left(  \operatorname*{Grasp}\nolimits_{0}\left(
I_{p}\mid C\right)  \right)  \left(  \left(  i,k\right)  \right)  $ is
well-defined if and only if \newline$\det\left(  \left(  I_{p}\mid C\right)
\left[  0+0:0+i\mid0+i+k:0+p+k\right]  \right)  \neq0$. Since
\begin{align}
&  \det\left(  \left(  I_{p}\mid C\right)  \left[  0+0:0+i\mid
0+i+k:0+p+k\right]  \right) \nonumber\\
&  =\det\left(  \left(  I_{p}\mid C\right)  \left[  0:i\mid i+k:p+k\right]
\right) \nonumber\\
&  =\mathfrak{D}_{\left(  i,k\right)  }\ \ \ \ \ \ \ \ \ \ \left(  \text{by
(\ref{lem.Grasp.generic.independency.Ddef})}\right) \nonumber\\
&  =\mathfrak{D}_{\mathbf{p}}\ \ \ \ \ \ \ \ \ \ \left(  \text{since }\left(
i,k\right)  =\mathbf{p}\right) \label{pf.Grasp.generic.independency.b.D}\\
&  \neq0\ \ \ \ \ \ \ \ \ \ \left(  \text{by (\ref{lem.Grasp.generic.calc.f.D}%
)}\right)  ,\nonumber
\end{align}
we can therefore conclude that the element $\left(  \operatorname*{Grasp}%
\nolimits_{0}\left(  I_{p}\mid C\right)  \right)  \left(  \left(  i,k\right)
\right)  $ is well-defined. Since $\left(  i,k\right)  =\mathbf{p}$, this
rewrites as follows: The element $\left(  \operatorname*{Grasp}\nolimits_{0}%
\left(  I_{p}\mid C\right)  \right)  \left(  \mathbf{p}\right)  $ is well-defined.

Moreover,
\begin{align}
&  \det\left(  \left(  I_{p}\mid C\right)  \left[  0+1:0+i\mid
0+i+k-1:0+p+k\right]  \right) \nonumber\\
&  =\det\left(  \left(  I_{p}\mid C\right)  \left[  1:i\mid i+k-1:p+k\right]
\right) \nonumber\\
&  =\mathfrak{N}_{\left(  i,k\right)  }\ \ \ \ \ \ \ \ \ \ \left(  \text{by
(\ref{lem.Grasp.generic.independency.Ndef})}\right) \nonumber\\
&  =\mathfrak{N}_{\mathbf{p}}\ \ \ \ \ \ \ \ \ \ \left(  \text{since }\left(
i,k\right)  =\mathbf{p}\right)  . \label{pf.Grasp.generic.independency.b.N}%
\end{align}

Now, $\mathbf{p}=\left(  i,k\right)  $, so that%
\begin{align*}
\left(  \operatorname*{Grasp}\nolimits_{0}\left(  I_{p}\mid C\right)  \right)
\left(  \mathbf{p}\right)   &  =\left(  \operatorname*{Grasp}\nolimits_{0}%
\left(  I_{p}\mid C\right)  \right)  \left(  \left(  i,k\right)  \right) \\
&  =\dfrac{\det\left(  \left(  I_{p}\mid C\right)  \left[  0+1:0+i\mid
0+i+k-1:0+p+k\right]  \right)  }{\det\left(  \left(  I_{p}\mid C\right)
\left[  0+0:0+i\mid0+i+k:0+p+k\right]  \right)  }\\
&  \ \ \ \ \ \ \ \ \ \ \left(  \text{by the definition of }\left(
\operatorname*{Grasp}\nolimits_{0}\left(  I_{p}\mid C\right)  \right)  \left(
\left(  i,k\right)  \right)  \right) \\
&  =\dfrac{\mathfrak{N}_{\mathbf{p}}}{\mathfrak{D}_{\mathbf{p}}}%
\ \ \ \ \ \ \ \ \ \ \left(  \text{by (\ref{pf.Grasp.generic.independency.b.N})
and (\ref{pf.Grasp.generic.independency.b.D})}\right)  .
\end{align*}
This completes the proof of Lemma \ref{lem.Grasp.generic.independency}
\textbf{(a)}.

\textbf{(b)} For every element $\mathbf{p}\in\mathbf{P}$, the element
$Q_{\mathbf{p}}$ of $\mathbb{F}\left(  x_{\mathbf{P}}\right)  $ is
well-defined\footnote{\textit{Proof.} Let $\mathbf{p}\in\mathbf{P}$. Then,
$Q_{\mathbf{p}}=\left(  \operatorname*{Grasp}\nolimits_{0}\left(  I_{p}\mid
C\right)  \right)  \left(  \mathbf{p}\right)  $ is well-defined (according to
Lemma \ref{lem.Grasp.generic.independency} \textbf{(a)}), qed.}. We are now
going to prove that the family $\left(  Q_{\mathbf{p}}\right)  _{\mathbf{p}%
\in\mathbf{P}}$ satisfies the algebraic triangularity condition\footnote{This
condition was described in Definition \ref{def.algebraic.triangularity.triang}%
.}.

Let $\mathbf{p}\in\mathbf{P}$. Then, Lemma \ref{lem.Grasp.generic.calc}
\textbf{(g)} says that there exist elements $\alpha_{\mathbf{p}}$,
$\beta_{\mathbf{p}}$, $\gamma_{\mathbf{p}}$ and $\delta_{\mathbf{p}}$ of
$\mathbb{F}\left[  x_{\mathbf{p}\Downarrow}\right]  $ satisfying
$\alpha_{\mathbf{p}}\delta_{\mathbf{p}}-\beta_{\mathbf{p}}\gamma_{\mathbf{p}%
}\neq0$, $\mathfrak{N}_{\mathbf{p}}=\alpha_{\mathbf{p}}x_{\mathbf{p}}%
+\beta_{\mathbf{p}}$ and $\mathfrak{D}_{\mathbf{p}}=\gamma_{\mathbf{p}%
}x_{\mathbf{p}}+\delta_{\mathbf{p}}$. Consider these elements $\alpha
_{\mathbf{p}}$, $\beta_{\mathbf{p}}$, $\gamma_{\mathbf{p}}$ and $\delta
_{\mathbf{p}}$. Clearly, the elements $\alpha_{\mathbf{p}}$, $\beta
_{\mathbf{p}}$, $\gamma_{\mathbf{p}}$ and $\delta_{\mathbf{p}}$ all belong to
$\mathbb{F}\left(  x_{\mathbf{p}\Downarrow}\right)  $ (because they belong to
$\mathbb{F}\left[  x_{\mathbf{p}\Downarrow}\right]  $, and because
$\mathbb{F}\left[  x_{\mathbf{p}\Downarrow}\right]  \subseteq\mathbb{F}\left(
x_{\mathbf{p}\Downarrow}\right)  $). Also,
\begin{align*}
Q_{\mathbf{p}}  &  =\left(  \operatorname*{Grasp}\nolimits_{0}\left(
I_{p}\mid C\right)  \right)  \left(  \mathbf{p}\right)  =\dfrac{\mathfrak{N}%
_{\mathbf{p}}}{\mathfrak{D}_{\mathbf{p}}}\ \ \ \ \ \ \ \ \ \ \left(  \text{by
(\ref{lem.Grasp.generic.independency.a})}\right) \\
&  =\dfrac{\alpha_{\mathbf{p}}x_{\mathbf{p}}+\beta_{\mathbf{p}}}%
{\gamma_{\mathbf{p}}x_{\mathbf{p}}+\delta_{\mathbf{p}}}%
\ \ \ \ \ \ \ \ \ \ \left(  \text{since }\mathfrak{N}_{\mathbf{p}}%
=\alpha_{\mathbf{p}}x_{\mathbf{p}}+\beta_{\mathbf{p}}\text{ and }%
\mathfrak{D}_{\mathbf{p}}=\gamma_{\mathbf{p}}x_{\mathbf{p}}+\delta
_{\mathbf{p}}\right)  .
\end{align*}
Forget now that we have chosen $\alpha_{\mathbf{p}}$, $\beta_{\mathbf{p}}$,
$\gamma_{\mathbf{p}}$ and $\delta_{\mathbf{p}}$. We have thus shown that there
exist elements $\alpha_{\mathbf{p}}$, $\beta_{\mathbf{p}}$, $\gamma
_{\mathbf{p}}$, $\delta_{\mathbf{p}}$ of $\mathbb{F}\left(  x_{\mathbf{p}%
\Downarrow}\right)  $ such that $\alpha_{\mathbf{p}}\delta_{\mathbf{p}}%
-\beta_{\mathbf{p}}\gamma_{\mathbf{p}}\neq0$ and $Q_{\mathbf{p}}=\dfrac
{\alpha_{\mathbf{p}}x_{\mathbf{p}}+\beta_{\mathbf{p}}}{\gamma_{\mathbf{p}%
}x_{\mathbf{p}}+\delta_{\mathbf{p}}}$.

Now, forget that we fixed $\mathbf{p}\in\mathbf{P}$. Thus, we have shown that
for every $\mathbf{p}\in\mathbf{P}$, there exist elements $\alpha_{\mathbf{p}%
}$, $\beta_{\mathbf{p}}$, $\gamma_{\mathbf{p}}$, $\delta_{\mathbf{p}}$ of
$\mathbb{F}\left(  x_{\mathbf{p}\Downarrow}\right)  $ such that $\alpha
_{\mathbf{p}}\delta_{\mathbf{p}}-\beta_{\mathbf{p}}\gamma_{\mathbf{p}}\neq0$
and $Q_{\mathbf{p}}=\dfrac{\alpha_{\mathbf{p}}x_{\mathbf{p}}+\beta
_{\mathbf{p}}}{\gamma_{\mathbf{p}}x_{\mathbf{p}}+\delta_{\mathbf{p}}}$. In
other words, the algebraic triangularity condition holds. In other words, the
family $\left(  Q_{\mathbf{p}}\right)  _{\mathbf{p}\in\mathbf{P}}\in\left(
\mathbb{F}\left(  x_{\mathbf{P}}\right)  \right)  ^{\mathbf{P}}$ is
$\mathbf{P}$-triangular (by the definition of ``$\mathbf{P}$-triangular'').
This proves Lemma \ref{lem.Grasp.generic.independency} \textbf{(b)}.

\textbf{(d)} Let $\mathbf{p}\in\mathbf{P}$. Since $\mathbf{p}\in
\mathbf{P}=\left\{  1,2,...,p\right\}  \times\left\{  1,2,...,q\right\}  $, we
can write $\mathbf{p}$ in the form $\mathbf{p}=\left(  i,k\right)  $ for some
$i\in\left\{  1,2,...,p\right\}  $ and $k\in\left\{  1,2,...,q\right\}  $.
Consider these$\ i$ and $k$.

The values $\mathfrak{N}_{\mathbf{p}}\left(  \left(  d_{\mathbf{r}}\right)
_{\mathbf{r}\in\mathbf{P}}\right)  $ and $\mathfrak{D}_{\mathbf{p}}\left(
\left(  d_{\mathbf{r}}\right)  _{\mathbf{r}\in\mathbf{P}}\right)  $ are
well-defined (since $\mathfrak{N}_{\mathbf{p}}$ and $\mathfrak{D}_{\mathbf{p}%
}$ are polynomials and thus have no denominators that could become zero upon
substitution of variables), and we have $\mathfrak{D}_{\mathbf{p}}\left(
\left(  d_{\mathbf{r}}\right)  _{\mathbf{r}\in\mathbf{P}}\right)  \neq0$.
Hence, the value $\dfrac{\mathfrak{N}_{\mathbf{p}}}{\mathfrak{D}_{\mathbf{p}}%
}\left(  \left(  d_{\mathbf{r}}\right)  _{\mathbf{r}\in\mathbf{P}}\right)  $
is well-defined. Since $Q_{\mathbf{p}}=\left(  \operatorname*{Grasp}%
\nolimits_{0}\left(  I_{p}\mid C\right)  \right)  \left(  \mathbf{p}\right)
=\dfrac{\mathfrak{N}_{\mathbf{p}}}{\mathfrak{D}_{\mathbf{p}}}$ (by
(\ref{lem.Grasp.generic.independency.a})), this rewrites as follows: The value
$Q_{\mathbf{p}}\left(  \left(  d_{\mathbf{r}}\right)  _{\mathbf{r}%
\in\mathbf{P}}\right)  $ is well-defined.

By the definition of $\left(  \operatorname*{Grasp}\nolimits_{0}\left(
I_{p}\mid D\right)  \right)  \left(  \left(  i,k\right)  \right)  $, we have%
\[
\left(  \operatorname*{Grasp}\nolimits_{0}\left(  I_{p}\mid D\right)  \right)
\left(  \left(  i,k\right)  \right)  =\dfrac{\det\left(  \left(  I_{p}\mid
D\right)  \left[  0+1:0+i\mid0+i+k-1:0+p+k\right]  \right)  }{\det\left(
\left(  I_{p}\mid D\right)  \left[  0+0:0+i\mid0+i+k:0+p+k\right]  \right)
}.
\]
Hence, the element $\left(  \operatorname*{Grasp}\nolimits_{0}\left(
I_{p}\mid D\right)  \right)  \left(  \left(  i,k\right)  \right)  $ is
well-defined if and only if \newline$\det\left(  \left(  I_{p}\mid D\right)
\left[  0+0:0+i\mid0+i+k:0+p+k\right]  \right)  \neq0$. Since%
\begin{align}
&  \det\left(  \left(  I_{p}\mid D\right)  \left[  0+0:0+i\mid
0+i+k:0+p+k\right]  \right) \nonumber\\
&  =\det\left(  \left(  I_{p}\mid D\right)  \left[  0:i\mid i+k:p+k\right]
\right) \nonumber\\
&  =\mathfrak{D}_{\left(  i,k\right)  }\left(  \left(  d_{\mathbf{r}}\right)
_{\mathbf{r}\in\mathbf{P}}\right)  \ \ \ \ \ \ \ \ \ \ \left(  \text{by
(\ref{lem.Grasp.generic.independency.cD})}\right) \nonumber\\
&  =\mathfrak{D}_{\mathbf{p}}\left(  \left(  d_{\mathbf{r}}\right)
_{\mathbf{r}\in\mathbf{P}}\right)  \ \ \ \ \ \ \ \ \ \ \left(  \text{since
}\left(  i,k\right)  =\mathbf{p}\right)
\label{pf.Grasp.generic.independency.e.D}\\
&  \neq0,\nonumber
\end{align}
we can therefore conclude that the element $\left(  \operatorname*{Grasp}%
\nolimits_{0}\left(  I_{p}\mid D\right)  \right)  \left(  \left(  i,k\right)
\right)  $ is well-defined. Since $\left(  i,k\right)  =\mathbf{p}$, this
rewrites as follows: The element $\left(  \operatorname*{Grasp}\nolimits_{0}%
\left(  I_{p}\mid D\right)  \right)  \left(  \mathbf{p}\right)  $ is
well-defined. Moreover,
\begin{align}
&  \det\left(  \left(  I_{p}\mid D\right)  \left[  0+1:0+i\mid
0+i+k-1:0+p+k\right]  \right) \nonumber\\
&  =\det\left(  \left(  I_{p}\mid D\right)  \left[  1:i\mid i+k-1:p+k\right]
\right) \nonumber\\
&  =\mathfrak{N}_{\left(  i,k\right)  }\left(  \left(  d_{\mathbf{r}}\right)
_{\mathbf{r}\in\mathbf{P}}\right)  \ \ \ \ \ \ \ \ \ \ \left(  \text{by
(\ref{lem.Grasp.generic.independency.cN})}\right) \nonumber\\
&  =\mathfrak{N}_{\mathbf{p}}\left(  \left(  d_{\mathbf{r}}\right)
_{\mathbf{r}\in\mathbf{P}}\right)  \ \ \ \ \ \ \ \ \ \ \left(  \text{since
}\left(  i,k\right)  =\mathbf{p}\right)  .
\label{pf.Grasp.generic.independency.e.N}%
\end{align}

Now, we have $Q_{\mathbf{p}}=\dfrac{\mathfrak{N}_{\mathbf{p}}}{\mathfrak{D}%
_{\mathbf{p}}}$, so that%
\[
Q_{\mathbf{p}}\left(  \left(  d_{\mathbf{r}}\right)  _{\mathbf{r}\in
\mathbf{P}}\right)  =\dfrac{\mathfrak{N}_{\mathbf{p}}}{\mathfrak{D}%
_{\mathbf{p}}}\left(  \left(  d_{\mathbf{r}}\right)  _{\mathbf{r}\in
\mathbf{P}}\right)  =\dfrac{\mathfrak{N}_{\mathbf{p}}\left(  \left(
d_{\mathbf{r}}\right)  _{\mathbf{r}\in\mathbf{P}}\right)  }{\mathfrak{D}%
_{\mathbf{p}}\left(  \left(  d_{\mathbf{r}}\right)  _{\mathbf{r}\in\mathbf{P}%
}\right)  }.
\]
Compared with%
\begin{align*}
\left(  \operatorname*{Grasp}\nolimits_{0}\left(  I_{p}\mid D\right)  \right)
\left(  \mathbf{p}\right)   &  =\left(  \operatorname*{Grasp}\nolimits_{0}%
\left(  I_{p}\mid D\right)  \right)  \left(  \left(  i,k\right)  \right)
\ \ \ \ \ \ \ \ \ \ \left(  \text{since }\mathbf{p}=\left(  i,k\right)
\right) \\
&  =\dfrac{\det\left(  \left(  I_{p}\mid D\right)  \left[  0+1:0+i\mid
0+i+k-1:0+p+k\right]  \right)  }{\det\left(  \left(  I_{p}\mid D\right)
\left[  0+0:0+i\mid0+i+k:0+p+k\right]  \right)  }\\
&  =\dfrac{\mathfrak{N}_{\mathbf{p}}\left(  \left(  d_{\mathbf{r}}\right)
_{\mathbf{r}\in\mathbf{P}}\right)  }{\mathfrak{D}_{\mathbf{p}}\left(  \left(
d_{\mathbf{r}}\right)  _{\mathbf{r}\in\mathbf{P}}\right)  }%
\ \ \ \ \ \ \ \ \ \ \left(  \text{by (\ref{pf.Grasp.generic.independency.e.N})
and (\ref{pf.Grasp.generic.independency.e.D})}\right)  ,
\end{align*}
this yields $Q_{\mathbf{p}}\left(  \left(  d_{\mathbf{r}}\right)
_{\mathbf{r}\in\mathbf{P}}\right)  =\left(  \operatorname*{Grasp}%
\nolimits_{0}\left(  I_{p}\mid D\right)  \right)  \left(  \mathbf{p}\right)
$. This completes the proof of Lemma \ref{lem.Grasp.generic.independency}
\textbf{(d)}.

Now, the proof of Lemma \ref{lem.Grasp.generic.independency} is finished at last.
\end{proof}

As explained above, Proposition \ref{prop.Grasp.generic} follows from Lemma
\ref{lem.Grasp.generic.independency}. Hence, Proposition
\ref{prop.Grasp.generic} is proven now.
\end{verlong}

\section{\label{sect.rect.finish}The rectangle: finishing the proofs}

As promised, we now use Propositions \ref{prop.Grasp.GraspR} and
\ref{prop.Grasp.generic} to derive our initially stated results on rectangles.
First, we formulate an easy consequence of Proposition \ref{prop.Grasp.GraspR}:

\begin{corollary}
\label{cor.Grasp.GraspR}Let $\mathbb{K}$ be a field. Let $p$ and $q$ be two
positive integers. Let $A\in\mathbb{K}^{p\times\left(  p+q\right)  }$ be a
matrix. Then, every $i \in\mathbb{N}$ satisfies%
\[
\operatorname*{Grasp}\nolimits_{-i}A=R_{\operatorname*{Rect}\left(
p,q\right)  }^{i}\left(  \operatorname*{Grasp}\nolimits_{0}A\right)
\]
(provided that $A$ is sufficiently generic in the sense of Zariski topology
that both sides of this equality are well-defined).
\end{corollary}

\begin{vershort}
\begin{proof}
[Proof of Corollary \ref{cor.Grasp.GraspR} (sketched).]We will prove Corollary
\ref{cor.Grasp.GraspR} by induction over $i$:

\textit{Induction base:} For $i = 0$, the claim of Corollary
\ref{cor.Grasp.GraspR} boils down to $\operatorname*{Grasp}\nolimits_{0}
A=R_{\operatorname*{Rect}\left(  p,q\right)  }^{0}\left(
\operatorname*{Grasp}\nolimits_{0}A\right)  $. This is obvious, and so the
induction base is complete.

\textit{Induction step:} Let $j\in\mathbb{N}$. Assume that Corollary
\ref{cor.Grasp.GraspR} holds for $i=j$. We need to prove that Corollary
\ref{cor.Grasp.GraspR} holds for $i=j+1$ as well.

Proposition \ref{prop.Grasp.GraspR} (applied to $-\left(  j+1\right)  $
instead of $j$) yields
\begin{align*}
\operatorname*{Grasp}\nolimits_{-\left(  j+1\right)  }A  &
=R_{\operatorname*{Rect}\left(  p,q\right)  }\left(  \operatorname*{Grasp}%
\nolimits_{-\left(  j+1\right)  +1}A\right)  =R_{\operatorname*{Rect}\left(
p,q\right)  }\left(  \underbrace{\operatorname*{Grasp}\nolimits_{-j}%
A}_{\substack{=R_{\operatorname*{Rect}\left(  p,q\right)  }^{j}\left(
\operatorname*{Grasp}\nolimits_{0}A\right)  \\\text{(since Corollary
\ref{cor.Grasp.GraspR}}\\\text{holds for }i=j\text{)}}}\right) \\
&  =R_{\operatorname*{Rect}\left(  p,q\right)  }\left(
R_{\operatorname*{Rect}\left(  p,q\right)  }^{j}\left(  \operatorname*{Grasp}%
\nolimits_{0}A\right)  \right)  =R_{\operatorname*{Rect}\left(  p,q\right)
}^{j+1}\left(  \operatorname*{Grasp}\nolimits_{0}A\right)  .
\end{align*}
In other words, Corollary \ref{cor.Grasp.GraspR} holds for $i=j+1$. This
completes the induction step. The induction proof of Corollary
\ref{cor.Grasp.GraspR} is thus finished.
\end{proof}
\end{vershort}

\begin{verlong}
\begin{proof}
[Proof of Corollary \ref{cor.Grasp.GraspR} (sketched).]We will prove Corollary
\ref{cor.Grasp.GraspR} by induction over $i$:

\textit{Induction base:} We have $R_{\operatorname*{Rect}\left(  p,q\right)
}^{0}=\operatorname*{id}$. Thus, $R_{\operatorname*{Rect}\left(  p,q\right)
}^{0}\left(  \operatorname*{Grasp}\nolimits_{0}A\right)
=\operatorname*{Grasp}\nolimits_{0}A=\operatorname*{Grasp}\nolimits_{-0}A$. In
other words, $\operatorname*{Grasp}\nolimits_{-0}A=R_{\operatorname*{Rect}%
\left(  p,q\right)  }^{0}\left(  \operatorname*{Grasp}\nolimits_{0}A\right)
$. In other words, Corollary \ref{cor.Grasp.GraspR} holds for $i=0$. This
completes the induction base.

\textit{Induction step:} Let $j\in\mathbb{N}$. Assume that Corollary
\ref{cor.Grasp.GraspR} holds for $i=j$. We need to prove that Corollary
\ref{cor.Grasp.GraspR} holds for $i=j+1$ as well.

Proposition \ref{prop.Grasp.GraspR} (applied to $-\left(  j+1\right)  $
instead of $j$) yields
\begin{align*}
\operatorname*{Grasp}\nolimits_{-\left(  j+1\right)  }A  &
=R_{\operatorname*{Rect}\left(  p,q\right)  }\left(  \operatorname*{Grasp}%
\nolimits_{-\left(  j+1\right)  +1}A\right)  =R_{\operatorname*{Rect}\left(
p,q\right)  }\left(  \underbrace{\operatorname*{Grasp}\nolimits_{-j}%
A}_{\substack{=R_{\operatorname*{Rect}\left(  p,q\right)  }^{j}\left(
\operatorname*{Grasp}\nolimits_{0}A\right)  \\\text{(since Corollary
\ref{cor.Grasp.GraspR}}\\\text{holds for }i=j\text{)}}}\right) \\
&  \ \ \ \ \ \ \ \ \ \ \left(  \text{since }-\left(  j+1\right)  +1=-j\right)
\\
&  =R_{\operatorname*{Rect}\left(  p,q\right)  }\left(
R_{\operatorname*{Rect}\left(  p,q\right)  }^{j}\left(  \operatorname*{Grasp}%
\nolimits_{0}A\right)  \right)  =\underbrace{\left(  R_{\operatorname*{Rect}%
\left(  p,q\right)  }\circ R_{\operatorname*{Rect}\left(  p,q\right)  }%
^{j}\right)  }_{=R_{\operatorname*{Rect}\left(  p,q\right)  }^{1+j}%
=R_{\operatorname*{Rect}\left(  p,q\right)  }^{j+1}}\left(
\operatorname*{Grasp}\nolimits_{0}A\right) \\
&  =R_{\operatorname*{Rect}\left(  p,q\right)  }^{j+1}\left(
\operatorname*{Grasp}\nolimits_{0}A\right)  .
\end{align*}
In other words, Corollary \ref{cor.Grasp.GraspR} holds for $i=j+1$. This
completes the induction step. The induction proof of Corollary
\ref{cor.Grasp.GraspR} is thus finished.
\end{proof}
\end{verlong}

\begin{proof}
[Proof of Theorem \ref{thm.rect.ord} (sketched).]We need to show that
$\operatorname*{ord}\left(  R_{\operatorname*{Rect}\left(  p,q\right)
}\right)  =p+q$. According to Proposition \ref{prop.rect.reduce}, it is enough
to prove that almost every (in the Zariski sense) reduced labelling
$f\in\mathbb{K}^{\widehat{\operatorname*{Rect}\left(  p,q\right)  }}$
satisfies $R_{\operatorname*{Rect}\left(  p,q\right)  }^{p+q}f=f$. So let
$f\in\mathbb{K}^{\widehat{\operatorname*{Rect}\left(  p,q\right)  }}$ be a
sufficiently generic reduced labelling. In other words, $f$ is a sufficiently
generic element of $\mathbb{K}^{\operatorname*{Rect}\left(  p,q\right)  }$
(because the reduced labellings $\mathbb{K}^{\widehat{\operatorname*{Rect}%
\left(  p,q\right)  }}$ are being identified with the elements of
$\mathbb{K}^{\operatorname*{Rect}\left(  p,q\right)  }$). Due to Proposition
\ref{prop.Grasp.generic}, there exists a matrix $A\in\mathbb{K}^{p\times
\left(  p+q\right)  }$ satisfying $f=\operatorname*{Grasp}\nolimits_{0}A$.
Consider this $A$. Due to Corollary \ref{cor.Grasp.GraspR} (applied to
$i=p+q$), we have%
\[
\operatorname*{Grasp}\nolimits_{-\left(  p+q\right)  }%
A=R_{\operatorname*{Rect}\left(  p,q\right)  }^{p+q}\left(
\underbrace{\operatorname*{Grasp}\nolimits_{0}A}_{=f}\right)
=R_{\operatorname*{Rect}\left(  p,q\right)  }^{p+q}f.
\]
But Proposition \ref{prop.Grasp.period} (applied to $j=-\left(  p+q\right)  $)
yields
\begin{align*}
\operatorname*{Grasp}\nolimits_{-\left(  p+q\right)  }A  &
=\operatorname*{Grasp}\nolimits_{p+q+\left(  -\left(  p+q\right)  \right)
}A=\operatorname*{Grasp}\nolimits_{0}A\ \ \ \ \ \ \ \ \ \ \left(  \text{since
}p+q+\left(  -\left(  p+q\right)  \right)  =0\right) \\
&  =f.
\end{align*}
Hence, $f=\operatorname*{Grasp}\nolimits_{-\left(  p+q\right)  }%
A=R_{\operatorname*{Rect}\left(  p,q\right)  }^{p+q}f$. In other words,
$R_{\operatorname*{Rect}\left(  p,q\right)  }^{p+q}f=f$. This (as we know)
proves Theorem \ref{thm.rect.ord}.
\end{proof}

\begin{proof}
[Proof of Theorem \ref{thm.rect.antip} (sketched).]Let us regard the reduced
labelling $f\in\mathbb{K}^{\widehat{\operatorname*{Rect}\left(  p,q\right)  }%
}$ as an element of $\mathbb{K}^{\operatorname*{Rect}\left(  p,q\right)  }$
(because we identify reduced labellings in $\mathbb{K}%
^{\widehat{\operatorname*{Rect}\left(  p,q\right)  }}$ with elements of
$\mathbb{K}^{\operatorname*{Rect}\left(  p,q\right)  }$). We assume WLOG that
this element $f\in\mathbb{K}^{\operatorname*{Rect}\left(  p,q\right)  }$ is
generic enough (among the reduced labellings) for Proposition
\ref{prop.Grasp.generic} to apply. By Proposition \ref{prop.Grasp.generic},
there exists a matrix $A\in\mathbb{K}^{p\times\left(  p+q\right)  }$
satisfying $f=\operatorname*{Grasp}\nolimits_{0}A$. Consider this $A$. Due to
Corollary \ref{cor.Grasp.GraspR} (applied to $i+k-1$ instead of $i$), we have%
\[
\operatorname*{Grasp}\nolimits_{-\left(  i+k-1\right)  }%
A=R_{\operatorname*{Rect}\left(  p,q\right)  }^{i+k-1}\left(
\underbrace{\operatorname*{Grasp}\nolimits_{0}A}_{=f}\right)
=R_{\operatorname*{Rect}\left(  p,q\right)  }^{i+k-1}f.
\]

But Proposition \ref{prop.Grasp.antipode} (applied to $j=-\left(
i+k-1\right)  $) yields%
\begin{align*}
\left(  \operatorname*{Grasp}\nolimits_{-\left(  i+k-1\right)  }A\right)
\left(  \left(  i,k\right)  \right)   &  =\dfrac{1}{\left(
\operatorname*{Grasp}\nolimits_{-\left(  i+k-1\right)  +i+k-1}A\right)
\left(  \left(  p+1-i,q+1-k\right)  \right)  }\\
&  =\dfrac{1}{f\left(  \left(  p+1-i,q+1-k\right)  \right)  }\\
&  \ \ \ \ \ \ \ \ \ \ \left(  \text{since }\operatorname*{Grasp}%
\nolimits_{-\left(  i+k-1\right)  +i+k-1}A=\operatorname*{Grasp}%
\nolimits_{0}A=f\right)  ,
\end{align*}
so that%
\[
f\left(  \left(  p+1-i,q+1-k\right)  \right)  =\dfrac{1}{\left(
\operatorname*{Grasp}\nolimits_{-\left(  i+k-1\right)  }A\right)  \left(
\left(  i,k\right)  \right)  }=\dfrac{1}{\left(  R_{\operatorname*{Rect}%
\left(  p,q\right)  }^{i+k-1}f\right)  \left(  \left(  i,k\right)  \right)  }%
\]
(since $\operatorname*{Grasp}\nolimits_{-\left(  i+k-1\right)  }%
A=R_{\operatorname*{Rect}\left(  p,q\right)  }^{i+k-1}f$). This proves Theorem
\ref{thm.rect.antip}.
\end{proof}

\begin{proof}
[Proof of Theorem \ref{thm.rect.antip.general} (sketched).]We will be using
the notation $\left(  a_{0},a_{1},...,a_{n+1}\right)  \flat f$ defined in
Definition \ref{def.bemol}.

Let $f\in\mathbb{K}^{\widehat{\operatorname*{Rect}\left(  p,q\right)  }}$ be
arbitrary. By genericity, we assume WLOG that $f\left(  0\right)  $ and
$f\left(  1\right)  $ are nonzero.

Let $n=p+q-1$. Then, $\operatorname*{Rect}\left(  p,q\right)  $ is an
$n$-graded poset. Also, $i+k-1\in\left\{  0,1,...,n\right\}  $. Moreover,
$1\leq n-i-k+2\leq n$.

Define an $\left(  n+2\right)  $-tuple $\left(  a_{0},a_{1},...,a_{n+1}%
\right)  \in\mathbb{K}^{n+2}$ by%
\[
a_{r}=\left\{
\begin{array}
[c]{c}%
\dfrac{1}{f\left(  0\right)  },\ \ \ \ \ \ \ \ \ \ \text{if }r=0;\\
1,\ \ \ \ \ \ \ \ \ \ \text{if }1\leq r\leq n;\\
\dfrac{1}{f\left(  1\right)  },\ \ \ \ \ \ \ \ \ \ \text{if }r=n+1
\end{array}
\right.  \ \ \ \ \ \ \ \ \ \ \text{for every }r\in\left\{
0,1,...,n+1\right\}  .
\]
Thus, $a_{n-i-k+2}=1$ (since $1\leq n-i-k+2\leq n$) and $a_{0}=\dfrac
{1}{f\left(  0\right)  }$ and $a_{n+1}=\dfrac{1}{f\left(  1\right)  }$.

Let $f^{\prime}=\left(  a_{0},a_{1},...,a_{n+1}\right)  \flat f$. Then, it is
easy to see from the definition of $\left(  a_{0},a_{1},...,a_{n+1}\right)
\flat f$ that $f^{\prime}\left(  0\right)  =1$ and $f^{\prime}\left(
1\right)  =1$. In other words, $f^{\prime}$ is a reduced $\mathbb{K}%
$-labelling. Hence, Theorem \ref{thm.rect.antip} (applied to $f^{\prime}$
instead of $f$) yields%
\begin{equation}
f^{\prime}\left(  \left(  p+1-i,q+1-k\right)  \right)  =\dfrac{1}{\left(
R_{\operatorname*{Rect}\left(  p,q\right)  }^{i+k-1}\left(  f^{\prime}\right)
\right)  \left(  \left(  i,k\right)  \right)  }.
\label{pf.rect.antip.general.1}%
\end{equation}

On the other hand, again from the definition of $f^{\prime}=\left(
a_{0},a_{1},...,a_{n+1}\right)  \flat f$, it is easy to see that $f^{\prime
}\left(  v\right)  =f\left(  v\right)  $ for every $v\in\operatorname*{Rect}%
\left(  p,q\right)  $. This yields, in particular, that $f^{\prime}\left(
\left(  p+1-i,q+1-k\right)  \right)  =f\left(  \left(  p+1-i,q+1-k\right)
\right)  $.

But let us define an element $\widehat{a}_{\kappa}^{\left(  \ell\right)  }%
\in\mathbb{K}^{\times}$ for every $\ell\in\left\{  0,1,...,n+1\right\}  $ and
$\kappa\in\left\{  0,1,...,n+1\right\}  $ as in Proposition
\ref{prop.Rl.scalmult}. Then, Proposition \ref{prop.Rl.scalmult} (applied to
$\ell=i+k-1$) yields%
\[
R_{\operatorname*{Rect}\left(  p,q\right)  }^{i+k-1}\left(  \left(
a_{0},a_{1},...,a_{n+1}\right)  \flat f\right)  =\left(  \widehat{a}%
_{0}^{\left(  i+k-1\right)  },\widehat{a}_{1}^{\left(  i+k-1\right)
},...,\widehat{a}_{n+1}^{\left(  i+k-1\right)  }\right)  \flat\left(
R_{\operatorname*{Rect}\left(  p,q\right)  }^{i+k-1}f\right)  .
\]
Since $\left(  a_{0},a_{1},...,a_{n+1}\right)  \flat f=f^{\prime}$, this
rewrites as
\[
R_{\operatorname*{Rect}\left(  p,q\right)  }^{i+k-1}\left(  f^{\prime}\right)
=\left(  \widehat{a}_{0}^{\left(  i+k-1\right)  },\widehat{a}_{1}^{\left(
i+k-1\right)  },...,\widehat{a}_{n+1}^{\left(  i+k-1\right)  }\right)
\flat\left(  R_{\operatorname*{Rect}\left(  p,q\right)  }^{i+k-1}f\right)  .
\]
Hence,%
\begin{align*}
&  \left(  R_{\operatorname*{Rect}\left(  p,q\right)  }^{i+k-1}\left(
f^{\prime}\right)  \right)  \left(  \left(  i,k\right)  \right) \\
&  =\left(  \left(  \widehat{a}_{0}^{\left(  i+k-1\right)  },\widehat{a}%
_{1}^{\left(  i+k-1\right)  },...,\widehat{a}_{n+1}^{\left(  i+k-1\right)
}\right)  \flat\left(  R_{\operatorname*{Rect}\left(  p,q\right)  }%
^{i+k-1}f\right)  \right)  \left(  \left(  i,k\right)  \right) \\
&  =\widehat{a}_{\deg\left(  \left(  i,k\right)  \right)  }^{\left(
i+k-1\right)  }\cdot\left(  R_{\operatorname*{Rect}\left(  p,q\right)
}^{i+k-1}f\right)  \left(  \left(  i,k\right)  \right) \\
&  \ \ \ \ \ \ \ \ \ \ \left(  \text{by the definition of }\left(
\widehat{a}_{0}^{\left(  i+k-1\right)  },\widehat{a}_{1}^{\left(
i+k-1\right)  },...,\widehat{a}_{n+1}^{\left(  i+k-1\right)  }\right)
\flat\left(  R_{\operatorname*{Rect}\left(  p,q\right)  }^{i+k-1}f\right)
\right) \\
&  =\widehat{a}_{i+k-1}^{\left(  i+k-1\right)  }\cdot\left(
R_{\operatorname*{Rect}\left(  p,q\right)  }^{i+k-1}f\right)  \left(  \left(
i,k\right)  \right)  \ \ \ \ \ \ \ \ \ \ \left(  \text{since }\deg\left(
\left(  i,k\right)  \right)  =i+k-1\right) \\
&  =\dfrac{1}{f\left(  0\right)  f\left(  1\right)  }\cdot\left(
R_{\operatorname*{Rect}\left(  p,q\right)  }^{i+k-1}f\right)  \left(  \left(
i,k\right)  \right)
\end{align*}
(since the definition of $\widehat{a}_{i+k-1}^{\left(  i+k-1\right)  }$ yields%
\begin{align*}
\widehat{a}_{i+k-1}^{\left(  i+k-1\right)  }  &  =\left\{
\begin{array}
[c]{c}%
\dfrac{a_{n+1}a_{\left(  i+k-1\right)  -\left(  i+k-1\right)  }}%
{a_{n+1-\left(  i+k-1\right)  }},\ \ \ \ \ \ \ \ \ \ \text{if }i+k-1\geq
i+k-1;\\
\dfrac{a_{n+1+\left(  i+k-1\right)  -\left(  i+k-1\right)  }a_{0}%
}{a_{n+1-\left(  i+k-1\right)  }},\ \ \ \ \ \ \ \ \ \ \text{if }i+k-1<i+k-1
\end{array}
\right. \\
&  =\dfrac{a_{n+1}a_{\left(  i+k-1\right)  -\left(  i+k-1\right)  }%
}{a_{n+1-\left(  i+k-1\right)  }}\ \ \ \ \ \ \ \ \ \ \left(  \text{since
}i+k-1\geq i+k-1\right) \\
&  =\dfrac{a_{n+1}a_{0}}{a_{n-i-k+2}}=\underbrace{a_{n+1}}_{=\dfrac
{1}{f\left(  1\right)  }}\underbrace{a_{0}}_{=\dfrac{1}{f\left(  0\right)  }%
}\ \ \ \ \ \ \ \ \ \ \left(  \text{since }a_{n-i-k+2}=1\right) \\
&  =\dfrac{1}{f\left(  0\right)  f\left(  1\right)  }%
\end{align*}
). Thus, (\ref{pf.rect.antip.general.1}) rewrites as%
\[
f^{\prime}\left(  \left(  p+1-i,q+1-k\right)  \right)  =\dfrac{1}{\dfrac
{1}{f\left(  0\right)  f\left(  1\right)  }\cdot\left(
R_{\operatorname*{Rect}\left(  p,q\right)  }^{i+k-1}f\right)  \left(  \left(
i,k\right)  \right)  }=\dfrac{f\left(  0\right)  f\left(  1\right)  }{\left(
R_{\operatorname*{Rect}\left(  p,q\right)  }^{i+k-1}f\right)  \left(  \left(
i,k\right)  \right)  }.
\]
This rewrites as%
\[
f\left(  \left(  p+1-i,q+1-k\right)  \right)  =\dfrac{f\left(  0\right)
f\left(  1\right)  }{\left(  R_{\operatorname*{Rect}\left(  p,q\right)
}^{i+k-1}f\right)  \left(  \left(  i,k\right)  \right)  }
\]
(since we know that $f^{\prime}\left(  \left(  p+1-i,q+1-k\right)  \right)
=f\left(  \left(  p+1-i,q+1-k\right)  \right)  $). This proves Theorem
\ref{thm.rect.antip.general}.
\end{proof}

\section{\texorpdfstring{The $\vartriangleright$ triangle}{The |> triangle}}

\label{sect.tria}

Having proven the main properties of birational rowmotion $R$ on the rectangle
$\operatorname*{Rect}\left(  p,q\right)  $ and on skeletal posets, we now turn
to other posets. We will spend the next three sections discussing the order of
birational rowmotion on certain triangle-shaped posets obtained as subsets of
the square $\operatorname*{Rect}\left(  p,p\right)  $. We start with the
easiest case:

\begin{definition}
\label{def.Leftri}Let $p$ be a positive integer. Define a subset
$\operatorname*{Tria}\left(  p\right)  $ of $\operatorname*{Rect}\left(
p,p\right)  $ by%
\[
\operatorname*{Tria}\left(  p\right)  =\left\{  \left(  i,k\right)
\in\left\{  1,2,...,p\right\}  ^{2}\ \mid\ i\leq k\right\}  .
\]
This subset $\operatorname*{Tria}\left(  p\right)  $ inherits a poset
structure from $\operatorname*{Rect}\left(  p,p\right)  $. In the following,
we will consider $\operatorname*{Tria}\left(  p\right)  $ as a poset using
this structure. This poset $\operatorname*{Tria}\left(  p\right)  $ is a
$\left(  2p-1\right)  $-graded poset. It has the form of a triangle (either of
$\vartriangleleft$ shape or of $\vartriangleright$ shape, depending on how you
draw the Hasse diagram).
\end{definition}

\Needspace{32\baselineskip}

\begin{example}
Here is the Hasse diagram of the poset $\operatorname*{Rect}\left(
4,4\right)  $, with the elements that belong to $\operatorname*{Tria}\left(
4\right)  $ marked by underlines:%
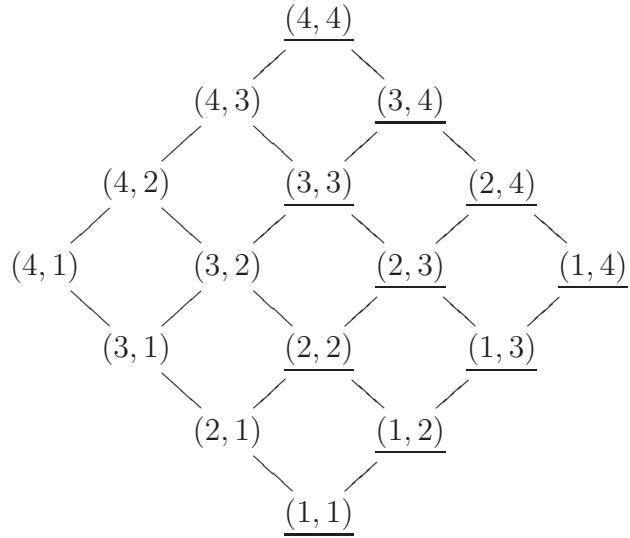
\[
\xymatrixrowsep{0.9pc}\xymatrixcolsep{0.20pc}\xymatrix{
& & & \underline{\left(4,4\right)} \ar@{-}[rd] \ar@{-}[ld] & & & \\
& & \left(4,3\right) \ar@{-}[rd] \ar@{-}[ld] & & \underline{\left(3,4\right)} \ar@{-}[rd] \ar@{-}[ld] & & \\
& \left(4,2\right) \ar@{-}[rd] \ar@{-}[ld] & & \underline{\left(3,3\right)} \ar@{-}[rd] \ar@{-}[ld] & & \underline{\left(2,4\right)} \ar@{-}[rd] \ar@{-}[ld] & \\
\left(4,1\right) \ar@{-}[rd] & & \left(3,2\right) \ar@{-}[rd] \ar@{-}[ld] & & \underline{\left(2,3\right)} \ar@{-}[rd] \ar@{-}[ld] & & \underline{\left(1,4\right)} \ar@{-}[ld] \\
& \left(3,1\right) \ar@{-}[rd] & & \underline{\left(2,2\right)} \ar@{-}[rd] \ar@{-}[ld] & & \underline{\left(1,3\right)} \ar@{-}[ld] & \\
& & \left(2,1\right) \ar@{-}[rd] & & \underline{\left(1,2\right)} \ar@{-}[ld] & & \\
& & & \underline{\left(1,1\right)} & & &
}.
\]
And here is the Hasse diagram of the poset $\operatorname*{Tria}\left(
4\right)  $ itself:%
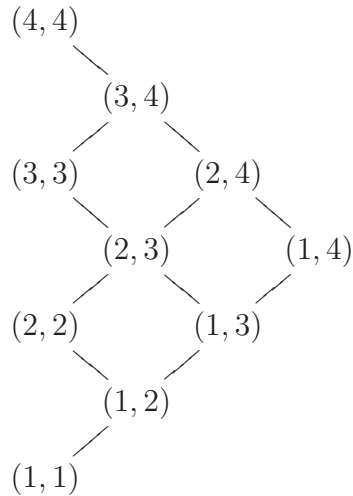
\[
\xymatrixrowsep{0.9pc}\xymatrixcolsep{0.20pc}\xymatrix{
& & & \left(4,4\right) \ar@{-}[rd] & & & \\
& & & & \left(3,4\right) \ar@{-}[rd] \ar@{-}[ld] & & \\
& & & \left(3,3\right) \ar@{-}[rd] & & \left(2,4\right) \ar@{-}[rd] \ar@{-}[ld] & \\
& & & & \left(2,3\right) \ar@{-}[rd] \ar@{-}[ld] & & \left(1,4\right) \ar@{-}[ld] \\
& & & \left(2,2\right) \ar@{-}[rd] & & \left(1,3\right) \ar@{-}[ld] & \\
& & & & \left(1,2\right) \ar@{-}[ld] & & \\
& & & \left(1,1\right) & & &
}.
\]

\end{example}

\begin{remark}
\label{rmk.Leftri.sw}Let $p$ be a positive integer. The poset
$\operatorname*{Tria}\left(  p\right)  $ appears in \cite[\S 6.2]%
{striker-williams} under the guise of the poset of order ideals (under
inclusion) of the rectangle $\operatorname*{Rect}\left(  2,p-1\right)  $. In
fact, it is easily checked that the poset of order ideals just mentioned
(denoted by $J\left(  \left[  2\right]  \times\left[  p-1\right]  \right)  $
in \cite{striker-williams}) is isomorphic to $\operatorname*{Tria}\left(
p\right)  $.
\end{remark}

We could also consider the subset $\left\{  \left(  i,k\right)  \in\left\{
1,2,...,p\right\}  ^{2}\ \mid\ i\geq k\right\}  $, but that would yield a
poset isomorphic to $\operatorname*{Tria}\left(  p\right)  $ and thus would
not be of any further interest.

\begin{theorem}
\label{thm.Leftri.ord}Let $p$ be a positive integer. Let $\mathbb{K}$ be a
field. Then, $\operatorname*{ord}\left(  R_{\operatorname*{Tria}\left(
p\right)  }\right)  =2p$.
\end{theorem}

This theorem yields $\operatorname*{ord}\left(  \overline{R}%
_{\operatorname*{Tria}\left(  p\right)  }\right)  \mid2p$. It can be shown
that actually $\operatorname*{ord}\left(  \overline{R}_{\operatorname*{Tria}%
\left(  p\right)  }\right)  =2p$ for $p>3$, while $\operatorname*{ord}\left(
\overline{R}_{\operatorname*{Tria}\left(  1\right)  }\right)  =1$,
$\operatorname*{ord}\left(  \overline{R}_{\operatorname*{Tria}\left(
2\right)  }\right)  =1$ and $\operatorname*{ord}\left(  \overline
{R}_{\operatorname*{Tria}\left(  3\right)  }\right)  =2$.

Again, Theorem \ref{thm.Leftri.ord} is the birational version of a known
result on classical rowmotion: From \cite[Theorem 6.2]{striker-williams} (and
our Remark \ref{rmk.Leftri.sw}), it follows that $\operatorname*{ord}\left(
\mathbf{r}_{\operatorname*{Tria}\left(  p\right)  }\right)  =2p$ (using the
notations of Definition \ref{def.classical.rm} and Definition
\ref{def.classical.conv.rP}). Theorem \ref{thm.Leftri.ord} thus shows that
birational rowmotion and classical rowmotion have the same order for
$\operatorname*{Tria}\left(  p\right)  $.

In order to prove Theorem \ref{thm.Leftri.ord}, we need a way to turn
labellings of $\operatorname*{Tria}\left(  p\right)  $ into labellings of
$\operatorname*{Rect}\left(  p,p\right)  $ in a rowmotion-equivariant way. It
turns out that the obvious \textquotedblleft unfolding\textquotedblright%
\ construction (with some fudge coefficients) works:

\begin{lemma}
\label{lem.Leftri.vrefl}Let $p$ be a positive integer. Let $\mathbb{K}$ be a
field of characteristic $\neq2$.

\textbf{(a)} Let $\operatorname*{vrefl}:\operatorname*{Rect}\left(
p,p\right)  \rightarrow\operatorname*{Rect}\left(  p,p\right)  $ be the map
sending every $\left(  i,k\right)  \in\operatorname*{Rect}\left(  p,p\right)
$ to $\left(  k,i\right)  $. This map $\operatorname*{vrefl}$ is an involutive
poset automorphism of $\operatorname*{Rect}\left(  p,p\right)  $. (In
intuitive terms, $\operatorname*{vrefl}$ is simply reflection across the
vertical axis.) We have $\operatorname*{vrefl}\left(  v\right)  \in
\operatorname*{Tria}\left(  p\right)  $ for every $v\in\operatorname*{Rect}%
\left(  p,p\right)  \setminus\operatorname*{Tria}\left(  p\right)  $.

We extend $\operatorname*{vrefl}$ to an involutive poset automorphism of
$\widehat{\operatorname*{Rect}\left(  p,p\right)  }$ by setting
$\operatorname*{vrefl}\left(  0\right)  =0$ and $\operatorname*{vrefl}\left(
1\right)  =1$.

\textbf{(b)} Define a map $\operatorname*{dble}:\mathbb{K}%
^{\widehat{\operatorname*{Tria}\left(  p\right)  }}\rightarrow\mathbb{K}%
^{\widehat{\operatorname*{Rect}\left(  p,p\right)  }}$ by setting%
\[
\left(  \operatorname*{dble}f\right)  \left(  v\right)  =\left\{
\begin{array}
[c]{l}%
\dfrac{1}{2}f\left(  1\right)  ,\ \ \ \ \ \ \ \ \ \ \text{if }v=1;\\
2f\left(  0\right)  ,\ \ \ \ \ \ \ \ \ \ \text{if }v=0;\\
f\left(  v\right)  ,\ \ \ \ \ \ \ \ \ \ \text{if }v\in\operatorname*{Tria}%
\left(  p\right)  ;\\
f\left(  \operatorname*{vrefl}\left(  v\right)  \right)
,\ \ \ \ \ \ \ \ \ \ \text{otherwise}%
\end{array}
\right.
\]
for all $v\in\widehat{\operatorname*{Rect}\left(  p,p\right)  }$ for all
$f\in\mathbb{K}^{\widehat{\operatorname*{Tria}\left(  p\right)  }}$. This is
well-defined. We have%
\begin{equation}
\left(  \operatorname*{dble}f\right)  \left(  v\right)  =f\left(  v\right)
\ \ \ \ \ \ \ \ \ \ \text{for every }v\in\operatorname*{Tria}\left(  p\right)
. \label{lem.Leftri.vrefl.b.1}%
\end{equation}
Also,%
\begin{equation}
\left(  \operatorname*{dble}f\right)  \left(  \operatorname*{vrefl}\left(
v\right)  \right)  =f\left(  v\right)  \ \ \ \ \ \ \ \ \ \ \text{for every
}v\in\operatorname*{Tria}\left(  p\right)  . \label{lem.Leftri.vrefl.b.2}%
\end{equation}

\textbf{(c)} We have%
\[
R_{\operatorname*{Rect}\left(  p,p\right)  }\circ\operatorname*{dble}%
=\operatorname*{dble}\circ R_{\operatorname*{Tria}\left(  p\right)  }.
\]

\end{lemma}

The coefficients $\dfrac{1}{2}$ and $2$ in the definition of
$\operatorname*{dble}$ ensure that the equality $R_{\operatorname*{Rect}%
\left(  p,p\right)  }\circ\operatorname*{dble}=\operatorname*{dble}\circ
R_{\operatorname*{Tria}\left(  p\right)  }$ in part \textbf{(c)} of the Lemma
holds on the level of labellings and not just up to homogeneous equivalence.

\begin{proof}
[Proof of Lemma \ref{lem.Leftri.vrefl} (sketched).]\textbf{(a)} Obvious.

\textbf{(b)} The well-definedness of $\operatorname*{dble}$ is pretty obvious.
The relation (\ref{lem.Leftri.vrefl.b.1}) follows from the definition of
$\operatorname*{dble}$. The relation (\ref{lem.Leftri.vrefl.b.2}) follows from
the fact that every $v\in\operatorname*{Tria}\left(  p\right)  $ satisfies
either $\operatorname*{vrefl}\left(  v\right)  \notin\operatorname*{Tria}%
\left(  p\right)  \cup\left\{  0,1\right\}  $ (in which case the definition of
$\operatorname*{dble}f$ yields $\left(  \operatorname*{dble}f\right)  \left(
\operatorname*{vrefl}\left(  v\right)  \right)  =f\left(
\underbrace{\operatorname*{vrefl}\left(  \operatorname*{vrefl}\left(
v\right)  \right)  }_{=v}\right)  =f\left(  v\right)  $) or
$\operatorname*{vrefl}\left(  v\right)  =v$ (in which case \newline$\left(
\operatorname*{dble}f\right)  \left(  \underbrace{\operatorname*{vrefl}\left(
v\right)  }_{=v}\right)  =\left(  \operatorname*{dble}f\right)  \left(
v\right)  =f\left(  v\right)  $ again by the definition of
$\operatorname*{dble}f$). This proves Lemma \ref{lem.Leftri.vrefl}
\textbf{(b)}.

\begin{vershort}
\textbf{(c)} We need to check that $\operatorname*{dble}\circ
R_{\operatorname*{Tria}\left(  p\right)  }=R_{\operatorname*{Rect}\left(
p,p\right)  }\circ\operatorname*{dble}$. In other words, we have to prove that
$\left(  \operatorname*{dble}\circ R_{\operatorname*{Tria}\left(  p\right)
}\right)  f=\left(  R_{\operatorname*{Rect}\left(  p,p\right)  }%
\circ\operatorname*{dble}\right)  f$ for every $f\in\mathbb{K}%
^{\widehat{\operatorname*{Tria}\left(  p\right)  }}$ for which
$R_{\operatorname*{Tria}\left(  p\right)  }\left(  f\right)  $ is
well-defined. So let $f\in\mathbb{K}^{\widehat{\operatorname*{Tria}\left(
p\right)  }}$ be such that $R_{\operatorname*{Tria}\left(  p\right)  }\left(
f\right)  $ is well-defined.
\end{vershort}

\begin{verlong}
\textbf{(c)} We need to check that $R_{\operatorname*{Rect}\left(  p,p\right)
}\circ\operatorname*{dble}=\operatorname*{dble}\circ R_{\operatorname*{Tria}%
\left(  p\right)  }$. In other words, we have to prove that $\left(
R_{\operatorname*{Rect}\left(  p,p\right)  }\circ\operatorname*{dble}\right)
f=\left(  \operatorname*{dble}\circ R_{\operatorname*{Tria}\left(  p\right)
}\right)  f$ for every $f\in\mathbb{K}^{\widehat{\operatorname*{Tria}\left(
p\right)  }}$ for which $R_{\operatorname*{Tria}\left(  p\right)  }\left(
f\right)  $ is well-defined.

So let $f\in\mathbb{K}^{\widehat{\operatorname*{Tria}\left(  p\right)  }}$ be
such that $R_{\operatorname*{Tria}\left(  p\right)  }\left(  f\right)  $ is
well-defined. We need to prove that $\left(  R_{\operatorname*{Rect}\left(
p,p\right)  }\circ\operatorname*{dble}\right)  f=\left(  \operatorname*{dble}%
\circ R_{\operatorname*{Tria}\left(  p\right)  }\right)  f$. In other words,
we need to prove that $\left(  \operatorname*{dble}\circ
R_{\operatorname*{Tria}\left(  p\right)  }\right)  f=\left(
R_{\operatorname*{Rect}\left(  p,p\right)  }\circ\operatorname*{dble}\right)
f$.
\end{verlong}

Set $f^{\prime}=\operatorname*{dble}f$ and $g=R_{\operatorname*{Tria}\left(
p\right)  }f$. Set $g^{\prime}=\operatorname*{dble}g$. Then,
\[
\left(  \operatorname*{dble}\circ R_{\operatorname*{Tria}\left(  p\right)
}\right)  f=\operatorname*{dble}\left(  \underbrace{R_{\operatorname*{Tria}%
\left(  p\right)  }f}_{=g}\right)  =\operatorname*{dble}g=g^{\prime}%
\]
and
\[
\left(  R_{\operatorname*{Rect}\left(  p,p\right)  }\circ\operatorname*{dble}%
\right)  f=R_{\operatorname*{Rect}\left(  p,p\right)  }\left(
\underbrace{\operatorname*{dble}f}_{=f^{\prime}}\right)
=R_{\operatorname*{Rect}\left(  p,p\right)  }f^{\prime}.
\]
Thus, our goal (namely, to prove that $\left(  \operatorname*{dble}\circ
R_{\operatorname*{Tria}\left(  p\right)  }\right)  f=\left(
R_{\operatorname*{Rect}\left(  p,p\right)  }\circ\operatorname*{dble}\right)
f$) is equivalent to showing that $g^{\prime}=R_{\operatorname*{Rect}\left(
p,p\right)  }f^{\prime}$.

So we need to prove that $g^{\prime}=R_{\operatorname*{Rect}\left(
p,p\right)  }f^{\prime}$. Since $f^{\prime}\left(  0\right)  =g^{\prime
}\left(  0\right)  $ (because the operation $\operatorname*{dble}$ multiplies
the label at $0$ with $2$, while the operation $R_{\operatorname*{Tria}\left(
p\right)  }$ leaves it unchanged) and $f^{\prime}\left(  1\right)  =g^{\prime
}\left(  1\right)  $ (for a similar reason), we know from Proposition
\ref{prop.R.implicit.converse} (applied to $\operatorname*{Rect}\left(
p,p\right)  $, $f^{\prime}$ and $g^{\prime}$ instead of $P$, $f$ and $g$) that
this will be done if we can show that%
\begin{equation}
g^{\prime}\left(  v\right)  =\dfrac{1}{f^{\prime}\left(  v\right)  }%
\cdot\dfrac{\sum\limits_{\substack{u\in\widehat{\operatorname*{Rect}\left(
p,p\right)  };\\u\lessdot v}}f^{\prime}\left(  u\right)  }{\sum
\limits_{\substack{u\in\widehat{\operatorname*{Rect}\left(  p,p\right)
};\\u\gtrdot v}}\dfrac{1}{g^{\prime}\left(  u\right)  }}%
\ \ \ \ \ \ \ \ \ \ \text{for every }v\in\operatorname*{Rect}\left(
p,p\right)  . \label{pf.Leftri.vrefl.goal}%
\end{equation}
Our goal is therefore to prove (\ref{pf.Leftri.vrefl.goal}).

\begin{verlong}
Let us first recall that $\operatorname*{vrefl}$ is a poset automorphism of
$\widehat{\operatorname*{Rect}\left(  p,p\right)  }$. Hence, if $u$ and $v$
are two elements of $\widehat{\operatorname*{Rect}\left(  p,p\right)  }$, then
we have the following equivalence of statements:%
\begin{equation}
\left(  u\lessdot v\right)  \Longleftrightarrow\left(  \operatorname*{vrefl}%
\left(  u\right)  \lessdot\operatorname*{vrefl}\left(  v\right)  \right)  .
\label{pf.Leftri.vrefl.equiv}%
\end{equation}
Similarly, if $u$ and $v$ are two elements of $\widehat{\operatorname*{Rect}%
\left(  p,p\right)  }$, then we have the following equivalence of statements:%
\begin{equation}
\left(  u\gtrdot v\right)  \Longleftrightarrow\left(  \operatorname*{vrefl}%
\left(  u\right)  \gtrdot\operatorname*{vrefl}\left(  v\right)  \right)  .
\label{pf.Leftri.vrefl.equiv2}%
\end{equation}
Also, $\operatorname*{vrefl}$ is a bijection $\widehat{\operatorname*{Rect}%
\left(  p,p\right)  }\rightarrow\widehat{\operatorname*{Rect}\left(
p,p\right)  }$ (since $\operatorname*{vrefl}$ is a poset automorphism).
\end{verlong}

But every $v\in\operatorname*{Tria}\left(  p\right)  $ satisfies%
\[
g\left(  v\right)  =\left(  R_{\operatorname*{Tria}\left(  p\right)
}f\right)  \left(  v\right)  =\dfrac{1}{f\left(  v\right)  }\cdot\dfrac
{\sum\limits_{\substack{u\in\widehat{\operatorname*{Tria}\left(  p\right)
};\\u \lessdot v}}f\left(  u\right)  }{\sum\limits_{\substack{u\in
\widehat{\operatorname*{Tria}\left(  p\right)  };\\u \gtrdot v}}\dfrac
{1}{\left(  R_{\operatorname*{Tria}\left(  p\right)  }f\right)  \left(
u\right)  }}%
\]
(by Proposition \ref{prop.R.implicit}, applied to $\operatorname*{Tria}\left(
p\right)  $ instead of $P$). Since $R_{\operatorname*{Tria}\left(  p\right)
}f=g$, this rewrites as%
\begin{equation}
g\left(  v\right)  =\dfrac{1}{f\left(  v\right)  }\cdot\dfrac{\sum
\limits_{\substack{u\in\widehat{\operatorname*{Tria}\left(  p\right)  };\\u
\lessdot v}}f\left(  u\right)  }{\sum\limits_{\substack{u\in
\widehat{\operatorname*{Tria}\left(  p\right)  };\\u \gtrdot v}}\dfrac
{1}{g\left(  u\right)  }}. \label{pf.Leftri.vrefl.2}%
\end{equation}

Now, let us prove (\ref{pf.Leftri.vrefl.goal}). So fix $v\in
\operatorname*{Rect}\left(  p,p\right)  $. Write $v$ in the form $v=\left(
i,k\right)  \in\left\{  1,2,...,p\right\}  ^{2}$. We distinguish between three cases:

\textit{Case 1:} We have $i<k$.

\textit{Case 2:} We have $i=k$.

\textit{Case 3:} We have $i>k$.

Let us first consider Case 1. In this case, $i<k$. As a consequence, every
$u\in\widehat{\operatorname*{Rect}\left(  p,p\right)  }$ satisfying $u\lessdot
v$ lies in $\operatorname*{Tria}\left(  p\right)  $. Hence, every
$u\in\widehat{\operatorname*{Rect}\left(  p,p\right)  }$ satisfying $u\lessdot
v$ satisfies
\begin{equation}
\underbrace{f^{\prime}}_{=\operatorname*{dble}f}\left(  u\right)  =\left(
\operatorname*{dble}f\right)  \left(  u\right)  =f\left(  u\right)
\label{pf.Leftri.vrefl.case1.1}%
\end{equation}
(by (\ref{lem.Leftri.vrefl.b.1})). Thus,%
\begin{equation}
\sum\limits_{\substack{u\in\widehat{\operatorname*{Rect}\left(  p,p\right)
};\\u\lessdot v}}\underbrace{f^{\prime}\left(  u\right)  }_{=f\left(
u\right)  }=\sum\limits_{\substack{u\in\widehat{\operatorname*{Rect}\left(
p,p\right)  };\\u\lessdot v}}f\left(  u\right)  =\sum\limits_{\substack{u\in
\widehat{\operatorname*{Tria}\left(  p\right)  };\\u\lessdot v}}f\left(
u\right)  \label{pf.Leftri.vrefl.case1.1a}%
\end{equation}
(since the elements $u\in\widehat{\operatorname*{Tria}\left(  p\right)  }$
such that $u\lessdot v$ are precisely the elements $u\in
\widehat{\operatorname*{Rect}\left(  p,p\right)  }$ such that $u\lessdot v$).

Also, every $u\in\widehat{\operatorname*{Rect}\left(  p,p\right)  }$
satisfying $u\gtrdot v$ lies in $\operatorname*{Tria}\left(  p\right)  $.
Hence, every $u\in\widehat{\operatorname*{Rect}\left(  p,p\right)  }$
satisfying $u\gtrdot v$ satisfies%
\begin{equation}
\underbrace{g^{\prime}}_{=\operatorname*{dble}g}\left(  u\right)  =\left(
\operatorname*{dble}g\right)  \left(  u\right)  =g\left(  u\right)
\label{pf.Leftri.vrefl.case1.2}%
\end{equation}
(by (\ref{lem.Leftri.vrefl.b.1})). Hence,%
\begin{equation}
\sum\limits_{\substack{u\in\widehat{\operatorname*{Rect}\left(  p,p\right)
};\\u\gtrdot v}}\dfrac{1}{g^{\prime}\left(  u\right)  }=\sum
\limits_{\substack{u\in\widehat{\operatorname*{Rect}\left(  p,p\right)
};\\u\gtrdot v}}\dfrac{1}{g\left(  u\right)  }=\sum\limits_{\substack{u\in
\widehat{\operatorname*{Tria}\left(  p\right)  };\\u\gtrdot v}}\dfrac
{1}{g\left(  u\right)  } \label{pf.Leftri.vrefl.case1.2a}%
\end{equation}
(because the elements $u\in\widehat{\operatorname*{Tria}\left(  p\right)  }$
such that $u\gtrdot v$ are precisely the elements $u\in
\widehat{\operatorname*{Rect}\left(  p,p\right)  }$ such that $u\gtrdot v$).

Finally, from $i<k$, we have $v\in\operatorname*{Tria}\left(  p\right)  $, so
that $\underbrace{f^{\prime}}_{=\operatorname*{dble}f}\left(  v\right)
=\left(  \operatorname*{dble}f\right)  \left(  v\right)  =f\left(  v\right)  $
(by (\ref{lem.Leftri.vrefl.b.1})) and similarly $g^{\prime}\left(  v\right)
=g\left(  v\right)  $.

Using the equalities (\ref{pf.Leftri.vrefl.case1.1a}) and
(\ref{pf.Leftri.vrefl.case1.2a}) as well as $f^{\prime}\left(  v\right)
=f\left(  v\right)  $ and $g^{\prime}\left(  v\right)  =g\left(  v\right)  $,
we can rewrite (\ref{pf.Leftri.vrefl.goal}) as%
\[
g\left(  v\right)  =\dfrac{1}{f\left(  v\right)  }\cdot\dfrac{\sum
\limits_{\substack{u\in\widehat{\operatorname*{Tria}\left(  p\right)
};\\u\lessdot v}}f\left(  u\right)  }{\sum\limits_{\substack{u\in
\widehat{\operatorname*{Tria}\left(  p\right)  };\\u\gtrdot v}}\dfrac
{1}{g\left(  u\right)  }}.
\]
But this follows from (\ref{pf.Leftri.vrefl.2}). Since
(\ref{pf.Leftri.vrefl.2}) is known to hold, we thus have proven
(\ref{pf.Leftri.vrefl.goal}) in Case 1.

\begin{vershort}
Let us next consider Case 3. It is very easy to check that every
$h\in\operatorname*{dble}\left(  \mathbb{K}^{\widehat{\operatorname*{Tria}%
\left(  p\right)  }}\right)  $ satisfies $h\left(  \operatorname*{vrefl}%
\left(  w\right)  \right)  =h\left(  w\right)  $ for every $w\in
\widehat{\operatorname*{Rect}\left(  p,p\right)  }$. Applied to $h=f^{\prime}$
(which belongs to $\operatorname*{dble}\left(  \mathbb{K}%
^{\widehat{\operatorname*{Tria}\left(  p\right)  }}\right)  $ because
$f^{\prime}=\operatorname*{dble}f$), this yields $f^{\prime}\left(
\operatorname*{vrefl}\left(  w\right)  \right)  =f^{\prime}\left(  w\right)  $
for every $w\in\widehat{\operatorname*{Rect}\left(  p,p\right)  }$. But
applied to $h=g^{\prime}$ (which belongs to $\operatorname*{dble}\left(
\mathbb{K}^{\widehat{\operatorname*{Tria}\left(  p\right)  }}\right)  $
because $g^{\prime}=\operatorname*{dble}g$), the same property yields
$g^{\prime}\left(  \operatorname*{vrefl}\left(  w\right)  \right)  =g^{\prime
}\left(  w\right)  $ for every $w\in\widehat{\operatorname*{Rect}\left(
p,p\right)  }$. We thus can rewrite the equality (\ref{pf.Leftri.vrefl.goal})
(which we desire to prove) by replacing each $g^{\prime}\left(  w\right)  $ by
$g^{\prime}\left(  \operatorname*{vrefl}\left(  w\right)  \right)  $ and by
replacing each $f^{\prime}\left(  w\right)  $ by $f^{\prime}\left(
\operatorname*{vrefl}\left(  w\right)  \right)  $. Additionally, we can
replace \textquotedblleft$u\lessdot v$\textquotedblright\ by \textquotedblleft%
$\operatorname*{vrefl}\left(  u\right)  \lessdot\operatorname*{vrefl}\left(
v\right)  $\textquotedblright, and replace \textquotedblleft$u\gtrdot
v$\textquotedblright\ by \textquotedblleft$\operatorname*{vrefl}\left(
u\right)  \gtrdot\operatorname*{vrefl}\left(  v\right)  $\textquotedblright.
Consequently, (\ref{pf.Leftri.vrefl.goal}) rewrites as%
\begin{equation}
g^{\prime}\left(  \operatorname*{vrefl}\left(  v\right)  \right)  =\dfrac
{1}{f^{\prime}\left(  \operatorname*{vrefl}\left(  v\right)  \right)  }%
\cdot\dfrac{\sum\limits_{\substack{u\in\widehat{\operatorname*{Rect}\left(
p,p\right)  };\\\operatorname*{vrefl}\left(  u\right)  \lessdot
\operatorname*{vrefl}\left(  v\right)  }}f^{\prime}\left(
\operatorname*{vrefl}\left(  u\right)  \right)  }{\sum\limits_{\substack{u\in
\widehat{\operatorname*{Rect}\left(  p,p\right)  };\\\operatorname*{vrefl}%
\left(  u\right)  \gtrdot\operatorname*{vrefl}\left(  v\right)  }}\dfrac
{1}{g^{\prime}\left(  \operatorname*{vrefl}\left(  u\right)  \right)  }}.
\label{pf.Leftri.vrefl.goal.short.1}%
\end{equation}
This equality can be simplified further by substituting $u$ for
$\operatorname*{vrefl}\left(  u\right)  $ on its right hand side:%
\begin{equation}
g^{\prime}\left(  \operatorname*{vrefl}\left(  v\right)  \right)  =\dfrac
{1}{f^{\prime}\left(  \operatorname*{vrefl}\left(  v\right)  \right)  }%
\cdot\dfrac{\sum\limits_{\substack{u\in\widehat{\operatorname*{Rect}\left(
p,p\right)  };\\u\lessdot\operatorname*{vrefl}\left(  v\right)  }}f^{\prime
}\left(  u\right)  }{\sum\limits_{\substack{u\in\widehat{\operatorname*{Rect}%
\left(  p,p\right)  };\\u\gtrdot\operatorname*{vrefl}\left(  v\right)
}}\dfrac{1}{g^{\prime}\left(  u\right)  }}.
\label{pf.Leftri.vrefl.goal.short.2}%
\end{equation}
This is precisely the statement of (\ref{pf.Leftri.vrefl.goal}) with
$\operatorname*{vrefl}\left(  v\right)  $ instead of $v$. But since we are in
Case 3 with our element $v$, we have $i>k$, so that $k<i$, and thus the
element $\operatorname*{vrefl}\left(  v\right)  =\left(  k,i\right)  $ of
$\operatorname*{Rect}\left(  p,p\right)  $ is in Case 1. Having already
verified (\ref{pf.Leftri.vrefl.goal}) in Case 1, we can thus apply
(\ref{pf.Leftri.vrefl.goal}) to $\operatorname*{vrefl}\left(  v\right)  $
instead of $v$, and conclude that (\ref{pf.Leftri.vrefl.goal.short.2}) holds.
This, as we know, is equivalent to (\ref{pf.Leftri.vrefl.goal}), and so
(\ref{pf.Leftri.vrefl.goal}) is proven in Case 3.
\end{vershort}

\begin{verlong}
Let us now consider Case 3 (the reason why we are doing the cases out of order
is that Case 3 is more similar to Case 1 than to Case 2). In this case, $i>k$.
As a consequence, $\operatorname*{vrefl}\left(  v\right)  \in
\operatorname*{Tria}\left(  p\right)  $.

Moreover, since $i>k$, it is clear that every $u\in
\widehat{\operatorname*{Rect}\left(  p,p\right)  }$ satisfying $u\lessdot v$
satisfies $\operatorname*{vrefl}\left(  u\right)  \in\operatorname*{Tria}%
\left(  p\right)  $. Hence, every $u\in\widehat{\operatorname*{Rect}\left(
p,p\right)  }$ satisfying $u\lessdot v$ satisfies
\[
\underbrace{f^{\prime}}_{=\operatorname*{dble}f}\left(  \underbrace{u}%
_{=\operatorname*{vrefl}\left(  \operatorname*{vrefl}\left(  u\right)
\right)  }\right)  =\left(  \operatorname*{dble}f\right)  \left(
\operatorname*{vrefl}\left(  \operatorname*{vrefl}\left(  u\right)  \right)
\right)  =f\left(  \operatorname*{vrefl}\left(  u\right)  \right)
\]
(by (\ref{lem.Leftri.vrefl.b.2}), since $\operatorname*{vrefl}\left(
u\right)  \in\operatorname*{Tria}\left(  p\right)  $). Thus,%
\begin{align}
\sum\limits_{\substack{u\in\widehat{\operatorname*{Rect}\left(  p,p\right)
};\\u\lessdot v}}\underbrace{f^{\prime}\left(  u\right)  }_{=f\left(
\operatorname*{vrefl}\left(  u\right)  \right)  }  &  =\sum
\limits_{\substack{u\in\widehat{\operatorname*{Rect}\left(  p,p\right)
};\\u\lessdot v}}f\left(  \operatorname*{vrefl}\left(  u\right)  \right)
=\sum\limits_{\substack{u\in\widehat{\operatorname*{Rect}\left(  p,p\right)
};\\\operatorname*{vrefl}\left(  u\right)  \lessdot\operatorname*{vrefl}%
\left(  v\right)  }}f\left(  \operatorname*{vrefl}\left(  u\right)  \right)
\nonumber\\
&  \ \ \ \ \ \ \ \ \ \ \left(
\begin{array}
[c]{c}%
\text{here, we replaced the summation sign \textquotedblleft}\sum
\limits_{\substack{u\in\widehat{\operatorname*{Rect}\left(  p,p\right)
};\\u\lessdot v}}\text{\textquotedblright}\\
\text{by a \textquotedblleft}\sum\limits_{\substack{u\in
\widehat{\operatorname*{Rect}\left(  p,p\right)  };\\\operatorname*{vrefl}%
\left(  u\right)  \lessdot\operatorname*{vrefl}\left(  v\right)
}}\text{\textquotedblright\ (because of the equivalence
(\ref{pf.Leftri.vrefl.equiv}))}%
\end{array}
\right) \nonumber\\
&  =\sum\limits_{\substack{u\in\widehat{\operatorname*{Rect}\left(
p,p\right)  };\\u\lessdot\operatorname*{vrefl}\left(  v\right)  }}f\left(
u\right) \nonumber\\
&  \ \ \ \ \ \ \ \ \ \ \left(
\begin{array}
[c]{c}%
\text{here, we substituted }u\text{ for }\operatorname*{vrefl}\left(
u\right)  \text{ in the sum, since}\\
\operatorname*{vrefl}\text{ is a bijection}%
\end{array}
\right) \nonumber\\
&  =\sum\limits_{\substack{u\in\widehat{\operatorname*{Tria}\left(  p\right)
};\\u\lessdot\operatorname*{vrefl}\left(  v\right)  }}f\left(  u\right)
\label{pf.Leftri.vrefl.case3.1a}%
\end{align}
(because the elements $u\in\widehat{\operatorname*{Tria}\left(  p\right)  }$
such that $u\lessdot\operatorname*{vrefl}\left(  v\right)  $ are precisely the
elements $u\in\widehat{\operatorname*{Rect}\left(  p,p\right)  }$ such that
$u\lessdot\operatorname*{vrefl}\left(  v\right)  $).

Also, every $u\in\widehat{\operatorname*{Rect}\left(  p,p\right)  }$
satisfying $u\gtrdot v$ satisfies $\operatorname*{vrefl}\left(  u\right)
\in\operatorname*{Tria}\left(  p\right)  $. Hence, every $u\in
\widehat{\operatorname*{Rect}\left(  p,p\right)  }$ satisfying $u\gtrdot v$
satisfies%
\[
\underbrace{g^{\prime}}_{=\operatorname*{dble}g}\left(  \underbrace{u}%
_{=\operatorname*{vrefl}\left(  \operatorname*{vrefl}\left(  u\right)
\right)  }\right)  =\left(  \operatorname*{dble}g\right)  \left(
\operatorname*{vrefl}\left(  \operatorname*{vrefl}\left(  u\right)  \right)
\right)  =g\left(  \operatorname*{vrefl}\left(  u\right)  \right)
\]
(by (\ref{lem.Leftri.vrefl.b.2}), since $\operatorname*{vrefl}\left(
u\right)  \in\operatorname*{Tria}\left(  p\right)  $). Therefore,%
\begin{align}
\sum\limits_{\substack{u\in\widehat{\operatorname*{Rect}\left(  p,p\right)
};\\u\gtrdot v}}\underbrace{\dfrac{1}{g^{\prime}\left(  u\right)  }}%
_{=\dfrac{1}{g\left(  \operatorname*{vrefl}\left(  u\right)  \right)  }}  &
=\sum\limits_{\substack{u\in\widehat{\operatorname*{Rect}\left(  p,p\right)
};\\u\gtrdot v}}\dfrac{1}{g\left(  \operatorname*{vrefl}\left(  u\right)
\right)  }=\sum\limits_{\substack{u\in\widehat{\operatorname*{Rect}\left(
p,p\right)  };\\\operatorname*{vrefl}\left(  u\right)  \gtrdot
\operatorname*{vrefl}\left(  v\right)  }}\dfrac{1}{g\left(
\operatorname*{vrefl}\left(  u\right)  \right)  }\nonumber\\
&  \ \ \ \ \ \ \ \ \ \ \left(
\begin{array}
[c]{c}%
\text{here, we replaced the summation sign \textquotedblleft}\sum
\limits_{\substack{u\in\widehat{\operatorname*{Rect}\left(  p,p\right)
};\\u\gtrdot v}}\text{\textquotedblright}\\
\text{by a \textquotedblleft}\dfrac{1}{g\left(  \operatorname*{vrefl}\left(
u\right)  \right)  }\text{\textquotedblright\ (because of the equivalence
(\ref{pf.Leftri.vrefl.equiv2}))}%
\end{array}
\right) \nonumber\\
&  =\sum\limits_{\substack{u\in\widehat{\operatorname*{Rect}\left(
p,p\right)  };\\u\gtrdot\operatorname*{vrefl}\left(  v\right)  }}\dfrac
{1}{g\left(  u\right)  }\nonumber\\
&  \ \ \ \ \ \ \ \ \ \ \left(
\begin{array}
[c]{c}%
\text{here, we substituted }u\text{ for }\operatorname*{vrefl}\left(
u\right)  \text{ in the sum, since}\\
\operatorname*{vrefl}\text{ is a bijection}%
\end{array}
\right) \nonumber\\
&  =\sum\limits_{\substack{u\in\widehat{\operatorname*{Tria}\left(  p\right)
};\\u\gtrdot\operatorname*{vrefl}\left(  v\right)  }}\dfrac{1}{g\left(
u\right)  } \label{pf.Leftri.vrefl.case3.2a}%
\end{align}
(because the elements $u\in\widehat{\operatorname*{Tria}\left(  p\right)  }$
such that $u\gtrdot\operatorname*{vrefl}\left(  v\right)  $ are precisely the
elements $u\in\widehat{\operatorname*{Rect}\left(  p,p\right)  }$ such that
$u\gtrdot\operatorname*{vrefl}\left(  v\right)  $).

Finally, from $i>k$, we have $\operatorname*{vrefl}\left(  v\right)
\in\operatorname*{Tria}\left(  p\right)  $, so that $\underbrace{f^{\prime}%
}_{=\operatorname*{dble}f}\left(  \underbrace{v}_{=\operatorname*{vrefl}%
\left(  \operatorname*{vrefl}\left(  v\right)  \right)  }\right)  =\left(
\operatorname*{dble}f\right)  \left(  \operatorname*{vrefl}\left(
\operatorname*{vrefl}\left(  v\right)  \right)  \right)  =f\left(
\operatorname*{vrefl}\left(  v\right)  \right)  $ (by
(\ref{lem.Leftri.vrefl.b.2}), since $\operatorname*{vrefl}\left(  v\right)
\in\operatorname*{Tria}\left(  p\right)  $) and similarly $g^{\prime}\left(
v\right)  =g\left(  \operatorname*{vrefl}\left(  v\right)  \right)  $.

Using the equalities (\ref{pf.Leftri.vrefl.case3.1a}) and
(\ref{pf.Leftri.vrefl.case3.2a}) as well as $f^{\prime}\left(  v\right)
=f\left(  \operatorname*{vrefl}\left(  v\right)  \right)  $ and $g^{\prime
}\left(  v\right)  =g\left(  \operatorname*{vrefl}\left(  v\right)  \right)
$, we can rewrite (\ref{pf.Leftri.vrefl.goal}) as%
\[
g\left(  \operatorname*{vrefl}\left(  v\right)  \right)  =\dfrac{1}{f\left(
\operatorname*{vrefl}\left(  v\right)  \right)  }\cdot\dfrac{\sum
\limits_{\substack{u\in\widehat{\operatorname*{Tria}\left(  p\right)
};\\u\lessdot\operatorname*{vrefl}\left(  v\right)  }}f\left(  u\right)
}{\sum\limits_{\substack{u\in\widehat{\operatorname*{Tria}\left(  p\right)
};\\u\gtrdot\operatorname*{vrefl}\left(  v\right)  }}\dfrac{1}{g\left(
u\right)  }}.
\]
But this follows from (\ref{pf.Leftri.vrefl.2}) (applied to
$\operatorname*{vrefl}\left(  v\right)  $ instead of $v$). Since
(\ref{pf.Leftri.vrefl.2}) is known to hold, we thus have proven
(\ref{pf.Leftri.vrefl.goal}) in Case 3.
\end{verlong}

Let us finally consider Case 2. In this case, $i=k$. Thus, $v=\left(
i,\underbrace{k}_{=i}\right)  =\left(  i,i\right)  $. Hence, $v\in
\widehat{\operatorname*{Tria}\left(  p\right)  }$. Thus,
$\underbrace{f^{\prime}}_{=\operatorname*{dble}f}\left(  v\right)  =\left(
\operatorname*{dble}f\right)  \left(  v\right)  =f\left(  v\right)  $ (by
(\ref{lem.Leftri.vrefl.b.1}), since $v\in\operatorname*{Tria}\left(  p\right)
$). Similarly, $g^{\prime}\left(  v\right)  =g\left(  v\right)  $.

We should now consider four subcases, depending on whether $i\notin\left\{
1,p\right\}  $ or $i=1\neq p$ or $i=p\neq1$ or $i=1=p$. But we are only going
to deal with the first of these subcases here, leaving the other three to the
reader. So let us consider the subcase when $i\notin\left\{  1,p\right\}  $.

We have $v=\left(  i,i\right)  $. Thus, the only element $u\in
\widehat{\operatorname*{Tria}\left(  p\right)  }$ such that $u\gtrdot v$ is
$\left(  i,i+1\right)  $, and the only element $u\in
\widehat{\operatorname*{Tria}\left(  p\right)  }$ such that $u\lessdot v$ is
$\left(  i-1,i\right)  $. Thus, (\ref{pf.Leftri.vrefl.2}) simplifies to%
\begin{equation}
g\left(  v\right)  =\dfrac{1}{f\left(  v\right)  }\cdot\dfrac{f\left(  \left(
i-1,i\right)  \right)  }{\left(  \dfrac{1}{g\left(  \left(  i,i+1\right)
\right)  }\right)  }. \label{pf.Leftri.vrefl.5}%
\end{equation}

\begin{vershort}
Now, recall that $g^{\prime}=\operatorname*{dble}g$. From the definition of
$\operatorname*{dble}g$, it therefore follows easily that $g^{\prime}\left(
\left(  i,i+1\right)  \right)  =g\left(  \left(  i,i+1\right)  \right)  $ and
$g^{\prime}\left(  \left(  i+1,i\right)  \right)  =g\left(  \left(
i,i+1\right)  \right)  $.

Also, $f^{\prime}=\operatorname*{dble}f$. From the definition of
$\operatorname*{dble}f$, we thus obtain $f^{\prime}\left(  \left(
i-1,i\right)  \right)  =f\left(  \left(  i-1,i\right)  \right)  $ and
$f^{\prime}\left(  \left(  i,i-1\right)  \right)  =f\left(  \left(
i-1,i\right)  \right)  $.
\end{vershort}

\begin{verlong}
Now, since $g^{\prime}=\operatorname*{dble}g$, we have $g^{\prime}\left(
\left(  i,i+1\right)  \right)  =\left(  \operatorname*{dble}g\right)  \left(
\left(  i,i+1\right)  \right)  =g\left(  \left(  i,i+1\right)  \right)  $ (by
(\ref{lem.Leftri.vrefl.b.1}), since $\left(  i,i+1\right)  \in
\operatorname*{Tria}\left(  p\right)  $) and \newline$g^{\prime}\left(
\operatorname*{vrefl}\left(  \left(  i,i+1\right)  \right)  \right)  =\left(
\operatorname*{dble}g\right)  \left(  \operatorname*{vrefl}\left(  \left(
i,i+1\right)  \right)  \right)  =g\left(  \left(  i,i+1\right)  \right)  $ (by
(\ref{lem.Leftri.vrefl.b.2}), since $\left(  i,i+1\right)  \in
\operatorname*{Tria}\left(  p\right)  $). The latter equality rewrites as
$g^{\prime}\left(  \left(  i+1,i\right)  \right)  =g\left(  \left(
i,i+1\right)  \right)  $ (since $\operatorname*{vrefl}\left(  \left(
i,i+1\right)  \right)  =\left(  i+1,i\right)  $).

Also, since $f^{\prime}=\operatorname*{dble}f$, we have $f^{\prime}\left(
\left(  i-1,i\right)  \right)  =\left(  \operatorname*{dble}f\right)  \left(
\left(  i-1,i\right)  \right)  =f\left(  \left(  i-1,i\right)  \right)  $ (by
(\ref{lem.Leftri.vrefl.b.1}), since $\left(  i-1,i\right)  \in
\operatorname*{Tria}\left(  p\right)  $) and \newline$f^{\prime}\left(
\operatorname*{vrefl}\left(  \left(  i-1,i\right)  \right)  \right)  =\left(
\operatorname*{dble}f\right)  \left(  \operatorname*{vrefl}\left(  \left(
i-1,i\right)  \right)  \right)  =f\left(  \left(  i-1,i\right)  \right)  $ (by
(\ref{lem.Leftri.vrefl.b.2}), since $\left(  i-1,i\right)  \in
\operatorname*{Tria}\left(  p\right)  $). The latter equality rewrites as
$f^{\prime}\left(  \left(  i,i-1\right)  \right)  =f\left(  \left(
i-1,i\right)  \right)  $ (since $\operatorname*{vrefl}\left(  \left(
i-1,i\right)  \right)  =\left(  i,i-1\right)  $).
\end{verlong}

But the elements $u\in\widehat{\operatorname*{Rect}\left(  p,p\right)  }$ such
that $u \gtrdot v$ are precisely $\left(  i+1,i\right)  $ and $\left(
i,i+1\right)  $, and the elements $u\in\widehat{\operatorname*{Rect}\left(
p,p\right)  }$ such that $u \lessdot v$ are precisely $\left(  i-1,i\right)  $
and $\left(  i,i-1\right)  $. Thus, the right hand side of
(\ref{pf.Leftri.vrefl.goal}) simplifies as follows:%
\begin{align*}
&  \dfrac{1}{f^{\prime}\left(  v\right)  }\cdot\dfrac{\sum
\limits_{\substack{u\in\widehat{\operatorname*{Rect}\left(  p,p\right)  };\\u
\lessdot v}}f^{\prime}\left(  u\right)  }{\sum\limits_{\substack{u\in
\widehat{\operatorname*{Rect}\left(  p,p\right)  };\\u \gtrdot v}}\dfrac
{1}{g^{\prime}\left(  u\right)  }}\\
&  =\dfrac{1}{f^{\prime}\left(  v\right)  }\cdot\dfrac{f^{\prime}\left(
\left(  i-1,i\right)  \right)  +f^{\prime}\left(  \left(  i,i-1\right)
\right)  }{\dfrac{1}{g^{\prime}\left(  \left(  i+1,i\right)  \right)  }%
+\dfrac{1}{g^{\prime}\left(  \left(  i,i+1\right)  \right)  }}=\dfrac
{1}{f\left(  v\right)  }\cdot\dfrac{f\left(  \left(  i-1,i\right)  \right)
+f\left(  \left(  i-1,i\right)  \right)  }{\dfrac{1}{g\left(  \left(
i,i+1\right)  \right)  }+\dfrac{1}{g\left(  \left(  i,i+1\right)  \right)  }%
}\\
&  \ \ \ \ \ \ \ \ \ \ \left(
\begin{array}
[c]{c}%
\text{since }f^{\prime}\left(  \left(  i-1,i\right)  \right)  =f\left(
\left(  i-1,i\right)  \right)  \text{, }f^{\prime}\left(  \left(
i,i-1\right)  \right)  =f\left(  \left(  i-1,i\right)  \right)  \text{,}\\
g^{\prime}\left(  \left(  i+1,i\right)  \right)  =g\left(  \left(
i,i+1\right)  \right)  \text{, }g^{\prime}\left(  \left(  i,i+1\right)
\right)  =g\left(  \left(  i,i+1\right)  \right) \\
\text{and }f^{\prime}\left(  v\right)  =f\left(  v\right)
\end{array}
\right) \\
&  =\dfrac{1}{f\left(  v\right)  }\cdot\dfrac{2\cdot f\left(  \left(
i-1,i\right)  \right)  }{2\cdot\dfrac{1}{g\left(  \left(  i,i+1\right)
\right)  }}=\dfrac{1}{f\left(  v\right)  }\cdot\dfrac{f\left(  \left(
i-1,i\right)  \right)  }{\left(  \dfrac{1}{g\left(  \left(  i,i+1\right)
\right)  }\right)  }=g\left(  v\right)  \ \ \ \ \ \ \ \ \ \ \left(  \text{by
(\ref{pf.Leftri.vrefl.5})}\right) \\
&  =g^{\prime}\left(  v\right)  .
\end{align*}
In other words, (\ref{pf.Leftri.vrefl.goal}) is proven in Case 2.

We have now proven (\ref{pf.Leftri.vrefl.goal}) in all three cases (not
counting the subcases which we left to the reader to ``enjoy''). Thus,
(\ref{pf.Leftri.vrefl.goal}) holds, and as we know this yields that
$g^{\prime}=R_{\operatorname*{Rect}\left(  p,p\right)  }f^{\prime}$. Lemma
\ref{lem.Leftri.vrefl} \textbf{(c)} is thus proven.
\end{proof}

\begin{proof}
[Proof of Theorem \ref{thm.Leftri.ord} (sketched).]Applying Proposition
\ref{prop.ord-projord} to $2p-1$ and $\operatorname*{Tria}\left(  p\right)  $
instead of $n$ and $P$, we obtain $\operatorname*{ord}\left(
R_{\operatorname*{Tria}\left(  p\right)  }\right)  =\operatorname{lcm}\left(
2p-1+1,\operatorname*{ord}\left(  \overline{R}_{\operatorname*{Tria}\left(
p\right)  }\right)  \right)  $. Hence, $\operatorname*{ord}\left(
R_{\operatorname*{Tria}\left(  p\right)  }\right)  $ is divisible by
$2p-1+1=2p$. Now, if we can prove that $\operatorname*{ord}\left(
R_{\operatorname*{Tria}\left(  p\right)  }\right)  \mid2p$, then we will
immediately obtain $\operatorname*{ord}\left(  R_{\operatorname*{Tria}\left(
p\right)  }\right)  =2p$, and Theorem \ref{thm.Leftri.ord} will be proven.

So let us show that $\operatorname*{ord}\left(  R_{\operatorname*{Tria}\left(
p\right)  }\right)  \mid2p$. This means showing that $R_{\operatorname*{Tria}%
\left(  p\right)  }^{2p}=\operatorname*{id}$. Since this statement boils down
to a collection of polynomial identities in the labels of an arbitrary
$\mathbb{K}$-labelling of $\operatorname*{Tria}\left(  p\right)  $, it is
clear that it is enough to prove it in the case when $\mathbb{K}$ is a field
of rational functions in finitely many variables over $\mathbb{Q}$. So let us
WLOG assume that $\mathbb{K}$ is a field of rational functions in finitely
many variables over $\mathbb{Q}$. Then, the characteristic of $\mathbb{K}$ is
$\neq2$ (it is $0$ indeed), so that we can apply Lemma \ref{lem.Leftri.vrefl}.

Let us use the notations of Lemma \ref{lem.Leftri.vrefl}. Lemma
\ref{lem.Leftri.vrefl} \textbf{(c)} yields%
\[
R_{\operatorname*{Rect}\left(  p,p\right)  }\circ\operatorname*{dble}%
=\operatorname*{dble}\circ R_{\operatorname*{Tria}\left(  p\right)  }.
\]
From this, it follows (by induction over $k$) that
\[
R_{\operatorname*{Rect}\left(  p,p\right)  }^{k}\circ\operatorname*{dble}%
=\operatorname*{dble}\circ R_{\operatorname*{Tria}\left(  p\right)  }^{k}%
\]
for every $k\in\mathbb{N}$. Applied to $k=2p$, this yields%
\begin{equation}
R_{\operatorname*{Rect}\left(  p,p\right)  }^{2p}\circ\operatorname*{dble}%
=\operatorname*{dble}\circ R_{\operatorname*{Tria}\left(  p\right)  }^{2p}.
\label{pf.Leftri.ord.1}%
\end{equation}
But Theorem \ref{thm.rect.ord} (applied to $q=p$) yields $\operatorname*{ord}%
\left(  R_{\operatorname*{Rect}\left(  p,p\right)  }\right)  =p+p=2p$, so that
$R_{\operatorname*{Rect}\left(  p,p\right)  }^{2p}=\operatorname*{id}$. Hence,
(\ref{pf.Leftri.ord.1}) simplifies to
\[
\operatorname*{dble}=\operatorname*{dble}\circ R_{\operatorname*{Tria}\left(
p\right)  }^{2p}.
\]
We can cancel $\operatorname*{dble}$ from this equation, because
$\operatorname*{dble}$ is an injective and therefore left-cancellable map. As
a consequence, we obtain $\operatorname*{id}=R_{\operatorname*{Tria}\left(
p\right)  }^{2p}$. In other words, $R_{\operatorname*{Tria}\left(  p\right)
}^{2p}=\operatorname*{id}$. This proves Theorem \ref{thm.Leftri.ord}.
\end{proof}

\section{\texorpdfstring{The $\Delta$ and $\nabla$ triangles}{The Delta and Nabla triangles}}%

\label{sect.DeltaNabla}

The next kind of triangle-shaped posets is more interesting.

\begin{definition}
\label{def.DeltaNabla}Let $p$ be a positive integer. Define three subsets
$\Delta\left(  p\right)  $, $\operatorname*{Eq}\left(  p\right)  $ and
$\nabla\left(  p\right)  $ of
$\operatorname*{Rect}\left(  p,p\right)
= \left\{ 1,2,...,p\right\} \times \left\{ 1,2,...,p\right\}
= \left\{ 1,2,...,p\right\}  ^{2}$ by%
\begin{align*}
\Delta\left(  p\right)   &  =\left\{  \left(  i,k\right)  \in\left\{
1,2,...,p\right\}  ^{2}\ \mid\ i+k>p+1\right\}  ;\\
\operatorname*{Eq}\left(  p\right)   &  =\left\{  \left(  i,k\right)
\in\left\{  1,2,...,p\right\}  ^{2}\ \mid\ i+k=p+1\right\}  ;\\
\nabla\left(  p\right)   &  =\left\{  \left(  i,k\right)  \in\left\{
1,2,...,p\right\}  ^{2}\ \mid\ i+k<p+1\right\}  .
\end{align*}
These subsets $\Delta\left(  p\right)  $, $\operatorname*{Eq}\left(  p\right)
$ and $\nabla\left(  p\right)  $ inherit a poset structure from
$\operatorname*{Rect}\left(  p,p\right)  $. In the following, we will consider
$\Delta\left(  p\right)  $, $\operatorname*{Eq}\left(  p\right)  $ and
$\nabla\left(  p\right)  $ as posets using this structure.

Clearly, $\operatorname*{Eq}\left(  p\right)  $ is an antichain with $p$
elements. (The name $\operatorname*{Eq}$ comes from ``equator''.) The posets
$\Delta\left(  p\right)  $ and $\nabla\left(  p\right)  $ are $\left(
p-1\right)  $-graded posets. They have the form of a ``Delta-shaped triangle''
and a ``Nabla-shaped triangle'', respectively (whence the names).
\end{definition}

\Needspace{42\baselineskip}

\begin{example}
Here is the Hasse diagram of the poset $\operatorname*{Rect}\left(
4,4\right)  $, where the elements belonging to $\Delta\left(  4\right)  $ have
been underlined and the elements belonging to $\operatorname*{Eq}\left(
4\right)  $ have been boxed:%
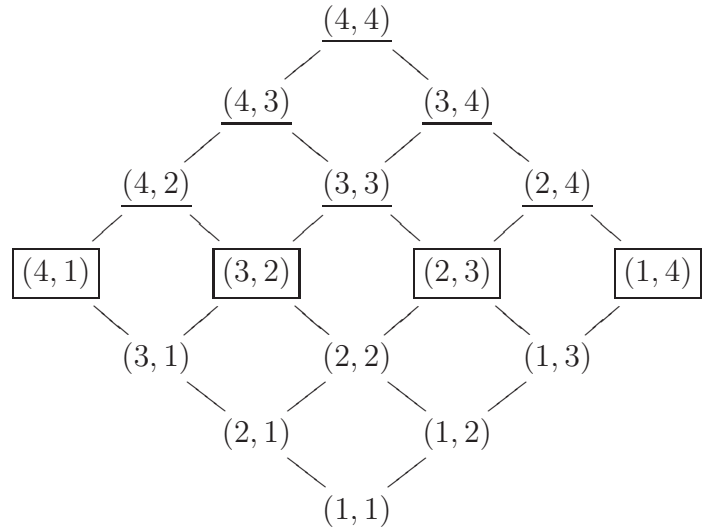
\[
\xymatrixrowsep{0.9pc}\xymatrixcolsep{0.20pc}\xymatrix{
& & & \underline{\left(4,4\right)} \ar@{-}[rd] \ar@{-}[ld] & & & \\
& & \underline{\left(4,3\right)} \ar@{-}[rd] \ar@{-}[ld] & & \underline{\left(3,4\right)} \ar@{-}[rd] \ar@{-}[ld] & & \\
& \underline{\left(4,2\right)} \ar@{-}[rd] \ar@{-}[ld] & & \underline{\left(3,3\right)} \ar@{-}[rd] \ar@{-}[ld] & & \underline{\left(2,4\right)} \ar@{-}[rd] \ar@{-}[ld] & \\
\fbox{$\left(4,1\right)$} \ar@{-}[rd] & & \fbox{$\left(3,2\right)$} \ar@{-}[rd] \ar@{-}[ld] & & \fbox{$\left(2,3\right)$} \ar@{-}[rd] \ar@{-}[ld] & & \fbox{$\left(1,4\right)$} \ar@{-}[ld] \\
& \left(3,1\right) \ar@{-}[rd] & & \left(2,2\right) \ar@{-}[rd] \ar@{-}[ld] & & \left(1,3\right) \ar@{-}[ld] & \\
& & \left(2,1\right) \ar@{-}[rd] & & \left(1,2\right) \ar@{-}[ld] & & \\
& & & \left(1,1\right) & & &
}.
\]
And here is the Hasse diagram of the poset $\Delta\left(  4\right)  $ itself:%
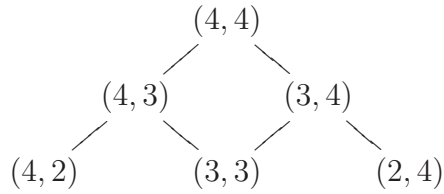
\[
\xymatrixrowsep{0.9pc}\xymatrixcolsep{0.20pc}\xymatrix{
& & & \left(4,4\right) \ar@{-}[rd] \ar@{-}[ld] & & & \\
& & \left(4,3\right) \ar@{-}[rd] \ar@{-}[ld] & & \left(3,4\right) \ar@{-}[rd] \ar@{-}[ld] & & \\
& \left(4,2\right) & & \left(3,3\right) & & \left(2,4\right) &
}.
\]
Here, on the other hand, is the Hasse diagram of the poset $\nabla\left(
4\right)  $:
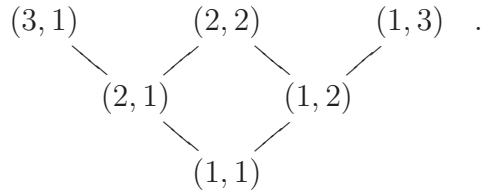
\[
\xymatrixrowsep{0.9pc}\xymatrixcolsep{0.20pc}\xymatrix{
& \left(3,1\right) \ar@{-}[rd] & & \left(2,2\right) \ar@{-}[rd] \ar@{-}[ld] & & \left(1,3\right) \ar@{-}[ld] & \\
& & \left(2,1\right) \ar@{-}[rd] & & \left(1,2\right) \ar@{-}[ld] & & \\
& & & \left(1,1\right) & & &
}.
\]

\end{example}

\begin{remark}
\label{rmk.Delta.sw}Let $p$ be a positive integer. The poset $\Delta\left(
p\right)  $ is isomorphic to the poset $\Phi^{+}\left(  A_{p-1}\right)  $ of
\cite[\S 3.2]{striker-williams}.
\end{remark}

\begin{remark}
\label{rmk.DeltaNabla}For every positive integer $p$, we have $\nabla\left(
p\right)  \cong\left(  \Delta\left(  p\right)  \right)  ^{\operatorname*{op}}$
as posets. This follows immediately from the poset antiautomorphism%
\begin{align*}
\operatorname*{hrefl}:\operatorname*{Rect}\left(  p,p\right)   &
\rightarrow\operatorname*{Rect}\left(  p,p\right)  ,\\
\left(  i,k\right)   &  \mapsto\left(  p+1-k,p+1-i\right)
\end{align*}
sending $\nabla\left(  p\right)  $ to $\Delta\left(  p\right)  $.
\end{remark}

Here we are using the following notions:

\begin{definition}
\textbf{(a)} If $P$ and $Q$ are two posets, then a map $f:P\rightarrow Q$ is
called a \textit{poset antihomomorphism} if and only if every $p_{1}\in P$ and
$p_{2}\in P$ satisfying $p_{1}\leq p_{2}$ in $P$ satisfy $f\left(
p_{1}\right)  \geq f\left(  p_{2}\right)  $ in $Q$. It is easy to see that the
poset antihomomorphisms $P\rightarrow Q$ are precisely the poset homomorphisms
$P\rightarrow Q^{\operatorname*{op}}$.

\textbf{(b)} If $P$ and $Q$ are two posets, then an invertible map
$f:P\rightarrow Q$ is called a \textit{poset antiisomorphism} if and only if
both $f$ and $f^{-1}$ are poset antihomomorphisms.

\textbf{(c)} If $P$ is a poset and $f:P\rightarrow P$ is an invertible map,
then $f$ is said to be a \textit{poset antiautomorphism} if $f$ is a poset antiisomorphism.
\end{definition}

We now state the main property of birational rowmotion $R$ on the posets
$\nabla\left(  p\right)  $ and $\Delta\left(  p\right)  $:

\begin{theorem}
\label{thm.Nabla.halfway}Let $p$ be an integer $\geq1$. Let $\mathbb{K}$ be a
field. For every $\left(  i,k\right)  \in\nabla\left(  p\right)  $ and every
$f\in\mathbb{K}^{\widehat{\nabla\left(  p\right)  }}$, we have%
\[
\left(  R_{\nabla\left(  p\right)  }^{p}f\right)  \left(  \left(  i,k\right)
\right)  =f\left(  \left(  k,i\right)  \right)  .
\]

\end{theorem}

\begin{theorem}
\label{thm.Delta.halfway}Let $p$ be an integer $\geq1$. Let $\mathbb{K}$ be a
field. For every $\left(  i,k\right)  \in\Delta\left(  p\right)  $ and every
$f\in\mathbb{K}^{\widehat{\Delta\left(  p\right)  }}$, we have%
\[
\left(  R_{\Delta\left(  p\right)  }^{p}f\right)  \left(  \left(  i,k\right)
\right)  =f\left(  \left(  k,i\right)  \right)  .
\]

\end{theorem}

The following two corollaries follow easily from these two theorems:

\begin{corollary}
\label{cor.Nabla.ord}Let $p$ be an integer $>1$. Let $\mathbb{K}$ be a field. Then:

\textbf{(a)} We have $\operatorname*{ord}\left(  R_{\nabla\left(  p\right)
}\right)  \mid2p$.

\textbf{(b)} If $p>2$, then $\operatorname*{ord}\left(  R_{\nabla\left(
p\right)  }\right)  =2p$.
\end{corollary}

\begin{corollary}
\label{cor.Delta.ord}Let $p$ be an integer $>1$. Let $\mathbb{K}$ be a field. Then:

\textbf{(a)} We have $\operatorname*{ord}\left(  R_{\Delta\left(  p\right)
}\right)  \mid2p$.

\textbf{(b)} If $p>2$, then $\operatorname*{ord}\left(  R_{\Delta\left(
p\right)  }\right)  =2p$.
\end{corollary}

Corollary \ref{cor.Delta.ord} is analogous to a known result for classical
rowmotion. In fact, from \cite[Conjecture 3.6]{striker-williams} (originally a
conjecture of Panyushev, then proven by Armstrong, Stump and Thomas) and our
Remark \ref{rmk.Delta.sw}, it can be seen that (using the notations of
Definition \ref{def.classical.rm} and Definition \ref{def.classical.conv.rP})
every integer $p>2$ satisfies $\operatorname*{ord}\left(  \mathbf{r}%
_{\Delta\left(  p\right)  }\right)  =2p$.

We now prepare for the proofs of Theorems \ref{thm.Nabla.halfway} and
\ref{thm.Delta.halfway}.

First of all, Corollary \ref{cor.Nabla.ord} is clearly equivalent to Corollary
\ref{cor.Delta.ord} (because of Remark \ref{rmk.DeltaNabla} and Proposition
\ref{prop.op.ord}). It is a bit more complicated to see that Theorem
\ref{thm.Nabla.halfway} is equivalent to Theorem \ref{thm.Delta.halfway}; we
will show this later. But let us first prove Theorem \ref{thm.Delta.halfway}.
The proof will use a mapping that transforms labellings of $\Delta\left(
p\right)  $ into labellings of $\operatorname*{Rect}\left(  p,p\right)  $ in a
way that is rowmotion-equivariant up to homogeneous equivalence. This mapping
is similar in its function to the mapping $\operatorname*{dble}$ of Lemma
\ref{lem.Leftri.vrefl}, but its definition is more intricate:\footnote{See
also Lemma \ref{lem.Delta.hrefl-general} further below for a generalization
of parts of this construction.}

\begin{lemma}
\label{lem.Delta.hrefl}Let $p$ be a positive integer. Clearly,
$\operatorname*{Rect}\left(  p,p\right)  $ is the disjoint union of the sets
$\Delta\left(  p\right)  $, $\nabla\left(  p\right)  $ and $\operatorname*{Eq}%
\left(  p\right)  $.

Let $\mathbb{K}$ be a field of characteristic $\neq2$.

\textbf{(a)} Let $\operatorname*{hrefl}:\operatorname*{Rect}\left(
p,p\right)  \rightarrow\operatorname*{Rect}\left(  p,p\right)  $ be the map
sending every $\left(  i,k\right)  \in\operatorname*{Rect}\left(  p,p\right)
$ to $\left(  p+1-k,p+1-i\right)  $. This map $\operatorname*{hrefl}$ is an
involution and a poset antiautomorphism of $\operatorname*{Rect}\left(
p,p\right)  $. (In intuitive terms, $\operatorname*{hrefl}$ is simply
reflection across the horizontal axis (i.e., the line $\operatorname*{Eq}%
\left(  p\right)  $).) We have $\operatorname*{hrefl}\mid_{\operatorname*{Eq}%
\left(  p\right)  }=\operatorname*{id}$ and $\operatorname*{hrefl}\left(
\Delta\left(  p\right)  \right)  =\nabla\left(  p\right)  $.

We extend $\operatorname*{hrefl}$ to an involutive poset antiautomorphism of
$\widehat{\operatorname*{Rect}\left(  p,p\right)  }$ by setting
$\operatorname*{hrefl}\left(  0\right)  =1$ and $\operatorname*{hrefl}\left(
1\right)  =0$.

\textbf{(b)} Define a rational map $\operatorname*{wing}:\mathbb{K}%
^{\widehat{\Delta\left(  p\right)  }}\dashrightarrow\mathbb{K}%
^{\widehat{\operatorname*{Rect}\left(  p,p\right)  }}$ by setting%
\[
\left(  \operatorname*{wing}f\right)  \left(  v\right)  =\left\{
\begin{array}
[c]{l}%
f\left(  v\right)  ,\ \ \ \ \ \ \ \ \ \ \text{if }v\in\Delta\left(  p\right)
\cup\left\{  1\right\}  ;\\
1,\ \ \ \ \ \ \ \ \ \ \text{if }v\in\operatorname*{Eq}\left(  p\right)  ;\\
\dfrac{1}{\left(  R_{\Delta\left(  p\right)  }^{p-\deg v}f\right)  \left(
\operatorname*{hrefl}v\right)  },\ \ \ \ \ \ \ \ \ \ \text{if }v\in
\nabla\left(  p\right)  \cup\left\{  0\right\}
\end{array}
\right.
\]
for all $v\in\widehat{\operatorname*{Rect}\left(  p,p\right)  }$ for all
$f\in\mathbb{K}^{\widehat{\Delta\left(  p\right)  }}$. This is well-defined.

\textbf{(c)} There exists a rational map $\overline{\operatorname*{wing}%
}:\overline{\mathbb{K}^{\widehat{\Delta\left(  p\right)  }}}\dashrightarrow
\overline{\mathbb{K}^{\widehat{\operatorname*{Rect}\left(  p,p\right)  }}}$
such that the diagram%
\begin{equation}
\xymatrixcolsep{5pc}
\xymatrix{ \mathbb{K}^{\widehat{\Delta\left(p\right)}} \ar@{-->}[d]_-{\pi} 
\ar@{-->}[r]^{\operatorname*{wing}} &
\mathbb{K}^{\widehat{\operatorname *{Rect}\left(p,p\right)}} \ar@{-->}[d]^-{\pi} \\ \overline{\mathbb{K}^{\widehat{\Delta\left(p\right)}}} \ar@{-->}[r]_{\overline{\operatorname*{wing}}} & \overline{\mathbb{K}^{\widehat {\operatorname*{Rect}\left(p,p\right)}}} }
\label{lem.Delta.hrefl.commut}%
\end{equation}
commutes.

\textbf{(d)} The rational map
$\overline{\operatorname*{wing}}$ defined in Lemma
\ref{lem.Delta.hrefl} \textbf{(c)} satisfies
\[
\overline{R}_{\operatorname*{Rect}\left(  p,p\right)  }\circ\overline
{\operatorname*{wing}}=\overline{\operatorname*{wing}}\circ\overline
{R}_{\Delta\left(  p\right)  }.
\]

\textbf{(e)} Consider the map $\operatorname*{vrefl}:\operatorname*{Rect}%
\left(  p,p\right)  \rightarrow\operatorname*{Rect}\left(  p,p\right)  $
defined in Lemma \ref{lem.Leftri.vrefl}. Define a map $\operatorname*{vrefl}%
\nolimits^{\ast}:\mathbb{K}^{\widehat{\operatorname*{Rect}\left(  p,p\right)
}}\rightarrow\mathbb{K}^{\widehat{\operatorname*{Rect}\left(  p,p\right)  }}$
by setting%
\[
\left(  \operatorname*{vrefl}\nolimits^{\ast}f\right)  \left(  v\right)
=f\left(  \operatorname*{vrefl}\left(  v\right)  \right)
\ \ \ \ \ \ \ \ \ \ \text{for all }v\in\widehat{\operatorname*{Rect}\left(
p,p\right)  }%
\]
for all $f\in\mathbb{K}^{\widehat{\operatorname*{Rect}\left(  p,p\right)  }}$.
Also, define a map $\operatorname*{vrefl}\nolimits^{\ast}:\mathbb{K}%
^{\widehat{\Delta\left(  p\right)  }}\rightarrow\mathbb{K}^{\widehat{\Delta
\left(  p\right)  }}$ by setting%
\[
\left(  \operatorname*{vrefl}\nolimits^{\ast}f\right)  \left(  v\right)
=f\left(  \operatorname*{vrefl}\left(  v\right)  \right)
\ \ \ \ \ \ \ \ \ \ \text{for all }v\in\widehat{\Delta\left(  p\right)  }%
\]
for all $f\in\mathbb{K}^{\widehat{\Delta\left(  p\right)  }}$. Then,%
\begin{equation}
\operatorname*{vrefl}\nolimits^{\ast}\circ R_{\Delta\left(  p\right)
}=R_{\Delta\left(  p\right)  }\circ\operatorname*{vrefl}\nolimits^{\ast}
\label{lem.Delta.hrefl.e.1}%
\end{equation}
(as rational maps $\mathbb{K}^{\widehat{\Delta\left(  p\right)  }%
}\dashrightarrow\mathbb{K}^{\widehat{\Delta\left(  p\right)  }}$).
Furthermore,%
\begin{equation}
\operatorname*{vrefl}\nolimits^{\ast}\circ R_{\operatorname*{Rect}\left(
p,p\right)  }=R_{\operatorname*{Rect}\left(  p,p\right)  }\circ
\operatorname*{vrefl}\nolimits^{\ast} \label{lem.Delta.hrefl.e.2}%
\end{equation}
(as rational maps $\mathbb{K}^{\widehat{\operatorname*{Rect}\left(
p,p\right)  }}\dashrightarrow\mathbb{K}^{\widehat{\operatorname*{Rect}\left(
p,p\right)  }}$). Finally,
\begin{equation}
\operatorname*{vrefl}\nolimits^{\ast}\circ\operatorname*{wing}%
=\operatorname*{wing}\circ\operatorname*{vrefl}\nolimits^{\ast}
\label{lem.Delta.hrefl.e.3}%
\end{equation}
(as rational maps $\mathbb{K}^{\widehat{\Delta\left(  p\right)  }%
}\dashrightarrow\mathbb{K}^{\widehat{\operatorname*{Rect}\left(  p,p\right)
}}$).

\textbf{(f)} Almost every (in the sense of Zariski topology) labelling
$f\in\mathbb{K}^{\widehat{\Delta\left(  p\right)  }}$ satisfying $f\left(
0\right)  =2$ satisfies%
\[
R_{\operatorname*{Rect}\left(  p,p\right)  }\left(  \operatorname*{wing}%
f\right)  =\operatorname*{wing}\left(  R_{\Delta\left(  p\right)  }f\right)
.
\]

\textbf{(g)} If $f$ and $g$ are two homogeneously equivalent zero-free
$\mathbb{K}$-labellings of $\Delta\left(  p\right)  $, then
$\operatorname*{vrefl}\nolimits^{\ast}f$ is homogeneously equivalent to
$\operatorname*{vrefl}\nolimits^{\ast}g$.
\end{lemma}

\begin{vershort}
\begin{proof}
[Proof of Lemma \ref{lem.Delta.hrefl} (sketched).]We will not delve into the
details of this tedious and yet straightforward proof. Let us merely make some
comments on the few interesting parts of it. Parts \textbf{(a)}, \textbf{(b)},
\textbf{(c)} and \textbf{(g)} are obvious. Part \textbf{(f)} can be verified
label-by-label using Propositions \ref{prop.R.implicit} and
\ref{prop.R.implicit.converse} and some nasty casework. Part \textbf{(d)}
won't be used in the following, but can easily be derived from part
\textbf{(f)}. Part \textbf{(e)} more or less follows from the fact that the
definitions of $R_{\Delta\left(  p\right)  }$, $R_{\operatorname*{Rect}\left(
p,p\right)  }$ and $\operatorname*{wing}$ are all \textquotedblleft
invariant\textquotedblright\ under the vertical reflection
$\operatorname*{vrefl}$; but proving part \textbf{(e)} in a pedestrian way
might be even more straightforward than formalizing this invariance
argument\footnote{Again, Propositions \ref{prop.R.implicit} and
\ref{prop.R.implicit.converse} come in handy for proving
(\ref{lem.Delta.hrefl.e.1}) and (\ref{lem.Delta.hrefl.e.2}). Then, one can
prove (by induction over $\ell$) that $\operatorname*{vrefl}\nolimits^{\ast
}\circ R_{\Delta\left(  p\right)  }^{\ell}=R_{\Delta\left(  p\right)  }^{\ell
}\circ\operatorname*{vrefl}\nolimits^{\ast}$ for all $\ell\in\mathbb{N}$.
Using this, (\ref{lem.Delta.hrefl.e.3}) is straightforward to check.}.
\end{proof}
\end{vershort}

\begin{verlong}
\begin{proof}
[Proof of Lemma \ref{lem.Delta.hrefl} (sketched).]The proof of this lemma is
tedious and computational (at some point requiring 7 cases to be checked), yet
very straightforward. We split it into several parts.

\bigskip

\underline{\textbf{Trivialities.}}

Let us first notice that the definitions of $\Delta\left(p\right)$,
$\operatorname{Eq}\left(p\right)$ and $\nabla\left(p\right)$ can
be rewritten as follows:
\begin{align*}
\Delta\left(  p\right)   &  =\left\{  v \in\left\{
1,2,...,p\right\}  ^{2}\ \mid\ \deg v > p\right\}  ;\\
\operatorname*{Eq}\left(  p\right)   &  =\left\{ v
\in\left\{  1,2,...,p\right\}  ^{2}\ \mid\ \deg v = p\right\}  ;\\
\nabla\left(  p\right)   &  =\left\{ v \in\left\{
1,2,...,p\right\}  ^{2}\ \mid\ \deg v < p\right\}  .
\end{align*}
Indeed, these equalities are clearly equivalent to the definitions
of $\Delta\left(p\right)$,
$\operatorname{Eq}\left(p\right)$ and $\nabla\left(p\right)$
because every $\left(i, k\right) \in \left\{ 1, 2, ..., p \right\}^2$
satisfies $\deg \left(\left(i,k\right)\right) = i+k-1$.

Let us now prove some easy parts of Lemma \ref{lem.Delta.hrefl}.

\textbf{(a)} All statements made in Lemma \ref{lem.Delta.hrefl} \textbf{(a)}
are obvious.

\textbf{(b)} Again, this is obvious.

\textbf{(c)} In order to prove Lemma \ref{lem.Delta.hrefl} \textbf{(c)}, it is
enough to show that if $f$ and $g$ are two homogeneously equivalent zero-free
$\mathbb{K}$-labellings in $\mathbb{K}^{\widehat{\Delta\left(  p\right)  }}$,
then $\operatorname*{wing}f$ and $\operatorname*{wing}g$ are homogeneously
equivalent. This is left to the reader (who is advised to remember Corollary
\ref{cor.hgRi}).

\bigskip

\underline{\textbf{Proof of (d): Exordium.}}

\textbf{(d)} This (and the proof of \textbf{(f)}) is the laborious part of the
proof of Lemma \ref{lem.Delta.hrefl}. Let us start by introducing some
notations and assumptions.

\bigskip

We WLOG assume that $p>2$, because otherwise the claim of Lemma
\ref{lem.Delta.hrefl} \textbf{(d)} is rather obvious.

Let us first recall that $\operatorname*{hrefl}$ is a poset antiautomorphism
of $\widehat{\operatorname*{Rect}\left(  p,p\right)  }$. Hence, if $u$ and $v$
are two elements of $\widehat{\operatorname*{Rect}\left(  p,p\right)  }$, then
we have the following equivalence of statements:%
\begin{equation}
\left(  u\lessdot v\right)  \Longleftrightarrow\left(  \operatorname*{hrefl}%
u\gtrdot\operatorname*{hrefl}v\right)  . \label{pf.Delta.hrefl.d.equiv}%
\end{equation}
Similarly, if $u$ and $v$ are two elements of $\widehat{\operatorname*{Rect}%
\left(  p,p\right)  }$, then we have the following equivalence of statements:%
\begin{equation}
\left(  u\gtrdot v\right)  \Longleftrightarrow\left(  \operatorname*{hrefl}%
u\lessdot\operatorname*{hrefl}v\right)  . \label{pf.Delta.hrefl.d.equiv2}%
\end{equation}
Also, $\operatorname*{hrefl}$ is a bijection $\widehat{\operatorname*{Rect}%
\left(  p,p\right)  }\rightarrow\widehat{\operatorname*{Rect}\left(
p,p\right)  }$ (since $\operatorname*{hrefl}$ is a poset automorphism).
Finally, every $v \in \widehat{\operatorname*{Rect}\left(  p,p\right)  }$
satisfies
\begin{equation}
\deg \left(\operatorname{hrefl} v\right) = 2p - \deg v .
\label{pf.Delta.hrefl.hrefl-inverts-deg}
\end{equation}
\footnote{\textit{Proof.} Let
$v \in \widehat{\operatorname*{Rect}\left(  p,p\right)  }$.
We need to prove (\ref{pf.Delta.hrefl.hrefl-inverts-deg}).

If $v = 0$, then (\ref{pf.Delta.hrefl.hrefl-inverts-deg}) is true
(because
$\deg \left(\underbrace{\operatorname{hrefl} 0}_{=1}\right) =
\deg 1 = 2p = 2p - \underbrace{0}_{=\deg 0} = 2p - \deg 0$).
Similarly, if $v = 1$, then (\ref{pf.Delta.hrefl.hrefl-inverts-deg})
is true. Hence, (\ref{pf.Delta.hrefl.hrefl-inverts-deg}) is proven
in the cases when $v = 0$ and when $v = 1$. Thus, for the rest of the
proof of (\ref{pf.Delta.hrefl.hrefl-inverts-deg}), we can WLOG assume
that $v$ is neither $0$ nor $1$. Assume this.

We have $v \in \widehat{\operatorname*{Rect}\left(  p,p\right)  }
= \operatorname*{Rect}\left(  p,p\right)  \cup \left\{0,1\right\}$,
so that $v \in \operatorname*{Rect}\left(  p,p\right)$ (since $v$ is
neither $0$ nor $1$). Thus,
$v \in \operatorname*{Rect}\left(  p,p\right)
= \left\{  1,2,...,p\right\}  \times\left\{  1,2,...,p\right\}$.
Hence, we can write $v$ in the form
$v=\left(  \mathfrak{i},\mathfrak{k}\right)  $
for some $\mathfrak{i}\in\left\{  1,2,...,p\right\}  $ and $\mathfrak{k}%
\in\left\{  1,2,...,p\right\}  $. Consider these $\mathfrak{i}$ and
$\mathfrak{k}$. Then,
$\deg\underbrace{v}_{=\left(  \mathfrak{i},\mathfrak{k}\right)  }=\deg\left(
\left(  \mathfrak{i},\mathfrak{k}\right)  \right)  =\mathfrak{i}%
+\mathfrak{k}-1$ (by Remark \ref{rmk.rect.cover} \textbf{(a)}, applied to
$\left(  \mathfrak{i},\mathfrak{k}\right)  $ instead of $\left(  i,k\right)
$). On the other hand,
\begin{align*}
\operatorname*{hrefl}\underbrace{v}_{=\left(  \mathfrak{i}%
,\mathfrak{k}\right)  }
= \operatorname*{hrefl}\left(  \mathfrak{i}%
,\mathfrak{k}\right)
=\left(  p+1-\mathfrak{k},p+1-\mathfrak{i}\right)
\ \ \ \ \ \ \ \ \ \ \left(  \text{by the definition of }\operatorname*{hrefl}%
\left(  \mathfrak{i},\mathfrak{k}\right)  \right)  ,
\end{align*}
so that
\begin{align*}
\deg \left(\operatorname*{hrefl} v\right)
 &  =\deg\left(  \left(  p+1-\mathfrak{k},p+1-\mathfrak{i}\right)
\right)  =\left(  p+1-\mathfrak{k}\right)  +\left(  p+1-\mathfrak{i}\right)
-1\\
&  \ \ \ \ \ \ \ \ \ \ \left(  \text{by Remark \ref{rmk.rect.cover}
\textbf{(a)}, applied to }\left(  p+1-\mathfrak{k},p+1-\mathfrak{i}\right)
\text{ instead of }\left(  i,k\right)  \right) \\
&  =2p-\underbrace{\left(  \mathfrak{i}+\mathfrak{k}-1\right)  }%
_{= \deg v} = 2p - \deg v.
\end{align*}
This proves (\ref{pf.Delta.hrefl.hrefl-inverts-deg}).}

Let $f$ be a zero-free $\mathbb{K}$-labelling in $\mathbb{K}^{\widehat{\Delta
\left(  p\right)  }}$ which is sufficiently generic for $R_{\Delta\left(
p\right)  }^{i}f$ to be well-defined for all $i\in\left\{  0,1,...,p\right\}
$. We are going to show that $\left(  \operatorname*{wing}\circ R_{\Delta
\left(  p\right)  }\right)  f$ and $\left(  R_{\operatorname*{Rect}\left(
p,p\right)  }\circ\operatorname*{wing}\right)  f$ are homogeneously equivalent
$\mathbb{K}$-labellings of $\operatorname*{Rect}\left(  p,p\right)  $.

We will use the notation introduced in Definition \ref{def.bemol}.

Notice that $f\left(  0\right)  \neq0$. Define a $\left(  2p+1\right)  $-tuple
$\left(  a_{0},a_{1},...,a_{2p}\right)  \in\left(  \mathbb{K}^{\times}\right)
^{2p+1}$ by%
\[
a_{i}=\left\{
\begin{array}
[c]{l}%
1,\ \ \ \ \ \ \ \ \ \ \text{if }i>p+1;\\
\dfrac{2}{f\left(  0\right)  },\ \ \ \ \ \ \ \ \ \ \text{if }p\leq i\leq
p+1;\\
1,\ \ \ \ \ \ \ \ \ \ \text{if }i<p
\end{array}
\right.  \ \ \ \ \ \ \ \ \ \ \text{for every }i\in\left\{  0,1,...,2p\right\}
.
\]

We are now going to prove that%
\begin{equation}
\left(  R_{\operatorname*{Rect}\left(  p,p\right)  }\circ\operatorname*{wing}%
\right)  f=\left(  a_{0},a_{1},...,a_{2p}\right)  \flat\left(  \left(
\operatorname*{wing}\circ R_{\Delta\left(  p\right)  }\right)  f\right)  .
\label{pf.Delta.hrefl.precise}%
\end{equation}

\textit{Proof of (\ref{pf.Delta.hrefl.precise}):} Let $g=\left(
\operatorname*{wing}\circ R_{\Delta\left(  p\right)  }\right)  f$ and
$f^{\prime}=\operatorname*{wing}f$. Let $g^{\prime}=\left(  a_{0}%
,a_{1},...,a_{2p}\right)  \flat g$. We are going to prove that $g^{\prime
}=R_{\operatorname*{Rect}\left(  p,p\right)  }f^{\prime}$.

\bigskip

\underline{\textbf{Values at }$0$ \textbf{and }$1$\textbf{.}}

First, we notice that $f^{\prime}=\operatorname*{wing}f$, so that%
\begin{align}
f^{\prime}\left(  0\right)   &  =\left(  \operatorname*{wing}f\right)  \left(
0\right)  =\dfrac{1}{\left(  R_{\Delta\left(  p\right)  }^{p-\deg0}f\right)
\left(  \operatorname*{hrefl}0\right)  }\ \ \ \ \ \ \ \ \ \ \left(  \text{by
the definition of }\operatorname*{wing}\right) \nonumber\\
&  =\dfrac{1}{\left(  R_{\Delta\left(  p\right)  }^{p-0}f\right)  \left(
1\right)  }\ \ \ \ \ \ \ \ \ \ \left(  \text{since }\operatorname*{hrefl}%
0=1\text{ and }\deg0=0\right) \nonumber\\
&  =\dfrac{1}{f\left(  1\right)  }\ \ \ \ \ \ \ \ \ \ \left(  \text{since
Corollary \ref{cor.R.implicit.01} yields }\left(  R_{\Delta\left(  p\right)
}^{p-0}f\right)  \left(  1\right)  =f\left(  1\right)  \right)  .
\label{pf.Delta.hrefl.f'0}%
\end{align}

On the other hand, $f^{\prime} = \operatorname*{wing} f$, so that
\[
f^{\prime}\left(1\right)
= \left( \operatorname*{wing} f\right)\left(1\right)
= f\left(1\right)
\]
(by the definition of $\operatorname*{wing}$, since
$1 \in \Delta\left(p\right) \cup \left\{1\right\}$).

But Proposition \ref{prop.R.implicit.01} yields $\left(  R_{\Delta\left(
p\right)  }f\right)  \left(  0\right)  =f\left(  0\right)  $ and $\left(
R_{\Delta\left(  p\right)  }f\right)  \left(  1\right)  =f\left(  1\right)  $.
But since $g=\left(  \operatorname*{wing}\circ R_{\Delta\left(  p\right)
}\right)  f$, we have%
\begin{align*}
g\left(  1\right)   &  =\left(  \left(  \operatorname*{wing}\circ
R_{\Delta\left(  p\right)  }\right)  f\right)  \left(  1\right)  =\left(
\operatorname*{wing}\left(  R_{\Delta\left(  p\right)  }f\right)  \right)
\left(  1\right)  =\left(  R_{\Delta\left(  p\right)  }f\right)  \left(
1\right) \\
&  \ \ \ \ \ \ \ \ \ \ \left(  \text{by the definition of }%
\operatorname*{wing}\right) \\
&  =f\left(  1\right)  \ \ \ \ \ \ \ \ \ \ \left(  \text{by Proposition
\ref{prop.R.implicit.01}, applied to }\Delta\left(  p\right)  \text{ instead
of }P\right) \\
& = f^{\prime}\left(1\right)
\ \ \ \ \ \ \ \ \ \ 
\left(\text{since } f^{\prime}\left(1\right) = f\left(1\right) \right)
\end{align*}
and%
\begin{align*}
g\left(  0\right)   &  =\left(  \left(  \operatorname*{wing}\circ
R_{\Delta\left(  p\right)  }\right)  f\right)  \left(  0\right)  =\left(
\operatorname*{wing}\left(  R_{\Delta\left(  p\right)  }f\right)  \right)
\left(  0\right) \\
&  =\dfrac{1}{\left(  R_{\Delta\left(  p\right)  }^{p-\deg0}\left(
R_{\Delta\left(  p\right)  }f\right)  \right)  \left(  \operatorname*{hrefl}%
0\right)  }\ \ \ \ \ \ \ \ \ \ \left(  \text{by the definition of
}\operatorname*{wing}\right) \\
&  =\dfrac{1}{\left(  R_{\Delta\left(  p\right)  }^{p}\left(  R_{\Delta\left(
p\right)  }f\right)  \right)  \left(  1\right)  }\ \ \ \ \ \ \ \ \ \ \left(
\text{since }\operatorname*{hrefl}0=1\text{ and }p-\underbrace{\deg0}%
_{=0}=p-0=p\right) \\
&  =\dfrac{1}{\left(  R_{\Delta\left(  p\right)  }^{p+1}f\right)  \left(
1\right)  }=\dfrac{1}{f\left(  1\right)  }\ \ \ \ \ \ \ \ \ \ \left(
\text{since Corollary \ref{cor.R.implicit.01} yields }\left(  R_{\Delta\left(
p\right)  }^{p+1}f\right)  \left(  1\right)  =f\left(  1\right)  \right)  .
\end{align*}

But since $g^{\prime}=\left(  a_{0},a_{1},...,a_{2p}\right)  \flat g$, we have%
\begin{align*}
g^{\prime}\left(  1\right)   &  =\left(  \left(  a_{0},a_{1},...,a_{2p}%
\right)  \flat g\right)  \left(  1\right)  =a_{\deg1}\cdot g\left(  1\right)
=\underbrace{a_{2p}}_{=1}\cdot g\left(  1\right)  \ \ \ \ \ \ \ \ \ \ \left(
\text{since }\deg1=2p\right) \\
&  =g\left(  1\right)  = f\left(  1\right)  =f^{\prime}\left(  1\right)
\end{align*}
and%
\begin{align*}
g^{\prime}\left(  0\right)   &  =\left(  \left(  a_{0},a_{1},...,a_{2p}%
\right)  \flat g\right)  \left(  0\right)  =a_{\deg0}\cdot g\left(  0\right)
=\underbrace{a_{0}}_{=1}\cdot g\left(  0\right)  \ \ \ \ \ \ \ \ \ \ \left(
\text{since }\deg0=0\right) \\
&  =g\left(  0\right)  =\dfrac{1}{f\left(  1\right)  }=f^{\prime}\left(
0\right)  \ \ \ \ \ \ \ \ \ \ \left(  \text{by (\ref{pf.Delta.hrefl.f'0}%
)}\right)  .
\end{align*}

\bigskip

\underline{\textbf{The plan of attack.}}

Recall that we need to prove that $g^{\prime}=R_{\operatorname*{Rect}\left(
p,p\right)  }f^{\prime}$. Since $f^{\prime}\left(  0\right)  =g^{\prime
}\left(  0\right)  $ and $f^{\prime}\left(  1\right)  =g^{\prime}\left(
1\right)  $, we know from Proposition \ref{prop.R.implicit.converse} (applied
to $\operatorname*{Rect}\left(  p,p\right)  $, $f^{\prime}$ and $g^{\prime}$
instead of $P$, $f$ and $g$) that this will be done if we can show that%
\begin{equation}
g^{\prime}\left(  v\right)  =\dfrac{1}{f^{\prime}\left(  v\right)  }%
\cdot\dfrac{\sum\limits_{\substack{u\in\widehat{\operatorname*{Rect}\left(
p,p\right)  };\\u \lessdot v}}f^{\prime}\left(  u\right)  }{\sum
\limits_{\substack{u\in\widehat{\operatorname*{Rect}\left(  p,p\right)  };\\u
\gtrdot v}}\dfrac{1}{g^{\prime}\left(  u\right)  }}%
\ \ \ \ \ \ \ \ \ \ \text{for every }v\in\operatorname*{Rect}\left(
p,p\right)  . \label{pf.Delta.hrefl.goal}%
\end{equation}

Now, let us prove (\ref{pf.Delta.hrefl.goal}). So fix $v\in
\operatorname*{Rect}\left(  p,p\right)  $. Since $g^{\prime}=\left(
a_{0},a_{1},...,a_{2p}\right)  \flat g$, we have%
\begin{align}
g^{\prime}\left(  v\right)   &  =\left(  \left(  a_{0},a_{1},...,a_{2p}%
\right)  \flat g\right)  \left(  v\right)  =a_{\deg v}\cdot g\left(  v\right)
\nonumber\\
&  =a_{\deg v}\cdot\left(  \operatorname*{wing}\left(  R_{\Delta\left(
p\right)  }f\right)  \right)  \left(  v\right)  \label{pf.Delta.hrefl.g'v}%
\end{align}
(since $g=\left(  \operatorname*{wing}\circ R_{\Delta\left(  p\right)
}\right)  f=\operatorname*{wing}\left(  R_{\Delta\left(  p\right)  }f\right)
$).

Moreover, every $u\in\widehat{\operatorname*{Rect}\left(  p,p\right)  }$ such
that $u \gtrdot v$ satisfies%
\begin{align}
g^{\prime}\left(  u\right)   &  =\left(  \left(  a_{0},a_{1},...,a_{2p}%
\right)  \flat g\right)  \left(  u\right)  \ \ \ \ \ \ \ \ \ \ \left(
\text{since }g^{\prime}=\left(  a_{0},a_{1},...,a_{2p}\right)  \flat g\right)
\nonumber\\
&  =a_{\deg u}\cdot g\left(  u\right)  =a_{\deg v+1}\cdot g\left(  u\right)
\nonumber\\
&  \ \ \ \ \ \ \ \ \ \ \left(  \text{since }u \gtrdot v\text{, so that }\deg
u=\deg v+1\right) \nonumber\\
&  =a_{\deg v+1}\cdot\left(  \operatorname*{wing}\left(  R_{\Delta\left(
p\right)  }f\right)  \right)  \left(  u\right)
\label{pf.Delta.hrefl.g'u.ucoversv}%
\end{align}
(since $g=\left(  \operatorname*{wing}\circ R_{\Delta\left(  p\right)
}\right)  f=\operatorname*{wing}\left(  R_{\Delta\left(  p\right)  }f\right)
$).

From $p > 2$, we obtain $2p - 2 > p$ and $2p-1 > p+1 > p > p-1 > 1$.

We distinguish between the following seven cases:

\textit{Case 1:} We have $\deg v=2p-1$.

\textit{Case 2:} We have $p+1<\deg v<2p-1$.

\textit{Case 3:} We have $\deg v=p+1$.

\textit{Case 4:} We have $\deg v=p$.

\textit{Case 5:} We have $\deg v=p-1$.

\textit{Case 6:} We have $1<\deg v<p-1$.

\textit{Case 7:} We have $\deg v=1$.

\bigskip

\underline{\textbf{Proof of (\ref{pf.Delta.hrefl.goal}) in Case 1.}}

Let us first consider Case 1. In this case, we have $\deg v=2p-1$. Thus, $\deg
v=2p-1>p+1$ (because $p>2$). As a consequence, both $v$ and all the elements
$u\in\widehat{\operatorname*{Rect}\left(  p,p\right)  }$ such that $u \lessdot
v$ belong to $\Delta\left(  p\right)  $. Hence, the elements $u\in
\widehat{\operatorname*{Rect}\left(  p,p\right)  }$ such that $u \lessdot v$
are precisely the elements $u\in\widehat{\Delta\left(  p\right)  }$ such that
$u \lessdot v$.

Since $\deg v=2p-1$, the equality (\ref{pf.Delta.hrefl.g'v}) rewrites as
\begin{equation}
g^{\prime}\left(  v\right)  =\underbrace{a_{2p-1}}_{=1}\cdot
\underbrace{\left(  \operatorname*{wing}\left(  R_{\Delta\left(  p\right)
}f\right)  \right)  \left(  v\right)  }_{\substack{=\left(  R_{\Delta\left(
p\right)  }f\right)  \left(  v\right)  \\\text{(by the definition of
}\operatorname*{wing}\text{,}\\\text{since }v\in\Delta\left(  p\right)
\cup\left\{  1\right\}  \text{)}}}=\left(  R_{\Delta\left(  p\right)
}f\right)  \left(  v\right)  . \label{pf.Delta.hrefl.g'v.c1.rewriter1}%
\end{equation}

Furthermore, since $f^{\prime}=\operatorname*{wing}f$, we have%
\begin{equation}
f^{\prime}\left(  v\right)  =\left(  \operatorname*{wing}f\right)  \left(
v\right)  =f\left(  v\right)  \label{pf.Delta.hrefl.g'v.c1.rewriter2}%
\end{equation}
(by the definition of $\operatorname*{wing}$, since $v\in\Delta\left(
p\right)  \cup\left\{  1\right\}  $). But every $u\in
\widehat{\operatorname*{Rect}\left(  p,p\right)  }$ such that $u \lessdot v$
belongs to $\Delta\left(  p\right)  \cup\left\{  1\right\}  $ (because $u
\lessdot v$, hence $\deg u=\underbrace{\deg v}_{=2p-1}-1=2p-1-1=2p-2>p$).
Hence, every $u\in\widehat{\operatorname*{Rect}\left(  p,p\right)  }$ such
that $u \lessdot v$ satisfies $\left(  \operatorname*{wing}f\right)  \left(
u\right)  =f\left(  u\right)  $ (by the definition of $\operatorname*{wing}$).
Since $\operatorname*{wing}f=f^{\prime}$, this rewrites as follows: Every
$u\in\widehat{\operatorname*{Rect}\left(  p,p\right)  }$ such that $u \lessdot
v$ satisfies $f^{\prime}\left(  u\right)  =f\left(  u\right)  $. Thus,
\begin{equation}
\sum\limits_{\substack{u\in\widehat{\operatorname*{Rect}\left(  p,p\right)
};\\u \lessdot v}}f^{\prime}\left(  u\right)  =\sum\limits_{\substack{u\in
\widehat{\operatorname*{Rect}\left(  p,p\right)  };\\u \lessdot v}}f\left(
u\right)  =\sum\limits_{\substack{u\in\widehat{\Delta\left(  p\right)  };\\u
\lessdot v}}f\left(  u\right)  \label{pf.Delta.hrefl.g'v.c1.rewriter3}%
\end{equation}
(since the elements $u\in\widehat{\operatorname*{Rect}\left(  p,p\right)  }$
such that $u \lessdot v$ are precisely the elements $u\in\widehat{\Delta
\left(  p\right)  }$ such that $u \lessdot v$).

Finally, there is only one $u\in\widehat{\operatorname*{Rect}\left(
p,p\right)  }$ such that $u \gtrdot v$, namely $u=1$ (because $\deg v=2p-1$).
Hence,%
\begin{equation}
\sum\limits_{\substack{u\in\widehat{\operatorname*{Rect}\left(  p,p\right)
};\\u \gtrdot v}}\dfrac{1}{g^{\prime}\left(  u\right)  }=\dfrac{1}{g^{\prime
}\left(  1\right)  }=\dfrac{1}{f\left(  1\right)  }
\label{pf.Delta.hrefl.g'v.c1.rewriter4}%
\end{equation}
(since $g^{\prime}\left(  1\right)  =g\left(  1\right)  =f\left(  1\right)  $).

Recall that we have to prove the equality (\ref{pf.Delta.hrefl.goal}). Upon
substitution of (\ref{pf.Delta.hrefl.g'v.c1.rewriter1}),
(\ref{pf.Delta.hrefl.g'v.c1.rewriter2}),
(\ref{pf.Delta.hrefl.g'v.c1.rewriter3}) and
(\ref{pf.Delta.hrefl.g'v.c1.rewriter4}), this equality
(\ref{pf.Delta.hrefl.goal}) becomes%
\begin{equation}
\left(  R_{\Delta\left(  p\right)  }f\right)  \left(  v\right)  =\dfrac
{1}{f\left(  v\right)  }\cdot\dfrac{\sum\limits_{\substack{u\in\widehat{\Delta
\left(  p\right)  };\\u \lessdot v}}f\left(  u\right)  }{\left(  \dfrac
{1}{f\left(  1\right)  }\right)  }. \label{pf.Delta.hrefl.g'v.c1.transgoal}%
\end{equation}
So (\ref{pf.Delta.hrefl.g'v.c1.transgoal}) is what we are going to prove now.

Applying Proposition \ref{prop.R.implicit} to $\Delta\left(  p\right)  $
instead of $P$, we find%
\begin{equation}
\left(  R_{\Delta\left(  p\right)  }f\right)  \left(  v\right)  =\dfrac
{1}{f\left(  v\right)  }\cdot\dfrac{\sum\limits_{\substack{u\in\widehat{\Delta
\left(  p\right)  };\\u \lessdot v}}f\left(  u\right)  }{\sum
\limits_{\substack{u\in\widehat{\Delta\left(  p\right)  };\\u \gtrdot
v}}\dfrac{1}{\left(  R_{\Delta\left(  p\right)  }f\right)  \left(  u\right)
}}. \label{pf.Delta.hrefl.g'v.c1.1}%
\end{equation}
But since $v$ is the maximal element of $\Delta\left(  p\right)  $ (because
$\deg v=2p-1$), there is only one element $u\in\widehat{\Delta\left(
p\right)  }$ such that $u \gtrdot v$, namely $1$. Thus,
\[
\sum\limits_{\substack{u\in\widehat{\Delta\left(  p\right)  };\\u\gtrdot
v}}\dfrac{1}{\left(  R_{\Delta\left(  p\right)  }f\right)  \left(  u\right)
}=\dfrac{1}{\left(  R_{\Delta\left(  p\right)  }f\right)  \left(  1\right)
}=\dfrac{1}{f\left(  1\right)  }%
\]
(since Proposition \ref{prop.R.implicit.01} yields $\left(  R_{\Delta\left(
p\right)  }f\right)  \left(  1\right)  =f\left(  1\right)  $). In view of
this, (\ref{pf.Delta.hrefl.g'v.c1.1}) becomes%
\[
\left(  R_{\Delta\left(  p\right)  }f\right)  \left(  v\right)  =\dfrac
{1}{f\left(  v\right)  }\cdot\dfrac{\sum\limits_{\substack{u\in\widehat{\Delta
\left(  p\right)  };\\u \lessdot v}}f\left(  u\right)  }{\left(  \dfrac
{1}{f\left(  1\right)  }\right)  }.
\]
But this is precisely (\ref{pf.Delta.hrefl.g'v.c1.transgoal}). Hence,
(\ref{pf.Delta.hrefl.g'v.c1.transgoal}) is proven. Since
(\ref{pf.Delta.hrefl.g'v.c1.transgoal}) was obtained by rewriting
(\ref{pf.Delta.hrefl.goal}), this yields that (\ref{pf.Delta.hrefl.goal}) is
proven in Case 1.

\bigskip

\underline{\textbf{Proof of (\ref{pf.Delta.hrefl.goal}) in Case 2.}}

Let us now consider Case 2. In this case, we have $p+1<\deg v<2p-1$. Thus,
$v\in\Delta\left(  p\right)  \subseteq\Delta\left(  p\right)  \cup\left\{
1\right\}  $. The equality (\ref{pf.Delta.hrefl.g'v}) rewrites as
\begin{equation}
g^{\prime}\left(  v\right)  =\underbrace{a_{\deg v}}%
_{\substack{=1\\\text{(since }\deg v>p+1\text{)}}}\cdot\underbrace{\left(
\operatorname*{wing}\left(  R_{\Delta\left(  p\right)  }f\right)  \right)
\left(  v\right)  }_{\substack{=\left(  R_{\Delta\left(  p\right)  }f\right)
\left(  v\right)  \\\text{(by the definition of }\operatorname*{wing}%
\text{,}\\\text{since }v\in\Delta\left(  p\right)  \cup\left\{  1\right\}
\text{)}}}=\left(  R_{\Delta\left(  p\right)  }f\right)  \left(  v\right)  .
\label{pf.Delta.hrefl.g'v.c2.rewriter1}%
\end{equation}

Furthermore, (\ref{pf.Delta.hrefl.g'v.c1.rewriter2}) and
(\ref{pf.Delta.hrefl.g'v.c1.rewriter3}) hold (this can be proven in the same
way as in Case 1).

Since $\deg v<2p-1$, we see that both $v$ and all the elements $u\in
\widehat{\operatorname*{Rect}\left(  p,p\right)  }$ such that $u \gtrdot v$
belong to $\Delta\left(  p\right)  $. Hence, the elements $u\in
\widehat{\operatorname*{Rect}\left(  p,p\right)  }$ such that $u \gtrdot v$
are precisely the elements $u\in\widehat{\Delta\left(  p\right)  }$ such that
$u \gtrdot v$.

But for every $u\in\widehat{\operatorname*{Rect}\left(  p,p\right)  }$ such
that $u \gtrdot v$, we have $u\in\Delta\left(  p\right)  \cup\left\{
1\right\}  $ (since $u \gtrdot v$, so that $\deg u=\underbrace{\deg
v}_{>p+1>p}+1>p+1$).

Now, every $u\in\widehat{\operatorname*{Rect}\left(  p,p\right)  }$ such that
$u \gtrdot v$ satisfies%
\begin{align*}
g^{\prime}\left(  u\right)   &  =\underbrace{a_{\deg v+1}}%
_{\substack{=1\\\text{(since }\deg v>p+1>p\text{,}\\\text{so that }\deg
v+1>p+1\text{)}}}\cdot\left(  \operatorname*{wing}\left(  R_{\Delta\left(
p\right)  }f\right)  \right)  \left(  u\right)  \ \ \ \ \ \ \ \ \ \ \left(
\text{by (\ref{pf.Delta.hrefl.g'u.ucoversv})}\right) \\
&  =\left(  \operatorname*{wing}\left(  R_{\Delta\left(  p\right)  }f\right)
\right)  \left(  u\right)  =\left(  R_{\Delta\left(  p\right)  }f\right)
\left(  u\right)
\end{align*}
(by the definition of $\operatorname*{wing}$, since $u\in\Delta\left(
p\right)  \cup\left\{  1\right\}  $). Thus,
\begin{equation}
\sum\limits_{\substack{u\in\widehat{\operatorname*{Rect}\left(  p,p\right)
};\\u \gtrdot v}}\dfrac{1}{g^{\prime}\left(  u\right)  }=\sum
\limits_{\substack{u\in\widehat{\operatorname*{Rect}\left(  p,p\right)  };\\u
\gtrdot v}}\dfrac{1}{\left(  R_{\Delta\left(  p\right)  }f\right)  \left(
u\right)  }=\sum\limits_{\substack{u\in\widehat{\Delta\left(  p\right)  };\\u
\gtrdot v}}\dfrac{1}{\left(  R_{\Delta\left(  p\right)  }f\right)  \left(
u\right)  } \label{pf.Delta.hrefl.g'v.c2.rewriter4}%
\end{equation}
(since the elements $u\in\widehat{\operatorname*{Rect}\left(  p,p\right)  }$
such that $u \gtrdot v$ are precisely the elements $u\in\widehat{\Delta\left(
p\right)  }$ such that $u \gtrdot v$).

Recall that we have to prove the equality (\ref{pf.Delta.hrefl.goal}). Upon
substitution of (\ref{pf.Delta.hrefl.g'v.c2.rewriter1}),
(\ref{pf.Delta.hrefl.g'v.c1.rewriter2}),
(\ref{pf.Delta.hrefl.g'v.c1.rewriter3}) and
(\ref{pf.Delta.hrefl.g'v.c2.rewriter4}), this equality
(\ref{pf.Delta.hrefl.goal}) becomes%
\begin{equation}
\left(  R_{\Delta\left(  p\right)  }f\right)  \left(  v\right)  =\dfrac
{1}{f\left(  v\right)  }\cdot\dfrac{\sum\limits_{\substack{u\in\widehat{\Delta
\left(  p\right)  };\\u \lessdot v}}f\left(  u\right)  }{\sum
\limits_{\substack{u\in\widehat{\Delta\left(  p\right)  };\\u\gtrdot v}%
}\dfrac{1}{\left(  R_{\Delta\left(  p\right)  }f\right)  \left(  u\right)  }}.
\label{pf.Delta.hrefl.g'v.c2.transgoal}%
\end{equation}
So (\ref{pf.Delta.hrefl.g'v.c2.transgoal}) is what we need to prove. But
(\ref{pf.Delta.hrefl.g'v.c2.transgoal}) follows directly from applying
Proposition \ref{prop.R.implicit} to $\Delta\left(  p\right)  $ instead of
$P$. Hence, (\ref{pf.Delta.hrefl.g'v.c2.transgoal}) is proven. Since
(\ref{pf.Delta.hrefl.g'v.c2.transgoal}) was obtained by rewriting
(\ref{pf.Delta.hrefl.goal}), this yields that (\ref{pf.Delta.hrefl.goal}) is
proven in Case 2.

\bigskip

\underline{\textbf{Proof of (\ref{pf.Delta.hrefl.goal}) in Case 3.}}

Next, let us consider Case 3. In this case, we have $\deg v=p+1$. Thus,
$v\in\Delta\left(  p\right)  \subseteq\Delta\left(  p\right)  \cup\left\{
1\right\}  $. Also, $\deg v=p+1<2p-1$ (since $p>2$).

The equalities (\ref{pf.Delta.hrefl.g'v.c1.rewriter2}) and
(\ref{pf.Delta.hrefl.g'v.c2.rewriter4}) are valid (this can be proven in the
same way as in Cases 1 and 2, respectively).

Since $\deg v=p+1$, we know that every $u\in\widehat{\operatorname*{Rect}%
\left(  p,p\right)  }$ such that $u \lessdot v$ satisfies $u\in
\operatorname*{Eq}\left(  p\right)  $. Hence, every $u\in
\widehat{\operatorname*{Rect}\left(  p,p\right)  }$ such that $u \lessdot v$
satisfies%
\begin{align*}
f^{\prime}\left(  u\right)   &  =\left(  \operatorname*{wing}f\right)  \left(
u\right)  \ \ \ \ \ \ \ \ \ \ \left(  \text{since }f^{\prime}%
=\operatorname*{wing}f\right) \\
&  =1\ \ \ \ \ \ \ \ \ \ \left(  \text{by the definition of }%
\operatorname*{wing}\text{, since }u\in\operatorname*{Eq}\left(  p\right)
\right)  .
\end{align*}
Thus,%
\begin{equation}
\sum\limits_{\substack{u\in\widehat{\operatorname*{Rect}\left(  p,p\right)
};\\u \lessdot v}}\underbrace{f^{\prime}\left(  u\right)  }_{=1}%
=\sum\limits_{\substack{u\in\widehat{\operatorname*{Rect}\left(  p,p\right)
};\\u \lessdot v}}1=2 \label{pf.Delta.hrefl.g'v.c3.rewriter3}%
\end{equation}
(since there exist exactly two $u\in\widehat{\operatorname*{Rect}\left(
p,p\right)  }$ such that $u \lessdot v$). Finally, the equality
(\ref{pf.Delta.hrefl.g'v}) rewrites as
\begin{equation}
g^{\prime}\left(  v\right)  =\underbrace{a_{\deg v}}_{\substack{=a_{p+1}%
\\\text{(since }\deg v=p+1\text{)}}}\cdot\underbrace{\left(
\operatorname*{wing}\left(  R_{\Delta\left(  p\right)  }f\right)  \right)
\left(  v\right)  }_{\substack{=\left(  R_{\Delta\left(  p\right)  }f\right)
\left(  v\right)  \\\text{(by the definition of }\operatorname*{wing}%
\text{,}\\\text{since }v\in\Delta\left(  p\right)  \cup\left\{  1\right\}
\text{)}}}=a_{p+1}\cdot\left(  R_{\Delta\left(  p\right)  }f\right)  \left(
v\right)  . \label{pf.Delta.hrefl.g'v.c3.rewriter1}%
\end{equation}

Recall that we have to prove the equality (\ref{pf.Delta.hrefl.goal}). Upon
substitution of (\ref{pf.Delta.hrefl.g'v.c3.rewriter1}),
(\ref{pf.Delta.hrefl.g'v.c1.rewriter2}),
(\ref{pf.Delta.hrefl.g'v.c3.rewriter3}) and
(\ref{pf.Delta.hrefl.g'v.c2.rewriter4}), this equality
(\ref{pf.Delta.hrefl.goal}) becomes%
\begin{equation}
a_{p+1}\cdot\left(  R_{\Delta\left(  p\right)  }f\right)  \left(  v\right)
=\dfrac{1}{f\left(  v\right)  }\cdot\dfrac{2}{\sum\limits_{\substack{u\in
\widehat{\Delta\left(  p\right)  };\\u \gtrdot v}}\dfrac{1}{\left(
R_{\Delta\left(  p\right)  }f\right)  \left(  u\right)  }}.
\label{pf.Delta.hrefl.g'v.c3.transgoal}%
\end{equation}
So (\ref{pf.Delta.hrefl.g'v.c3.transgoal}) is what we need to prove.

But (\ref{pf.Delta.hrefl.g'v.c1.1}) holds (this can be shown just as in Case
1). But since $v$ is a minimal element of $\Delta\left(  p\right)  $ (because
$\deg v=p+1$), there is only one element $u\in\widehat{\Delta\left(  p\right)
}$ such that $u \lessdot v$, namely $0$. Thus,
\[
\sum\limits_{\substack{u\in\widehat{\Delta\left(  p\right)  };\\u \lessdot
v}}f\left(  u\right)  =f\left(  0\right)  .
\]
In view of this, (\ref{pf.Delta.hrefl.g'v.c1.1}) becomes%
\[
\left(  R_{\Delta\left(  p\right)  }f\right)  \left(  v\right)  =\dfrac
{1}{f\left(  v\right)  }\cdot\dfrac{f\left(  0\right)  }{\sum
\limits_{\substack{u\in\widehat{\Delta\left(  p\right)  };\\u \gtrdot
v}}\dfrac{1}{\left(  R_{\Delta\left(  p\right)  }f\right)  \left(  u\right)
}}.
\]
Multiplying this with $a_{p+1}$, we obtain%
\[
a_{p+1}\cdot\left(  R_{\Delta\left(  p\right)  }f\right)  \left(  v\right)
=\underbrace{a_{p+1}}_{=\dfrac{2}{f\left(  0\right)  }}\cdot\dfrac{1}{f\left(
v\right)  }\cdot\dfrac{f\left(  0\right)  }{\sum\limits_{\substack{u\in
\widehat{\Delta\left(  p\right)  };\\u \gtrdot v}}\dfrac{1}{\left(
R_{\Delta\left(  p\right)  }f\right)  \left(  u\right)  }}=\dfrac{1}{f\left(
v\right)  }\cdot\dfrac{2}{\sum\limits_{\substack{u\in\widehat{\Delta\left(
p\right)  };\\u \gtrdot v}}\dfrac{1}{\left(  R_{\Delta\left(  p\right)
}f\right)  \left(  u\right)  }}.
\]
But this is precisely (\ref{pf.Delta.hrefl.g'v.c3.transgoal}). Hence,
(\ref{pf.Delta.hrefl.g'v.c3.transgoal}) is proven. Since
(\ref{pf.Delta.hrefl.g'v.c3.transgoal}) was obtained by rewriting
(\ref{pf.Delta.hrefl.goal}), this yields that (\ref{pf.Delta.hrefl.goal}) is
proven in Case 3.

\bigskip

\underline{\textbf{Proof of (\ref{pf.Delta.hrefl.goal}) in Case 4.}}

Next, let us consider Case 4. In this case, we have $\deg v=p$. Thus,
$v\in\operatorname*{Eq}\left(  p\right)  $. The equality
(\ref{pf.Delta.hrefl.g'v}) rewrites as
\begin{equation}
g^{\prime}\left(  v\right)  =\underbrace{a_{\deg v}}_{\substack{=a_{p}%
\\\text{(since }\deg v=p\text{)}}}\cdot\underbrace{\left(
\operatorname*{wing}\left(  R_{\Delta\left(  p\right)  }f\right)  \right)
\left(  v\right)  }_{\substack{=1\\\text{(by the definition of }%
\operatorname*{wing}\text{,}\\\text{since }v\in\operatorname*{Eq}\left(
p\right)  \text{)}}}=a_{p}=\dfrac{2}{f\left(  0\right)  }.
\label{pf.Delta.hrefl.g'v.c4.rewriter1}%
\end{equation}

But we have $v\in\operatorname*{Eq}\left(  p\right)  $. Thus, every element
$u\in\widehat{\operatorname*{Rect}\left(  p,p\right)  }$ such that $u\gtrdot
v$ satisfies $u\in\Delta\left(  p\right)  \subseteq\Delta\left(  p\right)
\cup\left\{  1\right\}  $. Hence, every element $u\in
\widehat{\operatorname*{Rect}\left(  p,p\right)  }$ such that $u\gtrdot v$
satisfies%
\begin{align*}
g^{\prime}\left(  u\right)   &  =\underbrace{a_{\deg v+1}}_{\substack{=a_{p+1}%
\\\text{(since }\deg v=p\text{)}}}\cdot\underbrace{\left(
\operatorname*{wing}\left(  R_{\Delta\left(  p\right)  }f\right)  \right)
\left(  u\right)  }_{\substack{=\left(  R_{\Delta\left(  p\right)  }f\right)
\left(  u\right)  \\\text{(by the definition of }\operatorname*{wing}%
\text{,}\\\text{since }u\in\Delta\left(  p\right)  \cup\left\{  1\right\}
\text{)}}}\ \ \ \ \ \ \ \ \ \ \left(  \text{by
(\ref{pf.Delta.hrefl.g'u.ucoversv})}\right) \\
&  =\underbrace{a_{p+1}}_{=\dfrac{2}{f\left(  0\right)  }}\cdot\left(
R_{\Delta\left(  p\right)  }f\right)  \left(  u\right)  =\dfrac{2}{f\left(
0\right)  }\left(  R_{\Delta\left(  p\right)  }f\right)  \left(  u\right)  .
\end{align*}
Thus,%
\begin{equation}
\sum\limits_{\substack{u\in\widehat{\operatorname*{Rect}\left(  p,p\right)
};\\u\gtrdot v}}\dfrac{1}{g^{\prime}\left(  u\right)  }=\sum
\limits_{\substack{u\in\widehat{\operatorname*{Rect}\left(  p,p\right)
};\\u\gtrdot v}}\dfrac{1}{\dfrac{2}{f\left(  0\right)  }\left(  R_{\Delta
\left(  p\right)  }f\right)  \left(  u\right)  }=\dfrac{f\left(  0\right)
}{2}\sum\limits_{\substack{u\in\widehat{\operatorname*{Rect}\left(
p,p\right)  };\\u\gtrdot v}}\dfrac{1}{\left(  R_{\Delta\left(  p\right)
}f\right)  \left(  u\right)  }. \label{pf.Delta.hrefl.g'v.c4.rewriter4}%
\end{equation}

On the other hand, recall that $v\in\operatorname*{Eq}\left(  p\right)  $.
Hence, every element $u\in\widehat{\operatorname*{Rect}\left(  p,p\right)  }$
such that $u\lessdot v$ satisfies $u\in\nabla\left(  p\right)  \subseteq
\nabla\left(  p\right)  \cup\left\{  0\right\}  $. Hence, every element
$u\in\widehat{\operatorname*{Rect}\left(  p,p\right)  }$ such that $u\lessdot
v$ satisfies
\begin{align*}
\underbrace{f^{\prime}}_{=\operatorname*{wing}f}\left(  u\right)   &  =\left(
\operatorname*{wing}f\right)  \left(  u\right) \\
&  =\dfrac{1}{\left(  R_{\Delta\left(  p\right)  }^{p-\deg u}f\right)  \left(
\operatorname*{hrefl}u\right)  }\ \ \ \ \ \ \ \ \ \ \left(  \text{by the
definition of }\operatorname*{wing}\text{, since }u\in\nabla\left(  p\right)
\cup\left\{  0\right\}  \right) \\
&  =\dfrac{1}{\left(  R_{\Delta\left(  p\right)  }f\right)  \left(
\operatorname*{hrefl}u\right)  }%
\end{align*}
(since $u\lessdot v$, so that $\deg u=\underbrace{\deg v}_{=p}-1=p-1$, so that
$p-\deg u=1$, so that $R_{\Delta\left(  p\right)  }^{p-\deg u}=R_{\Delta
\left(  p\right)  }^{1}=R_{\Delta\left(  p\right)  }$). Thus,%
\begin{align}
\sum\limits_{\substack{u\in\widehat{\operatorname*{Rect}\left(  p,p\right)
};\\u\lessdot v}}f^{\prime}\left(  u\right)   &  =\sum\limits_{\substack{u\in
\widehat{\operatorname*{Rect}\left(  p,p\right)  };\\u\lessdot v}}\dfrac
{1}{\left(  R_{\Delta\left(  p\right)  }f\right)  \left(
\operatorname*{hrefl}u\right)  }=\sum\limits_{\substack{u\in
\widehat{\operatorname*{Rect}\left(  p,p\right)  };\\\operatorname*{hrefl}%
u\lessdot v}}\underbrace{\dfrac{1}{\left(  R_{\Delta\left(  p\right)
}f\right)  \left(  \operatorname*{hrefl}\left(  \operatorname*{hrefl}u\right)
\right)  }}_{\substack{=\dfrac{1}{\left(  R_{\Delta\left(  p\right)
}f\right)  \left(  u\right)  }\\\text{(since }\operatorname*{hrefl}\left(
\operatorname*{hrefl}u\right)  =u\\\text{(because }\operatorname*{hrefl}\text{
is an involution))}}}\nonumber\\
&  \ \ \ \ \ \ \ \ \ \ \left(
\begin{array}
[c]{c}%
\text{here, we have substituted }\operatorname*{hrefl}u\text{ for }u\text{ in
the sum, since}\\
\operatorname*{hrefl}\text{ is a bijection}%
\end{array}
\right) \nonumber\\
&  =\sum\limits_{\substack{u\in\widehat{\operatorname*{Rect}\left(
p,p\right)  };\\\operatorname*{hrefl}u\lessdot v}}\dfrac{1}{\left(
R_{\Delta\left(  p\right)  }f\right)  \left(  u\right)  }.
\label{pf.Delta.hrefl.g'v.c4.rewriter3.pre}%
\end{align}
But let us notice that for every $u\in\widehat{\operatorname*{Rect}\left(
p,p\right)  }$, we have the following equivalence of statements:%
\begin{equation}
\left(  \operatorname*{hrefl}u\lessdot v\right)  \Longleftrightarrow\left(
u\gtrdot v\right)  \label{pf.Delta.hrefl.g'v.c4.rewriter3equiv}%
\end{equation}
\footnote{\textit{Proof of (\ref{pf.Delta.hrefl.g'v.c4.rewriter3equiv}):} Let
$u\in\widehat{\operatorname*{Rect}\left(  p,p\right)  }$. We know that
$\operatorname*{hrefl}$ is an antiautomorphism of the poset
$\widehat{\operatorname*{Rect}\left(  p,p\right)  }$. Hence, we have $u\gtrdot
v$ if and only if $\operatorname*{hrefl}u\lessdot\operatorname*{hrefl}v$.
Since $\operatorname*{hrefl}v=v$ (because $v\in\operatorname*{Eq}\left(
p\right)  $), this rewrites as follows: We have $u\gtrdot v$ if and only if
$\operatorname*{hrefl}u\lessdot v$. In other words, we have the equivalence of
statements $\left(  u\gtrdot v\right)  \Longleftrightarrow\left(
\operatorname*{hrefl}u\lessdot v\right)  $. This proves
(\ref{pf.Delta.hrefl.g'v.c4.rewriter3equiv}).}. Hence, we can replace the
summation sign \textquotedblleft$\sum\limits_{\substack{u\in
\widehat{\operatorname*{Rect}\left(  p,p\right)  };\\\operatorname*{hrefl}%
u\lessdot v}}$\textquotedblright\ by a \textquotedblleft$\sum
\limits_{\substack{u\in\widehat{\operatorname*{Rect}\left(  p,p\right)
};\\u\gtrdot v}}$\textquotedblright\ in
(\ref{pf.Delta.hrefl.g'v.c4.rewriter3.pre}). As a result, we obtain%
\begin{equation}
\sum\limits_{\substack{u\in\widehat{\operatorname*{Rect}\left(  p,p\right)
};\\u\lessdot v}}f^{\prime}\left(  u\right)  =\sum\limits_{\substack{u\in
\widehat{\operatorname*{Rect}\left(  p,p\right)  };\\u\gtrdot v}}\dfrac
{1}{\left(  R_{\Delta\left(  p\right)  }f\right)  \left(  u\right)  }.
\label{pf.Delta.hrefl.g'v.c4.rewriter3}%
\end{equation}

Since $f^{\prime}=\operatorname*{wing}f$, we have%
\begin{equation}
f^{\prime}\left(  v\right)  =\left(  \operatorname*{wing}f\right)  \left(
v\right)  =1 \label{pf.Delta.hrefl.g'v.c4.rewriter2}%
\end{equation}
(by the definition of $\operatorname*{wing}$, since $v\in\operatorname*{Eq}%
\left(  p\right)  $).

Recall that we have to prove the equality (\ref{pf.Delta.hrefl.goal}). Upon
substitution of (\ref{pf.Delta.hrefl.g'v.c4.rewriter1}),
(\ref{pf.Delta.hrefl.g'v.c4.rewriter2}),
(\ref{pf.Delta.hrefl.g'v.c4.rewriter3}) and
(\ref{pf.Delta.hrefl.g'v.c4.rewriter4}), this equality
(\ref{pf.Delta.hrefl.goal}) becomes%
\[
\dfrac{2}{f\left(  0\right)  }=\dfrac{1}{1}\cdot\dfrac{\sum
\limits_{\substack{u\in\widehat{\operatorname*{Rect}\left(  p,p\right)
};\\u\gtrdot v}}\dfrac{1}{\left(  R_{\Delta\left(  p\right)  }f\right)
\left(  u\right)  }}{\dfrac{f\left(  0\right)  }{2}\sum\limits_{\substack{u\in
\widehat{\operatorname*{Rect}\left(  p,p\right)  };\\u\gtrdot v}}\dfrac
{1}{\left(  R_{\Delta\left(  p\right)  }f\right)  \left(  u\right)  }}.
\]
But this equality is trivial (because after cancelling, it becomes $\dfrac
{2}{f\left(  0\right)  }=\dfrac{2}{f\left(  0\right)  }$). Hence,
(\ref{pf.Delta.hrefl.goal}) is proven in Case 4.

\bigskip

\underline{\textbf{Proof of (\ref{pf.Delta.hrefl.goal}) in Case 5.}}

Next, let us consider Case 5. In this case, we have $\deg v=p-1$. Thus,
$v\in\nabla\left(  p\right)  \subseteq\nabla\left(  p\right)  \cup\left\{
0\right\}  $ and $\operatorname*{hrefl}v\in\Delta\left(  p\right)
\subseteq\Delta\left(  p\right)  \cup\left\{  1\right\}  $.

Let $w=\operatorname*{hrefl}v$. Then, $w=\operatorname*{hrefl}v\in
\Delta\left(  p\right)  \subseteq\Delta\left(  p\right)  \cup\left\{
1\right\}  $.

Since $\deg v=p-1<p$, the equality (\ref{pf.Delta.hrefl.g'v}) rewrites as
\begin{align}
g^{\prime}\left(  v\right)   &  =\underbrace{a_{\deg v}}%
_{\substack{=1\\\text{(since }\deg v<p\text{)}}}\cdot\underbrace{\left(
\operatorname*{wing}\left(  R_{\Delta\left(  p\right)  }f\right)  \right)
\left(  v\right)  }_{\substack{=\dfrac{1}{\left(  R_{\Delta\left(  p\right)
}^{p-\deg v}\left(  R_{\Delta\left(  p\right)  }f\right)  \right)  \left(
\operatorname*{hrefl}v\right)  }\\\text{(by the definition of }%
\operatorname*{wing}\text{,}\\\text{since }v\in\nabla\left(  p\right)
\cup\left\{  0\right\}  \text{)}}}\nonumber\\
&  =\dfrac{1}{\left(  R_{\Delta\left(  p\right)  }^{p-\deg v}\left(
R_{\Delta\left(  p\right)  }f\right)  \right)  \left(  \operatorname*{hrefl}%
v\right)  }=\dfrac{1}{\left(  R_{\Delta\left(  p\right)  }^{p-\deg
v+1}f\right)  \left(  w\right)  } \label{pf.Delta.hrefl.g'v.c5.rewriter1}%
\end{align}
(since $R_{\Delta\left(  p\right)  }^{p-\deg v}\left(  R_{\Delta\left(
p\right)  }f\right)  =R_{\Delta\left(  p\right)  }^{p-\deg v+1}f$ and
$\operatorname*{hrefl}v=w$).

Also, since $f^{\prime}=\operatorname*{wing}f$, we have $f^{\prime}\left(
v\right)  =\left(  \operatorname*{wing}f\right)  \left(  v\right)  =\dfrac
{1}{\left(  R_{\Delta\left(  p\right)  }^{p-\deg v}f\right)  \left(
\operatorname*{hrefl}v\right)  }$ (by the definition of $\operatorname*{wing}%
$, since $v\in\nabla\left(  p\right)  \cup\left\{  0\right\}  $). Since
$\operatorname*{hrefl}v=w$, this rewrites as%
\begin{equation}
f^{\prime}\left(  v\right)  =\dfrac{1}{\left(  R_{\Delta\left(  p\right)
}^{p-\deg v}f\right)  \left(  w\right)  }.
\label{pf.Delta.hrefl.g'v.c5.rewriter2}%
\end{equation}

Since $\deg v=p-1$, we know that every $u\in\widehat{\operatorname*{Rect}%
\left(  p,p\right)  }$ such that $u \gtrdot v$ satisfies $u\in
\operatorname*{Eq}\left(  p\right)  $. Now, every $u\in
\widehat{\operatorname*{Rect}\left(  p,p\right)  }$ such that $u \gtrdot v$
satisfies%
\begin{align*}
g^{\prime}\left(  u\right)   &  =\underbrace{a_{\deg v+1}}_{\substack{=a_{p}%
\\\text{(since }\deg v=p-1\text{,}\\\text{so that }\deg v+1=p\text{)}}%
}\cdot\underbrace{\left(  \operatorname*{wing}\left(  R_{\Delta\left(
p\right)  }f\right)  \right)  \left(  u\right)  }_{\substack{=1\\\text{(by the
definition of }\operatorname*{wing}\text{,}\\\text{since }u\in
\operatorname*{Eq}\left(  p\right)  \text{)}}}\ \ \ \ \ \ \ \ \ \ \left(
\text{by (\ref{pf.Delta.hrefl.g'u.ucoversv})}\right) \\
&  =a_{p}=\dfrac{2}{f\left(  0\right)  }.
\end{align*}
Thus,%
\[
\sum\limits_{\substack{u\in\widehat{\operatorname*{Rect}\left(  p,p\right)
};\\u \gtrdot v}}\dfrac{1}{g^{\prime}\left(  u\right)  }=\sum
\limits_{\substack{u\in\widehat{\operatorname*{Rect}\left(  p,p\right)  };\\u
\gtrdot v}}\dfrac{1}{\left(  \dfrac{2}{f\left(  0\right)  }\right)  }%
=2\cdot\dfrac{1}{\left(  \dfrac{2}{f\left(  0\right)  }\right)  }%
\]
(since there exist exactly two $u\in\widehat{\operatorname*{Rect}\left(
p,p\right)  }$ such that $u \gtrdot v$). Thus,%
\begin{equation}
\sum\limits_{\substack{u\in\widehat{\operatorname*{Rect}\left(  p,p\right)
};\\u \gtrdot v}}\dfrac{1}{g^{\prime}\left(  u\right)  }=2\cdot\dfrac
{1}{\left(  \dfrac{2}{f\left(  0\right)  }\right)  }=f\left(  0\right)  .
\label{pf.Delta.hrefl.g'v.c5.rewriter4}%
\end{equation}

Furthermore, again because of $\deg v=p-1$, we have $p-1\geq\deg v>1$. Thus,
$v$ and the elements $u\in\widehat{\operatorname*{Rect}\left(  p,p\right)  }$
such that $u\lessdot v$ all belong to $\nabla\left(  p\right)
=\operatorname*{hrefl}\left(  \Delta\left(  p\right)  \right)  $. Hence, the
elements $u\in\widehat{\operatorname*{Rect}\left(  p,p\right)  }$ such that
$u\lessdot v$ are precisely the images under $\operatorname*{hrefl}$ of the
elements $u\in\widehat{\Delta\left(  p\right)  }$ such that $u\gtrdot
\operatorname*{hrefl}v$. \ \ \ \ \footnote{\textit{Proof.} We know that
the elements $u\in\widehat{\operatorname*{Rect}\left(  p,p\right)  }$
such that $u\lessdot v$ all belong to $\nabla\left(  p\right)$, and thus
all belong to $\widehat{\nabla\left(  p\right)}$. Hence,
\begin{align*}
&  \left\{ u \in \widehat{\operatorname*{Rect}\left(p, p\right)} \mid
           u \lessdot v \right\} \\
&= \left\{ u \in \widehat{\nabla\left(  p\right)} \mid
           u \lessdot v \right\} \\
&= \left\{ u \in \widehat{\nabla\left(  p\right)} \mid
           \operatorname*{hrefl}\left(u\right) \gtrdot
           \operatorname*{hrefl}\left(v\right) \right\}
\ \ \ \ \ \ \ \ \ \ \left(\text{because of (\ref{pf.Delta.hrefl.d.equiv})}\right) \\
&= \operatorname{hrefl}^{-1}\left(
     \left\{ u \in \widehat{\Delta\left(  p\right)} \mid
             u \gtrdot \operatorname*{hrefl}\left(v\right) \right\}
     \right)
\ \ \ \ \ \ \ \ \ \ \left(\text{since } \nabla\left(p\right)
                          = \operatorname{hrefl}^{-1}
                            \left(\Delta\left(p\right)\right) \right) \\
&= \operatorname*{hrefl}\left(
     \left\{ u \in \widehat{\Delta\left(  p\right)} \mid
             u \gtrdot \operatorname*{hrefl}\left(v\right) \right\}
     \right)
\end{align*}
(since $\operatorname{hrefl}^{-1} = \operatorname*{hrefl}$ (because
$\operatorname*{hrefl}$ is an involution)), qed.}

Notice that the elements $u\in\widehat{\operatorname*{Rect}\left(  p,p\right)
}$ such that $u\gtrdot w$ are precisely the elements $u\in\widehat{\Delta
\left(  p\right)  }$ such that $u\gtrdot w$. (This is because
$w=\operatorname*{hrefl}\underbrace{u}_{\in\nabla\left(  p\right)  }%
\in\operatorname*{hrefl}\left(  \nabla\left(  p\right)  \right)
=\Delta\left(  p\right)  $.)

But for every $u\in\widehat{\operatorname*{Rect}\left(  p,p\right)  }$ such
that $u\lessdot v$, we have%
\begin{align*}
\underbrace{f^{\prime}}_{=\operatorname*{wing}f}\left(  u\right)   &  =\left(
\operatorname*{wing}f\right)  \left(  u\right)  =\dfrac{1}{\left(
R_{\Delta\left(  p\right)  }^{p-\deg u}f\right)  \left(  \operatorname*{hrefl}%
u\right)  }\\
&  \ \ \ \ \ \ \ \ \ \ \left(  \text{by the definition of }%
\operatorname*{wing}\text{, since }u\in\nabla\left(  p\right)  \subseteq
\nabla\left(  p\right)  \cup\left\{  0\right\}  \right) \\
&  =\dfrac{1}{\left(  R_{\Delta\left(  p\right)  }^{p-\deg v+1}f\right)
\left(  \operatorname*{hrefl}u\right)  }\\
&  \ \ \ \ \ \ \ \ \ \ \left(
\begin{array}
[c]{c}%
\text{since }u\lessdot v\text{, so that }\deg u=\deg v-1\text{, so that}\\
p-\deg u=p-\left(  \deg v-1\right)  =p-\deg v+1
\end{array}
\right)  .
\end{align*}
Hence,%
\begin{align}
\sum\limits_{\substack{u\in\widehat{\operatorname*{Rect}\left(  p,p\right)
};\\u\lessdot v}}f^{\prime}\left(  u\right)   &  =\sum\limits_{\substack{u\in
\widehat{\operatorname*{Rect}\left(  p,p\right)  };\\u\lessdot v}}\dfrac
{1}{\left(  R_{\Delta\left(  p\right)  }^{p-\deg v+1}f\right)  \left(
\operatorname*{hrefl}u\right)  }=\sum\limits_{\substack{u\in
\widehat{\operatorname*{Rect}\left(  p,p\right)  };\\\operatorname*{hrefl}%
u\gtrdot\operatorname*{hrefl}v}}\dfrac{1}{\left(  R_{\Delta\left(  p\right)
}^{p-\deg v+1}f\right)  \left(  \operatorname*{hrefl}u\right)  }\nonumber\\
&  \ \ \ \ \ \ \ \ \ \ \left(
\begin{array}
[c]{c}%
\text{here, we replaced the summation sign \textquotedblleft}\sum
\limits_{\substack{u\in\widehat{\operatorname*{Rect}\left(  p,p\right)
};\\u\lessdot v}}\text{\textquotedblright}\\
\text{by a \textquotedblleft}\sum\limits_{\substack{u\in
\widehat{\operatorname*{Rect}\left(  p,p\right)  };\\\operatorname*{hrefl}%
u\gtrdot\operatorname*{hrefl}v}}\text{\textquotedblright\ (because of the
equivalence (\ref{pf.Delta.hrefl.d.equiv}))}%
\end{array}
\right) \nonumber\\
&  =\sum\limits_{\substack{u\in\widehat{\operatorname*{Rect}\left(
p,p\right)  };\\u\gtrdot\operatorname*{hrefl}v}}\dfrac{1}{\left(
R_{\Delta\left(  p\right)  }^{p-\deg v+1}f\right)  \left(  u\right)
}\nonumber\\
&  \ \ \ \ \ \ \ \ \ \ \left(
\begin{array}
[c]{c}%
\text{here, we have substituted }u\text{ for }\operatorname*{hrefl}u\text{ in
the sum}\\
\text{(since }\operatorname*{hrefl}\text{ is a bijection)}%
\end{array}
\right) \nonumber\\
&  =\sum\limits_{\substack{u\in\widehat{\operatorname*{Rect}\left(
p,p\right)  };\\u\gtrdot w}}\dfrac{1}{\left(  R_{\Delta\left(  p\right)
}^{p-\deg v+1}f\right)  \left(  u\right)  }\ \ \ \ \ \ \ \ \ \ \left(
\text{since }\operatorname*{hrefl}v=w\right) \nonumber\\
&  =\sum\limits_{\substack{u\in\widehat{\Delta\left(  p\right)  };\\u\gtrdot
w}}\dfrac{1}{\left(  R_{\Delta\left(  p\right)  }^{p-\deg v+1}f\right)
\left(  u\right)  } \label{pf.Delta.hrefl.g'v.c5.rewriter3}%
\end{align}
(since the elements $u\in\widehat{\operatorname*{Rect}\left(  p,p\right)  }$
such that $u\gtrdot w$ are precisely the elements $u\in\widehat{\Delta\left(
p\right)  }$ such that $u\gtrdot w$).

Recall that we have to prove the equality (\ref{pf.Delta.hrefl.goal}). Upon
substitution of (\ref{pf.Delta.hrefl.g'v.c5.rewriter1}),
(\ref{pf.Delta.hrefl.g'v.c5.rewriter2}),
(\ref{pf.Delta.hrefl.g'v.c5.rewriter3}) and
(\ref{pf.Delta.hrefl.g'v.c5.rewriter4}), this equality
(\ref{pf.Delta.hrefl.goal}) becomes%
\begin{equation}
\dfrac{1}{\left(  R_{\Delta\left(  p\right)  }^{p-\deg v+1}f\right)  \left(
w\right)  }=\dfrac{1}{\left(  \dfrac{1}{\left(  R_{\Delta\left(  p\right)
}^{p-\deg v}f\right)  \left(  w\right)  }\right)  }\cdot\dfrac{\sum
\limits_{\substack{u\in\widehat{\Delta\left(  p\right)  };\\u \gtrdot
w}}\dfrac{1}{\left(  R_{\Delta\left(  p\right)  }^{p-\deg v+1}f\right)
\left(  u\right)  }}{f\left(  0\right)  }.
\label{pf.Delta.hrefl.g'v.c5.transgoal}%
\end{equation}
So (\ref{pf.Delta.hrefl.g'v.c5.transgoal}) is what we need to prove.

But applying Proposition \ref{prop.R.implicit} to $\Delta\left(  p\right)  $,
$R_{\Delta\left(  p\right)  }^{p-\deg v}f$ and $w$ instead of $P$, $f$ and
$v$, we find%
\[
\left(  R_{\Delta\left(  p\right)  }\left(  R_{\Delta\left(  p\right)
}^{p-\deg v}f\right)  \right)  \left(  w\right)  =\dfrac{1}{\left(
R_{\Delta\left(  p\right)  }^{p-\deg v}f\right)  \left(  w\right)  }%
\cdot\dfrac{\sum\limits_{\substack{u\in\widehat{\Delta\left(  p\right)  };\\u
\lessdot w}}\left(  R_{\Delta\left(  p\right)  }^{p-\deg v}f\right)  \left(
u\right)  }{\sum\limits_{\substack{u\in\widehat{\Delta\left(  p\right)  };\\u
\gtrdot w}}\dfrac{1}{\left(  R_{\Delta\left(  p\right)  }\left(
R_{\Delta\left(  p\right)  }^{p-\deg v}f\right)  \right)  \left(  u\right)  }%
}.
\]
Since $R_{\Delta\left(  p\right)  }\left(  R_{\Delta\left(  p\right)
}^{p-\deg v}f\right)  =\underbrace{\left(  R_{\Delta\left(  p\right)  }\circ
R_{\Delta\left(  p\right)  }^{p-\deg v}\right)  }_{=R_{\Delta\left(  p\right)
}^{p-\deg v+1}}f=R_{\Delta\left(  p\right)  }^{p-\deg v+1}f$, this rewrites as%
\[
\left(  R_{\Delta\left(  p\right)  }^{p-\deg v+1}f\right)  \left(  w\right)
=\dfrac{1}{\left(  R_{\Delta\left(  p\right)  }^{p-\deg v}f\right)  \left(
w\right)  }\cdot\dfrac{\sum\limits_{\substack{u\in\widehat{\Delta\left(
p\right)  };\\u \lessdot w}}\left(  R_{\Delta\left(  p\right)  }^{p-\deg
v}f\right)  \left(  u\right)  }{\sum\limits_{\substack{u\in\widehat{\Delta
\left(  p\right)  };\\u \gtrdot w}}\dfrac{1}{\left(  R_{\Delta\left(
p\right)  }^{p-\deg v+1}f\right)  \left(  u\right)  }}.
\]
Taking the multiplicative inverses of both sides of this equality, we obtain%
\begin{equation}
\dfrac{1}{\left(  R_{\Delta\left(  p\right)  }^{p-\deg v+1}f\right)  \left(
w\right)  }=\left(  R_{\Delta\left(  p\right)  }^{p-\deg v}f\right)  \left(
w\right)  \cdot\dfrac{\sum\limits_{\substack{u\in\widehat{\Delta\left(
p\right)  };\\u \gtrdot w}}\dfrac{1}{\left(  R_{\Delta\left(  p\right)
}^{p-\deg v+1}f\right)  \left(  u\right)  }}{\sum\limits_{\substack{u\in
\widehat{\Delta\left(  p\right)  };\\u \lessdot w}}\left(  R_{\Delta\left(
p\right)  }^{p-\deg v}f\right)  \left(  u\right)  }.
\label{pf.Delta.hrefl.g'v.c5.1}%
\end{equation}
But since $w=\operatorname*{hrefl}\left(  v\right)  $ is a minimal element of
$\Delta\left(  p\right)  $ (because $\deg v=p-1$), there exists only one
element $u\in\widehat{\Delta\left(  p\right)  }$ such that $u \lessdot w$,
namely $0$. Hence,%
\[
\sum\limits_{\substack{u\in\widehat{\Delta\left(  p\right)  };\\u \lessdot
w}}\left(  R_{\Delta\left(  p\right)  }^{p-\deg v}f\right)  \left(  u\right)
=\left(  R_{\Delta\left(  p\right)  }^{p-\deg v}f\right)  \left(  0\right)
=f\left(  0\right)
\]
(by Corollary \ref{cor.R.implicit.01}). Thus, (\ref{pf.Delta.hrefl.g'v.c5.1})
simplifies to
\begin{align*}
\dfrac{1}{\left(  R_{\Delta\left(  p\right)  }^{p-\deg v+1}f\right)  \left(
w\right)  }  &  =\left(  R_{\Delta\left(  p\right)  }^{p-\deg v}f\right)
\left(  w\right)  \cdot\dfrac{\sum\limits_{\substack{u\in\widehat{\Delta
\left(  p\right)  };\\u \gtrdot w}}\dfrac{1}{\left(  R_{\Delta\left(
p\right)  }^{p-\deg v+1}f\right)  \left(  u\right)  }}{f\left(  0\right)  }\\
&  =\dfrac{1}{\left(  \dfrac{1}{\left(  R_{\Delta\left(  p\right)  }^{p-\deg
v}f\right)  \left(  w\right)  }\right)  }\cdot\dfrac{\sum
\limits_{\substack{u\in\widehat{\Delta\left(  p\right)  };\\u \gtrdot
w}}\dfrac{1}{\left(  R_{\Delta\left(  p\right)  }^{p-\deg v+1}f\right)
\left(  u\right)  }}{f\left(  0\right)  }.
\end{align*}
But this is precisely (\ref{pf.Delta.hrefl.g'v.c5.transgoal}). Hence,
(\ref{pf.Delta.hrefl.g'v.c5.transgoal}) is proven. Since
(\ref{pf.Delta.hrefl.g'v.c5.transgoal}) was obtained by rewriting
(\ref{pf.Delta.hrefl.goal}), this yields that (\ref{pf.Delta.hrefl.goal}) is
proven in Case 5.

\bigskip

\underline{\textbf{Proof of (\ref{pf.Delta.hrefl.goal}) in Case 6.}}

Let us now consider Case 6. In this case, we have $1<\deg v<p-1$. Thus,
$v\in\nabla\left(  p\right)  \subseteq\nabla\left(  p\right)  \cup\left\{
0\right\}  $ and $\operatorname*{hrefl}v\in\Delta\left(  p\right)
\subseteq\Delta\left(  p\right)  \cup\left\{  1\right\}  $.

Let $w=\operatorname*{hrefl}v$. Then, $w=\operatorname*{hrefl}v\in
\Delta\left(  p\right)  \subseteq\Delta\left(  p\right)  \cup\left\{
1\right\}  $.

The equalities (\ref{pf.Delta.hrefl.g'v.c5.rewriter1}),
(\ref{pf.Delta.hrefl.g'v.c5.rewriter2}) and
(\ref{pf.Delta.hrefl.g'v.c5.rewriter3}) hold (their proof proceeds in the same
way as in Case 5).

Since $1\leq\deg v<p-1$, we know that $v$ and the elements $u\in
\widehat{\operatorname*{Rect}\left(  p,p\right)  }$ such that $u\gtrdot v$ all
belong to $\nabla\left(  p\right)  =\operatorname*{hrefl}\left(  \Delta\left(
p\right)  \right)  $. Hence, the elements $u\in\widehat{\operatorname*{Rect}%
\left(  p,p\right)  }$ such that $u\gtrdot v$ are precisely the images under
$\operatorname*{hrefl}$ of the elements $u\in\widehat{\Delta\left(  p\right)
}$ such that $u\lessdot\operatorname*{hrefl}v$.

Notice that the elements $u\in\widehat{\operatorname*{Rect}\left(  p,p\right)
}$ such that $u\lessdot w$ are precisely the elements $u\in\widehat{\Delta
\left(  p\right)  }$ such that $u\lessdot w$.\ \ \ \ \footnote{\textit{Proof.}
We have $v \in \operatorname*{Rect}\left(  p,p\right)
\subseteq \widehat{\operatorname*{Rect}\left(  p,p\right)  }$. Hence,
(\ref{pf.Delta.hrefl.hrefl-inverts-deg}) yields
$\deg \left(\operatorname{hrefl} v\right) = 2p - \deg v$. Now,
$\deg \underbrace{w}_{= \operatorname{hrefl} v}
= \deg \left(\operatorname{hrefl} v\right) = 2p - \underbrace{\deg v}_{< p-1}
> 2p - \left(p-1\right) = p+1$.
Hence, the elements $u\in\widehat{\operatorname*{Rect}\left(  p,p\right)  }$
such that $u\lessdot w$ all belong to $\Delta\left(p\right)$ and thus also
to $\widehat{\Delta\left(p\right)}$.
Therefore, the elements $u\in\widehat{\operatorname*{Rect}\left(  p,p\right)  }$
such that $u\lessdot w$ are precisely the elements $u\in\widehat{\Delta\left(
p\right)  }$ such that $u\lessdot w$, qed.}

But for every $u\in\widehat{\operatorname*{Rect}\left(  p,p\right)  }$ such
that $u\gtrdot v$, we have%
\begin{align*}
g^{\prime}\left(  u\right)   &  =\underbrace{a_{\deg v+1}}%
_{\substack{=1\\\text{(since }\deg v<p-1\text{,}\\\text{so that }\deg
v+1<p\text{)}}}\cdot\left(  \operatorname*{wing}\left(  R_{\Delta\left(
p\right)  }f\right)  \right)  \left(  u\right)  \ \ \ \ \ \ \ \ \ \ \left(
\text{by (\ref{pf.Delta.hrefl.g'u.ucoversv})}\right) \\
&  =\left(  \operatorname*{wing}\left(  R_{\Delta\left(  p\right)  }f\right)
\right)  \left(  u\right)  =\dfrac{1}{\left(  R_{\Delta\left(  p\right)
}^{p-\deg u}\left(  R_{\Delta\left(  p\right)  }f\right)  \right)  \left(
\operatorname*{hrefl}u\right)  }\\
&  \ \ \ \ \ \ \ \ \ \ \left(  \text{by the definition of }%
\operatorname*{wing}\text{, since }u\in\nabla\left(  p\right)  \subseteq
\nabla\left(  p\right)  \cup\left\{  0\right\}  \right) \\
&  =\dfrac{1}{\left(  R_{\Delta\left(  p\right)  }^{p-\deg v}f\right)  \left(
\operatorname*{hrefl}u\right)  }%
\end{align*}
(since $u\gtrdot v$, so that $\deg u=\deg v+1$, so that $p-\deg u=p-\left(
\deg v+1\right)  =p-\deg v-1$, so that \newline$R_{\Delta\left(  p\right)
}^{p-\deg u}\left(  R_{\Delta\left(  p\right)  }f\right)  =R_{\Delta\left(
p\right)  }^{p-\deg v-1}\left(  R_{\Delta\left(  p\right)  }f\right)
=\underbrace{\left(  R_{\Delta\left(  p\right)  }^{p-\deg v-1}\circ
R_{\Delta\left(  p\right)  }\right)  }_{=R_{\Delta\left(  p\right)  }^{p-\deg
v}}f=R_{\Delta\left(  p\right)  }^{p-\deg v}f$). Hence, for every
$u\in\widehat{\operatorname*{Rect}\left(  p,p\right)  }$ such that $u\gtrdot
v$, we have $\dfrac{1}{g^{\prime}\left(  u\right)  }=\left(  R_{\Delta\left(
p\right)  }^{p-\deg v}f\right)  \left(  \operatorname*{hrefl}u\right)  $.
Thus,%
\begin{align}
\sum\limits_{\substack{u\in\widehat{\operatorname*{Rect}\left(  p,p\right)
};\\u\gtrdot v}}\dfrac{1}{g^{\prime}\left(  u\right)  }  &  =\sum
\limits_{\substack{u\in\widehat{\operatorname*{Rect}\left(  p,p\right)
};\\u\gtrdot v}}\left(  R_{\Delta\left(  p\right)  }^{p-\deg v}f\right)
\left(  \operatorname*{hrefl}u\right)  =\sum\limits_{\substack{u\in
\widehat{\operatorname*{Rect}\left(  p,p\right)  };\\\operatorname*{hrefl}%
u\lessdot\operatorname*{hrefl}v}}\left(  R_{\Delta\left(  p\right)  }^{p-\deg
v}f\right)  \left(  \operatorname*{hrefl}u\right) \nonumber\\
&  \ \ \ \ \ \ \ \ \ \ \left(
\begin{array}
[c]{c}%
\text{here, we replaced the summation sign \textquotedblleft}\sum
\limits_{\substack{u\in\widehat{\operatorname*{Rect}\left(  p,p\right)
};\\u\gtrdot v}}\text{\textquotedblright}\\
\text{by a \textquotedblleft}\sum\limits_{\substack{u\in
\widehat{\operatorname*{Rect}\left(  p,p\right)  };\\\operatorname*{hrefl}%
u\lessdot\operatorname*{hrefl}v}}\text{\textquotedblright\ (because of the
equivalence (\ref{pf.Delta.hrefl.d.equiv2}))}%
\end{array}
\right) \nonumber\\
&  =\sum\limits_{\substack{u\in\widehat{\operatorname*{Rect}\left(
p,p\right)  };\\u\lessdot\operatorname*{hrefl}v}}\left(  R_{\Delta\left(
p\right)  }^{p-\deg v}f\right)  \left(  u\right) \nonumber\\
&  \ \ \ \ \ \ \ \ \ \ \left(
\begin{array}
[c]{c}%
\text{here, we have substituted }u\text{ for }\operatorname*{hrefl}u\text{ in
the sum}\\
\text{(since }\operatorname*{hrefl}\text{ is a bijection)}%
\end{array}
\right) \nonumber\\
&  =\sum\limits_{\substack{u\in\widehat{\operatorname*{Rect}\left(
p,p\right)  };\\u\lessdot w}}\left(  R_{\Delta\left(  p\right)  }^{p-\deg
v}f\right)  \left(  u\right)  \ \ \ \ \ \ \ \ \ \ \left(  \text{since
}\operatorname*{hrefl}v=w\right) \nonumber\\
&  =\sum\limits_{\substack{u\in\widehat{\Delta\left(  p\right)  };\\u\lessdot
w}}\left(  R_{\Delta\left(  p\right)  }^{p-\deg v}f\right)  \left(  u\right)
\label{pf.Delta.hrefl.g'v.c6.rewriter4}%
\end{align}
(since the elements $u\in\widehat{\operatorname*{Rect}\left(  p,p\right)  }$
such that $u\lessdot w$ are precisely the elements $u\in\widehat{\Delta\left(
p\right)  }$ such that $u\lessdot w$).

Recall that we have to prove the equality (\ref{pf.Delta.hrefl.goal}). Upon
substitution of (\ref{pf.Delta.hrefl.g'v.c5.rewriter1}),
(\ref{pf.Delta.hrefl.g'v.c5.rewriter2}),
(\ref{pf.Delta.hrefl.g'v.c5.rewriter3}) and
(\ref{pf.Delta.hrefl.g'v.c6.rewriter4}), this equality
(\ref{pf.Delta.hrefl.goal}) becomes%
\begin{equation}
\dfrac{1}{\left(  R_{\Delta\left(  p\right)  }^{p-\deg v+1}f\right)  \left(
w\right)  }=\dfrac{1}{\left(  \dfrac{1}{\left(  R_{\Delta\left(  p\right)
}^{p-\deg v}f\right)  \left(  w\right)  }\right)  }\cdot\dfrac{\sum
\limits_{\substack{u\in\widehat{\Delta\left(  p\right)  };\\u \gtrdot
w}}\dfrac{1}{\left(  R_{\Delta\left(  p\right)  }^{p-\deg v+1}f\right)
\left(  u\right)  }}{\sum\limits_{\substack{u\in\widehat{\Delta\left(
p\right)  };\\u \lessdot w}}\left(  R_{\Delta\left(  p\right)  }^{p-\deg
v}f\right)  \left(  u\right)  }. \label{pf.Delta.hrefl.g'v.c6.transgoal}%
\end{equation}
So (\ref{pf.Delta.hrefl.g'v.c6.transgoal}) is what we need to prove.

But just as in Case 5, we can prove that (\ref{pf.Delta.hrefl.g'v.c5.1})
holds. Hence,%
\begin{align*}
\dfrac{1}{\left(  R_{\Delta\left(  p\right)  }^{p-\deg v+1}f\right)  \left(
w\right)  }  &  =\left(  R_{\Delta\left(  p\right)  }^{p-\deg v}f\right)
\left(  w\right)  \cdot\dfrac{\sum\limits_{\substack{u\in\widehat{\Delta
\left(  p\right)  };\\u \gtrdot w}}\dfrac{1}{\left(  R_{\Delta\left(
p\right)  }^{p-\deg v+1}f\right)  \left(  u\right)  }}{\sum
\limits_{\substack{u\in\widehat{\Delta\left(  p\right)  };\\u\lessdot
w}}\left(  R_{\Delta\left(  p\right)  }^{p-\deg v}f\right)  \left(  u\right)
}\\
&  =\dfrac{1}{\left(  \dfrac{1}{\left(  R_{\Delta\left(  p\right)  }^{p-\deg
v}f\right)  \left(  w\right)  }\right)  }\cdot\dfrac{\sum
\limits_{\substack{u\in\widehat{\Delta\left(  p\right)  };\\u \gtrdot
w}}\dfrac{1}{\left(  R_{\Delta\left(  p\right)  }^{p-\deg v+1}f\right)
\left(  u\right)  }}{\sum\limits_{\substack{u\in\widehat{\Delta\left(
p\right)  };\\u \lessdot w}}\left(  R_{\Delta\left(  p\right)  }^{p-\deg
v}f\right)  \left(  u\right)  }.
\end{align*}
But this is precisely (\ref{pf.Delta.hrefl.g'v.c6.transgoal}). Hence,
(\ref{pf.Delta.hrefl.g'v.c6.transgoal}) is proven. Since
(\ref{pf.Delta.hrefl.g'v.c6.transgoal}) was obtained by rewriting
(\ref{pf.Delta.hrefl.goal}), this yields that (\ref{pf.Delta.hrefl.goal}) is
proven in Case 6.

\bigskip

\underline{\textbf{Proof of (\ref{pf.Delta.hrefl.goal}) in Case 7.}}

Finally, let us consider Case 7. In this case, we have $\deg v=1$. Thus,
$v\in\nabla\left(  p\right)  \subseteq\nabla\left(  p\right)  \cup\left\{
0\right\}  $ and $\operatorname*{hrefl}v\in\Delta\left(  p\right)
\subseteq\Delta\left(  p\right)  \cup\left\{  1\right\}  $.

Let $w=\operatorname*{hrefl}v$. Then, $w=\operatorname*{hrefl}v\in
\Delta\left(  p\right)  \subseteq\Delta\left(  p\right)  \cup\left\{
1\right\}  $.

The equalities (\ref{pf.Delta.hrefl.g'v.c5.rewriter1}) and
(\ref{pf.Delta.hrefl.g'v.c5.rewriter2}) hold (their proof proceeds in the same
way as in Case 5). Moreover, (\ref{pf.Delta.hrefl.g'v.c6.rewriter4}) also
holds (this is proven in the same way as we did it in Case 6).

Since $\deg v=1$, we know that $v$ is a minimal element of
$\operatorname*{Rect}\left(  p,p\right)  $. Hence, there is only one element
$u\in\widehat{\operatorname*{Rect}\left(  p,p\right)  }$ such that $u \lessdot
v$, namely the element $0$. Thus,%
\begin{equation}
\sum\limits_{\substack{u\in\widehat{\operatorname*{Rect}\left(  p,p\right)
};\\u \lessdot v}}f^{\prime}\left(  u\right)  =f^{\prime}\left(  0\right)
=\dfrac{1}{f\left(  1\right)  }\ \ \ \ \ \ \ \ \ \ \left(  \text{by
(\ref{pf.Delta.hrefl.f'0})}\right)  . \label{pf.Delta.hrefl.g'v.c7.rewriter3}%
\end{equation}

Recall that we have to prove the equality (\ref{pf.Delta.hrefl.goal}). Upon
substitution of (\ref{pf.Delta.hrefl.g'v.c5.rewriter1}),
(\ref{pf.Delta.hrefl.g'v.c5.rewriter2}),
(\ref{pf.Delta.hrefl.g'v.c7.rewriter3}) and
(\ref{pf.Delta.hrefl.g'v.c6.rewriter4}), this equality
(\ref{pf.Delta.hrefl.goal}) becomes%
\begin{equation}
\dfrac{1}{\left(  R_{\Delta\left(  p\right)  }^{p-\deg v+1}f\right)  \left(
w\right)  }=\dfrac{1}{\left(  \dfrac{1}{\left(  R_{\Delta\left(  p\right)
}^{p-\deg v}f\right)  \left(  w\right)  }\right)  }\cdot\dfrac{\left(
\dfrac{1}{f\left(  1\right)  }\right)  }{\sum\limits_{\substack{u\in
\widehat{\Delta\left(  p\right)  };\\u \lessdot w}}\left(  R_{\Delta\left(
p\right)  }^{p-\deg v}f\right)  \left(  u\right)  }.
\label{pf.Delta.hrefl.g'v.c7.transgoal}%
\end{equation}
So (\ref{pf.Delta.hrefl.g'v.c7.transgoal}) is what we need to prove.

But we know that $v$ is a minimal element of $\operatorname*{Rect}\left(
p,p\right)  $, hence also a minimal element of $\nabla\left(  p\right)  $.
Thus, $\operatorname*{hrefl}v$ is a maximal element of $\Delta\left(
p\right)  $ (since $\operatorname*{hrefl}$ is a poset antiautomorphism and
sends $\nabla\left(  p\right)  $ to $\Delta\left(  p\right)  $). In other
words, $w$ is a maximal element of $\Delta\left(  p\right)  $ (since
$\operatorname*{hrefl}v=w$). Thus, there exists only one element
$u\in\widehat{\Delta\left(  p\right)  }$ such that $u \gtrdot w$, namely the
element $1$. Thus,%
\begin{equation}
\sum\limits_{\substack{u\in\widehat{\Delta\left(  p\right)  };\\u\gtrdot
w}}\dfrac{1}{\left(  R_{\Delta\left(  p\right)  }^{p-\deg v+1}f\right)
\left(  u\right)  }=\dfrac{1}{\left(  R_{\Delta\left(  p\right)  }^{p-\deg
v+1}f\right)  \left(  1\right)  }=\dfrac{1}{f\left(  1\right)  }
\label{pf.Delta.hrefl.g'v.c7.rewriter5}%
\end{equation}
(since Corollary \ref{cor.R.implicit.01} yields $\left(  R_{\Delta\left(
p\right)  }^{p-\deg v+1}f\right)  \left(  1\right)  =f\left(  1\right)  $).

But just as in Case 5, we can prove that (\ref{pf.Delta.hrefl.g'v.c5.1})
holds. Hence,%
\begin{align*}
\dfrac{1}{\left(  R_{\Delta\left(  p\right)  }^{p-\deg v+1}f\right)  \left(
w\right)  }  &  =\left(  R_{\Delta\left(  p\right)  }^{p-\deg v}f\right)
\left(  w\right)  \cdot\dfrac{\sum\limits_{\substack{u\in\widehat{\Delta
\left(  p\right)  };\\u\gtrdot w}}\dfrac{1}{\left(  R_{\Delta\left(  p\right)
}^{p-\deg v+1}f\right)  \left(  u\right)  }}{\sum\limits_{\substack{u\in
\widehat{\Delta\left(  p\right)  };\\u\lessdot w}}\left(  R_{\Delta\left(
p\right)  }^{p-\deg v}f\right)  \left(  u\right)  }\\
&  =\left(  R_{\Delta\left(  p\right)  }^{p-\deg v}f\right)  \left(  w\right)
\cdot\dfrac{\left(  \dfrac{1}{f\left(  1\right)  }\right)  }{\sum
\limits_{\substack{u\in\widehat{\Delta\left(  p\right)  };\\u\lessdot
w}}\left(  R_{\Delta\left(  p\right)  }^{p-\deg v}f\right)  \left(  u\right)
}\ \ \ \ \ \ \ \ \ \ \left(  \text{by (\ref{pf.Delta.hrefl.g'v.c7.rewriter5}%
)}\right) \\
&  =\dfrac{1}{\left(  \dfrac{1}{\left(  R_{\Delta\left(  p\right)  }^{p-\deg
v}f\right)  \left(  w\right)  }\right)  }\cdot\dfrac{\left(  \dfrac
{1}{f\left(  1\right)  }\right)  }{\sum\limits_{\substack{u\in\widehat{\Delta
\left(  p\right)  };\\u\lessdot w}}\left(  R_{\Delta\left(  p\right)
}^{p-\deg v}f\right)  \left(  u\right)  }.
\end{align*}
But this is precisely (\ref{pf.Delta.hrefl.g'v.c7.transgoal}). Hence,
(\ref{pf.Delta.hrefl.g'v.c7.transgoal}) is proven. Since
(\ref{pf.Delta.hrefl.g'v.c7.transgoal}) was obtained by rewriting
(\ref{pf.Delta.hrefl.goal}), this yields that (\ref{pf.Delta.hrefl.goal}) is
proven in Case 7.

\bigskip

\underline{\textbf{Finishing the proof of Lemma \ref{lem.Delta.hrefl} (d).}}

We have now proven (\ref{pf.Delta.hrefl.goal}) in each of the seven cases 1,
2, 3, 4, 5, 6, 7. Since these seven cases exhaust all possibilities, this
yields that (\ref{pf.Delta.hrefl.goal}) always holds. As we know, this shows
that $g^{\prime}=R_{\operatorname*{Rect}\left(  p,p\right)  }f^{\prime}$. Now,%
\begin{align*}
\left(  a_{0},a_{1},...,a_{2p}\right)  \flat\underbrace{\left(  \left(
\operatorname*{wing}\circ R_{\Delta\left(  p\right)  }\right)  f\right)
}_{=g}  &  =\left(  a_{0},a_{1},...,a_{2p}\right)  \flat g=g^{\prime
}=R_{\operatorname*{Rect}\left(  p,p\right)  }\underbrace{f^{\prime}%
}_{=\operatorname*{wing}f}\\
&  =R_{\operatorname*{Rect}\left(  p,p\right)  }\left(  \operatorname*{wing}%
f\right)  =\left(  R_{\operatorname*{Rect}\left(  p,p\right)  }\circ
\operatorname*{wing}\right)  f.
\end{align*}
In other words, (\ref{pf.Delta.hrefl.precise}) is proven. Hence, every
$x\in\widehat{\operatorname*{Rect}\left(  p,p\right)  }$ satisfies%
\[
\left(  \left(  R_{\operatorname*{Rect}\left(  p,p\right)  }\circ
\operatorname*{wing}\right)  f\right)  \left(  x\right)  =a_{\deg x}%
\cdot\left(  \left(  \operatorname*{wing}\circ R_{\Delta\left(  p\right)
}\right)  f\right)  \left(  x\right)  .
\]
Thus, the $\mathbb{K}$-labellings $\left(  \operatorname*{wing}\circ
R_{\Delta\left(  p\right)  }\right)  f$ and $\left(  R_{\operatorname*{Rect}%
\left(  p,p\right)  }\circ\operatorname*{wing}\right)  f$ are homogeneously equivalent.

Now, forget that we fixed $f$. We have thus shown that:%
\begin{equation}
\left(
\begin{array}
[c]{c}%
\text{if }f\text{ is a zero-free }\mathbb{K}\text{-labelling in }%
\mathbb{K}^{\widehat{\Delta\left(  p\right)  }}\text{ which is sufficiently
generic}\\
\text{for }R_{\Delta\left(  p\right)  }^{i}f\text{ to be well-defined for all
}i\in\left\{  0,1,...,p\right\}  \text{, then}\\
\left(  \operatorname*{wing}\circ R_{\Delta\left(  p\right)  }\right)  f\text{
and }\left(  R_{\operatorname*{Rect}\left(  p,p\right)  }\circ
\operatorname*{wing}\right)  f\text{ are homogeneously}\\
\text{equivalent }\mathbb{K}\text{-labellings of }\operatorname*{Rect}\left(
p,p\right)
\end{array}
\right)  . \label{pf.Delta.hrefl.almostthere}%
\end{equation}

Now, let $\widetilde{f}$ be any element of $\overline{\mathbb{K}%
^{\widehat{\Delta\left(  p\right)  }}}$ which is sufficiently generic for
$\overline{R}_{\Delta\left(  p\right)  }^{i}\widetilde{f}$ to be well-defined
for all $i\in\left\{  0,1,...,p\right\}  $. Then, there exists a zero-free
$\mathbb{K}$-labelling $f$ in $\mathbb{K}^{\widehat{\Delta\left(  p\right)  }%
}$ such that $\widetilde{f}=\pi\left(  f\right)  $, and such that $f$ is
sufficiently generic for $R_{\Delta\left(  p\right)  }^{i}f$ to be
well-defined for all $i\in\left\{  0,1,...,p\right\}  $. Consider this $f$.
According to (\ref{pf.Delta.hrefl.almostthere}), we know that $\left(
\operatorname*{wing}\circ R_{\Delta\left(  p\right)  }\right)  f$ and $\left(
R_{\operatorname*{Rect}\left(  p,p\right)  }\circ\operatorname*{wing}\right)
f$ are homogeneously equivalent $\mathbb{K}$-labellings of
$\operatorname*{Rect}\left(  p,p\right)  $. In other words,%
\begin{equation}
\pi\left(  \left(  \operatorname*{wing}\circ R_{\Delta\left(  p\right)
}\right)  f\right)  =\pi\left(  \left(  R_{\operatorname*{Rect}\left(
p,p\right)  }\circ\operatorname*{wing}\right)  f\right)  .
\label{pf.Delta.hrefl.equalpis}%
\end{equation}
But since $\widetilde{f}=\pi\left(  f\right)  $, we have
\begin{align*}
&  \left(  \overline{R}_{\operatorname*{Rect}\left(  p,p\right)  }%
\circ\overline{\operatorname*{wing}}\right)  \left(  \widetilde{f}\right) \\
&  =\left(  \overline{R}_{\operatorname*{Rect}\left(  p,p\right)  }%
\circ\overline{\operatorname*{wing}}\right)  \left(  \pi\left(  f\right)
\right)  =\left(  \overline{R}_{\operatorname*{Rect}\left(  p,p\right)  }%
\circ\underbrace{\overline{\operatorname*{wing}}\circ\pi}_{\substack{=\pi
\circ\operatorname*{wing}\\\text{(by the commutativity}\\\text{of the diagram
(\ref{lem.Delta.hrefl.commut}))}}}\right)  \left(  f\right) \\
&  =\left(  \underbrace{\overline{R}_{\operatorname*{Rect}\left(  p,p\right)
}\circ\pi}_{\substack{=\pi\circ R_{\operatorname*{Rect}\left(  p,p\right)
}\\\text{(by the commutativity}\\\text{of the diagram (\ref{def.hgR.commut}%
))}}}\circ\operatorname*{wing}\right)  \left(  f\right)  =\left(  \pi\circ
R_{\operatorname*{Rect}\left(  p,p\right)  }\circ\operatorname*{wing}\right)
\left(  f\right) \\
&  =\pi\left(  \left(  R_{\operatorname*{Rect}\left(  p,p\right)  }%
\circ\operatorname*{wing}\right)  f\right)  =\pi\left(  \left(
\operatorname*{wing}\circ R_{\Delta\left(  p\right)  }\right)  f\right)
\ \ \ \ \ \ \ \ \ \ \left(  \text{by (\ref{pf.Delta.hrefl.equalpis})}\right)
\\
&  =\left(  \underbrace{\pi\circ\operatorname*{wing}}_{\substack{=\overline
{\operatorname*{wing}}\circ\pi\\\text{(by the commutativity}\\\text{of the
diagram (\ref{lem.Delta.hrefl.commut}))}}}\circ R_{\Delta\left(  p\right)
}\right)  \left(  f\right)  =\left(  \overline{\operatorname*{wing}}%
\circ\underbrace{\pi\circ R_{\Delta\left(  p\right)  }}_{\substack{=\overline
{R}_{\Delta\left(  p\right)  }\circ\pi\\\text{(by the commutativity}\\\text{of
the diagram (\ref{def.hgR.commut}))}}}\right)  \left(  f\right) \\
&  =\left(  \overline{\operatorname*{wing}}\circ\overline{R}_{\Delta\left(
p\right)  }\circ\pi\right)  \left(  f\right)  =\left(  \overline
{\operatorname*{wing}}\circ\overline{R}_{\Delta\left(  p\right)  }\right)
\left(  \underbrace{\pi\left(  f\right)  }_{=\widetilde{f}}\right)  =\left(
\overline{\operatorname*{wing}}\circ\overline{R}_{\Delta\left(  p\right)
}\right)  \left(  \widetilde{f}\right)  .
\end{align*}

Now, forget that we fixed $\widetilde{f}$. We have thus proven that $\left(
\overline{R}_{\operatorname*{Rect}\left(  p,p\right)  }\circ\overline
{\operatorname*{wing}}\right)  \left(  \widetilde{f}\right)  =\left(
\overline{\operatorname*{wing}}\circ\overline{R}_{\Delta\left(  p\right)
}\right)  \left(  \widetilde{f}\right)  $ for every element $\widetilde{f}$ of
$\overline{\mathbb{K}^{\widehat{\Delta\left(  p\right)  }}}$ which is
sufficiently generic for $\overline{R}_{\Delta\left(  p\right)  }%
^{i}\widetilde{f}$ to be well-defined for all $i\in\left\{  0,1,...,p\right\}
$. In other words, $\overline{R}_{\operatorname*{Rect}\left(  p,p\right)
}\circ\overline{\operatorname*{wing}}=\overline{\operatorname*{wing}}%
\circ\overline{R}_{\Delta\left(  p\right)  }$ as rational maps. This proves
Lemma \ref{lem.Delta.hrefl} \textbf{(d)}.

\bigskip

\underline{\textbf{Proof of (f).}}

\textbf{(f)} As in the proof of Lemma \ref{lem.Delta.hrefl} \textbf{(d)}, we
WLOG assume that $p>2$.

Let $f$ be a zero-free $\mathbb{K}$-labelling in $\mathbb{K}^{\widehat{\Delta
\left(  p\right)  }}$ which satisfies $f\left(  0\right)  =2$ and is
sufficiently generic for $R_{\Delta\left(  p\right)  }^{i}f$ to be
well-defined for all $i\in\left\{  0,1,...,p\right\}  $. We must prove that
\[
R_{\operatorname*{Rect}\left(  p,p\right)  }\left(  \operatorname*{wing}%
f\right)  =\operatorname*{wing}\left(  R_{\Delta\left(  p\right)  }f\right)
.
\]

Notice that $f\left(  0\right)  \neq0$. Define a $\left(  2p+1\right)  $-tuple
$\left(  a_{0},a_{1},...,a_{2p}\right)  \in\left(  \mathbb{K}^{\times}\right)
^{2p+1}$ as in the proof of Lemma \ref{lem.Delta.hrefl} \textbf{(d)}. Then,
for every $i\in\left\{  0,1,...,2p\right\}  $, we have
\begin{align*}
a_{i}  &  =\left\{
\begin{array}
[c]{l}%
1,\ \ \ \ \ \ \ \ \ \ \text{if }i>p+1;\\
\dfrac{2}{f\left(  0\right)  },\ \ \ \ \ \ \ \ \ \ \text{if }p\leq i\leq
p+1;\\
1,\ \ \ \ \ \ \ \ \ \ \text{if }i<p
\end{array}
\right.  =\left\{
\begin{array}
[c]{l}%
1,\ \ \ \ \ \ \ \ \ \ \text{if }i>p+1;\\
1,\ \ \ \ \ \ \ \ \ \ \text{if }p\leq i\leq p+1;\\
1,\ \ \ \ \ \ \ \ \ \ \text{if }i<p
\end{array}
\right. \\
&  \ \ \ \ \ \ \ \ \ \ \left(  \text{since }\dfrac{2}{f\left(  0\right)
}=1\text{ (because }f\left(  0\right)  =2\text{)}\right) \\
&  =1.
\end{align*}
In other words, $\left(  a_{0},a_{1},...,a_{2p}\right)  =\left(
\underbrace{1,1,...,1}_{2p+1\text{ times}}\right)  $.

But as in the proof of Lemma \ref{lem.Delta.hrefl} \textbf{(d)}, we can show
that (\ref{pf.Delta.hrefl.precise}) holds. Hence,%
\begin{align*}
\left(  R_{\operatorname*{Rect}\left(  p,p\right)  }\circ\operatorname*{wing}%
\right)  f  &  =\underbrace{\left(  a_{0},a_{1},...,a_{2p}\right)  }_{=\left(
\underbrace{1,1,...,1}_{2p+1\text{ times}}\right)  }\flat\underbrace{\left(
\left(  \operatorname*{wing}\circ R_{\Delta\left(  p\right)  }\right)
f\right)  }_{=\operatorname*{wing}\left(  R_{\Delta\left(  p\right)
}f\right)  }\\
&  =\left(  \underbrace{1,1,...,1}_{2p+1\text{ times}}\right)  \flat\left(
\operatorname*{wing}\left(  R_{\Delta\left(  p\right)  }f\right)  \right)
=\operatorname*{wing}\left(  R_{\Delta\left(  p\right)  }f\right)  .
\end{align*}
Since $\left(  R_{\operatorname*{Rect}\left(  p,p\right)  }\circ
\operatorname*{wing}\right)  f=R_{\operatorname*{Rect}\left(  p,p\right)
}\left(  \operatorname*{wing}f\right)  $, this rewrites as \newline%
$R_{\operatorname*{Rect}\left(  p,p\right)  }\left(  \operatorname*{wing}%
f\right)  =\operatorname*{wing}\left(  R_{\Delta\left(  p\right)  }f\right)
$. This proves Lemma \ref{lem.Delta.hrefl} \textbf{(f)}.

\bigskip

\underline{\textbf{More trivialities.}}

\textbf{(e)} There is a quick way to prove the equalities
(\ref{lem.Delta.hrefl.e.1}) and (\ref{lem.Delta.hrefl.e.2}) by convincing
oneself that they are obvious.

Namely, let $\mathcal{P}_{\cong}$ be the category of all finite posets, where
the morphisms are given by poset isomorphisms. (So $\mathcal{P}_{\cong}$ is a
groupoid.) For every isomorphism $\phi:P\rightarrow Q$ between two finite
posets $P$ and $Q$, we can define a map $\phi^{\ast}:\mathbb{K}^{\widehat{Q}%
}\rightarrow\mathbb{K}^{\widehat{P}}$ by setting%
\[
\left(  \phi^{\ast}f\right)  \left(  v\right)  =f\left(  \phi\left(  v\right)
\right)  \ \ \ \ \ \ \ \ \ \ \text{for all }v\in\widehat{P}%
\]
(where $\phi$ is extended to send $1\in\widehat{P}$ to $1\in\widehat{Q}$ and
to send $0\in\widehat{P}$ to $0\in\widehat{Q}$). This definition does not
conflict with our formerly introduced notations $\operatorname*{vrefl}%
\nolimits^{\ast}$, because if we take $\phi$ to be either of the maps
$\operatorname*{vrefl}:\operatorname*{Rect}\left(  p,p\right)  \rightarrow
\operatorname*{Rect}\left(  p,p\right)  $ or $\operatorname*{vrefl}%
:\Delta\left(  p\right)  \rightarrow\Delta\left(  p\right)  $, then
$\phi^{\ast}$ is indeed one of the maps that we denoted by
$\operatorname*{vrefl}\nolimits^{\ast}$.

Let $\operatorname*{Birat}\nolimits_{\mathbb{K}}$ denote the category of
affine algebraic varieties over $\mathbb{K}$, with morphisms being invertible
dominant rational maps. Then, it is easy to see that we can define a
contravariant functor $\operatorname*{Labellings}:\mathcal{P}_{\cong%
}\rightarrow\operatorname*{Birat}\nolimits_{\mathbb{K}}$ by%
\begin{align*}
\operatorname*{Labellings}\left(  P\right)   &  =\mathbb{K}^{\widehat{P}%
}\ \ \ \ \ \ \ \ \ \ \text{for every finite poset }P;\\
\operatorname*{Labellings}\left(  \phi\right)   &  =\phi^{\ast}%
\ \ \ \ \ \ \ \ \ \ \text{for any isomorphism }\phi:P\rightarrow Q\text{
between two posets }P\text{ and }Q.
\end{align*}
Furthermore, it is easy to see that birational rowmotion $R$ is a natural
homomorphism $\operatorname*{Labellings}\rightarrow\operatorname*{Labellings}$
(this is more or less a formal way to state that the definition of $R$ only
depends on the structure of the poset $P$, not on how its elements were
labelled). That is, whenever $\phi:P\rightarrow Q$ is an isomorphism between
two finite posets $P$ and $Q$, we have $\phi^{\ast}\circ R_{Q}=R_{P}\circ
\phi^{\ast}$. Applying this to $P=\Delta\left(  p\right)  $, $Q=\Delta\left(
p\right)  $ and $\phi=\operatorname*{vrefl}$, we obtain $\operatorname*{vrefl}%
\nolimits^{\ast}\circ R_{\Delta\left(  p\right)  }=R_{\Delta\left(  p\right)
}\circ\operatorname*{vrefl}\nolimits^{\ast}$. Thus, (\ref{lem.Delta.hrefl.e.1}%
) is proven. Similarly, (\ref{lem.Delta.hrefl.e.2}) can be shown.

It remains to prove (\ref{lem.Delta.hrefl.e.3}). This isn't much harder, but
let us do this in a more pedestrian way. First let us notice that%
\begin{equation}
\operatorname*{vrefl}\nolimits^{\ast}\circ R_{\Delta\left(  p\right)  }^{\ell
}=R_{\Delta\left(  p\right)  }^{\ell}\circ\operatorname*{vrefl}\nolimits^{\ast
}\ \ \ \ \ \ \ \ \ \ \text{for every }\ell\in\mathbb{N}.
\label{pf.Delta.hrefl.e.4}%
\end{equation}
(Indeed, this follows easily by induction over $\ell$, using
(\ref{lem.Delta.hrefl.e.1}).)

Let $f\in\mathbb{K}^{\widehat{\Delta\left(  p\right)  }}$ be a $\mathbb{K}%
$-labelling which is sufficiently generic that both maps
$\operatorname*{vrefl}\nolimits^{\ast}\circ\operatorname*{wing}$ and
$\operatorname*{wing}\circ\operatorname*{vrefl}\nolimits^{\ast}$ can be
applied to $f$. For every $v\in\widehat{\operatorname*{Rect}\left(
p,p\right)  }$, we have%
\begin{align*}
&  \left(  \left(  \operatorname*{vrefl}\nolimits^{\ast}\circ
\operatorname*{wing}\right)  f\right)  \left(  v\right) \\
&  =\left(  \operatorname*{vrefl}\nolimits^{\ast}\left(  \operatorname*{wing}%
f\right)  \right)  \left(  v\right)  =\left(  \operatorname*{wing}f\right)
\left(  \operatorname*{vrefl}v\right)  \ \ \ \ \ \ \ \ \ \ \left(  \text{by
the definition of }\operatorname*{vrefl}\nolimits^{\ast}\right) \\
&  =\left\{
\begin{array}
[c]{l}%
f\left(  \operatorname*{vrefl}v\right)  ,\ \ \ \ \ \ \ \ \ \ \text{if
}\operatorname*{vrefl}v\in\Delta\left(  p\right)  \cup\left\{  1\right\}  ;\\
1,\ \ \ \ \ \ \ \ \ \ \text{if }\operatorname*{vrefl}v\in\operatorname*{Eq}%
\left(  p\right)  ;\\
\dfrac{1}{\left(  R_{\Delta\left(  p\right)  }^{p-\deg\left(
\operatorname*{vrefl}v\right)  }f\right)  \left(  \operatorname*{hrefl}\left(
\operatorname*{vrefl}v\right)  \right)  },\ \ \ \ \ \ \ \ \ \ \text{if
}\operatorname*{vrefl}v\in\nabla\left(  p\right)  \cup\left\{  0\right\}
\end{array}
\right. \\
&  \ \ \ \ \ \ \ \ \ \ \left(  \text{by the definition of }%
\operatorname*{wing}\right) \\
&  =\left\{
\begin{array}
[c]{l}%
f\left(  \operatorname*{vrefl}v\right)  ,\ \ \ \ \ \ \ \ \ \ \text{if }%
v\in\Delta\left(  p\right)  \cup\left\{  1\right\}  ;\\
1,\ \ \ \ \ \ \ \ \ \ \text{if }v\in\operatorname*{Eq}\left(  p\right)  ;\\
\dfrac{1}{\left(  R_{\Delta\left(  p\right)  }^{p-\deg\left(
\operatorname*{vrefl}v\right)  }f\right)  \left(  \operatorname*{hrefl}\left(
\operatorname*{vrefl}v\right)  \right)  },\ \ \ \ \ \ \ \ \ \ \text{if }%
v\in\nabla\left(  p\right)  \cup\left\{  0\right\}
\end{array}
\right. \\
&  \ \ \ \ \ \ \ \ \ \ \left(
\begin{array}
[c]{c}%
\text{since }\left(  \operatorname*{vrefl}v\in\Delta\left(  p\right)
\cup\left\{  1\right\}  \text{ is equivalent to }v\in\Delta\left(  p\right)
\cup\left\{  1\right\}  \right)  \text{,}\\
\left(  \operatorname*{vrefl}v\in\operatorname*{Eq}\left(  p\right)  \text{ is
equivalent to }v\in\operatorname*{Eq}\left(  p\right)  \right)  \text{ and}\\
\left(  \operatorname*{vrefl}v\in\nabla\left(  p\right)  \cup\left\{
0\right\}  \text{ is equivalent to }v\in\nabla\left(  p\right)  \cup\left\{
0\right\}  \right)
\end{array}
\right) \\
&  =\left\{
\begin{array}
[c]{l}%
f\left(  \operatorname*{vrefl}v\right)  ,\ \ \ \ \ \ \ \ \ \ \text{if }%
v\in\Delta\left(  p\right)  \cup\left\{  1\right\}  ;\\
1,\ \ \ \ \ \ \ \ \ \ \text{if }v\in\operatorname*{Eq}\left(  p\right)  ;\\
\dfrac{1}{\left(  R_{\Delta\left(  p\right)  }^{p-\deg v}f\right)  \left(
\operatorname*{vrefl}\left(  \operatorname*{hrefl}v\right)  \right)
},\ \ \ \ \ \ \ \ \ \ \text{if }v\in\nabla\left(  p\right)  \cup\left\{
0\right\}
\end{array}
\right. \\
&  \ \ \ \ \ \ \ \ \ \ \left(  \text{since }\deg\left(  \operatorname*{vrefl}%
v\right)  =\deg v\text{ and }\operatorname*{hrefl}\left(
\operatorname*{vrefl}v\right)  =\operatorname*{vrefl}\left(
\operatorname*{hrefl}v\right)  \right)
\end{align*}
and%
\begin{align*}
&  \left(  \left(  \operatorname*{wing}\circ\operatorname*{vrefl}%
\nolimits^{\ast}\right)  f\right)  \left(  v\right) \\
&  =\left(  \operatorname*{wing}\left(  \operatorname*{vrefl}\nolimits^{\ast
}f\right)  \right)  \left(  v\right) \\
&  =\left\{
\begin{array}
[c]{l}%
\left(  \operatorname*{vrefl}\nolimits^{\ast}f\right)  \left(  v\right)
,\ \ \ \ \ \ \ \ \ \ \text{if }v\in\Delta\left(  p\right)  \cup\left\{
1\right\}  ;\\
1,\ \ \ \ \ \ \ \ \ \ \text{if }v\in\operatorname*{Eq}\left(  p\right)  ;\\
\dfrac{1}{\left(  R_{\Delta\left(  p\right)  }^{p-\deg v}\left(
\operatorname*{vrefl}\nolimits^{\ast}f\right)  \right)  \left(
\operatorname*{hrefl}v\right)  },\ \ \ \ \ \ \ \ \ \ \text{if }v\in
\nabla\left(  p\right)  \cup\left\{  0\right\}
\end{array}
\right. \\
&  \ \ \ \ \ \ \ \ \ \ \left(  \text{by the definition of }%
\operatorname*{wing}\right) \\
&  =\left\{
\begin{array}
[c]{l}%
\left(  \operatorname*{vrefl}\nolimits^{\ast}f\right)  \left(  v\right)
,\ \ \ \ \ \ \ \ \ \ \text{if }v\in\Delta\left(  p\right)  \cup\left\{
1\right\}  ;\\
1,\ \ \ \ \ \ \ \ \ \ \text{if }v\in\operatorname*{Eq}\left(  p\right)  ;\\
\dfrac{1}{\left(  \operatorname*{vrefl}\nolimits^{\ast}\left(  R_{\Delta
\left(  p\right)  }^{p-\deg v}f\right)  \right)  \left(  \operatorname*{hrefl}%
v\right)  },\ \ \ \ \ \ \ \ \ \ \text{if }v\in\nabla\left(  p\right)
\cup\left\{  0\right\}
\end{array}
\right. \\
&  \ \ \ \ \ \ \ \ \ \ \left(
\begin{array}
[c]{c}%
\text{since (\ref{pf.Delta.hrefl.e.4}) (applied to }\ell=p-\deg v\text{)
yields}\\
\operatorname*{vrefl}\nolimits^{\ast}\circ R_{\Delta\left(  p\right)
}^{p-\deg v}=R_{\Delta\left(  p\right)  }^{p-\deg v}\circ\operatorname*{vrefl}%
\nolimits^{\ast}\text{, so that}\\
R_{\Delta\left(  p\right)  }^{p-\deg v}\left(  \operatorname*{vrefl}%
\nolimits^{\ast}f\right)  =\underbrace{\left(  R_{\Delta\left(  p\right)
}^{p-\deg v}\circ\operatorname*{vrefl}\nolimits^{\ast}\right)  }%
_{=\operatorname*{vrefl}\nolimits^{\ast}\circ R_{\Delta\left(  p\right)
}^{p-\deg v}}f\\
=\left(  \operatorname*{vrefl}\nolimits^{\ast}\circ R_{\Delta\left(  p\right)
}^{p-\deg v}\right)  f=\operatorname*{vrefl}\nolimits^{\ast}\left(
R_{\Delta\left(  p\right)  }^{p-\deg v}f\right)
\end{array}
\right) \\
&  =\left\{
\begin{array}
[c]{l}%
f\left(  \operatorname*{vrefl}v\right)  ,\ \ \ \ \ \ \ \ \ \ \text{if }%
v\in\Delta\left(  p\right)  \cup\left\{  1\right\}  ;\\
1,\ \ \ \ \ \ \ \ \ \ \text{if }v\in\operatorname*{Eq}\left(  p\right)  ;\\
\dfrac{1}{\left(  R_{\Delta\left(  p\right)  }^{p-\deg v}f\right)  \left(
\operatorname*{vrefl}\left(  \operatorname*{hrefl}v\right)  \right)
},\ \ \ \ \ \ \ \ \ \ \text{if }v\in\nabla\left(  p\right)  \cup\left\{
0\right\}
\end{array}
\right. \\
&  \ \ \ \ \ \ \ \ \ \ \left(
\begin{array}
[c]{c}%
\text{because the definition of }\operatorname*{vrefl}\nolimits^{\ast}\text{
yields }\left(  \operatorname*{vrefl}\nolimits^{\ast}f\right)  \left(
v\right)  =f\left(  \operatorname*{vrefl}v\right)  \text{ and}\\
\left(  \operatorname*{vrefl}\nolimits^{\ast}\left(  R_{\Delta\left(
p\right)  }^{p-\deg v}f\right)  \right)  \left(  \operatorname*{hrefl}%
v\right)  =\left(  R_{\Delta\left(  p\right)  }^{p-\deg v}f\right)  \left(
\operatorname*{vrefl}\left(  \operatorname*{hrefl}v\right)  \right)
\end{array}
\right)  .
\end{align*}
Hence, for every $v\in\widehat{\operatorname*{Rect}\left(  p,p\right)  }$, we
have
\begin{align*}
&  \left(  \left(  \operatorname*{vrefl}\nolimits^{\ast}\circ
\operatorname*{wing}\right)  f\right)  \left(  v\right) \\
&  =\left\{
\begin{array}
[c]{l}%
f\left(  \operatorname*{vrefl}v\right)  ,\ \ \ \ \ \ \ \ \ \ \text{if }%
v\in\Delta\left(  p\right)  \cup\left\{  1\right\}  ;\\
1,\ \ \ \ \ \ \ \ \ \ \text{if }v\in\operatorname*{Eq}\left(  p\right)  ;\\
\dfrac{1}{\left(  R_{\Delta\left(  p\right)  }^{p-\deg v}f\right)  \left(
\operatorname*{vrefl}\left(  \operatorname*{hrefl}v\right)  \right)
},\ \ \ \ \ \ \ \ \ \ \text{if }v\in\nabla\left(  p\right)  \cup\left\{
0\right\}
\end{array}
\right.  =\left(  \left(  \operatorname*{wing}\circ\operatorname*{vrefl}%
\nolimits^{\ast}\right)  f\right)  \left(  v\right)  .
\end{align*}
In other words, $\left(  \operatorname*{vrefl}\nolimits^{\ast}\circ
\operatorname*{wing}\right)  f=\left(  \operatorname*{wing}\circ
\operatorname*{vrefl}\nolimits^{\ast}\right)  f$. Now, forget that we fixed
$f$. We thus have shown that $\left(  \operatorname*{vrefl}\nolimits^{\ast
}\circ\operatorname*{wing}\right)  f=\left(  \operatorname*{wing}%
\circ\operatorname*{vrefl}\nolimits^{\ast}\right)  f$ for every $\mathbb{K}%
$-labelling $f\in\mathbb{K}^{\widehat{\Delta\left(  p\right)  }}$ which is
sufficiently generic that both maps $\operatorname*{vrefl}\nolimits^{\ast
}\circ\operatorname*{wing}$ and $\operatorname*{wing}\circ
\operatorname*{vrefl}\nolimits^{\ast}$ can be applied to $f$. In other words,
$\operatorname*{vrefl}\nolimits^{\ast}\circ\operatorname*{wing}%
=\operatorname*{wing}\circ\operatorname*{vrefl}\nolimits^{\ast}$. This proves
(\ref{lem.Delta.hrefl.e.3}). Thus, the proof of Lemma \ref{lem.Delta.hrefl}
\textbf{(e)} is complete.

\textbf{(g)} Lemma \ref{lem.Delta.hrefl} \textbf{(g)} follows easily from the
fact that $\operatorname*{vrefl}$ preserves the degree of an element of
$\Delta\left(  p\right)  $. This finally completes the proof of Lemma
\ref{lem.Delta.hrefl}.
\end{proof}
\end{verlong}

\bigskip

For easier reference, let us record a corollary of Lemma \ref{lem.Delta.hrefl}
\textbf{(f)}:

\begin{corollary}
\label{cor.Delta.hrefl.f}Let $p$ be a positive integer. Let $\mathbb{K}$ be a
field of characteristic $\neq2$. Consider the map $\operatorname*{wing}$
defined in Lemma \ref{lem.Delta.hrefl}. Let $\ell\in\mathbb{N}$.

Then, almost every (in the sense of Zariski topology) labelling $f\in
\mathbb{K}^{\widehat{\Delta\left(  p\right)  }}$ satisfying $f\left(
0\right)  =2$ satisfies%
\[
R_{\operatorname*{Rect}\left(  p,p\right)  }^{\ell}\left(
\operatorname*{wing}f\right)  =\operatorname*{wing}\left(  R_{\Delta\left(
p\right)  }^{\ell}f\right)  .
\]

\end{corollary}

\begin{vershort}
\begin{proof}
[Proof of Corollary \ref{cor.Delta.hrefl.f} (sketched).]The proof of Corollary
\ref{cor.Delta.hrefl.f} is an easy induction over $\ell$ (details left to the
reader), using Lemma \ref{lem.Delta.hrefl} \textbf{(f)} and the fact that
$R_{\Delta\left(  p\right)  }$ does not change the label at $1$.
\end{proof}
\end{vershort}

\begin{verlong}
\begin{proof}
[Proof of Corollary \ref{cor.Delta.hrefl.f} (sketched).]We will prove
Corollary \ref{cor.Delta.hrefl.f} by induction over $\ell$:

\textit{Induction base:} For $\ell=0$, Corollary \ref{cor.Delta.hrefl.f} is
obvious (because $R_{\operatorname*{Rect}\left(  p,p\right)  }^{0}%
=\operatorname*{id}$ and $R_{\Delta\left(  p\right)  }^{0}=\operatorname*{id}%
$). Hence, the induction base is complete.

\textit{Induction step:} Let $L$ be a positive integer. Assume that Corollary
\ref{cor.Delta.hrefl.f} holds for $\ell=L$. We need to show that Corollary
\ref{cor.Delta.hrefl.f} holds for $\ell=L+1$.

We know that Corollary \ref{cor.Delta.hrefl.f} holds for $\ell=L$. In other
words, almost every (in the sense of Zariski topology) labelling
$f\in\mathbb{K}^{\widehat{\Delta\left(  p\right)  }}$ satisfying $f\left(
0\right)  =2$ satisfies%
\begin{equation}
R_{\operatorname*{Rect}\left(  p,p\right)  }^{L}\left(  \operatorname*{wing}%
f\right)  =\operatorname*{wing}\left(  R_{\Delta\left(  p\right)  }%
^{L}f\right)  . \label{pf.Delta.hrefl.f.1}%
\end{equation}
Now, let $f$ be any labelling $f\in\mathbb{K}^{\widehat{\Delta\left(
p\right)  }}$ satisfying $f\left(  0\right)  =2$ which is sufficiently generic
for (\ref{pf.Delta.hrefl.f.1}) to hold and for $R_{\operatorname*{Rect}\left(
p,p\right)  }^{L+1}\left(  \operatorname*{wing}f\right)  $ and
$\operatorname*{wing}\left(  R_{\Delta\left(  p\right)  }^{L+1}f\right)  $ to
be well-defined. From Corollary \ref{cor.R.implicit.01}, we obtain $\left(
R_{\Delta\left(  p\right)  }^{L}f\right)  \left(  0\right)  =f\left(
0\right)  =2$. Moreover, $R_{\Delta\left(  p\right)  }^{L}f$ is also generic
(since $R_{\Delta\left(  p\right)  }^{L}$ is invertible). Hence, we can apply
Lemma \ref{lem.Delta.hrefl} \textbf{(f)} to $R_{\Delta\left(  p\right)  }%
^{L}f$ instead of $f$. We thus obtain%
\begin{equation}
R_{\operatorname*{Rect}\left(  p,p\right)  }\left(  \operatorname*{wing}%
\left(  R_{\Delta\left(  p\right)  }^{L}f\right)  \right)
=\operatorname*{wing}\left(  R_{\Delta\left(  p\right)  }\left(
R_{\Delta\left(  p\right)  }^{L}f\right)  \right)  .
\label{pf.Delta.hrefl.f.2}%
\end{equation}
Now,%
\begin{align*}
\underbrace{R_{\operatorname*{Rect}\left(  p,p\right)  }^{L+1}}%
_{=R_{\operatorname*{Rect}\left(  p,p\right)  }\circ R_{\operatorname*{Rect}%
\left(  p,p\right)  }^{L}}\left(  \operatorname*{wing}f\right)   &  =\left(
R_{\operatorname*{Rect}\left(  p,p\right)  }\circ R_{\operatorname*{Rect}%
\left(  p,p\right)  }^{L}\right)  \left(  \operatorname*{wing}f\right) \\
&  =R_{\operatorname*{Rect}\left(  p,p\right)  }\left(
\underbrace{R_{\operatorname*{Rect}\left(  p,p\right)  }^{L}\left(
\operatorname*{wing}f\right)  }_{\substack{=\operatorname*{wing}\left(
R_{\Delta\left(  p\right)  }^{L}f\right)  \\\text{(by
(\ref{pf.Delta.hrefl.f.1}))}}}\right) \\
&  =R_{\operatorname*{Rect}\left(  p,p\right)  }\left(  \operatorname*{wing}%
\left(  R_{\Delta\left(  p\right)  }^{L}f\right)  \right)
=\operatorname*{wing}\left(  \underbrace{R_{\Delta\left(  p\right)  }\left(
R_{\Delta\left(  p\right)  }^{L}f\right)  }_{=\left(  R_{\Delta\left(
p\right)  }\circ R_{\Delta\left(  p\right)  }^{L}\right)  f}\right) \\
&  \ \ \ \ \ \ \ \ \ \ \left(  \text{by (\ref{pf.Delta.hrefl.f.2})}\right) \\
&  =\operatorname*{wing}\left(  \underbrace{\left(  R_{\Delta\left(  p\right)
}\circ R_{\Delta\left(  p\right)  }^{L}\right)  }_{=R_{\Delta\left(  p\right)
}^{L+1}}f\right)  =\operatorname*{wing}\left(  R_{\Delta\left(  p\right)
}^{L+1}f\right)  .
\end{align*}

Now, forget that we fixed $f$. We thus have shown that almost every (in the
sense of Zariski topology) labelling $f\in\mathbb{K}^{\widehat{\Delta\left(
p\right)  }}$ satisfying $f\left(  0\right)  =2$ satisfies%
\[
R_{\operatorname*{Rect}\left(  p,p\right)  }^{L+1}\left(  \operatorname*{wing}%
f\right)  =\operatorname*{wing}\left(  R_{\Delta\left(  p\right)  }%
^{L+1}f\right)  .
\]
In other words, Corollary \ref{cor.Delta.hrefl.f} holds for $\ell=L+1$. This
completes the induction step. Thus, by induction, Corollary
\ref{cor.Delta.hrefl.f} is proven.
\end{proof}
\end{verlong}

We can now proceed to the proof of the theorems stated at the beginning of
this section:

\begin{vershort}
\begin{proof}
[Proof of Theorem \ref{thm.Delta.halfway} (sketched).]The result that we are
striving to prove is a collection of identities between rational functions,
hence boils down to a collection of polynomial identities in the labels of an
arbitrary $\mathbb{K}$-labelling of $\Delta\left(  p\right)  $. Therefore, it
is enough to prove it in the case when $\mathbb{K}$ is a field of rational
functions in finitely many variables over $\mathbb{Q}$. So let us WLOG assume
that we are in this case. Then, the characteristic of $\mathbb{K}$ is $\neq2$
(it is $0$ indeed), so that we can apply Lemma \ref{lem.Delta.hrefl} and
Corollary \ref{cor.Delta.hrefl.f}.

Consider the maps $\operatorname*{hrefl}$, $\operatorname*{wing}$,
$\operatorname*{vrefl}$ and $\operatorname*{vrefl}\nolimits^{\ast}$ defined in
Lemma \ref{lem.Delta.hrefl}. Clearly, it will be enough to prove that%
\[
R_{\Delta\left(  p\right)  }^{p}=\operatorname*{vrefl}\nolimits^{\ast}%
\]
as rational maps $\mathbb{K}^{\widehat{\Delta\left(  p\right)  }%
}\dashrightarrow\mathbb{K}^{\widehat{\Delta\left(  p\right)  }}$. In other
words, it will be enough to prove that $R_{\Delta\left(  p\right)  }%
^{p}g=\operatorname*{vrefl}\nolimits^{\ast}g$ for almost every $g\in
\mathbb{K}^{\widehat{\Delta\left(  p\right)  }}$.

So let $g\in\mathbb{K}^{\widehat{\Delta\left(  p\right)  }}$ be any
sufficiently generic zero-free labelling of $\Delta\left(  p\right)  $. We
need to show that $R_{\Delta\left(  p\right)  }^{p}g=\operatorname*{vrefl}%
\nolimits^{\ast}g$.

Let us use Definition \ref{def.bemol}. The poset $\Delta\left(  p\right)  $ is
$\left(  p-1\right)  $-graded. We can find a $\left(  p+1\right)  $-tuple
$\left(  a_{0},a_{1},...,a_{p}\right)  \in\left(  \mathbb{K}^{\times}\right)
^{p+1}$ such that $\left(  \left(  a_{0},a_{1},...,a_{p}\right)  \flat
g\right)  \left(  0\right)  =2$ (by setting $a_{0}=\dfrac{2}{g\left(
0\right)  }$, and choosing all other $a_{i}$ arbitrarily). Fix such a $\left(
p+1\right)  $-tuple, and set $f=\left(  a_{0},a_{1},...,a_{p}\right)  \flat
g$. Then, $f\left(  0\right)  =2$. We are going to prove that $R_{\Delta
\left(  p\right)  }^{p}f=\operatorname*{vrefl}\nolimits^{\ast}f$. Until we
have done this, we can forget about $g$; all we need to know is that $f$ is a
sufficiently generic $\mathbb{K}$-labelling of $\Delta\left(  p\right)  $
satisfying $f\left(  0\right)  =2$.

Let $\left(  i,k\right)  \in\Delta\left(  p\right)  $ be arbitrary. Then,
$i+k>p+1$ (since $\left(  i,k\right)  \in\Delta\left(  p\right)  $).
Consequently, $2p-\left(  i+k-1\right)  $ is a well-defined element of
$\left\{  1,2,...,p-1\right\}  $. Denote this element by $h$. Thus,
$h\in\left\{  1,2,...,p-1\right\}  $ and $i+k-1+h=2p$. Moreover, $\left(
k,i\right)  =\operatorname*{vrefl}v\in\Delta\left(  p\right)  $.

Let $v=\left(  p+1-k,p+1-i\right)  $. Then, $v=\operatorname*{hrefl}\left(
\left(  i,k\right)  \right)  \in\nabla\left(  p\right)  $ (since $\left(
i,k\right)  \in\Delta\left(  p\right)  $) and $\deg v=h$ (this follows by
simple computation). Moreover, $\operatorname*{hrefl}v=\left(  i,k\right)  $.

Applying Corollary \ref{cor.Delta.hrefl.f} to $\ell=h$, we obtain
$R_{\operatorname*{Rect}\left(  p,p\right)  }^{h}\left(  \operatorname*{wing}%
f\right)  =\operatorname*{wing}\left(  R_{\Delta\left(  p\right)  }%
^{h}f\right)  $, hence%
\begin{align}
&  \left(  R_{\operatorname*{Rect}\left(  p,p\right)  }^{h}\left(
\operatorname*{wing}f\right)  \right)  \left(  v\right) \nonumber\\
&  =\left(  \operatorname*{wing}\left(  R_{\Delta\left(  p\right)  }%
^{h}f\right)  \right)  \left(  v\right)  =\dfrac{1}{\left(  R_{\Delta\left(
p\right)  }^{p-\deg v}\left(  R_{\Delta\left(  p\right)  }^{h}f\right)
\right)  \left(  \operatorname*{hrefl}v\right)  }\nonumber\\
&  \ \ \ \ \ \ \ \ \ \ \left(  \text{by the definition of }%
\operatorname*{wing}\text{, since }v\in\nabla\left(  p\right)  \subseteq
\nabla\left(  p\right)  \cup\left\{  0\right\}  \right) \nonumber\\
&  =\dfrac{1}{\left(  R_{\Delta\left(  p\right)  }^{p-h}\left(  R_{\Delta
\left(  p\right)  }^{h}f\right)  \right)  \left(  \left(  i,k\right)  \right)
}\ \ \ \ \ \ \ \ \ \ \left(  \text{since }\deg v=h\text{ and }%
\operatorname*{hrefl}v=\left(  i,k\right)  \right) \nonumber\\
&  =\dfrac{1}{\left(  R_{\Delta\left(  p\right)  }^{p}f\right)  \left(
\left(  i,k\right)  \right)  }\ \ \ \ \ \ \ \ \ \ \left(  \text{since
}R_{\Delta\left(  p\right)  }^{p-h}\left(  R_{\Delta\left(  p\right)  }%
^{h}f\right)  =\underbrace{\left(  R_{\Delta\left(  p\right)  }^{p-h}\circ
R_{\Delta\left(  p\right)  }^{h}\right)  }_{=R_{\Delta\left(  p\right)  }^{p}%
}f=R_{\Delta\left(  p\right)  }^{p}f\right)  .
\label{pf.Delta.halfway.new.short.4}%
\end{align}
But Theorem \ref{thm.rect.antip.general} (applied to $p$,
$R_{\operatorname*{Rect}\left(  p,p\right)  }^{h}\left(  \operatorname*{wing}%
f\right)  $ and $\left(  k,i\right)  $ instead of $q$, $f$ and $\left(
i,k\right)  $) yields%
\begin{align*}
&  \left(  R_{\operatorname*{Rect}\left(  p,p\right)  }^{h}\left(
\operatorname*{wing}f\right)  \right)  \left(  \left(  p+1-k,p+1-i\right)
\right) \\
&  =\dfrac{\left(  R_{\operatorname*{Rect}\left(  p,p\right)  }^{h}\left(
\operatorname*{wing}f\right)  \right)  \left(  0\right)  \cdot\left(
R_{\operatorname*{Rect}\left(  p,p\right)  }^{h}\left(  \operatorname*{wing}%
f\right)  \right)  \left(  1\right)  }{\left(  R_{\operatorname*{Rect}\left(
p,p\right)  }^{i+k-1}\left(  R_{\operatorname*{Rect}\left(  p,p\right)  }%
^{h}\left(  \operatorname*{wing}f\right)  \right)  \right)  \left(  \left(
k,i\right)  \right)  }.
\end{align*}
Since $\left(  p+1-k,p+1-i\right)  =v$ and
\begin{align*}
R_{\operatorname*{Rect}\left(  p,p\right)  }^{i+k-1}\left(
R_{\operatorname*{Rect}\left(  p,p\right)  }^{h}\left(  \operatorname*{wing}%
f\right)  \right)   &  =\left(  \underbrace{R_{\operatorname*{Rect}\left(
p,p\right)  }^{i+k-1}\circ R_{\operatorname*{Rect}\left(  p,p\right)  }^{h}%
}_{\substack{=R_{\operatorname*{Rect}\left(  p,p\right)  }^{i+k-1+h}%
=R_{\operatorname*{Rect}\left(  p,p\right)  }^{2p}\\\text{(since
}i+k-1+h=2p\text{)}}}\right)  \left(  \operatorname*{wing}f\right) \\
&  =\underbrace{R_{\operatorname*{Rect}\left(  p,p\right)  }^{2p}%
}_{\substack{=\operatorname*{id}\\\text{(since Theorem \ref{thm.rect.ord}
(applied to }q=p\text{)}\\\text{yields }\operatorname*{ord}\left(
R_{\operatorname*{Rect}\left(  p,p\right)  }\right)  =p+p=2p\text{)}}}\left(
\operatorname*{wing}f\right)  =\operatorname*{wing}f,
\end{align*}
this equality rewrites as%
\[
\left(  R_{\operatorname*{Rect}\left(  p,p\right)  }^{h}\left(
\operatorname*{wing}f\right)  \right)  \left(  v\right)  =\dfrac{\left(
R_{\operatorname*{Rect}\left(  p,p\right)  }^{h}\left(  \operatorname*{wing}%
f\right)  \right)  \left(  0\right)  \cdot\left(  R_{\operatorname*{Rect}%
\left(  p,p\right)  }^{h}\left(  \operatorname*{wing}f\right)  \right)
\left(  1\right)  }{\left(  \operatorname*{wing}f\right)  \left(  \left(
k,i\right)  \right)  }.
\]
Since
\begin{align*}
&  \underbrace{\left(  R_{\operatorname*{Rect}\left(  p,p\right)  }^{h}\left(
\operatorname*{wing}f\right)  \right)  \left(  0\right)  }_{\substack{=\left(
\operatorname*{wing}f\right)  \left(  0\right)  \\\text{(by Corollary
\ref{cor.R.implicit.01})}}}\cdot\underbrace{\left(  R_{\operatorname*{Rect}%
\left(  p,p\right)  }^{h}\left(  \operatorname*{wing}f\right)  \right)
\left(  1\right)  }_{\substack{=\left(  \operatorname*{wing}f\right)  \left(
1\right)  \\\text{(by Corollary \ref{cor.R.implicit.01})}}}\\
&  =\underbrace{\left(  \operatorname*{wing}f\right)  \left(  0\right)
}_{\substack{=\dfrac{1}{\left(  R_{\Delta\left(  p\right)  }^{p-\deg
0}f\right)  \left(  \operatorname*{hrefl}0\right)  }\\\text{(by the definition
of }\operatorname*{wing}\text{)}}}\cdot\underbrace{\left(
\operatorname*{wing}f\right)  \left(  1\right)  }_{\substack{=f\left(
1\right)  \\\text{(by the definition of }\operatorname*{wing}\text{)}}}\\
&  =\dfrac{1}{\left(  R_{\Delta\left(  p\right)  }^{p-\deg0}f\right)  \left(
\operatorname*{hrefl}0\right)  }\cdot f\left(  1\right)  =1
\end{align*}
(since Corollary \ref{cor.R.implicit.01} yields $\left(  R_{\Delta\left(
p\right)  }^{p-\deg0}f\right)  \left(  \operatorname*{hrefl}0\right)
=f\left(  \operatorname*{hrefl}0\right)  =f\left(  1\right)  $), this
simplifies to
\[
\left(  R_{\operatorname*{Rect}\left(  p,p\right)  }^{h}\left(
\operatorname*{wing}f\right)  \right)  \left(  v\right)  =\dfrac{1}{\left(
\operatorname*{wing}f\right)  \left(  \left(  k,i\right)  \right)  }.
\]
Compared with (\ref{pf.Delta.halfway.new.short.4}), this yields $\dfrac
{1}{\left(  R_{\Delta\left(  p\right)  }^{p}f\right)  \left(  \left(
i,k\right)  \right)  }=\dfrac{1}{\left(  \operatorname*{wing}f\right)  \left(
\left(  k,i\right)  \right)  }$. Taking inverses in this equality, we get
\begin{align*}
\left(  R_{\Delta\left(  p\right)  }^{p}f\right)  \left(  \left(  i,k\right)
\right)   &  =\left(  \operatorname*{wing}f\right)  \left(  \left(
k,i\right)  \right)  =f\left(  \underbrace{\left(  k,i\right)  }%
_{=\operatorname*{vrefl}\left(  i,k\right)  }\right) \\
&  \ \ \ \ \ \ \ \ \ \ \left(  \text{by the definition of }%
\operatorname*{wing}\text{, since }\left(  k,i\right)  \in\Delta\left(
p\right)  \subseteq\Delta\left(  p\right)  \cup\left\{  1\right\}  \right) \\
&  =f\left(  \operatorname*{vrefl}\left(  i,k\right)  \right)  =\left(
\operatorname*{vrefl}\nolimits^{\ast}f\right)  \left(  \left(  i,k\right)
\right) \\
&  \ \ \ \ \ \ \ \ \ \ \left(  \text{since }\left(  \operatorname*{vrefl}%
\nolimits^{\ast}f\right)  \left(  \left(  i,k\right)  \right)  =f\left(
\operatorname*{vrefl}\left(  i,k\right)  \right)  \text{ by the definition of
}\operatorname*{vrefl}\nolimits^{\ast}\right)  .
\end{align*}

Now, we have shown this for \textbf{every} $\left(  i,k\right)  \in
\Delta\left(  p\right)  $. In other words, we have shown that $R_{\Delta
\left(  p\right)  }^{p}f=\operatorname*{vrefl}\nolimits^{\ast}f$.

Now, recall that $f=\left(  a_{0},a_{1},...,a_{p}\right)  \flat g$. Hence,%
\begin{equation}
R_{\Delta\left(  p\right)  }^{p}f=R_{\Delta\left(  p\right)  }^{p}\left(
\left(  a_{0},a_{1},...,a_{p}\right)  \flat g\right)  =\left(  a_{0}%
,a_{1},...,a_{p}\right)  \flat\left(  R_{\Delta\left(  p\right)  }%
^{p}g\right)  \label{pf.Delta.halfway.new.short.10}%
\end{equation}
(by Corollary \ref{cor.Rl.scalmult}, applied to $\Delta\left(  p\right)  $,
$p-1$ and $g$ instead of $P$, $n$ and $f$). On the other hand, $f=\left(
a_{0},a_{1},...,a_{p}\right)  \flat g$ yields%
\begin{equation}
\operatorname*{vrefl}\nolimits^{\ast}f=\operatorname*{vrefl}\nolimits^{\ast
}\left(  \left(  a_{0},a_{1},...,a_{p}\right)  \flat g\right)  =\left(
a_{0},a_{1},...,a_{p}\right)  \flat\left(  \operatorname*{vrefl}%
\nolimits^{\ast}g\right)  \label{pf.Delta.halfway.new.short.11}%
\end{equation}
(this is easy to check directly using the definitions of $\flat$ and
$\operatorname*{vrefl}\nolimits^{\ast}$, since $\operatorname*{vrefl}$
preserves degrees). In light of (\ref{pf.Delta.halfway.new.short.10}) and
(\ref{pf.Delta.halfway.new.short.11}), the equality $R_{\Delta\left(
p\right)  }^{p}f=\operatorname*{vrefl}\nolimits^{\ast}f$ becomes
\newline$\left(  a_{0},a_{1},...,a_{p}\right)  \flat\left(  R_{\Delta\left(
p\right)  }^{p}g\right)  =\left(  a_{0},a_{1},...,a_{p}\right)  \flat\left(
\operatorname*{vrefl}\nolimits^{\ast}g\right)  $. We can cancel the
\textquotedblleft$\left(  a_{0},a_{1},...,a_{p}\right)  \flat$%
\textquotedblright\ from both sides of this equation (since all $a_{i}$ are
nonzero), and thus obtain $R_{\Delta\left(  p\right)  }^{p}%
g=\operatorname*{vrefl}\nolimits^{\ast}g$. As we have seen, this is all we
need to prove Theorem \ref{thm.Delta.halfway}.
\end{proof}
\end{vershort}

\begin{verlong}
\begin{proof}
[Proof of Theorem \ref{thm.Delta.halfway} (sketched).]The result that we are
striving to prove is a collection of identities between rational functions,
hence boils down to a collection of polynomial identities in the labels of an
arbitrary $\mathbb{K}$-labelling of $\Delta\left(  p\right)  $. Therefore, it
is enough to prove it in the case when $\mathbb{K}$ is a field of rational
functions in finitely many variables over $\mathbb{Q}$. So let us WLOG assume
that $\mathbb{K}$ is a field of rational functions in finitely many variables
over $\mathbb{Q}$. Then, the characteristic of $\mathbb{K}$ is $\neq2$ (it is
$0$ indeed), so that we can apply Lemma \ref{lem.Delta.hrefl} and Corollary
\ref{cor.Delta.hrefl.f}.

Consider the maps $\operatorname*{hrefl}$, $\operatorname*{wing}$,
$\operatorname*{vrefl}$ and $\operatorname*{vrefl}\nolimits^{\ast}$ defined in
Lemma \ref{lem.Delta.hrefl}. We are going to prove that%
\[
R_{\Delta\left(  p\right)  }^{p}=\operatorname*{vrefl}\nolimits^{\ast}%
\]
as rational maps $\mathbb{K}^{\widehat{\Delta\left(  p\right)  }%
}\dashrightarrow\mathbb{K}^{\widehat{\Delta\left(  p\right)  }}$.

In fact, let $g\in\mathbb{K}^{\widehat{\Delta\left(  p\right)  }}$ be any
sufficiently generic zero-free labelling of $\Delta\left(  p\right)  $. Here,
\textquotedblleft sufficiently generic\textquotedblright\ means
\textquotedblleft generic enough for all terms appearing in the present proof
to be well-defined and for all genericity conditions to be
satisfied\textquotedblright.

We will use the notation introduced in Definition \ref{def.bemol}. Hence, a
$\mathbb{K}$-labelling $\left(  a_{0},a_{1},...,a_{\left(  p-1\right)
+1}\right)  \flat g$ of $\Delta\left(  p\right)  $ is well-defined for every
$\left(  \left(  p-1\right)  +2\right)  $-tuple $\left(  a_{0},a_{1}%
,...,a_{\left(  p-1\right)  +1}\right)  \in\left(  \mathbb{K}^{\times}\right)
^{\left(  p-1\right)  +2}$. Since $\left(  p-1\right)  +1=p$ and $\left(
p-1\right)  +2=p+1$, this rewrites as follows: A $\mathbb{K}$-labelling
$\left(  a_{0},a_{1},...,a_{p}\right)  \flat g$ of $\Delta\left(  p\right)  $
is well-defined for every $\left(  p+1\right)  $-tuple $\left(  a_{0}%
,a_{1},...,a_{p}\right)  \in\left(  \mathbb{K}^{\times}\right)  ^{p+1}$. It is
easy to see that every $\left(  a_{0},a_{1},...,a_{p}\right)  \in\left(
\mathbb{K}^{\times}\right)  ^{p+1}$ satisfies%
\begin{equation}
\operatorname*{vrefl}\nolimits^{\ast}\left(  \left(  a_{0},a_{1}%
,...,a_{p}\right)  \flat g\right)  =\left(  a_{0},a_{1},...,a_{p}\right)
\flat\left(  \operatorname*{vrefl}\nolimits^{\ast}g\right)
\label{pf.Delta.halfway.new.bemol-vrefl}%
\end{equation}
\footnote{\textit{Proof of (\ref{pf.Delta.halfway.new.bemol-vrefl}):} Let
$\left(  a_{0},a_{1},...,a_{p}\right)  \in\left(  \mathbb{K}^{\times}\right)
^{p+1}$. Let $v\in\widehat{\Delta\left(  p\right)  }$. Then, $\deg\left(
\operatorname*{vrefl}v\right)  =\deg v$ (this is easy to prove using just the
definition of the map $\operatorname*{vrefl}:\widehat{\Delta\left(  p\right)
}\rightarrow\widehat{\Delta\left(  p\right)  }$). Now,
\begin{align*}
&  \left(  \operatorname*{vrefl}\nolimits^{\ast}\left(  \left(  a_{0}%
,a_{1},...,a_{p}\right)  \flat g\right)  \right)  \left(  v\right) \\
&  =\left(  \left(  a_{0},a_{1},...,a_{p}\right)  \flat g\right)  \left(
\operatorname*{vrefl}v\right)  \ \ \ \ \ \ \ \ \ \ \left(  \text{by the
definition of }\operatorname*{vrefl}\nolimits^{\ast}\left(  \left(
a_{0},a_{1},...,a_{p}\right)  \flat g\right)  \right) \\
&  =a_{\deg\left(  \operatorname*{vrefl}v\right)  }\cdot g\left(
\operatorname*{vrefl}v\right)  \ \ \ \ \ \ \ \ \ \ \left(  \text{by the
definition of }\left(  a_{0},a_{1},...,a_{p}\right)  \flat g\right) \\
&  =a_{\deg v}\cdot g\left(  \operatorname*{vrefl}v\right)
\ \ \ \ \ \ \ \ \ \ \left(  \text{since }\deg\left(  \operatorname*{vrefl}%
v\right)  =\deg v\right)  .
\end{align*}
Compared with%
\begin{align*}
&  \left(  \left(  a_{0},a_{1},...,a_{p}\right)  \flat\left(
\operatorname*{vrefl}\nolimits^{\ast}g\right)  \right)  \left(  v\right) \\
&  =a_{\deg v}\cdot\underbrace{\left(  \operatorname*{vrefl}\nolimits^{\ast
}g\right)  \left(  v\right)  }_{\substack{=g\left(  \operatorname*{vrefl}%
v\right)  \\\text{(by the definition of }\operatorname*{vrefl}\nolimits^{\ast
}g\text{)}}}\ \ \ \ \ \ \ \ \ \ \left(  \text{by the definition of }\left(
a_{0},a_{1},...,a_{p}\right)  \flat\left(  \operatorname*{vrefl}%
\nolimits^{\ast}g\right)  \right) \\
&  =a_{\deg v}\cdot g\left(  \operatorname*{vrefl}v\right)  ,
\end{align*}
this yields%
\[
\left(  \operatorname*{vrefl}\nolimits^{\ast}\left(  \left(  a_{0}%
,a_{1},...,a_{p}\right)  \flat g\right)  \right)  \left(  v\right)  =\left(
\left(  a_{0},a_{1},...,a_{p}\right)  \flat\left(  \operatorname*{vrefl}%
\nolimits^{\ast}g\right)  \right)  \left(  v\right)  .
\]
\par
Now, forget that we fixed $v$. We thus have proven that $\left(
\operatorname*{vrefl}\nolimits^{\ast}\left(  \left(  a_{0},a_{1}%
,...,a_{p}\right)  \flat g\right)  \right)  \left(  v\right)  =\left(  \left(
a_{0},a_{1},...,a_{p}\right)  \flat\left(  \operatorname*{vrefl}%
\nolimits^{\ast}g\right)  \right)  \left(  v\right)  $ for every
$v\in\widehat{\Delta\left(  p\right)  }$. In other words,
\[
\operatorname*{vrefl}\nolimits^{\ast}\left(  \left(  a_{0},a_{1}%
,...,a_{p}\right)  \flat g\right)  =\left(  a_{0},a_{1},...,a_{p}\right)
\flat\left(  \operatorname*{vrefl}\nolimits^{\ast}g\right)  .
\]
This proves (\ref{pf.Delta.halfway.new.bemol-vrefl}).}.

Clearly, there exists a $\left(  p+1\right)  $-tuple $\left(  a_{0}%
,a_{1},...,a_{p}\right)  \in\left(  \mathbb{K}^{\times}\right)  ^{p+1}$ such
that $\left(  \left(  a_{0},a_{1},...,a_{p}\right)  \flat g\right)  \left(
0\right)  =2$ \ \ \ \ \footnote{\textit{Proof.} The labelling $g$ is
zero-free. Thus, there must exist some $v\in\widehat{P}_{0}$ satisfying
$g\left(  v\right)  \neq0$. This $v$ must be $0$ (since the only element of
$\widehat{P}_{0}$ is $0$). Hence, $0$ satisfies $g\left(  0\right)  \neq0$.
Thus, $\dfrac{2}{g\left(  0\right)  }$ is a well-defined element of
$\mathbb{K}$; it further lies in $\mathbb{K}^{\times}$ because the
characteristic of $\mathbb{K}$ is $\neq2$.
\par
Now, define a $\left(  p+1\right)  $-tuple $\left(  a_{0},a_{1},...,a_{p}%
\right)  \in\left(  \mathbb{K}^{\times}\right)  ^{p+1}$ by%
\[
\left(  a_{i}=\dfrac{2}{g\left(  0\right)  }\ \ \ \ \ \ \ \ \ \ \text{for
every }i\in\left\{  0,1,...,p\right\}  \right)  .
\]
Then, the definition of $\left(  a_{0},a_{1},...,a_{p}\right)  \flat g$
yields
\[
\left(  \left(  a_{0},a_{1},...,a_{p}\right)  \flat g\right)  \left(
0\right)  =\underbrace{a_{\deg0}}_{\substack{=\dfrac{2}{g\left(  0\right)
}\\\text{(by the definition of }a_{\deg0}\text{)}}}\cdot g\left(  0\right)
=\dfrac{2}{g\left(  0\right)  }\cdot g\left(  0\right)  =2.
\]
\par
Thus, the $\left(  p+1\right)  $-tuple $\left(  a_{0},a_{1},...,a_{p}\right)
\in\left(  \mathbb{K}^{\times}\right)  ^{p+1}$ we constructed satisfies
$\left(  \left(  a_{0},a_{1},...,a_{p}\right)  \flat g\right)  \left(
0\right)  =2$. Hence, we have shown that there exists a $\left(  p+1\right)
$-tuple $\left(  a_{0},a_{1},...,a_{p}\right)  \in\left(  \mathbb{K}^{\times
}\right)  ^{p+1}$ such that $\left(  \left(  a_{0},a_{1},...,a_{p}\right)
\flat g\right)  \left(  0\right)  =2$, qed.}. Consider such a $\left(
p+1\right)  $-tuple $\left(  a_{0},a_{1},...,a_{p}\right)  \in\left(
\mathbb{K}^{\times}\right)  ^{p+1}$.

Let $f$ denote the $\mathbb{K}$-labelling $\left(  a_{0},a_{1},...,a_{p}%
\right)  \flat g$. Then, $f=\left(  a_{0},a_{1},...,a_{p}\right)  \flat g$, so
that $f\left(  0\right)  =\left(  \left(  a_{0},a_{1},...,a_{p}\right)  \flat
g\right)  \left(  0\right)  =2$.

We are now going to prove $R_{\Delta\left(  p\right)  }^{p}%
f=\operatorname*{vrefl}\nolimits^{\ast}f$. While we will be proving
$R_{\Delta\left(  p\right)  }^{p}f=\operatorname*{vrefl}\nolimits^{\ast}f$, we
can forget about $g$; all we need to know is that $f$ is a sufficiently
generic $\mathbb{K}$-labelling of $\Delta\left(  p\right)  $ satisfying
$f\left(  0\right)  =2$.

It is easy to see that $\left(  R_{\Delta\left(  p\right)  }^{p}f\right)
\left(  1\right)  =\left(  \operatorname*{vrefl}\nolimits^{\ast}f\right)
\left(  1\right)  $ and $\left(  R_{\Delta\left(  p\right)  }^{p}f\right)
\left(  0\right)  =\left(  \operatorname*{vrefl}\nolimits^{\ast}f\right)
\left(  0\right)  $.\ \ \ \ \footnote{\textit{Proof.} Corollary
\ref{cor.R.implicit.01} yields $\left(  R_{\Delta\left(  p\right)  }%
^{p}f\right)  \left(  1\right)  =f\left(  1\right)  $ and $\left(
R_{\Delta\left(  p\right)  }^{p}f\right)  \left(  0\right)  =f\left(
0\right)  $. On the other hand, since $\operatorname*{vrefl}$ leaves the
elements $0$ and $1$ fixed, we have $\left(  \operatorname*{vrefl}%
\nolimits^{\ast}f\right)  \left(  1\right)  =f\left(  1\right)  $ and $\left(
\operatorname*{vrefl}\nolimits^{\ast}f\right)  \left(  0\right)  =f\left(
0\right)  $. Thus, $\left(  R_{\Delta\left(  p\right)  }^{p}f\right)  \left(
1\right)  =f\left(  1\right)  =\left(  \operatorname*{vrefl}\nolimits^{\ast
}f\right)  \left(  1\right)  $ and $\left(  R_{\Delta\left(  p\right)  }%
^{p}f\right)  \left(  0\right)  =f\left(  0\right)  =\left(
\operatorname*{vrefl}\nolimits^{\ast}f\right)  \left(  0\right)  $, qed.}
Also,%
\begin{equation}
\left(  \operatorname*{wing}f\right)  \left(  0\right)  \cdot\left(
\operatorname*{wing}f\right)  \left(  1\right)  =1.
\label{pf.Delta.halfway.new.0*1}%
\end{equation}
\footnote{\textit{Proof of (\ref{pf.Delta.halfway.new.0*1}):} By the
definition of $\operatorname*{wing}$, we have $\left(  \operatorname*{wing}%
f\right)  \left(  1\right)  =f\left(  1\right)  $ (since $1\in\Delta\left(
p\right)  \cup\left\{  1\right\}  $) and%
\begin{align*}
\left(  \operatorname*{wing}f\right)  \left(  0\right)   &  =\dfrac{1}{\left(
R_{\Delta\left(  p\right)  }^{p-\deg0}f\right)  \left(  \operatorname*{hrefl}%
0\right)  }\ \ \ \ \ \ \ \ \ \ \left(  \text{since }0\in\nabla\left(
p\right)  \cup\left\{  0\right\}  \right) \\
&  =\dfrac{1}{\left(  R_{\Delta\left(  p\right)  }^{p-0}f\right)  \left(
1\right)  }\ \ \ \ \ \ \ \ \ \ \left(  \text{since }\operatorname*{hrefl}%
0=1\text{ and }\deg0=0\right) \\
&  =\dfrac{1}{\left(  R_{\Delta\left(  p\right)  }^{p}f\right)  \left(
1\right)  }=\dfrac{1}{f\left(  1\right)  }\ \ \ \ \ \ \ \ \ \ \left(
\text{since Corollary \ref{cor.R.implicit.01} yields }\left(  R_{\Delta\left(
p\right)  }^{p}f\right)  \left(  1\right)  =f\left(  1\right)  \right)  .
\end{align*}
Hence,%
\[
\underbrace{\left(  \operatorname*{wing}f\right)  \left(  0\right)  }%
_{=\dfrac{1}{f\left(  1\right)  }}\cdot\underbrace{\left(
\operatorname*{wing}f\right)  \left(  1\right)  }_{=f\left(  1\right)
}=\dfrac{1}{f\left(  1\right)  }\cdot f\left(  1\right)  =1,
\]
and thus (\ref{pf.Delta.halfway.new.0*1}) is proven.}

Let $\left(  i,k\right)  \in\Delta\left(  p\right)  $ be arbitrary. Then,
$i+k>p+1$ (since $\left(  i,k\right)  \in\Delta\left(  p\right)  $).
Consequently, $2p-\left(  i+k-1\right)  $ is a well-defined element of
$\left\{  1,2,...,p-1\right\}  $. Denote this element by $h$. Thus,
$h=2p-\left(  i+k-1\right)  \in\left\{  1,2,...,p-1\right\}  $. Moreover, from
$h=2p-\left(  i+k-1\right)  $, we obtain $i+k-1+h=2p$. Also, since $\left(
i,k\right)  =v$, we have $\left(  k,i\right)  =\operatorname*{vrefl}v\in
\Delta\left(  p\right)  $ (since $v\in\Delta\left(  p\right)  $).

Let $v=\left(  p+1-k,p+1-i\right)  $. Then, $v=\operatorname*{hrefl}\left(
\left(  i,k\right)  \right)  \in\nabla\left(  p\right)  $ (since $\left(
i,k\right)  \in\Delta\left(  p\right)  $) and
\[
\deg v=\left(  p+1-k\right)  +\left(  p+1-i\right)  -1=2p-\left(
i+k-1\right)  =h.
\]
Furthermore, since $v=\operatorname*{hrefl}\left(  \left(  i,k\right)
\right)  $, we have $\operatorname*{hrefl}v=\left(  i,k\right)  $ (because
$\operatorname*{hrefl}$ is an involution).

Applying Corollary \ref{cor.Delta.hrefl.f} to $\ell=h$, we obtain
$R_{\operatorname*{Rect}\left(  p,p\right)  }^{h}\left(  \operatorname*{wing}%
f\right)  =\operatorname*{wing}\left(  R_{\Delta\left(  p\right)  }%
^{h}f\right)  $, hence%
\begin{align}
&  \left(  R_{\operatorname*{Rect}\left(  p,p\right)  }^{h}\left(
\operatorname*{wing}f\right)  \right)  \left(  v\right) \nonumber\\
&  =\left(  \operatorname*{wing}\left(  R_{\Delta\left(  p\right)  }%
^{h}f\right)  \right)  \left(  v\right)  =\dfrac{1}{\left(  R_{\Delta\left(
p\right)  }^{p-\deg v}\left(  R_{\Delta\left(  p\right)  }^{h}f\right)
\right)  \left(  \operatorname*{hrefl}v\right)  }\nonumber\\
&  \ \ \ \ \ \ \ \ \ \ \left(  \text{by the definition of }%
\operatorname*{wing}\text{, since }v\in\nabla\left(  p\right)  \subseteq
\nabla\left(  p\right)  \cup\left\{  0\right\}  \right) \nonumber\\
&  =\dfrac{1}{\left(  R_{\Delta\left(  p\right)  }^{p-h}\left(  R_{\Delta
\left(  p\right)  }^{h}f\right)  \right)  \left(  \left(  i,k\right)  \right)
}\ \ \ \ \ \ \ \ \ \ \left(  \text{since }\deg v=h\text{ and }%
\operatorname*{hrefl}v=\left(  i,k\right)  \right) \nonumber\\
&  =\dfrac{1}{\left(  R_{\Delta\left(  p\right)  }^{p}f\right)  \left(
\left(  i,k\right)  \right)  }\ \ \ \ \ \ \ \ \ \ \left(  \text{since
}R_{\Delta\left(  p\right)  }^{p-h}\left(  R_{\Delta\left(  p\right)  }%
^{h}f\right)  =\underbrace{\left(  R_{\Delta\left(  p\right)  }^{p-h}\circ
R_{\Delta\left(  p\right)  }^{h}\right)  }_{=R_{\Delta\left(  p\right)  }^{p}%
}f=R_{\Delta\left(  p\right)  }^{p}f\right)  . \label{pf.Delta.halfway.new.4}%
\end{align}
But Theorem \ref{thm.rect.antip.general} (applied to $p$,
$R_{\operatorname*{Rect}\left(  p,p\right)  }^{h}\left(  \operatorname*{wing}%
f\right)  $ and $\left(  k,i\right)  $ instead of $q$, $f$ and $\left(
i,k\right)  $) yields%
\begin{align*}
&  \left(  R_{\operatorname*{Rect}\left(  p,p\right)  }^{h}\left(
\operatorname*{wing}f\right)  \right)  \left(  \left(  p+1-k,p+1-i\right)
\right) \\
&  =\dfrac{\left(  R_{\operatorname*{Rect}\left(  p,p\right)  }^{h}\left(
\operatorname*{wing}f\right)  \right)  \left(  0\right)  \cdot\left(
R_{\operatorname*{Rect}\left(  p,p\right)  }^{h}\left(  \operatorname*{wing}%
f\right)  \right)  \left(  1\right)  }{\left(  R_{\operatorname*{Rect}\left(
p,p\right)  }^{i+k-1}\left(  R_{\operatorname*{Rect}\left(  p,p\right)  }%
^{h}\left(  \operatorname*{wing}f\right)  \right)  \right)  \left(  \left(
k,i\right)  \right)  }.
\end{align*}
Since $\left(  p+1-k,p+1-i\right)  =v$ and
\begin{align*}
R_{\operatorname*{Rect}\left(  p,p\right)  }^{i+k-1}\left(
R_{\operatorname*{Rect}\left(  p,p\right)  }^{h}\left(  \operatorname*{wing}%
f\right)  \right)   &  =\left(  \underbrace{R_{\operatorname*{Rect}\left(
p,p\right)  }^{i+k-1}\circ R_{\operatorname*{Rect}\left(  p,p\right)  }^{h}%
}_{\substack{=R_{\operatorname*{Rect}\left(  p,p\right)  }^{i+k-1+h}%
=R_{\operatorname*{Rect}\left(  p,p\right)  }^{2p}\\\text{(since
}i+k-1+h=2p\text{)}}}\right)  \left(  \operatorname*{wing}f\right) \\
&  =\underbrace{R_{\operatorname*{Rect}\left(  p,p\right)  }^{2p}%
}_{\substack{=\operatorname*{id}\\\text{(since Theorem \ref{thm.rect.ord}
(applied to }q=p\text{)}\\\text{yields }\operatorname*{ord}\left(
R_{\operatorname*{Rect}\left(  p,p\right)  }\right)  =p+p=2p\text{)}}}\left(
\operatorname*{wing}f\right) \\
&  =\operatorname*{id}\left(  \operatorname*{wing}f\right)
=\operatorname*{wing}f,
\end{align*}
this equality rewrites as%
\[
\left(  R_{\operatorname*{Rect}\left(  p,p\right)  }^{h}\left(
\operatorname*{wing}f\right)  \right)  \left(  v\right)  =\dfrac{\left(
R_{\operatorname*{Rect}\left(  p,p\right)  }^{h}\left(  \operatorname*{wing}%
f\right)  \right)  \left(  0\right)  \cdot\left(  R_{\operatorname*{Rect}%
\left(  p,p\right)  }^{h}\left(  \operatorname*{wing}f\right)  \right)
\left(  1\right)  }{\left(  \operatorname*{wing}f\right)  \left(  \left(
k,i\right)  \right)  }.
\]
Since
\[
\left(  R_{\operatorname*{Rect}\left(  p,p\right)  }^{h}\left(
\operatorname*{wing}f\right)  \right)  \left(  0\right)  =\left(
\operatorname*{wing}f\right)  \left(  0\right)  \ \ \ \ \ \ \ \ \ \ \left(
\text{by Corollary \ref{cor.R.implicit.01}}\right)
\]
and%
\[
\left(  R_{\operatorname*{Rect}\left(  p,p\right)  }^{h}\left(
\operatorname*{wing}f\right)  \right)  \left(  1\right)  =\left(
\operatorname*{wing}f\right)  \left(  1\right)  \ \ \ \ \ \ \ \ \ \ \left(
\text{by Corollary \ref{cor.R.implicit.01}}\right)  ,
\]
this equality rewrites as
\[
\left(  R_{\operatorname*{Rect}\left(  p,p\right)  }^{h}\left(
\operatorname*{wing}f\right)  \right)  \left(  v\right)  =\dfrac{\left(
\operatorname*{wing}f\right)  \left(  0\right)  \cdot\left(
\operatorname*{wing}f\right)  \left(  1\right)  }{\left(  \operatorname*{wing}%
f\right)  \left(  \left(  k,i\right)  \right)  }.
\]
Since $\left(  \operatorname*{wing}f\right)  \left(  0\right)  \cdot\left(
\operatorname*{wing}f\right)  \left(  1\right)  =1$ (by
(\ref{pf.Delta.halfway.new.0*1})), this rewrites as $\left(
R_{\operatorname*{Rect}\left(  p,p\right)  }^{h}\left(  \operatorname*{wing}%
f\right)  \right)  \left(  v\right)  =\dfrac{1}{\left(  \operatorname*{wing}%
f\right)  \left(  \left(  k,i\right)  \right)  }$. Compared with
(\ref{pf.Delta.halfway.new.4}), this yields $\dfrac{1}{\left(  R_{\Delta
\left(  p\right)  }^{p}f\right)  \left(  \left(  i,k\right)  \right)  }%
=\dfrac{1}{\left(  \operatorname*{wing}f\right)  \left(  \left(  k,i\right)
\right)  }$. Taking inverses in this equality, we get
\begin{align*}
\left(  R_{\Delta\left(  p\right)  }^{p}f\right)  \left(  \left(  i,k\right)
\right)   &  =\left(  \operatorname*{wing}f\right)  \left(  \left(
k,i\right)  \right)  =f\left(  \underbrace{\left(  k,i\right)  }%
_{=\operatorname*{vrefl}\left(  i,k\right)  }\right) \\
&  \ \ \ \ \ \ \ \ \ \ \left(  \text{by the definition of }%
\operatorname*{wing}\text{, since }\left(  k,i\right)  \in\Delta\left(
p\right)  \subseteq\Delta\left(  p\right)  \cup\left\{  1\right\}  \right) \\
&  =f\left(  \operatorname*{vrefl}\left(  i,k\right)  \right)  =\left(
\operatorname*{vrefl}\nolimits^{\ast}f\right)  \left(  \left(  i,k\right)
\right) \\
&  \ \ \ \ \ \ \ \ \ \ \left(  \text{since }\left(  \operatorname*{vrefl}%
\nolimits^{\ast}f\right)  \left(  \left(  i,k\right)  \right)  =f\left(
\operatorname*{vrefl}\left(  i,k\right)  \right)  \text{ by the definition of
}\operatorname*{vrefl}\nolimits^{\ast}\right)  .
\end{align*}

Now, forget that we fixed $\left(  i,k\right)  $. We thus have proven that
every $\left(  i,k\right)  \in\Delta\left(  p\right)  $ satisfies $\left(
R_{\Delta\left(  p\right)  }^{p}f\right)  \left(  \left(  i,k\right)  \right)
=\left(  \operatorname*{vrefl}\nolimits^{\ast}f\right)  \left(  \left(
i,k\right)  \right)  $. Hence, the $\mathbb{K}$-labellings $R_{\Delta\left(
p\right)  }^{p}f$ and $\operatorname*{vrefl}\nolimits^{\ast}f$ are equal on
each element of $\Delta\left(  p\right)  $. Since these two $\mathbb{K}%
$-labellings are also equal on each of the elements $0$ and $1$ (because
$\left(  R_{\Delta\left(  p\right)  }^{p}f\right)  \left(  1\right)  =\left(
\operatorname*{vrefl}\nolimits^{\ast}f\right)  \left(  1\right)  $ and
$\left(  R_{\Delta\left(  p\right)  }^{p}f\right)  \left(  0\right)  =\left(
\operatorname*{vrefl}\nolimits^{\ast}f\right)  \left(  0\right)  $), this
yields that the $\mathbb{K}$-labellings $R_{\Delta\left(  p\right)  }^{p}f$
and $\operatorname*{vrefl}\nolimits^{\ast}f$ are equal on each element of
$\Delta\left(  p\right)  \cup\left\{  0,1\right\}  =\widehat{\Delta\left(
p\right)  }$. In other words, $R_{\Delta\left(  p\right)  }^{p}%
f=\operatorname*{vrefl}\nolimits^{\ast}f$.

We now need to show that $R_{\Delta\left(  p\right)  }^{p}%
g=\operatorname*{vrefl}\nolimits^{\ast}g$.

Corollary \ref{cor.Rl.scalmult} (applied to $\Delta\left(  p\right)  $, $p-1$
and $g$ instead of $P$, $n$ and $f$) yields%
\[
R_{\Delta\left(  p\right)  }^{\left(  p-1\right)  +1}\left(  \left(
a_{0},a_{1},...,a_{\left(  p-1\right)  +1}\right)  \flat g\right)  =\left(
a_{0},a_{1},...,a_{\left(  p-1\right)  +1}\right)  \flat\left(  R_{\Delta
\left(  p\right)  }^{\left(  p-1\right)  +1}g\right)  .
\]
Hence,%
\begin{align}
&  \left(  a_{0},a_{1},...,a_{\left(  p-1\right)  +1}\right)  \flat\left(
R_{\Delta\left(  p\right)  }^{\left(  p-1\right)  +1}g\right) \nonumber\\
&  =R_{\Delta\left(  p\right)  }^{\left(  p-1\right)  +1}\left(  \left(
a_{0},a_{1},...,a_{\left(  p-1\right)  +1}\right)  \flat g\right) \nonumber\\
&  =R_{\Delta\left(  p\right)  }^{p}\left(  \underbrace{\left(  a_{0}%
,a_{1},...,a_{p}\right)  \flat g}_{=f}\right)  \ \ \ \ \ \ \ \ \ \ \left(
\text{since }\left(  p-1\right)  +1=p\right) \nonumber\\
&  =R_{\Delta\left(  p\right)  }^{p}f=\operatorname*{vrefl}\nolimits^{\ast}f
\label{pf.Delta.halfway.new.15}%
\end{align}
(by what we have already proven).

On the other hand,
\begin{align*}
&  \left(  a_{0},a_{1},...,a_{\left(  p-1\right)  +1}\right)  \flat\left(
\operatorname*{vrefl}\nolimits^{\ast}g\right) \\
&  =\left(  a_{0},a_{1},...,a_{p}\right)  \flat\left(  \operatorname*{vrefl}%
\nolimits^{\ast}g\right)  \ \ \ \ \ \ \ \ \ \ \left(  \text{since }\left(
p-1\right)  +1=p\right) \\
&  =\operatorname*{vrefl}\nolimits^{\ast}\left(  \underbrace{\left(
a_{0},a_{1},...,a_{p}\right)  \flat g}_{=f}\right)
\ \ \ \ \ \ \ \ \ \ \left(  \text{by (\ref{pf.Delta.halfway.new.bemol-vrefl}%
)}\right) \\
&  =\operatorname*{vrefl}\nolimits^{\ast}f=\left(  a_{0},a_{1},...,a_{\left(
p-1\right)  +1}\right)  \flat\left(  R_{\Delta\left(  p\right)  }^{\left(
p-1\right)  +1}g\right)  \ \ \ \ \ \ \ \ \ \ \left(  \text{by
(\ref{pf.Delta.halfway.new.15})}\right)  .
\end{align*}
Hence, Lemma \ref{lem.bemol.cancel} (applied to $\Delta\left(  p\right)  $,
$p-1$, $\operatorname*{vrefl}\nolimits^{\ast}g$ and $R_{\Delta\left(
p\right)  }^{\left(  p-1\right)  +1}g$ instead of $P$, $n$, $f$ and $g$)
yields $\operatorname*{vrefl}\nolimits^{\ast}g=R_{\Delta\left(  p\right)
}^{\left(  p-1\right)  +1}g=R_{\Delta\left(  p\right)  }^{p}g$. In other
words, $R_{\Delta\left(  p\right)  }^{p}g=\operatorname*{vrefl}\nolimits^{\ast
}g$.

Now, forget that we fixed $g$. We thus have shown that $R_{\Delta\left(
p\right)  }^{p}g=\operatorname*{vrefl}\nolimits^{\ast}g$ for every
sufficiently generic $g$. Hence, $R_{\Delta\left(  p\right)  }^{p}%
=\operatorname*{vrefl}\nolimits^{\ast}$. Thus, for every $\left(  i,k\right)
\in\Delta\left(  p\right)  $ and every $f\in\mathbb{K}^{\widehat{\Delta\left(
p\right)  }}$, we have%
\begin{align*}
\underbrace{\left(  R_{\Delta\left(  p\right)  }^{p}f\right)  }%
_{=\operatorname*{vrefl}\nolimits^{\ast}}\left(  \left(  i,k\right)  \right)
&  =\left(  \operatorname*{vrefl}\nolimits^{\ast}f\right)  \left(  \left(
i,k\right)  \right)  =f\left(  \underbrace{\operatorname*{vrefl}\left(
i,k\right)  }_{=\left(  k,i\right)  }\right)  \ \ \ \ \ \ \ \ \ \ \left(
\text{by the definition of }\operatorname*{vrefl}\nolimits^{\ast}\right) \\
&  =f\left(  \left(  k,i\right)  \right)  .
\end{align*}
This proves Theorem \ref{thm.Delta.halfway}.
\end{proof}
\end{verlong}

\begin{noncompile}
\textbf{[The following alternative proof of Theorem \ref{thm.Delta.halfway}
has been stuck into a ``noncompile'' environment because it is longer than the
proof given above without providing any additional insight. It is (to my
knowledge, at least) not wrong.]}

\begin{vershort}
\begin{proof}
[Second proof of Theorem \ref{thm.Delta.halfway} (sketched).]The result that
we are striving to prove is a collection of identities between rational
functions, hence boils down to a collection of polynomial identities in the
labels of an arbitrary $\mathbb{K}$-labelling of $\Delta\left(  p\right)  $.
Therefore, it is enough to prove it in the case when $\mathbb{K}$ is a field
of rational functions in finitely many variables over $\mathbb{Q}$. So let us
WLOG assume that $\mathbb{K}$ is a field of rational functions in finitely
many variables over $\mathbb{Q}$. Then, the characteristic of $\mathbb{K}$ is
$\neq2$ (it is $0$ indeed), so that we can apply Lemma \ref{lem.Delta.hrefl}
and Corollary \ref{cor.Delta.hrefl.f}.

Consider the maps $\operatorname*{hrefl}$, $\operatorname*{wing}$,
$\operatorname*{vrefl}$ and $\operatorname*{vrefl}\nolimits^{\ast}$ defined in
Lemma \ref{lem.Delta.hrefl}. Clearly, it will be enough to prove that%
\[
R_{\Delta\left(  p\right)  }^{p}=\operatorname*{vrefl}\nolimits^{\ast}%
\]
as rational maps $\mathbb{K}^{\widehat{\Delta\left(  p\right)  }%
}\dashrightarrow\mathbb{K}^{\widehat{\Delta\left(  p\right)  }}$. In other
words, it will be enough to prove that $R_{\Delta\left(  p\right)  }%
^{p}g=\operatorname*{vrefl}\nolimits^{\ast}g$ for almost every $g\in
\mathbb{K}^{\widehat{\Delta\left(  p\right)  }}$.

So let $g\in\mathbb{K}^{\widehat{\Delta\left(  p\right)  }}$ be any
sufficiently generic zero-free labelling of $\Delta\left(  p\right)  $. We
need to show that $R_{\Delta\left(  p\right)  }^{p}g=\operatorname*{vrefl}%
\nolimits^{\ast}g$. According to Proposition \ref{prop.reconstruct} (applied
to $p-1$, $\Delta\left(  p\right)  $, $R_{\Delta\left(  p\right)  }^{p}g$ and
$\operatorname*{vrefl}\nolimits^{\ast}g$ instead of $n$, $P$, $f$ and $g$),
this will be achieved if we can prove the following four claims:

\textit{Claim 1:} No $i\in\left\{  0,1,...,p-1\right\}  $ satisfies
$\mathbf{w}_{i}\left(  R_{\Delta\left(  p\right)  }^{p}g\right)  =0$.

\textit{Claim 2:} We have
\begin{align*}
&  \left(  \mathbf{w}_{0}\left(  R_{\Delta\left(  p\right)  }^{p}g\right)
,\mathbf{w}_{1}\left(  R_{\Delta\left(  p\right)  }^{p}g\right)
,...,\mathbf{w}_{p-1}\left(  R_{\Delta\left(  p\right)  }^{p}g\right)  \right)
\\
&  =\left(  \mathbf{w}_{0}\left(  \operatorname*{vrefl}\nolimits^{\ast
}g\right)  ,\mathbf{w}_{1}\left(  \operatorname*{vrefl}\nolimits^{\ast
}g\right)  ,...,\mathbf{w}_{p-1}\left(  \operatorname*{vrefl}\nolimits^{\ast
}g\right)  \right)  .
\end{align*}

\textit{Claim 3:} We have $\pi\left(  R_{\Delta\left(  p\right)  }%
^{p}g\right)  =\pi\left(  \operatorname*{vrefl}\nolimits^{\ast}g\right)  $.

\textit{Claim 4:} We have $\left(  R_{\Delta\left(  p\right)  }^{p}g\right)
\left(  0\right)  =\left(  \operatorname*{vrefl}\nolimits^{\ast}g\right)
\left(  0\right)  $.

Let us first check that Claims 2 and 4 are true. Indeed, Claim 4 obviously
holds. Claim 2 follows from%
\begin{align}
&  \left(  \mathbf{w}_{0}\left(  R_{\Delta\left(  p\right)  }^{p}g\right)
,\mathbf{w}_{1}\left(  R_{\Delta\left(  p\right)  }^{p}g\right)
,...,\mathbf{w}_{p-1}\left(  R_{\Delta\left(  p\right)  }^{p}g\right)  \right)
\nonumber\\
&  =\left(  \mathbf{w}_{0}\left(  g\right)  ,\mathbf{w}_{1}\left(  g\right)
,...,\mathbf{w}_{p-1}\left(  g\right)  \right)  \ \ \ \ \ \ \ \ \ \ \left(
\text{by (\ref{cor.wi.Rn+1.eq})}\right) \label{pf.Delta.halfway.short.1}\\
&  =\left(  \mathbf{w}_{0}\left(  \operatorname*{vrefl}\nolimits^{\ast
}g\right)  ,\mathbf{w}_{1}\left(  \operatorname*{vrefl}\nolimits^{\ast
}g\right)  ,...,\mathbf{w}_{p-1}\left(  \operatorname*{vrefl}\nolimits^{\ast
}g\right)  \right) \nonumber
\end{align}
(since it is very easy to see that every $i\in\left\{  0,1,...,p-1\right\}  $
satisfies $\mathbf{w}_{i}\left(  g\right)  =\mathbf{w}_{i}\left(
\operatorname*{vrefl}\nolimits^{\ast}g\right)  $).

Claim 1 does not follow from our assumptions per se. But since we are allowed
to restrict $g$ to a Zariski-dense open subset of $\mathbb{K}^{\widehat{\Delta
\left(  p\right)  }}$ (after all, we are proving an identity between rational
functions), we can just WLOG assume that $\mathbf{w}_{i}\left(  g\right)
\neq0$ for all $i\in\left\{  0,1,...,p-1\right\}  $. Assume this. Then, every
$i\in\left\{  0,1,...,p-1\right\}  $ satisfies
\begin{align*}
\mathbf{w}_{i}\left(  R_{\Delta\left(  p\right)  }^{p}g\right)   &
=\mathbf{w}_{i}\left(  g\right)  \ \ \ \ \ \ \ \ \ \ \left(  \text{by
(\ref{pf.Delta.halfway.short.1})}\right) \\
&  \neq0.
\end{align*}
In other words, Claim 1 is satisfied.

All that remains now is proving Claim 3. In other words, we need to prove that
$\pi\left(  R_{\Delta\left(  p\right)  }^{p}g\right)  =\pi\left(
\operatorname*{vrefl}\nolimits^{\ast}g\right)  $.

Clearly, there exists a zero-free labelling $f\in\mathbb{K}^{\widehat{\Delta
\left(  p\right)  }}$ such that $f$ is homogeneously equivalent to $g$ and
such that $f\left(  0\right)  =2$\ \ \ \ \footnote{Indeed, in order to
construct such an $f$, it is enough to replace the label at $0$ of $g$ by $2$,
while keeping all the other labels unchanged.}. Consider such an $f$. Since
$f$ and $g$ are homogeneously equivalent, so are $\operatorname*{vrefl}%
\nolimits^{\ast}f$ and $\operatorname*{vrefl}\nolimits^{\ast}g$ (by Lemma
\ref{lem.Delta.hrefl} \textbf{(g)}), so that $\pi\left(  \operatorname*{vrefl}%
\nolimits^{\ast}f\right)  =\pi\left(  \operatorname*{vrefl}\nolimits^{\ast
}g\right)  $. Of course, also $\pi\left(  f\right)  =\pi\left(  g\right)  $
since $f$ and $g$ are homogeneously equivalent.

But the commutativity of the diagram (\ref{def.hgR.commut}) yields
$\overline{R}_{\Delta\left(  p\right)  }\circ\pi=\pi\circ R_{\Delta\left(
p\right)  }$. Using this equality, it is easy to show that $\overline
{R}_{\Delta\left(  p\right)  }^{\ell}\circ\pi=\pi\circ R_{\Delta\left(
p\right)  }^{\ell}$ for every $\ell\in\mathbb{N}$. In particular,
$\overline{R}_{\Delta\left(  p\right)  }^{p}\circ\pi=\pi\circ R_{\Delta\left(
p\right)  }^{p}$. Since%
\[
\pi\left(  R_{\Delta\left(  p\right)  }^{p}g\right)  =\underbrace{\left(
\pi\circ R_{\Delta\left(  p\right)  }^{p}\right)  }_{=\overline{R}%
_{\Delta\left(  p\right)  }^{p}\circ\pi}\left(  g\right)  =\left(
\overline{R}_{\Delta\left(  p\right)  }^{p}\circ\pi\right)  \left(  g\right)
=\overline{R}_{\Delta\left(  p\right)  }^{p}\left(  \pi\left(  g\right)
\right)  .
\]
Similarly, $\pi\left(  R_{\Delta\left(  p\right)  }^{p}f\right)  =\overline
{R}_{\Delta\left(  p\right)  }^{p}\left(  \pi\left(  f\right)  \right)  $.
Thus,%
\[
\pi\left(  R_{\Delta\left(  p\right)  }^{p}f\right)  =\overline{R}%
_{\Delta\left(  p\right)  }^{p}\left(  \underbrace{\pi\left(  f\right)
}_{=\pi\left(  g\right)  }\right)  =\overline{R}_{\Delta\left(  p\right)
}^{p}\left(  \pi\left(  g\right)  \right)  =\pi\left(  R_{\Delta\left(
p\right)  }^{p}g\right)  .
\]
Hence, in order to prove that $\pi\left(  R_{\Delta\left(  p\right)  }%
^{p}g\right)  =\pi\left(  \operatorname*{vrefl}\nolimits^{\ast}g\right)  $, it
is enough to show that $\pi\left(  R_{\Delta\left(  p\right)  }^{p}f\right)
=\pi\left(  \operatorname*{vrefl}\nolimits^{\ast}g\right)  $. Since
$\pi\left(  \operatorname*{vrefl}\nolimits^{\ast}f\right)  =\pi\left(
\operatorname*{vrefl}\nolimits^{\ast}g\right)  $, this in turn boils down to
proving that $\pi\left(  R_{\Delta\left(  p\right)  }^{p}f\right)  =\pi\left(
\operatorname*{vrefl}\nolimits^{\ast}f\right)  $. In order to prove this, we
will verify the stronger assertion stating that $R_{\Delta\left(  p\right)
}^{p}f=\operatorname*{vrefl}\nolimits^{\ast}f$. Notice that $g$ does not
appear in this statement anymore, and we can forget about $g$; all we need to
know is that $f$ is a sufficiently generic $\mathbb{K}$-labelling of
$\Delta\left(  p\right)  $ satisfying $f\left(  0\right)  =2$.

Let $\left(  i,k\right)  \in\Delta\left(  p\right)  $ be arbitrary. Then,
$i+k>p+1$ (since $\left(  i,k\right)  \in\Delta\left(  p\right)  $).
Consequently, $2p-\left(  i+k-1\right)  $ is a well-defined element of
$\left\{  1,2,...,p-1\right\}  $. Denote this element by $h$. Thus,
$h\in\left\{  1,2,...,p-1\right\}  $ and $i+k-1+h=2p$. Moreover, $\left(
k,i\right)  =\operatorname*{vrefl}v\in\Delta\left(  p\right)  $.

Let $v=\left(  p+1-k,p+1-i\right)  $. Then, $v=\operatorname*{hrefl}\left(
\left(  i,k\right)  \right)  \in\nabla\left(  p\right)  $ (since $\left(
i,k\right)  \in\Delta\left(  p\right)  $) and $\deg v=h$ (this follows by
simple computation). Moreover, $\operatorname*{hrefl}v=\left(  i,k\right)  $.

Applying Corollary \ref{cor.Delta.hrefl.f} to $\ell=h$, we obtain
$R_{\operatorname*{Rect}\left(  p,p\right)  }^{h}\left(  \operatorname*{wing}%
f\right)  =\operatorname*{wing}\left(  R_{\Delta\left(  p\right)  }%
^{h}f\right)  $, hence%
\begin{align}
&  \left(  R_{\operatorname*{Rect}\left(  p,p\right)  }^{h}\left(
\operatorname*{wing}f\right)  \right)  \left(  v\right) \nonumber\\
&  =\left(  \operatorname*{wing}\left(  R_{\Delta\left(  p\right)  }%
^{h}f\right)  \right)  \left(  v\right)  =\dfrac{1}{\left(  R_{\Delta\left(
p\right)  }^{p-\deg v}\left(  R_{\Delta\left(  p\right)  }^{h}f\right)
\right)  \left(  \operatorname*{hrefl}v\right)  }\nonumber\\
&  \ \ \ \ \ \ \ \ \ \ \left(  \text{by the definition of }%
\operatorname*{wing}\text{, since }v\in\nabla\left(  p\right)  \subseteq
\nabla\left(  p\right)  \cup\left\{  0\right\}  \right) \nonumber\\
&  =\dfrac{1}{\left(  R_{\Delta\left(  p\right)  }^{p-h}\left(  R_{\Delta
\left(  p\right)  }^{h}f\right)  \right)  \left(  \left(  i,k\right)  \right)
}\ \ \ \ \ \ \ \ \ \ \left(  \text{since }\deg v=h\text{ and }%
\operatorname*{hrefl}v=\left(  i,k\right)  \right) \nonumber\\
&  =\dfrac{1}{\left(  R_{\Delta\left(  p\right)  }^{p}f\right)  \left(
\left(  i,k\right)  \right)  }\ \ \ \ \ \ \ \ \ \ \left(  \text{since
}R_{\Delta\left(  p\right)  }^{p-h}\left(  R_{\Delta\left(  p\right)  }%
^{h}f\right)  =\underbrace{\left(  R_{\Delta\left(  p\right)  }^{p-h}\circ
R_{\Delta\left(  p\right)  }^{h}\right)  }_{=R_{\Delta\left(  p\right)  }^{p}%
}f=R_{\Delta\left(  p\right)  }^{p}f\right)  .
\label{pf.Delta.halfway.short.4}%
\end{align}
But Theorem \ref{thm.rect.antip.general} (applied to $p$,
$R_{\operatorname*{Rect}\left(  p,p\right)  }^{h}\left(  \operatorname*{wing}%
f\right)  $ and $\left(  k,i\right)  $ instead of $q$, $f$ and $\left(
i,k\right)  $) yields%
\begin{align*}
&  \left(  R_{\operatorname*{Rect}\left(  p,p\right)  }^{h}\left(
\operatorname*{wing}f\right)  \right)  \left(  \left(  p+1-k,p+1-i\right)
\right) \\
&  =\dfrac{\left(  R_{\operatorname*{Rect}\left(  p,p\right)  }^{h}\left(
\operatorname*{wing}f\right)  \right)  \left(  0\right)  \cdot\left(
R_{\operatorname*{Rect}\left(  p,p\right)  }^{h}\left(  \operatorname*{wing}%
f\right)  \right)  \left(  1\right)  }{\left(  R_{\operatorname*{Rect}\left(
p,p\right)  }^{i+k-1}\left(  R_{\operatorname*{Rect}\left(  p,p\right)  }%
^{h}\left(  \operatorname*{wing}f\right)  \right)  \right)  \left(  \left(
k,i\right)  \right)  }.
\end{align*}
Since $\left(  p+1-k,p+1-i\right)  =v$ and
\begin{align*}
R_{\operatorname*{Rect}\left(  p,p\right)  }^{i+k-1}\left(
R_{\operatorname*{Rect}\left(  p,p\right)  }^{h}\left(  \operatorname*{wing}%
f\right)  \right)   &  =\left(  \underbrace{R_{\operatorname*{Rect}\left(
p,p\right)  }^{i+k-1}\circ R_{\operatorname*{Rect}\left(  p,p\right)  }^{h}%
}_{\substack{=R_{\operatorname*{Rect}\left(  p,p\right)  }^{i+k-1+h}%
=R_{\operatorname*{Rect}\left(  p,p\right)  }^{2p}\\\text{(since
}i+k-1+h=2p\text{)}}}\right)  \left(  \operatorname*{wing}f\right) \\
&  =\underbrace{R_{\operatorname*{Rect}\left(  p,p\right)  }^{2p}%
}_{\substack{=\operatorname*{id}\\\text{(since Theorem \ref{thm.rect.ord}
(applied to }q=p\text{)}\\\text{yields }\operatorname*{ord}\left(
R_{\operatorname*{Rect}\left(  p,p\right)  }\right)  =p+p=2p\text{)}}}\left(
\operatorname*{wing}f\right)  =\operatorname*{wing}f,
\end{align*}
this equality rewrites as%
\[
\left(  R_{\operatorname*{Rect}\left(  p,p\right)  }^{h}\left(
\operatorname*{wing}f\right)  \right)  \left(  v\right)  =\dfrac{\left(
R_{\operatorname*{Rect}\left(  p,p\right)  }^{h}\left(  \operatorname*{wing}%
f\right)  \right)  \left(  0\right)  \cdot\left(  R_{\operatorname*{Rect}%
\left(  p,p\right)  }^{h}\left(  \operatorname*{wing}f\right)  \right)
\left(  1\right)  }{\left(  \operatorname*{wing}f\right)  \left(  \left(
k,i\right)  \right)  }.
\]
Since
\begin{align*}
&  \underbrace{\left(  R_{\operatorname*{Rect}\left(  p,p\right)  }^{h}\left(
\operatorname*{wing}f\right)  \right)  \left(  0\right)  }_{\substack{=\left(
\operatorname*{wing}f\right)  \left(  0\right)  \\\text{(by Corollary
\ref{cor.R.implicit.01})}}}\cdot\underbrace{\left(  R_{\operatorname*{Rect}%
\left(  p,p\right)  }^{h}\left(  \operatorname*{wing}f\right)  \right)
\left(  1\right)  }_{\substack{=\left(  \operatorname*{wing}f\right)  \left(
1\right)  \\\text{(by Corollary \ref{cor.R.implicit.01})}}}\\
&  =\underbrace{\left(  \operatorname*{wing}f\right)  \left(  0\right)
}_{\substack{=\dfrac{1}{\left(  R_{\Delta\left(  p\right)  }^{p-\deg
0}f\right)  \left(  \operatorname*{hrefl}0\right)  }\\\text{(by the definition
of }\operatorname*{wing}\text{)}}}\cdot\underbrace{\left(
\operatorname*{wing}f\right)  \left(  1\right)  }_{\substack{=f\left(
1\right)  \\\text{(by the definition of }\operatorname*{wing}\text{)}}}\\
&  =\dfrac{1}{\left(  R_{\Delta\left(  p\right)  }^{p-\deg0}f\right)  \left(
\operatorname*{hrefl}0\right)  }\cdot f\left(  1\right)  =1
\end{align*}
(since Corollary \ref{cor.R.implicit.01} yields $\left(  R_{\Delta\left(
p\right)  }^{p-\deg0}f\right)  \left(  \operatorname*{hrefl}0\right)
=f\left(  \operatorname*{hrefl}0\right)  =f\left(  1\right)  $), this
simplifies to
\[
\left(  R_{\operatorname*{Rect}\left(  p,p\right)  }^{h}\left(
\operatorname*{wing}f\right)  \right)  \left(  v\right)  =\dfrac{1}{\left(
\operatorname*{wing}f\right)  \left(  \left(  k,i\right)  \right)  }.
\]
Compared with (\ref{pf.Delta.halfway.short.4}), this yields $\dfrac{1}{\left(
R_{\Delta\left(  p\right)  }^{p}f\right)  \left(  \left(  i,k\right)  \right)
}=\dfrac{1}{\left(  \operatorname*{wing}f\right)  \left(  \left(  k,i\right)
\right)  }$. Taking inverses in this equality, we get
\begin{align*}
\left(  R_{\Delta\left(  p\right)  }^{p}f\right)  \left(  \left(  i,k\right)
\right)   &  =\left(  \operatorname*{wing}f\right)  \left(  \left(
k,i\right)  \right)  =f\left(  \underbrace{\left(  k,i\right)  }%
_{=\operatorname*{vrefl}\left(  i,k\right)  }\right) \\
&  \ \ \ \ \ \ \ \ \ \ \left(  \text{by the definition of }%
\operatorname*{wing}\text{, since }\left(  k,i\right)  \in\Delta\left(
p\right)  \subseteq\Delta\left(  p\right)  \cup\left\{  1\right\}  \right) \\
&  =f\left(  \operatorname*{vrefl}\left(  i,k\right)  \right)  =\left(
\operatorname*{vrefl}\nolimits^{\ast}f\right)  \left(  \left(  i,k\right)
\right) \\
&  \ \ \ \ \ \ \ \ \ \ \left(  \text{since }\left(  \operatorname*{vrefl}%
\nolimits^{\ast}f\right)  \left(  \left(  i,k\right)  \right)  =f\left(
\operatorname*{vrefl}\left(  i,k\right)  \right)  \text{ by the definition of
}\operatorname*{vrefl}\nolimits^{\ast}\right)  .
\end{align*}

Now, recall that we have shown this for \textbf{every} $\left(  i,k\right)
\in\Delta\left(  p\right)  $. In other words, we have shown that
$R_{\Delta\left(  p\right)  }^{p}f=\operatorname*{vrefl}\nolimits^{\ast}f$. As
we have seen, this is all we need to prove Theorem \ref{thm.Delta.halfway}.
\end{proof}
\end{vershort}

\begin{verlong}
\begin{proof}
[Second proof of Theorem \ref{thm.Delta.halfway} (sketched).]The result that
we are striving to prove is a collection of identities between rational
functions, hence boils down to a collection of polynomial identities in the
labels of an arbitrary $\mathbb{K}$-labelling of $\Delta\left(  p\right)  $.
Therefore, it is enough to prove it in the case when $\mathbb{K}$ is a field
of rational functions in finitely many variables over $\mathbb{Q}$. So let us
WLOG assume that $\mathbb{K}$ is a field of rational functions in finitely
many variables over $\mathbb{Q}$. Then, the characteristic of $\mathbb{K}$ is
$\neq2$ (it is $0$ indeed), so that we can apply Lemma \ref{lem.Delta.hrefl}
and Corollary \ref{cor.Delta.hrefl.f}.

Consider the maps $\operatorname*{hrefl}$, $\operatorname*{wing}$,
$\operatorname*{vrefl}$ and $\operatorname*{vrefl}\nolimits^{\ast}$ defined in
Lemma \ref{lem.Delta.hrefl}. We are going to prove that%
\[
R_{\Delta\left(  p\right)  }^{p}=\operatorname*{vrefl}\nolimits^{\ast}%
\]
as rational maps $\mathbb{K}^{\widehat{\Delta\left(  p\right)  }%
}\dashrightarrow\mathbb{K}^{\widehat{\Delta\left(  p\right)  }}$.

In fact, let $g\in\mathbb{K}^{\widehat{\Delta\left(  p\right)  }}$ be any
sufficiently generic zero-free labelling of $\Delta\left(  p\right)  $. Here,
``sufficiently generic'' means ``generic enough for all terms appearing in the
present proof to be well-defined and for all genericity conditions to be satisfied''.

Clearly, there exists a zero-free labelling $f\in\mathbb{K}^{\widehat{\Delta
\left(  p\right)  }}$ such that $f$ is homogeneously equivalent to $g$ and
such that $f\left(  0\right)  =2$\ \ \ \ \footnote{Indeed, in order to
construct such an $f$, it is enough to replace the label at $0$ of $g$ by $2$,
while keeping all the other labels unchanged. The resulting labelling $f$ is
indeed homogeneously equivalent to $g$ (because $0$ is the only element of
$\widehat{\Delta\left(  p\right)  }$ which has degree $0$).}. Consider such an
$f$.

It is easy to see that $\left(  R_{\Delta\left(  p\right)  }^{p}f\right)
\left(  1\right)  =\left(  \operatorname*{vrefl}\nolimits^{\ast}f\right)
\left(  1\right)  $ and $\left(  R_{\Delta\left(  p\right)  }^{p}f\right)
\left(  0\right)  =\left(  \operatorname*{vrefl}\nolimits^{\ast}f\right)
\left(  0\right)  $.\ \ \ \ \footnote{\textit{Proof.} Corollary
\ref{cor.R.implicit.01} yields $\left(  R_{\Delta\left(  p\right)  }%
^{p}f\right)  \left(  1\right)  =f\left(  1\right)  $ and $\left(
R_{\Delta\left(  p\right)  }^{p}f\right)  \left(  0\right)  =f\left(
0\right)  $. On the other hand, since $\operatorname*{vrefl}$ leaves the
elements $0$ and $1$ fixed, we have $\left(  \operatorname*{vrefl}%
\nolimits^{\ast}f\right)  \left(  1\right)  =f\left(  1\right)  $ and $\left(
\operatorname*{vrefl}\nolimits^{\ast}f\right)  \left(  0\right)  =f\left(
0\right)  $. Thus, $\left(  R_{\Delta\left(  p\right)  }^{p}f\right)  \left(
1\right)  =f\left(  1\right)  =\left(  \operatorname*{vrefl}\nolimits^{\ast
}f\right)  \left(  1\right)  $ and $\left(  R_{\Delta\left(  p\right)  }%
^{p}f\right)  \left(  0\right)  =f\left(  0\right)  =\left(
\operatorname*{vrefl}\nolimits^{\ast}f\right)  \left(  0\right)  $, qed.}
Also,%
\begin{equation}
\left(  \operatorname*{wing}f\right)  \left(  0\right)  \cdot\left(
\operatorname*{wing}f\right)  \left(  1\right)  =1.
\label{pf.Delta.halfway.0*1}%
\end{equation}
\footnote{\textit{Proof of (\ref{pf.Delta.halfway.0*1}):} By the definition of
$\operatorname*{wing}$, we have $\left(  \operatorname*{wing}f\right)  \left(
1\right)  =f\left(  1\right)  $ (since $1\in\Delta\left(  p\right)
\cup\left\{  1\right\}  $) and%
\begin{align*}
\left(  \operatorname*{wing}f\right)  \left(  0\right)   &  =\dfrac{1}{\left(
R_{\Delta\left(  p\right)  }^{p-\deg0}f\right)  \left(  \operatorname*{hrefl}%
0\right)  }\ \ \ \ \ \ \ \ \ \ \left(  \text{since }0\in\nabla\left(
p\right)  \cup\left\{  0\right\}  \right) \\
&  =\dfrac{1}{\left(  R_{\Delta\left(  p\right)  }^{p-0}f\right)  \left(
1\right)  }\ \ \ \ \ \ \ \ \ \ \left(  \text{since }\operatorname*{hrefl}%
0=1\text{ and }\deg0=0\right) \\
&  =\dfrac{1}{\left(  R_{\Delta\left(  p\right)  }^{p}f\right)  \left(
1\right)  }=\dfrac{1}{f\left(  1\right)  }\ \ \ \ \ \ \ \ \ \ \left(
\text{since Corollary \ref{cor.R.implicit.01} yields }\left(  R_{\Delta\left(
p\right)  }^{p}f\right)  \left(  1\right)  =f\left(  1\right)  \right)  .
\end{align*}
Hence,%
\[
\underbrace{\left(  \operatorname*{wing}f\right)  \left(  0\right)  }%
_{=\dfrac{1}{f\left(  1\right)  }}\cdot\underbrace{\left(
\operatorname*{wing}f\right)  \left(  1\right)  }_{=f\left(  1\right)
}=\dfrac{1}{f\left(  1\right)  }\cdot f\left(  1\right)  =1,
\]
and thus (\ref{pf.Delta.halfway.0*1}) is proven.}

Let $\left(  i,k\right)  \in\Delta\left(  p\right)  $ be arbitrary. Then,
$i+k>p+1$ (since $\left(  i,k\right)  \in\Delta\left(  p\right)  $).
Consequently, $2p-\left(  i+k-1\right)  $ is a well-defined element of
$\left\{  1,2,...,p-1\right\}  $. Denote this element by $h$. Thus,
$h=2p-\left(  i+k-1\right)  \in\left\{  1,2,...,p-1\right\}  $. Moreover, from
$h=2p-\left(  i+k-1\right)  $, we obtain $i+k-1+h=2p$. Also, since $\left(
i,k\right)  =v$, we have $\left(  k,i\right)  =\operatorname*{vrefl}v\in
\Delta\left(  p\right)  $ (since $v\in\Delta\left(  p\right)  $).

Let $v=\left(  p+1-k,p+1-i\right)  $. Then, $v=\operatorname*{hrefl}\left(
\left(  i,k\right)  \right)  \in\nabla\left(  p\right)  $ (since $\left(
i,k\right)  \in\Delta\left(  p\right)  $) and
\[
\deg v=\left(  p+1-k\right)  +\left(  p+1-i\right)  -1=2p-\left(
i+k-1\right)  =h.
\]
Furthermore, since $v=\operatorname*{hrefl}\left(  \left(  i,k\right)
\right)  $, we have $\operatorname*{hrefl}v=\left(  i,k\right)  $ (because
$\operatorname*{hrefl}$ is an involution).

Applying Corollary \ref{cor.Delta.hrefl.f} to $\ell=h$, we obtain
$R_{\operatorname*{Rect}\left(  p,p\right)  }^{h}\left(  \operatorname*{wing}%
f\right)  =\operatorname*{wing}\left(  R_{\Delta\left(  p\right)  }%
^{h}f\right)  $, hence%
\begin{align}
&  \left(  R_{\operatorname*{Rect}\left(  p,p\right)  }^{h}\left(
\operatorname*{wing}f\right)  \right)  \left(  v\right) \nonumber\\
&  =\left(  \operatorname*{wing}\left(  R_{\Delta\left(  p\right)  }%
^{h}f\right)  \right)  \left(  v\right)  =\dfrac{1}{\left(  R_{\Delta\left(
p\right)  }^{p-\deg v}\left(  R_{\Delta\left(  p\right)  }^{h}f\right)
\right)  \left(  \operatorname*{hrefl}v\right)  }\nonumber\\
&  \ \ \ \ \ \ \ \ \ \ \left(  \text{by the definition of }%
\operatorname*{wing}\text{, since }v\in\nabla\left(  p\right)  \subseteq
\nabla\left(  p\right)  \cup\left\{  0\right\}  \right) \nonumber\\
&  =\dfrac{1}{\left(  R_{\Delta\left(  p\right)  }^{p-h}\left(  R_{\Delta
\left(  p\right)  }^{h}f\right)  \right)  \left(  \left(  i,k\right)  \right)
}\ \ \ \ \ \ \ \ \ \ \left(  \text{since }\deg v=h\text{ and }%
\operatorname*{hrefl}v=\left(  i,k\right)  \right) \nonumber\\
&  =\dfrac{1}{\left(  R_{\Delta\left(  p\right)  }^{p}f\right)  \left(
\left(  i,k\right)  \right)  }\ \ \ \ \ \ \ \ \ \ \left(  \text{since
}R_{\Delta\left(  p\right)  }^{p-h}\left(  R_{\Delta\left(  p\right)  }%
^{h}f\right)  =\underbrace{\left(  R_{\Delta\left(  p\right)  }^{p-h}\circ
R_{\Delta\left(  p\right)  }^{h}\right)  }_{=R_{\Delta\left(  p\right)  }^{p}%
}f=R_{\Delta\left(  p\right)  }^{p}f\right)  . \label{pf.Delta.halfway.4}%
\end{align}
But Theorem \ref{thm.rect.antip.general} (applied to $p$,
$R_{\operatorname*{Rect}\left(  p,p\right)  }^{h}\left(  \operatorname*{wing}%
f\right)  $ and $\left(  k,i\right)  $ instead of $q$, $f$ and $\left(
i,k\right)  $) yields%
\begin{align*}
&  \left(  R_{\operatorname*{Rect}\left(  p,p\right)  }^{h}\left(
\operatorname*{wing}f\right)  \right)  \left(  \left(  p+1-k,p+1-i\right)
\right) \\
&  =\dfrac{\left(  R_{\operatorname*{Rect}\left(  p,p\right)  }^{h}\left(
\operatorname*{wing}f\right)  \right)  \left(  0\right)  \cdot\left(
R_{\operatorname*{Rect}\left(  p,p\right)  }^{h}\left(  \operatorname*{wing}%
f\right)  \right)  \left(  1\right)  }{\left(  R_{\operatorname*{Rect}\left(
p,p\right)  }^{i+k-1}\left(  R_{\operatorname*{Rect}\left(  p,p\right)  }%
^{h}\left(  \operatorname*{wing}f\right)  \right)  \right)  \left(  \left(
k,i\right)  \right)  }.
\end{align*}
Since $\left(  p+1-k,p+1-i\right)  =v$ and
\begin{align*}
R_{\operatorname*{Rect}\left(  p,p\right)  }^{i+k-1}\left(
R_{\operatorname*{Rect}\left(  p,p\right)  }^{h}\left(  \operatorname*{wing}%
f\right)  \right)   &  =\left(  \underbrace{R_{\operatorname*{Rect}\left(
p,p\right)  }^{i+k-1}\circ R_{\operatorname*{Rect}\left(  p,p\right)  }^{h}%
}_{\substack{=R_{\operatorname*{Rect}\left(  p,p\right)  }^{i+k-1+h}%
=R_{\operatorname*{Rect}\left(  p,p\right)  }^{2p}\\\text{(since
}i+k-1+h=2p\text{)}}}\right)  \left(  \operatorname*{wing}f\right) \\
&  =\underbrace{R_{\operatorname*{Rect}\left(  p,p\right)  }^{2p}%
}_{\substack{=\operatorname*{id}\\\text{(since Theorem \ref{thm.rect.ord}
(applied to }q=p\text{)}\\\text{yields }\operatorname*{ord}\left(
R_{\operatorname*{Rect}\left(  p,p\right)  }\right)  =p+p=2p\text{)}}}\left(
\operatorname*{wing}f\right) \\
&  =\operatorname*{id}\left(  \operatorname*{wing}f\right)
=\operatorname*{wing}f,
\end{align*}
this equality rewrites as%
\[
\left(  R_{\operatorname*{Rect}\left(  p,p\right)  }^{h}\left(
\operatorname*{wing}f\right)  \right)  \left(  v\right)  =\dfrac{\left(
R_{\operatorname*{Rect}\left(  p,p\right)  }^{h}\left(  \operatorname*{wing}%
f\right)  \right)  \left(  0\right)  \cdot\left(  R_{\operatorname*{Rect}%
\left(  p,p\right)  }^{h}\left(  \operatorname*{wing}f\right)  \right)
\left(  1\right)  }{\left(  \operatorname*{wing}f\right)  \left(  \left(
k,i\right)  \right)  }.
\]
Since
\[
\left(  R_{\operatorname*{Rect}\left(  p,p\right)  }^{h}\left(
\operatorname*{wing}f\right)  \right)  \left(  0\right)  =\left(
\operatorname*{wing}f\right)  \left(  0\right)  \ \ \ \ \ \ \ \ \ \ \left(
\text{by Corollary \ref{cor.R.implicit.01}}\right)
\]
and%
\[
\left(  R_{\operatorname*{Rect}\left(  p,p\right)  }^{h}\left(
\operatorname*{wing}f\right)  \right)  \left(  1\right)  =\left(
\operatorname*{wing}f\right)  \left(  1\right)  \ \ \ \ \ \ \ \ \ \ \left(
\text{by Corollary \ref{cor.R.implicit.01}}\right)  ,
\]
this equality rewrites as
\[
\left(  R_{\operatorname*{Rect}\left(  p,p\right)  }^{h}\left(
\operatorname*{wing}f\right)  \right)  \left(  v\right)  =\dfrac{\left(
\operatorname*{wing}f\right)  \left(  0\right)  \cdot\left(
\operatorname*{wing}f\right)  \left(  1\right)  }{\left(  \operatorname*{wing}%
f\right)  \left(  \left(  k,i\right)  \right)  }.
\]
Since $\left(  \operatorname*{wing}f\right)  \left(  0\right)  \cdot\left(
\operatorname*{wing}f\right)  \left(  1\right)  =1$ (by
(\ref{pf.Delta.halfway.0*1})), this rewrites as $\left(
R_{\operatorname*{Rect}\left(  p,p\right)  }^{h}\left(  \operatorname*{wing}%
f\right)  \right)  \left(  v\right)  =\dfrac{1}{\left(  \operatorname*{wing}%
f\right)  \left(  \left(  k,i\right)  \right)  }$. Compared with
(\ref{pf.Delta.halfway.4}), this yields $\dfrac{1}{\left(  R_{\Delta\left(
p\right)  }^{p}f\right)  \left(  \left(  i,k\right)  \right)  }=\dfrac
{1}{\left(  \operatorname*{wing}f\right)  \left(  \left(  k,i\right)  \right)
}$. Taking inverses in this equality, we get
\begin{align*}
\left(  R_{\Delta\left(  p\right)  }^{p}f\right)  \left(  \left(  i,k\right)
\right)   &  =\left(  \operatorname*{wing}f\right)  \left(  \left(
k,i\right)  \right)  =f\left(  \underbrace{\left(  k,i\right)  }%
_{=\operatorname*{vrefl}\left(  i,k\right)  }\right) \\
&  \ \ \ \ \ \ \ \ \ \ \left(  \text{by the definition of }%
\operatorname*{wing}\text{, since }\left(  k,i\right)  \in\Delta\left(
p\right)  \subseteq\Delta\left(  p\right)  \cup\left\{  1\right\}  \right) \\
&  =f\left(  \operatorname*{vrefl}\left(  i,k\right)  \right)  =\left(
\operatorname*{vrefl}\nolimits^{\ast}f\right)  \left(  \left(  i,k\right)
\right) \\
&  \ \ \ \ \ \ \ \ \ \ \left(  \text{since }\left(  \operatorname*{vrefl}%
\nolimits^{\ast}f\right)  \left(  \left(  i,k\right)  \right)  =f\left(
\operatorname*{vrefl}\left(  i,k\right)  \right)  \text{ by the definition of
}\operatorname*{vrefl}\nolimits^{\ast}\right)  .
\end{align*}

Now, forget that we fixed $\left(  i,k\right)  $. We thus have proven that
every $\left(  i,k\right)  \in\Delta\left(  p\right)  $ satisfies $\left(
R_{\Delta\left(  p\right)  }^{p}f\right)  \left(  \left(  i,k\right)  \right)
=\left(  \operatorname*{vrefl}\nolimits^{\ast}f\right)  \left(  \left(
i,k\right)  \right)  $. Hence, the $\mathbb{K}$-labellings $R_{\Delta\left(
p\right)  }^{p}f$ and $\operatorname*{vrefl}\nolimits^{\ast}f$ are equal on
each element of $\Delta\left(  p\right)  $. Since these two $\mathbb{K}%
$-labellings are also equal on each of the elements $0$ and $1$ (because
$\left(  R_{\Delta\left(  p\right)  }^{p}f\right)  \left(  1\right)  =\left(
\operatorname*{vrefl}\nolimits^{\ast}f\right)  \left(  1\right)  $ and
$\left(  R_{\Delta\left(  p\right)  }^{p}f\right)  \left(  0\right)  =\left(
\operatorname*{vrefl}\nolimits^{\ast}f\right)  \left(  0\right)  $), this
yields that the $\mathbb{K}$-labellings $R_{\Delta\left(  p\right)  }^{p}f$
and $\operatorname*{vrefl}\nolimits^{\ast}f$ are equal on each element of
$\Delta\left(  p\right)  \cup\left\{  0,1\right\}  =\widehat{\Delta\left(
p\right)  }$. In other words, $R_{\Delta\left(  p\right)  }^{p}%
f=\operatorname*{vrefl}\nolimits^{\ast}f$.

Now, recall that $f$ is homogeneously equivalent to $g$. Thus, $\pi\left(
f\right)  =\pi\left(  g\right)  $. Also, $\operatorname*{vrefl}\nolimits^{\ast
}f$ is homogeneously equivalent to $\operatorname*{vrefl}\nolimits^{\ast}g$
(by Lemma \ref{lem.Delta.hrefl} \textbf{(g)}), so that $\pi\left(
\operatorname*{vrefl}\nolimits^{\ast}f\right)  =\pi\left(
\operatorname*{vrefl}\nolimits^{\ast}g\right)  $.

But the commutativity of the diagram (\ref{def.hgR.commut}) yields
$\overline{R}_{\Delta\left(  p\right)  }\circ\pi=\pi\circ R_{\Delta\left(
p\right)  }$. Using this equality, it is easy to show that%
\[
\overline{R}_{\Delta\left(  p\right)  }^{\ell}\circ\pi=\pi\circ R_{\Delta
\left(  p\right)  }^{\ell}\ \ \ \ \ \ \ \ \ \ \text{for every }\ell
\in\mathbb{N}%
\]
(in fact, this can be proven by induction over $\ell$). Applying this to
$\ell=p$, we obtain $\overline{R}_{\Delta\left(  p\right)  }^{p}\circ\pi
=\pi\circ R_{\Delta\left(  p\right)  }^{p}$. Hence,%
\[
\left(  \overline{R}_{\Delta\left(  p\right)  }^{p}\circ\pi\right)  \left(
f\right)  =\left(  \pi\circ R_{\Delta\left(  p\right)  }^{p}\right)  \left(
f\right)  =\pi\left(  \underbrace{R_{\Delta\left(  p\right)  }^{p}%
f}_{=\operatorname*{vrefl}\nolimits^{\ast}f}\right)  =\pi\left(
\operatorname*{vrefl}\nolimits^{\ast}f\right)  =\pi\left(
\operatorname*{vrefl}\nolimits^{\ast}g\right)  .
\]
Comparing this with
\begin{align*}
\left(  \overline{R}_{\Delta\left(  p\right)  }^{p}\circ\pi\right)  \left(
f\right)   &  =\overline{R}_{\Delta\left(  p\right)  }^{p}\left(
\underbrace{\pi\left(  f\right)  }_{=\pi\left(  g\right)  }\right)
=\overline{R}_{\Delta\left(  p\right)  }^{p}\left(  \pi\left(  g\right)
\right)  =\left(  \underbrace{\overline{R}_{\Delta\left(  p\right)  }^{p}%
\circ\pi}_{=\pi\circ R_{\Delta\left(  p\right)  }^{p}}\right)  \left(
g\right) \\
&  =\left(  \pi\circ R_{\Delta\left(  p\right)  }^{p}\right)  \left(
g\right)  =\pi\left(  R_{\Delta\left(  p\right)  }^{p}g\right)  ,
\end{align*}
we conclude that%
\begin{equation}
\pi\left(  R_{\Delta\left(  p\right)  }^{p}g\right)  =\pi\left(
\operatorname*{vrefl}\nolimits^{\ast}g\right)  .
\label{pf.Delta.halfway.sign2}%
\end{equation}

But $\Delta\left(  p\right)  $ is a $\left(  p-1\right)  $-graded poset.
Hence, we can apply (\ref{cor.wi.Rn+1.eq}) to $\Delta\left(  p\right)  $,
$p-1$ and $g$ instead of $P$, $n$ and $f$. We thus obtain%
\[
\left(  \mathbf{w}_{0}\left(  R_{\Delta\left(  p\right)  }^{\left(
p-1\right)  +1}g\right)  ,\mathbf{w}_{1}\left(  R_{\Delta\left(  p\right)
}^{\left(  p-1\right)  +1}g\right)  ,...,\mathbf{w}_{p-1}\left(
R_{\Delta\left(  p\right)  }^{\left(  p-1\right)  +1}g\right)  \right)
=\left(  \mathbf{w}_{0}\left(  g\right)  ,\mathbf{w}_{1}\left(  g\right)
,...,\mathbf{w}_{p-1}\left(  g\right)  \right)  .
\]
In other words, every $i\in\left\{  0,1,...,p-1\right\}  $ satisfies
$\mathbf{w}_{i}\left(  R_{\Delta\left(  p\right)  }^{\left(  p-1\right)
+1}g\right)  =\mathbf{w}_{i}\left(  g\right)  $. Since $\left(  p-1\right)
+1=p$, this rewrites as follows: Every $i\in\left\{  0,1,...,p-1\right\}  $
satisfies%
\begin{equation}
\mathbf{w}_{i}\left(  R_{\Delta\left(  p\right)  }^{p}g\right)  =\mathbf{w}%
_{i}\left(  g\right)  . \label{pf.Delta.halfway.5}%
\end{equation}

Since $g$ is sufficiently generic, we can WLOG assume that no $i\in\left\{
0,1,...,p-1\right\}  $ satisfies $\mathbf{w}_{i}\left(  g\right)  =0$. Assume this.

But it is easy to see that every $i\in\left\{  0,1,...,p-1\right\}  $
satisfies $\mathbf{w}_{i}\left(  \operatorname*{vrefl}\nolimits^{\ast
}g\right)  =\mathbf{w}_{i}\left(  g\right)  $\ \ \ \ \footnote{\textit{Proof.}
Denote by $P$ the poset $\Delta\left(  p\right)  $. Let $i\in\left\{
0,1,...,p-1\right\}  $. Then, the definition of $\mathbf{w}_{i}\left(
\operatorname*{vrefl}\nolimits^{\ast}g\right)  $ yields
\[
\mathbf{w}_{i}\left(  \operatorname*{vrefl}\nolimits^{\ast}g\right)
=\sum_{\substack{x\in\widehat{P}_{i};\ y\in\widehat{P}_{i+1};\\y\gtrdot
x}}\underbrace{\dfrac{\left(  \operatorname*{vrefl}\nolimits^{\ast}g\right)
\left(  x\right)  }{\left(  \operatorname*{vrefl}\nolimits^{\ast}g\right)
\left(  y\right)  }}_{\substack{=\dfrac{g\left(  \operatorname*{vrefl}%
x\right)  }{g\left(  \operatorname*{vrefl}y\right)  }\\\text{(because the
definition of }\operatorname*{vrefl}\nolimits^{\ast}\\\text{yields }\left(
\operatorname*{vrefl}\nolimits^{\ast}g\right)  \left(  x\right)  =g\left(
\operatorname*{vrefl}x\right)  \\\text{and }\left(  \operatorname*{vrefl}%
\nolimits^{\ast}g\right)  \left(  y\right)  =g\left(  \operatorname*{vrefl}%
y\right)  \text{)}}}=\sum_{\substack{x\in\widehat{P}_{i};\ y\in\widehat{P}%
_{i+1};\\y\gtrdot x}}\dfrac{g\left(  \operatorname*{vrefl}x\right)  }{g\left(
\operatorname*{vrefl}y\right)  }=\sum_{\substack{x\in\widehat{P}_{i}%
;\ y\in\widehat{P}_{i+1};\\y\gtrdot x}}\dfrac{g\left(  x\right)  }{g\left(
y\right)  }%
\]
(here, we have substituted $x$ and $y$ for $\operatorname*{vrefl}x$ and
$\operatorname*{vrefl}y$, because $\operatorname*{vrefl}$ is an involutive
poset automorphism of $P$). But the definition of $\mathbf{w}_{i}\left(
g\right)  $ yields $\mathbf{w}_{i}\left(  g\right)  =\sum
\limits_{\substack{x\in\widehat{P}_{i};\ y\in\widehat{P}_{i+1};\\y\gtrdot
x}}\dfrac{g\left(  x\right)  }{g\left(  y\right)  }$. Hence, $\mathbf{w}%
_{i}\left(  \operatorname*{vrefl}\nolimits^{\ast}g\right)  =\sum
\limits_{\substack{x\in\widehat{P}_{i};\ y\in\widehat{P}_{i+1};\\y\gtrdot
x}}\dfrac{g\left(  x\right)  }{g\left(  y\right)  }=\mathbf{w}_{i}\left(
g\right)  $, qed.}. Hence, every $i\in\left\{  0,1,...,p-1\right\}  $
satisfies $\mathbf{w}_{i}\left(  R_{\Delta\left(  p\right)  }^{p}g\right)
=\mathbf{w}_{i}\left(  g\right)  =\mathbf{w}_{i}\left(  \operatorname*{vrefl}%
\nolimits^{\ast}g\right)  $. Thus,%
\begin{align}
&  \left(  \mathbf{w}_{0}\left(  R_{\Delta\left(  p\right)  }^{p}g\right)
,\mathbf{w}_{1}\left(  R_{\Delta\left(  p\right)  }^{p}g\right)
,...,\mathbf{w}_{p-1}\left(  R_{\Delta\left(  p\right)  }^{p}g\right)  \right)
\nonumber\\
&  =\left(  \mathbf{w}_{0}\left(  \operatorname*{vrefl}\nolimits^{\ast
}g\right)  ,\mathbf{w}_{1}\left(  \operatorname*{vrefl}\nolimits^{\ast
}g\right)  ,...,\mathbf{w}_{p-1}\left(  \operatorname*{vrefl}\nolimits^{\ast
}g\right)  \right)  . \label{pf.Delta.halfway.sign1}%
\end{align}

Moreover, $\left(  R_{\Delta\left(  p\right)  }^{p}g\right)  \left(  0\right)
=g\left(  0\right)  $ (by Corollary \ref{cor.R.implicit.01}), but also
$\left(  \operatorname*{vrefl}\nolimits^{\ast}g\right)  \left(  0\right)
=g\left(  0\right)  $ (since $\operatorname*{vrefl}$ fixes $0$). Hence,%
\begin{equation}
\left(  R_{\Delta\left(  p\right)  }^{p}g\right)  \left(  0\right)  =g\left(
0\right)  =\left(  \operatorname*{vrefl}\nolimits^{\ast}g\right)  \left(
0\right)  . \label{pf.Delta.halfway.sign3}%
\end{equation}

Now, recall that no $i\in\left\{  0,1,...,p-1\right\}  $ satisfies
$\mathbf{w}_{i}\left(  g\right)  =0$. Because of (\ref{pf.Delta.halfway.5}),
this rewrites as follows: No $i\in\left\{  0,1,...,p-1\right\}  $ satisfies
$\mathbf{w}_{i}\left(  R_{\Delta\left(  p\right)  }^{p}g\right)  =0$. Thus,
because of (\ref{pf.Delta.halfway.sign1}), (\ref{pf.Delta.halfway.sign2}) and
(\ref{pf.Delta.halfway.sign3}), we can apply Proposition
\ref{prop.reconstruct} to $p-1$, $\Delta\left(  p\right)  $, $R_{\Delta\left(
p\right)  }^{p}g$ and $\operatorname*{vrefl}\nolimits^{\ast}g$ instead of $n$,
$P$, $f$ and $g$. Thus, we obtain $R_{\Delta\left(  p\right)  }^{p}%
g=\operatorname*{vrefl}\nolimits^{\ast}g$.

Now, forget that we fixed $g$. We thus have shown that $R_{\Delta\left(
p\right)  }^{p}g=\operatorname*{vrefl}\nolimits^{\ast}g$ for every
sufficiently generic $g$. Hence, $R_{\Delta\left(  p\right)  }^{p}%
=\operatorname*{vrefl}\nolimits^{\ast}$. Thus, for every $\left(  i,k\right)
\in\Delta\left(  p\right)  $ and every $f\in\mathbb{K}^{\widehat{\Delta\left(
p\right)  }}$, we have%
\begin{align*}
\underbrace{\left(  R_{\Delta\left(  p\right)  }^{p}f\right)  }%
_{=\operatorname*{vrefl}\nolimits^{\ast}}\left(  \left(  i,k\right)  \right)
&  =\left(  \operatorname*{vrefl}\nolimits^{\ast}f\right)  \left(  \left(
i,k\right)  \right)  =f\left(  \underbrace{\operatorname*{vrefl}\left(
i,k\right)  }_{=\left(  k,i\right)  }\right)  \ \ \ \ \ \ \ \ \ \ \left(
\text{by the definition of }\operatorname*{vrefl}\nolimits^{\ast}\right) \\
&  =f\left(  \left(  k,i\right)  \right)  .
\end{align*}
This proves Theorem \ref{thm.Delta.halfway}.
\end{proof}
\end{verlong}
\end{noncompile}

We can now obtain Theorem \ref{thm.Nabla.halfway} from Theorem
\ref{thm.Delta.halfway} using a construction from the proof of Proposition
\ref{prop.op.ord}:

\begin{vershort}
\begin{proof}
[Proof of Theorem \ref{thm.Nabla.halfway} (sketched).]The poset
antiautomorphism $\operatorname*{hrefl}$ of $\operatorname*{Rect}\left(
p,p\right)  $ defined in Remark \ref{rmk.DeltaNabla} restricts to a poset
antiisomorphism $\operatorname*{hrefl}:\nabla\left(  p\right)  \rightarrow
\Delta\left(  p\right)  $, that is, to a poset homomorphism
$\operatorname*{hrefl}:\nabla\left(  p\right)  \rightarrow\left(
\Delta\left(  p\right)  \right)  ^{\operatorname*{op}}$. We will use this
isomorphism to identify the poset $\nabla\left(  p\right)  $ with the opposite
poset $\left(  \Delta\left(  p\right)  \right)  ^{\operatorname*{op}}$ of
$\Delta\left(  p\right)  $.

Set $P=\Delta\left(  p\right)  $. Define a rational map $\kappa:\mathbb{K}%
^{\widehat{P}}\dashrightarrow\mathbb{K}^{\widehat{P^{\operatorname*{op}}}}$ as
in the proof of Proposition \ref{prop.op.ord}. Then, as in said proof, it can
be shown that the map $\kappa$ is a birational map and satisfies $\kappa\circ
R_{P}=R_{P^{\operatorname*{op}}}^{-1}\circ\kappa$. Since $P=\Delta\left(
p\right)  $ and $P^{\operatorname*{op}}=\left(  \Delta\left(  p\right)
\right)  ^{\operatorname*{op}}=\nabla\left(  p\right)  $, this rewrites as
$\kappa\circ R_{\Delta\left(  p\right)  }=R_{\nabla\left(  p\right)  }%
^{-1}\circ\kappa$. For the same reason, we know that $\kappa$ is a rational
map $\mathbb{K}^{\widehat{\Delta\left(  p\right)  }}\dashrightarrow
\mathbb{K}^{\widehat{\nabla\left(  p\right)  }}$.

From $\kappa\circ R_{\Delta\left(  p\right)  }=R_{\nabla\left(  p\right)
}^{-1}\circ\kappa$, we can easily obtain $\kappa\circ R_{\Delta\left(
p\right)  }^{m}=R_{\nabla\left(  p\right)  }^{-m}\circ\kappa$ for every
$m\in\mathbb{N}$. In particular, $\kappa\circ R_{\Delta\left(  p\right)  }%
^{p}=R_{\nabla\left(  p\right)  }^{-p}\circ\kappa$.

Now, consider the map $\operatorname*{vrefl}\nolimits^{\ast}:\mathbb{K}%
^{\widehat{\Delta\left(  p\right)  }}\rightarrow\mathbb{K}^{\widehat{\Delta
\left(  p\right)  }}$ defined in Lemma \ref{lem.Delta.hrefl} \textbf{(e)}, and
also consider the similarly defined map $\operatorname*{vrefl}\nolimits^{\ast
}:\mathbb{K}^{\widehat{\nabla\left(  p\right)  }}\rightarrow\mathbb{K}%
^{\widehat{\nabla\left(  p\right)  }}$. Both squares of the diagram%
\[
\xymatrixcolsep{5pc} \xymatrix{
\mathbb K^{\widehat{\Delta\left(p\right)}} \ar@{-->}[r]^{R_{\Delta
\left(p\right)}^p} \ar@{-->}[d]^{\kappa} & \mathbb K^{\widehat{\Delta
\left(p\right)}}
\ar[r]^{\operatorname*{vrefl}^{\ast}} \ar@{-->}[d]^{\kappa} & \mathbb
K^{\widehat{\Delta\left(p\right)}} \ar@{-->}[d]^{\kappa} \\
\mathbb K^{\widehat{\nabla\left(p\right)}} \ar@{-->}[r]_{R_{\nabla
\left(p\right)}^{-p}} & \mathbb K^{\widehat{\nabla\left(p\right)}}
\ar[r]_{\operatorname*{vrefl}^{\ast}} & \mathbb K^{\widehat{\nabla
\left(p\right)}}
}
\]
commute (the left square does so because of $\kappa\circ R_{\Delta\left(
p\right)  }^{p}=R_{\nabla\left(  p\right)  }^{-p}\circ\kappa$, and the
commutativity of the right square follows from a simple calculation), and so
the whole diagram commutes. In other words,
\begin{equation}
\kappa\circ\left(  \operatorname*{vrefl}\nolimits^{\ast}\circ R_{\Delta\left(
p\right)  }^{p}\right)  =\left(  \operatorname*{vrefl}\nolimits^{\ast}\circ
R_{\nabla\left(  p\right)  }^{-p}\right)  \circ\kappa.
\label{pf.Nabla.halfway.short.1}%
\end{equation}

But the statement of Theorem \ref{thm.Delta.halfway} can be rewritten as
$R_{\Delta\left(  p\right)  }^{p}=\operatorname*{vrefl}\nolimits^{\ast}$.
Since $\operatorname*{vrefl}\nolimits^{\ast}$ is an involution (this is clear
by inspection), we have $\operatorname*{vrefl}\nolimits^{\ast}=\left(
\operatorname*{vrefl}\nolimits^{\ast}\right)  ^{-1}$, so that
$\underbrace{\operatorname*{vrefl}\nolimits^{\ast}}_{=\left(
\operatorname*{vrefl}\nolimits^{\ast}\right)  ^{-1}}\circ\underbrace{R_{\Delta
\left(  p\right)  }^{p}}_{=\operatorname*{vrefl}\nolimits^{\ast}}=\left(
\operatorname*{vrefl}\nolimits^{\ast}\right)  ^{-1}\circ\operatorname*{vrefl}%
\nolimits^{\ast}=\operatorname*{id}$. Thus, (\ref{pf.Nabla.halfway.short.1})
simplifies to $\kappa\circ\operatorname*{id}=\left(  \operatorname*{vrefl}%
\nolimits^{\ast}\circ R_{\nabla\left(  p\right)  }^{-p}\right)  \circ\kappa$.
In other words, $\kappa=\left(  \operatorname*{vrefl}\nolimits^{\ast}\circ
R_{\nabla\left(  p\right)  }^{-p}\right)  \circ\kappa$. Since $\kappa$ is a
birational map, we can cancel $\kappa$ from this identity, obtaining
$\operatorname*{id}=\operatorname*{vrefl}\nolimits^{\ast}\circ R_{\nabla
\left(  p\right)  }^{-p}$. In other words, $R_{\nabla\left(  p\right)  }%
^{p}=\operatorname*{vrefl}\nolimits^{\ast}$. But this is precisely the
statement of Theorem \ref{thm.Nabla.halfway}.
\end{proof}
\end{vershort}

\begin{verlong}
\begin{proof}
[Proof of Theorem \ref{thm.Nabla.halfway} (sketched).]For every isomorphism
$\phi:P\rightarrow Q$ between two finite posets $P$ and $Q$, define a map
$\phi^{\ast}:\mathbb{K}^{\widehat{Q}}\rightarrow\mathbb{K}^{\widehat{P}}$ as
in the proof of Lemma \ref{lem.Delta.hrefl} \textbf{(e)}. Then, just in the
proof of Lemma \ref{lem.Delta.hrefl} \textbf{(e)}, we can show that
\begin{equation}
\phi^{\ast}\circ R_{Q}=R_{P}\circ\phi^{\ast} \label{pf.Nabla.halfway.functor}%
\end{equation}
whenever $\phi:P\rightarrow Q$ is an isomorphism between two finite posets $P$
and $Q$.

Set $P=\nabla\left(  p\right)  $. Define a rational map $\kappa:\mathbb{K}%
^{\widehat{P}}\dashrightarrow\mathbb{K}^{\widehat{P^{\operatorname*{op}}}}$ as
in the proof of Proposition \ref{prop.op.ord}. Then, as in the proof of
Proposition \ref{prop.op.ord}, we can see that this map $\kappa$ is a
birational map and satisfies $\kappa\circ R_{P}=R_{P^{\operatorname*{op}}%
}^{-1}\circ\kappa$. Thus,%
\begin{equation}
\kappa\circ R_{\nabla\left(  p\right)  }=R_{\nabla\left(  p\right)
^{\operatorname*{op}}}^{-1}\circ\kappa\label{pf.Nabla.halfway.commuteR}%
\end{equation}
(this is just a rewriting of $\kappa\circ R_{P}=R_{P^{\operatorname*{op}}%
}^{-1}\circ\kappa$, because $P=\nabla\left(  p\right)  $).

Moreover, from the definition of this map $\kappa$, we conclude that every
$f\in\mathbb{K}^{\widehat{P}}$ satisfies%
\begin{equation}
\left(  \kappa f\right)  \left(  w\right)  =\dfrac{1}{f\left(  w\right)
}\ \ \ \ \ \ \ \ \ \ \text{for every }w\in P. \label{pf.Nabla.halfway.kappa}%
\end{equation}

Notice that the antiautomorphism $\operatorname*{hrefl}$ restricts to an
antiisomorphism $\Delta\left(  p\right)  \rightarrow\nabla\left(  p\right)  $
of posets. We denote this restricted antiisomorphism by $\operatorname*{hrefl}%
$ again. Then, this $\operatorname*{hrefl}$ is an antiisomorphism
$\Delta\left(  p\right)  \rightarrow\nabla\left(  p\right)  $ of posets, hence
an isomorphism $\Delta\left(  p\right)  \rightarrow\nabla\left(  p\right)
^{\operatorname*{op}}$ of posets. Thus, it induces a map
$\operatorname*{hrefl}\nolimits^{\ast}:\mathbb{K}^{\widehat{\nabla\left(
p\right)  ^{\operatorname*{op}}}}\rightarrow\mathbb{K}^{\widehat{\Delta\left(
p\right)  }}$.

Applying (\ref{pf.Nabla.halfway.functor}) to $\Delta\left(  p\right)  $,
$\nabla\left(  p\right)  ^{\operatorname*{op}}$ and $\operatorname*{hrefl}$
instead of $P$, $Q$ and $\phi$, we obtain%
\begin{equation}
\operatorname*{hrefl}\nolimits^{\ast}\circ R_{\nabla\left(  p\right)
^{\operatorname*{op}}}=R_{\Delta\left(  p\right)  }\circ\operatorname*{hrefl}%
\nolimits^{\ast}. \label{pf.Nabla.halfway.commuteR.2}%
\end{equation}
Now,%
\begin{align}
&  \underbrace{\operatorname*{hrefl}\nolimits^{\ast}}_{=R_{\Delta\left(
p\right)  }^{-1}\circ R_{\Delta\left(  p\right)  }\circ\operatorname*{hrefl}%
\nolimits^{\ast}}\circ\underbrace{\kappa\circ R_{\nabla\left(  p\right)  }%
}_{\substack{=R_{\nabla\left(  p\right)  ^{\operatorname*{op}}}^{-1}%
\circ\kappa\\\text{(by (\ref{pf.Nabla.halfway.commuteR}))}}}\nonumber\\
&  =R_{\Delta\left(  p\right)  }^{-1}\circ\underbrace{R_{\Delta\left(
p\right)  }\circ\operatorname*{hrefl}\nolimits^{\ast}}%
_{\substack{=\operatorname*{hrefl}\nolimits^{\ast}\circ R_{\nabla\left(
p\right)  ^{\operatorname*{op}}}\\\text{(by (\ref{pf.Nabla.halfway.commuteR.2}%
))}}}\circ R_{\nabla\left(  p\right)  ^{\operatorname*{op}}}^{-1}\circ
\kappa\nonumber\\
&  =R_{\Delta\left(  p\right)  }^{-1}\circ\operatorname*{hrefl}\nolimits^{\ast
}\circ\underbrace{R_{\nabla\left(  p\right)  ^{\operatorname*{op}}}\circ
R_{\nabla\left(  p\right)  ^{\operatorname*{op}}}^{-1}}_{=\operatorname*{id}%
}\circ\kappa=R_{\Delta\left(  p\right)  }^{-1}\circ\operatorname*{hrefl}%
\nolimits^{\ast}\circ\kappa. \label{pf.Nabla.halfway.onestep}%
\end{align}

Using this, it is easy to see that%
\begin{equation}
\operatorname*{hrefl}\nolimits^{\ast}\circ\kappa\circ R_{\nabla\left(
p\right)  }^{\ell}=R_{\Delta\left(  p\right)  }^{-\ell}\circ
\operatorname*{hrefl}\nolimits^{\ast}\circ\kappa\ \ \ \ \ \ \ \ \ \ \text{for
every }\ell\in\mathbb{N}. \label{pf.Nabla.halfway.lstep}%
\end{equation}
\footnote{\textit{Proof of (\ref{pf.Nabla.halfway.lstep}):} We will prove
(\ref{pf.Nabla.halfway.lstep}) by induction over $\ell$:
\par
\textit{Induction base:} For $\ell=0$, the equality
(\ref{pf.Nabla.halfway.lstep}) is obvious (because for $\ell=0$, both
$R_{\nabla\left(  p\right)  }^{\ell}$ and $R_{\Delta\left(  p\right)  }%
^{-\ell}$ are identity maps). Hence, the induction base is complete.
\par
\textit{Induction step:} Let $L\in\mathbb{N}$. Assume that the equality
(\ref{pf.Nabla.halfway.lstep}) holds for $\ell=L$. We must then prove that the
equality (\ref{pf.Nabla.halfway.lstep}) holds for $\ell=L+1$.
\par
We have $\operatorname*{hrefl}\nolimits^{\ast}\circ\kappa\circ R_{\nabla
\left(  p\right)  }^{L}=R_{\Delta\left(  p\right)  }^{-L}\circ
\operatorname*{hrefl}\nolimits^{\ast}\circ\kappa$ (since the equality
(\ref{pf.Nabla.halfway.lstep}) holds for $\ell=L$). Now,%
\begin{align*}
\operatorname*{hrefl}\nolimits^{\ast}\circ\kappa\circ\underbrace{R_{\nabla
\left(  p\right)  }^{L+1}}_{=R_{\nabla\left(  p\right)  }^{L}\circ
R_{\nabla\left(  p\right)  }}  &  =\underbrace{\operatorname*{hrefl}%
\nolimits^{\ast}\circ\kappa\circ R_{\nabla\left(  p\right)  }^{L}}%
_{=R_{\Delta\left(  p\right)  }^{-L}\circ\operatorname*{hrefl}\nolimits^{\ast
}\circ\kappa}\circ R_{\nabla\left(  p\right)  }=R_{\Delta\left(  p\right)
}^{-L}\circ\underbrace{\operatorname*{hrefl}\nolimits^{\ast}\circ\kappa\circ
R_{\nabla\left(  p\right)  }}_{\substack{=R_{\Delta\left(  p\right)  }%
^{-1}\circ\operatorname*{hrefl}\nolimits^{\ast}\circ\kappa\\\text{(by
(\ref{pf.Nabla.halfway.onestep}))}}}\\
&  =\underbrace{R_{\Delta\left(  p\right)  }^{-L}\circ R_{\Delta\left(
p\right)  }^{-1}}_{=R_{\Delta\left(  p\right)  }^{-L-1}=R_{\Delta\left(
p\right)  }^{-\left(  L+1\right)  }}\circ\operatorname*{hrefl}\nolimits^{\ast
}\circ\kappa=R_{\Delta\left(  p\right)  }^{-\left(  L+1\right)  }%
\circ\operatorname*{hrefl}\nolimits^{\ast}\circ\kappa.
\end{align*}
In other words, the equality (\ref{pf.Nabla.halfway.lstep}) holds for
$\ell=L+1$. This completes the induction step. The induction proof of
(\ref{pf.Nabla.halfway.lstep}) is thus complete.}

But it is easy to see (using Theorem \ref{thm.Delta.halfway}) that%
\begin{equation}
\left(  R_{\Delta\left(  p\right)  }^{-p}f\right)  \left(  \left(  k,i\right)
\right)  =f\left(  \left(  i,k\right)  \right)  \ \ \ \ \ \ \ \ \ \ \text{for
every }\left(  i,k\right)  \in\Delta\left(  p\right)  \text{ and every }%
f\in\mathbb{K}^{\widehat{\Delta\left(  p\right)  }}
\label{pf.Nabla.halfway.halfway}%
\end{equation}
(where $f$ is assumed to be sufficiently generic for the left hand side of
this to be well-defined).\footnote{\textit{Proof of
(\ref{pf.Nabla.halfway.halfway}):} Let $f\in\mathbb{K}^{\widehat{\Delta\left(
p\right)  }}$ and $\left(  i,k\right)  \in\Delta\left(  p\right)  $. Applying
Theorem \ref{thm.Delta.halfway} to $R_{\Delta\left(  p\right)  }^{-p}f$
instead of $f$, we obtain%
\[
\left(  R_{\Delta\left(  p\right)  }^{p}\left(  R_{\Delta\left(  p\right)
}^{-p}f\right)  \right)  \left(  \left(  i,k\right)  \right)  =\left(
R_{\Delta\left(  p\right)  }^{-p}f\right)  \left(  \left(  k,i\right)
\right)  .
\]
Thus,
\begin{align*}
\left(  R_{\Delta\left(  p\right)  }^{-p}f\right)  \left(  \left(  k,i\right)
\right)   &  =\left(  R_{\Delta\left(  p\right)  }^{p}\left(  R_{\Delta\left(
p\right)  }^{-p}f\right)  \right)  \left(  \left(  i,k\right)  \right)
=\left(  \underbrace{\left(  R_{\Delta\left(  p\right)  }^{p}\circ
R_{\Delta\left(  p\right)  }^{-p}\right)  }_{=\operatorname*{id}}f\right)
\left(  \left(  i,k\right)  \right) \\
&  =\left(  \operatorname*{id}f\right)  \left(  \left(  i,k\right)  \right)
=f\left(  \left(  i,k\right)  \right)  .
\end{align*}
This proves (\ref{pf.Nabla.halfway.halfway}).}

Now, let $f\in\mathbb{K}^{\widehat{\nabla\left(  p\right)  }}$ be sufficiently
generic. Then, $f\in\mathbb{K}^{\widehat{\nabla\left(  p\right)  }}%
=\mathbb{K}^{\widehat{P}}$ (since $\nabla\left(  p\right)  =P$). Let $\left(
i,k\right)  \in\nabla\left(  p\right)  $. Then, $\left(  p+1-i,p+1-k\right)
\in\Delta\left(  p\right)  $ and $\left(  p+1-k,p+1-i\right)  \in\Delta\left(
p\right)  $. Thus, $\left(  \left(  \operatorname*{hrefl}\nolimits^{\ast}%
\circ\kappa\circ R_{\nabla\left(  p\right)  }^{\ell}\right)  f\right)  \left(
\left(  p+1-k,p+1-i\right)  \right)  $ is well-defined. We have
\begin{align*}
&  \left(  \left(  \operatorname*{hrefl}\nolimits^{\ast}\circ\kappa\circ
R_{\nabla\left(  p\right)  }^{p}\right)  f\right)  \left(  \left(
p+1-k,p+1-i\right)  \right) \\
&  =\left(  \operatorname*{hrefl}\nolimits^{\ast}\left(  \kappa\left(
R_{\nabla\left(  p\right)  }^{p}f\right)  \right)  \right)  \left(  \left(
p+1-k,p+1-i\right)  \right) \\
&  =\left(  \kappa\left(  R_{\nabla\left(  p\right)  }^{p}f\right)  \right)
\left(  \underbrace{\operatorname*{hrefl}\left(  p+1-k,p+1-i\right)
}_{\substack{=\left(  i,k\right)  \\\text{(this is easy to see)}}}\right)
\ \ \ \ \ \ \ \ \ \ \left(  \text{by the definition of }\operatorname*{hrefl}%
\nolimits^{\ast}\right) \\
&  =\left(  \kappa\left(  R_{\nabla\left(  p\right)  }^{p}f\right)  \right)
\left(  \left(  i,k\right)  \right)  =\dfrac{1}{\left(  R_{\nabla\left(
p\right)  }^{p}f\right)  \left(  \left(  i,k\right)  \right)  }\\
&  \ \ \ \ \ \ \ \ \ \ \left(
\begin{array}
[c]{c}%
\text{by (\ref{pf.Nabla.halfway.kappa}), applied to }R_{\nabla\left(
p\right)  }^{p}f\text{ and }\left(  i,k\right) \\
\text{instead of }f\text{ and }w\text{ (since }\left(  i,k\right)  \in
\nabla\left(  p\right)  =P\text{ and }R_{\nabla\left(  p\right)  }^{p}%
f\in\mathbb{K}^{\widehat{\nabla\left(  p\right)  }}=\mathbb{K}^{\widehat{P}%
}\text{)}%
\end{array}
\right)  .
\end{align*}
Compared to%
\begin{align*}
&  \left(  \underbrace{\left(  \operatorname*{hrefl}\nolimits^{\ast}%
\circ\kappa\circ R_{\nabla\left(  p\right)  }^{p}\right)  }%
_{\substack{=R_{\Delta\left(  p\right)  }^{-p}\circ\operatorname*{hrefl}%
\nolimits^{\ast}\circ\kappa\\\text{(by (\ref{pf.Nabla.halfway.lstep}), applied
to }\ell=p\text{)}}}f\right)  \left(  \left(  p+1-k,p+1-i\right)  \right) \\
&  =\left(  \left(  R_{\Delta\left(  p\right)  }^{-p}\circ
\operatorname*{hrefl}\nolimits^{\ast}\circ\kappa\right)  f\right)  \left(
\left(  p+1-k,p+1-i\right)  \right) \\
&  =\left(  R_{\Delta\left(  p\right)  }^{-p}\left(  \operatorname*{hrefl}%
\nolimits^{\ast}\left(  \kappa f\right)  \right)  \right)  \left(  \left(
p+1-k,p+1-i\right)  \right) \\
&  =\left(  \operatorname*{hrefl}\nolimits^{\ast}\left(  \kappa f\right)
\right)  \left(  \left(  p+1-i,p+1-k\right)  \right) \\
&  \ \ \ \ \ \ \ \ \ \ \left(
\begin{array}
[c]{c}%
\text{by (\ref{pf.Nabla.halfway.halfway}), applied to }\operatorname*{hrefl}%
\nolimits^{\ast}\left(  \kappa f\right)  \text{ and }\left(
p+1-i,p+1-k\right) \\
\text{instead of }f\text{ and }\left(  i,k\right)  \text{ (since }\left(
p+1-i,p+1-k\right)  \in\Delta\left(  p\right)  \text{)}%
\end{array}
\right) \\
&  =\left(  \kappa f\right)  \left(  \underbrace{\operatorname*{hrefl}\left(
p+1-i,p+1-k\right)  }_{\substack{=\left(  k,i\right)  \\\text{(this is easy to
prove)}}}\right)  \ \ \ \ \ \ \ \ \ \ \left(  \text{by the definition of
}\operatorname*{hrefl}\nolimits^{\ast}\right) \\
&  =\left(  \kappa f\right)  \left(  \left(  k,i\right)  \right)  =\dfrac
{1}{f\left(  \left(  k,i\right)  \right)  }\ \ \ \ \ \ \ \ \ \ \left(
\begin{array}
[c]{c}%
\text{by (\ref{pf.Nabla.halfway.kappa}), applied to }\left(  i,k\right)
\text{ instead of }\\
w\text{ (since }\left(  k,i\right)  \in\nabla\left(  p\right)  =P\text{ and
}f\in\mathbb{K}^{\widehat{P}}\text{)}%
\end{array}
\right)  ,
\end{align*}
this yields $\dfrac{1}{\left(  R_{\nabla\left(  p\right)  }^{p}f\right)
\left(  \left(  i,k\right)  \right)  }=\dfrac{1}{f\left(  \left(  k,i\right)
\right)  }$. Taking inverses of both sides, we obtain \newline$\left(
R_{\nabla\left(  p\right)  }^{p}f\right)  \left(  \left(  i,k\right)  \right)
=f\left(  \left(  k,i\right)  \right)  $. This proves Theorem
\ref{thm.Nabla.halfway}.
\end{proof}
\end{verlong}

\begin{proof}
[Proof of Corollary \ref{cor.Delta.ord} (sketched).]\textbf{(a)} Let
$f\in\mathbb{K}^{\widehat{\Delta\left(  p\right)  }}$ be sufficiently generic.
Then, every $\left(  i,k\right)  \in\Delta\left(  p\right)  $ satisfies%
\begin{align*}
&  \left(  \underbrace{R_{\Delta\left(  p\right)  }^{2p}}_{=R_{\Delta\left(
p\right)  }^{p}\circ R_{\Delta\left(  p\right)  }^{p}}f\right)  \left(
\left(  i,k\right)  \right) \\
&  =\left(  \left(  R_{\Delta\left(  p\right)  }^{p}\circ R_{\Delta\left(
p\right)  }^{p}\right)  f\right)  \left(  \left(  i,k\right)  \right)
=\left(  R_{\Delta\left(  p\right)  }^{p}\left(  R_{\Delta\left(  p\right)
}^{p}f\right)  \right)  \left(  \left(  i,k\right)  \right) \\
&  =\left(  R_{\Delta\left(  p\right)  }^{p}f\right)  \left(  \left(
k,i\right)  \right)  \ \ \ \ \ \ \ \ \ \ \left(  \text{by Theorem
\ref{thm.Delta.halfway}, applied to }R_{\Delta\left(  p\right)  }^{p}f\text{
instead of }f\right) \\
&  =f\left(  \left(  i,k\right)  \right)  \ \ \ \ \ \ \ \ \ \ \left(  \text{by
Theorem \ref{thm.Delta.halfway}, applied to }\left(  k,i\right)  \text{
instead of }\left(  i,k\right)  \right)  .
\end{align*}
Hence, the two labellings $R_{\Delta\left(  p\right)  }^{2p}f$ and $f$ are
equal on every element of $\Delta\left(  p\right)  $. Since these two
labellings are also equal on $0$ and $1$ (because Corollary
\ref{cor.R.implicit.01} yields $\left(  R_{\Delta\left(  p\right)  }%
^{2p}f\right)  \left(  0\right)  =f\left(  0\right)  $ and $\left(
R_{\Delta\left(  p\right)  }^{2p}f\right)  \left(  1\right)  =f\left(
1\right)  $), this yields that the two labellings $R_{\Delta\left(  p\right)
}^{2p}f$ and $f$ are equal on every element of $\Delta\left(  p\right)
\cup\left\{  0,1\right\}  =\widehat{\Delta\left(  p\right)  }$. Hence,
$R_{\Delta\left(  p\right)  }^{2p}f=f=\operatorname*{id}f$.

Now, forget that we fixed $f$. We thus have shown that $R_{\Delta\left(
p\right)  }^{2p}f=\operatorname*{id}f$ for every sufficiently generic
$f\in\mathbb{K}^{\widehat{\Delta\left(  p\right)  }}$. Hence, $R_{\Delta
\left(  p\right)  }^{2p}=\operatorname*{id}$. In other words,
$\operatorname*{ord}\left(  R_{\Delta\left(  p\right)  }\right)  \mid2p$. This
proves Corollary \ref{cor.Delta.ord} \textbf{(a)}.

\textbf{(b)} Proving Corollary \ref{cor.Delta.ord} \textbf{(b)} is left to the reader.
\end{proof}

\begin{proof}
[Proof of Corollary \ref{cor.Nabla.ord} (sketched).]Corollary
\ref{cor.Nabla.ord} can be deduced from Theorem \ref{thm.Nabla.halfway} in the
same way as Corollary \ref{cor.Delta.ord} is deduced from Theorem
\ref{thm.Delta.halfway}. We won't dwell on the details.
\end{proof}

Let us conclude this section by stating a generalization of parts
\textbf{(b)}, \textbf{(c)}, \textbf{(d)} and \textbf{(f)} of Lemma
\ref{lem.Delta.hrefl} that was pointed out by a referee. Rather than
restricting itself to $\operatorname*{Rect}\left(  p,p\right)  $, it is
concerned with an arbitrary $\left(  2p-1\right)  $-graded poset satisfying
certain axioms (which can be informally subsumed under the slogan
\textquotedblleft symmetric with respect to degree $p$ and regular near the
middle\textquotedblright):\footnote{We choose to label the parts of Lemma
\ref{lem.Delta.hrefl-general} by \textbf{(b)}, \textbf{(c)}, \textbf{(d)} and
\textbf{(f)}, since they generalize the parts \textbf{(b)}, \textbf{(c)},
\textbf{(d)} and \textbf{(f)} of Lemma \ref{lem.Delta.hrefl}, respectively.}

\Needspace{27\baselineskip}

\begin{lemma}
\label{lem.Delta.hrefl-general}Let $p$ be a positive integer. Let $P$ be a
$\left(  2p-1\right)  $-graded finite poset. Let $\operatorname*{hrefl}%
:P\rightarrow P$ be an involution such that $\operatorname*{hrefl}$ is a poset
antiautomorphism of $P$. (This $\operatorname*{hrefl}$ has nothing to do with
the $\operatorname*{hrefl}$ defined in Lemma \ref{lem.Delta.hrefl}, although
of course it is analogous to the latter.) We extend $\operatorname*{hrefl}$ to
an involutive poset antiautomorphism of $\widehat{P}$ by setting
$\operatorname*{hrefl}\left(  0\right)  =1$ and $\operatorname*{hrefl}\left(
1\right)  =0$.

Assume that every $v\in\widehat{P}$ satisfies
\begin{equation}
\deg\left(  \operatorname{hrefl}v\right)  =2p-\deg v.
\label{eq.lem.Delta.hrefl-general.deg}%
\end{equation}

Let $N$ be a positive integer. Assume that, for every $v\in P$ satisfying
$\deg v=p-1$, there exist precisely $N$ elements $u$ of $P$ satisfying
$u\gtrdot v$.

Define three subsets $\Delta$, $\operatorname*{Eq}$ and $\nabla$ of $P$ by%
\begin{align*}
\Delta &  =\left\{  v\in P\ \mid\ \deg v>p\right\}  ;\\
\operatorname*{Eq}  &  =\left\{  v\in P\ \mid\ \deg v=p\right\}  ;\\
\nabla &  =\left\{  v\in P\ \mid\ \deg v<p\right\}  .
\end{align*}

Clearly, $\Delta$, $\operatorname*{Eq}$ and $\nabla$ become subposets of $P$.
The poset $\operatorname*{Eq}$ is an antichain, while the posets $\Delta$ and
$\nabla$ are $\left(  p-1\right)  $-graded.

Assume that $\operatorname*{hrefl}\mid_{\operatorname*{Eq}}=\operatorname*{id}%
$. It is easy to see that $\operatorname*{hrefl}\left(  \Delta\right)
=\nabla$.

Let $\mathbb{K}$ be a field such that $N$ is invertible in
$\mathbb{K}$.

\textbf{(b)} Define a rational map $\operatorname*{wing}:\mathbb{K}%
^{\widehat{\Delta}}\dashrightarrow\mathbb{K}^{\widehat{P}}$ by setting%
\[
\left(  \operatorname*{wing}f\right)  \left(  v\right)  =\left\{
\begin{array}
[c]{l}%
f\left(  v\right)  ,\ \ \ \ \ \ \ \ \ \ \text{if }v\in\Delta\cup\left\{
1\right\}  ;\\
1,\ \ \ \ \ \ \ \ \ \ \text{if }v\in\operatorname*{Eq};\\
\dfrac{1}{\left(  R_{\Delta}^{p-\deg v}f\right)  \left(  \operatorname*{hrefl}%
v\right)  },\ \ \ \ \ \ \ \ \ \ \text{if }v\in\nabla\cup\left\{  0\right\}
\end{array}
\right.
\]
for all $v\in\widehat{P}$ for all $f\in\mathbb{K}^{\widehat{\Delta}}$. This is well-defined.

\textbf{(c)} There exists a rational map $\overline{\operatorname*{wing}%
}:\overline{\mathbb{K}^{\widehat{\Delta}}}\dashrightarrow\overline
{\mathbb{K}^{\widehat{P}}}$ such that the diagram%
\[
\xymatrixcolsep{5pc}\xymatrix{ \mathbb{K}^{\widehat{\Delta}} \ar@{-->}[d]_-{\pi} \ar@{-->}[r]^{\operatorname*{wing}} & \mathbb{K}^{\widehat{P}} \ar@{-->}[d]^-{\pi} \\ \overline{\mathbb{K}^{\widehat{\Delta}}} \ar@{-->}[r]_{\overline{\operatorname*{wing}}} & \overline{\mathbb{K}^{\widehat{P}}} }
\]
commutes.

\textbf{(d)} The rational map $\overline{\operatorname*{wing}}$ defined in
Lemma \ref{lem.Delta.hrefl-general} \textbf{(b)} satisfies
\[
\overline{R}_{P}\circ\overline{\operatorname*{wing}}=\overline
{\operatorname*{wing}}\circ\overline{R}_{\Delta}.
\]

\textbf{(f)} Almost every (in the sense of Zariski topology) labelling
$f\in\mathbb{K}^{\widehat{\Delta}}$ satisfying $f\left(  0\right)  =N$
satisfies%
\[
R_{P}\left(  \operatorname*{wing}f\right)  =\operatorname*{wing}\left(
R_{\Delta}f\right)  .
\]

\end{lemma}

Notice that the hypothesis (\ref{eq.lem.Delta.hrefl-general.deg}) is actually
redundant (it follows from the other requirements), but we have chosen to
state it because it is easily checked in practice and used in the proof.

\begin{example}
Let $P$ be a positive integer, and let $\mathbb{K}$ be a field of
characteristic $\neq2$. The hypotheses of Lemma \ref{lem.Delta.hrefl-general}
are satisfied if we set $P=\operatorname*{Rect}\left(  p,p\right)  $,
$\operatorname*{hrefl}=\operatorname*{hrefl}$ (by this, we mean that we define
$\operatorname*{hrefl}$ to be the map $\operatorname*{hrefl}$ defined in Lemma
\ref{lem.Delta.hrefl}) and $N=2$. In this case, the posets $\Delta$,
$\operatorname*{Eq}$ and $\nabla$ defined in Lemma
\ref{lem.Delta.hrefl-general} are precisely the posets $\Delta\left(
p\right)  $, $\operatorname*{Eq}\left(  p\right)  $ and $\nabla\left(
p\right)  $ introduced in Definition \ref{def.DeltaNabla}. Hence, Lemma
\ref{lem.Delta.hrefl-general} (when applied to this setting) yields the parts
\textbf{(b)}, \textbf{(c)}, \textbf{(d)} and \textbf{(f)} of Lemma
\ref{lem.Delta.hrefl}.
\end{example}

\Needspace{23\baselineskip}

\begin{example}
Here is another example of a situation in which Lemma
\ref{lem.Delta.hrefl-general} applies. Namely, the hypotheses of Lemma
\ref{lem.Delta.hrefl-general} are satisfied when $p=5$, $N=3$ and $P$ is the
poset with Hasse diagram%
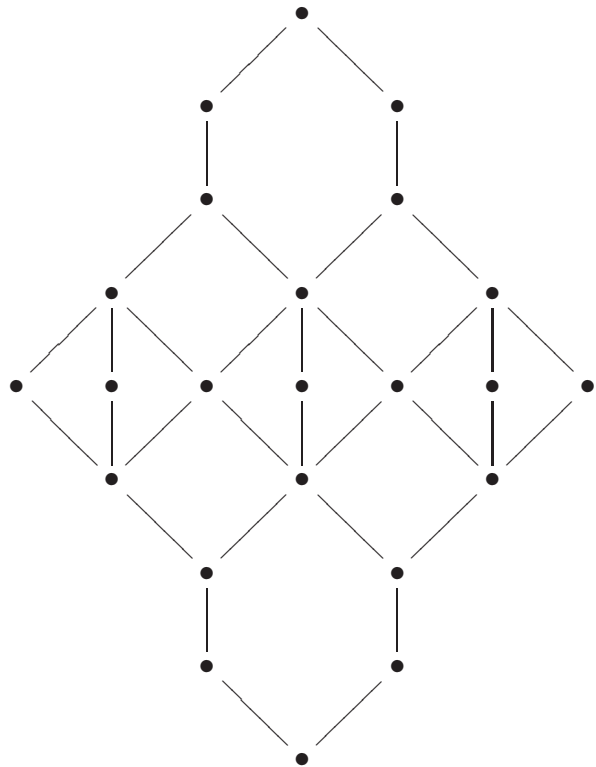
\[%
\xymatrix{
& & & \bullet\ar@{-}[dl] \ar@{-}[dr] \\
& & \bullet\ar@{-}[d] & & \bullet\ar@{-}[d] \\
& & \bullet\ar@{-}[dl] \ar@{-}[dr] & & \bullet\ar@{-}[dl] \ar@{-}[dr] \\
& \bullet\ar@{-}[dl] \ar@{-}[d] \ar@{-}[dr] & & \bullet\ar@{-}[dl] \ar@
{-}[d] \ar@{-}[dr] & & \bullet\ar@{-}[dl] \ar@{-}[d] \ar@{-}[dr] \\
\bullet\ar@{-}[dr] & \bullet\ar@{-}[d] & \bullet\ar@{-}[dl] \ar@
{-}[dr] & \bullet\ar@{-}[d] & \bullet\ar@{-}[dl] \ar@{-}[dr] & \bullet\ar@
{-}[d] & \bullet\ar@{-}[dl] \\
& \bullet\ar@{-}[dr] & & \bullet\ar@{-}[dl] \ar@{-}[dr] & & \bullet\ar@
{-}[dl] \\
& & \bullet\ar@{-}[d] & & \bullet\ar@{-}[d] \\
& & \bullet\ar@{-}[dr] & & \bullet\ar@{-}[dl] \\
& & & \bullet}%
\]
(with $\operatorname*{hrefl}:P\rightarrow P$ being the reflection with respect
to the horizontal axis of symmetry of this diagram).
\end{example}

\begin{proof}
[Proof of Lemma \ref{lem.Delta.hrefl-general} (sketched).] The proof of Lemma
\ref{lem.Delta.hrefl-general} is almost completely analogous to the proof of
parts \textbf{(b)}, \textbf{(c)}, \textbf{(d)} and \textbf{(f)} of Lemma
\ref{lem.Delta.hrefl}. Of course, several changes need to be made to the
latter proof to make it apply to Lemma \ref{lem.Delta.hrefl-general}: for instance,

\begin{itemize}
\item every appearance of $\operatorname*{Rect}\left(  p,p\right)  $,
$\Delta\left(  p\right)  $, $\nabla\left(  p\right)  $ or $\operatorname*{Eq}%
\left(  p\right)  $ must be replaced by $P$, $\Delta$, $\nabla$ or
$\operatorname*{Eq}$, respectively;

\item many (but not all) appearances of the number $2$ (such as its appearance
in the definition of $a_{i}$) have to be replaced by $N$;

\item various properties of $P$ now no longer follow from the definition of
$\operatorname*{Rect}\left(  p,p\right)  $ (because $P$ is no longer
$\operatorname*{Rect}\left(  p,p\right)  $), but instead have to be derived
from the hypotheses of Lemma \ref{lem.Delta.hrefl-general}\footnote{Most of
the time, this is obvious. For instance, the fact that $\operatorname*{hrefl}%
\left(  \Delta\right)  =\nabla$ follows from
(\ref{eq.lem.Delta.hrefl-general.deg}). The only fact that is not completely
trivial is that, for every $v\in P$ satisfying $\deg v=p+1$, there exist
precisely $N$ elements $u$ of $P$ satisfying $u\lessdot v$. Let us prove this
fact.
\par
We know that $\operatorname*{hrefl}$ is a poset antiautomorphism of
$\widehat{P}$. Hence, if $u$ and $v$ are two elements of $\widehat{P}$, then
we have the following equivalence of statements:%
\begin{equation}
\left(  u\lessdot v\right)  \Longleftrightarrow\left(  \operatorname*{hrefl}%
u\gtrdot\operatorname*{hrefl}v\right)  .\label{pf.Delta.hrefl-general.d.equiv}%
\end{equation}
\par
We also have assumed that, for every $v\in P$ satisfying $\deg v=p-1$, there
exist precisely $N$ elements $u$ of $P$ satisfying $u\gtrdot v$. In other
words, for every $v\in P$ satisfying $\deg v=p-1$, we have%
\begin{equation}
\left(  \text{the number of elements }u\text{ of }P\text{ satisfying }u\gtrdot
v\right)  =N.\label{pf.Delta.hrefl-general.almost-axiom.1}%
\end{equation}
\par
Now, let $v\in P$ be such that $\deg v=p+1$. We need to show that there exist
precisely $N$ elements $u$ of $P$ satisfying $u\lessdot v$.
\par
From (\ref{eq.lem.Delta.hrefl-general.deg}), we obtain $\deg\left(
\operatorname{hrefl}v\right)  =2p-\underbrace{\deg v}_{=p+1}=2p-\left(
p+1\right)  =p-1$. Hence, (\ref{pf.Delta.hrefl-general.almost-axiom.1})
(applied to $\operatorname*{hrefl}v$ instead of $v$) yields%
\begin{equation}
\left(  \text{the number of elements }u\text{ of }P\text{ satisfying }%
u\gtrdot\operatorname*{hrefl}v\right)
=N.\label{pf.Delta.hrefl-general.almost-axiom.2}%
\end{equation}
\par
But $\operatorname*{hrefl}:P\rightarrow P$ is a bijection (since
$\operatorname*{hrefl}$ is an involution). Thus, we can substitute
$\operatorname*{hrefl}u$ for $u$ in \textquotedblleft$\left(  \text{the number
of elements }u\text{ of }P\text{ satisfying }u\gtrdot\operatorname*{hrefl}%
v\right)  $\textquotedblright. We thus obtain%
\begin{align*}
& \left(  \text{the number of elements }u\text{ of }P\text{ satisfying
}u\gtrdot\operatorname*{hrefl}v\right)  \\
& =\left(  \text{the number of elements }u\text{ of }P\text{ satisfying
}\underbrace{\operatorname*{hrefl}u\gtrdot\operatorname*{hrefl}v}%
_{\substack{\text{this is equivalent to }\left(  u\lessdot v\right)
\\\text{(due to (\ref{pf.Delta.hrefl-general.d.equiv}))}}}\right)  \\
& =\left(  \text{the number of elements }u\text{ of }P\text{ satisfying
}u\lessdot v\right)  .
\end{align*}
Thus,%
\begin{align*}
& \left(  \text{the number of elements }u\text{ of }P\text{ satisfying
}u\lessdot v\right)  \\
& =\left(  \text{the number of elements }u\text{ of }P\text{ satisfying
}u\gtrdot\operatorname*{hrefl}v\right)  =N
\end{align*}
(by (\ref{pf.Delta.hrefl-general.almost-axiom.2})). In other words, there
exist precisely $N$ elements $u$ of $P$ satisfying $u\lessdot v$.
\par
This completes our proof of the fact that, for every $v\in P$ satisfying $\deg
v=p+1$, there exist precisely $N$ elements $u$ of $P$ satisfying $u\lessdot
v$.};

\item checking the case when $p\leq2$ is no longer trivial, but needs a bit
more work\footnote{The case when $p=1$ is still obvious (since $\Delta$ and
$\nabla$ are empty sets in this case). The case when $p=2$ can be handled by
the same arguments that were used to deal with the case when $p>2$ (in
particular, the same definition of the $\left(  2p+1\right)  $-tuple $\left(
a_{0},a_{1},...,a_{2p}\right)  $ applies), but the details are slightly
different (instead of the seven cases, there are now only three cases: $\deg
v=3$, $\deg v=2$ and $\deg v=1$).}.
\end{itemize}
\end{proof}

\section{The quarter-triangles}

\label{sect.quartertri}

We have now studied the order of birational rowmotion on all four triangles
(two of which are isomorphic as posets) which are obtained by cutting the
rectangle $\operatorname*{Rect}\left(  p,p\right)  $ along one of its
diagonals. But we can also cut $\operatorname*{Rect}\left(  p,p\right)  $
along \textbf{both} diagonals into four smaller triangles. These are
isomorphic in pairs, and we will analyze them now. The following definition is
an analogue of Definition \ref{def.DeltaNabla} but using $\operatorname*{Tria}%
\left(  p\right)  $ instead of $\operatorname*{Rect}\left(  p,p\right)  $:

\begin{definition}
\label{def.NEtri}Let $p$ be a positive integer. Define three subsets
$\operatorname*{NEtri}\left(  p\right)  $, $\operatorname*{Eqtri}\left(
p\right)  $ and $\operatorname*{SEtri}\left(  p\right)  $ of
$\operatorname*{Tria}\left(  p\right)  $ by%
\begin{align*}
\operatorname*{NEtri}\left(  p\right)   &  =\left\{  \left(  i,k\right)
\in\operatorname*{Tria}\left(  p\right)  \ \mid\ i+k>p+1\right\}  ;\\
\operatorname*{Eqtri}\left(  p\right)   &  =\left\{  \left(  i,k\right)
\in\operatorname*{Tria}\left(  p\right)  \ \mid\ i+k=p+1\right\}  ;\\
\operatorname*{SEtri}\left(  p\right)   &  =\left\{  \left(  i,k\right)
\in\operatorname*{Tria}\left(  p\right)  \ \mid\ i+k<p+1\right\}  .
\end{align*}
These subsets $\operatorname*{NEtri}\left(  p\right)  $,
$\operatorname*{Eqtri}\left(  p\right)  $ and $\operatorname*{SEtri}\left(
p\right)  $ inherit a poset structure from $\operatorname*{Tria}\left(
p\right)  $. In the following, we will consider $\operatorname*{NEtri}\left(
p\right)  $, $\operatorname*{Eqtri}\left(  p\right)  $ and
$\operatorname*{SEtri}\left(  p\right)  $ as posets using this structure.

Clearly, $\operatorname*{Eqtri}\left(  p\right)  $ is an antichain. The posets
$\operatorname*{NEtri}\left(  p\right)  $ and $\operatorname*{SEtri}\left(
p\right)  $ are $\left(  p-1\right)  $-graded posets having the form of
right-angled triangles.
\end{definition}

\Needspace{38\baselineskip}

\begin{example}
Here is the Hasse diagram of the poset $\operatorname*{Tria}\left(  4\right)
$, where the elements belonging to $\operatorname*{NEtri}\left(  4\right)  $
have been underlined and the elements belonging to $\operatorname*{Eqtri}%
\left(  4\right)  $ have been boxed:%
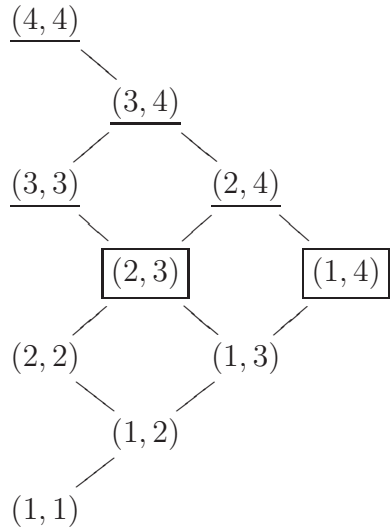
\[
\xymatrixrowsep{0.9pc}\xymatrixcolsep{0.20pc}\xymatrix{
& & & \underline{\left(4,4\right)} \ar@{-}[rd] & & & \\
& & & & \underline{\left(3,4\right)} \ar@{-}[rd] \ar@{-}[ld] & & \\
& & & \underline{\left(3,3\right)} \ar@{-}[rd] & & \underline{\left(2,4\right)} \ar@{-}[rd] \ar@{-}[ld] & \\
& & & & \fbox{$\left(2,3\right)$} \ar@{-}[rd] \ar@{-}[ld] & & \fbox{$\left(1,4\right)$} \ar@{-}[ld] \\
& & & \left(2,2\right) \ar@{-}[rd] & & \left(1,3\right) \ar@{-}[ld] & \\
& & & & \left(1,2\right) \ar@{-}[ld] & & \\
& & & \left(1,1\right) & & &
}.
\]
And here is the Hasse diagram of the poset $\operatorname*{NEtri}\left(
4\right)  $ itself:%
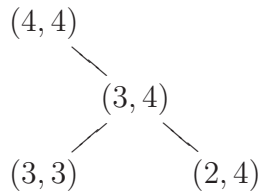
\[
\xymatrixrowsep{0.9pc}\xymatrixcolsep{0.20pc}\xymatrix{
& & & \left(4,4\right) \ar@{-}[rd] & & & \\
& & & & \left(3,4\right) \ar@{-}[rd] \ar@{-}[ld] & & \\
& & & \left(3,3\right) & & \left(2,4\right) &
}.
\]
Here, on the other hand, is the Hasse diagram of the poset
$\operatorname*{SEtri}\left(  4\right)  $:
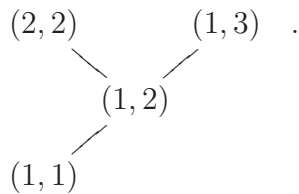
\[
\xymatrixrowsep{0.9pc}\xymatrixcolsep{0.20pc}\xymatrix{
& & & \left(2,2\right) \ar@{-}[rd] & & \left(1,3\right) \ar@{-}[ld] & \\
& & & & \left(1,2\right) \ar@{-}[ld] & & \\
& & & \left(1,1\right) & & &
}.
\]

\end{example}

\begin{remark}
\label{rmk.NEtri.SW}Let $p$ be an \textbf{even} positive integer. The poset
$\operatorname*{NEtri}\left(  p\right)  $ is isomorphic to the poset $\Phi
^{+}\left(  B_{p\diagup2}\right)  $ of \cite[\S 3.2]{striker-williams}. (For
odd $p$, the poset $\operatorname*{NEtri}\left(  p\right)  $ does not seem to
appear in \cite[\S 3.2]{striker-williams}.)
\end{remark}

The following conjectures have been verified using Sage for small values of
$p$:

\begin{conjecture}
\label{conj.SEtri.ord}Let $p$ be an integer $>1$. Then, $\operatorname*{ord}%
\left(  R_{\operatorname*{SEtri}\left(  p\right)  }\right)  =p$.
\end{conjecture}

\begin{conjecture}
\label{conj.NEtri.ord}Let $p$ be an integer $>1$. Then, $\operatorname*{ord}%
\left(  R_{\operatorname*{NEtri}\left(  p\right)  }\right)  =p$.
\end{conjecture}

The approach used to prove Theorem \ref{thm.Leftri.ord} allows proving these
two conjectures in the case of odd $p$, but in the even-$p$ case it fails
(although the order of \textbf{classical} rowmotion is again known to be $p$
in the even-$p$ case -- see \cite[Conjecture 3.6]{striker-williams}). Here is
how the proof proceeds in the case of odd $p$:

\begin{proposition}
\label{prop.SEtri.odd.ord}Let $p$ be an odd integer $>1$. Let $\mathbb{K}$ be
a field. Then, $\operatorname*{ord}\left(  R_{\operatorname*{SEtri}\left(
p\right)  }\right)  =p$.
\end{proposition}

\begin{proposition}
\label{prop.NEtri.odd.ord}Let $p$ be an odd integer $>1$. Let $\mathbb{K}$ be
a field. Then, $\operatorname*{ord}\left(  R_{\operatorname*{NEtri}\left(
p\right)  }\right)  =p$.
\end{proposition}

Our proof of Proposition \ref{prop.NEtri.odd.ord} rests upon the following fact:

\begin{lemma}
\label{lem.NEtri.vrefl}Let $\mathbb{K}$ be a field of characteristic $\neq2$.

Let $p$ be a positive integer.

\textbf{(a)} Let $\operatorname*{vrefl}:\Delta\left(  p\right)  \rightarrow
\Delta\left(  p\right)  $ be the map sending every $\left(  i,k\right)
\in\Delta\left(  p\right)  $ to $\left(  k,i\right)  $. This map
$\operatorname*{vrefl}$ is an involutive poset automorphism of $\Delta\left(
p\right)  $. (In intuitive terms, $\operatorname*{vrefl}$ is simply reflection
across the vertical axis.) We have $\operatorname*{vrefl}\left(  v\right)
\in\operatorname*{NEtri}\left(  p\right)  $ for every $v\in\Delta\left(
p\right)  \setminus\operatorname*{NEtri}\left(  p\right)  $.

We extend $\operatorname*{vrefl}$ to an involutive poset automorphism of
$\widehat{\Delta\left(  p\right)  }$ by setting $\operatorname*{vrefl}\left(
0\right)  =0$ and $\operatorname*{vrefl}\left(  1\right)  =1$.

\textbf{(b)} Define a map $\operatorname*{dble}:\mathbb{K}%
^{\widehat{\operatorname*{NEtri}\left(  p\right)  }}\rightarrow\mathbb{K}%
^{\widehat{\Delta\left(  p\right)  }}$ by setting%
\[
\left(  \operatorname*{dble}f\right)  \left(  v\right)  =\left\{
\begin{array}
[c]{l}%
\dfrac{1}{2}f\left(  1\right)  ,\ \ \ \ \ \ \ \ \ \ \text{if }v=1;\\
f\left(  0\right)  ,\ \ \ \ \ \ \ \ \ \ \text{if }v=0;\\
f\left(  v\right)  ,\ \ \ \ \ \ \ \ \ \ \text{if }v\in\operatorname*{NEtri}%
\left(  p\right)  ;\\
f\left(  \operatorname*{vrefl}\left(  v\right)  \right)
,\ \ \ \ \ \ \ \ \ \ \text{otherwise}%
\end{array}
\right.
\]
for all $v\in\widehat{\Delta\left(  p\right)  }$ for all $f\in\mathbb{K}%
^{\widehat{\operatorname*{NEtri}\left(  p\right)  }}$. This is well-defined.
We have%
\begin{equation}
\left(  \operatorname*{dble}f\right)  \left(  v\right)  =f\left(  v\right)
\ \ \ \ \ \ \ \ \ \ \text{for every }v\in\operatorname*{NEtri}\left(
p\right)  . \label{lem.NEtri.vrefl.b.doublev}%
\end{equation}
Also,%
\begin{equation}
\left(  \operatorname*{dble}f\right)  \left(  \operatorname*{vrefl}\left(
v\right)  \right)  =f\left(  v\right)  \ \ \ \ \ \ \ \ \ \ \text{for every
}v\in\operatorname*{NEtri}\left(  p\right)  .
\label{lem.NEtri.vrefl.b.doublevrefl}%
\end{equation}

\textbf{(c)} Assume that $p$ is odd. Then,%
\[
R_{\Delta\left(  p\right)  }\circ\operatorname*{dble}=\operatorname*{dble}%
\circ R_{\operatorname*{NEtri}\left(  p\right)  }.
\]

\end{lemma}

\begin{vershort}
\begin{proof}
[\nopunct]We omit the proofs of Lemma \ref{lem.NEtri.vrefl}, Proposition
\ref{prop.NEtri.odd.ord} and Proposition \ref{prop.SEtri.odd.ord} since
neither of them involves any new ideas. The first is analogous to that of
Lemma \ref{lem.Leftri.vrefl} (with $\Delta\left(  p\right)  $ and
$\operatorname*{NEtri}\left(  p\right)  $ taking the roles of
$\operatorname*{Rect}\left(  p,p\right)  $ and $\operatorname*{Tria}\left(
p\right)  $, respectively)\footnote{The only non-straightforward change that
must be made to the proof is the following: In Case 2 of the proof of Lemma
\ref{lem.Leftri.vrefl}, we used the (obvious) observation that $\left(
i-1,i\right)  $ and $\left(  i,i-1\right)  $ are elements of
$\operatorname*{Rect}\left(  p,p\right)  $ for every $\left(  i,i\right)
\in\operatorname*{Rect}\left(  p,p\right)  $ satisfying $i\neq1$. The
analogous observation that we need for proving Lemma \ref{lem.NEtri.vrefl} is
still true in the case of odd $p$, but a bit less obvious. In fact, it is the
observation that $\left(  i-1,i\right)  $ and $\left(  i,i-1\right)  $ are
elements of $\Delta\left(  p\right)  $ for every $\left(  i,i\right)
\in\Delta\left(  p\right)  $. This uses the oddness of $p$.}. The proof of
Proposition \ref{prop.NEtri.odd.ord} combines Lemma \ref{lem.NEtri.vrefl} with
Theorem \ref{thm.Delta.halfway}. Proposition \ref{prop.SEtri.odd.ord} is
derived from Proposition \ref{prop.NEtri.odd.ord} using Proposition
\ref{prop.op.ord}.
\end{proof}
\end{vershort}

\begin{verlong}
\begin{proof}
[Proof of Lemma \ref{lem.NEtri.vrefl} (sketched).]The proof of Lemma
\ref{lem.NEtri.vrefl} is analogous to that of Lemma \ref{lem.Leftri.vrefl}
(with $\Delta\left(  p\right)  $ and $\operatorname*{NEtri}\left(  p\right)  $
taking the roles of $\operatorname*{Rect}\left(  p,p\right)  $ and
$\operatorname*{Tria}\left(  p\right)  $, respectively). The only
non-straightforward change that must be made to the proof is the following: In
Case 2 of the proof of Lemma \ref{lem.Leftri.vrefl}, we used the (obvious)
observation that $\left(  i-1,i\right)  $ and $\left(  i,i-1\right)  $ are
elements of $\operatorname*{Rect}\left(  p,p\right)  $ for every $\left(
i,i\right)  \in\operatorname*{Rect}\left(  p,p\right)  $ satisfying $i\neq1$.
The analogous observation that we need for proving Lemma \ref{lem.NEtri.vrefl}
is still true when $p$ is odd; but it is a bit less obvious. In fact, it is
the observation that $\left(  i-1,i\right)  $ and $\left(  i,i-1\right)  $ are
elements of $\Delta\left(  p\right)  $ for every $\left(  i,i\right)
\in\Delta\left(  p\right)  $ when $p$ is odd.\footnote{Here is a
\textit{proof} of this observation: Let $\left(  i,i\right)  \in\Delta\left(
p\right)  $. Then, $i+i>p+1$ (since $\left(  i,i\right)  \in\Delta\left(
p\right)  $), so that $i+i-1\geq p+1$ (since $i+i$ and $p+1$ are integers).
But $i+i-1=2i-1$ is odd, while $p+1$ is even (since $p$ is odd). Hence,
$i+i-1\neq p+1$. Combined with $i+i-1\geq p+1$, this yields $i+i-1>p+1$. Thus,
$\left(  i,i-1\right)  \in\Delta\left(  p\right)  $, so that $\left(
i-1,i\right)  \in\Delta\left(  p\right)  $ as well, qed.} The reader can check
that, after this change and obvious replacements, the proof of Lemma
\ref{lem.Leftri.vrefl} becomes a proof of Lemma \ref{lem.NEtri.vrefl}.
\end{proof}

\begin{proof}
[Proof of Proposition \ref{prop.NEtri.odd.ord} (sketched).]From Proposition
\ref{prop.ord-projord} (applied to $p-1$ and $\operatorname*{NEtri}\left(
p\right)  $ instead of $n$ and $P$), we obtain $\operatorname*{ord}\left(
R_{\operatorname*{NEtri}\left(  p\right)  }\right)  =\operatorname{lcm}\left(
p-1+1,\operatorname*{ord}\left(  \overline{R}_{\operatorname*{NEtri}\left(
p\right)  }\right)  \right)  $. Hence, $\operatorname*{ord}\left(
R_{\operatorname*{NEtri}\left(  p\right)  }\right)  $ is divisible by
$p-1+1=p$. Thus, it remains to prove that $\operatorname*{ord}\left(
R_{\operatorname*{NEtri}\left(  p\right)  }\right)  \mid p$ (because then,
combining this with the fact that $\operatorname*{ord}\left(
R_{\operatorname*{NEtri}\left(  p\right)  }\right)  $ is divisible by $p$, we
will conclude that $\operatorname*{ord}\left(  R_{\operatorname*{NEtri}\left(
p\right)  }\right)  =p$, and this will prove Proposition
\ref{prop.NEtri.odd.ord}).

So let us show that $\operatorname*{ord}\left(  R_{\operatorname*{NEtri}%
\left(  p\right)  }\right)  \mid p$. This means showing that
$R_{\operatorname*{NEtri}\left(  p\right)  }^{p}=\operatorname*{id}$. Since
this statement boils down to a collection of polynomial identities in the
labels of an arbitrary $\mathbb{K}$-labelling of $\operatorname*{NEtri}\left(
p\right)  $, it is clear that it is enough to prove it in the case when
$\mathbb{K}$ is a field of rational functions in finitely many variables over
$\mathbb{Q}$. So let us WLOG assume that $\mathbb{K}$ is a field of rational
functions in finitely many variables over $\mathbb{Q}$. Then, the
characteristic of $\mathbb{K}$ is $\neq2$ (it is $0$ indeed), so that we can
apply Lemma \ref{lem.NEtri.vrefl}.

Define the maps $\operatorname*{vrefl}:\Delta\left(  p\right)  \rightarrow
\Delta\left(  p\right)  $ and $\operatorname*{vrefl}:\widehat{\Delta\left(
p\right)  }\rightarrow\widehat{\Delta\left(  p\right)  }$ as in Lemma
\ref{lem.NEtri.vrefl} \textbf{(a)}. Define the map $\operatorname*{dble}%
:\mathbb{K}^{\widehat{\operatorname*{NEtri}\left(  p\right)  }}\rightarrow
\mathbb{K}^{\widehat{\Delta\left(  p\right)  }}$ as in Lemma
\ref{lem.NEtri.vrefl} \textbf{(b)}.

Lemma \ref{lem.NEtri.vrefl} \textbf{(c)} yields $R_{\Delta\left(  p\right)
}\circ\operatorname*{dble}=\operatorname*{dble}\circ R_{\operatorname*{NEtri}%
\left(  p\right)  }$. Using this, we can prove that%
\begin{equation}
R_{\Delta\left(  p\right)  }^{\ell}\circ\operatorname*{dble}%
=\operatorname*{dble}\circ R_{\operatorname*{NEtri}\left(  p\right)  }^{\ell
}\ \ \ \ \ \ \ \ \ \ \text{for every }\ell\in\mathbb{N}.
\label{pf.NEtri.odd.ord.1}%
\end{equation}
(Indeed, (\ref{pf.NEtri.odd.ord.1}) easily follows from $R_{\Delta\left(
p\right)  }\circ\operatorname*{dble}=\operatorname*{dble}\circ
R_{\operatorname*{NEtri}\left(  p\right)  }$ by induction over $\ell$.)

Let $f\in\mathbb{K}^{\widehat{\operatorname*{NEtri}\left(  p\right)  }}$ be
sufficiently generic for $R_{\operatorname*{NEtri}\left(  p\right)  }^{\ell}f$
to be well-defined for all $\ell\in\left\{  0,1,...,p\right\}  $.

Applying (\ref{pf.NEtri.odd.ord.1}) to $\ell=p$, we obtain%
\[
R_{\Delta\left(  p\right)  }^{p}\circ\operatorname*{dble}=\operatorname*{dble}%
\circ R_{\operatorname*{NEtri}\left(  p\right)  }^{p}.
\]
Now,%
\begin{align}
\operatorname*{dble}\left(  R_{\operatorname*{NEtri}\left(  p\right)  }%
^{p}f\right)   &  =\underbrace{\left(  \operatorname*{dble}\circ
R_{\operatorname*{NEtri}\left(  p\right)  }^{p}\right)  }_{=R_{\Delta\left(
p\right)  }^{p}\circ\operatorname*{dble}}f=\left(  R_{\Delta\left(  p\right)
}^{p}\circ\operatorname*{dble}\right)  f\nonumber\\
&  =R_{\Delta\left(  p\right)  }^{p}\left(  \operatorname*{dble}f\right)  .
\label{pf.NEtri.odd.ord.2}%
\end{align}

Now, let $v\in\widehat{\operatorname*{NEtri}\left(  p\right)  }$ be arbitrary.
We are going to prove that $\left(  \operatorname*{dble}\left(
R_{\operatorname*{NEtri}\left(  p\right)  }^{p}f\right)  \right)  \left(
v\right)  =\left(  \operatorname*{dble}f\right)  \left(  v\right)  $.

We distinguish between the following three cases:

\textit{Case 1:} We have $v=0$.

\textit{Case 2:} We have $v=1$.

\textit{Case 3:} We have neither $v=0$ nor $v=1$.

Let us first consider Case 1. In this case, we have $v=0$. Thus,
\begin{align*}
\left(  \operatorname*{dble}\left(  R_{\operatorname*{NEtri}\left(  p\right)
}^{p}f\right)  \right)  \left(  \underbrace{v}_{=0}\right)   &  =\left(
\operatorname*{dble}\left(  R_{\operatorname*{NEtri}\left(  p\right)  }%
^{p}f\right)  \right)  \left(  0\right)  =2\underbrace{\left(
R_{\operatorname*{NEtri}\left(  p\right)  }^{p}f\right)  \left(  0\right)
}_{\substack{=f\left(  0\right)  \\\text{(by Corollary \ref{cor.R.implicit.01}%
)}}}\\
&  \ \ \ \ \ \ \ \ \ \ \left(  \text{by the definition of }%
\operatorname*{dble}\left(  R_{\operatorname*{NEtri}\left(  p\right)  }%
^{p}f\right)  \right) \\
&  =2f\left(  0\right)
\end{align*}
and%
\[
\left(  \operatorname*{dble}f\right)  \left(  \underbrace{v}_{=0}\right)
=\left(  \operatorname*{dble}f\right)  \left(  0\right)  =f\left(  0\right)
\ \ \ \ \ \ \ \ \ \ \left(  \text{by the definition of }\operatorname*{dble}%
f\right)  .
\]
Hence, $\left(  \operatorname*{dble}\left(  R_{\operatorname*{NEtri}\left(
p\right)  }^{p}f\right)  \right)  \left(  v\right)  =2f\left(  0\right)
=\left(  \operatorname*{dble}f\right)  \left(  v\right)  $. Thus, $\left(
\operatorname*{dble}\left(  R_{\operatorname*{NEtri}\left(  p\right)  }%
^{p}f\right)  \right)  \left(  v\right)  =\left(  \operatorname*{dble}%
f\right)  \left(  v\right)  $ is proven in Case 1.

Let us next consider Case 2. In this case, we have $v=1$. Thus,
\begin{align*}
\left(  \operatorname*{dble}\left(  R_{\operatorname*{NEtri}\left(  p\right)
}^{p}f\right)  \right)  \left(  \underbrace{v}_{=1}\right)   &  =\left(
\operatorname*{dble}\left(  R_{\operatorname*{NEtri}\left(  p\right)  }%
^{p}f\right)  \right)  \left(  1\right)  =\dfrac{1}{2}\underbrace{\left(
R_{\operatorname*{NEtri}\left(  p\right)  }^{p}f\right)  \left(  1\right)
}_{\substack{=f\left(  1\right)  \\\text{(by Corollary \ref{cor.R.implicit.01}%
)}}}\\
&  \ \ \ \ \ \ \ \ \ \ \left(  \text{by the definition of }%
\operatorname*{dble}\left(  R_{\operatorname*{NEtri}\left(  p\right)  }%
^{p}f\right)  \right) \\
&  =\dfrac{1}{2}f\left(  1\right)
\end{align*}
and%
\[
\left(  \operatorname*{dble}f\right)  \left(  \underbrace{v}_{=1}\right)
=\left(  \operatorname*{dble}f\right)  \left(  1\right)  =\dfrac{1}{2}f\left(
1\right)  \ \ \ \ \ \ \ \ \ \ \left(  \text{by the definition of
}\operatorname*{dble}f\right)  .
\]
Hence, $\left(  \operatorname*{dble}\left(  R_{\operatorname*{NEtri}\left(
p\right)  }^{p}f\right)  \right)  \left(  v\right)  =\dfrac{1}{2}f\left(
1\right)  =\left(  \operatorname*{dble}f\right)  \left(  v\right)  $. Thus,
$\left(  \operatorname*{dble}\left(  R_{\operatorname*{NEtri}\left(  p\right)
}^{p}f\right)  \right)  \left(  v\right)  =\left(  \operatorname*{dble}%
f\right)  \left(  v\right)  $ is proven in Case 2.

Finally, let us consider Case 3. In this case, we have neither $v=0$ nor
$v=1$. Consequently, $v\in\widehat{\operatorname*{NEtri}\left(  p\right)
}\setminus\left\{  0,1\right\}  =\operatorname*{NEtri}\left(  p\right)
\subseteq\Delta\left(  p\right)  \subseteq\operatorname*{Rect}\left(
p,p\right)  $. Thus, we can write $v$ in the form $\left(  i,k\right)  $ for
some $i\in\left\{  1,2,...,p\right\}  $ and $k\in\left\{  1,2,...,p\right\}
$. Consider these $i$ and $k$. Then,%
\begin{align*}
&  \underbrace{\left(  \operatorname*{dble}\left(  R_{\operatorname*{NEtri}%
\left(  p\right)  }^{p}f\right)  \right)  }_{\substack{=\left(  R_{\Delta
\left(  p\right)  }^{p}\left(  \operatorname*{dble}f\right)  \right)
\\\text{(by (\ref{pf.NEtri.odd.ord.2}))}}}\left(  \underbrace{v}_{=\left(
i,k\right)  }\right) \\
&  =\left(  R_{\Delta\left(  p\right)  }^{p}\left(  \operatorname*{dble}%
f\right)  \right)  \left(  \left(  i,k\right)  \right)  =\left(
\operatorname*{dble}f\right)  \left(  \underbrace{\left(  k,i\right)
}_{=\operatorname*{vrefl}\left(  i,k\right)  }\right) \\
&  \ \ \ \ \ \ \ \ \ \ \left(  \text{by Theorem \ref{thm.Delta.halfway},
applied to }\operatorname*{dble}f\text{ instead of }f\right) \\
&  =\left(  \operatorname*{dble}f\right)  \left(  \operatorname*{vrefl}%
\underbrace{\left(  i,k\right)  }_{=v}\right)  =\left(  \operatorname*{dble}%
f\right)  \left(  \operatorname*{vrefl}\left(  v\right)  \right)  =f\left(
v\right)  \ \ \ \ \ \ \ \ \ \ \left(  \text{by
(\ref{lem.NEtri.vrefl.b.doublevrefl})}\right) \\
&  =\left(  \operatorname*{dble}f\right)  \left(  v\right)
\ \ \ \ \ \ \ \ \ \ \left(  \text{by (\ref{lem.NEtri.vrefl.b.doublev}%
)}\right)  .
\end{align*}
Hence, $\left(  \operatorname*{dble}\left(  R_{\operatorname*{NEtri}\left(
p\right)  }^{p}f\right)  \right)  \left(  v\right)  =\left(
\operatorname*{dble}f\right)  \left(  v\right)  $ is proven in Case 3.

Now, $\left(  \operatorname*{dble}\left(  R_{\operatorname*{NEtri}\left(
p\right)  }^{p}f\right)  \right)  \left(  v\right)  =\left(
\operatorname*{dble}f\right)  \left(  v\right)  $ is proven in each of the
three Cases 1, 2 and 3. Since these three Cases cover all possibilities, this
shows that $\left(  \operatorname*{dble}\left(  R_{\operatorname*{NEtri}%
\left(  p\right)  }^{p}f\right)  \right)  \left(  v\right)  =\left(
\operatorname*{dble}f\right)  \left(  v\right)  $ always holds.

Now forget that we fixed $v$. We thus have shown that $\left(
\operatorname*{dble}\left(  R_{\operatorname*{NEtri}\left(  p\right)  }%
^{p}f\right)  \right)  \left(  v\right)  =\left(  \operatorname*{dble}%
f\right)  \left(  v\right)  $ for every $v\in\widehat{\operatorname*{NEtri}%
\left(  p\right)  }$. In other words, $\operatorname*{dble}\left(
R_{\operatorname*{NEtri}\left(  p\right)  }^{p}f\right)  =\operatorname*{dble}%
f$. Since the map $\operatorname*{dble}$ is injective, this yields
$R_{\operatorname*{NEtri}\left(  p\right)  }^{p}f=f$.

Now, forget that we fixed $f$. We thus have proven that
$R_{\operatorname*{NEtri}\left(  p\right)  }^{p}f=f$ for every $f\in
\mathbb{K}^{\widehat{\operatorname*{NEtri}\left(  p\right)  }}$ which is
sufficiently generic for $R_{\operatorname*{NEtri}\left(  p\right)  }^{\ell}f$
to be well-defined for all $\ell\in\left\{  0,1,...,p\right\}  $. In other
words, $R_{\operatorname*{NEtri}\left(  p\right)  }^{p}=\operatorname*{id}$ as
rational maps $\mathbb{K}^{\widehat{\operatorname*{NEtri}\left(  p\right)  }%
}\dashrightarrow\mathbb{K}^{\widehat{\operatorname*{NEtri}\left(  p\right)  }%
}$. Hence, $\operatorname*{ord}\left(  R_{\operatorname*{NEtri}\left(
p\right)  }\right)  \mid p$. As we know, this completes the proof of
Proposition \ref{prop.NEtri.odd.ord}.
\end{proof}

\begin{proof}
[Proof of Proposition \ref{prop.SEtri.odd.ord} (sketched).]Applying
Proposition \ref{prop.op.ord} to $P=\operatorname*{NEtri}\left(  p\right)  $,
we obtain $\operatorname*{ord}\left(  R_{\left(  \operatorname*{NEtri}\left(
p\right)  \right)  ^{\operatorname*{op}}}\right)  =\operatorname*{ord}\left(
R_{\operatorname*{NEtri}\left(  p\right)  }\right)  $ and $\operatorname*{ord}%
\left(  \overline{R}_{\left(  \operatorname*{NEtri}\left(  p\right)  \right)
^{\operatorname*{op}}}\right)  =\operatorname*{ord}\left(  \overline
{R}_{\operatorname*{NEtri}\left(  p\right)  }\right)  $.

Consider the poset antiautomorphism $\operatorname*{hrefl}$ of
$\operatorname*{Rect}\left(  p,p\right)  $ defined in Remark
\ref{rmk.DeltaNabla}. This antiautomorphism sends $\operatorname*{NEtri}%
\left(  p\right)  $ to $\operatorname*{SEtri}\left(  p\right)  $. Thus,
$\operatorname*{SEtri}\left(  p\right)  \cong\left(  \operatorname*{NEtri}%
\left(  p\right)  \right)  ^{\operatorname*{op}}$, so that%
\[
\operatorname*{ord}\left(  R_{\operatorname*{SEtri}\left(  p\right)  }\right)
=\operatorname*{ord}\left(  R_{\left(  \operatorname*{NEtri}\left(  p\right)
\right)  ^{\operatorname*{op}}}\right)  =\operatorname*{ord}\left(
R_{\operatorname*{NEtri}\left(  p\right)  }\right)  =p
\]
(by Proposition \ref{prop.NEtri.odd.ord}). This proves Proposition
\ref{prop.SEtri.odd.ord}.
\end{proof}
\end{verlong}

Nathan Williams suggested that the following generalization of Conjecture
\ref{conj.NEtri.ord} might hold:

\begin{conjecture}
\label{conj.NEtriminus.ord}Let $p$ be an integer $>1$. Let $s\in\mathbb{N}$.
Let $\operatorname*{NEtri}\nolimits^{\prime}\left(  p\right)  $ be the
subposet $\left\{  \left(  i,k\right)  \in\operatorname*{NEtri}\left(
p\right)  \ \mid\ k\geq s\right\}  $ of $\operatorname*{NEtri}\left(
p\right)  $. Then, $\operatorname*{ord}\left(  R_{\operatorname*{NEtri}%
\nolimits^{\prime}\left(  p\right)  }\right)  \mid p$.
\end{conjecture}

This conjecture has been verified using Sage for all $p\leq7$. Williams (based
on a philosophy from his thesis \cite{williams-cataland}) suspects there could
be a birational map between $\mathbb{K}^{\widehat{\operatorname*{NEtri}%
\nolimits^{\prime}\left(  p\right)  }}$ and $\mathbb{K}%
^{\widehat{\operatorname*{Rect}\left(  s-1,p-s+1\right)  }}$ which commutes
with the respective birational rowmotion operators for all $s>\dfrac{p}{2}$;
this, if shown, would obviously yield a proof of Conjecture
\ref{conj.NEtriminus.ord}. This already is an interesting question for
classical rowmotion; a bijection between the antichains (and thus between the
order ideals) of $\operatorname*{NEtri}\nolimits^{\prime}\left(  p\right)  $
and those of $\operatorname*{Rect}\left(  s-1,p-s+1\right)  $ was found by
Stembridge \cite[Theorem 5.4]{stembridge-trapezoid}, but does not commute with
classical rowmotion.

\section{\label{sect.negres}Negative results}

Generally, it is not true that if $P$ is an $n$-graded poset, then
$\operatorname*{ord}\left(  R_{P}\right)  $ is necessarily finite. When
$\operatorname*{char}\mathbb{K}=0$, the authors have proven the
following\footnote{See the ancillary files of the present arXiv preprint
(arXiv:1402.6178) for an outline of the (rather technical) proofs.}:

\begin{itemize}
\item If $P$ is the poset $\left\{  x_{1},x_{2},x_{3},x_{4},x_{5}\right\}  $
with relations $x_{1}<x_{3}$, $x_{1}<x_{4}$, $x_{1}<x_{5}$, $x_{2}<x_{4}$ and
$x_{2}<x_{5}$ (this is a $5$-element $2$-graded poset), then
$\operatorname*{ord}\left(  R_{P}\right)  =\infty$.

\item If $P$ is the \textquotedblleft chain-link fence\textquotedblright%
\ poset $/\backslash/\backslash/\backslash$ (that is, the subposet
\newline$\left\{  \left(  i,k\right)  \in\operatorname*{Rect}\left(
4,4\right)  \ \mid\ 5\leq i+k\leq6\right\}  $ of $\operatorname*{Rect}\left(
4,4\right)  $), then $\operatorname*{ord}\left(  R_{P}\right)  =\infty$.

\item If $P$ is the Boolean lattice $\left[  2\right]  \times\left[  2\right]
\times\left[  2\right]  $, then $\operatorname*{ord}\left(  R_{P}\right)
=\infty$.
\end{itemize}

The situation seems even more hopeless for non-graded posets.

\section{\label{sect.rootsys}The root system connection}

A question naturally suggesting itself is: What is it that makes certain
posets $P$ have finite $\operatorname*{ord}\left(  R_{P}\right)  $, while
others have not? Can we characterize the former posets? It might be too
optimistic to expect a full classification, given that our examples are
already rather diverse (skeletal posets, rectangles, triangles, posets like
that in Remark \ref{rmk.rect.ord.application}). As a first step (and inspired
by the general forms of the Zamolodchikov conjecture), we were tempted to
study posets arising from Dynkin diagrams. It appears that, unlike in the
Zamolodchikov conjecture, the interesting cases are not those having $P$ be a
product of Dynkin diagrams, but those having $P$ be a positive root poset of a
root system, or a parabolic quotient thereof. The idea is not new, as it was
already conjectured by Panyushev \cite[Conjecture 2.1]{panyushev} and proven
by Armstrong, Stump and Thomas \cite[Theorem 1.2]{armstrong-stump-thomas} that
if $W$ is a finite Weyl group with Coxeter number $h$, then classical
rowmotion on the set $J\left(  \Phi^{+}\left(  W\right)  \right)  $ (where
$\Phi^{+}\left(  W\right)  $ is the poset of positive roots of $W$) has order
$h$ or $2h$ (along with a few more properties, akin to our \textquotedblleft
reciprocity\textquotedblright\ statements)\footnote{Neither \cite{panyushev}
nor \cite{armstrong-stump-thomas} work directly with order ideals and
rowmotion, but instead they study antichains of the poset $\Phi^{+}\left(
W\right)  $ (which are called \textquotedblleft nonnesting
partitions\textquotedblright\ in \cite{armstrong-stump-thomas}) and an
operation on these antichains called Panyushev complementation. There is,
however, a simple bijection between the set of antichains of a poset $P$ and
the set $J\left(  P\right)  $, and the conjugate of Panyushev complementation
with respect to this bijection is precisely classical rowmotion.}.

In the case of birational rowmotion, the situation is less simple.
Specifically, the following can be said about positive root posets of
crystallographic root systems (as considered in \cite[\S 3.2]%
{striker-williams})\footnote{We refer to \cite[Definition 3.4]%
{striker-williams} for notations.}:

\begin{itemize}
\item If $P=\Phi^{+}\left(  A_{n}\right)  $ for $n\geq2$, then
$\operatorname*{ord}\left(  R_{P}\right)  =2\left(  n+1\right)  $. This is
just the assertion of Corollary \ref{cor.Delta.ord}. Note that for $n=1$, the
order $\operatorname*{ord}\left(  R_{P}\right)  $ is $2$ instead of $2\left(
1+1\right)  =4$.

\item If $P=\Phi^{+}\left(  B_{n}\right)  $ for $n\geq1$, then Conjecture
\ref{conj.SEtri.ord} claims that $\operatorname*{ord}\left(  R_{P}\right)
=2n$. Note that $\Phi^{+}\left(  B_{n}\right)  \cong\Phi^{+}\left(
C_{n}\right)  $.

\item We have $\operatorname*{ord}\left(  R_{P}\right)  =2$ for $P=\Phi
^{+}\left(  D_{2}\right)  $, and we have $\operatorname*{ord}\left(
R_{P}\right)  =8$ for $P=\Phi^{+}\left(  D_{3}\right)  $. However,
$\operatorname*{ord}\left(  R_{P}\right)  =\infty$ in the case when
$P=\Phi^{+}\left(  D_{4}\right)  $. This should not come as a surprise, since
$\Phi^{+}\left(  D_{4}\right)  $ has a property that none of the $\Phi
^{+}\left(  A_{n}\right)  $ or $\Phi^{+}\left(  B_{n}\right)  \cong\Phi
^{+}\left(  C_{n}\right)  $ have, namely an element covered by three other
elements. On the other hand, the finite orders in the $\Phi^{+}\left(
D_{2}\right)  $ and $\Phi^{+}\left(  D_{3}\right)  $ cases can be explained by
$\Phi^{+}\left(  D_{2}\right)  \cong\Phi^{+}\left(  A_{1}\times A_{1}\right)
\cong\left(  \text{two-element antichain}\right)  $ and $\Phi^{+}\left(
D_{3}\right)  \cong\Phi^{+}\left(  A_{3}\right)  $.
\end{itemize}

Nathan Williams has suggested that the behavior of $\Phi^{+}\left(
A_{n}\right)  $ and $\Phi^{+}\left(  B_{n}\right)  \cong\Phi^{+}\left(
C_{n}\right)  $ to have finite orders of $R_{P}$ could generalize to the
\textquotedblleft positive root posets\textquotedblright\ of the other
\textquotedblleft coincidental types\textquotedblright\ $H_{3}$ and
$I_{2}\left(  m\right)  $ (see, for example, Table 2.2 in
\cite{williams-cataland}). And indeed, computations in Sage have established
that $\operatorname*{ord}\left(  R_{P}\right)  =10$ for $P=\Phi^{+}\left(
H_{3}\right)  $, and we also have $\operatorname*{ord}\left(  R_{P}\right)
=\operatorname{lcm}\left(  2,m\right)  $ for $P=\Phi^{+}\left(  I_{2}\left(
m\right)  \right)  $ (this is a very easy consequence of Proposition
\ref{prop.ord-projord}).

It seems that minuscule heaps, as considered e.g. in \cite[\S 6]{rush-shi},
also lead to small $\operatorname*{ord}\left(  R_{P}\right)  $ values. Namely:

\begin{itemize}
\item The heap $P_{w_{0}^{J}}$ in \cite[Figure 8 (b)]{rush-shi} satisfies
$\operatorname*{ord}\left(  R_{P}\right)  =12$.

\item The heap $P_{w_{0}^{J}}$ in \cite[Figure 9 (b)]{rush-shi} seems to
satisfy $\operatorname*{ord}\left(  R_{P}\right)  =18$ (this was verified on
numerical examples, as the poset is too large for efficient general computations).
\end{itemize}

(These two posets also appear as posets corresponding to the \textquotedblleft
Cayley plane\textquotedblright\ and the \textquotedblleft Freudenthal
variety\textquotedblright\ in \cite[p. 2]{thomas-yong}.)

Various other families of posets related to root systems (minuscule posets,
d-complete posets, rc-posets, alternating sign matrix posets) remain to be studied.

\end{document}